\pdfoutput=1

\documentclass[10pt,reqno]{amsbook}

\usepackage[T1]{fontenc}
\usepackage[utf8x]{inputenc}
\usepackage[french,english]{babel}
\usepackage{amsmath,amsfonts,amssymb,amsthm,bm,shuffle}
\usepackage[math]{anttor}

\usepackage[top=3.5cm,bottom=3.5cm,left=3.4cm,right=3.4cm]{geometry}

\usepackage{xcolor}
\definecolor{ColBlack}{RGB}{0,0,0} 
\definecolor{ColWhite}{RGB}{255,255,255} 
\definecolor{Col1}{RGB}{133,6,6} 
\definecolor{Col2}{RGB}{198,8,0} 
\definecolor{Col3}{RGB}{174,74,52} 
\definecolor{Col4}{RGB}{103,113,121} 
\definecolor{Col5}{RGB}{90,94,107} 
\definecolor{Col6}{RGB}{70,63,50} 
\usepackage[hyperindex=true,frenchlinks=true,colorlinks=true,
citecolor=Col2,linkcolor=Col3,urlcolor=Col4,linktocpage,
pagebackref=true]{hyperref}

\usepackage{tikz}
\usetikzlibrary{shapes}
\usetikzlibrary{fit}
\usetikzlibrary{decorations.pathmorphing}
\usetikzlibrary{calc}
\usetikzlibrary{automata}

\usepackage{mathtools}
\usepackage{dsfont}
\usepackage{wasysym}
\usepackage{stmaryrd}
\usepackage{nameref}
\usepackage{youngtab}
\usepackage{cite}
\usepackage{xspace}
\usepackage{subfig}
\usepackage{multirow}
\usepackage{enumitem}
\usepackage{multicol}
\usepackage{chngcntr}

\tikzstyle{Centering}=[{baseline={([yshift=-0.5ex]current
    bounding box.center)}}]

\tikzstyle{Node}=[circle,draw=Col1!80,fill=Col1!8,inner sep=1pt,
minimum size=2mm,thick,font=\scriptsize]
\tikzstyle{Edge}=[draw=Col2!80,cap=round,thick]
\tikzstyle{Leaf}=[rectangle,draw=ColBlack!70,fill=ColBlack!16,
inner sep=0pt,minimum size=1mm,thick]
\tikzstyle{NodeClear}=[Node,fill=ColWhite!100]
\tikzstyle{NodeColorA}=[Node,draw=Col4!80,fill=Col4!8]
\tikzstyle{NodeColorB}=[Node,draw=Col6!80,fill=Col6!8]
\tikzstyle{NodeST}=[font=\footnotesize]
\tikzstyle{EdgeLabel}=[midway,inner sep=1pt,fill=ColWhite!0,
font=\scriptsize]
\tikzstyle{LeafLabel}=[font=\scriptsize,node distance=2mm]
\tikzstyle{EdgeValue}=[regular polygon,regular polygon sides=6,
draw=Col2!80,fill=Col2!2,inner sep=.2pt,line width=.75pt,
minimum size=2.8mm,midway,text=ColBlack,font=\tiny]
\tikzstyle{Subtree}=[regular polygon,regular polygon sides=3,
draw=Col1!80,fill=Col1!20,thick,minimum size=5mm,font=\scriptsize]

\tikzstyle{PathNode}=[circle,draw=Col1!90,fill=Col1!30,thick,
inner sep=0pt,minimum size=2.0mm]
\tikzstyle{PathNodeColorA}=[PathNode,draw=Col4!80,fill=Col4!18]
\tikzstyle{PathNodeColorB}=[PathNode,draw=Col6!80,fill=Col6!28]
\tikzstyle{PathStep}=[color=Col1!60,thick]

\tikzstyle{Box}=[rectangle,rounded corners=2pt,draw=Col1!80,fill=Col1!8,
inner sep=1pt,minimum size=.25cm,thick,font=\scriptsize]
\tikzstyle{BoxClear}=[Box,fill=ColWhite!100]
\tikzstyle{BoxColorA}=[Box,draw=Col4!80,fill=Col4!8]
\tikzstyle{BoxColorB}=[Box,draw=Col6!80,fill=Col6!8]
\tikzstyle{Arc}=[Edge,->,draw=Col2!80]
\tikzstyle{Grid} = [color=ColBlack!30]
\tikzstyle{EdgeRew}=[->,Col2!80,cap=round,thick]

\tikzstyle{Injection}=[ColBlack!100,draw,>->]
\tikzstyle{Surjection}=[ColBlack!100,draw,->>]

\tikzstyle{PolyEdgeGray}=[ColBlack!20,draw,cap=round]
\tikzstyle{PolyEdgeBlue}=[Col1!90,thick,draw,cap=round]
\tikzstyle{PolyEdgeRed}=[Col2!70,thick,draw,cap=round,dotted]

\tikzstyle{PosetVertex}=[circle,draw=Col1!80,fill=Col1!8,
inner sep=1pt,minimum size=2mm,thick,font=\scriptsize]

\tikzstyle{CliqueEdge}=[draw=Col2!90,thick]
\tikzstyle{CliqueEdgeColorA}=[CliqueEdge,draw=Col4!80,fill=Col4!8]
\tikzstyle{CliqueEmptyEdge}=[draw=Col4!90,thick,densely dashed]
\tikzstyle{CliqueLabel}=[midway,inner sep=1pt,fill=ColWhite!0,
font=\scriptsize]
\tikzstyle{CliquePoint}=[circle,inner sep=1pt,fill=Col2!25,
draw=Col2!70]
\tikzstyle{CliqueEdgeGray}=[ColBlack!30,draw,cap=round]
\tikzstyle{CliqueEdgeBlue}=[Col1!80,thick,draw,cap=round]
\tikzstyle{CliqueEdgeRed}=[Col2!80,thick,draw,cap=round,dotted]

\tikzstyle{SixVertexArrow}=[thick,Col1!90]
\tikzstyle{SixVertexGrid}=[thick,ColBlack!45]
\tikzstyle{SixVertexMarkedVertexA}=[circle,draw=Col1!80,
fill=Col1!20,inner sep=1pt,minimum size=2mm]
\tikzstyle{SixVertexMarkedVertexB}=[regular polygon,
regular polygon sides=5,draw=Col2!80,fill=Col2!20,
inner sep=1pt,minimum size=2mm]

\tikzstyle{Operator}=[rectangle,rounded corners,draw=Col1!100,
fill=Col1!10,minimum size=4mm,thick,font=\footnotesize]
\tikzstyle{OperatorColorA}=[Operator,draw=Col4!80,fill=Col4!8]
\tikzstyle{OperatorColorB}=[Operator,draw=Col6!80,fill=Col6!8]
\tikzstyle{OperatorColorC}=[Operator,draw=Col1!80,fill=Col1!8]

\tikzstyle{Domino}=[rectangle,draw=Col4!100,fill=Col4!20,
minimum height=.25cm,thick]
\tikzstyle{Domino3}=[Domino,minimum width=1.5cm]
\tikzstyle{Domino2}=[Domino,minimum width=1cm]
\tikzstyle{Domino1}=[Domino,minimum width=.5cm]
\tikzstyle{DominoFriable}=[rectangle,draw=Col2!100,
fill=Col2!20,minimum height=.25cm,thick,minimum width=.5cm]

\tikzstyle{MatchingVertex}=[circle,draw=Col1!80,fill=Col1!20,
minimum size=.85cm,inner sep=1pt,font=\footnotesize]
\tikzstyle{MatchingVertexColorA}=[MatchingVertex,draw=Col4!80,
fill=Col4!20]
\tikzstyle{MatchingEdge}=[draw=Col2!80,cap=round,thick]
\tikzstyle{DirectedMatchingEdge}=[draw=Col2!80,cap=round,thick,->]
\tikzstyle{PatternMatchingVertex}=[circle,draw=Col1!80,
fill=Col1!20,minimum size=.15cm,inner sep=1pt,font=\footnotesize]

\tikzstyle{AutState}=[circle,draw=Col2!80,fill=Col2!8,
inner sep=1pt,minimum size=8mm,thick,font=\footnotesize]
\tikzstyle{AutArc}=[draw=Col4!80,thick,->]
\tikzstyle{AutArcColorA}=[AutArc,draw=Col2!80]
\tikzstyle{AutArcColorB}=[AutArc,draw=Col4!80]

\tikzstyle{ChapterBox}=[rectangle,text width=2cm,text centered,
font=\footnotesize]

\newcommand{\Hide}[1]{\textcolor{Col4}{--- \tt hidden}}
\newcommand{\Def}[1]{\textcolor{Col3}{\em #1}}
\newcommand{\OEIS}[1]{\href{http://oeis.org/#1}{{\bf #1}}}

\DeclareRobustCommand{\gobblefive}[5]{}
\newcommand*{\SkipTocEntry}{\addtocontents{toc}{\gobblefive}}

\let\SavedCaption=\caption
\renewcommand*{\caption}[2][\shortcaption]{%
    \def\shortcaption{#2}
    \SavedCaption[\; #1]{#2}}

\let\SavedParagraph=\paragraph
\renewcommand{\paragraph}[1]{%
    \SavedParagraph{\it #1}}

\newcommand*\ClearToLeftPage{%
    \clearpage
    \ifodd\value{page}
        \hbox{}
        \vspace*{\fill}
        \thispagestyle{empty}
        \newpage
    \fi
}

\newcommand{\N}{\mathbb{N}}
\newcommand{\Z}{\mathbb{Z}}
\newcommand{\Q}{\mathbb{Q}}
\newcommand{\R}{\mathbb{R}}

\newcommand{\K}{\mathbb{K}}

\newcommand{\Aca}{\mathcal{A}}
\newcommand{\Bca}{\mathcal{B}}
\newcommand{\Cca}{\mathcal{C}}
\newcommand{\Dca}{\mathcal{D}}
\newcommand{\Eca}{\mathcal{E}}
\newcommand{\Fca}{\mathcal{F}}
\newcommand{\Gca}{\mathcal{G}}
\newcommand{\Hca}{\mathcal{H}}
\newcommand{\Ica}{\mathcal{I}}
\newcommand{\Lca}{\mathcal{L}}
\newcommand{\Mca}{\mathcal{M}}
\newcommand{\Nca}{\mathcal{N}}
\newcommand{\Oca}{\mathcal{O}}
\newcommand{\Pca}{\mathcal{P}}
\newcommand{\Qca}{\mathcal{Q}}
\newcommand{\Rca}{\mathcal{R}}
\newcommand{\Sca}{\mathcal{S}}
\newcommand{\Tca}{\mathcal{T}}
\newcommand{\Vca}{\mathcal{V}}
\newcommand{\Zca}{\mathcal{Z}}

\newcommand{\Afr}{\mathfrak{a}}
\newcommand{\Bfr}{\mathfrak{b}}
\newcommand{\Cfr}{\mathfrak{c}}
\newcommand{\Dfr}{\mathfrak{d}}
\newcommand{\Efr}{\mathfrak{e}}
\newcommand{\Ffr}{\mathfrak{f}}

\newcommand{\Lfr}{\mathfrak{l}}
\newcommand{\Sfr}{\mathfrak{s}}
\newcommand{\Rfr}{\mathfrak{r}}
\newcommand{\Pfr}{\mathfrak{p}}
\newcommand{\Qfr}{\mathfrak{q}}
\newcommand{\Tfr}{\mathfrak{t}}

\newcommand{\CFr}{\mathfrak{C}}
\newcommand{\GFr}{\mathfrak{G}}
\newcommand{\RFr}{\mathfrak{R}}

\newcommand{\Abb}{\mathbb{A}}
\newcommand{\Dbb}{\mathbb{D}}
\newcommand{\Ebb}{\mathbb{E}}

\newcommand{\Sbb}{\mathbb{S}}

\newcommand{\Ubb}{\mathbb{U}}
\newcommand{\Xbb}{\mathbb{X}}
\newcommand{\Ybb}{\mathbb{Y}}
\newcommand{\Zbb}{\mathbb{Z}}

\newcommand{\Fbf}{\mathbf{f}}
\newcommand{\Gbf}{\mathbf{g}}
\newcommand{\Hbf}{\mathbf{h}}
\newcommand{\Rbf}{\mathbf{r}}
\newcommand{\Ibf}{\mathbf{i}}
\newcommand{\Sbf}{\mathbf{s}}
\newcommand{\Tbf}{\mathbf{t}}
\newcommand{\Ubf}{\mathbf{u}}
\newcommand{\Xbf}{\mathbf{x}}

\newcommand{\Asf}{\mathsf{a}}
\newcommand{\Bsf}{\mathsf{b}}
\newcommand{\Csf}{\mathsf{c}}
\newcommand{\Dsf}{\mathsf{d}}
\newcommand{\Esf}{\mathsf{e}}

\newcommand{\VarX}{\mathsf{x}}
\newcommand{\VarY}{\mathsf{y}}
\newcommand{\VarZ}{\mathsf{z}}

\newcommand{\LambdaB}{\bm{\lambda}}
\newcommand{\MuB}{\bm{\mu}}

\newcommand{\Unit}{\mathds{1}}
\newcommand{\OneElement}{\epsilon}
\newcommand{\AtomElement}{\bullet}
\newcommand{\GenSeries}{\Gca}
\newcommand{\HilbSeries}{\Hca}
\newcommand{\Leaf}{\perp}
\newcommand{\Node}{\bullet}
\newcommand{\Alg}{\Aca}

\newcommand{\FreeAlg}{\Fca}
\newcommand{\PrefixOrder}{\leq_{\mathrm{pref}}}
\newcommand{\SuffixOrder}{\leq_{\mathrm{suff}}}
\newcommand{\FactorOrder}{\leq_{\mathrm{fact}}}
\newcommand{\NormalForms}{\Fca}
\newcommand{\RelationSpaceRewriteRule}{\Rca}
\newcommand{\RelationSpace}{\Rca}
\newcommand{\GeneratingSet}{\GFr}
\newcommand{\Comparable}{\Cca}
\newcommand{\Nmax}{\mathbb{M}}
\newcommand{\BubbleSeries}{\mathrm{B}}

\newcommand{\Congr}{\equiv}
\newcommand{\Hook}{\mathrm{hook}}

\newcommand{\Root}{\mathrm{root}}
\newcommand{\OrdDias}{\preccurlyeq}
\newcommand{\Equiv}{\leftrightarrow}
\newcommand{\Identity}{\mathrm{Id}}
\newcommand{\Op}{\star}
\newcommand{\OpA}{\square}
\newcommand{\OpB}{\triangle}
\newcommand{\OpDual}{\bar \Op}
\newcommand{\OpDualA}{\bar \OpA}
\newcommand{\OpDualB}{\bar \OpB}
\newcommand{\Ord}{\preccurlyeq}
\newcommand{\OrdStrict}{\prec}
\newcommand{\PosetEmpty}{\emptyset}
\newcommand{\UnitSet}{\Eca}
\newcommand{\Ideal}{\Ica}
\newcommand{\CorSchr}{\Cfr}

\newcommand{\OrdBE}{\preceq_{\mathrm{be}}}
\newcommand{\OrdD}{\preceq_{\mathrm{d}}}

\newcommand{\SeriesBubbles}{\mathrm{B}}
\newcommand{\SeriesElements}{\mathrm{F}}
\newcommand{\One}{1}
\newcommand{\Zero}{0}
\newcommand{\OrdPM}{\leq_{\mathrm{PM}}}
\newcommand{\OrdMQ}{\leq_{\mathrm{MQ}}}
\newcommand{\CoveringRelMQ}{\rightharpoonup}
\newcommand{\CoveringRelPM}{\rightarrow}

\newcommand{\Plus}{\mathbf{\text{+}}}
\newcommand{\Minus}{\mathbf{\text{-}}}
\newcommand{\NE}{\mathrm{ne}}
\newcommand{\SE}{\mathrm{se}}
\newcommand{\SW}{\mathrm{sw}}
\newcommand{\NW}{\mathrm{nw}}
\newcommand{\OI}{\mathrm{oi}}
\newcommand{\IO}{\mathrm{io}}
\newcommand{\AllStatsA}{\mathfrak{Z}}
\newcommand{\AllStatsB}{\mathfrak{N}}
\newcommand{\CongrMonoid}[1]{{\Equiv_{\mathrm{#1}}}}
\newcommand{\PropDPMShape}{\mathbf{PShape}}
\newcommand{\PropDPMValue}{\mathbf{PValue}}
\newcommand{\Lang}{\mathrm{L}}
\newcommand{\SyncLang}{\mathrm{L}_{\mathrm{S}}}
\newcommand{\GenS}{\Sbf}
\newcommand{\Synt}{\mathrm{synt}}
\newcommand{\Sync}{\mathrm{sync}}
\newcommand{\Dummy}{\lozenge}
\newcommand{\Perfect}{\mathrm{perf}}
\newcommand{\Height}{\mathrm{ht}}
\newcommand{\Deriv}{\to}
\newcommand{\SyncDeriv}{\leadsto}
\newcommand{\DerivGraph}{\mathrm{G}}
\newcommand{\SyncDerivGraph}{\mathrm{G}_{\mathrm{S}}}
\newcommand{\MultOutIn}{\chi}
\newcommand{\BMotz}{\Bca_{\mathrm{p}}}
\newcommand{\BBalTree}{\Bca_{\mathrm{bbt}}}
\newcommand{\ArityIn}{\uparrow}
\newcommand{\ArityOut}{\downarrow}

\newcommand{\List}{\mathrm{T}}
\newcommand{\Multiset}{\mathrm{S}}
\newcommand{\Set}{\mathrm{E}}
\newcommand{\Support}{\mathrm{Supp}}
\newcommand{\Charac}{\mathrm{ch}}
\newcommand{\Angle}[1]{\left\langle#1\right\rangle}
\newcommand{\AAngle}[1]{%
    \left\langle\left\langle#1\right\rangle\right\rangle}
\newcommand{\ExtProd}{\cdot}
\newcommand{\Suspension}{\mathrm{Sus}}
\newcommand{\Augmentation}{\mathrm{Aug}}
\newcommand{\Length}{\ell}
\newcommand{\Sector}{\wedge}
\newcommand{\LDup}{\ll}
\newcommand{\RDup}{\gg}
\newcommand{\Tensor}{\mathsf{T}}
\newcommand{\Symmetric}{\mathsf{S}}
\newcommand{\Exterior}{\mathsf{E}}
\newcommand{\Product}{\star}
\newcommand{\Coproduct}{\Delta}
\newcommand{\Dual}{\star}
\newcommand{\Rew}{\Rightarrow}
\newcommand{\RewRT}{\overset{*}{\Rew}}
\newcommand{\RewTrees}{\leadsto}
\newcommand{\RewTreesRT}{\overset{*}{\RewTrees}}

\newcommand{\RewSchr}{\rightarrowtriangle}
\newcommand{\RewTreesSchr}{\rightsquigarrow}
\newcommand{\Std}{\mathrm{std}}
\newcommand{\Cmp}{\mathrm{cmp}}
\newcommand{\Des}{\mathrm{Des}}
\newcommand{\CoInv}{\mathrm{Civ}}
\newcommand{\Conc}{\cdot}
\newcommand{\ProductTriangleRight}{\triangleright}
\newcommand{\RefinementOrder}{\preceq}
\newcommand{\Narayana}{\mathrm{nar}}
\newcommand{\Catalan}{\mathrm{cat}}
\newcommand{\T}{\mathsf{T}}
\newcommand{\TT}{\bar{\T}}
\newcommand{\PreLieProduct}{\curvearrowleft}
\newcommand{\Graft}{\curlywedge}
\newcommand{\Action}{\cdot}
\newcommand{\TreeLanguage}{\Nca}
\newcommand{\Arity}{\mathrm{ari}}
\newcommand{\TamariInvariant}{\mathrm{tam}}
\newcommand{\FreeOperad}{\mathbf{FO}}
\newcommand{\FreeColoredOperad}{\mathbf{FCO}}
\newcommand{\FreePro}{\mathbf{FP}}
\newcommand{\Corolla}[1]{\ocircle\left(#1\right)}
\newcommand{\CorollaNOARG}{\ocircle}
\newcommand{\Eval}{\mathrm{ev}}
\newcommand{\In}{\mathbf{in}}
\newcommand{\Out}{\mathbf{out}}
\newcommand{\Hull}{\mathbf{Hull}}
\newcommand{\AntiCol}{\mathbf{Anti}}
\newcommand{\Border}{\mathrm{bor}}
\newcommand{\Cpl}{\mathrm{cpl}}
\newcommand{\Ret}{\mathrm{ret}}
\newcommand{\OperadGen}[1]{\langle #1 \rangle}
\newcommand{\OperadColGen}[1]{\langle\langle #1 \rangle\rangle}

\newcommand{\LeftInvariant}{\mathrm{ln}}
\newcommand{\LDias}{\dashv}
\newcommand{\RDias}{\vdash}
\newcommand{\LDiasA}{\rotatebox[origin=c]{180}{$\Vdash$}}
\newcommand{\RDiasA}{\Vdash}
\newcommand{\LDendrA}{\leftharpoonup}
\newcommand{\RDendrA}{\rightharpoonup}
\newcommand{\LDendr}{\prec}
\newcommand{\RDendr}{\succ}
\newcommand{\MAs}{\star}
\newcommand{\MAsA}{\triangle}
\newcommand{\MDAsA}{\wasylozenge}
\newcommand{\MDAs}{\diamond}
\newcommand{\MTrias}{\perp}
\newcommand{\MTDendr}{\wedge}
\newcommand{\OpAsDendrA}{\bullet}
\newcommand{\OpAsDendr}{\odot}
\newcommand{\ProdZin}{\shuffle}
\newcommand{\Min}{\uparrow}
\newcommand{\Max}{\downarrow}
\newcommand{\Hamming}{\mathrm{ham}}
\newcommand{\Incr}{\mathrm{Incr}}
\newcommand{\Augm}{\mathrm{h}}
\newcommand{\Halo}{\mathrm{Halo}}
\newcommand{\MProjToPluri}{\mathrm{M}}
\newcommand{\NbInterv}{\mathrm{int}}
\newcommand{\ASchrMap}{\mathrm{s}}
\newcommand{\Cli}{\mathsf{C}}
\newcommand{\Frac}{\mathrm{F}}
\newcommand{\Skel}{\mathrm{skel}}
\newcommand{\Degr}{\mathrm{degr}}
\newcommand{\Cros}{\mathrm{cros}}
\newcommand{\Returned}{\mathrm{ret}}
\newcommand{\Del}{\mathrm{d}}
\newcommand{\BubbleTree}{\mathrm{bt}}
\newcommand{\OpAssoc}{\odot}
\newcommand{\Shift}{\mathrm{sh}}
\newcommand{\NullRows}{\mathrm{N}_{\mathrm{r}}}
\newcommand{\NullColumns}{\mathrm{N}_{\mathrm{c}}}
\newcommand{\LastRow}{\mathrm{last}_{\mathrm{r}}}
\newcommand{\LastColumn}{\mathrm{last}_{\mathrm{c}}}
\newcommand{\Over}{\diagup}
\newcommand{\Under}{\diagdown}
\newcommand{\Compr}{\mathrm{cp}}
\newcommand{\DecompC}{\circ}
\newcommand{\Lshuffle}{*}

\newcommand{\LCoDendr}{\Coproduct_\LDendr}
\newcommand{\RCoDendr}{\Coproduct_\RDendr}

\newcommand{\Statistics}{\mathrm{s}}
\newcommand{\Biproduct}{\square}
\newcommand{\Convolution}{*}
\newcommand{\Tit}{T}
\newcommand{\PvH}{\mathsf{H}}
\newcommand{\OvH}{H}
\newcommand{\OvP}{\mathsf{R}}
\newcommand{\MvP}{\mathsf{B}}
\newcommand{\Dec}{\mathrm{dec}}
\newcommand{\Reduced}{\mathrm{red}}
\newcommand{\Rev}{\mathrm{rev}}
\newcommand{\Sat}{\mathsf{S}}
\newcommand{\SuperShuffle}{\bullet}
\newcommand{\BinToPerm}{\mathrm{btp}}
\newcommand{\PermToBin}{\mathrm{ptb}}
\newcommand{\Bud}{\mathrm{Bud}}
\newcommand{\Compo}{\odot}
\newcommand{\Prune}{\mathfrak{pru}}
\newcommand{\Multi}{\mathrm{mu}}
\newcommand{\Mono}{\mathrm{mo}}
\newcommand{\Colors}{\mathfrak{col}}
\newcommand{\ColorTypes}{\mathfrak{colt}}
\newcommand{\Type}{\mathrm{type}}
\newcommand{\CFG}{\mathrm{CFG}}
\newcommand{\RTG}{\mathrm{RTG}}
\newcommand{\SG}{\mathrm{SG}}
\newcommand{\PrecompToOp}{\mathrm{PO}}
\newcommand{\Hadamard}{\wasylozenge}
\newcommand{\LangDescr}{\mathrm{Lang}}
\newcommand{\Algo}{\mathrm{alg}}
\newcommand{\HookLength}{\mathrm{hk}}
\newcommand{\Index}{\mathrm{ind}}
\newcommand{\DPlus}{\dot{+}}

\newcommand{\SymmetricGroup}{\mathfrak{S}}
\newcommand{\Compositions}{\mathrm{Comp}}
\newcommand{\IntPart}{\mathrm{Part}}
\newcommand{\BinaryTrees}{\mathrm{BT}}
\newcommand{\RootedTrees}{\mathrm{RT}}
\newcommand{\PlanarRootedTrees}{\mathrm{PRT}}
\newcommand{\ColoredPlanarRootedTrees}{\mathrm{CPRT}}
\newcommand{\Ladders}{\mathrm{Lad}}
\newcommand{\Corollas}{\mathrm{Cor}}
\newcommand{\Ary}{\mathrm{Ary}}
\newcommand{\Schroder}{\mathrm{Sch}}
\newcommand{\GenDias}{\GeneratingSet_{\Dias_\gamma}}
\newcommand{\GenDendr}{\GeneratingSet_{\Dendr_\gamma}}
\newcommand{\GenTDendr}{\GeneratingSet_{\TDendr_\gamma}}
\newcommand{\GenAs}{\GeneratingSet_{\As_\gamma}}
\newcommand{\GenDAs}{\GeneratingSet_{\DAs_\gamma}}
\newcommand{\GenDup}{\GeneratingSet_{\Dup_\gamma}}
\newcommand{\GenTrias}{\GeneratingSet_{\Trias_\gamma}}
\newcommand{\RelDias}{\RelationSpace_{\Dias_\gamma}}
\newcommand{\RelDendr}{\RelationSpace_{\Dendr_\gamma}}
\newcommand{\RelAs}{\RelationSpace_{\As_\gamma}}
\newcommand{\RelDAs}{\RelationSpace_{\DAs_\gamma}}
\newcommand{\RelDup}{\RelationSpace_{\Dup_\gamma}}
\newcommand{\RelTrias}{\RelationSpace_{\Trias_\gamma}}
\newcommand{\RelTDendr}{\RelationSpace_{\TDendr_\gamma}}
\newcommand{\SetASchr}{\Sca}
\newcommand{\ASchr}{\mathsf{ASchr}}
\newcommand{\Arcs}{\mathcal{A}}
\newcommand{\Diagonals}{\mathcal{D}}
\newcommand{\Edges}{\mathcal{E}}
\newcommand{\Cliques}{\mathcal{C}}
\newcommand{\Bubbles}{\mathcal{B}}
\newcommand{\Triangles}{\Tca}
\newcommand{\Primes}{\Pca}
\newcommand{\PackedMatrices}{\Pca}

\newcommand{\SeriesTotPrimitiveElements}{\Tca}
\newcommand{\UBPCollection}{\mathrm{UBP}}
\newcommand{\Matrices}{\mathrm{M}}
\newcommand{\AdmissibleCuts}{\mathrm{Adm}}
\newcommand{\Types}{\Tca}

\newcommand{\Preprographs}{\mathrm{PPrg}}
\newcommand{\Prographs}{\mathrm{Prg}}
\newcommand{\Wires}{\mathrm{Wir}}

\newcommand{\Sym}{\mathsf{Sym}}
\newcommand{\SymCom}{Sym}
\newcommand{\FdB}{FdB}
\newcommand{\FdBNC}{\mathsf{FdB}}
\newcommand{\FQSym}{\mathsf{FQSym}}
\newcommand{\FSym}{\mathsf{FSym}}
\newcommand{\PBT}{\mathsf{PBT}}
\newcommand{\Bell}{\mathsf{Bell}}
\newcommand{\Baxter}{\mathsf{Baxter}}
\newcommand{\WQSym}{\mathsf{WQSym}}
\newcommand{\PQSym}{\mathsf{PQSym}}
\newcommand{\Camb}{\mathsf{Camb}}
\newcommand{\CK}{\mathsf{CK}}
\newcommand{\As}{\mathsf{As}}
\newcommand{\Dup}{\mathsf{Dup}}
\newcommand{\Dendr}{\mathsf{Dendr}}
\newcommand{\Com}{\mathsf{Com}}
\newcommand{\PreLie}{\mathsf{PreLie}}

\newcommand{\Per}{\mathsf{Per}}
\newcommand{\Mag}{\mathsf{Mag}}
\newcommand{\PerZero}{\mathsf{Per}_0}
\newcommand{\End}{\mathsf{End}}
\newcommand{\PF}{\mathsf{PF}}
\newcommand{\PW}{\mathsf{PW}}
\newcommand{\PRT}{\mathsf{PRT}}
\newcommand{\FCat}[1]{\mathsf{FCat}^{(#1)}}
\newcommand{\Schr}{\mathsf{Schr}}
\newcommand{\Motz}{\mathsf{Motz}}
\newcommand{\Comp}{\mathsf{Comp}}
\newcommand{\DA}{\mathsf{DA}}
\newcommand{\SComp}{\mathsf{SComp}}
\newcommand{\Di}{\mathsf{Di}}
\newcommand{\Dias}{\mathsf{Dias}}
\newcommand{\Tr}{\mathsf{Tr}}
\newcommand{\Trias}{\mathsf{Trias}}
\newcommand{\TDendr}{\mathsf{TDendr}}

\newcommand{\NAP}{\mathsf{NAP}}
\newcommand{\BNC}{\mathsf{BNC}}
\newcommand{\Bubble}{{\mathsf{Bubble}}}
\newcommand{\FAs}{{\mathsf{FAs}}}
\newcommand{\LOp}{\mathsf{L}}
\newcommand{\NCT}{\mathsf{NCT}}
\newcommand{\NCP}{\mathsf{NCP}}
\newcommand{\AlgPos}{\mathsf{Pos}}
\newcommand{\AlgSets}{\mathsf{Sets}}
\newcommand{\AlgWords}{\mathsf{Words}}
\newcommand{\AlgComWords}{\mathsf{CWords}}
\newcommand{\AlgMarkWords}{\mathsf{MWords}}
\newcommand{\DAs}{\mathsf{DAs}}
\newcommand{\Zin}{\mathsf{Zin}}
\newcommand{\Leib}{\mathsf{Leib}}
\newcommand{\Lie}{\mathsf{Lie}}
\newcommand{\TwoAs}{\mathsf{2as}}
\newcommand{\DupDendr}{\mathsf{D}}
\newcommand{\MT}{\mathsf{MT}}
\newcommand{\DMT}{\mathsf{DMT}}
\newcommand{\Grav}{\mathsf{Grav}}
\newcommand{\Lab}{\mathsf{Lab}}
\newcommand{\Bub}{\mathsf{Bub}}
\newcommand{\Deg}{\mathsf{Deg}}
\newcommand{\Cro}{\mathsf{Cro}}
\newcommand{\Acy}{\mathsf{Acy}}
\newcommand{\Whi}{\mathrm{Whi}}
\newcommand{\NC}{\mathrm{NC}}
\newcommand{\Nes}{\mathrm{Nes}}
\newcommand{\Paths}{\mathrm{Pat}}
\newcommand{\Forests}{\mathrm{For}}
\newcommand{\Motzkin}{\mathrm{Mot}}
\newcommand{\Diss}{\mathrm{Dis}}
\newcommand{\WNC}{\mathrm{WNC}}
\newcommand{\Luc}{\mathsf{Luc}}
\newcommand{\FF}{\mathcal{F}\mathcal{F}}
\newcommand{\RatFct}{\mathsf{RatFct}}
\newcommand{\Mould}{\mathsf{Mould}}
\newcommand{\PM}[1]{{\mathsf{PM}_{#1}}}
\newcommand{\PMN}[1]{{\mathsf{PMN}_{#1}}}
\newcommand{\PML}[1]{{\mathsf{PML}_{#1}}}
\newcommand{\MQSym}{\mathsf{MQSym}}
\newcommand{\UBP}{\mathsf{UBP}}
\newcommand{\ASM}{\mathsf{ASM}}
\newcommand{\QASM}[1]{\ASM/_{#1}}
\newcommand{\AB}{\mathsf{AB}}
\newcommand{\Heap}{\mathsf{Heap}}
\newcommand{\FHeap}{\mathsf{FHeap}}
\newcommand{\FBT}{\mathsf{FBT}}
\newcommand{\BAs}{\mathsf{BAs}}
\newcommand{\PRF}{\mathsf{PRF}}
\newcommand{\Tree}{\mathsf{Tree}}
\newcommand{\HookSystem}{\mathrm{HS}}
\newcommand{\MonoidPrecomp}{\Mca_{\mathrm{p}}}
\newcommand{\MonoidPrecompBar}{\bar{\MonoidPrecomp}}

\newcommand{\InfMat}{\bar \Mca_\infty}
\newcommand{\InfMatP}{\bar \Pca_\infty}
\newcommand{\InfMatN}{\Mca_\infty}
\newcommand{\InfMatNP}{\Pca_\infty}
\newcommand{\Precomp}{\mathsf{Precomp}}
\newcommand{\Poset}{\mathsf{Poset}}
\newcommand{\Qoset}{\mathsf{Qoset}}
\newcommand{\Map}{\mathsf{Map}}
\newcommand{\NDMap}{\mathsf{NDMap}}
\newcommand{\CMag}{\mathrm{CMag}}

\newcommand{\BasisB}{\mathsf{B}}
\newcommand{\BasisE}{\mathsf{E}}
\newcommand{\BasisF}{\mathsf{F}}

\newcommand{\BasisH}{\mathsf{H}}
\newcommand{\BasisK}{\mathsf{K}}
\newcommand{\BasisM}{\mathsf{M}}
\newcommand{\BasisP}{\mathsf{P}}
\newcommand{\BasisR}{\mathsf{R}}
\newcommand{\BasisS}{\mathsf{S}}
\newcommand{\BasisT}{\mathsf{T}}
\newcommand{\BasisV}{\mathsf{V}}
\newcommand{\BasisW}{\mathsf{W}}

\newcommand{\ProblemSD}{\mathsf{SDP}}
\newcommand{\ProblemRC}{\mathsf{RCP}}

\newcommand{\PClass}{{\bf P}\xspace}
\newcommand{\NP}{{\bf NP}\xspace}
\newcommand{\NPC}{\mbox{\NP-complete\xspace}}

\newcommand{\LeafPic}{%
    \raisebox{.39em}{%
    \begin{tikzpicture}[yscale=.28,Centering]
        \node[Leaf](0)at(0,0){};
        \node(r)at(0,1){};
        \draw[Edge](r)--(0);
    \end{tikzpicture}}}

\newcommand{\NodePic}{%
    \raisebox{.39em}{%
    \begin{tikzpicture}[yscale=.25,Centering]
        \node[Node,minimum size=1.5mm](0)at(0,0){};
        \node(r)at(0,1){};
        \draw[Edge](r)--(0);
    \end{tikzpicture}}}

\newcommand{\EdgePic}{%
    \begin{tikzpicture}[xscale=.25,Centering]
        \draw[Edge](0,0)--(1,0);
    \end{tikzpicture}}

\newcommand{\BinaryNode}{
    \raisebox{.25em}{%
    \begin{tikzpicture}[xscale=.4,yscale=.3,Centering]
        \node[Leaf](0)at(0.50,-.75){};
        \node[Leaf](2)at(1.5,-.75){};
        \node[Node](1)at(1.00,0.00){};
        \draw[Edge](0)--(1);
        \draw[Edge](2)--(1);
        \node(r)at(1.00,.95){};
        \draw[Edge](r)--(1);
    \end{tikzpicture}}}

\newcommand{\CorollaTwo}[1]{
    \begin{tikzpicture}[xscale=.13,yscale=.12,Centering]
        \node[Leaf](0)at(0.00,-1.50){};
        \node[Leaf](2)at(2.00,-1.50){};
        \node[Node](1)at(1.00,0.00){\begin{math}#1\end{math}};
        \draw[Edge](0)--(1);
        \draw[Edge](2)--(1);
        \node(r)at(1.00,2.5){};
        \draw[Edge](r)--(1);
    \end{tikzpicture}}

\newcommand{\ColoredBinaryNode}[3]{
    \begin{tikzpicture}[xscale=.5,yscale=.4,Centering]
        \node[Leaf](0)at(0.50,-.75){};
        \node[Leaf](2)at(1.5,-.75){};
        \node[Node](1)at(1.00,0.00){};
        \draw[Edge](0)--(1);
        \draw[Edge](2)--(1);
        \node(r)at(1.00,.85){};
        \draw[Edge](r)--(1);
        \node[above of=r,node distance=1mm,font=\footnotesize]
            {\begin{math}#1\end{math}};
        \node[below of=0,node distance=3mm,font=\footnotesize]
            {\begin{math}#2\end{math}};
        \node[below of=2,node distance=3mm,font=\footnotesize]
            {\begin{math}#3\end{math}};
    \end{tikzpicture}}

\newcommand{\BinaryTreeLeft}{
    \begin{tikzpicture}[scale=.17,Centering]
        \node[Leaf](0)at(0.00,-3.33){};
        \node[Leaf](2)at(2.00,-3.33){};
        \node[Leaf](4)at(4.00,-1.67){};
        \node[Node](1)at(1.00,-1.67){};
        \node[Node](3)at(3.00,0.00){};
        \draw[Edge](0)--(1);
        \draw[Edge](1)--(3);
        \draw[Edge](2)--(1);
        \draw[Edge](4)--(3);
        \node(r)at(3.00,1.5){};
        \draw[Edge](r)--(3);
    \end{tikzpicture}}

\newcommand{\BinaryTreeRight}{
    \begin{tikzpicture}[scale=.17,Centering]
        \node[Leaf](0)at(0.00,-1.67){};
        \node[Leaf](2)at(2.00,-3.33){};
        \node[Leaf](4)at(4.00,-3.33){};
        \node[Node](1)at(1.00,0.00){};
        \node[Node](3)at(3.00,-1.67){};
        \draw[Edge](0)--(1);
        \draw[Edge](2)--(3);
        \draw[Edge](3)--(1);
        \draw[Edge](4)--(3);
        \node(r)at(1.00,1.5){};
        \draw[Edge](r)--(1);
    \end{tikzpicture}}

\newcommand{\BinaryRoot}[2]{%
    \begin{tikzpicture}[xscale=.4,yscale=.3,Centering]
        \node(0)at(-.50,-1.50){\begin{math}#1\end{math}};
        \node(2)at(2.50,-1.50){\begin{math}#2\end{math}};
        \node[Node](1)at(1.00,.50){};
        \draw[Edge](0)--(1);
        \draw[Edge](2)--(1);
        \node(r)at(1.00,1.75){};
        \draw[Edge](r)--(1);
    \end{tikzpicture}}

\newcommand{\LeafST}{%
    \begin{tikzpicture}[yscale=.35,Centering]
        \node(0)at(0,0){};
        \node(r)at(0,1){};
        \draw[Edge](r)--(0);
    \end{tikzpicture}}

\newcommand{\STLeft}[2]{%
    \begin{tikzpicture}[xscale=.27,yscale=.25,Centering]
        \node(0)at(0.00,-3.33){};
        \node(2)at(2.00,-3.33){};
        \node(4)at(4.00,-1.67){};
        \node[NodeST](1)at(1.00,-1.67){\begin{math}#1\end{math}};
        \node[NodeST](3)at(3.00,0.00){\begin{math}#2\end{math}};
        \draw[Edge](0)--(1);
        \draw[Edge](1)--(3);
        \draw[Edge](2)--(1);
        \draw[Edge](4)--(3);
        \node(r)at(3.00,1.75){};
        \draw[Edge](r)--(3);
    \end{tikzpicture}}

\newcommand{\STRight}[2]{%
    \begin{tikzpicture}[xscale=.27,yscale=.25,Centering]
        \node(0)at(0.00,-1.67){};
        \node(2)at(2.00,-3.33){};
        \node(4)at(4.00,-3.33){};
        \node[NodeST](1)at(1.00,0.00){\begin{math}#1\end{math}};
        \node[NodeST](3)at(3.00,-1.67){\begin{math}#2\end{math}};
        \draw[Edge](0)--(1);
        \draw[Edge](2)--(3);
        \draw[Edge](3)--(1);
        \draw[Edge](4)--(3);
        \node(r)at(1.00,1.5){};
        \draw[Edge](r)--(1);
    \end{tikzpicture}}

\newcommand{\AAA}{
\begin{tikzpicture}[scale=.175,Centering]
    \node[CliquePoint](C1)at(-0.87,-0.50){};
    \node[CliquePoint](C2)at(0.00,1.0){};
    \node[CliquePoint](C3)at(0.87,-0.50){};
    \draw[PolyEdgeGray](C1)--(C3);
    \draw[PolyEdgeGray](C1)--(C2);
    \draw[PolyEdgeGray](C2)--(C3);
    \draw[PolyEdgeBlue](C1)--(C2);
    \draw[PolyEdgeBlue](C1)--(C3);
    \draw[PolyEdgeBlue](C2)--(C3);
\end{tikzpicture}}

\newcommand{\AAB}{
\begin{tikzpicture}[scale=.175,Centering]
    \node[CliquePoint](C1)at(-0.87,-0.50){};
    \node[CliquePoint](C2)at(0.00,1.0){};
    \node[CliquePoint](C3)at(0.87,-0.50){};
    \draw[PolyEdgeGray](C1)--(C3);
    \draw[PolyEdgeGray](C1)--(C2);
    \draw[PolyEdgeGray](C2)--(C3);
    \draw[PolyEdgeBlue](C1)--(C2);
    \draw[PolyEdgeBlue](C1)--(C3);
\end{tikzpicture}}

\newcommand{\ABA}{
\begin{tikzpicture}[scale=.175,Centering]
    \node[CliquePoint](C1)at(-0.87,-0.50){};
    \node[CliquePoint](C2)at(0.00,1.0){};
    \node[CliquePoint](C3)at(0.87,-0.50){};
    \draw[PolyEdgeGray](C1)--(C3);
    \draw[PolyEdgeGray](C1)--(C2);
    \draw[PolyEdgeGray](C2)--(C3);
    \draw[PolyEdgeBlue](C1)--(C2);
    \draw[PolyEdgeBlue](C2)--(C3);
\end{tikzpicture}}

\newcommand{\ABB}{
\begin{tikzpicture}[scale=.175,Centering]
    \node[CliquePoint](C1)at(-0.87,-0.50){};
    \node[CliquePoint](C2)at(0.00,1.0){};
    \node[CliquePoint](C3)at(0.87,-0.50){};
    \draw[PolyEdgeGray](C1)--(C3);
    \draw[PolyEdgeGray](C1)--(C2);
    \draw[PolyEdgeGray](C2)--(C3);
    \draw[PolyEdgeBlue](C1)--(C2);
\end{tikzpicture}}

\newcommand{\BAA}{
\begin{tikzpicture}[scale=.175,Centering]
    \node[CliquePoint](C1)at(-0.87,-0.50){};
    \node[CliquePoint](C2)at(0.00,1.0){};
    \node[CliquePoint](C3)at(0.87,-0.50){};
    \draw[PolyEdgeGray](C1)--(C3);
    \draw[PolyEdgeGray](C1)--(C2);
    \draw[PolyEdgeGray](C2)--(C3);
    \draw[PolyEdgeBlue](C1)--(C3);
    \draw[PolyEdgeBlue](C2)--(C3);
\end{tikzpicture}}

\newcommand{\BAB}{
\begin{tikzpicture}[scale=.175,Centering]
    \node[CliquePoint](C1)at(-0.87,-0.50){};
    \node[CliquePoint](C2)at(0.00,1.0){};
    \node[CliquePoint](C3)at(0.87,-0.50){};
    \draw[PolyEdgeGray](C1)--(C3);
    \draw[PolyEdgeGray](C1)--(C2);
    \draw[PolyEdgeGray](C2)--(C3);
    \draw[PolyEdgeBlue](C1)--(C3);
\end{tikzpicture}}

\newcommand{\BBA}{
\begin{tikzpicture}[scale=.175,Centering]
    \node[CliquePoint](C1)at(-0.87,-0.50){};
    \node[CliquePoint](C2)at(0.00,1.0){};
    \node[CliquePoint](C3)at(0.87,-0.50){};
    \draw[PolyEdgeGray](C1)--(C3);
    \draw[PolyEdgeGray](C1)--(C2);
    \draw[PolyEdgeGray](C2)--(C3);
    \draw[PolyEdgeBlue](C2)--(C3);
\end{tikzpicture}}

\newcommand{\BBB}{
\begin{tikzpicture}[scale=.175,Centering]
    \node[CliquePoint](C1)at(-0.87,-0.50){};
    \node[CliquePoint](C2)at(0.00,1.0){};
    \node[CliquePoint](C3)at(0.87,-0.50){};
    \draw[PolyEdgeGray](C1)--(C3);
    \draw[PolyEdgeGray](C1)--(C2);
    \draw[PolyEdgeGray](C2)--(C3);
\end{tikzpicture}}

\newcommand{\BinTreeLTwo}[1]{
\begin{tikzpicture}[yscale=.25,xscale=.25,Centering]
    \node[Leaf](0)at(0.00,-3.33){};
    \node[Leaf](2)at(2.00,-3.33){};
    \node[Leaf](4)at(4.00,-1.67){};
    \node[Node](1)at(1.00,-1.67){};
    \node[Node](3)at(3.00,0.00){};
    \draw[Edge](0)--(1);
    \draw[Edge](1)edge[]node[EdgeValue]
            {\begin{math}#1\end{math}}(3);
    \draw[Edge](2)--(1);
    \draw[Edge](4)--(3);
    \node(r)at(3.00,1.25){};
    \draw[Edge](r)--(3);
\end{tikzpicture}}

\newcommand{\BinTreeRTwo}[1]{
\begin{tikzpicture}[yscale=.25,xscale=.25,Centering]
    \node[Leaf](0)at(0.00,-1.67){};
    \node[Leaf](2)at(2.00,-3.33){};
    \node[Leaf](4)at(4.00,-3.33){};
    \node[Node](1)at(1.00,0.00){};
    \node[Node](3)at(3.00,-1.67){};
    \draw[Edge](0)--(1);
    \draw[Edge](2)--(3);
    \draw[Edge](3)edge[]node[EdgeValue]
            {\begin{math}#1\end{math}}(1);
    \draw[Edge](4)--(3);
    \node(r)at(1.00,1.25){};
    \draw[Edge](r)--(1);
\end{tikzpicture}}

\newcommand{\ValuedBinTree}[4]{
\begin{tikzpicture}[xscale=.5,yscale=.4,Centering]
    \node(0)at(-.50,-1.50){\begin{math}#3\end{math}};
    \node(2)at(2.50,-1.50){\begin{math}#4\end{math}};
    \node[Node](1)at(1.00,.50){};
    \draw[Edge](0)edge[]node[EdgeValue]
        {\begin{math}#1\end{math}}(1);
    \draw[Edge](2)edge[]node[EdgeValue]
        {\begin{math}#2\end{math}}(1);
    \node(r)at(1.00,1.5){};
    \draw[Edge](r)--(1);
\end{tikzpicture}}

\newcommand{\FreeCorollaOne}[1]{
\begin{tikzpicture}[xscale=.13,yscale=.25,Centering]
    \node(0)at(1.00,-1.50){};
    \node[NodeST](1)at(1.00,0.00){\begin{math}#1\end{math}};
    \draw[Edge](0)--(1);
    \node(r)at(1.00,1.5){};
    \draw[Edge](r)--(1);
\end{tikzpicture}}

\newcommand{\AntiForestPattern}{
\begin{tikzpicture}[xscale=.2,yscale=.25,Centering]
    \node[PosetVertex](1)at(0,0){};
    \node[PosetVertex](2)at(-1,1){};
    \node[PosetVertex](3)at(1,1){};
    \draw[Edge](1)--(2);
    \draw[Edge](1)--(3);
\end{tikzpicture}}

\newcommand{\LinePattern}{
\begin{tikzpicture}[yscale=.28,Centering]
    \node[PosetVertex](1)at(0,0){};
    \node[PosetVertex](2)at(0,-1){};
    \draw[Edge](1)--(2);
\end{tikzpicture}}

\newcommand{\OnePoset}{
\begin{tikzpicture}[Centering]
    \node[PosetVertex](1)at(0,0){};
\end{tikzpicture}}

\newcommand{\UnitClique}{
\begin{tikzpicture}[scale=.6,Centering]
    \node[CliquePoint](1)at(0,0){};
    \node[CliquePoint](2)at(.75,0){};
    \draw[CliqueEmptyEdge](1)edge[]node[CliqueLabel]{}(2);
\end{tikzpicture}}

\newcommand{\SquareMotz}{
\begin{tikzpicture}[scale=.6,Centering]
    \node[CliquePoint](1)at(-0.71,0.71){};
    \node[CliquePoint](2)at(0.71,0.71){};
    \node[CliquePoint](3)at(0.71,-0.71){};
    \node[CliquePoint](4)at(-0.71,-0.71){};
    \draw[CliqueEmptyEdge](2)edge[]node[CliqueLabel]{}(3);
    \draw[CliqueEmptyEdge](1)edge[]node[CliqueLabel]{}(4);
    \draw[CliqueEmptyEdge](3)edge[]node[CliqueLabel]{}(4);
    \draw[CliqueEdge](1)edge[]node[CliqueLabel]
        {\begin{math}0\end{math}}(2);
\end{tikzpicture}}

\newcommand{\Triangle}[3]{
\begin{tikzpicture}[scale=.42,Centering]
    \node[CliquePoint](1)at(0,1){};
    \node[CliquePoint](2)at(0.87,-0.5){};
    \node[CliquePoint](3)at(-0.87,-0.5){};
    \draw[CliqueEdge](1)edge[]node[CliqueLabel]
        {\begin{math}#3\end{math}}(2);
    \draw[CliqueEdge](1)edge[]node[CliqueLabel]
        {\begin{math}#2\end{math}}(3);
    \draw[CliqueEdge](2)edge[]node[CliqueLabel]
        {\begin{math}#1\end{math}}(3);
\end{tikzpicture}}

\newcommand{\TriangleEXX}[3]{
\begin{tikzpicture}[scale=.42,Centering]
    \node[CliquePoint](1)at(0,1){};
    \node[CliquePoint](2)at(0.87,-0.5){};
    \node[CliquePoint](3)at(-0.87,-0.5){};
    \draw[CliqueEdge](1)edge[]node[CliqueLabel]
        {\begin{math}#3\end{math}}(2);
    \draw[CliqueEdge](1)edge[]node[CliqueLabel]
        {\begin{math}#2\end{math}}(3);
    \draw[CliqueEmptyEdge](2)edge[]node[CliqueLabel]{}(3);
\end{tikzpicture}}

\newcommand{\TriangleXEX}[3]{
\begin{tikzpicture}[scale=.42,Centering]
    \node[CliquePoint](1)at(0,1){};
    \node[CliquePoint](2)at(0.87,-0.5){};
    \node[CliquePoint](3)at(-0.87,-0.5){};
    \draw[CliqueEdge](1)edge[]node[CliqueLabel]
        {\begin{math}#3\end{math}}(2);
    \draw[CliqueEmptyEdge](1)edge[]node[CliqueLabel]{}(3);
    \draw[CliqueEdge](2)edge[]node[CliqueLabel]
        {\begin{math}#1\end{math}}(3);
\end{tikzpicture}}

\newcommand{\TriangleXXE}[3]{
\begin{tikzpicture}[scale=.42,Centering]
    \node[CliquePoint](1)at(0,1){};
    \node[CliquePoint](2)at(0.87,-0.5){};
    \node[CliquePoint](3)at(-0.87,-0.5){};
    \draw[CliqueEmptyEdge](1)edge[]node[CliqueLabel]{}(2);
    \draw[CliqueEdge](1)edge[]node[CliqueLabel]
        {\begin{math}#2\end{math}}(3);
    \draw[CliqueEdge](2)edge[]node[CliqueLabel]
        {\begin{math}#1\end{math}}(3);
\end{tikzpicture}}

\newcommand{\TriangleXEE}[3]{
\begin{tikzpicture}[scale=.42,Centering]
    \node[CliquePoint](1)at(0,1){};
    \node[CliquePoint](2)at(0.87,-0.5){};
    \node[CliquePoint](3)at(-0.87,-0.5){};
    \draw[CliqueEmptyEdge](1)edge[]node[CliqueLabel]{}(2);
    \draw[CliqueEmptyEdge](1)edge[]node[CliqueLabel]{}(3);
    \draw[CliqueEdge](2)edge[]node[CliqueLabel]
        {\begin{math}#1\end{math}}(3);
\end{tikzpicture}}

\newcommand{\TriangleEEX}[3]{
\begin{tikzpicture}[scale=.42,Centering]
    \node[CliquePoint](1)at(0,1){};
    \node[CliquePoint](2)at(0.87,-0.5){};
    \node[CliquePoint](3)at(-0.87,-0.5){};
    \draw[CliqueEdge](1)edge[]node[CliqueLabel]
        {\begin{math}#3\end{math}}(2);
    \draw[CliqueEmptyEdge](1)edge[]node[CliqueLabel]{}(3);
    \draw[CliqueEmptyEdge](2)edge[]node[CliqueLabel]{}(3);
\end{tikzpicture}}

\newcommand{\TriangleEXE}[3]{
\begin{tikzpicture}[scale=.42,Centering]
    \node[CliquePoint](1)at(0,1){};
    \node[CliquePoint](2)at(0.87,-0.5){};
    \node[CliquePoint](3)at(-0.87,-0.5){};
    \draw[CliqueEmptyEdge](1)edge[]node[CliqueLabel]{}(2);
    \draw[CliqueEdge](1)edge[]node[CliqueLabel]
        {\begin{math}#2\end{math}}(3);
    \draw[CliqueEmptyEdge](2)edge[]node[CliqueLabel]{}(3);
\end{tikzpicture}}

\newcommand{\TriangleEEE}[3]{
\begin{tikzpicture}[scale=.42,Centering]
    \node[CliquePoint](1)at(0,1){};
    \node[CliquePoint](2)at(0.87,-0.5){};
    \node[CliquePoint](3)at(-0.87,-0.5){};
    \draw[CliqueEmptyEdge](1)edge[]node[CliqueLabel]{}(2);
    \draw[CliqueEmptyEdge](1)edge[]node[CliqueLabel]{}(3);
    \draw[CliqueEmptyEdge](2)edge[]node[CliqueLabel]{}(3);
\end{tikzpicture}}

\newcommand{\TriangleOp}[3]{\;
\begin{tikzpicture}[scale=.35,Centering]
    \node[shape=coordinate](1)at(0,1){};
    \node[shape=coordinate](2)at(0.87,-0.5){};
    \node[shape=coordinate](3)at(-0.87,-0.5){};
    \draw[draw=Col6!90](1)edge[]node[CliqueLabel,font=\tiny]
        {\begin{math}#3\end{math}}(2);
    \draw[draw=Col6!90](1)edge[]node[CliqueLabel,font=\tiny]
        {\begin{math}#2\end{math}}(3);
    \draw[draw=Col6!90](2)edge[]node[CliqueLabel,font=\tiny]
        {\begin{math}#1\end{math}}(3);
\end{tikzpicture}\;}

\newcommand{\Matrix}[1]{
\left[
    \begin{smallmatrix}
        #1
    \end{smallmatrix}
\right]}

\newcommand{\SixVertexNE}{
\begin{tikzpicture}[scale=.8,Centering]
    \def \bi{(-.1,.1)--(0,-.1)--(.1,.1)};
    \def \ci{(-.1,-.1)--(0,.1)--(.1,-.1)};
    \def \di{(-.1,-.1)--(.1,0)--(-.1,.1)};
    \def \ei{(.1,-.1)--(-.1,0)--(.1,.1)};
    \draw[SixVertexGrid] (-.5,-.5) grid (.5,.5);
    \draw[SixVertexArrow,shift={(0,-.5)}] \ci;
    \draw[SixVertexArrow,shift={(0,0.5)}] \ci;
    \draw[SixVertexArrow,shift={(-.5,0)}] \di;
    \draw[SixVertexArrow,shift={(0.5,0)}] \di;
\end{tikzpicture}}

\newcommand{\SixVertexSW}{
\begin{tikzpicture}[scale=.8,Centering]
    \def \bi{(-.1,.1)--(0,-.1)--(.1,.1)};
    \def \ci{(-.1,-.1)--(0,.1)--(.1,-.1)};
    \def \di{(-.1,-.1)--(.1,0)--(-.1,.1)};
    \def \ei{(.1,-.1)--(-.1,0)--(.1,.1)};
    \draw[SixVertexGrid] (-.5,-.5) grid (.5,.5);
    \draw[SixVertexArrow,shift={(0,-.5)}] \bi;
    \draw[SixVertexArrow,shift={(0,0.5)}] \bi;
    \draw[SixVertexArrow,shift={(-.5,0)}] \ei;
    \draw[SixVertexArrow,shift={(0.5,0)}] \ei;
\end{tikzpicture}}

\newcommand{\SixVertexSE}{
\begin{tikzpicture}[scale=.8,Centering]
    \def \bi{(-.1,.1)--(0,-.1)--(.1,.1)};
    \def \ci{(-.1,-.1)--(0,.1)--(.1,-.1)};
    \def \di{(-.1,-.1)--(.1,0)--(-.1,.1)};
    \def \ei{(.1,-.1)--(-.1,0)--(.1,.1)};
    \draw[SixVertexGrid] (-.5,-.5) grid (.5,.5);
    \draw[SixVertexArrow,shift={(0,-.5)}] \bi;
    \draw[SixVertexArrow,shift={(0,0.5)}] \bi;
    \draw[SixVertexArrow,shift={(-.5,0)}] \di;
    \draw[SixVertexArrow,shift={(0.5,0)}] \di;
\end{tikzpicture}}

\newcommand{\SixVertexNW}{
\begin{tikzpicture}[scale=.8,Centering]
    \def \bi{(-.1,.1)--(0,-.1)--(.1,.1)};
    \def \ci{(-.1,-.1)--(0,.1)--(.1,-.1)};
    \def \di{(-.1,-.1)--(.1,0)--(-.1,.1)};
    \def \ei{(.1,-.1)--(-.1,0)--(.1,.1)};
    \draw[SixVertexGrid] (-.5,-.5) grid (.5,.5);
    \draw[SixVertexArrow,shift={(0,-.5)}] \ci;
    \draw[SixVertexArrow,shift={(0,0.5)}] \ci;
    \draw[SixVertexArrow,shift={(-.5,0)}] \ei;
    \draw[SixVertexArrow,shift={(0.5,0)}] \ei;
\end{tikzpicture}}

\newcommand{\SixVertexOI}{
\begin{tikzpicture}[scale=.8,Centering]
    \def \bi{(-.1,.1)--(0,-.1)--(.1,.1)};
    \def \ci{(-.1,-.1)--(0,.1)--(.1,-.1)};
    \def \di{(-.1,-.1)--(.1,0)--(-.1,.1)};
    \def \ei{(.1,-.1)--(-.1,0)--(.1,.1)};
    \draw[SixVertexGrid] (-.5,-.5) grid (.5,.5);
    \draw[SixVertexArrow,shift={(0,-.5)}] \bi;
    \draw[SixVertexArrow,shift={(0,0.5)}] \ci;
    \draw[SixVertexArrow,shift={(-.5,0)}] \di;
    \draw[SixVertexArrow,shift={(0.5,0)}] \ei;
\end{tikzpicture}}

\newcommand{\SixVertexIO}{
\begin{tikzpicture}[scale=.8,Centering]
    \def \bi{(-.1,.1)--(0,-.1)--(.1,.1)};
    \def \ci{(-.1,-.1)--(0,.1)--(.1,-.1)};
    \def \di{(-.1,-.1)--(.1,0)--(-.1,.1)};
    \def \ei{(.1,-.1)--(-.1,0)--(.1,.1)};
    \draw[SixVertexGrid] (-.5,-.5) grid (.5,.5);
    \draw[SixVertexArrow,shift={(0,-.5)}] \ci;
    \draw[SixVertexArrow,shift={(0,0.5)}] \bi;
    \draw[SixVertexArrow,shift={(-.5,0)}] \ei;
    \draw[SixVertexArrow,shift={(0.5,0)}] \di;
\end{tikzpicture}}

\newcommand{\SixVertexW}{
\begin{tikzpicture}[scale=.3,Centering]
    \def \bi{(-.1,.1)--(0,-.1)--(.1,.1)};
    \def \ci{(-.1,-.1)--(0,.1)--(.1,-.1)};
    \def \di{(-.1,-.1)--(.1,0)--(-.1,.1)};
    \def \ei{(.1,-.1)--(-.1,0)--(.1,.1)};
    \draw[SixVertexGrid] (-.5,0) grid (.5,0);
    \draw[SixVertexArrow,shift={(-.5,0)}] \ei;
    \draw[SixVertexArrow,shift={(0.5,0)}] \ei;
\end{tikzpicture}}

\newcommand{\SixVertexE}{
\begin{tikzpicture}[scale=.3,Centering]
    \def \bi{(-.1,.1)--(0,-.1)--(.1,.1)};
    \def \ci{(-.1,-.1)--(0,.1)--(.1,-.1)};
    \def \di{(-.1,-.1)--(.1,0)--(-.1,.1)};
    \def \ei{(.1,-.1)--(-.1,0)--(.1,.1)};
    \draw[SixVertexGrid] (-.5,0) grid (.5,0);
    \draw[SixVertexArrow,shift={(-.5,0)}] \di;
    \draw[SixVertexArrow,shift={(0.5,0)}] \di;
\end{tikzpicture}}

\newcommand{\SixVertexS}{
\begin{tikzpicture}[scale=.3,Centering]
    \def \bi{(-.1,.1)--(0,-.1)--(.1,.1)};
    \def \ci{(-.1,-.1)--(0,.1)--(.1,-.1)};
    \def \di{(-.1,-.1)--(.1,0)--(-.1,.1)};
    \def \ei{(.1,-.1)--(-.1,0)--(.1,.1)};
    \draw[SixVertexGrid] (0,-.5) grid (0,.5);
    \draw[SixVertexArrow,shift={(0,-.5)}] \bi;
    \draw[SixVertexArrow,shift={(0,.5)}] \bi;
\end{tikzpicture}}

\newcommand{\SixVertexN}{
\begin{tikzpicture}[scale=.3,Centering]
    \def \bi{(-.1,.1)--(0,-.1)--(.1,.1)};
    \def \ci{(-.1,-.1)--(0,.1)--(.1,-.1)};
    \def \di{(-.1,-.1)--(.1,0)--(-.1,.1)};
    \def \ei{(.1,-.1)--(-.1,0)--(.1,.1)};
    \draw[SixVertexGrid] (0,-.5) grid (0,.5);
    \draw[SixVertexArrow,shift={(0,-.5)}] \ci;
    \draw[SixVertexArrow,shift={(0,.5)}] \ci;
\end{tikzpicture}}

\newcommand{\CrossingLR}{%
  \begin{tikzpicture}[yscale=0.07,xscale=0.07,
      node distance=.25cm,line width=1pt,Centering]
    \node[PatternMatchingVertex](U1){};
    \node[PatternMatchingVertex][right of=U1](U2){};
    \node[PatternMatchingVertex][right of=U2](U3){};
    \node[PatternMatchingVertex][right of=U3](U4){};
    \draw[DirectedMatchingEdge] (U1.north) .. controls
        ($ (U1.north) + (0,6) $)
        and ($ (U3.north) + (0,6) $) .. (U3.north);
    \draw[DirectedMatchingEdge] (U2.north) .. controls
    ($ (U2.north) + (0,6) $)
    and ($ (U4.north) + (0,6) $) .. (U4.north);
\end{tikzpicture}}

\newcommand{\CrossingRL}{%
  \begin{tikzpicture}[yscale=0.07,xscale=0.07,
      node distance=.25cm,line width=1pt,Centering]
    \node[PatternMatchingVertex] (U1)               {};
    \node[PatternMatchingVertex] [right of=U1] (U2) {};
    \node[PatternMatchingVertex] [right of=U2] (U3) {};
    \node[PatternMatchingVertex] [right of=U3] (U4) {};
    \draw[DirectedMatchingEdge] (U1.north) .. controls
        ($ (U1.north) + (0,6) $)
        and ($ (U3.north) + (0,6) $) .. (U3.north);
    \draw[DirectedMatchingEdge] (U4.north) .. controls
        ($ (U4.north) + (0,6) $)
        and ($ (U2.north) + (0,6) $) .. (U2.north);
\end{tikzpicture}}

\newcommand{\InclusionLL}{%
  \begin{tikzpicture}[yscale=0.07,xscale=0.07,
      node distance=.25cm,line width=1pt,Centering]
    \node[PatternMatchingVertex] (U1)               {};
    \node[PatternMatchingVertex] [right of=U1] (U2) {};
    \node[PatternMatchingVertex] [right of=U2] (U3) {};
    \node[PatternMatchingVertex] [right of=U3] (U4) {};
    \draw[DirectedMatchingEdge] (U4.north) .. controls
        ($ (U4.north) + (0,6) $)
        and ($ (U1.north) + (0,6) $) .. (U1.north);
    \draw[DirectedMatchingEdge] (U3.north) .. controls
        ($ (U3.north) + (0,5) $)
        and ($ (U2.north) + (0,5) $) .. (U2.north);
\end{tikzpicture}}

\newcommand{\InclusionRR}{%
  \begin{tikzpicture}[yscale=0.07,xscale=0.07,
      node distance=.25cm,line width=1pt,Centering]
    \node[PatternMatchingVertex] (U1)               {};
    \node[PatternMatchingVertex] [right of=U1] (U2) {};
    \node[PatternMatchingVertex] [right of=U2] (U3) {};
    \node[PatternMatchingVertex] [right of=U3] (U4) {};
    \draw[DirectedMatchingEdge] (U1.north) .. controls
        ($ (U1.north) + (0,6) $)
        and ($ (U4.north) + (0,6) $) .. (U4.north);
    \draw[DirectedMatchingEdge] (U2.north) .. controls
    ($ (U2.north) + (0,5) $)
    and ($ (U3.north) + (0,5) $) .. (U3.north);
\end{tikzpicture}}

\newcommand{\InclusionRL}{%
  \begin{tikzpicture}[yscale=0.07,xscale=0.07,
      node distance=.25cm,line width=1pt,Centering]
    \node[PatternMatchingVertex] (U1)               {};
    \node[PatternMatchingVertex] [right of=U1] (U2) {};
    \node[PatternMatchingVertex] [right of=U2] (U3) {};
    \node[PatternMatchingVertex] [right of=U3] (U4) {};
    \draw[DirectedMatchingEdge] (U1.north) .. controls
        ($ (U1.north) + (0,6) $)
        and ($ (U4.north) + (0,6) $) .. (U4.north);
    \draw [DirectedMatchingEdge] (U3.north) .. controls
        ($ (U3.north) + (0,5) $)
        and ($ (U2.north) + (0,5) $) .. (U2.north);
\end{tikzpicture}}

\newcommand{\InclusionLR}{%
  \begin{tikzpicture}[yscale=0.07,xscale=0.07,
      node distance=.25cm,line width=1pt,Centering]
    \node[PatternMatchingVertex] (U1)               {};
    \node[PatternMatchingVertex] [right of=U1] (U2) {};
    \node[PatternMatchingVertex] [right of=U2] (U3) {};
    \node[PatternMatchingVertex] [right of=U3] (U4) {};
    \draw[DirectedMatchingEdge] (U4.north) .. controls
        ($ (U4.north) + (0,6) $)
        and ($ (U1.north) + (0,6) $) .. (U1.north);
    \draw[DirectedMatchingEdge] (U2.north) .. controls
        ($ (U2.north) + (0,5) $)
        and ($ (U3.north) + (0,5) $) .. (U3.north);
\end{tikzpicture}}

\newcommand{\MotzPoint}{
\begin{tikzpicture}[scale=.2,Centering]
    \draw[Grid](0,0) grid (0,0);
    \node[PathNode](1)at(0,0){};
\end{tikzpicture}}

\newcommand{\MotzHoriz}{
\begin{tikzpicture}[scale=.3,Centering]
    \draw[Grid](0,0) grid (1,0);
    \node[PathNode](Motz1)at(0,0){};
    \node[PathNode](Motz2)at(1,0){};
    \draw[PathStep](Motz1)--(Motz2);
\end{tikzpicture}}

\newcommand{\MotzPeak}{
\begin{tikzpicture}[scale=.22,Centering]
    \draw[Grid](0,0) grid (2,1);
    \node[PathNode](Motz1)at(0,0){};
    \node[PathNode](Motz2)at(1,1){};
    \node[PathNode](Motz3)at(2,0){};
    \draw[PathStep](Motz1)--(Motz2);
    \draw[PathStep](Motz2)--(Motz3);
\end{tikzpicture}}

\newcommand{\DyckUp}{
\begin{tikzpicture}[scale=.24,Centering]
    \draw[Grid](0,0) grid (1,1);
    \node[PathNode](Dyck1)at(0,0){};
    \node[PathNode](Dyck2)at(1,1){};
    \draw[PathStep](Dyck1)--(Dyck2);
\end{tikzpicture}}

\newcommand{\DyckDown}{
\begin{tikzpicture}[scale=.24,Centering]
    \draw[Grid](0,0) grid (1,1);
    \node[PathNode](Dyck1)at(0,1){};
    \node[PathNode](Dyck2)at(1,0){};
    \draw[PathStep](Dyck1)--(Dyck2);
\end{tikzpicture}}

\newcommand{\Cover}{
\begin{otherlanguage}{french}
\begin{normalfont}
\begin{Huge}
    {\sc Habilitation à diriger des recherches}
\end{Huge}
\vspace{1em}

\begin{large} \normalfont
    Spécialité Informatique
\end{large}
\vspace{1em}

\begin{large} \normalfont
    {\it Présentée le 4 décembre 2017 par}
\end{large}
\vspace{1em}

\begin{LARGE}
    {\bf Samuele Giraudo}
\end{LARGE}
\vspace{1em}

\begin{large} \normalfont
    Université Paris-Est Marne-la-Vallée
\end{large}
\vspace{3em}

\begin{Huge}
    {\bf Operads in algebraic combinatorics}
\end{Huge}
\begin{large}
    Opérades en combinatoire algébrique
\end{large}
\vspace{3em}

\begin{large} \normalfont
    {\it Devant le jury composé par}
\end{large}
\vspace{1em}

\begin{Large} \normalfont
    \begin{tabular}{ll}
        Frédéric Chapoton   \qquad & \qquad Président du jury \\[1em]
        Jean-Yves Thibon    \qquad & \qquad Garant d'habilitation
            \\[1em]
        Pierre-Louis Curien \qquad & \qquad Rapporteur \\
        Loïc Foissy         \qquad & \qquad Rapporteur \\
        Dominique Manchon   \qquad & \qquad Rapporteur
            \\[1em]
        Alessandra Frabetti \qquad & \qquad Examinateur \\
        Florent Hivert      \qquad & \qquad Examinateur \\
        Jean-Gabriel Luque  \qquad & \qquad Examinateur
    \end{tabular}
\end{Large}

\cleardoublepage

\end{normalfont}
\end{otherlanguage}}

\newcommand{\Keywords}{
    Algebraic combinatorics;
    Computer science;
    Formal power series;
    Tree;
    Rewrite system;
    Poset;
    Operad;
    Colored operad;
    Hopf bialgebra;
    Pro.
}

\newcommand{\Subjclass}{
    05A15, 
    05C05, 
    05E15, 
    16T05, 
    18D50, 
    68Q42. 
}

\newcommand{\Abstract}{
    \normalsize
    This habilitation thesis fits in the fields of algebraic and
    enumerative combinatorics, with connections with computer science.
    The main ideas developed in this work consist in endowing
    combinatorial objects (words, permutations, trees, integer
    partitions, Young tableaux, {\em etc.}) with operations in order to
    construct algebraic structures. This process allows, by studying
    algebraically the structures thus obtained (changes of bases,
    generating sets, presentations by generators and relations,
    morphisms, representations), to collect combinatorial information
    about the underlying objects. The algebraic structures the most
    encountered here are magmas, posets, associative algebras,
    dendriform algebras, Hopf bialgebras, operads, and pros.
    \smallbreak

    This work explores the aforementioned research direction and
    provides many (functorial or not) constructions having the
    particularity to build algebraic structures on combinatorial
    objects. We develop for instance a functor from nonsymmetric colored
    operads to nonsymmetric operads, from monoids to operads, from
    unitary magmas to nonsymmetric operads, from finite posets to
    nonsymmetric operads, from stiff pros to Hopf bialgebras, and from
    precompositions to nonsymmetric operads. These constructions bring
    alternative ways to describe already known structures and provide
    new ones, as for instance, some of the deformations of the
    noncommutative Faà di Bruno Hopf bialgebra of Foissy and a
    generalization of the dendriform operad of Loday.
    \smallbreak

    We also use algebraic structures to obtain enumerative results. In
    particular, nonsymmetric colored operads are promising devices to
    define formal series generalizing the usual ones. These series come
    with several products (for instance a pre-Lie product, an
    associative product, and their Kleene stars) enriching the usual
    ones on classical power series. This provides a framework and a
    toolbox to strike combinatorial questions in an original way.
    \smallbreak

    The text is organized as follows. The first two chapters pose the
    elementary notions of combinatorics and algebraic combinatorics
    used in the whole work. The last ten chapters contain our original
    research results fitting the context presented above.
    \medbreak
}

\newcommand{\BackCover}{%
\ClearToLeftPage

\thispagestyle{empty}

\begin{center} \LARGE
    {\bf Operads in algebraic combinatorics}
\end{center}
\begin{center}
    \begin{otherlanguage}{french}
        Opérades en combinatoire algébrique
    \end{otherlanguage}
\end{center}
\vspace{2em}

\begin{center} \Large
    {\bf Samuele Giraudo}
\end{center}
\vspace{2em}

{\it 2010 Mathematics Subject Classification.} \Subjclass
\bigbreak

{\it Key words and phrases.} \Keywords
\vspace{2em}

{\sc Abstract.} \Abstract}

\title[Operads in algebraic combinatorics]{\Cover}
\subjclass[2010]{\Subjclass}
\date{December 4, 2017}
\author{Samuele Giraudo}
\address{\scriptsize Université Paris-Est, LIGM (UMR $8049$), CNRS,
    ENPC, ESIEE Paris, UPEM, F-$77454$, Marne-la-Vallée, France}
\email{samuele.giraudo@u-pem.fr}

\numberwithin{equation}{subsection}

\counterwithin{figure}{chapter}
\counterwithin{table}{chapter}

\makeatletter
\def\l@part{\@tocline{-1}{20pt}{0pc}{5pc}{\Large \bf}}
\def\l@chapter{\@tocline{0}{10pt}{0pc}{5pc}{\bf}}
\def\l@section{\@tocline{1}{3pt}{1pc}{5pc}{}}
\def\l@subsection{\@tocline{2}{2pt}{2pc}{5pc}{}}
\makeatother

\renewcommand{\leq}{\leqslant}
\renewcommand{\geq}{\geqslant}

\newtheorem{Theorem}{Theorem}[subsection]
\newtheorem{Proposition}[Theorem]{Proposition}
\newtheorem{Lemma}[Theorem]{Lemma}

\begin{document}

\begin{abstract}
    \Abstract
\end{abstract}

\frontmatter

\maketitle

\begin{otherlanguage}{french}
\SkipTocEntry\chapter*{Remerciements}
Même si sur le papier je figure comme unique auteur de ce document, ce
travail n'a pu s'épanouir que grâce aux personnes que j'ai côtoyées.
Leurs apports, qu'ils soient scientifiques, culturels, spirituels ou
amicaux, sont tous à mes yeux d'une grande valeur. C'est à vous tous que
ces lignes sont dédiées.
\medbreak

Je souhaite tout d'abord remercier Jean-Yves Thibon d'avoir accepté ce
rôle de garant d'habilitation. Je le connais depuis que j'ai commencé
mon doctorat il y a huit ans. L'expérience montre qu'absolument tous
ses conseils et avis sont d'une importance capitale. L'étendue et la
profondeur de ses connaissances, ainsi que sa bienveillance, font que
l'avoir comme chef de l'équipe de combinatoire au LIGM est une chance
inestimable.
\medbreak

J'adresse maintenant mes plus sincères remerciements à l'ensemble des
rapporteurs. Je suis conscient du très lourd travail que représente
simplement la lecture de ce long manuscrit. Avoir un avis dessus et
l'évaluer est un travail encore plus exigeant. C'est donc à la fois un
honneur et une faveur immense qui m'est faite. J'ai rencontré
Pierre-Louis Curien en 2016 lors d'une conférence sur les opérades et
les systèmes de réécriture. Nous avons eu à cette occasion une très
agréable discussion à propos des opérades et des liens que cette
théorie entretient avec l'informatique théorique et les modèles de
calcul. Je suis très reconnaissant du travail de relecture attentif
de Pierre-Louis. Ses remarques ont beaucoup amélioré le texte. Les
travaux de Loïc Foissy sur les opérades, les structures dendriformes et
les bigèbres de toutes sortes m'ont beaucoup inspiré et vont sans
aucun doute le continuer pendant un très long moment. Je profite aussi
de cette opportunité pour le remercier pour ses réponses aux questions
que j'ai pu lui poser ces dernières années, que ce soit pendant les
conférences ou bien par courrier. J'ai assisté à un cours de Dominique
Manchon à Benasque en début 2010 sur les bigèbres de Hopf combinatoires
et les opérades. Je ne me doutais pas encore à ce moment là que ces
objets allaient être aussi importants pour moi. La pédagogie et le soin
avec lesquels il a présenté son cours ont beaucoup contribué à ma
passion alors naissante pour ces objets. Merci aussi à Dominique pour
ses remarques qui ont rendu ce document meilleur.
\medbreak

Je souhaite remercier maintenant l'ensemble des membres du jury. Vous
me faites un immense honneur par votre présence et le temps que vous
avez investi à cette fin. J'ai travaillé avec Frédéric Chapoton à
propos de constructions mêlant opérades et configurations de
diagonales dans les polygones. Cette expérience a été très bénéfique
puisque j'ai pu profiter de son point vue très précis sur ces
objets. Je profite de cette occasion pour remercier Frédéric pour ses
réponses patientes et complètes à mes questions régulières. Ses
contributions m'ont tellement inspiré que sans elles, ma recherche
serait aujourd'hui totalement différente. Les articles d'Alessandra
Frabetti m'ont appris beaucoup, surtout ceux traitant des séries sur
opérades ensemblistes. Je suis très content de la compter dans mon
jury, car, même si nous travaillons sur des domaines assez différents,
leur intersection contient au moins la théorie des opérades et des
bigèbres de Hopf. Florent Hivert a co-encadré ma thèse et m'a transmis
énormément de choses. Je l'ai connu alors que j'étais étudiant à
l'université de Rouen où j'ai été frappé par sa passion très vive et
communicante de la combinatoire. Je lui dois tant de choses que ces
quelques lignes ne pourraient le signifier. Que mon respect le plus
profond lui soit témoigné. J'ai collaboré avec Jean-Gabriel Luque au
sujet d'opérades émergeant d'opérateurs en théorie des langages. Sa
capacité à mélanger plusieurs concepts différents et de manière
cohérente est très enrichissante. Je le remercie aussi pour ses
invitations multiples et très agréables aux séminaires à Rouen.
\medbreak

Les collaborations en recherche sont des pratiques très épanouissantes.
Je remercie tous mes coauteurs et plus particulièrement (sans renommer
ceux déjà cités plus haut) Jean-Paul Bultel et Stéphane Vialette. J'ai
pris beaucoup de plaisir à travailler avec vous et j'espère bien que
cela va continuer.
\medbreak

Le LIGM est depuis quelques années dirigé de manière extraordinaire par
Cyril Nicaud. Je le remercie pour tout ce qu'il fait et la simplicité
qu'il y a installée. Je remercie les membres de l'équipe, Philippe
Biane, Nicolas Borie\footnote{Grand inventeur de la Borification et
du concept controversé du \og un Giraudal / des Giraudo\fg.}, Olivier
Bouillot, Cyrille Chenavier, Christophe Cordero\footnote{Mieux connu
sous le nom de\dots Ah non, je ne peux pas dire~!}, Matthieu
Josuat-Vergès, Arthur Nunge\footnote{Mieux connu sous le nom de Sire.}
et Luca Randazzo, Le groupe de travail du vendredi matin est
incontournable. Merci aussi à tout ceux qui ont fait partie de
l'équipe ou qui en sont plus ou moins proches~: Adrien Boussicault,
Grégory Châtel\footnote{Qui a, sans que personne ne sache comment, de
manière continue des thèmes et idées fascinants à partager.}, Zakaria
Chemli, Ali Chouria, Bérénice Delcroix-Oger, Gérard Duchamp, Valentin
Féray, Thibault Manneville, Gia-Thuy Pham\footnote{Qui me laisse
généreusement gagner beaucoup de parties d'échecs.}, Vincent Penelle,
Vincent Pilaud, Viviane Pons, Adrian Tanasa, Nicolas Thiéry, Christophe
Tollu et Vincent Vong. L'équipe administrative du laboratoire fait un
travail remarquable~: celui, entre autres, de simplifier de manière
miraculeuse les tâches qui paraissent insurmontables et complexes. Je
remercie à ce sujet Séverine Giboz, Patrice Hérault, Corinne
Palescandolo, Nathalie Rousseau et Pascale Souliez. Par ailleurs, on
peut facilement observer que l'équipe pédagogique du laboratoire tourne
bien, ce qui implique beaucoup de temps de gagné et du plaisir en plus.
Mes remerciements vont à ce sujet à tous ceux avec qui j'ai pu
collaborer pédagogiquement de près ou de loin, incluant Marie-Pierre
Béal, Claire David, Etienne Duris, Isabelle Fagnot, Rémi Forax, Xavier
Goaoc, Anthony Labarre, Stéphane Lohier, Antoine Meyer, Serge Midonnet,
Dominique Perrin, Carine Pivoteau, Dominique Revuz, Giuseppina Rindone,
Johan Thapper et Marc Zipstein.
\medbreak

Je saisis l'occasion pour témoigner toute ma gratitude envers
Jean-Christophe Novelli. Depuis que je le connais, il est toujours
disponible pour aider les autres, avec franchise, rigueur et diplomatie.
Il n'a cessé de me prodiguer des conseils précieux de tout ordre et de
me consacrer beaucoup de temps. Je lui témoigne ici mon respect le plus
sincère.
\medbreak

J'exprime tous mes remerciements envers mes amis (sans renommer ceux
déjà cités plus haut). En particulier, merci à Cédric, Jocelyn, Marie,
Christine et Simone.
\medbreak

Mes plus tendres remerciements sont destinés à Camille\footnote{Combe
mais sans \og s \fg.}. Merci de m'avoir soutenu depuis le
début\footnote{Et de m'avoir fait ---~bien~--- découvrir le concept de
pauses multiples en voyage.}. Je suis heureux, comblé et bien plus que
cela encore de t'avoir à mes côtés. Merci enfin à toute ma famille,
Papa, Maman, mes frères Carlo et Davide, et ma s\oe ur
Chiara\footnote{Se prononce \og Chiara \fg et non pas \og Chiara\fg.}.
Je sais que je peux compter sur vous en toutes
circonstances\footnote{Et réciproquement. Bip bip bip bip~! Je promets
d'essayer de progresser en informatique~!}. Je vous aime.
\medbreak

À tout ceux qui cherchent et ne voient pas leur nom~: mon insuffisance
est impardonnable. Si c'est le cas, ne cherchez plus, vous méritez d'y
être~; si vous n'êtes pas sous ma plume ce soir, vous êtes et resterez
dans mon c\oe ur.
\end{otherlanguage}

\tableofcontents
\listoffigures
\listoftables

\mainmatter


\chapter*{Introduction}
This dissertation contains the main research results developed since
our PhD, defended about six years ago. Our research fits in the fields
of combinatorics and algebraic combinatorics, with solid connections
with computer science. The purpose of this first part of the text is
to progressively contextualize the presented work and to provide a
preview of its main results.
\smallbreak

\SkipTocEntry\section*{Context}
Combinatorics is a subfield of both mathematics and computer science.
It is somewhat hard to provide a global and concise definition of this
field. For our part, we think that one of the best definitions of
combinatorics is that it is the science of the construction plans. A
construction plan is a list of rules expressed in a rigorous language,
whose goal is to define objects. Unlike construction plans of houses,
bridges, or space shuttles, a single construction plan in
combinatorics offers the possibility to build not only one object but
many similar ones. Indeed, a certain degree of freedom is contained in
such construction plans. All the objects thus described form a set,
named a combinatorial set. Given a construction plan, it is natural to
collect as many properties as possible of the objects of their
combinatorial set.
\smallbreak

One of the simplest examples of construction plans is the one describing
permutations. A permutation is a sequence of $n \in \N$ symbols taken in
the set $\{1, \dots, n\}$, each one appearing exactly once. From this
plan, the smallest objects are
\begin{equation}
    \epsilon, \quad
    1, \quad
    12, \enspace 21, \quad
    123, \enspace 132, \enspace 213, \enspace 231, \enspace 312,
    \enspace 321,
\end{equation}
where $\epsilon$ is the unique sequence of $n = 0$ symbols. A
slightly more elaborate example is the one of Motzkin paths. This
construction plan specifies that a Motzkin path is a possibly empty
sequence of steps of three kinds: a stationary step $\MotzHoriz$, a
rising step $\DyckUp$, or a descending step $\DyckDown$, with the
constraint that the path ends at the same level as its starting point
and never goes below its starting point. From this plan, we can build
among others the following Motzkin paths:
\begin{equation}
    \MotzPoint\,, \quad
\,,
\end{equation}
where $\MotzPoint$ is the unique sequence of $0$ steps.
\smallbreak

There are several flavors of combinatorics. Each of them depends on the
point of view on construction plans. Here follow the main ones related
with our work. To continue the  metaphor with construction plans of
buildings, let us express the point  of view of the architect, of the
workman, and of the electrician.
\smallbreak

\SkipTocEntry\subsection*{General combinatorics}
The role of the architect consists primarily in designing new
construction plans. The final aim is to use construction plans and
combinatorial sets as tools to solve precise problems or explain some
phenomena.
\smallbreak

For instance, we can study all the possible ways to bracket an
expression involving $n + 1$ occurrences of a variable $x$ and $n$
occurrences of a binary operation $\Product$ satisfying {\em a priori}
no relations. For instance,
\begin{equation}
    ((x \Product x) \Product ((x \Product x) \Product x))
\end{equation}
is one of these. A combinatorial modelization of this problem amounts
to seeing such expressions through their syntax trees. Since $\Product$
is binary, one can encode an expression with $n$ occurrences of
$\Product$ by a binary tree with $n$ internal nodes. The previous
expression is encoded in this way by the binary tree
\begin{equation} \label{equ:example_binary_tree}
\,.
\end{equation}
Now, the original problem is translated into a more combinatorial
language consisting in studying such binary trees. The construction
plan of binary trees is recursive: a binary tree is either a leaf
$\LeafPic$ or two binary trees attached to an internal node $\NodePic$.
The first ones are
\begin{equation}
    \LeafPic\,, \quad
    %
\,.
\end{equation}
From this translation, it is possible to enumerate the underlying
expressions of the trees for each size $n$. Moreover, this translation
helps to discover some properties of the expressions such as their
height, this statistics being the usual height of the binary trees.
\smallbreak

Another illustration of the work of the combinatorial architect consists
in designing combinatorial objects being the bases of some algebraic
structures. A classical example relies on free Lie algebras~\cite{Ret93}
and the description of their bases. Indeed, the combinatorial set of the
Lyndon words on a totally ordered alphabet $A$ is a basis of the free
Lie algebra generated by $A$. A Lyndon word on $A$ is a sequence $u$ of
$n \in \N$ symbols of $A$ such that all strict suffixes of $u$ are
greater than $u$ for the lexicographic order induced by the total order
on $A$. By knowing this property, the study of free Lie algebras
can be transfered on the combinatorial study of Lyndon words. A more
modern example in the same vein consists in describing the bases of free
pre-Lie algebras~\cite{Vin63,Ger63,Man11} generated by a set
$\GeneratingSet$. In this context, the right construction plan is the
one of the rooted trees on $\GeneratingSet$, that are connected acyclic
graphs whose vertices are labeled on $\GeneratingSet$ and
admitting a  distinguished vertex, the root~\cite{CL01}.
\smallbreak

The usual work in the field of general combinatorics consists hence in
modelizing a problem or a phenomenon coming from close domains such as
computer science, algebra, or physics, by combinatorial objects with the
hope of a better understanding.
\smallbreak

\SkipTocEntry\subsection*{Enumerative combinatorics}
The role of the workman, benefiting of the knowledge of a lot of
construction plans and of their internal functioning, is to understand
how does a construction plan work and to discover relations between a
bunch of them. The combinatorial workman asks in most cases the question
to count, given a construction plan, the combinatorial objects we can
build of a given size $n \in \N$. The notion of counting is primary in
enumerative combinatorics and is somewhat fuzzy. Counting may mean that
we expect a closed formula, a recurrence formula, a generating function,
or even a functional equation for a generating series. Generating series
are series of the form
\begin{equation}
    \GenSeries(t) = \sum_{n \in \N} a_n t^n
\end{equation}
where $a_n$ is the number of objects of size $n$ for each $n \in \N$.
They form a very important concept in enumerative combinatorics.
\smallbreak

For instance, by defining the size of a binary tree as its number $n$ of
internal nodes, it is possible to show that the generating series
$\GenSeries(t)$ of binary trees satisfies the algebraic equation
\begin{equation}
    \GenSeries(t) = 1 + t \GenSeries(t)^2,
\end{equation}
and expresses thus as a generating function by
\begin{equation}
    \GenSeries(t) = \frac{1 - \sqrt{1 - 4t}}{2t}.
\end{equation}
One can deduce from this that the number $a_n$ of binary trees with
$n \in \N$ internal nodes satisfies
\begin{equation} \label{equ:closed_formula_binary_trees}
    a_n = \frac{1}{n + 1} \binom{2n}{n}.
\end{equation}
\smallbreak

On the other hand, counting integer partitions is not so easy. An
integer partition of size $n \in \N$ is a multiset
$\lbag \lambda_1, \dots, \lambda_\ell\rbag$ of integers such that
$\lambda_1 + \dots + \lambda_\ell = n$. The  generating series
$\GenSeries(t)$ of these objects satisfies
\begin{equation}
    \GenSeries(t) = \prod_{k \in \N_{\geq 1}} \frac{1}{1 - t^k}.
\end{equation}
The situation here is less fruitful than in the case of binary trees
since there is no known closed formula for integer partitions similar
to~\eqref{equ:closed_formula_binary_trees}.
\smallbreak

Besides, as mentioned above, one of the roles of the combinatorial
workman consists in establishing links between different combinatorial
sets. Consider for instance the combinatorial set of Dyck paths, that
are Motzkin paths discussed before, but without horizontal steps
$\MotzHoriz$. Then, there is a bijection between the set of all binary
trees with $n$ internal nodes and the set of all Dyck paths having $n$
rising steps $\DyckUp$. This bijection can be computed by induction, but
there is a direct interpretation of it consisting in computing a Dyck
path in correspondence with a binary tree $\Tfr$ by performing a left to
right depth-first traversal of $\Tfr$ and outputting a step $\DyckUp$
when an internal node is visited and a step $\DyckDown$ when a leaf is
visited, without considering the last leaf. For instance, the Dyck path
in correspondence with the binary tree appearing
in~\eqref{equ:example_binary_tree} is
\begin{equation}
    \begin{tikzpicture}[scale=.28,Centering]
        \draw[Grid](0,0) grid (8,2);
        \node[PathNode](1)at(0,0){};
        \node[PathNode](2)at(1,1){};
        \node[PathNode](3)at(2,2){};
        \node[PathNode](4)at(3,1){};
        \node[PathNode](5)at(4,0){};
        \node[PathNode](6)at(5,1){};
        \node[PathNode](7)at(6,2){};
        \node[PathNode](8)at(7,1){};
        \node[PathNode](9)at(8,0){};
        \draw[PathStep](1)--(2);
        \draw[PathStep](2)--(3);
        \draw[PathStep](3)--(4);
        \draw[PathStep](4)--(5);
        \draw[PathStep](5)--(6);
        \draw[PathStep](6)--(7);
        \draw[PathStep](7)--(8);
        \draw[PathStep](8)--(9);
    \end{tikzpicture}\,.
\end{equation}
\smallbreak

Moreover, not only bijections between combinatorial sets are
interesting. Indeed, surjections or injections between combinatorial
sets are susceptible to establish interesting links between such sets.
For example, the algorithm of insertion of an element in a binary search
tree~\cite{Knu98} provides a surjection from the set of all
permutations of $n$ elements to the set of all binary trees with $n$
internal nodes. A binary search tree is a binary tree where all
internal nodes are labeled by integers with some extra conditions. The
insertion of a letter $a$ in a binary search tree $\Tfr$ consists in
following the path starting from the root of $\Tfr$ to one of its
leaves by going into the right subtree if $a$ is greater than the label
of the considered internal node and into the left one otherwise. For
instance, the insertion from left to right of the letters of the
permutation $\sigma := 451326$ gives the binary search tree
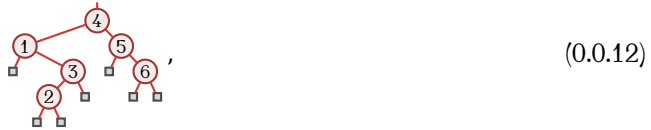
\begin{equation}
    \begin{tikzpicture}[xscale=.16,yscale=.13,Centering]
        \node[Leaf](0)at(0.00,-5.20){};
        \node[Leaf](10)at(10.00,-7.80){};
        \node[Leaf](12)at(12.00,-7.80){};
        \node[Leaf](2)at(2.00,-10.40){};
        \node[Leaf](4)at(4.00,-10.40){};
        \node[Leaf](6)at(6.00,-7.80){};
        \node[Leaf](8)at(8.00,-5.20){};
        \node[Node](1)at(1.00,-2.60){\begin{math}1\end{math}};
        \node[Node](11)at(11.00,-5.20){\begin{math}6\end{math}};
        \node[Node](3)at(3.00,-7.80){\begin{math}2\end{math}};
        \node[Node](5)at(5.00,-5.20){\begin{math}3\end{math}};
        \node[Node](7)at(7.00,0.00){\begin{math}4\end{math}};
        \node[Node](9)at(9.00,-2.60){\begin{math}5\end{math}};
        \draw[Edge](0)--(1);
        \draw[Edge](1)--(7);
        \draw[Edge](10)--(11);
        \draw[Edge](11)--(9);
        \draw[Edge](12)--(11);
        \draw[Edge](2)--(3);
        \draw[Edge](3)--(5);
        \draw[Edge](4)--(3);
        \draw[Edge](5)--(1);
        \draw[Edge](6)--(5);
        \draw[Edge](8)--(9);
        \draw[Edge](9)--(7);
        \node(r)at(7.00,2.75){};
        \draw[Edge](r)--(7);
    \end{tikzpicture}\,,
\end{equation}
and, by forgetting the labels of the nodes, we obtain the binary tree
image of $\sigma$. This surjection has also several algebraic
properties~\cite{LR98,HNT05}.
\smallbreak

A last important part of enumerative combinatorics includes algorithms
generating all the objects of a given size of a combinatorial
set~\cite{Rus03}. The efficiency of these algorithms is a highly
important feature, so that constant amortized time algorithms are the
most sought. One can cite in this context the algorithm of
Proskurowski and Ruskey for binary trees~\cite{PR85}, and the
Steinhaus--Johnson--Trotter algorithm for
permutations~\cite{Tro62,Joh63,Ste64}.
\smallbreak

\SkipTocEntry\subsection*{Algebraic combinatorics}
The role of the electrician consists in endowing an edifice with
a network of electric wires, making it capable to perform additional
functions. Given a construction plan, the combinatorial electrician
tries to define operations on its combinatorial objects. Operations on
combinatorial objects allow to assemble several of these to obtain
bigger ones, or, contrariwise, allow to disassemble a single object into
smaller pieces. In this last case, it is more accurate to speak of
co-operations. This point of view draws a bridge between combinatorics
and algebra, creating interactions in both ways between these two fields.
\smallbreak

For instance, operations on Motzkin paths offer an interesting way to
describe their generating series $\GenSeries(t)$, counting them with
respect to their number of steps. For this, consider the monoid
$(\Pca, \Conc)$ of all paths consisting in steps $\MotzHoriz$,
$\DyckUp$, and $\DyckDown$ (here we relax the conditions about the levels
of the starting and ending points of the paths), where $\Conc$ is the
concatenation of paths (obtained by superimposing the ending point of
the first path and the starting point of the second). For instance,
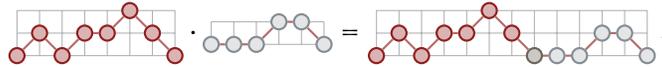
\begin{equation*}
\,.
\end{equation*}
Now, let $\Gbf$ be the formal series defined as the formal sum of all
Motzkin paths. Hence,
\begin{equation} \label{equ:series_Motzkin_monoid}
    \Gbf =
    \MotzPoint
    +
    %
    + \cdots.
\end{equation}
By nearly elementary properties of Motzkin paths about their unambiguous
decomposition, and by extending $\Conc$ linearly on series, $\Gbf$
can be expressed as
\begin{equation} \label{equ:series_Motzkin_paths_monoid}
    \Gbf =
    \MotzPoint
    +
    \MotzHoriz \Conc \Gbf
    +
    \DyckUp \Conc \Gbf \Conc \DyckDown \Conc \Gbf.
\end{equation}
By observing that $\GenSeries(t)$ is the series obtained by specializing
each Motzkin path of size $n$ by $t^n$ in $\Gbf$, we deduce
from~\eqref{equ:series_Motzkin_paths_monoid} that $\GenSeries(t)$
satisfies
\begin{equation}
    \GenSeries(t) = 1 + t \GenSeries(t) + t^2 \GenSeries(t)^2.
\end{equation}
From this algebraic equation, it is possible to obtain a generating
function of $\GenSeries(t)$ or a closed formula for its coefficients,
like in the case of binary trees presented above.
\smallbreak

On the other hand, endowing combinatorial sets with operations allows to
highlight some of their properties. Consider in this context the monoid
$(\SymmetricGroup, \Over)$ where $\SymmetricGroup$ is the combinatorial
set of all permutations and $\Over$ is the shifted concatenation of
permutations: given two permutations $\sigma$ and $\nu$,
$\sigma \Over \nu$ is the permutation obtained by concatenating $\sigma$
with the word obtained by incrementing each letter of $\nu$ by the size
of $\sigma$. For instance,
\begin{math}
    \textcolor{Col1}{312} \Over \textcolor{Col4}{21} =
    \textcolor{Col1}{312}\textcolor{Col4}{54}
\end{math}.
The minimal generating set of this monoid is the set of all connected
permutations, that are the nonempty permutations $\sigma$ having no
proper prefixes that are permutations~\cite{Com72}. For instance, the
first connected permutations are
\begin{multline}
    1, \quad 21, \quad
    231, \enspace 312, \enspace 321, \quad \\
    2341, \enspace 2413, \enspace 2431, \enspace 3142, \enspace
    3241, \enspace 3412, \enspace 3421, \enspace 4123, \enspace
    4132, \enspace 4213, \enspace 4231, \enspace 4312, \enspace 4321.
\end{multline}
From this very natural question about finding a minimal generating set
of an algebraic structure, we obtain the description of new natural
combinatorial objects. Moreover, since $(\SymmetricGroup, \Over)$ is
free as a monoid, the generating series $\GenSeries_\Cca(t)$ of
connected permutations and the generating series
$\GenSeries_\SymmetricGroup(t)$ of permutations are related by
\begin{equation}
    \GenSeries_\SymmetricGroup(t)
    =
    \frac{1}{1 - \GenSeries_\Cca(t)}
    =
    \sum_{n \in \N} n! t^n.
\end{equation}
Connected permutations have several properties. For instance, the Hopf
bialgebra of free quasi-symmetric functions (also known as the
Malvenuto-Reutenauer Hopf bialgebra~\cite{MR95}) is, as an associative
algebra, freely generated by the connected permutations~\cite{DHT02}.
\smallbreak

Let us consider another example consisting in a unary operation on
binary trees. As said before, binary trees are in one-to-one
correspondence with expressions involving variables $x$ and operations
$\Product$. Now, assume that $\Product$ is associative. This leads
to allow the relation
\begin{equation} \label{equ:rotation_expressions}
    \left(\dots \left(u_1 \Product u_2\right) \Product u_3 \dots\right)
    =
    \left(\dots u_1 \Product \left(u_2 \Product u_3\right) \dots\right),
\end{equation}
where $u_1$, $u_2$, and $u_3$ are expressions. In terms of binary
trees~\eqref{equ:rotation_expressions} translates as the identification
\begin{equation} \label{equ:rotation_binary_trees}
    \begin{tikzpicture}[xscale=.27,yscale=.3,Centering]
        \node(0)at(0.00,-3.33){\begin{math}\Tfr_1\end{math}};
        \node(2)at(2.00,-3.33){\begin{math}\Tfr_2\end{math}};
        \node(4)at(4.00,-1.67){\begin{math}\Tfr_3\end{math}};
        \node[Node](1)at(1.00,-1.67){};
        \node[Node](3)at(3.00,0.00){};
        \draw[Edge](0)--(1);
        \draw[Edge](1)--(3);
        \draw[Edge](2)--(1);
        \draw[Edge](4)--(3);
        \node(r)at(3.00,1){};
        \draw[Edge](r)--(3);
    \end{tikzpicture}
    =
    \begin{tikzpicture}[xscale=.27,yscale=.3,Centering]
        \node(0)at(0.00,-1.67){\begin{math}\Tfr_1\end{math}};
        \node(2)at(2.00,-3.33){\begin{math}\Tfr_2\end{math}};
        \node(4)at(4.00,-3.33){\begin{math}\Tfr_3\end{math}};
        \node[Node](1)at(1.00,0.00){};
        \node[Node](3)at(3.00,-1.67){};
        \draw[Edge](0)--(1);
        \draw[Edge](2)--(3);
        \draw[Edge](3)--(1);
        \draw[Edge](4)--(3);
        \node(r)at(1.00,1){};
        \draw[Edge](r)--(1);
    \end{tikzpicture},
\end{equation}
where $\Tfr_1$, $\Tfr_2$, and $\Tfr_3$ are binary trees. This
identification can be performed anywhere in the binary trees and not
only at their roots. One can see this identification as an operation
consisting in taking a binary tree and one of its edges oriented to
the left (like the left member of~\eqref{equ:rotation_binary_trees})
and changing it into an edge oriented to right  (like the right member
of~\eqref{equ:rotation_binary_trees}). This operation is known as a
right rotation~\cite{Knu98}. Now, it is possible to show by induction on
the number of internal nodes of the binary trees that all binary trees
with $n$ internal nodes can be identified with the right comb tree of
size $n$, that is the binary tree such that the left child of each
internal node is a leaf. This provides a (quite complicated) proof of
the well-known fact that all parenthesizings of an expression
involving $n$ occurrences of an associative operation $\Product$ are
equal.
\smallbreak

Finally, as mentioned before, interpreting an algebraic structure by
means of combinatorial objects endowed with operations---like free Lie
algebras in terms of Lyndon words, free pre-Lie algebras in terms of
rooted trees, free dendriform algebras~\cite{Lod01} in terms of binary
trees, or even free monoids in terms of words---brings a good
understanding of it. These interpretations of algebraic structures are
known as combinatorial realizations. In algebraic combinatorics, we
endow combinatorial sets, or more generally spaces whose bases are
indexed by combinatorial sets, with several algebraic structures. These
can be simply monoids or groups, but in some cases posets, lattices,
associative algebras, dendriform algebras, pre-Lie algebras,
duplicial algebras~\cite{Lod08}, {\em etc.} In this work, Hopf
bialgebras, operads, and pros are the structures encountered the most.
\smallbreak

\SkipTocEntry\subsection*{Other flavors of combinatorics}
In addition to enumerative and algebraic combinatorics, there are other
important flavors of combinatorics. Among these is analytic
combinatorics~\cite{FS09} wherein techniques coming from complex
analysis are employed at the level of generating series. This field
studies also the asymptotic behavior and the general form of
combinatorial objects. Probabilistic combinatorics is close to analytic
combinatorics. This domain uses methods coming from probability theory
to design algorithms randomly generating objects of a given
combinatorial set (see for instance~\cite{Rem85} for an algorithm
generating uniformly binary trees). Moreover, probabilistic
combinatorics is useful to show, within a given combinatorial set, that
there is at least one object satisfying a given property~\cite{AS08}.
This point of view was initiated by Erdős and has links with the
Ramsey theory~\cite{Soi10}. As a last flavor mentioned here, one can
cite geometric combinatorics wherein geometric realizations of
polytopes are designed, including the realizations of the
permutohedron and of the associahedron~\cite{CSZ15}.
\smallbreak

Let us now dive a little more deeply into algebraic combinatorics and
explain our point of view about it and the context of our contributions.
\smallbreak

\SkipTocEntry\section*{Point of view}
Historically, algebraic combinatorics was concerned with questions
related to representation theory~\cite{GL01}. This field consists in
studying algebraic structures (like monoids, groups, associative
algebras, Lie algebras, {\em etc.}) by regarding their elements as
linear maps. In an equivalent way, this amounts to letting the structure
act on a vector space in a reasonable way. One of the benefits of this
process rests upon the fact that algebraic problems are translated into
linear algebra questions. The underlying algorithmic of linear algebra
(like Gaussian elimination, matrix inversion, matrix reduction,
{\em etc.}) offers strategies coming from computer science and
combinatorics to explore these problems.
\smallbreak

In particular, representations of the symmetric groups
$\SymmetricGroup_n$ of permutations of size $n \in \N$ have a special
status in algebraic combinatorics. Indeed, the irreducible
representations of $\SymmetricGroup_n$ are indexed by integer
partitions of size $n$~\cite{FH91}. Schur functions are symmetric
functions indexed by integer partitions that appear in this context of
representation theory. They have the particularity to admit a lot of
very different but equivalent definitions~\cite{Las84,Sta99,Lot02}. The
set of all symmetric functions is naturally endowed with the structure
of an associative algebra $\SymCom$~\cite{Mac15} and the set of all
Schur functions is one of its bases. Many other bases of $\SymCom$
have been discovered, like the monomial, elementary, complete
homogeneous, and power sum functions. The changes of bases between
these different families of functions express most of the time by
simple and nice combinatorial algorithms.
\smallbreak

This work is distant from these classical questions about representation
theory and symmetric functions. Our point of view about algebraic
combinatorics is somewhat unrelated to these considerations but,
instead, related to the study of operations and algebraic structures on
combinatorial sets. Nevertheless, like in all these research areas, we
work most of the time with finite structures that can be encoded by
the computer. For this reason, we can use the computer to perform
large computations or to make experiments. These are very powerful tools
to establish conjectures and to collect as much information as possible
about a given research subject.
\smallbreak

\SkipTocEntry\subsection*{Objects, operations, and algebraic structures}
In accordance to what we have explained above, defining and studying
operations on combinatorial sets has several advantages. More
precisely, in this context, we try to progress in both of the
following axes:
\begin{enumerate}[label={\bf ({\Alph*})}]
    \item \label{item:defining_operations}
    Endowing combinatorial sets with algebraic structures by defining
    operations or co-operations;
    \smallbreak

    \item \label{item:realizing_structures}
    Given a type of algebraic structure, searching a realization of it
    in terms of combinatorial objects endowed with operations.
\end{enumerate}
Let us explain in more details these two directions.
\smallbreak

Point~\ref{item:defining_operations} consists, starting with a
combinatorial set $C$, in defining operations or co-operations on $C$.
In practice, we work rather on $\K \Angle{C}$, the linear span of $C$
where $\K$ is a field. To highlight some statistics on the objects,
$\K$ is often the field $\K(q_0, q_1, \dots)$ of rational functions
on the parameters $q_i$, $i \in \N$. The linear structure of
$\K \Angle{C}$ implies that we inherit the techniques coming from
linear algebra to  perform its study. When the (co)operations defined
on $\K \Angle{C}$ endow it with a certain algebraic structure (like
an associative algebra, a dendriform algebra, a pre-Lie algebra, or
even a coalgebra), we can ask all the algebraic questions related to
the structure and we can hope to harvest information about the
objects of~$C$.
\smallbreak

Among the classical questions, the first one consists in expressing new
bases of $\K \Angle{C}$ and observing how its operations behave on
these. Frequently, changes of bases are triangular and are defined
through partial orders on $C$ by considering sums of elements minored by
other ones. It is time to study an example. Let us endow
$\K \Angle{\SymmetricGroup}$ with the linear binary product $\cshuffle$,
where for any permutations $\sigma$ and $\nu$, $\sigma \cshuffle \nu$
is the sum of all the permutations that can be obtained by
interleaving the letters of $\sigma$ with the ones of the word obtained
by incrementing by the size of $\sigma$ the letters of $\nu$. For
instance,
\begin{equation}
    \textcolor{Col1}{12} \cshuffle \textcolor{Col4}{21} =
    \textcolor{Col1}{12}\textcolor{Col4}{43} +
    \textcolor{Col1}{1}\textcolor{Col4}{4}\textcolor{Col1}{2}
        \textcolor{Col4}{3} +
    \textcolor{Col1}{1}\textcolor{Col4}{43}\textcolor{Col1}{2} +
    \textcolor{Col4}{4}\textcolor{Col1}{12}\textcolor{Col4}{3} +
    \textcolor{Col4}{4}\textcolor{Col1}{1}\textcolor{Col4}{3}
        \textcolor{Col1}{2} +
    \textcolor{Col4}{43}\textcolor{Col1}{12}.
\end{equation}
This product is known as the shifted shuffle product of permutations.
Consider now the partial order $\Ord$ on $\SymmetricGroup$ being the
reflexive and transitive closure of the relation $\Rca$ such that
for any permutations $\sigma$ and $\nu$, $\sigma \,\Rca\, \nu$ if $\nu$
can be obtained from $\sigma$ by exchanging two adjacent letters
$\sigma(i)$ and $\sigma(i + 1)$ such that $\sigma(i) <\sigma(i + 1)$.
This order is known as the right weak order on
permutations~\cite{GR63,YO69}. Let the family
$\{\BasisE_\sigma : \sigma \in \SymmetricGroup\}$ of elements of
$\K \Angle{\SymmetricGroup}$ defined by
\begin{equation}
    \BasisE_\sigma :=
    \sum_{\substack{
        \nu \in \SymmetricGroup \\
        \sigma \Ord \nu
    }}
    \nu.
\end{equation}
For instance,
\begin{equation}
    \BasisE_{2341} = 2341 + 2431 + 3241 + 3421 + 4231 + 4321.
\end{equation}
By triangularity, this family forms a basis of
$\K \Angle{\SymmetricGroup}$ and it appears that the product $\cshuffle$
on this $\BasisE$-basis satisfies, for all permutations $\sigma$ and
$\nu$,
\begin{equation} \label{equ:product_E_basis_permutations}
    \BasisE_\sigma \cshuffle \BasisE_\nu =
    \BasisE_{\sigma \Over \nu},
\end{equation}
where $\Over$ is the shifted concatenation of permutations encountered
before. This provides an example of a rather complicated product when
considered in a given basis that becomes very simple in another one.
Moreover, proving that~\eqref{equ:product_E_basis_permutations} holds
provides an interesting proof of the associativity of~$\cshuffle$ since
$\Over$ is clearly associative.
\smallbreak

A second question almost as much immediate as the first one is to find
minimal generating sets of $\K \Angle{C}$, whose elements can be
interpreted as base blocks to build any object of $C$. This is even more
interesting when $\K \Angle{C}$ has some freeness properties; in this
case, any element decomposes in a unique way in a certain sense. As a
consequence, obtaining minimal generating sets of $\K \Angle{C}$ leads
to expressions for the Hilbert series
\begin{equation}
    \HilbSeries_{\K \Angle{C}}(t) =
    \sum_{n \in \N} \dim \K \Angle{C(n)} t^n
\end{equation}
of $\K \Angle{C}$, where $C(n)$ is the set of the objects of size
$n \in \N$ of $C$. Since as generating series
$\HilbSeries_{\K \Angle{C}}(t)$ is the generating series
$\GenSeries_C(t)$ of $C$, this may offer an alternative way to
enumerate the objects of $C$. To continue the example we started,
$(\K \Angle{\SymmetricGroup}, \cshuffle)$ admits as a minimal generating
set the set $\{\BasisE_\sigma : \sigma \in \Cca\}$ where $\Cca$ is the
set all connected permutations, and is freely generated by this set as
an associative algebra (for details see~\cite{DHT02}).
\smallbreak

Besides, to complete the study of $\K \Angle{C}$, it is natural to
study morphisms (with respect to the algebraic structure equipping
$\K \Angle{C}$) involving it. Automorphisms of $\K \Angle{C}$ lead
potentially to the discovery of more or less hidden symmetries between
the objects of $C$. Morphisms between $\K \Angle{C}$ and other  known
structures $\K \Angle{D}$ lead to establish connections between the
objects of $C$ and the ones of $D$. This also includes the study of
substructures and quotients of $\K \Angle{C}$. It is worth
observing that most of such morphisms use algorithms coming from
computer science in an unexpected way. For instance, the associative
algebra $(\K \Angle{\SymmetricGroup}, \cshuffle)$ admits several
substructures involving a large range of combinatorial objects. Some
of these can be constructed by considering a family
$\{\BasisP_x : x \in D\}$ of elements of $\K \Angle{\SymmetricGroup}$
defined by
\begin{equation}
    \BasisP_x :=
    \sum_{\substack{
        \sigma \in \SymmetricGroup \\
        \Algo(\sigma) = x
    }}
    \sigma,
\end{equation}
where $D$ is a certain combinatorial set and $\Algo$ is an algorithm
transforming a permutation into an object of $D$. When $\Algo$ satisfies
some precise properties (see~\cite{Hiv99,HN07,Gir11,NT14}), the
$\BasisP$-family spans an associative subalgebra of
$(\K \Angle{\SymmetricGroup}, \cshuffle)$. For instance, $D$ can be the
set of binary trees and $\Algo$, the algorithm of insertion in a binary
tree exposed above. In this case, one has for instance
\begin{equation}
    \BasisP_{
        \begin{tikzpicture}[xscale=.13,yscale=.11,Centering]
            \node[Leaf](0)at(0.00,-4.50){};
            \node[Leaf](2)at(2.00,-4.50){};
            \node[Leaf](4)at(4.00,-6.75){};
            \node[Leaf](6)at(6.00,-6.75){};
            \node[Leaf](8)at(8.00,-4.50){};
            \node[Node](1)at(1.00,-2.25){};
            \node[Node](3)at(3.00,0.00){};
            \node[Node](5)at(5.00,-4.50){};
            \node[Node](7)at(7.00,-2.25){};
            \draw[Edge](0)--(1);
            \draw[Edge](1)--(3);
            \draw[Edge](2)--(1);
            \draw[Edge](4)--(5);
            \draw[Edge](5)--(7);
            \draw[Edge](6)--(5);
            \draw[Edge](7)--(3);
            \draw[Edge](8)--(7);
            \node(r)at(3.00,2.5){};
            \draw[Edge](r)--(3);
        \end{tikzpicture}}
    =
    2143 + 2413 + 2431.
\end{equation}
Besides, $D$ can be the set of the standard Young tableaux and $\Algo$,
the algorithm consisting in inserting the letters of a permutation into
a standard Young tableau using the Schensted
algorithm~\cite{Sch61,Lot02}. In this case, one has for instance
\begin{equation}
    \BasisP_{\small \young(124,3)} = 1324 + 1342 + 3124.
\end{equation}
\smallbreak

These mechanisms, coming from algebraic combinatorics, can also be used
to conjecture properties and to obtain results in enumerative
combinatorics. Indeed, assume that $C$ and $D$ are two combinatorial
sets and that we look for a bijection between them. A tool to discover
a bijection consists in endowing $\K \Angle{C}$ and $\K \Angle{D}$
with algebraic structures satisfying similar properties. More
explicitly, when these structures admit minimal generating sets clearly
in bijection, and when the objects of $C$ and $D$ decompose in the
same way on the generators, one obtains a computable bijection between
$C$ and~$D$.
\smallbreak

In summary, Direction~\ref{item:defining_operations} uses algebra to
obtain results in combinatorics and in computer science.
\smallbreak

Conversely, Direction~\ref{item:realizing_structures} employs mechanisms
and techniques coming from combinatorics to solve algebraic questions.
Given a type of algebra, that is a set of (co)operation symbols together
with axioms they have to satisfy, the knowledge of the free structure
on a set $\GeneratingSet$ of generators brings a lot of information.
In many cases, the description of these structures is combinatorial, in
the sense that their bases are indexed by combinatorial objects labeled
in an adequate way by elements of $\GeneratingSet$. As mentioned before,
examples are abundant in the literature. They include pre-Lie algebras
using rooted trees and operations of grafting of trees~\cite{CL01},
Zinbiel algebras using words and half-shuffle operations~\cite{Lod95},
dendriform algebras using binary trees and operations of shuffling of
trees~\cite{Lod01}, operads using planar rooted trees and grafting
operations, and pros using prographs and operations of
compositions~\cite{Mar08}. Several of these combinatorial realizations
of algebraic structures can be established by orienting their axioms to
obtain rewrite rules~\cite{BN98}. When the obtained rewrite rules
satisfy some properties like termination and confluence, the normal
forms of the rewrite rules can be seen as the elements of the structure.
\smallbreak

Another classical example of use of combinatorial methods for algebra is
provided by the Littlewood-Richardson rule~\cite{LR34}. This rule offers
a way to compute the structure coefficients of the algebra of symmetric
functions $\SymCom$ in the basis of the Schur functions. A simple and
enlightening proof~\cite{DHT02,HNT05} of this rule is provided by the
combinatorics of Young tableaux and of the plactic
monoid~\cite{LS81,Lot02}.
\smallbreak

Let us now provide details about the main structures appearing here.
During our research, we work particularly with three types of algebraic
structures: Hopf bialgebras, operads, and pros. We now present some of
their features and why they are interesting and adapted structures in
the field of algebraic combinatorics.
\smallbreak

\SkipTocEntry\subsection*{Hopf bialgebras}
Hopf bialgebras are vector spaces endowed with an associative product
$\Product$ and a coassociative coproduct $\Coproduct$. These
(co)operations satisfy the relation
\begin{equation} \label{equ:Hopf_relation}
    \Coproduct(x \Product y) = \Coproduct(x) \Coproduct(y)
\end{equation}
for any elements $x$ and $y$. If we see $\Product$ as a product
assembling two elements to build another one, and $\Coproduct$ as a
coproduct breaking an element into two smaller parts,
Equation~\eqref{equ:Hopf_relation} says that assembling two elements and
then breaking the result is the same as assembling the results obtained
by breaking them before. This kind of commutation between $\Product$
and $\Coproduct$ is thus very natural. Hopf bialgebras $\K \Angle{C}$
where $C$ is a combinatorial set with exactly one element of size $0$
and where $\Product$ (resp. $\Coproduct$) is graded (resp. cograded)
are the most encountered ones in algebraic combinatorics. These
structures are known as combinatorial Hopf bialgebras. Main references
about these structures are~\cite{Car07} and~\cite{GR16}.
\smallbreak

The prototypal example of a Hopf bialgebra is the symmetric functions
$\SymCom$. Indeed, it is possible to add a coproduct on $\SymCom$ to
turn it into a Hopf bialgebra. Most of other Hopf bialgebras are
generalizations of $\SymCom$ in the sense that they contain it as a
quotient or as a Hopf sub-bialgebra. A famous full diagram of Hopf
bialgebras includes the Malvenuto-Reutenauer Hopf bialgebra~\cite{MR95},
also known as $\FQSym$~\cite{DHT02}. This structure is the space
$\K \Angle{\SymmetricGroup}$ endowed with the shifted shuffle product
and a deconcatenation coproduct of permutations. The
Malvenuto-Reutenauer Hopf bialgebra contains the Poirier-Reutenauer
Hopf bialgebra of tableaux~\cite{PR95}, also known as the Hopf
bialgebra of free symmetric functions $\FSym$~\cite{DHT02,HNT05} and
involves standard Young tableaux. The Loday-Ronco Hopf
bialgebra~\cite{LR98}, also known as the Hopf bialgebra of binary
search trees $\PBT$~\cite{HNT05} involves binary trees and is a Hopf
sub-bialgebra of $\FQSym$. Moreover, a noncommutative version
$\Sym$~\cite{GKLLRT94} of $\SymCom$ exists as a Hopf sub-bialgebra of
$\FQSym$ known as the Hopf bialgebra of noncommutative symmetric
functions. This structure involves integer compositions and provides
noncommutative versions of Schur functions. Furthermore, a lot of Hopf
bialgebras involving various sorts of trees with links with
renormalization theory like the Connes-Kreimer Hopf bialgebra
$\CK$~\cite{CK98} have been introduced. Several variations of this
structure exist~\cite{Foi02a,Foi02b} (see also~\cite{FNT14}).
\smallbreak

One of the main striking facts shared by most of these constructions is
that they establish links between combinatorial objects through
combinatorial algorithms, lead to the definition of monoids (like the
plactic~\cite{LS81,Lot02}, sylvester~\cite{HNT05}, and
hypoplactic~\cite{KT97,KT99} monoids), and use partial orders (the right
weak order on permutations~\cite{GR63,YO69}, the Tamari order on
binary trees~\cite{Tam62}, and the refinement order on integer
compositions).
\smallbreak

Besides Hopf bialgebras that are structures allowing, as explained, to
assemble or disassemble objects, operads are other ones manipulating
combinatorial objects. These last work by composing objects together
rather than assembling them.
\smallbreak

\SkipTocEntry\subsection*{Operads}
Operads are algebraic structures introduced in the context of algebraic
topology~\cite{May72,BV73}. These structures provide an abstraction of
the notion of operators (of any arities) and of their compositions. This
theory has somewhat been neglected during almost the first two decades
after its discovery. In the 1990s, the theory of operads enjoyed a
renaissance raised by Loday~\cite{Lod96} and, since the 2000s, many links
between the theory of operads and combinatorics have been developed.
A large survey of this theory can be found in~\cite{Mar08,LV12,Men15}.
\smallbreak

The modern treatment of operads in algebraic combinatorics consists in
regarding combinatorial objects like operators endowed with gluing
operations mimicking the composition (see for instance~\cite{Cha08}).
From an intuitive point of view, an operad is a set (or a space) of
abstract operators with several inputs and one output that can be
composed in many ways. More precisely, if $x$ is an operator with $n$
inputs and $y$ is an operator with $m$ inputs, $x \circ_i y$ denotes the
operator with $n + m - 1$ inputs obtained by gluing the output of $y$ to
the $i$th input of $x$. Pictorially,
\begin{equation}
\,.
\end{equation}
There is also an action $\Action$ of the symmetric group
$\SymmetricGroup_n$ on the elements of arity $n$ letting to permute
their inputs. Operads are algebraic structures related to trees in the
same way as monoids are algebraic structures related to words (by
their free objets). There are numerous variations and enrichments of
operads, like cyclic operads~\cite{GK95}, colored
operads~\cite{BV73,Yau16}, and nonsymmetric operads. In this
dissertation, we work mainly with nonsymmetric operads (also called ns
operads or pre-Lie systems).
\smallbreak

A large number of interactions between operads and combinatorics exist.
Let us explain four of these. Koszul duality of operads in an important
part of the theory. This kind of duality has been introduced by Ginzburg
and Kapranov~\cite{GK94} as an extension of the analogous duality for
quadratic associative algebras. An operad is by definition Koszul if
its Koszul complex is acyclic~\cite{GK94}. When $\Oca$ is a Koszul
operad, its Hilbert series $\HilbSeries_\Oca(t)$ and the one
$\HilbSeries_{\Oca^!}(t)$ of its Koszul dual $\Oca^!$ are related by
\begin{equation}
    \HilbSeries_\Oca\left(-\HilbSeries_{\Oca^!}(-t)\right) = t.
\end{equation}
Hence, from the knowledge of $\HilbSeries_\Oca(t)$, one can hope to
compute the coefficients of $\HilbSeries_{\Oca^!}(t)$. Moreover, the
Koszulity property for operads is strongly related to the theory of
rewrite rules on trees, this last theory providing a sufficient
combinatorial condition to prove the Koszulity of an
operad~\cite{Hof10,DK10,LV12}. Besides, another strategy to prove that
an operad is Koszul consists in constructing a family of posets from
an operad~\cite{MY91}, so that the Koszulity of the considered operad is
a consequence of a combinatorial property of these
posets~\cite{Val07}. Operads lead also to generalized versions of
generating series, enriching the usual techniques for enumeration. Given
an operad $\Oca$, one can consider formal series of the form
\begin{equation}
    \Fbf = \sum_{x \in \Oca} \lambda_x x,
\end{equation}
where the $\lambda_x$ are coefficients in $\K$. From the partial
compositions of $\Oca$, we can endow the set $\K \AAngle{\Oca}$ of all
series on $\Oca$ with a monoid structure. Chapoton studied some of
these for many usual operads~\cite{Cha02,Cha08,Cha09}. Some other
authors consider such series as {\em e.g.}, van der Laan~\cite{Vdl04},
Frabetti~\cite{Fra08}, and Loday and Nikolov~\cite{LN13}.
\smallbreak

Let us provide a simple example involving series on operads. The set of
all Motzkin paths forms a structure of a ns operad $\Motz$ where, given
two Motzkin paths $u$ and $v$, $u \circ_i v$ is the path obtained by
replacing the $i$th point of of $u$ (indexed from left to right) by $v$.
For example,
\begin{equation}
\,.
\end{equation}
In this operad, a Motzkin path of $n - 1$ steps is seen as an operator
of arity $n$. Now, by denoting by $\Gbf \in \K \AAngle{\Motz}$ the
formal sum of all the elements of $\Motz$, one obtains the relation
\begin{equation} \label{equ:series_Motzkin_operad}
    \Gbf =
    \MotzPoint + \MotzHoriz \circ \left[\MotzPoint, \Gbf\right]
    + \MotzPeak \circ \left[\MotzPoint, \Gbf, \Gbf\right],
\end{equation}
where $\circ$ is the complete composition map of $\Motz$ extended on
series on $\Motz$. Of course, this expression for $\Gbf$ is very similar
to the one provided by~\eqref{equ:series_Motzkin_monoid}
but~\eqref{equ:series_Motzkin_operad} admits at least two major
advantages. First, contrariwise to~\eqref{equ:series_Motzkin_monoid}
which relies on the monoid $(\Pca, \Conc)$ of all
paths (because $\DyckUp$ and $\DyckDown$ are not Motzkin paths),
\eqref{equ:series_Motzkin_operad} only uses elements of $\Motz$.
Second, the fact that~\eqref{equ:series_Motzkin_operad} holds is a
consequence of a presentation by generators and relations of $\Motz$.
Indeed, one can show that $\left\{\MotzHoriz, \MotzPeak\right\}$ is a
minimal generating set of $\Motz$ and that there is a convergent
orientation of the nontrivial relations between these generators so that
the normal forms are precisely the terms of the form $\MotzPoint$,
$\MotzHoriz \circ [\MotzPoint, u]$, or
$\MotzPeak \circ [\MotzPoint, u, v]$, where $u$ and $v$ are Motzkin
paths. This combinatorial property is a consequence of the Koszulity of
$\Motz$. All this provides another example of combinatorial properties
encapsulated into suitable algebraic structures.
\smallbreak

Let us provide now a little more elaborate example concerning the
enumeration of balanced binary trees. These trees were introduced in
an algorithmic context~\cite{AVL62} as efficient data structures to
represent dynamic finite sets. A binary tree $\Tfr$ is balanced if
for any of its internal node $u$, the height of the right and left
subtrees of $u$ differ by at most~$1$. For example,
\begin{equation}

\end{equation}
is a balanced binary tree. The generating series $\GenSeries(t)$ of
these trees, enumerating them with respect to their number of leaves,
satisfies $\GenSeries(t) = F(t, 0)$ where $F(x, y)$ is the bivariate
series satisfying the functional equation
\begin{equation}
    F(x, y) = x + F\left(x^2 + 2xy, x\right).
\end{equation}
The coefficients of $\GenSeries(t)$ can hence been computed by
iteration. This way of enumerating balanced binary trees is presented
in~\cite{BLL88,BLL98,Knu98}. By using colored operads, it is possible to
obtain a better description for the coefficients of $\GenSeries(t)$. For
this purpose, let $\CMag$ be the colored operad on the set of all binary
trees such that each leaf and the root has a color in $\{1, 2\}$. The
partial composition $\Sfr \circ_i \Tfr$ of two such trees is defined if
the output color of $\Tfr$ is the same as the color of the $i$th input
of $\Sfr$ and is the tree obtained by grafting the root of $\Tfr$ onto
the $i$th leaf of $\Sfr$. For example,
\begin{equation}
\,.
\end{equation}
In this operad, a binary tree having $n$ leaves is seen as an operator
of arity $n$. Now, Let $\Gbf \in \K \AAngle{\CMag}$ be the formal series
defined as the formal sum of all the balanced binary trees, seen as
elements of $\CMag$ where all colors are equal to $1$. Hence,
\begin{equation}
    \Gbf =
    %
    + \cdots.
\end{equation}
We obtain the relation
\begin{equation} \label{equ:series_CMag_operad}
    \Gbf =
    %
,
\end{equation}
where $\Compo$ is an associative product on series on $\CMag$ obtained
from its partial compositions maps and $\Compo_*$ is the Kleene star
of $\Compo$. One can deduce from~\eqref{equ:series_CMag_operad} and from
the properties of the operations $\Compo$ and $\Compo_*$ the
recurrence
\begin{equation}
    g(n, m) =
    \begin{cases}
        1 & \mbox{if } (n, m) = (1, 0), \\
        \sum\limits_{\substack{
            \ell_1, \ell_2 \in \N \\
            n = 2\ell_1 + \ell_2 + m
        }}
        \binom{\ell_1 + m}{\ell_1}
        2^m
        g(\ell_1 + m, \ell_2) & \mbox{otherwise}
    \end{cases}
\end{equation}
for the number $g(n, 0)$ of balanced binary trees with $n$ leaves.
\smallbreak

These two examples show that the formalization of combinatorial problems
in terms of operads offer tools for enumerative questions. Contrariwise,
operads on combinatorial objects may lead to algebraic observations.
Indeed, any operad $\Oca$ defines a category of algebras called
$\Oca$-algebras. Any $\Oca$-algebra can be seen as a representation of
$\Oca$ in the sense that $\Oca$ acts on any $\Oca$-algebra. For
instance, there is an operad $\Lie$ describing the category of all Lie
algebras, an operad $\As$ describing the category of all associative
algebras, and an operad $\Dendr$ describing the category of all
dendriform algebras~\cite{Lod01}. Morphisms $\phi : \Oca_1 \to \Oca_2$
between two operads $\Oca_1$ and $\Oca_2$ give rise to functors from the
category of $\Oca_2$-algebras to the one of $\Oca_1$-algebras. For
instance, $\Lie$ is a suboperad of $\As$ so that there is an injective
morphism $\phi : \Lie \to \As$. This morphism translates into the
well-known functor from associative algebras to Lie algebras consisting
in considering the commutator of an associative algebra as a Lie
bracket.
\smallbreak

\SkipTocEntry\subsection*{Pros}
A natural generalization of operads consists in authorizing multiple
outputs for its elements instead of a single one. This leads to the
theory of pros (this term is an abbreviation of product category). These
algebraic structures have been introduced by Mac Lane~\cite{McL65}.
Intuitively, a pro $\Pca$ is a set (or a space) of operators together
with two operations: an horizontal composition $*$ and a vertical
composition $\circ$. The first operation takes two operators $x$ and $y$
of $\Pca$ and builds a new one whose inputs (resp. outputs) are, from
left to right, those of $x$ and then those of $y$. The second operation
takes two operators $x$ and $y$ of $\Pca$ and produces a new one
obtained by plugging the outputs of $y$ onto the inputs of $x$. Basic
and modern references about pros are~\cite{Lei04} and~\cite{Mar08}.
\smallbreak

Like operads, pros can describe categories of algebras. Nevertheless, in
this case, pros can handle coproducts and can hence describe categories
of bialgebras. Consider the pro generated by the following three
operations:
\begin{equation}
\,,
\end{equation}
subjected to the following three relations:
\begin{equation}
    %
\,.
\end{equation}
The first (resp. second) one says that $\Product$ (resp. $\Coproduct$)
is associative (resp. coassociative). By seeing the operator $\omega$
as a map transposing its two inputs, the last one models
Relation~\eqref{equ:Hopf_relation}. Hence, this pro describes the
category of Hopf bialgebras.
\smallbreak

Another interaction between the theory of pros and combinatorics happens
when we consider presentations by generators and relations of pros
(see for instance~\cite{Laf11}). Recall that the symmetric group
$\SymmetricGroup_n$ is presented in the following way. It is generated
by symbols $\left\{s_i : 1 \leq i \leq n - 1\right\}$ whose elements are
called elementary transpositions. These generators are subjected to the
relations
\begin{subequations}
\begin{equation} \label{equ:symmetric_group_rel_1}
    s_i^2 = 1,
    \quad 1 \leq i \leq n - 1,
\end{equation}
\begin{equation} \label{equ:symmetric_group_rel_2}
    s_i s_j = s_j s_i,
    \quad 1 \leq i, j \leq n - 1 \mbox{ and } |i - j| \geq 2,
\end{equation}
\begin{equation} \label{equ:symmetric_group_rel_3}
    s_i s_{i + 1} s_i = s_i s_{i + 1} s_i,
    \quad 1 \leq i \leq n - 2.
\end{equation}
\end{subequations}
It is rather technical to show that $\SymmetricGroup_n$ admits the
stated presentation or, by going in the opposite direction, to show
that the group admitting the stated presentation is realized by
$\SymmetricGroup_n$. It is worth noting that there is a pro
$\K \Angle{\Per}$ of permutations offering a comfortable way to prove
these facts. Each permutation $\sigma$ of size $n \in \N$ is seen as an
operator with $n$ inputs and $n$ outputs, connecting each $i$th input
to the $\sigma(i)$th output. For instance, the permutation $42153$ is
seen as the element
\begin{equation}

\end{equation}
of $\K \Angle{\Per}$. The operations of pros, that are the horizontal
and vertical compositions, translate on permutations respectively as
the shifted concatenation $\Over$ of permutations and as the
composition $\circ$ of permutations. Therefore, $\K \Angle{\Per}$
contains all the symmetric groups $\SymmetricGroup_n$, $n \in \N$. Now,
one can ask about natural algebraic questions as finding a minimal
generating set of $\K \Angle{\Per}$. It is easy to show that the
singleton
\begin{equation}
    \GeneratingSet :=
    \left\{
    %
    \right\}
\end{equation}
is a minimal generating set of $\K \Angle{\Per}$. The unique element
$\Sfr$ of $\GeneratingSet$ encodes the permutation $21$. Now, by
finding the nontrivial relations satisfied by $\Sfr$~\cite{Laf03}, one
obtains the analogous relations of~\eqref{equ:symmetric_group_rel_1},
\eqref{equ:symmetric_group_rel_2},
and~\eqref{equ:symmetric_group_rel_3}, stated in the language of pros.
As a side remark, let us mention that the analogous relation
of~\eqref{equ:symmetric_group_rel_3} is axiomatic for pros. This
provides a nice strategy to establish the presentation of the symmetric
groups. Note that similar ideas work for establishing presentations
or realizations of other Coxeter groups.
\smallbreak

\SkipTocEntry\section*{Contributions}
Let us now present our contributions and the main results contained in
this dissertation. Before that, let us say a few words about the
organization of the text.
\smallbreak

\SkipTocEntry\subsection*{Global overview}
\begin{figure}[ht]
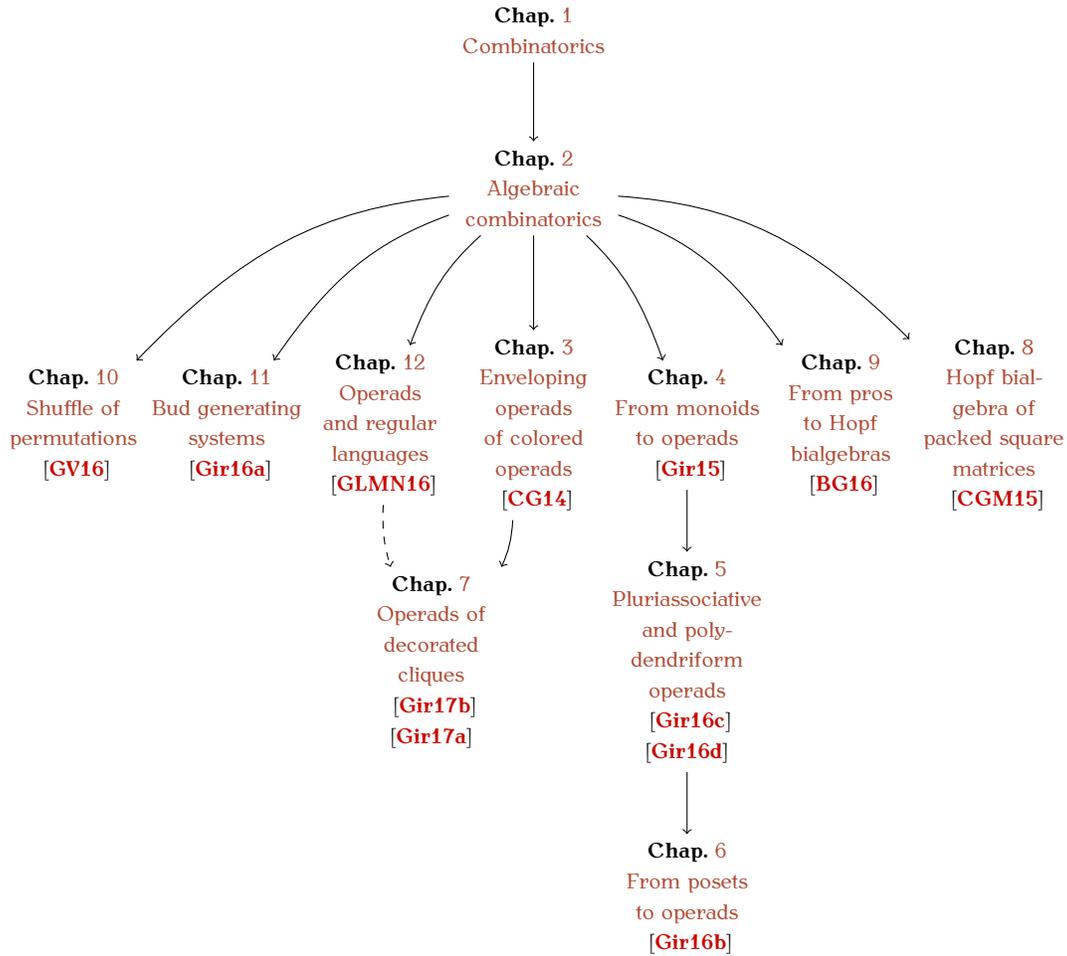

    \centering

    \caption[Diagram of the dependences between the chapters.]
    {Diagram of the dependences between the chapters. Each arrow
    $a \to b$ means that $b$ need some notions contained in~$a$.
    The dashed arrow means an optional dependence.}
    \label{fig:chapters}
\end{figure}

This text is divided into twelve chapters, the first two containing
preliminary notions, and the last ten containing original results
coming from published or submitted works. Figure~\ref{fig:chapters}
shows the diagram of dependences between the chapters and the references
to our work on which each chapter relies.
Our results fall into three categories: algebraic combinatorics,
enumerative combinatorics, and computer science.
\smallbreak

What follows is not a chapter-by-chapter summary. We follow the idea
to organize and present our contributions into the three categories
cited above. For this reason, a same chapter may appear several times
in the sequel.
\smallbreak

\SkipTocEntry\subsection*{Algebraic combinatorics}
Our main contributions in the field of algebraic combinatorics rely on
constructions, taking as input some algebraic structures, and outputting
other ones. Most of them are functorial and endow combinatorial sets
with (co)operations. We have presented above our point of view about the
advantages to endow objects with algebraic structures. Here, our
philosophy consists in designing general ways to achieve these goals.
For this reason, we create metatools (functorial constructions) whose
aim is to create tools (algebraic structures on combinatorial objects).
Let us list the main results, chapter by chapter, belonging to this
field.
\smallbreak

In {\bf Chapter~\ref{chap:enveloping}}, we introduce a tool to
facilitate the study of ns operads. This tools is a functor $\Hull$
from the category of ns colored operads to the category of ns
noncolored operads. It sends a ns colored operad $\Cca$ to the
smallest ns noncolored operad containing the elements of arities
greater than $1$ of $\Cca$. The ns operad $\Hull(\Cca)$ is realized
in terms of anticolored syntax trees labeled on $\Cca$, that are
particular syntax trees satisfying some conditions involving the
colors of $\Cca$. This construction is used to collect properties of
ns operads in the following way. Given a ns operad $\Oca$, finding a
ns colored operad $\Cca$ such that $\Hull(\Cca) = \Oca$ brings
information on $\Oca$. Indeed, some properties of $\Oca$ are implied
by properties of $\Cca$, such as the Hilbert series, the description
of suboperads and quotients, and presentations by generators and
relations. Due to the fact that a ns colored operad is a more
constrained structure than a noncolored one, it is in practice easier
to collect properties on $\Cca$ rather than on $\Oca$. These
techniques are illustrated to perform the study of the operad of
bicolored noncrossing configurations $\BNC$, an operad on some sorts
of noncrossing configurations~\cite{FN99}, introduced as a
generalization of the operads $\NCT$ of noncrossing trees and $\NCP$ of
noncrossing plants~\cite{Cha07}.
\smallbreak

The main contribution of {\bf Chapter~\ref{chap:monoids}} is a functor
$\T$ from the category of monoids to the category of operads. Given a
monoid $\Mca$, $\T \Mca$ is an operad of words on $\Mca$ seen as an
alphabet. The definitions of the partial compositions of $\T \Mca$
follow from the monoidal product of $\Mca$, and the symmetric groups act
on $\T \Mca$ by permuting the letters of the words. This functor is rich
from a combinatorial point of view since it leads to the construction of
several (ns) operads on combinatorial objects. Among others, $\T$ allows
to construct operads on endofunctions, parking functions, packed words,
permutations, and ns operads on planar rooted trees, $k$-ary trees (and
thus, $k$-Dyck paths, see~\cite{LPR15} for some structures on these),
Motzkin paths (the operad $\Motz$ appearing above is constructed in this
chapter), integer compositions, directed animals, and segmented
compositions. This construction $\T$ provides also alternative ways to
obtain the diassociative operad $\Dias$~\cite{Lod01} and the
triassociative operad $\Trias$~\cite{LR04}. By using rewriting
techniques, presentations of these operads are provided. We think that
there are a lot of other operads to be constructed through $\T$ on many
other families of combinatorial objects.
\smallbreak

We use in {\bf Chapter~\ref{chap:polydendr}} the functor $\T$ to define
generalizations of the diassociative operad depending on an integer
$\gamma \in \N$. These ns operads $\Dias_\gamma$ are realized in terms
of certain words on the alphabet $\{0, 1, \dots, \gamma\}$. Since the
diassociative operad is the Koszul dual of the dendriform operad
$\Dendr$~\cite{Lod01}, we obtain by Koszul duality the ns operads
$\Dendr_\gamma$, $\gamma \in \N$, each one being the Koszul dual of
$\Dias_\gamma$. These operads $\Dendr_\gamma$ are realized in terms of
binary trees with labeled edges endowed with tree shuffling operations.
The original motivation for this research direction is the following.
Dendriform algebras are algebraic structures consisting in two
operations $\LDendr$ and $\RDendr$ satisfying some axioms. As a
consequence of these axioms, the operation $\LDendr + \RDendr$ is
associative. This provides hence a framework to study associative
algebras by studying them as dendriform algebras through the definition
of two dendriform products such that their sum is the original product
of the algebra. The category of algebras described by $\Dendr_\gamma$
leads to a generalization of this method and forms a new device to study
associative algebras.
\smallbreak

In {\bf Chapter~\ref{chap:posets}}, we provide a functorial construction
associating with any finite poset $\Qca$ a ns operad $\As(\Qca)$. Under
some conditions on $\Qca$, $\As(\Qca)$ can be realized in terms of
Schröder trees labeled on $\Qca$ satisfying some conditions. The
original motivation for the introduction of this construction comes from
the previous chapter where two operads $\As_\gamma$ and $\DAs_\gamma$
were introduced. Both of these operads can be constructed as degenerate
cases of the construction $\As$ (by considering respectively trivial
orders and total orders). The question to study the operads $\As(\Qca)$
when $\Qca$ is nondegenerate is natural, and this leads to unexpected
results involving Koszul duality of ns operads. Indeed, the main result
here states that if $\Qca$ is a thin forest poset (a poset whose Hasse
diagram satisfies a certain property), the Koszul dual of $\As(\Qca)$ is
isomorphic to the ns operad $\As(\Qca^\perp)$ where $^\perp$ is an
involution on thin forest posets. The main contributions of this chapter
are of two kinds. First, new links between the theory of operads,
posets, and Koszul duality are developed (some different links
between these theories are included in~\cite{MY91,Val07,FFM16}).
Second, this work provides definitions of several algebraic
structures (as $\As(\Qca)$-algebras) having many associative
operations and generalizing associative algebras.
\smallbreak

Another functorial construction $\Cli$ is introduced in
{\bf Chapter~\ref{chap:cliques}}. It is defined from the category of
unitary magmas to the one of ns operads. Given a unitary magma $\Mca$,
$\Cli\Mca$ is a ns operad of regular polygons endowed with
configurations of arcs labeled on $\Mca$. The definitions of the partial
compositions of $\Cli\Mca$ follow from the magmatic product of $\Mca$.
This construction too is very rich from a combinatorial point of view
since it endows several families of configurations with ns operad
structures, as for instance, noncrossing configurations, Motzkin
configurations, Lucas configurations, and diagrams of involutions.
This leads to a new diagram of operads, each one realized in terms of
combinatorial objects. Moreover, the construction $\Cli$ allows to
define a large number of known operads in a unified way. For instance,
one can construct from $\Cli$ the operad $\BNC$ (and all its
suboperads as $\NCT$, $\NCP$, the dipterous operad~\cite{LR03,Zin12},
and the $2$-associative operad~\cite{LR06,Zin12}), the suboperad
$\FF_4$ of the operad of formal fractions $\FF$~\cite{CHN16}, the
operad of multi-tildes $\MT$~\cite{LMN13}, the operad of double
multi-tildes $\DMT$ (see Chapter~\ref{chap:languages}), and a ns
version of the gravity operad $\Grav$~\cite{Get94,AP15}.
\smallbreak

{\bf Chapter~\ref{chap:matrices}} is concerned with Hopf bialgebras and
more precisely with the definition of a Hopf bialgebra on packed square
matrices $\PM{k}$ depending on an integer parameter $k \in \N$. This
Hopf bialgebra is a generalization of $\FQSym$~\cite{DHT02} since the
subspace of $\PM{k}$ restricted to permutation matrices is isomorphic to
$\FQSym$. It contains also the Hopf bialgebra of uniform block
permutations $\UBP$~\cite{AO08}. Among the most notable facts, $\PM{k}$
contains a Hopf sub-bialgebra $\ASM$ involving alternating sign
matrices~\cite{MRR83}. This chapter provides also an algebraic point of
view of some statistics of alternating sign matrices by using the
associative algebra structure of~$\ASM$.
\smallbreak

Like the previous one, {\bf Chapter~\ref{chap:pros_bialgebras}} deals
with constructions of Hopf bialgebras. We introduce here a construction
$\PvH$ from stiff pros to Hopf bialgebras. A stiff pro is a quotient of
a free pro by a pro congruence satisfying some precise properties. This
construction is a generalization of the construction $\OvH$ associating
with an operad its natural Hopf bialgebra~\cite{Vdl04,CL07,ML14}.
Indeed, given an operad $\Oca$, one can construct a stiff pro
$\OvP(\Oca)$~\cite{Mar08} such that $\OvH(\Oca)$ is isomorphic to the
abelianization of $\PvH(\OvP(\Oca))$. In addition to creating a link
between the theories of pros and of Hopf bialgebras, the construction
$\PvH$ brings several new Hopf bialgebras on objects like forests of
trees with a fixed arity and heaps of pieces~\cite{Vie86}. It allows
also to construct some of the deformations of the noncommutative version
of the Fàa di Bruno Hopf bialgebra~\cite{BFK06} introduced by
Foissy~\cite{Foi08}.
\smallbreak

The main algebraic contribution of {\bf Chapter~\ref{chap:languages}}
concerns the introduction of a new category of algebraic objects,
the precompositions. These objects are used as inputs of a functorial
construction $\PrecompToOp$ producing ns operads. This construction is
used to provide alternative definitions of the operads $\MT$ and
$\Poset$ introduced in~\cite{LMN13} in the context of language theory
and the study of multi-tildes~\cite{CCM11}, and to construct new operads
$\DMT$ and $\Qoset$ as respective extensions of the last two.
\smallbreak

\SkipTocEntry\subsection*{Enumerative combinatorics}
As explained before, defining algebraic structures on combinatorial
objects leads to discover combinatorial properties on them. This form
a large part of our philosophy. We list here the main results, chapter
by chapter, in this context.
\smallbreak

As said above, {\bf Chapter~\ref{chap:enveloping}} presents a way to
study a ns operad $\Oca$ through a ns colored operad $\Cca$ satisfying
$\Hull(\Cca) = \Oca$. We highlight the fact that the Hilbert series of
$\Oca$ can be computed from the colored Hilbert series of $\Cca$ by
solving a system of equations, and, as the developed examples show, the
colored Hilbert series of $\Cca$ are simpler (rational) than the Hilbert
series of $\Oca$ (algebraic). Since the Hilbert series of a ns operad
and the generating series of its elements are the same series, this
offers a tool for enumeration. In this chapter, we enumerate bicolored
noncrossing configurations by using this technique.
\smallbreak

In {\bf Chapter~\ref{chap:buds}}, we work with formal power series on ns
colored operads to develop enumerative tools. We introduce a new kind of
formal grammar (see~\cite{Har78,HMU06}) generalizing both context-free
grammars of words and regular tree grammars~\cite{CDGJLLTT07}. These
grammars, called bud generating systems, allow to generate elements of a
ground ns colored operad. To enumerate the elements generated by a bud
generating system $\Bca$, we introduce three formal series on ns colored
operads: the hook generating series $\Hook(\Bca)$, the syntactic
generating series $\Synt(\Bca)$, and the synchronous generating series
$\Sync(\Bca)$. Each of these provide a different enumeration of the
elements generated by $\Bca$. For instance, $\Hook(\Bca)$ is a series
whose coefficients provide analogs of the hook-length statistics for
binary trees~\cite{Knu98}. Moreover, these series are defined through
operations on series: a pre-Lie product $\PreLieProduct$ and an
associative product $\Compo$. The example treated above about the
enumeration of balanced binary trees uses tools developed in this
chapter.
\smallbreak

\SkipTocEntry\subsection*{Computer science}
In this research some contributions to computer science and, more
precisely, to formal language theory and to computational complexity
theory have been developed. Let us summarize them.
\smallbreak

In {\bf Chapter~\ref{chap:shuffle}}, we consider the supershuffle
$\SuperShuffle$ of permutations introduced by Vargas~\cite{Var14}. This
operation is different from the shifted shuffle of permutations and can
be seen as an extension of the usual shuffle product $\shuffle$ on
words~\cite{EM53}. A classical question in algorithmic consists in
evaluating the complexity of recognizing words that are squares for this
operation. In other words, the problem amounts to decide if, given a
word $u$, there exists a word $v$ such that $u$ appears in
$v \shuffle v$. It is known from~\cite{RV13,BS14} that this problem is
$\NPC$. We ask in this chapter the analogous question for the
recognition of square permutations with respect to $\SuperShuffle$ and
show that this problem is also $\NPC$.
\smallbreak

We have explained that {\bf Chapter~\ref{chap:buds}} contains
enumerative results. In addition to this, the chapter introduces bud
generating systems as new kinds of grammars capable to generate any type
of combinatorial objects. Some elementary results about these grammars
are developed.
\smallbreak

As said before, {\bf Chapter~\ref{chap:languages}} provides algebraic
results. Nevertheless, the operads constructed here are intended to be
tools in formal language theory. A common research axis in this field
is to define a family of operations to express formal languages with the
smallest spatial complexity as possible. Multi-tildes~\cite{CCM11} have
been designed in this way. Here, we introduce an extension of
multi-tildes, namely the double multi-tildes, increasing their
expressive power. An operad $\DMT$ of double multi-tildes is constructed
and its action on languages is described. One of the main results of the
chapter is that every regular language can be expressed by the action of
a double multi-tilde seen as an operator of arity $n$ on $n$ languages
$\alpha_i$, $1 \leq i \leq n$, such that each $\alpha_i$ is empty or
contains one unique word of length $1$. However, this action is not
faithful, in the sense that there are different multi-tildes of $\DMT$
that act similarly on languages. For this reason, we introduce a
quotient $\Qoset$ of $\DMT$ such that elements are quasiorders. We show
that the action of $\Qoset$ on regular languages is faithful. This
establishes an unexpected link between quasiorders and regular
languages.
\medbreak

\part{Algebraic combinatorics}

\chapter{Combinatorics} \label{chap:combinatorics}
In combinatorics, counting how many objects a given family contains is
one of the most common, hard, and stimulating activities. Nevertheless,
even before trying to answer this kind of question, an important
preliminary and basic work consists in classifying the objects of a
family according to some of their particularities. The size of the
objects is, of course, one of these, but also, if we take the example of
permutations, the number of inversions and recoils are other features
which may be considered.
\smallbreak

Combinatorial collections are structures designed to work with such
structured sets of combinatorial objects. Roughly speaking, a
combinatorial collection is a set expressible by a disjoint union of
finite sets, indexed by a particular set $I$. Depending on $I$, these
collections are designed to represent various kinds of sets of objects.
For instance, when the indexing set $I$ is $\N$, one obtains graded
collections. These are sets endowed with a size function, a concept
fully developed in~\cite{FS09} under the name of combinatorial classes.
When the indexing set is $\N^2$, one obtains bigraded collections. These
collections provide a suitable framework to work with prographs and pros
(see Section~\ref{subsec:prographs} of this chapter and
Section~\ref{subsec:pros} of Chapter~\ref{chap:algebra}) since the
elements of these algebraic structures have an input and an output
arity. Moreover, when the indexing set is a set of words on a given
alphabet, one obtains colored collections. These collections provide a
suitable framework to work with colored syntax trees and colored operads
(see Section~\ref{subsubsec:colored_syntax_trees} of this chapter and
Section~\ref{subsubsec:colored_operads} of Chapter~\ref{chap:algebra})
since the elements of these algebraic structures have input and output
colors.
\smallbreak

There are other sensible tools to encode combinatorial collections.
One can cite species of structures introduced by Joyal~\cite{Joy81}
that allow to work with labeled objects. This theory has been
developed by the Quebec school of combinatorics~\cite{BLL98,BLL13}.
Species of structures are very good candidates to work with symmetric
operads~\cite{Men15} since the action of the symmetric group of a
symmetric operad is encapsulated into the action of the symmetric
group on an underlying species of structure. Another interesting way
to describe combinatorial objects passes through polynomial
functors~\cite{Koc09}.
\smallbreak

This machinery of combinatorial collections is applied in this work
mainly to rigorously define several families of trees. Let us remark
that this concept of tree encompasses a large range of quite different
combinatorial objects. For instance, in graph theory, trees are
connected acyclic graphs while in combinatorics, one encounters mostly
rooted trees. Among rooted trees, some of these can be planar (the order
of the children of a node is relevant) or not. In addition to this, the
internal nodes, the leaves, or the edges of the trees can be labeled,
and some conditions for the arities of their nodes can be imposed. One
of the first occurrences of the concept of tree came from the work of
Cayley~\cite{Cay1857}. Nowadays, trees appear among other in computer
science as data structures~\cite{Knu98,CLRS09}, in combinatorics in
relation with enumeration questions and Lagrange
inversion~\cite{Lab81,FS09}, and in algebraic combinatorics, where
several families of trees are endowed with algebraic
structures~\cite{LR98,HNT05,Cha08}. In our context, the most important
families of trees are the  syntax trees, which are kind of labeled
planar rooted trees. These trees are central objects in the study of
operads.
\smallbreak

This chapter is devoted to set the main definitions and notations about
combinatorics. In Section~\ref{sec:collections}, we introduce
combinatorial collections and structured combinatorial collections,
including the notion of posets and rewrite systems. Trees and syntax
trees are considered in Section~\ref{sec:trees}. We describe here
several families of trees and rewrite systems on syntax trees. Finally,
Section~\ref{sec:secondary_objects} contains a list of definitions of
combinatorial objects met afterwards in this dissertation.
\medbreak

\section{Structured collections} \label{sec:collections}
We introduce here the general notion of collection. Then, we consider
very usual concepts as magmas, monoids, posets, rewrite systems under
this context of collections.
\medbreak

\subsection{Collections}
After defining collections, combinatorial collections, (multi)graded
collections, and colored collections, a bunch of operations on graded
collections are reviewed.
\medbreak

\subsubsection{General collections}
Let $I$ be a nonempty set called \Def{index}. An \Def{$I$-collection} is
a set $C$ expressible as a disjoint union
\begin{equation} \label{equ:collection}
    C = \bigsqcup_{i \in I} C(i)
\end{equation}
where all $C(i)$, $i \in I$, are sets. All the elements of $C$ (resp.
$C(i)$ for an $i \in I$) are called \Def{objects} (resp.
\Def{$i$-objects}) of $C$. If $x$ is an $i$-object of $C$, we say that
the \Def{index} $\Index(x)$ of $x$ is $i$. When for all $i \in I$, all
$C(i)$ are finite sets, $C$ is \Def{combinatorial}. Besides, $C$ is
\Def{finite} if $C$ is finite as a set. The \Def{empty $I$-collection}
is the set $\emptyset$. When $I$ is a singleton, $C$ is \Def{simple}.
Any set can thus be seen as a simple collection and conversely.
\medbreak

A \Def{relation} on $C$ is a binary relation $\Rca$ on $C$ such that for
any objects $x$ and $y$ of $C$ satisfying $x \,\Rca\, y$, $x$ and $y$
have the same index. Let $C_1$ and $C_2$ be two $I$-collections. A map
$\phi : C_1 \to C_2$ is an \Def{$I$-collection morphism} if, for all
$x \in C_1$, $\Index(x) = \Index(\phi(x))$. We express by
$C_1 \simeq C_2$ the fact that there exists an isomorphism between $C_1$
and $C_2$. Besides, if for all $i \in I$, $C_1(i) \subseteq C_2(i)$,
$C_1$ is a \Def{subcollection} of $C_2$. For any $i \in I$, we can
regard each $C(i)$ as a subcollection of $C$ consisting in all its
$i$-objects.
\medbreak

Let us now consider particular $I$-collections for precise sets $I$.
Table~\ref{tab:collections} contains an overview of the properties that
such collections can satisfy.
\begin{table}[ht]
    \centering
    \begin{tabular}{|c|c|} \hline
        \multicolumn{2}{|c|}{{\bf Collections}} \\
        \multicolumn{2}{|c|}{
            {\it Combinatorial} \qquad
            {\it Finite} \qquad
            {\it Simple}}
            \\ \hline
        {\bf $k$-graded} & \multirow{2}{*}{{\bf Colored}} \\ \cline{1-1}
        {\bf $1$-graded} & \multirow{2}{*}{{\it Monochrome}} \\
        {\it Connected} \quad
        {\it Augmented} \quad
        {\it Monatomic} & \\ \hline
    \end{tabular}
    \bigbreak

    \caption[Properties of $I$-collections.]
    {The most common $I$-collections (in bold) and the properties
    (in italic) they can satisfy. The inclusions relations between
    these collections read from bottom to top. For instance,
    $1$-graded collections are particular $k$-graded collections which
    are particular collections.}
    \label{tab:collections}
\end{table}
\medbreak

\subsubsection{Graded collections} \label{subsubsec:graded_collections}
An $\N$-collection is called a \Def{graded collection}. If $C$ is a
graded collection, for any object $x$ of $C$, the \Def{size} $|x|$ of
$x$ is the integer $\Index(x)$. The map $|-| : C \to \N$ is the
\Def{size function} of~$C$.
\medbreak

Let us from now on assume that $C$ is a combinatorial graded collection.
The \Def{generating series} of $C$ is the series
\begin{equation} \label{equ:generating_series_graded_collection}
    \GenSeries_C(t) := \sum_{n \in \N} \# C(n) t^n,
\end{equation}
where $\# S$ denotes the cardinality of any finite set $S$. This formal
power series encodes the \Def{sequence of integers} associated with $C$,
that is the sequence $\left(\# C(n)\right)_{n \in \N}$. Observe that if
$C_1$ and $C_2$ are two combinatorial graded collections,
$C_1 \simeq C_2$ holds if and only if
$\GenSeries_{C_1}(t) = \GenSeries_{C_2}(t)$.
\medbreak

We say that $C$ is \Def{connected} if $C(0)$ is a singleton, and that
$C$ is \Def{augmented} if $C(0) = \emptyset$. Moreover, $C$ is
\Def{monatomic} if it is augmented and $C(1)$ is a singleton. We denote
by $\{\OneElement\}$ the graded collection such that $\OneElement$ is an
object satisfying $|\OneElement| = 0$. This collection is called the
\Def{unit collection}. Observe that $\{\OneElement\}$ is connected, and
that $C$ is connected if and only if there is a unique collection
morphism from $\{\OneElement\}$ to $C$. We denote by $\{\AtomElement\}$
the collection such that $\AtomElement$ is an \Def{atom}, that is an
object satisfying $|\AtomElement| = 1$. This collection is called the
\Def{neutral collection}. Observe that $\{\AtomElement\}$ is monatomic,
and that $C$ is monatomic if and only if $C$ is augmented and there is a
unique collection morphism from $\{\AtomElement\}$ to~$C$.
\medbreak

\subsubsection{Statistics and multigraded collections}
\label{subsubsec:multigraded_collections}
Let $C$ be a collection. A \Def{statistics} on $C$ is a map
$\Statistics : C \to \N$, associating a nonnegative integer value with
any object of $C$. A \Def{$k$-graded collection} (also called
\Def{multigraded collection}) is an $\N^k$-collection for an integer
$k \geq 1$. To not overload the notation, we denote by
$C(n_1, \dots, n_k)$ the subset $C((n_1, \dots, n_k))$ of any
$k$-graded collection $C$. These collections are useful to work with
objects endowed with many statistics. Indeed, if $x$ is an
$(n_1, \dots, n_k)$-object, one sets $\Statistics_i(x) := n_i$ for each
$1 \leq i \leq k$. This defines in this way $k$ statistics
$\Statistics_i : C \to \N$, $1 \leq i \leq k$. Besides, the
\Def{generating series} of a combinatorial $k$-graded $\N^k$-collection
$C$ is the series
\begin{equation} \label{equ:generating_series_multigraded_collection}
    \GenSeries_C(t_1, \dots, t_k)
    := \sum_{(n_1, \dots, n_k) \in \N^k}
    \# C(n_1, \dots, n_k) \; t_1^{n_1} \dots t_k^{n_k}.
\end{equation}
Of course, \eqref{equ:generating_series_graded_collection} is a
particular case of~\eqref{equ:generating_series_multigraded_collection}
when~$k = 1$.
\medbreak

\subsubsection{Colored collections}
\label{subsubsec:colored_collections}
Let $\CFr$ be a finite set, called \Def{set of colors}. A
\Def{$\CFr$-colored collection} $C$ is an $I$-collection such that
\begin{equation}
    I := \left\{
        (a, u) : a \in \CFr
        \mbox{ and } u \in \CFr^\ell \mbox{ for an } \ell \geq 1
    \right\}.
\end{equation}
In other terms, any object $x$ of $C$ has an index $(a, u)$. By setting
that the \Def{size} of $x$ is the length $|u|$ of $u$ (that is the
integer $\ell$ such that $u \in \CFr^\ell$), we can see $C$ as an
augmented graded collection. Moreover, the \Def{output color} of $x$ is
$\Out(x) := a$, and the \Def{word of input colors} of $x$ is
$\In(x) := u$. The \Def{$i$th input color} of $x$ is the $i$th letter
of $\In(x)$, denoted by $\In_i(x)$. We say that $C$ is \Def{monochrome}
if $\CFr$ is a singleton. For any nonnegative integer $k$, a
\Def{$k$-colored collection} is a $\CFr$-colored collection where
$\CFr$ is the set of integers $\{1, \dots, k\}$.  Assume now that
$\CFr = \{a_1, \dots, a_k\}$ and let
$\Xbb_\CFr := \left\{\VarX_{a_1}, \dots, \VarX_{a_k}\right\}$ and
$\Ybb_\CFr := \left\{\VarY_{a_1}, \dots, \VarY_{a_k}\right\}$ be two
alphabets of commutative letters. The \Def{generating series} of $C$ is
the series
\begin{equation}
    \GenSeries_C\left(\VarX_{a_1}, \dots, \VarX_{a_k},
    \VarY_{a_1}, \dots, \VarY_{a_k}\right) :=
    \sum_{x \in C} \VarX_{\Out(x)}
    \prod_{1 \leq i \leq |\In(x)|} \VarY_{\In_i(x)}.
\end{equation}
Observe that when $C$ is monochrome, the specialization
$\GenSeries_C(1, t)$ is the generating series of $C$ seen as a graded
collection.
\medbreak

\subsubsection{Products in collections}
\label{subsubsec:products_collections}
Let $C$ be an $I$-collection. A \Def{product} on $C$ is a partial map
\begin{equation}
    \Product : C^p \to C
\end{equation}
where $p \in \N \setminus \{0\}$. The \Def{arity} of $\Product$ is $p$.
Any product of arity $p$ can be seen as an operation taking $p$ elements
of $C$ as input and outputting one element of $C$. When, for any
$i \in I$ and any objects $x_1$, \dots, $x_p$ of $C(i)$ such that
$\Product(x_1, \dots, x_p)$ is defined,
\begin{math}
    \Product(x_1, \dots, x_p) \in C(i)
\end{math}
holds, we say that $\Product$ is \Def{internal}.
\medbreak

When $I$ is endowed with an associative binary product $\DPlus$, if for
any objects $x_1$, \dots, $x_p$ of $C$ such that
$\Product(x_1, \dots, x_p)$ is defined,
\begin{equation}
    \Index(\Product(x_1, \dots, x_p))
    =
    \Index(x_1) \DPlus \cdots \DPlus \Index(x_p),
\end{equation}
we say that $\Product$ is \Def{$\DPlus$-compatible}. In the particular
case where $C$ is a graded collection and $\Product$ is $+$-compatible,
$\Product$ is \Def{graded}.
\medbreak

Let us now assume that $C$ is a simple collection endowed with a binary
and total product $\Product$. In this case, $C$ is a \Def{magma}. We say
that $C$ is \Def{right cancellable} if for any $x, y, z \in C$, the
relation $y \Product x = z \Product x$ implies $y = z$. An element
$\Unit$ of $C$ is a \Def{unit} for $\Product$ if for all $x \in C$,
$x \Product \Unit = x = \Unit \Product x$. When $\Product$ admits a
unit, $C$ is a \Def{unitary magma}. If, additionally, the product
$\Product$ is associative, $C$ is a \Def{monoid}.
\medbreak

\subsubsection{Operations over graded collections}
\label{subsubsec:operations_graded_comb_collections}
We list here the most important operations that take as input graded
collections and output new ones. Most of these are binary or unary, and
under some precise conditions, they produce combinatorial collections.
In what follows, $C$, $C_1$, $C_2$, and $C_3$ are four graded
collections.
\medbreak

\paragraph{Suspension and augmentation}
For any $k \in \Z$, the \Def{$k$-suspension} of $C$ is the graded
collection $\Suspension_k(C)$ defined for all $n \in \N$ by
\begin{equation} \label{equ:suspension_operation}
    \left(\Suspension_k(C)\right)(n) :=
    \begin{cases}
        C(n - k) & \mbox{if } n - k \geq 0, \\
        \emptyset & \mbox{otherwise}.
    \end{cases}
\end{equation}
Observe that $\Suspension_1(\Suspension_{-1}(C))$ is the subcollection
$C \setminus C(0)$ of $C$, that is the augmented collection having the
objects of $C$ without its objects of size $0$. We call this collection
the \Def{augmentation} of~$C$ and we denote it by~$\Augmentation(C)$.
\medbreak

\paragraph{Sum}
The \Def{sum} of $C_1$ and $C_2$ is the graded collection
$C_1 + C_2$ such that, for all $n \in \N$,
\begin{equation} \label{equ:sum_operation}
    (C_1 + C_2)(n) := C_1(n) \sqcup C_2(n).
\end{equation}
In other words, each object of size $n$ of $C_1 + C_2$ is either an
object of size $n$ of $C_1$ or an object of size $n$ of $C_2$. Since the
sum operation~\eqref{equ:sum_operation} is defined through a disjoint
union, when the sets $C_1(n)$ and $C_2(n)$ are not disjoint, there are
in $(C_1 + C_2)(n)$ two copies of each element belonging to the
intersection $C_1(n) \cap C_2(n)$, one coming from $C_1(n)$, the other
from $C_2(n)$. Moreover, when $C_1$ and $C_2$ are combinatorial,
$C_1 + C_2$ is also combinatorial and its generating series satisfies
\begin{equation} \label{equ:generating_series_sum}
    \GenSeries_{C_1 + C_2}(t)
    = \GenSeries_{C_1}(t) + \GenSeries_{C_2}(t).
\end{equation}
The iterated version of the operation $+$ is denoted by $\bigsqcup$ in
the sequel.
\medbreak

\paragraph{Product}
The \Def{product} of $C_1$ and $C_2$ is the graded collection
$C_1 \times C_2$ such that, for all $n \in \N$,
\begin{equation} \label{equ:product_operation}
    (C_1 \times C_2)(n) :=
    \left\{
        (x_1, x_2) : x_1 \in C_1, x_2 \in C_2,
        \mbox{ and } |x_1| + |x_2| = n
    \right\}.
\end{equation}
In other words, each object of size $n$ of $C_1 \times C_2$ is an
ordered pair $(x_1, x_2)$ such that $x_1$ (resp. $x_2$) is an object of
$C_1$ (resp. $C_2$) and the sum of the sizes of $x_1$ and $x_2$ is~$n$.
Moreover, when $C_1$ and $C_2$ are combinatorial, $C_1 \times C_2$ is
also combinatorial and its generating series satisfies
\begin{equation} \label{equ:generating_series_product}
    \GenSeries_{C_1 \times C_2}(t)
    = \GenSeries_{C_1}(t) \GenSeries_{C_2}(t).
\end{equation}
\medbreak

\paragraph{Hadamard product}
The \Def{Hadamard product} of $C_1$ and $C_2$ is the graded collection
$C_1 \Hadamard C_2$ such that, for all $n \in \N$,
\begin{equation} \label{equ:product_operation}
    (C_1 \Hadamard C_2)(n) := C_1(n) \times C_2(n).
\end{equation}
In other words, each object of size $n$ of $C_1 \Hadamard C_2$ is an
ordered pair $(x_1, x_2)$ such that $x_1$ (resp. $x_2$) is an object of
size $n$ of $C_1$ (resp. $C_2$).
Moreover, when $C_1$ and $C_2$ are combinatorial, $C_1 \Hadamard C_2$ is
also combinatorial and its generating series satisfies
\begin{equation} \label{equ:generating_series_Hadamard_product}
    \GenSeries_{C_1 \Hadamard C_2}(t)
    = \GenSeries_{C_1}(t) \Hadamard \GenSeries_{C_2}(t)
    = \sum_{n \in \N} \# C_1(n) \# C_2(n) t^n.
\end{equation}
\medbreak

\paragraph{List operation}
For any $k \geq 0$, the \Def{$k$-list operation} applied to $C$ produces
the graded collection $\List_k(C)$ such that, for all $n \in \N$,
\begin{equation} \label{equ:k_list_operation}
    \left(\List_k(C)\right)(n) :=
    \left\{
        (x_1, \dots, x_k) : x_1, \dots, x_k \in C,
        \mbox{ and } |x_1| + \dots + |x_k| = n
    \right\}.
\end{equation}
In other words, each object of size $n$ of $\List_k(C)$ is a tuple
$(x_1, \dots, x_k)$ of objects of $C$ such that the sum of the sizes of
$x_1$, \dots, $x_k$ is $n$. When $C$ is combinatorial, $\List_k(C)$ is
also combinatorial and its generating series satisfies
\begin{equation}
    \GenSeries_{\List_k(C)}(t) = \GenSeries_C(t)^k.
\end{equation}
The \Def{list operation} applied to $C$ produces the graded collection
$\List(C)$ defined by
\begin{equation} \label{equ:list_operation}
    \List(C) := \bigsqcup_{k \in \N} \List_k(C).
\end{equation}
Moreover, when $C$ is combinatorial and augmented, $\List(C)$ is also
combinatorial (but not augmented) and its generating series satisfies
\begin{equation} \label{equ:generating_series_list}
    \GenSeries_{\List(C)}(t) =
    \frac{1}{1 - \GenSeries_C(t)}.
\end{equation}
Besides, for any $k \in \N$, we denote by $\List_{\geq k}(C)$ the graded
collection defined by
\begin{equation} \label{equ:list_operation}
    \List_{\geq k}(C) :=
    \bigsqcup_{\substack{
        \ell \in \N \\
        \ell \geq k
    }}
    \List_\ell(C).
\end{equation}
This notation ``$\List$'' comes from tensor algebras (see
Section~\ref{subsubsec:tensor_algebras} of Chapter~\ref{chap:algebra}).
\medbreak

\paragraph{Multiset operation}
For any $k \geq 0$, the \Def{$k$-multiset operation} applied to $C$
produces the graded collection $\Multiset_k(C)$ such that, for all
$n \in \N$,
\begin{equation} \label{equ:k_multiset_operation}
    \left(\Multiset_k(C)\right)(n) :=
    \left\{
        \lbag x_1, \dots, x_k \rbag : x_1, \dots, x_k \in C,
        \mbox{ and } |x_1| + \dots + |x_k| = n
    \right\}.
\end{equation}
In other words, each object of size $n$ of $\Multiset_k(C)$ is a
multiset $\lbag x_1, \dots, x_k \rbag$ of objects of $C$ such that the
sum of the sizes of $x_1$, \dots, $x_k$ is $n$. The
\Def{multiset operation} applied to $C$ produces the graded collection
$\Multiset(C)$ defined by
\begin{equation} \label{equ:multiset_operation}
    \Multiset(C) := \bigsqcup_{k \in \N} \Multiset_k(C).
\end{equation}
Moreover, when $C$ is combinatorial and augmented, $\Multiset(C)$ is
also combinatorial (but not augmented) and its generating series
satisfies
\begin{equation}
    \label{equ:generating_series_multiset}
    \GenSeries_{\Multiset(C)}(t) =
    \prod_{n \in \N \setminus \{0\}}
    \left(\frac{1}{1 - t^n}\right)^{\# C(n)}.
\end{equation}
This notation ``$\Multiset$'' comes from symmetric algebras (see
Section~\ref{subsubsec:symmetric_algebras} of
Chapter~\ref{chap:algebra}).
\medbreak

\paragraph{Set operation}
For any $k \geq 0$, the \Def{$k$-set operation} applied to $C$ produces
the graded collection $\Set_k(C)$ such that, for all $n \in \N$,
\begin{equation} \label{equ:k_set_operation}
    \left(\Set_k(C)\right)(n) :=
    \left\{
        \{x_1, \dots, x_k\} \subseteq C :
        |x_1| + \dots + |x_k| = n
    \right\}.
\end{equation}
In other words, each object of size $n$ of $\Set_k(C)$ is a set
$\{x_1, \dots, x_k\}$ of objects of $C$ such that the sum of the sizes
of $x_1$, \dots, $x_k$ is $n$. The \Def{set operation} applied to $C$
produces the graded collection $\Set(C)$ defined by
\begin{equation} \label{equ:set_operation}
    \Set(C) := \bigsqcup_{k \in \N} \Set_k(C).
\end{equation}
Moreover, when $C$ is combinatorial, $\Set(C)$ is also combinatorial and
its generating series satisfies
\begin{equation} \label{equ:generating_series_set}
    \GenSeries_{\Set(C)}(t) =
    \prod_{n \in \N \setminus \{0\}}
    \left(1 + t^n\right)^{\# C(n)}.
\end{equation}
Unlike the cases of the list and multiset operations, $\Set(C)$ is a
combinatorial collection without requiring that $C$ is augmented.
This notation ``$\Set$'' comes from exterior algebras (see
Section~\ref{subsubsec:exterior_algebras} of
Chapter~\ref{chap:algebra}).
\medbreak

\paragraph{Composition product}
For any $k \geq 0$ and graded collections $C_1$, \dots, $C_k$, the
\Def{homogeneous composition} of $C$ with $C_1$, \dots, $C_k$ is the
graded collection $C \circ \left[C_1, \dots, C_k\right]$ such that, for
all $n \in \N$,
\begin{equation} \label{equ:homogeneous_composition_operation}
    \left(C \circ \left[C_1, \dots, C_k\right]\right)(n) :=
    \bigsqcup_{x \in C(k)}
    \{(x, (y_1, \dots, y_k)) : y_i \in C_i, 1 \leq i \leq k,
    \mbox{ and } |y_1| + \dots + |y_k| = n\}.
\end{equation}
In other words, each object of size $n$ of $C \circ [C_1, \dots, C_k]$
is an ordered pair $(x, (y_1, \dots, y_k))$ where $x$ is an object of
$C$ of size $k$, and $(y_1, \dots, y_k)$ is a tuple such that each $y_i$
is an object of $C_i$, $1 \leq i \leq k$, and the sum of the sizes of
these objects $y_i$ is~$n$. The \Def{composition} of $C_1$ and $C_2$ is
the graded collection $C_1 \circ C_2$ such that, for all $n \in \N$,
\begin{equation} \label{equ:composition_operation}
    (C_1 \circ C_2)(n) :=
    \bigsqcup_{k \in \N}
    C_1 \circ
    \underbrace{\left[C_2, \dots, C_2\right]}
    _{k \mbox{ \footnotesize terms}}.
\end{equation}
Moreover, when $C_1$ and $C_2$ are both combinatorial and $C_2$ is
augmented, $C_1 \circ C_2$ is also combinatorial (but not necessarily
augmented) and its generating series satisfies
\begin{equation} \label{equ:generating_series_composition}
    \GenSeries_{C_1 \circ C_2}(t) =
    \GenSeries_{C_1}(\GenSeries_{C_2}(t)).
\end{equation}
\medbreak

\subsection{Main collections}
We define, in some cases by using the operations of
Section~\ref{subsubsec:operations_graded_comb_collections}, some usual
graded combinatorial collections. At the same time, we set here our main
notations and definitions about their objects.
\medbreak

\subsubsection{Natural numbers} \label{subsubsec:natural_numbers}
We can regard the set $\N$ as the graded collection satisfying
$\N(n) := \{n\}$ for all $n \in \N$. Hence,
$\List(\{\AtomElement\}) \simeq \N$ for an atom $\AtomElement$.
Moreover, for any $k \in \N$, let $\N_{\geq k}$ be the graded collection
defined by
\begin{equation}
    \N_{\geq k} := \Suspension_k\left(\Suspension_{-k}(\N)\right).
\end{equation}
By definition of the suspension operation over graded collections,
$\N_{\geq k}$ is the set of all integers greater than or equal to $k$.
Observe that $\N_{\geq 1} = \Augmentation(\N)$. The generating series of
$\N_{\geq k}$ satisfies
\begin{equation} \label{equ:generating_series_natural_numbers}
    \GenSeries_{\N_{\geq k}}(t) = \frac{t^k}{1 - t}
    = t^k + t^{k + 1} + t^{k + 2} + \cdots.
\end{equation}
Observe also that the list operation over graded collections can be
expressed as a composition involving $\N$ since
\begin{equation}
    \List(C) \simeq \N \circ C
\end{equation}
for any augmented combinatorial graded collection $C$. Let also, for any
$x, z \in \N$, the subcollection
$[x, z] := \{y \in \N : x \leq y \leq z\}$, and $[x] := [1, x]$. These
examples of graded collections are among the simplest nontrivial ones.
\medbreak

It is time to provide some notations about natural numbers. For any
multiset $S := \lbag s_1, \dots, s_n \rbag$ of elements of
$\N$, we denote by $\sum S$ the sum $s_1 + \dots + s_n$ of its elements.
We moreover denote by $S!$ the \Def{multinomial coefficient}
\begin{equation}
    S! :=
    \binom{\sum S}{s_1, \dots, s_n} =
    \frac{(s_1 + \dots + s_n)!}{s_1! \dots s_n!}.
\end{equation}
\medbreak

\subsubsection{Words} \label{subsubsec:words}
Let $A$ be an \Def{alphabet}, that is a set whose elements are called
\Def{letters}. One can see $A$ as a graded collection wherein all
letters are atoms. In this case, we denote by $A^*$ the graded
collection $\List(A)$. By definition, the objects of $A^*$ are finite
sequences of elements of $A$. We call \Def{words} on $A$ these objects.
When $A$ is finite, $A^*$ is combinatorial and it follows
from~\eqref{equ:generating_series_list} that the generating series of
$A^*$ is
\begin{equation}
    \GenSeries_{A^*}(t)
    = \sum_{n \in \N} m^n t^n
    = 1 + m t + m^2 t^2 + m^3 t^3 + \cdots,
\end{equation}
where $m := \# A$. If $u := (a_1, \dots, a_n)$ is a word on $A$, it
follows from the definition of $A^*$ that the size $|u|$ of $u$ is $n$.
The \Def{$i$th letter} of $u$ is $a_i$ and is denoted by $u(i)$ (and
also denoted by $u_i$ in some contexts). For any letter $\Bsf \in A$,
the \Def{number of occurrences} $|u|_\Bsf$ of $\Bsf$ in $u$ is the
cardinality of the set $\{i \in [|u|] : u(i) = \Bsf\}$. The unique word
on $A$ of size $0$ is denoted by $\epsilon$ and is called
\Def{empty word}. The subcollection $A^+ := \Augmentation(A^*)$ of $A^*$
contains all nonempty words on $A$. For any $n \in \N$, $A^n$ denote the
subcollection $A^*(n)$ of $A^*$. When $A$ is endowed with a total order
$\Ord$ and $u$ is nonempty, $\max_\Ord(u)$ is the greatest letter
appearing in $u$ with respect to $\Ord$. Moreover, an \Def{inversion} of
$u$ is a pair $(i, j)$ such that $i < j$, $u(i) \ne u(j)$, and
$u(j) \Ord u(i)$. Given two words $u$ and $v$ on $A$, the
\Def{concatenation} of $u$ and $v$ is the word $u \Conc v$ containing
from left to right the letters of $u$ and then the ones of $v$. If $u$
can be expressed as $u = u_1 \Conc u_2 \Conc u_3$ where
$u_1, u_2, u_3 \in A^*$, we say that $u_1$ (resp. $u_3$, $u_2$) is a
\Def{prefix} (resp. \Def{suffix}, \Def{factor}) of $u$. We denote by
$u \PrefixOrder v$ (resp. $u \SuffixOrder v$, $u \FactorOrder v$) the
fact that $u$ is a prefix (resp. suffix, factor) of $v$. For any subset
$P := \{p_1 \leq \cdots \leq p_k\}$ of $[|u|]$, $u_{|P}$ is the word
$u(p_1) \dots u(p_k)$. Moreover, when $v$ is a word such that there
exists $P \subseteq [|u|]$ satisfying $v = u_{|P}$, $v$ is a
\Def{subword} of $u$. The \Def{commutative image} of $u$ is the multiset
$\lbag u(i) : i \in [|u|]\rbag$. Given two words $u$ and $v$ of the same
size $n$, the \Def{Hamming distance} $\Hamming(u, v)$ between $u$ and
$v$ is the number of integers $i \in [n]$ such that $u(i) \ne v(i)$. A
\Def{language} on $A$ is subcollection of $A^*$. A language $\Lca$ on
$A$ is \Def{prefix} if for all $u \in \Lca$ and $v \in A^*$,
$v \PrefixOrder u$ implies~$v \in \Lca$.
\medbreak

\subsubsection{Integer compositions}
\label{subsubsec:integer_compositions}
By regarding the set $\N$ as a graded collection as explained in
Section~\ref{subsubsec:natural_numbers}, let $\Compositions$ be the
combinatorial graded collection $\List\left(\N_{\geq 1}\right)$. It
follows from~\eqref{equ:generating_series_list}
and~\eqref{equ:generating_series_natural_numbers} that the generating
series of $\Compositions$ is
\begin{equation}
    \GenSeries_{\Compositions}(t)
    = \frac{1 - t}{1 - 2t}
    = 1 + \sum_{n \geq 1} 2^{n - 1} t^n
    = 1 + t + 2t^2 + 4t^3 + 8t^4 + 16t^5 + 32t^6 + 64t^7
    + \cdots.
\end{equation}
By definition, the objects of $\Compositions$ are finite sequences of
positive numbers. We call \Def{integer compositions} (or, for short,
\Def{compositions}) these objects. If
$\LambdaB := (\LambdaB_1, \dots, \LambdaB_k)$ is a composition, it
follows from the definition of $\Compositions$ that the size
$|\LambdaB|$ of $\LambdaB$ is $\LambdaB_1 + \dots + \LambdaB_k$. The
\Def{length} $\Length(\LambdaB)$ of $\LambdaB$ is $k$, and for any
$i \in [\Length(\LambdaB)]$, the \Def{$i$th part} of $\LambdaB$ is
$\LambdaB_i$. The unique composition of size $0$ is denoted by
$\epsilon$ and is called \Def{empty composition} (even if $\epsilon$ is
already used to express the empty word, this overloading of notation
is not a problem in practice).
\medbreak

The \Def{descents set} of $\LambdaB$ is the set
\begin{equation}
    \Des(\LambdaB) :=
    \{\LambdaB_1, \LambdaB_1 + \LambdaB_2, \dots,
    \LambdaB_1 + \LambdaB_2 + \dots + \LambdaB_{k - 1}\}.
\end{equation}
For instance, $\Des(4131) = \{4, 5, 8\}$. Moreover, for any word $u$
defined on an alphabet $A$ equipped with a total order $\Ord$, the
\Def{composition} $\Cmp(u)$ of $u$ is the composition of size $|u|$
defined by
\begin{equation}
    \Cmp(u) := \left(|u_1|, \dots, |u_k|\right),
\end{equation}
where $u = u_1 \Conc \dots \Conc u_k$ is the factorization of $u$ in
longest nondecreasing factors (with respect to the order $\Ord$). For
instance, if
\begin{math}
    u :=
    \textcolor{Col1}{a_2 a_2 a_3}
    \textcolor{Col4}{a_1 a_3}
    \textcolor{Col1}{a_2}
    \textcolor{Col4}{a_1 a_2}
\end{math}
is a word on the alphabet $A := \{a_1, a_2, a_3\}$ ordered by
$a_1 \Ord a_2 \Ord a_3$, $\Cmp(u) = 3212$. When $\# A \geq 2$, this map
$\Cmp$ is a surjective collection morphism from $A^*$ to
$\Compositions$.
\medbreak

Integer compositions are drawn as \Def{ribbon diagrams} in the following
way. For each part $\LambdaB_i$ of $\LambdaB$, we draw a horizontal line
of $\LambdaB_i$ boxes. These lines are organized so that the line for
the first part of $\LambdaB$ is the uppermost, and the first box of the
line of the part $\LambdaB_{i + 1}$ is glued below the last box of the
line of the part $\LambdaB_i$, for all $i \in [\Length(\LambdaB) - 1]$.
For instance, the ribbon diagram of the composition $4131$ is
\begin{equation}
    \begin{tikzpicture}[Centering]
        \node[Box](1)at(0,0){};
        \node[Box,right of=1,node distance=0.25cm](2){};
        \node[Box,right of=2,node distance=0.25cm](3){};
        \node[Box,right of=3,node distance=0.25cm](4){};
        \node[Box,below of=4,node distance=0.25cm](5){};
        \node[Box,below of=5,node distance=0.25cm](6){};
        \node[Box,right of=6,node distance=0.25cm](7){};
        \node[Box,right of=7,node distance=0.25cm](8){};
        \node[Box,below of=8,node distance=0.25cm](9){};
    \end{tikzpicture}\,.
\end{equation}
\medbreak

\subsubsection{Integer partitions} \label{subsubsec:integer_partitions}
Again by regarding the set $\N$ as a graded collection as considered in
Section~\ref{subsubsec:natural_numbers}, let $\IntPart$ be the graded
combinatorial collection $\Multiset\left(\N_{\geq 1}\right)$. Since
$\# \N_{\geq 1}(n) = 1$ for all $n \geq 1$, it follows
from~\eqref{equ:generating_series_multiset} that the generating series
of $\IntPart$ is
\begin{equation}
    \GenSeries_{\IntPart}(t)
    = \prod_{n \geq 1} \frac{1}{1 - t^n}
    = 1 + t + 2t^2 + 3t^3 + 5t^4 + 7t^5 + 11t^6 + 15t^7 + 22t^8
    + \cdots.
\end{equation}
By definition, the objects of $\IntPart$ are finite multisets of
positive integers. We call \Def{integer partitions} (or, for short,
\Def{partitions}) these objects. As a consequence of the definition of
$\IntPart$, the size $|\lambda|$ of any partition $\lambda$ is the sum
of the integers appearing in the multiset $\lambda$. Due to the
definition of partitions as multisets, we can present a partition as an
ordered sequence of positive integers with respect to any total order on
$\N_{\geq 1}$. For this reason, we denote any partition $\lambda$ by a
nondecreasing sequence $(\lambda_1, \dots, \lambda_k)$ of positive
integers (that is, $\lambda_i \geq \lambda_{i + 1}$ for all
$i \in [k - 1]$). Under this convention, the \Def{length}
$\Length(\lambda)$ of $\lambda$ is $k$, and for any
$i \in [\Length(\lambda)]$, the \Def{i$th$ part} of $\lambda$
is~$\lambda_i$.
\medbreak

\subsubsection{Permutations and colored permutations}
\label{subsubsec:permutations}
A \Def{permutation} of \Def{size} $n$ is a bijection $\sigma$ from $[n]$
to $[n]$. The combinatorial graded collection of all permutations is
denoted by $\SymmetricGroup$. The generating series of
$\SymmetricGroup$ is
\begin{equation}
    \GenSeries_\SymmetricGroup(t)
    = \sum_{n \in \N} n! t^n
    = 1 + t + 2t^2 + 6t^3 + 24t^4 + 120t^5 + 720t^6 + 5040t^7 +
    40320t^8 + \cdots.
\end{equation}
Any permutation $\sigma$ of $\SymmetricGroup(n)$ is denoted as a word
$\sigma(1) \dots \sigma(n)$ on $\N_{\geq 1}$. Under this convention, a
permutation of size $n$ is a word on the alphabet $[n]$ with exactly one
occurrence of each letter of~$[n]$. The composition operation $\circ$ of
maps forms a binary internal operation on~$\SymmetricGroup$.
\medbreak

A \Def{descent} of $\sigma$ is a position $i \in [|\sigma| - 1]$ such
that $\sigma(i) > \sigma(i + 1)$. The set of all descents of $\sigma$ is
denoted by $\Des(\sigma)$. A \Def{coinversion} of $\sigma$ is an ordered
pair of letters $(\Asf, \Bsf)$ occurring in $\sigma$ such that
$\Asf < \Bsf$ and the position of $\Asf$ is greater than the position of
$\Bsf$ in $\sigma$. The set of all coinversions of $\sigma$ is denoted
by $\CoInv(\sigma)$. For any word $u$ defined on an alphabet $A$
equipped with a total order $\Ord$, the \Def{standardized} $\Std(u)$ of
$u$ is the permutation of size $|u|$ having the same inversions as the
ones of $u$. In other terms $\Std(u)$ has its letters in the same
relative order as those of $u$, with respect to $\Ord$, where equal
letters of $u$ are ordered from left to right as the smallest to the
greatest. For example, by considering the alphabet $\N$ equipped with
the natural order of integers, $\Std(211241) = 412563$. This map $\Std$
is a surjective collection morphism from $\N^*$ to~$\SymmetricGroup$.
\medbreak

This collection $\SymmetricGroup$ admits the following straightforward
generalization. For any $\ell \geq 1$, let $\SymmetricGroup^{(\ell)}$ be
the set of all pairs $(\sigma, u)$ where $\sigma$ is a permutation and
$u$ is a word of $[\ell]^{|\sigma|}$. We call this object an
\Def{$\ell$-colored permutation}. The \Def{size} of $(\sigma, u)$ in
$\SymmetricGroup^{(\ell)}$ is the size of $\sigma$ in $\SymmetricGroup$.
\medbreak

\subsubsection{Binary trees} \label{subsubsec:binary_trees}
Let $\BinaryTrees_\Leaf$ be the combinatorial graded collection
satisfying the relation
\begin{equation} \label{equ:algebraic_relation_binary_trees_leaves}
    \BinaryTrees_\Leaf = \{\Leaf\} +
    \{\Node\} \times {\BinaryTrees_\Leaf}^2,
\end{equation}
where $\Leaf$ is an atomic object called \Def{leaf} and $\Node$ is an
object of size $0$ called \Def{internal node}. We call \Def{binary tree}
each object of $\BinaryTrees_\Leaf$. By definition, a binary tree $\Tfr$
is either the leaf $\Leaf$ or an ordered pair
$(\Node, (\Tfr_1, \Tfr_2))$ where $\Tfr_1$ and $\Tfr_2$ are binary
trees. Observe that this description of binary trees is recursive. For
instance,
\begin{equation} \label{equ:examples_binary_trees}
    \Leaf, \quad
    (\Node, (\Leaf, \Leaf)), \quad
    (\Node, ((\Node, (\Leaf, \Leaf)), \Leaf)), \quad
    (\Node, (\Leaf, (\Node, (\Leaf, \Leaf)))), \quad
    (\Node, ((\Node, (\Leaf, \Leaf)), (\Node, (\Leaf, \Leaf)))),
\end{equation}
are binary trees.  If $\Tfr$ is a binary tree different from the leaf,
by definition, $\Tfr$ can be expressed as
$\Tfr = (\Node, (\Tfr_1, \Tfr_2))$ where $\Tfr_1$ and $\Tfr_2$ are two
binary trees. In this case, $\Tfr_1$ (resp. $\Tfr_2$) is the
\Def{left subtree} (resp. \Def{right subtree}) of~$\Tfr$. By drawing
each leaf by $\LeafPic$ and each binary tree with at least one internal
node by an internal node $\NodePic$ attached below it, from left to
right, to its left and right subtrees by means of edges $\EdgePic$, the
binary trees of~\eqref{equ:examples_binary_trees} are depicted by
\begin{equation}
    \LeafPic, \quad
    \begin{tikzpicture}[xscale=.2,yscale=.17,Centering]
        \node[Leaf](0)at(0.00,-1.50){};
        \node[Leaf](2)at(2.00,-1.50){};
        \node[Node](1)at(1.00,0.00){};
        \draw[Edge](0)--(1);
        \draw[Edge](2)--(1);
        \node(r)at(1.00,1.75){};
        \draw[Edge](r)--(1);
    \end{tikzpicture}\,, \quad
    \begin{tikzpicture}[xscale=.16,yscale=.14,Centering]
        \node[Leaf](0)at(0.00,-3.33){};
        \node[Leaf](2)at(2.00,-3.33){};
        \node[Leaf](4)at(4.00,-1.67){};
        \node[Node](1)at(1.00,-1.67){};
        \node[Node](3)at(3.00,0.00){};
        \draw[Edge](0)--(1);
        \draw[Edge](1)--(3);
        \draw[Edge](2)--(1);
        \draw[Edge](4)--(3);
        \node(r)at(3.00,2){};
        \draw[Edge](r)--(3);
    \end{tikzpicture}\,, \quad
    \begin{tikzpicture}[xscale=.16,yscale=.14,Centering]
        \node[Leaf](0)at(0.00,-1.67){};
        \node[Leaf](2)at(2.00,-3.33){};
        \node[Leaf](4)at(4.00,-3.33){};
        \node[Node](1)at(1.00,0.00){};
        \node[Node](3)at(3.00,-1.67){};
        \draw[Edge](0)--(1);
        \draw[Edge](2)--(3);
        \draw[Edge](3)--(1);
        \draw[Edge](4)--(3);
        \node(r)at(1.00,2){};
        \draw[Edge](r)--(1);
    \end{tikzpicture}\,, \quad
    \begin{tikzpicture}[xscale=.15,yscale=.12,Centering]
        \node[Leaf](0)at(0.00,-4.67){};
        \node[Leaf](2)at(2.00,-4.67){};
        \node[Leaf](4)at(4.00,-4.67){};
        \node[Leaf](6)at(6.00,-4.67){};
        \node[Node](1)at(1.00,-2.33){};
        \node[Node](3)at(3.00,0.00){};
        \node[Node](5)at(5.00,-2.33){};
        \draw[Edge](0)--(1);
        \draw[Edge](1)--(3);
        \draw[Edge](2)--(1);
        \draw[Edge](4)--(5);
        \draw[Edge](5)--(3);
        \draw[Edge](6)--(5);
        \node(r)at(3.00,2.5){};
        \draw[Edge](r)--(3);
    \end{tikzpicture}\,.
\end{equation}
By definition of the sum and the product operations over graded
collections, the size of a binary tree $\Tfr$ satisfies
\begin{equation} \label{equ:size_binary_tree_leaves}
    |\Tfr| =
    \begin{cases}
        1 & \mbox{if } \Tfr = \Leaf, \\
        |\Tfr_1| + |\Tfr_2|
            & \mbox{otherwise (} \Tfr = (\Node, (\Tfr_1, \Tfr_2))
                \mbox{)}.
    \end{cases}
\end{equation}
In other words, the size of $\Tfr$ is the number of occurrences of
$\Leaf$ it contains. Since $\GenSeries_{\{\Leaf\}}(t) = t$ and
$\GenSeries_{\{\Node\}}(t) = 1$, it follows
from~\eqref{equ:generating_series_sum}
and~\eqref{equ:generating_series_product} that the generating series of
$\BinaryTrees_\Leaf$ satisfies the quadratic algebraic equation
\begin{equation} \label{equ:algebraic_relation_binary_trees_leaves}
    t - \GenSeries_{\BinaryTrees_\Leaf}(t)
    + \GenSeries_{\BinaryTrees_\Leaf}(t)^2 = 0.
\end{equation}
The unique solution having a combinatorial meaning
of~\eqref{equ:algebraic_relation_binary_trees_leaves} is
\begin{equation}
    \GenSeries_{\BinaryTrees_\Leaf}(t) =
    \frac{1 - \sqrt{1 - 4t}}{2} =
    \sum_{n \in \N_{\geq 1}}
    \frac{1}{n} \binom{2n - 2}{n - 1} t^n
\end{equation}
The sequence of integers associated with $\BinaryTrees_\Leaf$ begins by
\begin{equation}
    1, 1, 2, 5, 14, 42, 132, 429,
\end{equation}
and is Sequence~\OEIS{A000108} of~\cite{Slo}. These numbers are known
as \Def{Catalan numbers}.
\medbreak

\subsection{Posets on collections} \label{subsec:posets}
We consider now collections endowed with partial order relations
compatible with their indexations. Such structures are important in
combinatorics since they lead for instance to the construction of
alternative bases of combinatorial spaces (see
Section~\ref{subsec:change_basis_posets} of Chapter~\ref{chap:algebra}).
We provide general definitions about posets and consider as examples
three important ones: the cube, Tamari, and right weak order posets.
\medbreak

\subsubsection{Elementary definitions}
\label{subsubsec:definitions_posets}
An \Def{$I$-poset} is a pair $(\Qca, \Ord_\Qca)$ where $\Qca$ is an
$I$-collection and $\Ord_\Qca$ is both a relation on $\Qca$ (recall
that relations on collections preserve the indices) and a partial order
relation. For any property $P$ of collections, we say that
$(\Qca, \Ord)$ \Def{satisfies the property $P$} if, as a collection,
$\Qca$ satisfies $P$. Observe in particular that our terminology
concerning graded posets differs from the classical one~\cite{Sta11}
(where a poset is graded when all its maximal chains have the same
length). Moreover, simple posets are usual posets (that are sets
endowed with partial order relations, without extra structure).
\medbreak

The \Def{strict order relation} of $\Ord$ is the binary relation
$\OrdStrict$ on $\Qca$ satisfying, for all $x, y \in \Qca$,
$x \OrdStrict y$ if $x \Ord y$ and $x \ne y$. The \Def{interval}
between two objects $x$ and $z$ of $\Qca$ is the set
$[x, z] := \{y \in \Qca : x \Ord_\Qca y \Ord_\Qca z\}$. When all
intervals of $\Qca$ are finite, $\Qca$ is \Def{locally finite}. Observe
that when $\Qca$ is combinatorial, $\Qca$ is locally finite. When
$\Qca$ is finite, the number of intervals of $\Qca$ is finite and is
denoted by $\NbInterv(\Qca)$. For any $i \in I$, an object $x$ of
$\Qca(i)$ is a \Def{greatest} (resp. \Def{least}) \Def{element} if for
all $y \in \Qca(i)$, $y \Ord_\Qca x$ (resp. $x \Ord_\Qca y$). Moreover,
for any $i \in I$, an object $x$ of $\Qca(i)$ is a \Def{maximal} (resp.
\Def{minimal}) \Def{element} if for all $y \in \Qca(i)$, $x \Ord_\Qca y$
(resp. $y \Ord_\Qca x$) implies $x = y$. The partial binary operation
$\min$ (resp. $\max$) with respect to the order $\Ord_\Qca$ is denoted
by $\Min_\Qca$ (resp. $\Max_\Qca$). If $x$ and $y$ are two objects of
$\Qca$, $y$ \Def{covers} $x$ if $x \Ord_\Qca y$ and $[x, y] = \{x, y\}$.
Two objects $x$ and $y$ are \Def{comparable} (resp. \Def{incomparable})
in $\Qca$ if $x \Ord_\Qca y$ or $y \Ord_\Qca x$ (resp. neither
$x \Ord_\Qca y$ nor $y \Ord_\Qca x$ holds). If for any $i \in I$
and any $i$-objects $x$ and $y$ of $\Qca$, $x$ and $y$ are comparable,
$\Qca$ is a \Def{total order}. A \Def{chain} of $\Qca$ is a sequence
$(x_1, \dots, x_k)$ such that $x_j \Ord_\Qca x_{j + 1}$ for all
$j \in [k - 1]$. An \Def{antichain} of $\Qca$ is a subset of pairwise
incomparable elements of $\Qca$. A \Def{linear extension} of $\Qca$ is
an $I$-poset $(\Qca', \Ord_\Qca')$ being a total order and such that
$\Ord_\Qca'$ contains $\Ord_\Qca$ as a relation. An \Def{order filter}
of $\Qca$ is a subset $\Fca$ of $\Qca$ such that for all $x \in \Fca$
and all $y \in \Qca$ satisfying $x \Ord_\Qca y$, $y$ is in~$\Fca$. For
any $i \in I$, the \Def{$i$-subposet} of $\Qca$ is the poset obtained by
restricting $\Ord_\Qca$ on $\Qca(i)$. The \Def{Hasse diagram} of
$(\Qca, \Ord_\Qca)$ is the directed graph having $\Qca$ as set of
vertices and all the pairs $(x, y)$ where $y$ covers $x$ as set of arcs.
\medbreak

We shall define posets $\Qca$ by drawing Hasse diagrams, where minimal
elements are drawn uppermost and vertices are labeled by the elements
of $\Qca$. For instance, the Hasse diagram
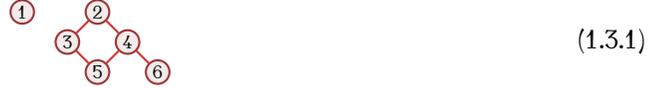
\begin{equation}
    \begin{tikzpicture}[xscale=.4,yscale=.4,Centering]
        \node[PosetVertex](1)at(-1.5,1){\begin{math}1\end{math}};
        \node[PosetVertex](2)at(1,1){\begin{math}2\end{math}};
        \node[PosetVertex](3)at(0,0){\begin{math}3\end{math}};
        \node[PosetVertex](4)at(2,0){\begin{math}4\end{math}};
        \node[PosetVertex](5)at(1,-1){\begin{math}5\end{math}};
        \node[PosetVertex](6)at(3,-1){\begin{math}6\end{math}};
        \draw[Edge](2)--(3);
        \draw[Edge](2)--(4);
        \draw[Edge](3)--(5);
        \draw[Edge](4)--(5);
        \draw[Edge](4)--(6);
    \end{tikzpicture}
\end{equation}
denotes the simple poset $([6], \Ord)$ satisfying among others
$3 \Ord 5$ and~$2 \Ord 6$.
\medbreak

The \Def{dual} of $\Qca$ is the poset $(\Qca, \bar{\Ord_\Qca})$ such
that $x \bar{\Ord_\Qca} y$ holds whenever $y \Ord_\Qca x$ for any
$x, y \in \Qca$. Besides, if $(\Qca_1, \Ord_{\Qca_1})$ and
$(\Qca_2, \Ord_{\Qca_2})$ are two posets, a map
$\phi : \Qca_1 \to \Qca_2$ is a \Def{poset morphism} if $\phi$ is a
collection morphism and for all $x, y \in \Qca_1$ such that
$x \Ord_{\Qca_1} y$, $\phi(x) \Ord_{\Qca_2} \phi(y)$. Besides, $\Qca_2$
is a \Def{subposet} of $\Qca_1$ if $\Qca_2$ is a subcollection of
$\Qca_1$ and $\Ord_{\Qca_2}$ is the restriction of $\Ord_{\Qca_1}$
on~$\Qca_2$.
\medbreak

Let us state the following easy lemma, used for instance in
Chapter~\ref{chap:posets}.
\medbreak

\begin{Lemma} \label{lem:morphism_posets_min}
    Let $\Qca_1$ and $\Qca_2$ be two posets and
    $\phi : \Qca_1 \to \Qca_2$ be a morphism of posets. Then, for all
    comparable objects $x$ and $y$ of $\Qca_1$,
    \begin{equation}
        \phi\left(x \Min_{\Qca_1} y\right) =
        \phi(x) \Min_{\Qca_2} \phi(y).
    \end{equation}
\end{Lemma}
\medbreak

\subsubsection{Patterns} \label{subsubsec:poset_patterns}
Let $(\Qca_1, \Ord_{\Qca_1})$ and $(\Qca_2, \Ord_{\Qca_1})$ be two
posets. We say that $\Qca_1$ admits an \Def{occurrence} of (the
\Def{pattern}) $\Qca_2$ if there is an isomorphism of posets
$\phi : \Qca_1' \to \Qca_2$ where $\Qca_1'$ is a subposet of $\Qca_1$.
Conversely, we say that $\Qca_1$ \Def{avoids} $\Qca_2$ if there is no
occurrence of $\Qca_2$ in~$\Qca_1$. Since only the isomorphism class of
a pattern is important to decide if a poset admits an occurrence of it,
we shall draw unlabeled Hasse diagrams to specify patterns. For
instance, the simple poset
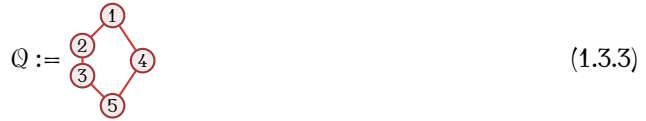
\begin{equation}
    \Qca :=
    \begin{tikzpicture}[xscale=.4,yscale=.4,Centering]
        \node[PosetVertex](1)at(0,0){\begin{math}1\end{math}};
        \node[PosetVertex](2)at(-1,-1){\begin{math}2\end{math}};
        \node[PosetVertex](3)at(-1,-2){\begin{math}3\end{math}};
        \node[PosetVertex](4)at(1,-1.5){\begin{math}4\end{math}};
        \node[PosetVertex](5)at(0,-3){\begin{math}5\end{math}};
        \draw[Edge](1)--(2);
        \draw[Edge](2)--(3);
        \draw[Edge](3)--(5);
        \draw[Edge](1)--(4);
        \draw[Edge](4)--(5);
    \end{tikzpicture}
\end{equation}
admits two occurrences of the simple pattern
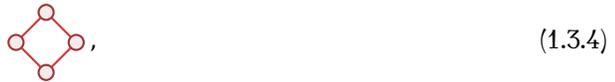
\begin{equation}
    \begin{tikzpicture}[xscale=.4,yscale=.4,Centering]
        \node[PosetVertex](1)at(1,1){};
        \node[PosetVertex](2)at(0,0){};
        \node[PosetVertex](3)at(2,0){};
        \node[PosetVertex](4)at(1,-1){};
        \draw[Edge](1)--(2);
        \draw[Edge](1)--(3);
        \draw[Edge](2)--(4);
        \draw[Edge](3)--(4);
    \end{tikzpicture}\,,
\end{equation}
a first one since $1 \Ord_\Qca 2$, $1 \Ord_\Qca 4$, $2 \Ord_\Qca 5$, and
$4 \Ord_\Qca 5$, and a second one since $1 \Ord_\Qca 3$,
$1 \Ord_\Qca 4$, $3 \Ord_\Qca 5$, and $4 \Ord_\Qca 5$. Moreover, $\Qca$
avoids the simple pattern
\begin{equation}
    \begin{tikzpicture}[xscale=.4,yscale=.4,Centering]
        \node[PosetVertex](1)at(0,0){};
        \node[PosetVertex](2)at(1,0){};
        \node[PosetVertex](3)at(2,0){};
    \end{tikzpicture}
\end{equation}
since $\Qca$ has no antichain of cardinality~$3$.
\medbreak

We call \Def{forest poset} any finite simple poset avoiding the pattern
$\AntiForestPattern$. In other words, a forest poset is a poset for
which its Hasse diagram is a forest of rooted trees (where roots are the
minimal elements).
\medbreak

\begin{Lemma} \label{lem:forest_poset_max_three_elements}
    Let $\Qca$ be a forest poset and $x$, $y$, and $z$ be three elements
    of $\Qca$ such that $x$ and $y$ are comparable and $y$ and $z$ are
    comparable. Then, $x \Min_\Qca y \Min_\Qca z$ is a well-defined
    element of~$\Qca$.
\end{Lemma}
\medbreak

\subsubsection{Examples} \label{subsubsec:examples_posets}
We consider here three well-known combinatorial posets.
\medbreak

\begin{figure}[ht]
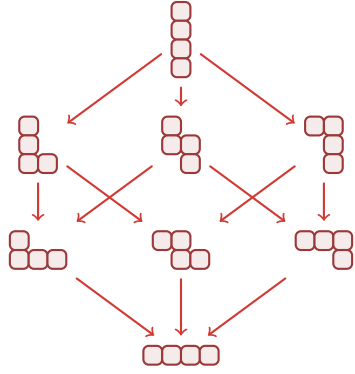

    \centering
    \begin{equation*}

            };
            \draw[Arc](1111)--(112);
            \draw[Arc](1111)--(121);
            \draw[Arc](1111)--(211);
            \draw[Arc](112)--(13);
            \draw[Arc](121)--(13);
            \draw[Arc](112)--(22);
            \draw[Arc](211)--(22);
            \draw[Arc](211)--(31);
            \draw[Arc](121)--(31);
            \draw[Arc](13)--(4);
            \draw[Arc](22)--(4);
            \draw[Arc](31)--(4);
        \end{tikzpicture}
    \end{equation*}
    \caption[The cube poset of order $4$.]
    {The Hasse diagram of the refinement order of compositions of
    size $4$, where each composition is represented through its ribbon
    diagram.}
    \label{fig:cube_order_4}
\end{figure}

\paragraph{The cube poset}
Let $\RefinementOrder$ be the partial order relation on the
combinatorial collection $\Compositions$ of compositions generated by
the covering relation $\Rca$ defined, for any composition $\LambdaB$ of
length $k$, by
\begin{equation}
    \left(\LambdaB_1, \dots, \LambdaB_{i - 1},
    \LambdaB_i, \LambdaB_{i + 1},
    \LambdaB_{i + 2}, \dots, \LambdaB_k\right)
    \, \Rca \,
    \left(\LambdaB_1, \dots, \LambdaB_{i - 1},
    \LambdaB_i + \LambdaB_{i + 1},
    \LambdaB_{i + 2}, \dots, \LambdaB_k\right).
\end{equation}
For instance, $2123 \RefinementOrder 215$ and $2123 \RefinementOrder 8$.
This order is the \Def{refinement order} of compositions. The Hasse
diagram of $(\Compositions, \RefinementOrder)$ restricted on
$\Compositions(4)$ is shown in Figure~\ref{fig:cube_order_4}.
\medbreak

Observe that for all compositions $\LambdaB$ and $\MuB$,
$\LambdaB \RefinementOrder \MuB$ if and only if
$\Des(\MuB) \subseteq \Des(\LambdaB)$. Each $n$-subposet of the
refinement order of compositions is known as the \Def{cube poset} of
dimension~$n - 1$. Moreover, the cube poset of dimension $n - 1$ is
isomorphic to the dual of the poset of all subsets of $[n - 1]$ ordered
by set inclusion. An isomorphism is provided by the map $\Des$ sending a
composition of size $n$ to a subset of~$[n - 1]$.
\medbreak

\begin{figure}[ht]
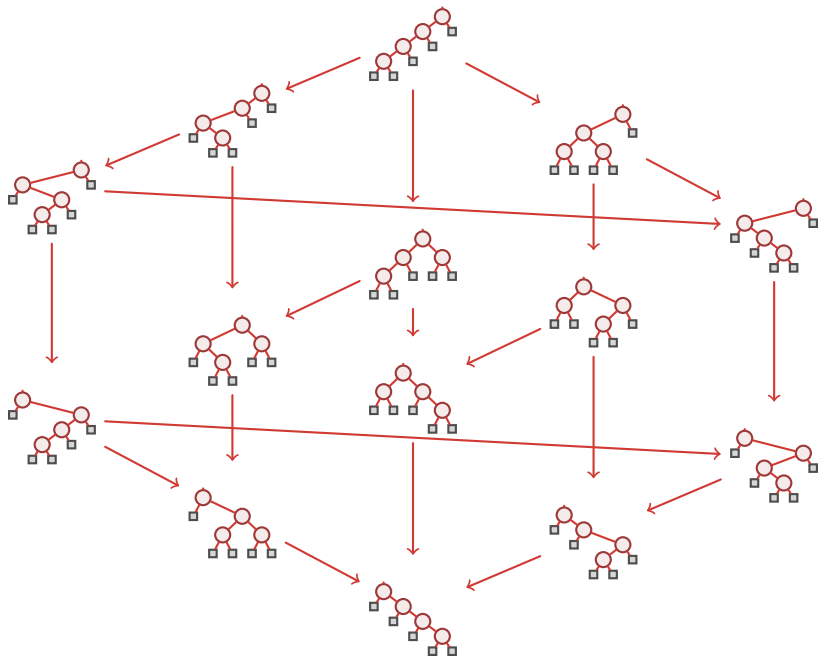

    \centering
    \begin{equation*}

            };
            \draw[Arc](222200000)--(222020000);
            \draw[Arc](222200000)--(222002000);
            \draw[Arc](222020000)--(220220000);
            \draw[Arc](222002000)--(220202000);
            \draw[Arc](220220000)--(220202000);
            \draw[Arc](222200000)--(222000200);
            \draw[Arc](222020000)--(220200200);
            \draw[Arc](222000200)--(220200200);
            \draw[Arc](222002000)--(220022000);
            \draw[Arc](222000200)--(220020200);
            \draw[Arc](220022000)--(220020200);
            \draw[Arc](220220000)--(202220000);
            \draw[Arc](220202000)--(202202000);
            \draw[Arc](202220000)--(202202000);
            \draw[Arc](202220000)--(202200200);
            \draw[Arc](220200200)--(202200200);
            \draw[Arc](220022000)--(202022000);
            \draw[Arc](202202000)--(202022000);
            \draw[Arc](202200200)--(202020200);
            \draw[Arc](220020200)--(202020200);
            \draw[Arc](202022000)--(202020200);
        \end{tikzpicture}
    \end{equation*}
    \caption[The Tamari poset of order $5$.]
    {The Hasse diagram of the Tamari poset of binary trees of
    size~$5$.}
    \label{fig:tamari_order_4}
\end{figure}

\paragraph{The Tamari order on binary trees}
Let $\Ord$ be the partial order relation on the combinatorial collection
$\BinaryTrees_\Leaf$ of binary trees generated by the covering relation
$\Rca$ defined by
\begin{equation} \label{equ:right_rotation_Tamari}
    (\dots(\Node, ((\Node, (\Rfr_1, \Rfr_2)), \Rfr_3))\dots)
    \, \Rca \,
    (\dots(\Node, (\Rfr_1, (\Node, (\Rfr_2, \Rfr_3))))\dots),
\end{equation}
where $\Rfr_1$, $\Rfr_2$, and $\Rfr_3$ are any binary trees. We call
$\Rca$ the \Def{right rotation} relation. At this moment, the definition
of this relation on binary trees is informal, but, in
Section~\ref{subsec:rewrite_rules_syntax_trees}, we shall develop
precise tools to define and handle such operations on binary trees and
more generally on syntax trees. The order $\Ord$ is the
\Def{Tamari order} on binary trees. The Hasse diagram of
$(\BinaryTrees_\Leaf, \Ord)$ restricted on $\BinaryTrees_\Leaf(5)$ is
shown in Figure~\ref{fig:tamari_order_4}.
\medbreak

The Tamari poset is a combinatorial poset on binary trees introduced in
the study of nonassociative operations~\cite{Tam62}. Indeed, the
covering relation generating this poset can be thought as a way to move
brackets in expressions where a nonassociative product intervenes.
Moreover, seen on binary trees, this operation translates as a right
rotation, a fundamental operation on binary search trees, used in an
algorithmic context~\cite{Knu98}. This operation is used to maintain
binary trees with a small height in order to access efficiently, from
the roots, to their internal nodes. Some of these trees are known as
balanced binary trees~\cite{AVL62} and form efficient structures to
represent dynamic sets (sets supporting the addition and the suppression
of elements). A lot of properties of the Tamari poset are known,
like the number of intervals of each of its $n$-subposets~\cite{Cha06}
(equivalently, this is the number of pairs of comparable trees
enumerated by their size), and the fact that these posets are
lattices~\cite{HT72}, for all $n \in \N_{\geq 1}$. Generalizations of
this poset have been introduced by Bergeron and
Préville-Ratelle~\cite{BP12} under the name of $m$-Tamari poset. This
poset is defined on the combinatorial collection of all
$m\! +\! 1$-ary trees (see Section~\ref{subsubsec:k_ary_trees}). The
number of intervals of each of its $n$-subposets, and the fact that
these posets are lattices are known from~\cite{BFP11}, for all
$n \in \N_{\geq 1}$.
\medbreak

\begin{figure}[ht]
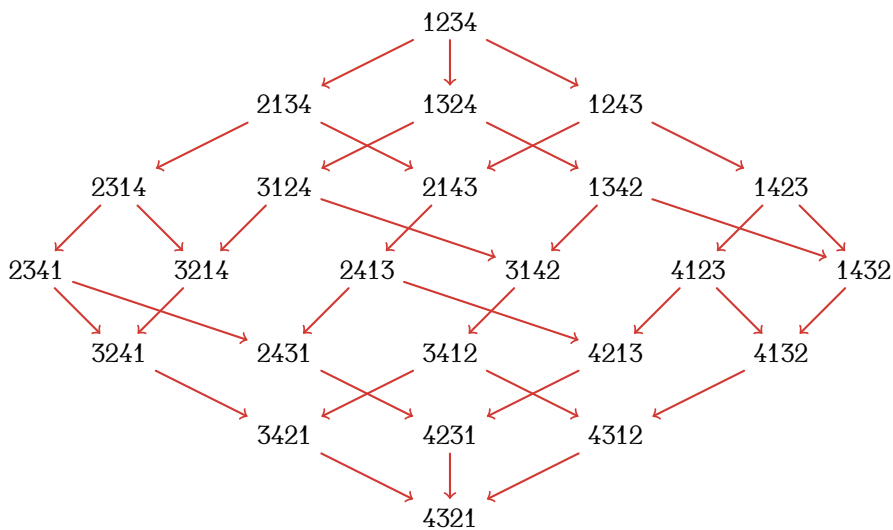

    \centering
    \begin{equation*}

    \end{equation*}
    \caption[The right weak poset of order $4$.]
    {The Hasse diagram of the right weak poset of permutations of
    size~$4$.}
    \label{fig:right_weak_order_order_4}
\end{figure}

\paragraph{The right weak order on permutations}
Let $\Ord$ be the partial order relation on the combinatorial collection
$\SymmetricGroup$ of permutations generated by the covering relation
$\Rca$ defined by
\begin{equation}
    \textcolor{Col1}{u}
    \textcolor{Col4}{\Asf \Bsf}
    \textcolor{Col1}{v}
    \, \Rca \,
    \textcolor{Col1}{u}
    \textcolor{Col4}{\Bsf \Asf}
    \textcolor{Col1}{v},
\end{equation}
where $u$ and $v$ are words on $\N_{\geq 1}$, and $\Asf$ and $\Bsf$ are
letters such that $\Asf < \Bsf$. This order is the \Def{right weak
order} of permutations. The Hasse diagram of
$(\SymmetricGroup, \Ord)$ restricted on $\SymmetricGroup(4)$ is shown in
Figure~\ref{fig:right_weak_order_order_4}.
\medbreak

The right weak poset of permutations is also a lattice~\cite{GR63,YO69}.
In a surprising way, despite its apparent simplicity, there is no known
description of the number of intervals of each $n$-subposet,
$n \in \N$, of the right weak poset. Some other combinatorial poset
structures exist on $\SymmetricGroup$ like the Bruhat order, whose
generating relation is similar to the one of the right weak poset. The
definition of the Bruhat order on permutations comes from the general
notion of Bruhat order~\cite{Bjo84} in Coxeter groups~\cite{Cox34}.
\medbreak

As a last noteworthy fact, the cube, the Tamari, and the right weak
posets are linked through surjective morphisms of combinatorial
posets~\cite{LR02}. Indeed, a map between the right weak poset to
the Tamari poset is based upon the binary search tree insertion
algorithm~\cite{Knu98,HNT05}. This algorithm consists in inserting the
letters of a permutation to form step by step a binary tree. Moreover, a
map between the Tamari poset to the cube poset uses the
canopies~\cite{LR98} of the binary trees. The canopy of a binary tree is
a binary word encoding the orientations (to the left or to the right) of
its leaves.
\medbreak

\subsection{Rewrite systems on collections} \label{subsec:rewrite_rules}
A rewrite rule describes a process whose goal is to transform
iteratively a combinatorial object into another one. We consider rewrite
rules on $I$-collections, so that an $i$-object, $i \in I$, can be
transformed only into $i$-objects. As we shall see, rewrite rules and
posets have some close connections because it is possible, in some cases,
to construct posets from rewrite systems. A general reference about
rewrite rules and rewrite systems is~\cite{BN98}.
\medbreak

Two properties of rewrite systems are fundamental: the termination and
the confluence. We provide strategies to prove that a given rewrite
system satisfies one or the other.
\medbreak

\subsubsection{Elementary definitions}
Let $C$ be an $I$-collection. An \Def{$I$-rewrite system} is a pair
$(C, \Rew)$ where $C$ is an $I$-collection and $\Rew$ is a relation on
$C$. We call $\Rew$ a \Def{rewrite rule}. For any property $P$ of
collections, we say that $(C, \Rew)$ \Def{satisfies the property $P$}
if, as a collection, $C$ satisfies~$P$. If $x$, $y_1$, \dots, $y_k$, and
$x'$ are objects of $C$ such that $k \in \N$ and
\begin{equation}
    x \Rew y_1 \Rew \cdots \Rew y_k \Rew x',
\end{equation}
we say that $x$ is \Def{rewritable} by $\Rew$ into $x'$ in $k + 1$
steps. The reflexive and transitive closure of $\Rew$ is denoted
by~$\RewRT$. The directed graph $(C, \Rew)$ consisting in $C$ as set
of vertices and $\Rew$ as set of arcs is the \Def{rewriting graph}
of~$(C, \Rew)$.
\medbreak

\subsubsection{Termination} \label{subsubsec:termination}
When there is no infinite chain
\begin{equation} \label{equ:infinite_rewrite_chain}
    x_1 \Rew x_2 \Rew x_3 \Rew \cdots
\end{equation}
where all $x_j \in C$, $j \in \N_{\geq 1}$, $(C, \Rew)$ is
\Def{terminating}. Observe that, if $C$ is combinatorial, due to the
fact that for any $i \in I$, each set $C(i)$ is finite and the fact that
the rewriting relation preserves the indices, if an infinite
chain~\eqref{equ:infinite_rewrite_chain} exists, then it is of the form
\begin{equation} \label{equ:infinite_rewrite_chain_repetition}
    x_1 \Rew \cdots \Rew x_j \Rew \cdots \Rew x_j \Rew \cdots,
\end{equation}
for a $j \in \N_{\geq 1}$. A \Def{normal form} of $(C, \Rew)$ is an
object $x$ of $C$ such that for all $x' \in C$, $x \RewRT x'$ imply
$x' = x$. In other words, a normal form of $(C, \Rew)$ is an object
which is not rewritable by $\Rew$. This set of objects, which is a
subcollection of $C$, is denoted by $\NormalForms_{(C, \Rew)}$. The
following result provides a tool in the aim to show that a rewrite
system is terminating.
\medbreak

\begin{Lemma} \label{lem:terminating_rewrite_rules_posets}
    Let $(C, \Rew)$ be a combinatorial rewrite system.
    Then, $(C, \Rew)$ is terminating if and only if the binary relation
    $\RewRT$ is an order relation and endows~$C$ with a structure of a
    combinatorial poset.
\end{Lemma}
\medbreak

In practice, Lemma~\ref{lem:terminating_rewrite_rules_posets} is used
as follows. To show that a combinatorial rewrite system $(C, \Rew)$ is
terminating, we construct a map $\theta : C \to \Qca$ where
$(\Qca, \Ord)$ is an $I$-poset such that for any $x, x' \in C$,
$x \Rew x'$ implies $\theta(x) \OrdStrict \theta(x')$. Such a map
$\theta$ is a \Def{termination invariant}. Indeed, since each $C(i)$,
$i \in I$, is finite, this property leads to the fact that there is no
infinite chain of the
form~\eqref{equ:infinite_rewrite_chain_repetition}. In most cases,
$\Qca$ is a set of tuples of integers of a fixed length, and $\Ord$ is
the lexicographic order on these tuples.
\medbreak

In~\cite{BN98}, a general method using maps called measure functions to
show that (not necessarily combinatorial) rewrite systems are
terminating is presented.
\medbreak

When $C$ is combinatorial and $\Rew$ is terminating, by
Lemma~\ref{lem:terminating_rewrite_rules_posets}, $(C, \RewRT)$ is a
combinatorial poset and we call it the \Def{poset generated} by~$\Rew$.
\medbreak

\subsubsection{Confluence}
When for any objects $x$, $y_1$, and $y_2$ of $C$ such that
$x \RewRT y_1$ and $x \RewRT y_2$, there exists an object $x'$ of $C$
such that $y_1 \RewRT x'$ and $y_2 \RewRT x'$, the rewrite system
$(C, \Rew)$ is \Def{confluent}. An object $x$ of $C$ is a
\Def{branching} object if there exist two different objects $y_1$ and
$y_2$ satisfying $x \Rew y_1$ and $x \Rew y_2$. In this case, the
pair $\{y_1, y_2\}$ is a \Def{branching pair} for~$x$. We say that a
branching pair $\{y_1, y_2\}$ is \Def{joinable} is there exists an
object $z$ of $C$ such that $y_1 \RewRT z$ and $y_2 \RewRT z$. In
practice, showing that a terminating rewrite system is confluent is
made simple thank to the following result, known as the
\Def{diamond lemma}.
\medbreak

\begin{Lemma} \label{lem:diamond_lemma}
    Let $(C, \Rew)$ be a rewrite system. If $(C, \Rew)$ is terminating
    and all of its branching pairs are joinable, $(C, \Rew)$ is
    confluent.
\end{Lemma}
\medbreak

Lemma~\ref{lem:diamond_lemma} is a highly important result in the theory
of rewrite systems and is due to Newman~\cite{New42}. There are some
additional useful tools in this theory like the Knuth-Bendix completion
algorithm~\cite{KB70}. This semi-algorithm takes as input a
non-confluent rewrite system and outputs, if possible, a confluent one
having the same reflexive, symmetric, and transitive closures.
\medbreak

When $\Rew$ is both terminating and confluent, $\Rew$ is
\Def{convergent}.
\medbreak

\subsubsection{Closures}
Let $(C, \Rew)$ be an $I$-rewrite system and assume that $C$ is endowed
with a set $\Pca$ of $\DPlus$-compatible products, where $\DPlus$ is an
associative binary product on $I$. Then, let $(C, \Rew_\Pca)$ be the
rewrite system such that $\Rew_\Pca$ contains $\Rew$ (as a binary
relation) and satisfies
\begin{equation} \label{equ:closure_rewrite_rule}
    \Product\left(x_1, \dots, x_{j - 1}, y,
    x_{j + 1}, \dots, x_p\right)
    \Rew_\Pca
    \Product\left(x_1, \dots, x_{j - 1}, y',
    x_{j + 1}, \dots, x_p\right)
\end{equation}
for any product $\Product$ of arity $p$ of $\Pca$, $j \in [p]$,
$x_\ell \in C$, $\ell \in [p] \setminus \{j\}$, $y, y' \in C$ such that
$y \Rew y'$, and when both members
of~\eqref{equ:closure_rewrite_rule} are defined (because the products of
$\Pca$ can be partial, see
Section~\ref{subsubsec:products_collections}). The fact that all
products $\Product$ of $\Pca$ are $\DPlus$-compatible ensures that
$(C, \Rew_\Pca)$ is a rewrite system. We call~$(C, \Rew_\Pca)$ the
\Def{$\Pca$-closure} of $(C, \Rew)$. Such closures provide convenient
and concise ways to define rewrite systems.
\medbreak

\subsubsection{Examples} \label{subsubsec:examples_rewrite_rules}
Let us review some examples of rewrite systems on various combinatorial
sets.
\medbreak

\paragraph{A rewrite rule on words}
Let $A := \{\Asf, \Bsf\}$ be an alphabet, and consider the rewrite
system $(A^*, \Rew$) defined by
\begin{equation}
    \textcolor{Col1}{u(1) \dots u(n - 1)} \textcolor{Col4}{u(n)}
    \Rew
    \textcolor{Col4}{u(n)} \textcolor{Col1}{u(1) \dots u(n - 1)}
\end{equation}
for any $u \in A^n$ and $n \in \N_{\geq 2}$. We have for instance
\begin{equation}
    \Asf \Asf \Bsf \Asf \Rew \Asf \Asf \Asf \Bsf
    \Rew \Bsf \Asf \Asf \Asf \Rew \Asf \Bsf \Asf \Asf
    \Rew \Asf \Asf \Bsf \Asf.
\end{equation}
This rewrite system is not terminating but, since for each word
$u \in A^*$ there is at most a word $v \in A^*$ satisfying $u \Rew v$,
$(A^*, \Rew)$ is confluent.
\medbreak

Let also be the rewrite system $(A^*, \Rew)$ defined by
$\Asf \Bsf \Asf \Rew \Bsf \Asf \Bsf$. Consider the ternary product
$\Product$ on $A^*$ defined by $\Product(u, v, w) := u \Conc v \Conc w$
where $\Conc$ is the concatenation product of words.
Let $(A^*, \Rew_\Pca)$ be the $\Pca$-closure of $(A^*, \Rew)$ where
$\Pca := \{\Product\}$. By definition of closures, $\Rew_\Pca$ satisfies
\begin{subequations}
\begin{equation}
    \Asf \Bsf \Asf \Rew_\Pca \Bsf \Asf \Bsf,
\end{equation}
\begin{equation}
    \Asf \Bsf \Asf \Conc v \Conc w
    \Rew_\Pca
    \Bsf \Asf \Bsf \Conc v \Conc w,
\end{equation}
\begin{equation}
    u \Conc \Asf \Bsf \Asf \Conc w
    \Rew_\Pca
    u \Conc \Bsf \Asf \Bsf \Conc w,
\end{equation}
\begin{equation}
    u \Conc v \Conc \Asf \Bsf \Asf
    \Rew_\Pca
    u \Conc v \Conc \Bsf \Asf \Bsf,
\end{equation}
\end{subequations}
for any words $u$, $v$, and $w$ on $A$. All this is equivalent to the
fact that $\Rew_\Pca$ is the rewrite rule satisfying
\begin{equation}
    u \Conc \Asf \Bsf \Asf \Conc w
    \Rew_\Pca
    u \Conc \Bsf \Asf \Bsf \Conc w,
\end{equation}
for any words $u$ and $w$ on $A$. The rewrite system $(A^*, \Rew_\Pca)$
is terminating since, for any words $u$ and $v$ on $A$, if
$u \Rew_\Pca v$, then $|v|_\Bsf  = |u|_\Bsf + 1$. Hence, the map
$\theta : A^n \to [0, n]$ defined for any $n \in \N$ and $u \in A^n$ by
$\theta(u) := |u|_\Bsf$ is a termination invariant. The normal forms of
$(A^*, \Rew_\Pca)$ are the words that do not admit $\Asf \Bsf \Asf$ as
factor. Moreover, $(A^*, \Rew_\Pca)$ is not confluent since
\begin{math}
    \Asf \Bsf \Asf \Bsf \Asf \Rew_\Pca \Bsf \Asf \Bsf \Bsf \Asf
\end{math}
and
\begin{math}
    \Asf \Bsf \Asf \Bsf \Asf \Rew_\Pca \Asf \Bsf \Bsf \Asf \Bsf
\end{math},
and $\{\Bsf \Asf \Bsf \Bsf \Asf, \Asf \Bsf \Bsf \Asf \Bsf\}$ is a
non-joinable branching pair for $\Asf \Bsf \Asf \Bsf \Asf$ (since these
two elements are normal forms).
\medbreak

\section{Collections of trees} \label{sec:trees}
This section is devoted mainly to set all basic definitions about
trees used in this work. We define here the collection of planar
rooted trees and present some of its properties. We then consider
enrichments of planar rooted trees, namely the syntax trees. These are
one of the most important objects in this work since bases of free
operads are indexed by syntax trees.  Moreover, rewrite systems on
syntax trees are reviewed. These rewrite systems are a major tool to
study operads  since they allow to establish presentation by generators
and relations, or the Koszulity of an operad.
\medbreak

\subsection{Planar rooted trees} \label{subsec:planar_rooted_trees}
The graded combinatorial collection of the planar rooted trees can be
defined concisely in a recursive way by using some operations over
graded combinatorial collections (see
Section~\ref{subsubsec:operations_graded_comb_collections}). However, to
define rigorously the usual notions of internal node, leaf, child,
father, path, subtree, {\em etc.}, we need the notion of language
associated with a tree. Indeed, a planar rooted tree is in fact a finite
language satisfying some properties. Therefore, in this section, we
shall adopt the point of view of defining most of the properties of a
planar rooted tree through its language.
\medbreak

\subsubsection{Collection of planar rooted trees}
\label{subsubsec:comb_collection_planar_rooted_trees}
Let $\PlanarRootedTrees$ be the graded combinatorial collection
satisfying the relation
\begin{equation} \label{equ:combinatorial_set_planar_rooted_trees}
    \PlanarRootedTrees = \{\Node\} \times \List(\PlanarRootedTrees)
\end{equation}
where $\Node$ is an atomic object called \Def{node}. We call
\Def{planar rooted tree} each object of $\PlanarRootedTrees$. By
definition, a planar rooted tree $\Tfr$ is an ordered pair
$(\Node, (\Tfr_1, \dots, \Tfr_k))$ where $(\Tfr_1, \dots, \Tfr_k)$ is a
(possibly empty) tuple of planar rooted trees. This definition is
recursive. By convention, the planar rooted tree $(\Node, ())$ is
denoted by $\Leaf$ and is called the \Def{leaf}. Observe that the leaf
is of size $1$. For instance,
\begin{equation} \label{equ:examples_planar_rooted_trees}
    \Leaf, \quad
    (\Node, (\Leaf)), \quad
    (\Node, (\Leaf, \Leaf)), \quad
    (\Node, (\Leaf, (\Node, (\Leaf)))), \quad
    (\Node, ((\Node, ((\Node, (\Leaf, \Leaf)))),
        \Leaf, (\Node, (\Leaf, \Leaf))))
\end{equation}
are planar rooted trees. The \Def{root arity} of a planar rooted tree
$\Tfr := (\Node, (\Tfr_1, \dots, \Tfr_k))$ is $k$. If $\Tfr$ is a planar
rooted tree different from the leaf, by definition, $\Tfr$ can be
expressed as $\Tfr = (\Node, (\Tfr_1, \dots, \Tfr_k))$ where $k \geq 1$
and all $\Tfr_i$, $i \in [k]$, are planar rooted trees. In this case,
for any $i \in [k]$, $\Tfr_i$ is the \Def{$i$th suffix subtree} of
$\Tfr$. Planar rooted trees are depicted by drawing each leaf by
$\LeafPic$ and each planar rooted tree different from the leaf by a node
$\NodePic$ attached below it, from left to right, to its suffix subtrees
$\Tfr_1$, \dots, $\Tfr_k$ by means of edges $\EdgePic$. For instance,
the planar rooted trees of~\eqref{equ:examples_planar_rooted_trees}
are depicted by
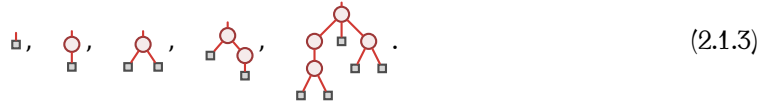
\begin{equation}
    \LeafPic, \quad
    \begin{tikzpicture}[scale=.3,Centering]
        \node[Leaf](1)at(0.00,-1.00){};
        \node[Node](0)at(0.00,0.00){};
        \draw[Edge](1)--(0);
        \node(r)at(0.00,1){};
        \draw[Edge](r)--(0);
    \end{tikzpicture}\,, \quad
    \begin{tikzpicture}[xscale=.2,yscale=.2,Centering]
        \node[Leaf](0)at(0.00,-1.50){};
        \node[Leaf](2)at(2.00,-1.50){};
        \node[Node](1)at(1.00,0.00){};
        \draw[Edge](0)--(1);
        \draw[Edge](2)--(1);
        \node(r)at(1.00,1.5){};
        \draw[Edge](r)--(1);
    \end{tikzpicture}\,, \quad
    \begin{tikzpicture}[xscale=.23,yscale=.2,Centering]
        \node[Leaf](0)at(0.00,-1.33){};
        \node[Leaf](3)at(2.00,-2.67){};
        \node[Node](1)at(1.00,0.00){};
        \node[Node](2)at(2.00,-1.33){};
        \draw[Edge](0)--(1);
        \draw[Edge](2)--(1);
        \draw[Edge](3)--(2);
        \node(r)at(1.00,1.5){};
        \draw[Edge](r)--(1);
    \end{tikzpicture}\,, \quad
    \begin{tikzpicture}[xscale=.18,yscale=.16,Centering]
        \node[Leaf](1)at(0.00,-6.75){};
        \node[Leaf](3)at(2.00,-6.75){};
        \node[Leaf](5)at(3.00,-2.25){};
        \node[Leaf](6)at(4.00,-4.50){};
        \node[Leaf](8)at(6.00,-4.50){};
        \node[Node](0)at(1.00,-2.25){};
        \node[Node](2)at(1.00,-4.50){};
        \node[Node](4)at(3.00,0.00){};
        \node[Node](7)at(5.00,-2.25){};
        \draw[Edge](0)--(4);
        \draw[Edge](1)--(2);
        \draw[Edge](2)--(0);
        \draw[Edge](3)--(2);
        \draw[Edge](5)--(4);
        \draw[Edge](6)--(7);
        \draw[Edge](7)--(4);
        \draw[Edge](8)--(7);
        \node(r)at(3.00,1.75){};
        \draw[Edge](r)--(4);
    \end{tikzpicture}\,.
\end{equation}
\medbreak

By definition of the product and the list operations over graded
collections (see
Section~\ref{subsubsec:operations_graded_comb_collections}), the size of
a planar rooted tree $\Tfr$ having a root arity of $k$ satisfies
\begin{equation} \label{equ:size_planar_rooted_trees}
    |\Tfr| = 1 + \sum\limits_{i \in [k]} |\Tfr_k|.
\end{equation}
In other words, the size of $\Tfr$ is the number of occurrences of
$\Node$ it contains.  We also deduce
from~\eqref{equ:combinatorial_set_planar_rooted_trees} that the
generating series of $\PlanarRootedTrees$ satisfies
\begin{equation}
    \GenSeries_\PlanarRootedTrees(t) =
    \frac{t}{1 - \GenSeries_\PlanarRootedTrees(t)}
\end{equation}
so that it satisfies the quadratic algebraic equation
\begin{equation} \label{equ:generating_series_planar_rooted_trees}
    t - \GenSeries_\PlanarRootedTrees(t)
    + \GenSeries_\PlanarRootedTrees(t)^2 = 0.
\end{equation}
\medbreak

\subsubsection{Induction and structural induction}
\label{subsubsec:structural_induction}
One among the most obvious techniques to prove that all the planar
rooted trees of a subcollection $C$ of $\PlanarRootedTrees$ satisfy a
predicate $P$ consists in performing a proof by induction on the size of
the trees of~$C$.
\medbreak

There is another method which is in some cases much more elegant than
this approach, called \Def{structural induction} on trees. A
subcollection $C$ of $\PlanarRootedTrees$ is \Def{inductive} if
$C$ is nonempty and, if $\Tfr \in C$, all suffix subtrees $\Tfr_i$ of
$\Tfr$ belong to~$C$. Observe in particular that $\Leaf$ belongs to any
inductive subcollection of~$\PlanarRootedTrees$.
\medbreak

\begin{Theorem} \label{thm:structural_induction}
    Let $C$ be an inductive subcollection of $\PlanarRootedTrees$
    and $P$ be a predicate on~$C$. If
    \begin{enumerate}[label={(\it\roman*)}]
        \item \label{item:structural_induction_1}
        the leaf $\Leaf$ satisfies $P$;
        \item \label{item:structural_induction_2}
        for any $\Tfr_1, \dots, \Tfr_k \in C$ such that
        $\Tfr := (\Node, (\Tfr_1, \dots, \Tfr_k))$ belongs to $C$, the
        fact that all $P(\Tfr_i)$, $i \in [k]$, hold implies that
        $P(\Tfr)$ holds;
    \end{enumerate}
    then, all objects of $C$ satisfy~$P$.
\end{Theorem}
\medbreak

Theorem~\ref{thm:structural_induction} provides a powerful tool to prove
properties $P$ of planar rooted trees belonging to inductive
combinatorial subsets $C$. In practice, to perform a structural
induction in order to show that all objects $\Tfr$ of $C$ satisfy $P$,
we check that $C$ is inductive and that
Properties~\ref{item:structural_induction_1}
and~\ref{item:structural_induction_2} of
Theorem~\ref{thm:structural_induction} hold.
\medbreak

\subsubsection{Links with binary trees}
\label{subsubsec:links_binary_trees}
As a consequence of~\eqref{equ:generating_series_planar_rooted_trees},
we observe that the generating series of $\PlanarRootedTrees$ satisfies
the same algebraic relation as the one of the graded collection
$\BinaryTrees_\Leaf$ of binary trees where the size of a binary tree is
its number of leaves (defined in Section~\ref{subsubsec:binary_trees}).
Therefore, $\PlanarRootedTrees$ and $\BinaryTrees_\Leaf$ are isomorphic
as graded collections. Let us describe an explicit isomorphism between
these two collections. Let
$\phi : \PlanarRootedTrees \to \BinaryTrees_\Leaf$ be the map
recursively defined, for any planar rooted tree $\Tfr$, by
\begin{equation} \label{equ:bijection_planar_rooted_trees_binary_trees}
    \phi(\Tfr) :=
    \begin{cases}
        \Leaf \in \BinaryTrees_\Leaf & \mbox{if } \Tfr = \Leaf, \\
        \left(\Node,
        \left(\phi(\Tfr_1),
            \phi
            \left(\left(\Node,
            \left(\Tfr_2, \dots, \Tfr_k\right)
            \right)\right)
        \right)
        \right)
            & \mbox{otherwise (}
            \Tfr = (\Node, (\Tfr_1, \Tfr_2, \dots, \Tfr_k))
            \mbox{ with } k \geq 1 \mbox{)} .
    \end{cases}
\end{equation}
One has for instance
\begin{subequations}
\begin{equation}
    \phi\left(
\,.
\end{equation}
\end{subequations}
\medbreak

\begin{Proposition}
\label{prop:bijection_planar_rooted_trees_binary_trees}
    The graded combinatorial collections $\PlanarRootedTrees$ and
    $\BinaryTrees_\Leaf$ are isomorphic. The map~$\phi$ defined
    by~\eqref{equ:bijection_planar_rooted_trees_binary_trees} is an
    isomorphism between these two collections.
\end{Proposition}
\medbreak

This bijection is known as the rotation correspondence and is due to
Knuth~\cite{Knu97}. It offers a means to encode a planar rooted tree by
a binary tree and admits applications in algebraic
combinatorics~\cite{NT13,EM14}.
\medbreak

\subsubsection{Tree languages}
To rigorously specify nodes in planar rooted trees, we shall use a
useful interpretation of planar rooted trees as special languages on the
alphabet $\N_{\geq 1}$. Recall that a right monoid action of a monoid
$A^*$ of words (endowed with the concatenation product) on a set $S$ is
a map $\Action : S \times A^* \to S$ satisfying $x \Action \epsilon = x$
and $x \Action u a = (x \Action u) \Action a$, for all $x \in S$,
$u \in A^*$, and $a \in A$. Let
\begin{equation}
    \Action : \PlanarRootedTrees \times \N_{\geq 1}^* \to
    \PlanarRootedTrees
\end{equation}
be the right partial monoid action defined recursively by
\begin{equation} \label{equ:action_monoid_trees}
    (\Node, (\Tfr_1, \dots, \Tfr_k)) \Action u :=
    \begin{cases}
        (\Node, (\Tfr_1, \dots, \Tfr_k)) & \mbox{if } u = \epsilon, \\
        \Tfr_i \Action v & \mbox{otherwise (}
            u = i v \mbox{ where } v \in \N_{\geq 1}^*
            \mbox{ and } i \in \N_{\geq 1} \mbox{)},
    \end{cases}
\end{equation}
for any $(\Node, (\Tfr_1, \dots, \Tfr_k)) \in \PlanarRootedTrees$ and
$u \in \N_{\geq 1}^*$. Observe that this action is partial since each
$\Tfr_i$ in~\eqref{equ:action_monoid_trees} is well-defined only if $i$
is no greater than the root arity of~$\Tfr$. The \Def{tree language}
$\TreeLanguage(\Tfr)$ of $\Tfr$ is the finite language on $\N_{\geq 1}$
of all the words $u$ such that $\Tfr \Action u$ is a well-defined planar
rooted tree.
\medbreak

For instance, by setting
\begin{equation} \label{equ:example_planar_rooted_tree}
    \Tfr :=
    \begin{tikzpicture}[scale=.16,Centering]
        \node[Leaf](0)at(0.00,-2.40){};
        \node[Leaf](10)at(7.00,-7.20){};
        \node[Leaf](11)at(8.00,-2.40){};
        \node[Leaf](3)at(1.00,-9.60){};
        \node[Leaf](5)at(3.00,-9.60){};
        \node[Leaf](7)at(4.00,-4.80){};
        \node[Leaf](8)at(5.00,-7.20){};
        \node[Node](1)at(4.00,0.00){};
        \node[Node](2)at(2.00,-4.80){};
        \node[Node](4)at(2.00,-7.20){};
        \node[Node](6)at(4.00,-2.40){};
        \node[Node](9)at(6.00,-4.80){};
        \draw[Edge](0)--(1);
        \draw[Edge](10)--(9);
        \draw[Edge](11)--(1);
        \draw[Edge](2)--(6);
        \draw[Edge](3)--(4);
        \draw[Edge](4)--(2);
        \draw[Edge](5)--(4);
        \draw[Edge](6)--(1);
        \draw[Edge](7)--(6);
        \draw[Edge](8)--(9);
        \draw[Edge](9)--(6);
        \node(r)at(4.00,1.80){};
        \draw[Edge](r)--(1);
    \end{tikzpicture}\,,
\end{equation}
we have
\begin{equation}
    \Tfr \Action 1 = \LeafPic, \quad
    \Tfr \Action 231 = \LeafPic, \quad
    \Tfr \Action 3 = \LeafPic, \quad
    \Tfr \Action 21 =
    \begin{tikzpicture}[scale=.2,Centering]
        \node[Leaf](1)at(0.00,-2.67){};
        \node[Leaf](3)at(2.00,-2.67){};
        \node[Node](0)at(1.00,0.00){};
        \node[Node](2)at(1.00,-1.33){};
        \draw[Edge](1)--(2);
        \draw[Edge](2)--(0);
        \draw[Edge](3)--(2);
        \node(r)at(1.00,1.5){};
        \draw[Edge](r)--(0);
    \end{tikzpicture}\,, \quad
    \Tfr \Action 23 =
    \begin{tikzpicture}[scale=.17,Centering]
        \node[Leaf](0)at(0.00,-1.50){};
        \node[Leaf](2)at(2.00,-1.50){};
        \node[Node](1)at(1.00,0.00){};
        \draw[Edge](0)--(1);
        \draw[Edge](2)--(1);
        \node(r)at(1.00,1.75){};
        \draw[Edge](r)--(1);
    \end{tikzpicture}\,,
\end{equation}
and, among others, the actions of the words $11$, $24$, and $2321$ on
$\Tfr$ are all undefined. Moreover, the tree language of $\Tfr$ is
\begin{equation} \label{equ:example_tree_language}
    \TreeLanguage(\Tfr) =
    \{\epsilon, 1, 2, 21, 211, 2111, 2112, 22, 23, 231, 232, 3\}.
\end{equation}
\medbreak

Let $\Lca_\PlanarRootedTrees$ be the graded combinatorial collection of
all finite and nonempty prefix languages $\Lca$ on $\N_{\geq 1}$ such
that if $ui \in \Lca$ where $u \in \N_{\geq 1}^*$ and
$i \in \N_{\geq 2}$, $u i' \in \Lca$ where $i' := i - 1$. The size of
such a language is its cardinality. For instance, the set
$\TreeLanguage(\Tfr)$ of~\eqref{equ:example_tree_language} is an object
of size $12$ of~$\Lca_\PlanarRootedTrees$, and
$\{\epsilon, 1, 11, 12, 2\}$ is an object of size~$5$.
\medbreak

\begin{Proposition}
\label{prop:bijection_planar_rooted_trees_prefix_languages}
    The graded combinatorial collections $\PlanarRootedTrees$ and
    $\Lca_\PlanarRootedTrees$ are isomorphic. Seen as a morphism of
    combinatorial collections
    \begin{math}
        \TreeLanguage :
        \PlanarRootedTrees \to \Lca_\PlanarRootedTrees
    \end{math},
    $\TreeLanguage$ is an isomorphism between these two collections.
\end{Proposition}
\medbreak

Proposition~\ref{prop:bijection_planar_rooted_trees_prefix_languages} is
used in practice to define planar rooted trees through their languages.
This will be useful later when operations on planar rooted trees will be
described.
\medbreak

\subsubsection{Additional definitions}
\label{subsubsec:definitions_trees}
Let $\Tfr$ be a planar rooted tree. We say that each word of
$\TreeLanguage(\Tfr)$ is a \Def{node} of $\Tfr$. A node $u$ of $\Tfr$
is an \Def{internal node} if there is an $i \in \N_{\geq 1}$ such that
$ui$ is a node of $\Tfr$. A node $u$ of $\Tfr$ which is not an internal
node is a \Def{leaf}. The set of all internal nodes (resp. leaves) of
$\Tfr$ is denoted by $\TreeLanguage_\Node(\Tfr)$ (resp.
$\TreeLanguage_\Leaf(\Tfr)$). The \Def{root} of $\Tfr$ is the node
$\epsilon$ (which can be either an internal node or a leaf). The
\Def{degree} $\deg(\Tfr)$ of $\Tfr$ is $\# \TreeLanguage_\Node(\Tfr)$
and the \Def{arity} $\Arity(\Tfr)$ of $\Tfr$ is
$\# \TreeLanguage_\Leaf(\Tfr)$. A node $u$ of $\Tfr$ is an
\Def{ancestor} of a node $v$ of $\Tfr$ if $u \ne v$ and
$u \PrefixOrder v$. Moreover, $v$ is the \Def{$i$th child} of $u$ if
$v = u i$ for an $i \in \N_{\geq 1}$. In this case, $u$ is the (unique)
\Def{father} of $v$. The \Def{arity} of a node is the number of children
it has. Two nodes $v$ and $v'$ of $\Tfr$ are \Def{brothers} if there
exist a node $u$ of $\Tfr$ and $i \ne i' \in \N_{\geq 1}$ such that $v$
is the $i$th child of $u$ and $v'$ is the $i'$th child of $u$.
The lexicographic order on the words of $\TreeLanguage(\Tfr)$ induces a
total order on the nodes of $\Tfr$ called \Def{depth-first order}. The
\Def{$i$th leaf} of $\Tfr$ is the $i$th leaf encountered by considering
the nodes of $\Tfr$ according to the depth-first order. A \Def{sector}
of $\Tfr$ is an ordered pair $(u_i, u_{i + 1})$ of leaves of $\Tfr$ such
that $u_i$ (resp. $u_{i + 1}$) is the $i$th (resp. $i\! +\! 1$st) leaf
of $\Tfr$. The \Def{number of sectors} of $\Tfr$ is denoted by
$\Sector(\Tfr)$ and is equal to $\Arity(\Tfr) - 1$. A \Def{path} in
$\Tfr$ is a sequence $(u_1, \dots, u_k)$ of nodes of $\Tfr$ such that
for any $j \in [k - 1]$, $u_j$ is the father of $u_{j + 1}$. Such a path
is \Def{maximal} if $u_1$ is the root of $\Tfr$ and $u_k$ is a leaf. The
\Def{length} of a path is the number of nodes it contains. The
\Def{height} $\Height(\Tfr)$ of $\Tfr$ is the maximal length of its
maximal paths minus $1$. This is also the length of a longest word of
$\TreeLanguage(\Tfr)$ minus~$1$. When all maximal paths of $\Tfr$ have
the same length, $\Tfr$ is \Def{perfect}. For any node $u$ of $\Tfr$,
the planar rooted tree $\Tfr \Action u$ is the \Def{suffix subtree} of
$\Tfr$ rooted at $u$. By extension, the \Def{$i$th suffix subtree} of
$u$ is the planar rooted tree $\Tfr \Action ui$ when $i$ is no greater
than the arity of $u$. A planar rooted tree $\Sfr$ is a
\Def{prefix subtree} of  $\Tfr$ if
$\TreeLanguage(\Sfr) \subseteq \TreeLanguage(\Tfr)$. A planar rooted
tree $\Sfr$ is a \Def{factor subtree} of $\Tfr$ rooted at a node $u$ if
$\Sfr$ is a prefix subtree of a suffix subtree of $\Tfr$ rooted at~$u$.
The \Def{poset} induced by $\Tfr$ is the poset $(\Qca_\Tfr, \Ord_\Tfr)$
where $\Qca_\Tfr := \TreeLanguage_\Node(\Tfr)$ and $\Ord_\Tfr$ is
the prefix order relation $\PrefixOrder$ on words. In other terms, the
poset induced by $\Tfr$ is a poset on the internal nodes of $\Tfr$
where $\Tfr$ is its Hasse diagram wherein the root is the least element.
\medbreak

Let us provide some examples for these notions. Consider the planar
rooted tree $\Tfr$ of~\eqref{equ:example_planar_rooted_tree}. Then,
\begin{subequations}
\begin{equation}
    \TreeLanguage_\Node(\Tfr) =
    \{\epsilon, 2, 21, 211, 23\},
\end{equation}
\begin{equation}
    \TreeLanguage_\Leaf(\Tfr) =
    \{1, 2111, 2112, 22, 231, 232, 3\},
\end{equation}
\end{subequations}
so that $\deg(\Tfr) = 5$ and $\Arity(\Tfr) = 7$. Besides, the sequences
$(\epsilon, 2, 21)$ and $(\epsilon, 2, 23)$ are nonmaximal paths in
$\Tfr$, and on the contrary, the paths $(\epsilon, 1)$,
$(\epsilon, 2, 21, 211, 2112)$, and $(\epsilon, 2, 22)$ are maximal. The
maximal path $(\epsilon, 2, 21, 211, 2112)$ have a maximal length among
all maximal paths of $\Tfr$ and hence, the height of $\Tfr$ is~$4$.
In the poset induced by $\Tfr$, one has
$\epsilon \Ord_\Tfr 2 \Ord_\Tfr 21 \Ord 211$ and
$\epsilon \Ord_\Tfr 2 \Ord_\Tfr 23$. Finally, the planar rooted tree
\begin{equation}
    \Sfr :=
    \begin{tikzpicture}[xscale=.25,yscale=.18,Centering]
        \node[Leaf](0)at(0.00,-2.00){};
        \node[Leaf](3)at(1.00,-6.00){};
        \node[Leaf](5)at(2.00,-4.00){};
        \node[Leaf](6)at(3.00,-4.00){};
        \node[Leaf](7)at(4.00,-2.00){};
        \node[Node](1)at(2.00,0.00){};
        \node[Node](2)at(1.00,-4.00){};
        \node[Node](4)at(2.00,-2.00){};
        \draw[Edge](0)--(1);
        \draw[Edge](2)--(4);
        \draw[Edge](3)--(2);
        \draw[Edge](4)--(1);
        \draw[Edge](5)--(4);
        \draw[Edge](6)--(4);
        \draw[Edge](7)--(1);
        \node(r)at(2.00,1.50){};
        \draw[Edge](r)--(1);
    \end{tikzpicture}
\end{equation}
is a prefix subtree of $\Tfr$, and, the planar rooted tree
\begin{equation}
    \Rfr :=
    \begin{tikzpicture}[xscale=.25,yscale=.18,Centering]
        \node[Leaf](1)at(0.00,-3.33){};
        \node[Leaf](3)at(1.00,-1.67){};
        \node[Leaf](4)at(2.00,-1.67){};
        \node[Node](0)at(0.00,-1.67){};
        \node[Node](2)at(1.00,0.00){};
        \draw[Edge](0)--(2);
        \draw[Edge](1)--(0);
        \draw[Edge](3)--(2);
        \draw[Edge](4)--(2);
        \node(r)at(1.00,1.5){};
        \draw[Edge](r)--(2);
    \end{tikzpicture}\,,
\end{equation}
being a suffix subtree of $\Sfr$ rooted at the node $2$, is a factor
subtree of $\Tfr$ rooted at~$2$.
\medbreak

\subsection{Subcollections of planar rooted trees}
\label{subsec:families_planar_rooted_trees}
By basically restraining the possible arities of the internal nodes of
planar rooted trees, we obtain several subcollections of
$\PlanarRootedTrees$. We review here the families formed by ladders,
corollas, $k$-ary trees, and Schröder trees. Besides, among these
families, some admit alternative size functions.
\medbreak

\subsubsection{Ladders and corollas}
A \Def{ladder} is a planar rooted tree of arity $1$. The first ladders
are
\begin{equation}
    \LeafPic, \quad
\,.
\end{equation}
This set of ladders forms a subcollection $\Ladders$ of
$\PlanarRootedTrees$. Besides, a \Def{corolla} is a planar rooted tree
of degree $1$. The first corollas are
\begin{equation}
    %
\,.
\end{equation}
This set of corollas forms a subcollection $\Corollas$ of
$\PlanarRootedTrees$. Observe that $(\Node, (\Leaf))$ is the only
planar rooted that is both a ladder and a corolla.
\medbreak

\subsubsection{$k$-ary trees} \label{subsubsec:k_ary_trees}
Let $k \in \N_{\geq 1}$. A \Def{$k$-ary} tree is a planar rooted tree
$\Tfr$ such that all internal nodes are of arity $k$. For instance, the
first $3$-ary trees are
\begin{equation}
    \LeafPic, \quad
    %
\,.
\end{equation}
This set of $k$-ary trees forms a subcollection $\Ary^{(k)}$ of
$\PlanarRootedTrees$ expressing recursively as
\begin{equation}
    \Ary^{(k)} = \{\Leaf\} + \{\Node\} \times {\Ary^{(k)}}^k,
\end{equation}
where $\Leaf$ and $\Node$ are both atomic. One can immediately observe
that $\Ary^{(1)} = \Ladders$.
\medbreak

By structural induction (see Theorem~\ref{thm:structural_induction}) on
$\Ary^{(k)}$ (which is an inductive subcollection of
$\PlanarRootedTrees$), it follows that for any $k$-ary tree $\Tfr$, the
arity and the degree of $\Tfr$ are related by
\begin{equation} \label{equ:k-ary_trees_arity_degree}
    \Arity(\Tfr) - \deg(\Tfr) (k - 1) = 1.
\end{equation}
This implies that a $k$-ary tree of a given arity has an imposed degree
and conversely, a $k$-ary tree of a given degree has an imposed arity.
Hence, since the size of a $k$-ary tree $\Tfr$ is
\begin{math}
    \Arity(\Tfr) + \deg(\Tfr)
\end{math}
and there are finitely many planar rooted trees of a fixed size, there
are finitely many $k$-ary trees of a fixed arity, and there are finitely
many $k$-ary trees of a fixed degree. As a consequence, the graded
collections $\Ary^{(k)}_\Leaf$ and $\Ary^{(k)}_\Node$ of all $k$-ary
trees such that the size of a tree of $\Ary^{(k)}_\Leaf$ is its arity
and the size of a tree of $\Ary^{(k)}_\Node$ is its degree are
combinatorial. Observe that $\Ary^{(2)}_\Leaf \simeq \BinaryTrees_\Leaf$
where $\BinaryTrees_\Leaf$ is defined in
Section~\ref{subsubsec:binary_trees}. Moreover, the generating series of
$\Ary^{(k)}_\Node$ satisfies the algebraic equation
\begin{equation}
    1 - \GenSeries_{\Ary^{(k)}_\Node}(t)
    + t {\GenSeries_{\Ary^{(k)}_\Node}(t)}^k
    = 0.
\end{equation}
and it is known~\cite{DM47} that
\begin{equation} \label{equ:number_k-ary_trees}
    \# \Ary^{(k)}_\Node(n) =
    \frac{1}{(k - 1) n + 1} \binom{kn}{n}.
\end{equation}
For instance, the sequences of integers associated with
$\Ary^{(k)}_\Node$ begin with
\begin{subequations}
\begin{equation}
    1, 1, 1, 1, 1, 1, 1, 1,
    \qquad k = 1,
\end{equation}
\begin{equation}
    1, 1, 2, 5, 14, 42, 132, 429, 1430,
    \qquad k = 2,
\end{equation}
\begin{equation}
    1, 1, 3, 12, 55, 273, 1428, 7752, 43263,
    \qquad k = 3,
\end{equation}
\begin{equation}
    1, 1, 4, 22, 140, 969, 7084, 53820, 420732,
    \qquad k = 4.
\end{equation}
\end{subequations}
The second, third, and fourth sequences are respectively
Sequences~\OEIS{A000108}, \OEIS{A001764}, and~\OEIS{A002293}
of~\cite{Slo}. These are known as the Fuss-Catalan numbers.
\medbreak

From now on, we call \Def{binary tree} any $2$-ary tree. If $\Tfr$ is a
binary tree and $u$ is an internal node of $\Tfr$, $u1$ and $u2$ are
nodes of $\Tfr$. We call $u1$ (resp. $u2$) the \Def{left} (resp.
\Def{right}) \Def{child} of $u$, and $\Tfr \Action u1$ (resp.
$\Tfr \Action u2$) the \Def{left} (resp. \Def{right}) \Def{subtree} of
$u$ in $\Tfr$. The left (resp. right) subtree of $\Tfr$ is the
\Def{left} (resp. \Def{right}) \Def{subtree} of the root of $\Tfr$.
Besides, a \Def{left} (resp. \Def{right}) \Def{comb tree} is a binary
tree $\Tfr$ such that for all internal nodes $u$ of $\Tfr$, all right
(resp. left) subtrees of $u$ are leaves. The \Def{infix order} induced
by $\Tfr$ is the total order on the set of its internal nodes defined
recursively by setting that all the internal nodes of $\Tfr \Action 1$
are smaller than the root of $\Tfr$, and that the root of $\Tfr$ is
smaller that all the internal nodes of~$\Tfr \Action 2$.
\medbreak

\subsubsection{Schröder trees} \label{subsubsec:schroder_trees}
A \Def{Schröder tree} is a planar rooted tree such that all internal
nodes are of arities $2$ or more. Some among the first Schröder trees
are
\begin{equation}
    \LeafPic, \quad
\,.
\end{equation}
This set of Schröder trees forms a subcollection $\Schroder$ of
$\PlanarRootedTrees$ expressing recursively as
\begin{equation} \label{equ:combinatorial_set_schroder_trees}
    \Schroder = \{\Leaf\} +
    \{\Node\} \times \List_{\geq 2}\left(\Schroder\right),
\end{equation}
where $\Leaf$ and $\Node$ are both atomic.
\medbreak

By structural induction on $\Schroder$ (which is
an inductive subcollection of $\PlanarRootedTrees$), it follows that
there are finitely many Schröder trees of a given arity $n$. For this
reason, the graded collection $\Schroder_\Leaf$ of all the Schröder
trees such that the size of a tree of $\Schroder_\Leaf$ is its arity is
combinatorial. Conversely, considering the degrees of the trees for
their sizes does not form a combinatorial graded collection since there
are infinitely many Schröder trees of degree $1$ (the corollas). The
generating series of $\Schroder_\Leaf$ satisfies the algebraic
quadratic equation
\begin{equation}
    t - (1 + t) \GenSeries_{\Schroder_\Leaf}(t)
    + 2 {\GenSeries_{\Schroder_\Leaf}(t)}^2 = 0.
\end{equation}
Let $\Narayana(n, k)$ be the number of binary trees of arity $n$ having
exactly $k$ internal nodes having an internal node as a left child.
Then, for all $0 \leq k \leq n - 2$, it is known~\cite{Nar55} that
\begin{equation} \label{equ:narayana_numbers}
    \Narayana(n, k) =
    \frac{1}{k + 1} \binom{n - 2}{k} \binom{n - 1}{k}.
\end{equation}
These are \Def{Narayana numbers}. The cardinalities of the sets
$\Schroder_\Leaf(n)$ hence express by
\begin{equation} \label{equ:enumeration_schroder_trees}
    \# \Schroder_\Leaf(n) =
    \sum_{k \in [0, n - 2]}
    2^k \; \Narayana(n, k),
\end{equation}
for all $n \in \N_{\geq 2}$. The sequence of integers associated with
$\Schroder_\Leaf$ begins by
\begin{equation}
    1, 1, 3, 11, 45, 197, 903, 4279,
\end{equation}
and forms Sequence~\OEIS{A001003} of~\cite{Slo}.
\medbreak

\subsection{Syntax trees} \label{subsec:syntax_trees}
We are now in position to introduce syntax trees. Such trees are,
roughly speaking, planar rooted trees where internal nodes are labeled
by objects of a fixed graded collection. These trees can be endowed with
two size functions (where the size is the degree or the arity), leading
to the definition of two graded collections of syntax trees.
\medbreak

\subsubsection{Collections of syntax trees}
\label{subsubsec:graded_collection_syntax_trees}
Let $C$ be an augmented graded collection. A \Def{syntax tree} on $C$
(or, for short, a \Def{$C$-syntax tree}) is a planar rooted tree $\Tfr$
endowed with a map $\omega_\Tfr : \TreeLanguage_\Node(\Tfr) \to C$
sending each internal node $u$ of $\Tfr$ of arity $k$ to an element of
size $k$ of $C$. This map $\omega_\Tfr$ is the \Def{labeling map} of
$\Tfr$. We say that an internal node $u$ of $\Tfr$ is \Def{labeled} by
$x \in C$ if $\omega_\Tfr(u) = x$. The collection $C$ is the
\Def{labeling collection} of $\Tfr$. The \Def{underlying planar rooted
tree} of $\Tfr$ is the planar rooted tree obtained by forgetting the map
$\omega_\Tfr$. For any $x \in C$, the \Def{corolla} labeled by $x$ is
the $C$-syntax tree $\Corolla{x}$ having exactly one internal node
labeled by $x$ and with $|x|$ leaves as children. All the notions about
planar rooted trees defined in Sections~\ref{subsec:planar_rooted_trees}
and~\ref{subsec:families_planar_rooted_trees} apply to $C$-syntax trees
as well. More precisely, for any property $P$ on planar rooted trees, we
say that $\Tfr$ satisfies the property $P$ if the underlying planar
rooted tree of $\Tfr$ satisfies~$P$. Moreover, the notions of suffix,
prefix, and factor subtrees of planar rooted trees naturally extend on
$C$-syntax trees by taking into account the labeling maps. In graphical
representations of a $C$-syntax tree $\Tfr$, instead of drawing each
internal node $u$ of $\Tfr$ by $\NodePic$, we draw $u$ by its
label~$\omega_\Tfr(u)$.
\medbreak

For instance, consider the labeling collection
$C := C(1) \sqcup C(2) \sqcup C(3)$ where $C(1) := \{\Asf, \Bsf\}$,
$C(2) := \{\Csf\}$, and $C(3) := \{\Dsf, \Esf\}$, and the planar rooted
tree
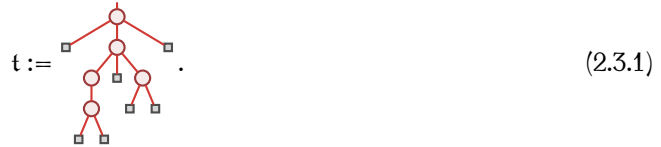
\begin{equation} \label{equ:example_syntax_tree}
    \Tfr :=
\,.
\end{equation}
By endowing $\Tfr$ with the labeling map defined by
$\omega_\Tfr(\epsilon) := \Esf$, $\omega_\Tfr(2) := \Dsf$,
$\omega_\Tfr(21) := \Asf$, $\omega_\Tfr(211) := \Csf$, and
$\omega_\Tfr(23) := \Csf$, $\Tfr$ is a $C$-syntax tree. This $C$-syntax
tree is depicted more concisely as
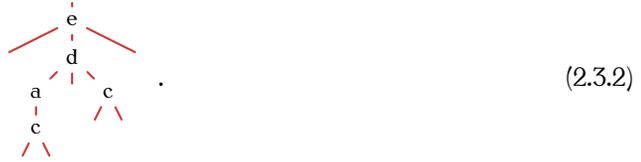
\begin{equation}
    %
\,.
\end{equation}
\medbreak

We denote by $\PlanarRootedTrees^C$ the graded collection of all the
$C$-syntax trees, where the size of a $C$-syntax tree $\Tfr$ is the size
of its underlying planar rooted tree in $\PlanarRootedTrees$. When $C$
is additionally combinatorial, by structural induction on planar rooted
trees, it follows that for any $\Tfr \in \PlanarRootedTrees$, there are
finitely many labeling maps $\omega_\Tfr$ for $\Tfr$. For this reason,
$\PlanarRootedTrees^C$ is in this case combinatorial. Besides, let
$\Ladders^C$, $\Corollas^C$, $\Ary^{(k), C}$, and $\Schroder^C$ be
respectively the subcollections of $\PlanarRootedTrees^C$ consisting in
the $C$-syntax trees whose underlying planar rooted trees are ladders,
corollas, $k$-ary trees, and Schröder trees. The concepts of inductive
subcollections of $\PlanarRootedTrees^C$ and of structural induction
presented in Section~\ref{subsubsec:structural_induction} extend
obviously on $C$-syntax trees.
\medbreak

\subsubsection{Alternative definition and generating series}
The graded collection $\PlanarRootedTrees^C$ can be described as
follows. Let $\Sca^C$ be the graded collection satisfying the relation
\begin{equation} \label{equ:expression_combinatorial_set_syntax_trees}
    \Sca^C = \{\Leaf\} +
    \{\Node\} \times \left(C \circ \Sca^C\right)
\end{equation}
where both $\Leaf$ and $\Node$ are atomic, and $\circ$ is the
composition product over combinatorial collections defined in
Section~\ref{subsubsec:products_collections}. Then, the combinatorial
collections $\PlanarRootedTrees^C$ and $\Sca^C$ are isomorphic through
the map $\phi : \PlanarRootedTrees^C \to \Sca^C$ of combinatorial
collections recursively defined, for any
$\Tfr \in \PlanarRootedTrees^C$ of root arity $k$, by
\begin{equation}
    \phi(\Tfr) :=
    \begin{cases}
        \Leaf \in \Sca^C & \mbox{if } \Tfr = \Leaf, \\
        \left(\Node,
            \left(
                \omega_\Tfr(\epsilon),
                (\phi(\Tfr_1), \dots, \phi(\Tfr_k))
            \right)
        \right)
        & \mbox{otherwise}.
    \end{cases}
\end{equation}
From this equivalence
and~\eqref{equ:expression_combinatorial_set_syntax_trees}, we obtain,
when $C$ is combinatorial, that the generating series of
$\PlanarRootedTrees^C$ satisfies
\begin{equation}
    \GenSeries_{\PlanarRootedTrees^C}(t) =
    t +
    t \; \GenSeries_C\left(\GenSeries_{\PlanarRootedTrees^C}(t)\right),
\end{equation}
where $\GenSeries_C(t)$ is the generating series of $C$. For instance,
by considering the combinatorial collection $C$ defined above, we have
$\GenSeries_C(t) = 2t + t^2 + 2t^3$, so that
\begin{equation}
    t + (2t - 1) \GenSeries_{\PlanarRootedTrees^C}(t)
    + t {\GenSeries_{\PlanarRootedTrees^C}(t)}^2
    + 2t {\GenSeries_{\PlanarRootedTrees^C}(t)}^3 = 0.
\end{equation}
\medbreak

\subsubsection{Subcollections of syntax trees}
For well-chosen augmented combinatorial collections $C$, it is
possible to recover a large part of the families of planar rooted
trees described in Section~\ref{subsec:families_planar_rooted_trees}.
Indeed, one has
\begin{math}
    \PlanarRootedTrees^{\N_{\geq 1}}
    \simeq \PlanarRootedTrees
\end{math},
\begin{math}
    \PlanarRootedTrees^{\N_{\geq 2}}
    \simeq \Schroder
\end{math},
and, when $\Node_k$ is an object of size $k \in \N_{\geq 1}$,
\begin{math}
    \PlanarRootedTrees^{\{\Node_k\}}
    \simeq \Ary^{(k)}
\end{math}.
\medbreak

\subsubsection{Alternative sizes}
Let $\PlanarRootedTrees^C_\Leaf$ be the graded collection of all the
$C$-syntax trees such that the size of a tree is its arity. One has
$\PlanarRootedTrees_\Leaf^C \simeq \Sca^C$ where $\Sca^C$ is the graded
collection defined
in~\eqref{equ:expression_combinatorial_set_syntax_trees} wherein $\Leaf$
is atomic and $\Node$ is of size $0$. When $C$ is combinatorial,
augmented, and has no object of size $1$, we can show by structural
induction on $\PlanarRootedTrees_\Leaf^C$ that there are finitely many
$C$-syntax trees of a given arity $n$. For this reason,
$\PlanarRootedTrees^C_\Leaf$ is combinatorial. In this case, the
generating series of $\PlanarRootedTrees_\Leaf^C$ satisfies
\begin{equation}
    \GenSeries_{\PlanarRootedTrees_\Leaf^C}(t) =
    t +
    \GenSeries_C
    \left(\GenSeries_{\PlanarRootedTrees_\Leaf^C}(t)\right).
\end{equation}
\medbreak

Let also $\PlanarRootedTrees_\Node^C$ be the graded collection of
all the $C$-syntax trees such that the size of a tree is its degree. One
has $\PlanarRootedTrees_\Node^C \simeq \Sca^C$ where $\Sca^C$ is the
graded collection defined
in~\eqref{equ:expression_combinatorial_set_syntax_trees} wherein $\Node$
is atomic and $\Leaf$ is of size $0$. When $C$ is augmented and is
finite, we can show by structural induction on
$\PlanarRootedTrees_\Node^C$ that there are finitely many $C$-syntax
trees of a given degree $n$. For this reason,
$\PlanarRootedTrees_\Node^C$ is combinatorial. In this case, the
generating series of $\PlanarRootedTrees_\Node^C$ satisfies
\begin{equation}
    \GenSeries_{\PlanarRootedTrees_\Node^C}(t) =
    1 +
    t\; \GenSeries_C
    \left(\GenSeries_{\PlanarRootedTrees_\Node^C}(t)\right).
\end{equation}
Observe that $\PlanarRootedTrees_\Node^C$ is not an augmented graded
collection.
\medbreak

\subsection{Syntax tree patterns and rewrite systems}
\label{subsec:rewrite_rules_syntax_trees}
We focus now the theory of rewrite systems on the particular case of
syntax trees. Intuitively, a rewrite rule on syntax trees works by
replacing factor subtrees in a syntax tree by other ones. We explain
techniques to prove termination and confluence of these particular
rewrite systems.
\medbreak

\subsubsection{Occurrence and avoidance of patterns}
Let $C$ be an augmented graded collection, and $\Sfr$ and $\Tfr$ be two
$C$-syntax trees. For any node $u$ of $\Tfr$, $\Sfr$ \Def{occurs} at
position $u$ in $\Tfr$ if $\Sfr$ is a factor subtree of $\Tfr$ rooted at
$u$. In this case, we say that $\Tfr$  \Def{admits an occurrence} of the
\Def{pattern} $\Sfr$. Conversely,  $\Tfr$ \Def{avoids} $\Sfr$ if there
is no occurrence of $\Sfr$ in~$\Tfr$. By extension, $\Tfr$ avoids a set
$P$ of $C$-syntax trees if $\Tfr$ avoids all the patterns of $P$. For
instance, consider the combinatorial collection
$C := C(2) \sqcup C(3)$ where $C(2) := \{\Asf, \Bsf\}$ and
$C(3) := \{\Csf\}$, and the $C$-syntax tree
\begin{equation}
    \Tfr :=
\,.
\end{equation}
Then, $\Tfr$ admits an occurrence of
\begin{equation}
    %
\end{equation}
at position $1$ and two occurrences of
\begin{equation}
    %
\end{equation}
at positions $11$ and~$21$.
\medbreak

\subsubsection{Grafting of syntax trees}
\label{subsubsec:grafting_syntax_trees}
Let $\Tfr$ be a $C$-syntax tree of arity $n$, $i \in [n]$, and $\Sfr$
be a $C$-syntax tree. The \Def{grafting} of $\Sfr$ onto the $i$th leaf
$u$ of $\Tfr$ is the $C$-syntax tree $\Rfr := \Tfr \circ_i \Sfr$ defined
as follows. The underlying planar rooted tree of $\Rfr$ admits the tree
language
\begin{equation} \label{equ:tree_language_grafting}
    \TreeLanguage(\Rfr)
    :=
    \left(\TreeLanguage(\Tfr) \setminus \{u\}\right)
    \cup
    \{u v : v \in \TreeLanguage(\Sfr)\},
\end{equation}
and the labeling map of $\Rfr$ satisfies, for any
$w \in \TreeLanguage_\Node(\Rfr)$,
\begin{equation}
    \omega_\Rfr(w) :=
    \begin{cases}
            \omega_\Tfr(w)
                & \mbox{if } w \in \TreeLanguage_\Node(\Tfr), \\
            \omega_\Sfr(v)
                & \mbox{otherwise (}
                w = u v \mbox{ and } v \in \TreeLanguage_\Node(\Sfr)
                \mbox{)}.
    \end{cases}
\end{equation}
Observe that by
Proposition~\ref{prop:bijection_planar_rooted_trees_prefix_languages},
$\Rfr$ is wholly specified by its tree language $\TreeLanguage(\Rfr)$
defined in~\eqref{equ:tree_language_grafting}. In more intuitive terms,
the tree $\Rfr$ is obtained by connecting the root of $\Sfr$ onto the
$i$th leaf of $\Tfr$. For instance, by considering the same labeling
collection $C$ as above,
\begin{equation}
\,.
\end{equation}
\medbreak

The operations $\circ_i$ thus defined are binary products
\begin{equation}
    \circ_i :
    \PlanarRootedTrees^C_\Leaf
    \times \PlanarRootedTrees^C_\Leaf
    \to \PlanarRootedTrees^C_\Leaf
\end{equation}
on $\PlanarRootedTrees^C_\Leaf$, in the sense of
Section~\ref{subsubsec:products_collections}. We call each $\circ_i$ a
\Def{grafting operation}. Since, for any $C$-syntax trees $\Tfr$ and
$\Sfr$, and $i \in [\Arity(\Tfr)]$,
\begin{equation}
    \Arity(\Tfr \circ_i \Sfr) = \Arity(\Tfr) + \Arity(\Sfr) - 1,
\end{equation}
these operations are $\DPlus$-compatible for the product $\DPlus$
defined by $n \DPlus m := n + m - 1$ for all $n, m \in \N_{\geq 1}$.
\medbreak

\subsubsection{Complete grafting of syntax trees}
\label{subsubsec:complete_grafting_syntax_trees}
Let $\Tfr$ be a $C$-syntax tree of arity $n$, and let $\Sfr_1$, \dots,
$\Sfr_n$ be $C$-syntax trees. The \Def{complete grafting} of $\Sfr_1$,
\dots, $\Sfr_n$ onto $\Tfr$ is the $C$-syntax tree
$\circ^{(n)}(\Tfr, \Sfr_1, \dots, \Sfr_n)$ defined by
\begin{equation} \label{equ:complete_grafting_syntax_trees}
    \circ^{(n)}(\Tfr, \Sfr_1, \dots, \Sfr_n)
    :=
    \left(\dots \left(\left(\Tfr \circ_n \Sfr_n\right) \circ_{n - 1}
    \Sfr_{n - 1}\right)\dots\right)
    \circ_1 \Sfr_1.
\end{equation}
In more intuitive terms, the tree
\begin{math}
    \circ^{(n)}(\Tfr, \Sfr_1, \dots, \Sfr_n)
\end{math}
is obtained by connecting the root of each $\Sfr_i$ onto the $i$th leaf
of~$\Tfr$. For instance, by considering the same labeling collection $C$
as before,
\begin{equation}
    \circ^{(4)}
    \left(
\,.
\end{equation}
\medbreak

The operations $\circ^{(n)}$ thus defined are products
\begin{equation}
    \circ^{(n)} :
    \PlanarRootedTrees^C_\Leaf
    \times
    \left(\PlanarRootedTrees^C_\Leaf\right)^n
    \to
    \PlanarRootedTrees^C_\Leaf
\end{equation}
of arity $n + 1$ on $\PlanarRootedTrees^C_\Leaf$. We call each
$\circ^{(n)}$ a \Def{complete grafting operation}. Since, for any
$C$-syntax trees $\Tfr$, $\Sfr_1$, \dots, $\Sfr_n$ such that
$n = \Arity(\Tfr)$,
\begin{equation}
    \Arity\left(\circ^{(n)}(\Tfr, \Sfr_1, \dots, \Sfr_n)\right) =
    \Arity(\Tfr) + \Arity(\Sfr_1) + \dots + \Arity(\Sfr_n) - n,
\end{equation}
these operations are $\DPlus$-compatible for the product $\DPlus$
defined in Section~\ref{subsubsec:grafting_syntax_trees}. Moreover, to
gain concision, we shall denote by $\Tfr \circ [\Sfr_1, \dots, \Sfr_n]$
the $C$-syntax tree
\begin{math}
    \circ^{(n)}(\Tfr, \Sfr_1, \dots, \Sfr_n)
\end{math}.
\medbreak

\subsubsection{Rewrite systems}
Let for any $m, i \in \N_{\geq 1}$ the operation
$\circleddash^{(m)}_i$ defined as follows.
For any $C$-syntax trees $\Tfr$, $\Rfr$, $\Sfr_1$, \dots, $\Sfr_m$
where $\Tfr$ is of arity $n \geq i$, we set
\begin{equation}
    \circleddash^{(m)}_i\left(\Tfr, \Rfr, \Sfr_1, \dots, \Sfr_m\right)
    :=
    \Tfr \circ_i
        \left(\Rfr \circ \left[\Sfr_1, \dots, \Sfr_m\right]\right).
\end{equation}
These operations $\circleddash^{(m)}_i$ thus defined are products
\begin{equation}
    \circleddash^{(m)}_i :
    \PlanarRootedTrees^C_\Leaf
    \times
    \PlanarRootedTrees^C_\Leaf
    \times
    \left(\PlanarRootedTrees^C_\Leaf\right)^m
    \to
    \PlanarRootedTrees^C_\Leaf
\end{equation}
of arity $m + 2$ on $\PlanarRootedTrees^C_\Leaf$. It is easy to see
that $\circleddash^{(m)}_i$ is $\DPlus$-compatible for the product
$\DPlus$ considered in Sections~\ref{subsubsec:grafting_syntax_trees}
and~\ref{subsubsec:complete_grafting_syntax_trees}.
\medbreak

Let $(\PlanarRootedTrees^C_\Leaf, \Rew)$ be a rewrite system. We
denote by $(\PlanarRootedTrees^C_\Leaf, \RewTrees)$ the
$\Pca$-closure of $(\PlanarRootedTrees^C_\Leaf, \Rew)$, where
\begin{equation}
    \Pca := \left\{\circleddash^{(m)}_i : m, i \in \N_{\geq 1}\right\}.
\end{equation}
We call $(\PlanarRootedTrees^C_\Leaf, \RewTrees)$ the \Def{closure} of
$(\PlanarRootedTrees^C_\Leaf, \Rew)$. In other terms, $\RewTrees$ is the
rewrite rule satisfying
\begin{equation} \label{equ:closure_rewrite_rule_syntax_trees}
    \Tfr \circ_i
        \left(\Rfr \circ \left[\Sfr_1, \dots, \Sfr_m\right]\right)
    \RewTrees
    \Tfr \circ_i
         \left(\Rfr' \circ \left[\Sfr_1, \dots, \Sfr_m\right]\right)
\end{equation}
for any $C$-syntax trees $\Tfr$, $\Rfr$, $\Rfr'$, $\Sfr_1$, \dots,
$\Sfr_m$ where $\Tfr$ of arity $n$, $i \in [n]$, and $\Rfr \Rew \Rfr'$.
In intuitive terms, one has
$\Qfr \RewTrees \Qfr'$ for two $C$-syntax trees $\Qfr$ and $\Qfr'$ if
there are two $C$-syntax trees $\Rfr$ and $\Rfr'$ such that
$\Rfr \Rew \Rfr'$ and, by replacing an occurrence of $\Rfr$ by $\Rfr'$
in $\Qfr$, we obtain~$\Qfr'$. For instance, by considering the same
labeling set $C$ as before, let $(\PlanarRootedTrees^C_\Leaf, \Rew)$
be the rewrite system defined by
\begin{equation} \label{equ:example_rewrite_rule_syntax_trees}
\,.
\end{equation}
One has the following chain of rewritings
\begin{equation}
    %
\,.
\end{equation}
\medbreak

Observe by the way that the right rotation operation on binary trees
considered in Section~\ref{subsubsec:examples_posets}
(see~\eqref{equ:right_rotation_Tamari}) can be expressed as the closure
of the rewrite system $(\PlanarRootedTrees^B_\Leaf, \Rew)$ such that
$B := B(2) := \{\Bsf\}$ defined by
\begin{equation}
    %
\,.
\end{equation}
\medbreak

In this dissertation, we shall mainly consider rewrite systems
$(\PlanarRootedTrees^C_\Leaf, \RewTrees)$ defined as closures of rewrite
systems $(\PlanarRootedTrees^C_\Leaf, \Rew)$ such that the number of
pairs $(\Tfr, \Tfr')$ satisfying $\Tfr \Rew \Tfr'$ is finite. We say in
this case that $(\PlanarRootedTrees^C_\Leaf, \RewTrees)$ is of
\Def{finite type}. In this context, the \Def{degree} of
$(\PlanarRootedTrees^C_\Leaf, \RewTrees)$ is the maximal degree among
the $C$-syntax trees appearing as left members of $\Rew$. The
\Def{arity} of $(\PlanarRootedTrees^C_\Leaf, \RewTrees)$ is the maximal
arity among the $C$-syntax trees appearing as left (or right) members
of~$\Rew$.
\medbreak

\subsubsection{Proving termination}
\label{subsubsec:proving_termination}
We have observed in Section~\ref{subsubsec:termination} that termination
invariants provide tools to show that a combinatorial rewrites ystem is
terminating. This idea extends on rewrite systems on syntax trees
defined as closures of other ones in the following way.
\medbreak

Let $(\PlanarRootedTrees^C_\Leaf, \Rew)$ be a combinatorial rewrite
system and $(\PlanarRootedTrees^C_\Leaf, \RewTrees)$ its closure.
Assume that $\theta : \PlanarRootedTrees^C_\Leaf \to \Qca$ is a
termination invariant for $(\PlanarRootedTrees^C_\Leaf, \Rew)$, where
$(\Qca, \Ord)$ is a poset. We say that $\theta$ is \Def{compatible with
the closure} if, for any $C$-syntax trees $\Rfr$ and $\Rfr'$ such that
$\Rfr \Rew \Rfr'$, the inequality
\begin{equation} \label{equ:compatibility_termination_invariant}
    \theta\left(
    \Sfr \circ_i \left(\Rfr \circ
    \left[\Qfr_1, \dots, \Qfr_k\right]\right)
    \right)
    \OrdStrict
    \theta\left(
    \Sfr \circ_i \left(\Rfr' \circ
    \left[\Qfr_1, \dots, \Qfr_k\right]\right)
    \right)
\end{equation}
holds for all $C$-syntax trees $\Sfr$, $\Qfr_1$, \dots, $\Qfr_k$, where
$k := \Arity(\Rfr) = \Arity(\Rfr')$. Now, as a consequence
of~\eqref{equ:closure_rewrite_rule_syntax_trees} and
Lemma~\ref{lem:terminating_rewrite_rules_posets}, one has the following
result.
\medbreak

\begin{Proposition} \label{prop:compatible_terminating_invariant}
    Let $C$ be an augmented combinatorial collection without object of
    size~$1$, $(\PlanarRootedTrees^C_\Leaf, \Rew)$ be a rewrite
    system, and $(\PlanarRootedTrees^C_\Leaf, \RewTrees)$ be the
    closure of $(\PlanarRootedTrees^C_\Leaf, \Rew)$. If
    $\theta : \PlanarRootedTrees^C_\Leaf \to \Qca$ is a termination
    invariant for $(\PlanarRootedTrees^C_\Leaf, \Rew)$ and $\theta$ is
    compatible with the closure,
    $(\PlanarRootedTrees^C_\Leaf, \RewTrees)$ is terminating.
\end{Proposition}
\medbreak

Consider for instance the rewrite rule
$(\PlanarRootedTrees^C_\Leaf, \Rew)$ defined
by~\eqref{equ:example_rewrite_rule_syntax_trees}. By setting
$\Qca := \N^2$ and $\Ord$ as the lexicographic order on $\N^2$, let us
define the map $\theta : \PlanarRootedTrees^C_\Leaf \to \Qca$, for any
$C$-syntax tree $\Tfr$, by
\begin{math}
    \theta(\Tfr) := \left(\deg(\Tfr), \TamariInvariant(\Tfr)\right)
\end{math},
where
\begin{equation} \label{equ:Tamari_invariant}
    \TamariInvariant(\Tfr) :=
    \sum_{\substack{
        u \in \TreeLanguage_\Node(\Tfr) \\
        u \mbox{ \footnotesize of arity } 2
    }}
    \deg(\Tfr \Action u2).
\end{equation}
In other words, $\TamariInvariant(\Tfr)$ is the sum, for all internal
and binary nodes $u$ of $\Tfr$, of the number of internal nodes
appearing in the $2$nd suffix subtrees of $u$. One can check that
$\theta(\Tfr) \OrdStrict \theta(\Tfr')$ for all the $C$-syntax trees
$\Rfr$ and $\Rfr'$ such that $\Rfr \Rew \Rfr'$. Indeed,
\begin{equation}
    \theta\left(

    \right).
\end{equation}
Moreover, the fact that $\theta$ is compatible with the closure is a
straightforward verification. Therefore, the closure
$(\PlanarRootedTrees^C_\Leaf, \RewTrees)$ of
$(\PlanarRootedTrees^C_\Leaf, \Rew)$ is terminating.
\medbreak

\subsubsection{Proving confluence}
In the same way as the tool to show that a rewrite system on $C$-syntax
trees is terminating presented in
Section~\ref{subsubsec:proving_termination}, we present here a tool to
prove that rewrite systems on syntax trees defined as closures of other
ones are confluent. This criterion requires now some precise properties.
\medbreak

\begin{Proposition} \label{prop:confluence_rewrite_rule_syntax_trees}
    Let $C$ be an augmented combinatorial collection without object of
    size~$1$, $(\PlanarRootedTrees^C_\Leaf, \Rew)$ be a rewrite
    system, and $(\PlanarRootedTrees^C_\Leaf, \RewTrees)$ be the
    closure of $(\PlanarRootedTrees^C_\Leaf, \Rew)$. If
    $(\PlanarRootedTrees^C_\Leaf, \RewTrees)$ is terminating, is of
    finite type, has $\ell \geq 0$ as degree, and all branching pairs of
    $(\PlanarRootedTrees^C_\Leaf, \RewTrees)$ consisting in trees with
    $2 \ell - 1$ internal nodes or less are joinable, then
    $(\PlanarRootedTrees^C_\Leaf, \RewTrees)$ is confluent.
\end{Proposition}
\begin{proof}
    Assume that all the hypotheses of the statement hold. Let $\Tfr$
    be a branching tree of $(\PlanarRootedTrees^C_\Leaf, \RewTrees)$ and
    $\{\Rfr_1, \Rfr_2\}$ be a branching pair for $\Tfr$. We have thus
    $\Tfr \RewTrees \Rfr_1$ and $\Tfr \RewTrees \Rfr_2$. By definition
    of $\RewTrees$, there are four $C$-syntax trees $\Tfr'_1$,
    $\Rfr'_1$, $\Tfr'_2$, and $\Rfr'_2$ such that
    $\Tfr'_1 \Rew \Rfr'_1$, $\Tfr'_2 \Rew \Rfr'_2$, and $\Rfr_1$ (resp.
    $\Rfr_2$) is obtained by replacing an occurrence of $\Tfr'_1$
    (resp. $\Tfr'_2$) rooted at a node $u_1$ (resp. $u_2$) by $\Rfr'_1$
    (resp. $\Rfr'_2$) in $\Tfr$. We have now two cases to consider,
    depending on the positions of the nodes $u_1$ and $u_2$ in~$\Tfr$.
    \begin{enumerate}[fullwidth,label={\it Case \arabic*.}]
        \item Assume first that the occurrences of $\Tfr'_1$ and
        $\Tfr'_2$ at positions $u_1$ and $u_2$ in $\Tfr$ do not share
        any internal node of $\Tfr$. Then, $\Rfr_1$ (resp. $\Rfr_2$)
        admits an occurrence of $\Tfr'_2$ (resp. $\Tfr'_1$) at a
        position $u'_2$ (resp. $u'_1$). These positions $u'_1$ and
        $u'_2$ are obtained from the original positions $u_1$ and $u_2$
        of the occurrences of $\Tfr'_1$ and $\Tfr'_2$ in $\Tfr$. Now,
        let $\Tfr'$ be the tree obtained by replacing the occurrence of
        $\Tfr'_2$ at position $u'_2$ by $\Rfr'_2$ in $\Rfr_1$.
        Equivalently, due to the above assumption, $\Tfr'$ is the tree
        obtained by replacing the occurrence of $\Tfr'_1$ at position
        $u'_1$ by $\Rfr'_1$ in~$\Rfr_2$. Thereby, we have
        $\Rfr_1 \RewTrees \Tfr'$ and $\Rfr_2 \RewTrees \Tfr'$, showing
        that the branching pair $\{\Rfr_1, \Rfr_2\}$ is joinable.
        \smallbreak

        \item Otherwise, the occurrences of $\Tfr'_1$ and $\Tfr'_2$ at
        positions $u_1$ and $u_2$ in $\Tfr$ share at least one internal
        node of $\Tfr$. Denote by $\Sfr$ the factor subtree of $\Tfr$,
        rooted at a node $v$, of the smallest degree and whose internal
        nodes come from the internal nodes of $\Tfr$ involved in the
        occurrences of $\Tfr'_1$ and $\Tfr'_2$ at positions $u_1$ and
        $u_2$. Observe that $v = u_1$ or $v = u_2$. Let us denote by
        $v_1$ (resp. $v_2$) the position of the occurrence of $\Tfr'_1$
        (resp. $\Tfr'_2$) in $\Sfr$. Let $\Sfr_1$ (resp. $\Sfr_2$) be
        the tree obtained by replacing the occurrence of $\Tfr'_1$
        (resp. $\Tfr'_2$) at position $v_1$ (resp. $v_2$) by $\Rfr'_1$
        (resp. $\Rfr'_2$) in $\Sfr$. Now, due to the above assumption,
        by the minimality of the degree of $\Sfr$, $\Sfr$ has at most
        $2 \ell - 1$ internal nodes. Now, since we have
        $\Sfr \RewTrees \Sfr_1$ and $\Sfr \RewTrees \Sfr_2$, $\Sfr$ is
        a branching tree for $\RewTrees$ and $\{\Sfr_1, \Sfr_2\}$ is a
        branching pair for $\Sfr$. By hypothesis, $\{\Sfr_1, \Sfr_2\}$
        is joinable, so that there is a tree $\Sfr'$ such that
        $\Sfr_1 \RewTreesRT \Sfr'$ and $\Sfr_2 \RewTreesRT \Sfr'$. By
        setting $\Tfr'$ as the tree obtained by replacing the
        occurrence of $\Sfr$ at position $v$ by $\Sfr'$ in $\Tfr$, we
        finally have $\Rfr_1 \RewTreesRT \Tfr'$ and
        $\Rfr_2 \RewTreesRT \Tfr'$, showing that $\{\Rfr_1, \Rfr_2\}$ is
        joinable.
    \end{enumerate}
    We have shown that all branching pairs of $\RewTrees$ are joinable.
    Since $(\PlanarRootedTrees^C_\Leaf, \RewTrees)$ is terminating,
    this implies by the diamond lemma (Lemma~\ref{lem:diamond_lemma})
    that $(\PlanarRootedTrees^C_\Leaf, \RewTrees)$ is confluent.
\end{proof}
\medbreak

Proposition~\ref{prop:confluence_rewrite_rule_syntax_trees} leads to an
algorithmic way to check if a terminating rewrite system
$(\PlanarRootedTrees^C_\Leaf, \RewTrees)$ defined as the closure of an
other one $(\PlanarRootedTrees^C_\Leaf, \Rew)$ is confluent by
enumerating all the $C$-syntax trees $\Tfr$ of degrees at most
$2\ell - 1$ (where $\ell$ is the degree of
$(\PlanarRootedTrees^C_\Leaf, \RewTrees)$) and by computing the parts
$G_\Tfr$ of the rewriting graphs of
$(\PlanarRootedTrees^C_\Leaf, \RewTrees)$ consisting in the trees
reachable from~$\Tfr$. If each $G_\Tfr$ contains exactly one normal form
(which correspond to a vertex with no outgoing edge in $G_\Tfr$),
$(\PlanarRootedTrees^C_\Leaf, \RewTrees)$ is confluent.
\medbreak

For instance, by considering the same labeling set $C$ as above, let
$(\PlanarRootedTrees^C_\Leaf, \Rew)$ be the rewrite system defined by
\begin{equation}
\,.
\end{equation}
The degree of the closure $(\PlanarRootedTrees^C_\Leaf, \RewTrees)$ of
$(\PlanarRootedTrees^C_\Leaf, \Rew)$ is $\ell := 2$ and it is possible
to show that $(\PlanarRootedTrees^C_\Leaf, \RewTrees)$ is terminating.
Consider
\begin{equation}
    \Tfr :=
    %
\,,
\end{equation}
which is a $C$-syntax tree of degree $2\ell - 1 = 3$. The graph $G_\Tfr$
associated with $\Tfr$ is of the form
\begin{equation}
    %
\,,
\end{equation}
and shows that $(\PlanarRootedTrees^C_\Leaf, \RewTrees)$ is not
confluent. Indeed, $\Tfr$ is a non-joinable branching tree. On the other
hand, consider the rewrite system $(\PlanarRootedTrees^C_\Leaf, \Rew)$
defined by
\begin{equation}
    %
\,.
\end{equation}
The degree of the closure $(\PlanarRootedTrees^C_\Leaf, \RewTrees)$ of
$(\PlanarRootedTrees^C_\Leaf, \Rew)$ is $\ell := 2$ and here also, it
is possible to show that $(\PlanarRootedTrees^C_\Leaf, \RewTrees)$ is
terminating. Consider
\begin{equation}
    \Tfr :=
    %
\,,
\end{equation}
a $C$-syntax tree of degree $2\ell - 1 = 3$. The graph $G_\Tfr$
associated with $\Tfr$ is of the form
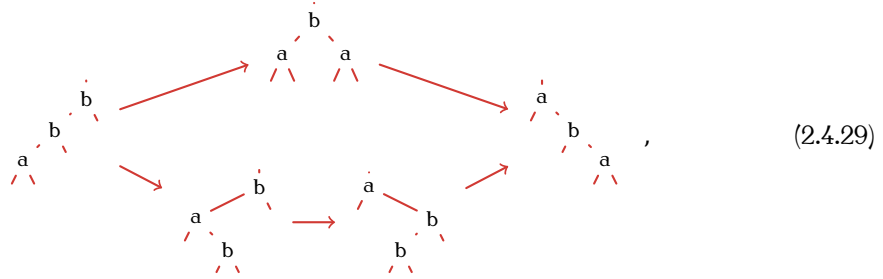
\begin{equation}
    %
        };
        \draw[Arc](2220000_bba)--(2200200_baa);
        \draw[Arc](2220000_bba)--(2202000_bab);
        \draw[Arc](2202000_bab)--(2022000_abb);
        \draw[Arc](2022000_abb)--(2020200_aba);
        \draw[Arc](2200200_baa)--(2020200_aba);
    \end{tikzpicture}\,,
\end{equation}
This graph satisfies the required property stated above, and, as a
systematic study of cases shows, all other graphs $G_\Sfr$ where $\Sfr$
is a $C$-syntax tree of degree $3$ or less, also. For this reason,
$(\PlanarRootedTrees^C_\Leaf, \RewTrees)$ is confluent.
\medbreak

\section{Combinatorial objects} \label{sec:secondary_objects}
This last section of the chapter contains a list of definitions about
combinatorial objects appearing in some next chapters.
\medbreak

\subsection{Other kinds of trees}
Let us set some definitions about two other kinds to trees: rooted trees
and colored syntax trees.
\medbreak

\subsubsection{Rooted trees} \label{subsubsec:rooted_trees}
Let $\RootedTrees$ be the graded collection satisfying the relation
\begin{equation}
    \RootedTrees = \{\Node\} \times \Multiset(\RootedTrees).
\end{equation}
where $\Node$ is an atomic object called \Def{node}. We call
\Def{rooted tree} each object of $\RootedTrees$. By definition, a rooted
tree $\Tfr$ is an ordered pair
$(\Node, \lbag \Tfr_1, \dots, \Tfr_k \rbag)$ where
$\lbag \Tfr_1, \dots, \Tfr_k \rbag$ is a multiset of rooted trees. Like
the case of planar rooted trees, this definition is  recursive. For
instance,
\begin{equation} \label{equ:examples_rooted_trees}
    (\Node, \emptyset), \quad
    (\Node, \lbag (\Node, \emptyset) \rbag), \quad
    (\Node, \lbag (\Node, \emptyset), (\Node, \emptyset) \rbag), \quad
    (\Node, \lbag (\Node, \emptyset), (\Node, \emptyset),
        (\Node, \emptyset) \rbag), \quad
    (\Node, \lbag (\Node, \lbag (\Node, \emptyset), (\Node, \emptyset)
        \rbag) \rbag),
\end{equation}
are rooted trees. If $\Tfr = (\Node, \lbag \Tfr_1, \dots, \Tfr_k \rbag)$
is a rooted tree, each $\Tfr_i$, $i \in [k]$, is a \Def{suffix subtree}
of~$\Tfr$.
\medbreak

Rooted trees are different kinds of trees than planar rooted trees
presented in Section~\ref{sec:trees}. The difference is due to the fact
that rooted trees are defined by using multisets of rooted trees, while
planar rooted trees are defined by using lists of planar rooted trees.
Hence, the order of the suffix subtrees of a rooted tree is not
significant.
\medbreak

By drawing each rooted tree by a node $\NodePic$ attached below it to
its subtrees by means of edges $\EdgePic$, the rooted trees
of~\eqref{equ:examples_rooted_trees} are depicted by
\begin{equation}
    \NodePic \quad,
    \begin{tikzpicture}[scale=.3,Centering]
        \node[Node](1)at(0.00,-1.00){};
        \node[Node](0)at(0.00,0.00){};
        \draw[Edge](1)--(0);
        \node(r)at(0.00,1){};
        \draw[Edge](r)--(0);
    \end{tikzpicture}\,, \quad
    \begin{tikzpicture}[scale=.2,Centering]
        \node[Node](0)at(0.00,-1.50){};
        \node[Node](2)at(2.00,-1.50){};
        \node[Node](1)at(1.00,0.00){};
        \draw[Edge](0)--(1);
        \draw[Edge](2)--(1);
        \node(r)at(1.00,1.5){};
        \draw[Edge](r)--(1);
    \end{tikzpicture}\,, \quad
    \begin{tikzpicture}[xscale=.3,yscale=.2,Centering]
        \node[Node](0)at(0.00,-2.00){};
        \node[Node](2)at(1.00,-2.00){};
        \node[Node](3)at(2.00,-2.00){};
        \node[Node](1)at(1.00,0.00){};
        \draw[Edge](0)--(1);
        \draw[Edge](2)--(1);
        \draw[Edge](3)--(1);
        \node(r)at(1.00,1.50){};
        \draw[Edge](r)--(1);
    \end{tikzpicture}\,, \quad
    \begin{tikzpicture}[scale=.2,Centering]
        \node[Node](1)at(0.00,-2.67){};
        \node[Node](3)at(2.00,-2.67){};
        \node[Node](0)at(1.00,0.00){};
        \node[Node](2)at(1.00,-1.33){};
        \draw[Edge](1)--(2);
        \draw[Edge](2)--(0);
        \draw[Edge](3)--(2);
        \node(r)at(1.00,1.5){};
        \draw[Edge](r)--(0);
    \end{tikzpicture}\,.
\end{equation}
By definition of the product and multiset operations over combinatorial
collections, the size of a rooted tree $\Tfr$ satisfies
\begin{equation} \label{equ:size_rooted_tree}
    |\Tfr| := 1 + \sum_{i \in [k]} |\Tfr_i|.
\end{equation}
The sequence of integers associated with $\RootedTrees$ begins by
\begin{equation}
    1, 1, 2, 4, 9, 20, 48, 115,
\end{equation}
and forms Sequence~\OEIS{A000081} of~\cite{Slo}.
\medbreak

\subsubsection{Colored syntax trees}
\label{subsubsec:colored_syntax_trees}
Let $\CFr$ be a set of colors and $C$ be a $\CFr$-colored collection
(see Section~\ref{subsubsec:colored_collections}). A
\Def{$\CFr$-colored $C$-syntax tree} is a triple $(a, \Tfr, u)$ where
$\Tfr$ is a $C$-syntax tree, $a \in \CFr$, $u \in \CFr^{\Arity(\Tfr)}$,
and for any internal nodes $u$ and $v$ of $\Tfr$ such that $v$ is the
$i$th child of $u$, $\Out(y) = \In_i(x)$ where $x$ (resp. $y$) is the
label of $u$ (resp. $v$). The set of all $\CFr$-colored $C$-syntax trees
is denoted by $\ColoredPlanarRootedTrees^C$. This set is a
$\CFr$-colored collection by setting that $\Out(((a, \Tfr, u))) := a$
and $\In((a, \Tfr, u)) := u$ for all
$(a, \Tfr, u) \in \ColoredPlanarRootedTrees^C$. By a slight abuse of
notation, if $u$ is an internal node of $\Tfr$, we denote by $\Out(u)$
(resp. $\In(u)$) the color $\Out(x)$ (resp. word of colors $\In(x)$)
where $x$ is the label of $u$. We say that a $\CFr$-colored $C$-syntax
tree $\Tfr$ is \Def{monochrome} if $C$ is a monochrome colored
collection. In graphical representations of a $\CFr$-colored $C$-syntax
tree $(a, \Tfr, u)$, we draw $\Tfr$ together with its output color above
its root and its input color $u(i)$ below its $i$th leaf for
any~$i \in [|u|]$.
\medbreak

For instance, consider the set of colors $\CFr := \{1, 2\}$ and the
$\CFr$-colored collection $C$ defined by
$C := C(2) \sqcup C(3)$ with $C(2) := \{\Asf, \Bsf\}$,
$C(3) := \{\Csf\}$, $\Out(\Asf) := 1$, $\Out(\Bsf) := 2$,
$\Out(\Csf) := 1$, $\In(\Asf) := 11$, $\In(\Bsf) := 21$, and
$\In(\Csf) := 221$. The tree
\begin{equation}
    \begin{tikzpicture}[xscale=.37,yscale=.15,Centering]
        \node(0)at(0.00,-6.50){};
        \node(10)at(8.00,-9.75){};
        \node(12)at(10.00,-9.75){};
        \node(2)at(2.00,-6.50){};
        \node(4)at(3.00,-3.25){};
        \node(5)at(4.00,-9.75){};
        \node(7)at(5.00,-9.75){};
        \node(8)at(6.00,-9.75){};
        \node[NodeST](1)at(1.00,-3.25){\begin{math}\Bsf\end{math}};
        \node[NodeST](11)at(9.00,-6.50){\begin{math}\Asf\end{math}};
        \node[NodeST](3)at(3.00,0.00){\begin{math}\Csf\end{math}};
        \node[NodeST](6)at(5.00,-6.50){\begin{math}\Csf\end{math}};
        \node[NodeST](9)at(7.00,-3.25){\begin{math}\Asf\end{math}};
        \node(r)at(3.00,2.75){};
        \draw[Edge](0)--(1);
        \draw[Edge](1)--(3);
        \draw[Edge](10)--(11);
        \draw[Edge](11)--(9);
        \draw[Edge](12)--(11);
        \draw[Edge](2)--(1);
        \draw[Edge](4)--(3);
        \draw[Edge](5)--(6);
        \draw[Edge](6)--(9);
        \draw[Edge](7)--(6);
        \draw[Edge](8)--(6);
        \draw[Edge](9)--(3);
        \draw[Edge](r)--(3);
        \node[LeafLabel,above of=r]{\begin{math}1\end{math}};
        \node[LeafLabel,below of=0]{\begin{math}2\end{math}};
        \node[LeafLabel,below of=2]{\begin{math}1\end{math}};
        \node[LeafLabel,below of=4]{\begin{math}2\end{math}};
        \node[LeafLabel,below of=5]{\begin{math}2\end{math}};
        \node[LeafLabel,below of=7]{\begin{math}2\end{math}};
        \node[LeafLabel,below of=8]{\begin{math}1\end{math}};
        \node[LeafLabel,below of=10]{\begin{math}1\end{math}};
        \node[LeafLabel,below of=12]{\begin{math}1\end{math}};
    \end{tikzpicture}
\end{equation}
is a $\CFr$-colored $C$-syntax tree. Its degree is $5$, its arity is
$8$, and its height is $3$. Moreover, its output color is $1$ and its
word of input colors is $21222111$. Besides, $(1, \Leaf, 1)$ and
$(1, \Leaf, 2)$ are two $\CFr$-colored $C$-syntax trees of degree $0$
and arity~$1$.
\medbreak

Let $(a, \Tfr, u)$ and $(b, \Sfr, v)$ be two $\CFr$-colored $C$-syntax
trees and  $i \in [|\Arity(\Tfr)|]$. If $b = u(i)$, the \Def{grafting}
on $\Sfr$ onto the $i$th leaf of $\Tfr$ is defined by
\begin{equation} \label{equ:grafting_colored_syntax_trees}
    (a, \Tfr, u) \circ_i (b, \Sfr, v) :=
    (a, \Tfr \circ_i \Sfr, u \mapsfrom_i v),
\end{equation}
where $u \mapsfrom_i v$ is the word obtained by replacing the $i$th
letter of $u$ by $v$, and the second occurrence of $\circ_i$
in~\eqref{equ:grafting_colored_syntax_trees} is the grafting of
syntax trees defined in Section~\ref{subsubsec:grafting_syntax_trees}.
For instance, by considering the same labeling $\CFr$-colored collection
as above,
\begin{equation}
\,.
\end{equation}
\medbreak

Let $\ColoredPlanarRootedTrees^C_\Leaf$ be the $\CFr$-colored collection
of the $\CFr$-colored $C$-syntax trees. The operations $\circ_i$ thus
defined are binary products
\begin{equation}
    \circ_i :
    \ColoredPlanarRootedTrees^C_\Leaf \times
    \ColoredPlanarRootedTrees^C_\Leaf \to
    \ColoredPlanarRootedTrees^C_\Leaf
\end{equation}
on $\ColoredPlanarRootedTrees^C_\Leaf$, in the sense of
Section~\ref{subsubsec:products_collections}. We call each $\circ_i$ a
\Def{grafting operation}. By seeing $\ColoredPlanarRootedTrees^C_\Leaf$
as a graded collection (see
Section~\ref{subsubsec:colored_collections}), the $\circ_i$,
$i \in \N_{\geq 1}$, are $\DPlus$-compatible products, where $\DPlus$
is the operation considered in
Section~\ref{subsubsec:grafting_syntax_trees}. Observe also that, due
to the condition on the colors between the two operands to the
operation, the $\circ_i$ are partial products.
\medbreak

Most of the notions exposed in
Section~\ref{subsec:rewrite_rules_syntax_trees} about syntax trees
and rewrite systems on syntax trees naturally extend on colored syntax
trees like, among others, the notions of occurrences of patterns, the
complete grafting operations, and the criteria offered by
Propositions~\ref{prop:compatible_terminating_invariant}
and~\ref{prop:confluence_rewrite_rule_syntax_trees} to respectively
prove the termination and the confluence of rewrite system on syntax
trees.
\medbreak

\subsection{Configurations of chords} \label{subsubsec:configurations}
Configurations of chords are very classical combinatorial objects
defined as collections of diagonals and edges in regular polygons. The
literature abounds of studies of various kinds of configurations. One
can cite for instance~\cite{DRS10} about triangulations, \cite{FN99}
about noncrossing configurations, and~\cite{CP92} about
multi-trian\-gulations. We provide here definitions about them and
consider a generalization of configurations wherein the edges and
diagonals are labeled on a set.
\medbreak

\subsubsection{Polygons}
A \Def{polygon} of \Def{size} $n \geq 1$ is a directed graph $\Pfr$ on
the set of vertices $[n + 1]$. An \Def{arc} of $\Pfr$ is a pair of
integers $(x, y)$ with $1 \leq x < y \leq n + 1$, a \Def{diagonal} is
an arc $(x, y)$ different from $(x, x + 1)$ and $(1, n + 1)$, and an
\Def{edge} is an arc of the form $(x, x + 1)$ and different from
$(1, n + 1)$. We denote by $\Arcs_\Pfr$ (resp. $\Diagonals_\Pfr$,
$\Edges_\Pfr$) the set of all arcs (resp. diagonals, edges) of $\Pfr$.
For any $i \in [n]$, the \Def{$i$th edge} of $\Pfr$ is the edge
$(i, i + 1)$, and the arc $(1, n + 1)$ is the \Def{base} of~$\Pfr$.
\medbreak

In our graphical representations, each polygon is depicted so that its
base is the bottommost segment, vertices are implicitly numbered from
$1$ to $n + 1$ in the clockwise direction, and the diagonals are not
drawn. For example,
\begin{equation}
    \Pfr :=
    \begin{tikzpicture}[scale=.85,Centering]
        \node[CliquePoint](1)at(-0.50,-0.87){};
        \node[CliquePoint](2)at(-1.00,-0.00){};
        \node[CliquePoint](3)at(-0.50,0.87){};
        \node[CliquePoint](4)at(0.50,0.87){};
        \node[CliquePoint](5)at(1.00,0.00){};
        \node[CliquePoint](6)at(0.50,-0.87){};
        \draw[CliqueEmptyEdge](1)edge[]node[CliqueLabel]{}(2);
        \draw[CliqueEmptyEdge](1)edge[]node[CliqueLabel]{}(6);
        \draw[CliqueEmptyEdge](2)edge[]node[CliqueLabel]{}(3);
        \draw[CliqueEmptyEdge](3)edge[]node[CliqueLabel]{}(4);
        \draw[CliqueEmptyEdge](4)edge[]node[CliqueLabel]{}(5);
        \draw[CliqueEmptyEdge](5)edge[]node[CliqueLabel]{}(6);
        \node[left of=1,node distance=3mm,font=\scriptsize]
            {\begin{math}1\end{math}};
        \node[left of=2,node distance=3mm,font=\scriptsize]
            {\begin{math}2\end{math}};
        \node[above of=3,node distance=3mm,font=\scriptsize]
            {\begin{math}3\end{math}};
        \node[above of=4,node distance=3mm,font=\scriptsize]
            {\begin{math}4\end{math}};
        \node[right of=5,node distance=3mm,font=\scriptsize]
            {\begin{math}5\end{math}};
        \node[right of=6,node distance=3mm,font=\scriptsize]
            {\begin{math}6\end{math}};
    \end{tikzpicture}
\end{equation}
is a polygon of size $5$. Its set of all diagonals is
\begin{equation}
    \Diagonals_\Pfr =
    \{(1, 3), (1, 4), (1, 5), (2, 4), (2, 5), (2, 6),
    (3, 5), (3, 6), (4, 6)\},
\end{equation}
its set of all edges is
\begin{equation}
    \Edges_\Pfr = \{(1, 2), (2, 3), (3, 4), (4, 5), (5, 6)\},
\end{equation}
and its set of all arcs is
\begin{equation}
    \Arcs_\Pfr = \Diagonals_\Pfr \sqcup \Edges_\Pfr \sqcup \{(1, 6)\}.
\end{equation}
\medbreak

\subsubsection{Configurations}
For any set $S$, an \Def{$S$-configuration} (or a \Def{configuration}
when $S$ is known without ambiguity) is a polygon $\Cfr$ endowed with a
partial function
\begin{equation}
    \phi_\Cfr : \Arcs_\Cfr \to S.
\end{equation}
When $\phi_\Cfr((x, y))$ is defined, we say that the arc $(x, y)$ is
\Def{labeled} and we denote it by $\Cfr(x, y)$. When the base of $\Cfr$
is labeled, we denote it by $\Cfr_0$, and when the $i$th edge of $\Cfr$
is labeled, we denote it by $\Cfr_i$.
\medbreak

In our graphical representations, we shall represent any
$S$-configuration $\Cfr$ by drawing a polygon of the same size as the
one of $\Cfr$ following the conventions explained before, and by
labeling its arcs accordingly. For instance
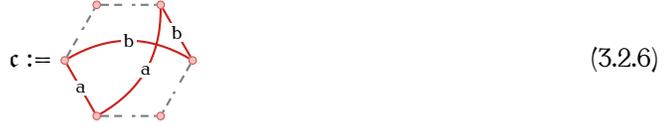
\begin{equation}
    \Cfr :=
    \begin{tikzpicture}[scale=.85,Centering]
        \node[CliquePoint](1)at(-0.50,-0.87){};
        \node[CliquePoint](2)at(-1.00,-0.00){};
        \node[CliquePoint](3)at(-0.50,0.87){};
        \node[CliquePoint](4)at(0.50,0.87){};
        \node[CliquePoint](5)at(1.00,0.00){};
        \node[CliquePoint](6)at(0.50,-0.87){};
        \draw[CliqueEdge](1)edge[]node[CliqueLabel]
            {\begin{math}\Asf\end{math}}(2);
        \draw[CliqueEmptyEdge](1)edge[]node[CliqueLabel]{}(6);
        \draw[CliqueEmptyEdge](2)edge[]node[CliqueLabel]{}(3);
        \draw[CliqueEmptyEdge](3)edge[]node[CliqueLabel]{}(4);
        \draw[CliqueEdge](4)edge[]node[CliqueLabel]
            {\begin{math}\Bsf\end{math}}(5);
        \draw[CliqueEmptyEdge](5)edge[]node[CliqueLabel]{}(6);
        \draw[CliqueEdge](1)edge[bend right=30]node[CliqueLabel]
            {\begin{math}\Asf\end{math}}(4);
        \draw[CliqueEdge](2)edge[bend left=30]node[CliqueLabel]
            {\begin{math}\Bsf\end{math}}(5);
    \end{tikzpicture}
\end{equation}
is an $\{\Asf, \Bsf, \Cfr\}$-configuration. The arcs $(1, 2)$ and
$(1, 4)$ of $\Cfr$ are labeled by $\Asf$, the arcs $(2, 5)$ and
$(4, 5)$ are labeled by $\Bsf$, and the other arcs are unlabeled.
\medbreak

\subsubsection{Additional definitions}
Let us now provide some definitions and statistics on configurations.
Let $\Cfr$ be a configuration of size $n$. The \Def{skeleton} of $\Cfr$
is the undirected graph $\Skel(\Cfr)$ on the set of vertices $[n + 1]$
and such that for any $x < y \in [n + 1]$, there is an arc $\{x, y\}$
in $\Skel(\Cfr)$ if $(x, y)$ is labeled in $\Cfr$. The \Def{degree} of
a vertex $x$ of $\Cfr$ is the number of vertices adjacent to $x$ in
$\Skel(\Cfr)$. The \Def{degree} $\Degr(\Cfr)$ of $\Cfr$ is the maximal
degree among its vertices. Two (non-necessarily labeled) diagonals
$(x, y)$ and $(x', y')$ of $\Cfr$ are \Def{crossing} if
$x < x' < y < y'$ or $x' < x < y' < y$. The \Def{crossing} of a labeled
diagonal $(x, y)$ of $\Cfr$ is the number of labeled diagonals
$(x', y')$ such that $(x, y)$ and $(x', y')$ are crossing. The
\Def{crossing} $\Cros(\Cfr)$ of $\Cfr$ is the maximal crossing among
its labeled diagonals. When $\Cros(\Cfr) = 0$, there are no crossing
diagonals in $\Cfr$ and in this case, $\Cfr$ is \Def{noncrossing}. A
(non-necessarily labeled) arc $(x', y')$ is \Def{nested} in a
(non-necessarily labeled) arc $(x, y)$ of $\Cfr$ if
$x \leq x' < y' \leq y$. We say that $\Cfr$ is \Def{nesting-free} if for
any labeled arcs $(x, y)$ and $(x', y')$ of $\Cfr$ such that $(x', y')$
is nested in $(x, y)$, $(x, y) = (x', y')$. Besides, $\Cfr$ is
\Def{acyclic} if $\Skel(\Cfr)$ is acyclic. When $\Cfr$ has no labeled
edges nor labeled base, $\Cfr$ is \Def{white}. If $\Cfr$ has no labeled
diagonals, $\Cfr$ is a \Def{bubble}. A \Def{triangle} is a configuration
of size~$2$. Obviously, all triangles are bubbles, and all bubbles are
noncrossing.
\medbreak

\subsection{Prographs} \label{subsec:prographs}
We present here prographs, that are combinatorial objects modeling
operations with several inputs and several outputs. These objects are
elements of free pros (see Section~\ref{subsec:pros} of
Chapter~\ref{chap:algebra}) and admit many different definitions. A
first one consists in defining prographs (called in this context
diagrams) through an equivalence relation~\cite{Laf11}. A second one
consists in defining prographs (called in this context directed
$(m, n)$-graphs) by graphs satisfying some conditions~\cite{Mar08}. We
chose to define these objects by using the tools offered by the theory
of bigraded collections. Our approach is however similar to the one
of~\cite{Laf11}.
\medbreak

For all bigraded collections $C$ considered in this section, if
$x$ is an object of index $(p, q) \in \N^2$, the \Def{input} (resp.
\Def{output}) \Def{arity} of $x$ is $|x|_\ArityIn := p$ (resp.
$|x|_\ArityOut := q$).
\medbreak

\subsubsection{Sequences of wires} \label{subsubsec:wires}
Let $\Wires$ be the bigraded collection satisfying
\begin{equation}
    \Wires(p, q) :=
    \begin{cases}
        \left\{\Unit_p\right\} & \mbox{if } p = q, \\
        \emptyset & \mbox{otherwise}.
    \end{cases}
\end{equation}
We call $\Unit_1$ the \Def{wire} and each element $\Unit_p$ the
\Def{sequence of wires} of arity $p$. Each $\Unit_p$ is depicted by $p$
vertical lines. For instance, $\Unit_5$ is depicted as
\begin{equation}
    \begin{tikzpicture}[xscale=.17,yscale=.15,Centering]
        \node[Leaf](S1)at(0,0){};
        \node[Leaf](S2)at(2,0){};
        \node[Leaf](S3)at(4,0){};
        \node[Leaf](S4)at(6,0){};
        \node[Leaf](S5)at(8,0){};
        \node[Leaf](E1)at(0,-5){};
        \node[Leaf](E2)at(2,-5){};;
        \node[Leaf](E3)at(4,-5){};
        \node[Leaf](E4)at(6,-5){};
        \node[Leaf](E5)at(8,-5){};
        \draw[Edge](S1)--(E1);
        \draw[Edge](S2)--(E2);
        \draw[Edge](S3)--(E3);
        \draw[Edge](S4)--(E4);
        \draw[Edge](S5)--(E5);
    \end{tikzpicture}\,.
\end{equation}
\medbreak

\subsubsection{Preprographs}
From now on, $C$ is a bigraded collection such that
$C(p, q) = \emptyset$ if $p = 0$ or $q = 0$.
\medbreak

An \Def{elementary prograph} $x$ on $C$ (or, for short, an
\Def{elementary $C$-prograph}) is an object $\Asf$ of $C(p, q)$. We
represent $x$ as a rectangle labeled by $\Asf$ with $p$ incoming edges
(below the rectangle) and $q$ outgoing edges (above the rectangle). For
instance, if $\Asf \in C(2, 3)$, the elementary prograph $\Asf$ is
depicted as
\begin{equation}
\,.
\end{equation}
\medbreak

An \Def{enriched elementary prograph} $x$ on $C$ (or, for short, an
\Def{enriched elementary $C$-prograph}) is a triple $(k, \Asf, \ell)$
where $k, \ell \in \N$, and $\Asf$ is an elementary $C$-prograph. These
objects form a bigraded collection where the index of $(k, \Asf, \ell)$
is $(k + |\Asf|_\ArityIn + \ell, k + |\Asf|_\ArityOut + \ell)$.
We represent $x$ by drawing from left to right $\Unit_k$, $\Asf$, and
$\Unit_\ell$. For instance, the enriched elementary prograph
$(1, \Asf, 2)$ where $\Asf \in C(2, 3)$ is depicted as
\begin{equation}
    %
\,.
\end{equation}
\medbreak

A \Def{preprograph} $x$ on $C$ (or, for short, a \Def{$C$-preprograph})
is a sequence $(x_1, \dots, x_k)$ of enriched elementary $C$-prographs
such that $|x_i|_\ArityIn = |x_{i + 1}|_\ArityOut$ for any
$i \in [k - 1]$. These objects form a bigraded collection
$\Preprographs^C$ where the index of $x$ is
$(|x_k|_\ArityIn, |x_1|_\ArityOut)$. We represent $x$ by drawing each
enriched elementary prograph $x_i$, $i \in [k]$, vertically, where
$x_1$ is at the top and $x_k$ is at the bottom. We moreover draw dashed
lines for all $i \in [k - 1]$ between $x_i$ and $x_{i + 1}$. For
instance, the $C$-preprograph
\begin{equation}
    x := ((1, \Asf, 5), (5, \Bsf, 2), (7, \Asf, 1))
\end{equation}
where $\Asf \in C(2, 2)$ and $\Bsf \in C(3, 1)$ is depicted as
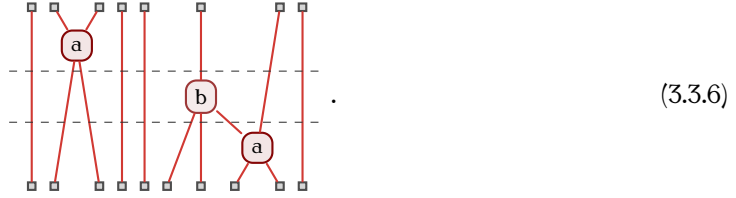
\begin{equation} \label{equ:example_preprograph}
    \begin{tikzpicture}[xscale=.3,yscale=.17,Centering]
        \node[Leaf](S0)at(-1,0){};
        \node[Leaf](S1)at(0,0){};
        \node[Leaf](S2)at(2,0){};
        \node[Leaf](S22)at(3,0){};
        \node[Leaf](S222)at(4,0){};
        \node[Leaf](S3)at(6.5,0){};
        \node[Leaf](S4)at(10,0){};
        \node[Leaf](S5)at(11,0){};
        \node[Operator](N1)at(1,-3){\begin{math}\Asf\end{math}};
        \node[OperatorColorC](N2)at(6.5,-7){\begin{math}\Bsf\end{math}};
        \node[Operator](N3)at(9,-11){\begin{math}\Asf\end{math}};
        \node[Leaf](E0)at(-1,-14){};
        \node[Leaf](E1)at(0,-14){};
        \node[Leaf](E2)at(2,-14){};
        \node[Leaf](E22)at(3,-14){};
        \node[Leaf](E222)at(4,-14){};
        \node[Leaf](E3)at(5,-14){};
        \node[Leaf](E4)at(6.5,-14){};
        \node[Leaf](E5)at(8,-14){};
        \node[Leaf](E6)at(10,-14){};
        \node[Leaf](E7)at(11,-14){};
        \draw[Edge](S0)--(E0);
        \draw[Edge](N1)--(S1);
        \draw[Edge](N1)--(S2);
        \draw[Edge](N1)--(E1);
        \draw[Edge](N1)--(E2);
        \draw[Edge](N2)--(S3);
        \draw[Edge](N2)--(E3);
        \draw[Edge](N2)--(E4);
        \draw[Edge](N2)--(N3);
        \draw[Edge](N3)--(S4);
        \draw[Edge](N3)--(E5);
        \draw[Edge](N3)--(E6);
        \draw[Edge](S5)--(E7);
        \draw[Edge](S22)--(E22);
        \draw[Edge](S222)--(E222);
        \draw[draw=ColBlack!80,dashed](-2,-5)--(12,-5);
        \draw[draw=ColBlack!80,dashed](-2,-9)--(12,-9);
    \end{tikzpicture}\,.
\end{equation}
\medbreak

\subsubsection{Prographs} \label{subsubsec:prographs}
Let $(\Preprographs^C, \Rew)$ be the rewrite system satisfying
\begin{equation}
    \left(\left(k_1, \Asf, k_2 + |\Bsf|_\ArityOut + k_3\right),
    \left(k_1 + |\Asf|_\ArityIn + k_2, \Bsf, k_3\right)\right)
    \Rew
    \left(\left(k_1 + |\Asf|_\ArityOut + k_2, \Bsf, k_3\right),
    \left(k_1, \Asf, k_2 + |\Bsf|_\ArityIn + k_3\right)\right)
\end{equation}
where $\Asf$ and $\Bsf$ are elementary $C$-prographs and
$k_1, k_2, k_3 \geq 0$. Pictorially,
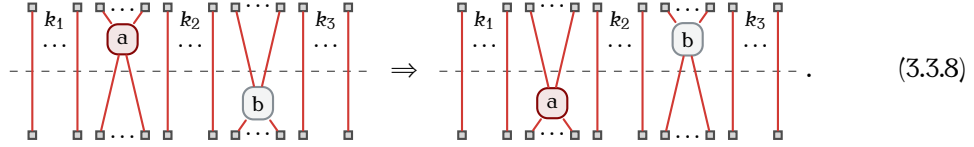
\begin{equation}
\,.
\end{equation}
\medbreak

If $x := (x_1, \dots, x_k)$ and $y := (y_1, \dots, y_\ell) $ are two
$C$-preprographs such that $|x|_\ArityIn = |y|_\ArityOut$, we denote by
$x \circ y$ the $C$-preprograph $(x_1, \dots, x_k, y_1, \dots, y_\ell)$.
Let also $\Product$ be the ternary product on $\Preprographs^C$ defined
by $\Product(x, y, z) := x \circ y \circ z$ where $x$, $y$, and $z$ are
three $C$-preprographs satisfying $|x|_\ArityIn = |y|_\ArityOut$
and $|y|_\ArityIn = |z|_\ArityOut$.
\medbreak

Let $(\Preprographs^C, \Rew_\Product)$ be the $\{\Product\}$-closure of
$(\Preprographs^C, \Rew)$ (it is possible to define a product $\DPlus$
on $\N^2$ so that $\Product$ is $\DPlus$-compatible) and $\Equiv$ be
the reflexive, symmetric, and transitive closure of $\Rew_\Product$. Let
the bigraded collection $\Prographs^C$ defined by
\begin{equation}
    \Prographs^C := \Preprographs^C/_\Equiv + \Wires.
\end{equation}
We call \Def{prograph} on $C$ (or, for short, a \Def{$C$-prograph})
any element $x$ of $\Prographs^C$. When $x$ is an object of
$\Preprographs^C/_\Equiv$, we represent $x$ by considering the drawing
of any preprograph of the $\Equiv$-equivalence class of $x$ and by
letting the elementary prographs constituting $x$ to move vertically
along the edges. When $x$ is an object of $\Wires$, we represent $x$ as
explained in Section~\ref{subsubsec:wires}. For instance, the prograph
having the preprograph of~\eqref{equ:example_preprograph} in its
$\Equiv$-equivalence class is depicted as
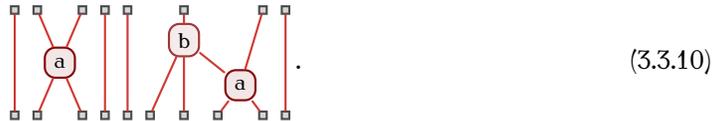
\begin{equation} \label{equ:example_prograph}
    \begin{tikzpicture}[xscale=.3,yscale=.2,Centering]
        \node[Leaf](S0)at(-1,0){};
        \node[Leaf](S1)at(0,0){};
        \node[Leaf](S2)at(2,0){};
        \node[Leaf](S22)at(3,0){};
        \node[Leaf](S222)at(4,0){};
        \node[Leaf](S3)at(6.5,0){};
        \node[Leaf](S4)at(10,0){};
        \node[Leaf](S5)at(11,0){};
        \node[Operator](N1)at(1,-3.5){\begin{math}\Asf\end{math}};
        \node[OperatorColorC](N2)at(6.5,-2){\begin{math}\Bsf\end{math}};
        \node[Operator](N3)at(9,-5){\begin{math}\Asf\end{math}};
        \node[Leaf](E0)at(-1,-7){};
        \node[Leaf](E1)at(0,-7){};
        \node[Leaf](E2)at(2,-7){};
        \node[Leaf](E22)at(3,-7){};
        \node[Leaf](E222)at(4,-7){};
        \node[Leaf](E3)at(5,-7){};
        \node[Leaf](E4)at(6.5,-7){};
        \node[Leaf](E5)at(8,-7){};
        \node[Leaf](E6)at(10,-7){};
        \node[Leaf](E7)at(11,-7){};
        \draw[Edge](S0)--(E0);
        \draw[Edge](N1)--(S1);
        \draw[Edge](N1)--(S2);
        \draw[Edge](N1)--(E1);
        \draw[Edge](N1)--(E2);
        \draw[Edge](N2)--(S3);
        \draw[Edge](N2)--(E3);
        \draw[Edge](N2)--(E4);
        \draw[Edge](N2)--(N3);
        \draw[Edge](N3)--(S4);
        \draw[Edge](N3)--(E5);
        \draw[Edge](N3)--(E6);
        \draw[Edge](S5)--(E7);
        \draw[Edge](S22)--(E22);
        \draw[Edge](S222)--(E222);
    \end{tikzpicture}\,.
\end{equation}
\medbreak

The \Def{degree} $\deg(x)$ of a prograph $x$ is defined in the following
way. When $x$ is an object of $\Preprographs^C/_\Equiv$, $\deg(x)$ is
the length of $x'$ where $x'$ is any preprograph of the
$\Equiv$-equivalence class $x$. In other terms, $\deg(x)$ is the number
of elementary prographs constituting $x$. When $x$ is an object of
$\Wires$, $\deg(x) := 0$. For instance, the prograph
of~\eqref{equ:example_prograph} has $3$ as degree, and each sequence of
wires $\Unit_p$, $p \in \N$, has $0$ as degree.
\medbreak

\subsubsection{Operations on prographs}
\label{subsubsec:operations_prographs}
Let $x$ and $y$ be two $C$-prographs such that
$|x|_\ArityIn = |y|_\ArityOut$. The \Def{vertical composition}
$x \circ y$ of $x$ and $y$ is defined as follows. When $x$ (resp. $y$)
is a sequence of wires, $x \circ y$ is equal to $y$ (resp. $x$).
Otherwise, $x$ and $y$ are not sequences of wires, and $x \circ y$ is
the $C$-prograph $[x' \circ y']_\Equiv$ where $x'$ and $y'$ are
respectively any elements of the $\Equiv$-equivalence classes $x$ and
$y$, and $\circ$ is the operation on preprographs defined in
Section~\ref{subsubsec:prographs}. For instance,
\begin{equation}
\,.
\end{equation}
\medbreak

Let $x$ and $y$ be two $C$-prographs. The \Def{horizontal composition}
$x * y$ of $x$ and $y$ is defined as follows. If $x$ and $y$ are both
the sequences of wires $\Unit_p$ and $\Unit_q$ for some $p, q \in \N$,
\begin{equation}
    \Unit_p * \Unit_q := \Unit_{p + q}.
\end{equation}
If $x$ is the sequence of wires $\Unit_p$ for
a $p \in \N$ and $y$ is an object of $\Preprographs^C/_\Equiv$,
\begin{equation}
    \Unit_p * y :=
    \left[\left(\left(p + k_1, \Asf_1, \ell_1\right), \dots,
        \left(p + k_r, \Asf_r, \ell_r\right)\right)\right]_\Equiv,
\end{equation}
where
\begin{math}
    \left(\left(k_1, \Asf_1, \ell_1\right),
    \dots,
    \left(k_r, \Asf_r, \ell_r\right)\right)
\end{math}
is any preprograph in the $\Equiv$-equivalence class $y$. Similarly, if
$x$ is an object of $\Preprographs^C/_\Equiv$ and $y$ is the sequence
of wires $\Unit_p$, $p \in \N$,
\begin{equation}
    x * \Unit_p :=
    \left[\left(\left(k_1, \Asf_1, \ell_1 + p\right), \dots,
        \left(k_r, \Asf_r, \ell_r + p\right)\right)\right]_\Equiv,
\end{equation}
where
\begin{math}
    \left(\left(k_1, \Asf_1, \ell_1\right),
    \dots,
    \left(k_r, \Asf_r, \ell_r\right)\right)
\end{math}
is any preprograph in the $\Equiv$-equivalence class $x$. Finally, when
$x$ and $y$ are both objects of $\Preprographs^C/_\Equiv$, $x * y$ is
defined, by using the particular cases for the horizontal composition
explained above and the vertical composition, by
\begin{equation}
    x * y :=
    \left(x * \Unit_{|y|_\ArityOut}\right)
    \circ \left(\Unit_{|x|_\ArityIn} * y\right).
\end{equation}
For instance,
\begin{equation}
\,.
\end{equation}
\medbreak

\subsection{Alternating sign matrices} \label{subsec:ASMs}
We recall here some definitions about alternating sign matrices and
usual statistics on them.
\medbreak

\subsubsection{Alternating sign matrices and six-vertex configurations}
An \Def{alternating sign matrix}~\cite{MRR83}, or an \Def{ASM} for
short, of size $n$ is a square matrix of order $n$ with entries in the
alphabet $\{\Zero, \Plus, \Minus\}$ such that every row and column
starts and ends by $0$ or by $\Plus$ and in every row and column, the
$\Plus$ and the $\Minus$ alternate. For instance,
\begin{equation}
    \delta :=
    \Matrix{
        \Zero & \Plus & \Zero & \Zero & \Zero \\
        \Zero & \Zero & \Plus & \Zero & \Zero \\
        \Plus & \Minus & \Zero & \Zero & \Plus \\
        \Zero & \Plus & \Minus & \Plus & \Zero \\
        \Zero & \Zero & \Plus & \Zero & \Zero}
\end{equation}
is an ASM of size $5$. The number $a_n$ of these objects of size $n$
satisfies
\begin{equation}
    a_n = \prod_{0 \leq i \leq n - 1} \frac{(3i + 1)!}{(n + i)!},
\end{equation}
a formula conjectured in~\cite{MRR83} and proven independently by
Zeilberger~\cite{Zei96} and Kuperberg~\cite{Kup96}.
\medbreak

\begin{figure}[ht]
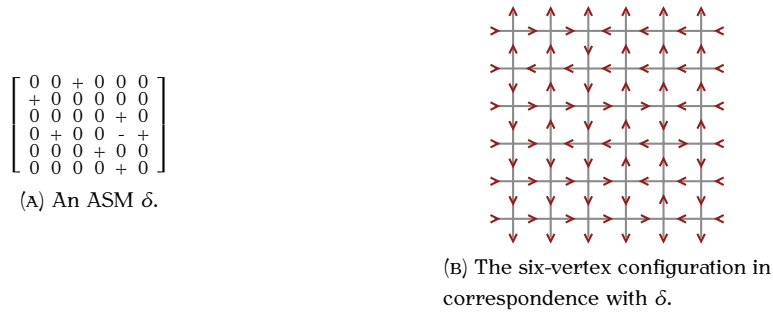

    \centering
    \subfloat[][An ASM $\delta$.]{
    \begin{minipage}[c]{.45\textwidth}
    \centering
    \begin{math}
    \Matrix{
        0 & 0 & \Plus & 0 & 0 & 0 \\
        \Plus & 0 & 0 & 0 & 0 & 0 \\
        0 & 0 & 0 & 0 & \Plus & 0 \\
        0 & \Plus & 0 & 0 & \Minus & \Plus \\
        0 & 0 & 0 & \Plus & 0 & 0 \\
        0 & 0 & 0 & 0 & \Plus & 0
    }
    \end{math}
    \end{minipage}
    \label{subfig:objects_bijection_ASMs_ASM}}
    \subfloat[][The six-vertex configuration in correspondence
    with $\delta$.]{
    \begin{minipage}[c]{.45\textwidth}
    \centering
    \begin{math}

    \end{math}
    \end{minipage}
    \label{subfig:objects_bijection_ASMs_six_vertex}}
    \caption{An ASM and a six-vertex configuration in correspondence.}
    \label{fig:objects_bijection_ASMs}
\end{figure}

A \Def{six-vertex configuration} of size $n$ is an $n \times n$ square
grid with oriented edges so that each vertex has two incoming and two
outgoing edges. There are six possible configurations for each vertex,
whence the name. A six-vertex configuration satisfies the \Def{domain
wall boundary condition} if all its horizontal (resp. vertical) edges
on the boundary are oriented inwardly (resp. outwardly).
Figure~\ref{subfig:objects_bijection_ASMs_six_vertex} shows an example
of such an object. In what follows, we shall exclusively and implicitly
consider six-vertex configurations satisfying the domain wall boundary
condition.
\medbreak

Six-vertex configurations of size $n$ are in one-to-one correspondence
with ASMs of the same size. To compute the ASM in correspondence
with a six-vertex configuration, we replace each of its vertices by a
symbol $0$, $\Plus$, or $\Minus$ according to the rules described in
Table~\ref{tab:correspondence_ASM}.
\begin{table}[ht]
    \centering
    \begin{tabular}{c||c|c|c|c|c|c}
        Statistics & $\NE$ & $\SW$ & $\SE$ & $\NW$ & $\OI$ & $\IO$
            \\ \hline
        ASM entry & $0$ & $0$ & $0$ & $0$ & $\Plus$ & $\Minus$
            \\ \hline
        Six-vertex-configuration &
        \begin{math}\SixVertexNE\end{math} &
        \begin{math}\SixVertexSW\end{math} &
        \begin{math}\SixVertexSE\end{math} &
        \begin{math}\SixVertexNW\end{math} &
        \begin{math}\SixVertexOI\end{math} &
        \begin{math}\SixVertexIO\end{math}
    \end{tabular}
    \bigbreak

    \caption[Six-vertex configurations and statistics on ASMs.]
    {Correspondence between entries of ASMs, vertices of six-vertex
    configurations, and statistics on six-vertex-configurations.}
    \label{tab:correspondence_ASM}
\end{table}
Reciprocally, to recover a six-vertex configuration from an ASM
$\delta$, we first replace each nonzero entry of $\delta$ by the
corresponding vertex configuration (see the last two columns of
Table~\ref{tab:correspondence_ASM}). Then, for each zero entry of
$\delta$, we look at the sum $\ell$ (resp. $a$) of the entries (a
$\Plus$ counts as $1$ and a $\Minus$ counts as $-1$) to the left (resp.
above) of it and in the same row (resp. column). By the alternating
property of the ASMs, $\ell$ and $a$ belong to $\{0, 1\}$. Now, set in
$\delta$ the configuration $\SixVertexW$ (resp. $\SixVertexE$) if
$\ell = 1$ (resp. $\ell = 0$) together with the configuration
$\SixVertexS$ (resp. $\SixVertexN$) if $a = 1$ (resp. $a = 0$).
Figure~\ref{fig:objects_bijection_ASMs} shows an example.
\medbreak

\subsubsection{Statistics on alternating sign matrices}
It is possible to define several statistics on ASMs by counting how
many entries of an ASM play a special role, seen as vertices of the
six-vertex configurations in correspondence.
\medbreak

Let us denote by $\NE(\delta)$ (resp. $\SW(\delta)$, $\SE(\delta)$,
$\NW(\delta)$, $\OI(\delta)$, $\IO(\delta)$) the number of vertices
$\NE$ (resp. $\SW$, $\SE$, $\NW$, $\OI$, $\IO$) in the six-vertex
configuration in bijection with the ASM $\delta$ (see
Table~\ref{tab:correspondence_ASM}). Let
$\AllStatsA := \{\SE,\NW,\SW,\NE\}$ be the set of the statistics
counting the four configurations of $0$ and $\AllStatsB := \{\IO,\OI\}$
be the set of the statistics counting the two nonzero configurations.
\medbreak

Let us end this section on ASMs by stating the following result
establishing some symmetries satisfied by these statistics.
\medbreak

\begin{Proposition} \label{prop:symmetries_statistics_ASMs}
    Let $\delta$ be an ASM of size $n$. Then,
    \begin{equation}
        \SE(\delta) = \NW(\delta), \quad
        \NE(\delta) = \SW(\delta), \quad
        \OI(\delta) = \IO(\delta) + n.
    \end{equation}
\end{Proposition}
\medbreak


\chapter{Algebraic combinatorics} \label{chap:algebra}
One of the main activities in algebraic combinatorics consists in
developing interactions between combinatorics (enumerative combinatorics
and even computer science) and algebra. The benefits are twofold: one
obtains combinatorial properties and results by seeing combinatorial
objects under an algebraic framework, and studying algebraic structures
with the help of combinatorics leads to general algebraic results.
\smallbreak

The first direction consists in endowing collections of combinatorial
objects with operations. This provides a framework to collect
combinatorial and enumerative properties on the objects by exploring
natural and usual algebraic questions on the obtained algebraic
structures. To be a little more precise, let $C$ be a collection and
$\K \Angle{C}$ be the linear span of $C$ where $\K$ is any field. To
get a better understanding of properties of the objects of $C$, we
endow $\K \Angle{C}$ with operations or co-operations. In this way, we
can ask about the behavior of these operations under different bases of
$\K \Angle{C}$ (leading to discovering links between different products,
for instance, the shifted shuffle product and the shifted concatenation
products of permutations are the same ones~\cite{DHT02,DHNT11}), minimal
generating sets of $\K \Angle{C}$ (leading to describe the objects of
$C$ as assemblies of elementary building blocks), morphisms involving
$\K \Angle{C}$ and other linear spans of collections (leading to
discover symmetries of $C$---useful for enumeration problems---, or
establishing links between $C$ and other
collections~\cite{LR98,DHT02,HNT05}).
\smallbreak

The second direction consists, on the contrary, in seeing abstract
algebraic structures as linear spans of combinatorial objects endowed
with operations. This process is known as a combinatorial realization
of an algebraic structure. To be more concrete, given a category of
algebras defined by the relations their (co)operations have to satisfy,
the problem consists in understanding the free object on a set
$\GeneratingSet$ of generators. This reinforces the understanding of the
category since all other ones are, in most cases, quotients or
substructures of free ones. The literature contains a lot of such
constructions. For instance, free Lie algebras are realized in terms of
Lyndon words and concatenation operations~\cite{Ret93}, free pre-Lie
algebras in terms of rooted trees and grafting operations~\cite{CL01},
free dendriform algebras in terms of binary trees and shuffling
operations~\cite{Lod01}, free duplicial algebras in terms of binary
trees and over and under operations~\cite{Lod08}, and free Zinbiel
algebras using words and half-shuffle operations~\cite{Lod95}.
\smallbreak

Additionally to very classical algebraic structures like magmas,
monoids, groups, posets, and associative algebras, in our work we
consider Hopf bialgebras~\cite{Car07,GR16},
operads~\cite{Mar08,LV12,Men15}, and pros~\cite{Lei04,Mar08}. The
purpose of this chapter is to present a unified approach to work with
these structures. We introduce in this way the notion of polynomial
spaces and of biproducts, that are operations working with several
inputs and several outputs. All the aforementioned algebraic structures
can be seen as particular cases of these objects.
\smallbreak

This chapter begins in Section~\ref{sec:polynomial_spaces} by defining
polynomial and series spaces on collections. Then, in
Section~\ref{sec:bialgebras}, we introduce biproducts, bialgebras, and
list some examples of such algebraic structures. Finally, in
Sections~\ref{sec:hopf_bialgebras}, \ref{sec:operads},
and~\ref{sec:pros}, we provide the mains definitions and properties of
Hopf bialgebras, operads, and pros used in the next chapters.
\medbreak

\section{Polynomial spaces} \label{sec:polynomial_spaces}
We introduce here the notion of polynomial spaces and series spaces. All
the algebraic structures considered in this dissertation are polynomial
or series spaces endowed with some operations or co-operations. A set of
operations, analogous to the operations on graded collections of
Section~\ref{subsubsec:products_collections} of
Chapter~\ref{chap:combinatorics}, over graded polynomial spaces are
considered. We also review some links between changes of bases of
polynomial spaces, posets, and incidence algebras.
\medbreak

\subsection{Series and polynomials on collections}
Intuitively, a series (resp. polynomial) on a collection $C$ is a formal
sum (resp. finite formal sum) of objects of the $C$ with coefficients in
a field $\K$. In what follows, $\K$ can be any field of
characteristic~$0$.
\medbreak

\subsubsection{Rational functions} \label{subsubsec:rational_functions}
In a combinatorial context, it is nevertheless convenient to set $\K$ as
the space $\Q(q_0, q_1, \dots)$ of \Def{rational functions} on the
formal parameters $q_i$, $i \in \N$. Let us recall some classical
notations. For any $i \in \N$,
\begin{subequations}
\begin{equation} \label{equ:q_analogs_integers}
    (n)_{q_i} := 1 + q_i + q_i^2 + \dots + q_i^{n - 1},
    \qquad n \in \N_{\geq 1},
\end{equation}
\begin{equation} \label{equ:q_factorials}
    (n)_{q_i}! :=
    \begin{cases}
        1 & \mbox{if } n = 0, \\
        (n)_{q_i} \left((n - 1)_{q_i}\right)!
            & \mbox{otherwise (} n \in \N_{\geq 1} \mbox{)},
    \end{cases}
\end{equation}
\begin{equation} \label{equ:q_binomials}
    \binom{n_1 + n_2}{n_1, n_2}_{q_i}
    :=
    \frac{(n_1 + n_2)_{{q_i}!}} {(n_1)_{q_i}! \; (n_2)_{q_i}!},
    \qquad n_1, n_2 \in \N.
\end{equation}
\end{subequations}
Elements~\eqref{equ:q_analogs_integers} are known as \Def{$q$-analogs}
of integers. Indeed, the specialization $q_i := 1$ in $(n)_{q_i}$ is
equal to~$n$. Elements~\eqref{equ:q_factorials} are \Def{$q$-factorials}
and~\eqref{equ:q_binomials} are \Def{$q$-binomials}.
\medbreak

\subsubsection{Series and polynomials} \label{subsubsec:C_polynomials}
Let $C$ be an $I$-collection. A \Def{series} on $C$ (or, for short, a
\Def{$C$-series}) is a map $f : C \to \K$. The
\Def{coefficient} $f(x)$ of $x \in C$ in $f$ is denoted by
$\Angle{x, f}$. The \Def{support} of $f$ is the set
\begin{equation} \label{equ:support_polynomial}
    \Support(f) := \{x \in C : \Angle{x, f} \ne 0\},
\end{equation}
where the symbol $0$ of~\eqref{equ:support_polynomial} is the zero of
$\K$. A \Def{polynomial} on $C$ (or, for short, a \Def{$C$-polynomial})
is a $C$-series having a finite support. A $C$-series $f$ is a
\Def{$C$-monomial} if $\Support(f)$ is a singleton. We say that $f$ is
\Def{homogeneous} if there is an $i \in I$ such that
$\Support(f) \subseteq C(i)$. For any subset $X$ of $C$, the
\Def{characteristic series} of $X$ is the $C$-series $\Charac(X)$
defined, for any $x \in C$, by
\begin{equation}
    \Angle{x, \Charac(X)} :=
    \begin{cases}
        1 \in \K & \mbox{if } x \in X, \\
        0 \in \K & \mbox{otherwise}.
    \end{cases}
\end{equation}
Given two $C$-series $f$ and $g$, the \Def{scalar product} of $f$ and
$g$ is the scalar
\begin{equation} \label{equ:scalar_product_C_series}
    \Angle{f, g} := \sum_{x \in C} \Angle{x, f} \Angle{x, g}
\end{equation}
of $\K$. Of course, when $f$ or $g$ are $C$-polynomials, $\Angle{f, g}$
is well-defined. When $f$ and $g$ are both $C$-series, $\Angle{f, g}$
may not. This notation for the scalar product of $C$-series is
consistent with the notation $\Angle{x, f}$ for the coefficient of $x$
in $f$ because by~\eqref{equ:scalar_product_C_series}, the coefficient
$\Angle{x, f}$ and the scalar product $\Angle{\Charac(\{x\}), f}$ are
equal.
\medbreak

When $C$ is a graded collection and $f$ is a $C$-polynomial, the
\Def{degree} $\deg(f)$ of $f$ is undefined if $\Support(f) = \emptyset$
and is otherwise the greatest size of an object appearing
in~$\Support(f)$.
\medbreak

\subsubsection{Polynomial spaces} \label{subsubsec:polynomial_spaces}
The \Def{$C$-polynomial space} is the set $\K \Angle{C}$ of all the
$C$-poly\-nomials. We say that $C$ is the \Def{underlying collection}
of $\K \Angle{C}$. For any property $P$ of collections (see
Section~\ref{sec:collections} of Chapter~\ref{chap:combinatorics}), we
say that $\K \Angle{C}$ \Def{satisfies the property $P$} if $C$
satisfies $P$. This set $\K \Angle{C}$ is endowed with the following two
operations. First, the \Def{addition}
\begin{equation}
    + : \K \Angle{C} \times \K \Angle{C} \to \K \Angle{C}
\end{equation}
is defined, for any $f_1, f_2 \in \K \Angle{C}$ and $x \in C$, by
\begin{equation}
    \Angle{x, f_1 + f_2} := \Angle{x, f_1} + \Angle{x, f_2}.
\end{equation}
Second, the \Def{multiplication by a scalar}
\begin{equation}
    \ExtProd : \K \times \K \Angle{C} \to \K \Angle{C}
\end{equation}
is defined, for any $f \in \K \Angle{C}$, $\lambda \in \K$, and
$x \in C$, by
\begin{equation}
    \Angle{x, \lambda \ExtProd f} = \lambda \Angle{x, f}.
\end{equation}
Endowed with these two operations, $\K \Angle{C}$ is a $\K$-vector
space. Moreover, $\K \Angle{C}$ decomposes as a direct sum
\begin{equation} \label{equ:decomposition_space_polynomials}
    \K \Angle{C} = \bigoplus_{i \in I} \K \Angle{C(i)}.
\end{equation}
We call each $\K \Angle{C(i)}$ the \Def{$i$-homogeneous component} of
$\K \Angle{C}$. In the sequel, we shall also write
$\K \Angle{C}(i)$ for $\K \Angle{C(i)}$.
\medbreak

Besides, by using the linear structure of $\K \Angle{C}$, any
$C$-polynomial $f$ can be expressed as the finite sum of $C$-monomials
\begin{equation}
    f = \sum_{x \in C} \Angle{x, f} \ExtProd \Charac(\{x\}),
\end{equation}
which is denoted, by a slight abuse of notation, by
\begin{equation} \label{equ:sum_notation_C_polynomials}
    f = \sum_{x \in C} \Angle{x, f} x.
\end{equation}
The notation~\eqref{equ:sum_notation_C_polynomials} for $f$ as a linear
combination of objects of $C$ is the \Def{sum notation} of
$C$-polynomials. By using this notation, it appears that the set
$\{\Charac(\{x\}) : x \in C\}$ forms a basis of $\K \Angle{C}$. This
basis is called \Def{fundamental basis} of $\K \Angle{C}$, and, by a
slight but convenient abuse of notations, each basis element
$\Charac(\{x\})$, $x \in C$, is simply denoted by~$x$.
\medbreak

We would like to emphasize the fact a polynomial space $\K \Angle{C}$
is always seen through its explicit basis $C$ (contrarily to working
with a vector space $\Vca$ without explicit basis). In the sequel, we
shall define (co-)operations on $C$ which extend by linearity on
$\K \Angle{C}$. Properties of such (co-)operations (like associativity
or commutativity) can be defined and checked only on~$C$.
\medbreak

Besides, we are sometimes led to consider several bases of
$\K \Angle{C}$ and work with many of them at the same time. In this
case, to distinguish elements expressed on different bases, we denote
them by putting elements of $C$ as indices of a letter naming the
basis. For instance, the elements of the $\BasisB$-basis of
$\K \Angle{C}$ are denoted by $\BasisB_x$, $x \in C$.
\medbreak

Let $\K \Angle{C_1}$ and $\K \Angle{C_2}$ be two polynomial spaces such
that $C_1$ and $C_2$ are both $I$-collections. A linear map
\begin{equation}
    \phi : \K \Angle{C_1} \to \K \Angle{C_2}
\end{equation}
is a \Def{polynomial space morphism} if for all $i \in I$ and all
$x \in C_1(i)$, $\phi(x) \in \K \Angle{C_2(i)}$. Observe that any
combinatorial collection morphism $\psi : C_1 \to C_2$ gives rise to a
polynomial space morphism $\phi : \K \Angle{C_1} \to \K \Angle{C_2}$
obtained by extending $\psi$ linearly. Besides, we say that
$\K \Angle{C_2}$ is a \Def{subspace} of $\K \Angle{C_1}$ if there is an
injective polynomial space morphism from $\K \Angle{C_2}$
to~$\K \Angle{C_1}$.
\medbreak

\subsubsection{Graded combinatorial polynomial spaces}
When $C$ is a graded combinatorial collection, as a particular
case of~\eqref{equ:decomposition_space_polynomials}, $\K \Angle{C}$
decomposes as a direct sum
\begin{equation}
    \K \Angle{C} = \bigoplus_{n \in \N} \K \Angle{C}(n).
\end{equation}
Moreover, since $C$ is combinatorial, each $\K \Angle{C(n)}$,
$n \in \N$, is finite dimensional. For this reason, the
\Def{Hilbert series} of $\K \Angle{C}$, defined by
\begin{equation}
    \HilbSeries_{\K \Angle{C}}(t)
    = \sum_{n \in \N} \dim \K \Angle{C}(n) \; t^n,
\end{equation}
is a well-defined series. We can observe that the Hilbert series
$\HilbSeries_{\K \Angle{C}}(t)$ of $\K \Angle{C}$ and the generating
series $\GenSeries_C(t)$ of $C$ are the same power series.
\medbreak

\subsubsection{Duality} \label{subsubsec:duality_polynomial_spaces}
The \Def{dual} of $\K \Angle{C}$ is the $\K$-vector space
$\K \Angle{C}^\Dual$ defined by
\begin{equation}
    \K \Angle{C}^\Dual :=
    \bigoplus_{i \in I} \K \Angle{C}(i)^\Dual,
\end{equation}
where for any $i \in I$, $\K \Angle{C}(i)^\Dual$ is the dual space of
$\K \Angle{C}(i)$. If $C$ is combinatorial, all $\K \Angle{C}(i)$ are
finite dimensional spaces, so that
$\K \Angle{C}(i)^\Dual \simeq \K \Angle{C}(i)$, and hence,
\begin{equation}
    \K \Angle{C}^\Dual \simeq \K \Angle{C}.
\end{equation}
For this reason, we shall identify $\K \Angle{C}$ and
$\K \Angle{C}^\Dual$ in this work once $C$ is combinatorial. The
\Def{duality bracket} between $\K \Angle{C}$ and $\K \Angle{C}^\Dual$ is
the linear map
\begin{equation} \label{equ:duality_bracket}
    \Angle{-, -} :
    \K \Angle{C} \otimes \K \Angle{C}^\Dual \to \K
\end{equation}
defined, for all $x, x' \in C$, by
\begin{equation}
    \Angle{x, x'} :=
    \begin{cases}
        1 & \mbox{if } x = x', \\
        0 & \mbox{otherwise}.
    \end{cases}
\end{equation}
If $\Vca$ is a $\K$-vector space, $\Vca^{\otimes k}$ denotes the space
of all tensors on $\Vca$ of order $k \in \N$. The duality bracket
extends for any $k \in \N$ on
\begin{math}
    \K \Angle{C}^{\otimes k} \otimes {\K \Angle{C}^\Dual}^{\otimes k}
\end{math}
linearly by
\begin{equation}
    \Angle{x_1 \otimes \dots \otimes x_k,
    x'_1 \otimes \dots \otimes x'_k}
    :=
    \prod_{i \in [k]} \Angle{x_i, x'_i}
\end{equation}
for any $x_i, x'_i \in C$, $i \in [k]$.
\medbreak

\subsubsection{Rewrite systems and quotient spaces}
Any rewrite system $(C, \Rew)$ gives rise to a subspace
$\RelationSpaceRewriteRule_{(C, \Rew)}$ of $\K \Angle{C}$ generated by
all the $C$-polynomials $x' - x$ whenever $x$ and $x'$ are two objects
of $C$ such that $x \Rew x'$. We call
$\RelationSpaceRewriteRule_{(C, \Rew)}$ the \Def{space induced} by
$(C, \Rew)$. Conversely, when $\RelationSpace$ is a subspace of
$\K \Angle{C}$ such that there exists a rewrite system $(C, \Rew)$ such
that $\RelationSpace$ and $\RelationSpaceRewriteRule_{(C, \Rew)}$ are
isomorphic, we say that $(C, \Rew)$ is an \Def{orientation} of
$\RelationSpace$. When $(C, \Rew)$ is convergent, one has a concrete
description of the quotient space
$\K \Angle{C}/_{\RelationSpaceRewriteRule_{(C, \Rew)}}$ provided by the
following result.
\medbreak

\begin{Proposition}
\label{prop:space_induced_space_quotient_rewrite_rule}
    Let $(C, \Rew)$ be a convergent rewrite system. Then, as spaces
    \begin{equation}
        \K \Angle{C} /_{\RelationSpaceRewriteRule_{(C, \Rew)}}
        \simeq
        \K \Angle{\NormalForms_{(C, \Rew)}}.
    \end{equation}
\end{Proposition}
\begin{proof}
    First, observe that since $(C, \Rew)$ is convergent, $(C, \Rew)$ is
    terminating and admits a set $\NormalForms_{(C, \Rew)}$ of normal
    forms. Moreover, since $(C, \Rew)$ is confluent, for any object $x$
    of $C$, there is a unique normal form $\eta(x)$ such that
    $x \RewRT \eta(x)$. Let $\Vca$ be the subspace of $\K \Angle{C}$
    generated by all the $C$-polynomials $\eta(x) - x$ such that
    $x \in C$, and let us show that
    $\RelationSpaceRewriteRule_{(C, \Rew)}$ is isomorphic to $\Vca$. Let
    $y - x$ be an element of $\RelationSpaceRewriteRule_{(C, \Rew)}$
    such that $x, y \in C$ and $x \Rew y$. Then, $\eta(y) - y$ and
    $\eta(x) - x$ are elements of $\Vca$. Now, since $x \Rew y$ and
    $(C, \Rew)$ is confluent, the two normal forms $\eta(x)$ and
    $\eta(y)$ are equal. This implies that
    $\eta(x) - x - (\eta(y) - y) = y - x$, showing that
    $y - x \in \Vca$. By linearity, this implies that
    $\RelationSpaceRewriteRule_{(C, \Rew)}$ is a subspace of $\Vca$.
    Conversely, let $\eta(x) - x$ be an nonzero element of $\Vca$ where
    $x \in C$. By definition, there is a chain
    \begin{math}
        x \Rew y_1 \Rew \cdots \Rew y_k \Rew \eta(x)
    \end{math}
    with $k \in \N$ and $y_i \in C$, $i \in [k]$. Hence, the
    $C$-polynomials $\eta(x) - y_k$, $y_k - y_{k - 1}$, \dots,
    $y_2 - y_1$, $y_1 - x$ belong to
    $\RelationSpaceRewriteRule_{(C, \Rew)}$. By summing all these
    $C$-polynomials, we obtain that $\eta(x) - x$ belongs to
    $\RelationSpaceRewriteRule_{(C, \Rew)}$. By linearity, this implies
    that $\Vca$ is a subspace of
    $\RelationSpaceRewriteRule_{(C, \Rew)}$. Now, let $B$ be the family
    of all the $C$-polynomials of the form $\eta(x) - x$ where $x$ is an
    object of $C$ which is not a normal form for $(C, \Rew)$ (otherwise,
    the polynomial would be zero). Since all the elements of $B$ are
    linearly independent, $B$ is a basis of $\Vca$. Hence, the linear
    map $\phi : \Vca \to \Vca$ satisfying $\phi(\eta(x) - x) := x$ for
    all $x \in C \setminus \NormalForms_{(C, \Rew)}$ is an isomorphism.
    This shows that $\RelationSpaceRewriteRule_{(C, \Rew)}$ is
    isomorphic to $\K \Angle{C \setminus \NormalForms_{(C, \Rew)}}$. The
    statement of the proposition follows.
\end{proof}
\medbreak

\subsubsection{Series spaces} \label{subsubsec:C_series}
The \Def{$C$-series space} is the set $\K \AAngle{C}$ of all $C$-series.
Most of the definitions concerning $C$-polynomials and $C$-polynomial
spaces developed in Section~\ref{subsubsec:polynomial_spaces} remain
valid in the context of $C$-series. For instance, $\K \AAngle{C}$ is a
vector space (for the similar operations of addition and multiplication
by a scalar as the ones of $C$-polynomials) and each element of
$\K \AAngle{C}$ can be expressed by a sum
notation~\eqref{equ:sum_notation_C_polynomials}, which is possibly
infinite. One of the main differences of features between $\K \Angle{C}$
and $\K \AAngle{C}$ is that the first one admits $C$ as a basis, while
the second does not.
\medbreak

Observe that a generating series of a combinatorial graded collection is
an element of $\K \AAngle{\Multiset(\{t\})}$, where $\{t\}$ is the
graded collection wherein $t$ is an atomic object and $\Multiset$ is
multiset operation over graded collections (see
Sections~\ref{subsubsec:graded_collections}
and~\ref{subsubsec:operations_graded_comb_collections} of
Chapter~\ref{chap:combinatorics}). Since the introduction of formal
power series, a lot of generalizations were proposed in order to extend
the range of enumerative problems they can help to solve.
\medbreak

The most obvious ones are multivariate series allowing to count objects
not only with respect to their sizes but also with respect to various
other statistics. This encompasses the case of the generating series of
combinatorial multigraded collections (see
Section~\ref{subsubsec:multigraded_collections} of
Chapter~\ref{chap:combinatorics}). Such series are elements of
$\K \AAngle{\Multiset(\{t_1, \dots, t_k\})}$ where all the $t_i$,
$i \in [k]$, are atomic objects. Another one consists in considering
noncommutative series on words~\cite{Eil74,SS78,BR10} (and hence,
elements of $\K \AAngle{A^*}$, where $A$ is an alphabet), or even,
pushing the generalization one step further, on elements of a
monoid~\cite{Sak09} (and hence, elements of $\K \AAngle{\Mca}$ where
$\Mca$ is a monoid). Besides, as another generalization, series on trees
have been considered~\cite{BR82,Boz01}. Series on operads (see
Section~\ref{subsec:operads} about these algebraic structures) increase
the list of these generalizations. Chapoton is the first to have
considered such series on operads~\cite{Cha02,Cha08,Cha09}. Several
authors have contributed to this field by considering slight variations
in the definitions of these series. Among these, one can cite van der
Laan~\cite{Vdl04}, Frabetti~\cite{Fra08}, and Loday and
Nikolov~\cite{LN13}.
\medbreak

\subsection{Operations over graded polynomial spaces}
In the same way as operations over graded collections allow to create
new graded collections from already existing ones, there exist
operations over graded polynomial spaces. Some of these are consequences
of the definitions of operations over graded collections. We present
here the main ones. In what follows, $\K \Angle{C}$, $\K \Angle{C_1}$
and $\K \Angle{C_2}$ are three graded polynomial spaces.
\medbreak

\subsubsection{Direct sum}
The sum of two graded collections translates as the direct sum of the
associated graded polynomial spaces. Indeed,
\begin{equation} \label{equ:direct_sum_polynomial_spaces}
    \K \Angle{C_1 + C_2} \simeq \K \Angle{C_1} \oplus \K \Angle{C_2}.
\end{equation}
An isomorphism between the two spaces
of~\eqref{equ:direct_sum_polynomial_spaces} is provided by the map
\begin{math}
    \phi : \K \Angle{C_1 + C_2}
    \to \K \Angle{C_1} \oplus \K \Angle{C_2}
\end{math},
linearly defined for any $x \in C_1 + C_2$ by
\begin{equation}
    \phi(x) :=
    \begin{cases}
        x \in \K \Angle{C_1} & \mbox{if } x \in C_1, \\
        x \in \K \Angle{C_2} & \mbox{otherwise (} x \in C_2 \mbox{)}.
    \end{cases}
\end{equation}
\medbreak

\subsubsection{Tensor product}
The product of two graded collections translates as the tensor product
of the associated graded polynomial spaces. Indeed
\begin{equation} \label{equ:tensor_product_polynomial_spaces}
    \K \Angle{C_1 \times C_2}
    \simeq \K \Angle{C_1} \otimes \K \Angle{C_2}.
\end{equation}
An isomorphism between the two spaces
of~\eqref{equ:tensor_product_polynomial_spaces} is provided by the map
\begin{math}
    \phi : \K \Angle{C_1 \times C_2}
    \to \K \Angle{C_1} \otimes \K \Angle{C_2}
\end{math},
linearly defined for any $(x_1, x_2) \in C_1 \times C_2$ by
\begin{equation}
    \phi((x_1, x_2)) := x_1 \otimes x_2.
\end{equation}
\medbreak

\subsubsection{Tensor algebras} \label{subsubsec:tensor_algebras}
If $\Vca$ is a $\K$-vector space, the \Def{tensor algebra} of $\Vca$ is
the space $\Tensor \Vca$ defined by
\begin{equation}
    \Tensor \Vca :=
    \bigoplus_{k \in \N}
    \Vca^{\otimes k}.
\end{equation}
A basis of $\Tensor \Vca$ is formed by all tensors on any basis of
$\Vca$. The list operation applied to a graded collection translates as
the tensor algebra of the associated graded polynomial space. Indeed,
for any~$k \geq 0$,
\begin{equation}
    \K \Angle{\List_k(C)} \simeq \K \Angle{C}^{\otimes k}
\end{equation}
and
\begin{equation} \label{equ:tensor_algebra_polynomial_spaces}
    \K \Angle{\List(C)} \simeq \Tensor \K \Angle{C}.
\end{equation}
An isomorphism between the two spaces
of~\eqref{equ:tensor_algebra_polynomial_spaces} is provided by the map
\begin{math}
    \phi : \K \Angle{\List(C)}
    \to \Tensor \K \Angle{C}
\end{math},
linearly defined for any $(x_1, \dots, x_k) \in \List(C)$ by
\begin{equation}
    \phi((x_1, \dots, x_k)) :=
    x_1 \otimes \dots \otimes x_k.
\end{equation}
\medbreak

\subsubsection{Symmetric algebras} \label{subsubsec:symmetric_algebras}
If $\Vca$ is a $\K$-vector space, the \Def{symmetric algebra} of $\Vca$
is the space $\Symmetric \Vca$ defined by
\begin{equation}
    \Symmetric \Vca :=
    \Tensor \Vca/_{\Vca_\Symmetric},
\end{equation}
where $\Vca_\Symmetric$ is the subspace of $\Tensor \K \Angle{C}$
consisting in all the tensors
\begin{equation}
    u \otimes x_1 \otimes x_2 \otimes v
    -
    u \otimes x_2 \otimes x_1 \otimes v,
\end{equation}
where $u, v \in \Tensor \Vca$ and $x_1, x_2 \in \Vca$. A basis of
$\Symmetric \Vca$ is formed by all monomials on any basis of $\Vca$.
The multiset operation applied to a graded collection translates as
the symmetric algebra of the associated graded polynomial space. Indeed,
\begin{equation} \label{equ:symmetric_algebra_polynomial_spaces}
    \K \Angle{\Multiset(C)} \simeq \Symmetric \K \Angle{C}.
\end{equation}
An isomorphism between the two spaces
of~\eqref{equ:symmetric_algebra_polynomial_spaces} is provided by the
map
\begin{math}
    \phi : \K \Angle{\Multiset(C)}
    \to \Symmetric \K \Angle{C}
\end{math},
linearly defined for any $\lbag x_1, \dots, x_k \rbag \in \Multiset(C)$
by
\begin{equation}
    \phi(\lbag x_1, \dots, x_k \rbag)
    := y_1^{\alpha_1} \dots y_\ell^{\alpha_\ell},
\end{equation}
where $\ell$ is the number of distinct elements of
$\lbag x_1, \dots, x_k \rbag$ and each $\alpha_i$, $i \in [\ell]$,
denotes the multiplicity of $y_i$ in $\lbag x_1, \dots, x_k \rbag$.
\medbreak

\subsubsection{Exterior algebras} \label{subsubsec:exterior_algebras}
If $\Vca$ is a $\K$-vector space, the \Def{exterior algebra} of $\Vca$
is the space $\Exterior \Vca$ defined by
\begin{equation}
    \Exterior \Vca := \Tensor \Vca/_{\Vca_\Exterior},
\end{equation}
where $\Vca_\Exterior$ is the subspace of $\Tensor \Vca$ consisting in
all the tensors
\begin{equation}
    u \otimes x_1 \otimes x_2 \otimes v
    +
    u \otimes x_2 \otimes x_1 \otimes v,
\end{equation}
where $u, v \in \Tensor \Vca$ and $x_1, x_2 \in \Vca$. A basis of
$\Exterior \Vca$ is formed by all monomials on a basis of $\Vca$ without
repeated variables. The set operation applied to a graded collection
translates as the exterior algebra of the associated graded
polynomial space. Indeed,
\begin{equation} \label{equ:exterior_algebra_polynomial_spaces}
    \K \Angle{\Set(C)} \simeq \Exterior \K \Angle{C}.
\end{equation}
An isomorphism between the two spaces
of~\eqref{equ:exterior_algebra_polynomial_spaces} is provided by the map
\begin{math}
    \phi : \K \Angle{\Set(C)}
    \to \Exterior \K \Angle{C}
\end{math},
linearly defined for any $\{x_1, \dots, x_k\} \in \Set(C)$ by
\begin{equation}
    \phi(\{x_1, \dots, x_k\}) := x_1 \dots x_k.
\end{equation}
\medbreak

\subsection{Changes of basis and posets}
\label{subsec:change_basis_posets}
It is very usual, given a polynomial space $\K \Angle{C}$, to consider
a poset structure on $C$ to define new bases of $\K \Angle{C}$. Indeed,
such new bases are defined by considering sums of elements greater (or
smaller) than other ones. In this context, incidence algebras of posets
and their Möbius functions play an important role. We expose here these
concepts.
\medbreak

\subsubsection{Incidence algebras}
One of the first apparitions of incidence algebras in combinatorics is
due to Rota~\cite{Rot64}. These structures, associated with any locally
finite poset, provide an abstraction of the principle of
inclusion-exclusion~\cite{Sta11} through their Möbius functions. Indeed,
the usual inclusion-exclusion principle comes from the Möbius function
of the cube poset.
\medbreak

Let $(\Qca, \Ord)$ be a locally finite $I$-poset and $\Comparable^\Qca$
be the $I$-collection of all the ordered pairs $(x, y)$ of objects of
$\Qca$ such that $x \Ord y$, called \Def{pairs of comparable objects}.
The index of $(x, y)$ is the index of $x$ in $\Qca$ (or equivalently,
the index of $y$ in $\Qca$). The \Def{incidence algebra} of
$(\Qca, \Ord)$ is the polynomial space $\K \Angle{\Comparable^\Qca}$
endowed with the linear binary product $\Product$ (the notion of
products in polynomial spaces is presented in the following
Section~\ref{sec:bialgebras} but here, only elementary notions about
these are needed) defined, for any objects $(x, y)$ and $(x', y')$ of
$\Comparable^\Qca$ by
\begin{equation}
    (x, y) \Product \left(x', y'\right) :=
    \begin{cases}
        \left(x, y'\right) & \mbox{if } y = x', \\
        0 & \mbox{otherwise}.
    \end{cases}
\end{equation}
This product is obviously associative. Moreover, on each $i$-homogeneous
component of $\K \Angle{\Comparable^\Qca}$,
$i \in I$, the $\Comparable^\Qca$-polynomial
\begin{equation}
    \Unit_i := \sum_{x \in C(i)} (x, x)
\end{equation}
plays the role of a unit, that is,
\begin{math}
    f \Product \Unit_i = f = \Unit_i \Product f
\end{math}
for all $f \in \K \Angle{\Comparable^C}(i)$. Let for any $i \in I$ the
$\Comparable^\Qca$-polynomial $\zeta_i$, called
\Def{$i$-zeta polynomial} of $(\Qca, \Ord)$, defined by
\begin{equation}
    \zeta_i :=
    \sum_{\substack{
        x, y \in C(i) \\
        x \Ord y
    }}
    (x, y).
\end{equation}
This $\Comparable^\Qca$-polynomial encodes some properties of the order
$\Ord$. For instance, the coefficient in $\zeta_i \Product \zeta_i$ of
each $(x, y) \in \Comparable^\Qca(i)$ is the cardinality of the interval
$[x, y]$ in $(\Qca, \Ord)$. The \Def{$i$-Möbius polynomial} of
$(\Qca, \Ord)$ is the $\Comparable^\Qca$-polynomial $\mu_i$ satisfying
\begin{equation}
    \mu_i \Product \zeta_i
    = \Unit_i
    = \zeta_i \Product \mu_i.
\end{equation}
In other words, $\mu_i$ is the inverse of $\zeta_i$ with respect to the
product~$\Product$. Recall that, as exposed in
Section~\ref{subsubsec:C_polynomials}, polynomials on collections are
functions associating a coefficient with any object. For this reason,
$\zeta_i$ and $\mu_i$ are functions associating a coefficient with any
pair of comparable objects of $\Qca$. We have presented the elements of
incidence algebras as polynomials of pairs of comparable elements, but
in the literature~\cite{Sta11}, it is most common to see these elements
as maps associating a coefficient with each pair of comparable
elements. These two points of view are therefore equivalent but the
definition of the product of incidence algebras in terms of polynomials
is simpler.
\medbreak

\begin{Theorem} \label{thm:mobius_polynomial_coefficients}
    Let $(\Qca, \Ord)$ be a locally finite $I$-poset. Then, the
    $i$-Möbius polynomial~$\mu_i$, $i \in I$, of $(\Qca, \Ord)$ is a
    well-defined element of $\K \Angle{\Comparable^\Qca}$ and its
    coefficients satisfy $\Angle{(x, x), \mu_i} = 1$ for all
    $x \in \Qca(i)$, and
    \begin{equation} \label{equ:mobius_polynomial_coefficients}
        \Angle{(x, z), \mu_i}
        =
        -
        \sum_{\substack{
            y \in \Qca(i) \\
            x \Ord y \OrdStrict z
        }}
        \Angle{(x, y), \mu_i}
    \end{equation}
    for all $x, z \in C(i)$ such that $x \ne z$.
\end{Theorem}
\begin{proof}
    Let $f$ be a $\Comparable^\Qca$-polynomial satisfying
    $f \Product \zeta_i = \Unit_i$. By using the definitions of
    $\Product$ and of $\zeta_i$, we obtain
    \begin{equation}
        \sum_{\substack{
            x, y, z \in \Qca(i) \\
            x \Ord y \Ord z
        }}
        \Angle{(x, y), f} (x, z)
        =
        \sum_{x \in \Qca(i)} (x, x).
    \end{equation}
    This leads to the fact that $\Angle{(x, x), f} = 1$ for all
    $x \in \Qca(i)$, and, for all $x, z \in \Qca(i)$ such
    that~$x \ne z$,
    \begin{equation} \label{equ:mobius_polynomial_coefficients_proof_1}
        \sum_{\substack{
            y \in \Qca(i) \\
            x \Ord y \Ord z
        }}
        \Angle{(x, y), f}
        = 0.
    \end{equation}
    Then, \eqref{equ:mobius_polynomial_coefficients_proof_1} rewrites as
    \begin{equation} \label{equ:mobius_polynomial_coefficients_proof_2}
        \Angle{(x, z), f}
        +
        \sum_{\substack{
            y \in \Qca(i) \\
            x \Ord y \OrdStrict z
        }}
        \Angle{(x, y), f}
        = 0.
    \end{equation}
    Moreover, in the same way, one can prove that if $f'$ is a
    $\Comparable^\Qca$-polynomial satisfying
    $\zeta_i \Product f' = \Unit_i$, the coefficients of $f'$ satisfy
    the same relations as the ones of $f$. Recall now that in any
    algebraic structure endowed with a unitary and associative product,
    if an element has an inverse, it is unique. For this reason,
    $\mu_i = f$. Hence,
    \eqref{equ:mobius_polynomial_coefficients_proof_2} implies that the
    coefficients of $\mu_i$
    satisfy~\eqref{equ:mobius_polynomial_coefficients}. Finally, since
    all coefficients appearing in $\mu_i$ are well-defined, $\mu_i$ is a
    well-defined $\Comparable^\Qca$-polynomial.
\end{proof}
\medbreak

Theorem~\ref{thm:mobius_polynomial_coefficients} provides a recursive
way to compute the coefficients of $\mu_i$, $i \in I$, as a consequence
of the finiteness of each interval of $\Qca(i)$.
\medbreak

\subsubsection{Changes of basis}
\label{subsubsec:polynomial_spaces_change_basis}
Let $C$ be a combinatorial $I$-collection and $\Ord$ be a partial order
relation on $C$ such that $(C, \Ord)$ is an $I$-poset. Consider the
family
\begin{equation} \label{equ:family_basis_change_order}
    \left\{\BasisB_x^\Ord, x \in C\right\}
\end{equation}
of elements of $\K \Angle{C}$ defined, from the fundamental basis of
$\K \Angle{C}$, by
\begin{equation} \label{equ:partial_order_basis_definition}
    \BasisB_x^\Ord :=
    \sum_{\substack{
        y \in C \\
        x \Ord y
    }}
    y.
\end{equation}
Observe that since $C$ is combinatorial and $\Ord$ preserves the indices
of the objects of $C$, each $\BasisB^\Ord_x$ is a homogeneous
$C$-polynomial. We call the family~\eqref{equ:family_basis_change_order}
the \Def{$\BasisB^\Ord$-family} of~$\K \Angle{C}$.
\medbreak

\begin{Proposition} \label{prop:partial_order_bases}
    Let $(C, \Ord)$ be a combinatorial $I$-poset. The
    $\BasisB^\Ord$-family forms a basis of $\K \Angle{C}$ and
    \begin{equation} \label{equ:partial_order_bases}
        x =
        \sum_{\substack{
            y \in C \\
            x \Ord y
        }}
        \Angle{(x, y), \mu_i}
        \BasisB_y^\Ord
    \end{equation}
    for all $x \in C(i)$, $i \in I$, where $\mu_i$ is the $i$-Möbius
    polynomial of~$(C, \Ord)$.
\end{Proposition}
\begin{proof}
    Let us compute the right member of~\eqref{equ:partial_order_bases}
    by using~\eqref{equ:partial_order_basis_definition}. Then, for any
    $x \in C(i)$, $i \in I$, by using the relations satisfied by the
    coefficients of $\mu_i$ provided by
    Theorem~\ref{thm:mobius_polynomial_coefficients}, we obtain
    \begin{equation}\begin{split}
        \sum_{\substack{
            y \in C \\
            x \Ord y
        }}
        \Angle{(x, y), \mu_i}
        \BasisB_y^\Ord
        & =
        \sum_{\substack{
            y \in C \\
            x \Ord y
        }}
        \Angle{(x, y), \mu_i}
        \sum_{\substack{
            z \in C \\
            y \Ord z
        }}
        z \\
        & =
        \sum_{\substack{
            z \in C \\
            x \Ord z
        }}
        \left(
        \sum_{\substack{
            y \in C \\
            x \Ord y \Ord z
        }}
        \Angle{(x, y), \mu_i}
        \right)
        z \\
        & =
        \Angle{(x, x), \mu_i} x
        +
        \sum_{\substack{
            z \in C \\
            x < z
        }}
        \left(
        \sum_{\substack{
            y \in C \\
            x \Ord y \Ord z
        }}
        \Angle{(x, y), \mu_i}
        \right)
        z \\
        & =
        x + 0.
    \end{split}\end{equation}
    Therefore, \eqref{equ:partial_order_bases} holds. Finally, since $C$
    is a basis of $\K \Angle{C}$ and, as we have shown, any $x \in C$
    can be expressed as a linear combination of elements of the
    $\BasisB^\Ord$-family, this family is a basis of~$\K \Angle{C}$.
\end{proof}
\medbreak

\section{Bialgebras} \label{sec:bialgebras}
Bialgebras are polynomial spaces endowed with operations. These
operations are very general in the sense that they can have several
inputs and outputs. These structures encompass all the algebraic
structures seen in this work.
\medbreak

\subsection{Biproducts and duality}
Polynomial spaces are rather poor algebraic structures. It is usual in
combinatorics to handle spaces endowed with several products. When the
products are compatible with the sizes of the underlying combinatorial
objects, all this form a graded algebra. This notion is detailed here,
as well as the concepts of coproduct, duality, and coalgebras and
bialgebras.
\medbreak

\subsubsection{Biproducts} \label{subsubsec:biproducts}
Let $C$ be an $I$-collection. A \Def{biproduct} on $\K \Angle{C}$ is a
linear map
\begin{equation}
    \Biproduct :
    \K \Angle{C}^{\otimes p} \to \K \Angle{C}^{\otimes q}
\end{equation}
where $p, q \in \N$. The \Def{arity} (resp. \Def{coarity}) of
$\Biproduct$ is $p$ (resp. $q$). Any biproduct of arity $p$ and coarity
$q$ can be seen as an operation taking $p$ elements of $\K \Angle{C}$
as input and outputting $q$ elements of $\K \Angle{C}$. This biproduct
is depicted by a rectangle labeled by its name, with $p$ incoming edges
(below the rectangle) and $q$ outgoing edges (above the rectangle) as
\begin{equation}
    \begin{tikzpicture}[xscale=.45,yscale=.2,Centering]
        \node(S1)at(0,0){};
        \node(Sq)at(2,0){};
        \node[Operator](N1)at(1,-3.5){\begin{math}\Biproduct\end{math}};
        \node(E1)at(0,-7){};
        \node(Ep)at(2,-7){};
        \draw[Edge](S1)--(N1);
        \draw[Edge](Sq)--(N1);
        \draw[Edge](N1)--(E1);
        \draw[Edge](N1)--(Ep);
        \node[above of=N1,node distance=6mm,font=\footnotesize]
            {\begin{math}\dots\end{math}};
        \node[below of=N1,node distance=6mm,font=\footnotesize]
            {\begin{math}\dots\end{math}};
        \node[above of=N1,node distance=8mm,font=\footnotesize]
            {\begin{math}q\end{math}};
        \node[below of=N1,node distance=8mm,font=\footnotesize]
            {\begin{math}p\end{math}};
    \end{tikzpicture}\,.
\end{equation}
\medbreak

Under a concrete point of view, given any objects $x_1$, \dots, $x_p$ of
$C$,
\begin{equation}
    \Biproduct(x_1 \otimes \dots \otimes x_p)
    =
    \sum_{y_1, \dots, y_q \in C}
    \xi^{(x_1 \otimes \dots \otimes x_p,
    y_1 \otimes \dots \otimes y_q)}_\Biproduct
    y_1 \otimes \dots \otimes y_q,
\end{equation}
where the
\begin{math}
    \xi^{(x_1 \otimes \dots \otimes x_p,
    y_1 \otimes \dots \otimes y_q)}_\Biproduct
\end{math}
are coefficients of $\K$. These coefficients are called
\Def{structure coefficients} of $\Biproduct$ and wholly determine the
behavior of~$\Biproduct$. We say that $\Biproduct$ is \Def{degenerate}
if all its structure coefficients are zero.
\medbreak

The set of all the biproducts of arity $p$ and coarity $q$ on
$\K \Angle{C}$ has a structure of a $\K$-vector space. Indeed, if
$\Biproduct_1$ and $\Biproduct_2$ are two such biproducts, the
\Def{addition} of $\Biproduct_1$ and $\Biproduct_2$ is the biproduct
$\Biproduct_1 + \Biproduct_2$ defined by
\begin{equation}
    \left(\Biproduct_1 + \Biproduct_2\right)
        (x_1 \otimes \dots \otimes x_p)
    :=
    \Biproduct_1(x_1 \otimes \dots \otimes x_p)
    +
    \Biproduct_2(x_1 \otimes \dots \otimes x_p)
\end{equation}
for any objects $x_1$, \dots, $x_p$ of $C$. Moreover, for any
coefficient $\lambda \in \K$, if $\Biproduct$ is such a biproduct, the
\Def{multiplication by a scalar} of $\Biproduct$ by $\lambda$ is the
biproduct $\lambda \Biproduct$ defined by
\begin{equation}
    (\lambda \Biproduct)(x_1 \otimes \dots \otimes x_p)
    :=
    \lambda \Biproduct(x_1 \otimes \dots \otimes x_p)
\end{equation}
for any objects $x_1$, \dots, $x_p$ of $C$.
\medbreak

\subsubsection{Dual biproducts}
Assume that $C$ is combinatorial so that we can identify $\K \Angle{C}$
with its dual $\K \Angle{C}^\Dual$ as mentioned in
Section~\ref{subsubsec:duality_polynomial_spaces}. Given a biproduct
$\Biproduct$ on $\K \Angle{C}$ of arity $p$ and coarity $q$, the
\Def{dual biproduct} of $\Biproduct$ is the biproduct
\begin{equation}
    \Biproduct^\Dual :
    {\K \Angle{C}^\Dual}^{\otimes q} \to
    {\K \Angle{C}^\Dual}^{\otimes p}
\end{equation}
of arity $q$ and coarity $p$, linearly defined, for all
$y_1, \dots, y_k \in C$, by
\begin{equation} \label{equ:dual_biproduct}
    \Biproduct_\Dual(y_1 \otimes \cdots \otimes y_q)
    :=
    \sum_{x_1, \dots, x_p \in C}
    \Angle{\Biproduct(x_1 \otimes \dots \otimes x_p),
        y_1 \otimes \dots \otimes y_q}
    x_1 \otimes \dots \otimes x_p.
\end{equation}
Let us observe that in~\eqref{equ:dual_biproduct}, the coefficient
\begin{math}
    \Angle{\Biproduct(x_1 \otimes \dots \otimes x_p),
        y_1 \otimes \dots \otimes y_q}
\end{math}
is the structure coefficient
\begin{math}
    \xi^{(x_1 \otimes \dots \otimes x_p,
    y_1 \otimes \dots \otimes y_q)}_\Biproduct
\end{math}
of $\Biproduct$. Hence, if one sees the structure coefficients of
$\Biproduct$ as a matrix whose rows are indexed by the
$x_1 \otimes \dots \otimes x_p$ and the columns by the
$y_1 \otimes \dots \otimes y_q$, the structure coefficients of
$\Biproduct^\Dual$ is the transpose of this matrix.
\medbreak

\subsubsection{Products} \label{subsubsec:products}
A \Def{product} is a biproduct of coarity $1$. In this section,
$\Product$ is a product of arity $p \in \N$.
\medbreak

For any $\ell \geq 1$, let $\Tensor_\ell(\Product)$ be the biproduct
\begin{equation}
    \Tensor_\ell(\Product) :
    {\left({\K \Angle{C}}^{\otimes \ell}\right)}^{\otimes p}
    \simeq
    {\K \Angle{C}}^{\otimes \ell p}
    \to {\K \Angle{C}}^{\otimes \ell}
\end{equation}
defined linearly by
\begin{multline} \label{equ:product_on_tensors}
    \Tensor_\ell(\Product)\left(x_{1, 1} \otimes \dots \otimes
        x_{\ell, 1}
    \otimes x_{1, 2} \otimes \dots \otimes x_{\ell, 2} \otimes
    \dots
    \otimes x_{1, p} \otimes \dots \otimes x_{\ell, p}\right) \\
    :=
    \Product\left(x_{1, 1} \otimes \dots \otimes x_{1, p}\right)
    \otimes \Product\left(x_{2, 1} \otimes \dots \otimes x_{2, p}\right)
    \otimes \dots \otimes
    \Product\left(x_{\ell, 1} \otimes \dots \otimes x_{\ell, p}\right),
\end{multline}
for any
\begin{math}
    x_{1, 1}, \dots, x_{1, p}, x_{2, 1}, \dots, x_{2, p},
    \dots, x_{\ell, 1}, \dots, x_{\ell, p} \in C
\end{math}.
Graphically, $\Tensor_\ell(\Product)$ is the biproduct
\begin{equation}
\,.
\end{equation}
This product $\Tensor_\ell(\Product)$ can be seen as the $\ell$th-tensor
power of $\Product$ seen as a linear map. For this reason,
$\Tensor_\ell(\Product)$ is called the \Def{$\ell$th tensor power} of
$\Product$. For instance, when $\Product$ is a binary product, one has
\begin{equation}
    \Tensor_2(\Product) :
    \K \Angle{C}^{\otimes 2 \times 2} \to
    \K \Angle{C}^{\otimes 2}
\end{equation}
and
\begin{equation} \label{equ:tensor_power_binary_product}
    (x_{1, 1} \otimes x_{2, 1}) \; \Tensor_2(\Product) \;
    (x_{1, 2} \otimes x_{2, 2}) =
    (x_{1, 1} \Product x_{1, 2}) \otimes (x_{2, 1s} \Product x_{2, 2})
\end{equation}
for all $x_{1, 1}, x_{2, 1}, x_{1, 2}, x_{2, 2} \in C$.
In~\eqref{equ:tensor_power_binary_product}, since $\Product$ and
$\Tensor_2(\Product)$ are binary products, we denote them in infix way.
We follow this convention in all this text. Graphically,
$\Tensor_2(\Product)$ is the biproduct
\begin{equation}
    \begin{tikzpicture}[xscale=.3,yscale=.15,Centering]
        \node(S1)at(1,0){};
        \node(S2)at(5,0){};
        \node[Operator](N1)at(1,-3.5){\begin{math}\Product\end{math}};
        \node[Operator](N2)at(5,-3.5){\begin{math}\Product\end{math}};
        \node(E11)at(0,-7){};
        \node(E12)at(2,-7){};
        \node(E21)at(4,-7){};
        \node(E22)at(6,-7){};
        \draw[Edge](S1)--(N1);
        \draw[Edge](S2)--(N2);
        \draw[Edge](N1)--(E11);
        \draw[Edge](N1)--(E12);
        \draw[Edge](N2)--(E21);
        \draw[Edge](N2)--(E22);
        \node(I11)at(0,-12){};
        \node(I12)at(2,-12){};
        \node(I21)at(4,-12){};
        \node(I22)at(6,-12){};
        \draw[Edge](I11)--(E11);
        \draw[Edge](I22)--(E22);
        \draw[Edge](I12)edge[out=90,in=-90]node[]{}(E21);
        \draw[Edge](I21)edge[out=90,in=-90]node[]{}(E12);
        \node[below of=I11,node distance=3mm,font=\footnotesize]
            {\begin{math}x_{1, 1}\end{math}};
        \node[below of=I12,node distance=3mm,font=\footnotesize]
            {\begin{math}x_{2, 1}\end{math}};
        \node[below of=I21,node distance=3mm,font=\footnotesize]
            {\begin{math}x_{1, 2}\end{math}};
        \node[below of=I22,node distance=3mm,font=\footnotesize]
            {\begin{math}x_{2, 2}\end{math}};
    \end{tikzpicture}
\end{equation}
\medbreak

Let us now list some properties a product can satisfy.
\medbreak

When $I$ is endowed with an associative binary product $\DPlus$, if for
any $x_1, \dots, x_p \in C$ one has
\begin{equation}
    \Product(x_1 \otimes \dots \otimes x_p)
    \in \K \Angle{C}(\Index(x_1) \DPlus \cdots \DPlus \Index(x_p)),
\end{equation}
we say that $\Product$ is \Def{$\DPlus$-compatible}. In the particular
case where $C$ is a graded collection and $\Product$ is $+$-compatible,
$\Product$ is \Def{graded}. If $\left\{\BasisB_x : x \in C\right\}$ is
a basis of $\K \Angle{C}$ such that, for any objects $x_1$, \dots, $x_p$
of $C$ there is an $x \in C$ satisfying
\begin{equation}
    \Product
    \left(\BasisB_{x_1} \otimes \dots \otimes \BasisB_{x_p}\right)
    = \BasisB_x,
\end{equation}
we say that the $\BasisB$-basis of $\K \Angle{C}$ is a \Def{set-basis}
(or a \Def{multiplicative basis}) with respect to~$\Product$.
\medbreak

Assume now that $\Product$ is of arity $2$. The \Def{associator} of
$\Product$ is the ternary product
\begin{equation}
    (-, -, -)_\Product :
    \K \Angle{C} \otimes \K \Angle{C} \otimes \K \Angle{C}
    \to \K \Angle{C}
\end{equation}
defined linearly for all $x_1, x_2, x_3 \in C$ by
\begin{equation}
    (x_1, x_2, x_3)_\Product := (x_1 \Product x_2) \Product x_3
    - x_1 \Product (x_2 \Product x_3).
\end{equation}
When for  all $x_1, x_2, x_3 \in C$,
\begin{equation}
    (x_1, x_2, x_3)_\Product = 0,
\end{equation}
we say that $\Product$ is \Def{associative}. The \Def{commutator} of
$\Product$ is the binary product
\begin{equation}
    [-, -]_\Product :
    \K \Angle{C} \otimes \K \Angle{C}
    \to \K \Angle{C}
\end{equation}
defined linearly for all $x_1, x_2 \in C$ by
\begin{equation}
    [x_1, x_2]_\Product := x_1 \Product x_2 - x_2 \Product x_1.
\end{equation}
When for all $x_1, x_2 \in C$,
\begin{equation}
    [x_1, x_2]_\Product = 0,
\end{equation}
the product $\Product$ is \Def{commutative}. When there is an element
$\Unit_\Product$ of $\K \Angle{C}$ such that, for all $x \in C$,
\begin{equation}
    x \Product \Unit_\Product
    = x
    = \Unit_\Product \Product x,
\end{equation}
we say that $\Product$ is \Def{unitary} and that $\Unit_\Product$ is the
\Def{unit} of~$\Product$. This element $\Unit_\Product$ of
$\K \Angle{C}$ can be seen as a product of arity $0$, that is
$\Unit_\Product : ({\K \Angle{C}}^{\otimes 0} = \K) \to \K \Angle{C}$ is
the map sending linearly the multiplicative unit of $\K$ to the element
$\Unit_\Product$ of $\K \Angle{C}$. Observe that if $\Unit_\Product$ is
a graded product, $\Unit_\Product$ is necessarily of degree~$0$.
\medbreak

Finally, any product $\Product$ on $C$ of arity $p$ gives rise to a
product $\bar{\Product}$ on $\K \Angle{C}$ defined linearly, for any
objects $x_1$, \dots, $x_p$ of $C$, by
\begin{equation}
    \bar{\Product}(x_1 \otimes \dots \otimes x_p) :=
    \begin{cases}
        \Product(x_1, \dots, x_p)
            & \mbox{if } \Product(x_1, \dots, x_p)
            \mbox{ is well-defined}, \\
        0 & \mbox{otherwise}.
    \end{cases}
\end{equation}
This product $\bar{\Product}$ is the \Def{linearization} of~$\Product$.
\medbreak

\subsubsection{Coproducts}
A \Def{coproduct} is a biproduct of arity $1$. Let $\Coproduct$ be a
product of coarity $q \in \N$. Observe that the dual of a coproduct is a
product and conversely.
\medbreak

All the properties of products defined in
Section~\ref{subsubsec:products} hold for coproducts in the following
way. For any property $P$ on products, we say that $\Coproduct$
\Def{satisfies the property ``co$P$''} if the dual product
$\Product_\Coproduct$ of $\Coproduct$ satisfies $P$. For instance,
$\Coproduct$ is \Def{cograded} if $\Product_\Coproduct$ is graded, and
$\Coproduct$ is \Def{coassociative} if $\Product_\Coproduct$ is
associative.
\medbreak

\subsection{Polynomial bialgebras} \label{subsec:polynomial_bialgebras}
We now consider polynomial spaces endowed with a set of biproducts. The
main definitions and properties of these structures are listed.
\medbreak

\subsubsection{Elementary definitions}
A \Def{polynomial bialgebra} is a pair $(\K \Angle{C}, \Bca)$ were
$\K \Angle{C}$ is a polynomial space endowed with a (possibly infinite)
set $\Bca$ of biproducts. Let $(\K \Angle{C_1}, \Bca_1)$ and
$(\K \Angle{C_2}, \Bca_2)$ be two polynomial bialgebras. These algebras
are \Def{$\mu$-compatible} if there exists a bijective map
$\mu : \Bca_1 \to \Bca_2$ that sends any biproduct of $\Bca_1$ to a
biproduct of $\Bca_2$ of the same arity and coarity. When
$(\K \Angle{C_1}, \Bca_1)$ and $(\K \Angle{C_2}, \Bca_2)$ are
$\mu$-compatible, a \Def{$\mu$-polynomial bialgebra morphism} (or simply
a \Def{polynomial bialgebra morphism} when there is no ambiguity) from
$\K \Angle{C_1}$ to $\K \Angle{C_2}$ is a polynomial space morphism
\begin{math}
    \phi : \K \Angle{C_1} \to \K \Angle{C_2}
\end{math}
such that
\begin{equation} \label{equ:morphism_polynomial_bialgebras}
    \left(\phi^{\otimes q}\right)
    \left(\Biproduct(x_1 \otimes \dots \otimes x_p)\right)
    =
    (\mu(\Biproduct))
        \left(\phi(x_1) \otimes \dots \otimes \phi(x_p)\right)
\end{equation}
for all biproducts $\Biproduct$ or arity $p$ and coarity $q$ of
$\Bca_1$, and $x_1, \dots, x_p \in C_1$, where $\phi^{\otimes q}$ is the
$q$th tensor power $\Tensor_q(\phi)$ of $\phi$. Graphically,
\eqref{equ:morphism_polynomial_bialgebras} reads as
\begin{equation}
\,.
\end{equation}
\medbreak

Besides, when $(\K \Angle{C_1}, \Bca_1)$ and
$(\K \Angle{C_2}, \Bca_2)$ are $\mu$-compatible,
$(\K \Angle{C_2}, \Bca_2)$ is a \Def{sub-bialgebra} of
$(\K \Angle{C_1}, \Bca_1)$ if there is an injective $\mu$-polynomial
bialgebra morphism from $\K \Angle{C_2}$ to $\K \Angle{C_1}$. Let
$(\K \Angle{C}, \Bca)$ be a polynomial bialgebra. For any subset
$\GeneratingSet$ of $\K \Angle{C}$, the \Def{bialgebra generated} by
$\GeneratingSet$ is the smallest sub-bialgebra
$\K \Angle{C}^\GeneratingSet$ of $\K \Angle{C}$ containing
$\GeneratingSet$. When $\K \Angle{C}^\GeneratingSet = \K \Angle{C}$ and
$\GeneratingSet$ is minimal with respect to the inclusion among the
subsets of $\GeneratingSet$ satisfying this property, $\GeneratingSet$
is a \Def{minimal generating set} of $\K \Angle{C}$. A \Def{polynomial
bialgebra ideal} of $\K \Angle{C}$ is a subspace $\Vca$ of
$\K \Angle{C}$ such that
\begin{equation}
    \Biproduct\left(x_1 \otimes \dots \otimes x_{i - 1}
    \otimes f \otimes x_{i + 1} \otimes \dots \otimes x_p\right)
    \in
    \bigoplus_{j \in [q]}
    \K \Angle{C}^{\otimes j - 1}
    \otimes \Vca \otimes
    \K \Angle{C}^{\otimes q - j}
\end{equation}
for all biproducts $\Biproduct$ of $\Bca$ of arity $p$ and coarity $q$,
$i \in [p]$, $f \in \Vca$, and $x_r \in C$ where
$r \in [p] \setminus \{i\}$. Given a polynomial bialgebra ideal $\Vca$
of $\K \Angle{C}$, the \Def{quotient bialgebra} $\K \Angle{C}/_\Vca$ of
$\K \Angle{C}$ by $\Vca$ is defined in the usual way.
\medbreak

When $\Bca$ contains only products (resp. coproducts),
$(\K \Angle{C}, \Bca)$ is a \Def{polynomial algebra} (resp.
\Def{polynomial coalgebra}).
\medbreak

\subsubsection{Combinatorial polynomial bialgebras}
\label{subsubsec:combinatorial_polynomial_bialgebras}
In practice, and even more so in this dissertation, most of the
encountered polynomial bialgebras are of the form $(\K \Angle{C}, \Bca)$
where $C$ is a combinatorial $I$-collection and $\Bca$ contains only
products and coproducts. When all products (resp. coproducts) of $\Bca$
are $\DPlus$-compatible (resp. $\DPlus$-cocompatible) for some
associative binary products $\DPlus$ on $I$, we say that $\K \Angle{C}$
is a \Def{combinatorial bialgebra}. In most practical cases, $C$ is a
graded, a bigraded, or a colored combinatorial collection.
\medbreak

Let us assume that $(\K \Angle{C}, \Bca)$  is a combinatorial bialgebra.
The \Def{dual bialgebra} of $\K \Angle{C}$ is the bialgebra
$(\K \Angle{C}^\Dual, \Bca^\Dual)$ where $\Bca^\Dual$ is the set of the
dual biproducts of the biproducts of~$\Bca$.
\medbreak

It is very common, given a combinatorial bialgebra
$(\K \Angle{C}, \Bca)$, to endow $C$ with a structure of a combinatorial
poset $(C, \Ord)$ in order to construct $\BasisB^\Ord$-families (see
Section~\ref{subsubsec:polynomial_spaces_change_basis}). For instance,
when a biproduct $\Biproduct$ has complicated structure coefficients,
considering an adequate partial order relation $\Ord$ on $C$ such that
the $\BasisB^\Ord$-family is a set-basis with respect to~$\Biproduct$
allows to infer properties of $\Biproduct$ (such as generating sets of
$\K \Angle{C}$, a description of the nontrivial relations satisfied by
these generators, or even freeness properties).
\medbreak

\subsubsection{Set-theoretic algebras} \label{subsubsec:set_algebras}
When $(\K \Angle{C}, \Pca)$ is a polynomial algebra such that its
fundamental basis is a set-basis with respect to all the products of
$\Pca$, each product $\bar{\Product}$ of $\Pca$ is the linearization of
a product $\Product$ on $C$. In this case, it is possible to forget the
linear structure of $\K \Angle{C}$ and work only with $C$ and its set of
products $\Pca' := \{\Product : \bar{\Product} \in \Pca\}$. We say in
this case that $C$ is a \Def{set-theoretic algebra}.
\medbreak

A large part of the concepts presented above about bialgebras work for
the particular case of set-theoretic algebras with some adjustments. For
instance, to define quotients of a set-theoretic algebra $(C, \Pca')$,
we do not work with polynomial algebra ideals but with congruences of
set-theoretic algebras. To be a little more precise, a
\Def{set-theoretic algebra congruence} is a relation $\Congr$ on $C$
which is an equivalence relation satisfying
\begin{equation}
    \Product
        \left(x_1, \dots, x_{i - 1}, x_i, x_{i + 1}, \dots, x_p\right)
    \Congr
    \Product
        \left(x_1, \dots, x_{i - 1}, x_i', x_{i + 1}, \dots, x_p\right)
\end{equation}
for all products $\Product$ of arities $p$, $i \in [p]$,
$x_i, x_i' \in C$, $x_j \in C$, $j \in [p] \setminus \{i\}$,
whenever $x_i \Congr x_i'$.
\medbreak

In the sequel, if ``$N$'' is the name of an algebraic structure, we call
``set-$N$'' the corresponding set-theoretic structure. For instance, a
set-theoretic unitary associative algebra is a monoid. We shall further
encounter in this way set-operads, colored set-operads, and set-pros.
\medbreak

\subsection{Types of polynomial bialgebras}
A \Def{type of polynomial bialgebra} is specified by biproduct symbols
together with their arities and coarities, and the possible relations
between them (like, for instance, associativity, commutativity,
cocommutativity, or distributivity). In this section, we list some of
the very ordinary types of polynomial bialgebras in combinatorics, and
give concrete examples for each of them. Hopf bialgebras are other
very important types of polynomial bialgebras and are presented in
Section~\ref{sec:hopf_bialgebras}.
\medbreak

\subsubsection{Associative algebras}
\label{subsubsec:associative_algebras}
An \Def{associative algebra} is a polynomial space endowed with an
associative binary product. An associative algebra is \Def{unitary} if
its product is unitary. Besides, an associative algebra is
\Def{commutative} if its product is commutative. To perfectly fit to the
definition of types of bialgebras given above, the type of unitary
associative and commutative algebras is made of a product symbol
$\Product$ of arity $2$ and a product symbol $\Unit$ of arity $0$
together with the relations $(f_1, f_2, f_3)_\Product = 0$,
$[f_1, f_2]_\Product = 0$,
\begin{math}
    f \Product \Unit(\lambda) = \lambda f = \Unit(\lambda) \Product f
\end{math},
where $\lambda$ is any coefficient of $\K$, and $f_1$, $f_2$, and $f_3$
are any elements of the space.
\medbreak

\paragraph{Concatenation algebra}
Let $A := \{a_1, \dots, a_\ell\}$ be an alphabet. The
\Def{concatenation product} is the binary product $\Conc$ on
$\K \Angle{A^*}$ defined as the linearization of the concatenation
product on $A^*$. Since $\Conc$ is graded and all $\K \Angle{A^n}$ are
finite dimensional for all $n \in \N$, $(\K \Angle{A^*}, \Conc)$ is a
combinatorial algebra. Moreover, $\Conc$ is associative, noncommutative,
and admits the empty word $\epsilon$ as unit so that
$(\K \Angle{A^*}, \Conc)$ is a unitary noncommutative associative
algebra.
\medbreak

\paragraph{Shuffle algebra}
The \Def{shuffle product} is the binary product $\shuffle$ on
$\K \Angle{A^*}$ linearly and recursively defined by
\begin{subequations}
\begin{equation} \label{equ:shuffle_words_1}
    u \shuffle \epsilon := u =: \epsilon \shuffle u,
\end{equation}
\begin{equation} \label{equ:shuffle_words_2}
    ua \shuffle vb :=
    (u \shuffle vb) \Conc a + (ua \shuffle v) \Conc b
\end{equation}
\end{subequations}
for any $u, v \in A^*$ and $a, b \in A$, where $\Conc$ is the
concatenation product of words. Intuitively, $\shuffle$ consists in
summing in all the ways of interlacing the two operands. For instance,
\begin{equation}\begin{split}
    \textcolor{Col1}{a_1 a_2}
    \shuffle \textcolor{Col4}{a_2 a_1 a_1} & =
    \textcolor{Col1}{a_1} \textcolor{Col1}{a_2}
    \textcolor{Col4}{a_2} \textcolor{Col4}{a_1}
    \textcolor{Col4}{a_1}
    +
    \textcolor{Col1}{a_1} \textcolor{Col4}{a_2}
    \textcolor{Col1}{a_2} \textcolor{Col4}{a_1}
    \textcolor{Col4}{a_1}
    +
    \textcolor{Col1}{a_1} \textcolor{Col4}{a_2}
    \textcolor{Col4}{a_1} \textcolor{Col1}{a_2}
    \textcolor{Col4}{a_1}
    +
    \textcolor{Col1}{a_1} \textcolor{Col4}{a_2}
    \textcolor{Col4}{a_1} \textcolor{Col4}{a_1}
    \textcolor{Col1}{a_2} \\
    & \quad +
    \textcolor{Col4}{a_2} \textcolor{Col1}{a_1}
    \textcolor{Col1}{a_2} \textcolor{Col4}{a_1}
    \textcolor{Col4}{a_1}
    +
    \textcolor{Col4}{a_2} \textcolor{Col1}{a_1}
    \textcolor{Col4}{a_1} \textcolor{Col1}{a_2}
    \textcolor{Col4}{a_1}
    +
    \textcolor{Col4}{a_2} \textcolor{Col1}{a_1}
    \textcolor{Col4}{a_1} \textcolor{Col4}{a_1}
    \textcolor{Col1}{a_2}
    +
    \textcolor{Col4}{a_2} \textcolor{Col4}{a_1}
    \textcolor{Col1}{a_1} \textcolor{Col1}{a_2}
    \textcolor{Col4}{a_1} \\
    & \quad +
    \textcolor{Col4}{a_2} \textcolor{Col4}{a_1}
    \textcolor{Col1}{a_1} \textcolor{Col4}{a_1}
    \textcolor{Col1}{a_2}
    +
    \textcolor{Col4}{a_2} \textcolor{Col4}{a_1}
    \textcolor{Col4}{a_1} \textcolor{Col1}{a_1}
    \textcolor{Col1}{a_2} \\
    & =
    2 a_1 a_2 a_2 a_1 a_1 + a_1 a_2 a_1 a_2 a_1 +
    a_1 a_2 a_1 a_1 a_2 + a_2 a_1 a_2 a_1 a_1 \\
    & \quad +
    2 a_2 a_1 a_1 a_2 a_1 + 3 a_2 a_1 a_1 a_1 a_2.
\end{split}\end{equation}
Since $\shuffle$ is graded and all $\K \Angle{A^n}$ are finite
dimensional for all $n \in \N$, $(\K \Angle{A^*}, \shuffle)$ is a
combinatorial algebra. Moreover, $\shuffle$ is associative, commutative,
and admits $\epsilon$ as unit so that $(\K \Angle{A^*}, \shuffle)$ is a
unitary commutative associative algebra.
\medbreak

\subsubsection{Coassociative coalgebras}
\label{subsubsec:coassociative_coalgebras}
A \Def{coassociative coalgebra} is a polynomial space endowed with a
coassociative coproduct. A coassociative coalgebra is \Def{counitary} if
its coproduct is counitary. Besides, a coassociative coalgebra is
\Def{cocommutative} if its coproduct is cocommutative.
\medbreak

\paragraph{Deconcatenation coalgebra}
let $\Coproduct_\Conc$ be the dual coproduct of the concatenation
product~$\Conc$ of $\K \Angle{A^*}$ considered in
Section~\ref{subsubsec:associative_algebras}.
By~\eqref{equ:dual_biproduct}, for all $u \in A^*$,
\begin{equation}
    \Coproduct_\Conc(u) =
    \sum_{v, w \in A^*} \Angle{v \Conc w, u} v \otimes w
    =
    \sum_{\substack{
        v, w \in A^* \\
        v \Conc w = u
    }}
    v \otimes w.
\end{equation}
For instance,
\begin{equation}
    \Coproduct_\Conc(a_1 a_1 a_2) =
    \epsilon \otimes a_1 a_1 a_2 +
    a_1 \otimes a_1 a_2 +
    a_1 a_1 \otimes a_2 +
    a_1 a_1 a_2 \otimes \epsilon.
\end{equation}
This coproduct is known as the \Def{deconcatenation coproduct} and
endows $\K \Angle{A^*}$ with a structure of a counitary coassociative
noncocommutative coalgebra.
\medbreak

\paragraph{Unshuffle coalgebra}
Let $\Coproduct_\shuffle$ be the dual coproduct of the shuffle product
$\shuffle$ of $\K \Angle{A^*}$. Again by~\eqref{equ:dual_biproduct}, for
all $u \in A^*$,
\begin{equation}
    \Coproduct_\shuffle(u) =
    \sum_{v, w \in A^*}
    \Angle{v \shuffle w, u}
    v \otimes w.
\end{equation}
The coefficient $\Angle{v \shuffle w, u}$ counts the number of ways to
decompose $u$ as two disjoint subwords $v$ and $w$, and thus,
\begin{equation} \label{equ:coproduct_dual_shuffle}
    \Coproduct_\shuffle(u) =
    \sum_{\substack{
        P_1, P_2 \subseteq [|u|] \\
        P_1 \sqcup P_2 = [|u|]
    }}
    u_{|P_1} \otimes u_{|P_2}.
\end{equation}
This coproduct can also be expressed by
\begin{equation}
    \Coproduct_\shuffle(a) = \epsilon \otimes a + a \otimes \epsilon
\end{equation}
for any $a \in A$, and
\begin{equation} \label{equ:coproduct_dual_shuffle_iterative}
    \Coproduct_\shuffle(u) =
    \prod_{i \in [|u|]} \Coproduct(u_i)
\end{equation}
for any $u \in A^*$, where the product
of~\eqref{equ:coproduct_dual_shuffle_iterative} denotes the iterated
version of the $2$nd tensor power $\Tensor_2(\Conc)$ of the
concatenation product $\Conc$. This product $\Tensor_2(\Conc)$ is
associative due to the fact that $\Conc$ is associative, and hence, its
iterated version is well-defined. For instance,
\begin{equation}\begin{split}
    \Coproduct_\shuffle(a_1 a_1 a_2) & =
    (\epsilon \otimes a_1 + a_1 \otimes \epsilon) \;\Tensor_2(\Conc)\;
    (\epsilon \otimes a_1 + a_1 \otimes \epsilon) \; \Tensor_2(\Conc) \;
    (\epsilon \otimes a_2 + a_2 \otimes \epsilon) \\
    & =
    \epsilon \otimes a_1 a_1 a_2 +
    a_2 \otimes a_1 a_1 +
    a_1 \otimes a_1 a_2 +
    a_1 a_2 \otimes a_1 \\
    & \quad +
    a_1 \otimes a_1 a_2 +
    a_1 a_2 \otimes a_1 +
    a_1 a_1 \otimes a_2 +
    a_1 a_1 a_2 \otimes \epsilon \\
    & =
    \epsilon \otimes a_1 a_1 a_2 +
    a_2 \otimes a_1 a_1 +
    2 a_1 \otimes a_1 a_2 \\
    & \quad +
    2 a_1 a_2 \otimes a_1 +
    a_1 a_1 \otimes a_2 +
    a_1 a_1 a_2 \otimes \epsilon.
\end{split}\end{equation}
This coproduct is known as the \Def{unshuffling coproduct} and endows
$\K \Angle{A^*}$ with a structure of a counitary coassociative
cocommutative coalgebra.
\medbreak

\subsubsection{Dendriform algebras}
\label{subsubsec:dendriform_algebras}
A \Def{dendriform algebra}~\cite{Lod01} is a polynomial space
$\K \Angle{C}$ endowed with two binary products $\LDendr$ and $\RDendr$
satisfying
\begin{subequations}
    \begin{equation} \label{equ:relation_dendr_1}
        (f_1 \LDendr f_2) \LDendr f_3 =
        f_1 \LDendr (f_2 \LDendr f_3) + f_1 \LDendr (f_2 \RDendr f_3),
    \end{equation}
    \begin{equation} \label{equ:relation_dendr_2}
        (f_1 \RDendr f_2) \LDendr f_3 = f_1 \RDendr (f_2 \LDendr f_3),
    \end{equation}
    \begin{equation} \label{equ:relation_dendr_3}
        (f_1 \LDendr f_2) \RDendr f_3 + (f_1 \RDendr f_2) \RDendr f_3
        = f_1 \RDendr (f_2 \RDendr f_3),
    \end{equation}
\end{subequations}
for all $f_1, f_2, f_3 \in \K \Angle{C}$.
\medbreak

\paragraph{Dendriform algebra structure}
A polynomial algebra $(\K \Angle{C}, \Product)$, where $\Product$ is a
binary product, admits a \Def{dendriform algebra structure} if its
product can be split into two operations
\begin{equation}
    \Product = \LDendr + \RDendr,
\end{equation}
where $\LDendr$ and $\RDendr$ are two non-degenerate binary products
such that $(\K \Angle{C}, \LDendr, \RDendr)$ is a dendriform algebra.
Observe that if $(\K \Angle{C}, \Product)$ admits a dendriform algebra
structure, $\Product$ is associative. The associativity of
$\LDendr + \RDendr$ is a consequence of
Relations~\eqref{equ:relation_dendr_1}, \eqref{equ:relation_dendr_2},
and~\eqref{equ:relation_dendr_3} of dendriform algebras.
\medbreak

\paragraph{Codendriform coalgebra structure}
By dualizing the notion of dendriform algebra structure, one obtains
the notion of \Def{codendriform coalgebras}~\cite{Foi07}. More
precisely, a codendriform coalgebra is a polynomial space $\K \Angle{C}$
endowed with two binary coproducts $\LCoDendr$ and $\RCoDendr$ such that
the dual products $\LDendr$ and $\RDendr$ of respectively $\LCoDendr$
and $\RCoDendr$ endow $\K \Angle{C}$ with a dendriform algebra
structure.
\medbreak

In the same way as above, we say that a polynomial coalgebra
$(\K \Angle{C}, \Coproduct)$, where $\Coproduct$ is a binary coproduct,
admits a \Def{codendriform algebra structure} if its coproduct can be
split into two operations
\begin{equation}
    \Coproduct = \LCoDendr + \RCoDendr,
\end{equation}
where $\LCoDendr$ and $\RCoDendr$ are two non-degenerate binary
coproducts such that $(\K \Angle{C}, \LCoDendr, \RCoDendr)$ is a
codendriform colalgebra.
\medbreak

\paragraph{Bidendriform bialgebra structure}
A bialgebra $(\K \Angle{C}, \Product, \Coproduct)$, where $\Product$ is
a binary product and $\Coproduct$ is a binary coproduct, admits a
\Def{bidendriform bialgebra structure}~\cite{Foi07} if $\K \Angle{C}$
admits both a dendriform algebra $(\K \Angle{C}, \LDendr, \RDendr)$ and
a codendriform coalgebra $(\K \Angle{C}, \LCoDendr, \RCoDendr)$
structure with some extra compatibility relations between the products
$\LDendr$ and $\RDendr$ and the coproducts $\LCoDendr$ and $\RCoDendr$.
\medbreak

One among the main benefits of showing that $\K \Angle{C}$ admits a
bidendriform bialgebra structure relies a rigidity theorem~\cite{Foi07}
implying several properties of $\K \Angle{C}$. For instance, when
$\K \Angle{C}$ is a Hopf bialgebra (see
Section~\ref{sec:hopf_bialgebras}), the fact that $\K \Angle{C}$ admits
a bidendriform bialgebra structure implies its self-duality, its
freeness as an associative algebra, and its freeness as a coassociative
coalgebra.
\medbreak

\paragraph{Remarks and generalizations}
We invite the reader to take a look
at~\cite{LR98,Agu00,Lod02,Foi07,EMP08,EM09,LV12} for a supplementary
review of properties of dendriform algebras.
\medbreak

Besides, in the recent years, a lot of generalizations of dendriform
algebras and their dual notions were introduced, each of them splitting
an associative product in different ways and in more that two pieces.
Tridendriform algebras~\cite{LR04}, quadri-algebras~\cite{AL04},
ennea-algebras~\cite{Ler04}, $m$-dendriform algebras of
Leroux~\cite{Ler07}, $m$-dendriform algebras of Novelli~\cite{Nov14},
and polydendriform algebras (see Chapter~\ref{chap:polydendr}) are
examples of such structures.
\medbreak

\paragraph{Shuffle dendriform algebra}
Consider on $\K \Angle{A^*}$ the binary products $\LDendr$ and $\RDendr$
defined linearly and recursively by
\begin{subequations}
\begin{equation} \label{equ:LRDendr_words_epsilson_1}
    u \LDendr \epsilon := u =: \epsilon \RDendr u,
\end{equation}
\begin{equation} \label{equ:LRDendr_words_epsilson_2}
    w \RDendr \epsilon =: 0 := \epsilon \LDendr w,
\end{equation}
\begin{equation} \label{equ:LDendr_words}
    ua \LDendr v := (u \shuffle v) \Conc a,
\end{equation}
\begin{equation} \label{equ:RDendr_words}
    u \RDendr vb := (u \shuffle v) \Conc b
\end{equation}
\end{subequations}
for any $u, v \in A^*$, $w \in A^+$, and $a, b \in A$, where $\Conc$ is
the concatenation product of words. In other words, $u \LDendr v$ (resp.
$u \RDendr v$) is the sum of all the words $w$ obtained by shuffling $u$
and $v$ such that the last letter of $w$ comes from $u$ (resp. $v$). For
example,
\begin{subequations}
\begin{equation}\begin{split}
    \textcolor{Col1}{a_1 a_2}
    \LDendr \textcolor{Col4}{a_2 a_1 a_1} & =
    \textcolor{Col1}{a_1} \textcolor{Col4}{a_2}
    \textcolor{Col4}{a_1} \textcolor{Col4}{a_1}
    \textcolor{Col1}{a_2}
    +
    \textcolor{Col4}{a_2} \textcolor{Col1}{a_1}
    \textcolor{Col4}{a_1} \textcolor{Col4}{a_1}
    \textcolor{Col1}{a_2}
    +
    \textcolor{Col4}{a_2} \textcolor{Col4}{a_1}
    \textcolor{Col1}{a_1} \textcolor{Col4}{a_1}
    \textcolor{Col1}{a_2}
    +
    \textcolor{Col4}{a_2} \textcolor{Col4}{a_1}
    \textcolor{Col4}{a_1} \textcolor{Col1}{a_1}
    \textcolor{Col1}{a_2} \\
    & =
    a_1 a_2 a_1 a_1 a_2 + 3 a_2 a_1 a_1 a_1 a_2,
\end{split}\end{equation}
\begin{equation}\begin{split}
    \textcolor{Col1}{a_1 a_2}
    \RDendr \textcolor{Col4}{a_2 a_1 a_1} & =
    \textcolor{Col1}{a_1} \textcolor{Col1}{a_2}
    \textcolor{Col4}{a_2} \textcolor{Col4}{a_1}
    \textcolor{Col4}{a_1}
    +
    \textcolor{Col1}{a_1} \textcolor{Col4}{a_2}
    \textcolor{Col1}{a_2} \textcolor{Col4}{a_1}
    \textcolor{Col4}{a_1}
    +
    \textcolor{Col1}{a_1} \textcolor{Col4}{a_2}
    \textcolor{Col4}{a_1} \textcolor{Col1}{a_2}
    \textcolor{Col4}{a_1}
    +
    \textcolor{Col4}{a_2} \textcolor{Col1}{a_1}
    \textcolor{Col1}{a_2} \textcolor{Col4}{a_1}
    \textcolor{Col4}{a_1} \\
    & \quad +
    \textcolor{Col4}{a_2} \textcolor{Col1}{a_1}
    \textcolor{Col4}{a_1} \textcolor{Col1}{a_2}
    \textcolor{Col4}{a_1}
    +
    \textcolor{Col4}{a_2} \textcolor{Col4}{a_1}
    \textcolor{Col1}{a_1} \textcolor{Col1}{a_2}
    \textcolor{Col4}{a_1} \\
    & =
    2 a_1 a_2 a_2 a_1 a_1 + a_1 a_2 a_1 a_2 a_1
    + a_2 a_1 a_2 a_1 a_1 + 2 a_2 a_1 a_1 a_2 a_1.
\end{split}\end{equation}
\end{subequations}
These two products endow $\K \Angle{A^*}$ with a structure of a
dendriform algebra. Moreover, the products $\LDendr$ and $\RDendr$
divide the shuffle product into two parts in the sense that
\begin{equation} \label{equ:splitting_associative_product}
    u \shuffle v = u \LDendr v + u \RDendr v
\end{equation}
for all $u, v \in A^*$. This shows that $(\K \Angle{A^*}, \shuffle)$
admits a dendriform algebra structure and offers a way to recover the
recursive definition (see~\eqref{equ:shuffle_words_1}
and~\eqref{equ:shuffle_words_2}) of $\shuffle$. This recursive
description of the shuffle product was known since Ree~\cite{Ree58}.
\medbreak

\paragraph{Max dendriform algebra}
Assume here that $A$ is a totally ordered alphabet by $a_i \leq a_j$
if $i \leq j$. Consider on $\K \Angle{A^+}$ the binary products
$\LDendr$ and $\RDendr$ defined linearly by
\begin{subequations}
\begin{equation}
    u \LDendr v :=
    \begin{cases}
        u \Conc v & \mbox{if } \max_\leq(u) \geq \max_\leq(v) \\
        0 & \mbox{otherwise},
    \end{cases}
\end{equation}
\begin{equation}
    u \RDendr v :=
    \begin{cases}
        u \Conc v & \mbox{if } \max_\leq(u) < \max_\leq(v) \\
        0 & \mbox{otherwise},
    \end{cases}
\end{equation}
\end{subequations}
for all $u, v \in A^+$, where $\Conc$ is the concatenation product of
words. These two products endow $\K \Angle{A^+}$ with a structure of a
dendriform algebra. Moreover, we have here $\Conc = \LDendr + \RDendr$
where $\Conc$ is the associative algebra product of concatenation
of~$\K \Angle{A^+}$.
\medbreak

\subsubsection{Pre-Lie algebras} \label{subsubsec:pre_lie_algebras}
A \Def{pre-Lie algebra} is a polynomial space $\K \Angle{C}$ endowed
with a binary product $\PreLieProduct$ satisfying
\begin{equation} \label{equ:relation_pre_Lie}
    (f_1 \PreLieProduct f_2) \PreLieProduct f_3 -
    f_1 \PreLieProduct (f_2 \PreLieProduct f_3)
    =
    (f_1 \PreLieProduct f_3) \PreLieProduct f_2 -
    f_1 \PreLieProduct (f_3 \PreLieProduct f_2)
\end{equation}
for all $f_1, f_2, f_3 \in \K \Angle{C}$. This
relation~\eqref{equ:relation_pre_Lie} of pre-Lie algebras says that the
associator $(-, -, -)_\PreLieProduct$ is symmetric in its two last
entries.
\medbreak

Pre-Lie algebras were introduced by Vinberg~\cite{Vin63} and
Gerstenhaber~\cite{Ger63} independently. These structures appear under
different names in the literature, for instance as Vinberg algebras,
left-symmetric algebras, or chronological algebras. The appellation
pre-Lie algebra is now very natural since, given a pre-Lie algebra
$(\K \Angle{C}, \PreLieProduct)$, the commutator of $\PreLieProduct$
endows $\K \Angle{C}$ with a structure of a Lie algebra. In the context
of combinatorics, several pre-Lie products are defined on combinatorial
spaces by summing over all the ways to compose (in a certain sense) two
combinatorial objects. For this reason, in an intuitive way, pre-Lie
algebras encode the combinatorics of the composition of combinatorial
objects in all possible ways~\cite{Cha08}. For more details on
pre-Lie algebras, see~\cite{Man11}.
\medbreak

\paragraph{Pre-Lie algebras from associative algebras}
When $(\K \Angle{C}, \Product)$ is a an associative algebra, $\Product$
satisfies in particular~\eqref{equ:relation_pre_Lie} since both left and
right members are equal to zero. For this reason,
$(\K \Angle{C}, \Product)$ is a pre-Lie algebra.
\medbreak

\paragraph{Pre-Lie algebra of rooted trees}
Recall that $\RootedTrees$ is the graded combinatorial collection of all
rooted trees (see Section~\ref{subsubsec:rooted_trees} of
Chapter~\ref{chap:combinatorics}). Consider now on
$\K \Angle{\RootedTrees}$ the products
\begin{math}
    \Graft^{(k)} : \K \Angle{\RootedTrees}^{\otimes k}
    \to \K \Angle{\RootedTrees}
\end{math}
defined linearly for all $k \geq 1$ and all rooted trees $\Tfr_1$,
\dots, $\Tfr_k$ by
\begin{equation}
    \Graft^{(k)}(\Tfr_1 \otimes \dots \otimes \Tfr_k) :=
    (\Node, \lbag \Tfr_1, \dots, \Tfr_k\rbag).
\end{equation}
Intuitively, $\Graft^{(k)}$ consists in grafting all the trees $\Tfr_1$,
\dots, $\Tfr_k$ onto a common root. This product is symmetric with
respect to all its inputs. Now, let $\PreLieProduct$ be the binary
product on $\K \Angle{\RootedTrees}$ defined linearly and recursively by
\begin{equation}
    \Sfr \PreLieProduct \Tfr :=
    \Graft^{(k + 1)}
        \left(\Sfr_1 \otimes \dots \otimes \Sfr_k \otimes \Tfr\right)
    +
    \sum_{i \in [k]}
    \Graft^{(k)}(\Sfr_1 \otimes \dots \otimes \Sfr_{i - 1}
    \otimes \left(\Sfr_i \PreLieProduct \Tfr\right)
    \otimes \Sfr_{i + 1} \otimes \dots \otimes \Sfr_k)
\end{equation}
for any $\Sfr, \Tfr \in \RootedTrees$ where
$\Sfr = \left(\Node, \lbag \Sfr_1, \dots, \Sfr_k \rbag\right)$.
Intuitively, $\PreLieProduct$ consists in summing all the ways of
connecting the root of the second operand on a node of the first. For
example,
\begin{equation}
\,.
\end{equation}
This product endows $\K \Angle{\RootedTrees}$ with a structure of a
pre-Lie algebra.
\medbreak

The free objects in the category of pre-Lie algebras have been described
by Chapoton and Livernet~\cite{CL01}. They have shown that the free
pre-Lie algebra generated by a set $\GeneratingSet$ is the
combinatorial space of all rooted trees whose nodes are labeled on
$\GeneratingSet$, and the product of two such rooted trees is the sum
of all the ways to connect the root of the second tree to a node of
the first. Thereby, the pre-Lie algebra
$(\K \Angle{\RootedTrees}, \PreLieProduct)$ is the free pre-Lie algebra
generated by a singleton.
\medbreak

\subsubsection{About bialgebras}
In the field of algebraic combinatorics, many types of bialgebras have
emerged recently. As previously mentioned, bidendriform
bialgebras~\cite{Foi07} are one of these. In~\cite{Lod08}, Loday defined
the notion of triples of operads, leading to the constructions of
various kinds of bialgebras and analogs of the Poincaré-Birkhoff-Witt
and Cartier-Milnor-Moore theorems (see also~\cite{Cha02}). He defined,
among others, infinitesimal bialgebras, forming an example of
combinatorial bialgebras having one associative binary product and one
coassociative binary coproduct satisfying a compatibility relation.
Moreover, in~\cite{Foi12}, Foissy considered algebraic structures, named
$\Dup$-$\Dendr$ bialgebras, having two binary products satisfying the
duplicial relations~\cite{BF03,Lod08}, two binary coproducts such that
their dual products satisfy the dendriform relations, and such that
these four (co)products satisfy several compatibility relations. These
structures lead to rigidity theorems, in the sense that any
$\Dup$-$\Dendr$ bialgebra is free as a duplicial algebra. In the same
way, Foissy introduced also in~\cite{Foi15} structures named
$\Com$-$\PreLie$ bialgebras, that are spaces with an associative and
commutative binary product, a pre-Lie product, and a binary coproduct
that satisfy compatibility relations.
\medbreak

\section{Hopf bialgebras in combinatorics} \label{sec:hopf_bialgebras}
Hopf bialgebras are polynomial spaces endowed with an associative
product $\Product$ and a coassociative coproduct $\Delta$ satisfying a
kind of commutativity relation very natural in combinatorics. We list
the basic concepts related with these structures and provide some
examples.
\medbreak

\subsection{Hopf bialgebras}
A \Def{Hopf bialgebra} is a polynomial space $\K \Angle{C}$ endowed with
a binary product $\Product$ and a binary coproduct $\Coproduct$ such
that $(\K \Angle{C}, \Product)$ is a unitary associative algebra,
$(\K \Angle{C}, \Coproduct)$ is a counitary coassociative coalgebra,
and, for all $f_1, f_2 \in \K \Angle{C}$,
\begin{equation} \label{equ:hopf_compatibility}
    \Coproduct(f_1 \Product f_2) =
    \Coproduct(f_1) \; \Tensor_2(\Product) \; \Coproduct(f_2).
\end{equation}
The dual bialgebra of a Hopf bialgebra is still a Hopf bialgebra.
\medbreak

Let us now provide some classical definitions about Hopf bialgebras.
\medbreak

\subsubsection{Primitive and group-like elements}
An element $f$ of $\K \Angle{C}$ is \Def{primitive} if
$\Coproduct(f) = 1 \otimes f + f \otimes 1$. The set
$\Pca_{\K \Angle{C}}$ of all primitive elements of $\K \Angle{C}$ forms
a subspace of $\K \Angle{C}$ and the commutator $[-, -]_\Product$ endows
$\Pca_{\K \Angle{C}}$ with a structure of a Lie algebra. Besides, an
element $f$ of $\K \Angle{C}$ is \Def{group-like} if
$\Coproduct(f) = f \otimes f$.
\medbreak

\subsubsection{Convolution product and antipode}
Given two Hopf bialgebras $(\K \Angle{C_1}, \Product_1, \Coproduct_1)$
and $(\K \Angle{C_2}, \Product_2, \Coproduct_2)$, if $\omega$ and
$\omega'$ are two Hopf bialgebra morphisms from $\K \Angle{C_1}$ to
$\K \Angle{C_2}$, the \Def{convolution} of $\omega$ and $\omega'$ is the
map
\begin{equation}
    \omega \Convolution \omega' : \K \Angle{C_1} \to \K \Angle{C_2}
\end{equation}
defined linearly, for any $x \in C_1$, by
\begin{equation}
    (\omega \Convolution \omega')(x)
    :=
    \sum_{y_1, y_2 \in C_1}
    \xi^{(x, y_1 \otimes y_2)}_{\Coproduct_1}
    \omega(y_1) \Product_2 \omega'(y_2),
\end{equation}
where the $\xi^{(-, -)}_{\Coproduct_1}$ are the structure coefficients
of $\Coproduct_1$. This convolution product is associative, as a
consequence of the fact that $\Coproduct_1$ is coassociative and
$\Product_2$ is associative.
\medbreak

Now, let $(\K \Angle{C}, \Product, \Coproduct)$ be a Hopf bialgebra. Let
$\nu : \K \Angle{C} \to \K \Angle{C}$ be the linear map defined as the
inverse of the identity map $\Identity_{\K \Angle{C}}$ on
$\K \Angle{C}$. This map $\nu$ is the \Def{antipode} of $\K \Angle{C}$
and it can be undefined in certain cases.
\medbreak

\subsubsection{Combinatorial Hopf bialgebras}
In algebraic combinatorics, one encounters very particular Hopf
bialgebras. A \Def{combinatorial Hopf bialgebra} is a Hopf bialgebra
$(\K \Angle{C}, \Product, \Coproduct)$ which is graded and combinatorial
(that is $\K \Angle{C}$ is a graded combinatorial space, and $\Product$
and $\Coproduct$ are respectively graded and cograded) and such that $C$
is connected as a graded collection (as a consequence, $\K \Angle{C}(0)$
is of dimension $1$ and can be identified with $\K$).
\medbreak

All combinatorial Hopf bialgebras admit a unique well-defined antipode
$\nu$. Indeed, consider a combinatorial Hopf bialgebra
$(\K \Angle{C}, \Product, \Coproduct)$ and let us denote by $\Unit$ its
unique element of $C(0)$. We can consider, without loss of generality
that $\Unit$ is the unit of $\Product$. The antipode $\nu$ must satisfy
\begin{equation}
    \left(\nu \Convolution \Identity_{\K \Angle{C}}\right)(x) =
    \begin{cases}
        \Unit & \mbox{if } x = \Unit, \\
        0 & \mbox{otherwise}.
    \end{cases}
\end{equation}
Hence, we obtain
\begin{equation}
    \sum_{y_1, y_2 \in C}
    \xi^{(x, y_1 \otimes y_2)}_{\Coproduct} \nu(y_1) \Product y_2
    =
    \begin{cases}
        \Unit & \mbox{if } x = \Unit, \\
        0 & \mbox{otherwise}.
    \end{cases}
\end{equation}
Now, by using the fact that $\Delta$ is counitary and cograded,
we obtain
\begin{equation} \label{equ:antipode_formula_1}
    \nu(\Unit) = \Unit,
\end{equation}
and, for any $x \in \Augmentation(C)$,
\begin{equation} \label{equ:antipode_formula_2}
    \nu(x)
    =
    -\sum_{\substack{
        y_1, y_2 \in C \\
        y_2 \ne \Unit
    }}
    \xi^{(x, y_1 \otimes y_2)}_{\Coproduct} \nu(y_1) \Product y_2.
\end{equation}
Therefore, \eqref{equ:antipode_formula_1}
and~\eqref{equ:antipode_formula_2} imply that the antipode of
$\K \Angle{C}$ is  well-defined can be computed by induction.
\medbreak

\subsection{Main Hopf bialgebras in combinatorics}
Hopf bialgebras are a heavily studied subject. In the last years, many
Hopf bialgebras have been introduced involving a very wide range of
combinatorial spaces. Let us review the main examples.
\medbreak

\subsubsection{Shuffle deconcatenation Hopf bialgebra}
Let $A := \{a_1, \dots, a_\ell\}$ be an alphabet. The concatenation
product $\Conc$ and the unshuffling coproduct $\Coproduct_\shuffle$ (see
Section~\ref{subsubsec:coassociative_coalgebras}) endow $\K \Angle{A^*}$
with a structure of a combinatorial Hopf bialgebra
$\left(\K \Angle{A^*}, \Conc, \Coproduct_\shuffle\right)$. Its dual
bialgebra is the Hopf bialgebra
$\left(\K \Angle{A^*}, \shuffle, \Coproduct_\Conc\right)$ where
$\shuffle$ is the shuffle product and $\Coproduct_\Conc$ is the
deconcatenation coproduct (see again
Section~\ref{subsubsec:coassociative_coalgebras}).
\medbreak

\subsubsection{Noncommutative symmetric functions}
Consider the graded combinatorial polynomial space
$\Sym := \K \Angle{\Compositions}$ of the compositions. Let
\begin{math}
    \left\{\BasisS_{\bm{\lambda}} :
    \bm{\lambda} \in \Compositions\right\}
\end{math}
be the basis of the \Def{complete noncommutative symmetric functions} of
$\Sym$ and $\Product$ be the binary product defined linearly, for any
$\bm{\lambda}, \bm{\mu} \in \Compositions$, by
\begin{equation}
    \BasisS_{\bm{\lambda}} \Product \BasisS_{\bm{\mu}} :=
    \BasisS_{\bm{\lambda} \Conc \bm{\mu}},
\end{equation}
where $\bm{\lambda} \Conc \bm{\mu}$ is the concatenation of the
compositions (seen as words of integers). Moreover, let $\Coproduct$ be
the binary coproduct defined linearly, for any
$\bm{\lambda} \in \Compositions$, by
\begin{equation} \label{equ:coproduct_sym_S_basis}
    \Coproduct\left(\BasisS_{\bm{\lambda}}\right) :=
    \prod_{j \in [\Length(\bm{\lambda})]}
    \left(
    \sum_{\substack{
        n, m \in \N \\
        n + m = \bm{\lambda}_j}
    }
    \BasisS_{(n)} \otimes \BasisS_{(m)}
    \right),
\end{equation}
where the product of~\eqref{equ:coproduct_sym_S_basis} denotes the
iterated version of $2$nd tensor power $\Tensor_2(\Product)$
of~$\Product$, and for any $n \geq 1$, $\BasisS_{(n)}$ is the basis
element indexed by the composition of length $1$ whose only part is $n$,
and $\BasisS_{(0)}$ is identified with the unit $1$ of $\K$. For
instance,
\begin{equation}\begin{split}
    \Coproduct\left(\BasisS_{121}\right) & =
    (1 \otimes \BasisS_1 + \BasisS_1 \otimes 1)
    \; \Tensor_2(\Product) \;
    (1 \otimes \BasisS_2 + \BasisS_1 \otimes \BasisS_1
        + \BasisS_2 \otimes 1)
    \; \Tensor_2(\Product) \;
    (1 \otimes \BasisS_1 + \BasisS_1 \otimes 1) \\
    & =
    1 \otimes \BasisS_{121} +
    \BasisS_1 \otimes \BasisS_{111} +
    \BasisS_1 \otimes \BasisS_{12} +
    \BasisS_1 \otimes \BasisS_{21} +
    2 \BasisS_{11} \otimes \BasisS_{11} +
    \BasisS_{11} \otimes \BasisS_2 \\
    & \quad +
    \BasisS_2 \otimes \BasisS_{11} +
    \BasisS_{111} \otimes \BasisS_1 +
    \BasisS_{12} \otimes \BasisS_1 +
    \BasisS_{21} \otimes \BasisS_1 +
    \BasisS_{121} \otimes 1.
\end{split}\end{equation}
The product $\Product$ and the coproduct $\Coproduct$ endow $\Sym$ with
a structure of a combinatorial Hopf bialgebra.
\medbreak

Moreover, let
\begin{math}
    \left\{\BasisR_{\bm{\lambda}} :
    \bm{\lambda} \in \Compositions\right\}
\end{math}
be the family defined by
\begin{equation}
    \BasisR_{\bm{\lambda}} :=
    \sum_{\substack{
        \bm{\mu} \in \Compositions \\
        \bm{\lambda} \preceq \bm{\mu}}}
    (-1)^{\Length(\bm{\lambda}) - \Length(\bm{\mu})}
    \BasisS_{\bm{\mu}},
\end{equation}
where $\preceq$ is the refinement order of compositions. For instance,
\begin{equation}
    \BasisR_{212} = \BasisS_{212} - \BasisS_{23} - \BasisS_{32}
    + \BasisS_{5}.
\end{equation}
By triangularity, this family forms a basis of $\Sym$ and is known as
the basis of \Def{ribbon noncommutative symmetric functions}. On this
basis, one has, for any $\bm{\lambda}, \bm{\mu} \in \Compositions$,
\begin{equation}
    \BasisR_{\bm{\lambda}} \Product \BasisR_{\bm{\mu}} :=
    \BasisR_{\bm{\lambda} \Product \bm{\mu}} +
    \BasisR_{\bm{\lambda} \ProductTriangleRight \bm{\mu}},
\end{equation}
for any $\bm{\lambda}, \bm{\mu} \in \Compositions$, where
$\bm{\lambda} \Conc \bm{\mu}$ is the concatenation of the compositions
and
\begin{equation}
    \bm{\lambda} \ProductTriangleRight \bm{\mu} :=
    \left(\bm{\lambda}_1, \dots,
    \bm{\lambda}_{\Length(\bm{\lambda}) - 1},
    \bm{\lambda}_{\Length(\bm{\lambda})} + \bm{\mu}_1,
    \bm{\mu}_2, \dots, \bm{\mu}_{\Length(\bm{\mu})}\right).
\end{equation}
For instance,
\begin{equation}
    \BasisR_{\textcolor{Col1}{3112}} \Product
    \BasisR_{\textcolor{Col4}{142}} =
    \BasisR_{\textcolor{Col1}{3112}\textcolor{Col4}{142}} +
    \BasisR_{\textcolor{Col1}{311}
        \textcolor{Col2}{3}\textcolor{Col4}{42}}.
\end{equation}
\medbreak

This Hopf bialgebra $\Sym$ is usually known as the \Def{Hopf bialgebra
of noncommutative symmetric functions}. To explain this name, consider a
totally ordered alphabet $A := \{a_1, a_2, \dots\}$ where
$1 \leq i \leq j$ implies $a_i \Ord a_j$. Now, let the series
\begin{equation}
    \BasisR_{\bm{\lambda}}(A) :=
        \sum_{\substack{
            u \in A^* \\
            \Cmp(u) = \bm{\lambda}
        }}
        u,
\end{equation}
of $\K \AAngle{A^*}$ defined for all $\bm{\lambda} \in \Compositions$,
where $\Cmp$ is defined in Section~\ref{subsubsec:integer_compositions}
of Chapter~\ref{chap:combinatorics}. Observe that all
$\BasisR_{\bm{\lambda}}(A)$ are polynomials when $A$ is finite, but are
series in the other case. For instance,
\begin{subequations}
\begin{equation}
    \BasisR_{\textcolor{Col1}{3}\textcolor{Col4}{1}}(\{a_1, a_2\}) =
    \textcolor{Col1}{a_1 a_1 a_2} \textcolor{Col4}{a_1} +
    \textcolor{Col1}{a_1 a_2 a_2} \textcolor{Col4}{a_1} +
    \textcolor{Col1}{a_2 a_2 a_2} \textcolor{Col4}{a_1},
\end{equation}
\begin{equation}\begin{split} \label{equ:example_realization_ribbon_sym}
    \BasisR_{\textcolor{Col1}{2}\textcolor{Col4}{1}}
        (\{a_1, a_2, a_3\}) & =
    \textcolor{Col1}{a_1 a_2} \textcolor{Col4}{a_1} +
    \textcolor{Col1}{a_1 a_3} \textcolor{Col4}{a_1} +
    \textcolor{Col1}{a_1 a_3} \textcolor{Col4}{a_2} +
    \textcolor{Col1}{a_2 a_2} \textcolor{Col4}{a_1} \\
    & \quad +
    \textcolor{Col1}{a_2 a_3} \textcolor{Col4}{a_1} +
    \textcolor{Col1}{a_2 a_3} \textcolor{Col4}{a_2} +
    \textcolor{Col1}{a_3 a_3} \textcolor{Col4}{a_1} +
    \textcolor{Col1}{a_3 a_3} \textcolor{Col4}{a_2},
\end{split}\end{equation}
 \begin{equation}\begin{split}
    \BasisR_{\textcolor{Col1}{1}\textcolor{Col4}{2}
        \textcolor{Col2}{1}}(\{a_1, a_2, a_3\}) & =
    \textcolor{Col1}{a_2} \textcolor{Col4}{a_1 a_2}
        \textcolor{Col2}{a_1} +
    \textcolor{Col1}{a_2} \textcolor{Col4}{a_1 a_3}
        \textcolor{Col2}{a_1} +
    \textcolor{Col1}{a_2} \textcolor{Col4}{a_1 a_3}
        \textcolor{Col2}{a_2} +
    \textcolor{Col1}{a_3} \textcolor{Col4}{a_1 a_2}
        \textcolor{Col2}{a_1} +
    \textcolor{Col1}{a_3} \textcolor{Col4}{a_1 a_3}
        \textcolor{Col2}{a_1} \\
    & \quad +
    \textcolor{Col1}{a_3} \textcolor{Col4}{a_1 a_3}
        \textcolor{Col2}{a_2} +
    \textcolor{Col1}{a_3} \textcolor{Col4}{a_2 a_2}
        \textcolor{Col2}{a_1} +
    \textcolor{Col1}{a_3} \textcolor{Col4}{a_2 a_3}
        \textcolor{Col2}{a_1} +
    \textcolor{Col1}{a_3} \textcolor{Col4}{a_2 a_3}
        \textcolor{Col2}{a_2}.
\end{split}\end{equation}
\end{subequations}
The linear span of all the $\BasisR_{\bm{\lambda}}(A)$,
$\bm{\lambda} \in \Compositions$, is the space of noncommutative
symmetric functions on $A$. The associative algebra structure of $\Sym$
is compatible with these series in the sense that
\begin{equation} \label{equ:product_sym_series}
    \BasisR_{\bm{\lambda}}(A) \Conc \BasisR_{\bm{\mu}}(A) =
    \left(\BasisR_{\bm{\lambda}} \Product \BasisR_{\bm{\mu}}\right)(A)
\end{equation}
for all $\bm{\lambda}, \bm{\mu} \in \Compositions$, where the product
$\Conc$ of the left member of~\eqref{equ:product_sym_series} is the
usual product of noncommutative series of~$\K \AAngle{A^*}$.
\medbreak

This Hopf bialgebra has been introduced in~\cite{GKLLRT94} as a
generalization of the usual symmetric functions~\cite{Mac15}. This
generalization is a consequence of the fact that there is a surjective
morphism from $\Sym$ to the algebra of symmetric functions.
\medbreak

\subsubsection{Free quasi-symmetric noncommutative symmetric functions}
\label{subsubsec:FQSym}
Consider the gra\-ded combinatorial polynomial
space $\FQSym := \K \Angle{\SymmetricGroup}$ of the permutations. Let
$\{\BasisF_\sigma : \sigma \in \SymmetricGroup\}$ be the basis of the
\Def{fundamental free quasi-symmetric functions} of $\FQSym$ and
$\Product$ be the binary product defined linearly, for any
$\sigma, \nu \in \SymmetricGroup$, by
\begin{equation} \label{equ:product_FQSym}
    \BasisF_\sigma \Product \BasisF_\nu :=
    \sum_{\pi \in \SymmetricGroup}
    \Angle{\pi, \sigma \shuffle \bar{\nu}} \BasisF_\pi,
\end{equation}
where $\bar{\nu}$ is the word obtained by incrementing each letter of
$\nu$ by $|\sigma|$, and $\shuffle$ is the shuffle product of words
defined in Section~\ref{subsubsec:associative_algebras}. For instance
\begin{equation}
    \BasisF_{\textcolor{Col1}{21}}
    \Product
    \BasisF_{\textcolor{Col4}{12}} =
    \BasisF_{\textcolor{Col1}{21}\textcolor{Col4}{34}}
    +
    \BasisF_{\textcolor{Col1}{2}\textcolor{Col4}{3}
    \textcolor{Col1}{1}\textcolor{Col4}{4}}
    +
    \BasisF_{\textcolor{Col1}{2}\textcolor{Col4}{3}
    \textcolor{Col4}{4}\textcolor{Col1}{1}}
    +
    \BasisF_{\textcolor{Col4}{3}\textcolor{Col1}{2}
    \textcolor{Col1}{1}\textcolor{Col4}{4}}
    +
    \BasisF_{\textcolor{Col4}{3}\textcolor{Col1}{2}
    \textcolor{Col4}{4}\textcolor{Col1}{1}}
    +
    \BasisF_{\textcolor{Col4}{34}\textcolor{Col1}{21}}.
\end{equation}
This product is known as the \Def{shifted shuffle product}. Let also
$\Coproduct$ be the binary coproduct defined linearly, for any
$\pi \in \SymmetricGroup$, by
\begin{equation} \label{equ:coproduct_FQSym}
    \Coproduct(\BasisF_\pi) :=
    \sum_{0 \leq i \leq |\pi|}
    \BasisF_{\Std(\pi(1) \dots \pi(i))}
    \otimes
    \BasisF_{\Std\left(\pi(i + 1) \dots \pi(|\pi|)\right)},
\end{equation}
where $\Std$ is defined in Section~\ref{subsubsec:permutations} of
Chapter~\ref{chap:combinatorics}. For instance
\begin{equation}
    \Coproduct(\BasisF_{42513}) =
    1 \otimes \BasisF_{42513} +
    \BasisF_{1} \otimes \BasisF_{2413} +
    \BasisF_{21} \otimes \BasisF_{312} +
    \BasisF_{213} \otimes \BasisF_{12} +
    \BasisF_{3241} \otimes \BasisF_{1} +
    \BasisF_{42513} \otimes 1.
\end{equation}
The product $\Product$ and the coproduct $\Coproduct$ endow $\FQSym$
with a structure of a combinatorial Hopf bialgebra.
\medbreak

This Hopf bialgebra $\FQSym$ is usually known as the
\Def{Hopf bialgebra of free quasi-symmetric functions}. Indeed, as for
$\Sym$, there is a way to see the elements of $\FQSym$ as noncommutative
series. For this, consider a totally ordered alphabet
$A := \{a_1, a_2, \dots\}$ where $1 \leq i \leq j$ implies
$a_i \Ord a_j$. Let the series
\begin{equation}
    \BasisF_\sigma(A) :=
    \sum_{\substack{
        u \in A^* \\
        \Std(u) = \sigma^{-1}
    }}
    u,
\end{equation}
of $\K \AAngle{A^*}$ defined  for all $\sigma \in \SymmetricGroup$. For
instance
\begin{subequations}
\begin{equation} \label{equ:example_realization_fqsym_1}
    \BasisF_{312}(\{a_1, a_2, a_3\}) =
    a_2 a_2 a_1 + a_2 a_3 a_1 +
    a_3 a_3 a_1 + a_3 a_3 a_2,
\end{equation}
\begin{equation} \label{equ:example_realization_fqsym_2}
    \BasisF_{132}(\{a_1, a_2, a_3\}) =
    a_1 a_2 a_1 + a_1 a_3 a_1 +
    a_1 a_3 a_2 + a_2 a_3 a_2.
\end{equation}
\end{subequations}
\medbreak

Furthermore, the Hopf bialgebras $\FQSym$ and $\Sym$ are related through
the injective morphism of Hopf bialgebras $\phi : \Sym \to \FQSym$
defined linearly by
\begin{equation}
    \phi\left(\BasisR_{\bm{\lambda}}\right) :=
    \sum_{\substack{
        \sigma \in \SymmetricGroup \\
        \Des\left(\sigma^{-1}\right) = \Des(\bm{\lambda})
    }}
    \BasisF_\sigma
\end{equation}
for all $\bm{\lambda} \in \Compositions$. For instance,
\begin{equation} \label{equ:example_morphism_from_sym_to_fqsym}
    \phi(\BasisR_{21}) = \BasisF_{312} + \BasisF_{132}.
\end{equation}
Observe, with the help of~\eqref{equ:example_realization_ribbon_sym},
\eqref{equ:example_realization_fqsym_1},
and~\eqref{equ:example_realization_fqsym_2}, in particular
that~\eqref{equ:example_morphism_from_sym_to_fqsym} holds on the
noncommutative series associated with the elements of $\Sym$ and
$\FQSym$, that is,
$\BasisR_{21}(A) = \BasisF_{312}(A) + \BasisF_{132}(A)$.
\medbreak

This Hopf bialgebra has been introduced by Malvenuto and
Reutenauer~\cite{MR95} and is sometimes called the
Malvenuto-Reutenauer algebra. Due to its interpretation~\cite{DHT02} as
an algebra of noncommutative series $\BasisF_\sigma(A)$, it is also
called the algebra of free quasi-symmetric functions. Other classical
examples include the Poirier-Reutenauer Hopf bialgebra of
tableaux~\cite{PR95}, also known as the Hopf bialgebra of free symmetric
functions $\FSym$~\cite{DHT02,HNT05}. This Hopf bialgebra is defined on
the combinatorial space of all standard Young tableaux. The Loday-Ronco
Hopf bialgebra~\cite{LR98}, also known as the Hopf bialgebra of binary
search trees $\PBT$~\cite{HNT05} is defined on the combinatorial space
of all binary trees. As other modern examples of combinatorial spaces
endowed with a Hopf bialgebra structure, one can cite
$\WQSym$~\cite{Hiv99} involving packed words, $\PQSym$~\cite{NT07}
involving parking functions, $\Bell$~\cite{Rey07} involving set
partitions, $\Baxter$~\cite{LR12,Gir12a} involving ordered pairs of twin
binary trees, and $\Camb$~\cite{CP17} involving Cambrian trees. The
study of all these structures uses a large set of tools. Indeed, it
relies on algorithms transforming words into combinatorial objects,
congruences of free monoids, partials orders structures and lattices,
and polytopes and their geometric realizations. Besides, a polynomial
realization of a combinatorial Hopf bialgebra $\K \Angle{C}$ consists in
seeing $\K \Angle{C}$ as an algebra of noncommutative series so that its
product is the usual product of series and its coproduct is obtained by
alphabet doubling (see for instance~\cite{Hiv03}). In this text, only
the polynomial realizations of $\Sym$ and $\FQSym$ have been detailed,
but all the Hopf bialgebras discussed here have polynomial realizations.
\medbreak

\subsubsection{Congruences and Hopf sub-bialgebras of $\FQSym$}
\label{subsubsec:subbialgebras_FQSym_congruences}
It is worth to note that some of the structures discussed above (and
many other ones) can be constructed through congruences of the free
monoid $A^*$ where $A := \{a_1, a_2, \dots\}$. Indeed, if $\Congr$ is a
congruence of $A^*$, one can construct a family
\begin{math}
    \left\{\BasisP_{[\sigma]_\Congr} :
        \sigma \in \SymmetricGroup\right\}
\end{math}
where $[\sigma]_\Congr$ is the $\Congr$-equivalence class of the
permutation $\sigma$ seen as a word on $A$ by identifying each letter
$i$ of $\sigma$ with the letter $a_i$ of $A$, and, for any
$\sigma \in \SymmetricGroup$,
\begin{equation} \label{equ:elements_congruence_FQSym}
    \BasisP_{[\sigma]_\Congr} :=
    \sum_{\sigma \in [\sigma]_\Congr}
    \BasisF_\sigma.
\end{equation}
\medbreak

Of course, the elements~\eqref{equ:elements_congruence_FQSym} do not
form a Hopf sub-bialgebra of $\FQSym$ without precise properties on
$\Congr$. Let us state them. First, we consider that $A$ is totally
ordered by the relation $\Ord$ satisfying $a_i \Ord a_j$ if $i \leq j$.
For any interval $J$ of $A$ and any word $u$ on $A$, we denote by
$u_{|J}$ the subword of $u$ consisting in the letters belonging to $J$.
We say that $\Congr$ is \Def{compatible with the restriction of alphabet
intervals} if, for any interval $J$ of $A$ and any $u, v \in A^*$,
$u \Congr v$ implies $u_{|J} \Congr v_{|J}$. We say that $\Congr$ is
\Def{compatible with the destandardization process} if, for any
$u, v \in A^*$, $u \Congr v$ if and only if $\Std(u) \Congr \Std(v)$ and
$u$ and $v$ have the same commutative image.
\medbreak

\begin{Theorem} \label{thm:congruences_FQSym}
    Let $\Congr$ be a monoid congruence of $A^*$ compatible with the
    restriction of alphabet intervals and with the destandardization
    process. Then, the elements~\eqref{equ:elements_congruence_FQSym}
    form a combinatorial Hopf sub-bialgebra of $\FQSym$ whose bases
    are index by the $\Congr$-equivalence classes of permutations.
\end{Theorem}
\medbreak

This way to construct combinatorial Hopf sub-bialgebras of $\FQSym$ has
been introduced in~\cite{HN07}. One can see also~\cite{Hiv99,Gir11,NT14}
where properties of this construction are studied.
\medbreak

Let us now provide some examples of congruences satisfying the
requirements of Theorem~\ref{thm:congruences_FQSym}.
\medbreak

\paragraph{Sylvester congruence}
The \Def{sylvester congruence}~\cite{HNT05} is the finest monoid
congruence $\Congr$ of $A^*$ satisfying, for any $u \in A^*$ and
$\Asf, \Bsf, \Csf \in A$,
\begin{equation}
    \Asf \Csf u \Bsf \Congr \Csf \Asf u \Bsf,
    \qquad \Asf \Ord \Bsf \OrdStrict \Csf.
\end{equation}
For example, the $\Congr$-equivalence class of the permutation $15423$
(see Figure \ref{fig:sylvester_equiv_class}) is
\begin{equation}
    \{12543, 15243, 15423, 51243, 51423, 54123\}.
\end{equation}
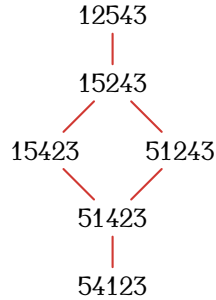
\begin{figure}[ht]
    \begin{tikzpicture}[scale=.45]
        \pgfmathsetmacro\x{2}
        \pgfmathsetmacro\y{2}
        \node(sig1)at(0, 0){\begin{math}12543\end{math}};
        \node(sig2)at(0, -\y){\begin{math}15243\end{math}};
        \node(sig3)at(-\x, -2*\y){\begin{math}15423\end{math}};
        \node(sig4)at(\x, -2*\y){\begin{math}51243\end{math}};
        \node(sig5)at(0, -3*\y){\begin{math}51423\end{math}};
        \node(sig6)at(0, -4*\y){\begin{math}54123\end{math}};
        \draw[Edge]   (sig1)--(sig2);
        \draw[Edge]   (sig2)--(sig3);
        \draw[Edge]   (sig2)--(sig4);
        \draw[Edge]   (sig3)--(sig5);
        \draw[Edge]   (sig4)--(sig5);
        \draw[Edge]   (sig5)--(sig6);
    \end{tikzpicture}
    \caption[The sylvester equivalence class of the permutation
        $15423$.]
    {The sylvester equivalence class of the permutation $15423$.}
    \label{fig:sylvester_equiv_class}
\end{figure}
The set of all $\Congr$-equivalence classes of permutations of size $n$
are in one-to-one correspondence with the set of all binary trees with
$n$ internal nodes. A possible bijection between these two sets is
furnished by the binary search tree insertion algorithm~\cite{Knu98}.
\medbreak

\paragraph{Plactic congruence}
The \Def{plactic congruence}~\cite{LS81,Lot02} is the finest monoid
congruence $\Congr$ of $A^*$ satisfying, for any $\Asf, \Bsf, \Csf \in A$,
\begin{subequations}
\begin{equation}
    \Asf \Csf \Bsf \Congr \Csf \Asf \Bsf,
    \qquad \Asf \Ord \Bsf \OrdStrict \Csf,
\end{equation}
\begin{equation}
    \Bsf \Asf \Csf \Congr \Bsf \Csf \Asf,
    \qquad \Asf \OrdStrict \Bsf \Ord \Csf.
\end{equation}
\end{subequations}
The set of all $\Congr$-equivalence classes of permutations of size $n$
are in one-to-one correspondence with the set of all standard Young
tableaux~\cite{LS81,Ful97}. A possible bijection between these two
sets is furnished by the Robinson-Schensted correspondence~\cite{Sch61}.
\medbreak

\paragraph{Baxter congruence}
The \Def{Baxter congruence}~\cite{Gir12d} (see also~\cite{Rea05,LR12})
is the finest monoid congruence $\Congr$ of $A^*$ satisfying, for any
$u, v \in A^*$ and $\Asf, \Bsf, \Csf, \Dsf \in A$,
\begin{subequations}
\begin{equation}
    \Csf u \Asf\Dsf v \Bsf \Congr \Csf u \Dsf\Asf v \Bsf,
    \qquad \Asf \Ord \Bsf \OrdStrict \Csf \Ord \Dsf,
\end{equation}
\begin{equation}
    \Bsf u \Dsf\Asf v \Csf \Congr \Bsf u \Asf\Dsf v \Csf,
    \qquad \Asf \OrdStrict \Bsf \Ord \Csf \OrdStrict \Dsf.
\end{equation}
\end{subequations}
The set of all $\Congr$-equivalence classes of permutations of size $n$
are in one-to-one correspondence with the set of all ordered pairs of
twin binary trees, objects introduced in~\cite{DG94}. A possible
bijection between these two sets uses the classical binary search tree
insertion algorithm together with a variant of it where the last
inserted node becomes the root of the tree~\cite{Gir12d}.
\medbreak

\paragraph{Bell congruence}
The \Def{Bell congruence}~\cite{Rey07} is the finest monoid congruence
$\Congr$ of $A^*$ satisfying, for any $u \in A^*$ and
$\Asf, \Bsf, \Csf \in A$,
\begin{equation}
    \Asf \Csf u \Bsf \Congr \Csf \Asf u \Bsf,
    \qquad \Asf \Ord \Bsf \OrdStrict \Csf
    \mbox{ and for all letters } \Dsf \mbox{ of } u, \Csf \Ord \Dsf.
\end{equation}
The set of all $\Congr$-equivalence classes of permutations of size $n$
are in one-to-one correspondence with the set of all set partitions of
$[n]$. A possible bijection between these two sets uses a variant of
the patience sorting algorithm~\cite{Rey07}.
\medbreak

\paragraph{Hypoplactic congruence}
The \Def{hypoplactic congruence}~\cite{KT97,KT99} is the finest monoid
congruence $\Congr$ of $A^*$ satisfying, for any $u \in A^*$ and
$\Asf, \Bsf, \Csf \in A$,
\begin{subequations}
\begin{equation}
    \Asf\Csf u \Bsf \Congr \Csf\Asf u \Bsf,
    \qquad \Asf \Ord \Bsf \OrdStrict \Csf,
\end{equation}
\begin{equation}
    \Bsf u \Csf\Asf \Congr \Bsf u \Asf\Csf,
    \qquad \Asf \OrdStrict \Bsf \Ord \Csf.
\end{equation}
\end{subequations}
The set of all $\Congr$-equivalence classes of permutations of size $n$
are in one-to-one correspondence with the set of all compositions of
size $n$.
\medbreak

\paragraph{Total congruence}
The \Def{total congruence} is the monoid congruence $\Congr$ satisfying
$u \Congr v$ if $u$ and $v$ have the same commutative image. There is
exactly one $\Congr$-equivalence class of permutations of size $n$.
\medbreak

\subsubsection{Hopf bialgebra of colored permutations}
\label{subsubsec:FQSym_colored}
Let, for any $\ell \geq 1$, the graded combinatorial polynomial space
$\FQSym^{(\ell)} := \K \Angle{\SymmetricGroup^{(\ell)}}$ of the
$\ell$-colored permutations. Let
\begin{math}
    \{\BasisF_{(\sigma, u)} : (\sigma, u) \in \SymmetricGroup^{(\ell)}\}
\end{math}
be the basis of the \Def{fundamental $\ell$-free quasi-symmetric
functions} of $\FQSym^{(\ell)}$. The space $\FQSym^{(\ell)}$ is endowed
with a binary product $\Product$ similar to the product of $\FQSym$
(see~\eqref{equ:product_FQSym}) wherein the letters of the permutations
and their colors are shuffled. For instance, in $\FQSym^{(5)}$,
\begin{equation}
    \BasisF_{(12, 43)} \Product \BasisF_{(1, 5)}
    =
    \BasisF_{(123, 435)} + \BasisF_{(132, 453)} +
    \BasisF_{(312, 543)}.
\end{equation}
Let also $\Delta$ be the binary coproduct defined in $\FQSym^{(\ell)}$
in a similar way as the coproduct of $\FQSym$
(see~\eqref{equ:coproduct_FQSym}). Again in this case, the colors follow
the letters of the permutations. For instance, in $\FQSym^{(4)}$,
\begin{equation}
    \BasisF_{(312, 411)} =
    1 \otimes \BasisF_{(312, 411)}
    + \BasisF_{(1, 4)} \otimes \BasisF_{(12, 11)}
    + \BasisF_{(21, 41)} \otimes \BasisF_{(1, 1)}
    + \BasisF_{(312, 411)} \otimes 1.
\end{equation}
The product $\Product$ and the coproduct $\Delta$ endow
$\FQSym^{(\ell)}$ with a structure of a combinatorial Hopf bialgebra.
\medbreak

These Hopf bialgebras have been introduced in~\cite{NT10}. Obviously,
they provide a generalization of $\FQSym$ since since
$\FQSym = \FQSym^{(1)}$ and, for any $\ell \geq 1$, $\FQSym^{(\ell)}$ is
a Hopf sub-bialgebra of $\FQSym^{(\ell + 1)}$.
\medbreak

\subsubsection{Hopf bialgebra of uniform block permutations}
\label{subsubsec:UBP}
A \Def{uniform block permutation} (or a \Def{UBP} for short) of size $n$
is a bijection $\pi : \pi^d \to \pi^c$ where $\pi^d$ and $\pi^c$ are set
partitions of $[n]$, and, for any $e \in \pi^d$, $\# e = \# \pi(e)$.
These objects are obvious generalizations of permutations since a
permutation is a UBP where $\pi^d$ and $\pi^c$ are sets of singletons.
For instance, the map $\pi$ defined by
\begin{equation}
    \pi(\{1, 4, 5\}) := \{2, 5, 6\}, \quad
    \pi(\{2\}) := \{1\}, \quad
    \pi(\{3, 6\}) := \{3, 4\}
\end{equation}
is a UBP of size $6$. We denote by $\UBPCollection$ the graded
combinatorial collection of all UBPs. The sequence of integers
associated with $\UBPCollection$ starts by
\begin{equation}
    1, 1, 3, 16, 131, 1496, 22482, 426833,
\end{equation}
and is Sequence~\OEIS{A023998} of~\cite{Slo}.
\medbreak

The graded combinatorial polynomial space
$\UBP := \K \Angle{\UBPCollection}$ admits a combinatorial Hopf
bialgebra structure defined through its basis
$\{\BasisF_\pi : \pi \in \UBPCollection\}$ (see~\cite{AO08}). This Hopf
bialgebra contains $\FQSym$.
\medbreak

\subsubsection{Hopf bialgebra of matrix quasi-symmetric functions}
A \Def{packed matrix} is a matrix with entries in $\N$ such that each
row and each column contains at least one nonzero entry. We denote by
$\Matrices$ the graded combinatorial collection of all packed matrices,
where the size of a packed matrix is the sum of its entries.
\medbreak

The graded combinatorial polynomial space
$\MQSym := \K \Angle{\Matrices}$ admits a combinatorial Hopf bialgebra
structure defined through its basis
$\left\{\BasisM_M : M \in \Matrices\right\}$ of the
\Def{quasi-multiword functions}. Let $\Product$ be the binary product
defined linearly, for any $M_1, M_2 \in \Matrices$ in the following way.
The product $\BasisM_{M_1} \Product \BasisM_{M_2}$ is the sum of all the
$\BasisM_M$ such that the packed matrix $M$ is obtained by horizontally
concatenating $N_1$ and $N_2$ where $N_1$ (resp. $N_2$) is obtained from
$M_1$ (resp. $M_2$) by inserting some null rows, and so that $N_1$ and
$N_2$ have both a same number of rows. For example,
\begin{equation}
    \BasisM_{\Matrix{\textcolor{Col1}{2} & \textcolor{Col1}{1} \\
    \textcolor{Col1}{0} & \textcolor{Col1}{1}}} \Product
    \BasisM_{\Matrix{
        \textcolor{Col4}{1} & \textcolor{Col4}{3}
    }} =
    \BasisM_{\Matrix{
        \textcolor{Col1}{2} & \textcolor{Col1}{1} & 0 & 0 \\
        \textcolor{Col1}{0} & \textcolor{Col1}{1} & 0 & 0 \\
        0 & 0 & \textcolor{Col4}{1} & \textcolor{Col4}{3}
    }} +
    \BasisM_{\Matrix{
        \textcolor{Col1}{2} & \textcolor{Col1}{1} & 0 & 0 \\
        \textcolor{Col1}{0} & \textcolor{Col1}{1} &
            \textcolor{Col4}{1} & \textcolor{Col4}{3} }} +
    \BasisM_{\Matrix{
        \textcolor{Col1}{2} & \textcolor{Col1}{1} & 0 & 0 \\
        0 & 0 & \textcolor{Col4}{1} & \textcolor{Col4}{3} \\
        \textcolor{Col1}{0} & \textcolor{Col1}{1} & 0 & 0
    }} +
    \BasisM_{\Matrix{
        \textcolor{Col1}{2} & \textcolor{Col1}{1} &
            \textcolor{Col4}{1} & \textcolor{Col4}{3} \\
            \textcolor{Col1}{0} & \textcolor{Col1}{1} & 0 & 0
    }} +
    \BasisM_{\Matrix{
        0 & 0 & \textcolor{Col4}{1} & \textcolor{Col4}{3} \\
        \textcolor{Col1}{2} & \textcolor{Col1}{1} & 0 & 0 \\
        \textcolor{Col1}{0} & \textcolor{Col1}{1} & 0 & 0
    }}.
\end{equation}
There is also a binary coproduct $\Delta$ defined on $\MQSym$ that, on
the basis of quasi-multiword functions, splits the packed matrices
horizontally and delete the columns of zeros.
\medbreak

This Hopf bialgebra has been introduced in~\cite{Hiv99} (see
also~\cite{DHT02}).
\medbreak

\section{Operads in combinatorics} \label{sec:operads}
We regard here operads as polynomial algebras and provide the main
definitions used in the following chapters. We also give examples of
some usual operads.
\medbreak

\subsection{Operads} \label{subsec:operads}
Operads have been introduced in the field of algebraic
topology~\cite{May72,BV73}. Here we see operads under a combinatorial
point of view. The notions of operads, nonsymmetric operads,  free
operads, presentations by generators and relations, Koszul duality, and
algebras over operads are reviewed.
\medbreak

\subsubsection{Nonsymmetric operads} \label{subsubsec:ns_operads}
A \Def{nonsymmetric operad} (or a \Def{ns operad} for short) is
a graded augmented polynomial space $\K \Angle{C}$ endowed with a set of
binary linear products $\left\{\circ_i : i \in \N_{\geq 1}\right\}$.
These products have to satisfy several relations. First, when
$x \in C(n)$ and $i \geq n + 1$, for any $y \in C$,
\begin{equation}
    x \circ_i y = 0.
\end{equation}
Moreover, the products $\circ_i$, $i \in \N_{\geq 1}$, satisfy
\begin{equation} \label{equ:graduation_operads}
    \circ_i : \K \Angle{C}(n) \otimes \K \Angle{C}(m) \to
    \K \Angle{C}(n + m - 1),
    \qquad
    n, m \in \N_{\geq 1}, i \in [n].
\end{equation}
This is equivalent to the fact that if $x \in C(n)$, $y \in C(m)$, and
$i \in [n]$,
\begin{math}
    x \circ_i y \in \K \Angle{C}(n + m - 1).
\end{math}
For any $x \in C(n)$, $y \in C(m)$, and $z \in C(k)$, one must have
\begin{subequations}
\begin{equation} \label{equ:operad_axiom_1}
    (x \circ_i y) \circ_{i + j - 1} z = x \circ_i (y \circ_j z),
    \qquad i \in [n], j \in [m],
\end{equation}
\begin{equation} \label{equ:operad_axiom_2}
    (x \circ_i y) \circ_{j + m - 1} z = (x \circ_j z) \circ_i y,
    \qquad i < j \in [n],
\end{equation}
\end{subequations}
Finally, we demand the existence of an element $\Unit$ of
$\K \Angle{C}(1)$ satisfying, for any $x \in C(n)$,
\begin{equation} \label{equ:operad_axiom_3}
    \Unit \circ_1 x = x = x \circ_i \Unit,
    \qquad i \in [n].
\end{equation}
\medbreak

Let us now provide an intuitive meaning of these relations. Any element
$f$ of $\K \Angle{C}(n)$ is seen as a product of arity $n$, depicted as
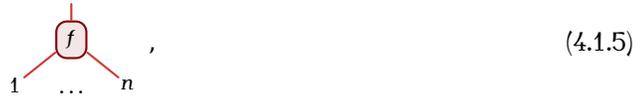
\begin{equation}
\,,
\end{equation}
and where the inputs are indexed from $1$ to $n$ from left to right. The
$\circ_i$ act by composing these products: for any
$f \in \K \Angle{C}(n)$, $g \in \K \Angle{C}(m)$, and $i \in [n]$,
$f \circ_i g$ is the product obtained by plugging the output of $g$ onto
the $i$th input of $f$. This is depicted as
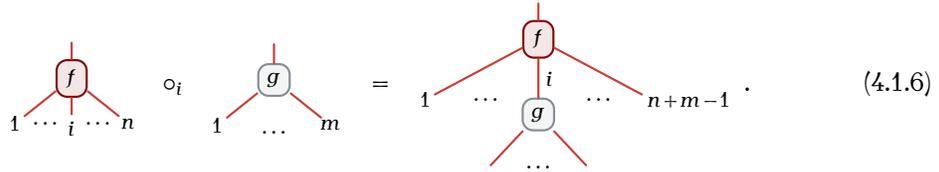
\begin{equation}
    %
\,.
\end{equation}
Under this formalism, Relation~\eqref{equ:operad_axiom_1} says that the
product
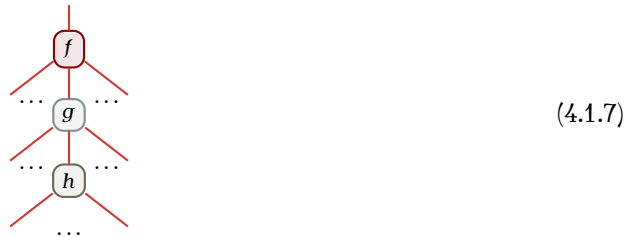
\begin{equation}
    %
\end{equation}
can be formed by two ways: by starting by composing $f$ and $g$, and the
result with $h$, or by starting by composing $g$ and $h$, and the result
with $f$. Relation~\eqref{equ:operad_axiom_2} says that the product
\begin{equation}
    %
\end{equation}
can be formed by two ways: by starting by composing $f$ and $g$, and the
result with $h$, or by starting by composing $f$ and $h$, and the result
with $g$. Finally, Relation~\eqref{equ:operad_axiom_3} says
that~$\Unit$ behaves as an identity product, so that
\begin{equation}
    %
\,.
\end{equation}
\medbreak

Let us fix some vocabulary. Each element $f$ of $\K \Angle{C}(n)$ is
of \Def{arity} $n$. The arity of $f$ is denoted by $|f|$. The maps
$\circ_i$, $i \geq 1$, are \Def{partial composition maps}.
Relation~\eqref{equ:operad_axiom_1} is the \Def{series associativity
relation}, while~\eqref{equ:operad_axiom_2} is the \Def{parallel
associativity relation}. The element $\Unit$ of arity $1$
satisfying~\eqref{equ:operad_axiom_3} is the \Def{unit} of
$\K \Angle{C}$. This element is unique.
\medbreak

Since a ns operad is a particular polynomial algebra, all the
properties and definitions about polynomial algebras exposed in
Section~\ref{subsec:polynomial_bialgebras} remain valid for ns operads
(like ns operad morphisms, ns suboperads, generating  sets, operad
ideals and quotients, {\em etc.}). Observe also,
from~\eqref{equ:graduation_operads} that the $\circ_i$,
$i \in \N_{\geq 1}$, are not graded products. Nevertheless, these
partial composition maps are $\DPlus$-compatible products where $\DPlus$
is the operation satisfying $n \DPlus m := n + m - 1$ for any
$n, m \in \N_{\geq 1}$. As a side remark, since $C$ augmented, one can
see the partial composition maps as graded products on the space
$\K \Angle{\Suspension_{-1}(C)}$. In this way, when $C$ is
combinatorial, one can see $\K \Angle{C}$ as a particular combinatorial
algebra (see
Section~\ref{subsubsec:combinatorial_polynomial_bialgebras}).
\medbreak

\subsubsection{Additional definitions}
Given a ns operad $\K \Angle{C}$, the \Def{complete composition maps}
of $\K \Angle{C}$ are, the linear maps $\circ^{(n)}$,
$n \in \N_{\geq 1}$, satisfying
\begin{equation}
    \circ^{(n)} :
    \K \Angle{C}(n) \otimes
    \K \Angle{C}(m_1) \otimes \dots \otimes \K \Angle{C}(m_n)
    \to \K \Angle{C}(m_1 + \dots + m_n)
\end{equation}
defined linearly, for any $x \in C(n)$ and $y_1, \dots, y_n \in C$, by
\begin{equation}
    \circ^{(n)}(x \otimes y_1 \otimes \dots \otimes y_n)
    := (\dots ((x \circ_n y_n) \circ_{n - 1} y_{n - 1}) \dots)
    \circ_1 y_1.
\end{equation}
To gain concision, we shall denote by $f \circ [g_1, \dots, g_n]$ the
element $\circ^{(n)}(f \otimes g_1 \otimes \dots \otimes g_n)$ of
$\K \Angle{C}$, for all $f \in \K \Angle{C}(n)$ and
$g_1, \dots, g_n \in \K \Angle{C}$. Under this notation, we call $\circ$
the \Def{complete composition map} of $\K \Angle{C}$.
\medbreak

Let us now provide some particular definitions about ns operads that do
not come from the general ones of polynomial algebras of
Section~\ref{subsec:polynomial_bialgebras}.
\medbreak

An element $f$ of arity $2$ of $\K \Angle{C}$ is \Def{associative} if
$f \circ_1 f = f \circ_2 f$. If $\K \Angle{C_1}$ and $\K \Angle{C_2}$
are two ns operads, a \Def{ns operad antimorphism} is a graded
polynomial space morphism $\phi : \K \Angle{C_1} \to \K \Angle{C_2}$
such that $\phi(f \circ_i g) = \phi(f) \circ_{n - i + 1} \phi(g)$ for
any element $f$ of arity $n$ of $\K \Angle{C_1}$, any
$g \in \K \Angle{C_1}$, and $i \in [n]$. A \Def{symmetry} of
$\K \Angle{C}$ is either a ns operad automorphism or a ns operad
antiautomorphism of $\K \Angle{C}$. The set of all symmetries of
$\K \Angle{C}$ forms a group for the map composition, called
\Def{group of symmetries} of $\K \Angle{C}$.
\medbreak

Given two ns operads $\K \Angle{C_1}$ and $\K \Angle{C_2}$, the
\Def{Hadamard product} of $\K \Angle{C_1}$ and $\K \Angle{C_2}$ is the
ns operad denoted by $\K \Angle{C_1} \Hadamard \K \Angle{C_2}$ and
defined on the polynomial space $\K \Angle{C_1 \Hadamard C_1}$ where
$\Hadamard$ is the Hadamard product of graded collections (see
Section~\ref{subsubsec:products_collections} of
Chapter~\ref{chap:combinatorics}). The partial compositions $\circ_i$,
$i \in \N_{\geq 1}$, are defined linearly by
\begin{equation} \label{equ:partial_composition_map_Hadamard}
    (x_1, x_2) \circ_i (y_1, y_2) :=
    (x_1 \circ_i y_1, x_2 \circ_i y_2),
\end{equation}
for any $(x_1, x_2) \in (C_1 \Hadamard C_2)(n)$,
$(y_1, y_2) \in C_1 \Hadamard C_2$, $i \in [n]$, where the second (resp.
third) occurrence of $\circ_i$
in~\eqref{equ:partial_composition_map_Hadamard} is a partial composition
map of $\K \Angle{C_1}$ (resp. $\K \Angle{C_2}$).
\medbreak

\subsubsection{Free ns operads}
Let $\GeneratingSet$ be an augmented graded collection. The
\Def{free ns operad} over $\GeneratingSet$ is the ns operad
\begin{equation}
    \FreeOperad(\GeneratingSet)
    :=
    \K \Angle{\PlanarRootedTrees_\Leaf^\GeneratingSet},
\end{equation}
where $\PlanarRootedTrees_\Leaf^\GeneratingSet$ is the graded collection
of all the $\GeneratingSet$-syntax trees (see
Section~\ref{subsec:syntax_trees} of Chapter~\ref{chap:combinatorics}).
The space $\FreeOperad(\GeneratingSet)$ is endowed with the
linearizations of the grafting operations $\circ_i$,
$i \in \N_{\geq 1}$, defined in
Section~\ref{subsubsec:grafting_syntax_trees} of
Chapter~\ref{chap:combinatorics}. The unit of
$\FreeOperad(\GeneratingSet)$ is the only $\GeneratingSet$-syntax tree
$\Leaf$ of arity~$1$.
\medbreak

Let also
\begin{equation}
    \CorollaNOARG : \GeneratingSet \to \FreeOperad(\GeneratingSet)
\end{equation}
be the \Def{inclusion map}, that is the map sending any
$x \in \GeneratingSet$ to the corolla $\Corolla{x}$ (see
Section~\ref{subsubsec:graded_collection_syntax_trees} of
Chapter~\ref{chap:combinatorics}).
\medbreak

\subsubsection{Evaluations and treelike expressions}
Let now $\K \Angle{C}$ be a ns operad. Since $C$ is a graded augmented
collection, one can consider the free ns operad $\FreeOperad(C)$ of the
$C$-syntax trees. The \Def{evaluation map} of $\K \Angle{C}$ is the map
\begin{equation}
    \Eval : \FreeOperad(C) \to \K \Angle{C}
\end{equation}
defined linearly by induction, for any $C$-syntax tree $\Tfr$, by
\begin{equation}
    \Eval(\Tfr) :=
    \begin{cases}
        \Unit \in \K \Angle{C}
            & \mbox{if } \Tfr = \Leaf, \\
        \omega_\Tfr(\epsilon) \circ
        \left[\Eval\left(\Tfr_1\right), \dots,
            \Eval\left(\Tfr_k\right)\right]
            & \mbox{otherwise},
    \end{cases}
\end{equation}
where $\circ$ is the complete composition map of $\K \Angle{C}$,
$\omega_\Tfr(\epsilon)$ is the label of the root of $\Tfr$, and $k$ is
the root arity of $\Tfr$. This map $\Eval$ is the unique surjective ns
operad morphism from $\FreeOperad(C)$ to $\K \Angle{C}$ satisfying
$\Eval(\Corolla{x}) = x$ for all $x \in C$.
\medbreak

For any subset $S$ of $C$, an \Def{$S$-treelike expression} of an
element $f$ of $\K \Angle{C}$ is an element $g$ of $\FreeOperad(S)$ such
that $\Eval(g) = f$. A treelike expression can be thought as a
factorization in a ns operad.
\medbreak

\subsubsection{Presentations by generators and relations}
A \Def{presentation} of a ns operad $\K \Angle{C}$ consists in a pair
$(\GeneratingSet, \RelationSpace)$ such that $\GeneratingSet$ is an
augmented graded collection, $\RelationSpace$ is a subspace of
$\FreeOperad(\GeneratingSet)$ and
\begin{equation}
    \K \Angle{C} \simeq
    \FreeOperad(\GeneratingSet)/_{\Angle{\RelationSpace}}
\end{equation}
where $\Angle{\RelationSpace}$ is the ns operad ideal of
$\FreeOperad(\GeneratingSet)$ generated by $\RelationSpace$. We call
$\GeneratingSet$ the \Def{set of generators} and $\RelationSpace$ the
\Def{space of relations} of $\K \Angle{C}$.
\medbreak

We say that a presentation $(\GeneratingSet, \RelationSpace)$ of
$\K \Angle{C}$ is \Def{quadratic} if $\RelationSpace$ is a homogeneous
subspace of $\FreeOperad(\GeneratingSet)$ consisting in syntax trees of
degree $2$. Besides, we say that $(\GeneratingSet, \RelationSpace)$  is
\Def{binary} if $\GeneratingSet$ has only elements of size~$2$. By
extension, we say also that $\K \Angle{C}$ is \Def{quadratic} (resp.
\Def{binary}) if it admits a quadratic (resp. binary) presentation.
\medbreak

In practice, to establish presentations of ns operads, we use rewrite
systems on syntax trees (see
Section~\ref{subsec:rewrite_rules_syntax_trees} of
Chapter~\ref{chap:combinatorics}).
\medbreak

\begin{Theorem} \label{thm:presentation_operads}
    Let $\K \Angle{C}$ be a ns operad, $\GeneratingSet$ be a
    subcollection of $C$, and $\RelationSpace$ be a subspace of
    $\FreeOperad(\GeneratingSet)$ of syntax trees of degrees $2$ or
    more. If
    \begin{enumerate}[label={(\it\roman*)}]
        \item \label{item:presentation_operads_1}
        the collection $\GeneratingSet$ is a generating set of
        $\K \Angle{C}$;
        \item \label{item:presentation_operads_2}
        for any $f \in \RelationSpace$, $\Eval(f) = 0$;
        \item \label{item:presentation_operads_3}
        there exists a rewrite system
        $(\PlanarRootedTrees^\GeneratingSet_\Leaf, \Rew)$ being an
        orientation of $\RelationSpace$, such that
        its closure
        $(\PlanarRootedTrees^\GeneratingSet_\Leaf, \RewTrees)$ is
        convergent, and its set of normal forms
        \begin{math}
            \NormalForms_{
                (\PlanarRootedTrees^\GeneratingSet_\Leaf, \RewTrees)}
        \end{math}
        is isomorphic to $C$,
    \end{enumerate}
    then $(\GeneratingSet, \RelationSpace)$ is a presentation of
    $\K \Angle{C}$.
\end{Theorem}
\begin{proof}
    By definition of the evaluation map
    $\Eval : \FreeOperad(C) \to \K \Angle{C}$ and
    by~\ref{item:presentation_operads_2}, the map
    \begin{equation}
        \phi : \FreeOperad(\GeneratingSet)/_{\Angle{\RelationSpace}}
        \to
        \K \Angle{C}
    \end{equation}
    defined linearly for any $x \in C$ by $\phi([x]) := \Eval(x)$, where
    $[x]$ is the image of $x$ through the canonical surjection from
    $\FreeOperad(C)$ to $\FreeOperad(C)/_{\Angle{\RelationSpace}}$, is a
    ns operad morphism. Moreover, by~\ref{item:presentation_operads_1},
    and since $\RelationSpace$ has no element of degree $0$ or $1$,
    $\phi([g]) = g$ for all $g \in \GeneratingSet$. This implies that
    $\phi$ is surjective.
    \smallbreak

    Besides, \ref{item:presentation_operads_3} and
    Proposition~\ref{prop:space_induced_space_quotient_rewrite_rule}
    imply, as spaces, the isomorphisms
    \begin{equation}
        \FreeOperad(\GeneratingSet)/_{\Angle{\RelationSpace}}
        \simeq
        \K \Angle{\NormalForms_{
            (\PlanarRootedTrees^\GeneratingSet_\Leaf, \RewTrees)}}
        \simeq
        \K \Angle{C}
    \end{equation}
    This, together with the fact that $\phi$ is surjective implies that
    $\phi$ is a ns operad isomorphism. Hence, $\K \Angle{C}$ admits the
    claimed presentation.
\end{proof}
\medbreak

In practice, there are at least two ways to use
Theorem~\ref{thm:presentation_operads} to establish a presentation of
a ns operad $\K \Angle{C}$. The first one is the most obvious: it
consists first in finding a generating set $\GeneratingSet$ of
$\K \Angle{C}$, then conjecturing a space of relations $\RelationSpace$
and a rewrite system $(\PlanarRootedTrees^\GeneratingSet_\Leaf, \Rew)$
such that all conditions~\ref{item:presentation_operads_1},
\ref{item:presentation_operads_2}, and~\ref{item:presentation_operads_3}
are satisfied. This can be technical (especially to prove that the
closure $(\PlanarRootedTrees^\GeneratingSet_\Leaf, \RewTrees)$ is
convergent), and relies heavily on computer exploration. The second way
requires as a prerequisite that $\K \Angle{C}$ is combinatorial (and
hence, all its homogeneous components are finite dimensional). In this
case, we need here also to find a generating set $\GeneratingSet$ of
$\K \Angle{C}$, a space of relations $\RelationSpace$ and a rewrite
system $(\PlanarRootedTrees^\GeneratingSet_\Leaf, \Rew)$ such
that~\ref{item:presentation_operads_1},
\ref{item:presentation_operads_2} hold. The difference with the first
way occurs for~\ref{item:presentation_operads_3}: it is now sufficient
to prove $(\PlanarRootedTrees^\GeneratingSet_\Leaf, \RewTrees)$ is
terminating (and not necessarily convergent). Indeed, if
$(\PlanarRootedTrees^\GeneratingSet_\Leaf, \RewTrees)$ is terminating,
since $\K \Angle{C}$ is combinatorial,
\begin{equation} \label{equ:inequality_dim_presentation_operad}
    \dim \K \Angle{C(n)} =
    \# \NormalForms_{
        (\PlanarRootedTrees^\GeneratingSet_\Leaf, \RewTrees)}
    (n)
    \geq
    \dim \FreeOperad(\GeneratingSet(n))/_{\Angle{\RelationSpace}}
\end{equation}
for all $n \geq 1$. The inequality
of~\eqref{equ:inequality_dim_presentation_operad} comes from the fact
that, since we do not know if
$(\PlanarRootedTrees^\GeneratingSet_\Leaf, \RewTrees)$ is confluent,
it can have more normal forms of arity $n$ than the dimension of
$\FreeOperad(\GeneratingSet)/_{\Angle{\RelationSpace}}$ in arity $n$. It
follows from~\eqref{equ:inequality_dim_presentation_operad}, by using
similar arguments as the ones used in the proof of
Theorem~\ref{thm:presentation_operads}, that there is a ns operad
isomorphism from $\FreeOperad(\GeneratingSet)/_{\Angle{\RelationSpace}}$
to $\K \Angle{C}$.
\medbreak

\subsubsection{Koszulity and Koszul duality}
In~\cite{GK94}, Ginzburg and Kapranov extended the notion of Koszul
duality of quadratic associative algebras to quadratic ns operads.
Starting with a ns operad $\K \Angle{C}$ admitting a binary and
quadratic presentation $(\GeneratingSet, \RelationSpace)$ where
$\GeneratingSet$ is finite, the \Def{Koszul dual} of $\K \Angle{C}$ is
the ns operad $\K \Angle{C}^!$, isomorphic to the ns operad admitting
the presentation $\left(\GeneratingSet, \RelationSpace^\perp\right)$
where $\RelationSpace^\perp$ is the annihilator of $\RelationSpace$ in
$\FreeOperad(\GeneratingSet)$ with respect to the scalar product
\begin{equation}
    \langle -, - \rangle :
    \FreeOperad(\GeneratingSet)(3) \otimes
        \FreeOperad(\GeneratingSet)(3)
    \to \K
\end{equation}
linearly defined, for all $x, x', y, y' \in \GeneratingSet(2)$, by
\begin{equation} \label{equ:scalar_product_koszul}
    \left\langle \Corolla{x} \circ_i \Corolla{y},
    \Corolla{x'} \circ_{i'} \Corolla{y'} \right\rangle :=
    \begin{cases}
        1 & \mbox{if }
            x = x', y = y', \mbox{ and } i = i' = 1, \\
        -1 & \mbox{if }
            x = x', y = y', \mbox{ and } i = i' = 2, \\
        0 & \mbox{otherwise}.
    \end{cases}
\end{equation}
Then, with knowledge of a presentation of $\K \Angle{C}$, one can
compute a presentation of~$\K \Angle{C}^!$.
\medbreak

A quadratic ns operad $\K \Angle{C}$ is \Def{Koszul} if its Koszul
complex is acyclic~\cite{GK94,LV12}. Furthermore, when $\K \Angle{C}$ is
Koszul and admits an Hilbert series, the Hilbert series of
$\K \Angle{C}$ and of its Koszul dual $\K \Angle{C}^!$ are
related~\cite{GK94} by
\begin{equation} \label{equ:Hilbert_series_Koszul_operads}
    \HilbSeries_{\K \Angle{C}}
        \left(-\HilbSeries_{\K \Angle{C}^!}(-t)\right)
    = t.
\end{equation}
Relation~\eqref{equ:Hilbert_series_Koszul_operads} can be used either to
prove that a ns operad is not Koszul (it is the case when the
coefficients of the hypothetical Hilbert series of the Koszul dual
admits coefficients that are not nonnegative integers) or to compute the
Hilbert series of the Koszul dual of a Koszul operad.
\medbreak

In all this work, to prove the Koszulity of a ns operad $\K \Angle{C}$,
we shall make use of a tool introduced by Dotsenko and
Khoroshkin~\cite{DK10} in the context of Gröbner bases for operads,
which reformulates in our context, by using rewrite rules on syntax
trees, in the following way.
\medbreak

\begin{Lemma} \label{lem:koszulity_criterion_pbw}
    Let $\K \Angle{C}$ be a ns operad admitting a quadratic presentation
    $(\GeneratingSet, \RelationSpace)$. If there exists an orientation
    $(\PlanarRootedTrees^\GeneratingSet_\Leaf, \Rew)$ of
    $\RelationSpace$ such that its closure
    $(\PlanarRootedTrees^\GeneratingSet_\Leaf, \RewTrees)$ is a
    convergent rewrite system, then $\K \Angle{C}$ is Koszul.
\end{Lemma}
\medbreak

When $(\PlanarRootedTrees^\GeneratingSet_\Leaf, \RewTrees)$ satisfies
the conditions contained in the statement of
Lemma~\ref{lem:koszulity_criterion_pbw}, the set of
$\GeneratingSet$-syntax trees that are normal forms
\begin{math}
    \NormalForms_{(\PlanarRootedTrees^\GeneratingSet_\Leaf, \RewTrees)}
\end{math}
forms a basis of
\begin{math}
    \FreeOperad(\GeneratingSet)/_{\Angle{\RelationSpace}}
\end{math},
called \Def{Poincaré-Birkhoff-Witt basis}. These bases arise from the
work of Hoffbeck~\cite{Hof10} (see also~\cite{LV12}).
\medbreak

\subsubsection{Algebras over ns operads}
Any ns operad $\K \Angle{C}$ encodes a type of graded polynomial
algebras, called \Def{algebras over $\K \Angle{C}$} (or, for short,
\Def{$\K \Angle{C}$-algebras}). A $\K \Angle{C}$-algebra
is a graded polynomial space $\K \Angle{D}$, where $D$ is a graded
collection, and endowed with a linear left action
\begin{equation}
    \Action : \K \Angle{C}(n) \otimes
    {\K \Angle{D}}^{\otimes n} \to \K \Angle{D},
    \qquad n \geq 1,
\end{equation}
satisfying the relations imposed by the structure of a ns operad of
$\K \Angle{C}$, that are
\begin{multline} \label{equ:algebra_over_operad}
    (f \circ_i g) \Action
    \left(a_1 \otimes \dots \otimes a_{n + m - 1}\right) = \\
    f \Action \left(a_1 \otimes \dots
        \otimes a_{i - 1} \otimes
        g \Action \left(a_i \otimes \dots \otimes a_{i + m - 1}\right)
        \otimes a_{i + m} \otimes
        \dots \otimes a_{n + m - 1}\right),
\end{multline}
for all $f \in \K \Angle{C}(n)$, $g \in \K \Angle{C}(m)$, $i \in [n]$,
and
\begin{math}
    a_1 \otimes \dots \otimes a_{n + m - 1}
    \in {\K \Angle{D}}^{\otimes {n + m - 1}}.
\end{math}
In other words, any element $f$ of $\K \Angle{C}$ of arity $n$ plays the
role of a linear operation
\begin{equation}
    f : {\K \Angle{D}}^{\otimes n} \to \K \Angle{D},
\end{equation}
taking $n$ elements of $\K \Angle{D}$ as inputs and computing an element
of $\K \Angle{D}$. Under this point of view,
Relation~\eqref{equ:algebra_over_operad} reads as
\begin{equation}
\,.
\end{equation}
\medbreak

Notice that, by~\eqref{equ:algebra_over_operad}, if $\GeneratingSet$ is
a generating set of $\K \Angle{C}$, it is enough to define the action of
each $x \in \GeneratingSet$ on ${\K \Angle{D}}^{\otimes |x|}$ to wholly
define~$\Action$.
\medbreak

By a slight but convenient abuse of notation, for any
$f \in \K \Angle{C}(n)$, we shall denote by $f(a_1, \dots, a_n)$, or by
$a_1 \, f \, a_2$ if $f$ has arity $2$, the element
$f \Action (a_1 \otimes \dots \otimes a_n)$ of $\K \Angle{D}$, for any
$a_1 \otimes \dots \otimes a_n \in {\K \Angle{D}}^{\otimes n}$. Observe
that by~\eqref{equ:algebra_over_operad}, any associative element of
$\K \Angle{C}$ gives rise to an associative operation on~$\K \Angle{D}$.
\medbreak

The class of all the $\K \Angle{C}$-algebras forms a category, called
\Def{category of $\K \Angle{C}$-algebras}, wherein morphisms
\begin{equation}
    \phi : \K \Angle{D_1} \to \K \Angle{D_2}
\end{equation}
between two $\K \Angle{C}$-algebras $\K \Angle{D_1}$ and
$\K \Angle{D_2}$ are \Def{$\K \Angle{C}$-algebra morphisms}, that are
graded polynomial algebra morphisms satisfying
\begin{equation}
    \phi\left(f\left(a_1, \dots, a_n\right)\right)
    =
    f\left(\phi\left(a_1\right), \dots, \phi\left(a_n\right)\right)
\end{equation}
for all $a_1, \dots, a_n \in \K \Angle{D_1}$ and
$f \in \K \Angle{C}(n)$, $n \in \N_{\geq 1}$.
\medbreak

\subsubsection{Set-operads}
Following Section~\ref{subsubsec:set_algebras}, in a ns set-operad
$\K \Angle{C}$, any partial composition of two elements of $C$ belongs
to $C$. We say in this case that the fundamental basis of $\K \Angle{C}$
is a \Def{set-operad basis}. Besides, by extension, if $\K \Angle{C}$ is
a ns operad admitting a basis $C'$ which is a set-operad basis, we say
that $\K \Angle{C}$ is a \Def{set-operad}. To study a set-operad
$\K \Angle{C}$, it is in some cases convenient to forget about its
linear structure and see it as a graded collection $C$ endowed with
partial composition maps $ \circ_i$, $i \in \N_{\geq 1}$. We will follow
this idea multiple times in the sequel.
In the context of set-operads, a $\K \Angle{C}$-algebra is
called a \Def{$C$-monoid}.
\medbreak

We now state some useful lemmas and notions about ns set-operads.
\medbreak

\begin{Lemma} \label{lem:generation_elements_operads}
    Let $C$ be a ns set-operad generated by a set $\GeneratingSet$ of
    generators. Then any object $x$ of $C$ different from the unit can
    be written as
    \begin{math}
        x = y \circ_i g
    \end{math},
    where $y \in C(n)$, $n \in \N_{\geq 1}$, $g \in \GeneratingSet$,
    and~$i \in [n]$.
\end{Lemma}
\medbreak

Lemma~\ref{lem:generation_elements_operads} is a consequence of the fact
that, since $\GeneratingSet$ is a generating set of $C$, any object of
$C$ admits a treelike expression being a $\GeneratingSet$-syntax tree.
\medbreak

Now, let $C$ be a ns set-operad. Given a subset $S$ of $C$, the
\Def{$S$-degree} of an object $x$ of $C$ is defined by
\begin{equation} \label{equ:S_degree_operad}
    \deg_S(x) := \max \left\{\deg(\Tfr) : \Tfr \in \FreeOperad(S)
        \mbox{ and } \Eval(\Tfr) = x\right\}.
\end{equation}
Of course, the set appearing in~\eqref{equ:S_degree_operad} could be
empty or infinite, so that some elements of $C$ could have no
$S$-degree.
\medbreak

\subsubsection{Left expressions in ns set-operads and hook-length
    formula}
\label{subsubsec:left_expressions_hook_length}
Let $C$ be a ns set-operad, $S$ be a subset of $C$, and $x$ be an object
of $C$. An \Def{$S$-left expression} of $x$ is an expression for $x$ of
the form
\begin{equation}
    x = \left(\dots \left(s_1 \circ_{i_1} s_2\right)
    \circ_{i_2} \dots\right) \circ_{i_{\ell - 1}} s_\ell
\end{equation}
where $s_1, \dots, s_\ell \in S$ and
$i_1, \dots, i_{\ell - 1} \in \N_{\geq 1}$.
\medbreak

A \Def{linear extension} of a syntax tree $\Tfr$ is a linear extension
of the poset $\Qca_\Tfr$ induced by $\Tfr$ (see
Section~\ref{subsubsec:definitions_trees} of
Chapter~\ref{chap:combinatorics}). Left expressions and linear
extensions of treelike expressions in ns set-operads are related, as
shown by the following lemma.
\medbreak

\begin{Lemma} \label{lem:left_expressions_linear_extensions}
    Let $C$ be a combinatorial ns operad and $S$ be a subset
    of $C$. Then, for any object $x$ of $C$, the set of all the
    $S$-left expressions of $x$ is in one-to-one correspondence with the
    set of all pairs $(\Tfr, u)$ where $\Tfr$ is an $S$-treelike
    expression of $x$ and $u$ is a linear extension of~$\Tfr$.
\end{Lemma}
\medbreak

A famous result of Knuth~\cite{Knu98}, known as the
\Def{hook-length formula for trees}, stated here in our setting, says
that given a syntax tree $\Tfr$, the number of linear
extensions of the poset induced by $\Tfr$ is
\begin{equation} \label{equ:hook_length_syntax_tree}
    \HookLength(\Tfr) :=
    \frac{\deg(\Tfr)!}{\prod_{v \in \TreeLanguage_\Node(\Tfr)}
    \deg(\Tfr_v)}.
\end{equation}
\medbreak

A subset $S$ of $C$ \Def{finitely factorizes} $C(1)$ if any element of
$C(1)$ admits finitely many factorizations on $S$ with respect to the
operation~$\circ_1$.
\medbreak

When $S$ finitely factorizes $C(1)$, the number of $S$-treelike
expressions for any object $x$ of $C$ is finite. Hence, in this case, we
deduce from Lemma~\ref{lem:left_expressions_linear_extensions}
and~\eqref{equ:hook_length_syntax_tree} that the number of $S$-left
expressions of $x$ is
\begin{equation} \label{equ:number_left_expressions}
    \sum_{\substack{
        \Tfr \in \FreeOperad(S) \\
        \Eval(\Tfr) = x
    }}
    \HookLength(\Tfr).
\end{equation}
\medbreak

\subsubsection{Colored operads} \label{subsubsec:colored_operads}
Let $\CFr$ be a set of colors. A \Def{nonsymmetric $\CFr$-colored
operad} (or a \Def{ns $\CFr$-colored operad} for short) is a polynomial
space $\K \Angle{C}$ where $C$ is a $\CFr$-colored collection
and $\K \Angle{C}$ is endowed with a set of partially defined binary
linear products $\left\{\circ_i : i \in \N_{\geq 1}\right\}$ of the
form~\eqref{equ:graduation_operads}. The following conditions have to
hold. First, the partial composition $x \circ_i y$ is defined if and
only $\Out(y) = \In_i(y)$ for any $x \in C(n)$, $y \in C(m)$, and
$i \in [n]$. Moreover, when they are well-defined,
Relations~\eqref{equ:graduation_operads}, \eqref{equ:operad_axiom_1},
and~\eqref{equ:operad_axiom_2} have to hold. Finally, we demand the
existence of a set of elements $\{\Unit_a : a \in \CFr\}$ of arity $1$
satisfying
\begin{subequations}
\begin{equation}
    \Out(\Unit_a) = a = \In(\Unit_a),
    \qquad a \in \CFr,
\end{equation}
\begin{equation}
    \Unit_a \circ_1 x = x = x \circ_i \Unit_b,
    \qquad
    x \in C(n), a, b \in \CFr, i \in [n], n \in \N_{\geq 1},
\end{equation}
whenever $\Out(x) = a$ and $\In_i(x) = b$.
\end{subequations}
We call each $\Unit_a$, $a \in \CFr$, the \Def{unit of color $a$}. For
any nonnegative integer $k$, a \Def{ns $k$-colored operad} is a ns
$\CFr$-colored operad where $\CFr$ is a $k$-colored collection. A
\Def{ns monochrome operad} is a ns $\CFr$-colored operad where $\CFr$ is
monochrome.
\medbreak

To describe free ns colored operads, we need the notion of
$\CFr$-colored syntax trees (see
Section~\ref{subsubsec:colored_syntax_trees} of
Chapter~\ref{chap:combinatorics}). Let $\GeneratingSet$ be a
$\CFr$-colored graded collection. The \Def{free ns
$\CFr$-colored operad} over $\GeneratingSet$ is the ns $\CFr$-colored
operad
\begin{equation}
    \FreeColoredOperad(\GeneratingSet) :=
    \K \Angle{\ColoredPlanarRootedTrees^\GeneratingSet},
\end{equation}
where $\ColoredPlanarRootedTrees^\GeneratingSet$ is the $\CFr$-colored
collection of all the $\CFr$-colored $\GeneratingSet$-syntax trees. The
space $\FreeColoredOperad(\GeneratingSet)$ is endowed with the
linearizations of the grafting operations $\circ_i$,
$i \in \N_{\geq 1}$, defined in
Section~\ref{subsubsec:colored_syntax_trees} of
Chapter~\ref{chap:combinatorics}. The units of color $a$, $a \in \CFr$,
are the trees of degree $0$ and arity $1$ with output and input colors
equal to $a$.
\medbreak

\begin{Lemma}
    \label{lem:finitely_factorizing_finite_treelike_expressions}
    Let $C$ be a combinatorial ns $\CFr$-colored set-operad and $S$ be
    a subset of $C$ such that $S$ finitely factorizes $C(1)$. Then, any
    element of $C$ admits finitely many $S$-treelike expressions.
\end{Lemma}
\medbreak

All the notions developed in the above sections about ns operads extends
on ns colored operads by taking colors into account. Classical
references about colored operads are~\cite{BV73,Yau16}.
\medbreak

\subsubsection{Categorical point of view}
In the same way as a monoid $\Mca$ can be seen as a category with
exactly one object $x$ (the elements of $\Mca$ are interpreted as
morphisms $\phi : x \to x$), a ns operad $\Oca$ can be seen as a
multicategory with exactly one object $x$. In this case, the elements
of $\Oca$ of arity $n \in \N_{\geq 1}$ are interpreted as
multimorphisms $\phi : x^n \to x$. The complete composition maps of
$\Oca$ translate as the composition of multimorphisms.
\medbreak

In a similar way, a ns $\CFr$-colored operad $\Cca$ can be seen as a
multicategory having $\CFr$ as set of objects (see~\cite{Cur12}). In
this case, the elements of $\Cca$ having $u_1 \dots u_n \in \CFr^n$ as
word of input colors and $a \in \CFr$ as output color are interpreted as
multimorphisms $\phi :u_1 \times \dots \times u_n \to a$. The complete
composition maps of $\Cca$ translate as the composition of
multimorphisms, where the constraints imposed by the colors in $\Cca$
become constraints imposed by the domains and codomains of
multimorphisms.
\medbreak

\subsubsection{Enrichments}
Ns operads can have some additional structure. Let us now describe some
usual enrichments of operads. In what follows, $\K \Angle{C}$ is a ns
operad.
\medbreak

\paragraph{Basic ns set-operads}
When $\K \Angle{C}$ is a ns set-operad, let for any $y \in C(m)$ the
maps
\begin{equation}
    \circ_i^y : C(n) \to C(n + m - 1),
    \qquad n \in \N_{\geq 1}, i \in [n],
\end{equation}
defined by
\begin{equation}
    \circ_i^y(x) := x \circ_i y
\end{equation}
for all $x \in C$. When all the maps $\circ_i^y$, $y \in C$, are
injective, $\K \Angle{C}$ is a \Def{basic ns set-operad} and $C$ is a
\Def{basic ns set-operad basis} of $\K \Angle{C}$. This notion is a
slightly modified version of the original notion of basic set-operads
introduced by Vallette~\cite{Val07}.
\medbreak

\paragraph{Rooted ns operads}
The ns operad $\K \Angle{C}$ is \Def{rooted} if there is a map
\begin{equation}
    \Root : C(n) \to [n],
    \qquad n \in \N_{\geq 1},
\end{equation}
satisfying, for all $x \in C(n)$, $y \in C(m)$, and $i \in [n]$,
\begin{equation} \label{equ:rooted_operad}
    \Root(x \circ_i y) =
    \begin{cases}
        \Root(x) + m - 1 & \mbox{if } i \leq \Root(x) - 1, \\
        \Root(x) + \Root(y) - 1 & \mbox{if } i = \Root(x), \\
        \Root(x) & \mbox{otherwise (} i \geq \Root(x) + 1 \mbox{)}.
    \end{cases}
\end{equation}
We call such a map a \Def{root map}. More intuitively, the root map of a
rooted ns operad associates a particular input with any of its basis
elements and this input is preserved by partial compositions. It is
immediate that any ns operad $\K \Angle{C}$ is rooted for both the
root maps $\Root_{\mathrm{L}}$ and $\Root_{\mathrm{R}}$ which send
respectively all objects $x$ of $C$ of arity $n$ to $1$ or to $n$. For
this reason, we say that $\K \Angle{C}$ is \Def{nontrivially rooted}
if it can be endowed with a root map different from
$\Root_{\mathrm{L}}$ and $\Root_{\mathrm{R}}$. These notions are slight
variations of ones introduced first in~\cite{Cha14}.
\medbreak

\paragraph{Cyclic ns operads}
Finally, $\K \Angle{C}$ is \Def{cyclic} if there is a map
\begin{equation} \label{equ:rotation_map_0}
    \rho : \K \Angle{C}(n) \to \K \Angle{C}(n),
    \qquad n \in \N_{\geq 1},
\end{equation}
satisfying, for all $f \in \K \Angle{C}(n)$, $g \in \K \Angle{C}(m)$,
and $i \in [n]$,
\begin{subequations}
\begin{equation} \label{equ:rotation_map_1}
    \rho(\Unit) = \Unit,
\end{equation}
\begin{equation} \label{equ:rotation_map_2}
    \rho^{n + 1}(f) = f,
\end{equation}
\begin{equation} \label{equ:rotation_map_3}
    \rho(f \circ_i g) =
    \begin{cases}
        \rho(g) \circ_m \rho(f) & \mbox{if } i = 1, \\
        \rho(f) \circ_{i - 1} g & \mbox{otherwise}.
    \end{cases}
\end{equation}
\end{subequations}
We call such a map $\rho$ a \Def{rotation map}. Intuitively, a rotation
map in a ns operad acts by transforming its $1$st input of an element
$f$ in an output, its $2$nd input in a $1$st input, its $3$rd input in
a $2$nd input, and so on,  and its output in a last input. This notion
has been introduced in\cite{GK95}.
\medbreak

\subsubsection{Symmetric operads} \label{subsubsec:symmetric_operads}
Let $\Per$ be the ns operad defined by
$\Per := \K \Angle{\Augmentation(\SymmetricGroup)}$ with the partial
compositions defined as follow. For all $\sigma \in \SymmetricGroup(n)$,
$\nu \in \SymmetricGroup(m)$, and $i \in [n]$,
\begin{equation}
    \sigma \circ_i \nu :=
    \sigma'(1) \dots \sigma'(i - 1)
    \nu'(1) \dots \nu'(m)
    \sigma'(i + 1) \dots \sigma'(n),
\end{equation}
where, for any $j \in [m]$,
\begin{equation}
    \nu'(j) := \nu(j) + \sigma(i) - 1,
\end{equation}
and, or any $j \in [n]$,
\begin{equation}
    \sigma'_j :=
    \begin{cases}
        \sigma_j & \mbox{if } \sigma_j < \sigma_i, \\
        \sigma_j + m - 1 & \mbox{otherwise}.
    \end{cases}
\end{equation}
For instance, here are two examples of compositions in $\Per$:
\begin{subequations}
\begin{equation}
    \textcolor{Col1}{1}{\bf 2}\textcolor{Col1}{3}
    \circ_2
    \textcolor{Col4}{12}
    = \textcolor{Col1}{1}\textcolor{Col4}{23}\textcolor{Col1}{4},
\end{equation}
\begin{equation}
    \textcolor{Col1}{741}{\bf 5}\textcolor{Col1}{623}
    \circ_4
    \textcolor{Col4}{231}
    = \textcolor{Col1}{941}\textcolor{Col4}{675}\textcolor{Col1}{823}.
\end{equation}
\end{subequations}
This ns operad is known as the \Def{associative noncommutative operad}
or more prosaically, the \Def{operad of permutations}.
\medbreak

A \Def{symmetric operad}, or an \Def{operad} for short, is a ns operad
$\K \Angle{C}$ together with linear maps
\begin{equation}
    \Action : \K \Angle{C}(n) \otimes
    \Per(n) \to \K \Angle{C}(n), \qquad n \geq 1,
\end{equation}
satisfying, for any $x \in C(n)$, $y \in C(m)$,
$\sigma \in \SymmetricGroup(n)$,
$\nu \in \SymmetricGroup(m)$, and $i \in [n]$,
\begin{equation} \label{equ:equivariance}
    (x \Action \sigma) \circ_i (y \Action \nu) =
    \left(x \circ_{\sigma_i} y\right) \Action
    \left(\sigma \circ_i \nu \right),
\end{equation}
and in such a way that $\Action$ also is a symmetric group action. Note
that any operad $\K \Angle{C}$ is also (and thus can be seen as) a ns
operad by forgetting its action of $\Per$.
\medbreak

A \Def{simple permutation} is a permutation $\sigma$ such that for all
factors $u$ of $\sigma$, if the letters of $u$ form an interval of
$\N$ then $|u| = 1$ or $|u| = |\sigma|$. For instance, the permutation
$624135$ is not simple since the letters of the factor $u := 2413$
form an interval of $\N$. On the other hand, the permutation
$5137462$ is simple.
\medbreak

\begin{Proposition} \label{prop:generating_set_per}
    As a ns operad, $\Per$ is minimally generated by the set of
    all simple permutations of sizes $2$ or more.
\end{Proposition}
\medbreak

One of the simplest examples of symmetric operads is the
\Def{commutative associative operad} $\Com$. This operad is defined
as $\Com := \K \Angle{C}$ where $C$ is the augmented graded collection
satisfying $C(n) := \{\Afr_n\}$ for all $n \in \N_{\geq 1}$, its partial
compositions satisfy
\begin{equation}
    \Afr_n \circ_i \Afr_m := \Afr_{n + m - 1},
\end{equation}
for any $n, m \in \N_{\geq 1}$ and $i \in [n]$, and $\Per$ acts
trivially on $\Com$.
\medbreak

There is also a notion of algebras over symmetric operads. An
\Def{algebra over $\K \Angle{C}$} (or, short, a $\K \Angle{C}$-algebra)
is an algebra $\K \Angle{D}$ over $\K \Angle{C}$ seen as a ns operad.
We ask additionally that the relation
\begin{equation}
    (f \Action \sigma)(a_1, \dots, a_n) =
    f\left(a_{\sigma^{-1}(1)}, \dots a_{\sigma^{-1}(n)}\right)
\end{equation}
holds for any $f \in \K \Angle{C}(n)$, $\sigma \in \SymmetricGroup(n)$,
$a_1, \dots, a_n \in \K \Angle{D}$, $n \in \N_{\geq 1}$.
\medbreak

\subsection{Main operads in combinatorics}
This section contains a list of common ns operads studied or encountered
in the sequel. A classical text containing a list of definitions of
operads is~\cite{Zin12}.
\medbreak

\subsubsection{Associative operad} \label{subsubsec:associative_operad}
The \Def{associative operad} $\As$ is defined as the operad $\Com$ seen
as a nonsymmetric one (see Section~\ref{subsubsec:symmetric_operads}).
\medbreak

The ns operad $\As$ is a set-operad and its Hilbert series satisfies
\begin{equation}
    \HilbSeries_\As(t) = \frac{t}{1 - t}.
\end{equation}
Moreover, $\As$ admits the presentation
$(\GeneratingSet, \RelationSpace)$ where $\GeneratingSet := \{\Afr_2\}$
and $\RelationSpace$ is the space generated by
\begin{equation}
    \Corolla{\Afr_2} \circ_1 \Corolla{\Afr_2} -
    \Corolla{\Afr_2} \circ_2 \Corolla{\Afr_2}.
\end{equation}
Any algebra over $\As$ is a space $\K \Angle{D}$ endowed with a binary
associative operation.
\medbreak

\subsubsection{Magmatic operad} \label{subsubsec:magmatic_operad}
Let $\Mag := \K \Angle{\Ary^{(2)}_\Leaf}$ be the ns operad where
for any binary trees $\Tfr$ and $\Sfr$,
the partial composition $\Tfr \circ_i \Sfr$ is the grafting of $\Sfr$
onto the $i$th leaf of $\Tfr$, seen as syntax trees. In other terms,
$\Mag$ is the operad $\FreeOperad(C)$ where $C := C(2) := \{\Asf\}$.
This ns operad is the \Def{magmatic operad}.
\medbreak

The ns operad $\Mag$ is a set-operad and its Hilbert series satisfies
\begin{equation}
    \HilbSeries_\Mag(t) = \frac{1 - \sqrt{1 - 4t}}{2}.
\end{equation}
Moreover, $\Mag$ admits the presentation
$(\GeneratingSet, \RelationSpace)$ where
\begin{equation}
    \GeneratingSet := \left\{\BinaryNode\right\}
\end{equation}
and $\RelationSpace$ is the trivial space. Any algebra over $\Mag$ is a
space $\K \Angle{D}$ endowed with a binary operation which do satisfy
any required relation.
\medbreak

\subsubsection{Duplicial operad} \label{subsubsec:duplicial_operad}
Let $\Dup := \K \Angle{\Augmentation(\Ary^{(2)}_\Node)}$ be the ns
operad where for any nonempty binary trees $\Tfr$ and $\Sfr$, the
partial composition $\Tfr \circ_i \Sfr$ consists in replacing the $i$th
(with respect to the infix order) internal node $u$ of $\Tfr$ by a copy
of $\Sfr$, and by grafting the left subtree of $u$ to the first leaf of
the copy, and the right subtree of $u$ to the last leaf of the copy.
For instance,
\begin{equation}
\,.
\end{equation}
This ns operad is the \Def{duplicial operad}~\cite{Lod08}.
\medbreak

The ns operad $\Dup$ is a set-operad and its Hilbert series satisfies
\begin{equation}
    \HilbSeries_{\Dup}(t) = \frac{1 - 2t - \sqrt{1 - 4t}}{2t}.
\end{equation}
Moreover, $\Dup$ admits the presentation $(\GeneratingSet, \RelationSpace)$
where
\begin{equation} \label{equ:generating_set_dup}
    \GeneratingSet :=
    \left\{\BinaryTreeLeft\,, \BinaryTreeRight\right\}
\end{equation}
and $\RelationSpace$ is the space generated by,
by denoting by $\LDup$ (resp. $\RDup$) the first (resp. second) tree
of~\eqref{equ:generating_set_dup},
\begin{subequations}
\begin{equation} \label{equ:operad_dup_relation_1}
    \Corolla{\LDup} \circ_1 \Corolla{\LDup} -
    \Corolla{\LDup} \circ_2 \Corolla{\LDup}\,,
\end{equation}
\begin{equation} \label{equ:operad_dup_relation_2}
    \Corolla{\RDup} \circ_1 \Corolla{\LDup} -
    \Corolla{\LDup} \circ_2 \Corolla{\RDup}\,,
\end{equation}
\begin{equation} \label{equ:operad_dup_relation_3}
    \Corolla{\RDup} \circ_1 \Corolla{\RDup} -
    \Corolla{\RDup} \circ_2 \Corolla{\RDup}\,.
\end{equation}
\end{subequations}
Any algebra over $\Dup$ is a space $\K \Angle{D}$ endowed with two
binary relations $\LDup$ and $\RDup$ such that both $\LDup$ and
$\RDup$ are associative (as consequences
of~\eqref{equ:operad_dup_relation_1}
and~\eqref{equ:operad_dup_relation_2}), and, for any $x, y, z \in D$,
\begin{equation}
    (x \LDup y) \RDup z = x \LDup (y \RDup z).
\end{equation}
These structures are called \Def{duplicial algebras}.
\medbreak

\subsubsection{Operad of rational functions}
\label{subsubsec:rational_functions_operad}
The graded vector space of all commutative rational functions
$\K(\Ubb)$, where $\Ubb$ is the infinite commutative alphabet
$\{u_1, u_2, \dots\}$, has the structure of a ns operad $\RatFct$
introduced by Loday~\cite{Lod10} and defined as follows. Let
$\RatFct(n)$ be the subspace $\K(u_1, \dots, u_n)$ of $\K(\Ubb)$ and
\begin{equation}
    \RatFct := \bigoplus_{n \geq 1} \RatFct(n).
\end{equation}
Observe that since $\RatFct$ is a graded space, each rational function
has an arity. Hence, by setting $f_1(u_1) := 1$ and $f_2(u_1, u_2) := 1$,
$f_1$ is of arity $1$ while $f_2$ is of arity $2$, so that $f_1$ and
$f_2$ are considered as different rational functions. The partial
composition of two rational functions $f \in \RatFct(n)$ and
$g \in \RatFct(m)$ satisfies, for any $i \in [n]$,
\begin{equation} \label{equ:partial_composition_RatFct}
    f \circ_i g :=
    f\left(u_1, \dots, u_{i - 1}, u_i + \dots + u_{i + m - 1},
        u_{i + m}, \dots, u_{n + m - 1}\right)
    \;
    g\left(u_i, \dots, u_{i + m - 1}\right).
\end{equation}
This ns operad is the \Def{operad of rational functions}.
The rational function $f$ of $\RatFct(1)$ defined by $f(u_1) := 1$ is
the unit of~$\RatFct$. As shown by Loday, this operad is (nontrivially)
isomorphic to the operad $\Mould$ introduced by Chapoton~\cite{Cha07}.
\medbreak

\subsubsection{Diassociative operad}
\label{subsubsec:diassociative_operad}
Let the operad $\Dias := \K \Angle{C}$ where $C$ is the augmented
graded collection satisfying
$C(n) := \left\{\Efr_{n, k} : k \in [n]\right\}$
for all $n \in \N_{\geq 1}$. The partial compositions of $\Dias$
are defined by
\begin{equation} \label{equ:partial_composition_Dias}
    \Efr_{n, k} \circ_i \Efr_{m, \ell} =
    \begin{cases}
        \Efr_{n + m - 1, k + m - 1}
            & \mbox{if } i < k, \\
        \Efr_{n + m - 1, k + \ell - 1}
            & \mbox{if } i = k, \\
        \Efr_{n + m - 1, k}
            & \mbox{otherwise (} i > k \mbox{)},
    \end{cases}
\end{equation}
for all $n, m \in \N_{\geq 1}$, $k \in [n]$, $\ell \in [m]$, and
$i \in [n]$. This operad is the \Def{diassociative operad}.
\medbreak

The ns operad $\Dias$ is a set-operad and its Hilbert series satisfies
\begin{equation}
    \HilbSeries_\Dias(t) = \frac{t}{(1 - t)^2}.
\end{equation}
Moreover, $\Dias$ admits the presentation
$(\GeneratingSet, \RelationSpace)$ where
\begin{equation} \label{equ:generating_set_dias}
    \GeneratingSet :=
    \left\{\Efr_{2, 1}, \Efr_{2, 2}\right\}
\end{equation}
and $\RelationSpace$ is the space generated by, by denoting by
$\LDias$ (resp. $\RDias$) the first (resp. second) element
of~\eqref{equ:generating_set_dias},
\begin{subequations}
\begin{equation} \label{equ:operad_dias_relation_1}
    \Corolla{\LDias} \circ_1 \Corolla{\LDias} -
    \Corolla{\LDias} \circ_2 \Corolla{\LDias},
    \quad
    \Corolla{\LDias} \circ_1 \Corolla{\LDias} -
    \Corolla{\LDias} \circ_2 \Corolla{\RDias},
\end{equation}
\begin{equation} \label{equ:operad_dias_relation_2}
    \Corolla{\LDias} \circ_1 \Corolla{\RDias} -
    \Corolla{\RDias} \circ_2 \Corolla{\LDias},
\end{equation}
\begin{equation} \label{equ:operad_dias_relation_3}
    \Corolla{\RDias} \circ_1 \Corolla{\LDias} -
    \Corolla{\RDias} \circ_2 \Corolla{\RDias},
    \quad
    \Corolla{\RDias} \circ_1 \Corolla{\RDias} -
    \Corolla{\RDias} \circ_2 \Corolla{\RDias}.
\end{equation}
\end{subequations}
This operad, by its presentation by generators and relations, has
been introduced in~\cite{Lod01}. Its realization in terms of the
elements $\Efr_{n, k}$ and the partial compositions
maps~\eqref{equ:partial_composition_Dias} appears in~\cite{Cha05}.
Any algebra over $\Dias$ is a space $\K \Angle{D}$ endowed with two
binary relations $\LDias$ and $\RDias$ such that both $\LDias$ and
$\RDias$ are associative (as particular consequences
of~\eqref{equ:operad_dias_relation_1}
and~\eqref{equ:operad_dias_relation_3}), and, for any $x, y, z \in D$,
\begin{subequations}
\begin{equation}
    (x \LDias y) \LDias z = x \LDias (y \RDias z),
\end{equation}
\begin{equation}
    (x \RDias y) \LDias z = x \RDias (y \LDias z),
\end{equation}
\begin{equation}
    (x \LDias y) \RDias z = x \RDias (y \RDias z).
\end{equation}
\end{subequations}
These structures are called \Def{diassociative algebras}.
\medbreak

\subsubsection{Dendriform operad} \label{subsubsec:dendriform_operad}
The \Def{dendriform operad} $\Dendr$ is the ns suboperad of $\RatFct$
generated by the set
$\left\{\frac{1}{u_1}, \frac{1}{u_2}\right\}$~\cite{Lod10}. This operad
admits the presentation $(\GeneratingSet, \RelationSpace)$ where
\begin{equation}
    \GeneratingSet := \GeneratingSet(2) := \{\LDendr, \RDendr\}
\end{equation}
and $\RelationSpace$ is the space generated by
\begin{subequations}
\begin{equation} \label{equ:operad_dendr_relation_1}
    \Corolla{\LDendr} \circ_1 \Corolla{\LDendr} -
    \Corolla{\LDendr} \circ_2 \Corolla{\LDendr} -
    \Corolla{\LDendr} \circ_2 \Corolla{\RDendr},
\end{equation}
\begin{equation} \label{equ:operad_dendr_relation_2}
    \Corolla{\LDendr} \circ_1 \Corolla{\RDendr} -
    \Corolla{\RDendr} \circ_2 \Corolla{\LDendr},
\end{equation}
\begin{equation} \label{equ:operad_dendr_relation_3}
    \Corolla{\RDendr} \circ_1 \Corolla{\LDendr} +
    \Corolla{\RDendr} \circ_1 \Corolla{\RDendr} -
    \Corolla{\RDendr} \circ_2 \Corolla{\RDendr}.
\end{equation}
\end{subequations}
This operad, by its presentation by generators and relations, has
been introduced in~\cite{Lod01}. It is shown here that $\Dendr$ is the
Koszul dual of $\Dias$ and that these operads are Koszul operads.
Hence, the Hilbert series $\HilbSeries_\Dendr(t)$ of $\Dendr$
satisfies, by~\eqref{equ:Hilbert_series_Koszul_operads},
\begin{equation}
    \HilbSeries_\Dendr(t) = \frac{1 - 2t - \sqrt{1 - 4t}}{2t}.
\end{equation}
This shows that $\Dendr$ is, as a combinatorial polynomial space,
the space $\K \Angle{\Augmentation(\Ary^{(2)}_\Node)}$. Moreover,
the operads $\Dias$ and $\Dendr$ are Koszul dual one of the other.
Finally, any algebra of $\Dendr$ is a dendriform algebra
(see Section~\ref{subsubsec:dendriform_algebras}).
\medbreak

The free dendriform algebra over one generator is the space $\Dendr$,
that is the linear span of all nonempty binary trees, endowed with
the linear operations
\begin{equation}
    \LDendr, \RDendr :
    \Dendr \otimes \Dendr \to \Dendr,
\end{equation}
defined recursively, for any nonempty  tree $\Sfr$, and binary trees
$\Tfr_1$ and $\Tfr_2$ by
\begin{subequations}
\begin{equation}
    \Sfr \LDendr \LeafPic
    := \Sfr =:
    \LeafPic \RDendr \Sfr,
\end{equation}
\begin{equation}
    \LeafPic \LDendr \Sfr := 0 =: \Sfr \RDendr \LeafPic,
\end{equation}
\begin{equation}
    \BinaryRoot{\Tfr_1}{\Tfr_2} \LDendr \Sfr :=
    \BinaryRoot{\Tfr_1}{\Tfr_2 \LDendr \Sfr}
    + \BinaryRoot{\Tfr_1}{\Tfr_2 \RDendr \Sfr}\,,
\end{equation}
\begin{equation}
    \begin{split} \Sfr \RDendr \end{split}
    \BinaryRoot{\Tfr_1}{\Tfr_2} :=
    \BinaryRoot{\Sfr \RDendr \Tfr_1}{\Tfr_2}
    + \BinaryRoot{\Sfr \LDendr \Tfr_1}{\Tfr_2}.
\end{equation}
\end{subequations}
Note that neither $\LeafPic \LDendr \LeafPic$ nor
$\LeafPic \RDendr \LeafPic$ are defined.
We have for instance,
\begin{subequations}
\begin{equation}
\,.
\end{equation}
\end{subequations}
\medbreak

\subsubsection{Operad of gravity chord diagrams}
\label{subsubsec:operad_Grav}
The \Def{operad of gravity chord diagrams} $\Grav$ is an operad defined
in~\cite{AP15}. This operad is the nonsymmetric version of the gravity
operad, a symmetric operad introduced by Getzler~\cite{Get94}. Let us
describe this operad.
\medbreak

A \Def{gravity chord diagram} is a configuration $\Cfr$ (see
Section~\ref{subsubsec:configurations} of
Chapter~\ref{chap:combinatorics}) where each arc can be
\Def{blue} (drawn as a thick line) and that satisfies the following
conditions. By denoting by $n$ the size of $\Cfr$, all the edges and the
base of $\Cfr$ are blue, and if $(x, y)$ and $(x', y')$ are two blue
crossing diagonals of $\Cfr$ such that $x < x'$, the arc $(x', y)$ is
not labeled. In other words, the quadrilateral formed by the vertices
$x$, $x'$, $y$, and $y'$ of $\Cfr$ is such that its side $(x', y)$ is
not labeled. For instance,
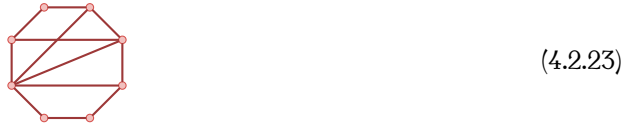
\begin{equation}
    \begin{tikzpicture}[scale=.8,Centering]
        \node[CliquePoint](1)at(-0.38,-0.92){};
        \node[CliquePoint](2)at(-0.92,-0.38){};
        \node[CliquePoint](3)at(-0.92,0.38){};
        \node[CliquePoint](4)at(-0.38,0.92){};
        \node[CliquePoint](5)at(0.38,0.92){};
        \node[CliquePoint](6)at(0.92,0.38){};
        \node[CliquePoint](7)at(0.92,-0.38){};
        \node[CliquePoint](8)at(0.38,-0.92){};
        \draw[CliqueEdgeBlue](1)--(2);
        \draw[CliqueEdgeBlue](1)--(8);
        \draw[CliqueEdgeBlue](2)--(3);
        \draw[CliqueEdgeBlue](2)--(5);
        \draw[CliqueEdgeBlue](2)--(6);
        \draw[CliqueEdgeBlue](2)--(7);
        \draw[CliqueEdgeBlue](3)--(4);
        \draw[CliqueEdgeBlue](3)--(6);
        \draw[CliqueEdgeBlue](4)--(5);
        \draw[CliqueEdgeBlue](5)--(6);
        \draw[CliqueEdgeBlue](6)--(7);
        \draw[CliqueEdgeBlue](7)--(8);
    \end{tikzpicture}
\end{equation}
is a gravity chord diagram of arity $7$ having four blue diagonals
(observe in particular that, as required, the arc $(3, 5)$ is not
labeled). For any $n \geq 2$, $\Grav(n)$ is the linear span of all
gravity chord diagrams of size $n$. Moreover, $\Grav(1)$ is the linear
span of the singleton containing the only polygon of arity $1$ where
its only arc is not labeled. The partial composition of $\Grav$ is
defined, in a geometric way, as follows. For any gravity chord diagrams
$\Cfr$ and $\Dfr$ of respective arities $n$ and $m$, and $i \in [n]$,
the gravity chord diagram $\Cfr \circ_i \Dfr$ is obtained by gluing the
base of $\Dfr$ onto the $i$th edge of $\Cfr$, so that the arc
$(i, i + m)$ of $\Cfr \circ_i \Dfr$ is blue. For example,
\begin{equation}
\,.
\end{equation}
\medbreak

\section{Pros in combinatorics} \label{sec:pros}
We regard here pros as polynomial algebras and provide the main
definitions used in the following chapters. We also give examples
of some usual pros.
\medbreak

\subsection{Pros} \label{subsec:pros}
Pros can be thought as variations of operads allowing multiple outputs
for some elements. Surprisingly, pros appeared earlier than operads in
the work of Mac Lane~\cite{McL65}. Intuitively, a pro is a space
$\K \Angle{C}$ wherein elements are biproducts (see
Section~\ref{subsubsec:biproducts}) and is endowed with two
operations: an horizontal composition and a vertical composition. The
first operation takes two operators $x$ and $y$ of $\K \Angle{C}$ and
builds a new one whose inputs (resp. outputs) are, from left to right,
those of $x$ and then those of $y$. The second operation takes two
operators $x$ and $y$ of $\K \Angle{C}$ and produces a new one obtained
by plugging the outputs of $y$ onto the inputs of $x$. Basic and modern
references about pros are~\cite{Lei04} and~\cite{Mar08}.
\medbreak

\subsubsection{Categorical definition}
A \Def{product category} (or, for short, a \Def{pro}) is a category
$\Pca$ endowed with a associative bifunctor
$* : \Pca \times \Pca \to \Pca$ such that the objects of $\Pca$ are
the elements of $\N$ and $x * y := x + y$ for all $x, y \in \N$.
\medbreak

This formal definition of pros is not combinatorial. Let us provide in
the next section a more concrete one.
\medbreak

\subsubsection{Axioms of pros}
A \Def{pro} is a bigraded polynomial space $\K \Angle{C}$ endowed with
a binary linear product
\begin{equation}
    * : \K \Angle{C}(p, q) \times \K \Angle{C}(p', q')
    \to \K \Angle{C}(p + p', q + q'),
    \qquad p, p', q, q' \geq 0,
\end{equation}
called \Def{horizontal composition} and
a binary linear product
$\circ$ is a map of the form
\begin{equation}
    \circ : \K \Angle{C}(q, r) \times \K \Angle{C}(p, q)
    \to \K \Angle{C}(p, r),
    \qquad p, q, r \geq 0,
\end{equation}
called \Def{vertical composition}.
We demand, for any $p \in \N$, the existence of an element $\Unit_p$ of
$\K \Angle{C}(p, p)$ called \Def{unit of arity $p$}. The \Def{input}
(resp. \Def{output}) \Def{arity} of $f \in \K \Angle{C}(p, q)$
is $|f|_\ArityIn := p$ (resp. $|x|_\ArityOut := q$).
\medbreak

These data have to satisfy for all $x, y, z, x', y' \in C$ the six relations
\begin{equation} \label{equ:assoc_compo_h}
    (x * y) * z = x * (y * z),
\end{equation}
\begin{equation} \label{equ:assoc_compo_v}
    (x \circ y) \circ z = x \circ (y \circ z),
    \qquad
    |x|_\ArityIn = |y|_\ArityOut,
    |y|_\ArityIn = |z|_\ArityOut,
\end{equation}
\begin{equation} \label{equ:compo_h_v}
    (x \circ y) * (x' \circ y') = (x * x') \circ (y * y'),
    \qquad \mbox{ if }
    |x|_\ArityIn = |y|_\ArityOut,
    |x'|_\ArityIn = |y'|_\ArityOut,
\end{equation}
\begin{equation} \label{equ:relation_unit_h}
    \Unit_p * \Unit_q = \Unit_{p + q},
    \qquad p, q \geq 0,
\end{equation}
\begin{equation} \label{equ:relation_unit_zero}
    x * \Unit_0 = x = \Unit_0 * x,
\end{equation}
\begin{equation} \label{equ:relation_unit_v}
    x \circ \Unit_p = x = \Unit_q \circ x,
    \qquad p, q \geq 0, \mbox{ if }
    |x|_\ArityIn = p,
    |x|_\ArityOut = q.
\end{equation}
\medbreak

Since a pro is a particular polynomial algebra, all the properties and
definitions about polynomial algebras exposed in
Section~\ref{subsec:polynomial_bialgebras} remain valid for pros (like
pros morphisms, sub-pros, generating  sets, pro ideals and quotients,
{\em etc.}).
\medbreak

\subsubsection{Free pros} \label{subsubsec:free_pros}
Let $\GeneratingSet$ be a bigraded collection such that
$\GeneratingSet(p, q) = \emptyset$ if $p = 0$ or $q = 0$.
The \Def{free pro} over $\GeneratingSet$ is the pro
\begin{equation}
    \FreePro(\GeneratingSet) :=
    \K \Angle{\Prographs^\GeneratingSet},
\end{equation}
where $\Prographs^\GeneratingSet$ is the bigraded collection of all
the $\GeneratingSet$-prographs (see Section~\ref{subsec:prographs} of
Chapter~\ref{chap:combinatorics}). The space $\FreePro(\GeneratingSet)$
is endowed with linearizations of the horizontal composition
of prographs and of the vertical composition of prographs
(see Section~\ref{subsubsec:operations_prographs} of
Chapter~\ref{chap:combinatorics}). The unit of arity $p$, $p \geq 0$,
is the sequence of wires $\Unit_p$. Notice that by the above assumption
on $\GeneratingSet$, there is no elementary $\GeneratingSet$-prograph in
$\FreePro(G)$ with a null input or output arity. Therefore, $\Unit_0$ is
the only element of $\FreePro(G)$ with a null input (resp. output)
arity. In this dissertation, we consider only free pros satisfying this
property.
\medbreak

Let us now state some definitions and properties about free pros.
\medbreak

Let $\FreePro(\GeneratingSet)$ be a free pro. Since
$\FreePro(\GeneratingSet)$ is free (and, by convention,
$\GeneratingSet$ has no generator of input or output arity $0$), any
$\GeneratingSet$-prograph $x$  can be uniquely written as
\begin{equation}
    x = x_1 * \dots * x_\ell
\end{equation}
where the $x_i$ are $\GeneratingSet$-prographs different from $\Unit_0$,
and $\ell \geq 0$ is maximal. We call the word
$\Dec(x) := (x_1, \dots, x_\ell)$ the \Def{maximal decomposition} of $x$
and the $x_i$ the \Def{factors} of $x$. Notice that the maximal
decomposition of $\Unit_0$ is the empty word. We have in
$\FreePro(\GeneratingSet)$, for instance, by setting $\GeneratingSet$
as the bigraded collection defined by
$\GeneratingSet := \GeneratingSet(2, 2) \sqcup \GeneratingSet(3, 1)$
where $\GeneratingSet(2, 2) := \{\Asf\}$ and
$\GeneratingSet(3, 1) := \{\Bsf\}$,
\begin{equation}
    \Dec\left(\;

    \;\right).
\end{equation}
\medbreak

An $\GeneratingSet$-prograph $x$ is \Def{reduced} if all its factors are
different from $\Unit_1$. For any $\GeneratingSet$-prograph $x$, we
denote by $\Reduced(x)$ the reduced $\GeneratingSet$-prograph admitting
as maximal decomposition the longest subword of $\Dec(x)$ consisting in
factors different from $\Unit_1$. We have in $\FreePro(\GeneratingSet)$,
for instance,
\begin{equation}
    \Reduced\left(\;
    %
\,.
\end{equation}
By extension, we denote by $\Reduced(\FreePro(\GeneratingSet))$ the set
of all the reduced $\GeneratingSet$-prographs. Note that $\Unit_0$
belongs to~$\Reduced(\FreePro(\GeneratingSet))$.
\medbreak

Besides, we say that a $\GeneratingSet$-prograph is \Def{indecomposable}
if its maximal decomposition consists in exactly one factor. Note that
$\Unit_0$ is not indecomposable while $\Unit_1$ is.
\medbreak

\begin{Lemma} \label{lem:free_pro_element_reduced}
    Let $x$ and $y$ be two $\GeneratingSet$-prographs such that
    $x = \Reduced(y)$. Then, by denoting by $(x_1, \dots, x_\ell)$ the
    maximal decomposition of $x$, there exists a unique sequence of
    nonnegative integers $p_1, \dots, p_\ell, p_{\ell + 1}$ such that
    \begin{equation}
        y = \Unit_{p_1} * x_1 * \Unit_{p_2} * x_2 *
            \dots * x_\ell * \Unit_{p_{\ell + 1}}.
    \end{equation}
\end{Lemma}
\medbreak

\begin{Lemma} \label{lem:square_rule_free_pros}
    Let $x$, $y$, $z$, and $t$ be four $\GeneratingSet$-prographs such
    that $x * y = z \circ t$. Then, there exist four unique elements
    $x_1$, $x_2$, $y_1$, $y_2$ of $\Pca$ such that $x = x_1 \circ x_2$,
    $y = y_1 \circ y_2$, $z = x_1 * y_1$, and $t = x_2 * y_2$.
\end{Lemma}
\medbreak

\subsection{Main pros in combinatorics}
This section contains a short list of common pros
(see~\cite{Laf03,Laf11} for many examples of pros).
\medbreak

\subsubsection{Pro of maps}
Let $\Map$ be the bigraded collection wherein for any $p, q \in \N$,
$\Map(p, q)$ is the set of maps from $[p]$ to $[q]$. We endow
$\K \Angle{\Map}$ with an horizontal composition $*$ defined linearly,
for any maps $f \in \Map(p_1, q_1)$, $g \in \Map(p_2, q_2)$, and
$x \in [p_1 + q_1]$, by
\begin{equation}
    (f * g)(x) :=
    \begin{cases}
        f(x) & \mbox{if } x \in [|f|_\ArityIn], \\
        g(x) + |f|_\ArityOut & \mbox{otherwise}.
    \end{cases}
\end{equation}
We endow the polynomial space $\K \Angle{\Map}$ with a vertical
composition $\circ$ defined linearly, for any maps
$f \in \Map(p_1, q_1)$, $g \in \Map(p_2, p_1)$, and $x \in [p_2]$, by
\begin{equation}
    (f \circ g)(x) = f(g(x)).
\end{equation}
The unit $\Unit_p$ of arity $p \in \N$ is the identity map on $[p]$.
\medbreak

For instance, by denoting maps of $\Map$ by words (the $i$th letters of
the words being the images of $i$), one has
\begin{subequations}
\begin{equation} \label{equ:example_pro_Map_1}
    \textcolor{Col1}{3115} * \textcolor{Col4}{133} =
    \textcolor{Col1}{3115} \textcolor{Col4}{799},
    \qquad
    \textcolor{Col1}{3115} \in \Map(4, 6),
    \textcolor{Col4}{133} \in \Map(3, 5),
\end{equation}
\begin{equation} \label{equ:example_pro_Map_2}
    \textcolor{Col1}{1224244} \circ \textcolor{Col4}{3312} =
    2212,
    \qquad
    \textcolor{Col1}{1224244} \in \Map(7, 6),
    \textcolor{Col4}{3312} \in \Map(4, 7).
\end{equation}
\end{subequations}
By seeing each map $f \in \Map(p, q)$ as an operator with
$p$ inputs and $q$ outputs where each $i$th input is connected to the
$f(i)$th output, \eqref{equ:example_pro_Map_1}
and~\eqref{equ:example_pro_Map_2} read respectively as
\begin{subequations}
\begin{equation}
\,.
\end{equation}
\end{subequations}
\medbreak

Observe hence that the horizontal composition of $\K \Angle{\Map}$ is a
shifted concatenation of words and that the vertical composition of
$\K \Angle{\Map}$ is a functional composition. We call
$\K \Angle{\Map}$ the \Def{pro of maps}.
\medbreak

\subsubsection{Pro of increasing maps}
Let $\NDMap$ be the bigraded collection wherein for any $p, q \in \N$,
$\NDMap(p, q)$ is the set of nondecreasing maps from $[p]$ to $[q]$,
that is the maps $f$ such that $i \leq j$ implies $f(i) \leq f(j)$.
Since the operations $*$ and $\circ$ of $\K \Angle{\Map}$ are stable in
$\K \Angle{\NDMap}$ and the identity maps $\Unit_p$, $p \in \N$,
are nondecreasing maps, $\K \Angle{\NDMap}$ is a sub-pro of
$\K \Angle{\Map}$. In particular, we call $\K \Angle{\NDMap}$ the
\Def{pro of nondecreasing maps}.
\medbreak

\subsubsection{Pro of permutations and props}
Let $\Per$ be the bigraded subcollection of $\Map$ consisting in all
bijective maps. Hence, $\Per(p, q) = \emptyset$ when $p \ne q$. Since
the operations $*$ and $\circ$ of $\K \Angle{\Map}$ are stable in
$\K \Angle{\Per}$ and the identity maps $\Unit_p$, $p \in \N$,
are bijections, $\K \Angle{\Per}$ is a sub-pro of $\K \Angle{\Map}$.
It is moreover possible to show that the singleton
$\GeneratingSet := \{21\} \subseteq \Per(2, 2)$ is a minimal generating
set of $\K \Angle{\Per}$. We call $\K \Angle{\Per}$ the
\Def{pro of permutations}.
\medbreak

A \Def{prop} is a pro $\K \Angle{C}$ containing $\K \Angle{\Per}$ as a
sub-pro.
\medbreak

\part{Combinatorial operads}

\chapter{Enveloping operads of colored operads} \label{chap:enveloping}
The content of this chapter comes from~\cite{CG14} and is a joint work
with Frédéric Chapoton. We include here some new results that do not
appear in the aforementioned publication like the Koszulity of some of
the constructed operads.
\medbreak

\section*{Introduction}
In~\cite{Cha07}, Chapoton considered a ns operad structure on the
objects called noncrossing trees and noncrossing plants. These objects
can be depicted as simple graphs inside regular polygons, and are some
kinds of noncrossing configurations that are well-known combinatorial
objects~\cite{FN99,FS09}. The partial compositions of these ns operads
have very simple graphical descriptions and it is tempting and easy to
generalize this composition as much as possible, by removing some
constraints on the objects. This leads to a very big ns operad of
noncrossing configurations. This research initially started as a study
of this ns operad, with possible aim the description of its suboperads.
\smallbreak

This study has led us to the following results. First, we introduce a
general functorial construction from ns colored operads to ns operads,
which is called the enveloping operad. This can be compared to the
amalgamated product of groups, in the sense that it takes a compound
object to build a unified object in the simplest possible way, by
imposing as few relations as possible. The main interest of this
construction relies on the fact that a lot of properties of an
enveloping operad (as {\em e.g.}, its Hilbert series and a presentation
by generators and relations) can be obtained from its underlying ns
colored operad.
\smallbreak

Next, we consider the ns operad $\BNC$ of bicolored noncrossing
configurations, defined by a simple graphical composition, and show
that it admits a description as the enveloping operad of a very simple
ns colored operad on two colors called $\Bubble$. We also obtain a
presentation by generators and relations of the ns colored
operad~$\Bubble$.
\smallbreak

Then this understanding of the operad $\BNC$ is used to describe in
details some of its suboperads, namely those generated by two chosen
generators among the binary generators of $\BNC$. This already gives an
interesting family of operads, where one can recognize some known ones:
the operad of based noncrossing trees $\NCT$~\cite{Cha07}, the operad of
noncrossing plants $\NCP$~\cite{Cha07}, the dipterous
operad~\cite{LR03,Zin12}, and the $2$-associative
operad~\cite{LR06,Zin12}. Our main results here are a presentation by
generators and relations for all these suboperads except one, and also
the description of all the generating series. It should be noted that
the presentations are obtained in a case-by-case fashion, using similar
techniques involving rewrite rules on syntax trees (see
Section~\ref{subsec:rewrite_rules_syntax_trees} of
Chapter~\ref{chap:combinatorics}).
\smallbreak

This chapter is organized as follows. In
Section~\ref{sec:enveloping_operads}, the general construction of
enveloping operads is given and its properties described. Next, in
Section~\ref{sec:operad_BNC}, we introduce the operad $\BNC$ and prove
that this operad is isomorphic to an enveloping operad. Finally, in
Section~\ref{sec:suboperads_BNC}, several suboperads of $\BNC$ are
considered, in a more or less detailed way.
\medbreak

\subsubsection*{Note}
This chapter deals only with ns set-operads and ns colored set-operads.
For this reason, ``operad'' means ``ns set-operad''. Moreover, we
consider only colored operads $\Cca$ such that $\Cca(1)$ is trivial,
that is $\Cca(1) = \{\Unit_c : c \in [k]\}$. Moreover, all considered
colored operads are $k$-colored operads (see
Section~\ref{subsubsec:colored_operads} of Chapter~\ref{chap:algebra}).
\medbreak

\section{Enveloping operads of colored operads}
\label{sec:enveloping_operads}
The aim of this section is twofold. We begin by introducing the main
object of this chapter: the construction which associates an operad with
a colored one, namely its enveloping operad. We finally justify the
benefits of seeing an operad $\Oca$ as an enveloping operad of a colored
one $\Cca$ by reviewing some properties of $\Oca$ that can be deduced
from the ones of~$\Cca$.
\medbreak

\subsection{The construction}
Let us now introduce the construction associating a (noncolored) operad
with a colored one. We begin by giving the formal definition of what
enveloping operads of colored operads are, and then, give a
combinatorial interpretation of the construction in terms of
anticolored syntax trees.
\medbreak

\subsubsection{Enveloping operads}
Let $\Cca$ be a $k$-colored operad. Recall that $\Augmentation(\Cca)$
is the set $\Cca \setminus \Cca(1)$. The \Def{enveloping operad}
$\Hull(\Cca)$ of $\Cca$ is the smallest (noncolored) operad containing
$\Augmentation(\Cca)$. In other terms,
\begin{equation} \label{equ:enveloping_operad}
    \Hull(\Cca) :=
    \FreeOperad\left(\Augmentation(\Cca)\right)/_\equiv,
\end{equation}
where $\equiv$ is the smallest operad congruence of
$\FreeOperad(\Augmentation(\Cca))$ satisfying
\begin{equation}
    \Corolla{x} \circ_i \Corolla{y} \equiv \Corolla{x \circ_i y},
\end{equation}
for all $x, y \in \Augmentation(\Cca)$ such that $x \circ_i y$ are
well-defined in~$\Cca$. Observe that in~\eqref{equ:enveloping_operad},
$\FreeOperad(\Augmentation(\Cca))$ is the free noncolored operad
generated by $\Augmentation(\Cca)$, where $\Augmentation(\Cca)$ is here
a combinatorial graded collection whose input and output colors are
forgotten.
\medbreak

\subsubsection{Reductions}
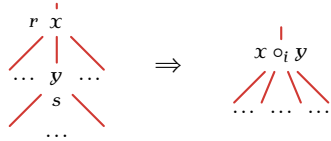
\begin{figure}[ht]
    \begin{equation*}
        \begin{tikzpicture}[scale=.25,Centering]
            \node[NodeST](1)at(0,0){\begin{math}x\end{math}};
            \node[NodeST](2)at(0,-3){\begin{math}y\end{math}};
            \node(1G)at(-3,-3){};
            \node(1D)at(3,-3){};
            \node(2G)at(-3,-6){};
            \node(2D)at(3,-6){};
            \node(v)at(0,1.5){};
            \draw[Edge](1)--(2);
            \draw[Edge](1)--(1G);
            \draw[Edge](1)--(1D);
            \draw[Edge](2)--(2G);
            \draw[Edge](2)--(2D);
            \draw[Edge](1)--(v);
            \node[font=\scriptsize]at(-1.75,-3)
                {\begin{math}\dots\end{math}};
            \node[font=\scriptsize]at(1.75,-3)
                {\begin{math}\dots\end{math}};
            \node[font=\scriptsize]at(0,-6)
                {\begin{math}\dots\end{math}};
            \node[left of=1,node distance=3mm,font=\scriptsize]
                {\begin{math}r\end{math}};
            \node[below of=2,node distance=3mm,font=\scriptsize]
                {\begin{math}s\end{math}};
        \end{tikzpicture}
    \quad \Rew \quad
    \begin{tikzpicture}[scale=.25,Centering]
        \node[NodeST](1)at(0,0){\begin{math}x \circ_i y\end{math}};
        \node(v)at(0,2){};
        \node(1G)at(-3,-3){};
        \node(1D)at(3,-3){};
        \node(2G)at(-1.25,-3){};
        \node(2D)at(1.25,-3){};
        \draw[Edge](1)--(v);
        \draw[Edge](1)--(1G);
        \draw[Edge](1)--(1D);
        \draw[Edge](1)--(2G);
        \draw[Edge](1)--(2D);
        \node[font=\scriptsize]at(-2,-3){\begin{math}\dots\end{math}};
        \node[font=\scriptsize]at(2,-3){\begin{math}\dots\end{math}};
        \node[font=\scriptsize]at(0,-3){\begin{math}\dots\end{math}};
    \end{tikzpicture}
    \end{equation*}
    \caption[Reduction of syntax trees on colored operads.]
    {The reduction of syntax trees. The internal node $s$ is the $i$th
    child of $r$.}
    \label{fig:reduction_1_colored_syntax_trees}
\end{figure}

Let $\Tfr$ be a syntax tree of $\FreeOperad(\Augmentation(\Cca))$
and $e$ be an edge of $\Tfr$ connecting two internal nodes $r$ and $s$
respectively labeled by $x$ and $y$, such that $s$ is the $i$th child of
$r$ and, as elements of $\Cca$, $\In_i(x) = \Out(y)$. Then, the
\Def{reduction} of $\Tfr$ with respect to $e$ is the tree obtained by
replacing $r$ and $s$ by an internal node labeled by $x \circ_i y$ (see
Figure~\ref{fig:reduction_1_colored_syntax_trees}). This tree is an
element of $\FreeOperad(\Augmentation(\Cca))$.
\medbreak

\subsubsection{Anticolored syntax trees}
Let $C$ be a $k$-colored collection. A \Def{$k$-anticolored syntax tree}
on $C$ (or for short, a \Def{$k$-anticolored $C$-syntax tree}) is a
(noncolored) $C$-syntax tree $\Tfr$ such that for any internal nodes $r$
and $s$ of $\Tfr$ such that $s$ is the $i$th child of $r$, we have
$\In_i(x) \ne \Out(y)$  where $x$ (resp. $y$) is the label of $r$ (resp.
$s$). The set of all $k$-anticolored $C$-syntax trees is denoted by
$\AntiCol(C)$. Observe that the leaf $\Leaf$ is a $k$-anticolored syntax
tree. Since anticolored syntax trees are particular syntax trees, the
usual terminology and tools about them applies (see
Section~\ref{sec:trees} of Chapter~\ref{chap:combinatorics}).
\medbreak

\subsubsection{The operad of anticolored syntax trees}
\begin{figure}[ht]
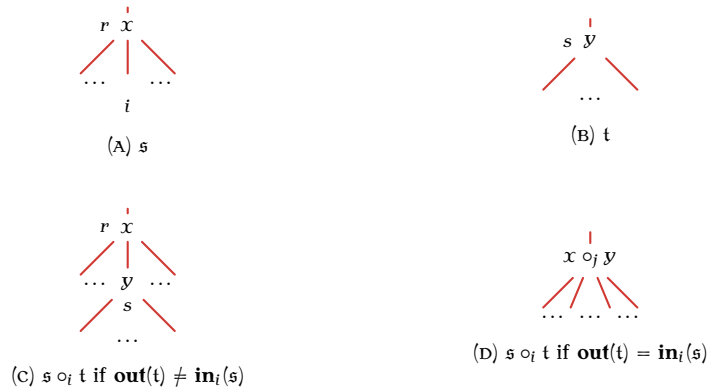

    \subfloat[][$\Sfr$]{
    \begin{minipage}[c]{.4\textwidth}
    \centering

    \end{minipage}
    \label{subfig:anticolored_same_colors}}
    \caption[Partial composition of anticolored syntax trees.]
    {The two cases for the partial composition of two anticolored trees
    $\Sfr$ and $\Tfr$.
    In~\protect\subref{subfig:anticolored_same_colors}, $j$ is
    the index of the $i$th leaf of $\Sfr$ among the children of the
    internal node~$r$.}
    \label{fig:composition_anticolored_trees}
\end{figure}

For any $k$-colored operad $\Cca$, the set
$\AntiCol(\Augmentation(\Cca))$ is endowed with an operad structure for
the partial composition defined as follows. Let $\Sfr$ and $\Tfr$ be two
anticolored syntax trees on $\Augmentation(\Cca)$. If
$\Out(\Tfr) \ne \In_i(\Sfr)$, $\Sfr \circ_i \Tfr$ is the anticolored
syntax tree obtained by grafting the root of $\Tfr$ on the $i$th leaf
of $\Sfr$. Otherwise, when $\Out(\Tfr) = \In_i(\Sfr)$,
$\Sfr \circ_i \Tfr$ is the anticolored syntax tree obtained by grafting
the root of $\Tfr$ on the $i$th leaf of $\Sfr$ and then, by reducing
the obtained tree with respect to the edge connecting the nodes $r$ and
$s$, where $r$ is the parent of the $i$th leaf of $\Sfr$ and $s$ is the
root of $\Tfr$ (see Figure~\ref{fig:composition_anticolored_trees}).
\medbreak

\begin{Proposition} \label{prop:interpretation_enveloping_operad}
    For any colored operad $\Cca$, the operads $\Hull(\Cca)$ and
    $\AntiCol(\Augmentation(\Cca))$ are isomorphic.
\end{Proposition}
\begin{proof}
    Let $\phi : \Hull(\Cca) \to \AntiCol(\Augmentation(\Cca))$ be the
    map associating with any $\equiv$-equivalence class of syntax trees
    on $\Augmentation(\Cca)$, the only anticolored
    $\Augmentation(\Cca)$-syntax tree on  belonging to it. To prove
    the statement, let us show that $\phi$ is a well-defined operad
    isomorphism.
    \smallbreak

    For that, let $\RewTrees$ be the closure of the rewrite relation
    $\Rew$ of reduction with respect to the partial compositions
    operations of trees. The axioms of operads ensure that $\RewTrees$
    is confluent, and since any rewriting decreases the degrees of the
    trees, $\RewTrees$ is terminating. The normal forms of $\RewTrees$
    are the trees that cannot be reduced, and thus, are anticolored
    $\Augmentation(\Cca)$-syntax trees. Since by definition of $\equiv$,
    $\Sfr \RewTrees \Tfr$ implies $\Sfr \equiv \Tfr$, the
    application $\phi$ is well-defined and is a bijection.
    \smallbreak

    Finally, let $[\Sfr]_\equiv, [\Tfr]_\equiv \in \Hull(\Cca)$,
    $\Sfr := \phi([\Sfr]_\equiv)$, and $\Tfr := \phi([\Tfr]_\equiv)$.
    The only anticolored syntax tree in $[\Sfr \circ_i \Tfr]_\equiv$ is
    obtained by grafting $\Sfr$ and $\Tfr$ together and performing, if
    possible, a reduction with respect to the edge linking them. Since
    the obtained tree is also the anticolored syntax tree
    $\Sfr \circ_i \Tfr$ of $\AntiCol(\Augmentation(\Cca))$, $\phi$ is
    an operad morphism.
\end{proof}
\medbreak

Proposition~\ref{prop:interpretation_enveloping_operad} implies that
the elements of $\Hull(\Cca)$ can be regarded as anticolored trees,
endowed with their partial composition defined above. We shall maintain
this point of view in the rest of this chapter by setting
$\Hull(\Cca) := \AntiCol(\Augmentation(\Cca))$.
\medbreak

\subsubsection{Functoriality}
Let $\Cca_1$ and $\Cca_2$ be two $k$-colored operads. Given
$\phi : \Cca_1 \to \Cca_2$ a colored operad morphism, let the map
$\Hull(\phi) : \Hull(\Cca_1) \to \Hull(\Cca_2)$ be the unique operad
morphism satisfying
\begin{equation}
    \Hull(\phi)\left(\Corolla{x}\right) := \Corolla{\phi(x)}
\end{equation}
for any $x \in \Cca_1$.
\medbreak

\begin{Theorem} \label{thm:functoriality_enveloping_operad}
    The construction $\Hull$ is a functor from the category of colored
    operads to the category of operads that preserves injections and
    surjections.
\end{Theorem}
\begin{proof}
    For any colored operad $\Cca$, $\Hull(\Cca)$ is by definition an
    operad on anticolored syntax trees on $\Augmentation(\Cca)$.
    Moreover, by induction on the number of internal nodes of the
    anticolored syntax trees, it follows that for any colored operad
    morphism $\phi$, $\Hull(\phi)$ is a well-defined operad morphism.
    \smallbreak

    Since $\Hull$ is compatible with map composition and sends the
    identity colored operad morphism to the identity operad morphism,
    $\Hull$ is a functor. It is moreover plain that if $\phi$ is an
    injective (resp. surjective) colored operad morphism, then
    $\Hull(\phi)$ is an injective (resp. surjective) operad morphism.
\end{proof}
\medbreak

Theorem~\ref{thm:functoriality_enveloping_operad} is rich in
consequences: Propositions~\ref{prop:suboperad_enveloping},
\ref{prop:generating_set_enveloping}, \ref{prop:presentation},
\ref{prop:group_symmetries_enveloping} of next section directly rely on
it.
\medbreak

Notice that $\Hull$ is a surjective functor. Indeed, since an
anticolored syntax tree on a $1$-colored collection is necessarily a
corolla, for any operad $\Oca$, $\Hull(\Oca)$ contains only corollas
labeled on $\Augmentation(\Oca)$ and it is therefore isomorphic
to~$\Oca$.
\medbreak

Notice also that $\Hull$ is not an injective functor. Let us exhibit two
$2$-colored operads not themselves isomorphic that produce by $\Hull$
two isomorphic operads. Let $\Cca$ be the $2$-colored operad where
$\Cca(2) := \{\Asf_2\}$ with $\Out(\Asf_2) := 1$ and
$\In_1(\Asf_2) := \In_2(\Asf_2) := 2$, and for all $n \geq 3$,
$\Cca(n) := \emptyset$. Due to the output and input colors of $\Asf_2$,
there is no nontrivial partial composition in $\Cca$. On the other
hand, let $\FAs$ be the $2$-colored operad where, for all $n \geq 2$,
$\FAs(n) := \{\Bsf_n\}$ with $\Out(\Bsf_n) := 1$, $\In_1(\Bsf_n) := 1$,
and $\In_i(\Bsf_n) := 2$ for all $2 \leq i \leq n$. Nontrivial partial
compositions of $\FAs$ are only defined for the first position by
\begin{math}
    \Bsf_n \circ_1 \Bsf_m := \Bsf_{n + m - 1},
\end{math}
for any $n, m \geq 2$. One observes that $\Hull(\Cca)$ and $\Hull(\FAs)$
are both the free operad generated by one element of arity $2$ with no
nontrivial relations, and hence, are isomorphic. The isomorphism
between $\Hull(\Cca)$ and $\Hull(\FAs)$ can be described by a left-child
right-sibling bijection~\cite{CLRS09} between binary trees and planar
rooted trees.
\medbreak

\subsubsection{Example}
Consider the $2$-colored operad $\FAs$ defined in the previous section.
The elements of $\Hull(\FAs)$ are anticolored syntax trees on
$\Augmentation(\FAs)$. Because of the output and input colors of the
elements of $\FAs$, $\Hull(\FAs)$ contains trees where all internal
nodes have no child in the first position. For instance,
\begin{equation}

\end{equation}
is a partial composition in $\Hull(\FAs)$ which does not require any
reduction. On the other hand,
\begin{equation}
    %
\end{equation}
is a partial composition requiring a reduction.
\medbreak

\subsection{Bubble decompositions of operads and consequences}
\label{subsec:consequences_bubble_decomposition}
Let $\Oca$ be an operad. We say that $\Cca$ is a
\Def{$k$-bubble decomposition} of $\Oca$ if $\Cca$ is a $k$-colored
operad such that $\Hull(\Cca)$ and $\Oca$ are isomorphic. In this case,
we say that the elements of $\Cca$ are \Def{bubbles}. As we shall show,
since a bubble decomposition $\Cca$ of an operad $\Oca$ contains a lot
of information about $\Oca$, the study of $\Oca$ can be confined to the
study of $\Cca$. Then, the main interest of the construction $\Hull$ is
here: the study of an operad $\Oca$ is confined to the study of one of
its bubble decompositions. Since colored operads are more constrained
structures than operads, this study is in most cases simpler than the
direct study of the operad itself.
\medbreak

\subsubsection{Hilbert series}
Let $c \in [k]$ be a color. The \Def{$c$-colored Hilbert series}
of $\Cca$ is the series
\begin{equation}
    \BubbleSeries_c(\VarZ_1, \dots, \VarZ_k) :=
    \VarX_c^{-1} \;
    \GenSeries_{\Augmentation(\Cca)}(0, \dots, 0, \VarX_c, 0, \dots, 0,
    \VarZ_1, \dots, \VarZ_k),
\end{equation}
where $\GenSeries_{\Augmentation(\Cca)}$ is the generating series of the
colored collection $\Augmentation(\Cca)$ (see
Section~\ref{subsubsec:colored_collections} of
Chapter~\ref{chap:combinatorics}). In more concrete terms,
$\BubbleSeries_c(\VarZ_1, \dots, \VarZ_k)$ is the series wherein the
coefficient of $\VarZ_1^{\alpha_1} \dots \VarZ_k^{\alpha_k}$ counts the
nontrivial elements of $\Cca$ having $c$ as output color and $\alpha_a$
inputs of color $a$ for all $a \in [k]$.
\medbreak

\begin{Proposition} \label{prop:hilbert_series_enveloping}
    Let $\Cca$ be a $k$-colored operad. Then, the Hilbert series
    $\HilbSeries(t)$ of the enveloping operad of $\Cca$ satisfies
    \begin{equation}
        \HilbSeries(t) =
        t + \HilbSeries_1(t) + \dots + \HilbSeries_k(t),
    \end{equation}
    where for all $c \in [k]$, the series $\HilbSeries_c(t)$ satisfy
    \begin{equation*}
        \HilbSeries_c(t) =
        \BubbleSeries_c(\HilbSeries(t) - \HilbSeries_1(t), \dots,
            \HilbSeries(t) - \HilbSeries_k(t)).
    \end{equation*}
\end{Proposition}
\medbreak

Note that Proposition~\ref{prop:hilbert_series_enveloping} implies
that, if the colored Hilbert series of $\Cca$ is algebraic, the Hilbert
series of $\Hull(\Cca)$ also is. Nevertheless, as we shall see,
rationality is not preserved.
\medbreak

\subsubsection{Suboperads and quotients}
\begin{Proposition} \label{prop:suboperad_enveloping}
    Let $\Cca$ be a colored operad and $\Cca'$ be one of its colored
    suboperads (resp. quotients). Then, the enveloping operad of $\Cca'$
    is a suboperad (resp. quotient) of the enveloping operad of $\Cca$.
\end{Proposition}
\medbreak

\subsubsection{Generating sets}
\begin{Proposition} \label{prop:generating_set_enveloping}
    Let $\Cca$ be a colored operad admitting $\GeneratingSet$
    as a generating set. Then, the enveloping operad of $\Cca$ is
    generated by
    \begin{equation}
        \Hull(\GeneratingSet) :=
        \{\Corolla{g} : g \in \GeneratingSet\}.
    \end{equation}
\end{Proposition}
\medbreak

\subsubsection{Symmetries}
\begin{Proposition} \label{prop:group_symmetries_enveloping}
    Let $\Cca$ be a colored operad and $\Gca$ its group of symmetries.
    Then, the group of symmetries of the enveloping operad of $\Cca$ is
    $\Hull(\Gca)$ where
    \begin{equation}
        \Hull(\Gca) := \{\Hull(\phi) : \phi \in \Gca\}.
    \end{equation}
\end{Proposition}
\medbreak

\subsubsection{Presentations by generators and relations}
\begin{Proposition} \label{prop:presentation}
    Let $\Cca$ be a colored operad admitting a presentation
    $(\GeneratingSet, \Equiv)$.
    Then, the enveloping operad of $\Cca$ admits the presentation
    $(\Hull(\GeneratingSet), \Equiv')$, where
    $\Equiv'$ is the equivalence relation satisfying
    \begin{equation}
        \Sfr' \Equiv' \Tfr'
        \quad \mbox{if and only if} \quad
        \Sfr \Equiv \Tfr,
    \end{equation}
    where $\Sfr'$ (resp. $\Tfr'$) is the
    $\Hull(\GeneratingSet)$-colored syntax tree obtained by replacing
    any node labeled by $x$ of $\Sfr$ (resp. $\Tfr$) by~$\Corolla{x}$.
\end{Proposition}
\medbreak

\section{The operad of bicolored noncrossing configurations}
\label{sec:operad_BNC}
In this section, we shall define an operad over a new kind of
noncrossing configurations. In order to study it and apply the results
of Section~\ref{sec:enveloping_operads}, we shall see this operad as an
enveloping operad of a colored one.
\medbreak

\subsection{Bicolored noncrossing configurations}
Let us start by introducing our new combinatorial object, some of its
properties, and its operad structure.
\medbreak

\subsubsection{Bicolored noncrossing configurations}
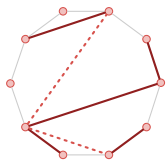
\begin{figure}[ht]
    \centering
    \begin{tikzpicture}[Centering]
        \node[CliquePoint](0)at(-0.3,-0.95){};
        \node[CliquePoint](1)at(-0.8,-0.58){};
        \node[CliquePoint](2)at(-1.,-0.){};
        \node[CliquePoint](3)at(-0.8,0.59){};
        \node[CliquePoint](4)at(-0.3,0.96){};
        \node[CliquePoint](5)at(0.31,0.96){};
        \node[CliquePoint](6)at(0.81,0.59){};
        \node[CliquePoint](7)at(1.,0.01){};
        \node[CliquePoint](8)at(0.81,-0.58){};
        \node[CliquePoint](9)at(0.31,-0.95){};
        \draw[PolyEdgeGray](0)--(1);
        \draw[PolyEdgeGray](1)--(2);
        \draw[PolyEdgeGray](2)--(3);
        \draw[PolyEdgeGray](3)--(4);
        \draw[PolyEdgeGray](4)--(5);
        \draw[PolyEdgeGray](5)--(6);
        \draw[PolyEdgeGray](6)--(7);
        \draw[PolyEdgeGray](7)--(8);
        \draw[PolyEdgeGray](8)--(9);
        \draw[PolyEdgeGray](9)--(0);
        \draw[PolyEdgeBlue](0)--(1);
        \draw[PolyEdgeBlue](1)--(7);
        \draw[PolyEdgeBlue](3)--(5);
        \draw[PolyEdgeBlue](6)--(7);
        \draw[PolyEdgeBlue](8)--(9);
        \draw[PolyEdgeRed](1)--(5);
        \draw[PolyEdgeRed](1)--(9);
    \end{tikzpicture}
    \caption[A bicolored noncrossing configuration.]
    {A nonbased BNC of size $9$. Blue arcs are $(1, 2)$, $(2, 8)$,
    $(4, 6)$, $(7, 8)$, and $(9, 10)$, and red arcs are $(2, 6)$ and
    $(2, 10)$. All other arcs are uncolored. The border of this
    BNC is $211111212$.}
    \label{fig:example_BNC}
\end{figure}

A \Def{bicolored noncrossing configuration} (or, for short, a \Def{BNC})
is a noncrossing configuration (see
Section~\ref{subsubsec:configurations} of
Chapter~\ref{chap:combinatorics}) where each arc can be either
\Def{blue} (drawn as a thick line) or \Def{red} (drawn as a dotted line)
and such that all red arcs are diagonals. We say that $\Cfr$ is
\Def{based} if its base is blue and \Def{nonbased} otherwise. Besides,
we impose by definition that there is only one BNC of size $1$: the
segment consisting in one blue arc. Figure~\ref{fig:example_BNC} shows a
BNC.
\medbreak

When the size of $\Cfr$ is not smaller than $2$, the \Def{border} of
$\Cfr$ is the word $\Border(\Cfr)$ of length~$n$ such that, for any
$i \in [n]$, $\Border(\Cfr)_i := 1$ if the $i$th edge of $\Cfr$ is
uncolored and $\Border(\Cfr)_i := 2$ otherwise. See
Figure~\ref{fig:example_BNC} for an example.
\medbreak

\subsubsection{Operad structure}
From now on, the \Def{arity} $|\Cfr|$ of a BNC $\Cfr$ is its size. Let
$\Cfr$ and $\Dfr$ be two BNCs of respective arities $n$ and $m$, and
$i \in [n]$. The partial composition $\Cfr \circ_i \Dfr =: \Efr$ is
obtained by gluing the base of $\Dfr$ onto the $i$th edge of $\Cfr$,
and then,
\begin{enumerate}
    \item if the base of $\Dfr$ and the $i$th edge of $\Cfr$ are both
    uncolored, the arc $(i, i + m)$ of $\Efr$ becomes red;
    \item if the base of $\Dfr$ and the $i$th edge of $\Cfr$ are both
    blue, the arc $(i, i + m)$ of $\Efr$ becomes blue;
    \item otherwise, the base of $\Dfr$ and the $i$th edge of $\Cfr$
    have different colors; in this case, the arc $(i, i + m)$ of $\Efr$
    is uncolored.
\end{enumerate}
For aesthetic reasons, the resulting shape is reshaped to form a regular
polygon. For instance,
\begin{subequations}
\begin{equation}
\,,
\end{equation}
\end{subequations}
are partial compositions in~$\BNC$.
\medbreak

\begin{Proposition} \label{prop:operad_BNC}
    The set of all the BNCs, together with the partial composition maps
    $\circ_i$ and the BNC of arity $1$ as unit form an operad, denoted
    by~$\BNC$.
\end{Proposition}
\medbreak

\subsection{The colored operad of bubbles} \label{subsec:operade_Bulle}
We now define a colored operad involving particular BNCs and perform
a complete study of it.
\medbreak

\subsubsection{Bubbles}
A \Def{bubble} is a BNC of size no smaller than $2$ with no diagonal
(hence the name). For instance,
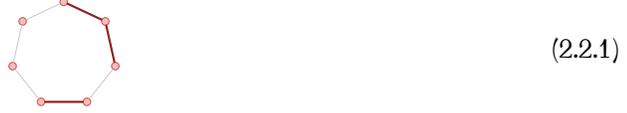
\begin{equation}
    \begin{tikzpicture}[scale=.7,Centering]
        \node[CliquePoint](0)at(-0.43,-0.9){};
        \node[CliquePoint](1)at(-0.97,-0.22){};
        \node[CliquePoint](2)at(-0.78,0.63){};
        \node[CliquePoint](3)at(-0.,1.){};
        \node[CliquePoint](4)at(0.79,0.63){};
        \node[CliquePoint](5)at(0.98,-0.22){};
        \node[CliquePoint](6)at(0.44,-0.9){};
        \draw[PolyEdgeGray](0)--(1);
        \draw[PolyEdgeGray](1)--(2);
        \draw[PolyEdgeGray](2)--(3);
        \draw[PolyEdgeGray](3)--(4);
        \draw[PolyEdgeGray](4)--(5);
        \draw[PolyEdgeGray](5)--(6);
        \draw[PolyEdgeGray](6)--(0);
        \draw[PolyEdgeBlue](3)--(4);
        \draw[PolyEdgeBlue](4)--(5);
        \draw[PolyEdgeBlue](0)--(6);
    \end{tikzpicture}
\end{equation}
is a bubble of size $6$ and whose border is $111221$.
\medbreak

\subsubsection{Colored operad structure}
Let $\Bfr$ be a bubble of arity $n$. Let us assign input and output
colors to $\Bfr$ in the following way. The output color $\Out(\Bfr)$ of
$\Bfr$ is $1$ if $\Bfr$ is based and $2$ otherwise, and the color
$\In_i(\Bfr)$ of the $i$th input of $\Bfr$ is the $i$th letter of the
border of~$\Bfr$.
\medbreak

Let us denote by $\Unit_1$ and $\Unit_2$ two virtual bubbles of arity
$1$ such that $\Out(\Unit_1) := \In_1(\Unit_1) := 1$ and
$\Out(\Unit_2) := \In_1(\Unit_2) := 2$.
\medbreak

\begin{Proposition} \label{prop:operad_Bubble}
    The set of all the bubbles, together with the partial composition
    map $\circ_i$ of $\BNC$ and the units $\Unit_1$ and $\Unit_2$ form a
    $2$-colored operad, denoted by~$\Bubble$.
\end{Proposition}
\medbreak

Notice that any bubble $\Bfr$ is wholly encoded by the pair
\begin{math}
    \left(\Out(\Bfr), \left(\In_i(\Bfr)\right)_{i \in |\Bfr|}\right)
\end{math}.
Therefore, $\Bubble$ is a very simple colored operad: for any $n$, the
set of elements of arity $n \geq 2$ is $[2] \times [2]^n$ and the
partial composition, when defined, is a substitution in words.
For instance, the partial composition
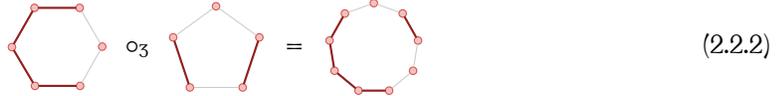
\begin{equation}

\end{equation}
of $\Bubble$ can be expressed concisely as
\begin{equation}
    (1, 22\textcolor{Col1}{2}11)
    \circ_3
    (\textcolor{Col1}{2}, \textcolor{Col4}{{\bf 2112}})
    =
    (1, 22\textcolor{Col4}{{\bf 2112}}11).
\end{equation}
\medbreak

\subsubsection{Colored Hilbert series} \label{sec:series_bubbles_Bubble}
Since $\Bubble$ contains by definition all the bubbles, the colored
Hilbert series of $\Bubble$ satisfies
\begin{equation} \label{equ:series_bubbles_Bubble}
    \BubbleSeries_1(\VarZ_1, \VarZ_2) =
    \BubbleSeries_2(\VarZ_1, \VarZ_2) =
    \sum_{n \geq 2} (\VarZ_1 + \VarZ_2)^n =
    \frac{(\VarZ_1 + \VarZ_2)^2}{1 - \VarZ_1 - \VarZ_2}.
\end{equation}
\medbreak

\subsubsection{Generating set}
\begin{Proposition} \label{prop:generators_Bubble}
    The set
    \begin{equation}
        \GeneratingSet_\Bubble :=
        \left\{\AAA, \AAB, \BAA, \BAB, \ABA, \ABB, \BBA, \BBB\right\}
    \end{equation}
    of bubbles of arity $2$ is the unique minimal generating set
    of~$\Bubble$.
\end{Proposition}
\medbreak

\subsubsection{Symmetries} \label{subsubsec:symmetries_Bubble}
The \Def{complementary} $\Cpl(\Bfr)$ of a bubble $\Bfr$ is the bubble
obtained by swapping the colors of the edges of $\Bfr$ (blue edges
become uncolored and conversely). The \Def{returned} $\Ret(\Bfr)$ of
$\Bfr$ is the bubble obtained by applying to $\Bfr$ the reflection
through the vertical line passing by its base.
Figure~\ref{fig:examples_symmetries} shows examples of these symmetries.
\begin{figure}[ht]
    \subfloat[][Complementary.]{
    \begin{minipage}[c]{.4\textwidth}
    \begin{equation*}

    \end{equation*}
    \end{minipage}
    \label{subfig:returned}}
    \caption{The complementary and the returned of a bubble.}
    \label{fig:examples_symmetries}
\end{figure}
\medbreak

\begin{Proposition} \label{prop:symmetries_Bubble}
    The group of symmetries of $\Bubble$ is generated by $\Cpl$ and
    $\Ret$ and satisfies the relations
    \begin{equation} \label{equ:relations_group_symmetries_Bubble}
        \Ret = \Ret^{-1}, \quad
        \Cpl = \Cpl^{-1}, \quad
        \Ret\, \Cpl = \Cpl\, \Ret.
    \end{equation}
\end{Proposition}
\medbreak

\subsubsection{Presentation by generators and relations}

\begin{Theorem} \label{thm:presentation_Bubble}
    The $2$-colored operad $\Bubble$ admits the presentation
    $\left(\GeneratingSet_\Bubble, \Equiv\right)$ where
    $\Equiv$ is the equivalence relation satisfying
    \begin{subequations}
    \begin{equation} \label{equ:presentation_Bubble_1}
        \Corolla{\BBB} \circ_1 \Corolla{\BAB}
        \enspace \Equiv \enspace
        \Corolla{\ABB} \circ_1 \Corolla{\BBB}
        \enspace \Equiv \enspace
        \Corolla{\BBB} \circ_2 \Corolla{\BAB}
        \enspace \Equiv \enspace
        \Corolla{\BBA} \circ_2 \Corolla{\BBB},
    \end{equation}
    \begin{equation}
        \Corolla{\BBB} \circ_1 \Corolla{\BAA}
        \enspace \Equiv \enspace
        \Corolla{\ABB} \circ_1 \Corolla{\BBA}
        \enspace \Equiv \enspace
        \Corolla{\BBB} \circ_2 \Corolla{\AAB}
        \enspace \Equiv \enspace
        \Corolla{\BBA} \circ_2 \Corolla{\ABB},
    \end{equation}
    \begin{equation}
        \Corolla{\BBA} \circ_1 \Corolla{\BAA}
        \enspace \Equiv \enspace
        \Corolla{\ABA} \circ_1 \Corolla{\BBA}
        \enspace \Equiv \enspace
        \Corolla{\BBA} \circ_2 \Corolla{\ABA}
        \enspace \Equiv \enspace
        \Corolla{\BBB} \circ_2 \Corolla{\AAA},
    \end{equation}
    \begin{equation}
        \Corolla{\BBA} \circ_1 \Corolla{\BAB}
        \enspace \Equiv \enspace
        \Corolla{\ABA} \circ_1 \Corolla{\BBB}
        \enspace \Equiv \enspace
        \Corolla{\BBA} \circ_2 \Corolla{\BBA}
        \enspace \Equiv \enspace
        \Corolla{\BBB} \circ_2 \Corolla{\BAA},
    \end{equation}
    \begin{equation}
        \Corolla{\ABB} \circ_1 \Corolla{\ABB}
        \enspace \Equiv \enspace
        \Corolla{\BBB} \circ_1 \Corolla{\AAB}
        \enspace \Equiv \enspace
        \Corolla{\ABB} \circ_2 \Corolla{\BAB}
        \enspace \Equiv \enspace
        \Corolla{\ABA} \circ_2 \Corolla{\BBB},
    \end{equation}
    \begin{equation}
        \Corolla{\ABB} \circ_1 \Corolla{\ABA}
        \enspace \Equiv \enspace
        \Corolla{\BBB} \circ_1 \Corolla{\AAA}
        \enspace \Equiv \enspace
        \Corolla{\ABB} \circ_2 \Corolla{\AAB}
        \enspace \Equiv \enspace
        \Corolla{\ABA} \circ_2 \Corolla{\ABB},
    \end{equation}
    \begin{equation}
        \Corolla{\ABA} \circ_1 \Corolla{\ABA}
        \enspace \Equiv \enspace
        \Corolla{\BBA} \circ_1 \Corolla{\AAA}
        \enspace \Equiv \enspace
        \Corolla{\ABB} \circ_2 \Corolla{\AAA}
        \enspace \Equiv \enspace
        \Corolla{\ABA} \circ_2 \Corolla{\ABA},
    \end{equation}
    \begin{equation}
        \Corolla{\ABA} \circ_1 \Corolla{\ABB}
        \enspace \Equiv \enspace
        \Corolla{\BBA} \circ_1 \Corolla{\AAB}
        \enspace \Equiv \enspace
        \Corolla{\ABB} \circ_2 \Corolla{\BAA}
        \enspace \Equiv \enspace
        \Corolla{\ABA} \circ_2 \Corolla{\BBA},
    \end{equation}
    \begin{equation}
        \Corolla{\BAB} \circ_1 \Corolla{\BAB}
        \enspace \Equiv \enspace
        \Corolla{\AAB} \circ_1 \Corolla{\BBB}
        \enspace \Equiv \enspace
        \Corolla{\BAB} \circ_2 \Corolla{\BAB}
        \enspace \Equiv \enspace
        \Corolla{\BAA} \circ_2 \Corolla{\BBB},
    \end{equation}
    \begin{equation}
        \Corolla{\BAB} \circ_1 \Corolla{\BAA}
        \enspace \Equiv \enspace
        \Corolla{\AAB} \circ_1 \Corolla{\BBA}
        \enspace \Equiv \enspace
        \Corolla{\BAB} \circ_2 \Corolla{\AAB}
        \enspace \Equiv \enspace
        \Corolla{\BAA} \circ_2 \Corolla{\ABB},
    \end{equation}
    \begin{equation}
        \Corolla{\BAA} \circ_1 \Corolla{\BAA}
        \enspace \Equiv \enspace
        \Corolla{\AAA} \circ_1 \Corolla{\BBA}
        \enspace \Equiv \enspace
        \Corolla{\BAA} \circ_2 \Corolla{\ABA}
        \enspace \Equiv \enspace
        \Corolla{\BAB} \circ_2 \Corolla{\AAA},
    \end{equation}
    \begin{equation}
        \Corolla{\BAA} \circ_1 \Corolla{\BAB}
        \enspace \Equiv \enspace
        \Corolla{\AAA} \circ_1 \Corolla{\BBB}
        \enspace \Equiv \enspace
        \Corolla{\BAA} \circ_2 \Corolla{\BBA}
        \enspace \Equiv \enspace
        \Corolla{\BAB} \circ_2 \Corolla{\BAA},
    \end{equation}
    \begin{equation}
        \Corolla{\AAB} \circ_1 \Corolla{\ABB}
        \enspace \Equiv \enspace
        \Corolla{\BAB} \circ_1 \Corolla{\AAB}
        \enspace \Equiv \enspace
        \Corolla{\AAB} \circ_2 \Corolla{\BAB}
        \enspace \Equiv \enspace
        \Corolla{\AAA} \circ_2 \Corolla{\BBB},
    \end{equation}
    \begin{equation}
        \Corolla{\AAB} \circ_1 \Corolla{\ABA}
        \enspace \Equiv \enspace
        \Corolla{\BAB} \circ_1 \Corolla{\AAA}
        \enspace \Equiv \enspace
        \Corolla{\AAB} \circ_2 \Corolla{\AAB}
        \enspace \Equiv \enspace
        \Corolla{\AAA} \circ_2 \Corolla{\ABB},
    \end{equation}
    \begin{equation}
        \Corolla{\AAA} \circ_1 \Corolla{\ABA}
        \enspace \Equiv \enspace
        \Corolla{\BAA} \circ_1 \Corolla{\AAA}
        \enspace \Equiv \enspace
        \Corolla{\AAB} \circ_2 \Corolla{\AAA}
        \enspace \Equiv \enspace
        \Corolla{\AAA} \circ_2 \Corolla{\ABA},
    \end{equation}
    \begin{equation} \label{equ:presentation_Bubble_16}
        \Corolla{\AAA} \circ_1 \Corolla{\ABB}
        \enspace \Equiv \enspace
        \Corolla{\BAA} \circ_1 \Corolla{\AAB}
        \enspace \Equiv \enspace
        \Corolla{\AAB} \circ_2 \Corolla{\BAA}
        \enspace \Equiv \enspace
        \Corolla{\AAA} \circ_2 \Corolla{\BBA}.
    \end{equation}
    \end{subequations}
    Moreover, the set of the colored $\GeneratingSet_\Bubble$-syntax
    trees avoiding the trees appearing as second, third or fourth
    members of Relations~\eqref{equ:presentation_Bubble_1}---%
    \eqref{equ:presentation_Bubble_16}
    is a Poincaré-Birkhoff-Witt basis of~$\Bubble$.
\end{Theorem}
\begin{proof}
    To prove the presentation of the statement, we shall show that there
    exists a colored operad isomorphism
    \begin{math}
        \phi :
        \FreeColoredOperad\left(\GeneratingSet_\Bubble\right)/_\equiv
        \to \Bubble
    \end{math}
    where $\equiv$ is the operad congruence generated by
    $\Equiv$.
    \smallbreak

    Let us set $\phi\left([\Corolla{g}]_\equiv\right) := g$ for any $g$
    of $\GeneratingSet_\Bubble$. We observe that for any relation
    \begin{math}
        \Corolla{x} \circ_i \Corolla{y}
        \Equiv
        \Corolla{z} \circ_j \Corolla{t}
    \end{math}
    of the statement, we have $x \circ_i y = z \circ_j t$. It then
    follows that $\phi$ can be uniquely extended into a colored operad
    morphism. Moreover, since the image of $\phi$ contains all the
    generators of $\Bubble$, $\phi$ is surjective.
    \smallbreak

    Let us now prove that $\phi$ is a bijection. For that, let us orient
    the relation $\Equiv$ by means of the rewrite rule $\Rew$
    on the colored syntax trees on $\GeneratingSet_\Bubble$ satisfying
    $\Sfr \Rew \Tfr$ if $\Sfr \Equiv \Tfr$ and $\Tfr$ is one
    of the  following sixteen target trees
    \begin{equation} \label{equ:target_trees_Bubble}
        \begin{split}
            \Corolla{\BBB} \circ_1 \Corolla{\BAB}, \enspace
            \Corolla{\BBB} \circ_1 \Corolla{\BAA}, \enspace
            \Corolla{\BBA} \circ_1 \Corolla{\BAA}, \enspace
            \Corolla{\BBA} \circ_1 \Corolla{\BAB}, \\
            \Corolla{\ABB} \circ_1 \Corolla{\ABB}, \enspace
            \Corolla{\ABB} \circ_1 \Corolla{\ABA}, \enspace
            \Corolla{\ABA} \circ_1 \Corolla{\ABA}, \enspace
            \Corolla{\ABA} \circ_1 \Corolla{\ABB}, \\
            \Corolla{\BAB} \circ_1 \Corolla{\BAB}, \enspace
            \Corolla{\BAB} \circ_1 \Corolla{\BAA}, \enspace
            \Corolla{\BAA} \circ_1 \Corolla{\BAA}, \enspace
            \Corolla{\BAA} \circ_1 \Corolla{\BAB}, \\
            \Corolla{\AAB} \circ_1 \Corolla{\ABB}, \enspace
            \Corolla{\AAB} \circ_1 \Corolla{\ABA}, \enspace
            \Corolla{\AAA} \circ_1 \Corolla{\ABA}, \enspace
            \Corolla{\AAA} \circ_1 \Corolla{\ABB}.
        \end{split}
    \end{equation}
    The target trees of $\Rew$ are the only left comb trees appearing
    in each $\Equiv$-equivalence class of the statement such
    that the color of the first input of the root is the same
    as the color of the first input of its child.
    \smallbreak

    Let $\RewTrees$ be the closure of $\Rew$ and let us prove that
    $\RewTrees$ is terminating. Let $\psi$ be the map associating the
    pair $(\TamariInvariant(\Tfr), \LeftInvariant(\Tfr))$ with a colored
    syntax tree $\Tfr$, where $\TamariInvariant$ is defined
    by~\eqref{equ:Tamari_invariant} in Chapter~\ref{chap:combinatorics},
    and $\LeftInvariant(\Tfr)$ is the number of internal nodes labeled
    by $x$ of $\Tfr$ having an internal node labeled by $y$ as left
    child such that $\In_1(x) \ne \In_1(y)$. We observe that, for any
    trees $\Tfr_0$ and $\Tfr_1$ such that $\Tfr_0 \RewTrees \Tfr_1$,
    $\psi(\Tfr_1)$ is lexicographically smaller than $\psi(\Tfr_0)$.
    Hence, $\RewTrees$ is terminating.
    \smallbreak

    The normal forms of $\RewTrees$ are the colored
    $\GeneratingSet_\Bubble$-syntax trees avoiding the trees $\Sfr$
    appearing as a left members of $\Rew$. These are left comb trees
    $\Tfr$ such that for all internal nodes $x$ and $y$ of $\Tfr$,
    $\In_1(x) = \In_1(y)$. Pictorially, $\Tfr$ is of the form
    \begin{equation}
        \Tfr = \enspace
        \begin{tikzpicture}[xscale=.45,yscale=.6,Centering]
            \node(v)at(0,1.25){};
            \node[NodeST](0)at(0,0){\begin{math}x_{n-1}\end{math}};
            \node[NodeST](1)at(-2,-2){\begin{math}x_1\end{math}};
            \node[Leaf](2)at(-3.5,-3.5){};
            \node[Leaf](3)at(-.5,-3.5){};
            \node[Leaf](4)at(1.5,-1.5){};
            \draw[Edge](v)edge node[EdgeLabel]
                {\begin{math}c\end{math}}(0);
            \draw[Edge,dashed](0)edge node[EdgeLabel]
                {\begin{math}d_1\end{math}}(1);
            \draw[Edge](1)edge node[EdgeLabel]
                {\begin{math}d_1\end{math}}(2);
            \draw[Edge](1)edge node[EdgeLabel,right]
                {\begin{math}d_2\end{math}}(3);
            \draw[Edge](0)edge node[EdgeLabel,right]
                {\begin{math}d_n\end{math}}(4);
        \end{tikzpicture}\,,
    \end{equation}
    where $c \in [2]$, $d_i \in [2]$ for all $i \in [n]$, and
    $x_j \in \GeneratingSet_\Bubble$ for all $j \in [n - 1]$. Since
    $\Tfr$ is a colored syntax tree, given $c$ and the $d_i$, there is
    exactly one possibility for all the $x_j$. Therefore, there are
    $f_c(n) := 2^n$ normal forms of $\RewTrees$ of arity $n$ with $c$
    as output color. This imply that
    $\FreeColoredOperad\left(\GeneratingSet_\Bubble\right)/_\equiv$
    contains at most $f_c(n)$ elements of arity $n$ and $c$ as output
    color. Then, since $f_c(n)$ is also the number of elements of
    $\Bubble$ with arity $n$ and $c$ as output color (see
    Section~\ref{sec:series_bubbles_Bubble}), $\phi$ is a bijection.
    \smallbreak

    Finally, since $\Rew$ is an orientation of $\Equiv$ and
    the normal forms of $\RewTrees$ are the colored
    $\GeneratingSet_\Bubble$-syntax trees avoiding the trees appearing
    as second, third or fourth members of
    Relations~\eqref{equ:presentation_Bubble_1}---%
    \eqref{equ:presentation_Bubble_16}
    (see~\eqref{equ:target_trees_Bubble}), the last part of the
    statement follows.
\end{proof}
\medbreak

\subsection{Properties of the operad of bicolored noncrossing
    configurations}
Let us come back to the study of the operad $\BNC$. We show here that
$\BNC$ is the enveloping operad of $\Bubble$ and then, by using the
results of Section~\ref{subsec:operade_Bulle} together with the ones of
Section~\ref{subsec:consequences_bubble_decomposition}, give some of its
properties.
\medbreak

\subsubsection{Bubble decomposition}
Let $\Cfr$ be a BNC. An \Def{area} of $\Cfr$ is a maximal component of
$\Cfr$ without colored diagonals and bounded by colored arcs or by
uncolored edges. Any area $a$ of $\Cfr$ defines a bubble $\Bfr$
consisting in the edges of $a$. The base of $\Bfr$ is the only edge of
$a$ that splits $\Cfr$ in two parts where one contains the base of
$\Cfr$ and the other contains $a$. Blue edges of $a$ remain blue edges
in $\Bfr$ and red edges of $a$ become uncolored edges in~$\Bfr$.
\medbreak

\begin{figure}[ht]
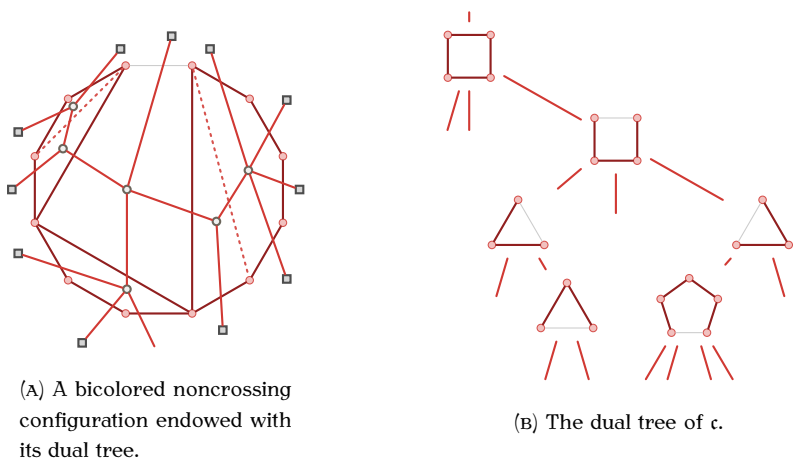

    \centering
    \subfloat[]
    [A bicolored noncrossing configuration endowed with its dual tree.]{
    \begin{minipage}[c]{.4\textwidth}
    \begin{equation*}

    \end{equation*}
    \end{minipage}
    \label{subfig:bicolored_noncrossing_configuration}}
    \subfloat[][The dual tree of $\Cfr$.]{
    \begin{minipage}[c]{.4\textwidth}
    \begin{equation*}
        %
};
            \node(9)at(17.50,-17.00){};
            \draw[Edge](0)--(1);
            \draw[Edge](10)--(12);
            \draw[Edge](11)--(12);
            \draw[Edge](12)--(15);
            \draw[Edge](13)--(12);
            \draw[Edge](14)--(12);
            \draw[Edge](15)--(8);
            \draw[Edge](16)--(15);
            \draw[Edge](2)--(1);
            \draw[Edge](3)--(4);
            \draw[Edge](4)--(8);
            \draw[Edge](5)--(6);
            \draw[Edge](6)--(4);
            \draw[Edge](7)--(6);
            \draw[Edge](8)--(1);
            \draw[Edge](9)--(8);
            \draw[Edge](1)--(2.5,4.5);
        \end{tikzpicture}
    \end{equation*}
    \end{minipage}
    \label{subfig:dual_tree}}
    \caption{A bicolored noncrossing configuration and its dual
    tree.}
    \label{fig:example_dual_tree_BNC}
\end{figure}

The \Def{dual tree} of $\Cfr$ is the planar rooted tree labeled by
bubbles defined as follows. If $\Cfr$ is of size $1$, its dual tree is
the leaf. Otherwise, put an internal node in each area of $\Cfr$ and
connect any pair of nodes that are in adjacent areas. Put also leaves
outside $\Cfr$, one for each edge, except the base, and connect these
with the internal nodes of their adjacent areas. This forms a tree
rooted at the node of the area containing the base of $\Cfr$. Finally,
label each internal node of the tree by the bubble associated with the
area containing it. Figure~\ref{fig:example_dual_tree_BNC}
shows an example of a BNC and its dual tree.
\medbreak

\begin{Theorem} \label{thm:bubble_decomposition_BNC}
    The $2$-colored operad $\Bubble$ is a $2$-bubble decomposition of
    the operad $\BNC$.
\end{Theorem}
\medbreak

\subsubsection{Enumeration of the bicolored noncrossing configurations}
By using the fact that, by Theorem~\ref{thm:bubble_decomposition_BNC},
$\Bubble$ is a $2$-bubble decomposition of $\BNC$, together with
Proposition~\ref{prop:hilbert_series_enveloping} and the colored Hilbert
series~\eqref{equ:series_bubbles_Bubble} of $\Bubble$, we obtain the
following algebraic equation for the generating series of the BNCs.
\medbreak

\begin{Proposition} \label{prop:hilbert_series_BNC}
    The Hilbert series $\HilbSeries(t)$ of $\BNC$ satisfies
    \begin{equation}
        -t - t^2 + (1 - 4t) \HilbSeries(t) - 3 \HilbSeries(t)^2 = 0.
    \end{equation}
\end{Proposition}
\medbreak

First numbers of BNCs by size are
\begin{equation}
    1, 8, 80, 992, 13760, 204416, 3180800, 51176960, 844467200.
\end{equation}
This forms Sequence~\OEIS{A234596} of~\cite{Slo}.
\medbreak

\subsubsection{Other consequences}
Since $\Bubble$ is, by Theorem~\ref{thm:bubble_decomposition_BNC}, a
$2$-bubble decomposition of $\BNC$, we can use the results of
Section~\ref{subsec:consequences_bubble_decomposition} to obtain the
generating set, the group of symmetries, and the presentation by
generators and relations of $\BNC$. Thus, by
Propositions~\ref{prop:generating_set_enveloping}
and~\ref{prop:generators_Bubble}, the generating set of $\BNC$ is the
set of the eight BNCs of arity $2$. By
Propositions~\ref{prop:group_symmetries_enveloping}
and~\ref{prop:symmetries_Bubble}, the group of symmetries of $\BNC$ is
generated by the maps $\Cpl' := \Hull(\Cpl)$ and $\Ret' := \Hull(\Ret)$.
For any BNC $\Cfr$, $\Cpl'(\Cfr)$ is the BNC obtained by swapping the
colors of the red and blue diagonals of $\Cfr$, and by swapping the
colors of the edges of $\Cfr$. Moreover, for any BNC $\Cfr$,
$\Ret'(\Cfr)$ is the BNC obtained by applying to $\Cfr$ the reflection
through the vertical line passing by its base. Finally, by
Proposition~\ref{prop:presentation} and
Theorem~\ref{thm:presentation_Bubble}, $\BNC$ admits the presentation
by generators and relations of the statement of
Theorem~\ref{thm:presentation_Bubble}.
\medbreak

\section{Suboperads of the operad of bicolored noncrossing
    configurations} \label{sec:suboperads_BNC}
We now study some of the suboperads of $\BNC$ generated by various sets
of BNCs. We shall mainly focus on the suboperads generated by sets of
two BNCs of arity~$2$.
\medbreak

\subsection{Overview of the obtained suboperads}
In what follows, we denote by $\OperadGen{\GeneratingSet}$ the
suboperad of $\BNC$ generated by a set $\GeneratingSet$ of BNCs, and
when $\GeneratingSet$ is a set of bubbles, by
$\OperadColGen{\GeneratingSet}$ the colored suboperad of $\Bubble$
generated by~$\GeneratingSet$.
\medbreak

\subsubsection{Orbits of suboperads}
There are $2^8 = 256$ suboperads of $\BNC$ generated by elements of
arity $2$. The symmetries provided by the group of symmetries of $\BNC$
(see Proposition~\ref{prop:symmetries_Bubble}) allow to gather some of
these together. Indeed, if $\GeneratingSet_1$ and $\GeneratingSet_2$
are two sets of BNCs and $\phi$ is a map of the group of symmetries of
$\BNC$ such that $\phi(\GeneratingSet_1) = \GeneratingSet_2$, the
suboperads $\OperadGen{\GeneratingSet_1}$ and
$\OperadGen{\GeneratingSet_2}$ would be isomorphic or antiisomorphic.
We say in this case that these two operads are \Def{equivalent}. There
are in this way only $88$ orbits of suboperads that are pairwise
nonequivalent.
\medbreak

\subsubsection{Suboperads on one generator}
There are three orbits of suboperads of $\BNC$ generated by one
generator of arity $2$. The first contains $\OperadGen{\AAA}$. By
induction on the arity, one can show that this operad contains all the
triangulations and that it is free. The second one contains
$\OperadGen{\AAB}$. By using similar arguments, one can show that this
operad is also free and isomorphic to the latter. The third orbit
contains $\OperadGen{\BAB}$. This operad contains exactly one element of
any arity, and hence, is the associative operad.
\medbreak

\subsubsection{Operads of noncrossing trees and plants}
Chapoton defined in~\cite{Cha07} an operad $\NCT$ involving based
noncrossing trees and an operad $\NCP$ involving noncrossing plants. As
follows directly from the definition, these operads are the suboperads
$\OperadGen{\AAB, \BAA}$ and $\OperadGen{\AAB, \ABA ,\BAA}$ of $\BNC$
respectively. The operad $\NCT$ governs $\LOp$-algebras, a sort of
algebras introduced by Leroux~\cite{Ler11}.
\medbreak

\subsubsection{Suboperads on two generators}
The $\binom{8}{2} = 28$ suboperads of $\BNC$ generated by two BNCs of
arity $2$ form eleven orbits. Table~\ref{tab:suboperads_BNC}
summarizes some information about these.
\begin{table}[ht]
    \centering
    \begin{tabular}{c|c|c}
        Operad & Dimensions & Presentation \\ \hline \hline
        $\OperadGen{\AAA, \AAB}$ &
            1, 2, 8, 40, 224, 1344, 8448, 54912 & free \\ \hline
        $\OperadGen{\AAA, \BBB}$ &
            1, 2, 8, 40, 216, 1246, 7516, 46838
            & quartic or more \\ \hline
        $\OperadGen{\AAA, \ABB}$ &
            \multirow{3}{*}{1, 2, 8, 38, 200, 1124, 6608, 40142} &
            \multirow{3}{*}{cubic} \\
        $\OperadGen{\AAB, \ABB}$ & & \\
        $\OperadGen{\AAB, \BBA}$ & & \\ \hline
        $\OperadGen{\AAA, \BAB}$ &
            1, 2, 7, 31, 154, 820, 4575, 26398 & quadratic \\ \hline
        $\OperadGen{\AAB, \BAA}$ &
            1, 2, 7, 30, 143, 728, 3876, 21318 & quadratic \\ \hline
        $\OperadGen{\AAB, \BAB}$ &
            \multirow{4}{*}{1, 2, 6, 22, 90, 394, 1806, 8558} &
            \multirow{4}{*}{quadratic} \\
        $\OperadGen{\AAA, \ABA}$ & & \\
        $\OperadGen{\BAA, \ABA}$ & & \\
        $\OperadGen{\BAB, \ABA}$ & &
    \end{tabular}
    \bigbreak

    \caption[Data about the binary suboperads of the ns operad~$\BNC$.]
    {The eleven orbits of suboperads of $\BNC$ generated by two
    generators of arity $2$, their dimensions and the degrees of
    nontrivial relations between their generators.}
    \label{tab:suboperads_BNC}
\end{table}
Some of these operads are well-known operads: the free operad
$\OperadGen{\AAA, \AAB}$ on two generators of arity $2$, the operad of
noncrossing trees~\cite{Cha07,Ler11} $\OperadGen{\AAB, \BAA}$, the
dipterous operad~\cite{LR03,Zin12} $\OperadGen{\BAA, \ABA}$, and the
$2$-associative operad~\cite{LR06,Zin12} $\OperadGen{\BAB, \ABA}$. All
the Hilbert series of the eleven operads are algebraic, with the genus
of the associated algebraic curve being $0$.  The sequences of the
dimensions of the operads of Table~\ref{tab:suboperads_BNC} are
respectively Sequences~\OEIS{A052701}, \OEIS{A234938}, \OEIS{A234939},
\OEIS{A007863}, \OEIS{A006013}, and~\OEIS{A006318} of~\cite{Slo}.
\medbreak

\subsubsection{Suboperads on more than two generators}
Some suboperads of $\BNC$ generated by more than two generators are very
complicated to study. For instance, the operad
$\langle \AAA, \BAB, \BBB\rangle$ has two equivalence classes of
nontrivial relations in degree $2$, three in degree $3$, ten in degree
$4$ and seems to have no nontrivial relations in higher degree (this
has been checked until degree $6$). The operad
$\langle \AAA, \ABA, \BAB, \BBB\rangle$ is also complicated since it
has four equivalence classes of nontrivial relations in degree $2$,
sixteen in degree $3$ and seems to have no nontrivial relations in
higher degree (this has been checked until degree $6$).
\medbreak

\subsection{Suboperads generated by two elements of arity 2}
For any of the eleven nonequivalent suboperads of $\BNC$ generated by
two elements of arity $2$, we compute its dimensions and provide a
presentation by generators and relations by passing through a bubble
decomposition of it.
\medbreak

\subsubsection{Outline of the study}
Let $\OperadGen{\GeneratingSet}$ be one of these suboperads. Since, by
Theorem~\ref{thm:bubble_decomposition_BNC}, $\Bubble$ is a bubble
decomposition of $\BNC$ and $\OperadGen{\GeneratingSet}$ is generated
by bubbles, $\OperadColGen{\GeneratingSet}$ is a bubble decomposition
of $\OperadGen{\GeneratingSet}$. We shall compute the dimensions and
establish the presentation by generators and relations of
$\OperadColGen{\GeneratingSet}$ to obtain in return, by
Propositions~\ref{prop:hilbert_series_enveloping}
and~\ref{prop:presentation}, the dimensions and the presentation by
generators and relations of $\OperadGen{\GeneratingSet}$. To compute the
dimensions of $\OperadColGen{G}$, we shall furnish a description of its
elements and then deduce from the description its colored Hilbert
series. Table~\ref{tab:suboperads_bubbles_BNC} shows the first
coefficients of the colored Hilbert series of the eleven colored
suboperads. All of these series are rational.
\begin{table}[ht]
    \centering
    \begin{tabular}{c|c|c}
        colored operad & Based bubbles & Nonbased bubbles \\
            \hline \hline
        $\OperadColGen{\AAA, \AAB}$ &
            2, 2, 2, 2, 2, 2, 2 &
            0, 0, 0, 0, 0, 0, 0 \\ \hline
        $\OperadColGen{\AAA, \BBB}$ &
            1, 2, 5, 10, 21, 42, 85 & 1, 2, 5, 10, 21, 42, 85 \\ \hline
        $\OperadColGen{\AAA, \ABB}$ &
            \multirow{3}{*}{1, 2, 4, 8, 16, 32, 64} &
            \multirow{3}{*}{1, 2, 4, 8, 16, 32, 64} \\
        $\OperadColGen{\AAB, \ABB}$ & & \\
        $\OperadColGen{\AAB, \BBA}$ & & \\ \hline
        $\OperadColGen{\AAA, \BAB}$ &
            2, 3, 5, 8, 13, 21, 34 &
            0, 0, 0, 0, 0, 0, 0 \\ \hline
        $\OperadColGen{\AAB, \BAA}$ &
            2, 3, 4, 5, 6, 7, 8 &
            0, 0, 0, 0, 0, 0, 0 \\ \hline
        $\OperadColGen{\AAB, \BAB}$ &
            2, 4, 8, 16, 32, 64, 128 &
            0, 0, 0, 0, 0, 0, 0 \\ \hline
        $\OperadColGen{\AAA, \ABA}$ &
            \multirow{3}{*}{1, 1, 1, 1, 1, 1, 1} &
            \multirow{3}{*}{1, 1, 1, 1, 1, 1, 1} \\
        $\OperadColGen{\BAA, \ABA}$ & & \\
        $\OperadColGen{\BAB, \ABA}$ & &
    \end{tabular}
    \bigbreak

    \caption[Data about the binary suboperads of the ns colored
    operad~$\Bubble$.]
    {The eleven orbits of $2$-colored suboperads of $\Bubble$
    generated by two generators of arity $2$ and the number of their
    bubbles, based and nonbased.}
    \label{tab:suboperads_bubbles_BNC}
\end{table}
To establish the presentation of $\OperadColGen{\GeneratingSet}$, we
shall use the same strategy as the one used for the proof of the
presentation of $\Bubble$ (see the proof of
Theorem~\ref{thm:presentation_Bubble}). Recall that this consists
in exhibiting an orientation $\Rew$ of the presentation we want to prove
such that its closure $\RewTrees$ is a terminating rewrite rule on
colored syntax trees and its normal forms are in bijection with the
elements of~$\OperadColGen{\GeneratingSet}$.
\medbreak

\subsubsection{First orbit}
This orbit consists of the operads $\OperadGen{\AAA, \AAB}$,
$\OperadGen{\AAA, \BAA}$, $\OperadGen{\BBB, \ABB}$, and
$\OperadGen{\BBB, \BBA}$. We choose $\OperadGen{\AAA, \AAB}$ as a
representative of the orbit.
\medbreak

\begin{Proposition} \label{prop:bubbles_aaa_aab}
    The set of bubbles of $\OperadColGen{\AAA, \AAB}$ is the set of
    based bubbles such that all edges of the border except possibly the
    last one are blue. Moreover, the colored Hilbert series of
    $\OperadColGen{\AAA, \AAB}$ satisfy
    \begin{equation}
        \BubbleSeries_1(\VarZ_1, \VarZ_2) =
        \frac{\VarZ_1\VarZ_2 + \VarZ_2^2}{1 - \VarZ_2}
        \qquad \mbox{ and } \qquad
        \BubbleSeries_2(\VarZ_1, \VarZ_2) = 0.
    \end{equation}
\end{Proposition}
\medbreak

\begin{Proposition} \label{prop:hilbert_series_aaa_aab}
    The Hilbert series $\HilbSeries(t)$ of $\OperadGen{\AAA, \AAB}$
    satisfies
    \begin{equation} \label{equ:hilbert_series_aaa_aab}
        t - \HilbSeries(t) + 2 \HilbSeries(t)^2 = 0.
    \end{equation}
\end{Proposition}
\medbreak

\begin{Theorem} \label{thm:presentation_aaa_aab}
    The operad $\OperadGen{\AAA, \AAB}$ is the free operad generated by
    two generators of arity $2$.
\end{Theorem}
\medbreak

\subsubsection{Second orbit}
This orbit consists of the operad $\OperadGen{\AAA, \BBB}$.
\medbreak

\begin{Proposition} \label{prop:bubbles_aaa_bbb}
    The set of based (resp. nonbased) bubbles of
    $\OperadColGen{\AAA, \BBB}$ of arity $n$ is the set of based (resp.
    nonbased) bubbles having at least two consecutive edges of the
    border of a same color and the number of blue (resp. uncolored)
    edges of the border is congruent to $1 - n$ modulo $3$. Moreover,
    the colored Hilbert series of $\OperadColGen{\AAA, \BBB}$ satisfy
    \begin{subequations}
    \begin{equation}
        \BubbleSeries_1(\VarZ_1, \VarZ_2) =
        \frac{\VarZ_1 + \VarZ_2^2}
        {1 - 3\VarZ_1\VarZ_2 - \VarZ_1^3 - \VarZ_2^3} -
        \frac{\VarZ_1}{1 - \VarZ_1\VarZ_2}
    \end{equation}
    and
    \begin{equation}
        \BubbleSeries_2(\VarZ_1, \VarZ_2) =
        \frac{\VarZ_2 + \VarZ_1^2}
        {1 - 3\VarZ_1\VarZ_2 - \VarZ_1^3 - \VarZ_2^3} -
        \frac{\VarZ_2}{1 - \VarZ_1\VarZ_2}.
    \end{equation}
    \end{subequations}
\end{Proposition}
\medbreak

\begin{Proposition} \label{prop:hilbert_series_aaa_bbb}
    The Hilbert series $\HilbSeries(t)$ of $\OperadGen{\AAA, \BBB}$
    satisfies
    \begin{equation}
        4t - 2t^2 - t^3 + t^4 + (-4 + 4t - t^2 + 2t^3)\HilbSeries(t) +
        (6 + t)\HilbSeries(t)^2 + (1 - 2t)\HilbSeries(t)^3
        - \HilbSeries(t)^4
        = 0.
    \end{equation}
\end{Proposition}
\medbreak

\begin{Proposition} \label{prop:relations_aaa_bbb}
    The operad $\OperadGen{\AAA, \BBB}$ does not admit nontrivial
    relations between its generators in degree two, three, five and
    six. It admits the following non trivial relations between its
    generators in degree four:
    \begin{subequations}
    \begin{equation}
        ((\AAA \circ_2 \BBB) \circ_3 \AAA) \circ_3 \BBB
        \enspace = \enspace
        ((\AAA \circ_1 \BBB) \circ_1 \AAA) \circ_2 \BBB,
    \end{equation}
    \begin{equation}
        ((\AAA \circ_2 \BBB) \circ_2 \AAA) \circ_4 \AAA
        \enspace = \enspace
        ((\AAA \circ_1 \BBB) \circ_1 \AAA) \circ_3 \AAA,
    \end{equation}
    \begin{equation}
        ((\AAA \circ_2 \BBB) \circ_2 \AAA) \circ_3 \BBB
        \enspace = \enspace
        ((\AAA \circ_1 \BBB) \circ_1 \AAA) \circ_4 \BBB,
    \end{equation}
    \begin{equation}
        ((\AAA \circ_1 \BBB) \circ_3 \BBB) \circ_4 \AAA
        \enspace = \enspace
        ((\AAA \circ_1 \BBB) \circ_2 \AAA) \circ_2 \BBB,
    \end{equation}
    \begin{equation}
        ((\BBB \circ_2 \AAA) \circ_3 \BBB) \circ_3 \AAA
        \enspace = \enspace
        ((\BBB \circ_1 \AAA) \circ_1 \BBB) \circ_2 \AAA,
    \end{equation}
    \begin{equation}
        ((\BBB \circ_2 \AAA) \circ_2 \BBB) \circ_4 \BBB
        \enspace = \enspace
        ((\BBB \circ_1 \AAA) \circ_1 \BBB) \circ_3 \BBB,
    \end{equation}
    \begin{equation}
        ((\BBB \circ_2 \AAA) \circ_2 \BBB) \circ_3 \AAA
        \enspace = \enspace
        ((\BBB \circ_1 \AAA) \circ_1 \BBB) \circ_4 \AAA,
    \end{equation}
    \begin{equation}
        ((\BBB \circ_1 \AAA) \circ_3 \AAA) \circ_4 \BBB
        \enspace = \enspace
        ((\BBB \circ_1 \AAA) \circ_2 \BBB) \circ_2 \AAA.
    \end{equation}
    \end{subequations}
\end{Proposition}
\begin{proof}
    This statement is proven with the help of the computer. All partial
    compositions between the generators $\AAA$ and $\BBB$ are computed
    up to degree six and relations thus established.
\end{proof}
\medbreak

Proposition~\ref{prop:relations_aaa_bbb} does not provide a presentation
by generators and relations of $\OperadGen{\AAA, \BBB}$. The methods
employed in this chapter fail to establish the presentation of
$\OperadColGen{\AAA, \BBB}$ because it is not possible to define an
orientation $\Rew$ of the relations of the statement of
Proposition~\ref{prop:relations_aaa_bbb}. Indeed, in degree six, all
the closures $\RewTrees$ of $\Rew$ have no less than $7518$ normal forms
whereas they should be $7516$. Nevertheless, these relations seem to be
the only nontrivial ones; this may be proved by using the Knuth-Bendix
completion algorithm (see~\cite{KB70,BN98}) over an appropriate
orientation of the relations.
\medbreak

\subsubsection{Third orbit}
This orbit consists of the operads $\OperadGen{\AAA, \ABB}$,
$\OperadGen{\AAA, \BBA}$, $\OperadGen{\AAB, \BBB}$, and
$\OperadGen{\BAA, \BBB}$. We choose $\OperadGen{\AAA, \ABB}$ as a
representative of the orbit.
\medbreak

\begin{Proposition} \label{prop:bubbles_aaa_abb}
    The set of based (resp. nonbased) bubbles of
    $\OperadColGen{\AAA, \ABB}$ of arity $n$ is the set of based (resp.
    nonbased) bubbles such that first edge is blue and the number of
    uncolored edges of the border is congruent to $n$ (resp. $n + 1$)
    modulo $2$. Moreover, the colored Hilbert series of
    $\OperadColGen{\AAA, \ABB}$ satisfy
    \begin{equation}
        \BubbleSeries_1(\VarZ_1, \VarZ_2) =
        \frac{\VarZ_2^2}{1 - 2\VarZ_1 + \VarZ_1^2 - \VarZ_2^2}
        \qquad \mbox{ and } \qquad
        \BubbleSeries_2(\VarZ_1, \VarZ_2) =
        \frac{\VarZ_1\VarZ_2 - \VarZ_1^2\VarZ_2 + \VarZ_2^3}
        {1 - 2\VarZ_1 + \VarZ_1^2 - \VarZ_2^2}.
    \end{equation}
\end{Proposition}
\medbreak

\begin{Proposition} \label{prop:hilbert_series_aaa_abb}
    The Hilbert series $\HilbSeries(t)$ of $\OperadGen{\AAA, \ABB}$
    satisfies
    \begin{equation}
        2t - t^2 + (2t - 2)\HilbSeries(t) + 3\HilbSeries(t)^2 = 0.
    \end{equation}
\end{Proposition}
\medbreak

\begin{Theorem} \label{thm:presentation_aaa_abb}
    The operad $\OperadGen{\AAA, \ABB}$ admits the presentation
    $(\{\AAA, \ABB\}, \Equiv)$ where $\Equiv$ is the
    equivalence relation satisfying
    \begin{subequations}
    \begin{equation}
        \left(\Corolla{\ABB} \circ_2 \Corolla{\AAA}\right)
        \circ_3 \Corolla{\ABB}
        \enspace \Equiv \enspace
        \left(\Corolla{\ABB} \circ_1 \Corolla{\ABB}\right)
        \circ_2 \Corolla{\AAA},
    \end{equation}
    \begin{equation}
        \left(\Corolla{\AAA} \circ_2 \Corolla{\ABB}\right)
        \circ_3 \Corolla{\AAA}
        \enspace \Equiv \enspace
        \left(\Corolla{\AAA} \circ_1 \Corolla{\ABB}\right)
        \circ_2 \Corolla{\AAA}.
    \end{equation}
    \end{subequations}
    Moreover, the set of the $\{\AAA, \ABB\}$-syntax trees avoiding
    the trees
    \begin{math}
        \left(\Corolla{\ABB} \circ_1 \Corolla{\ABB}\right)
        \circ_2 \Corolla{\AAA}
    \end{math}
    and
    \begin{math}
        \left(\Corolla{\AAA} \circ_2 \Corolla{\ABB}\right)
        \circ_3 \Corolla{\AAA}
    \end{math}
    is a Poincaré-Birkhoff-Witt basis of~$\OperadGen{\AAA, \ABB}$.
\end{Theorem}
\medbreak

\subsubsection{Fourth orbit}
This orbit consists of the operads $\OperadGen{\AAB, \ABB}$ and
$\OperadGen{\BAA, \BBA}$. We choose $\OperadGen{\AAB, \ABB}$ as a
representative of the orbit.
\medbreak

\begin{Proposition} \label{prop:bubbles_aab_abb}
    The set of bubbles of $\OperadColGen{\AAB, \ABB}$ is the set of
    bubbles such that first edge is blue and last edge uncolored.
    Moreover, the colored Hilbert series of $\OperadColGen{\AAB, \ABB}$
    satisfy
    \begin{equation}
        \BubbleSeries_1(\VarZ_1, \VarZ_2) =
        \frac{\VarZ_1\VarZ_2}{1 - \VarZ_1 - \VarZ_2}
        \qquad \mbox{ and } \qquad
        \BubbleSeries_2(\VarZ_1, \VarZ_2) =
        \frac{\VarZ_1\VarZ_2}{1 - \VarZ_1 - \VarZ_2}.
    \end{equation}
\end{Proposition}
\medbreak

\begin{Proposition} \label{prop:hilbert_series_aab_abb}
    The Hilbert series $\HilbSeries(t)$ of $\OperadGen{\AAB, \ABB}$
    satisfies
    \begin{equation}
        2t - t^2 + (2t - 2)\HilbSeries(t) + 3\HilbSeries(t)^2 = 0.
    \end{equation}
\end{Proposition}
\medbreak

\begin{Theorem} \label{thm:presentation_aab_abb}
    The operad $\OperadGen{\AAB, \ABB}$ admits the presentation
    $(\{\AAB, \ABB\}, \Equiv)$ where $\Equiv$ is the
    equivalence relation satisfying
    \begin{subequations}
    \begin{equation}
        \left(\Corolla{\ABB} \circ_2 \Corolla{\AAB}\right)
        \circ_2 \Corolla{\ABB}
        \enspace \Equiv \enspace
        \left(\Corolla{\ABB} \circ_1 \Corolla{\ABB}\right)
        \circ_2 \Corolla{\AAB},
    \end{equation}
    \begin{equation}
        \left(\Corolla{\AAB} \circ_2 \Corolla{\AAB}\right)
        \circ_2 \Corolla{\ABB}
        \enspace \Equiv \enspace
        \left(\Corolla{\AAB} \circ_1 \Corolla{\ABB}\right)
        \circ_2 \Corolla{\AAB}.
    \end{equation}
    \end{subequations}
    Moreover, the set of the $\{\AAB, \ABB\}$-syntax trees avoiding
    the trees
    \begin{math}
        \left(\Corolla{\ABB} \circ_2 \Corolla{\AAB}\right)
        \circ_2 \Corolla{\ABB}
    \end{math}
    and
    \begin{math}
        \left(\Corolla{\AAB} \circ_2 \Corolla{\AAB}\right)
        \circ_2 \Corolla{\ABB}
    \end{math}
    is a Poincaré-Birkhoff-Witt basis of~$\OperadGen{\AAB, \ABB}$.
\end{Theorem}
\medbreak

\subsubsection{Fifth orbit}
This orbit consists of the operads $\OperadGen{\AAB, \BBA}$ and
$\OperadGen{\BAA, \ABB}$. We choose $\OperadGen{\AAB, \BBA}$ as a
representative of the orbit.
\medbreak

\begin{Proposition} \label{prop:bubbles_aab_bba}
    The set of based (resp. nonbased) bubbles of
    $\OperadColGen{\AAB, \BBA}$ is the set of based (resp. nonbased)
    bubbles such that penultimate edge is blue (resp. uncolored) and
    the last edge is uncolored (resp. blue). Moreover, the colored
    Hilbert series of $\OperadColGen{\AAB, \BBA}$ satisfy
    \begin{equation}
        \BubbleSeries_1(\VarZ_1, \VarZ_2) =
        \frac{\VarZ_1\VarZ_2}{1 - \VarZ_1 - \VarZ_2}
        \qquad \mbox{ and } \qquad
        \BubbleSeries_2(\VarZ_1, \VarZ_2) =
        \frac{\VarZ_1\VarZ_2}{1 - \VarZ_1 - \VarZ_2}.
    \end{equation}
\end{Proposition}
\medbreak

\begin{Proposition} \label{prop:hilbert_series_aab_bba}
    The Hilbert series $\HilbSeries(t)$ of $\OperadGen{\AAB, \BBA}$
    satisfies
    \begin{equation}
        2t - t^2 + (2t - 2)\HilbSeries(t) + 3\HilbSeries(t)^2 = 0.
    \end{equation}
\end{Proposition}
\medbreak

\begin{Theorem} \label{thm:presentation_aab_bba}
    The operad $\OperadGen{\AAB, \BBA}$ admits the presentation
    $(\{\AAB, \AAB\}, \Equiv)$ where $\Equiv$ is the
    equivalence relation satisfying
    \begin{equation}
        \left(\Corolla{\AAB} \circ_2 \Corolla{\AAB}\right)
        \circ_2 \Corolla{\BBA}
        \enspace \Equiv \enspace
        \left(\Corolla{\AAB} \circ_1 \Corolla{\BBA}\right)
        \circ_1 \Corolla{\AAB},
    \end{equation}
    \begin{equation}
        \left(\Corolla{\BBA} \circ_2 \Corolla{\BBA}\right)
        \circ_2 \Corolla{\AAB}
        \enspace \Equiv \enspace
        \left(\Corolla{\BBA} \circ_1 \Corolla{\AAB}\right)
        \circ_1 \Corolla{\BBA}.
    \end{equation}
    Moreover, the set of the $\{\AAB, \BBA\}$-syntax trees avoiding
    the trees
    \begin{math}
        \left(\Corolla{\AAB} \circ_2 \Corolla{\AAB}\right)
        \circ_2 \Corolla{\BBA}
    \end{math}
    and
    \begin{math}
        \left(\Corolla{\BBA} \circ_2 \Corolla{\BBA}\right)
        \circ_2 \Corolla{\AAB}
    \end{math}
    is a Poincaré-Birkhoff-Witt basis of~$\OperadGen{\AAB, \BBA}$.
\end{Theorem}
\medbreak

\subsubsection{Sixth orbit}
This orbit consists of the operads $\OperadGen{\AAA, \BAB}$ and
$\OperadGen{\ABA, \BBB}$. We choose $\OperadGen{\AAA, \BAB}$ as a
representative of the orbit.
\medbreak

\begin{Proposition} \label{prop:bubbles_aaa_bab}
    The set of bubbles of $\OperadColGen{\AAA, \BAB}$ is the set of
    based bubbles such that maximal sequences of blues edges of the
    border have even length.Moreover, the colored Hilbert series of
    $\OperadColGen{\AAA, \BAB}$ satisfy
    \begin{equation}
        \BubbleSeries_1(\VarZ_1, \VarZ_2) =
        \frac{\VarZ_1^2 + \VarZ_2^2 + \VarZ_1\VarZ_2^2}
        {1 - \VarZ_1 - \VarZ_2^2}
        \qquad \mbox{ and } \qquad
        \BubbleSeries_2(\VarZ_1, \VarZ_2) = 0.
    \end{equation}
\end{Proposition}
\medbreak

\begin{Proposition} \label{prop:hilbert_series_aaa_bab}
    The Hilbert series $\HilbSeries(t)$ of $\OperadGen{\AAA, \BAB}$
    satisfies
    \begin{equation}
        t + (t - 1) \HilbSeries(t) + \HilbSeries(t)^2 +
        \HilbSeries(t)^3 = 0.
    \end{equation}
\end{Proposition}
\medbreak

\begin{Theorem} \label{thm:presentation_aaa_bab}
    The operad $\OperadGen{\AAA, \BAB}$ admits the presentation
    $(\{\AAA, \BAB\}, \Equiv)$ where $\Equiv$ is the
    equivalence relation satisfying
    \begin{equation}
        \Corolla{\BAB} \circ_1 \Corolla{\BAB}
        \enspace \Equiv \enspace
        \Corolla{\BAB} \circ_2 \Corolla{\BAB}.
    \end{equation}
    Moreover, $\OperadGen{\AAA, \BAB}$ is Koszul and the set of the
    $\{\AAA, \BAB\}$-syntax trees avoiding the tree
    \begin{math}
        \Corolla{\BAB} \circ_2 \Corolla{\BAB}
    \end{math}
    is a Poincaré-Birkhoff-Witt basis of~$\OperadGen{\AAA, \BAB}$.
\end{Theorem}
\medbreak

\subsubsection{Seventh orbit} \label{subsubsec:construction_NCT}
This orbit consists of the operads $\OperadGen{\AAB, \BAA}$ and
$\OperadGen{\ABB, \BBA}$. We choose $\OperadGen{\AAB, \BAA}$ as a
representative of the orbit.
\medbreak

\begin{Proposition} \label{prop:bubbles_aab_baa}
    The set of bubbles of $\OperadColGen{\AAB, \BAA}$ is the set of
    based bubbles having exactly one uncolored edge in the border.
    Moreover, the colored Hilbert series of $\OperadColGen{\AAB, \BAA}$
    satisfy
    \begin{equation}
        \BubbleSeries_1(\VarZ_1, \VarZ_2) =
        \frac{2\VarZ_1\VarZ_2 - \VarZ_1\VarZ_2^2}{(1 - \VarZ_2)^2}
        \qquad \mbox{ and } \qquad
        \BubbleSeries_2(\VarZ_1, \VarZ_2) = 0.
    \end{equation}
\end{Proposition}
\medbreak

\begin{Proposition} \label{prop:hilbert_series_aab_baa}
    The Hilbert series $\HilbSeries(t)$ of $\OperadGen{\AAB, \BAA}$
    satisfies
    \begin{equation}
        t - \HilbSeries(t) + 2 \HilbSeries(t)^2 - \HilbSeries(t)^3 = 0.
    \end{equation}
\end{Proposition}
\medbreak

\begin{Theorem} \label{thm:presentation_aab_baa}
    The operad $\OperadGen{\AAB, \BAA}$ admits the presentation
    $(\{\AAB, \BAA\}, \Equiv)$ where $\Equiv$ is the
    equivalence relation satisfying
    \begin{equation}
        \Corolla{\BAA} \circ_1 \Corolla{\AAB}
        \enspace \Equiv \enspace
        \Corolla{\AAB} \circ_2 \Corolla{\BAA}.
    \end{equation}
    Moreover, $\OperadGen{\AAB, \BAA}$ is Koszul and the set of the
    $\{\AAB, \BAA\}$-syntax trees avoiding the tree
    \begin{math}
        \Corolla{\AAB} \circ_2 \Corolla{\BAA}
    \end{math}
    is a Poincaré-Birkhoff-Witt basis of~$\OperadGen{\AAB, \BAA}$.
\end{Theorem}
\medbreak

The above presentation shows that $\OperadGen{\AAB, \BAA}$ is the
operad of based noncrossing trees~$\NCT$.
\medbreak

\subsubsection{Eighth orbit}
This orbit consists of the operads $\OperadGen{\AAB, \BAB}$,
$\OperadGen{\BAA, \BAB}$, $\OperadGen{\ABA, \ABB}$, and
$\OperadGen{\ABA, \BBA}$. We choose $\OperadGen{\AAB, \BAB}$ as a
representative of the orbit.
\medbreak

\begin{Proposition} \label{prop:bubbles_aab_bab}
    The set of bubbles of $\OperadColGen{\AAB, \BAB}$ is the set of
    based bubbles such that last edge is uncolored. Moreover, the
    colored Hilbert series of $\OperadColGen{\AAB, \BAB}$ satisfy
    \begin{equation}
        \BubbleSeries_1(\VarZ_1, \VarZ_2) =
        \frac{\VarZ_1^2 + \VarZ_1\VarZ_2}{1 - \VarZ_1 - \VarZ_2}
        \qquad \mbox{ and } \qquad
        \BubbleSeries_2(\VarZ_1, \VarZ_2) = 0.
    \end{equation}
\end{Proposition}
\medbreak

\begin{Proposition} \label{prop:hilbert_series_aab_bab}
    The Hilbert series $\HilbSeries(t)$ of $\OperadGen{\AAB, \BAB}$
    satisfies
    \begin{equation}
        t - (1 - t) \HilbSeries(t) + \HilbSeries(t)^2 = 0.
    \end{equation}
\end{Proposition}
\medbreak

\begin{Theorem} \label{thm:presentation_aab_bab}
    The operad $\OperadGen{\AAB, \BAB}$ admits the presentation
    $(\{\AAB, \BAB\}, \Equiv)$ where $\Equiv$ is the
    equivalence relation satisfying
    \begin{equation}
        \Corolla{\BAB} \circ_1 \Corolla{\BAB}
        \enspace \Equiv \enspace
        \Corolla{\BAB} \circ_2 \Corolla{\BAB},
    \end{equation}
    \begin{equation}
        \Corolla{\BAB} \circ_1 \Corolla{\AAB}
        \enspace \Equiv \enspace
        \Corolla{\AAB} \circ_2 \Corolla{\BAB}.
    \end{equation}
    Moreover, $\OperadGen{\AAB, \BAB}$ is Koszul and the set of the
    $\{\AAB, \BAB\}$-syntax trees avoiding the trees
    \begin{math}
        \Corolla{\BAB} \circ_1 \Corolla{\BAB}
    \end{math}
    and
    \begin{math}
        \Corolla{\BAB} \circ_1 \Corolla{\AAB}
    \end{math}
    is a Poincaré-Birkhoff-Witt basis of~$\OperadGen{\AAB, \BAB}$.
\end{Theorem}
\medbreak

This operad $\OperadGen{\AAB, \BAB}$ is similar to the duplicial
operad~\cite{Lod08,Zin12} with the difference that $\AAB$ is not
associative in~$\OperadGen{\AAB, \BAB}$.
\medbreak

\subsubsection{Ninth orbit}
This orbit consists of the operads $\OperadGen{\AAA, \ABA}$ and
$\OperadGen{\BAB, \BBB}$. We choose $\OperadGen{\AAA, \ABA}$ as a
representative of the orbit.
\medbreak

\begin{Proposition} \label{prop:bubbles_aaa_aba}
    The set of bubbles of $\OperadColGen{\AAA, \ABA}$ is the set of
    bubbles such that all edges of the border are blue. Moreover, the
    colored Hilbert series of $\OperadColGen{\AAA, \ABA}$ satisfy
    \begin{equation}
        \BubbleSeries_1(\VarZ_1, \VarZ_2) =
        \frac{\VarZ_2^2}{1 - \VarZ_2}
        \qquad \mbox{ and } \qquad
        \BubbleSeries_2(\VarZ_1, \VarZ_2) =
        \frac{\VarZ_2^2}{1 - \VarZ_2}.
    \end{equation}
\end{Proposition}
\medbreak

\begin{Proposition} \label{prop:hilbert_series_aaa_aba}
    The Hilbert series $\HilbSeries(t)$ of $\OperadGen{\AAA, \ABA}$
    satisfies
    \begin{equation}
        t - (1 - t) \HilbSeries(t) + \HilbSeries(t)^2 = 0.
    \end{equation}
\end{Proposition}
\medbreak

\begin{Theorem} \label{thm:presentation_aaa_aba}
    The operad $\OperadGen{\AAA, \ABA}$ admits the presentation
    $(\{\AAA, \ABA\}, \Equiv)$ where $\Equiv$ is the
    equivalence relation satisfying
    \begin{equation}
        \Corolla{\ABA} \circ_1 \Corolla{\ABA}
        \enspace \Equiv \enspace
        \Corolla{\ABA} \circ_2 \Corolla{\ABA},
    \end{equation}
    \begin{equation}
        \Corolla{\AAA} \circ_1 \Corolla{\ABA}
        \enspace \Equiv \enspace
        \Corolla{\AAA} \circ_2 \Corolla{\ABA}.
    \end{equation}
    Moreover, $\OperadGen{\AAA, \ABA}$ is Koszul and the set of the
    $\{\AAA, \ABA\}$-syntax trees avoiding the trees
    \begin{math}
        \Corolla{\ABA} \circ_2 \Corolla{\ABA}
    \end{math}
    and
    \begin{math}
        \Corolla{\AAA} \circ_2 \Corolla{\ABA}
    \end{math}
    is a Poincaré-Birkhoff-Witt basis of~$\OperadGen{\AAA, \ABA}$.
\end{Theorem}
\medbreak

\subsubsection{Tenth orbit}
This orbit consists of the operads $\OperadGen{\AAB, \ABA}$,
$\OperadGen{\BAA, \ABA}$, $\OperadGen{\BAB, \ABB}$, and
$\OperadGen{\BAB, \BBA}$. We choose $\OperadGen{\BAA, \ABA}$ as a
representative of the orbit.
\medbreak

\begin{Proposition} \label{prop:bubbles_aab_aba}
    The set of based (resp. nonbased) bubbles of
    $\OperadColGen{\BAA, \ABA}$ is the set of based (resp. nonbased)
    bubbles such that first edge is uncolored (resp. blue) and the
    other edges of the border are blue. Moreover, the colored Hilbert
    series of $\OperadColGen{\BAA, \ABA}$ satisfy
    \begin{equation}
        \BubbleSeries_1(\VarZ_1, \VarZ_2) =
        \frac{\VarZ_1\VarZ_2}{1 - \VarZ_2}
        \qquad \mbox{ and } \qquad
        \BubbleSeries_2(\VarZ_1, \VarZ_2) =
        \frac{\VarZ_2^2}{1 - \VarZ_2}.
    \end{equation}
\end{Proposition}
\medbreak

\begin{Proposition} \label{prop:hilbert_series_aab_aba}
    The Hilbert series $\HilbSeries(t)$ of $\OperadGen{\BAA, \ABA}$
    satisfies
    \begin{equation}
        t - (1 - t) \HilbSeries(t) + \HilbSeries(t)^2 = 0.
    \end{equation}
\end{Proposition}
\medbreak

\begin{Theorem} \label{thm:presentation_aab_aba}
    The operad $\OperadGen{\BAA, \ABA}$ admits the presentation
    $(\{\BAA, \ABA\}, \Equiv)$ where $\Equiv$ is the
    equivalence relation satisfying
    \begin{equation}
        \Corolla{\ABA} \circ_1 \Corolla{\ABA}
        \enspace \Equiv \enspace
        \Corolla{\ABA} \circ_2 \Corolla{\ABA},
    \end{equation}
    \begin{equation}
        \Corolla{\BAA} \circ_1 \Corolla{\BAA}
        \enspace \Equiv \enspace
        \Corolla{\BAA} \circ_2 \Corolla{\ABA}.
    \end{equation}
    Moreover, $\OperadGen{\BAA, \ABA}$ is Koszul and the set of the
    $\{\BAA, \ABA\}$-syntax trees avoiding the trees
    \begin{math}
        \Corolla{\ABA} \circ_1 \Corolla{\ABA}
    \end{math}
    and
    \begin{math}
        \Corolla{\BAA} \circ_1 \Corolla{\BAA}
    \end{math}
    is a Poincaré-Birkhoff-Witt basis of~$\OperadGen{\BAA, \ABA}$.
\end{Theorem}
\medbreak

This operad $\OperadGen{\BAA, \ABA}$ is hence the operad governing
dipterous algebras~\cite{LR03,Zin12}.
\medbreak

\subsubsection{Eleventh orbit}
This orbit consists of the operad $\OperadGen{\BAB, \ABA}$.
\medbreak

\begin{Proposition} \label{prop:bubbles_bab_aba}
    The set of based (resp. nonbased) bubbles of
    $\OperadColGen{\BAB, \ABA}$ is the set of based (resp. nonbased)
    bubbles such that all edges of the border are uncolored (resp.
    blue). Moreover, the colored Hilbert series of
    $\OperadColGen{\BAB, \ABA}$ satisfy
    \begin{equation}
        \BubbleSeries_1(\VarZ_1, \VarZ_2) =
        \frac{\VarZ_1^2}{1 - \VarZ_1}
        \qquad \mbox{ and } \qquad
        \BubbleSeries_2(\VarZ_1, \VarZ_2) =
        \frac{\VarZ_2^2}{1 - \VarZ_2}.
    \end{equation}
\end{Proposition}
\medbreak

\begin{Proposition} \label{prop:hilbert_series_bab_aba}
    The Hilbert series $\HilbSeries(t)$ of $\OperadGen{\BAB, \ABA}$
    satisfies
    \begin{equation}
        t - (1 - t) \HilbSeries(t) + \HilbSeries(t)^2 = 0.
    \end{equation}
\end{Proposition}
\medbreak

\begin{Theorem} \label{thm:presentation_bab_aba}
    The operad $\OperadGen{\BAB, \ABA}$ admits the presentation
    $(\{\BAB, \ABA\}, \Equiv)$ where $\Equiv$ is the
    equivalence relation satisfying
    \begin{equation}
        \Corolla{\ABA} \circ_1 \Corolla{\ABA}
        \enspace \Equiv \enspace
        \Corolla{\ABA} \circ_2 \Corolla{\ABA},
    \end{equation}
    \begin{equation}
        \Corolla{\BAB} \circ_1 \Corolla{\BAB}
        \enspace \Equiv \enspace
        \Corolla{\BAB} \circ_2 \Corolla{\BAB}.
    \end{equation}
    Moreover, $\OperadGen{\BAB, \ABA}$ is Koszul and the set of the
    $\{\BAB, \ABA\}$-syntax trees avoiding the trees
    \begin{math}
        \Corolla{\ABA} \circ_2 \Corolla{\ABA}
    \end{math}
    and
    \begin{math}
        \Corolla{\BAB} \circ_2 \Corolla{\BAB}
    \end{math}
    is a Poincaré-Birkhoff-Witt basis of~$\OperadGen{\BAB, \ABA}$.
\end{Theorem}
\medbreak

This operad $\OperadGen{\BAA, \ABA}$ is hence the operad governing
two-associative algebras~\cite{LR06,Zin12}.
\medbreak

\section*{Concluding remarks}
We have developed in this chapter a tool to study a noncolored operad
$\Oca$ by considering one of its bubble decompositions $\Cca$. As
explained, most of the properties of $\Oca$ come from properties of
$\Cca$. This, together with the fact that a colored operad is a more
constrained structure than a noncolored one, leads to easier proofs
for most of their properties (for instance to establish a presentation
by generators and relations, a bubble decomposition reduces the number
of orientations of relations to consider). This framework has been
applied to the operad $\BNC$ of bicolored noncrossing configurations
and on some of its suboperads. As an additional remark, this chapter
considers only colored or uncolored set-operads but the notion of
enveloping operads and bubble decompositions also work for linear
operads.
\medbreak


\chapter{From monoids to operads} \label{chap:monoids}
The content of this chapter comes from~\cite{Gir12b,Gir12c,Gir15a}. A
very preliminary version of the results presented here was sketched
in~\cite{Gir11} while the work was still in progress. We include here
some new results that do not appear in the aforementioned publications
like the Koszulity and presentations of the Koszul duals of some of the
constructed ns operads.
\medbreak

\section*{Introduction}
We propose here a new generic method to build combinatorial operads. The
starting point is to pick a monoid $\Mca$. We then consider the set of
words whose letters are elements of $\Mca$. The arity of such words are
their length, the composition of two words is expressed from the
product of $\Mca$, and permutations act on words by permuting letters.
In this way, we associate an operad denoted by $\T \Mca$ with any
monoid $\Mca$. This construction is rich from a combinatorial point of
view since it allows us, by considering suboperads and quotients of
$\T \Mca$, to get new (symmetric or not) operads on various
combinatorial objects. Our construction is related to two previous ones.
\smallbreak

The first one is a construction of Méndez and Nava~\cite{MN93} emerging
from the context of the species theory~\cite{Joy81}. Roughly speaking, a
species is a combinatorial construction $U$ which takes an underlying
finite set $E$ as input and produces a set $U[E]$ of objects by adding
some structure on the elements of $E$ (see~\cite{BLL98}). This theory
has many links with the theory of operads since an operad is a monoid
with respect to the operation of substitution of species.
In~\cite{MN93}, the authors defined the plethystic species, that are
species taking as input sets where any element has a color picked from a
fixed monoid $\Mca$. This monoid has to satisfy some precise conditions
(as to be left cancellable and without proper divisor of the unit, and
such that any element has finitely many factorizations). It appears
that the elements of the so-called uniform plethystic species can be
seen as words of colors and hence, as elements of $\T \Mca$. Moreover,
the composition of this operad is the one of $\T \Mca$. The main
difference between the construction of Méndez and Nava and ours lies in
the fact that the construction $\T$ can be applied to any monoid.
\smallbreak

The second one, introduced by Berger and Moerdijk~\cite{BM03}, is a
construction which allows to obtain, from a commutative bialgebra
$\Bca$, a cooperad $\TT \Bca$. Our construction $\T$ and the
construction $\TT$ of these two authors are different but coincide in
many cases. For instance, when $(\Mca, \Product)$ is a monoid such that
for any $x \in \Mca$, the set of pairs $(y, z) \in \Mca^2$ satisfying
$y \Product z = x$ is finite, the operad $\T \Mca$ is the dual of the
cooperad $\TT \Bca$ where $\Bca$ is the dual bialgebra of
$\K \Angle{\Mca}$ endowed with the diagonal coproduct. On the other
hand, there are operads that we can build by the construction $\T$ but
not by the construction $\TT$, and conversely. For  example, the operad
$\T \Z$, where $\Z$ is the additive monoid of integers, cannot be
obtained as the dual of a cooperad built by the construction $\TT$ of
Berger and Moerdijk.
\smallbreak

Furthermore, the operads $\T \Mca$ are defined directly on set-theoretic
bases. Hence, these operads are well-defined in the category of sets and
the computations are explicit. It is therefore possible given a monoid
$\Mca$, to make experiments on the operad $\T \Mca$, using if necessary
a computer. In this chapter, we study many applications of the
construction $\T$ focusing on its combinatorial aspect. More precisely,
we define, by starting from very simple monoids like the additive or
$\max$ monoids of integers, or cyclic monoids, various
nonsymmetric operads involving well-known combinatorial objects.
\smallbreak

This chapter is organized as follows. We begin, in
Section~\ref{sec:construction_T}, by defining the construction $\T$
associating an operad with a monoid and establishing its first
properties. We show that this construction is a functor from the
category of monoids to the category of operads which preserves
injections and surjections. We then apply this construction in
Section~\ref{sec:construction_T_concrete_constructions} to various
monoids and obtain several new (symmetric or not) operads. We construct
in this way some operads on combinatorial objects which were not
provided with such a structure: planar rooted trees with a fixed arity,
Motzkin words, integer compositions, directed animals, and segmented
integer compositions. We also obtain new operads on objects which are
already provided with such a structure: endofunctions, parking
functions, packed words, permutations, planar rooted trees, and Schröder
trees. By using the construction $\T$, we also give an alternative
construction for the diassociative operad~\cite{Lod01} and for the
triassociative operad~\cite{LR04}.
\medbreak

\subsubsection*{Note}
In this chapter, ``operad'' means ``symmetric operad''. To refer to a
nonsymmetric operad, we shall write ``ns operad''.
\medbreak

\section{A functor from monoids to operads} \label{sec:construction_T}
We describe in this section the main ingredient of this chapter, namely
the \Def{construction $\T$}.
\medbreak

\subsection{The operad of a monoid} \label{subsec:operad_of_monoid}
We explain here how the construction $\T$ associates an operad
$\T \Mca$ with any monoid $\Mca$ and an operad morphism
$\T \theta : \T \Mca_1 \to \T \Mca_2$ with any monoid morphism
$\theta : \Mca_1 \to \Mca_2$. We also review some of the main properties
of $\T$.
\medbreak

Let $\Mca$ be a monoid with an associative product $\Product$ admitting
a unit $\Unit_\Mca$. We denote by $\T \Mca$ the space
\begin{equation}
    \T \Mca := \bigoplus_{n \geq 1} \T \Mca(n)
\end{equation}
where for all $n \geq 1$,
\begin{math}
    \T \Mca(n) := \K \Angle{\Mca^n}.
\end{math}
The set $\Mca^+$ forms hence a basis of $\T \Mca$ called
\Def{fundamental basis}. We endow $\T \Mca$ with the partial composition
maps
\begin{equation}
    \circ_i : \T \Mca(n) \otimes \T \Mca(m)
    \to \T \Mca(n + m - 1),
    \qquad n, m \geq 1, i \in [n],
\end{equation}
defined linearly, over the fundamental basis, for any words
$u \in \Mca^n$ and $v \in \Mca^m$ by
\begin{equation}
    u \circ_i v := u_1 \dots u_{i - 1}
    \, (u_i \Product v_1) \, \dots \, (u_i \Product v_m) \,
    u_{i + 1} \dots u_n,
    \qquad i \in [n].
\end{equation}
Moreover, we endow $\T \Mca$ with right actions
\begin{equation}
    \Action : \T \Mca(n) \otimes \Per(n) \to \T \Mca(n),
    \qquad n \geq 1,
\end{equation}
defined linearly, for any permutation $\sigma \in \SymmetricGroup(n)$
and word $u \in \Mca^n$ by
\begin{equation}
    u \Action \sigma := u_{\sigma_1} \dots u_{\sigma_n}.
\end{equation}
In other words, $\T \Mca$ is the vector space of the words on $\Mca$
seen as an alphabet, the partial composition returns to insert a word
$v$ onto the $i$th letter $u_i$ of a word $u$ together with a left
multiplication by $u_i$, and permutations act by permuting the letters
of the words. The arity $|u|$ of an element $u$ of $\T M(n)$ is~$n$.
\medbreak

Now, let $\Mca_1$ and $\Mca_2$ be two monoids and
$\theta : \Mca_1 \to \Mca_2$ be a monoid morphism. Let us denote by
$\T \theta$ the map
\begin{math}
    \T \theta : \T \Mca_1 \to \T \Mca_2
\end{math},
defined for all $u_1 \dots u_n \in \T \Mca(n)$ by
\begin{math}
    \T \theta\left(u_1 \dots u_n\right) :=
    \theta(u_1) \dots \theta(u_n)
\end{math}.
\medbreak

\begin{Theorem} \label{thm:construction_T_functor}
    The construction $\T$ is a functor from the category of monoids with
    monoid morphisms to the category of operads with operad morphisms.
    Moreover, $\T$ preserves injections and surjections.
\end{Theorem}
\medbreak

Observe that the word $\Unit_\Mca \in \T \Mca(1)$ is the unit
of~$\T \Mca$.
\medbreak

Let us consider an example. Let $\Mca := \{\Asf, \Bsf\}^*$ be a free
monoid. Then, $\T \Mca$ is the space of all words whose letters are
words on $\{\Asf, \Bsf\}$. We call such element \Def{multiwords}. For
instance,
\begin{math}
    (\Asf \Asf, \Bsf \Asf, \Bsf, \epsilon, \Asf)
\end{math}
is an element of arity $5$ of $\T \Mca$ and
\begin{equation}
    (\textcolor{Col1}{\Asf \Asf},
    \textcolor{Col1}{\Bsf \Asf}, \mathbf{\Bsf},
    \textcolor{Col1}{\epsilon},
    \textcolor{Col1}{\Asf})
    \circ_3
    (\textcolor{Col4}{\Asf \Bsf},
    \textcolor{Col4}{\epsilon},
    \textcolor{Col4}{\Asf})
    =
    (\textcolor{Col1}{\Asf \Asf},
    \textcolor{Col1}{\Bsf \Asf},
    \mathbf{\Bsf} \textcolor{Col4}{\Asf \Bsf},
    \mathbf{\Bsf}, \mathbf{\Bsf} \textcolor{Col4}{\Asf},
    \textcolor{Col1}{\epsilon}, \textcolor{Col1}{\Asf})
\end{equation}
and
\begin{equation}
    (\Asf \Asf, \Bsf \Asf, \Bsf, \epsilon, \Asf)
    \Action 41352
    =
    (\epsilon, \Asf \Asf, \Bsf, \Asf, \Bsf \Asf).
\end{equation}
Moreover, if $\theta : \Mca \to \N$ is the monoid morphism
defined by $\theta(u) := |u|$, where $\N$ is the additive monoid of
natural numbers, one has
\begin{equation}
    \theta((\Asf \Asf, \Bsf \Asf, \Bsf, \epsilon, \Asf))
    = 22101.
\end{equation}
\medbreak

\subsection{Main properties of the construction}
\label{subsec:main_properties_T}
Let us review the main properties of the construction $\T$.
\medbreak

\begin{Proposition} \label{prop:construction_T_basic}
    Let $\Mca$ be a monoid. The fundamental basis of $\T \Mca$ is a
    basic set-operad basis if and only if $\Mca$ is a right cancellable
    monoid.
\end{Proposition}
\medbreak

\begin{Proposition} \label{prop:construction_T_generating_set}
    Let $\Mca$ be a monoid. Then, the set
    \begin{math}
        \GeneratingSet_\Mca :=
        \Mca \sqcup \left\{\Afr_2\right\}
    \end{math},
    is a generating set of $\T \Mca$ as a ns operad, where
    $\Afr_2 := \Unit_\Mca \Unit_\Mca$ and the elements
    of $\Mca$ are seen as elements of arity $1$ of~$\T \Mca$.
\end{Proposition}
\medbreak

\begin{Proposition} \label{prop:presentation_T}
    Let $\Mca$ be a monoid. The ns operad $\T \Mca$ admits the
    presentation
    $\left(\GeneratingSet_\Mca, \RelationSpace_\Mca\right)$ where
    $\RelationSpace_\Mca$ is the subspace of
    $\FreeOperad\left(\GeneratingSet_\Mca\right)$ generated by the
    elements
    \begin{subequations}
    \begin{equation} \label{equ:relation_T_1}
        \Corolla{\Afr_2} \circ_1 \Corolla{\Afr_2}
        -
        \Corolla{\Afr_2} \circ_2 \Corolla{\Afr_2},
    \end{equation}
    \begin{equation} \label{equ:relation_T_2}
        \Corolla{x} \circ_1 \Corolla{y}
        -
        \Corolla{x \Product y},
        \qquad x, y \in \Mca,
    \end{equation}
    \begin{equation} \label{equ:relation_T_3}
        \Corolla{\Afr_2} \circ
        \left[\Corolla{x} \otimes \Corolla{x}\right]
        -
        \Corolla{x} \circ_1 \Corolla{\Afr_2},
        \qquad x \in \Mca.
    \end{equation}
    \end{subequations}
\end{Proposition}
\medbreak

The proof of Proposition~\ref{prop:presentation_T} relies on a
orientation $\Rew$ of $\RelationSpace_\Mca$ satisfying
\begin{subequations}
\begin{equation}
    \STLeft{\Afr_2}{\Afr_2}
    \Rew
    \STRight{\Afr_2}{\Afr_2},
\end{equation}
\begin{equation}
    \begin{tikzpicture}[scale=.5,Centering]
        \node(2)at(0.00,-2.00){};
        \node[NodeST](0)at(0.00,0.00){\begin{math}x\end{math}};
        \node[NodeST](1)at(0.00,-1.00){\begin{math}y\end{math}};
        \draw[Edge](1)--(0);
        \draw[Edge](2)--(1);
        \node(r)at(0.00,0.75){};
        \draw[Edge](r)--(0);
    \end{tikzpicture}
    \Rew
    \begin{tikzpicture}[scale=.5,Centering]
        \node(1)at(0.00,-1.00){};
        \node[NodeST](0)at(0.00,0.00)
            {\begin{math}x \Product  y\end{math}};
        \draw[Edge](1)--(0);
        \node(r)at(0.00,1){};
        \draw[Edge](r)--(0);
    \end{tikzpicture}\,,
    \qquad x, y \in \Mca,
\end{equation}
\begin{equation}
    \begin{tikzpicture}[scale=.4,Centering]
        \node(1)at(0.00,-2.67){};
        \node(3)at(2.00,-2.67){};
        \node[NodeST](0)at(1.00,0.00){\begin{math}x\end{math}};
        \node[NodeST](2)at(1.00,-1.33)
            {\begin{math}\Afr_2\end{math}};
        \draw[Edge](1)--(2);
        \draw[Edge](2)--(0);
        \draw[Edge](3)--(2);
        \node(r)at(1.00,1.00){};
        \draw[Edge](r)--(0);
    \end{tikzpicture}
    \Rew
    \begin{tikzpicture}[scale=.3,Centering]
        \node(1)at(0.00,-3.33){};
        \node(4)at(2.00,-3.33){};
        \node[NodeST](0)at(0.00,-1.67){\begin{math}x\end{math}};
        \node[NodeST](2)at(1.00,0.00){\begin{math}\Afr_2\end{math}};
        \node[NodeST](3)at(2.00,-1.67){\begin{math}x\end{math}};
        \draw[Edge](0)--(2);
        \draw[Edge](1)--(0);
        \draw[Edge](3)--(2);
        \draw[Edge](4)--(3);
        \node(r)at(1.00,1.5){};
        \draw[Edge](r)--(2);
    \end{tikzpicture}\,,
    \qquad x \in \Mca.
\end{equation}
\end{subequations}
The closure $\RewTrees$ of $\Rew$ is a convergent rewrite rule and its
normal forms of arity $n$ are in one-to-one correspondence with the set
of the words of arity $n$ of~$\T \Mca$.
\medbreak

Let $\Alg$ be an associative algebra with associative product denoted
by $\Conc$, and
\begin{equation}
    \uparrow_x : \Alg \to \Alg, \qquad x \in \Mca,
\end{equation}
be a family of associative algebra morphisms satisfying
\begin{equation} \label{equ:compatible_algebra_morphisms_monoids}
    \uparrow_x \circ \uparrow_y = \uparrow_{x \Product y},
\end{equation}
for all $x, y \in \Mca$. Observe
that~\eqref{equ:compatible_algebra_morphisms_monoids} implies in
particular that $\uparrow_{\Unit_\Mca} = \Identity_\Alg$ where
$\Identity_\Alg$ is the identity map on $\Alg$. We call
\Def{$\Mca$-compatible algebra} such an algebra.
\medbreak

\begin{Theorem} \label{thm:construction_T_algebras}
    Let $\Mca$ be a monoid and $\Alg$ be an $\Mca$-compatible algebra.
    Then, $\Alg$ is an algebra on~$\T \Mca$.
\end{Theorem}
\medbreak

Theorem~\ref{thm:construction_T_algebras} is a consequence of the
presentation of $\T \Mca$ provided by
Proposition~\ref{prop:presentation_T}. Indeed, the associative
product $\Conc$ comes from the generator $\Afr_2$ of $\T \Mca$,
and each map $\uparrow_x$, $x \in \Mca$, comes from the generator
$x \in \Mca$ of $\T \Mca$.
\medbreak

For instance, by Theorem~\ref{thm:construction_T_algebras}, the space
$\K \Angle{\Mca^*}$ of noncommutative polynomials on $\Mca$, endowed
with the associative product~$\Conc$ of concatenation of words of
$\Mca^*$ and with the maps $\uparrow_x$, $x \in \Mca$, defined
linearly for all words $u$ on $\Mca$ by
\begin{equation}
    \uparrow_x(u) :=
    (x \Product u_1) \dots \left(x \Product u_{|u|}\right)
\end{equation}
is an algebra on~$\T \Mca$.
\medbreak

\section{Concrete constructions}
\label{sec:construction_T_concrete_constructions}
Through this section, we consider examples of applications of the
functor $\T$. We shall mainly consider, given a monoid $\Mca$, some
suboperads of $\T \Mca$, symmetric or not, which have for all $n \geq 1$
finitely many elements of arity $n$. For the most part of the
constructed operads $\Oca$, we shall establish isomorphisms of
combinatorial spaces $\phi : \Oca \to \K \Angle{C}$ where the $C$ are
well-chosen combinatorial sets. To this aim, we shall consider
bijections between the basis elements of $\Oca(n)$ and the elements of
size $n$ of $C$, for all $n \geq 1$. Interpreting the partial
compositions of $\Oca$ on $\K \Angle{C}$ amounts to endow
$\K \Angle{C}$ with the structure of an operad, and thus to the
construction of an operad on the objects of $C$. Moreover, we shall
also establish presentations by generators and relations of the
constructed ns operads by using tools from rewrite theory on syntax
trees (see Section~\ref{subsec:rewrite_rules_syntax_trees} of
Chapter~\ref{chap:combinatorics}).
\medbreak

\subsection{Operads from the additive monoid}
\begin{figure}[ht]
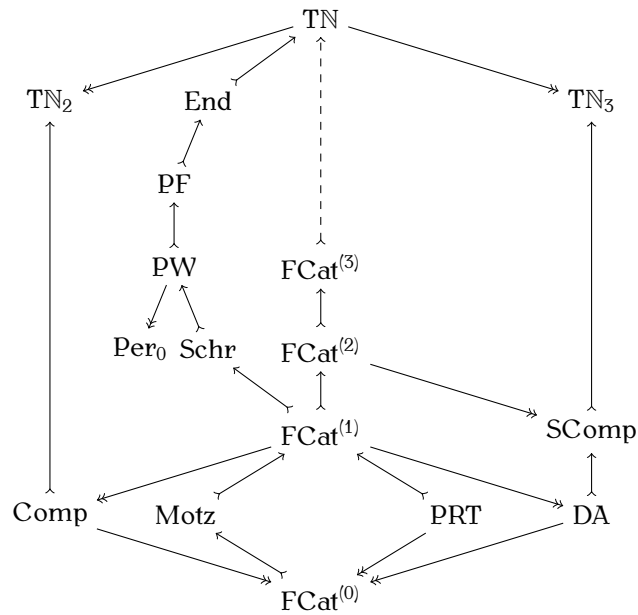

    \centering

    \caption[The diagram of ns suboperads and quotients of the
    operad~$\T \N$.]
    {The diagram of ns suboperads and quotients of $\T \N$. Arrows
    $\rightarrowtail$ (resp. $\twoheadrightarrow$) are injective (resp.
    surjective) ns operad morphisms.}
    \label{fig:construction_TN_diagram}
\end{figure}

Let us denote by $\N$ the additive monoid of natural numbers, and for
all $\ell \geq 1$, by $\N_\ell$ the cyclic additive monoid on
$\Z/_{\ell \Z}$. Note that since, by
Theorem~\ref{thm:construction_T_functor}, $\T$ is a functor which
preserves surjective maps, $\T \N_\ell$ is a quotient operad of $\T \N$.
Besides, since the monoids $\N$ and $\N_\ell$ are right cancellable, by
Proposition~\ref{prop:construction_T_basic}, the fundamental bases of
the operads $\T \N$ and $\T \N_\ell$ are basic set-operad bases. As a
consequence, the fundamental bases of all the suboperads of $\T \N$
and $\T \N_\ell$ constructed in this section are basic set-operad bases.
All these operads fit into the diagram of ns operads represented by
Figure~\ref{fig:construction_TN_diagram}.
Table~\ref{tab:construction_TN_table} summarizes some information about
these ns operads.
\begin{table}[ht]
    \centering
    \begin{tabular}{c|c|c|c|c}
        Monoid & Ns operad & Generators & First dimensions
            & Combinatorial objects \\ \hline \hline
        \multirow{8}{*}{\begin{math}\N\end{math}} & $\End$ & ---
            & $1, 4, 27, 256, 3125$ & Endofunctions \\
        & $\PF$ & --- & $1, 3, 16, 125, 1296$ & Parking functions \\
        & $\PW$ & --- & $1, 3, 13, 75, 541$ & Packed words \\
        & $\PerZero$ & --- & $1, 2, 6, 24, 120$ & Permutations \\
        & $\PRT$ & $01$  & $1, 1, 2, 5, 14, 42$ & Planar rooted trees \\
        & $\FCat{k}$ & $00$, $01$, \dots, $0k$ & Fu\ss-Catalan numbers
            & $k$-ary trees \\
        & $\Schr$ & $00$, $01$, $10$ & $1, 3, 11, 45, 197$
            & Schröder trees \\
        & $\Motz$ & $00$, $010$ & $1, 1, 2, 4, 9, 21, 51$
            & Motzkin words \\ \hline
        $\N_2$ & $\Comp$ & $00$, $01$ & $1, 2, 4, 8, 16, 32$
            & Compositions \\ \hline
        \multirow{2}{*}{\begin{math}\N_3\end{math}} & $\DA$ & $00$, $01$
            & $1, 2, 5, 13, 35, 96$ & Directed animals \\
        & $\SComp$ & $00$, $01$, $02$ & $1, 3, 27, 81, 243$
            & Seg. compositions
    \end{tabular}
    \bigbreak

    \caption[Data about ns suboperads and quotients of the
    operad~$\T \N$.]
    {Ground monoids, generators, first dimensions, and combinatorial
    objects involved in the ns suboperads and quotients of $\T \N$.}
    \label{tab:construction_TN_table}
\end{table}
\medbreak

\subsubsection{Operads on endofunctions, parking functions, packed
words, and permutations}
Recall that an \Def{endofunction} of size $n$ is a word $u$ of length
$n$ on the alphabet $\{1, \dots, n\}$. A \Def{parking function} of size
$n$ is an endofunction $u$ of size $n$ such that the nondecreasing
rearrangement $v$ of $u$ satisfies $v_i \leq i$ for all $i \in [n]$. A
\Def{packed word} of size $n$ is an endofunction $u$ of size $n$ such
that for any letter $u_i \geq 2$ of $u$, there is in $u$ a letter
$u_j = u_i - 1$. Note that neither the set of endofunctions, of parking
functions, of packed words, nor of permutations are suboperads of
$\T \N$. Indeed, one has the following counterexample:
\begin{equation}
    \textcolor{Col1}{1}{\bf 2} \circ_2 \textcolor{Col4}{12} =
    \textcolor{Col1}{1}\textcolor{Col4}{34},
\end{equation}
showing that, even if $12$ is a permutation, $134$ is not an
endofunction.
\medbreak

As a consequence of this observation, let us call a word $u$ a
\Def{twisted} endofunction (resp. parking function, packed word,
permutation) if the word $(u_1 + 1) (u_2 + 1) \dots (u_n + 1)$ is an
endofunction (resp. parking function, packed word, permutation). For
example, the word $2300$ is a twisted endofunction since $3411$ is an
endofunction. Let us denote by $\End$ (resp. $\PF$, $\PW$, $\PerZero$)
the linear span of all twisted endofunctions (resp. parking functions,
packed words, permutations). Under this reformulation, one has the
following result:
\begin{Proposition} \label{prop:construction_T_End_FP_MT}
    The spaces $\End$, $\PF$, and $\PW$ form suboperads of $\T \N$.
\end{Proposition}
\medbreak

For example, we have in $\End$ the partial composition
\begin{equation}
    \textcolor{Col1}{2}{\bf 1} \textcolor{Col1}{23}
    \circ_2
    \textcolor{Col4}{30313} =
    \textcolor{Col1}{2}\textcolor{Col4}{41424}\textcolor{Col1}{23},
\end{equation}
and the action $\Action$ of a permutation
\begin{equation}
    11210 \Action 23514 = 12011.
\end{equation}

Note that $\End$ is not a finitely generated operad. Indeed, the twisted
endofunctions $u$ of size $n$ satisfying $u_i := n - 1$ for all
$i \in [n]$ cannot be obtained by partial compositions involving
elements of $\End$ of arity smaller than $n$. Similarly, $\PF$ is not a
finitely generated operad since the twisted parking functions $u$ of
size $n$ satisfying $u_i := 0$ for all $i \in [n - 1]$ and
$u_n := n - 1$ cannot be obtained by partial compositions involving
elements of $\PF$ of arity smaller than $n$. However, the operad $\PW$
is a finitely generated operad:
\begin{Proposition}  \label{prop:generation_MT}
    The operad $\PW$ is the suboperad of $\T \N$ generated
    by~$\{00, 01\}$.
\end{Proposition}
\medbreak

Let $\Vca$ be the linear span of all twisted packed words having
multiple occurrences of a same letter.
\medbreak

\begin{Proposition} \label{prop:quotient_MT_PerZero}
    The space $\Vca$ is an operad ideal of $\PW$. Moreover, the quotient
    operad $\PW/_\Vca$ is the space $\PerZero$ of twisted
    permutations. Finally, for all twisted permutations $u$ and $v$,
    the partial composition map in $\PerZero$ can be expressed as
    \begin{equation} \label{equ:quotient_MT_PerZero}
        u \circ_i v =
        \begin{cases}
            u \circ_i v & \mbox{if } u_i = |u| - 1, \\
            0_\K & \mbox{otherwise},
        \end{cases}
    \end{equation}
    where $0_\K$ is the null vector of $\PerZero$ and the partial
    composition map $\circ_i$ in the right member
    of~\eqref{equ:quotient_MT_PerZero} is the partial composition map
    of~$\PW$.
\end{Proposition}
\medbreak

Here are two examples of compositions in $\PerZero$:
\begin{subequations}
\begin{equation}
    {\bf 2}\textcolor{Col1}{0431} \circ_1 \textcolor{Col4}{102}
    = 0_\K,
\end{equation}
\begin{equation}
    \textcolor{Col1}{20}{\bf 4}\textcolor{Col1}{31}
    \circ_3
    \textcolor{Col4}{102}
    = \textcolor{Col1}{20}\textcolor{Col4}{546}\textcolor{Col1}{31}.
\end{equation}
\end{subequations}
\medbreak

Let us recall that any minimal generating set of $\Per$ (seen as a ns
operad) has no element of arity $3$ (see
Section~\ref{subsubsec:symmetric_operads} of Chapter~\ref{chap:algebra}).
Moreover, since the homogeneous component of arity $3$ of
$\PerZero^{\{01, 10\}}$ is only of dimension $4$, any minimal
generating set of $\PerZero$ has two elements of arity $3$. Therefore,
$\PerZero$ and $\Per$ are not isomorphic as ns operads.
\medbreak

\subsubsection{A ns operad on planar rooted trees}
Let $\PRT$ be the ns suboperad of $\T \N$ generated by
$\GeneratingSet_\PRT := \{01\}$.
\medbreak

\begin{Proposition} \label{prop:basis_PRT}
    The fundamental basis of $\PRT$ is the set of all the words
    $u$ on the alphabet $\N$ satisfying $u_1 = 0$ and
    $1 \leq u_{i + 1} \leq u_i + 1$ for all $i \in [|u| - 1]$.
\end{Proposition}
\medbreak

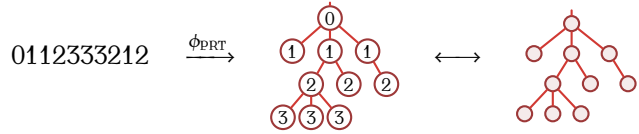
\begin{figure}[ht]
    \centering
    \begin{equation*}
        0112333212
        \quad \xrightarrow{\phi_\PRT} \quad
        \begin{tikzpicture}[xscale=.25,yscale=.22,Centering]
            \node[NodeClear](1)at(0,0){\begin{math}0\end{math}};
            \node[NodeClear](2)at(-2,-2){\begin{math}1\end{math}};
            \node[NodeClear](3)at(0,-2){\begin{math}1\end{math}};
            \node[NodeClear](4)at(-1,-4){\begin{math}2\end{math}};
            \node[NodeClear](5)at(-2.5,-6){\begin{math}3\end{math}};
            \node[NodeClear](6)at(-1,-6){\begin{math}3\end{math}};
            \node[NodeClear](7)at(0.5,-6){\begin{math}3\end{math}};
            \node[NodeClear](8)at(1,-4){\begin{math}2\end{math}};
            \node[NodeClear](9)at(2,-2){\begin{math}1\end{math}};
            \node[NodeClear](10)at(3,-4){\begin{math}2\end{math}};
            \draw[Edge](1)--(2);
            \draw[Edge](1)--(3);
            \draw[Edge](3)--(4);
            \draw[Edge](4)--(5);
            \draw[Edge](4)--(6);
            \draw[Edge](4)--(7);
            \draw[Edge](3)--(8);
            \draw[Edge](1)--(9);
            \draw[Edge](9)--(10);
            \node(r)at(0,1.5){};
            \draw[Edge](1)--(r);
        \end{tikzpicture}
        \quad \longleftrightarrow \quad
        \begin{tikzpicture}[xscale=.25,yscale=.2,Centering]
            \node[Node](1)at(0,0){};
            \node[Node](2)at(-2,-2){};
            \node[Node](3)at(0,-2){};
            \node[Node](4)at(-1,-4){};
            \node[Node](5)at(-2.5,-6){};
            \node[Node](6)at(-1,-6){};
            \node[Node](7)at(0.5,-6){};
            \node[Node](8)at(1,-4){};
            \node[Node](9)at(2,-2){};
            \node[Node](10)at(3,-4){};
            \draw[Edge](1)--(2);
            \draw[Edge](1)--(3);
            \draw[Edge](3)--(4);
            \draw[Edge](4)--(5);
            \draw[Edge](4)--(6);
            \draw[Edge](4)--(7);
            \draw[Edge](3)--(8);
            \draw[Edge](1)--(9);
            \draw[Edge](9)--(10);
            \node(r)at(0,1.5){};
            \draw[Edge](1)--(r);
        \end{tikzpicture}
    \end{equation*}
    \caption[A planar rooted tree of the ns operad $\PRT$.]
    {Interpretation of a word of the ns operad $\PRT$ in terms of a
    planar rooted tree via the bijection $\phi_\PRT$. The nodes of the
    planar rooted tree in the middle are labeled by their depth.}
    \label{fig:interpretation_elements_PRT}
\end{figure}

Let us consider the combinatorial graded collection of all planar rooted
trees where the size of such a tree is its number of nodes (this is the
collection $\PlanarRootedTrees$ defined in
Section~\ref{subsubsec:comb_collection_planar_rooted_trees} of
Chapter~\ref{chap:combinatorics}). There is a bijection $\phi_\PRT$
between the words of $\PRT$ of arity $n$ and planar rooted trees of size
$n$. To compute $\phi_\PRT(u)$ where $u$ is a word of $\PRT$,
iteratively insert the letters of $u$ from left to right according to
the following procedure. If $|u| = 1$, then $u = 0$ and $\phi_\PRT(0)$
is the only planar rooted tree with one node. Otherwise, the insertion
of a letter $\Asf \geq 1$ into a planar rooted tree $\Tfr$ consists in
grafting in $\Tfr$ a new node as the rightmost child of the last node of
depth $\Asf - 1$ for the depth-first traversal of $\Tfr$. The inverse
bijection is computed as follows. Given a planar rooted tree $\Tfr$ of
size $n$, one computes a word of $\PRT$ of arity $n$ by labeling each
node of $\Tfr$ by its depth and then, by reading its labels following a
depth-first traversal of $\Tfr$. Since the words of $\PRT$ satisfy
Proposition~\ref{prop:basis_PRT}, $\phi_\PRT$ is well-defined. Hence, we
can regard the words of arity $n$ of $\PRT$ as planar rooted trees with
$n$ nodes. Figure~\ref{fig:interpretation_elements_PRT} shows an example
of this bijection.
\medbreak

Hence, the Hilbert series of $\PRT$ satisfies
\begin{equation}
    \HilbSeries_\PRT(t) =
    \sum_{n \geq 1}
    \frac{1}{n} \binom{2n - 2}{n - 1} t^n,
\end{equation}
so that its dimensions are the Catalan numbers.
\medbreak

In terms of planar rooted trees, the partial composition of
$\PRT$ can be expressed as follows:
\begin{Proposition} \label{prop:partial_composition_PRT}
    Let $\Sfr$ and $\Tfr$ be two planar rooted trees and $s$ be the
    $i$th node for the depth-first traversal of $\Sfr$. The composition
    $\Sfr \circ_i \Tfr$ in $\PRT$ amounts to replace $s$ by the root of
    $\Tfr$ and graft the children of $s$ as rightmost sons of the root
    of $\Tfr$.
\end{Proposition}
\medbreak

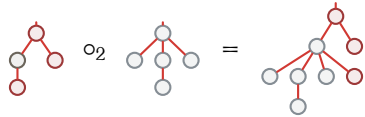
\begin{figure}[ht]
    \centering
    \begin{equation*}
        \begin{tikzpicture}[xscale=.25,yscale=.18,Centering]
            \node[Node](1)at(0,0){};
            \node[NodeColorB](2)at(-1,-2){};
            \node[Node](3)at(-1,-4){};
            \node[Node](4)at(1,-2){};
            \draw[Edge](1)--(2);
            \draw[Edge](2)--(3);
            \draw[Edge](1)--(4);
            \node(r)at(0,1.5){};
            \draw[Edge](1)--(r);
        \end{tikzpicture}
        \enspace \circ_2 \enspace
        \begin{tikzpicture}[xscale=.25,yscale=.18,Centering]
            \node[NodeColorA](1)at(0,0){};
            \node[NodeColorA](2)at(-1.5,-2){};
            \node[NodeColorA](3)at(0,-2){};
            \node[NodeColorA](4)at(0,-4){};
            \node[NodeColorA](5)at(1.5,-2){};
            \draw[Edge](1)--(2);
            \draw[Edge](1)--(3);
            \draw[Edge](3)--(4);
            \draw[Edge](1)--(5);
            \node(r)at(0,1.5){};
            \draw[Edge](1)--(r);
        \end{tikzpicture}
        \enspace = \enspace
        \begin{tikzpicture}[xscale=.25,yscale=.2,Centering]
            \node[Node](1)at(0,0){};
            \node[NodeColorA](2)at(-1,-2){};
            \node[NodeColorA](3)at(-3.5,-4){};
            \node[NodeColorA](4)at(-2,-4){};
            \node[NodeColorA](5)at(-2,-6){};
            \node[NodeColorA](6)at(-.5,-4){};
            \node[Node](7)at(1,-4){};
            \node[Node](8)at(1,-2){};
            \draw[Edge](1)--(2);
            \draw[Edge](2)--(3);
            \draw[Edge](2)--(4);
            \draw[Edge](4)--(5);
            \draw[Edge](2)--(6);
            \draw[Edge](2)--(7);
            \draw[Edge](1)--(8);
            \node(r)at(0,1.5){};
            \draw[Edge](1)--(r);
        \end{tikzpicture}
    \end{equation*}
    \caption[A partial composition of two planar rooted trees of
    $\PRT$.]
    {Interpretation of a partial composition in the ns operad $\PRT$ in
    terms of planar rooted trees.}
    \label{fig:partial_composition_PRT}
\end{figure}

Figure~\ref{fig:partial_composition_PRT} shows an example of a partial
composition in~$\PRT$.
\medbreak

\begin{Proposition} \label{prop:presentation_PRT}
    The ns operad $\PRT$ is isomorphic to the free ns operad generated
    by one element of arity $2$.
\end{Proposition}
\medbreak

Proposition~\ref{prop:presentation_PRT} also says that $\PRT$ is
isomorphic to the ns magmatic operad and hence, that $\PRT$ is a
realization of the magmatic operad. This result is already known since
in~\cite{MY91}, Méndez and Yang point out that the species of
parenthesizations (binary trees) and the species of planar rooted trees
are isomorphic. This isomorphism implies that these species are also
isomorphic as ns operads. Moreover, $\PRT$ can be seen as a planar
version of the \Def{non-associative permutative operad}
$\NAP$~\cite{MY91} (see also~\cite{Liv06}) seen as a ns operad, which
is an operad involving labeled non-planar rooted trees.
\medbreak

\subsubsection{A ns operad on $k$-ary trees}
Let $k \geq 0$ be an integer and $\FCat{k}$ be the ns suboperad
of $\T \N$ generated by
\begin{math}
    \GeneratingSet_{\FCat{k}} := \{00, 01, \dots, 0k\}
\end{math}.
It is immediate from the definition of $\FCat{k}$ that for any
$k \geq 0$, $\FCat{k}$ is a ns suboperad of $\FCat{k + 1}$. Hence, the
ns operads $\FCat{k}$ form an increasing sequence (for the inclusion)
of ns operads. Note that $\FCat{0}$ is isomorphic to the ns associative
operad $\As$. Note also that since $\FCat{1}$ is generated by $00$ and
$01$ and since $\PRT$ is generated by $01$, $\PRT$ is a ns suboperad
of~$\FCat{1}$.
\medbreak

\begin{Proposition} \label{prop:basis_FCat}
    For any $k \geq 0$, the fundamental basis of $\FCat{k}$ is the set
    of all the words $u$ on the alphabet $\N$ satisfying $u_1 = 0$ and
    $0 \leq u_{i + 1} \leq u_i + k$ for all $i \in [|u| - 1]$.
\end{Proposition}
\medbreak

\begin{figure}[ht]
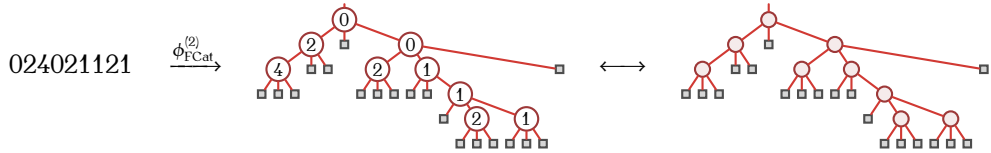

    \centering
    \begin{equation*}
        024021121
        \quad \xrightarrow{\phi_\FCat{2}} \quad

    \end{equation*}
    \caption[A $3$-ary tree of the ns operad $\FCat{2}$.]
    {Interpretation of a word of the ns operad $\FCat{2}$ in terms of
    $3$-ary trees via the bijection $\phi_\FCat{2}$. The $3$-ary tree
    in the middle is well-labeled.}
    \label{fig:interpretation_elements_FCat}
\end{figure}

Let us consider the combinatorial graded collection of the
$k\! +\! 1$-ary trees where the size of such a tree is its number of
internal nodes (this is the collection $\Ary^{(k + 1)}_\Node$ defined
in Section~\ref{subsubsec:k_ary_trees} of
Chapter~\ref{chap:combinatorics}). Let
$\Tfr$ be a $k\! +\! 1$-ary tree of size $n$. We say that an internal
node $x$ of $\Tfr$ is \Def{smaller} than an internal node $y$ if, in
the depth-first traversal of $\Tfr$, $x$ appears before $y$. We also
say that a $k\! +\! 1$-ary tree $\Tfr$ is \Def{well-labeled} if its
root is labeled by $0$, and, for each internal node $x$ of $\Tfr$
labeled by $\Asf$, the children of $x$ which are not leaves are labeled,
from left to right, by $\Asf + k$, \dots, $\Asf + 1$, $\Asf$. There is a
unique way to label a $k\! +\! 1$-ary tree so that it is well-labeled.
There is a bijection $\phi_\FCat{k}$ between the words of $\FCat{k}$ of
arity $n$ and well-labeled $k\! +\! 1$-ary trees of size $n$. To
compute $\phi_\FCat{k}(u)$ where $u$ is a word of $\FCat{k}$,
iteratively insert the letters of $u$ from left to right according to
the following procedure. If $|u| = 1$, then $u = 0$ and
$\phi_\FCat{k}(u)$ is the only well-labeled $k\! +\! 1$-ary tree of size
$1$. Otherwise, the insertion of a letter $\Asf \geq 0$ into a
well-labeled $k\! +\! 1$-ary tree $\Tfr$ consists in replacing a leaf of
$\Tfr$ by the $k\! +\! 1$-ary tree $\Sfr$ of size $1$ labeled by $\Asf$
so that $\Sfr$ is the child of the greatest internal node such that the
obtained $k\! +\! 1$-ary tree is still well-labeled. The inverse
bijection is computed as follows. Given a well-labeled $k\! +\! 1$-ary
tree $\Tfr$, one computes a word of $\FCat{k}$ of arity $n$ by reading
its labels following a depth-first traversal of $\Tfr$. Since the words
of $\FCat{k}$ satisfy Proposition~\ref{prop:basis_FCat}, $\phi_\FCat{k}$
is well-defined. Hence, we can regard the words of arity $n$ of
$\FCat{k}$ as $k\! +\! 1$-ary trees of size $n$.
Figure~\ref{fig:interpretation_elements_FCat} shows an example of this
bijection.
\medbreak

Hence, the Hilbert series of $\FCat{k}$ satisfies the algebraic
relation
\begin{equation}
    1 - \HilbSeries_{\FCat{k}}(t) + t \HilbSeries_{\FCat{k}}(t)^{k + 1}
    = 0,
\end{equation}
so that
\begin{equation}
    \HilbSeries_{\FCat{k}}(t) =
    \sum_{n \geq 1}
    \frac{1}{kn + 1}\binom{kn + n}{n} t^n.
\end{equation}
\medbreak

In terms of $k\! +\! 1$-ary trees, the partial composition of $\FCat{k}$
can be expressed as follows:
\begin{Proposition} \label{prop:partial_composition_FCat}
    Let $\Sfr$ and $\Tfr$ be two $k\! +\! 1$-ary trees and $s$ be the
    $i$th internal node for the depth-first traversal of $\Sfr$. The
    composition $\Sfr \circ_i \Tfr$ in $\FCat{k}$ amounts to replace $s$
    by the root of $\Tfr$ and graft the children of $s$ from right to
    left on the rightmost leaves of~$\Tfr$.
\end{Proposition}
\medbreak

\begin{figure}[ht]
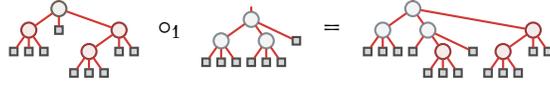

    \centering
    \begin{equation*}

    \end{equation*}
    \caption[A partial composition of two $3$-ary trees of $\FCat{2}$.]
    {Interpretation of the partial composition map of the ns operad
    $\FCat{2}$ in terms of $3$-ary trees.}
    \label{fig:partial_composition_FCat}
\end{figure}

Figure~\ref{fig:partial_composition_FCat} shows an example of
composition in $\FCat{2}$.
\medbreak

\begin{Theorem} \label{thm:presentation_FCat}
    For any $k \geq 0$, the ns operad $\FCat{k}$ admits the presentation
    \begin{math}
        \left(\GeneratingSet_{\FCat{k}}, \RelationSpace\right)
    \end{math}
    where $\RelationSpace$ is the subspace of
    $\FreeOperad\left(\GeneratingSet_{\FCat{k}}\right)$ generated by the
    elements
    \begin{equation}
        \Corolla{0 (\Asf + \Bsf)} \circ_1 \Corolla{0 \Asf}
        -
        \Corolla{0 \Asf} \circ_2 \Corolla{0 \Bsf},
        \qquad \Asf, \Bsf \geq 0, \Asf + \Bsf \leq k.
    \end{equation}
    Moreover, $\FCat{k}$ is a Koszul operad and the set of the
    $\GeneratingSet_{\FCat{k}}$-syntax trees avoiding the trees
    $\Corolla{0 (\Asf + \Bsf)} \circ_1 \Corolla{0 \Asf}$ for all
    $\Asf, \Bsf \geq 0$ and $\Asf + \Bsf \leq k$ is a
    Poincaré-Birkhoff-Witt basis of~$\FCat{k}$.
\end{Theorem}
\medbreak

It is now natural to ask if $\FCat{1}$ is isomorphic to the dendriform
operad $\Dendr$ of to the duplicial operad $\Dup$ since all these ns
operads share the same dimensions. Let us show that $\FCat{1}$ is not
isomorphic to $\Dup$ nor to $\Dendr$. Assume first that
$\phi : \Dup \to \FCat{1}$ is an operad morphism so that
$\phi(\LDup) = \lambda_1 00 + \lambda_2 01$ and
$\phi(\RDup) = \lambda_3 00 + \lambda_4 01$ for some coefficients
$\lambda_1$, $\lambda_2$, $\lambda_3$, and $\lambda_4$ of $\K$. Now,
Relation~\eqref{equ:operad_dup_relation_1} of
Section~\ref{subsubsec:duplicial_operad} of Chapter~\ref{chap:algebra}
of the presentation of $\Dup$ leads to
\begin{equation}
    \phi(\LDup) \circ_1 \phi(\LDup) -
    \phi(\LDup) \circ_2 \phi(\LDup) = 0,
\end{equation}
so that
\begin{equation}
    \lambda_1^2 000 +
    \lambda_1\lambda_2 010 +
    \lambda_1\lambda_2 001 +
    \lambda_2^2 011
    -
    \lambda_1^2 000 -
    \lambda_1\lambda_2 001 -
    \lambda_1\lambda_2 011 -
    \lambda_2^2 012 = 0,
\end{equation}
implying that $\lambda_2 = 0$. Similarly,
Relation~\eqref{equ:operad_dup_relation_3} of the
presentation of $\Dup$ implies that $\lambda_4 = 0$. Therefore, we
obtain $\phi(\LDup) = \lambda_1 00$ and $\phi(\RDup) = \lambda_3 00$,
showing that $\phi$ is not an isomorphism. Assume now that
$\phi : \Dendr \to \FCat{0}$ is an operad morphism so that
$\phi(\LDup) = \lambda_1 00 + \lambda_2 01$ and
$\phi(\RDup) = \lambda_3 00 + \lambda_4 01$ for some coefficients
$\lambda_1$, $\lambda_2$, $\lambda_3$, and $\lambda_4$ of $\K$.
Relation~\eqref{equ:operad_dendr_relation_2} of
Section~\ref{subsubsec:dendriform_operad} of Chapter~\ref{chap:algebra}
of the presentation of $\Dendr$ leads to
\begin{equation}
    \phi(\RDendr) \circ_1 \phi(\LDendr) -
    \phi(\LDendr) \circ_2 \phi(\RDendr) = 0,
\end{equation}
so that
\begin{equation}
    \lambda_1\lambda_3 000 +
    \lambda_2\lambda_3 010 +
    \lambda_1\lambda_4 001 +
    \lambda_2\lambda_4 011
    -
    \lambda_1\lambda_3 000 -
    \lambda_1\lambda_4 001 -
    \lambda_2\lambda_3 011 -
    \lambda_2\lambda_4 012
    = 0,
\end{equation}
implying that $\lambda_2 = 0$, or $\lambda_3 = 0 = \lambda_4$. When
$\lambda_2 = 0$, one has $\phi(\LDendr) = \lambda_1 00$ and, since
$\lambda_1 00$ is associative in $\FCat{1}$ but $\LDendr$ is not
associative in $\Dendr$, $\phi$ cannot be an isomorphism. Moreover,
when $\lambda_3 = 0 = \lambda_4$, the kernel of $\phi$ is nontrivial,
showing that $\phi$ is not an isomorphism.
\medbreak

Since by Theorem~\ref{thm:presentation_FCat}, $\FCat{k}$ is binary and
quadratic, this ns operad admits a Koszul dual. Let ${\FCat{k}}^!$ be the
Koszul dual of $\FCat{k}$.
\medbreak

\begin{Proposition} \label{prop:presentation_FCat_dual}
    For any $k \geq 0$, the ns operad ${\FCat{k}}^!$ admits the
    presentation
    \begin{math}
        \left(\GeneratingSet_{\FCat{k}}, \RelationSpace\right)
    \end{math}
    where $\RelationSpace$ is the subspace of
    $\FreeOperad\left(\GeneratingSet_{\FCat{k}}\right)$ generated by the
    elements
    \begin{subequations}
    \begin{equation}
        \Corolla{0 (\Asf + \Bsf)} \circ_1 \Corolla{0 \Asf}
        -
        \Corolla{0 \Asf} \circ_2 \Corolla{0 \Bsf},
        \qquad \Asf, \Bsf \geq 0, \Asf + \Bsf \leq k,
    \end{equation}
    \begin{equation}
        \Corolla{0 \Asf} \circ_1 \Corolla{0 (\Asf + \Bsf + 1)},
        \qquad 0 \leq \Asf \leq k, 0 \leq \Bsf \leq k - 1,
        \Asf + \Bsf + 1 \leq k,
    \end{equation}
    \begin{equation}
        \Corolla{0 \Asf} \circ_2 \Corolla{0 \Bsf},
        \qquad 0 \leq \Asf \leq k, 0 \leq \Bsf \leq k,
        \Asf + \Bsf \geq k + 1.
    \end{equation}
    \end{subequations}
\end{Proposition}
\medbreak

\begin{Proposition} \label{prop:Hilbert_series_FCat_dual}
    For any $k \geq 0$, the Hilbert series of the ns operad
    ${\FCat{k}}^!$ can be expressed as
    \begin{equation}
        \HilbSeries_{{\FCat{k}}^!}(t) =
        \frac{t}{(1 - t)^{k + 1}}.
    \end{equation}
\end{Proposition}
\medbreak

We deduce from Proposition~\ref{prop:Hilbert_series_FCat_dual} that
\begin{equation}
    \HilbSeries_{{\FCat{k}}^!}(t) =
    \sum_{n \geq 1} \binom{n + k - 1}{k} t^n.
\end{equation}
\medbreak

\subsubsection{A ns operad on Schröder trees}
Let $\Schr$ be the ns suboperad of $\T \N$ generated by
$\GeneratingSet_\Schr := \{00, 01, 10\}$. Since
$\GeneratingSet_\FCat{1} \subseteq \GeneratingSet_\Schr$,
$\FCat{1}$ is a ns suboperad of $\Schr$. Moreover, since $\PW$ is, by
Proposition~\ref{prop:generation_MT}, generated as an operad by
$\{00, 01\}$, $\Schr$ is a ns suboperad of~$\PW$.
\medbreak

\begin{Proposition} \label{prop:basis_Schr}
    The fundamental basis of $\Schr$ is the set of all the words $u$ on
    the alphabet $\N$ having at least one occurrence of $0$ and, for
    all letter $\Bsf \geq 1$ of $u$, there exists a letter
    $\Asf := \Bsf - 1$ such that $u$ has a factor $\Asf v \Bsf$ or
    $\Bsf v \Asf$ where $v$ is a word consisting in letters $\Csf$
    satisfying~$\Csf \geq \Bsf$.
\end{Proposition}
\medbreak

\begin{figure}[ht]
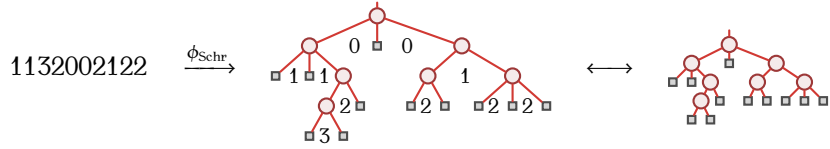

    \centering
    \begin{equation*}
        1132002122
        \quad \xrightarrow{\phi_\Schr} \quad

    \end{equation*}
    \caption[A Schröder tree of the ns operad $\Schr$.]
    {Interpretation of a word of the ns operad $\Schr$ in terms of
    Schröder trees via the bijection $\phi_\Schr$.}
    \label{fig:interpretation_elements_Schr}
\end{figure}

Let us consider the combinatorial graded collection of the Schröder
trees where the size of such a tree is its number of sectors (this is
the collection $\Suspension_{-1}(\Schroder_\Leaf)$ where
$\Schroder_\Leaf$ is the collection defined in
Section~\ref{subsubsec:schroder_trees} of
Chapter~\ref{chap:combinatorics}). There is a bijection $\phi_\Schr$
between the words of $\Schr$ of arity $n$ and Schröder trees of size
$n$. To compute $\phi_\Schr(u)$ where $u$ is a word of $\Schr$,
factorize $u$ as
\begin{math}
    u = u^{(1)} \Asf \dots \Asf u^{(\ell)}
\end{math}
where $\Asf$ is the smallest letter occurring in $u$ and the $u^{(i)}$,
$i \in [\ell]$, are factors of $u$ without $\Asf$. Then, set
\begin{equation}
    \phi_\Schr(u) :=
    \begin{cases}
        \LeafPic & \mbox{if } u = \epsilon, \\
        \Graft\left(\phi_\Schr\left(u^{(1)}\right), \dots,
        \phi_\Schr\left(u^{(\ell)}\right)\right) & \mbox{otherwise},
    \end{cases}
\end{equation}
where $\epsilon$ denotes the empty word and
$\Graft(\Tfr_1, \dots, \Tfr_\ell)$ is the Schröder tree consisting in a
root that has $\Tfr_1$, \dots, $\Tfr_\ell$ as subtrees from left to
right. The inverse bijection is computed as follows. Given a Schröder
tree $\Tfr$, one computes a word of $\Schr$ by assigning to each sector
$(x_i, x_{i + 1})$ of $\Sfr$ the maximal depth of the common ancestors
to the leaves $x_i$ and $x_{i + 1}$. The word of $\Schr$ is obtained by
reading the labels from left to right. Since the words of $\Schr$
satisfy Proposition~\ref{prop:basis_Schr}, $\phi_\Schr$ is well-defined.
Figure~\ref{fig:interpretation_elements_Schr} shows an example of this
bijection.
\medbreak

Hence, the Hilbert series of $\Schr$ satisfies the algebraic relation
\begin{equation}
    t + (3t - 1) \HilbSeries_\Schr(t) + 2t \HilbSeries_\Schr(t)^2 = 0
\end{equation}
so that its first dimensions are
\begin{equation}
    1, 3, 11, 45, 197, 903, 4279, 20793,
\end{equation}
forming Sequence~\OEIS{A001003} of~\cite{Slo}.
\medbreak

\begin{figure}[ht]
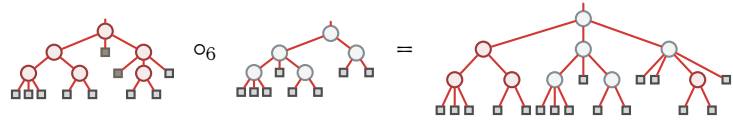

    \centering
    \begin{equation*}

    \end{equation*}
    \caption[A partial composition of two Schröder trees of $\Schr$.]
    {Interpretation of the partial composition map of the ns operad
    $\Schr$ in terms of Schröder trees.}
    \label{fig:partial_composition_Schr}
\end{figure}

It is possible to use the bijection $\phi_\Schr$ to express the partial
composition of $\Schr$ in terms of Schröder trees. We shall not
describe it here but Figure~\ref{fig:partial_composition_Schr} shows an
example of such a composition.
\medbreak

\begin{Theorem} \label{thm:presentation_Schr}
    The ns operad $\Schr$ admits the presentation
    $\left(\GeneratingSet_\Schr, \RelationSpace\right)$ where
    $\RelationSpace$ is the subspace of
    $\FreeOperad\left(\GeneratingSet_\Schr\right)$ generated by the
    elements
    \begin{subequations}
    \begin{equation}
        \Corolla{00} \circ_1 \Corolla{00}
        -
        \Corolla{00} \circ_2 \Corolla{00},
    \end{equation}
    \begin{equation}
        \Corolla{01} \circ_1 \Corolla{10}
        -
        \Corolla{10} \circ_2 \Corolla{01},
    \end{equation}
    \begin{equation}
        \Corolla{00} \circ_1 \Corolla{01}
        -
        \Corolla{00} \circ_2 \Corolla{10},
    \end{equation}
    \begin{equation}
        \Corolla{01} \circ_1 \Corolla{00}
        -
        \Corolla{00} \circ_2 \Corolla{01},
    \end{equation}
    \begin{equation}
        \Corolla{00} \circ_1 \Corolla{10}
        -
        \Corolla{10} \circ_2 \Corolla{00},
    \end{equation}
    \begin{equation}
        \Corolla{01} \circ_1 \Corolla{01}
        -
        \Corolla{01} \circ_2 \Corolla{00},
    \end{equation}
    \begin{equation}
        \Corolla{10} \circ_1 \Corolla{00}
        -
        \Corolla{10} \circ_2 \Corolla{10}.
    \end{equation}
    \end{subequations}
    Moreover, $\Schr$ is a Koszul operad and the set of the
    $\GeneratingSet_\Schr$-syntax trees avoiding the trees
    $\Corolla{00} \circ_1 \Corolla{00}$,
    $\Corolla{01} \circ_1 \Corolla{10}$,
    $\Corolla{00} \circ_1 \Corolla{01}$,
    $\Corolla{01} \circ_1 \Corolla{00}$,
    $\Corolla{00} \circ_1 \Corolla{10}$,
    $\Corolla{01} \circ_1 \Corolla{01}$, and
    $\Corolla{10} \circ_2 \Corolla{10}$ is a Poincaré-Birkhoff-Witt
    basis of~$\Schr$.
\end{Theorem}
\medbreak

Since by Theorem~\ref{thm:presentation_Schr}, $\Schr$ is binary and
quadratic, this ns operad admits a Koszul dual. Let $\Schr^!$ be the
Koszul dual of $\Schr$.
\medbreak

\begin{Proposition} \label{prop:presentation_Schr_dual}
    The ns operad $\Schr^!$ admits the presentation
    \begin{math}
        \left(\GeneratingSet_{\Schr}, \RelationSpace\right)
    \end{math}
    where $\RelationSpace$ is the subspace of
    $\FreeOperad\left(\GeneratingSet_{\Schr}\right)$ generated by the
    elements
    \begin{subequations}
    \begin{equation}
        \Corolla{00} \circ_1 \Corolla{00}
        -
        \Corolla{00} \circ_2 \Corolla{00},
    \end{equation}
    \begin{equation}
        \Corolla{01} \circ_1 \Corolla{10}
        -
        \Corolla{10} \circ_2 \Corolla{01},
    \end{equation}
    \begin{equation}
        \Corolla{00} \circ_1 \Corolla{01}
        -
        \Corolla{00} \circ_2 \Corolla{10},
    \end{equation}
    \begin{equation}
        \Corolla{01} \circ_1 \Corolla{00}
        -
        \Corolla{00} \circ_2 \Corolla{01},
    \end{equation}
    \begin{equation}
        \Corolla{00} \circ_1 \Corolla{10}
        -
        \Corolla{10} \circ_2 \Corolla{00},
    \end{equation}
    \begin{equation}
        \Corolla{01} \circ_1 \Corolla{01}
        -
        \Corolla{01} \circ_2 \Corolla{00},
    \end{equation}
    \begin{equation}
        \Corolla{10} \circ_1 \Corolla{00}
        -
        \Corolla{10} \circ_2 \Corolla{10},
    \end{equation}
    \begin{equation}
        \Corolla{10} \circ_1 \Corolla{01},
    \end{equation}
    \begin{equation}
        \Corolla{10} \circ_1 \Corolla{10},
    \end{equation}
    \begin{equation}
        \Corolla{01} \circ_2 \Corolla{01},
    \end{equation}
    \begin{equation}
        \Corolla{01} \circ_2 \Corolla{10}.
    \end{equation}
    \end{subequations}
\end{Proposition}
\medbreak

\begin{Proposition} \label{prop:Hilbert_series_Schr_dual}
    For any $k \geq 0$, the Hilbert series of the ns operad
    $\Schr^!$ can be expressed as
    \begin{equation}
        \HilbSeries_{\Schr^!}(t) =
        \frac{t}{(1 - t)(1 - 2t)}.
    \end{equation}
\end{Proposition}
\medbreak

We deduce from Proposition~\ref{prop:Hilbert_series_Schr_dual} that
\begin{equation}
    \HilbSeries_{\Schr^!}(t) = \sum_{n \geq 1} (2^n - 1) t^n.
\end{equation}
\medbreak

\subsubsection{A ns operad on Motzkin words}
\label{subsubsec:operad_Motz}
Let $\Motz$ be the ns suboperad of $\T \N$ generated by
$\GeneratingSet_\Motz := \{00, 010\}$. Since $00$ and $010$ are
elements of $\FCat{1}$, $\Motz$ is a ns suboperad of $\FCat{1}$.
Moreover, since
$\GeneratingSet_\FCat{0} \subseteq \GeneratingSet_\Motz$, $\FCat{0}$ is
a ns suboperad of $\Motz$.
\medbreak

\begin{Proposition} \label{prop:basis_Motz}
    The fundamental basis of $\Motz$ is the set of all the words $u$ on
    the alphabet $\N$ beginning and starting by $0$ and such that
    $|u_i - u_{i + 1}| \leq 1$ for all $i \in [|u| - 1]$.
\end{Proposition}
\medbreak

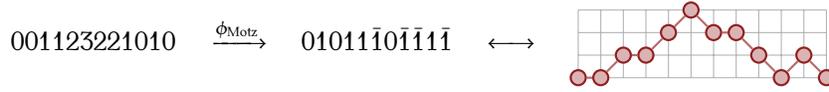
\begin{figure}[ht]
    \centering
    \begin{equation*}
        001123221010
        \quad \xrightarrow{\phi_\Motz} \quad
        0 1 0 1 1 \bar1 0 \bar1 \bar1 1 \bar1
        \quad \longleftrightarrow \quad
        \begin{tikzpicture}[scale=.3,Centering]
            \draw[Grid] (0,0) grid (11,3);
            \node[PathNode](0)at(0,0){};
            \node[PathNode](1)at(1,0){};
            \node[PathNode](2)at(2,1){};
            \node[PathNode](3)at(3,1){};
            \node[PathNode](4)at(4,2){};
            \node[PathNode](5)at(5,3){};
            \node[PathNode](6)at(6,2){};
            \node[PathNode](7)at(7,2){};
            \node[PathNode](8)at(8,1){};
            \node[PathNode](9)at(9,0){};
            \node[PathNode](10)at(10,1){};
            \node[PathNode](11)at(11,0){};
            \draw[PathStep](0)--(1);
            \draw[PathStep](1)--(2);
            \draw[PathStep](2)--(3);
            \draw[PathStep](3)--(4);
            \draw[PathStep](4)--(5);
            \draw[PathStep](5)--(6);
            \draw[PathStep](6)--(7);
            \draw[PathStep](7)--(8);
            \draw[PathStep](8)--(9);
            \draw[PathStep](9)--(10);
            \draw[PathStep](10)--(11);
        \end{tikzpicture}
    \end{equation*}
    \caption[A Motzkin path of the ns operad $\Motz$.]
    {Interpretation of a word of the ns operad $\Motz$ in terms of
    Motzkin words and Motzkin paths via the bijection $\phi_\Motz$.}
    \label{fig:interpretation_elements_Motz}
\end{figure}

A \Def{Motzkin word} is a word $u$ on the alphabet
$\left\{-1, 0, 1\right\}$ such that the sum of all letters of $u$ is $0$
and, for any prefix $u'$ of $u$, the sum of all letters of $u'$ is a
nonnegative integer. The size $|u|$ of a Motzkin word $u$ is its length
plus $1$. In the sequel, we shall denote by $\bar 1$ the letter $-1$.
We can represent a Motzkin word $u$ graphically by a \Def{Motzkin path}
that is the path in $\N^2$ connecting the points $(0, 0)$ and $(n, 0)$
obtained by drawing a step $(1, -1)$ (resp. $(1, 0)$, $(1, 1)$) for
each letter $\bar 1$ (resp. $0$, $1$) of $u$. There is a bijection
$\phi_\Motz$ between the words of $\Motz$ of arity $n$ and Motzkin
words of size $n$. To compute $\phi_\Motz(u)$ where $u$ is an word of
$\Motz(n)$, build the word $v$ of length $n - 1$ satisfying
$v_i := u_{i + 1} - u_i$ for all $i \in [n - 1]$. The inverse bijection
is computed as follows. The word of $\Motz$ in bijection with a Motzkin
word $v$ is the word $u$ such that $u_i$ is the sum of the letters of
the prefix $v_1 \dots v_{i - 1}$  of $v$, for all $i \in [n]$. Since
the words of $\Motz$ satisfy Proposition~\ref{prop:basis_Motz},
$\phi_\Motz$ is well-defined.
Figure~\ref{fig:interpretation_elements_Motz} shows an example of this
bijection.
\medbreak

Hence, the Hilbert series of $\Motz$ satisfies the algebraic relation
\begin{equation}
    t + (t - 1)\HilbSeries_\Motz(t) + t\HilbSeries_\Motz(t)^2 = 0
\end{equation}
so that its first dimensions are
\begin{equation}
    1, 1, 2, 4, 9, 21, 51, 127,
\end{equation}
forming Sequence~\OEIS{A001006}.
\medbreak

In terms of Motzkin words, the partial composition of $\Motz$
can be expressed as follows:
\begin{Proposition} \label{prop:partial_composition_Motz}
    Let $u$ and $v$ be two Motzkin words where $u$ is of size $n$,
    and $i \in [n]$ be an integer. Then the composition $u \circ_i v$
    in $\Motz$ amounts to insert $v$ at the $i$th position into~$u$.
\end{Proposition}
\medbreak

\begin{figure}[ht]
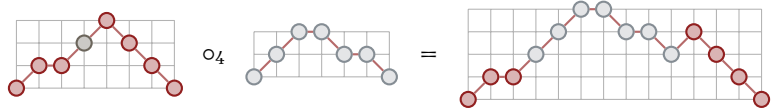

    \centering
    \begin{equation*}

    \end{equation*}
    \caption[A partial composition of two Motzkin paths of $\Motz$.]
    {Interpretation of the partial composition map of the ns operad
    $\Motz$ in terms of Motzkin paths.}
    \label{fig:partial_composition_Motz}
\end{figure}

Figure~\ref{fig:partial_composition_Motz} shows an example of
composition in $\Motz$.
\medbreak

\begin{Theorem} \label{thm:presentation_Motz}
    The ns operad $\Motz$ admits the presentation
    $\left(\GeneratingSet_\Motz, \RelationSpace\right)$ where
    $\RelationSpace$ is the subspace of
    $\FreeOperad\left(\GeneratingSet_\Motz\right)$ generated by the
    elements
    \begin{subequations}
    \begin{equation}
        \Corolla{00} \circ_1 \Corolla{00}
        -
        \Corolla{00} \circ_2 \Corolla{00},
    \end{equation}
    \begin{equation}
        \Corolla{010} \circ_1 \Corolla{00}
        -
        \Corolla{00} \circ_2 \Corolla{010},
    \end{equation}
    \begin{equation}
        \Corolla{00} \circ_1 \Corolla{010}
        -
        \Corolla{010} \circ_3 \Corolla{00},
    \end{equation}
    \begin{equation}
        \Corolla{010} \circ_1 \Corolla{010}
        -
        \Corolla{010} \circ_3 \Corolla{010}.
    \end{equation}
    \end{subequations}
    Moreover, $\Motz$ is a Koszul operad and the set of the
    $\GeneratingSet_\Motz$-syntax trees avoiding the trees
    $\Corolla{00} \circ_1 \Corolla{00}$,
    $\Corolla{010} \circ_1 \Corolla{00}$,
    $\Corolla{00} \circ_1 \Corolla{010}$,
    and $\Corolla{010} \circ_1 \Corolla{010}$ is a
    Poincaré-Birkhoff-Witt basis of $\Motz$.
\end{Theorem}
\medbreak

\subsubsection{A ns operad on compositions}
Let $\Comp$ be the ns suboperad of $\T \N_2$ generated by
$\GeneratingSet_\Comp := \{00, 01\}$. Since $\FCat{1}$ is the ns
suboperad of $\T \N$ generated by
$\GeneratingSet_\FCat{1} = \{00, 01\}$, and since $\T \N_2$ is a
quotient of $\T \N$, $\Comp$ is a quotient of~$\FCat{1}$.
\medbreak

\begin{Proposition} \label{prop:basis_Comp}
    The fundamental basis of $\Comp$ is the set of all the  words on
    the alphabet $\{0, 1\}$ beginning by~$0$.
\end{Proposition}
\begin{proof}
    It is immediate, from the definition of $\Comp$ and
    Lemma~\ref{lem:generation_elements_operads} of
    Chapter~\ref{chap:algebra}, that any element of this ns operad
    starts by $0$ since its generators $00$ and $01$ all start by~$0$.
    \smallbreak

    Let us now show by induction on the length of the words that $\Comp$
    contains any word $u$ satisfying the statement. This is true when
    $|u| = 1$. When $n := |u| \geq 2$, let us observe that if $u$ only
    consists in letters $0$, $\Comp$ contains $u$ because $u$ can be
    obtained by composing the generator $00$ with itself. Otherwise, $u$
    has at least one occurrence of $1$. Since its first letter is $0$,
    there is in $u$ a factor $u_i u_{i + 1} = 01$. By setting
    \begin{math}
        v := u_1 \dots u_i u_{i + 2} \dots u_n
    \end{math},
    we have $u = v \circ_i 01$. Since $v$ satisfies the statement, by
    induction hypothesis $\Comp$ contains $v$. Hence, $\Comp$ also
    contains~$u$.
\end{proof}
\medbreak

Let us consider the combinatorial graded collection of the compositions
where the size of such a composition is the sum of its parts (this
is the collection $\Compositions$ defined in
Section~\ref{subsubsec:integer_compositions} of
Chapter~\ref{chap:combinatorics}. The \Def{$i$th box} of a ribbon
diagram of a composition $\bm{\lambda}$ is the $i$th encountered box by
traversing $\bm{\lambda}$ column by column from left to right and from
top to bottom. The \Def{transpose} of $\bm{\lambda}$ is the ribbon
diagram obtained by applying to $\bm{\lambda}$ the reflection through
the line passing by its first and its last boxes. There is a bijection
$\phi_\Comp$ between the words of $\Comp$ of arity $n$ and ribbon
diagrams of compositions of size $n$. To compute $\phi_\Comp(u)$ where
$u$ is a word of $\Comp$, iteratively insert the letters of $u$ from
left to right according to the following procedure. If $|u| = 1$, then
$u = 0$ and $\phi_\Comp(0)$ is the only ribbon diagram consisting in one
box. Otherwise, the insertion of a letter $\Asf$ into $\bm{\lambda}$
consists in adding a new box below (resp. to the right of) the right
bottommost box of $\bm{\lambda}$ if $\Asf = 1$ (resp. $\Asf = 0$). The
inverse bijection is computed as follows. Given a ribbon diagram
$\bm{\lambda}$ of a composition of size $n$, one computes a word of
$\Comp$ of arity $n$ by labeling the first box of $\bm{\lambda}$ by $0$
and the $i$th box $b$ by $0$ if the $(i - 1)$st box is on the left of
$b$ or by $1$ otherwise, for any $1 \leq i \leq n$. The corresponding
word of $\Comp$ is obtained by reading the labels of $\bm{\lambda}$ from
top to bottom and left to right. Since the words of $\Comp$ satisfy
Proposition~\ref{prop:basis_Comp}, $\phi_\Comp$ is well-defined. Hence,
we can regard the words of arity $n$ of $\Comp$ as ribbon diagrams with
$n$ boxes. Figure~\ref{fig:interpretation_elements_Comp} shows an
example of this bijection.
\begin{figure}[ht]
    \centering
    \begin{equation*}
        0100001011011010
        \quad \xrightarrow{\phi_\Comp} \quad

    \end{equation*}
    \caption[A ribbon diagram of the ns operad $\Comp$.]
    {Interpretation of a word of the ns operad $\Comp$ in terms of
    compositions via the bijection $\phi_\Comp$. Boxes of the ribbon
    diagram in the middle are labeled.}
    \label{fig:interpretation_elements_Comp}
\end{figure}
\medbreak

Hence, the Hilbert series of $\Comp$ satisfies
\begin{equation}
    \HilbSeries_\Comp(t) = \sum_{n \geq 1} 2^{n - 1} t^n.
\end{equation}
\medbreak

In terms of ribbon diagrams, the partial composition of $\Comp$
can be expressed as follows:
\begin{Proposition} \label{prop:partial_composition_Comp}
    Let $\bm{\lambda}$ and $\bm{\mu}$ be two ribbon diagrams, $i$ be an
    integer, and $c$ be the $i$th box of $\bm{\lambda}$. Then, the
    composition $\bm{\lambda} \circ_i \bm{\mu}$ in $\Comp$ amounts to
    replace $c$ by $\bm{\mu}$ if $c$ is the upper box of its column, or
    to replace $c$ by the transpose ribbon diagram of $\bm{\mu}$
    otherwise.
\end{Proposition}
\medbreak

Figure~\ref{fig:partial_composition_Comp} shows two examples
of compositions in $\Comp$.
\begin{figure}[ht]
    \centering
    \subfloat[]
    [Composition in an upper box of a column.]{
    \begin{minipage}[c]{.45\textwidth}
    \begin{equation*}

    \end{equation*}
    \end{minipage}
    \label{subfig:partial_composition_Comp_1}}
    \qquad
    \subfloat[]
    [Composition in a non-upper box of a column.]{
    \begin{minipage}[c]{.45\textwidth}
    \begin{equation*}
        %
    \end{equation*}
    \end{minipage}
    \label{subfig:partial_composition_Comp_2}}
    \caption[Two partial compositions of ribbon diagrams of $\Comp$.]
    {Interpretation of the partial composition map of the ns operad
    $\Comp$ in terms of ribbon diagrams.}
    \label{fig:partial_composition_Comp}
\end{figure}

\begin{Theorem} \label{thm:presentation_Comp}
    The ns operad $\Comp$ admits the presentation
    $\left(\GeneratingSet_\Comp, \RelationSpace\right)$ where
    $\RelationSpace$ is the subspace of
    $\FreeOperad\left(\GeneratingSet_\Comp\right)$ generated by the
    elements
    \begin{subequations}
    \begin{equation} \label{equ:relation_Comp_1}
        \Corolla{00} \circ_1 \Corolla{00}
        -
        \Corolla{00} \circ_2 \Corolla{00},
    \end{equation}
    \begin{equation}
        \Corolla{01} \circ_1 \Corolla{00}
        -
        \Corolla{00} \circ_2 \Corolla{01},
    \end{equation}
    \begin{equation}
        \Corolla{01} \circ_1 \Corolla{01}
        -
        \Corolla{01} \circ_2 \Corolla{00},
    \end{equation}
    \begin{equation} \label{equ:relation_Comp_4}
        \Corolla{00} \circ_1 \Corolla{01}
        -
        \Corolla{01} \circ_2 \Corolla{01}.
    \end{equation}
    \end{subequations}
    Moreover, $\Comp$ is a Koszul operad and the set of the
    $\GeneratingSet_\Comp$-syntax trees avoiding the trees
    $\Corolla{00} \circ_1 \Corolla{00}$,
    $\Corolla{01} \circ_1 \Corolla{00}$,
    $\Corolla{01} \circ_1 \Corolla{01}$,
    and $\Corolla{00} \circ_1 \Corolla{01}$ is a Poincaré-Birkhoff-Witt
    basis of~$\Comp$.
\end{Theorem}
\begin{proof}
    Observe first that since the evaluations of all the
    elements~\eqref{equ:relation_Comp_1}---\eqref{equ:relation_Comp_4}
    are $0$, for all element $x$ of $\RelationSpace$, $\Eval(x) = 0$.
    Let now $\Rew$ be the rewrite rule, being an orientation of
    $\RelationSpace$, defined by
    \begin{subequations}
    \begin{equation}
        \STLeft{00}{00} \Rew \STRight{00}{00},
    \end{equation}
    \begin{equation}
        \STLeft{00}{01} \Rew \STRight{00}{01},
    \end{equation}
    \begin{equation}
        \STLeft{01}{01} \Rew \STRight{01}{00},
    \end{equation}
    \begin{equation}
        \STLeft{01}{00} \Rew \STRight{01}{01}.
    \end{equation}
    \end{subequations}
    Let $\RewTrees$ be the closure of $\Rew$. It is immediate that the
    map $\TamariInvariant$ (see~\eqref{equ:Tamari_invariant} of
    Chapter~\ref{chap:combinatorics}) is a termination invariant for
    $\RewTrees$ so that $\RewTrees$ is terminating. Moreover, the normal
    forms of $\RewTrees$ are the $\GeneratingSet_\Comp$-syntax trees
    which have no internal node with an internal node as left child.
    Hence, the generating series
    $\GenSeries_{\NormalForms_\RewTrees}(t)$ of the normal forms of
    $\RewTrees$ satisfies
    \begin{equation}
        \GenSeries_{\NormalForms_\RewTrees}(t) =
        \sum_{n \geq 1} 2^{n - 1} t^n.
    \end{equation}
    By Proposition~\ref{prop:basis_Comp},
    $\GenSeries_{\NormalForms_\RewTrees}(t)$ also is the Hilbert series
    of $\Comp$. Hence, by Theorem~\ref{thm:presentation_operads} of
    Chapter~\ref{chap:algebra}, $\Comp$ admits the claimed presentation.
\end{proof}
\medbreak

Since by Theorem~\ref{thm:presentation_Comp}, $\Comp$ is binary and
quadratic, this ns operad admits a Koszul dual. Let $\Comp^!$ be the
Koszul dual of~$\Comp$.
\medbreak

\begin{Proposition} \label{prop:presentation_Comp_dual}
    The ns operad $\Comp$ is self-dual, that is $\Comp \simeq \Comp^!$.
\end{Proposition}
\medbreak

\subsubsection{A ns operad on directed animals}
Let $\DA$ be the ns suboperad of $\T \N_3$ generated by
$\GeneratingSet_\DA := \{00, 01\}$. Since $\FCat{1}$ is the ns
suboperad of $\T \N$ generated by
$\GeneratingSet_\FCat{1} = \{00, 01\}$, and since $\T \N_3$ is a
quotient of $\T \N$, $\DA$ is a quotient of $\FCat{1}$. We denote here
by $\bar 1$ the representative of the equivalence class of $2$
in~$\N_3$.
\medbreak

\begin{Proposition} \label{prop:bijection_words_DA_motzkin_prefixes}
    Let, for any $n \geq 1$, $\Bca(n)$ be the set of the words of
    forming the fundamental basis of $\DA(n)$ and let
    $\phi_\DA : \Bca(n) \to \left\{\bar 1, 0, 1\right\}^{n - 1}$ be
    the map defined for any word $u$ of arity $n$ of $\DA$ by
    \begin{equation}
        \phi_\DA(u) := (u_1 * u_2) \; (u_2 * u_3) \dots
        \left(u_{n - 1} * u_n\right),
    \end{equation}
    where $u_i * u_{i + 1} := u_{i + 1} - u_i \mod 3$ for all
    $i \in [n - 1]$. Then, $\phi_\DA$ is a bijection between the words
    of arity $n$ of $\DA$ and prefixes of Motzkin words of
    length~$n - 1$.
\end{Proposition}
\medbreak

Here are two examples of images by $\phi_\DA$ of words of $\DA$:
\begin{subequations}
\begin{equation}
    \phi_\DA\left(011\bar 1\bar 10\bar 101\right) = 10101\bar111,
\end{equation}
\begin{equation}
    \phi_\DA\left(010010101\bar 11\right) = 1\bar101\bar11\bar111\bar1.
\end{equation}
\end{subequations}
\medbreak

\begin{figure}[ht]
    \centering
    \begin{tikzpicture}[scale=.3,Centering]
        \draw[Grid] (0,0) grid (7,6);
        \node[PathNode](0)at(0,0){};
        \node[PathNode](1)at(1,0){};
        \node[PathNode](2)at(1,1){};
        \node[PathNode](3)at(1,2){};
        \node[PathNode](4)at(1,3){};
        \node[PathNode](5)at(1,4){};
        \node[PathNode](6)at(1,5){};
        \node[PathNode](7)at(1,6){};
        \node[PathNode](8)at(2,1){};
        \node[PathNode](9)at(2,2){};
        \node[PathNode](10)at(2,4){};
        \node[PathNode](11)at(3,2){};
        \node[PathNode](12)at(3,3){};
        \node[PathNode](13)at(3,4){};
        \node[PathNode](14)at(4,3){};
        \node[PathNode](15)at(5,3){};
        \node[PathNode](16)at(5,4){};
        \node[PathNode](17)at(5,5){};
        \node[PathNode](18)at(6,3){};
        \node[PathNode](19)at(6,5){};
        \node[PathNode](20)at(7,5){};
        \draw[PathStep](0)--(1);
        \draw[PathStep](1)--(2);
        \draw[PathStep](2)--(3);
        \draw[PathStep](3)--(4);
        \draw[PathStep](4)--(5);
        \draw[PathStep](5)--(6);
        \draw[PathStep](6)--(7);
        \draw[PathStep](2)--(8);
        \draw[PathStep](5)--(10);
        \draw[PathStep](10)--(13);
        \draw[PathStep](8)--(9);
        \draw[PathStep](9)--(11);
        \draw[PathStep](11)--(12);
        \draw[PathStep](12)--(13);
        \draw[PathStep](12)--(14);
        \draw[PathStep](14)--(15);
        \draw[PathStep](15)--(18);
        \draw[PathStep](15)--(16);
        \draw[PathStep](16)--(17);
        \draw[PathStep](17)--(19);
        \draw[PathStep](19)--(20);
        \draw[PathStep](3)--(9);
    \end{tikzpicture}
    \caption[A directed animal.]
    {A directed animal of size $21$. The point $(0, 0)$ is the lowest
    and leftmost point.}
    \label{fig:directed_animal}
\end{figure}
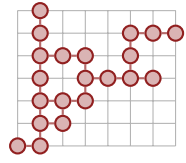

A \Def{directed animal} is a subset $A$ of $\N^2$ such that
$(0, 0) \in A$ and $(i, j) \in A$ with $i \geq 1$ or $j \geq 1$ implies
$(i - 1, j) \in A$ or $(i, j - 1) \in A$. The size of a directed animal
$A$ is its cardinality. Figure~\ref{fig:directed_animal} shows a
directed animal.
\medbreak

According to~\cite{GV88}, there is a bijection $\alpha$ between the set
of prefixes of Motzkin words of length $n - 1$ and the set of directed
animals of size $n$. Hence, by
Proposition~\ref{prop:bijection_words_DA_motzkin_prefixes}, the map
$\alpha \circ \phi_\DA$ is a bijection between the words of $\DA$ of
arity $n$ and directed animals of size $n$. Therefore, $\DA$ can be seen
as a ns operad on directed animals.
\medbreak

Hence, the Hilbert series of $\DA$ satisfies the algebraic relation
\begin{equation}
    t + (3t - 1) \HilbSeries_\DA(t) + (3t - 1) \HilbSeries_\DA(t)^2 = 0
\end{equation}
so that its first dimensions are
\begin{equation}
    1, 2, 5, 13, 35, 96, 267, 750, 2123,
\end{equation}
forming Sequence~\OEIS{A005773} of~\cite{Slo}.
\medbreak

\begin{Theorem} \label{thm:presentation_DA}
    The ns operad $\DA$ admits the presentation
    $\left(\GeneratingSet_\DA, \RelationSpace\right)$ where
    $\RelationSpace$ is the subspace of
    $\FreeOperad\left(\GeneratingSet_\DA\right)$ generated by the
    elements
    \begin{subequations}
    \begin{equation}
        \Corolla{00} \circ_1 \Corolla{00}
        -
        \Corolla{00} \circ_2 \Corolla{00},
    \end{equation}
    \begin{equation}
        \Corolla{01} \circ_1 \Corolla{00}
        -
        \Corolla{00} \circ_2 \Corolla{01},
    \end{equation}
    \begin{equation}
        \Corolla{01} \circ_1 \Corolla{01}
        -
        \Corolla{01} \circ_2 \Corolla{00},
    \end{equation}
    \begin{equation} \label{equ:relation_DA_3}
        (\Corolla{00} \circ_1 \Corolla{01}) \circ_2 \Corolla{01}
        -
        (\Corolla{01} \circ_2 \Corolla{01}) \circ_3 \Corolla{01}.
    \end{equation}
    \end{subequations}
    Moreover, the set of the $\GeneratingSet_\DA$-syntax trees avoiding
    the trees $\Corolla{00} \circ_1 \Corolla{00}$,
    $\Corolla{01} \circ_1 \Corolla{00}$,
    $\Corolla{01} \circ_2 \Corolla{00}$,
    $(\Corolla{01} \circ_2 \Corolla{01}) \circ_3 \Corolla{01}$ is
    a Poincaré-Birkhoff-Witt basis of~$\DA$.
\end{Theorem}
\medbreak

Since the nontrivial relation~\eqref{equ:relation_DA_3} has degree $3$,
the presentation of $\DA$ exhibited by
Theorem~\ref{thm:presentation_DA} is not quadratic.
\medbreak

\subsubsection{A ns operad on segmented compositions}
Let $\SComp$ be the ns suboperad of $\T \N_3$ generated by
$\GeneratingSet_\SComp := \{00, 01, 02\}$. Since $\FCat{2}$ is the ns
suboperad of $\T \N$ generated by
$\GeneratingSet_\FCat{2} = \{00, 01, 02\}$, and since $\T \N_3$ is a
quotient of $\T \N$, $\SComp$ is a quotient of $\FCat{2}$. Moreover,
since $\DA$ is generated by
$\GeneratingSet_\DA \subseteq \GeneratingSet_\SComp$, $\DA$ is a ns
suboperad of $\SComp$.
\medbreak

\begin{Proposition} \label{prop:basis_SComp}
    The fundamental basis of $\SComp$ is the set of all the words on
    the alphabet $\{0, 1, 2\}$ beginning by~$0$.
\end{Proposition}
\medbreak

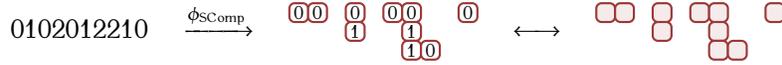
\begin{figure}[ht]
    \centering
    \begin{equation*}
        0102012210
        \quad \xrightarrow{\phi_\SComp} \quad
        \begin{tikzpicture}[Centering]
            \node[BoxClear](1)at(0,0){\begin{math}0\end{math}};
            \node[BoxClear,right of=1,node distance=0.25cm](2)
                {\begin{math}0\end{math}};
            \node[BoxClear,right of=0,node distance=.75cm](3)
                {\begin{math}0\end{math}};
            \node[BoxClear,below of=3,node distance=0.25cm](4)
                {\begin{math}1\end{math}};
            \node[BoxClear,right of=0,node distance=1.25cm](5)
                {\begin{math}0\end{math}};
            \node[BoxClear,right of=5,node distance=0.25cm](6)
                {\begin{math}0\end{math}};
            \node[BoxClear,below of=6,node distance=0.25cm](7)
                {\begin{math}1\end{math}};
            \node[BoxClear,below of=7,node distance=0.25cm](8)
                {\begin{math}1\end{math}};
            \node[BoxClear,right of=8,node distance=0.25cm](9)
                {\begin{math}0\end{math}};
            \node[BoxClear,right of=1,node distance=2.25cm](10)
                {\begin{math}0\end{math}};
        \end{tikzpicture}
        \quad \longleftrightarrow \quad
        \begin{tikzpicture}[Centering]
            \node[Box](1)at(0,0){};
            \node[Box,right of=1,node distance=0.25cm](2){};
            \node[Box,right of=0,node distance=.75cm](3){};
            \node[Box,below of=3,node distance=0.25cm](4){};
            \node[Box,right of=0,node distance=1.25cm](5){};
            \node[Box,right of=5,node distance=0.25cm](6){};
            \node[Box,below of=6,node distance=0.25cm](7){};
            \node[Box,below of=7,node distance=0.25cm](8){};
            \node[Box,right of=8,node distance=0.25cm](9){};
            \node[Box,right of=1,node distance=2.25cm](10){};
        \end{tikzpicture}
    \end{equation*}
    \caption[A ribbon diagram of the ns operad $\SComp$.]
    {Interpretation of a word of the operad $\SComp$ in terms of a
    segmented composition via the bijection $\phi_\SComp$. Boxes of the
    ribbon diagram in the middle are labeled.}
    \label{fig:interpretation_elements_SComp}
\end{figure}

A \Def{segmented composition} is a sequence
\begin{math}
    (\bm{\lambda}_1, \dots, \bm{\lambda}_\ell)
\end{math}
of compositions $\bm{\lambda}_i$, $i \in [\ell]$. The size of a
segmented composition is the sum of the sizes of the compositions
constituting it. We shall represent a segmented composition
$\bm{\lambda}$ by a \Def{ribbon diagram}, that is the diagram
consisting in the sequence of the ribbon diagrams of the compositions
that constitute $\bm{\lambda}$. There is a bijection between the words
of $\SComp$ of arity $n$ and ribbon diagrams of segmented compositions
of size $n$. To compute $\phi_\SComp(u)$ where $u$ is a word of
$\SComp$, factorize $u$ as $x = 0 u^{(1)} \dots 0 u^{(\ell)}$ such that
for any $i \in [\ell]$, the factor $u^{(i)}$ has no occurrence of $0$,
and compute the sequence
\begin{math}
    \left(\phi_\Comp\left(0\bar u^{(1)}\right), \dots,
    \phi_\Comp\left(0\bar u^{(\ell)}\right)\right)
\end{math},
where for any $i \in [\ell]$, $\bar u^{(i)}$ is the word obtained from
$u^{(i)}$ by decreasing all letters. The inverse bijection is computed
as follows. Given a ribbon diagram
\begin{math}
    \bm{\lambda} := (\bm{\lambda}_1, \dots, \bm{\lambda}_\ell)
\end{math}
of a segmented composition of size $n$, one computes a word of $\SComp$
of arity $n$ by computing the sequence
$\left(u^{(1)}, \dots, u^{(\ell)}\right)$ where for any
$i \in [|\ell|]$, $u^{(i)}$ is the word of $\Comp$ obtained by applying
the inverse bijection of $\phi_\Comp$ to $\bm{\lambda}_i$, then by
incrementing in each $u^{(i)}$ all letters, excepted the first one, and
finally by concatenating the words of the sequence together. Since the
words of $\SComp$ satisfy Proposition~\ref{prop:basis_SComp},
$\phi_\SComp$ is well-defined.
Figure~\ref{fig:interpretation_elements_SComp} shows an example of this
bijection.
\medbreak

Hence, the Hilbert of $\SComp$ satisfies
\begin{equation}
    \HilbSeries_\SComp(t) =
    \sum_{n \geq 1} 3^{n - 1} t^n.
\end{equation}
\medbreak

In terms of ribbon diagrams, the partial composition of $\SComp$
can be expressed as follows:
\begin{Theorem} \label{thm:presentation_SComp}
    The ns operad $\SComp$ admits the presentation
    $\left(\GeneratingSet_\SComp, \RelationSpace\right)$ where
    $\RelationSpace$ is the subspace of
    $\FreeOperad\left(\GeneratingSet_\SComp\right)$ generated by the
    elements
    \begin{subequations}
    \begin{equation}
        \Corolla{00} \circ_1 \Corolla{00}
        -
        \Corolla{00} \circ_2 \Corolla{00},
    \end{equation}
    \begin{equation}
        \Corolla{01} \circ_1 \Corolla{00}
        -
        \Corolla{00} \circ_2 \Corolla{01},
    \end{equation}
    \begin{equation}
        \Corolla{01} \circ_1 \Corolla{01}
        -
        \Corolla{01} \circ_2 \Corolla{00},
    \end{equation}
    \begin{equation}
        \Corolla{00} \circ_1 \Corolla{01}
        -
        \Corolla{01} \circ_2 \Corolla{02},
    \end{equation}
    \begin{equation}
        \Corolla{01} \circ_1 \Corolla{02}
        -
        \Corolla{02} \circ_2 \Corolla{02},
    \end{equation}
    \begin{equation}
        \Corolla{00} \circ_1 \Corolla{02}
        -
        \Corolla{02} \circ_2 \Corolla{01},
    \end{equation}
    \begin{equation}
        \Corolla{02} \circ_1 \Corolla{00}
        -
        \Corolla{00} \circ_2 \Corolla{02},
    \end{equation}
    \begin{equation}
        \Corolla{02} \circ_1 \Corolla{01}
        -
        \Corolla{01} \circ_2 \Corolla{01},
    \end{equation}
    \begin{equation}
        \Corolla{02} \circ_1 \Corolla{02}
        -
        \Corolla{02} \circ_2 \Corolla{00}.
    \end{equation}
    \end{subequations}
    Moreover, $\SComp$ is a Koszul operad and the set of the
    $\GeneratingSet_\SComp$-syntax trees avoiding the trees
    $\Corolla{00} \circ_1 \Corolla{00}$,
    $\Corolla{01} \circ_1 \Corolla{00}$,
    $\Corolla{01} \circ_1 \Corolla{01}$,
    $\Corolla{00} \circ_1 \Corolla{01}$,
    $\Corolla{01} \circ_1 \Corolla{02}$,
    $\Corolla{00} \circ_1 \Corolla{02}$,
    $\Corolla{02} \circ_1 \Corolla{00}$,
    $\Corolla{02} \circ_1 \Corolla{01}$,
    and $\Corolla{02} \circ_1 \Corolla{02}$ is a Poincaré-Birkhoff-Witt
    basis of~$\SComp$.
\end{Theorem}
\medbreak

Since by Theorem~\ref{thm:presentation_SComp}, $\SComp$ is binary and
quadratic, this ns operad admits a Koszul dual. Let $\SComp^!$ be the
Koszul dual of~$\SComp$.
\medbreak

\begin{Proposition} \label{prop:presentation_SComp_dual}
    The ns operad $\SComp$ is self-dual, that is
    $\SComp \simeq \SComp^!$.
\end{Proposition}
\medbreak

\subsection{Operads from the $\max$ monoid}
We shall denote by $\Nmax$ the monoid $\N$ with the binary operation
$\max$ as product. Note that the ns suboperad of $\T \Nmax$ generated
by $\{\Asf \Asf\}$ where $\Asf \in \Nmax$ are all isomorphic to the ns
associative operad $\As$. The operads constructed in this section fit
into the diagram of ns operads represented by
Figure~\ref{fig:construction_TM_diagram}.
Table~\ref{tab:construction_TM_table} summarizes some information about
these ns operads.
\begin{figure}[ht]
    \centering
    \begin{tikzpicture}[scale=.5]
        \node(TN)at(0,0){\begin{math}\T \Nmax\end{math}};
        \node(Tr)at(0,-2){\begin{math}\Tr\end{math}};
        \node(D)at(0,-4){\begin{math}\Di\end{math}};
        \node(As)at(0,-6){\begin{math}\As\end{math}};
        \draw[Injection](Tr)--(TN);
        \draw[Injection](D)--(Tr);
        \draw[Surjection](D)--(As);
    \end{tikzpicture}
    \caption[The diagram of ns suboperads and quotients of the
    operad~$\T \Nmax$.]
    {The diagram of ns suboperads and quotients of $\T \Nmax$. Arrows
    $\rightarrowtail$ (resp. $\twoheadrightarrow$) are injective (resp.
    surjective) ns operad morphisms.}
    \label{fig:construction_TM_diagram}
\end{figure}
\begin{table}[ht]
    \centering
    \begin{tabular}{c|c|c|c|c}
        Monoid & Ns operad & Generators & First dimensions
            & Combinatorial objects \\ \hline \hline
        \multirow{2}{*}{\begin{math}\Nmax\end{math}}
            & $\Di$ & $01$, $10$ & $1, 2, 3, 4, 5$
            & Binary words with exactly one $0$ \\
            & $\Tr$ & $00$, $01$, $10$ & $1, 3, 7, 15, 31$
            & Binary words with at least one $0$
    \end{tabular}
    \bigbreak

    \caption[Data about ns suboperads of the ns operad~$\T \Nmax$.]
    {Ground monoids, generators, first dimensions, and combinatorial
    objects involved in the ns suboperads of~$\T \Nmax$.}
    \label{tab:construction_TM_table}
\end{table}
\medbreak

\subsubsection{The diassociative operad}
Let $\Di$ be the ns suboperad of $\T \Nmax$ generated by
$\GeneratingSet_\Di := \{01, 10\}$.
\medbreak

\begin{Proposition} \label{prop:basis_Di}
    The fundamental basis of $\Di$ is the set of all the words on the
    alphabet $\{0, 1\}$ containing exactly one~$0$.
\end{Proposition}
\medbreak

The diassociative operad $\Dias$~\cite{Lod01} is a ns operad whose
definition is recalled in Section~\ref{subsubsec:diassociative_operad}
of Chapter~\ref{chap:algebra}.
\medbreak

\begin{Proposition} \label{prop:isomorphism_Dias_Di}
    The ns operad $\Di$ and the ns diassociative operad $\Dias$ are
    isomorphic and the map
    \begin{equation}
        \phi : \Dias \to \Di
    \end{equation}
    satisfying $\phi(\dashv) = 01$ and $\phi(\vdash) = 10$ is an
    isomorphism of operads.
\end{Proposition}
\medbreak

Proposition~\ref{prop:isomorphism_Dias_Di} also shows that $\Di$ is a
realization of the ns diassociative operad.
\medbreak

\subsubsection{The triassociative operad}
Let $\Tr$ be the ns suboperad of $\T \Nmax$ generated by
$\GeneratingSet_\Tr := \{00, 01, 10\}$. Since
$\GeneratingSet_\Di \subseteq \GeneratingSet_\Tr$, $\Di$ is a ns
suboperad of $\Tr$.
\medbreak

\begin{Proposition} \label{prop:basis_Tr}
    The fundamental basis of $\Tr$ is the set of all the words on the
    alphabet $\{0, 1\}$ containing at least one~$0$.
\end{Proposition}
\medbreak

The triassociative operad $\Trias$ is a ns operad introduced
in~\cite{LR04}.
\medbreak

\begin{Proposition} \label{prop:isomorphism_Trias_Tr}
    The ns operad $\Tr$ and the ns triassociative operad $\Trias$ are
    isomorphic and the map
    \begin{equation}
        \phi : \Trias \to \Tr
    \end{equation}
    satisfying $\phi(\dashv) = 01$, $\phi(\vdash) = 10$, and
    $\phi(\perp) = 00$ is an isomorphism.
\end{Proposition}
\medbreak

Proposition~\ref{prop:isomorphism_Trias_Tr} also shows that $\Tr$ is a
realization of the ns triassociative operad.
\medbreak

\section*{Concluding remarks}
We have presented here the functorial construction $\T$ producing an
operad given a monoid. As we have seen, this construction is very rich
from a combinatorial point of view since most of the obtained operads
coming from usual monoids involve a wide range of combinatorial
objects. There are various ways to continue this work. Let us address
here the main directions.
\smallbreak

In the first place, it appears that we have somewhat neglected the fact
that $\T$ is a functor to operads and not only to ns ones. Indeed,
except for the operads $\End$, $\PF$, $\PW$, and $\PerZero$, we only
have regarded the obtained operads as ns ones. Computer experiments let
us think that the dimensions of the operads $\PRT$, $\FCat{2}$, $\Motz$,
$\DA$ and $\SComp$ seen as symmetric ones are, respectively, Sequences
\OEIS{A052882}, \OEIS{A050351}, \OEIS{A032181}, \OEIS{A101052}, and
\OEIS{A001047} of~\cite{Slo}. Bijections between elements of these
operads and combinatorial objects enumerated by these sequences,
together with presentations by generators and relations in this
symmetric context, would be worth studying.
\smallbreak

Another line of research is the following. It is well-known that the
Koszul dual of the operads $\Dias$ and $\Trias$ are respectively the
dendriform $\Dendr$~\cite{Lod01} and the tridendriform
$\TDendr$~\cite{LR04} operads. Since the operads $\Di$ and $\Tr$,
obtained from the $\T$ construction, are respectively isomorphic to the
operads $\Dias$ and $\Trias$, we can ask if there are generalizations
of $\Di$ and $\Tr$ so that their Koszul duals provide generalizations
of the operads $\Dendr$ and $\TDendr$. It turns out that it is the case
and Chapter~\ref{chap:polydendr} contains our research and results
about this subject.
\medbreak


\chapter{Pluriassociative and polydendriform operads}
\label{chap:polydendr}
The content of this chapter comes from~\cite{Gir16a} and~\cite{Gir16b}.
\medbreak

\section*{Introduction}
Associative algebras play an obvious and primary role in algebraic
combinatorics. Let us cite for instance the algebra of symmetric
functions~\cite{Mac15} involving integer partitions, the algebra of
noncommutative symmetric functions~\cite{GKLLRT94} involving integer
compositions, the Malvenuto-Reutenauer algebra of free quasi-symmetric
functions~\cite{MR95} (see also~\cite{DHT02}) involving permutations,
the Loday-Ronco Hopf algebra of binary trees~\cite{LR98} (see
also~\cite{HNT05}), and the Connes-Kreimer Hopf algebra of forests of
rooted trees~\cite{CK98}. There are several ways to understand and to
gather information about such structures and their associative
operations. A very fruitful strategy consists in splitting an
associative product $\Product$ into two separate operations $\LDendr$
and $\RDendr$ in such a way that $\Product$ turns to be the sum of
$\LDendr$ and $\RDendr$ (see Section~\ref{subsubsec:biproducts} of
Chapter~\ref{chap:algebra} about the sum of operations). One of the most
obvious example occurs by considering the shuffle product on words (see
Section~\ref{subsubsec:associative_algebras} of
Chapter~\ref{chap:algebra}). Indeed, this product can be separated into
two operations according to the origin (first or second operand) of the
last letter of the words appearing in the result~\cite{Ree58}. Other
main examples include the split of the shifted shuffle product of
permutations of the Malvenuto-Reutenauer Hopf algebra and of the product
of binary trees of the Loday-Ronco Hopf algebra~\cite{Foi07}. The
original formalization and the germs of generalization of these notions,
due to Loday~\cite{Lod01}, lead to the introduction of dendriform
algebras. Dendriform algebras are vector spaces endowed with two
operations $\LDendr$ and $\RDendr$ so that $\LDendr + \RDendr$ is
associative and satisfy some few other relations. Since any dendriform
algebra is a quotient of a certain free dendriform algebra, the study of
free dendriform algebras is worth considering. Besides, the
description of free dendriform algebras has a nice combinatorial
interpretation involving binary trees and shuffle of binary trees.
\medbreak

In recent years, several generalizations of dendriform algebras were
introduced and studied. Among these, one can cite tridendriform
algebras~\cite{LR04}, quadri-algebras~\cite{AL04},
ennea-algebras~\cite{Ler04}, $m$-dendriform algebras of
Leroux~\cite{Ler07}, and $m$-dendriform algebras of
Novelli~\cite{Nov14}, all providing new ways to split associative
products into more than two pieces. Besides, free objects in the
corresponding categories of these algebras can be described by
relatively complex combinatorial objects and more or less tricky
operations on these. For instance, free tridendriform algebras involve
Schröder trees, free quadri-algebras involve noncrossing connected
graphs on a circle, and free $m$-dendriform algebras of Leroux and free
$m$-dendriform algebras of Novelli involve planar rooted trees where
internal nodes have a constant number of children.
\medbreak

The first goal of this chapter is to define and justify a new
generalization of dendriform algebras, using the point of view offered
by the theory of operads. Our long term primary objective is to develop
new implements to split associative products in smaller pieces. We use
the approach consisting in considering the diassociative operad
$\Dias$~\cite{Lod01}, the Koszul dual of the dendriform operad $\Dendr$,
rather that focusing on $\Dendr$. Since $\Dias$ admits a description far
simpler than $\Dendr$, starting by constructing a generalization of
$\Dias$ to obtain a generalization of $\Dendr$ by Koszul duality is a
convenient path to explore. To obtain a generalization of the
diassociative operad, we exploit the general functorial construction
$\T$ producing an operad from any monoid (see
Chapter~\ref{chap:monoids}). We show here that this functor $\T$
provides an original construction for the diassociative operad. In the
present chapter, we rely on $\T$ to construct the operads $\Dias_\gamma$,
where $\gamma$ is a nonnegative integer, in such a way that
$\Dias_1 = \Dias$.
\medbreak

The operads $\Dias_\gamma$, called $\gamma$-pluriassociative operads,
are operads defined on the linear span of some words on the alphabet
$\{0\} \sqcup [\gamma]$. By computing the Koszul dual of $\Dias_\gamma$,
we obtain the operads $\Dendr_\gamma$, satisfying $\Dendr_1 = \Dendr$.
The operads $\Dendr_\gamma$ govern the category of the so-called
$\gamma$-polydendriform algebras, that are algebras with
$2\gamma$ operations $\LDendrA_a$, $\RDendrA_a$, $a \in [\gamma]$,
satisfying some relations. Free algebras in these categories involve
binary trees where all edges connecting two internal nodes are labeled
on $[\gamma]$. These algebras are endowed with $2 \gamma$ products
described by induction and can be seen as kinds of shuffle of trees,
generalizing the shuffle of trees introduced by Loday~\cite{Lod01}
intervening in the construction of free dendriform algebras. Moreover,
the introduction of $\gamma$-polydendriform algebras offers to split an
associative product $\Product$ by
\begin{equation}
    \Product = \LDendrA_1 + \RDendrA_1 + \dots + \LDendrA_\gamma
    + \RDendrA_\gamma,
\end{equation}
with, among others, the stiffening conditions that all partial sums
\begin{equation}
    \LDendrA_1 + \RDendrA_1 + \dots + \LDendrA_a + \RDendrA_a
\end{equation}
are associative for all $a \in [\gamma]$. Besides, this work naturally
leads to the consideration and the definition of numerous operads.
Table~\ref{tab:gamma_operads} summarizes some information about these.
\begin{table}
    \centering
    \begin{tabular}{c|c|c|c}
        Operad & Bases & Dimensions & Symmetric \\ \hline \hline
        $\Dias_\gamma$ & Some words on $\{0\} \sqcup [\gamma]$
            & $n \gamma^{n - 1}$ & No \\ \hline
        $\Dendr_\gamma$ & $\gamma$-edge valued binary trees
            & $\gamma^{n - 1}\, \Catalan(n)$ & No \\ \hline
        $\As_\gamma$ & $\gamma$-corollas & $\gamma$ & No \\ \hline
        $\DAs_\gamma$ & $\gamma$-alternating Schröder trees
            & $\sum\limits_{k = 0}^{n - 2} \gamma^{k + 1}
            (\gamma - 1)^{n - k - 2}\, \Narayana(n, k)$ & No \\ \hline
        $\Dup_\gamma$ & $\gamma$-edge valued binary trees
            & $\gamma^{n - 1}\, \Catalan(n)$ & No \\ \hline
        $\Trias_\gamma$ & Some words on $\{0\} \sqcup [\gamma]$
            & $(\gamma + 1)^n - \gamma^n$ & No \\ \hline
        $\TDendr_\gamma$ & $\gamma$-edge valued Schröder trees
            & $\sum\limits_{k = 0}^{n - 1} (\gamma + 1)^k
            \gamma^{n - k - 1}\,\Narayana(n + 1, k)$ & No \\ \hline
        $\Com_\gamma$ & --- & --- & Yes \\ \hline
        $\Zin_\gamma$ & --- & --- & Yes
    \end{tabular}
    \bigbreak

    \caption[Data about operads related to the ns operads $\Dias$ and
    $\Dendr$.]
    {The main operads constructed in this chapter. All these depend on
    a nonnegative integer parameter $\gamma$. The shown dimensions are
    the ones of their homogeneous components of arities~$n \geq 2$.}
    \label{tab:gamma_operads}
\end{table}
\medbreak

This work is organized as follows. Section~\ref{sec:dias_gamma} is
devoted to the introduction and the complete study of the operads
$\Dias_\gamma$, and in Section~\ref{sec:dias_gamma_algebras}, algebras
over $\Dias_\gamma$ are studied. In Section~\ref{sec:pluritriass} we
presents an analogous generalization $\Trias_\gamma$ of the
triassociative operad~\cite{LR04}. The study of $\Dendr_\gamma$ is
performed in Section~\ref{sec:dendr_gamma}, where we provide several
presentations of this operad and a construction of free
$\gamma$-polydendriform algebras. Section~\ref{sec:as_gamma} extends a
part of the operadic butterfly~\cite{Lod01,Lod06}. This extension
contains the operads $\Dias_\gamma$, $\Dendr_\gamma$, and two
generalizations $\As_\gamma$ and $\DAs_\gamma$ of the associative
operad, Koszul duals one of the other. Finally, in
Section~\ref{sec:further_generalizations_gamma}, we sustain our previous
ideas to propose still new generalizations of some more operads like the
operad $\Dup_\gamma$ generalizing the duplicial operad~\cite{Lod08}
and the operad $\TDendr_\gamma$ generalizing the tridendriform operad.
We also then define the operads $\Com_\gamma$, $\Lie_\gamma$,
$\Zin_\gamma$, and $\Leib_\gamma$, that are respective generalizations
of the commutative operad, the Lie operad, the Zinbiel
operad~\cite{Lod95} and the Leibniz operad~\cite{Lod93}.
\medbreak

\subsubsection*{Note}
This chapter deals mostly with ns operads. For this reason, ``operad''
means ``ns operad''.
\medbreak

\section{Pluriassociative operads} \label{sec:dias_gamma}
We define here one of the main object of this chapter: a generalization
on a nonnegative integer parameter $\gamma$ of the diassociative operad
(see~\cite{Lod01} or Section~\ref{subsubsec:dendriform_operad} of
Chapter~\ref{chap:algebra}). We provide a complete study of this new
operad.
\medbreak

\subsection{Construction and first properties}
Our generalization of the diassociative operad passes through the
functor $\T$ (see Section~\ref{subsec:operad_of_monoid} of
Chapter~\ref{chap:monoids}). We begin here by describing a basis and by
establishing the Hilbert series of our generalization.
\medbreak

\subsubsection{Construction} \label{subsubsec:construction_dias_gamma}
For any integer $\gamma \geq 0$, let $\Nmax_\gamma$ be the monoid
$\{0\} \sqcup [\gamma]$ with the binary operation $\max$ as product,
denoted by $\Max$. We define the \Def{$\gamma$-pluriassociative operad}
$\Dias_\gamma$ as the suboperad of $\T \Nmax_\gamma$ generated by
\begin{equation} \label{equ:generators_dias_gamma}
    \GenDias := \{0a, a0 : a \in [\gamma]\}.
\end{equation}
By definition, $\Dias_\gamma$ is the linear span of all the words that
can be obtained by partial compositions of words of $\GenDias$. We have,
for instance,
\begin{subequations}
\begin{equation}
    \Dias_2(1) = \K \Angle{\{0\}},
\end{equation}
\begin{equation}
    \Dias_2(2)
    = \K \Angle{\{01, 02, 10, 20\}},
\end{equation}
\begin{equation}
    \Dias_2(3)
    = \K \Angle{\{011, 012, 021, 022, 101, 102,
        201, 202, 110, 120, 210, 220\}},
\end{equation}
\end{subequations}
and, as examples of partial compositions in $\Dias_3$,
\begin{subequations}
\begin{equation}
    \textcolor{Col1}{211} {\bf 2} \textcolor{Col1}{01}
        \circ_4 \textcolor{Col4}{31103}
    = \textcolor{Col1}{211} \textcolor{Col4}{32223}
    \textcolor{Col1}{01},
\end{equation}
\begin{equation}
    \textcolor{Col1}{11} {\bf 1} \textcolor{Col1}{101}
        \circ_3 \textcolor{Col4}{20}
    = \textcolor{Col1}{11} \textcolor{Col4}{21}
    \textcolor{Col1}{101},
\end{equation}
\begin{equation}
    \textcolor{Col1}{1} {\bf 0} \textcolor{Col1}{13}
        \circ_2 \textcolor{Col4}{210}
    = \textcolor{Col1}{1} \textcolor{Col4}{210}
    \textcolor{Col1}{13}.
\end{equation}
\end{subequations}
\medbreak

\subsubsection{First properties}
In the first place, observe that $\Dias_1$ is the operad $\Di$ defined
in Chapter~\ref{chap:monoids}. For this reason, $\Dias_1$ is the
diassociative operad $\Dias$. Moreover, observe that $\Dias_0$ is the
trivial operad and that $\Dias_\gamma$ is a suboperad of
$\Dias_{\gamma + 1}$. Then, for all integers $\gamma \geq 0$, the
operads $\Dias_\gamma$ are generalizations of the diassociative operad.
Besides, it follows immediately from the definition of $\Dias_\gamma$ as
a suboperad of $\T \Nmax_\gamma$ that its fundamental basis is a
set-operad basis. Indeed, any partial composition of two basis elements
of $\Dias_\gamma$ gives rises to exactly one basis element.
\medbreak

\subsubsection{Elements and dimensions}
\begin{Proposition} \label{prop:basis_dias_gamma}
    For any integer $\gamma \geq 0$, the fundamental basis of
    $\Dias_\gamma$ is the set of all the words on the alphabet
    $\{0\} \sqcup [\gamma]$ containing exactly one occurrence of~$0$.
\end{Proposition}
\medbreak

We deduce from Proposition~\ref{prop:basis_dias_gamma} that the
Hilbert series of $\Dias_\gamma$ satisfies
\begin{equation} \label{equ:serie_hilbert_dias_gamma}
    \Hca_{\Dias_\gamma}(t) = \frac{t}{(1 - \gamma t)^2}
\end{equation}
and that for all $n \geq 1$, $\dim \Dias_\gamma(n) = n \gamma^{n - 1}$.
For instance, the first dimensions of $\Dias_1$, $\Dias_2$, $\Dias_3$,
and $\Dias_4$ are respectively
\begin{subequations}
\begin{equation}
    1, 2, 3, 4, 5, 6, 7, 8, 9, 10, 11,
\end{equation}
\begin{equation}
    1, 4, 12, 32, 80, 192, 448, 1024, 2304, 5120, 11264,
\end{equation}
\begin{equation}
    1, 6, 27, 108, 405, 1458, 5103, 17496, 59049, 196830, 649539,
\end{equation}
\begin{equation}
    1, 8, 48, 256, 1280, 6144, 28672, 131072, 589824, 2621440, 11534336.
\end{equation}
\end{subequations}
The second one is Sequence~\OEIS{A001787}, the third one is
Sequence~\OEIS{A027471}, and the last one is Sequence~\OEIS{A002697}
of~\cite{Slo}.
\medbreak

\subsection{Additional properties}\label{subsec:presentation_dias_gamma}
We exhibit here, among others, two presentations of $\Dias_\gamma$ and
establish the  fact that it is a Koszul operad.
\medbreak

\subsubsection{Presentation by generators and relations}
\label{subsubsec:presentation_dias_gamma}
For any $a \in [\gamma]$, let us denote by $\LDias_a$ (resp. $\RDias_a$)
the generator $0a$ (resp. $a0$) of $\Dias_\gamma$.
\medbreak

\begin{Theorem} \label{thm:presentation_dias_gamma}
    For any integer $\gamma \geq 0$, the operad $\Dias_\gamma$ admits
    the presentation $\left(\GenDias, \RelDias\right)$ where $\RelDias$
    is the space induced by the equivalence relation $\Equiv_\gamma$
    satisfying
    \begin{subequations}
    \begin{equation} \label{equ:relation_dias_gamma_1}
        \Corolla{\LDias_a} \circ_1 \Corolla{\RDias_{a'}}
         \Equiv_\gamma
        \Corolla{\RDias_{a'}} \circ_2 \Corolla{\LDias_a},
        \qquad a, a' \in [\gamma],
    \end{equation}
    \begin{equation} \label{equ:relation_dias_gamma_2}
        \Corolla{\LDias_a} \circ_1 \Corolla{\LDias_b}
         \Equiv_\gamma
        \Corolla{\LDias_a} \circ_2 \Corolla{\RDias_b},
        \qquad a < b \in [\gamma],
    \end{equation}
    \begin{equation} \label{equ:relation_dias_gamma_3}
        \Corolla{\RDias_a} \circ_1 \Corolla{\LDias_b}
         \Equiv_\gamma
        \Corolla{\RDias_a} \circ_2 \Corolla{\RDias_b},
        \qquad a < b \in [\gamma],
    \end{equation}
    \begin{equation} \label{equ:relation_dias_gamma_4}
        \Corolla{\LDias_b} \circ_1 \Corolla{\LDias_a}
         \Equiv_\gamma
        \Corolla{\LDias_a} \circ_2 \Corolla{\LDias_b},
        \qquad a < b \in [\gamma],
    \end{equation}
    \begin{equation} \label{equ:relation_dias_gamma_5}
        \Corolla{\RDias_a} \circ_1 \Corolla{\RDias_b}
         \Equiv_\gamma
        \Corolla{\RDias_b} \circ_2 \Corolla{\RDias_a},
        \qquad a < b \in [\gamma],
    \end{equation}
    \begin{equation} \label{equ:relation_dias_gamma_6}
        \Corolla{\LDias_d} \circ_1 \Corolla{\LDias_d}
         \Equiv_\gamma
        \Corolla{\LDias_d} \circ_2 \Corolla{\LDias_c}
         \Equiv_\gamma
        \Corolla{\LDias_d} \circ_2 \Corolla{\RDias_c},
        \qquad c \leq d \in [\gamma],
    \end{equation}
    \begin{equation} \label{equ:relation_dias_gamma_7}
        \Corolla{\RDias_d} \circ_1 \Corolla{\LDias_c}
         \Equiv_\gamma
        \Corolla{\RDias_d} \circ_1 \Corolla{\RDias_c}
         \Equiv_\gamma
        \Corolla{\RDias_d} \circ_2 \Corolla{\RDias_d},
        \qquad c \leq d \in [\gamma].
    \end{equation}
    \end{subequations}
\end{Theorem}
\medbreak

Our proof of Theorem~\ref{thm:presentation_dias_gamma} does not follow
the usual technique consisting in providing a convergent orientation
$\Rew_\gamma$ of $\Equiv_\gamma$ and proving that its closure
$\RewTrees_\gamma$ admits as many normal forms of arity $n$ as basis
words of $\Dias_\gamma(n)$ (as for instance in
Chapter~\ref{chap:monoids}). Instead, we consider the evaluation
morphism
\begin{equation}
    \Eval :
    \FreeOperad\left(\GenDias\right) \to \Dias_\gamma,
\end{equation}
and show that its kernel is generated by $\RelDias$. This strategy uses
the noteworthy fact that the image of a $\GenDias$-syntax tree $\Tfr$
can be computed as follows. We say that an integer
$a \in \{0\} \sqcup [\gamma]$ is \Def{eligible} for a leaf $x$ of $\Tfr$
if $a = 0$ or there is an ancestor $y$ of $x$ labeled by $\LDias_a$
(resp. $\RDias_a$) and $x$ is in the right (resp. left) subtree of $y$.
The \Def{image} of $x$ is its greatest eligible integer. Now,
$\Eval(\Tfr)$ is the word obtained by considering, from left to right,
the images of the leaves of $\Tfr$ (see
Figure~\ref{fig:example_eval_dias_gamma}).
\begin{figure}[ht]
    \centering

    \caption[A syntax tree on $\GenDias$ and its evaluation.]
    {A $\GenDias$-syntax tree $\Tfr$ where the images of its leaves
    are shown. This tree satisfies
    $\Eval(\Tfr) = \textcolor{Col1}{340122332242}$.}
    \label{fig:example_eval_dias_gamma}
\end{figure}
\medbreak

The space of relations $\RelDias$ of $\Dias_\gamma$ exhibited by
Theorem~\ref{thm:presentation_dias_gamma} can be rephrased in a more
compact way as the space generated by
\begin{subequations}
\begin{equation} \label{equ:relation_dias_gamma_1_concise}
    \Corolla{\LDias_a} \circ_1 \Corolla{\RDias_{a'}}
    -
    \Corolla{\RDias_{a'}} \circ_2 \Corolla{\LDias_a},
    \qquad a, a' \in [\gamma],
\end{equation}
\begin{equation} \label{equ:relation_dias_gamma_2_concise}
    \Corolla{\LDias_a} \circ_1 \Corolla{\LDias_{a \Max a'}}
    -
    \Corolla{\LDias_a} \circ_2 \Corolla{\RDias_{a'}},
    \qquad a, a' \in [\gamma],
\end{equation}
\begin{equation} \label{equ:relation_dias_gamma_3_concise}
    \Corolla{\RDias_a} \circ_1 \Corolla{\LDias_{a'}}
    -
    \Corolla{\RDias_a} \circ_2 \Corolla{\RDias_{a \Max a'}},
    \qquad a, a' \in [\gamma],
\end{equation}
\begin{equation} \label{equ:relation_dias_gamma_4_concise}
    \Corolla{\LDias_{a \Max a'}} \circ_1 \Corolla{\LDias_a}
    -
    \Corolla{\LDias_a} \circ_2 \Corolla{\LDias_{a'}},
    \qquad a, a' \in [\gamma],
\end{equation}
\begin{equation} \label{equ:relation_dias_gamma_5_concise}
    \Corolla{\RDias_a} \circ_1 \Corolla{\RDias_{a'}}
    -
    \Corolla{\RDias_{a \Max a'}} \circ_2 \Corolla{\RDias_a},
    \qquad a, a' \in [\gamma].
\end{equation}
\end{subequations}
\medbreak

\subsubsection{Koszulity}
\begin{Theorem} \label{thm:koszulity_dias_gamma}
    For any integer $\gamma \geq 0$, $\Dias_\gamma$ is a Koszul operad
    and the set of the $\GenDias$-syntax trees avoiding the trees
    \begin{subequations}
    \begin{equation}
        \Corolla{\RDias_{a'}} \circ_2 \Corolla{\LDias_a},
        \qquad a, a' \in [\gamma],
    \end{equation}
    \begin{equation}
        \Corolla{\LDias_a} \circ_2 \Corolla{\RDias_{a'}},
        \qquad a, a' \in [\gamma],
    \end{equation}
    \begin{equation}
        \Corolla{\RDias_a} \circ_1 \Corolla{\LDias_{a'}},
        \qquad a, a' \in [\gamma],
    \end{equation}
    \begin{equation}
        \Corolla{\LDias_a} \circ_2 \Corolla{\LDias_{a'}},
        \qquad a, a' \in [\gamma],
    \end{equation}
    \begin{equation}
        \Corolla{\RDias_a} \circ_1 \Corolla{\RDias_{a'}},
        \qquad a, a' \in [\gamma],
    \end{equation}
    \end{subequations}
    is a Poincaré-Birkhoff-Witt basis of~$\Dias_\gamma$.
\end{Theorem}
\medbreak

\subsubsection{Miscellaneous properties}
We list some secondary properties of $\Dias_\gamma$. The definitions of
these properties can be found in Section~\ref{subsec:operads} of
Chapter~\ref{chap:algebra}.
\medbreak

\begin{Proposition} \label{prop:symmetries_dias_gamma}
    For any integer $\gamma \geq 0$, the group of symmetries of
    $\Dias_\gamma$ contains the linear map sending any word of
    $\Dias_\gamma$ to its mirror image.
\end{Proposition}
\medbreak

\begin{Proposition} \label{prop:dias_gamma_basic}
    For any integer $\gamma \geq 0$, the fundamental basis of
    $\Dias_\gamma$ is a basic set-operad basis.
\end{Proposition}
\medbreak

\begin{Proposition} \label{prop:dias_gamma_rooted}
    For any integer $\gamma \geq 0$, $\Dias_\gamma$ is a nontrivially
    rooted operad for the root map sending any word of $\Dias_\gamma$ to
    the position of its $0$.
\end{Proposition}
\medbreak

\subsubsection{Alternative basis} \label{subsubsec:basis_K_dias_gamma}
Let $\OrdDias_\gamma$ be the order relation on the underlying set of
$\Dias_\gamma(n)$, $n \geq 1$, where for all words $x$ and $y$ of
$\Dias_\gamma$ of a same arity $n$, we have
\begin{equation}
    x \OrdDias_\gamma y
    \qquad \mbox{ if } x_i \leq y_i \mbox{ for all } i \in [n].
\end{equation}
This order relation allows to define
for all words $x$ of $\Dias_\gamma$ the elements
\begin{equation} \label{equ:basis_K_to_fundamental_dias_gamma}
    \BasisK^{(\gamma)}_x :=
    \sum_{x \OrdDias_\gamma x'} \mu_\gamma(x, x') \, x',
\end{equation}
where $\mu_\gamma$ is the Möbius function of the poset defined by
$\OrdDias_\gamma$. For instance,
\begin{subequations}
\begin{equation}
    \BasisK^{(2)}_{102} = 102 - 202,
\end{equation}
\begin{equation}
    \BasisK^{(3)}_{102} = \BasisK^{(4)}_{102} = 102 - 103 - 202 + 203,
\end{equation}
\begin{equation}
    \BasisK^{(3)}_{23102}
    = 23102 - 23103 - 23202 + 23203 - 33102 + 33103 + 33202 - 33203.
\end{equation}
\end{subequations}
\medbreak

Since, by Möbius inversion (see
Proposition~\ref{prop:partial_order_bases} of
Chapter~\ref{chap:algebra}), for any word $x$ of
$\Dias_\gamma$ one has
\begin{equation} \label{equ:fundamental_to_basis_K_dias_gamma}
    x = \sum_{x \OrdDias_\gamma x'} \BasisK^{(\gamma)}_{x'},
\end{equation}
the family of all $\BasisK^{(\gamma)}_x$, where the $x$ are words of
$\Dias_\gamma$, forms by triangularity a basis of $\Dias_\gamma$, called
the \Def{$\BasisK$-basis}.
\medbreak

Recall that $\Hamming(u, v)$ denotes the Hamming distance between the
words $u$ and $v$ of the same length. For any word $x$ of $\Dias_\gamma$
of length $n$, we denote by $\Incr_\gamma(x)$ the set of all words
obtained by incrementing by $1$ some letters of $x$ smaller than
$\gamma$ and greater than $0$. Proposition~\ref{prop:basis_dias_gamma}
ensures that $\Incr_\gamma(x)$ is a set of words of~$\Dias_\gamma$.
\medbreak

\begin{Proposition}
    \label{prop:basis_K_to_fundamental_direct_dias_gamma}
    For any integer $\gamma \geq 0$ and any word $x$ of $\Dias_\gamma$,
    \begin{equation}
    \label{equ:basis_K_to_fundamental_direct_dias_gamma}
        \BasisK^{(\gamma)}_x =
        \sum_{x' \in \Incr_\gamma(x)} (-1)^{\Hamming(x, x')} \, x'.
    \end{equation}
\end{Proposition}
\medbreak

To compute a direct expression for the partial composition of
$\Dias_\gamma$ over the $\BasisK$-basis, we have to introduce two
notations. If $x$ is a word of $\Dias_\gamma$ of length nons maller than
$2$, we denote by $\min(x)$ the smallest letter of $x$ among its letters
different from $0$. Proposition~\ref{prop:basis_dias_gamma} ensures that
$\min(x)$ is well-defined. Moreover, for all words $x$ and $y$ of
$\Dias_\gamma$, a position $i$ such that $x_i \ne 0$, and
$a \in [\gamma]$, we denote by $x \circ_{a, i} y$ the word $x \circ_i y$
in which the $0$ coming from $y$ is replaced by $a$ instead of~$x_i$.
\medbreak

\begin{Theorem} \label{thm:composition_basis_K_dias_gamma}
    For any integer $\gamma \geq 0$, the partial composition of
    $\Dias_\gamma$ over the $\BasisK$-basis satisfies, for all words $x$
    and $y$ of $\Dias_\gamma$ of arities non smaller than $2$,
    \begin{equation}
        \BasisK^{(\gamma)}_x \circ_i \BasisK^{(\gamma)}_y =
        \begin{cases}
            \BasisK^{(\gamma)}_{x \circ_i y}
                & \mbox{if } \min(y) > x_i, \\
            \sum_{a \in [x_i, \gamma]}
                \BasisK^{(\gamma)}_{x \circ_{a, i} y}
                & \mbox{if } \min(y) = x_i, \\
            0 & \mbox{otherwise (} \min(y) < x_i \mbox{)}.
        \end{cases}
    \end{equation}
\end{Theorem}
\medbreak

We have for instance
\begin{subequations}
\begin{equation}
    \BasisK^{(5)}_{20413} \circ_1 \BasisK^{(5)}_{304}
    = \BasisK^{(5)}_{3240413},
\end{equation}
\begin{equation}
    \BasisK^{(5)}_{20413} \circ_2 \BasisK^{(5)}_{304}
    = \BasisK^{(5)}_{2304413},
\end{equation}
\begin{equation}
    \BasisK^{(5)}_{20413} \circ_3 \BasisK^{(5)}_{304} = 0,
\end{equation}
\begin{equation}
    \BasisK^{(5)}_{20413} \circ_4 \BasisK^{(5)}_{304}
    = \BasisK^{(5)}_{2043143},
\end{equation}
\begin{equation}
    \BasisK^{(5)}_{20413} \circ_5 \BasisK^{(5)}_{304}
    = \BasisK^{(5)}_{2041334} + \BasisK^{(5)}_{2041344}
        + \BasisK^{(5)}_{2041354}.
\end{equation}
\end{subequations}
\medbreak

Theorem~\ref{thm:composition_basis_K_dias_gamma} implies in particular
that the structure coefficients of the partial composition of
$\Dias_\gamma$ over the $\BasisK$-basis are $0$ or $1$. It is possible
to define another bases of $\Dias_\gamma$ by reversing
in~\eqref{equ:basis_K_to_fundamental_dias_gamma} the relation
$\OrdDias_\gamma$ and by suppressing or keeping the Möbius function
$\mu_\gamma$. This gives obviously rise to three other bases. It is
worth to note that, as small computations reveal, over all these
additional bases, the structure coefficients of the partial
composition of $\Dias_\gamma$ can be negative or different from $1$.
This observation makes the $\BasisK$-basis even more particular and
interesting. It has some other properties, as next section will show.
\medbreak

\subsubsection{Alternative presentation}
\label{subsubsec:presentation_dias_gamma_alternative}
The $\BasisK$-basis introduced in the previous section leads to state a
new presentation for $\Dias_\gamma$ in the following way. For any
$a \in [\gamma]$, let us denote by $\LDiasA_a$ (resp. $\RDiasA_a$) the
element $\BasisK_{0a}$ (resp. $\BasisK_{a0}$) of $\Dias_\gamma$. Then,
for all $a \in [\gamma]$ we have
\begin{subequations}
\begin{equation}
    \LDias_a = \sum_{a \leq b \in [\gamma]} \LDiasA_b
\end{equation}
and
\begin{equation}
    \RDias_a = \sum_{a \leq b \in [\gamma]} \RDiasA_b,
\end{equation}
\end{subequations}
and by triangularity, the family
\begin{equation}
    \GenDias' := \left\{\LDiasA_a, \RDiasA_a : a \in [\gamma]\right\}
\end{equation}
is a generating set of~$\Dias_\gamma$.
\medbreak

\begin{Proposition} \label{prop:presentation_dias_gamma_alternative}
    For any integer $\gamma \geq 0$, the operad $\Dias_\gamma$ admits
    the presentation
    $\left(\GenDias', \RelationSpace'_{\Dias_\gamma}\right)$
    where $\RelationSpace'_{\Dias_\gamma}$ is the space generated by
    \begin{subequations}
    \begin{equation} \label{equ:relation_dias_gamma_1_alternative}
        \Corolla{\LDiasA_a} \circ_1 \Corolla{\RDiasA_{a'}} -
        \Corolla{\RDiasA_{a'}} \circ_2 \Corolla{\LDiasA_a},
        \qquad a, a' \in [\gamma],
    \end{equation}
    \begin{equation} \label{equ:relation_dias_gamma_2_alternative}
        \Corolla{\RDiasA_b} \circ_1 \Corolla{\RDiasA_a},
        \qquad a < b \in [\gamma],
    \end{equation}
    \begin{equation} \label{equ:relation_dias_gamma_3_alternative}
        \Corolla{\LDiasA_b} \circ_2 \Corolla{\LDiasA_a},
        \qquad a < b \in [\gamma],
    \end{equation}
    \begin{equation} \label{equ:relation_dias_gamma_4_alternative}
        \Corolla{\RDiasA_b} \circ_1 \Corolla{\LDiasA_a},
        \qquad a < b \in [\gamma],
    \end{equation}
    \begin{equation} \label{equ:relation_dias_gamma_5_alternative}
        \Corolla{\LDiasA_b} \circ_2 \Corolla{\RDiasA_a},
        \qquad a < b \in [\gamma],
    \end{equation}
    \begin{equation} \label{equ:relation_dias_gamma_6_alternative}
        \Corolla{\RDiasA_a} \circ_1 \Corolla{\RDiasA_b}
        -
        \Corolla{\RDiasA_b} \circ_2 \Corolla{\RDiasA_a},
        \qquad a < b \in [\gamma],
    \end{equation}
    \begin{equation} \label{equ:relation_dias_gamma_7_alternative}
        \Corolla{\LDiasA_b} \circ_1 \Corolla{\LDiasA_a}
        -
        \Corolla{\LDiasA_a} \circ_2 \Corolla{\LDiasA_b},
        \qquad a < b \in [\gamma],
    \end{equation}
    \begin{equation} \label{equ:relation_dias_gamma_8_alternative}
        \Corolla{\RDiasA_a} \circ_1 \Corolla{\LDiasA_b}
        -
        \Corolla{\RDiasA_a} \circ_2 \Corolla{\RDiasA_b},
        \qquad a < b \in [\gamma],
    \end{equation}
    \begin{equation} \label{equ:relation_dias_gamma_9_alternative}
        \Corolla{\LDiasA_a} \circ_1 \Corolla{\LDiasA_b}
        -
        \Corolla{\LDiasA_a} \circ_2 \Corolla{\RDiasA_b},
        \qquad a < b \in [\gamma],
    \end{equation}
    \begin{equation} \label{equ:relation_dias_gamma_10_alternative}
        \Corolla{\RDiasA_a} \circ_1 \Corolla{\RDiasA_a}
        -
        \left(\sum_{a \leq b \in [\gamma]}
            \Corolla{\RDiasA_a} \circ_2 \Corolla{\RDiasA_b}
        \right),
        \qquad a \in [\gamma],
    \end{equation}
    \begin{equation} \label{equ:relation_dias_gamma_11_alternative}
        \left(\sum_{a \leq b \in [\gamma]}
            \Corolla{\LDiasA_a} \circ_1 \Corolla{\LDiasA_b}
        \right)
        -
        \Corolla{\LDiasA_a} \circ_2 \Corolla{\LDiasA_a},
        \qquad a \in [\gamma],
    \end{equation}
    \begin{equation} \label{equ:relation_dias_gamma_12_alternative}
        \Corolla{\RDiasA_a} \circ_1 \Corolla{\LDiasA_a}
        -
        \left(\sum_{a \leq b \in [\gamma]}
            \Corolla{\RDiasA_b} \circ_2 \Corolla{\RDiasA_a}
        \right),
        \qquad a \in [\gamma],
    \end{equation}
    \begin{equation} \label{equ:relation_dias_gamma_13_alternative}
        \left(\sum_{a \leq b \in [\gamma]}
            \Corolla{\LDiasA_b} \circ_1 \Corolla{\LDiasA_a}
        \right)
        -
        \Corolla{\LDiasA_a} \circ_2 \Corolla{\RDiasA_a},
        \qquad a \in [\gamma].
    \end{equation}
    \end{subequations}
\end{Proposition}
\medbreak

Despite the apparent complexity of the presentation of $\Dias_\gamma$
exhibited by Proposition~\ref{prop:presentation_dias_gamma_alternative},
as we will see in Section~\ref{sec:dendr_gamma}, the Koszul dual of
$\Dias_\gamma$ computed from this presentation has a very simple and
manageable expression.
\medbreak

\section{Pluriassociative algebras} \label{sec:dias_gamma_algebras}
We now focus on algebras over $\gamma$-pluriassociative operads.
For this purpose, we construct free $\Dias_\gamma$-algebras over one
generator, and define and study two notions of units for
$\Dias_\gamma$-algebras. We end this section by introducing a convenient
way to define $\Dias_\gamma$-algebras and give several examples of such
algebras.
\medbreak

\subsection{Category of pluriassociative algebras and free objects}
Let us study the category of $\Dias_\gamma$-algebras and the units for
algebras in this category.
\medbreak

\subsubsection{Pluriassociative algebras}
We call \Def{$\gamma$-pluriassociative algebra} any
$\Dias_\gamma$-algebra. From the presentation of $\Dias_\gamma$ provided
by Theorem~\ref{thm:presentation_dias_gamma}, any
$\gamma$-pluriassociative algebra is a vector space endowed with linear
operations $\LDias_a, \RDias_a$, $a \in [\gamma]$, satisfying
the relations encoded by~\eqref{equ:relation_dias_gamma_1_concise}---%
\eqref{equ:relation_dias_gamma_5_concise}.
\medbreak

\subsubsection{General definitions}
Let $\Pca$ be a $\gamma$-pluriassociative algebra. We say that $\Pca$
is \Def{commutative} if for all $x, y \in \Pca$ and $a \in [\gamma]$,
$x \LDias_a y = y \RDias_a x$. Besides, $\Pca$ is \Def{pure} for all
$a, a' \in [\gamma]$, $a \ne a'$ implies $\LDias_a \ne \LDias_{a'}$ and
$\RDias_a \ne \RDias_{a'}$.
\medbreak

Given a subset $C$ of $[\gamma]$, one can keep on the vector space
$\Pca$ only the operations $\LDias_a$ and $\RDias_a$ such that
$a \in C$. By renumbering the indices of these operations from $1$ to
$\# C$ by respecting their former relative numbering, we obtain a
$\# C$-pluriassociative algebra. We call it the
\Def{$\# C$-pluriassociative subalgebra induced by} $C$ of~$\Pca$.
\medbreak

\subsubsection{Free pluriassociative algebras}
\label{subsubsec:free_dias_gamma_algebras}
Recall that $\FreeAlg_{\Dias_\gamma}$ denotes the free
$\Dias_\gamma$-algebra over one generator. By definition,
$\FreeAlg_{\Dias_\gamma}$ is the linear span of the set of the words on
$\{0\} \sqcup [\gamma]$ with exactly one occurrence of~$0$. Let us endow
this space with the linear operations
\begin{equation}
    \LDias_a, \RDias_a :
    \FreeAlg_{\Dias_\gamma} \otimes \FreeAlg_{\Dias_\gamma}
    \to \FreeAlg_{\Dias_\gamma},
    \qquad
    a \in [\gamma],
\end{equation}
satisfying, for any such words $u$ and $v$,
\begin{subequations}
\begin{equation}
    u \LDias_a v := u \; \Augm_a(v)
\end{equation}
and
\begin{equation}
    u \RDias_a v := \Augm_a(u) \; v,
\end{equation}
\end{subequations}
where $\Augm_a(u)$ (resp. $\Augm_a(v)$) is the word obtained by replacing
in $u$ (resp. $v$) any occurrence of a letter smaller than $a$ by $a$.
\medbreak

\begin{Proposition} \label{prop:free_dias_gamma_algebras}
    For any integer $\gamma \geq 0$, the vector space
    $\FreeAlg_{\Dias_\gamma}$ of all nonempty words on
    $\{0\} \sqcup [\gamma]$ containing exactly one occurrence of $0$
    endowed with the operations $\LDias_a$, $\RDias_a$,
    $a \in [\gamma]$, is the free $\gamma$-pluriassociative algebra
    over one generator.
\end{Proposition}
\medbreak

One has for instance in $\FreeAlg_{\Dias_4}$,
\begin{subequations}
\begin{equation}
    \textcolor{Col1}{101241} \LDias_2 \textcolor{Col4}{203} =
    \textcolor{Col1}{101241}\textcolor{Col4}{2{\bf 2}3},
\end{equation}
\begin{equation}
    \textcolor{Col1}{101241} \RDias_3 \textcolor{Col4}{203} =
    \textcolor{Col1}{{\bf 3333}4{\bf 3}}\textcolor{Col4}{203}.
\end{equation}
\end{subequations}
\medbreak

\subsection{Bar and wire-units}
Loday has defined in~\cite{Lod01} some notions of units in diassociative
algebras. We generalize here these definitions to the context of
$\gamma$-pluriassociative algebras.
\medbreak

\subsubsection{Bar-units}
Let $\Pca$ be a $\gamma$-pluriassociative algebra and $a \in [\gamma]$.
We say that an element $e$ of $\Pca$ is an \Def{$a$-bar-unit}, or
simply a \Def{bar-unit} when taking into account the value of $a$
is not necessary, of $\Pca$ if for all $x \in \Pca$,
\begin{equation}
    x \LDias_a e = x = e \RDias_a x.
\end{equation}
As we shall see below, a $\gamma$-pluriassociative algebra can have,
for a given $a \in [\gamma]$, several $a$-bar-units. The \Def{$a$-halo}
of $\Pca$, denoted by $\Halo_a(\Pca)$, is the set of the $a$-bar-units
of~$\Pca$.
\medbreak

\subsubsection{Wire-units}
Let $\Pca$ be a $\gamma$-pluriassociative algebra and $a \in [\gamma]$.
We say that an element $e$ of $\Pca$ is an \Def{$a$-wire-unit}, or
simply a \Def{wire-unit} when taking into account the value of $a$
is not necessary, of $\Pca$ if for all $x \in \Pca$,
\begin{equation}
    e \LDias_a x = x = x \RDias_a e.
\end{equation}
As the following proposition shows, the presence of a wire-unit in
$\Pca$ has some implications.
\medbreak

\begin{Proposition} \label{prop:wire_unit_dias_gamma}
    Let $\gamma \geq 0$ be an integer and $\Pca$ be a
    $\gamma$-pluriassociative algebra admitting a $b$-wire-unit $e$
    for a $b \in [\gamma]$. Then
    \begin{enumerate}[label={\it (\roman*)}]
        \item for all $a \in [b]$, the operations $\LDias_a$,
        $\LDias_b$, $\RDias_a$, and $\RDias_b$ of $\Pca$ are equal;
        \item $e$ is also an $a$-wire-unit for all $a \in [b]$;
        \item $e$ is the only wire-unit of $\Pca$;
        \item if $e'$ is an $a$-bar unit for a $a \in [b]$,
        then~$e' = e$.
    \end{enumerate}
\end{Proposition}
\medbreak

Relying on Proposition~\ref{prop:wire_unit_dias_gamma}, we define the
\Def{height} of a $\gamma$-pluriassociative algebra $\Pca$ as $0$ if
$\Pca$ has no wire-unit, otherwise as the greatest integer
$h \in [\gamma]$ such that the unique wire-unit $e$ of $\Pca$ is a
$h$-wire-unit. Observe that any pure $\gamma$-pluriassociative algebra
has height $0$ or~$1$.
\medbreak

\subsection{Construction of pluriassociative algebras}
We now present a general way to construct $\gamma$-pluriassociative
algebras. Our construction is a natural generalization of some
constructions introduced by Loday~\cite{Lod01} in the context of
diassociative algebras. In this section, we introduce new algebraic
structures, the so-called $\gamma$-multiprojection algebras, which are
inputs of our construction.
\medbreak

\subsubsection{Multiassociative algebras}
\label{subsubsec:multiassociative_algebras}
For any integer $\gamma \geq 0$, a
\Def{$\gamma$-multiassociative algebra} is a vector space $\Mca$ endowed
with linear operations
\begin{equation}
    \MAs_a : \Mca \otimes \Mca \to \Mca,
    \qquad a \in [\gamma],
\end{equation}
satisfying, for all $x, y, z \in \Mca$, the relations
\begin{equation} \label{equ:relations_multiassociative_algebras}
    (x \MAs_a y) \MAs_b z =
    (x \MAs_b y) \MAs_{a'} z =
    x \MAs_{a''} (y \MAs_b z) =
    x \MAs_b (y \MAs_{a'''} z),
    \qquad a, a', a'', a''' \leq b \in [\gamma].
\end{equation}
These algebras are obvious generalizations of associative algebras since
all of its operations are associative. Observe that
by~\eqref{equ:relations_multiassociative_algebras}, all bracketings of
an expression involving elements of a $\gamma$-multiassociative algebra
and some of its operations are equal. Then, since the bracketings of
such expressions are not significant, we shall denote these without
parenthesis. In upcoming Section~\ref{sec:as_gamma}, we will study the
underlying operads of the category of $\gamma$-multiassociative
algebras, called $\As_\gamma$, for a very specific purpose.
\medbreak

If $\Mca_1$ and $\Mca_2$ are two $\gamma$-multiassociative algebras,
a linear map $\phi : \Mca_1 \to \Mca_2$ is a
\Def{$\gamma$-multiassociative algebra morphism} if it commutes with
the operations of $\Mca_1$ and $\Mca_2$. We say that $\Mca$ is
\Def{commutative} when all operations of $\Mca$ are commutative.
Besides, for an $a \in [\gamma]$, an element $\Unit$ of $\Mca$ is an
\Def{$a$-unit}, or simply a \Def{unit} when taking into account the
value of $a$ is not necessary, of $\Mca$ if for all $x \in \Mca$,
$\Unit \MAs_a x = x = x \MAs_a \Unit$. When $\Mca$ admits a unit, we say
that $\Mca$ is \Def{unital}. As the following proposition shows, the
presence of a unit in $\Mca$ has some implications.
\medbreak

\begin{Proposition} \label{prop:units_multiassociative_algebras}
    Let $\gamma \geq 0$ be an integer and $\Mca$ be a
    $\gamma$-multiassociative algebra admitting a $b$-unit $\Unit$
    for a $b \in [\gamma]$. Then
    \begin{enumerate}[label={\it (\roman*)}]
        \item for all $a \in [b]$, the operations $\MAs_a$ and $\MAs_b$
        of $\Mca$ are equal;
        \item $\Unit$ is also an $a$-unit for all $a \in [b]$;
        \item $\Unit$ is the only unit of~$\Mca$.
    \end{enumerate}
\end{Proposition}
\medbreak

Relying on Proposition~\ref{prop:units_multiassociative_algebras},
similarly to the case of $\gamma$-pluriassociative algebras, we define
the \Def{height} of a $\gamma$-multiassociative algebra $\Mca$ as zero
if $\Mca$ has no unit, otherwise as the greatest integer $h \in [\gamma]$
such that the unit $\Unit$ of $\Mca$ is an $h$-unit.
\medbreak

\subsubsection{Multiprojection algebras}
A \Def{$\gamma$-multiprojection algebra} is a $\gamma$-multiassociative
algebra $\Mca$ endowed with endomorphisms
\begin{equation}
    \pi_a : \Mca \to \Mca,
    \qquad a \in [\gamma],
\end{equation}
satisfying
\begin{equation} \label{equ:relation_algebre_multiproj}
    \pi_a \circ \pi_{a'} = \pi_{a \Max a'},
    \qquad a, a' \in [\gamma].
\end{equation}
\medbreak

By extension, the \Def{height} of $\Mca$ is its height as a
$\gamma$-multiassociative algebra. We say that $\Mca$ is \Def{unital} as
a $\gamma$-multiprojection algebra if $\Mca$ is unital as a
$\gamma$-multiassociative algebra and its only, by
Proposition~\ref{prop:units_multiassociative_algebras}, unit $\Unit$
satisfies $\pi_a(\Unit) = \Unit$ for all $a \in [h]$ where $h$ is the
height of~$\Mca$.
\medbreak

\subsubsection{From multiprojection algebras to pluriassociative
algebras}
The next result describes how to construct $\gamma$-pluriassociative
algebras from $\gamma$-multiprojection algebras.
\medbreak

\begin{Theorem}
\label{thm:multiprojection_algebra_to_dias_gamma_algebra}
    For any integer $\gamma \geq 0$ and any $\gamma$-multiprojection
    algebra $\Mca$, the vector space $\Mca$ endowed with binary linear
    operations $\LDias_a$, $\RDias_a$, $a \in [\gamma]$, defined for all
    $x, y \in \Mca$ by
    \begin{subequations}
    \begin{equation}
        x \LDias_a y := x \MAs_a \pi_a(y)
    \end{equation}
    and
    \begin{equation}
        x \RDias_a y := \pi_a(x) \MAs_a y,
    \end{equation}
    \end{subequations}
    where the $\MAs_a$, $a \in [\gamma]$, are the operations of $\Mca$ and
    the $\pi_a$, $a \in [\gamma]$, are its endomorphisms, is a
    $\gamma$-pluriassociative algebra, denoted by $\MProjToPluri(\Mca)$.
\end{Theorem}
\begin{proof}
    This is a verification of the relations of $\gamma$-pluriassociative
    algebras in $\MProjToPluri(\Mca)$. Let $x$, $y$, and $z$ be three
    elements of $\MProjToPluri(\Mca)$ and  $a, a' \in [\gamma]$.
    \smallbreak

    By~\eqref{equ:relations_multiassociative_algebras}, we have
    \begin{equation}
        (x \RDias_{a'} y) \LDias_a z =
        \pi_{a'}(x) \MAs_{a'} y \MAs_a \pi_a(z) =
        x \RDias_{a'} (y \LDias_a z),
    \end{equation}
    showing that~\eqref{equ:relation_dias_gamma_1_concise} is satisfied
    in $\MProjToPluri(\Mca)$.
    \smallbreak

    Moreover, by~\eqref{equ:relations_multiassociative_algebras}
    and~\eqref{equ:relation_algebre_multiproj}, we have
    \begin{equation} \begin{split}
        x \LDias_a (y \RDias_{a'} z)
            & = x \MAs_a \pi_a(\pi_{a'}(y) \MAs_{a'} z) \\
            & = x \MAs_a \pi_{a \Max a'}(y) \MAs_{a'} \pi_a(z) \\
            & = x \MAs_{a \Max a'} \pi_{a \Max a'}(y) \MAs_a \pi_a(z) \\
            & = (x \LDias_{a \Max a'} y) \LDias_a z,
    \end{split} \end{equation}
    so that~\eqref{equ:relation_dias_gamma_2_concise}, and for the same
    reasons~\eqref{equ:relation_dias_gamma_3_concise}, check out in
    $\MProjToPluri(\Mca)$.
    \smallbreak

    Finally, again by~\eqref{equ:relations_multiassociative_algebras}
    and~\eqref{equ:relation_algebre_multiproj}, we have
    \begin{equation} \begin{split}
        x \LDias_a (y \LDias_{a'} z)
            & = x \MAs_a \pi_a(y \MAs_{a'} \pi_{a'}(z)) \\
            & = x \MAs_a \pi_a(y) \MAs_{a'} \pi_{a \Max a'}(z) \\
            & = x \MAs_a \pi_a(y) \MAs_{a \Max a'} \pi_{a \Max a'}(z) \\
            & = (x \LDias_a y) \LDias_{a \Max a'} z,
    \end{split} \end{equation}
    showing that~\eqref{equ:relation_dias_gamma_4_concise},
    and for the same reasons~\eqref{equ:relation_dias_gamma_5_concise},
    are satisfied in $\MProjToPluri(\Mca)$.
\end{proof}
\medbreak

When $\Mca$ is commutative, since for all $x, y \in \MProjToPluri(\Mca)$
and $a \in [\gamma]$,
\begin{equation}
    x \LDias_a y = x \MAs_a \pi_a(y) =
    \pi_a(y) \MAs_a x = y \RDias_a x,
\end{equation}
it appears that $\MProjToPluri(\Mca)$ is a commutative
$\gamma$-pluriassociative algebra.
\medbreak

When $\Mca$ is unital, $\MProjToPluri(\Mca)$ has several properties,
summarized in the next proposition.
\medbreak

\begin{Proposition}
\label{prop:multiprojection_algebra_to_dias_gamma_algebra_properties}
    Let $\gamma \geq 0$ be an integer, $\Mca$ be a unital
    $\gamma$-multiprojection algebra of height~$h$.
    Then, by denoting by $\Unit$ the unit of $\Mca$ and by $\pi_a$,
    $a \in [\gamma]$, its endomorphisms,
    \begin{enumerate}[label={\it (\roman*)}]
        \item for any $a \in [h]$, $\Unit$ is an $a$-bar-unit of
        $\MProjToPluri(\Mca)$;
        \item
\label{item:multiprojection_algebra_to_dias_gamma_algebra_properties_2}
        for any $a \leq b \in [h]$, $\Halo_a(\MProjToPluri(\Mca))$
        is a subset of $\Halo_b(\MProjToPluri(\Mca))$;
        \item for any $a \in [h]$, the linear span of
        $\Halo_a(\MProjToPluri(\Mca))$ forms an
        $h\!-\!a\!+\!1$-pluriassociative subalgebra of the
        $h\!-\!a\!+\!1$-pluriassociative subalgebra of
        $\MProjToPluri(\Mca)$ induced by $[a, h]$;
        \item for any $a \in [h]$, $\pi_a$ is the identity map if and
        only if $\Unit$ is an $a$-wire-unit of~$\MProjToPluri(\Mca)$.
    \end{enumerate}
\end{Proposition}
\medbreak

\subsubsection{Examples of constructions of pluriassociative algebras}
The construction $\MProjToPluri$ of
Theorem~\ref{thm:multiprojection_algebra_to_dias_gamma_algebra} allows
to build several $\gamma$-pluriassociative algebras. A few examples
follow.
\medbreak

\paragraph{The $\gamma$-pluriassociative algebra of positive integers}
Let $\gamma \geq 1$ be an integer and consider the vector space
$\AlgPos$ spanned by positive integers, endowed with the operations
$\MAs_a$,$a \in [\gamma]$, all equal to the operation $\Max$ extended by
linearity and with the endomorphisms $\pi_a$, $a \in [\gamma]$,
linearly defined for any positive integer $x$ by $\pi_a(x) := a \Max x$.
Then, $\AlgPos$ is a non-unital $\gamma$-multiprojection algebra. By
Theorem~\ref{thm:multiprojection_algebra_to_dias_gamma_algebra},
$\MProjToPluri(\AlgPos)$ is a $\gamma$-pluriassociative algebra. We have
for instance
\begin{subequations}
\begin{equation}
    \textcolor{Col1}{2} \LDias_3 \textcolor{Col4}{5} =
    \textcolor{Col4}{5},
\end{equation}
\begin{equation}
    \textcolor{Col1}{1} \RDias_3 \textcolor{Col4}{2} = 3.
\end{equation}
\end{subequations}
We can observe that $\MProjToPluri(\AlgPos)$ is commutative, pure, and
its $1$-halo is $\{1\}$. Moreover, when
$\gamma \geq 2$, $\MProjToPluri(\AlgPos)$ has no wire-unit and no
$a$-bar-unit for $a \geq 2 \in [\gamma]$. This example is important
because it provides a counterexample
for~\ref{item:multiprojection_algebra_to_dias_gamma_algebra_properties_2}
of Proposition~\ref%
{prop:multiprojection_algebra_to_dias_gamma_algebra_properties}
in the case when the construction $\MProjToPluri$ is applied
to a non-unital $\gamma$-multiprojection algebra.
\medbreak

\paragraph{The $\gamma$-pluriassociative algebra of finite sets}
Let $\gamma \geq 1$ be an integer and consider the vector space
$\AlgSets$ of finite sets of positive integers, endowed with the
operations $\MAs_a$, $a \in [\gamma]$, all equal to the union operation
$\cup$ extended by linearity and with the endomorphisms $\pi_a$,
$a \in [\gamma]$, linearly defined for any finite set of positive
integers $x$ by $\pi_a(x) := x \cap [a, \gamma]$. Then, $\AlgSets$ is a
$\gamma$-multiprojection algebra. By
Theorem~\ref{thm:multiprojection_algebra_to_dias_gamma_algebra},
$\MProjToPluri(\AlgSets)$ is a $\gamma$-pluriassociative algebra. We
have for instance
\begin{subequations}
\begin{equation}
    \{\textcolor{Col1}{2}, \textcolor{Col1}{4}\}
    \LDias_3
    \{\textcolor{Col4}{1}, \textcolor{Col4}{3},
        \textcolor{Col4}{5}\}
    = \{\textcolor{Col1}{2}, \textcolor{Col4}{3},
    \textcolor{Col1}{4}, \textcolor{Col4}{5}\},
\end{equation}
\begin{equation}
    \{\textcolor{Col1}{1}, \textcolor{Col1}{2},
    \textcolor{Col1}{4}\}
    \RDias_3
    \{\textcolor{Col4}{1}, \textcolor{Col4}{3},
        \textcolor{Col4}{5}\}
    = \{\textcolor{Col4}{1}, \textcolor{Col4}{3},
    \textcolor{Col1}{4}, \textcolor{Col4}{5}\}.
\end{equation}
\end{subequations}
We can observe that $\MProjToPluri(\AlgSets)$ is commutative and pure.
Moreover, $\emptyset$ is a $1$-wire-unit of $\MProjToPluri(\AlgSets)$
and, by Proposition~\ref{prop:wire_unit_dias_gamma}, it is its only
wire-unit. Therefore, $\MProjToPluri(\AlgSets)$ has height $1$. Observe
that for any $a \in [\gamma]$, the $a$-halo of $\MProjToPluri(\AlgSets)$
consists in the subsets of $[a - 1]$. Besides, since $\AlgSets$ is a
unital $\gamma$-multiprojection algebra, $\MProjToPluri(\AlgSets)$
satisfies all properties exhibited by
Proposition~\ref%
{prop:multiprojection_algebra_to_dias_gamma_algebra_properties}.
\medbreak

\paragraph{The $\gamma$-pluriassociative algebra of words}
Let $\gamma \geq 1$ be an integer and consider the vector space
$\AlgWords$ of the words of positive integers. Let us endow $\AlgWords$
with the operations $\MAs_a$, $a \in [\gamma]$, all equal to the
concatenation operation extended by linearity and with the
endomorphisms $\pi_a$, $a \in [\gamma]$, where for any word $x$ of
positive integers, $\pi_a(x)$ is the longest subword of $x$ consisting
in letters greater than or equal to $a$. Then, $\AlgWords$ is a
$\gamma$-multiprojection algebra. By
Theorem~\ref{thm:multiprojection_algebra_to_dias_gamma_algebra},
$\MProjToPluri(\AlgWords)$ is a $\gamma$-pluriassociative algebra. We
have for instance
\begin{subequations}
\begin{equation}
    \textcolor{Col1}{412} \LDias_3
    \textcolor{Col4}{14231} =
    \textcolor{Col1}{412}\textcolor{Col4}{43},
\end{equation}
\begin{equation}
    \textcolor{Col1}{11} \RDias_2 \textcolor{Col4}{323} =
    \textcolor{Col4}{323}.
\end{equation}
\end{subequations}
We can observe that $\MProjToPluri(\AlgWords)$ is not commutative and
is pure. Moreover, $\epsilon$ is a $1$-wire-unit of
$\MProjToPluri(\AlgWords)$ and by
Proposition~\ref{prop:wire_unit_dias_gamma}, it is its only wire-unit.
Therefore, $\MProjToPluri(\AlgWords)$ has height $1$. Observe that for
any $a \in [\gamma]$, the $a$-halo of $\MProjToPluri(\AlgWords)$
consists in the words on the alphabet $[a - 1]$. Besides, since
$\AlgWords$ is a unital $\gamma$-multiprojection algebra,
$\MProjToPluri(\AlgWords)$ satisfies all properties exhibited by
Proposition~\ref%
{prop:multiprojection_algebra_to_dias_gamma_algebra_properties}.
\medbreak

The $\gamma$-pluriassociative algebras $\MProjToPluri(\AlgSets)$ and
$\MProjToPluri(\AlgWords)$ are related in the following way. Let
$I_{\rm com}$ be the subspace of $\MProjToPluri(\AlgWords)$ generated
by the $x - x'$ where $x$ and $x'$ are words of positive integers and
have the same commutative image. Since $I_{\rm com}$ is a
$\gamma$-pluriassociative algebra ideal of $\MProjToPluri(\AlgWords)$,
one can consider the quotient $\gamma$-pluriassociative algebra
$\AlgComWords := \MProjToPluri(\AlgWords)/_{I_{\rm com}}$. Its elements
can be seen as commutative words of positive integers.
\medbreak

Moreover, let $I_{\rm occ}$ be the subspace of
$\MProjToPluri(\AlgComWords)$ generated by the $x - x'$ where $x$ and
$x'$ are commutative words of positive integers and for any letter
$a \in [\gamma]$, $a$ appears in $x$ if and only if $a$ appears in $x'$.
Since $I_{\rm occ}$ is a $\gamma$-pluriassociative algebra ideal of
$\MProjToPluri(\AlgComWords)$, one can consider the quotient
$\gamma$-pluriassociative algebra
$\MProjToPluri(\AlgComWords)/_{I_{\rm occ}}$. Its elements can be seen
as finite subsets of positive integers and we observe that
$\MProjToPluri(\AlgComWords)/_{I_{\rm occ}} = \MProjToPluri(\AlgSets)$.
\medbreak

\paragraph{The $\gamma$-pluriassociative algebra of marked words}
Let $\gamma \geq 1$ be an integer and consider the vector space
$\AlgMarkWords$ of the words of positive integers where letters can be
marked or not, with at least one occurrence of a marked letter. We
denote by $\bar a$ any \Def{marked letter} $a$ and we say that the
\Def{value} of $\bar a$ is $a$. Let us endow $\AlgMarkWords$ with the
linear operations $\MAs_a$, $a \in [\gamma]$, where for all words $u$
and $v$ of $\AlgMarkWords$, $u \MAs_a v$ is obtained by concatenating
$u$ and $v$, and by replacing therein all marked letters by $\bar c$
where $c := \max(u) \Max a \Max \max(v)$ where $\max(u)$ (resp.
$\max(v)$) denotes the greatest value among the marked letters of $u$
(resp. $v$). For instance,
\begin{subequations}
\begin{equation}
    \textcolor{Col1}{2 \bar 1 3 1 \bar 3}
    \MAs_2
    \textcolor{Col4}{3 \bar 4 \bar 1 2 1}
    =
    \textcolor{Col1}{2 {\bf \bar 4} 3 1 {\bf \bar 4}}
    \textcolor{Col4}{3 \bar 4 {\bf \bar 4} 2 1},
\end{equation}
\begin{equation}
    \textcolor{Col1}{\bar 2 1 1 \bar 1}
    \MAs_3
    \textcolor{Col4}{3 4 \bar 2} =
    \textcolor{Col1}{{\bf \bar 3} 1 1 {\bf \bar 3}}
    \textcolor{Col4}{3 4 {\bf \bar 3}}.
\end{equation}
\end{subequations}
We also endow $\AlgMarkWords$ with the endomorphisms $\pi_a$,
$a \in [\gamma]$, where for any word $u$ of $\AlgMarkWords$, $\pi_a(u)$
is obtained by replacing in $u$ any occurrence of a nonmarked letter
smaller than $a$ by $a$. For instance,
\begin{equation}
    \pi_3\left(\textcolor{Col4}{2} \textcolor{Col1}{\bar 2}
        \textcolor{Col4}{14} \textcolor{Col1}{\bar 4}
        \textcolor{Col4}{3} \textcolor{Col1}{\bar 5}\right)
    = \textcolor{Col4}{\bf 3} \textcolor{Col1}{\bar 2}
        \textcolor{Col4}{{\bf 3} 4} \textcolor{Col1}{\bar 4}
        \textcolor{Col4}{\bf 3} \textcolor{Col1}{\bar 5}.
\end{equation}
One can show without difficulty that $\AlgMarkWords$ is a
$\gamma$-multiprojection algebra. By
Theorem~\ref{thm:multiprojection_algebra_to_dias_gamma_algebra},
$\MProjToPluri(\AlgMarkWords)$ is a $\gamma$-pluriassociative algebra.
We have for instance
\begin{subequations}
\begin{equation}
    \textcolor{Col1}{3 \bar 2 5}
    \LDias_3
    \textcolor{Col4}{4 \bar 4 1}
    = \textcolor{Col1}{3 {\bf \bar 4} 5}
        \textcolor{Col4}{4 \bar 4 {\bf 3}},
\end{equation}
\begin{equation}
    \textcolor{Col1}{1 \bar 3 4 \bar 1 3}
    \RDias_2
    \textcolor{Col4}{3 1 \bar 2 3 \bar 1 1}
    = \textcolor{Col1}{{\bf 2} \bar 3 4 {\bf \bar 3} 3}
      \textcolor{Col4}{3 1 {\bf \bar 3} 3 {\bf \bar 3} 1}.
\end{equation}
\end{subequations}
We can observe that $\MProjToPluri(\AlgMarkWords)$ is not commutative,
pure, and has no wire-units neither bar-units.
\medbreak

\paragraph{The free $\gamma$-pluriassociative algebra over one
generator}
Let $\gamma \geq 0$ be an integer. We give here a construction of the
free $\gamma$-pluriassociative algebra $\FreeAlg_{\Dias_\gamma}$ over
one generator described in
Section~\ref{subsubsec:free_dias_gamma_algebras}
passing through the following $\gamma$-multiprojection algebra and the
construction $\MProjToPluri$. Consider the vector space of
nonempty words on the alphabet $\{0\} \sqcup [\gamma]$ with exactly one
occurrence of $0$, endowed with the operations $\MAs_a$,
$a \in [\gamma]$, all equal to the concatenation operation extended by
linearity and with the endomorphisms $\Augm_a$, $a \in [\gamma]$,
defined in Section~\ref{subsubsec:free_dias_gamma_algebras}. This vector
space is a $\gamma$-multiprojection algebra. Therefore, by
Theorem~\ref{thm:multiprojection_algebra_to_dias_gamma_algebra}, it
gives rise by the construction $\MProjToPluri$ to a
$\gamma$-pluriassociative algebra and it appears that it is
$\FreeAlg_{\Dias_\gamma}$. Besides, we can now observe that
$\FreeAlg_{\Dias_\gamma}$ is not commutative, pure, and has no
wire-units neither bar-units.
\medbreak

\section{Pluritriassociative operads} \label{sec:pluritriass}
We describe in this section a generalization on a nonnegative integer
parameter $\gamma$ of the triassociative operad~\cite{LR04}.
\medbreak

\subsection{Construction and first properties}\label{subsec:trias_gamma}
Our original idea of using the $\T$ construction (see
Section~\ref{subsubsec:construction_dias_gamma}) to obtain a
generalization of the diassociative operad admits an analogue in the
context of the triassociative operad. Let us describe it.
\medbreak

\subsubsection{Construction}
For any integer $\gamma \geq 0$, we define the
\Def{$\gamma$-pluritriassociative operad} $\Trias_\gamma$ as the
suboperad of $\T \Nmax_\gamma$ generated by
\begin{equation}
   \GenTrias := \left\{0a, 00, a0 : a \in [\gamma]\right\}.
\end{equation}
By definition, $\Trias_\gamma$ is the vector space of words that can be
obtained by partial compositions of words
of $\GenTrias$. We have, for instance,
\begin{subequations}
\begin{equation}
    \Trias_2(1) = \K \Angle{\{0\}},
\end{equation}
\begin{equation}
    \Trias_2(2) = \K \Angle{\{00, 01, 02, 10, 20\}},
\end{equation}
\begin{multline}
    \Trias_2(3)
    = \K \Angle{\{000, 001, 002, 010, 011, 012, 020, 021, \right. \\
        \left. 022, 100, 101, 102, 110, 120, 200, 201, 202, 210, 220\}},
\end{multline}
\end{subequations}
\medbreak

\subsubsection{First properties}
In the first place, observe that $\Trias_1$ is the operad $\Tr$ defined
in Chapter~\ref{chap:monoids}. For this reason, $\Trias_1$ is the
triassociative operad $\Trias$. Moreover, observe that $\Trias_0$ is
the trivial operad and that $\Trias_\gamma$ is a suboperad of
$\Trias_{\gamma + 1}$. Then, for all integers $\gamma \geq 0$, the
operads $\Trias_\gamma$ are generalizations of the triassociative
operad. Observe that since $\GenTrias = \GenDias \sqcup \{00\}$,
$\Dias_\gamma$ is a suboperad of $\Trias_\gamma$. Finally, remark that
the fundamental basis of $\Trias_\gamma$ is a set-operad basis.
\medbreak

\subsubsection{Elements and dimensions}
\begin{Proposition} \label{prop:basis_trias_gamma}
    For any integer $\gamma \geq 0$, the fundamental basis of
    $\Trias_\gamma$ is the set of all the words on the alphabet
    $\{0\} \sqcup [\gamma]$ containing at least one occurrence of~$0$.
\end{Proposition}
\medbreak

We deduce from Proposition~\ref{prop:basis_trias_gamma} that the
Hilbert series of $\Trias_\gamma$ satisfies
\begin{equation}
    \Hca_{\Trias_\gamma}(t) =
    \frac{t}{(1 - \gamma t)(1 - \gamma t - t)}
\end{equation}
and that for all $n \geq 1$,
$\dim \Trias_\gamma(n) = (\gamma + 1)^n - \gamma^n$. For instance, the
first dimensions of  $\Trias_1$, $\Trias_2$, $\Trias_3$, and $\Trias_4$
are respectively
\begin{subequations}
\begin{equation}
    1, 3, 7, 15, 31, 63, 127, 255, 511, 1023, 2047,
\end{equation}
\begin{equation}
    1, 5, 19, 65, 211, 665, 2059, 6305, 19171, 58025, 175099,
\end{equation}
\begin{equation}
    1, 7, 37, 175, 781, 3367, 14197, 58975, 242461, 989527, 4017157,
\end{equation}
\begin{equation}
    1, 9, 61, 369, 2101, 11529, 61741, 325089, 1690981, 8717049,
    44633821.
\end{equation}
\end{subequations}
These sequences are respectively Sequences~\OEIS{A000225},
\OEIS{A001047}, \OEIS{A005061}, and \OEIS{A005060} of~\cite{Slo}.
\medbreak

\subsection{Additional properties}
We exhibit here a presentation of $\Trias_\gamma$ and establish the fact
that it is a Koszul operad.
\medbreak

\subsubsection{Presentation by generators and relations}
For any $a \in [\gamma]$, let us denote by $\LDias_a$ (resp. $\RDias_a$,
$\MTrias$) the generator $0a$ (resp. $a0$, $00$) of~$\Trias_\gamma$.
\medbreak

\begin{Theorem} \label{thm:presentation_trias_gamma}
    For any integer $\gamma \geq 0$, the operad $\Trias_\gamma$ admits
    the presentation $\left(\GenTrias, \RelTrias\right)$ where
    $\RelTrias$ is the space induced by the
    equivalence relation $\Equiv_\gamma$ satisfying
    \begin{subequations}
    \begin{equation}\label{equ:relation_presentation_trias_gamma_1}
        \Corolla{\MTrias} \circ_1 \Corolla{\MTrias}
         \Equiv_\gamma
        \Corolla{\MTrias} \circ_2 \Corolla{\MTrias},
    \end{equation}
    \begin{equation}\label{equ:relation_presentation_trias_gamma_2}
        \Corolla{\LDias_a} \circ_1 \Corolla{\MTrias}
         \Equiv_\gamma
        \Corolla{\MTrias} \circ_2 \Corolla{\LDias_a},
        \qquad a \in [\gamma],
    \end{equation}
    \begin{equation}\label{equ:relation_presentation_trias_gamma_3}
        \Corolla{\MTrias} \circ_1 \Corolla{\RDias_a}
         \Equiv_\gamma
        \Corolla{\RDias_a} \circ_2 \Corolla{\MTrias},
        \qquad a \in [\gamma],
    \end{equation}
    \begin{equation}\label{equ:relation_presentation_trias_gamma_4}
        \Corolla{\MTrias} \circ_1 \Corolla{\LDias_a}
         \Equiv_\gamma
        \Corolla{\MTrias} \circ_2 \Corolla{\RDias_a},
        \qquad a \in [\gamma],
    \end{equation}
    \begin{equation}\label{equ:relation_presentation_trias_gamma_5}
        \Corolla{\LDias_a} \circ_1 \Corolla{\RDias_{a'}}
         \Equiv_\gamma
        \Corolla{\RDias_{a'}} \circ_2 \Corolla{\LDias_a},
        \qquad a, a' \in [\gamma],
    \end{equation}
    \begin{equation}\label{equ:relation_presentation_trias_gamma_6}
        \Corolla{\LDias_a} \circ_1 \Corolla{\LDias_b}
         \Equiv_\gamma
        \Corolla{\LDias_a} \circ_2 \Corolla{\RDias_b},
        \qquad a < b \in [\gamma],
    \end{equation}
    \begin{equation}\label{equ:relation_presentation_trias_gamma_7}
        \Corolla{\RDias_a} \circ_1 \Corolla{\LDias_b}
         \Equiv_\gamma
        \Corolla{\RDias_a} \circ_2 \Corolla{\RDias_b},
        \qquad a < b \in [\gamma],
    \end{equation}
    \begin{equation}\label{equ:relation_presentation_trias_gamma_8}
        \Corolla{\LDias_b} \circ_1 \Corolla{\LDias_a}
         \Equiv_\gamma
        \Corolla{\LDias_a} \circ_2 \Corolla{\LDias_b},
        \qquad a < b \in [\gamma],
    \end{equation}
    \begin{equation}\label{equ:relation_presentation_trias_gamma_9}
        \Corolla{\RDias_a} \circ_1 \Corolla{\RDias_b}
         \Equiv_\gamma
        \Corolla{\RDias_b} \circ_2 \Corolla{\RDias_a},
        \qquad a < b \in [\gamma],
    \end{equation}
    \begin{equation}\label{equ:relation_presentation_trias_gamma_10}
        \Corolla{\LDias_d} \circ_1 \Corolla{\LDias_d}
         \Equiv_\gamma
        \Corolla{\LDias_d} \circ_2 \Corolla{\MTrias}
         \Equiv_\gamma
        \Corolla{\LDias_d} \circ_2 \Corolla{\LDias_c}
         \Equiv_\gamma
        \Corolla{\LDias_d} \circ_2 \Corolla{\RDias_c},
        \qquad c \leq d \in [\gamma],
    \end{equation}
    \begin{equation}\label{equ:relation_presentation_trias_gamma_11}
        \Corolla{\RDias_d} \circ_1 \Corolla{\LDias_c}
         \Equiv_\gamma
        \Corolla{\RDias_d} \circ_1 \Corolla{\RDias_c}
         \Equiv_\gamma
        \Corolla{\RDias_d} \circ_1 \Corolla{\MTrias}
         \Equiv_\gamma
        \Corolla{\RDias_d} \circ_2 \Corolla{\RDias_d},
        \qquad c \leq d \in [\gamma].
    \end{equation}
    \end{subequations}
\end{Theorem}
\medbreak

In the same fashion as we have done for
Theorem~\ref{thm:presentation_dias_gamma}, our proof of
Theorem~\ref{thm:presentation_trias_gamma} is based upon the computation
of the kernel of the evaluation morphism
\begin{equation}
    \Eval : \FreeOperad\left(\GenTrias\right) \to \Trias_\gamma.
\end{equation}
In this case, the image of a $\GenTrias$-syntax tree $\Tfr$ can be
computed in the same way as in the case of $\GenDias$-syntax trees
(see Section~\ref{subsubsec:presentation_dias_gamma}). The internal
nodes of $\Tfr$ labeled by $\MTrias$ do not play any role in this
computation (see Figure~\ref{fig:example_eval_trias_gamma}).
\begin{figure}[ht]
    \centering

    \caption[A syntax tree on $\GenTrias$ and its evaluation.]
    {A $\GenTrias$-syntax tree $\Tfr$ where the images of its leaves
    are shown. This tree satisfies
    $\Eval(\Tfr) = \textcolor{Col1}{332440433201}$.}
    \label{fig:example_eval_trias_gamma}
\end{figure}
\medbreak

The space of relations $\RelTrias$ of $\Trias_\gamma$ exhibited by
Theorem~\ref{thm:presentation_trias_gamma} can be rephrased a bit more
concisely as the space generated by
\begin{subequations}
\begin{equation} \label{equ:relation_trias_gamma_1_concise}
    \Corolla{\MTrias} \circ_1 \Corolla{\MTrias}
    -
    \Corolla{\MTrias} \circ_2 \Corolla{\MTrias},
\end{equation}
\begin{equation} \label{equ:relation_trias_gamma_2_concise}
    \Corolla{\LDias_a} \circ_1 \Corolla{\MTrias}
    -
    \Corolla{\MTrias} \circ_2 \Corolla{\LDias_a},
    \qquad a \in [\gamma],
\end{equation}
\begin{equation} \label{equ:relation_trias_gamma_3_concise}
    \Corolla{\MTrias} \circ_1 \Corolla{\RDias_a}
    -
    \Corolla{\RDias_a} \circ_2 \Corolla{\MTrias},
    \qquad a \in [\gamma],
\end{equation}
\begin{equation} \label{equ:relation_trias_gamma_4_concise}
    \Corolla{\MTrias} \circ_1 \Corolla{\LDias_a}
    -
    \Corolla{\MTrias} \circ_2 \Corolla{\RDias_a},
    \qquad a \in [\gamma],
\end{equation}
\begin{equation} \label{equ:relation_trias_gamma_5_concise}
    \Corolla{\LDias_a} \circ_1 \Corolla{\RDias_{a'}}
    -
    \Corolla{\RDias_{a'}} \circ_2 \Corolla{\LDias_a},
    \qquad a, a' \in [\gamma],
\end{equation}
\begin{equation} \label{equ:relation_trias_gamma_6_concise}
    \Corolla{\LDias_a} \circ_1 \Corolla{\LDias_{a \Max a'}}
    -
    \Corolla{\LDias_a} \circ_2 \Corolla{\RDias_{a'}},
    \qquad a, a' \in [\gamma],
\end{equation}
\begin{equation} \label{equ:relation_trias_gamma_7_concise}
    \Corolla{\RDias_a} \circ_1 \Corolla{\LDias_{a'}}
    -
    \Corolla{\RDias_a} \circ_2 \Corolla{\RDias_{a \Max a'}},
    \qquad a, a' \in [\gamma],
\end{equation}
\begin{equation} \label{equ:relation_trias_gamma_8_concise}
    \Corolla{\LDias_{a \Max a'}} \circ_1 \Corolla{\LDias_a}
    -
    \Corolla{\LDias_a} \circ_2 \Corolla{\LDias_{a'}},
    \qquad a, a' \in [\gamma],
\end{equation}
\begin{equation} \label{equ:relation_trias_gamma_9_concise}
    \Corolla{\RDias_a} \circ_1 \Corolla{\RDias_{a'}}
    -
    \Corolla{\RDias_{a \Max a'}} \circ_2 \Corolla{\RDias_a},
    \qquad a, a' \in [\gamma],
\end{equation}
\begin{equation} \label{equ:relation_trias_gamma_10_concise}
    \Corolla{\LDias_a} \circ_1 \Corolla{\LDias_a}
    -
    \Corolla{\LDias_a} \circ_2 \Corolla{\MTrias},
    \qquad a \in [\gamma],
\end{equation}
\begin{equation} \label{equ:relation_trias_gamma_11_concise}
    \Corolla{\RDias_a} \circ_2 \Corolla{\RDias_a}
    -
    \Corolla{\RDias_a} \circ_1 \Corolla{\MTrias}.
    \qquad a \in [\gamma].
\end{equation}
\end{subequations}
\medbreak

\subsubsection{Koszulity}
\begin{Theorem} \label{thm:koszulity_trias_gamma}
    For any integer $\gamma \geq 0$, $\Trias_\gamma$ is a Koszul operad
    and the set of the $\GenTrias$-syntax trees avoiding the
    trees
    \begin{subequations}
    \begin{equation}
        \Corolla{\MTrias} \circ_2 \Corolla{\MTrias},
    \end{equation}
    \begin{equation}
        \Corolla{\LDias_a} \circ_1 \Corolla{\MTrias}
        \qquad a \in [\gamma],
    \end{equation}
    \begin{equation}
        \Corolla{\RDias_a} \circ_2 \Corolla{\MTrias},
        \qquad a \in [\gamma],
    \end{equation}
    \begin{equation}
        \Corolla{\MTrias} \circ_2 \Corolla{\RDias_a},
        \qquad a \in [\gamma],
    \end{equation}
    \begin{equation}
        \Corolla{\RDias_{a'}} \circ_2 \Corolla{\LDias_a},
        \qquad a, a' \in [\gamma],
    \end{equation}
    \begin{equation}
        \Corolla{\LDias_a} \circ_2 \Corolla{\RDias_{a'}},
        \qquad a, a' \in [\gamma],
    \end{equation}
    \begin{equation}
        \Corolla{\RDias_a} \circ_1 \Corolla{\LDias_{a'}},
        \qquad a, a' \in [\gamma],
    \end{equation}
    \begin{equation}
        \Corolla{\LDias_a} \circ_2 \Corolla{\LDias_{a'}},
        \qquad a, a' \in [\gamma],
    \end{equation}
    \begin{equation}
        \Corolla{\RDias_a} \circ_1 \Corolla{\RDias_{a'}},
        \qquad a, a' \in [\gamma],
    \end{equation}
    \begin{equation}
        \Corolla{\LDias_a} \circ_2 \Corolla{\MTrias},
        \qquad a \in [\gamma],
    \end{equation}
    \begin{equation}
        \Corolla{\RDias_a} \circ_1 \Corolla{\MTrias},
        \qquad a \in [\gamma].
    \end{equation}
    \end{subequations}
    is a Poincaré-Birkhoff-Witt basis of~$\Trias_\gamma$.
\end{Theorem}
\medbreak

\section{Polydendriform operads} \label{sec:dendr_gamma}
We introduce at this point our generalization on a nonnegative integer
parameter $\gamma$ of the dendriform operad and dendriform algebras. We
first construct this operad, compute its dimensions, and give then two
presentations by generators and relations. This section ends by a
description of free algebras over one generator in the category encoded
by our generalization.
\medbreak

\subsection{Construction and properties}
\label{subsec:construcion_dendr_gamma}
Theorem \ref{thm:presentation_dias_gamma}, by exhibiting a presentation
of $\Dias_\gamma$, shows that this operad is binary and quadratic. It
then admits a Koszul dual, denoted by $\Dendr_\gamma$ and called
\Def{$\gamma$-polydendriform operad}.
\medbreak

\subsubsection{Definition and presentation}
\label{subsubsec:presentation_dendr_gamma_alternative}
A description of $\Dendr_\gamma$ is provided by the following presentation
by generators and relations.
\medbreak

\begin{Theorem} \label{thm:presentation_dendr_gamma}
    For any integer $\gamma \geq 0$, the operad $\Dendr_\gamma$
    admits the presentation $\left(\GenDendr, \RelDendr\right)$ where
    \begin{math}
        \GenDendr := \GenDendr(2):=
        \left\{\LDendrA_a, \RDendrA_a : a \in [\gamma]\right\}
    \end{math}
    and $\RelDendr$ is the space generated by
    \begin{subequations}
    \begin{equation} \label{equ:relation_dendr_gamma_1_alternative}
        \Corolla{\LDendrA_a} \circ_1 \Corolla{\RDendrA_{a'}}
        -
        \Corolla{\RDendrA_{a'}} \circ_2 \Corolla{\LDendrA_a},
        \qquad a, a' \in [\gamma],
    \end{equation}
    \begin{equation} \label{equ:relation_dendr_gamma_2_alternative}
        \Corolla{\LDendrA_a} \circ_1 \Corolla{\LDendrA_b}
        -
        \Corolla{\LDendrA_a} \circ_2 \Corolla{\RDendrA_b},
        \qquad a < b \in [\gamma],
    \end{equation}
    \begin{equation} \label{equ:relation_dendr_gamma_3_alternative}
        \Corolla{\RDendrA_a} \circ_1 \Corolla{\LDendrA_b}
        -
        \Corolla{\RDendrA_a} \circ_2 \Corolla{\RDendrA_b},
        \qquad a < b \in [\gamma],
    \end{equation}
    \begin{equation} \label{equ:relation_dendr_gamma_4_alternative}
        \Corolla{\LDendrA_a} \circ_1 \Corolla{\LDendrA_b}
        -
        \Corolla{\LDendrA_a} \circ_2 \Corolla{\LDendrA_b},
        \qquad a < b \in [\gamma],
    \end{equation}
    \begin{equation} \label{equ:relation_dendr_gamma_5_alternative}
        \Corolla{\RDendrA_a} \circ_1 \Corolla{\RDendrA_b}
        -
        \Corolla{\RDendrA_a} \circ_2 \Corolla{\RDendrA_b},
        \qquad a < b \in [\gamma],
    \end{equation}
    \begin{equation} \label{equ:relation_dendr_gamma_6_alternative}
        \Corolla{\LDendrA_d} \circ_1 \Corolla{\LDendrA_d}
        -
        \left(\sum_{c \in [d]}
            \Corolla{\LDendrA_d} \circ_2 \Corolla{\LDendrA_c}
            +
            \Corolla{\LDendrA_d} \circ_2 \Corolla{\RDendrA_c}
        \right),
        \qquad d \in [\gamma],
    \end{equation}
    \begin{equation} \label{equ:relation_dendr_gamma_7_alternative}
        \left(\sum_{c \in [d]}
            \Corolla{\RDendrA_d} \circ_1 \Corolla{\RDendrA_c}
            +
            \Corolla{\RDendrA_d} \circ_1 \Corolla{\LDendrA_c}
        \right)
        -
        \Corolla{\RDendrA_d} \circ_2 \Corolla{\RDendrA_d},
        \qquad d \in [\gamma].
    \end{equation}
    \end{subequations}
\end{Theorem}
\medbreak

Theorem \ref{thm:presentation_dendr_gamma} provides a quite complicated
presentation of $\Dendr_\gamma$. We shall define below a more convenient
basis for the space of relations of~$\Dendr_\gamma$.
\medbreak

\subsubsection{Elements and dimensions}
\begin{Proposition} \label{prop:hilbert_series_dendr_gamma}
    For any integer $\gamma \geq 0$, the Hilbert series
    $\Hca_{\Dendr_\gamma}(t)$ of the operad $\Dendr_\gamma$ satisfies
    \begin{equation} \label{equ:hilbert_series_dendr_gamma}
        t + (2\gamma t - 1)\, \Hca_{\Dendr_\gamma}(t)
        + \gamma^2 t \, \Hca_{\Dendr_\gamma}(t)^2 = 0.
    \end{equation}
\end{Proposition}
\begin{proof}
    Let $G(t)$ be the generating series such that $G(-t)$
    satisfies~\eqref{equ:hilbert_series_dendr_gamma}. Therefore, $G(t)$
    satisfies
    \begin{equation}
        t = \frac{-G(t)}
                 {\left(1 + \gamma \, G(t)\right)^2}.
    \end{equation}
    Moreover, by setting $F(t) := \Hca_{\Dias_\gamma}(-t)$, where
    $\Hca_{\Dias_\gamma}(t)$ is the Hilbert series of $\Dias_\gamma$
    defined by~\eqref{equ:serie_hilbert_dias_gamma}, we have
    \begin{equation} \label{equ:hilbert_series_dendr_gamma_demo}
        F\left(G(t)\right)
            = \frac{-G(t)}
                  {\left(1 + \gamma \, G(t)\right)^2}
            = t,
    \end{equation}
    showing that $F(t)$ and $G(t)$ are the inverses for each other for
    series composition.
    \smallbreak

    Now, since  by Theorem~\ref{thm:koszulity_dias_gamma} and
    Proposition~\ref{prop:basis_dias_gamma}, $\Dias_\gamma$ is a Koszul
    operad and its Hilbert series is $\Hca_{\Dias_\gamma}(t)$, and
    since $\Dendr_\gamma$ is by definition the Koszul dual of
    $\Dias_\gamma$, the Hilbert series of these two operads
    satisfy Relation~\eqref{equ:Hilbert_series_Koszul_operads} of
    Chapter~\ref{chap:algebra}. Therefore,
    \eqref{equ:hilbert_series_dendr_gamma_demo} implies that
    the Hilbert series of $\Dendr_\gamma$ is $\Hca_{\Dendr_\gamma}(t)$.
\end{proof}
\medbreak

\begin{figure}[ht]
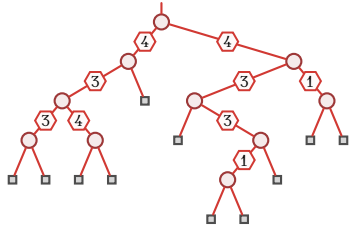

    \centering
    \begin{equation*}

    \end{equation*}
    \caption[A $4$-edge valued binary tree of the ns
    operad~$\Dendr_4(10)$.]
    {A $4$-edge valued binary tree of arity $10$. This tree is
    a basis element of $\Dendr_4(10)$.}
    \label{fig:example_edge_valued_binary_tree}
\end{figure}

By examining the expression for $\Hca_{\Dendr_\gamma}(t)$ of the
statement of Proposition~\ref{prop:hilbert_series_dendr_gamma}, we
observe that for any $n \geq 1$, $\Dendr_\gamma(n)$ can be seen as the
vector space $\FreeAlg_{\Dendr_\gamma}(n)$ of all binary trees with $n$
internal nodes wherein its $n - 1$ edges connecting two internal nodes
are labeled on $[\gamma]$. We call these trees \Def{$\gamma$-edge valued
binary trees}. In our graphical representations of $\gamma$-edge valued
binary trees, any edge label is drawn into a hexagon located half the
edge (see Figure~\ref{fig:example_edge_valued_binary_tree}).
\medbreak

We deduce from Proposition~\ref{prop:hilbert_series_dendr_gamma} that the
Hilbert series of $\Dendr_\gamma$ satisfies
\begin{equation}
    \Hca_{\Dendr_\gamma}(t) =
    \frac{1 - \sqrt{1 - 4 \gamma t} - 2 \gamma t}{2\gamma^2 t},
\end{equation}
and we also obtain that for all $n \geq 1$,
$\dim \Dendr_\gamma(n) = \gamma^{n - 1} \Catalan(n)$ where $\Catalan(n)$
is the number $\frac{1}{n + 1}\binom{2n}{n}$ of binary trees with $n$
internal nodes. For instance, the first dimensions of $\Dendr_1$,
$\Dendr_2$, $\Dendr_3$, and $\Dendr_4$ are respectively
\begin{subequations}
\begin{equation}
    1, 2, 5, 14, 42, 132, 429, 1430, 4862, 16796, 58786,
\end{equation}
\begin{equation}
    1, 4, 20, 112, 672, 4224, 27456, 183040, 1244672, 8599552, 60196864,
\end{equation}
\begin{equation}
    1, 6, 45, 378, 3402, 32076, 312741, 3127410, 31899582, 330595668,
    3471254514,
\end{equation}
\begin{equation}
    1, 8, 80, 896, 10752, 135168, 1757184, 23429120, 318636032,
    4402970624, 61641588736.
\end{equation}
\end{subequations}
These sequences are respectively Sequences~\OEIS{A000108},
\OEIS{A003645}, \OEIS{A101600}, and \OEIS{A269796} of~\cite{Slo}.
\medbreak

\subsubsection{Associative operations}
In the same manner as in the dendriform operad the sum of its two
operations produces an associative operation, in the $\gamma$-dendriform
operad there is a way to build associative operations, as the next
statement shows.
\medbreak

\begin{Proposition} \label{prop:associative_operation_dendr_gamma_other}
    For any integers $\gamma \geq 0$ and $b \in [\gamma]$, the element
    \begin{equation}
        \OpAsDendrA_b := \sum_{a \in [b]} \LDendrA_a + \RDendrA_a
    \end{equation}
    of $\Dendr_\gamma$ is associative.
\end{Proposition}
\medbreak

\subsubsection{Alternative presentation}
For any integer $\gamma \geq 0$, let $\LDendr_b$ and $\RDendr_b$,
$b \in [\gamma]$, be the elements of $\Dendr_\gamma$ defined by
\begin{subequations}
\begin{equation} \label{equ:definition_operation_dendr_gamma_left}
    \LDendr_b := \sum_{a \in [b]} \LDendrA_a,
\end{equation}
and
\begin{equation} \label{equ:definition_operation_dendr_gamma_right}
    \RDendr_b := \sum_{a \in [b]} \RDendrA_a.
\end{equation}
\end{subequations}
Then, since for all $b \in [\gamma]$ we have
\begin{subequations}
\begin{equation}
    \LDendrA_b =
    \begin{cases}
        \LDendr_1 & \mbox{if } b = 1, \\
        \LDendr_b - \LDendr_{b - 1} & \mbox{otherwise},
    \end{cases}
\end{equation}
and
\begin{equation}
    \RDendrA_b =
    \begin{cases}
        \RDendr_1 & \mbox{if } b = 1, \\
        \RDendr_b - \RDendr_{b - 1} & \mbox{otherwise},
    \end{cases}
\end{equation}
\end{subequations}
by triangularity, the family
\begin{equation}
    \GenDendr' := \left\{\LDendr_b, \RDendr_b : b \in [\gamma]\right\}
\end{equation}
is a generating set of $\Dendr_\gamma$. Remark that this change of basis
from is similar to the change of basis of $\Dias_\gamma$ considered in
Section~\ref{subsubsec:presentation_dias_gamma_alternative}. Let us now
express a presentation of $\Dendr_\gamma$ through the
family~$\GenDendr'$.
\medbreak

\begin{Theorem} \label{thm:autre_presentation_dendr_gamma}
    For any integer $\gamma \geq 0$, the operad $\Dendr_\gamma$
    admits the presentation $\left(\GenDendr', \RelDendr'\right)$
    where $\RelDendr'$ is the space generated by
    \begin{subequations}
    \begin{equation} \label{equ:relation_dendr_gamma_1}
        \Corolla{\LDendr_a} \circ_1 \Corolla{\RDendr_{a'}}
        -
        \Corolla{\RDendr_{a'}} \circ_2 \Corolla{\LDendr_a},
        \qquad a, a' \in [\gamma],
    \end{equation}
    \begin{equation} \label{equ:relation_dendr_gamma_2}
        \Corolla{\LDendr_a} \circ_1 \Corolla{\LDendr_b}
        -
        \Corolla{\LDendr_a} \circ_2 \Corolla{\RDendr_b}
        -
        \Corolla{\LDendr_a} \circ_2 \Corolla{\LDendr_a},
        \qquad a < b \in [\gamma],
    \end{equation}
    \begin{equation} \label{equ:relation_dendr_gamma_3}
        \Corolla{\RDendr_a} \circ_1 \Corolla{\RDendr_a}
        +
        \Corolla{\RDendr_a} \circ_1 \Corolla{\LDendr_b}
        -
        \Corolla{\RDendr_a} \circ_2 \Corolla{\RDendr_b},
        \qquad a < b \in [\gamma],
    \end{equation}
    \begin{equation} \label{equ:relation_dendr_gamma_4}
        \Corolla{\LDendr_b} \circ_1 \Corolla{\LDendr_a}
        -
        \Corolla{\LDendr_a} \circ_2 \Corolla{\LDendr_b}
        -
        \Corolla{\LDendr_a} \circ_2 \Corolla{\RDendr_a},
        \qquad a < b \in [\gamma],
    \end{equation}
    \begin{equation} \label{equ:relation_dendr_gamma_5}
        \Corolla{\RDendr_a} \circ_1 \Corolla{\LDendr_a}
        +
        \Corolla{\RDendr_a} \circ_1 \Corolla{\RDendr_b}
        -
        \Corolla{\RDendr_b} \circ_2 \Corolla{\RDendr_a},
        \qquad a < b \in [\gamma],
    \end{equation}
    \begin{equation} \label{equ:relation_dendr_gamma_6}
        \Corolla{\LDendr_a} \circ_1 \Corolla{\LDendr_a}
        -
        \Corolla{\LDendr_a} \circ_2 \Corolla{\RDendr_a}
        -
        \Corolla{\LDendr_a} \circ_2 \Corolla{\LDendr_a},
        \qquad a \in [\gamma],
    \end{equation}
    \begin{equation} \label{equ:relation_dendr_gamma_7}
        \Corolla{\RDendr_a} \circ_1 \Corolla{\RDendr_a}
        +
        \Corolla{\RDendr_a} \circ_1 \Corolla{\LDendr_a}
        -
        \Corolla{\RDendr_a} \circ_2 \Corolla{\RDendr_a},
        \qquad a \in [\gamma].
    \end{equation}
    \end{subequations}
\end{Theorem}
\begin{proof}
    Let us show that $\RelDendr'$ is equal to the space of relations
    $\RelDendr$ of $\Dendr_\gamma$ defined in the statement of
    Theorem~\ref{thm:presentation_dendr_gamma}. By this last theorem,
    for any $x \in \FreeOperad\left(\GenDendr\right)(3)$, $x$ is in
    $\RelDendr$ if and only if $\Eval(x) = 0$. By straightforward
    computations, by expanding any element $x$
    of~\eqref{equ:relation_dendr_gamma_1}---%
    \eqref{equ:relation_dendr_gamma_7} over the elements
    $\LDendrA_a$, $\RDendrA_a$, $a \in [\gamma]$, by
    using~\eqref{equ:definition_operation_dendr_gamma_left}
    and~\eqref{equ:definition_operation_dendr_gamma_right} we obtain
    that $x$ can be expressed as a sum of elements of $\RelDendr$. This
    implies that $\Eval(x) = 0$ and hence that $\RelDendr'$ is a
    subspace of $\RelDendr$. Now, one can observe that
    elements~\eqref{equ:relation_dendr_gamma_1}---%
    \eqref{equ:relation_dendr_gamma_6} are linearly independent.
    Then, $\RelDendr'$ has dimension $3\gamma^2$ which is also, by
    Theorem~\ref{thm:presentation_dendr_gamma}, the dimension of
    $\RelDendr$. The statement of the theorem follows.
\end{proof}
\medbreak

The presentation of $\Dendr_\gamma$ provided by
Theorem~\ref{thm:autre_presentation_dendr_gamma} is easier to handle
than the one provided by Theorem \ref{thm:presentation_dendr_gamma}. The
main reason is that
Relations~\eqref{equ:relation_dendr_gamma_6_alternative}
and~\eqref{equ:relation_dendr_gamma_7_alternative} of the first
presentation involve a nonconstant number of terms, while all relations
of this second presentation always involve only two or three terms. As a
very remarkable fact, it is worthwhile to note that the presentation of
$\Dendr_\gamma$ provided by
Theorem~\ref{thm:autre_presentation_dendr_gamma} can be directly
obtained by considering the Koszul dual of $\Dias_\gamma$ over the
$\BasisK$-basis (see Sections~\ref{subsubsec:basis_K_dias_gamma}
and~\ref{subsubsec:presentation_dias_gamma_alternative}). Therefore, an
alternative way to establish this presentation consists in computing the
Koszul dual of $\Dias_\gamma$ seen through the presentation having
$\RelDendr'$ as space of relations, which is made of the relations of
$\Dias_\gamma$ expressed over the $\BasisK$-basis (see
Proposition~\ref{prop:presentation_dias_gamma_alternative}).
\medbreak

From now on, $\Min$ denotes the operation $\min$ on integers. Using this
notation, the space of relations $\RelDendr'$ of $\Dendr_\gamma$
exhibited by Theorem~\ref{thm:autre_presentation_dendr_gamma}
can be rephrased in a more compact way as the space generated by
\begin{subequations}
\begin{equation} \label{equ:relation_dendr_gamma_1_concise}
    \Corolla{\LDendr_a} \circ_1 \Corolla{\RDendr_{a'}}
    -
    \Corolla{\RDendr_{a'}} \circ_2 \Corolla{\LDendr_a},
    \qquad a, a' \in [\gamma],
\end{equation}
\begin{equation} \label{equ:relation_dendr_gamma_2_concise}
    \Corolla{\LDendr_a} \circ_1 \Corolla{\LDendr_{a'}}
    -
    \Corolla{\LDendr_{a \Min a'}} \circ_2 \Corolla{\LDendr_a}
    -
    \Corolla{\LDendr_{a \Min a'}} \circ_2 \Corolla{\RDendr_{a'}},
    \qquad a, a' \in [\gamma],
\end{equation}
\begin{equation} \label{equ:relation_dendr_gamma_3_concise}
    \Corolla{\RDendr_{a \Min a'}} \circ_1 \Corolla{\LDendr_{a'}}
    +
    \Corolla{\RDendr_{a \Min a'}} \circ_1 \Corolla{\RDendr_a}
    -
    \Corolla{\RDendr_a} \circ_2 \Corolla{\RDendr_{a'}},
    \qquad a, a' \in [\gamma].
\end{equation}
\end{subequations}
\medbreak

Over the family $\GenDendr'$, one can build associative operations in
$\Dendr_\gamma$ in the following way.
\medbreak

\begin{Proposition} \label{prop:associative_operation_dendr_gamma}
    For any integers $\gamma \geq 0$ and $b \in [\gamma]$, the element
    \begin{equation}
        \OpAsDendr_b := \LDendr_b + \RDendr_b
    \end{equation}
    of $\Dendr_\gamma$ is associative. Moreover, any associative
    element of $\Dendr_\gamma$ is proportional to $\OpAsDendr_b$ for
    a~$b \in [\gamma]$.
\end{Proposition}
\medbreak

\subsection{Category of polydendriform algebras and free objects}
The aim of this section is to describe the category of
$\Dendr_\gamma$-algebras and more particularly the free
$\Dendr_\gamma$-algebra over one generator.
\medbreak

\subsubsection{Polydendriform algebras}
We call \Def{$\gamma$-polydendriform algebra} any
$\Dendr_\gamma$-algebra. From the presentation of $\Dendr_\gamma$
provided by Theorem~\ref{thm:presentation_dendr_gamma}, any
$\gamma$-polydendriform algebra is a vector space endowed with linear
operations $\LDendrA_a, \RDendrA_a$, $a \in [\gamma]$, satisfying the
relations encoded by~\eqref{equ:relation_dendr_gamma_1_alternative}---%
\eqref{equ:relation_dendr_gamma_7_alternative}. By considering the
presentation of $\Dendr_\gamma$ exhibited by
Theorem~\ref{thm:autre_presentation_dendr_gamma}, any
$\gamma$-polydendriform algebra is a vector space endowed with linear
operations $\LDendr_a, \RDendr_a$, $a \in [\gamma]$, satisfying the
relations encoded by~\eqref{equ:relation_dendr_gamma_1_concise}---%
\eqref{equ:relation_dendr_gamma_3_concise}.
\medbreak

\subsubsection{Two ways to split associativity}
Like dendriform algebras, which offer a way to split an associative
operation into two parts, $\gamma$-polydendriform algebras propose two
ways to split associativity depending on its chosen presentation.
\medbreak

On the one hand, in a $\gamma$-polydendriform algebra $\Dca$ over the
operations $\LDendrA_a$, $\RDendrA_a$, $a \in [\gamma]$, by
Proposition~\ref{prop:associative_operation_dendr_gamma_other}, an
associative operation $\OpAsDendrA$ is split into the $2\gamma$
operations $\LDendrA_a$, $\RDendrA_a$, $a \in [\gamma]$, so that for all
$x, y \in \Dca$,
\begin{equation}
    x \OpAsDendrA y =
    \sum_{a \in [\gamma]}  x \LDendrA_a y + x \RDendrA_a y,
\end{equation}
and all partial sums operations $\OpAsDendrA_b$, $b \in [\gamma]$,
satisfying
\begin{equation}
    x \OpAsDendrA_b y =
    \sum_{a \in [b]} x \LDendrA_a y + x \RDendrA_a x,
\end{equation}
also are associative.
\medbreak

On the other hand, in a $\gamma$-polydendriform algebra over the
operations $\LDendr_a$, $\RDendr_a$, $a \in [\gamma]$, by
Proposition~\ref{prop:associative_operation_dendr_gamma}, several
associative operations $\OpAsDendr_a$, $a \in [\gamma]$, are each split
into two operations $\LDendr_a$, $\RDendr_a$, $a \in [\gamma]$, so that
for all $x, y \in \Dca$,
\begin{equation}
    x \OpAsDendr_a y = x \LDendr_a y + x \RDendr_a y.
\end{equation}
\medbreak

Therefore, we can observe that $\gamma$-polydendriform algebras over the
operations $\LDendrA_a$, $\RDendrA_a$, $a \in [\gamma]$, are adapted to
study associative algebras (by splitting its single product in the way
we have described above) while $\gamma$-polydendriform algebras over the
operations $\LDendr_a$, $\RDendr_a$, $a \in [\gamma]$, are adapted to
study vectors spaces endowed with several associative products (by
splitting each one in the way we have described above). Algebras with
several associative products will be studied in
Section~\ref{sec:as_gamma}.
\medbreak

\subsubsection{Free polydendriform algebras}
From now on, in order to simplify and make the next definitions
uniform, we consider that in any $\gamma$-edge valued binary tree
$\Tfr$, all edges connecting internal nodes of $\Tfr$ with leaves are
labeled by $\infty$. By convention, for all $a \in [\gamma]$, we have
$a \Min \infty = a = \infty \Min a$.
\medbreak

Let us endow the vector space $\FreeAlg_{\Dendr_\gamma}$ of
$\gamma$-edge valued binary trees with linear operations
\begin{equation}
    \LDendr_a, \RDendr_a :
    \FreeAlg_{\Dendr_\gamma} \otimes \FreeAlg_{\Dendr_\gamma}
    \to \FreeAlg_{\Dendr_\gamma},
    \qquad
    a \in [\gamma],
\end{equation}
recursively defined, for any $\gamma$-edge valued binary tree $\Sfr$
and any $\gamma$-edge valued binary trees or leaves $\Tfr_1$ and
$\Tfr_2$ by
\begin{subequations}
\begin{equation}
    \Sfr \LDendr_a \LeafPic
    := \Sfr =:
    \LeafPic \RDendr_a \Sfr,
\end{equation}
\begin{equation}
    \LeafPic \LDendr_a \Sfr := 0 =: \Sfr \RDendr_a \LeafPic,
\end{equation}
\begin{equation}
    \ValuedBinTree{x}{y}{\Tfr_1}{\Tfr_2}
    \LDendr_a \Sfr :=
    \ValuedBinTree{x}{z}{\Tfr_1}{\Tfr_2 \LDendr_a \Sfr}
    +
    \ValuedBinTree{x}{z}{\Tfr_1}{\Tfr_2 \RDendr_y \Sfr}\,,
    \qquad
    z := a \Min y,
\end{equation}
\begin{equation}
    \Sfr \RDendr_a
    \ValuedBinTree{x}{y}{\Tfr_1}{\Tfr_2}
    :=
    \ValuedBinTree{z}{y}{\Sfr \RDendr_a \Tfr_1}{\Tfr_2}
    +
    \ValuedBinTree{z}{y}{\Sfr \LDendr_x \Tfr_1}{\Tfr_2}\,,
    \qquad
    z := a \Min x.
\end{equation}
\end{subequations}
Note that neither $\LeafPic \LDendr_a \LeafPic$ nor
$\LeafPic \RDendr_a \LeafPic$ are defined.
\medbreak

For example, we have
\begin{subequations}
\begin{multline}
\,.
\end{multline}
\end{subequations}
\medbreak

\begin{Theorem} \label{thm:free_dendr_gamma_algebra}
    For any integer $\gamma \geq 0$, the vector space
    $\FreeAlg_{\Dendr_\gamma}$ of all $\gamma$-edge valued binary trees
    endowed with the operations $\LDendr_a$, $\RDendr_a$,
    $a \in [\gamma]$, is the free $\gamma$-polydendriform algebra over
    one generator.
\end{Theorem}
\medbreak

\section{Multiassociative operads} \label{sec:as_gamma}
There is a well-known diagram, whose definition is recalled below,
gathering the diassociative, associative, and  dendriform operads.
The main goal of this section is to define a generalization on a
nonnegative integer parameter of the associative operad to obtain a new
version of this diagram, suited to the context of pluriassociative and
polydendriform operads.
\medbreak

\subsection{Two generalizations of the associative operad}
The associative operad is generated by one binary element. This operad
admits two different generalizations generated by $\gamma$ binary
elements with the particularity that one is the Koszul dual of the
other. In this section, we introduce and study these two operads.
\medbreak

\subsubsection{Multiassociative operads} \label{subsubsec:as_gamma}
For any integer $\gamma \geq 0$, we define the
\Def{$\gamma$-multiassociative operad} $\As_\gamma$ as the operad
admitting the presentation $\left(\GenAs, \RelAs\right)$, where
\begin{equation}
    \GeneratingSet_{\As_\gamma}
    := \GeneratingSet_{\As_\gamma}(2)
    := \{\MAs_a : a \in [\gamma]\}
\end{equation}
and $\RelationSpace_{\As_\gamma}$ is generated by
\begin{subequations}
\begin{equation} \label{equ:relation_as_gamma_1}
    \Corolla{\MAs_a} \circ_1 \Corolla{\MAs_b}
    -
    \Corolla{\MAs_b} \circ_2 \Corolla{\MAs_b},
    \qquad a \leq b \in [\gamma],
\end{equation}
\begin{equation} \label{equ:relation_as_gamma_2}
    \Corolla{\MAs_b} \circ_1 \Corolla{\MAs_a}
    -
    \Corolla{\MAs_b} \circ_2 \Corolla{\MAs_b},
    \qquad a < b \in [\gamma],
\end{equation}
\begin{equation} \label{equ:relation_as_gamma_3}
    \Corolla{\MAs_a} \circ_2 \Corolla{\MAs_b}
    -
    \Corolla{\MAs_b} \circ_2 \Corolla{\MAs_b},
    \qquad a < b \in [\gamma],
\end{equation}
\begin{equation} \label{equ:relation_as_gamma_4}
    \Corolla{\MAs_b} \circ_2 \Corolla{\MAs_a}
    -
    \Corolla{\MAs_b} \circ_2 \Corolla{\MAs_b},
    \qquad a < b \in [\gamma].
\end{equation}
\end{subequations}
This space of relations can be rephrased in a more compact way as the
space generated by
\begin{subequations}
\begin{equation} \label{equ:relation_as_gamma_1_concise}
    \Corolla{\MAs_a} \circ_1 \Corolla{\MAs_{a'}}
    -
    \Corolla{\MAs_{a \Max a'}} \circ_2 \Corolla{\MAs_{a \Max a'}},
    \qquad a, a' \in [\gamma],
\end{equation}
\begin{equation} \label{equ:relation_as_gamma_1_concise}
    \Corolla{\MAs_a} \circ_2 \Corolla{\MAs_{a'}}
    -
    \Corolla{\MAs_{a \Max a'}} \circ_2 \Corolla{\MAs_{a \Max a'}},
    \qquad a, a' \in [\gamma].
\end{equation}
\end{subequations}
\medbreak

It follows immediately that $\As_\gamma$ is well-defined as a
set-operad. Moreover, since $\As_1$ is isomorphic to the associative
operad $\As$ and $\As_\gamma$ is a suboperad of $\As_{\gamma + 1}$,
for all integers $\gamma \geq 0$, the operads  $\As_\gamma$ are
generalizations of the associative operad. Observe that the algebras
over $\As_\gamma$ are the $\gamma$-multiassociative algebras
introduced in Section~\ref{subsubsec:multiassociative_algebras}.
\medbreak

Let us now provide a realization of $\As_\gamma$. A
\Def{$\gamma$-corolla} is a rooted tree with at most one internal node
labeled on $[\gamma]$. Denote by $\FreeAlg_{\As_\gamma}$ the graded
vector space of all $\gamma$-corollas where the arity of a
$\gamma$-corolla is its arity, and let
\begin{equation}
    \MAs : \FreeAlg_{\As_\gamma} \otimes \FreeAlg_{\As_\gamma}
    \to \FreeAlg_{\As_\gamma}
\end{equation}
be the linear operation where, for any $\gamma$-corollas $\Cfr_1$ and
$\Cfr_2$, $\Cfr_1 \MAs \Cfr_2$ is the $\gamma$-corolla
with $n + m - 1$ leaves and labeled by $a \Max a'$ where $n$ (resp. $m$)
is the number of leaves of $\Cfr_1$ (resp. $\Cfr_2$) and $a$
(resp. $a'$) is the label of $\Cfr_1$ (resp. $\Cfr_2$).
\medbreak

\begin{Proposition} \label{prop:realization_koszulity_as_gamma}
    For any integer $\gamma \geq 0$, the operad $\As_\gamma$ is the
    vector space $\FreeAlg_{\As_\gamma}$ of $\gamma$-corollas and its
    partial compositions satisfy, for any $\gamma$-corollas $\Cfr_1$ and
    $\Cfr_2$, $\Cfr_1 \circ_i \Cfr_2 = \Cfr_1 \MAs \Cfr_2$ for all valid
    integer $i$. Besides, $\As_\gamma$ is a Koszul operad and the set of
    right comb $\GenAs$-syntax trees where all internal nodes have the
    same label forms a Poincaré-Birkhoff-Witt basis of~$\As_\gamma$.
\end{Proposition}
\medbreak

We have for instance in $\As_3$,
\begin{subequations}
\begin{equation}
\,.
\end{equation}
\end{subequations}
\medbreak

We deduce from Proposition~\ref{prop:realization_koszulity_as_gamma}
that the Hilbert series of $\As_\gamma$ satisfies
\begin{equation} \label{equ:hilbert_series_as_gamma}
    \Hca_{\As_\gamma}(t) = \frac{t + (\gamma - 1)t^2}{1 - t}.
\end{equation}
and that for all $n \geq 2$, $\dim \As_\gamma(n) = \gamma$.
\medbreak

\subsubsection{Dual multiassociative operads}
\label{subsubsec:das_gamma}
Since $\As_\gamma$ is a binary and quadratic operad, its admits a Koszul
dual, denoted by $\DAs_\gamma$ and called
\Def{$\gamma$-dual multiassociative operad}. The presentation of this
operad is provided by the next result.
\medbreak

\begin{Proposition} \label{prop:presentation_das_gamma}
    For any integer $\gamma \geq 0$, the operad $\DAs_\gamma$ admits the
    following presentation $\left(\GenDAs, \RelDAs\right)$ where
    \begin{math}
        \GenDAs := \GenDAs(2) := \{\MDAsA_a : a \in [\gamma]\}
    \end{math}
    and $\RelDAs$ is the space generated by
    \begin{multline} \label{equ:relation_das_gamma_alternative}
        \left(\sum_{a < b}
            \Corolla{\MDAsA_a} \circ_1 \Corolla{\MDAsA_b}
            +
            \Corolla{\MDAsA_b} \circ_1 \Corolla{\MDAsA_a}
            -
            \Corolla{\MDAsA_a} \circ_2 \Corolla{\MDAsA_b}
            -
            \Corolla{\MDAsA_b} \circ_2 \Corolla{\MDAsA_a}
        \right) \\
        +
        \Corolla{\MDAsA_b} \circ_1 \Corolla{\MDAsA_b}
        -
        \Corolla{\MDAsA_b} \circ_2 \Corolla{\MDAsA_b},
        \qquad b \in [\gamma].
    \end{multline}
\end{Proposition}
\medbreak

For any integer $\gamma \geq 0$, let $\MDAs_b$, $b \in [\gamma]$, the
elements of $\DAs$ defined by
\begin{equation} \label{equ:definition_operation_das_gamma}
    \MDAs_b := \sum_{a \in [b]} \MDAsA_a.
\end{equation}
Then, since for all $b \in [\gamma]$ we have
\begin{equation}
    \MDAsA_b =
    \begin{cases}
        \MDAs_1 & \mbox{if } b = 1, \\
        \MDAs_b - \MDAs_{b - 1} & \mbox{otherwise},
    \end{cases}
\end{equation}
by triangularity, the family
\begin{equation}
    \GenDAs' := \{\MDAs_b : b \in [\gamma]\}
\end{equation}
is a generating set of $\DAs_\gamma$. Let us now express a presentation
of $\DAs_\gamma$ through the family~$\GenDAs'$.
\medbreak

\begin{Proposition} \label{prop:presentation_das_gamma_alternative}
    For any integer $\gamma \geq 0$, the operad $\DAs_\gamma$ admits the
    presentation $\left(\GenDAs', \RelDAs'\right)$ where $\RelDAs'$ is
    the space generated by
    \begin{equation} \label{equ:relation_das_gamma}
        \Corolla{\MDAs_a} \circ_1 \Corolla{\MDAs_a}
        -
        \Corolla{\MDAs_a} \circ_2 \Corolla{\MDAs_a},
        \qquad a \in [\gamma].
    \end{equation}
\end{Proposition}
\medbreak

Observe, from the presentation provided by
Proposition~\ref{prop:presentation_das_gamma_alternative} of
$\DAs_\gamma$, that $\DAs_2$ is the operad denoted by $\TwoAs$
in~\cite{LR06}.
\medbreak

Notice that the presentation of the Koszul dual of $\DAs_\gamma$
computed from the presentation $\left(\GenDAs', \RelDAs'\right)$ of
Proposition~\ref{prop:presentation_das_gamma_alternative} gives rise to
the following presentation for $\As_\gamma$. This last operad
admits the presentation $\left(\GenAs', \RelAs'\right)$ where
\begin{equation}
    \GenAs' := \GenAs'(2) := \{\MAsA_a : a \in [\gamma]\}
\end{equation}
and $\RelAs'$ is the space generated by
\begin{subequations}
\begin{equation}
    \Corolla{\MAsA_a} \circ_1 \Corolla{\MAsA_{a'}},
    \qquad a \ne a' \in [\gamma],
\end{equation}
\begin{equation}
    \Corolla{\MAsA_a} \circ_2 \Corolla{\MAsA_{a'}},
    \qquad a \ne a' \in [\gamma],
\end{equation}
\begin{equation}
    \Corolla{\MAsA_a} \circ_1 \Corolla{\MAsA_a}
    -
    \Corolla{\MAsA_a} \circ_2 \Corolla{\MAsA_a},
    \qquad a \in [\gamma].
\end{equation}
\end{subequations}
Indeed, $\RelAs'$ is the space $\RelAs$ through the identification
\begin{equation}
    \MAsA_a =
    \begin{cases}
        \MAs_\gamma & \mbox{if } a = \gamma, \\
        \MAs_a - \MAs_{a + 1} & \mbox{otherwise}.
    \end{cases}
\end{equation}
\medbreak

\begin{Proposition} \label{prop:hilbert_series_das_gamma}
    For any integer $\gamma \geq 0$, the Hilbert series
    $\Hca_{\DAs_\gamma}(t)$ of the operad $\DAs_\gamma$ satisfies
    \begin{equation} \label{equ:serie_hilbert_das_gamma}
        t + (t - 1)\, \Hca_{\DAs_\gamma}(t) +
        (\gamma - 1) \, \Hca_{\DAs_\gamma}(t)^2
        = 0.
    \end{equation}
\end{Proposition}
\medbreak

\begin{figure}[ht]
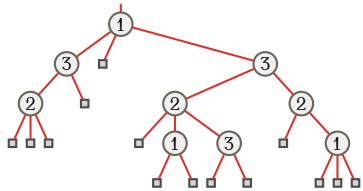

    \centering
    \begin{equation*}

    \end{equation*}
    \caption[A $3$-alternating Schröder tree of the ns
    operad~$\DAs_3(14)$.]
    {A $3$-alternating Schröder tree of size $14$. This tree is a basis
    element of $\DAs_3(14)$.}
    \label{fig:example_alternating_schroder_tree}
\end{figure}

By examining the expression for $\Hca_{\DAs_\gamma}(t)$ of the
statement of Proposition~\ref{prop:hilbert_series_das_gamma}, we
observe that for any $n \geq 1$, $\DAs_\gamma(n)$ can be seen as the
vector space $\FreeAlg_{\DAs_\gamma}(n)$ of all Schröder trees of arity
$n$, all labeled on $[\gamma]$ such that the label of an internal node
is different from the labels of its children that are internal nodes
(see Figure~\ref{fig:example_alternating_schroder_tree}).
We call these trees \Def{$\gamma$-alternating Schröder trees}. Let us
also denote by $\FreeAlg_{\DAs_\gamma}$ the graded vector space of all
$\gamma$-alternating Schröder trees.
\medbreak

We deduce also from Proposition~\ref{prop:hilbert_series_das_gamma} that
\begin{equation}
    \Hca_{\DAs_\gamma}(t) =
    \frac{1 - \sqrt{1 - (4\gamma - 2)t + t^2} - t}{2(\gamma - 1)}.
\end{equation}
\medbreak

\begin{Proposition} \label{prop:dimensions_DAs_gamma}
    For any integer $\gamma \geq 0$, the dimensions of the operad
    $\DAs_\gamma$ satisfy, for all $n \geq 2$,
    \begin{equation}
        \dim \DAs_\gamma(n) =
        \sum_{k = 0}^{n - 2}
        \gamma^{k + 1} (\gamma - 1)^{n - k - 2} \, \Narayana(n, k).
    \end{equation}
\end{Proposition}
\medbreak

In the statement of Proposition~\ref{prop:dimensions_DAs_gamma},
$\Narayana(n, k)$ is a Narayana number whose definition is recalled in
Section~\ref{equ:narayana_numbers} of Chapter~\ref{chap:combinatorics}.
For instance, the first dimensions of  $\DAs_1$, $\DAs_2$, $\DAs_3$, and
$\DAs_4$ are respectively
\begin{subequations}
\begin{equation}
    1, 1, 1, 1, 1, 1, 1, 1, 1, 1, 1,
\end{equation}
\begin{equation}
    1, 2, 6, 22, 90, 394, 1806, 8558, 41586, 206098, 1037718,
\end{equation}
\begin{equation}
    1, 3, 15, 93, 645, 4791, 37275, 299865, 2474025, 20819307,
    178003815,
\end{equation}
\begin{equation}
    1, 4, 28, 244, 2380, 24868, 272188, 3080596, 35758828, 423373636,
    5092965724.
\end{equation}
\end{subequations}
The second one is Sequence~\OEIS{A006318}, the third one is
Sequence~\OEIS{A103210}, and the last one is Sequence~\OEIS{A103211}
of~\cite{Slo}.
\medbreak

Let us now establish a realization of $\DAs_\gamma$.
\medbreak

\begin{Proposition} \label{prop:realization_das_gamma}
    For any nonnegative integer $\gamma$, the operad $\DAs_\gamma$ is
    the vector space $\FreeAlg_{\DAs_\gamma}$ of $\gamma$-alternating
    Schröder trees. Moreover, for any $\gamma$-alternating Schröder
    trees $\Sfr$ and $\Tfr$, $\Sfr \circ_i \Tfr$ is the
    $\gamma$-alternating Schröder tree obtained by grafting the root
    of $\Tfr$ on the $i$th leaf $x$ of $\Sfr$ and then, if the father
    $y$ of $x$ and the root $z$ of $\Tfr$ have a same label, by
    contracting the edge connecting $y$ and~$z$.
\end{Proposition}
\medbreak

We have for instance in $\DAs_3$,
\begin{subequations}
\begin{equation}
\,.
\end{equation}
\end{subequations}
\medbreak

\subsection{A diagram of operads}
\label{subsec:diagramme_dias_as_dendr_gamma}
We now define morphisms between the operads $\Dias_\gamma$,
$\As_\gamma$, $\DAs_\gamma$, and $\Dendr_\gamma$ to obtain a
generalization of a classical diagram involving the diassociative,
associative, and dendriform operads.
\medbreak

\subsubsection{Relating the diassociative and dendriform operads}
The diagram
\begin{equation} \label{equ:diagram_dendr_as_dias}
    \begin{tikzpicture}[yscale=.65,Centering]
        \node(Dendr)at(0,0){\begin{math} \Dendr \end{math}};
        \node(As)at(2,0){\begin{math} \As \end{math}};
        \node(Dias)at(4,0){\begin{math} \Dias \end{math}};
        \draw[->](Dias)edge node[anchor=south]
            {\begin{math} \eta \end{math}}(As);
        \draw[->](As)edge node[anchor=south]
            {\begin{math} \zeta \end{math}}(Dendr);
        \draw[<->,dotted,loop above,looseness=13](As)
            edge node[anchor=south]{\begin{math} ! \end{math}}(As);
        \draw[<->,dotted,loop above,looseness=1.5](Dendr)
            edge node[anchor=south]{\begin{math} ! \end{math}}(Dias);
    \end{tikzpicture}
\end{equation}
is a well-known diagram of operads, being a part of the so-called
operadic butterfly~\cite{Lod01,Lod06} and summarizing in a nice
way the links between the dendriform, associative, and diassociative
operads. The operad $\As$, being at the center of the diagram, is its
own Koszul dual, while $\Dias$ and $\Dendr$ are Koszul dual one of the
other.
\medbreak

The operad morphisms $\eta : \Dias \to \As$ and $\zeta : \As \to \Dendr$
are linearly defined through the realizations of $\Dias$ and $\Dendr$
recalled respectively in Sections~\ref{subsubsec:diassociative_operad}
and\ref{subsubsec:dendriform_operad} of Chapter~\ref{chap:algebra} by
\begin{equation}
    \eta(\Efr_{2, 1}) :=
\,.
\end{equation}
Since $\Dias$ is generated by $\Efr_{2, 1}$ and $\Efr_{2, 2}$, and since
$\As$ is generated by
\begin{math}
%
\end{math},
$\eta$ and $\zeta$ are wholly defined.
\medbreak

\subsubsection{Relating the pluriassociative and polydendriform operads}
\begin{Proposition} \label{prop:morphism_dias_gamma_to_as_gamma}
    For any integer $\gamma \geq 0$, the map
    $\eta_\gamma : \Dias_\gamma \to \As_\gamma$ satisfying
    \begin{equation}
        \eta_\gamma(0a) =
        %
        = \eta_\gamma(a0),
        \qquad a \in [\gamma],
    \end{equation}
    extends in a unique way into an operad morphism. Moreover, this
    morphism is surjective.
\end{Proposition}
\medbreak

By Proposition~\ref{prop:morphism_dias_gamma_to_as_gamma}, the map
$\eta_\gamma$, whose definition is only given in arity $2$, defines an
operad morphism. Nevertheless, by induction on the arity, one can prove
that for any word $x$ of $\Dias_\gamma$, $\eta_\gamma(x)$ is the
$\gamma$-corolla of arity $|x|$ labeled by the greatest letter of~$x$.
\medbreak

\begin{Proposition} \label{prop:morphism_das_gamma_to_dendr_gamma}
    For any integer $\gamma \geq 0$, the map
    $\zeta_\gamma : \DAs_\gamma \to \Dendr_\gamma$ satisfying
    \begin{equation}
        \zeta_\gamma\left(

        \right)
        = \BinTreeLTwo{a} + \BinTreeRTwo{a},
        \qquad a \in [\gamma],
    \end{equation}
    extends in a unique way into an operad morphism.
\end{Proposition}
\medbreak

We have to observe that the morphism $\zeta_\gamma$ defined in the
statement of Proposition~\ref{prop:morphism_das_gamma_to_dendr_gamma}
is injective only for $\gamma \leq 1$. Indeed, when $\gamma \geq 2$,
we have the relation
\begin{equation}
    \zeta_2\left(
    %
    \right)\,.
\end{equation}
\medbreak

\begin{Theorem} \label{thm:diagram_dias_as_das_dendr_gamma}
    For any integer $\gamma \geq 0$, the operads $\Dias_\gamma$,
    $\Dendr_\gamma$, $\As_\gamma$, and $\DAs_\gamma$ fit into the
    diagram
    \begin{equation} \label{equ:diagram_dias_as_das_dendr_gamma}
        %
\,,
    \end{equation}
    where $\eta_\gamma$ is the surjection defined in the statement
    of Proposition~\ref{prop:morphism_dias_gamma_to_as_gamma}
    and $\zeta_\gamma$ is the operad morphism defined in the statement
    of Proposition~\ref{prop:morphism_das_gamma_to_dendr_gamma}.
\end{Theorem}
\medbreak

Diagram~\eqref{equ:diagram_dias_as_das_dendr_gamma} is a
generalization of \eqref{equ:diagram_dendr_as_dias} in which the
associative operad splits into operads $\As_\gamma$ and~$\DAs_\gamma$.
\medbreak

\section{Further generalizations}
\label{sec:further_generalizations_gamma}
In this last section of this chapter, we propose some generalizations
on a nonnegative integer parameter of well-known operads. For this, we
use similar tools as the ones used in the first sections of the chapter.
\medbreak

\subsection{Duplicial operad}
We construct here a generalization on a nonnegative integer parameter of
the duplicial operad and describe the free algebras over one generator
in the category encoded by this generalization.
\medbreak

\subsubsection{Multiplicial operads} \label{subsub:dup_gamma}
It is well-known~\cite{LV12} that the dendriform operad and the
duplicial operad $\Dup$~\cite{Lod08} are both specializations of a
same operad $\DupDendr_q$ with one parameter $q \in \K$. This operad
admits the presentation
\begin{math}
    \left(
        \GeneratingSet_{\DupDendr_q}, \RelationSpace_{\DupDendr_q}
    \right)
\end{math}, where
$\GeneratingSet_{\DupDendr_q} := \GeneratingSet_\Dendr$ and
$\RelationSpace_{\DupDendr_q}$ is the space generated by
\begin{subequations}
\begin{equation}
    \Corolla{\LDendr} \circ_1 \Corolla{\RDendr}
    -
    \Corolla{\RDendr} \circ_2 \Corolla{\LDendr},
\end{equation}
\begin{equation}
    \Corolla{\LDendr} \circ_1 \Corolla{\LDendr}
    -
    \Corolla{\LDendr} \circ_2 \Corolla{\LDendr}
    - q
    \Corolla{\LDendr} \circ_2 \Corolla{\RDendr},
\end{equation}
\begin{equation}
    q \Corolla{\RDendr} \circ_1 \Corolla{\LDendr}
    +
    \Corolla{\RDendr} \circ_1 \Corolla{\RDendr}
    -
    \Corolla{\RDendr} \circ_2 \Corolla{\RDendr}.
\end{equation}
\end{subequations}
One can observe that $\DupDendr_1$ is the dendriform operad and that
$\DupDendr_0$ is the duplicial operad.
\medbreak

On the basis of this observation, from the presentation of
$\Dendr_\gamma$ provided by
Theorem~\ref{thm:autre_presentation_dendr_gamma} and its concise form
provided by Relations~\eqref{equ:relation_dendr_gamma_1_concise},
\eqref{equ:relation_dendr_gamma_2_concise},
and~\eqref{equ:relation_dendr_gamma_3_concise} for its space of
relations, we define the operad $\DupDendr_{q, \gamma}$ with two
parameters, an integer $\gamma \geq 0$ and $q \in \K$, in the
following way. We set $\DupDendr_{q, \gamma}$ as the operad admitting
the presentation
\begin{math}
    \left(
    \GeneratingSet_{\DupDendr_{q, \gamma}},
    \RelationSpace_{\DupDendr_{q, \gamma}}
    \right)
\end{math},
where $\GeneratingSet_{\DupDendr_{q, \gamma}} := \GenDendr'$ and
$\RelationSpace_{\DupDendr_{q, \gamma}}$ is the space generated by
\begin{subequations}
\begin{equation} \label{equ:dupdendr_gamma_1}
    \Corolla{\LDendr_a} \circ_1 \Corolla{\RDendr_{a'}}
    -
    \Corolla{\RDendr_{a'}} \circ_2 \Corolla{\LDendr_a},
    \qquad a, a' \in [\gamma],
\end{equation}
\begin{equation} \label{equ:dupdendr_gamma_2}
    \Corolla{\LDendr_a} \circ_1 \Corolla{\LDendr_{a'}}
    -
    \Corolla{\LDendr_{a \Min a'}} \circ_2 \Corolla{\LDendr_a}
    - q
    \Corolla{\LDendr_{a \Min a'}} \circ_2 \Corolla{\RDendr_{a'}},
    \qquad a, a' \in [\gamma],
\end{equation}
\begin{equation} \label{equ:dupdendr_gamma_3}
    q \Corolla{\RDendr_{a \Min a'}} \circ_1 \Corolla{\LDendr_{a'}}
    +
    \Corolla{\RDendr_{a \Min a'}} \circ_1 \Corolla{\RDendr_a}
    -
    \Corolla{\RDendr_a} \circ_2 \Corolla{\RDendr_{a'}},
    \qquad a, a' \in [\gamma].
\end{equation}
\end{subequations}
One can observe that $\DupDendr_{1, \gamma}$ is the
operad~$\Dendr_\gamma$.
\medbreak

Let us define the operad the \Def{$\gamma$-multiplicial operad}
$\Dup_\gamma$ as the operad $\DupDendr_{0, \gamma}$.
By using respectively the symbols $\LDup_a$ and $\RDup_a$ instead of
$\LDendr_a$ and $\RDendr_a$ for all $a \in [\gamma]$, we obtain that the
space of relations $\RelDup$ of $\Dup_\gamma$ is generated by
\begin{subequations}
\begin{equation} \label{equ:relation_dup_gamma_1}
    \LDup_a \circ_1 \RDup_{a'} - \RDup_{a'} \circ_2 \LDup_a,
    \qquad a, a' \in [\gamma],
\end{equation}
\begin{equation} \label{equ:relation_dup_gamma_2}
    \LDup_a \circ_1 \LDup_{a'} - \LDup_{a \Min a'} \circ_2 \LDup_a,
    \qquad a, a' \in [\gamma],
\end{equation}
\begin{equation} \label{equ:relation_dup_gamma_3}
    \RDup_{a \Min a'} \circ_1 \RDup_a - \RDup_a \circ_2 \RDup_{a'},
    \qquad a, a' \in [\gamma].
\end{equation}
\end{subequations}
We denote by $\GenDup$ the generating set
$\{\LDup_a, \RDup_a : a \in [\gamma] \}$ of $\Dup_\gamma$.
\medbreak

\begin{Proposition} \label{prop:properties_dup_gamma}
    For any integer $\gamma \geq 0$, the operad $\Dup_\gamma$ is Koszul
    and for any integer $n \geq 1$, $\Dup_\gamma(n)$ is the vector space
    of $\gamma$-edge valued binary trees with $n$ internal nodes.
\end{Proposition}
\medbreak

Since Proposition~\ref{prop:properties_dup_gamma} shows that the operads
$\Dup_\gamma$ and $\Dendr_\gamma$ have the same underlying vector space,
asking if these two operads are isomorphic is natural. The next result
implies that this is not the case.
\medbreak

\begin{Proposition} \label{prop:associative_operations_dup_gamma}
    For any integer $\gamma \geq 0$, any associative element of
    $\Dup_\gamma$ is proportional to $\LDup_a$ or to $\RDup_a$ for
    an~$a \in [\gamma]$.
\end{Proposition}
\medbreak

By Proposition~\ref{prop:associative_operations_dup_gamma} there are
exactly $2\gamma$ nonproportional associative operations in
$\Dup_\gamma$ while, by
Proposition~\ref{prop:associative_operation_dendr_gamma}
there are exactly $\gamma$ such operations in $\Dendr_\gamma$.
Therefore, $\Dup_\gamma$ and $\Dendr_\gamma$ are not isomorphic.
\medbreak

\subsubsection{Free multiplicial algebras}
A \Def{$\gamma$-multiplicial algebra} is a $\Dup_\gamma$-algebra. From
the definition of $\Dup_\gamma$, any $\gamma$-multiplicial algebra is a
vector space endowed with linear operations $\LDup_a, \RDup_a$,
$a \in [\gamma]$, satisfying the relations encoded
by~\eqref{equ:relation_dup_gamma_1}---\eqref{equ:relation_dup_gamma_3}.
\medbreak

In order the simplify and make uniform next definitions, we consider
that in any $\gamma$-edge valued binary tree $\Tfr$, all edges
connecting internal nodes of $\Tfr$ with leaves are labeled by $\infty$.
By convention, for all $a \in [\gamma]$, we have
$a \Min \infty = a = \infty \Min a$.
\medbreak

Let us endow the vector space $\FreeAlg_{\Dup_\gamma}$ of
$\gamma$-edge valued binary trees with linear operations
\begin{equation}
    \LDup_a, \RDup_a :
    \FreeAlg_{\Dup_\gamma} \otimes \FreeAlg_{\Dup_\gamma}
    \to \FreeAlg_{\Dup_\gamma},
    \qquad a \in [\gamma],
\end{equation}
recursively defined, for any $\gamma$-edge valued binary tree $\Sfr$
and any $\gamma$-edge valued binary trees or leaves $\Tfr_1$ and
$\Tfr_2$ by
\begin{subequations}
\begin{equation}
    \Sfr \LDup_a \LeafPic
    := \Sfr =:
    \LeafPic \RDup_a \Sfr,
\end{equation}
\begin{equation}
    \LeafPic \LDup_a \Sfr := 0 =: \Sfr \RDup_a \LeafPic,
\end{equation}
\begin{equation}
    \ValuedBinTree{x}{y}{\Tfr_1}{\Tfr_2}
    \LDup_a \Sfr :=
    \ValuedBinTree{x}{z}{\Tfr_1}{\Tfr_2 \LDup_a \Sfr}\,,
    \qquad
    z := a \Min y,
\end{equation}
\begin{equation}
    \Sfr \RDup_a
    \ValuedBinTree{x}{y}{\Tfr_1}{\Tfr_2}
    :=
    \ValuedBinTree{z}{y}{\Sfr \RDup_a \Tfr_1}{\Tfr_2}\,,
    \qquad
    z := a \Min x.
\end{equation}
\end{subequations}
Note that neither $\LeafPic \LDup_a \LeafPic$ nor
$\LeafPic \RDup_a \LeafPic$ are defined.
\medbreak

These recursive definitions for the operations $\LDup_a$, $\RDup_a$,
$a \in [\gamma]$, lead to the following direct reformulations. If $\Sfr$
and $\Tfr$ are two $\gamma$-edge valued binary trees,
$\Tfr \LDup_a \Sfr$ (resp. $\Sfr \RDup_a \Tfr$) is obtained by
replacing each label $y$ (resp. $x$) of any edge in the rightmost (resp.
leftmost) path of $\Tfr$ by $a \Min y$ (resp. $a \Min x$) to obtain a
tree $\Tfr'$, and by grafting the root of $\Sfr$ on the rightmost (resp.
leftmost) leaf of $\Tfr'$. These two operations are respective
generalizations of the operations under and over on binary trees
introduced by Loday and Ronco~\cite{LR02}.
\medbreak

For example, we have
\begin{subequations}
\begin{equation}
\,.
\end{equation}
\end{subequations}
\medbreak

\begin{Theorem} \label{thm:free_dup_gamma_algebra}
    For any integer $\gamma \geq 0$, the vector space
    $\FreeAlg_{\Dup_\gamma}$ of all $\gamma$-edge valued binary trees
    endowed with the operations $\LDup_a$, $\RDup_a$, $a \in [\gamma]$,
    is the free $\gamma$-multiplicial algebra over one generator.
\end{Theorem}
\medbreak

\subsection{Polytridendriform operads} \label{subsec:polytdendr}
We propose here a generalization $\TDendr_\gamma$ on a nonnegative
integer parameter $\gamma$ of the tridendriform operad~\cite{LR04}.
This last operad is the Koszul dual of the triassociative operad. We
proceed by using an analogous strategy as the one used to define the
operads $\Dendr_\gamma$ as Koszul duals of $\Dias_\gamma$. Indeed, we
define $\TDendr_\gamma$ as the Koszul dual of the operad
$\Trias_\gamma$, called $\gamma$-pluritriassociative operad, a
generalization of the triassociative operad defined in
Section~\ref{sec:pluritriass}.
\medbreak

Theorem~\ref{thm:presentation_trias_gamma}, by exhibiting
a presentation of $\Trias_\gamma$, shows that this operad is binary and
quadratic. It then admits a Koszul dual, denoted by $\TDendr_\gamma$ and
called \Def{$\gamma$-polytridendriform operad}.
\medbreak

\begin{Theorem} \label{thm:presentation_tdendr_gamma}
    For any integer $\gamma \geq 0$, the operad $\TDendr_\gamma$ admits
    the presentation $\left(\GenTDendr, \RelTDendr\right)$ where
    \begin{math}
        \GenTDendr := \GenTDendr(2) :=
        \left\{\LDendrA_a, \MTDendr, \RDendrA_a : a \in [\gamma]\right\}
    \end{math}
    and $\RelTDendr$ is the space generated by
    \begin{subequations}
    \begin{equation}\label{equ:relation_presentation_tdendr_gamma_1}
        \Corolla{\MTDendr} \circ_1 \Corolla{\MTDendr}
        -
        \Corolla{\MTDendr} \circ_2 \Corolla{\MTDendr},
    \end{equation}
    \begin{equation}\label{equ:relation_presentation_tdendr_gamma_2}
        \Corolla{\LDendrA_a} \circ_1 \Corolla{\MTDendr}
        -
        \Corolla{\MTDendr} \circ_2 \Corolla{\LDendrA_a},
        \qquad a \in [\gamma],
    \end{equation}
    \begin{equation}\label{equ:relation_presentation_tdendr_gamma_3}
        \Corolla{\MTDendr} \circ_1 \Corolla{\RDendrA_a}
        -
        \Corolla{\RDendrA_a} \circ_2 \Corolla{\MTDendr},
        \qquad a \in [\gamma],
    \end{equation}
    \begin{equation}\label{equ:relation_presentation_tdendr_gamma_4}
        \Corolla{\MTDendr} \circ_1 \Corolla{\LDendrA_a}
        -
        \Corolla{\MTDendr} \circ_2 \Corolla{\RDendrA_a},
        \qquad a \in [\gamma],
    \end{equation}
    \begin{equation}\label{equ:relation_presentation_tdendr_gamma_5}
        \Corolla{\LDendrA_a} \circ_1 \Corolla{\RDendrA_{a'}}
        -
        \Corolla{\RDendrA_{a'}} \circ_2 \Corolla{\LDendrA_a},
        \qquad a, a' \in [\gamma],
    \end{equation}
    \begin{equation}\label{equ:relation_presentation_tdendr_gamma_6}
        \Corolla{\LDendrA_a} \circ_1 \Corolla{\LDendrA_b}
        -
        \Corolla{\LDendrA_a} \circ_2 \Corolla{\RDendrA_b},
        \qquad a < b \in [\gamma],
    \end{equation}
    \begin{equation}\label{equ:relation_presentation_tdendr_gamma_7}
        \Corolla{\RDendrA_a} \circ_1 \Corolla{\LDendrA_b}
        -
        \Corolla{\RDendrA_a} \circ_2 \Corolla{\RDendrA_b},
        \qquad a < b \in [\gamma],
    \end{equation}
    \begin{equation}\label{equ:relation_presentation_tdendr_gamma_8}
        \Corolla{\LDendrA_b} \circ_1 \Corolla{\LDendrA_a}
        -
        \Corolla{\LDendrA_a} \circ_2 \Corolla{\LDendrA_b},
        \qquad a < b \in [\gamma],
    \end{equation}
    \begin{equation}\label{equ:relation_presentation_tdendr_gamma_9}
        \Corolla{\RDendrA_a} \circ_1 \Corolla{\RDendrA_b}
        -
        \Corolla{\RDendrA_b} \circ_2 \Corolla{\RDendrA_a},
        \qquad a < b \in [\gamma],
    \end{equation}
    \begin{equation}\label{equ:relation_presentation_tdendr_gamma_10}
        \Corolla{\LDendrA_d} \circ_1 \Corolla{\LDendrA_d}
        -
        \Corolla{\LDendrA_d} \circ_2 \Corolla{\MTDendr}
        -
        \left(\sum_{c \in [d]}
            \Corolla{\LDendrA_d} \circ_2 \Corolla{\LDendrA_c}
            +
            \Corolla{\LDendrA_d} \circ_2 \Corolla{\RDendrA_c}
        \right),
        \qquad d \in [\gamma],
    \end{equation}
    \begin{equation}\label{equ:relation_presentation_tdendr_gamma_11}
        \left(\sum_{c \in [d]}
            \Corolla{\RDendrA_d} \circ_1 \Corolla{\LDendrA_c}
            +
            \Corolla{\RDendrA_d} \circ_1 \Corolla{\RDendrA_c}
        \right)
        +
        \Corolla{\RDendrA_d} \circ_1 \Corolla{\MTDendr}
        -
        \Corolla{\RDendrA_d} \circ_2 \Corolla{\RDendrA_d},
        \qquad d \in [\gamma].
    \end{equation}
    \end{subequations}
\end{Theorem}
\medbreak

\begin{Proposition} \label{prop:hilbert_series_tdendr_gamma}
    For any integer $\gamma \geq 0$, the Hilbert series
    $\Hca_{\TDendr_\gamma}(t)$ of the operad $\TDendr_\gamma$ satisfies
    \begin{equation}
        t + ((2\gamma + 1) t - 1) \, \Hca_{\TDendr_\gamma}(t) +
        \gamma (\gamma + 1) t \, \Hca_{\TDendr_\gamma}(t)^2
        = 0.
    \end{equation}
\end{Proposition}
\medbreak

\begin{figure}[ht]
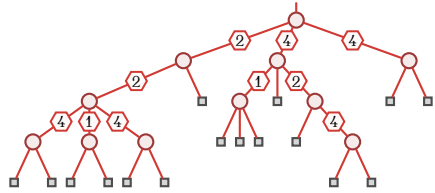

    \centering
    \begin{equation*}

    \end{equation*}
    \caption[A $4$-edge valued Schröder of the ns
    operad~$\TDendr_4(16)$.]
    {A $4$-edge valued Schröder tree of arity $16$. This tree is
    a basis element of $\TDendr_4(16)$.}
    \label{fig:example_edge_valued_schroder_tree}
\end{figure}

By examining the expression for $\Hca_{\TDendr_\gamma}(t)$ of the
statement of Proposition~\ref{prop:hilbert_series_tdendr_gamma}, we
observe that for any $n \geq 1$, $\TDendr(n)$ can be seen as the vector
space $\FreeAlg_{\TDendr_\gamma}(n)$ of Schröder trees with $n$
sectors wherein its edges connecting two internal nodes are labeled on
$[\gamma]$. We call these trees \Def{$\gamma$-edge valued Schröder
trees}. In our graphical representations of $\gamma$-edge valued
Schröder trees, any edge label is drawn into a hexagon located half
the edge (see Figure~\ref{fig:example_edge_valued_schroder_tree}).
\medbreak

We deduce from Proposition~\ref{prop:hilbert_series_tdendr_gamma} that
\begin{equation}
    \Hca_{\TDendr_\gamma}(t) =
    \frac{1 - \sqrt{1 - (4\gamma + 2)t + t^2} - (2\gamma + 1)t}
    {2(\gamma + \gamma^2)t}.
\end{equation}
\medbreak

\begin{Proposition} \label{prop:dimensions_TDendr_gamma}
    For any integer $\gamma \geq 0$, the dimensions of the operad
    $\TDendr_\gamma$ satisfy, for all $n \geq 1$,
    \begin{equation}
        \dim \TDendr_\gamma(n) =
        \sum_{k = 0}^{n - 1} (\gamma + 1)^k \gamma^{n - k - 1} \,
        \Narayana(n + 1, k).
    \end{equation}
\end{Proposition}
\medbreak

For instance, the first dimensions of $\TDendr_1$, $\TDendr_2$,
$\TDendr_3$, and $\TDendr_4$ are respectively
\begin{subequations}
\begin{equation}
    1, 3, 11, 45, 197, 903, 4279, 20793, 103049, 518859, 2646723,
\end{equation}
\begin{equation}
    1, 5, 31, 215, 1597, 12425, 99955, 824675, 6939769, 59334605,
    513972967,
\end{equation}
\begin{equation}
    1, 7, 61, 595, 6217, 68047, 770149, 8939707, 105843409,
    1273241431, 15517824973,
\end{equation}
\begin{equation}
    1, 9, 101, 1269, 17081, 240849, 3511741, 52515549, 801029681,
    12414177369, 194922521301.
\end{equation}
\end{subequations}
These sequences are respectively Sequences~\OEIS{A001003},
\OEIS{A269730}, \OEIS{A269731}, and \OEIS{A269732} of~\cite{Slo}.
\medbreak

\subsection{Operads of the operadic butterfly}
In what follows, we shall work with algebraic structures satisfying
relations involving possibly permutations of some inputs. For
simplicity, instead of working with symmetric operads, we shall just
work with types of algebras (see
Section~\ref{subsubsec:symmetric_operads} of
Chapter~\ref{chap:algebra}).
\medbreak

\subsubsection{A generalization of the operadic butterfly}
Let us consider the diagram of symmetric operads
\begin{equation} \label{equ:butterfly_diagram_gamma}
    \begin{tikzpicture}[xscale=.7,yscale=.65,Centering]
        \node(Dendr)at(-4,2){\begin{math}\Dendr_\gamma\end{math}};
        \node(As)at(1.5,0){\begin{math}\As_\gamma\end{math}};
        \node(DAs)at(-1.5,0){\begin{math}\DAs_\gamma\end{math}};
        \node(Dias)at(4,2){\begin{math}\Dias_\gamma\end{math}};
        \node(Com)at(-4,-2){\begin{math}\Com_\gamma\end{math}};
        \node(Lie)at(4,-2){\begin{math}\Lie_\gamma\end{math}};
        \node(Zin)at(-6,0){\begin{math}\Zin_\gamma\end{math}};
        \node(Leib)at(6,0){\begin{math}\Leib_\gamma\end{math}};
        \draw[->](Dias)--(As);
        \draw[->](DAs)--(Dendr);
        \draw[->](Dendr)--(Zin);
        \draw[->](Leib)--(Dias);
        \draw[->](Com)--(Zin);
        \draw[->](DAs)--(Com);
        \draw[->](Lie)--(As);
        \draw[->](Leib)--(Lie);
        \draw[<->,dotted](As)
            edge node[anchor=south]{\begin{math} ! \end{math}}(DAs);
        \draw[<->,dotted](Dendr)
            edge node[anchor=south]{\begin{math} ! \end{math}}(Dias);
        \draw[<->,dotted](Com)
            edge node[anchor=south]{\begin{math} ! \end{math}}(Lie);
        \draw[<->,dotted,loop below,looseness=.2](Zin)
            edge node[anchor=south]{\begin{math} ! \end{math}}(Leib);
    \end{tikzpicture}
\end{equation}
where $\DAs_\gamma$ is the $\gamma$-dual multiassociative operad defined
in Section~\ref{subsubsec:das_gamma} and $\Com_\gamma$, $\Lie_\gamma$,
$\Zin_\gamma$, and $\Leib_\gamma$, respectively are generalizations on a
nonnegative integer parameter $\gamma$ of the operads $\Com$, $\Lie$,
$\Zin$, and $\Leib$. Let us now define these operads.
\medbreak

\subsubsection{Commutative and Lie operads} \label{subsubsec:com_gamma}
Observe that the commutative operad $\Com$ is a commutative version
of~$\As = \DAs_1$ (see Section~\ref{subsubsec:symmetric_operads} of
Chapter~\ref{chap:algebra}). We define the symmetric operad
$\Com_\gamma$ by using the same idea of being a commutative version of
$\DAs_\gamma$. Therefore, $\Com_\gamma$ is the symmetric operad
describing the category of algebras $\Cca$ with binary operations
$\MDAs_a$, $a \in [\gamma]$, subjected for any elements $x$, $y$, and
$z$ of $\Cca$ to the two sorts of relations
\begin{subequations}
\begin{equation} \label{equ:relation_com_gamma_1}
    x \MDAs_a y = y \MDAs_a x,
    \qquad a \in [\gamma],
\end{equation}
\begin{equation} \label{equ:relation_com_gamma_2}
    (x \MDAs_a y) \MDAs_a z = x \MDAs_a (y \MDAs_a z),
    \qquad a \in [\gamma].
\end{equation}
\end{subequations}
Moreover, we define the symmetric operad $\Lie_\gamma$ as the Koszul
dual of~$\Com_\gamma$.
\medbreak

\subsubsection{Zinbiel and Leibniz operads} \label{subsubsec:zin_gamma}
It is well-known that the Zinbiel operad $\Zin$~\cite{Lod95} is a
commutative version of~$\Dendr = \Dendr_1$~\cite{Lod01}. We define the
symmetric operad $\Zin_\gamma$ by using the same idea of having the
property to be a commutative version of $\Dendr_\gamma$. Therefore,
$\Zin_\gamma$ is the symmetric operad describing the category of
algebras $\Zca$ with binary operations $\ProdZin_a$, $a \in [\gamma]$,
subjected for any elements $x$, $y$, and $z$ of $\Zca$ to the relation
\begin{equation} \label{equ:relation_zinbiel_gamma}
    (x \ProdZin_{a'} y) \ProdZin_a z =
    x \ProdZin_{a \Min a'} (y \ProdZin_a z) +
    x \ProdZin_{a \Min a'} (z \ProdZin_{a'} y),
    \qquad a, a' \in [\gamma].
\end{equation}
Relation~\eqref{equ:relation_zinbiel_gamma} is obtained
from Relations~\eqref{equ:relation_dendr_gamma_1_concise},
\eqref{equ:relation_dendr_gamma_2_concise},
and~\eqref{equ:relation_dendr_gamma_3_concise} of
$\gamma$-polydendri\-form algebras with the condition that for any
elements $x$ and $y$ and $a \in [\gamma]$,
$x \LDendr_a y = y \RDendr_a x$, and by setting
$x \ProdZin_a y := x \LDendr_a y$. Moreover, we define the symmetric
operad $\Leib_\gamma$ as the Koszul dual of $\Zin_\gamma$.
\medbreak

\begin{Proposition} \label{prop:morphism_com_zin_gamma}
    For any integer $\gamma \geq 0$ and any $\Zin_\gamma$-algebra
    $\Zca$, the binary operations $\MDAs_a$, $a \in [\gamma]$, defined
    for all elements $x$ and $y$ of $\Zca$ by
    \begin{equation}
        x \MDAs_a y := x \ProdZin_a y + y \ProdZin_a x,
        \qquad a \in [\gamma],
    \end{equation}
    endow $\Zca$ with a $\Com_\gamma$-algebra structure.
\end{Proposition}
\begin{proof}
    Since for all $a \in [\gamma]$ and all elements $x$ and $y$ of
    $\Zca$, by~\eqref{equ:relation_zinbiel_gamma}, we have
    \begin{equation}
        x \MDAs_a y - y \MDAs_a x =
        x \ProdZin_a y + y \ProdZin_a x -
        y \ProdZin_a x - x \ProdZin_a y
        = 0,
    \end{equation}
    the operations $\MDAs_a$ satisfy
    Relation~\eqref{equ:relation_com_gamma_1} of
    $\Com_\gamma$-algebras. Moreover, since for all $a \in [\gamma]$ and
    all elements $x$, $y$, and $z$ of $\Zca$,
    by~\eqref{equ:relation_zinbiel_gamma}, we have
    \begin{equation}\begin{split}
        (x \MDAs_a y) \MDAs_a z
        & - x \MDAs_a (y \MDAs_a z) \\
        & =
        (x \ProdZin_a y + y \ProdZin_a x) \ProdZin_a z +
        z \ProdZin_a (x \ProdZin_a y + y \ProdZin_a x) \\
        & \qquad - x \ProdZin_a (y \ProdZin_a z + z \ProdZin_a y)
        - (y \ProdZin_a z + z \ProdZin_a y) \ProdZin_a x \\
        & =
        (x \ProdZin_a y) \ProdZin_a z +
        (y \ProdZin_a x) \ProdZin_a z +
        z \ProdZin_a (x \ProdZin_a y) +
        z \ProdZin_a (y \ProdZin_a x) \\
        & \qquad - x \ProdZin_a (y \ProdZin_a z) -
        x \ProdZin_a (z \ProdZin_a y) -
        (y \ProdZin_a z) \ProdZin_a x -
        (z \ProdZin_a y) \ProdZin_a x \\
        & =
        (y \ProdZin_a x) \ProdZin_a z -
        (y \ProdZin_a z) \ProdZin_a x \\
        & =
        y \ProdZin_a (x \ProdZin_a z) +
        y \ProdZin_a (z \ProdZin_a x) -
        y \ProdZin_a (z \ProdZin_a x) -
        y \ProdZin_a (x \ProdZin_a z) \\
        & = 0,
    \end{split}\end{equation}
    the operations $\MDAs_a$ satisfy
    Relation~\eqref{equ:relation_com_gamma_2} of $\Com_\gamma$-algebras.
    Hence, $\Zca$ is a $\Com_\gamma$-algebra.
\end{proof}
\medbreak

\begin{Proposition} \label{prop:morphism_dendr_zin_gamma}
    For any integer $\gamma \geq 0$, and any $\Zin_\gamma$-algebra
    $\Zca$, the binary operations $\LDendr_a$, $\RDendr_a$,
    $a \in [\gamma]$ defined for all elements $x$ and $y$ of $\Zca$ by
    \begin{equation}
        x \LDendr_a y := x \ProdZin_a y,
        \qquad a \in [\gamma],
    \end{equation}
    and
    \begin{equation}
        x \RDendr_a y := y \ProdZin_a x,
        \qquad a \in [\gamma],
    \end{equation}
    endow $\Zca$ with a $\gamma$-polydendriform algebra structure.
\end{Proposition}

The constructions stated by
Propositions~\ref{prop:morphism_com_zin_gamma}
and~\ref{prop:morphism_dendr_zin_gamma} producing from a
$\Zin_\gamma$-algebra respectively a $\Com_\gamma$-algebra and a
$\gamma$-polydendriform algebra are functors from the category of
$\Zin_\gamma$-algebras respectively to the category of
$\Com_\gamma$-algebras and the category of $\gamma$-polydendriform
algebras. These functors respectively translate into symmetric operad
morphisms from $\Com_\gamma$ to $\Zin_\gamma$ and from $\Dendr_\gamma$
to $\Zin_\gamma$. These morphisms are generalizations of known morphisms
between $\Com$, $\Dendr$, and $\Zin$ of the operadic butterfly
(see~\cite{Lod01,Lod06,Zin12}).
\medbreak

\section*{Concluding remarks}
In this chapter, we have defined a new generalization $\Dendr_\gamma$
of the dendriform operad and also several ones of related operads.
Among its most important features, $\Dendr_\gamma$ encodes the notion
of splitting an associative product in several pieces. Moreover, as
illustrated, the underlying combinatorics of this operad involves a new
kind of combinatorial objects, which are binary trees with labeled
edges. A natural question about these trees consists in investigating
whether some particular subfamilies of these form suboperads of
$\Dendr_\gamma$. Moreover, like the dendriform operad which admits
a realization in term of rational functions~\cite{Lod10} (see
also Section~\ref{subsubsec:rational_functions_operad} of
Chapter~\ref{chap:algebra}), we can ask whether $\Dendr_\gamma$ admits a
similar realization.
\smallbreak

Besides, a complete study of the operads $\Com_\gamma$, $\Lie_\gamma$,
$\Zin_\gamma$, and $\Leib_\gamma$ (like computing their presentations
and providing realizations), and suitable definitions for all the
morphisms intervening in our generalization of the operadic
butterfly~\eqref{equ:butterfly_diagram_gamma} is worth to interest for
future works.
\smallbreak

Finally, one of our generalizations of the associative operad, namely
the multiassociative operad $\As_\gamma$, admits a direct
generalization $\As(\Qca)$ wherein its presentation is parametrized by
a finite poset $\Qca$. These operads and their Koszul duals have nice
combinatorial properties and will be studied in
Chapter~\ref{chap:posets}.
\medbreak


\chapter{From posets to operads} \label{chap:posets}
The content of this chapter comes from~\cite{Gir16c}.
\medbreak

\section*{Introduction}
This chapter is devoted to enrich connections between operads and
combinatorics by establishing a new link between posets and operads by
means of a construction associating a operad $\As(\Qca)$, called
$\Qca$-associative operad, with any finite poset~$\Qca$. This
construction is a functor $\As$ from the category of finite posets to
the category of binary and quadratic operads. The will to generalize two
families of operads Koszul dual to each other, constructed in
Chapter~\ref{chap:polydendr}, is the first impetus of this work. The
operads of these families are the multiassociative operads $\As_\gamma$
and the dual multiassociative operads $\DAs_\gamma$ (see
Section~\ref{sec:as_gamma} of Chapter~\ref{chap:polydendr}). In this
present work, we retrieve $\As_\gamma$ by applying the construction
$\As$ to the total order on a set of $\gamma$ elements and we retrieve
$\DAs_\gamma$ by applying the construction $\As$ to the trivial order
on the same set. Note that different constructions of operads involving
posets~\cite{FFM16}, and not directly related constructions involving
posets and operads~\cite{MY91,Val07} have been considered in the
literature.
\medbreak

Let us describe some main properties of $\As$. First, each operad
obtained by our construction provides a generalization of the
associative operad since all its generating operations are associative.
Besides, many combinatorial properties of the starting poset $\Qca$ lead
to algebraic properties for $\As(\Qca)$ (see
Table~\ref{tab:properties_as_poset}).
\begin{table}[ht]
    \centering
    \begin{tabular}{c|c|c}
        Properties of the poset $\Qca$
            & Properties of the operad $\As(\Qca)$
            & Statement \\ \hline \hline
        None & Binary and quadratic
                & Definition,
                Section~\ref {subsubsec:definition_as_poset}
                \\ \hline
        Forest & Koszul & Theorem~\ref{thm:koszulity_as_poset} \\ \hline
        Thin forest & Closed under Koszul duality
                & Theorem~\ref{thm:isomorphism_thin_forests_dual_operad}
                \\ \hline
        Trivial & Basic set-operad basis
            & Proposition~\ref{prop:basic_operad_as_poset}
    \end{tabular}
    \bigbreak

    \caption[Properties of $\Qca$-associative operads.]
    {Summary of the properties satisfied by a poset $\Qca$ implying
    properties for the operad $\As(\Qca)$. Note that any trivial poset
    is also a thin forest poset, and that a thin forest poset is also a
    forest poset. In particular, if $\Qca$ is a trivial poset,
    $\As(\Qca)$ has all properties mentioned in the middle column.}
    \label{tab:properties_as_poset}
\end{table}
For instance, when $\Qca$ is a forest (with the meaning that no element
of $\Qca$ covers two different elements), $\As(\Qca)$ is a Koszul
operad. Moreover, when $\Qca$ is not a trivial poset, the fundamental
basis of $\As(\Qca)$ is not a basic set-operad basis. This last
property seems to be interesting since almost all set-operad bases of
common operads are basic, such as the associative operad or the
diassociative operad~\cite{Lod01} (see additionally~\cite{Zin12}). This
gives to our construction a very unique flavor.
\medbreak

The further study of the operads obtained by the construction $\As$ is
driven by computer exploration. Indeed, computer experiments bring us
the observation that some operads obtained by the construction $\As$ are
Koszul duals to each other. This observation raises several questions.
The first one consists in describing a family of posets, called thin
forest posets, such that the construction $\As$ restricted to this
family is closed under Koszul duality. The second one consists in
defining an operation $^\perp$ on this family of posets such that for
any of these posets $\Qca$, $\As(\Qca^\perp)$ is isomorphic to the
Koszul dual $\As(\Qca)^!$ of $\As(\Qca)$. The last one relies on an
expression of an explicit isomorphism between $\As(\Qca)^!$ and
$\As(\Qca^\perp)$. We answer all these questions in this work, forming
its main results. As additional results, we provide a complete study
of the operads $\As(\Qca)$, including, when $\Qca$ satisfies some
precise properties, an expression for its Hilbert series and a
realization involving labeled Schröder trees.
\medbreak

This chapter is organized as follows.
Section~\ref{sec:posets_to_operads} is concerned with the description
of the construction $\As$ and the first general properties of the
obtained operads. In Section~\ref{sec:forest_posets}, we focus on the
case where the poset $\Qca$ at input of the construction $\As$ is a
forest poset. We show that in this case, $\As(\Qca)$ is a Koszul operad
and derive some consequences. This chapter ends by introducing in
Section~\ref{sec:thin_forest_posets} the class of thin forest posets.
The construction $\As$ restricted to this class of posets has the
property to be closed under Koszul duality.
\medbreak

\subsubsection*{Note}
In this chapter all posets are finite. For this reason, ``poset'' means
``finite poset''. Moreover, since this chapter deals only with ns
operads, ``operad'' means ``ns operad''. If $\Op$ is a generator of an
operad $\Oca$, we denote by $\OpDual$ the associated generator in the
Koszul dual of~$\Oca$.
\medbreak

\section{From posets to operads} \label{sec:posets_to_operads}
This section is devoted to the introduction of our construction
producing an operad from a poset. We also establish here some of its
first general properties. We end this section by presenting algebras
over our operads and some of their properties.
\medbreak

\subsection{Construction}
Let us describe the construction $\As$, associating with any poset a
binary and quadratic operad presentation, and prove that it is
functorial.
\medbreak

\subsubsection{Operad presentations from posets}
\label{subsubsec:definition_as_poset}
For any poset $(\Qca, \Ord_\Qca)$, we define the
\Def{$\Qca$-associa\-tive operad} $\As(\Qca)$ as the operad admitting the
presentation
\begin{math}
    \left(\GeneratingSet_\Qca^\Op, \RelationSpace_\Qca^\Op\right)
\end{math}
where $\GeneratingSet_\Qca^\Op$ is the set of generators
\begin{equation}
    \GeneratingSet_\Qca^\Op :=
    \GeneratingSet_\Qca^\Op(2) :=
    \left\{\Op_a : a \in \Qca\right\},
\end{equation}
and $\RelationSpace_\Qca^\Op$ is the space of relations generated by
\begin{subequations}
\begin{equation} \label{equ:relation_as_poset_1}
    \Corolla{\Op_a} \circ_1 \Corolla{\Op_b}
    -
    \Corolla{\Op_{a \Min_\Qca b}} \circ_2 \Corolla{\Op_{a \Min_\Qca b}},
    \qquad
    a, b \in \Qca
    \mbox{ and }
    (a \Ord_\Qca b \mbox{ or } b \Ord_\Qca a),
\end{equation}
\begin{equation} \label{equ:relation_as_poset_2}
    \Corolla{\Op_{a \Min_\Qca b}} \circ_1 \Corolla{\Op_{a \Min_\Qca b}}
    -
    \Corolla{\Op_a} \circ_2 \Corolla{\Op_b},
    \qquad
    a, b \in \Qca
    \mbox{ and }
    (a \Ord_\Qca b \mbox{ or } b \Ord_\Qca a).
\end{equation}
\end{subequations}
By definition, $\As(\Qca)$ is a binary and quadratic operad. The
\Def{fundamental basis} of $\As(\Qca)$ is the basis induced by
$\GeneratingSet_\Qca^\Op$. Moreover, since
Relations~\eqref{equ:relation_as_poset_1}
and~\eqref{equ:relation_as_poset_2} are of the form $\Sfr - \Tfr$ where
$\Sfr$ and $\Tfr$ are $\GeneratingSet_\Qca^\Op$-syntax trees,
$\As(\Qca)$ is well-defined in the category of sets.
\medbreak

\begin{Lemma} \label{lem:relations_as_poset}
    Let $\Qca$ be a poset. For any $a \in \Qca$, let $R_a$ be the set
    \begin{multline}
        R_a := \left\{
        \Corolla{\Op_a} \circ_1 \Corolla{\Op_b}, \;
        \Corolla{\Op_b} \circ_1 \Corolla{\Op_a}, \right. \\
        \left. \Corolla{\Op_b} \circ_2 \Corolla{\Op_a}, \;
        \Corolla{\Op_a} \circ_2 \Corolla{\Op_b}
        : b \in \Qca \mbox{ and } a \Ord_\Qca b\right\}.
    \end{multline}
    Then, for all $\Sfr, \Tfr \in R_a$, $\Sfr - \Tfr$ is an element of
    the space of relations $\RelationSpace_\Qca^\Op$ of $\As(\Qca)$.
\end{Lemma}
\medbreak

By Lemma~\ref{lem:relations_as_poset}, we observe that all $\Op_a$,
$a \in \Qca$, are associative. For this reason, $\As(\Qca)$ is a
generalization of the associative operad on several binary generating
operations. As we will see in the sequel, this very simple way to
produce operads has many combinatorial and algebraic properties.
\medbreak

\subsubsection{Functoriality}
For any morphism of posets $\phi : \Qca_1 \to \Qca_2$, we denote by
$\As(\phi)$ the map
\begin{equation}
    \As(\phi) : \As(\Qca_1)(2) \to \As(\Qca_2)(2)
\end{equation}
defined by
\begin{equation}
    \As(\phi)\left(\Op_x\right) := \Op_{\phi(x)}
\end{equation}
for all $x \in \Qca_1$.
\medbreak

\begin{Lemma} \label{lem:maps_as_poset}
    Let $\Qca_1$ and $\Qca_2$ be two posets and
    $\phi : \Qca_1 \to \Qca_2$ be a morphism of posets. Then, the map
    $\As(\phi)$ uniquely extends into an operad morphism from
    $\As(\Qca_1)$ to $\As(\Qca_2)$.
\end{Lemma}

\begin{Theorem} \label{thm:as_poset_functor}
    The construction $\As$ is a functor from the category of posets to
    the category of binary and quadratic operads.
\end{Theorem}
\medbreak

\subsubsection{First examples} \label{subsubsec:first_examples_as_poset}
Let us use Theorem~\ref{thm:as_poset_functor} to exhibit some examples
of constructions of binary an quadratic operads from posets.
\medbreak

From the poset
\begin{equation}
    \Qca :=
    \begin{tikzpicture}[xscale=.4,yscale=.4,Centering]
        \node[PosetVertex](1)at(0,0){\begin{math}1\end{math}};
        \node[PosetVertex](2)at(2,0){\begin{math}2\end{math}};
        \node[PosetVertex](3)at(1,-1){\begin{math}3\end{math}};
        \node[PosetVertex](4)at(3,-1){\begin{math}4\end{math}};
        \draw[Edge](1)--(3);
        \draw[Edge](2)--(3);
        \draw[Edge](2)--(4);
    \end{tikzpicture}\,,
\end{equation}
the operad $\As(\Qca)$ is binary and quadratic, generated by the set
$\GeneratingSet_\Qca^\Op = \{\Op_1, \Op_2, \Op_3, \Op_4\}$, and,
by Lemma~\ref{lem:relations_as_poset}, subjected to the relations
\begin{subequations}
\begin{equation}
    \Op_1 \circ_1 \Op_1
    = \Op_1 \circ_1 \Op_3 = \Op_3 \circ_1 \Op_1
    = \Op_3 \circ_2 \Op_1 = \Op_1 \circ_2 \Op_3
    = \Op_1 \circ_2 \Op_1,
\end{equation}
\begin{equation}\begin{split}
    \Op_2 \circ_1 \Op_2
    = \Op_2 \circ_1 \Op_3 & = \Op_2 \circ_1 \Op_4
    = \Op_3 \circ_1 \Op_2 = \Op_4 \circ_1 \Op_2
    \\ & = \Op_4 \circ_2 \Op_2 = \Op_3 \circ_2 \Op_2
    = \Op_2 \circ_2 \Op_4 = \Op_2 \circ_2 \Op_3
    = \Op_2 \circ_2 \Op_2,
\end{split}\end{equation}
\begin{equation}
    \Op_3 \circ_1 \Op_3 = \Op_3 \circ_2 \Op_3,
\end{equation}
\begin{equation}
    \Op_4 \circ_1 \Op_4 = \Op_4 \circ_2 \Op_4.
\end{equation}
\end{subequations}
\medbreak

Besides, when $\Qca$ is the trivial poset on the set $[\ell]$,
$\ell \geq 0$, the operad $\As(\Qca)$ is generated by the set
$\GeneratingSet_\Qca^\Op = \{\Op_1, \dots, \Op_\ell\}$ and, by
Lemma~\ref{lem:relations_as_poset}, subjected to the relations
\begin{equation}
    \Op_a \circ_1 \Op_a = \Op_a \circ_2 \Op_a,
    \qquad a \in [\ell].
\end{equation}
In particular, when $\ell = 0$, $\As(\Qca)$ is the trivial operad, when
$\ell = 1$, $\As(\Qca)$ is the associative operad, and when $\ell = 2$,
$\As(\Qca)$ is the operad $\TwoAs$~\cite{LR06}. These operads, for a
generic $\ell \geq 0$, are the dual multiassociative operads
$\DAs_\ell$, introduced in Section~\ref{sec:as_gamma} of
Chapter~\ref{chap:polydendr}. These operads can be realized using
Schröder trees endowed with labels satisfying some conditions. In
Section~\ref{subsubsec:realization_as_poset_forest}, we shall describe
a generalized version of this realization.
\medbreak

Furthermore, when $\Qca$ is the total order on the set $[\ell]$,
$\ell \geq 0$, the operad $\As(\Qca)$ is generated by the set
$\GeneratingSet_\Qca^\Op = \{\Op_1, \dots, \Op_\ell\}$ and, by
Lemma~\ref{lem:relations_as_poset}, subjected to the relations
\begin{equation}
    \Op_a \circ_1 \Op_a =
    \Op_a \circ_1 \Op_b = \Op_b \circ_1 \Op_a =
    \Op_b \circ_2 \Op_a = \Op_a \circ_2 \Op_b =
    \Op_a \circ_2 \Op_a,
    \qquad a \leq b \in [\ell].
\end{equation}
In particular, when $\ell = 0$, $\As(\Qca)$ is the trivial operad and
when $\ell = 1$, $\As(\Qca)$ is the associative operad. These operads,
for a generic $\ell \geq 0$, are the multiassociative operads
$\As_\ell$, introduced in Section~\ref{sec:as_gamma} of
Chapter~\ref{chap:polydendr}. They have the particularity to have
stationary dimensions since $\dim \As(\Qca)(1) = 1$ and
$\dim \As(\Qca)(n) = \ell$ for all~$n \geq 2$.
\medbreak

\subsection{General properties}
Let us now list some general properties of the operad $\As(\Qca)$ where
$\Qca$ is a poset without particular requirements. We provide the
dimension of the space of relations of $\As(\Qca)$, describe its
associative elements, and give a necessary and sufficient condition for
the fact that its fundamental basis is a basic set-operad basis.
\medbreak

\subsubsection{Space of relations dimensions}

\begin{Proposition} \label{prop:dimension_relations_as_poset}
    Let $\Qca$ be a poset. Then, the dimension of the space
    $\RelationSpace_\Qca^\Op$ of relations of $\As(\Qca)$ satisfies
    \begin{equation} \label{equ:dimension_relations_as_poset}
        \dim \RelationSpace_\Qca^\Op =
        4 \, \NbInterv(\Qca) - 3 \, \# \Qca.
    \end{equation}
\end{Proposition}
\medbreak

Recall that $\NbInterv(\Qca)$ denotes the number of intervals of
$\Qca$ (see Section~\ref{subsubsec:definitions_posets}
of~Chapter~\ref{chap:combinatorics}).
\medbreak

\subsubsection{Associative elements}

\begin{Proposition} \label{prop:associative_elements_as_poset}
    Let $\Qca$ be a poset and
    $C := \{c_1 \OrdStrict_\Qca \dots \OrdStrict_\Qca c_\ell\}$ be a
    chain of $\Qca$. Then, $\K \Angle{C}$ contains only associative
    elements of $\As(\Qca)$. Conversely, any associative element of
    $\As(\Qca)$ is an element of $\K \Angle{C}$ for a chain $C$
    of~$\Qca$.
\end{Proposition}
\medbreak

\subsubsection{Basicity}
\begin{Proposition} \label{prop:basic_operad_as_poset}
    Let $\Qca$ be a poset. The fundamental basis of $\As(\Qca)$ is a
    basic set-operad basis if and only if $\Qca$ is a trivial poset.
\end{Proposition}
\medbreak

\subsection{Algebras over poset associative operads}
Let $\Qca$ be a poset. From the presentation
$\left(\GeneratingSet_\Qca^\Op, \RelationSpace_\Qca^\Op\right)$ of the
operad $\As(\Qca)$ provided by its definition in
Section~\ref{subsubsec:definition_as_poset}, an $\As(\Qca)$-algebra is
a vector space $\Alg$ endowed with linear operations
\begin{equation}
    \Op_a : \Alg \otimes \Alg \to \Alg,
    \qquad a \in \Qca,
\end{equation}
satisfying, for all $x, y, z \in \Alg$, the relations
\begin{equation} \label{equ:relations_algebras_as_poset}
    (x \Op_a y) \Op_b z
    = x \Op_a (y \Op_b z)
    = (x \Op_c y) \Op_a z
    = x \Op_c (y \Op_a z),
    \qquad a, b, c \in \Qca
    \mbox{ and } a \Ord_\Qca b \mbox{ and } a \Ord_\Qca c.
\end{equation}
We call \Def{$\Qca$-associative algebra} any $\As(\Qca)$-algebra.
\medbreak

We shall exhibit two examples of $\Qca$-associative algebras in the
sequel: in Section~\ref{subsubsec:free_as_forest_poset_algebras}, free
$\Qca$-associative algebras over one generator when $\Qca$ is a forest
poset and in Section~\ref{subsubsec:antichains_algebra},
$\Qca$-associative algebras involving the antichains of the
poset~$\Qca$.
\medbreak

\subsubsection{Units}
Let $\Qca$ be a poset and $\Alg$ be a $\Qca$-associative algebra. An
\Def{$a$-unit}, $a \in \Qca$, of $\Alg$ is an element $\Unit_a$ of
$\Alg$ satisfying
\begin{equation}
    \Unit_a \Op_a x = x = x \Op_a \Unit_a
\end{equation}
for all $x \in \Alg$. Obviously, for any $a \in \Qca$ there is at most
one $a$-unit in $\Alg$.
\medbreak

Besides, for any element $x$ of $\Alg$, we denote by $\UnitSet_\Alg(x)$
the set of elements $a$ of $\Qca$ such that $x$ is an $a$-unit of
$\Alg$. Obviously, if $\Unit_a$ is an $a$-unit of $\Alg$,
$a \in \UnitSet_\Alg(\Unit_a)$.
\medbreak

\begin{Proposition} \label{prop:units_as_gamma}
    Let $\Qca$ be a poset and $\Alg$ be a $\Qca$-associative algebra.
    Then:
    \begin{enumerate}[label={\it (\roman*)}]
        \item \label{item:units_as_gamma_1}
        for any element $x$ of $\Alg$, $\UnitSet_\Alg(x)$ is an order
        filter of $\Qca$;
        \item \label{item:units_as_gamma_2}
        for all elements $x$ and $y$ of $\Alg$ such that $x \ne y$, the
        sets $\UnitSet_\Alg(x)$ and $\UnitSet_\Alg(y)$ are disjoint.
    \end{enumerate}
\end{Proposition}
\medbreak

Proposition~\ref{prop:units_as_gamma} implies that the sets
$\UnitSet_\Alg(x)$, $x \in \Alg$, form a partition of an order filter
of~$\Qca$ where each part is itself an order filter of~$\Qca$.
\medbreak

\subsubsection{Antichains algebra} \label{subsubsec:antichains_algebra}
Let $\Qca$ be a poset and set $\Xbb_\Qca := \{x_a : a \in \Qca\}$ as
a set of commutative parameters and consider the commutative and
associative polynomial algebra $\K[\Xbb_\Qca]/_{\Ideal_\Qca}$, where
$\Ideal_\Qca$ is the ideal of $\K[\Xbb_\Qca]$ generated by
\begin{equation}
    x_a x_b - x_a, \qquad a \Ord_\Qca b \in \Qca.
\end{equation}
Then, one observes that $x_{a_1} \dots x_{a_k}$ is a reduced monomial of
$\K[\Xbb_\Qca]/_{\Ideal_\Qca}$ if and only if the set
$\{a_1, \dots, a_k\}$ is an antichain of $\Qca$ of size $k$.
\medbreak

We endow $\K[\Xbb_\Qca]/_{\Ideal_\Qca}$ with linear operations
\begin{equation}
    \Op_a :
    \K[\Xbb_\Qca]/_{\Ideal_\Qca} \otimes
    \K[\Xbb_\Qca]/_{\Ideal_\Qca} \to
    \K[\Xbb_\Qca]/_{\Ideal_\Qca},
    \qquad a \in \Qca,
\end{equation}
defined, for all reduced monomials $x_{b_1} \dots x_{b_k}$ and
$x_{c_1} \dots x_{c_\ell}$ of $\K[\Xbb]/_{\Ideal_\Qca}$, by
\begin{equation}
    x_{b_1} \dots x_{b_k} \Op_a x_{c_1} \dots x_{c_\ell} :=
    \pi(x_{b_1} \dots x_{b_k} x_a x_{c_1} \dots x_{c_\ell}),
\end{equation}
where $\pi : \K[\Xbb_\Qca] \to \K[\Xbb_\Qca]/_{\Ideal_\Qca}$ is the
canonical projection. These operations $\Op_a$, $a \in \Qca$, endow
$\K[\Xbb_\Qca]/_{\Ideal_\Qca}$ with a structure of a $\Qca$-associative
algebra.
\medbreak

Consider for instance the poset
\begin{equation}
    \Qca :=
    \begin{tikzpicture}[xscale=.45,yscale=.45,Centering]
        \node[PosetVertex](1)at(0,0){\begin{math}1\end{math}};
        \node[PosetVertex](2)at(-1,-1){\begin{math}2\end{math}};
        \node[PosetVertex](3)at(0,-1){\begin{math}3\end{math}};
        \node[PosetVertex](4)at(1,-1){\begin{math}4\end{math}};
        \node[PosetVertex](5)at(0,-2){\begin{math}5\end{math}};
        \draw[Edge](1)--(2);
        \draw[Edge](1)--(3);
        \draw[Edge](1)--(4);
        \draw[Edge](2)--(5);
        \draw[Edge](4)--(5);
    \end{tikzpicture}\,.
\end{equation}
The space $\K[\Xbb_\Qca]/_{\Ideal_\Qca}$ is the linear span of the
reduced monomials
\begin{equation}
    x_1, \; x_2, \; x_3, \; x_4, \; x_5, \;
    x_2 x_3, \; x_2 x_4, \; x_3 x_4, \; x_3 x_5, \;
    x_2 x_3 x_4,
\end{equation}
and one has for instance
\begin{subequations}
\begin{equation}
    x_2 \Op_3 x_4 = x_2 x_3 x_4,
\end{equation}
\begin{equation}
    x_2 x_3 \Op_1 x_4 = x_1,
\end{equation}
\begin{equation}
    x_2 x_3 \Op_5 x_4 = x_2 x_3 x_4.
\end{equation}
\end{subequations}
\medbreak

\section{Forest posets, Koszul duality, and Koszulity}
\label{sec:forest_posets}
Here, we focus on the construction $\As$ when the input poset $\Qca$
of the construction is a forest poset. In this case, we show that
$\As(\Qca)$ is Koszul, we provide a realization of $\As(\Qca)$, and
we obtain a functional equation for its Hilbert series. We end this
section by computing presentations of the Koszul dual of $\As(\Qca)$.
\medbreak

\subsection{Koszulity and Poincaré-Birkhoff-Witt bases}
We prove here that when $\Qca$ is a forest poset, $\As(\Qca)$ is Koszul.
For that, we consider an orientation of the space of relations
$\RelationSpace_\Qca^\Op$ of $\As(\Qca)$ and show that this orientation
is a convergent rewrite rule. As a consequence, the Koszulity of
$\As(\Qca)$ follows (see Lemma~\ref{lem:koszulity_criterion_pbw}
of Chapter~\ref{chap:algebra}).
\medbreak

\subsubsection{Forest posets}
We call \Def{forest poset} any poset avoiding the pattern
$\AntiForestPattern$ (see Section~\ref{subsubsec:poset_patterns} of
Chapter~\ref{chap:combinatorics} for the definition of pattern avoidance
in posets). In other words, a forest poset is a poset for which its
Hasse diagram is a forest of rooted trees (where roots are minimal
elements). Figure~\ref{fig:not_thin_forest_poset} shows an example of a
forest poset.
\medbreak

\subsubsection{Orientation of the space of relations}
Let $\Qca$ be a poset (not necessarily a forest poset just now) and
$\Rew_\Qca$ be the rewrite rule on $\GeneratingSet_\Qca^\Op$-syntax
trees satisfying
\begin{subequations}
\begin{equation} \label{equ:rewrite_as_poset_1}
    \Corolla{\Op_a} \circ_1 \Corolla{\Op_b}
    \Rew_\Qca
    \Corolla{\Op_{a \Min_\Qca b}} \circ_2 \Corolla{\Op_{a \Min_\Qca b}},
    \qquad a, b \in \Qca
    \mbox{ and }
    (a \Ord_\Qca b \mbox{ or } b \Ord_\Qca a),
\end{equation}
\begin{equation} \label{equ:rewrite_as_poset_2}
    \Corolla{\Op_a} \circ_2 \Corolla{\Op_b}
    \Rew_\Qca
    \Corolla{\Op_{a \Min_\Qca b}} \circ_2 \Corolla{\Op_{a \Min_\Qca b}},
    \qquad a, b \in \Qca
    \mbox{ and }
    (a \OrdStrict_\Qca b \mbox{ or } b \OrdStrict_\Qca a).
\end{equation}
\end{subequations}
Let also $\RewTrees_\Qca$ be the closure of $\Rew_\Qca$.
\medbreak

\subsubsection{Convergent rewrite rule}

\begin{Lemma} \label{lem:rewrite_rule_space_relations_as_poset}
    Let $\Qca$ be a poset. Then, $\Rew_\Qca$ is an orientation of the
    space of relations~$\RelationSpace_\Qca^\Op$ of~$\As(\Qca)$.
\end{Lemma}
\medbreak

\begin{Lemma} \label{lem:terminating_rewrite_rule}
    Let $\Qca$ be a poset. Then, $\RewTrees_\Qca$ is a terminating
    rewrite rule.
\end{Lemma}
\medbreak

\begin{Lemma} \label{lem:normal_forms_as_poset}
    Let $\Qca$ be a poset. Then, the set of the normal forms of
    $\RewTrees_\Qca$ is the set of the
    $\GeneratingSet_\Qca^\Op$-syntax trees $\Tfr$ such that for any
    internal node of $\Tfr$ labeled by $\Op_a$ having a left (resp.
    right) child labeled by $\Op_b$, $a$ and $b$ are incomparable
    (resp. are equal or are incomparable) in~$\Qca$.
\end{Lemma}
\medbreak

Let us denote by $\NormalForms(\Qca)$ the set of the normal forms of
$\RewTrees_\Qca$, described in the statement of
Lemma~\ref{lem:normal_forms_as_poset}. Moreover, we denote by
$\NormalForms(\Qca)(n)$, $n \geq 1$, the set $\NormalForms(\Qca)$
restricted to syntax trees with exactly~$n$ leaves. From their
description provided by Lemma~\ref{lem:normal_forms_as_poset}, any
tree $\Tfr$ of $\NormalForms(\Qca)$ different from the leaf is of the
recursive unique general form
\begin{equation} \label{equ:normal_forms_as_poset}
    \Tfr =
    \begin{tikzpicture}[xscale=.5,yscale=.45,Centering]
        \node(0)at(0.00,-1.67){\begin{math}\Sfr_1\end{math}};
        \node(2)at(2.00,-3.33){\begin{math}\Sfr_{\ell - 1}\end{math}};
        \node(4)at(4.00,-3.33){\begin{math}\Sfr_\ell\end{math}};
        \node(1)at(1.00,0.00){\begin{math}\Op_a\end{math}};
        \node(3)at(3.00,-1.67){\begin{math}\Op_a\end{math}};
        \draw[Edge](0)--(1);
        \draw[Edge](2)--(3);
        \draw[Edge,densely dashed](3)--(1);
        \draw[Edge](4)--(3);
        \node(r)at(1.00,1){};
        \draw[Edge](r)--(1);
    \end{tikzpicture}\,,
\end{equation}
where $a \in \Qca$ and the dashed edge denotes a right comb tree wherein
internal nodes are labeled by $\Op_a$, and for any $i \in [\ell]$,
$\Sfr_i$ is a tree of $\NormalForms(\Qca)$ such that $\Sfr_i$ is the
leaf or its root is labeled by a~$\Op_b$, $b \in\Qca$, so that $a$ and
$b$ are incomparable in~$\Qca$.
\medbreak

\begin{Lemma} \label{lem:confluent_rewrite_rule_as_poset}
    Let $\Qca$ be a forest poset. Then, $\RewTrees_\Qca$ is a confluent
    rewrite rule.
\end{Lemma}
\medbreak

In Lemma~\ref{lem:confluent_rewrite_rule_as_poset}, the condition on
$\Qca$ to be a forest poset is a necessary condition. Indeed, by setting
\begin{equation}
    \Qca :=
\,,
\end{equation}
the branching tree
\begin{equation}
    %
\end{equation}
of $\RewTrees_\Qca$ admits the branching pair consisting in the two
syntax trees
\begin{equation}
    %
\,.
\end{equation}
Since these two trees are normal forms of $\RewTrees_\Qca$, this
branching pair is not joinable, hence showing that $\RewTrees_\Qca$ is
not confluent.
\medbreak

\subsubsection{Koszulity}
Lemmas~\ref{lem:rewrite_rule_space_relations_as_poset},
\ref{lem:terminating_rewrite_rule},
\ref{lem:normal_forms_as_poset},
and~\ref{lem:confluent_rewrite_rule_as_poset} imply the following
result.
\begin{Theorem} \label{thm:koszulity_as_poset}
    Let $\Qca$ be a forest poset. Then, the operad $\As(\Qca)$ is Koszul
    and the set $\NormalForms(\Qca)$ forms a Poincaré-Birkhoff-Witt
    basis of~$\As(\Qca)$.
\end{Theorem}
\medbreak

\subsection{Dimensions and realization}
The Koszulity, and more specifically the existence of a
Poincaré-Birkhoff-Witt basis $\NormalForms(\Qca)$ highlighted by
Theorem~\ref{thm:koszulity_as_poset} for $\As(\Qca)$ when $\Qca$ is a
forest poset, lead to a combinatorial realization of $\As(\Qca)$. Before
describing this realization, we shall provide a functional equation for
the Hilbert series of $\As(\Qca)$.
\medbreak

\subsubsection{Dimensions}

\begin{Proposition} \label{prop:hilbert_series_as_poset}
    Let $\Qca$ be a forest poset. Then, the Hilbert series
    $\HilbSeries_\Qca(t)$ of $\As(\Qca)$ satisfies
    \begin{equation} \label{equ:hilbert_series_as_poset}
        \HilbSeries_\Qca(t) =
        t + \sum_{a \in \Qca} \HilbSeries_\Qca^a(t),
    \end{equation}
    where for all $a \in \Qca$, the $\HilbSeries_\Qca^a(t)$ satisfy
    \begin{equation} \label{equ:hilbert_series_h_a_as_poset}
        \HilbSeries_\Qca^a(t) =
        \left(t + \bar \HilbSeries_\Qca^a(t)\right)
        \left(t + \bar \HilbSeries_\Qca^a(t)
            + \HilbSeries_\Qca^a(t)\right),
    \end{equation}
    and for all $a \in \Qca$, the $\bar \HilbSeries_\Qca^a(t)$ satisfy
    \begin{equation} \label{equ:hilbert_series_h_a_incomp_as_poset}
        \bar \HilbSeries_\Qca^a(t) =
        \sum_{\substack{b \in \Qca \\ a \not \Ord_\Qca b \\
                b \not \Ord_\Qca a}}
        \HilbSeries_\Qca^b(t).
    \end{equation}
\end{Proposition}
\medbreak

For instance, let us use
Proposition~\ref{prop:hilbert_series_as_poset} for the operad
$\As(\Qca)$ when $\Qca$ is the total order on the set $[\ell]$,
$\ell \geq 0$. This operad is the multiassociative operad, whose
definition is recalled in
Section~\ref{subsubsec:first_examples_as_poset}.
By~\eqref{equ:hilbert_series_h_a_incomp_as_poset}, we have
\begin{equation}
    \bar \HilbSeries_\Qca^a(t) = 0,
    \qquad a \in [\ell],
\end{equation}
and hence, by~\eqref{equ:hilbert_series_h_a_as_poset},
\begin{equation}
    \HilbSeries_\Qca^a(t) = \frac{t^2}{1 - t},
    \qquad a \in [\ell].
\end{equation}
Then, by~\eqref{equ:hilbert_series_as_poset}, the Hilbert series of
$\As(\Qca)$ satisfies
\begin{equation}
    \HilbSeries_\Qca(t) = t + \frac{\ell t^2}{1 - t},
    \qquad \ell \geq 0.
\end{equation}
\medbreak

Let us use Proposition~\ref{prop:hilbert_series_as_poset} for the operad
$\As(\Qca)$ when $\Qca$ is the trivial poset on the set $[\ell]$,
$\ell \geq 0$. This operad is the dual multiassociative
operad, whose definition is recalled in
Section~\ref{subsubsec:first_examples_as_poset}.
By~\eqref{equ:hilbert_series_h_a_incomp_as_poset}, one has
\begin{equation}
    \bar \HilbSeries_\Qca^a(t) =
    \sum_{\substack{b \in [\ell] \\ b \ne a}}
    \bar \HilbSeries_\Qca^b(t),
    \qquad a \in [\ell],
\end{equation}
implying, by~\eqref{equ:hilbert_series_as_poset}, that
\begin{equation}
    \bar \HilbSeries_\Qca^a(t)
    = \HilbSeries_\Qca(t) - t - \HilbSeries_\Qca^a(t),
    \qquad a \in [\ell].
\end{equation}
Now, by~\eqref{equ:hilbert_series_h_a_as_poset}, we obtain
\begin{equation}
    \HilbSeries_\Qca^a(t) =
    \frac{\HilbSeries_\Qca(t)^2}{1 + \HilbSeries_\Qca(t)},
    \qquad a \in [\ell].
\end{equation}
Therefore, by~\eqref{equ:hilbert_series_as_poset}, the
Hilbert series of $\As(\Qca)$ satisfies the quadratic functional
equation
\begin{equation}
    t + (t - 1)\HilbSeries_\Qca(t) + (\ell - 1) \HilbSeries_\Qca(t)^2
    = 0,
    \qquad \ell \geq 0,
\end{equation}
and can be expressed as
\begin{equation}
    \HilbSeries_\Qca(t) =
    \frac{1 - t - \sqrt{1 + (2 - 4 \ell) t + t^2}}
    {2 (\ell - 1)},
    \qquad \ell = 0 \mbox{ or } \ell \geq 2.
\end{equation}
The dimensions of the first homogeneous components of $\As(\Qca)$ are
\begin{subequations}
\begin{equation}
    1, 2, 6, 22, 90, 394, 1806, 8558, 41586, 206098,
    \qquad \ell = 2,
\end{equation}
\begin{equation}
    1, 3, 15, 93, 645, 4791, 37275, 299865, 2474025, 20819307,
    \qquad \ell = 3,
\end{equation}
\begin{equation}
    1, 4, 28, 244, 2380, 24868, 272188, 3080596, 35758828, 423373636,
    \qquad \ell = 4,
\end{equation}
\begin{equation}
    1, 5, 45, 505, 6345, 85405, 1204245, 17558705, 262577745,
    4005148405,
    \qquad \ell = 5.
\end{equation}
\end{subequations}
These sequences are respectively Sequences~\OEIS{A006318},
\OEIS{A103210}, \OEIS{A103211}, and \OEIS{A133305} of~\cite{Slo}.
\medbreak

Finally, let us use Proposition~\ref{prop:hilbert_series_as_poset} for
the operad $\As(\Qca)$ when $\Qca$ is the forest poset
\begin{equation}
    \Qca :=
    \begin{tikzpicture}[xscale=.5,yscale=.5,Centering]
        \node[PosetVertex](1)at(0,0){\begin{math}1\end{math}};
        \node[PosetVertex](2)at(0,-1){\begin{math}2\end{math}};
        \node[PosetVertex](3)at(1,0){\begin{math}3\end{math}};
        \node[PosetVertex](4)at(1,-1){\begin{math}4\end{math}};
        \draw[Edge](1)--(2);
        \draw[Edge](3)--(4);
    \end{tikzpicture}\,.
\end{equation}
By~\eqref{equ:hilbert_series_h_a_incomp_as_poset}, one has
\begin{subequations}
\begin{equation}
    \bar \HilbSeries_\Qca^1(t) = \bar \HilbSeries_\Qca^2(t) =
    \HilbSeries_\Qca^3(t) + \HilbSeries_\Qca^4(t),
\end{equation}
\begin{equation}
    \bar \HilbSeries_\Qca^3(t) = \bar \HilbSeries_\Qca^4(t) =
    \HilbSeries_\Qca^1(t) + \HilbSeries_\Qca^2(t),
\end{equation}
\end{subequations}
and, by~\eqref{equ:hilbert_series_h_a_as_poset} and straightforward
computations, we obtain that
\begin{equation}
    \HilbSeries_\Qca^1(t) = \HilbSeries_\Qca^2(t) =
    \HilbSeries_\Qca^3(t) = \HilbSeries_\Qca^4(t),
\end{equation}
so that the Hilbert series of $\As(\Qca)$ satisfies the quadratic
functional equation
\begin{equation}
    \frac{1}{2}t + \frac{1}{4}t^2 +
    \left(t - \frac{1}{2}\right)\HilbSeries_\Qca(t) +
    \frac{3}{4}\HilbSeries_\Qca(t)^2 = 0.
\end{equation}
This Hilbert series can expressed as
\begin{equation}
    \HilbSeries_\Qca(t) =
    \frac{1 - 2t - \sqrt{1 - 10t + t^2}}{3},
\end{equation}
and the dimensions of the first homogeneous components of $\As(\Qca)$
are
\begin{equation}
    1, 4, 20, 124, 860, 6388, 49700, 399820, 3298700, 27759076.
\end{equation}
Terms of this sequence are the ones of Sequence~\OEIS{A107841}
of~\cite{Slo} multiplied by $2$.
\medbreak

\subsubsection{Realization}
\label{subsubsec:realization_as_poset_forest}
Let us describe a combinatorial realization of $\As(\Qca)$ when $\Qca$
is a forest poset in terms of Schröder trees (see
Section~\ref{subsubsec:schroder_trees} of
Chapter~\ref{chap:combinatorics}) with a certain labeling and through an
algorithm to compute their partial composition.
\medbreak

If $\Qca$ is a poset (not necessarily a forest poset just now), a
\Def{$\Qca$-Schröder tree} is a Schröder tree where internal nodes are
labeled on $\Qca$. For any element $a$ of $\Qca$ and any $n \geq 2$, we
denote by $\CorSchr_a^n$ the $\Qca$-Schröder tree consisting in a single
internal node labeled by $a$ attached to $n$ leaves. We call these trees
\Def{$\Qca$-corollas}. A \Def{$\Qca$-alternating Schröder tree} is a
$\Qca$-Schröder tree $\Tfr$ such that for any internal node $y$ of
$\Tfr$ having a father $x$, the labels of $x$ and $y$ are incomparable
in $\Qca$. We denote by $\SetASchr(\Qca)$ the set of all
$\Qca$-alternating Schröder trees and by $\SetASchr(\Qca)(n)$,
$n \geq 1$, the set $\SetASchr(\Qca)$ restricted to trees with exactly
$n$ leaves. Any tree $\Tfr$ of $\SetASchr(\Qca)$ different from the
leaf is of the recursive unique general form
\begin{equation} \label{equ:description_aschr}
    \Tfr =
    \begin{tikzpicture}[xscale=.6,yscale=.3,Centering]
        \node(0)at(0.00,-2.00)
            {\footnotesize \begin{math}\Sfr_1\end{math}};
        \node(3)at(2.00,-2.00)
            {\footnotesize \begin{math}\Sfr_\ell\end{math}};
        \node[Node](1)at(1.00,0.00){\begin{math}a\end{math}};
        \draw[Edge](0)--(1);
        \draw[Edge](3)--(1);
        \node(r)at(1.00,1.25){};
        \draw[Edge](r)--(1);
        \node[below of=1,node distance=.6cm]
            {\footnotesize\begin{math}\dots\end{math}};
    \end{tikzpicture}\,,
\end{equation}
where $a \in \Qca$ and for any $i \in [\ell]$, $\Sfr_i$ is a tree of
$\SetASchr(\Qca)$ such that $\Sfr_i$ is a leaf or its root is labeled by
a $b \in \Qca$ and $a$ and $b$ and incomparable in~$\Qca$.
\medbreak

Relying on the description of the elements of $\NormalForms(\Qca)$
provided by Lemma~\ref{lem:normal_forms_as_poset} and on their recursive
general form provided by~\eqref{equ:normal_forms_as_poset}, let us
consider the map
\begin{equation}
    \ASchrMap_\Qca : \NormalForms(\Qca)(n) \to \SetASchr(\Qca)(n),
    \qquad n \geq 1,
\end{equation}
defined recursively by sending the leaf to the leaf and, for any tree
$\Tfr$ of $\NormalForms(\Qca)$ different from the leaf, by
\begin{equation} \label{equ:recursive_description_aschr_map}
    \ASchrMap_\Qca(\Tfr) =
    \ASchrMap_\Qca
    \left(
    \begin{tikzpicture}[xscale=.4,yscale=.35,Centering]
        \node[NodeST](0)at(0.00,-1.67)
            {\begin{math}\Sfr_1\end{math}};
        \node[NodeST](2)at(2.00,-3.33)
            {\begin{math}\Sfr_{\ell - 1}\end{math}};
        \node[NodeST](4)at(4.00,-3.33)
            {\begin{math}\Sfr_\ell\end{math}};
        \node[NodeST](1)at(1.00,0.00)
            {\begin{math}\Op_a\end{math}};
        \node[NodeST](3)at(3.00,-1.67)
            {\begin{math}\Op_a\end{math}};
        \draw[Edge](0)--(1);
        \draw[Edge](2)--(3);
        \draw[Edge,densely dashed](3)--(1);
        \draw[Edge](4)--(3);
        \node(r)at(1.00,1.25){};
        \draw[Edge](r)--(1);
    \end{tikzpicture}
    \right) :=
    \begin{tikzpicture}[xscale=.9,yscale=.3,Centering]
        \node(0)at(0.00,-2.00)
        {\footnotesize \begin{math}\ASchrMap_\Qca(\Sfr_1)\end{math}};
        \node(3)at(2.00,-2.00)
        {\footnotesize \begin{math}\ASchrMap_\Qca(\Sfr_\ell)\end{math}};
        \node[Node](1)at(1.00,0.00){\begin{math}a\end{math}};
        \draw[Edge](0)--(1);
        \draw[Edge](3)--(1);
        \node(r)at(1.00,1.25){};
        \draw[Edge](r)--(1);
        \node[below of=1,node distance=.6cm]
            {\footnotesize\begin{math}\dots\end{math}};
    \end{tikzpicture}\,,
\end{equation}
where $a \in \Qca$ and, in the syntax tree
of~\eqref{equ:recursive_description_aschr_map}, the dashed edge denotes
a right comb tree wherein internal nodes are labeled by $\Op_a$, and for
any $i \in [\ell]$, $\Sfr_i$ is a tree of $\NormalForms(\Qca)$ such that
$\Sfr_i$ is the leaf or its root is labeled by $\Op_b$, $b \in \Qca$,
and $a$ and $b$ are incomparable in $\Qca$. It is immediate that
$\ASchrMap_\Qca(\Tfr)$ is a $\Qca$-alternating Schröder tree, so that
$\ASchrMap_\Qca$ is a well-defined map.
\medbreak

\begin{Lemma} \label{lem:bijection_pbw_basis_aschr}
    Let $\Qca$ be a poset. Then, for any $n \geq 1$, the map
    $\ASchrMap_\Qca$ is a bijection between the set of syntax trees of
    $\NormalForms(\Qca)(n)$ with $n$ leaves and the set
    $\SetASchr(\Qca)(n)$ of $\Qca$-alternating Schröder trees with $n$
    leaves.
\end{Lemma}
\medbreak

In order to define a partial composition for $\Qca$-alternating Schröder
trees, we introduce the following rewrite rule. When $\Qca$ is a forest
poset, consider the rewrite rule $\RewSchr_\Qca$ on $\Qca$-Schröder
trees (not necessarily $\Qca$-alternating Schröder trees) satisfying
\begin{equation} \label{equ:rewrite_schroder}
    \begin{tikzpicture}[xscale=.35,yscale=.25,Centering]
        \node[Leaf](0)at(0.00,-2.00){};
        \node[Leaf](2)at(1.00,-4.00){};
        \node[Leaf](4)at(3.00,-4.00){};
        \node[Leaf](5)at(4.00,-2.00){};
        \node[Node](1)at(2.00,0.00){\begin{math}b\end{math}};
        \node[Node](3)at(2.00,-2.00){\begin{math}a\end{math}};
        \draw[Edge](0)--(1);
        \draw[Edge](2)--(3);
        \draw[Edge](3)--(1);
        \draw[Edge](4)--(3);
        \draw[Edge](5)--(1);
        \node(r)at(2.00,1.50){};
        \draw[Edge](r)--(1);
        \node[right of=0,node distance=.32cm]
            {\footnotesize\begin{math}\dots\end{math}};
        \node[left of=5,node distance=.32cm]
            {\footnotesize\begin{math}\dots\end{math}};
        \node[right of=2,node distance=.4cm]
            {\footnotesize\begin{math}\dots\end{math}};
    \end{tikzpicture}
    \enspace \RewSchr_\Qca \enspace
    \begin{tikzpicture}[xscale=.45,yscale=.5,Centering]
        \node[Leaf](0)at(0.00,-1.50){};
        \node[Leaf](2)at(2.00,-1.50){};
        \node[Node](1)at(1.00,0.00){\begin{math}a \Min_\Qca b\end{math}};
        \draw[Edge](0)--(1);
        \draw[Edge](2)--(1);
        \node(r)at(1.00,1.5){};
        \draw[Edge](r)--(1);
        \node[below of=1,node distance=.8cm]
            {\footnotesize\begin{math}\dots\end{math}};
    \end{tikzpicture},
    \qquad
    a, b \in \Qca \mbox{ and }
    (a \Ord_\Qca b \mbox{ or } b \Ord_\Qca a).
\end{equation}
Let also $\RewTreesSchr_\Qca$ be the closure of $\RewSchr_\Qca$.
Equation~\eqref{equ:example_schroder_rewriting_steps} shows examples
of steps of rewritings by $\RewTreesSchr_\Qca$ for the poset $\Qca$
defined in~\eqref{equ:example_forest_poset}.
\medbreak

\begin{Lemma} \label{lem:rewrite_rule_schroder_trees}
    Let $\Qca$ be a poset. Then, $\RewTreesSchr_\Qca$ is a terminating
    rewrite rule and the set of its normal forms is the set of all
    $\Qca$-alternating Schröder trees. Moreover, when~$\Qca$ is a
    forest poset, $\RewTreesSchr_\Qca$ is confluent.
\end{Lemma}
\medbreak

In Lemma~\ref{lem:rewrite_rule_schroder_trees}, the condition on $\Qca$
to be a forest poset is a necessary condition for the confluence of
$\RewTreesSchr_\Qca$. Indeed, by setting
\begin{equation}
    \Qca :=
\,,
\end{equation}
the branching tree
\begin{equation}
    %
\end{equation}
of $\RewTreesSchr_\Qca$ admits the branching pair consisting in the two
trees
\begin{equation}
    %
\,.
\end{equation}
Since these two trees are normal forms of $\RewTreesSchr_\Qca$, this
branching pair is not joinable and hence, $\RewTreesSchr_\Qca$ is not
confluent.
\medbreak

We define the partial composition $\Sfr \circ_i \Tfr$ of two
$\Qca$-alternating Schröder trees $\Sfr$ and $\Tfr$ as the
$\Qca$-alternating Schröder tree being the normal form by
$\RewTreesSchr_\Qca$ of the $\Qca$-Schröder tree obtained by grafting
the root of $\Tfr$ on the $i$th leaf of $\Sfr$. We denote by
$\ASchr(\Qca)$ the linear span of the set of the $\Qca$-alternating
Schröder trees endowed with the partial composition described above
and extended by linearity. Consider for instance the forest poset
\begin{equation} \label{equ:example_forest_poset}
    \Qca :=
\,.
\end{equation}
Then, we have in $\ASchr(\Qca)$ the partial composition
\begin{equation}
    %
\end{equation}
is a sequence of rewritings steps by $\RewTreesSchr_\Qca$, where the
leftmost tree of~\eqref{equ:example_schroder_rewriting_steps} is
obtained by grafting the root of the second tree
of~\eqref{equ:example_schroder_partial_composition} onto the first leaf
of the first tree of~\eqref{equ:example_schroder_partial_composition}.
\medbreak

\begin{Proposition} \label{prop:operad_aschr}
    Let $\Qca$ be a forest poset. Then, $\ASchr(\Qca)$ is an operad
    graded by the number of the leaves of the trees. Moreover, as an
    operad, $\ASchr(\Qca)$ is generated by the set of $\Qca$-corollas of
    arity two.
\end{Proposition}
\medbreak

\begin{Theorem} \label{thm:realization_as_poset}
    Let $\Qca$ be a forest poset. Then, the operads $\As(\Qca)$ and
    $\ASchr(\Qca)$ are isomorphic.
\end{Theorem}
\begin{proof}
    First, by Proposition~\ref{prop:operad_aschr}, $\ASchr(\Qca)$ is an
    operad wherein for any $n \geq 1$, its graded component of arity
    $n$ has bases indexed by $\Qca$-alternating Schröder trees with $n$
    leaves. By Lemma~\ref{lem:bijection_pbw_basis_aschr}, these trees
    are in bijection with the elements of the Poincaré-Birkhoff-Witt
    basis $\NormalForms(\Qca)$ of $\As(\Qca)$ provided by
    Theorem~\ref{thm:koszulity_as_poset}. By~\cite{Hof10}, this shows
    that $\ASchr(\Qca)$ and $\As(\Qca)$ are isomorphic as graded
    vector spaces.
    \smallbreak

    The generators of $\ASchr(\Qca)$, that are by
    Proposition~\ref{prop:operad_aschr} $\Qca$-corollas of arity two,
    satisfy at least the nontrivial relations
    \begin{subequations}
    \begin{equation} \label{equ:realization_as_poset_relation_1}
        \CorSchr_a^2 \circ_1 \CorSchr_b^2
        -
        \CorSchr_{a \Min_\Qca b}^2 \circ_2 \CorSchr_{a \Min_\Qca b}^2
        = 0,
        \qquad a, b \in \Qca \mbox{ and }
        (a \Ord_\Qca b \mbox{ or } b \Ord_\Qca a),
    \end{equation}
    \begin{equation} \label{equ:realization_as_poset_relation_2}
        \CorSchr_{a \Min_\Qca b}^2 \circ_1 \CorSchr_{a \Min_\Qca b}^2
        -
        \CorSchr_a^2 \circ_2 \CorSchr_b^2 = 0,
        \qquad a, b \in \Qca \mbox{ and }
        (a \Ord_\Qca b \mbox{ or } b \Ord_\Qca a),
    \end{equation}
    \end{subequations}
    obtained by a direct computation in $\ASchr(\Qca)$. By using the
    same reasoning as the one used to establish
    Proposition~\ref{prop:dimension_relations_as_poset}, we obtain that
    there are as many elements of the
    form~\eqref{equ:realization_as_poset_relation_1}
    or~\eqref{equ:realization_as_poset_relation_2} as generating
    relations (see~\eqref{equ:relation_as_poset_1}
    and~\eqref{equ:relation_as_poset_2}) for the space of relations
    $\RelationSpace_\Qca^\Op$ of $\As(\Qca)$ . Therefore, as
    $\ASchr(\Qca)$ and $\As(\Qca)$ are isomorphic as graded vector
    spaces, it cannot be more nontrivial relations in $\ASchr(\Qca)$
    than Relations~\eqref{equ:realization_as_poset_relation_1}
    and~\eqref{equ:realization_as_poset_relation_2}.
    \smallbreak

    Finally, by identifying all symbols $\CorSchr_a^2$, $a \in \Qca$,
    with $\Op_a$, we observe that $\As(\Qca)$ and $\ASchr(\Qca)$ admit
    the same presentation. This implies that $\As(\Qca)$ and
    $\ASchr(\Qca)$ are isomorphic operads.
\end{proof}
\medbreak

As announced, Theorem~\ref{thm:realization_as_poset} provides a
combinatorial realization $\ASchr(\Qca)$ of $\As(\Qca)$ when $\Qca$ is
a forest poset.
\medbreak

\subsubsection{Free forest poset associative algebras over one
generator} \label{subsubsec:free_as_forest_poset_algebras}
The realization of $\As(\Qca)$, when $\Qca$ is a forest poset, provided
by Theorem~\ref{thm:realization_as_poset} in terms of
$\Qca$-alternating Schröder trees leads to the following description.
The free $\Qca$-associative algebra over one generator, where $\Qca$
is a forest poset, has $\ASchr(\Qca)$ as underlying vector space and
is endowed with linear operations
\begin{equation}
    \Op_a : \ASchr(\Qca) \otimes \ASchr(\Qca) \to \ASchr(\Qca),
    \qquad a \in \Qca,
\end{equation}
satisfying for all $\Qca$-alternating Schröder trees $\Sfr$ and $\Tfr$,
\begin{equation}
    \Sfr \Op_a \Tfr
    = \left(\CorSchr_a^2 \circ_2 \Tfr\right) \circ_1 \Sfr.
\end{equation}
In an alternative way, $\Sfr \Op_a \Tfr$ is the $\Qca$-alternating
Schröder obtained by considering the normal form by
$\RewTreesSchr_\Qca$ of the tree obtained by grafting $\Sfr$ and $\Tfr$
respectively as left and right child of a binary corolla labeled by~$a$.
\medbreak

Let us provide examples of computations in the free $\Qca$-associative
algebra over one generator where $\Qca$ is the forest poset
\begin{equation}
    \Qca :=
\,.
\end{equation}
\end{subequations}
\medbreak

\subsection{Koszul dual}
We now establish a first presentation for the Koszul dual $\As(\Qca)^!$
of $\As(\Qca)$ where $\Qca$ is a poset (and not necessarily a forest
poset) and provide moreover a second presentation of $\As(\Qca)^!$ when
$\Qca$ is a forest poset. This second presentation of $\As(\Qca)^!$ is
simpler than the first one and it shall be considered in the next
section.
\medbreak

\subsubsection{Presentation by generators and relations}

\begin{Proposition} \label{prop:presentation_koszul_dual_as_poset}
    Let $\Qca$ be a poset. Then, the Koszul dual $\As(\Qca)^!$ of
    $\As(\Qca)$ admits the presentation
    \begin{math}
        \left(\GeneratingSet_\Qca^{\OpDual},
        \RelationSpace_\Qca^{\OpDual} \right)
    \end{math}
    where
    \begin{equation}
        \GeneratingSet_\Qca^{\OpDual}
        := \GeneratingSet_\Qca^{\OpDual}(2)
        := \left\{\OpDual_a : a \in \Qca\right\},
    \end{equation}
    and $\RelationSpace_\Qca^{\OpDual}$ is the subspace of
    $\FreeOperad\left(\GeneratingSet_\Qca^{\OpDual}\right)$ generated by
    \begin{subequations}
    \begin{multline} \label{equ:relation_dual_as_poset_1}
        \Corolla{\OpDual_a} \circ_1 \Corolla{\OpDual_a}
        -
        \Corolla{\OpDual_a} \circ_2 \Corolla{\OpDual_a} \\
        +
        \sum_{\substack{b \in \Qca \\ a \OrdStrict_\Qca b}}
        \left(\Corolla{\OpDual_b} \circ_1 \Corolla{\OpDual_a}
        +
        \Corolla{\OpDual_a} \circ_1 \Corolla{\OpDual_b}
        -
        \Corolla{\OpDual_b} \circ_2 \Corolla{\OpDual_a}
        -
        \Corolla{\OpDual_a} \circ_2 \Corolla{\OpDual_b}\right),
        \qquad a \in \Qca,
    \end{multline}
    \begin{equation} \label{equ:relation_dual_as_poset_2}
        \Corolla{\OpDual_c} \circ_1 \Corolla{\OpDual_d},
        \qquad c, d \in \Qca
        \mbox{ and } c \not \Ord_\Qca d
        \mbox{ and } d \not \Ord_\Qca c,
    \end{equation}
    \begin{equation} \label{equ:relation_dual_as_poset_3}
        \Corolla{\OpDual_c} \circ_2 \Corolla{\OpDual_d},
        \qquad c, d \in \Qca
        \mbox{ and } c \not \Ord_\Qca d
        \mbox{ and } d \not \Ord_\Qca c.
    \end{equation}
    \end{subequations}
\end{Proposition}
\medbreak

\begin{Proposition} \label{prop:dimension_relations_dual_as_poset}
    Let $\Qca$ be a poset. Then, the dimension of the space
    $\RelationSpace_\Qca^{\OpDual}$ of relations of $\As(\Qca)^!$
    satisfies
    \begin{equation}
        \dim \RelationSpace_\Qca^{\OpDual} =
        2\, (\# \Qca)^2 + 3\, \# \Qca - 4\, \NbInterv(\Qca).
    \end{equation}
\end{Proposition}
\medbreak

Observe that, by Propositions~\ref{prop:dimension_relations_as_poset}
and~\ref{prop:dimension_relations_dual_as_poset}, we have
\begin{equation}\begin{split}
    \dim \RelationSpace_\Qca^\Op + \dim \RelationSpace_\Qca^{\OpDual}
    & = 4 \, \NbInterv(\Qca) - 3 \, \# \Qca +
        2\, (\# \Qca)^2 + 3\, \# \Qca - 4\, \NbInterv(\Qca) \\
    & = 2\, (\# \Qca)^2 \\
    & = \dim \FreeOperad(\GeneratingSet_\Qca^\Op)(3),
\end{split}\end{equation}
as expected by Koszul duality.
\medbreak

\subsubsection{Alternative presentation}
\label{subsubsec:alternative_presentation_koszul_dual_as_poset}
For any element $a$ of a poset $\Qca$ (not necessarily a forest poset
just now), let $\OpDualB_a$ be the element of $\As(\Qca)^!(2)$ defined
by
\begin{equation} \label{equ:basis_change_dual_as_poset}
    \OpDualB_a :=
    \sum_{\substack{b \in \Qca \\ a \Ord_\Qca b}}
    \OpDual_b.
\end{equation}
We denote by $\GeneratingSet_\Qca^{\OpDualB}$ the set of all
$\OpDualB_a$, $a \in \Qca$. By triangularity, the family
$\GeneratingSet_{\Qca}^{\OpDualB}$ forms a basis of
$\As(\Qca)^!(2)$ and hence, generates $\As(\Qca)^!$. Consider for
instance the poset
\begin{equation}
    \Qca :=
    \begin{tikzpicture}[Centering,xscale=.4,yscale=.42]
        \node[PosetVertex](1)at(0,0){\begin{math}1\end{math}};
        \node[PosetVertex](2)at(2,0){\begin{math}2\end{math}};
        \node[PosetVertex](3)at(1,-1){\begin{math}3\end{math}};
        \node[PosetVertex](4)at(3,-1){\begin{math}4\end{math}};
        \node[PosetVertex](5)at(2,-2){\begin{math}5\end{math}};
        \draw[Edge](1)--(3);
        \draw[Edge](2)--(3);
        \draw[Edge](2)--(4);
        \draw[Edge](4)--(5);
    \end{tikzpicture}\,.
\end{equation}
The elements of $\GeneratingSet_\Qca^{\OpDualB}$ then express as
\begin{subequations}
\begin{equation}
    \OpDualB_1 = \OpDual_1 + \OpDual_3,
\end{equation}
\begin{equation}
    \OpDualB_2 = \OpDual_2 + \OpDual_3 + \OpDual_4 + \OpDual_5,
\end{equation}
\begin{equation}
    \OpDualB_3 = \OpDual_3,
\end{equation}
\begin{equation}
    \OpDualB_4 = \OpDual_4 + \OpDual_5,
\end{equation}
\begin{equation}
    \OpDualB_5 = \OpDual_5.
\end{equation}
\end{subequations}
\medbreak

\begin{Proposition}
    \label{prop:alternative_presentation_koszul_dual_as_poset}
    Let $\Qca$ be a forest poset. Then, the operad $\As(\Qca)^!$ admits
    the presentation
    \begin{math}
        \left(\GeneratingSet_{\Qca}^{\OpDualB},
        \RelationSpace_{\Qca}^{\OpDualB}\right)
    \end{math}
    where $\RelationSpace_{\Qca}^{\OpDualB}$ is the subspace of
    $\FreeOperad\left(\GeneratingSet_{\Qca}^{\OpDualB}\right)$
    generated by
    \begin{subequations}
    \begin{equation} \label{equ:alternative_relation_dual_as_poset_1}
        \Corolla{\OpDualB_a} \circ_1 \Corolla{\OpDualB_a}
        -
        \Corolla{\OpDualB_a} \circ_2 \Corolla{\OpDualB_a},
        \qquad a \in \Qca,
    \end{equation}
    \begin{equation} \label{equ:alternative_relation_dual_as_poset_2}
        \Corolla{\OpDualB_c} \circ_1 \Corolla{\OpDualB_d},
        \qquad c, d \in \Qca
        \mbox{ and } c \not \Ord_\Qca d
        \mbox{ and } d \not \Ord_\Qca c,
    \end{equation}
    \begin{equation} \label{equ:alternative_relation_dual_as_poset_3}
        \Corolla{\OpDualB_c} \circ_2 \Corolla{\OpDualB_d},
        \qquad c, d \in \Qca
        \mbox{ and } c \not \Ord_\Qca d
        \mbox{ and } d \not \Ord_\Qca c.
    \end{equation}
    \end{subequations}
\end{Proposition}
\medbreak

By considering the presentation of $\As(\Qca)^!$ furnished by
Proposition~\ref{prop:alternative_presentation_koszul_dual_as_poset}
when $\Qca$ is a forest poset, we obtain by Koszul duality a new
presentation
\begin{math}
    \left(\GeneratingSet_\Qca^{\OpB},
    \RelationSpace_\Qca^{\OpB}\right)
\end{math}
for $\As(\Qca)$ where the set of generators
$\GeneratingSet_\Qca^{\OpB}$ is defined by
\begin{equation}
    \GeneratingSet_\Qca^{\OpB}
    := \GeneratingSet_\Qca^{\OpB}(2)
    := \left\{\OpB_a : a \in \Qca\right\},
\end{equation}
and the space of relations $\RelationSpace_\Qca^{\OpB}$ is generated by
\begin{subequations}
\begin{equation}
    \Corolla{\OpB_a} \circ_1 \Corolla{\OpB_a}
    -
    \Corolla{\OpB_a} \circ_2 \Corolla{\OpB_a},
    \qquad a \in \Qca,
\end{equation}
\begin{equation}
    \Corolla{\OpB_a} \circ_1 \Corolla{\OpB_b},
    \qquad a, b \in \Qca
    \mbox{ and (} a \OrdStrict_\Qca b
    \mbox{ or } b \OrdStrict_\Qca a \mbox{)},
\end{equation}
\begin{equation}
    \Corolla{\OpB_a} \circ_2 \Corolla{\OpB_b},
    \qquad a, b \in \Qca
    \mbox{ and (} a \OrdStrict_\Qca b
    \mbox{ or } b \OrdStrict_\Qca a \mbox{)}.
\end{equation}
\end{subequations}
\medbreak

\subsubsection{Example}
To end this section, let us give a complete example of the spaces of
relations $\RelationSpace_\Qca^\Op$, $\RelationSpace_\Qca^{\OpDual}$,
$\RelationSpace_\Qca^{\OpDualB}$, and $\RelationSpace_\Qca^{\OpB}$ of
the operads $\As(\Qca)$ and $\As(\Qca)^!$ where $\Qca$ is the forest
poset
\begin{equation}
    \Qca :=
    \begin{tikzpicture}[Centering,xscale=.35,yscale=.42]
        \node[PosetVertex](1)at(0,0){\begin{math}1\end{math}};
        \node[PosetVertex](2)at(-1,-1){\begin{math}2\end{math}};
        \node[PosetVertex](3)at(1,-1){\begin{math}3\end{math}};
        \draw[Edge](1)--(2);
        \draw[Edge](1)--(3);
    \end{tikzpicture}\,.
\end{equation}
\medbreak

First, by definition of $\As$ and by
Lemma~\ref{lem:relations_as_poset}, the generators of
$\GeneratingSet_\Qca^\Op$ are subjected to the relations
\begin{subequations}
\begin{equation}\begin{split}
    \Op_1 \circ_1 \Op_1
    = \Op_1 \circ_1 \Op_2
    & = \Op_2 \circ_1 \Op_1
    = \Op_1 \circ_1 \Op_3
    = \Op_3 \circ_1 \Op_1 \\
    & = \Op_3 \circ_2 \Op_1
    = \Op_1 \circ_2 \Op_3
    = \Op_2 \circ_2 \Op_1
    = \Op_1 \circ_2 \Op_2
    = \Op_1 \circ_2 \Op_1,
\end{split}
\end{equation}
\begin{equation}
    \Op_2 \circ_1 \Op_2 = \Op_2 \circ_2 \Op_2,
\end{equation}
\begin{equation}
    \Op_3 \circ_1 \Op_3 = \Op_3 \circ_2 \Op_3.
\end{equation}
\end{subequations}
This describes $\RelationSpace_\Qca^\Op$.
\medbreak

By Proposition~\ref{prop:presentation_koszul_dual_as_poset}, the
generators of $\GeneratingSet_\Qca^{\OpDual}$ are subjected to the
relations
\begin{subequations}
\begin{equation}\begin{split}
    \OpDual_1 \circ_1 \OpDual_1
    + \OpDual_1 \circ_1 \OpDual_2
    & + \OpDual_2 \circ_1 \OpDual_1
    + \OpDual_1 \circ_1 \OpDual_3
    + \OpDual_3 \circ_1 \OpDual_1 \\
    & = \OpDual_3 \circ_2 \OpDual_1
    + \OpDual_1 \circ_2 \OpDual_3
    + \OpDual_2 \circ_2 \OpDual_1
    + \OpDual_1 \circ_2 \OpDual_2
    + \OpDual_1 \circ_2 \OpDual_1,
\end{split}
\end{equation}
\begin{equation}
    \OpDual_2 \circ_1 \OpDual_2 = \OpDual_2 \circ_2 \OpDual_2,
\end{equation}
\begin{equation}
    \OpDual_3 \circ_1 \OpDual_3 = \OpDual_3 \circ_2 \OpDual_3,
\end{equation}
\begin{equation}
    \OpDual_2 \circ_1 \OpDual_3
    = \OpDual_3 \circ_1 \OpDual_2
    = \OpDual_3 \circ_2 \OpDual_2
    = \OpDual_2 \circ_2 \OpDual_3 = 0.
\end{equation}
\end{subequations}
This describes $\RelationSpace_\Qca^{\OpDual}$.
\medbreak

By
Proposition~\ref{prop:alternative_presentation_koszul_dual_as_poset},
the generators of $\GeneratingSet_\Qca^{\OpDualB}$ are subjected to the
relations
\begin{subequations}
\begin{equation}
    \OpDualB_1 \circ_1 \OpDualB_1 = \OpDualB_1 \circ_2 \OpDualB_1,
\end{equation}
\begin{equation}
    \OpDualB_2 \circ_1 \OpDualB_2 = \OpDualB_2 \circ_2 \OpDualB_2,
\end{equation}
\begin{equation}
    \OpDualB_3 \circ_1 \OpDualB_3 = \OpDualB_3 \circ_2 \OpDualB_3,
\end{equation}
\begin{equation}
    \OpDualB_2 \circ_1 \OpDualB_3
    = \OpDualB_3 \circ_1 \OpDualB_2
    = \OpDualB_3 \circ_2 \OpDualB_2
    = \OpDualB_2 \circ_2 \OpDualB_3 = 0.
\end{equation}
\end{subequations}
This describes $\RelationSpace_\Qca^{\OpDualB}$.
\medbreak

Finally, by the observation established at the end of
Section~\ref{subsubsec:alternative_presentation_koszul_dual_as_poset},
the generators of $\GeneratingSet_\Qca^{\OpB}$ are subjected to the
relations
\begin{subequations}
\begin{equation}
    \OpB_1 \circ_1 \OpB_1 = \OpB_1 \circ_2 \OpB_1,
\end{equation}
\begin{equation}
    \OpB_2 \circ_1 \OpB_2 = \OpB_2 \circ_2 \OpB_2,
\end{equation}
\begin{equation}
    \OpB_3 \circ_1 \OpB_3 = \OpB_3 \circ_2 \OpB_3,
\end{equation}
\begin{equation}\begin{split}
    \OpB_1 \circ_1 \OpB_2
    = \OpB_2 \circ_1 \OpB_1
    & = \OpB_1 \circ_1 \OpB_3
    = \OpB_3 \circ_1 \OpB_1 \\
    & = \OpB_3 \circ_2 \OpB_1
    = \OpB_1 \circ_2 \OpB_3
    = \OpB_2 \circ_2 \OpB_1
    = \OpB_1 \circ_2 \OpB_2 = 0.
\end{split}\end{equation}
\end{subequations}
This describes $\RelationSpace_\Qca^{\OpB}$.
\medbreak

\section{Thin forest posets and Koszul duality}
\label{sec:thin_forest_posets}
As we have seen in Section~\ref{sec:forest_posets}, certain properties
satisfied by the poset $\Qca$ imply properties for the operad
$\As(\Qca)$. In this section, we show that when $\Qca$ is a forest
poset with an extra condition, the Koszul dual $\As(\Qca)^!$ of
$\As(\Qca)$ can be constructed via the construction~$\As$.
\medbreak

\subsection{Thin forest posets}
A subclass of the class of forest posets, whose elements are called thin
forest posets, is described here. We also define an involution on these
posets that is linked, as we shall see later, to Koszul duality of
the concerned operads.
\medbreak

\subsubsection{Description}
\begin{figure}[ht]
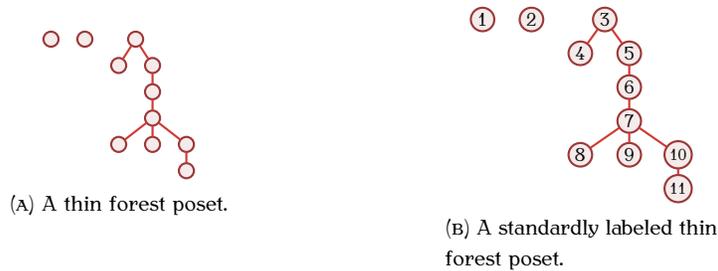

    \subfloat[][A thin forest poset.]{
    \begin{minipage}[c]{.4\textwidth}
    \centering

    \end{minipage}
    \label{subfig:thin_forest_poset}}
    \subfloat[][A standardly labeled thin forest poset.]{
    \begin{minipage}[c]{.4\textwidth}
    \centering
    %
    \end{minipage}
    \label{subfig:standard_thin_forest_poset}}
    \caption[A thin forest poset and a standard thin forest poset.]
    {Hasse diagrams of a thin forest poset and a standardly labeled
    version.}
    \label{fig:thin_forest_poset}
\end{figure}

\begin{figure}[ht]
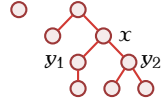

    \centering
    %
    \caption[A forest poset.]
    {The Hasse diagram of a forest poset which is not a thin
    forest poset. Indeed, the node $x$ has two children $y_1$ and $y_2$
    such that the suffix subtrees rooted at $y_1$ and $y_2$ have both
    two nodes or more.}
    \label{fig:not_thin_forest_poset}
\end{figure}

A \Def{thin forest poset} is a forest poset avoiding the pattern
$\LinePattern\; \LinePattern$ (see
Section~\ref{subsubsec:poset_patterns} of
Chapter~\ref{chap:combinatorics} for the definition of pattern avoidance
in posets). In other words, a thin forest poset is a poset so that the
nonplanar rooted tree $\Tfr$ obtained by adding a (new) root to the
Hasse diagram of $\Qca$ has the following property. Any node $x$ of
$\Tfr$ has at most one child $y$ such that the suffix subtree of $\Tfr$
rooted at $y$ has two nodes or more. For instance,
Figure~\ref{subfig:thin_forest_poset} shows a thin forest poset, while
Figure~\ref{fig:not_thin_forest_poset} shows a forest poset that does
not satisfies the described property.
\medbreak

A \Def{standard labeling} of a thin forest poset $\Qca$ consists in
labeling the vertices of the Hasse diagram of $\Qca$ from $1$ to
$\# \Qca$ in the order they appear in a depth first traversal, by always
visiting in a same sibling the node with the biggest subtree as last.
For instance, a standard labeling of the poset of
Figure~\ref{subfig:thin_forest_poset} is the poset shown in
Figure~\ref{subfig:standard_thin_forest_poset}. In what follows, we
shall consider only standardly labeled thin forest posets and we shall
identify any element $x$ of a thin forest posets $\Qca$ as the label
of $x$ in a standard labeling of $\Qca$. Moreover, we shall see thin
forest posets as forests of nonplanar rooted trees, obtained by
considering the Hasse diagrams of these posets.
\medbreak

Thin forest posets admit the following recursive description. If $\Qca$
is a thin forest poset, then $\Qca$ is the empty forest $\PosetEmpty$,
or it is the forest
\begin{equation}
    \OnePoset \; \Qca'
\end{equation}
consisting in the tree of one node $\OnePoset$ (labeled by $1$) and a
thin forest poset $\Qca'$, or it is the forest
\begin{equation}
    \begin{tikzpicture}[Centering,yscale=.6]
        \node[PosetVertex](1)at(0,0){};
        \node(2)at(0,-1){\begin{math}\Qca'\end{math}};
        \draw[Edge](1)--(2);
    \end{tikzpicture}
\end{equation}
consisting in one root (labeled by $1$) attached to the roots of the
trees of the thin forest poset~$\Qca'$. Therefore, there are $2^{n - 1}$
thin forest posets of size $n \geq 1$.
\medbreak

\subsubsection{Duality}
Given a thin forest poset $\Qca$, the \Def{dual} of $\Qca$ is the poset
$\Qca^\perp$ such that, for all $a, b \in \Qca^\perp$,
$a \Ord_{\Qca^\perp} b$ if and only if $a = b$ or $a$ and $b$ are
incomparable in $\Qca$ and $a < b$ (for the natural order on the labels
$a$ and $b$ that are integers). For instance, consider the poset
\begin{equation}
    \Qca :=
    \begin{tikzpicture}[Centering,xscale=.45,yscale=.42]
        \node[PosetVertex](1)at(0,2){\begin{math}1\end{math}};
        \node[PosetVertex](2)at(2,2){\begin{math}2\end{math}};
        \node[PosetVertex](3)at(1,1){\begin{math}3\end{math}};
        \node[PosetVertex](4)at(2,1){\begin{math}4\end{math}};
        \node[PosetVertex](5)at(3,1){\begin{math}5\end{math}};
        \node[PosetVertex](6)at(3,0){\begin{math}6\end{math}};
        \draw[Edge](2)--(3);
        \draw[Edge](2)--(4);
        \draw[Edge](2)--(5);
        \draw[Edge](5)--(6);
    \end{tikzpicture}\,.
\end{equation}
Since $1 \not \Ord_\Qca 2$, $1 \not \Ord_\Qca 3$, $1 \not \Ord_\Qca 4$,
$1 \not \Ord_\Qca 5$, $1 \not \Ord_\Qca 6$, $3 \not \Ord_\Qca 4$,
$3 \not \Ord_\Qca 5$, $3 \not \Ord_\Qca 6$, $4 \not \Ord_\Qca 5$,
$4 \not \Ord_\Qca 6$, in the dual $\Qca^\perp$ of $\Qca$ we have
$1 \Ord_{\Qca^\perp} 2$, $1 \Ord_{\Qca^\perp} 3$,
$1 \Ord_{\Qca^\perp} 4$, $1 \Ord_{\Qca^\perp} 5$,
$1 \Ord_{\Qca^\perp} 6$, $3 \Ord_{\Qca^\perp} 4$,
$3 \Ord_{\Qca^\perp} 5$, $3 \Ord_{\Qca^\perp} 6$,
$4 \Ord_{\Qca^\perp} 5$, $4 \Ord_{\Qca^\perp} 6$ and hence,
\begin{equation}
    \Qca^\perp =
    \begin{tikzpicture}[Centering,xscale=.45,yscale=.42]
        \node[PosetVertex](1)at(.25,3){\begin{math}1\end{math}};
        \node[PosetVertex](2)at(-.5,2){\begin{math}2\end{math}};
        \node[PosetVertex](3)at(1,2){\begin{math}3\end{math}};
        \node[PosetVertex](4)at(1,1){\begin{math}4\end{math}};
        \node[PosetVertex](5)at(0.5,0){\begin{math}5\end{math}};
        \node[PosetVertex](6)at(1.5,0){\begin{math}6\end{math}};
        \draw[Edge](1)--(2);
        \draw[Edge](1)--(3);
        \draw[Edge](3)--(4);
        \draw[Edge](4)--(5);
        \draw[Edge](4)--(6);
    \end{tikzpicture}\,.
\end{equation}
Observe that this operation $^\perp$ is an involution on thin forest
posets.
\medbreak

We now state two lemmas about thin forest posets and the
operation~$^\perp$.
\medbreak

\begin{Lemma} \label{lem:recursive_description_dual_poset}
    Let $\Qca$ be a thin forest poset. The dual $\Qca^\perp$ of $\Qca$
    admits the following recursive expression:
    \begin{enumerate}[label={\it (\roman*)}]
        \item if $\Qca$ is the empty forest $\PosetEmpty$, then
        \begin{equation}
            \PosetEmpty^\perp = \PosetEmpty;
        \end{equation}
        \item if $\Qca$ is of the form
        $\Qca = \OnePoset \; \Qca'$ where $\Qca'$ is
        a thin forest poset, then
        \begin{equation}
            \left(\OnePoset \; \Qca'\right)^\perp =
            \begin{tikzpicture}[Centering,yscale=.6]
                \node[PosetVertex](1)at(0,0){};
                \node(2)at(0,-1)
                    {\begin{math}{\Qca'}^\perp\end{math}};
                \draw[Edge](1)--(2);
            \end{tikzpicture};
        \end{equation}
        \item otherwise, $\Qca$ is of the form
        \begin{math}
            \Qca =
            \begin{tikzpicture}[Centering,yscale=.4]
                \node[PosetVertex](1)at(0,0){};
                \node(2)at(0,-1)
                    {\footnotesize \begin{math}\Qca'\end{math}};
                \draw[Edge](1)--(2);
            \end{tikzpicture}
        \end{math}
        where $\Qca'$ is a thin forest poset, and then
        \begin{equation}
            \left(
            \begin{tikzpicture}
                [Centering,yscale=.6]
                \node[PosetVertex](1)at(0,0){};
                \node(2)at(0,-1){\begin{math}\Qca'\end{math}};
                \draw[Edge](1)--(2);
            \end{tikzpicture}
            \right)^\perp
            = \OnePoset \; {\Qca'}^\perp.
        \end{equation}
    \end{enumerate}
\end{Lemma}
\medbreak

\begin{Lemma} \label{lem:relation_intervals_poset_dual}
    Let $\Qca$ be a thin forest poset. Then, the number of intervals of
    $\Qca$ and the number of intervals of its dual are related by
    \begin{equation}
        \NbInterv(\Qca) + \NbInterv\left(\Qca^\perp\right) =
        \frac{(\# \Qca)^2 + 3\, \# \Qca}{2}.
    \end{equation}
\end{Lemma}
\medbreak

\subsection{Koszul duality and poset duality}
By defining here an alternative basis for $\As(\Qca)$ when $\Qca$ is a
thin forest poset, we show that the construction $\As$ is closed under
Koszul duality on thin forest posets. More precisely, we show that
$\As(\Qca)^!$ and $\As(\Qca^\perp)$ are two isomorphic operads.
\medbreak

\subsubsection{Alternative basis}
Let $\Qca$ be a thin forest poset. For any element $b$ of $\Qca$, let
$\OpDualA_b$ be the element of $\As(\Qca)^!(2)$ defined by
\begin{equation} \label{equ:basis_change_dual_thin_forest_poset}
    \OpDualA_b :=
    \sum_{\substack{a \in \Qca^\perp \\ a \Ord_{\Qca^\perp} b}}
    \OpDualB_a.
\end{equation}
We denote by $\GeneratingSet_\Qca^{\OpDualA}$ the set of all
$\OpDualA_b$, $b \in \Qca$. By triangularity, the family
$\GeneratingSet_\Qca^{\OpDualA}$ forms a basis of
$\As(\Qca)^!(2)$ and hence, generates $\As(\Qca)^!$.
Consider for instance the thin forest poset
\begin{equation} \label{equ:thin_forest_poset_example}
    \Qca :=
    \begin{tikzpicture}[Centering,xscale=.45,yscale=.42]
        \node[PosetVertex](1)at(0,0){\begin{math}1\end{math}};
        \node[PosetVertex](2)at(1,0){\begin{math}2\end{math}};
        \node[PosetVertex](3)at(2.75,0){\begin{math}3\end{math}};
        \node[PosetVertex](4)at(2.25,-1){\begin{math}4\end{math}};
        \node[PosetVertex](5)at(3.25,-1){\begin{math}5\end{math}};
        \node[PosetVertex](6)at(3.25,-2){\begin{math}6\end{math}};
        \draw[Edge](3)--(4);
        \draw[Edge](3)--(5);
        \draw[Edge](5)--(6);
    \end{tikzpicture}\,.
\end{equation}
The dual poset of $\Qca$ is
\begin{equation} \label{equ:thin_forest_poset_dual_example}
    \Qca^\perp =
    \begin{tikzpicture}[Centering,xscale=.45,yscale=.42]
        \node[PosetVertex](1)at(0,0){\begin{math}1\end{math}};
        \node[PosetVertex](2)at(0,-1){\begin{math}2\end{math}};
        \node[PosetVertex](3)at(-1,-2){\begin{math}3\end{math}};
        \node[PosetVertex](4)at(1,-2){\begin{math}4\end{math}};
        \node[PosetVertex](5)at(.5,-3){\begin{math}5\end{math}};
        \node[PosetVertex](6)at(1.5,-3){\begin{math}6\end{math}};
        \draw[Edge](1)--(2);
        \draw[Edge](2)--(3);
        \draw[Edge](2)--(4);
        \draw[Edge](4)--(5);
        \draw[Edge](4)--(6);
    \end{tikzpicture}
\end{equation}
and hence, the elements of $\GeneratingSet_\Qca^{\OpDualA}$ express as
\begin{subequations}
\begin{equation}
    \OpDualA_1 = \OpDualB_1,
\end{equation}
\begin{equation}
    \OpDualA_2 = \OpDualB_1 + \OpDualB_2,
\end{equation}
\begin{equation}
    \OpDualA_3 = \OpDualB_1 + \OpDualB_2 + \OpDualB_3,
\end{equation}
\begin{equation}
    \OpDualA_4 = \OpDualB_1 + \OpDualB_2 + \OpDualB_4,
\end{equation}
\begin{equation}
    \OpDualA_5 = \OpDualB_1 + \OpDualB_2 + \OpDualB_4 + \OpDualB_5,
\end{equation}
\begin{equation}
    \OpDualA_6 = \OpDualB_1 + \OpDualB_2 + \OpDualB_4 + \OpDualB_6.
\end{equation}
\end{subequations}
\medbreak

\begin{Lemma} \label{lem:dimensions_relations_thin_forest_poset}
    Let $\Qca$ be a thin forest poset. Then, the dimension of the space
    $\RelationSpace_{\Qca^\perp}^\Op$ of relations of $\As(\Qca^\perp)$
    and the dimension of the space $\RelationSpace_\Qca^{\OpDual}$ of
    relations of $\As(\Qca)^!$ are related by
    \begin{equation} \label{equ:dimensions_relations_thin_forest_poset}
        \dim \RelationSpace_{\Qca^\perp}^\Op =
        4 \, \NbInterv\left(\Qca^\perp\right) - 3 \, \# \Qca =
        \dim \RelationSpace_\Qca^{\OpDual}.
    \end{equation}
\end{Lemma}
\medbreak

\subsubsection{Isomorphism}

\begin{Theorem} \label{thm:isomorphism_thin_forests_dual_operad}
    Let $\Qca$ be a thin forest poset. Then, the map
    $\phi : \As(\Qca^\perp) \to \As(\Qca)^!$ defined for any
    $a \in \Qca^\perp$ by $\phi(\Op_a) := \OpDualA_a$ extends in a
    unique way to an isomorphism of operads.
\end{Theorem}
\begin{proof}
    Let us denote by $\RelationSpace_\Qca^{\OpDualA}$ the space of
    relations of $\As(\Qca)^!$, expressed on the generating family
    $\GeneratingSet_\Qca^{\OpDualA}$. This space is the same as the
    space $\RelationSpace_\Qca^{\OpDualB}$, described by
    Proposition~%
    \ref{prop:alternative_presentation_koszul_dual_as_poset}.
    Let us exhibit a generating family of
    $\RelationSpace_\Qca^{\OpDualA}$ as a vector space. For this, let
    $a$ and $b$ be two elements of $\Qca^\perp$ such that
    $a \Ord_{\Qca^\perp} b$. We have, by
    using~\eqref{equ:basis_change_dual_thin_forest_poset},
    \begin{equation}\begin{split}
    \label{equ:isomorphism_thin_forests_dual_operad_computation_1}
        \OpDualA_a \circ_1 \OpDualA_b
                - \OpDualA_a \circ_2 \OpDualA_a & =
        \sum_{\substack{a', b' \in \Qca^\perp \\
            a' \Ord_{\Qca^\perp} a \\
            b' \Ord_{\Qca^\perp} b}}
        \left(\OpDualB_{a'} \circ_1 \OpDualB_{b'}\right)
        -
        \sum_{\substack{a', a'' \in \Qca^\perp \\
            a' \Ord_{\Qca^\perp} a \\
            a'' \Ord_{\Qca^\perp} a}}
        \left(\OpDualB_{a'} \circ_2 \OpDualB_{a''}\right) \\
        & =
        \sum_{\substack{a' \in \Qca^\perp \\
            a' \Ord_{\Qca^\perp} a}}
        \left(\OpDualB_{a'} \circ_1 \OpDualB_{a'}\right)
        -
        \sum_{\substack{a' \in \Qca^\perp \\
            a' \Ord_{\Qca^\perp} a}}
        \left(\OpDualB_{a'} \circ_2 \OpDualB_{a'}\right) \\
        & = 0.
    \end{split}\end{equation}
    Indeed the second equality
    of~\eqref{equ:isomorphism_thin_forests_dual_operad_computation_1}
    comes, by
    Proposition~%
    \ref{prop:alternative_presentation_koszul_dual_as_poset}, from the
    presence of the
    elements~\eqref{equ:alternative_relation_dual_as_poset_2}
    and~\eqref{equ:alternative_relation_dual_as_poset_3} in
    $\RelationSpace_\Qca^{\OpDualB}$, together with the fact that for
    all comparable elements $a'$ and $b'$ in $\Qca$, the fact that
    $a' \Ord_{\Qca^\perp} a$, $b' \Ord_{\Qca^\perp} b$, and
    $a \Ord_{\Qca^\perp} b$ implies that $a' = b'$. Besides, the last
    equality
    of~\eqref{equ:isomorphism_thin_forests_dual_operad_computation_1}
    comes, by
    Proposition~%
    \ref{prop:alternative_presentation_koszul_dual_as_poset}, from the
    presence of the
    elements~\eqref{equ:alternative_relation_dual_as_poset_1} in
    $\RelationSpace_\Qca^{\OpDualB}$. Similar arguments show that
    \begin{subequations}
    \begin{equation}
        \OpDualA_b \circ_1 \OpDualA_a
            - \OpDualA_a \circ_2 \OpDualA_a = 0,
    \end{equation}
    \begin{equation}
        \OpDualA_a \circ_1 \OpDualA_a
            - \OpDualA_b \circ_2 \OpDualA_a = 0,
    \end{equation}
    \begin{equation}
        \OpDualA_a \circ_1 \OpDualA_a
            - \OpDualA_a \circ_2 \OpDualA_b = 0.
    \end{equation}
    \end{subequations}
    \smallbreak

    We then have shown that the elements
    \begin{subequations}
    \begin{equation}
        \label{equ:isomorphism_thin_forests_dual_operad_generator_1}
        \OpDualA_a \circ_1 \OpDualA_b
        -
        \OpDualA_{a \Min_{\Qca^\perp} b} \circ_2
                \OpDualA_{a \Min_{\Qca^\perp} b},
        \qquad a, b \in {\Qca^\perp} \mbox{ and }
            (a \Ord_{\Qca^\perp} b \mbox{ or } b \Ord_{\Qca^\perp} a),
    \end{equation}
    \begin{equation}
        \label{equ:isomorphism_thin_forests_dual_operad_generator_2}
        \OpDualA_{a \Min_{\Qca^\perp} b} \circ_1
                \OpDualA_{a \Min_{\Qca^\perp} b}
        -
        \OpDualA_a \circ_2 \OpDualA_b,
        \qquad  a, b \in {\Qca^\perp} \mbox{ and }
            (a \Ord_{\Qca^\perp} b \mbox{ or } b \Ord_{\Qca^\perp} a).
    \end{equation}
    \end{subequations}
    are in $\RelationSpace_\Qca^{\OpDualB}$. It is immediate that the
    family consisting in the
    elements~%
    \eqref{equ:isomorphism_thin_forests_dual_operad_generator_1}
    and~\eqref{equ:isomorphism_thin_forests_dual_operad_generator_2} is
    free. We denote by $\RelationSpace$ the vector space generated by
    this family. By using the same arguments as the ones used in the
    proof of Proposition~\ref{prop:dimension_relations_as_poset}, we
    obtain that the dimension of $\RelationSpace$ is
    \begin{equation}
        \dim \RelationSpace =
        4 \, \NbInterv\left(\Qca^\perp\right) - 3 \, \# \Qca^\perp.
    \end{equation}
    Now, by Lemma~\ref{lem:dimensions_relations_thin_forest_poset},
    we deduce that
    \begin{math}
        \dim \RelationSpace =
        \dim \RelationSpace_\Qca^{\OpDual} =
        \dim \RelationSpace_\Qca^{\OpDualB},
    \end{math}
    implying that $\RelationSpace$ and $\RelationSpace_\Qca^{\OpDualB}$
    are equal.
    \smallbreak

    Therefore, the family of the $\GeneratingSet_\Qca^{\OpDualA}$
    generating $\As(\Qca)^!$ is subjected to the same relations as the
    family of the $\GeneratingSet_\Qca^\Op$ generating $\As(\Qca^\perp)$
    (compare~\eqref{equ:isomorphism_thin_forests_dual_operad_generator_1}
    with~\eqref{equ:relation_as_poset_1}
    and~\eqref{equ:isomorphism_thin_forests_dual_operad_generator_2}
    with~\eqref{equ:relation_as_poset_2}). Whence the statement of the
    theorem.
\end{proof}
\medbreak

The isomorphism $\phi$ between $\As(\Qca^\perp)$ and $\As(\Qca)^!$
provided by Theorem~\ref{thm:isomorphism_thin_forests_dual_operad} can
be expressed from the generating family
$\GeneratingSet_{\Qca^\perp}^\Op$ of $\As(\Qca^\perp)$ to the
generating family $\GeneratingSet_\Qca^{\OpDual}$ of $\As(\Qca)^!$,
for any $b \in \Qca^\perp$, as
\begin{equation}
    \phi(\Op_b) =
    \sum_{\substack{a \in \Qca^\perp \\ a \Ord_{\Qca^\perp} b}} \;
    \sum_{\substack{c \in \Qca \\ a \Ord_\Qca c}} \OpDual_c.
\end{equation}
\medbreak

For instance, by considering the pair of thin forest posets in duality
\begin{equation}
    \left(\Qca, \Qca^\perp\right) =
    \left(
    \begin{tikzpicture}[Centering,xscale=.45,yscale=.42]
        \node[PosetVertex](1)at(0,0){\begin{math}1\end{math}};
        \node[PosetVertex](2)at(1,0){\begin{math}2\end{math}};
        \node[PosetVertex](3)at(2.75,0){\begin{math}3\end{math}};
        \node[PosetVertex](4)at(2.25,-1){\begin{math}4\end{math}};
        \node[PosetVertex](5)at(3.25,-1){\begin{math}5\end{math}};
        \node[PosetVertex](6)at(3.25,-2){\begin{math}6\end{math}};
        \draw[Edge](3)--(4);
        \draw[Edge](3)--(5);
        \draw[Edge](5)--(6);
    \end{tikzpicture}\,, \enspace
    \begin{tikzpicture}[Centering,xscale=.45,yscale=.42]
        \node[PosetVertex](1)at(0,0){\begin{math}1\end{math}};
        \node[PosetVertex](2)at(0,-1){\begin{math}2\end{math}};
        \node[PosetVertex](3)at(-1,-2){\begin{math}3\end{math}};
        \node[PosetVertex](4)at(1,-2){\begin{math}4\end{math}};
        \node[PosetVertex](5)at(.5,-3){\begin{math}5\end{math}};
        \node[PosetVertex](6)at(1.5,-3){\begin{math}6\end{math}};
        \draw[Edge](1)--(2);
        \draw[Edge](2)--(3);
        \draw[Edge](2)--(4);
        \draw[Edge](4)--(5);
        \draw[Edge](4)--(6);
    \end{tikzpicture}\,\right),
\end{equation}
the map $\phi : \As(\Qca^\perp) \to \As(\Qca)^!$ defined in the
statement of Theorem~\ref{thm:isomorphism_thin_forests_dual_operad}
satisfies
\begin{subequations}
\begin{equation}
    \phi(\Op_1) = \OpDual_1,
\end{equation}
\begin{equation}
    \phi(\Op_2) = \OpDual_1 + \OpDual_2,
\end{equation}
\begin{equation}
    \phi(\Op_3) = \OpDual_1 + \OpDual_2 + \OpDual_3 + \OpDual_4
        + \OpDual_5 + \OpDual_6,
\end{equation}
\begin{equation}
    \phi(\Op_4) = \OpDual_1 + \OpDual_2 + \OpDual_4,
\end{equation}
\begin{equation}
    \phi(\Op_5) = \OpDual_1 + \OpDual_2 + \OpDual_4 + \OpDual_5
        + \OpDual_6,
\end{equation}
\begin{equation}
    \phi(\Op_6) = \OpDual_1 + \OpDual_2 + \OpDual_4 + \OpDual_6.
\end{equation}
\end{subequations}
\medbreak

Moreover, by considering the opposite pair of thin posets forests in
duality
\begin{equation}
    \left(\Qca, \Qca^\perp\right) =
    \left(
    \begin{tikzpicture}[Centering,xscale=.45,yscale=.42]
        \node[PosetVertex](1)at(0,0){\begin{math}1\end{math}};
        \node[PosetVertex](2)at(0,-1){\begin{math}2\end{math}};
        \node[PosetVertex](3)at(-1,-2){\begin{math}3\end{math}};
        \node[PosetVertex](4)at(1,-2){\begin{math}4\end{math}};
        \node[PosetVertex](5)at(.5,-3){\begin{math}5\end{math}};
        \node[PosetVertex](6)at(1.5,-3){\begin{math}6\end{math}};
        \draw[Edge](1)--(2);
        \draw[Edge](2)--(3);
        \draw[Edge](2)--(4);
        \draw[Edge](4)--(5);
        \draw[Edge](4)--(6);
    \end{tikzpicture}\,, \enspace
    \begin{tikzpicture}[Centering,xscale=.45,yscale=.42]
        \node[PosetVertex](1)at(0,0){\begin{math}1\end{math}};
        \node[PosetVertex](2)at(1,0){\begin{math}2\end{math}};
        \node[PosetVertex](3)at(2.75,0){\begin{math}3\end{math}};
        \node[PosetVertex](4)at(2.25,-1){\begin{math}4\end{math}};
        \node[PosetVertex](5)at(3.25,-1){\begin{math}5\end{math}};
        \node[PosetVertex](6)at(3.25,-2){\begin{math}6\end{math}};
        \draw[Edge](3)--(4);
        \draw[Edge](3)--(5);
        \draw[Edge](5)--(6);
    \end{tikzpicture}\,\right),
\end{equation}
the map $\phi : \As(\Qca^\perp) \to \As(\Qca)^!$ defined in the
statement of Theorem~\ref{thm:isomorphism_thin_forests_dual_operad}
satisfies
\begin{subequations}
\begin{equation}
    \phi(\Op_1) = \OpDual_1 + \OpDual_2 + \OpDual_3 + \OpDual_4
        + \OpDual_5 + \OpDual_6,
\end{equation}
\begin{equation}
    \phi(\Op_2) = \OpDual_2 + \OpDual_3 + \OpDual_4
        + \OpDual_5 + \OpDual_6,
\end{equation}
\begin{equation}
    \phi(\Op_3) = \OpDual_3,
\end{equation}
\begin{equation}
    \phi(\Op_4) = \OpDual_3 + \OpDual_4 + \OpDual_5 + \OpDual_6,
\end{equation}
\begin{equation}
    \phi(\Op_5) = \OpDual_3 + \OpDual_5,
\end{equation}
\begin{equation}
    \phi(\Op_6) = \OpDual_3 + \OpDual_5 + \OpDual_6.
\end{equation}
\end{subequations}
\medbreak

Notice also that since the dual of the total order $\Qca$ on a set of
$\ell \geq 0$ elements is the trivial order $\Qca^\perp$ on the same
set, by Theorem~\ref{thm:isomorphism_thin_forests_dual_operad},
$\As(\Qca)$ is the Koszul dual of $\As(\Qca^\perp)$. This is coherent
with the results of Section~\ref{sec:as_gamma} of
Chapter~\ref{chap:polydendr} about the multiassociative operad
(equal to $\As(\Qca)$) and the dual multiassociative operad (equal
to~$\As(\Qca^\perp)$).
\medbreak

\section*{Concluding remarks}
Through this chapter, we have presented a functorial construction $\As$
from posets to operads establishing a link between the two underlying
categories. The operads obtained through this construction generalize
the (dual) multiassociative operads. As we
have seen, some combinatorial properties of the starting posets $\Qca$
imply properties on the obtained operads $\As(\Qca)$ as, among others,
basicity and Koszulity.
\smallbreak

This work raises several questions. We have presented two classes of
$\Qca$-associative algebras: the free $\Qca$-associative algebras over
one generator where $\Qca$ are forest posets and a polynomial algebra
involving the antichains of a poset $\Qca$. The question to characterize
free $\Qca$-associative algebras over one generator with no assumption
on $\Qca$ is open. Also, the question to define some other interesting
$\Qca$-associative algebras has not been considered in this work and
deserves to be addressed.
\smallbreak

Besides, we have shown that when $\Qca$ is a forest poset, $\As(\Qca)$
is Koszul. The property of being a forest poset for $\Qca$ is only a
sufficient condition for the Koszulity of $\As(\Qca)$ and the question
to find a necessary condition is worthwhile. Notice that the strategy to
prove the Koszulity of an operad by the partition poset
method~\cite{MY91,Val07} (see also~\cite{LV12}) cannot be applied to our
context. Indeed, this strategy applies only to basic operads and we have
shown that almost all operads $\As(\Qca)$ are not basic.
\medbreak


\chapter{Operads of decorated cliques} \label{chap:cliques}
The content of this chapter comes from~\cite{Gir17a,Gir17b}.
\medbreak

\section*{Introduction}
Regular polygons endowed with configurations of diagonals are very
classical combinatorial objects. Up to some restrictions or enrichments,
sets defined on these polygons can be put in bijection with several
combinatorial families. For instance, it is well-known that
triangulations~\cite{DRS10}, forming a particular subset of the set
of all polygons, are in one-to-one correspondence with binary trees,
and a lot of structures and operations on binary trees translate nicely
on triangulations. Indeed, among others, the rotation operation on
binary trees~\cite{Knu98} is the covering relation of the Tamari
order~\cite{Tam62,HT72} (see also
Section~\ref{subsubsec:examples_posets} of
Chapter~\ref{chap:combinatorics}) and this operation translates as a
diagonal flip in triangulations. Also, noncrossing
configurations~\cite{FN99} form another interesting subfamily of such
polygons. Natural generalizations of noncrossing configurations consist
in allowing, with more or less restrictions, some crossing diagonals.
One of these families is formed by the multi-triangulations~\cite{CP92},
that are polygons wherein the number of mutually crossing diagonal is
bounded. Besides, let us emphasize that the class of combinatorial
objects in bijection with sets of polygons with configurations of
diagonals is large enough in order to contain, among others, dissections
of polygons, noncrossing partitions, permutations, and involutions.
\smallbreak

The purpose of this work is twofold. First, we are concerned in endowing
the linear span of the polygons with configurations of arcs with a
structure of an operad. This is justified by the preliminary observation
that most of the subfamilies of polygons endowed with configurations
of diagonals discussed above are stable for several natural composition
operations. Even better, some of these can be described as the closure
with respect to these composition operations of small sets  of polygons.
For this reason, operads are very promising candidates, among the
modern algebraic structures, to study such objects under an algebraic
and combinatorial flavor. This leads to see these objects under a new
light, stressing some of their combinatorial and algebraic properties.
Second, we would provide a general construction of operads of polygons
rich enough so that it includes some already known operads. As a
consequence, we obtain alternative definitions of existing operads and
new interpretations of these.
\smallbreak

For this aim, we work here with $\Mca$-decorated cliques (or
$\Mca$-cliques for short), that are complete graphs whose arcs are
labeled on $\Mca$, where $\Mca$ is a unitary magma. These objects are
natural generalizations of polygons with configurations of arcs since
the arcs of any $\Mca$-clique labeled by the unit of $\Mca$ are
considered as missing. The elements of $\Mca$ different from the unit
allow moreover to handle polygons with arcs of different colors. For
instance, each usual noncrossing configuration $\Cfr$ can be encoded
by an $\N_2$-clique $\Pfr$, where $\N_2$ is the cyclic additive unitary
magma $\Z/_{2\Z}$, wherein each arc labeled by $1 \in \N_2$ in
$\Pfr$ denotes the presence of the same arc in $\Cfr$, and each arc
labeled by $0 \in \N_2$ in $\Pfr$ denotes its absence in~$\Cfr$. Our
construction is materialized by a functor $\Cli$ from the category of
unitary magmas to the category of operads. It builds, from any unitary
magma $\Mca$, an operad $\Cli\Mca$ on $\Mca$-cliques. The partial
composition $\Pfr \circ_i \Qfr$ of two $\Mca$-cliques $\Pfr$ and $\Qfr$
of $\Cli\Mca$ consists in gluing the $i$th edge of $\Pfr$ (with respect
to a precise indexation) and a special arc of $\Qfr$, called the base,
together to form a new $\Mca$-clique. The magmatic operation of $\Mca$
explains how to relabel the two overlapping arcs.
\smallbreak

This operad $\Cli\Mca$ has a lot of properties, which can be apprehended
both under a combinatorial and an algebraic point of view. First, many
families of particular polygons with configurations of arcs form
quotients or suboperads of $\Cli\Mca$. We can for instance control the
degrees of the vertices or the crossings between diagonals to obtain new
operads. We can also forbid all diagonals, or some labels for the
diagonals or the edges, or all nestings of diagonals, or even all
cycles formed by arcs. All these combinatorial particularities and
restrictions on $\Mca$-cliques behave well algebraically. Moreover, by
using the fact that the direct sum of two ideals of an operad $\Oca$ is
still an ideal of $\Oca$, these constructions can be mixed to get even
more operads. For instance, it is well-known that Motzkin
configurations, that are polygons with disjoint noncrossing diagonals,
are enumerated by Motzkin numbers~\cite{Mot48}. Since a Motzkin
configuration can be encoded by an $\Mca$-clique where all vertices
are of degrees at most $1$ and no diagonal crosses another one, we
obtain an operad $\Motzkin\Mca$ on colored Motzkin configurations
which is both a quotient of $\Deg_1\Mca$, the quotient of $\Cli\Mca$
consisting in all $\Mca$-cliques such that all vertices are of degrees
at most $1$, and of $\NC\Mca$, the quotient (and suboperad) of
$\Cli\Mca$ consisting in all noncrossing $\Mca$-cliques. We also get
quotients of $\Cli\Mca$ involving, among others, Schröder trees,
forests of paths, forests of trees, dissections of polygons, Lucas
configurations, with colored versions for each of these. This leads to
a new hierarchy of operads, wherein links between its components
appear as surjective or injective operad morphisms.
Table~\ref{tab:constructed_operads_C_M} lists the main operads
constructed in this work and gathers some information about these.
\begin{table}[ht]
    \centering
    \begin{tabular}{c|c|c|c}
        Operad & Objects & Status with respect to $\Cli\Mca$
            & Place \\ \hline \hline
        $\Cli\Mca$ & $\Mca$-cliques & --- &
        Section~\ref{sec:construction_Cli} \\ \hline
        $\Lab_{B, E, D}\Mca$ & $\Mca$-cliques with restricted labels
            & Suboperad
            & Section~\ref{subsubsec:suboperad_Cli_M_labels}\\
        $\Whi\Mca$ & White $\Mca$-cliques & Suboperad
            & Section~\ref{subsubsec:suboperad_Cli_M_white} \\
        $\Cro_k\Mca$ & $\Mca$-cliques of crossings at most $k$
            & Suboperad and quotient
            & Section~\ref{subsubsec:quotient_Cli_M_crossings} \\
        $\Bub\Mca$ & $\Mca$-bubbles
            & Quotient
            & Section~\ref{subsubsec:quotient_Cli_M_bubbles} \\
        $\Deg_k\Mca$ & $\Mca$-cliques of degrees at most $k$
            & Quotient
            & Section~\ref{subsubsec:quotient_Cli_M_degrees} \\
        $\Nes\Mca$ & Nesting-free $\Mca$-cliques
            & Quotient
            & Section~\ref{subsubsec:quotient_Cli_M_Nes} \\
        $\Acy\Mca$ & Acyclic $\Mca$-cliques
            & Quotient
            & Section~\ref{subsubsec:quotient_Cli_M_acyclic} \\ \hline
        $\NC\Mca$ & noncrossing $\Mca$-cliques & Suboperad and quotient
            & Section~\ref{sec:operad_noncrossing} \\
    \end{tabular}
    \bigbreak

    \caption[Data about operads related to the ns operad~$\Cli\Mca$.]
    {The main operads defined in this work. All these operads depend
    on a unitary magma $\Mca$ which has, in some cases, to satisfy
    some precise conditions. Some of these operads depend also on a
    nonnegative integer~$k$ or subsets $B$, $E$, and $D$ of~$\Mca$.}
    \label{tab:constructed_operads_C_M}
\end{table}
\smallbreak

One of the most notable of these substructures is built by considering
the $\Dbb_0$-cliques that have vertices of degrees at most $1$, where
$\Dbb_0$ is the multiplicative unitary magma on $\{0, 1\}$. This is in
fact the quotient $\Deg_1\Dbb_0$ of $\Cli\Dbb_0$ and involves
involutions (or equivalently, standard Young tableaux by the
Robinson-Schensted correspondence~\cite{Sch61,Lot02}). To the best of
our knowledge, $\Deg_1\Dbb_0$ is the first nontrivial operad on these
objects. As an important remark at this stage, let us highlight that
when $\Mca$ is nontrivial, $\Cli\Mca$ is not a binary operad. Indeed,
all its minimal generating sets are infinite and its generators have
arbitrary high arities. Nevertheless, the biggest binary suboperad of
$\Cli\Mca$ is the operad $\NC\Mca$ of noncrossing configurations and
this operad is quadratic and Koszul, regardless of $\Mca$. Furthermore,
the construction $\Cli$ maintains some links with the operad $\RatFct$
of rational functions introduced by Loday~\cite{Lod10} (see also
Section~\ref{subsubsec:rational_functions_operad} of
Chapter~\ref{chap:algebra}). In fact, provided that $\Mca$ satisfies
some conditions, each $\Mca$-clique encodes a rational function. This
defines an operad morphism from $\Cli\Mca$ to $\RatFct$. Moreover, the
construction $\Cli$ allows to construct already known operads in
original ways. For instance, for well-chosen unitary magmas $\Mca$, the
operads $\Cli\Mca$ contain $\FF_4$, a suboperad of the operad of formal
fractions $\FF$~\cite{CHN16}, the operad $\NCT$ of based noncrossing
trees~\cite{Cha07,Ler11}, and $\MT$ and $\DMT$, two operads respectively
defined in~\cite{LMN13} and in Chapter~\ref{chap:languages} that involve
multi-tildes and double multi-tildes, which are operators coming from
formal language theory~\cite{CCM11}. Moreover, $\Cli$ provides a
construction of $\BNC$, the operad of bicolored noncrossing
configurations (see Chapter~\ref{chap:enveloping}). For this reason, in
particular, all the suboperads of $\BNC$ can be obtained from the
construction $\Cli$. This includes for example the dipterous
operad~\cite{LR03,Zin12}. The operads $\Cli\Mca$ also contains $\Grav$,
the gravity operad, a symmetric operad introduced by
Getzler~\cite{Get94}, seen here as a nonsymmetric one~\cite{AP15}.
\smallbreak

This chapter is organized as follows. In
Section~\ref{sec:construction_Cli}, we introduce $\Mca$-cliques, the
construction $\Cli$, and study some of its properties. Then,
Section~\ref{sec:quotients_suboperads_C_M} is devoted to define several
suboperads and quotients of $\Cli\Mca$. This leads to a bunch of new
operads on particular $\Mca$-cliques. We focus next, in
Section~\ref{sec:operad_noncrossing}, on the study of the suboperad
$\NC\Mca$ of $\Cli\Mca$ on the noncrossing $\Mca$-cliques. Among
others, we provide a presentation by generators and relations of
$\NC\Mca$ and of its Koszul dual. Finally, in
Section~\ref{sec:concrete_constructions}, we use the construction
$\Cli$ to provide alternative definitions of some known operads.
\medbreak

\subsubsection*{Note}
This chapter deals onlys with ns operads. For this reason, ``operad''
means ``ns operad''.
\medbreak

\section{From unitary magmas to operads}\label{sec:construction_Cli}
We describe in this section our construction from unitary magmas to
operads and study its main algebraic and combinatorial properties.
\medbreak

\subsection{Unitary magmas, decorated cliques, and operads}
\label{subsec:decorated_cliques}
We present here our main combinatorial objects, the decorated cliques.
The construction $\Cli$, which takes a unitary magma as input and
produces an operad, is defined.
\medbreak

\subsubsection{Unitary magmas}
Recall first that a unitary magma is a set endowed with a binary
operation $\Op$ admitting a left and right unit $\Unit_\Mca$. For
convenience, we denote by $\bar{\Mca}$ the set
$\Mca \setminus \{\Unit_\Mca\}$. To explore some examples in this
chapter, we shall mostly consider four sorts of unitary magmas: the
additive unitary magma on all integers denoted by $\Z$, the cyclic
additive unitary magma on $\Z/_{\ell \Z}$ denoted by $\N_\ell$, the
unitary magma
\begin{equation}
    \Dbb_\ell := \{\Unit, 0, \Dsf_1, \dots, \Dsf_\ell\}
\end{equation}
where $\Unit$ is the unit of $\Dbb_\ell$, $0$ is absorbing, and
$\Dsf_i \Op \Dsf_j = 0$ for all $i, j \in [\ell]$, and the unitary magma
\begin{equation}
    \Ebb_\ell := \{\Unit, \Esf_1, \dots, \Esf_\ell\}
\end{equation}
where $\Unit$ is the unit of $\Ebb_\ell$ and $\Esf_i \Op \Esf_j = \Unit$
for all $i, j \in [\ell]$. Observe that since
\begin{equation}
    \Esf_1 \Op (\Esf_1 \Op \Esf_2) = \Esf_1 \Op \Unit = \Esf_1
    \ne
    \Esf_2 = \Unit \Op \Esf_2 = (\Esf_1 \Op \Esf_1) \Op \Esf_2,
\end{equation}
all unitary magmas $\Ebb_\ell$, $\ell \geq 2$, are not monoids.
\medbreak

\subsubsection{Decorated cliques}
An \Def{$\Mca$-decorated clique} (or an \Def{$\Mca$-clique} for short)
is an $\Mca$-configuration $\Pfr$ (see
Section~\ref{subsubsec:configurations} of
Chapter~\ref{chap:combinatorics}) such that each arc of $\Pfr$ has a
label. When the arc $(x, y)$ of $\Pfr$ is labeled by an element
different from $\Unit_\Mca$, we say that the arc $(x, y)$ is
\Def{solid}. By convention, we require that the $\Mca$-clique
$\UnitClique$ of size $1$ having its base labeled by $\Unit_\Mca$ is the
only such object of size $1$. The set of all $\Mca$-cliques is denoted
by~$\Cliques_\Mca$.
\medbreak

In our graphical representations, we shall represent any $\Mca$-clique
$\Pfr$ by following the drawing conventions of configurations explained
in Section~\ref{subsubsec:configurations} of
Chapter~\ref{chap:combinatorics} with the difference that non-solid
diagonals are not drawn. For instance,
\begin{equation}
    \Pfr :=
    \begin{tikzpicture}[scale=.85,Centering]
        \node[CliquePoint](1)at(-0.43,-0.90){};
        \node[CliquePoint](2)at(-0.97,-0.22){};
        \node[CliquePoint](3)at(-0.78,0.62){};
        \node[CliquePoint](4)at(-0.00,1.00){};
        \node[CliquePoint](5)at(0.78,0.62){};
        \node[CliquePoint](6)at(0.97,-0.22){};
        \node[CliquePoint](7)at(0.43,-0.90){};
        \draw[CliqueEdge](1)edge[]node[CliqueLabel]
            {\begin{math}-1\end{math}}(2);
        \draw[CliqueEdge](1)
            edge[bend left=30]node[CliqueLabel,near start]
            {\begin{math}2\end{math}}(5);
        \draw[CliqueEdge](1)edge[]node[CliqueLabel]
            {\begin{math}1\end{math}}(7);
        \draw[CliqueEmptyEdge](2)edge[]node[CliqueLabel]{}(3);
        \draw[CliqueEmptyEdge](3)edge[]node[CliqueLabel]{}(4);
        \draw[CliqueEdge](3)
            edge[bend left=30]node[CliqueLabel,near start]
            {\begin{math}-1\end{math}}(7);
        \draw[CliqueEdge](4)edge[]node[CliqueLabel]
            {\begin{math}3\end{math}}(5);
        \draw[CliqueEdge](5)edge[]node[CliqueLabel]
            {\begin{math}2\end{math}}(6);
        \draw[CliqueEdge](5)edge[]node[CliqueLabel]
            {\begin{math}1\end{math}}(7);
        \draw[CliqueEmptyEdge](6)edge[]node[CliqueLabel]{}(7);
    \end{tikzpicture}
\end{equation}
is a $\Z$-clique such that, among others $\Pfr(1, 2) = -1$,
$\Pfr(1, 5) = 2$, $\Pfr(3, 7) = -1$, $\Pfr(5, 7) = 1$, $\Pfr(2, 3) = 0$
(because $0$ is the unit of $\Z$), and $\Pfr(2, 6) = 0$ (for the same
reason).
\medbreak

Let us now provide some definitions and statistics on $\Mca$-cliques.
The \Def{underlying configuration} of $\Pfr$ is the
$\bar{\Mca}$-configuration $\bar{\Pfr}$ of the same size as the one of
$\Pfr$ and such $\bar{\Pfr}(x, y) := \Pfr(x, y)$ for all solid arcs
$(x, y)$ of $\Pfr$, and all other arcs of $\bar{\Pfr}$ are unlabeled.
The \Def{skeleton}, (resp. \Def{degree}, \Def{crossing}) of $\Pfr$ is
the skeleton (resp. the degree, the crossing) of $\bar{\Pfr}$. Moreover,
$\Pfr$ is \Def{nesting-free}, (resp. \Def{acyclic}, \Def{white}, an
\Def{$\Mca$-bubble}, an \Def{$\Mca$-triangle}), if $\bar{\Pfr}$ is
nesting-free (resp. acyclic, white, a bubble, a triangle). The set of
all $\Mca$-bubbles (resp. $\Mca$-triangles) is denoted by
$\Bubbles_\Mca$ (resp. $\Triangles_\Mca$). Finally, the \Def{border} of
$\Pfr$ is the word $\Border(\Pfr)$ of length $n$ such that for any
$i \in [n]$, $\Border(\Pfr)_i = \Pfr_i$.
\medbreak

\subsubsection{Partial composition of $\Mca$-cliques}
\begin{figure}[ht]
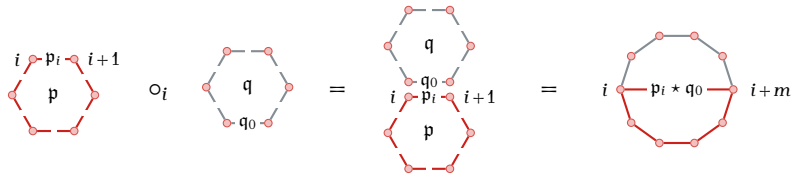

    \centering
    \begin{equation*}

    \end{equation*}
    \caption[The partial composition of $\Mca$-cliques.]
    {The partial composition of $\Cli\Mca$, described in geometric
    terms. Here, $\Pfr$ and $\Qfr$ are two $\Mca$-cliques. The arity
    of $\Qfr$ is~$m$ and $i$ is an integer between $1$ and~$|\Pfr|$.}
    \label{fig:composition_Cli_M}
\end{figure}

From now on, the \Def{arity} of an $\Mca$-clique $\Pfr$ is its size and
is denoted by $|\Pfr|$. For any unitary magma $\Mca$, we define the
vector space
\begin{equation}
    \Cli\Mca := \bigoplus_{n \geq 1} \Cli\Mca(n),
\end{equation}
where $\Cli\Mca(n)$ is the linear span of all $\Mca$-cliques of arity
$n$, $n \geq 1$. The set $\Cliques_\Mca$ forms hence a basis of
$\Cli\Mca$ called \Def{fundamental basis}. Observe that the space
$\Cli\Mca(1)$ has dimension $1$ since it is the linear span of the
$\Mca$-clique $\UnitClique$. We endow $\Cli\Mca$ with partial
composition maps
\begin{equation}
    \circ_i : \Cli\Mca(n) \otimes \Cli\Mca(m) \to
    \Cli\Mca(n + m - 1),
    \qquad n, m \geq 1, i \in [n],
\end{equation}
defined linearly, in the fundamental basis, in the following way. Let
$\Pfr$ and $\Qfr$ be two $\Mca$-cliques of respective arities $n$ and
$m$, and $i \in [n]$ be an integer. We set $\Pfr \circ_i \Qfr$ as the
$\Mca$-clique of arity $n + m - 1$ such that, for any arc $(x, y)$ where
$1 \leq x < y \leq n + m$,
\begin{equation} \label{equ:partial_composition_Cli_M}
    (\Pfr \circ_i \Qfr)(x, y) :=
    \begin{cases}
        \Pfr(x, y)
            & \mbox{if } y \leq i, \\
        \Pfr(x, y - m + 1)
            & \mbox{if } x \leq i < i + m \leq y
            \mbox{ and } (x, y) \ne (i, i + m), \\
        \Pfr(x - m + 1, y - m + 1)
            & \mbox{if } i + m \leq x, \\
        \Qfr(x - i + 1, y - i + 1)
            & \mbox{if } i \leq x < y \leq i + m
              \mbox{ and } (x, y) \ne (i, i + m), \\
        \Pfr_i \Op \Qfr_0
            & \mbox{if } (x, y) = (i, i + m), \\
        \Unit_\Mca
            & \mbox{otherwise}.
    \end{cases}
\end{equation}
We recall that $\Op$ denotes the operation of $\Mca$ and $\Unit_\Mca$
its unit. In a geometric way, $\Pfr \circ_i \Qfr$ is obtained by gluing
the base of $\Qfr$ onto the $i$th edge of $\Pfr$, by relabeling the
common arcs between $\Pfr$ and $\Qfr$, respectively the arcs
$(i, i + 1)$ and $(1, m + 1)$, by $\Pfr_i \Op \Qfr_0$, and by adding all
required non solid diagonals on the graph thus obtained to become a
clique (see Figure~\ref{fig:composition_Cli_M}).
For example, in $\Cli\Z$, one has the two partial compositions
\begin{subequations}
\begin{equation}
    %
\,.
\end{equation}
\end{subequations}
\medbreak

\subsubsection{Functorial construction from unitary magmas to operads}
If $\Mca_1$ and $\Mca_2$ are two unitary magmas and
$\theta : \Mca_1 \to \Mca_2$ is a unitary magma morphism, we define
\begin{equation}
    \Cli\theta : \Cli\Mca_1 \to \Cli\Mca_2
\end{equation}
as the linear map sending any $\Mca_1$-clique $\Pfr$ of arity $n$ to the
$\Mca_2$-clique $(\Cli\theta)(\Pfr)$ of the same arity such that, for
any arc $(x, y)$ where $1 \leq x < y \leq  n + 1$,
\begin{equation} \label{equ:morphism_Cli_M}
    ((\Cli\theta)(\Pfr))(x, y) := \theta(\Pfr(x, y)).
\end{equation}
In a geometric way, $(\Cli\theta)(\Pfr)$ is the $\Mca_2$-clique obtained
by relabeling each arc of $\Pfr$ by the image of its label by~$\theta$.
\medbreak

\begin{Theorem} \label{thm:clique_construction}
    The construction $\Cli$ is a functor from the category of unitary
    magmas to the category of operads. Moreover, $\Cli$ respects
    injections and surjections.
\end{Theorem}
\medbreak

We name the construction $\Cli$ as the \Def{clique construction} and
$\Cli\Mca$ as the \Def{$\Mca$-clique operad}. Observe that the
fundamental basis of $\Cli\Mca$ is a set-operad basis of $\Cli\Mca$.
Besides, when $\Mca$ is the trivial unitary magma $\{\Unit_\Mca\}$,
$\Cli\Mca$ is the linear span of all decorated cliques having only
non-solid arcs. Thus, each space $\Cli\Mca(n)$, $n \geq 1$, is of
dimension $1$ and it follows from the definition of the partial
composition of $\Cli\Mca$ that this operad is isomorphic to the
associative operad $\As$. The next result shows that the clique
construction is compatible with the Cartesian product of unitary magmas.
\medbreak

\begin{Proposition} \label{prop:Cli_M_Cartesian_product}
    Let $\Mca_1$ and $\Mca_2$ be two unitary magmas. Then, the operads
    $(\Cli\Mca_1) \Hadamard (\Cli\Mca_2)$ and
    $\Cli(\Mca_1 \times \Mca_2)$ are isomorphic.
\end{Proposition}
\medbreak

\subsection{General properties}
We investigate here some properties of clique operads, as their
dimensions, their minimal generating sets, the fact that they admit a
cyclic operad structure, and describe their partial compositions over
two alternative bases.
\medbreak

\subsubsection{Dimensions and minimal generating set}
\label{subsubsec:dimensions_clique_operads}

\begin{Proposition} \label{prop:dimensions_Cli_M}
    Let $\Mca$ be a finite unitary magma. For all $n \geq 2$,
    \begin{equation} \label{equ:dimensions_Cli_M}
        \dim \Cli\Mca(n) = m^{\binom{n + 1}{2}},
    \end{equation}
    where $m := \# \Mca$.
\end{Proposition}
\medbreak

From Proposition~\ref{prop:dimensions_Cli_M}, the first dimensions of
$\Cli\Mca$ depending on $m := \# \Mca$ are
\begin{subequations}
\begin{equation}
    1, 1, 1, 1, 1, 1, 1, 1,
    \qquad m = 1,
\end{equation}
\begin{equation}
    1, 8, 64, 1024, 32768, 2097152, 268435456, 68719476736,
    \qquad m = 2,
\end{equation}
\begin{multline}
    1, 27, 729, 59049, 14348907, 10460353203, 22876792454961, \\
    150094635296999121,
    \qquad m = 3,
\end{multline}
\begin{multline}
    1, 64, 4096, 1048576, 1073741824, 4398046511104,
    72057594037927936, \\
    4722366482869645213696,
    \qquad m = 4.
\end{multline}
\end{subequations}
Except for the first terms, the second one forms
Sequence~\OEIS{A006125}, the third one forms Sequence~\OEIS{A047656},
and the last one forms Sequence~\OEIS{A053763} of~\cite{Slo}.
\medbreak

Let $\Primes_\Mca$ be the set of all $\Mca$-cliques $\Pfr$ or arity
$n \geq 2$ such that, for any (non-necessarily solid) diagonal $(x, y)$
of $\Pfr$, there is at least one solid diagonal $(x', y')$ of $\Pfr$
such that $(x, y)$ and $(x', y')$ are crossing. We call $\Primes_\Mca$
the set of all \Def{prime $\Mca$-cliques}. Observe that, according to
this description, all $\Mca$-triangles are prime.
\medbreak

\begin{Proposition} \label{prop:generating_set_Cli_M}
    Let $\Mca$ be a unitary magma. The set $\Primes_\Mca$ is a minimal
    generating set of~$\Cli\Mca$.
\end{Proposition}
\medbreak

Computer experiments tell us that, when $m := \# \Mca = 2$, the first
numbers of prime $\Mca$-cliques are, size by size,
\begin{equation} \label{equ:prime_cliques_numbers_2}
    0, 8, 16, 352, 16448, 1380224.
\end{equation}
Moreover, remark that each $n$th term of this sequence is divisible
by $m^{n + 1}$ since the labels of the base and the edges of an
$\Mca$-clique $\Pfr$ have no influence on the fact that $\Pfr$ is prime.
This gives the sequence
\begin{equation} \label{equ:white_prime_cliques_numbers}
    0, 1, 1, 11, 257, 10783.
\end{equation}
None of these sequences appear in~\cite{Slo} at this time.
\medbreak

\subsubsection{Associative elements}

\begin{Proposition} \label{prop:associative_elements_Cli_M}
    Let $\Mca$ be a unitary magma and $f$ be an element of $\Cli\Mca(2)$
    of the form
    \begin{equation} \label{equ:associative_elements_Cli_M_0}
        f :=
        \sum_{\Pfr \in \Triangles_\Mca}
        \lambda_\Pfr \Pfr,
    \end{equation}
    where the $\lambda_\Pfr$, $\Pfr \in \Triangles_\Mca$, are
    coefficients of $\K$. Then, $f$ is associative if and only if
    \begin{subequations}
    \begin{equation} \label{equ:associative_elements_Cli_M_1}
        \sum_{\substack{
            \Pfr_1, \Qfr_0 \in \Mca \\
            \delta =  \Pfr_1 \Op \Qfr_0
        }}
        \lambda_{\Triangle{\Pfr_0}{\Pfr_1}{\Pfr_2}}
        \lambda_{\Triangle{\Qfr_0}{\Qfr_1}{\Qfr_2}}
        = 0,
        \qquad
        \Pfr_0, \Pfr_2, \Qfr_1, \Qfr_2 \in \Mca,
        \delta \in \bar{\Mca},
    \end{equation}
    \begin{equation} \label{equ:associative_elements_Cli_M_2}
        \sum_{\substack{
            \Pfr_1, \Qfr_0 \in \Mca \\
            \Pfr_1 \Op \Qfr_0 = \Unit_\Mca
        }}
        \lambda_{\Triangle{\Pfr_0}{\Pfr_1}{\Pfr_2}}
        \lambda_{\Triangle{\Qfr_0}{\Qfr_1}{\Qfr_2}}
        -
        \lambda_{\Triangle{\Pfr_0}{\Qfr_1}{\Pfr_1}}
        \lambda_{\Triangle{\Qfr_0}{\Qfr_2}{\Pfr_2}}
        = 0,
        \qquad
        \Pfr_0, \Pfr_2, \Qfr_1, \Qfr_2 \in \Mca,
    \end{equation}
    \begin{equation} \label{equ:associative_elements_Cli_M_3}
        \sum_{\substack{
            \Pfr_2, \Qfr_0 \in \Mca \\
            \delta =  \Pfr_2 \Op \Qfr_0
        }}
        \lambda_{\Triangle{\Pfr_0}{\Pfr_1}{\Pfr_2}}
        \lambda_{\Triangle{\Qfr_0}{\Qfr_1}{\Qfr_2}}
        = 0,
        \qquad
        \Pfr_0, \Pfr_1, \Qfr_1, \Qfr_2 \in \Mca,
        \delta \in \bar{\Mca}.
    \end{equation}
    \end{subequations}
\end{Proposition}
\medbreak

For instance, by Proposition~\ref{prop:associative_elements_Cli_M}, the
binary elements
\begin{subequations}
\begin{equation}
    \Triangle{1}{1}{1},
\end{equation}
\begin{equation}
    \TriangleEEE{}{}{}
    +
    \TriangleEXE{}{1}{}
    -
    \TriangleXEE{1}{}{}
    +
    \TriangleEEX{}{}{1}
    -
    \TriangleXXE{1}{1}{}
    +
    \TriangleEXX{}{1}{1}
    -
    \TriangleXEX{1}{}{1}
    -
    \Triangle{1}{1}{1}
\end{equation}
\end{subequations}
of $\Cli\N_2$ are associative.
\medbreak

\subsubsection{Symmetries}
Let $\Returned : \Cli\Mca \to \Cli\Mca$ be the linear map sending any
$\Mca$-clique $\Pfr$ of arity $n$ to the $\Mca$-clique
$\Returned(\Pfr)$ of the same arity such that, for any arc $(x, y)$
where $1 \leq x < y \leq n + 1$,
\begin{equation} \label{equ:returned_map_Cli_M}
    \left(\Returned(\Pfr)\right)(x, y) := \Pfr(n - y + 2, n - x + 2).
\end{equation}
In a geometric way, $\Returned(\Pfr)$ is the $\Mca$-clique obtained by
applying on $\Pfr$ a reflection trough the vertical line passing by its
base. For instance, one has in $\Cli\Z$,
\begin{equation}
    \Returned\left(
\,.
\end{equation}
\medbreak

\begin{Proposition} \label{prop:symmetries_Cli_M}
    Let $\Mca$ be a unitary magma. Then, the group of symmetries of
    $\Cli\Mca$ contains the map $\Returned$ and all the maps
    $\Cli\theta$ where $\theta$ are unitary magma automorphisms
    of~$\Mca$.
\end{Proposition}
\medbreak

\subsubsection{Basic set-operad basis}

\begin{Proposition} \label{prop:basic_Cli_M}
    Let $\Mca$ be a unitary magma. The fundamental basis of $\Cli\Mca$
    is a basic set-operad basis if and only if $\Mca$ is right
    cancellable.
\end{Proposition}
\medbreak

\subsubsection{Cyclic operad structure}
Let $\rho : \Cli\Mca \to \Cli\Mca$ be the linear map sending any
$\Mca$-clique $\Pfr$ of arity $n$ to the $\Mca$-clique $\rho(\Pfr)$ of
the same arity such that, for any arc $(x, y)$ where
$1 \leq x < y \leq n + 1$,
\begin{equation} \label{equ:rotation_map_Cli_M}
    (\rho(\Pfr))(x, y) :=
    \begin{cases}
        \Pfr(x + 1, y + 1) & \mbox{if } y \leq n, \\
        \Pfr(1, x + 1) & \mbox{otherwise (} y = n + 1 \mbox{)}.
    \end{cases}
\end{equation}
In a geometric way, $\rho(\Pfr)$ is the $\Mca$-clique obtained by
applying a rotation of one step of $\Pfr$ in the counterclockwise
direction. For instance, one has in $\Cli\Z$,
\begin{equation}
    \rho\left(
\,.
\end{equation}
\medbreak

\begin{Proposition} \label{prop:cyclic_Cli_M}
    Let $\Mca$ be a unitary magma. The map $\rho$ is a rotation map
    of $\Cli\Mca$, endowing this operad with a cyclic operad structure.
\end{Proposition}
\medbreak

\subsubsection{Alternative bases}
If $\Pfr$ and $\Qfr$ are two $\Mca$-cliques of the same arity, the
\Def{Hamming distance} $\Hamming(\Pfr, \Qfr)$ between $\Pfr$ and $\Qfr$
is the number of arcs $(x, y)$ such that $\Pfr(x, y) \ne \Qfr(x, y)$.
Let $\OrdBE$ be the partial order relation on the set of all
$\Mca$-cliques, where, for any $\Mca$-cliques $\Pfr$ and $\Qfr$, one
has $\Pfr \OrdBE \Qfr$ if $\Qfr$ can be obtained from $\Pfr$ by
replacing some labels $\Unit_\Mca$ of its edges or its base by other
labels of $\Mca$. In the same way, let $\OrdD$ be the partial order
on the same set where $\Pfr \OrdD \Qfr$ if $\Qfr$ can be obtained from
$\Pfr$ by replacing some labels $\Unit_\Mca$ of its diagonals by other
labels of $\Mca$.
\medbreak

For all $\Mca$-cliques $\Pfr$, let the elements of $\Cli\Mca$ defined by
\begin{subequations}
\begin{equation}
    \BasisH_\Pfr :=
    \sum_{\substack{
        \Pfr' \in \Cliques_\Mca \\
        \Pfr' \OrdBE \Pfr
    }}
    \Pfr',
\end{equation}
and
\begin{equation}
    \BasisK_\Pfr :=
    \sum_{\substack{
        \Pfr' \in \Cliques_\Mca \\
        \Pfr' \OrdD \Pfr
    }}
    (-1)^{\Hamming(\Pfr', \Pfr)}
    \Pfr'.
\end{equation}
\end{subequations}
For instance, in~$\Cli\Z$,
\begin{subequations}
\begin{equation}
    \BasisH_{
\,.
\end{equation}
\end{subequations}
\medbreak

Since by Möbius inversion (see
Proposition~\ref{prop:partial_order_bases} of
Chapter~\ref{chap:algebra}), one has for any $\Mca$-clique $\Pfr$,
\begin{equation}
    \sum_{\substack{
        \Pfr' \in \Cliques_\Mca \\
        \Pfr' \OrdBE \Pfr
    }}
    (-1)^{\Hamming(\Pfr', \Pfr)}
    \BasisH_{\Pfr'}
    =
    \Pfr
    =
    \sum_{\substack{
        \Pfr' \in \Cliques_\Mca \\
        \Pfr' \OrdD \Pfr
    }}
    \BasisK_{\Pfr'},
\end{equation}
by triangularity, the family of all the $\BasisH_\Pfr$ (resp.
$\BasisK_\Pfr$) forms a  basis of $\Cli\Mca$ called the
\Def{$\BasisH$-basis} (resp. the \Def{$\BasisK$-basis}).
\medbreak

If $\Pfr$ is an $\Mca$-clique, $\Del_0(\Pfr)$ (resp. $\Del_i(\Pfr)$) is
the $\Mca$-clique obtained by replacing the label of the base
(resp. $i$th edge) of $\Pfr$ by $\Unit_\Mca$.
\medbreak

\begin{Proposition} \label{prop:composition_Cli_M_basis_H}
    Let $\Mca$ be a unitary magma. The partial composition of $\Cli\Mca$
    can be expressed over the $\BasisH$-basis, for any $\Mca$-cliques
    $\Pfr$ and $\Qfr$ different from $\UnitClique$ and any valid
    integer $i$, as
    \begin{equation}
        \BasisH_\Pfr \circ_i \BasisH_\Qfr
        =
        \begin{cases}
            \BasisH_{\Pfr \circ_i \Qfr}
            + \BasisH_{\Del_i(\Pfr) \circ_i \Qfr}
            + \BasisH_{\Pfr \circ_i \Del_0(\Qfr)}
            + \BasisH_{\Del_i(\Pfr) \circ_i \Del_0(\Qfr)}
                & \mbox{if } \Pfr_i \ne \Unit_\Mca \mbox{ and }
                    \Qfr_0 \ne \Unit_\Mca, \\
            \BasisH_{\Pfr \circ_i \Qfr}
            + \BasisH_{\Del_i(\Pfr) \circ_i \Qfr}
                & \mbox{if } \Pfr_i \ne \Unit_\Mca, \\
            \BasisH_{\Pfr \circ_i \Qfr}
            + \BasisH_{\Pfr \circ_i \Del_0(\Qfr)}
                & \mbox{if } \Qfr_0 \ne \Unit_\Mca, \\
            \BasisH_{\Pfr \circ_i \Qfr} & \mbox{otherwise}.
        \end{cases}
    \end{equation}
\end{Proposition}
\medbreak

\begin{Proposition} \label{prop:composition_Cli_M_basis_K}
    Let $\Mca$ be a unitary magma. The partial composition of $\Cli\Mca$
    can be expressed over the $\BasisK$-basis, for any $\Mca$-cliques
    $\Pfr$ and $\Qfr$ different from $\UnitClique$ and any valid
    integer $i$, as
    \begin{equation}
        \BasisK_\Pfr \circ_i \BasisK_\Qfr
        =
        \begin{cases}
            \BasisK_{\Pfr \circ_i \Qfr}
                & \mbox{if }
                \Pfr_i \Op \Qfr_0 = \Unit_\Mca, \\
            \BasisK_{\Pfr \circ_i \Qfr} +
            \BasisK_{\Del_i(\Pfr) \circ_i \Del_0(\Qfr)}
                & \mbox{otherwise}.
        \end{cases}
    \end{equation}
\end{Proposition}
\medbreak

For instance, in $\Cli\Z$,
\begin{subequations}
\begin{equation}
    \BasisH_{
}\,.
\end{equation}
\end{subequations}
\medbreak

\subsubsection{Rational functions} \label{subsubsec:rational_functions}
We develop here a link between $\Cli\Mca$ and the operad $\RatFct$ of
rational functions introduced by Loday~\cite{Lod10} (see also
Section~\ref{subsubsec:rational_functions_operad} of
Chapter~\ref{chap:algebra}).
\medbreak

Let us assume that $\Mca$ is a \Def{$\Z$-graded unitary magma}, that
is a unitary magma such that there exists a unitary magma morphism
$\theta : \Mca \to \Z$. We say that $\theta$ is a \Def{rank function}
of $\Mca$. In this context, let
\begin{equation}
    \Frac_\theta : \Cli\Mca \to \RatFct
\end{equation}
be the linear map defined, for any $\Mca$-clique $\Pfr$, by
\begin{equation} \label{equ:definition_frac_clique}
    \Frac_\theta(\Pfr) :=
    \prod_{(x, y) \in \Arcs_\Pfr}
    \left(u_x + \dots + u_{y - 1}\right)^{\theta(\Pfr(x, y))}.
\end{equation}
For instance, by considering the unitary magma $\Z$ together with its
identity map $\Identity$ as rank function, one has
\begin{equation}
    \Frac_\Identity\left(
    \begin{tikzpicture}[scale=.85,Centering]
        \node[CliquePoint](1)at(-0.43,-0.90){};
        \node[CliquePoint](2)at(-0.97,-0.22){};
        \node[CliquePoint](3)at(-0.78,0.62){};
        \node[CliquePoint](4)at(-0.00,1.00){};
        \node[CliquePoint](5)at(0.78,0.62){};
        \node[CliquePoint](6)at(0.97,-0.22){};
        \node[CliquePoint](7)at(0.43,-0.90){};
        \draw[CliqueEdge](1)edge[]node[CliqueLabel]
            {\begin{math}-1\end{math}}(2);
        \draw[CliqueEdge](1)
            edge[bend left=30]node[CliqueLabel,near start]
            {\begin{math}2\end{math}}(5);
        \draw[CliqueEdge](1)edge[]node[CliqueLabel]
            {\begin{math}1\end{math}}(7);
        \draw[CliqueEmptyEdge](2)edge[]node[CliqueLabel]{}(3);
        \draw[CliqueEmptyEdge](3)edge[]node[CliqueLabel]{}(4);
        \draw[CliqueEdge](3)
            edge[bend left=30]node[CliqueLabel,near start]
            {\begin{math}-2\end{math}}(7);
        \draw[CliqueEdge](4)edge[]node[CliqueLabel]
            {\begin{math}3\end{math}}(5);
        \draw[CliqueEmptyEdge](5)edge[]node[CliqueLabel]{}(6);
        \draw[CliqueEmptyEdge](6)edge[]node[CliqueLabel]{}(7);
        \draw[CliqueEdge](5)edge[]node[CliqueLabel]
            {\begin{math}-1\end{math}}(7);
    \end{tikzpicture}\right)
    =
    \frac{
        \left(u_1 + u_2 + u_3 + u_4\right)^2
        \left(u_1 + u_2 + u_3 + u_4 + u_5 + u_6\right)
        u_4^3
    }{u_1 \left(u_3 + u_4 + u_5 + u_6\right)^2 \left(u_5 + u_6\right)}.
\end{equation}
\medbreak

\begin{Theorem} \label{thm:rat_fct_cliques}
    Let $\Mca$ be a $\Z$-graded unitary magma and $\theta$ be a rank
    function of $\Mca$. The map $\Frac_\theta$ is an operad morphism
    from $\Cli\Mca$ to~$\RatFct$.
\end{Theorem}
\medbreak

The operad morphism $\Frac_\theta$ is not injective. Indeed, by
considering the magma $\Z$ together with its identity map $\Identity$ as
rank function, one one has for instance
\begin{subequations}
\begin{equation} \label{equ:frac_not_injective_1}
    \Frac_\Identity\left(
    \TriangleXEE{1}{}{}
    - \TriangleEXE{}{1}{}
    - \TriangleEEX{}{}{1}
    \right)
    = (u_1 + u_2) - u_1 - u_2
    = 0,
\end{equation}
\begin{equation} \label{equ:frac_not_injective_2}
    \Frac_\Identity\left(
    \begin{tikzpicture}[scale=0.6,Centering]
        \node[CliquePoint](1)at(-0.71,-0.71){};
        \node[CliquePoint](2)at(-0.71,0.71){};
        \node[CliquePoint](3)at(0.71,0.71){};
        \node[CliquePoint](4)at(0.71,-0.71){};
        \draw[CliqueEmptyEdge](1)edge[]node[]{}(2);
        \draw[CliqueEmptyEdge](1)edge[]node[]{}(4);
        \draw[CliqueEdge](2)edge[]node[CliqueLabel]
            {\begin{math}-1\end{math}}(3);
        \draw[CliqueEdge](3)edge[]node[CliqueLabel]
            {\begin{math}-1\end{math}}(4);
    \end{tikzpicture}
    -
    \begin{tikzpicture}[scale=0.6,Centering]
        \node[CliquePoint](1)at(-0.71,-0.71){};
        \node[CliquePoint](2)at(-0.71,0.71){};
        \node[CliquePoint](3)at(0.71,0.71){};
        \node[CliquePoint](4)at(0.71,-0.71){};
        \draw[CliqueEmptyEdge](1)edge[]node[]{}(2);
        \draw[CliqueEmptyEdge](1)edge[]node[]{}(4);
        \draw[CliqueEmptyEdge](2)edge[]node[]{}(3);
        \draw[CliqueEdge](2)edge[]node[CliqueLabel]
            {\begin{math}-1\end{math}}(4);
        \draw[CliqueEdge](3)edge[]node[CliqueLabel]
            {\begin{math}-1\end{math}}(4);
    \end{tikzpicture}
    -
    \begin{tikzpicture}[scale=0.6,Centering]
        \node[CliquePoint](1)at(-0.71,-0.71){};
        \node[CliquePoint](2)at(-0.71,0.71){};
        \node[CliquePoint](3)at(0.71,0.71){};
        \node[CliquePoint](4)at(0.71,-0.71){};
        \draw[CliqueEmptyEdge](1)edge[]node[]{}(2);
        \draw[CliqueEmptyEdge](1)edge[]node[]{}(4);
        \draw[CliqueEmptyEdge](3)edge[]node[]{}(4);
        \draw[CliqueEdge](2)edge[]node[CliqueLabel]
            {\begin{math}-1\end{math}}(3);
        \draw[CliqueEdge](2)edge[]node[CliqueLabel]
            {\begin{math}-1\end{math}}(4);
    \end{tikzpicture}
    \right)
    =
    \frac{1}{u_2 u_3}
    -
    \frac{1}{(u_2 + u_3) u_3}
    -
    \frac{1}{u_2 (u_2 + u_3)}
    = 0.
\end{equation}
\end{subequations}
\medbreak

\begin{Proposition} \label{prop:rat_fct_cliques_map_Laurent_polynomials}
    The subspace of $\RatFct$ of all Laurent polynomials on~$\Ubb$ is
    the image by $\Frac_\Identity : \Cli\Z \to \RatFct$ of the subspace
    of $\Cli\Z$ consisting in the linear span of all $\Z$-bubbles.
\end{Proposition}
\medbreak

On each homogeneous subspace $\Cli\Mca(n)$ of the elements of arity
$n \geq 1$ of $\Cli\Mca$, let the product
\begin{equation}
    \Op : \Cli\Mca(n) \otimes \Cli\Mca(n) \to \Cli\Mca(n)
\end{equation}
defined linearly, for each $\Mca$-cliques $\Pfr$ and $\Qfr$ of
$\Cli\Mca(n)$, by
\begin{equation}
    (\Pfr \Op \Qfr)(x, y) := \Pfr(x, y) \Op \Qfr(x, y),
\end{equation}
where $(x, y)$ is any arc such that $1 \leq x < y \leq n + 1$. For
instance, in $\Cli\Z$,
\begin{equation}
\,.
\end{equation}
\medbreak

\begin{Proposition} \label{prop:rat_fct_cliques_product}
    Let $\Mca$ be a $\Z$-graded unitary magma and $\theta$ be a rank
    function of~$\Mca$. For any homogeneous elements $f$ and $g$ of
    $\Cli\Mca$ of the same arity,
    \begin{equation} \label{equ:rat_fct_cliques_product}
        \Frac_\theta(f) \Frac_\theta(g) = \Frac_\theta(f \Op g).
    \end{equation}
\end{Proposition}
\medbreak

\begin{Proposition} \label{prop:rat_fct_cliques_inverse}
    Let $\Pfr$ be an $\Mca$-clique of $\Cli\Z$. Then,
    \begin{equation}
        \frac{1}{\Frac_\Identity(\Pfr)}
        = \Frac_\Identity((\Cli \eta)(\Pfr)),
    \end{equation}
    where $\eta : \Z \to \Z$ is the unitary magma morphism defined by
    $\eta(x) := -x$ for all $x \in \Z$.
\end{Proposition}
\medbreak

\section{Quotients and suboperads}\label{sec:quotients_suboperads_C_M}
We define here quotients and suboperads of $\Cli\Mca$, leading to the
construction of some new operads involving various combinatorial objects
which are, basically, $\Mca$-cliques with some restrictions.
\medbreak

\subsection{Main substructures} \label{subsec:main_substructures_C_M}
Most of the natural subfamilies of $\Mca$-cliques that can be described
by simple combinatorial properties as $\Mca$-cliques with restrained
labels for the bases, edges, and diagonals, white $\Mca$-cliques,
$\Mca$-cliques with a fixed maximal value for their crossings,
$\Mca$-bubbles, $\Mca$-cliques with a fixed maximal value for their
degrees, nesting-free $\Mca$-cliques, and acyclic $\Mca$-cliques
inherit the algebraic structure of operad of $\Cli\Mca$ and form
quotients and suboperads of $\Cli\Mca$. We construct and briefly study
here these main substructures of~$\Cli\Mca$.
\medbreak

\subsubsection{Restricting the labels}
\label{subsubsec:suboperad_Cli_M_labels}
In what follows, if $X$ and $Y$ are two subsets of $\Mca$, $X \Op Y$
denotes the set $\{x \Op y : x \in X \mbox{ and } y \in Y\}$.
\medbreak

Let $B$, $E$, and $D$ be three subsets of $\Mca$ and
$\Lab_{B, E, D}\Mca$ be the subspace of $\Cli\Mca$ generated by all
$\Mca$-cliques $\Pfr$ such that the bases of $\Pfr$ are labeled on $B$,
all edges of $\Pfr$ are labeled on $E$, and all diagonals of $\Pfr$ are
labeled on~$D$.
\medbreak

\begin{Proposition} \label{prop:suboperad_Cli_M_labels}
    Let $\Mca$ be a unitary magma and $B$, $E$, and $D$ be three subsets
    of $\Mca$. If $\Unit_\Mca \in B$, $\Unit_\Mca \in D$, and
    $E \Op B \subseteq D$, $\Lab_{B, E, D}\Mca$ is a suboperad
    of~$\Cli\Mca$.
\end{Proposition}
\medbreak

\begin{Proposition} \label{prop:suboperad_Cli_M_labels_dimensions}
    Let $\Mca$ be a unitary magma and $B$, $E$, and $D$ be three
    finite subsets of $\Mca$. For all $n \geq 2$,
    \begin{equation} \label{equ:suboperad_Cli_M_labels_dimensions}
        \dim \Lab_{B, E, D}\Mca(n) =
         b e^n d^{(n + 1)(n - 2) / 2},
    \end{equation}
    where $b := \# B$, $e := \# E$, and $d := \# D$.
\end{Proposition}
\medbreak

\subsubsection{White cliques} \label{subsubsec:suboperad_Cli_M_white}
Let $\Whi\Mca$ be the subspace of $\Cli\Mca$ generated by all white
$\Mca$-cliques. Since, by definition of white $\Mca$-cliques,
\begin{equation}
    \Whi\Mca = \Lab_{\{\Unit_\Mca\}, \{\Unit_\Mca\}, \Mca}\Mca,
\end{equation}
by Proposition~\ref{prop:suboperad_Cli_M_labels}, $\Whi\Mca$ is a
suboperad of $\Cli\Mca$. It follows from
Proposition~\ref{prop:suboperad_Cli_M_labels_dimensions} that when
$\Mca$ is finite, the dimensions of $\Whi\Mca$ satisfy, for any
$n \geq 2$,
\begin{equation}
    \dim \Whi\Mca(n) =
    m^{(n + 1)(n - 2) / 2},
\end{equation}
where $m := \# \Mca$.
\medbreak

\subsubsection{Restricting the crossings}
\label{subsubsec:quotient_Cli_M_crossings}
Let $k \geq 0$ be an integer and $\RelationSpace_{\Cro_k\Mca}$ be the
subspace of $\Cli\Mca$ generated by all $\Mca$-cliques $\Pfr$ such that
$\Cros(\Pfr) \geq k + 1$. As a quotient of graded vector spaces,
\begin{equation}
    \Cro_k\Mca := \Cli\Mca/_{\RelationSpace_{\Cro_k\Mca}}
\end{equation}
is the linear span of all $\Mca$-cliques $\Pfr$ such that
$\Cros(\Pfr) \leq k$.
\medbreak

\begin{Proposition} \label{prop:quotient_Cli_M_crossings}
    Let $\Mca$ be a unitary magma and $k \geq 0$ be an integer. Then,
    the space $\Cro_k\Mca$ is both a quotient and a suboperad
    of~$\Cli\Mca$.
\end{Proposition}
\medbreak

For instance, in the operad $\Cro_2\Z$, we have
\begin{equation}
\,.
\end{equation}
\medbreak

When $0 \leq k' \leq k$ are integers, by
Proposition~\ref{prop:quotient_Cli_M_crossings}, $\Cro_k\Mca$ and
$\Cro_{k'}\Mca$ are both quotients and suboperads of $\Cli\Mca$. First,
since any $\Mca$-clique of $\Cro_{k'}\Mca$ is also an $\Mca$-clique of
$\Cro_k\Mca$, $\Cro_{k'}\Mca$ is a suboperad of~$\Cro_k\Mca$. Second,
since $\RelationSpace_{\Cro_k\Mca}$ is a subspace of
$\RelationSpace_{\Cro_{k'}\Mca}$, $\Cro_{k'}\Mca$ is a quotient
of~$\Cro_k\Mca$.
\medbreak

Remark that $\Cro_0\Mca$ is the linear span of all noncrossing
$\Mca$-cliques. We can see these objects as noncrossing
configurations~\cite{FN99} where the edges and bases are colored by
elements of $\Mca$ and the diagonals, by elements of $\bar{\Mca}$. The
operad $\Cro_0\Mca$ has a lot of properties and will be studied in
details in Section~\ref{sec:operad_noncrossing}.
\medbreak

\subsubsection{Bubbles}
\label{subsubsec:quotient_Cli_M_bubbles}
Let $\RelationSpace_{\Bub\Mca}$ be the subspace of $\Cli\Mca$
generated by all $\Mca$-cliques that are not bubbles. As a quotient of
graded vector spaces,
\begin{equation}
    \Bub\Mca := \Cli\Mca/_{\RelationSpace_{\Bub\Mca}}
\end{equation}
is the linear span of all $\Mca$-bubbles.
\medbreak

\begin{Proposition} \label{prop:quotient_Cli_M_bubbles}
    Let $\Mca$ be a unitary magma. Then, the space $\Bub_\Mca$ is a
    quotient operad of $\Cli\Mca$.
\end{Proposition}
\medbreak

For instance, in the operad $\Bub\Z$, we have
\vspace{-1.75em}
\begin{multicols}{2}
\begin{subequations}
\begin{equation}

    = 0.
\end{equation}
\end{subequations}
\end{multicols}
\medbreak

When $\Mca$ is finite, the dimensions of $\Bub\Mca$ satisfy, for any
$n \geq 2$,
\begin{equation}
    \dim \Bub\Mca(n) = m^{n + 1},
\end{equation}
where $m := \# \Mca$.
\medbreak

\subsubsection{Restricting the degrees}
\label{subsubsec:quotient_Cli_M_degrees}
Let $k \geq 0$ be an integer and $\RelationSpace_{\Deg_k\Mca}$ be the
subspace of $\Cli\Mca$ generated by all $\Mca$-cliques $\Pfr$ such that
$\Degr(\Pfr) \geq k + 1$. As a quotient of graded vector spaces,
\begin{equation}
    \Deg_k\Mca := \Cli\Mca/_{\RelationSpace_{\Deg_k\Mca}}
\end{equation}
is the linear span of all $\Mca$-cliques $\Pfr$ such that
$\Degr(\Pfr) \leq k$.
\medbreak

\begin{Proposition} \label{prop:quotient_Cli_M_degrees}
    Let $\Mca$ be a unitary magma without nontrivial unit divisors and
    $k \geq 0$ be an integer. Then, the space $\Deg_k\Mca$ is a quotient
    operad of $\Cli\Mca$.
\end{Proposition}
\medbreak

For instance, in the operad $\Deg_3\Dbb_2$ (observe that $\Dbb_2$ is
a unitary magma without nontrivial unit divisors), we have
\vspace{-1.75em}
\begin{multicols}{2}
\begin{subequations}
\begin{equation}

    = 0.
\end{equation}
\end{subequations}
\end{multicols}
\medbreak

When $0 \leq k' \leq k$ are integers, by
Proposition~\ref{prop:quotient_Cli_M_degrees}, $\Deg_k\Mca$ and
$\Deg_{k'}\Mca$ are both quotients operads of $\Cli\Mca$. Moreover,
since $\RelationSpace_{\Deg_k\Mca}$ is a subspace of
$\RelationSpace_{\Deg_{k'}\Mca}$, $\Deg_{k'}\Mca$ is a quotient operad
of~$\Deg_k\Mca$.
\medbreak

Observe that $\Deg_0\Mca$ is the linear span of all $\Mca$-cliques
without solid arcs. If $\Pfr$ and $\Qfr$ are such $\Mca$-cliques, all
partial compositions $\Pfr \circ_i \Qfr$ are equal to the unique
$\Mca$-clique without solid arcs of arity $|\Pfr| + |\Qfr| - 1$. For
this reason, $\Deg_0\Mca$ is the associative operad~$\As$.
\medbreak

Any skeleton of an $\Mca$-clique of arity $n$ of $\Deg_1\Mca$ can be
seen as a partition of the set $[n + 1]$ in singletons or pairs.
Therefore, $\Deg_1\Mca$ can be seen as an operad on such colored
partitions, where each pair of the partitions have one color among the
set $\bar{\Mca}$. In the operad $\Deg_1\Dbb_0$ (observe that $\Dbb_0$ is
the only unitary magma without nontrivial unit divisors on two
elements), one has for instance
\vspace{-1.75em}
\begin{multicols}{2}
\begin{subequations}
\begin{equation} \label{equ:example_involutions_1}

    = 0.
\end{equation}
\end{subequations}
\end{multicols}
\medbreak

By seeing each solid arc $(x, y)$ of an $\Mca$-clique $\Pfr$ of
$\Deg_1\Dbb_0$ of arity $n$ as the transposition exchanging the letter
$x$ and the letter $y$, we can interpret $\Pfr$ as an involution of
$\mathfrak{S}_{n + 1}$ made of the product of these transpositions.
Hence, $\Deg_1\Dbb_0$ can be seen as an operad on involutions. Under
this point of view, the partial
compositions~\eqref{equ:example_involutions_1}
and~\eqref{equ:example_involutions_2} translate on permutations as
\vspace{-1.75em}
\begin{multicols}{2}
\begin{subequations}
\begin{equation}
    42315 \circ_2 3412 = 6452317,
\end{equation}

\begin{equation}
    42315 \circ_3 3412 = 0.
\end{equation}
\end{subequations}
\end{multicols}
\noindent Equivalently, by the Robinson-Schensted correspondence (see
for instance~\cite{Sch61,Lot02}), $\Deg_1\Dbb_0$ is an operad of
standard Young tableaux. The dimensions of $\Deg_1\Dbb_0$ operad begin
by
\begin{equation}
    1, 4, 10, 26, 76, 232, 764, 2620,
\end{equation}
and form, except for the first terms, Sequence~\OEIS{A000085}
of~\cite{Slo}. Moreover, when $\# \Mca = 3$, the dimensions of
$\Deg_1\Mca$ begin by
\begin{equation}
    1, 7, 25, 81, 331, 1303, 5937, 26785,
\end{equation}
and form, except for the first terms, Sequence~\OEIS{A047974}
of~\cite{Slo}.
\medbreak

Besides, any skeleton of an $\Mca$-clique of $\Deg_2\Mca$ can be seen as
a \Def{thunderstorm graph}, {\em i.e.}, a graph where connected
components are cycles or paths. Therefore, $\Deg_2\Mca$ can be seen as
an operad on such colored graphs, where the arcs of the graphs have one
color among the set $\bar{\Mca}$. When $\# \Mca = 2$, the
dimensions of this operad begin by
\begin{equation}
    1, 8, 41, 253, 1858, 15796, 152219, 1638323,
\end{equation}
and form, except for the first terms, Sequence~\OEIS{A136281}
of~\cite{Slo}.
\medbreak

\subsubsection{Nesting-free cliques}\label{subsubsec:quotient_Cli_M_Nes}
Let $\RelationSpace_{\Nes\Mca}$ be the subspace of $\Cli\Mca$
generated by all $\Mca$-cliques that are not nesting-free. As a
quotient of graded vector spaces,
\begin{equation}
    \Nes\Mca := \Cli\Mca/_{\RelationSpace_{\Nes\Mca}}
\end{equation}
is the linear span of all nesting-free $\Mca$-cliques.
\medbreak

\begin{Proposition} \label{prop:quotient_Cli_M_nesting_free}
    Let $\Mca$ be a unitary magma without nontrivial unit divisors.
    Then, the space $\Nes\Mca$ is a quotient operad of $\Cli\Mca$.
\end{Proposition}
\medbreak

For instance, in the operad $\Nes\Dbb_2$,
\vspace{-1.75em}
\begin{multicols}{2}
\begin{subequations}
\begin{equation}

    = 0.
\end{equation}
\end{subequations}
\end{multicols}
\medbreak

Remark that in the same way as considering $\Mca$-cliques of crossings
no greater than $k$ leads to quotients $\Cro_k\Mca$ of $\Cli\Mca$ (see
Section~\ref{subsubsec:quotient_Cli_M_crossings}), it is possible to
define analogous quotients $\Nes_k\Mca$ spanned by $\Mca$-cliques having
solid arcs that nest at most $k$ other ones.
\medbreak

\begin{Proposition} \label{prop:dimensions_Nes_M}
    Let $\Mca$ be a finite unitary magma without nontrivial unit
    divisors. For all $n \geq 2$,
    \begin{equation} \label{equ:dimensions_Nes_M}
        \dim \Nes\Mca(n) =
        \sum_{0 \leq k \leq n}
        (m - 1)^k \; \Narayana(n + 2, k),
    \end{equation}
    where $m := \# \Mca$.
\end{Proposition}
\medbreak

In the statement of Proposition~\ref{prop:dimensions_Nes_M},
$\Narayana(n, k)$ is a Narayana number whose definition is
recalled in Section~\ref{subsubsec:schroder_trees} of
Chapter~\ref{chap:combinatorics}.
\medbreak

The skeletons of the $\Mca$-cliques of $\Nes\Mca$ of arities greater
than $1$ are the graphs such that, if $\{x, y\}$ and $\{x', y'\}$ are
two arcs such that $x \leq x' < y' \leq y$, then $x = x'$ and $y = y'$.
Therefore, $\Nes\Mca$ can be seen as an operad on such colored graphs,
where the arcs of the graphs have one color among the set $\bar{\Mca}$.
\medbreak

By Proposition~\ref{prop:dimensions_Nes_M}, when $\# \Mca = 2$, the
dimensions of $\Nes\Mca$ begin by
\begin{equation}
    1, 5, 14, 42, 132, 429, 1430, 4862,
\end{equation}
and form, except for the first terms, Sequence~\OEIS{A000108}
of~\cite{Slo}. When $\# \Mca = 3$, the dimensions of
$\Nes\Mca$ begin by
\begin{equation}
    1, 11, 45, 197, 903, 4279, 20793, 103049,
\end{equation}
and form, except for the first terms, Sequence~\OEIS{A001003}
of~\cite{Slo}. When $\# \Mca = 4$, the dimensions of
$\Nes\Mca$ begin by
\begin{equation}
    1, 19, 100, 562, 3304, 20071, 124996, 793774,
\end{equation}
and form, except for the first terms, Sequence~\OEIS{A007564}
of~\cite{Slo}.
\medbreak

\subsubsection{Acyclic decorated cliques}
\label{subsubsec:quotient_Cli_M_acyclic}
Let $\RelationSpace_{\Acy\Mca}$ be the subspace of $\Cli\Mca$ generated
by all $\Mca$-cliques that are not acyclic. As a quotient of graded
vector spaces,
\begin{equation}
    \Acy\Mca := \Cli\Mca/_{\RelationSpace_{\Acy\Mca}}
\end{equation}
is the linear span of all acyclic $\Mca$-cliques.
\medbreak

\begin{Proposition} \label{prop:quotient_Cli_M_acyclic}
    Let $\Mca$ be a unitary magma without nontrivial unit divisors.
    Then, the space $\Acy\Mca$ is a quotient operad of $\Cli\Mca$.
\end{Proposition}
\medbreak

For instance, in the operad $\Acy\Dbb_2$,
\vspace{-1.75em}
\begin{multicols}{2}
\begin{subequations}
\begin{equation}

    = 0.
\end{equation}
\end{subequations}
\end{multicols}
\medbreak

The skeletons of the $\Mca$-cliques of $\Acy\Mca$ of arities greater
than $1$ are acyclic graphs or equivalently, forest of non-rooted trees.
Therefore, $\Acy\Mca$ can be seen as an operad on colored forests of
trees, where the edges of the trees of the forests have one color among
the set $\bar{\Mca}$. When $\# \Mca = 2$, the dimensions of $\Acy\Mca$
begin by
\begin{equation}
    1, 7, 38, 291, 2932, 36961, 561948, 10026505,
\end{equation}
and form, except for the first terms, Sequence~\OEIS{A001858}
of~\cite{Slo}.
\medbreak

\subsection{Secondary substructures}
\label{subsec:secondary_substructures_C_M}
Some more substructures of $\Cli\Mca$ are constructed and briefly
studied here. They are constructed by mixing some of the constructions
of the seven main substructures of $\Cli\Mca$ defined in
Section~\ref{subsec:main_substructures_C_M} in the following sense.
\medbreak

For any operad $\Oca$ and operad ideals $\RelationSpace_1$ and
$\RelationSpace_2$ of $\Oca$, the space
$\RelationSpace_1 + \RelationSpace_2$ is still an operad ideal of
$\Oca$, and $\Oca/_{\RelationSpace_1 + \RelationSpace_2}$ is a quotient
of both $\Oca/_{\RelationSpace_1}$ and $\Oca/_{\RelationSpace_2}$.
Moreover, if $\Oca'$ is a suboperad of $\Oca$ and $\RelationSpace$ is an
operad ideal of $\Oca$, the space $\RelationSpace \cap \Oca'$ is an
operad ideal of $\Oca'$, and $\Oca'/_{\RelationSpace \cap \Oca'}$ is a
quotient of $\Oca'$ and a suboperad of $\Oca/_\RelationSpace$. For these
reasons (straightforwardly provable), we can combine the constructions
of the previous section to build a bunch of new suboperads and quotients
of~$\Cli\Mca$.
\medbreak

\subsubsection{Colored white noncrossing configurations}
When $\Mca$ is a unitary magma, let
\begin{equation}
    \WNC\Mca := \Whi\Mca/_{\RelationSpace_{\Cro_0\Mca} \cap \Whi\Mca}.
\end{equation}
The $\Mca$-cliques of $\WNC\Mca$ are white noncrossing $\Mca$-cliques.
\medbreak

\begin{Proposition} \label{prop:dimensions_WNC_M}
    Let $\Mca$ be a finite unitary magma. For all $n \geq 2$,
    \begin{equation}
        \dim \WNC\Mca(n) =
        \sum_{0 \leq k \leq n - 2}
        m^k (m - 1)^{n - k - 2} \; \Narayana(n, k),
    \end{equation}
    where $m := \# \Mca$.
\end{Proposition}
\medbreak

When $\# \Mca = 2$, the dimensions of $\WNC\Mca$ begin by
\begin{equation}
    1, 1, 3, 11, 45, 197, 903, 4279,
\end{equation}
and form Sequence~\OEIS{A001003} of~\cite{Slo}. When $\# \Mca = 3$, the
dimensions of $\WNC\Mca$ begin by
\begin{equation}
    1, 1, 5, 31, 215, 1597, 12425, 99955,
\end{equation}
and form Sequence~\OEIS{A269730} of~\cite{Slo}. Observe that these
dimensions are shifted versions the ones of the
$\gamma$-polytridendriform operads $\TDendr_\gamma$ (see
Section~\ref{subsec:polytdendr} of Chapter~\ref{chap:polydendr})
with~$\gamma := \# \Mca - 1$.
\medbreak

\subsubsection{Colored forests of paths}
When $\Mca$ is a unitary magma without nontrivial unit divisors, let
\begin{equation}
    \Paths\Mca :=
    \Cli\Mca/_{\RelationSpace_{\Deg_2\Mca} + \RelationSpace_{\Acy\Mca}}.
\end{equation}
The skeletons of the $\Mca$-cliques of $\Paths\Mca$ are forests of
non-rooted trees that are paths. Therefore, $\Paths\Mca$ can be seen as
an operad on colored such graphs, where the arcs of the graphs have one
color among the set~$\bar{\Mca}$.
\medbreak

When $\# \Mca = 2$, the dimensions of $\Paths\Mca$ begin by
\begin{equation}
    1, 7, 34, 206, 1486, 12412, 117692, 1248004,
\end{equation}
an form, except for the first terms, Sequence~\OEIS{A011800}
of~\cite{Slo}.
\medbreak

\subsubsection{Colored forests}
When $\Mca$ is a unitary magma without nontrivial unit divisors, let
\begin{equation}
    \Forests\Mca :=
    \Cli\Mca/_{\RelationSpace_{\Cro_0\Mca} + \RelationSpace_{\Acy\Mca}}.
\end{equation}
The skeletons of the $\Mca$-cliques of $\Forests\Mca$ are forests of
rooted trees having no arcs $\{x, y\}$ and $\{x', y'\}$ satisfying
$x < x' < y < y'$. Therefore, $\Forests\Mca$ can be seen as an operad
on such colored forests, where the edges of the forests have one color
among the set $\bar{\Mca}$. When $\# \Mca = 2$, the dimensions of
$\Forests\Mca$ begin by
\begin{equation}
    1, 7, 33, 81, 1083, 6854, 45111, 305629,
\end{equation}
and form, except for the first terms, Sequence~\OEIS{A054727},
of~\cite{Slo}.
\medbreak

\subsubsection{Colored Motzkin configurations}
\label{subsubsec:Motzkin_configurations}
When $\Mca$ is a unitary magma without nontrivial unit divisors, let
\begin{equation}
    \Motzkin\Mca :=
    \Cli\Mca/_%
    {\RelationSpace_{\Cro_0\Mca} + \RelationSpace_{\Deg_1\Mca}}.
\end{equation}
The skeletons of the $\Mca$-cliques of $\Motzkin\Mca$ are configurations
of non-intersecting chords on a circle. Equivalently, these objects are
graphs of involutions (see
Section~\ref{subsubsec:quotient_Cli_M_degrees}) having no arcs
$\{x, y\}$ and $\{x', y'\}$ satisfying $x < x' < y < y'$. These objects
are enumerated by Motzkin numbers~\cite{Mot48}. Therefore,
$\Motzkin\Mca$ can be seen as an operad on such colored graphs, where
the arcs of the graphs have one color among the set $\bar{\Mca}$. When
$\# \Mca = 2$, the dimensions of $\Motzkin\Mca$ begin by
\begin{equation}
    1, 4, 9, 21, 51, 127, 323, 835,
\end{equation}
and form, except for the first terms, Sequence~\OEIS{A001006},
of~\cite{Slo}.
\medbreak

\subsubsection{Colored dissections of polygons}
When $\Mca$ is a unitary magma without nontrivial unit divisors, let
\begin{equation}
    \Diss\Mca :=
    \Whi\Mca/_%
    {(\RelationSpace_{\Cro_0\Mca} + \RelationSpace_{\Deg_1\Mca})
    \cap \Whi\Mca}.
\end{equation}
The skeletons of the $\Mca$-cliques of $\Diss\Mca$ are \Def{strict
dissections of polygons}, that are graphs of Motzkin configurations
with no arcs of the form $\{x, x + 1\}$ or $\{1, n + 1\}$, where $n + 1$
is the number of vertices of the graphs. Therefore, $\Diss\Mca$ can be
seen as an operad on such colored graphs, where the arcs of the graphs
have one color among the set $\bar{\Mca}$. When $\# \Mca = 2$, the
dimensions of $\Diss\Mca$ begin by
\begin{equation}
    1, 1, 3, 6, 13, 29, 65, 148,
\end{equation}
and form, except for the first terms, Sequence~\OEIS{A093128}
of~\cite{Slo}.
\medbreak

\subsubsection{Colored Lucas configurations}
When $\Mca$ is a unitary magma without nontrivial unit divisors, let
\begin{equation}
    \Luc\Mca :=
    \Cli\Mca/_{\RelationSpace_{\Bub\Mca} + \RelationSpace_{\Deg_1\Mca}}.
\end{equation}
The skeletons of the $\Mca$-cliques of $\Luc\Mca$ are graphs such that
all vertices are of degrees at most $1$ and all arcs are of the form
$\{x, x + 1\}$ or $\{1, n + 1\}$, where $n + 1$ is the number of
vertices of the graphs. Therefore, $\Luc\Mca$ can be seen as an operad
on such colored graphs, where the arcs of the graphs have one color
among the set $\bar{\Mca}$. When $\# \Mca = 2$, the dimensions of
$\Luc\Mca$ begin by
\begin{equation}
    1, 4, 7, 11, 18, 29, 47, 76,
\end{equation}
and form, except for the first terms, Sequence~\OEIS{A000032}
of~\cite{Slo}.
\medbreak

\subsection{Relations between substructures}
The suboperads and quotients of $\Cli\Mca$ construc\-ted in
Sections~\ref{subsec:main_substructures_C_M}
and~\ref{subsec:secondary_substructures_C_M} are linked by injective or
surjective operad morphisms. To establish these, we rely on the
following lemma.
\medbreak

\begin{Lemma} \label{lem:inclusion_families_cliques}
    Let $\Mca$ be a unitary magma. Then,
    \begin{enumerate}[fullwidth,label={(\it\roman*)}]
        \item \label{item:inclusion_families_cliques_1}
        the space $\RelationSpace_{\Acy\Mca}$ is a subspace of
        $\RelationSpace_{\Deg_1\Mca}$;
        \item \label{item:inclusion_families_cliques_2}
        the spaces $\RelationSpace_{\Nes\Mca}$ and
        $\RelationSpace_{\Bub\Mca}$ are subspaces of
        $\RelationSpace_{\Deg_0\Mca}$;
        \item \label{item:inclusion_families_cliques_3}
        the spaces $\RelationSpace_{\Cro_0\Mca}$ and
        $\RelationSpace_{\Deg_2\Mca}$ are subspaces of
        $\RelationSpace_{\Bub\Mca}$;
        \item \label{item:inclusion_families_cliques_4}
        the spaces $\RelationSpace_{\Deg_2\Mca}$ and
        $\RelationSpace_{\Acy\Mca}$ are subspaces of
        $\RelationSpace_{\Nes\Mca}$.
    \end{enumerate}
\end{Lemma}
\medbreak

\subsubsection{Relations between the main substructures}
\begin{figure}[ht]
    \centering
    \scalebox{.73}{
    \begin{tikzpicture}[xscale=1.2,yscale=1.1,Centering]
        \node[text=Col2](CliM)at(8,10)
            {\begin{math}\Cli\Mca\end{math}};
        \node[text=Col4](AcyM)at(4,8)
            {\begin{math}\Acy\Mca\end{math}};
        \node[text=Col4](DegkM)at(6,8)
            {\begin{math}\Deg_k\Mca\end{math}};
        \node[text=Col1](CrokM)at(10,8)
            {\begin{math}\Cro_k\Mca\end{math}};
        \node[text=Col1](LabM)at(12,8)
            {\begin{math}\Lab_{B, E, D}\Mca\end{math}};
        \node[text=Col4](Deg2M)at(6,6)
            {\begin{math}\Deg_2\Mca\end{math}};
        \node[text=Col1](Cro0M)at(10,6)
            {\begin{math}\Cro_0\Mca\end{math}};
        \node[text=Col4](NesM)at(4,4)
            {\begin{math}\Nes\Mca\end{math}};
        \node[text=Col4](Deg1M)at(6,4)
            {\begin{math}\Deg_1\Mca\end{math}};
        \node[text=Col4](BubM)at(8,4)
            {\begin{math}\Bub\Mca\end{math}};
        \node[text=Col1](WhiM)at(12,4)
            {\begin{math}\Whi\Mca\end{math}};
        \node[text=Col4](Deg0M)at(8,2)
            {\begin{math}\Deg_0\Mca\end{math}};
        \draw[Surjection](CliM)--(AcyM);
        \draw[Surjection](CliM)--(DegkM);
        \draw[Surjection](CliM)to[bend right=15](CrokM);
        \draw[Injection](CrokM)to[bend right=15](CliM);
        \draw[Injection](LabM)--(CliM);
        \draw[Surjection](AcyM)--(NesM);
        \draw[Surjection](AcyM)--(Deg1M);
        \draw[Surjection](DegkM)--(Deg2M);
        \draw[Surjection](CrokM)to[bend right=15](Cro0M);
        \draw[Injection](Cro0M)to[bend right=15](CrokM);
        \draw[Injection](WhiM)--(LabM);
        \draw[Surjection](Deg2M)--(NesM);
        \draw[Surjection](Deg2M)--(Deg1M);
        \draw[Surjection](Deg2M)--(BubM);
        \draw[Surjection](Cro0M)--(BubM);
        \draw[Surjection](NesM)--(Deg0M);
        \draw[Surjection](Deg1M)--(Deg0M);
        \draw[Surjection](BubM)--(Deg0M);
        \draw[Surjection](WhiM)--(Deg0M);
    \end{tikzpicture}}
    \caption[The diagram of the main suboperads and quotients
    of the ns operad~$\Cli\Mca$.]
    {The diagram of the main suboperads and quotients of $\Cli\Mca$.
    Arrows~$\rightarrowtail$ (resp.~$\twoheadrightarrow$) are injective
    (resp. surjective) operad morphisms. Here, $\Mca$ is a unitary magma
    without nontrivial unit divisors, $k$ is a positive integer, and
    $B$, $E$, and $D$ are subsets of $\Mca$ such that
    $\Unit_\Mca \in B$, $\Unit_\Mca \in E$, and $E \Op B \subseteq D$.}
    \label{fig:diagram_main_operads_C_M}
\end{figure}
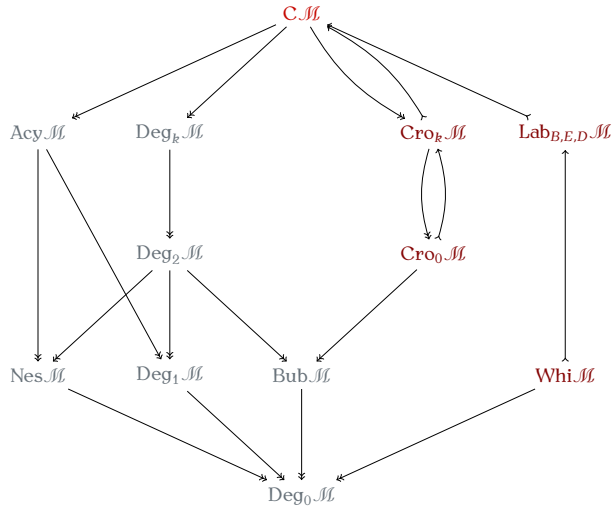

Let us list and explain the morphisms between the main substructures of
$\Cli\Mca$. First, Lemma~\ref{lem:inclusion_families_cliques} implies
that there are surjective operad morphisms from $\Acy\Mca$ to
$\Deg_1\Mca$, from $\Nes\Mca$ to $\Deg_0\Mca$, from $\Bub\Mca$ to
$\Deg_0\Mca$, from $\Cro_0\Mca$ to $\Bub\Mca$, from $\Deg_2\Mca$ to
$\Bub\Mca$, from $\Deg_2\Mca$ to $\Nes\Mca$, and from $\Acy\Mca$
to~$\Nes\Mca$. Second, when $B$, $E$, and $D$ are subsets of $\Mca$ such
that $\Unit_\Mca \in B$, $\Unit_\Mca \in E$, and $E \Op B \subseteq D$,
$\Whi\Mca$ is a suboperad of $\Lab_{B,E,D}\Mca$. Finally, there is a
surjective operad morphism from $\Whi\Mca$ to the associative operad
$\As$ sending any $\Mca$-clique $\Pfr$ of $\Whi\Mca$ to the unique basis
element of $\As$ of the same arity as the one of~$\Pfr$. The relations
between the main suboperads and quotients of $\Cli\Mca$ built here are
summarized in the diagram of operad morphisms of
Figure~\ref{fig:diagram_main_operads_C_M}.
\medbreak

\subsubsection{Relations between the secondary and main substructures}
\begin{figure}[ht]
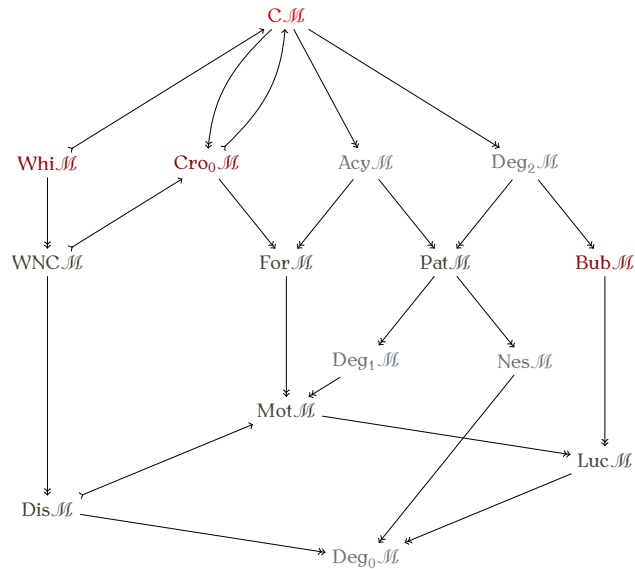

    \centering
    \scalebox{.73}{
}
    \caption[The diagram of the secondary suboperads and quotients
    of~$\Cli\Mca$.]
    {The diagram of the secondary suboperads and quotients of $\Cli\Mca$
    together with some of their related main suboperads and quotients of
    $\Cli\Mca$. Arrows~$\rightarrowtail$ (resp.~$\twoheadrightarrow$)
    are injective (resp. surjective) operad morphisms. Here, $\Mca$ is a
    unitary magma without nontrival unit divisors.}
    \label{fig:diagram_secondary_operads_C_M}
\end{figure}

Let us now list and explain the morphisms between the secondary and main
substructures of $\Cli\Mca$. First, immediately from their definitions,
$\WNC\Mca$ is a suboperad of $\Cro_0\Mca$ and a quotient of~$\Whi\Mca$,
$\Paths\Mca$ is both a quotient of $\Deg_2\Mca$ and $\Acy\Mca$,
$\Forests\Mca$ is both a quotient of $\Cro_0\Mca$ and $\Acy\Mca$,
$\Motzkin\Mca$ is both a quotient of $\Cro_0\Mca$ and $\Deg_1\Mca$,
$\Diss\Mca$ is a suboperad of $\Motzkin\Mca$ and a quotient of
$\WNC\Mca$, and $\Luc\Mca$ is both a quotient of $\Bub\Mca$ and
$\Deg_1\Mca$. Moreover, since by
Lemma~\ref{lem:inclusion_families_cliques},
$\RelationSpace_{\Acy\Mca}$ is a subspace of
$\RelationSpace_{\Deg_1\Mca}$, $\RelationSpace_{\Deg_2\Mca}$ and
$\RelationSpace_{\Acy\Mca}$ are subspaces of
$\RelationSpace_{\Nes\Mca}$, and $\RelationSpace_{\Cro_0\Mca}$ is a
subspace of $\RelationSpace_{\Bub\Mca}$, we respectively have that
$\RelationSpace_{\Deg_2\Mca} + \RelationSpace_{\Acy\Mca}$ is a
subspace of both $\RelationSpace_{\Deg_1\Mca}$ and
$\RelationSpace_{\Nes\Mca}$,
$\RelationSpace_{\Cro_0\Mca} + \RelationSpace_{\Acy\Mca}$ is a subspace
of $\RelationSpace_{\Cro_0\Mca} + \RelationSpace_{\Deg_1\Mca}$, and
$\RelationSpace_{\Cro_0\Mca} + \RelationSpace_{\Deg_1\Mca}$ is a
subspace of $\RelationSpace_{\Bub\Mca} + \RelationSpace_{\Deg_1\Mca}$.
For these reasons, there are surjective operad morphisms from
$\Paths\Mca$ to $\Deg_1\Mca$, from $\Paths\Mca$ to~$\Nes\Mca$, from
$\Forests\Mca$ to~$\Motzkin\Mca$, and from $\Motzkin\Mca$ to~$\Luc\Mca$.
The relations between the secondary suboperads and quotients of
$\Cli\Mca$ built here are summarized in the diagram of operad
morphisms of Figure~\ref{fig:diagram_secondary_operads_C_M}.
\medbreak

\section{Operads of noncrossing decorated cliques}
\label{sec:operad_noncrossing}
We perform here a complete study of the suboperad $\Cro_0\Mca$ of
noncrossing $\Mca$-cliques defined in
Section~\ref{subsubsec:quotient_Cli_M_crossings}. For simplicity, this
operad is denoted in the sequel as $\NC\Mca$ and named as the
\Def{noncrossing $\Mca$-clique operad}. The process giving from any
unitary magma $\Mca$ the operad $\NC\Mca$ is called the
\Def{noncrossing clique construction}.
\medbreak

\subsection{General properties}
To study $\NC\Mca$, we begin by establishing the fact that $\NC\Mca$
inherits some properties of~$\Cli\Mca$. Then, we shall describe a
realization of $\NC\Mca$ in terms of decorated Schröder trees, compute a
minimal generating set of $\NC\Mca$, and compute its dimensions.
\medbreak

First of all, we call \Def{fundamental basis} of $\NC\Mca$ the
fundamental basis of $\Cli\Mca$ restricted on noncrossing
$\Mca$-cliques. By definition of $\NC\Mca$ and by
Proposition~\ref{prop:quotient_Cli_M_crossings}, the partial composition
$\Pfr \circ_i \Qfr$ of two noncrossing $\Mca$-cliques $\Pfr$ and $\Qfr$
in $\NC\Mca$ is equal to the partial composition $\Pfr \circ_i \Qfr$ in
$\Cli\Mca$. Therefore, the fundamental basis of $\NC\Mca$ is a
set-operad basis.
\medbreak

\subsubsection{First properties}

\begin{Proposition} \label{prop:inherited_properties_NC_M}
    Let $\Mca$ be a unitary magma. Then,
    \begin{enumerate}[fullwidth,label={(\it\roman*)}]
        \item \label{item:inherited_properties_NC_M_1}
        the associative elements of $\NC\Mca$ are the ones of $\Cli\Mca$;
        \item \label{item:inherited_properties_NC_M_2}
        the group of symmetries of $\NC\Mca$ contains the map
        $\Returned$ (defined by~\eqref{equ:returned_map_Cli_M}) and all
        the maps $\Cli\theta$ where $\theta$ are unitary magma
        automorphisms of $\Mca$;
        \item \label{item:inherited_properties_NC_M_3}
        the fundamental basis of $\NC\Mca$ is a basic set-operad basis
        if and only if $\Mca$ is right cancellable;
        \item \label{item:inherited_properties_NC_M_4}
        the map $\rho$ (defined by~\eqref{equ:rotation_map_Cli_M}) is a
        rotation map of $\NC\Mca$ endowing it with a cyclic operad
        structure.
    \end{enumerate}
\end{Proposition}
\medbreak

\subsubsection{Treelike expressions on bubbles}
\label{subsubsec:treelike_bubbles}
Let $\Pfr$ be a noncrossing $\Mca$-clique or arity $n \geq 2$, and
$(x, y)$ be a diagonal or the base of $\Pfr$. Consider the path
$(x = z_1, z_2, \dots, z_k, z_{k + 1} = y)$ in $\Pfr$ such that
$k \geq 2$, for all $i \in [k + 1]$, $x \leq z_i \leq y$, and for all
$i \in [k]$, $z_{i + 1}$ is the greatest vertex of $\Pfr$ so that
$(z_i, z_{i + 1})$ is a solid diagonal or a (non-necessarily solid) edge
of $\Pfr$. The \Def{area} of $\Pfr$ adjacent to $(x, y)$ is the
$\Mca$-bubble $\Qfr$ of arity $k$ whose base is labeled by $\Pfr(x, y)$
and $\Qfr_i = \Pfr(z_i, z_{i + 1})$ for all $i \in [k]$. From a
geometric point of view, $\Qfr$ is the unique maximal component of
$\Pfr$ adjacent to the arc $(x, y)$, without solid diagonals, and
bounded by solid diagonals or edges of~$\Pfr$. For instance, for the
noncrossing $\Z$-clique
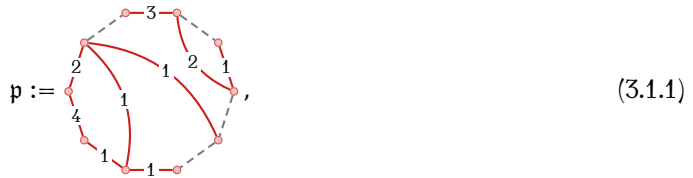
\begin{equation}
    \Pfr :=
\,,
\end{equation}
the path associated with the diagonal $(4, 9)$ of $\Pfr$ is
$(4, 5, 6, 8, 9)$. For this reason, the area of $\Pfr$ adjacent to
$(4, 9)$ is the $\Z$-bubble
\begin{equation}
    %
\,.
\end{equation}
\medbreak

\begin{Proposition} \label{prop:unique_decomposition_NC_M}
    Let $\Mca$ be a unitary magma and $\Pfr$ be a noncrossing
    $\Mca$-clique of arity greater than $1$. Then, there is a unique
    $\Mca$-bubble $\Qfr$ with a maximal arity $k \geq 2$ such that
    $\Pfr = \Qfr \circ [\Rfr_1, \dots, \Rfr_k]$, where each $\Rfr_i$,
    $i \in [k]$, is a noncrossing $\Mca$-clique with a base labeled
    by~$\Unit_\Mca$.
\end{Proposition}
\medbreak

Consider the map
\begin{equation}
    \BubbleTree :
    \NC\Mca \to \FreeOperad\left(\Bubbles_\Mca\right)
\end{equation}
defined linearly an recursively by $\BubbleTree(\UnitClique) := \Leaf$
and, for any noncrossing $\Mca$-clique $\Pfr$ of arity greater than~$1$,
by
\begin{equation}
    \BubbleTree(\Pfr) :=
    \Corolla{\Qfr} \circ
    \left[\BubbleTree(\Rfr_1), \dots, \BubbleTree(\Rfr_k)\right],
\end{equation}
where $\Pfr = \Qfr \circ [\Rfr_1, \dots, \Rfr_k]$ is the unique
decomposition of $\Pfr$ stated in
Proposition~\ref{prop:unique_decomposition_NC_M}. We call
$\BubbleTree(\Pfr)$ the \Def{bubble tree} of~$\Pfr$. For instance,
in~$\NC\Z$,
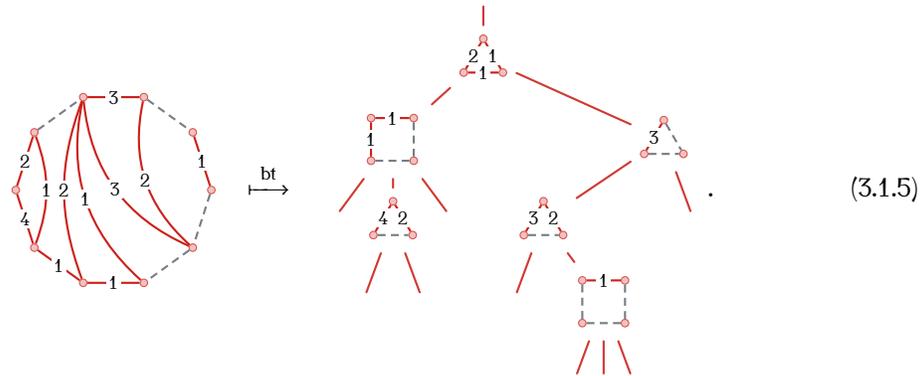
\begin{equation} \label{equ:bubble_tree_example}
};
        \draw[Edge](0)--(1);
        \draw[Edge](1)--(6);
        \draw[Edge](10)--(8);
        \draw[Edge](11)--(10);
        \draw[Edge](12)--(10);
        \draw[Edge](13)--(6);
        \draw[Edge](14)--(13);
        \draw[Edge](2)--(3);
        \draw[Edge](3)--(1);
        \draw[Edge](4)--(3);
        \draw[Edge](5)--(1);
        \draw[Edge](7)--(8);
        \draw[Edge](8)--(13);
        \draw[Edge](9)--(10);
        \node(r)at(5.00,2.25){};
        \draw[Edge](r)--(6);
    \end{tikzpicture}\,.
\end{equation}
\medbreak

\begin{Lemma} \label{lem:map_NC_M_bubble_tree_treelike_expression}
    Let $\Mca$ be a unitary magma. For any noncrossing $\Mca$-clique
    $\Pfr$, $\BubbleTree(\Pfr)$ is a treelike expression on
    $\Bubbles_\Mca$ of~$\Pfr$.
\end{Lemma}
\medbreak

\begin{Proposition} \label{prop:map_NC_M_bubble_tree}
    Let $\Mca$ be a unitary magma. Then, the map $\BubbleTree$ is
    injective and the image of $\BubbleTree$ is the linear span of all
    syntax trees $\Tfr$ on $\Bubbles_\Mca$ such that
    \begin{enumerate}[fullwidth,label={(\it\roman*)}]
        \item \label{item:map_NC_M_bubble_tree_1}
        the root of $\Tfr$ is labeled by an $\Mca$-bubble;
        \item \label{item:map_NC_M_bubble_tree_2}
        the internal nodes of $\Tfr$ different from the root are
        labeled by $\Mca$-bubbles whose bases are labeled
        by~$\Unit_\Mca$;
        \item \label{item:map_NC_M_bubble_tree_3}
        if $x$ and $y$ are two internal nodes of $\Tfr$ such that $y$ is
        the $i$th child of $x$, the $i$th edge of the bubble labeling
        $x$ is solid.
    \end{enumerate}
\end{Proposition}
\medbreak

Observe that $\BubbleTree$ is not an operad morphism. Indeed,
\begin{equation} \label{equ:bubble_tree_not_morphism}
    \BubbleTree\left(

    \right).
\end{equation}
Observe moreover that~\eqref{equ:bubble_tree_not_morphism} holds for all
unitary magmas $\Mca$ since $\Unit_\Mca$ is always idempotent.
\medbreak

\subsubsection{Realization in terms of decorated Schröder trees}
\label{subsubsec:M_Schroder_trees}
An \Def{$\Mca$-Schröder tree} $\Tfr$ is a Schröder tree (see
Section~\ref{subsubsec:schroder_trees} of
Chapter~\ref{chap:combinatorics}) such that each edge connecting two
internal nodes is labeled on $\bar{\Mca}$, each edge connecting an
internal node an a leaf is labeled on $\Mca$, and the outgoing edge from
the root of $\Tfr$ is labeled on $\Mca$
(see~\eqref{equ:example_M_Schroder_tree} for an example of a
$\Z$-Schröder tree).
\medbreak

From the description of the image of the map $\BubbleTree$ provided by
Proposition~\ref{prop:map_NC_M_bubble_tree}, any bubble tree $\Tfr$
of a noncrossing $\Mca$-clique $\Pfr$ of arity $n$ can be encoded by an
$\Mca$-Schröder tree $\Sfr$ with $n$ leaves. Indeed, this
$\Mca$-Schröder tree is obtained by considering each internal node $x$
of $\Tfr$ and by labeling the edge connecting $x$ and its $i$th child by
the label of the $i$th edge of the $\Mca$-bubble labeling $x$. The
outgoing edge from the root of $\Sfr$ is labeled by the label of the
base of the $\Mca$-bubble labeling the root of $\Tfr$. For instance, the
bubble tree of~\eqref{equ:bubble_tree_example} is encoded by the
$\Z$-Schröder tree
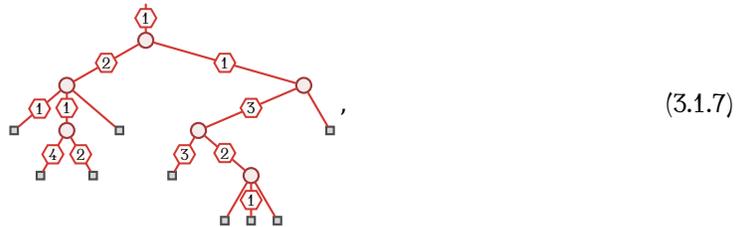
\begin{equation} \label{equ:example_M_Schroder_tree}
    \begin{tikzpicture}[xscale=.35,yscale=.2,Centering]
        \node[Leaf](0)at(0.00,-6.00){};
        \node[Leaf](11)at(9.00,-12.00){};
        \node[Leaf](12)at(10.00,-12.00){};
        \node[Leaf](14)at(12.00,-6.00){};
        \node[Leaf](2)at(1.00,-9.00){};
        \node[Leaf](4)at(3.00,-9.00){};
        \node[Leaf](5)at(4.00,-6.00){};
        \node[Leaf](7)at(6.00,-9.00){};
        \node[Leaf](9)at(8.00,-12.00){};
        \node[Node](1)at(2.00,-3.00){};
        \node[Node](10)at(9.00,-9.00){};
        \node[Node](13)at(11.00,-3.00){};
        \node[Node](3)at(2.00,-6.00){};
        \node[Node](6)at(5.00,0.00){};
        \node[Node](8)at(7.00,-6.00){};
        \draw[Edge](0)edge[]node[EdgeValue]{\begin{math}1\end{math}}(1);
        \draw[Edge](1)edge[]node[EdgeValue]{\begin{math}2\end{math}}(6);
        \draw[Edge](10)edge[]node[EdgeValue]{\begin{math}2\end{math}}(8);
        \draw[Edge](11)edge[]node[EdgeValue]
            {\begin{math}1\end{math}}(10);
        \draw[Edge](12)--(10);
        \draw[Edge](13)edge[]node[EdgeValue]
            {\begin{math}1\end{math}}(6);
        \draw[Edge](14)--(13);
        \draw[Edge](2)edge[]node[EdgeValue]{\begin{math}4\end{math}}(3);
        \draw[Edge](3)edge[]node[EdgeValue]{\begin{math}1\end{math}}(1);
        \draw[Edge](4)edge[]node[EdgeValue]{\begin{math}2\end{math}}(3);
        \draw[Edge](5)--(1);
        \draw[Edge](7)edge[]node[EdgeValue]{\begin{math}3\end{math}}(8);
        \draw[Edge](8)edge[]node[EdgeValue]
            {\begin{math}3\end{math}}(13);
        \draw[Edge](9)--(10);
        \node(r)at(5.00,3){};
        \draw[Edge](r)edge[]node[EdgeValue]{\begin{math}1\end{math}}(6);
    \end{tikzpicture}\,,
\end{equation}
where the labels of the edges are drawn in the hexagons and where
unlabeled edges are implicitly labeled by $\Unit_\Mca$. We shall use
these drawing conventions in the sequel. As a side remark, observe that
the $\Mca$-Schröder tree encoding a noncrossing $\Mca$-clique $\Pfr$ and
the dual tree of $\Pfr$ (in the usual meaning) have the same underlying
unlabeled tree.
\medbreak

This encoding of noncrossing $\Mca$-cliques by bubble trees is
reversible and hence, one can interpret $\NC\Mca$ as an operad on the
linear span of all $\Mca$-Schröder trees. Hence, through this
interpretation, if $\Sfr$ and $\Tfr$ are two $\Mca$-Schröder trees and
$i$ is a valid integer, the tree $\Sfr \circ_i \Tfr$ is computed by
grafting the root of $\Tfr$ to the $i$th leaf of $\Sfr$. Then, by
denoting by $b$ the label of the edge adjacent to the root of $\Tfr$ and
by $a$ the label of the edge adjacent to the $i$th leaf of $\Sfr$, we
have two cases to consider, depending on the value of $c := a \Op b$. If
$c \ne \Unit_\Mca$, we label the edge connecting $\Sfr$ and $\Tfr$ by
$c$. Otherwise, when $c = \Unit_\Mca$, we contract the edge connecting
$\Sfr$ and $\Tfr$ by merging the root of $\Tfr$ and the father of the
$i$th leaf of $\Sfr$ (see
Figure~\ref{fig:composition_NC_M_Schroder_trees}).
\begin{figure}[ht]
    \centering
    \subfloat[][The expression $\Sfr \circ_i \Tfr$ to compute. The
    displayed leaf is the $i$th one of~$\Sfr$.]{
    \begin{minipage}[c]{.8\textwidth}
        \centering
        \begin{equation*}

        \end{equation*}
    \end{minipage}
    \label{subfig:composition_NC_M_Schroder_trees_1}}

    \subfloat[][The resulting tree when $a \Op b \ne \Unit_\Mca$.]{
    \begin{minipage}[c]{.4\textwidth}
        \centering
        %
    \end{minipage}
    \label{subfig:composition_NC_M_Schroder_trees_2}}
    \subfloat[][The resulting tree when $a \Op b = \Unit_\Mca$.]{
    \begin{minipage}[c]{.4\textwidth}
        \centering
        %
    \end{minipage}
    \label{subfig:composition_NC_M_Schroder_trees_3}}
    \caption[The partial composition of the ns operad $\NC\Mca$ in
    terms of $\Mca$-Schröder trees.]
    {The partial composition of $\NC\Mca$ realized on $\Mca$-Schröder
    trees. Here, the two
    cases~\protect\subref{subfig:composition_NC_M_Schroder_trees_2}
    and~\protect\subref{subfig:composition_NC_M_Schroder_trees_3}
    for the computation of $\Sfr \circ_i \Tfr$
    are shown, where $\Sfr$ and $\Tfr$ are two $\Mca$-Schröder trees.
    In these drawings, the triangles denote subtrees.}
    \label{fig:composition_NC_M_Schroder_trees}
\end{figure}
For instance, in $\NC\N_3$, one has the two partial compositions
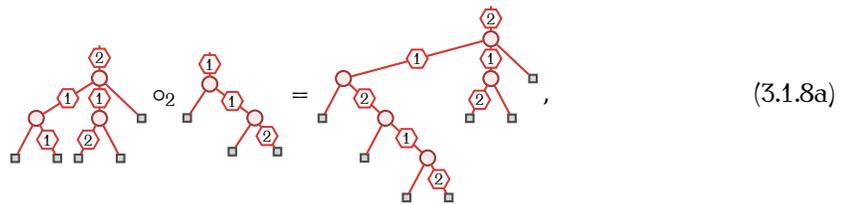
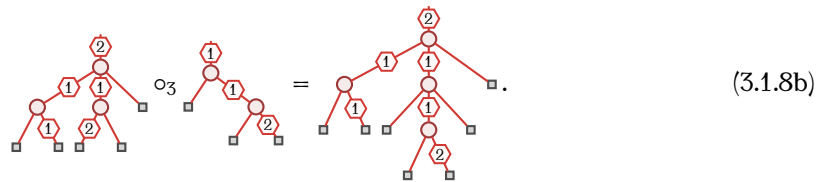
\begin{subequations}
\begin{equation}
\,.
\end{equation}
\end{subequations}
\medbreak

In the sequel, we shall indifferently see $\NC\Mca$ as an operad on
noncrossing $\Mca$-cliques or on $\Mca$-Schröder trees.
\medbreak

\subsubsection{Minimal generating set}

\begin{Proposition} \label{prop:generating_set_NC_M}
    Let $\Mca$ be a unitary magma. The set $\Triangles_\Mca$ of all
    $\Mca$-triangles is a minimal generating set of~$\NC\Mca$.
\end{Proposition}
\medbreak

Proposition~\ref{prop:generating_set_NC_M} also says that $\NC\Mca$ is
the smallest suboperad of $\Cli\Mca$ that contains all $\Mca$-triangles
and that $\NC\Mca$ is the biggest binary suboperad of~$\Cli\Mca$.
\medbreak

\subsubsection{Dimensions}
We now use the notion of bubble trees introduced in
Section~\ref{subsubsec:treelike_bubbles} to compute the dimensions
of~$\NC\Mca$.
\medbreak

\begin{Proposition} \label{prop:Hilbert_series_NC_M}
    Let $\Mca$ be a finite unitary magma. The Hilbert series
    $\HilbSeries_{\NC\Mca}(t)$ of $\NC\Mca$ satisfies
    \begin{equation} \label{equ:Hilbert_series_NC_M}
        t + \left(m^3 - 2m^2 + 2m - 1\right)t^2
        + \left(2m^2t - 3mt + 2t - 1\right) \HilbSeries_{\NC\Mca}(t)
        + \left(m - 1\right) \HilbSeries_{\NC\Mca}(t)^2
        = 0,
    \end{equation}
    where $m := \# \Mca$.
\end{Proposition}
\medbreak

We deduce from Proposition~\ref{prop:Hilbert_series_NC_M} that the
Hilbert series of $\NC\Mca$ satisfies
\begin{equation} \label{equ:Hilbert_series_NC_M_function}
    \HilbSeries_{\NC\Mca}(t) =
    \frac{1 - (2m^2 - 3m + 2)t - \sqrt{1 - 2(2m^2 - m)t + m^2t^2}}
    {2(m - 1)},
\end{equation}
where $m := \# \Mca \ne 1$.
\medbreak

By using Narayana numbers, whose definition is recalled in
Section~\ref{subsubsec:quotient_Cli_M_Nes}, one can state the following
result.
\medbreak

\begin{Proposition} \label{prop:dimensions_NC_M}
    Let $\Mca$ be a finite unitary magma. For all $n \geq 2$,
    \begin{equation} \label{equ:dimensions_NC_M}
        \dim \NC\Mca(n) =
        \sum_{0 \leq k \leq n - 2}
            m^{n + k + 1} (m - 1)^{n - k - 2} \;
            \Narayana(n, k),
    \end{equation}
    where $m := \# \Mca$.
\end{Proposition}
\medbreak

We can use Proposition~\ref{prop:dimensions_NC_M} to compute the
first dimensions of $\NC\Mca$. For instance, depending on
$m := \# \Mca$, we have the following sequences of dimensions:
\begin{subequations}
\begin{equation}
    1, 1, 1, 1, 1, 1, 1, 1,
    \qquad m = 1,
\end{equation}
\begin{equation}
    1, 8, 48, 352, 2880, 25216, 231168, 2190848,
    \qquad m = 2,
\end{equation}
\begin{equation}
    1, 27, 405, 7533, 156735, 349263, 81520425, 1967414265,
    \qquad m = 3,
\end{equation}
\begin{equation}
    1, 64, 1792, 62464, 2437120, 101859328, 4459528192, 201889939456.
    \qquad m = 4,
\end{equation}
\end{subequations}
The second one forms, except for the first terms,
Sequence~\OEIS{A054726} of~\cite{Slo}. The last two sequences are not
listed in~\cite{Slo} at this time.
\medbreak

\subsection{Presentation and Koszulity}
The aim of this section is to establish a presentation by generators and
relations of $\NC\Mca$. For this, we will define an adequate rewrite
rule on the set of the syntax trees on $\Triangles_\Mca$ and prove that
it admits the required properties.
\medbreak

\subsubsection{Space of relations}
\label{subsubsec:space_of_relations_NC_M}
Let $\RelationSpace_{\NC\Mca}$ be the subspace of
$\FreeOperad\left(\Triangles_\Mca\right)(3)$ generated by the elements
\begin{subequations}
\begin{equation} \label{equ:relation_1_NC_M}
    \Corolla{\Triangle{\Pfr_0}{\Pfr_1}{\Pfr_2}}
    \circ_1
    \Corolla{\Triangle{\Qfr_0}{\Qfr_1}{\Qfr_2}}
    -
    \Corolla{\Triangle{\Pfr_0}{\Rfr_1}{\Pfr_2}}
    \circ_1
    \Corolla{\Triangle{\Rfr_0}{\Qfr_1}{\Qfr_2}},
    \qquad
    \mbox{if } \Pfr_1 \Op \Qfr_0 = \Rfr_1 \Op \Rfr_0 \ne \Unit_\Mca,
\end{equation}
\begin{equation} \label{equ:relation_2_NC_M}
    \Corolla{\Triangle{\Pfr_0}{\Pfr_1}{\Pfr_2}}
    \circ_1
    \Corolla{\Triangle{\Qfr_0}{\Qfr_1}{\Qfr_2}}
    -
    \Corolla{\Triangle{\Pfr_0}{\Qfr_1}{\Rfr_2}}
    \circ_2
    \Corolla{\Triangle{\Rfr_0}{\Qfr_2}{\Pfr_2}},
    \qquad
    \mbox{if } \Pfr_1 \Op \Qfr_0 = \Rfr_2 \Op \Rfr_0 = \Unit_\Mca,
\end{equation}
\begin{equation} \label{equ:relation_3_NC_M}
    \Corolla{\Triangle{\Pfr_0}{\Pfr_1}{\Pfr_2}}
    \circ_2
    \Corolla{\Triangle{\Qfr_0}{\Qfr_1}{\Qfr_2}}
    -
    \Corolla{\Triangle{\Pfr_0}{\Pfr_1}{\Rfr_2}}
    \circ_2
    \Corolla{\Triangle{\Rfr_0}{\Qfr_1}{\Qfr_2}},
    \qquad
    \mbox{if } \Pfr_2 \Op \Qfr_0 = \Rfr_2 \Op \Rfr_0 \ne \Unit_\Mca,
\end{equation}
\end{subequations}
where $\Pfr$, $\Qfr$, and $\Rfr$ are $\Mca$-triangles.
\medbreak

\begin{Lemma} \label{lem:quadratic_relations_NC_M}
    Let $\Mca$ be a unitary magma, and $\Sfr$ and $\Tfr$ be two syntax
    trees of arity $3$ on $\Triangles_\Mca$. Then, $\Sfr - \Tfr$ belongs
    to $\RelationSpace_{\NC\Mca}$ if and only
    if~$\Eval(\Sfr) = \Eval(\Tfr)$.
\end{Lemma}
\medbreak

\begin{Proposition} \label{prop:dimensions_relations_NC_M}
    Let $\Mca$ be a finite unitary magma. Then, the dimension of the
    space $\RelationSpace_{\NC\Mca}$ satisfies
    \begin{equation}
        \dim \RelationSpace_{\NC\Mca} = 2m^6 - 2m^5 + m^4,
    \end{equation}
    where $m := \# \Mca$.
\end{Proposition}
\medbreak

Observe that, by Proposition~\ref{prop:dimensions_relations_NC_M}, the
dimension of $\RelationSpace_{\NC\Mca}$ only depends on the
cardinality of $\Mca$ and not on its operation~$\Op$.
\medbreak

\subsubsection{Rewrite rule} \label{subsubsec:rewrite_rule_NC_M}
Let $\Rew_\Mca$ be the rewrite rule on the set of the
$\Triangles_\Mca$-syntax trees on satisfying
\begin{subequations}
\begin{equation} \label{equ:rewrite_1_NC_M}
    \Corolla{\Triangle{\Pfr_0}{\Pfr_1}{\Pfr_2}}
    \circ_1
    \Corolla{\Triangle{\Qfr_0}{\Qfr_1}{\Qfr_2}}
    \Rew_\Mca
    \Corolla{\Triangle{\Pfr_0}{\delta}{\Pfr_2}}
    \circ_1
    \Corolla{\TriangleEXX{\Unit_\Mca}{\Qfr_1}{\Qfr_2}},
    \qquad
    \mbox{if } \Qfr_0 \ne \Unit_\Mca,
    \mbox{ where } \delta := \Pfr_1 \Op \Qfr_0,
\end{equation}
\begin{equation} \label{equ:rewrite_2_NC_M}
    \Corolla{\Triangle{\Pfr_0}{\Pfr_1}{\Pfr_2}}
    \circ_1
    \Corolla{\Triangle{\Qfr_0}{\Qfr_1}{\Qfr_2}}
    \Rew_\Mca
    \Corolla{\TriangleXXE{\Pfr_0}{\Qfr_1}{\Unit_\Mca}}
    \circ_2
    \Corolla{\TriangleEXX{\Unit_\Mca}{\Qfr_2}{\Pfr_2}},
    \qquad
    \mbox{if } \Pfr_1 \Op \Qfr_0 = \Unit_\Mca,
\end{equation}
\begin{equation} \label{equ:rewrite_3_NC_M}
    \Corolla{\Triangle{\Pfr_0}{\Pfr_1}{\Pfr_2}}
    \circ_2
    \Corolla{\Triangle{\Qfr_0}{\Qfr_1}{\Qfr_2}}
    \Rew_\Mca
    \Corolla{\Triangle{\Pfr_0}{\Pfr_1}{\delta}}
    \circ_2
    \Corolla{\TriangleEXX{\Unit_\Mca}{\Qfr_1}{\Qfr_2}},
    \qquad
    \mbox{if } \Qfr_0 \ne \Unit_\Mca,
    \mbox{ where } \delta := \Pfr_2 \Op \Qfr_0,
\end{equation}
\end{subequations}
where $\Pfr$ and $\Qfr$ are $\Mca$-triangles. Let also $\RewTrees_\Mca$
be the closure of $\Rew_\Mca$.
\medbreak

\begin{Proposition} \label{prop:rewrite_rule_NC_M}
    Let $\Mca$ be a finite unitary magma. Then, $\RewTrees_\Mca$ is a
    convergent rewrite rule and an orientation
    of~$\RelationSpace_{\NC\Mca}$.
\end{Proposition}
\medbreak

\begin{Lemma} \label{lem:rewrite_rule_NC_M_normal_forms}
    Let $\Mca$ be a unitary magma. The set of the normal forms of
    $\RewTrees_\Mca$ is the set of the $\Triangles_\Mca$-syntax trees
    $\Tfr$ such that, for any internal nodes $x$ and $y$ of $\Tfr$
    where $y$ is a child of $x$,
    \begin{enumerate}[fullwidth,label={(\it\roman*)}]
        \item \label{item:rewrite_rule_NC_M_normal_forms_1}
        the base of the $\Mca$-triangle labeling $y$ is labeled
        by~$\Unit_\Mca$;
        \item \label{item:rewrite_rule_NC_M_normal_forms_2}
        if $y$ is a left child of $x$, the first edge of the
        $\Mca$-triangle labeling $x$ is not labeled by~$\Unit_\Mca$.
    \end{enumerate}
    Moreover, when $\Mca$ is finite, the generating series of the
    normal forms of $\RewTrees_\Mca$ is the Hilbert series
    $\HilbSeries_{\NC\Mca(t)}$ of~$\NC\Mca$.
\end{Lemma}
\medbreak

\subsubsection{Presentation and Koszulity}
The results of Sections~\ref{subsubsec:space_of_relations_NC_M}
and~\ref{subsubsec:rewrite_rule_NC_M} are finally used here to state a
presentation of $\NC\Mca$ and the fact that $\NC\Mca$ is a Koszul
operad.
\medbreak

\begin{Theorem} \label{thm:presentation_NC_M}
    Let $\Mca$ be a finite unitary magma. Then, $\NC\Mca$ admits the
    presentation~%
    \begin{math}
        \left(\Triangles_\Mca, \RelationSpace_{\NC\Mca}\right)
    \end{math}.
\end{Theorem}
\begin{proof}
    First, Proposition~\ref{prop:rewrite_rule_NC_M} implies that
    we can regard the underlying space of the quotient operad
    \begin{equation}
        \Oca :=
        \FreeOperad\left(\Triangles_\Mca\right)/_%
        {\langle\RelationSpace_{\NC\Mca}\rangle}
    \end{equation}
    as the linear span of all normal forms of $\RewTrees_\Mca$.
    Moreover, as a consequence of
    Lemma~\ref{lem:quadratic_relations_NC_M}, the map
    \begin{math}
        \phi : \Oca \to \NC\Mca
    \end{math}
    defined linearly for any normal form $\Tfr$ of $\RewTrees_\Mca$ by
    $\phi(\Tfr) := \Eval(\Tfr)$ is an operad morphism. Now, by
    Proposition~\ref{prop:generating_set_NC_M}, $\phi$ is surjective.
    Moreover, by Lemma~\ref{lem:rewrite_rule_NC_M_normal_forms}, we
    obtain that the dimensions of the spaces $\Oca(n)$, $n \geq 1$, are
    the ones of $\NC\Mca(n)$. Hence, $\phi$ is an operad isomorphism
    and the statement of the theorem follows.
\end{proof}
\medbreak

By Theorem~\ref{thm:presentation_NC_M}, the operad $\NC\N_2$ is
generated by
\begin{equation}
    \Triangles_{\N_2} =
    \left\{
    \TriangleEEE{}{}{},
    \TriangleEXE{}{1}{},
    \TriangleEEX{}{}{1},
    \TriangleEXX{}{1}{1},
    \TriangleXEE{1}{}{},
    \TriangleXXE{1}{1}{},
    \TriangleXEX{1}{}{1},
    \Triangle{1}{1}{1}
    \right\},
\end{equation}
and these generators are subjected exactly to the nontrivial relations
\begin{subequations}
\begin{equation}
    \TriangleXEX{a}{}{b_3}
    \circ_1
    \Triangle{1}{b_1}{b_2}
    =
    \Triangle{a}{1}{b_3}
    \circ_1
    \TriangleEXX{}{b_1}{b_2},
    \qquad
    a, b_1, b_2, b_3 \in \N_2,
\end{equation}
\begin{equation}
    \Triangle{a}{1}{b_3}
    \circ_1
    \Triangle{1}{b_1}{b_2}
    =
    \TriangleXEX{a}{}{b_3}
    \circ_1
    \TriangleEXX{}{b_1}{b_2}
    =
    \TriangleXXE{a}{b_1}{}
    \circ_2
    \TriangleEXX{}{b_2}{b_3}
    =
    \Triangle{a}{b_1}{1}
    \circ_2
    \Triangle{1}{b_2}{b_3},
    \qquad
    a, b_1, b_2, b_3 \in \N_2,
\end{equation}
\begin{equation}
    \TriangleXXE{a}{b_1}{}
    \circ_2
    \Triangle{1}{b_2}{b_3}
    =
    \Triangle{a}{b_1}{1}
    \circ_2
    \TriangleEXX{}{b_2}{b_3},
    \qquad
    a, b_1, b_2, b_3 \in \N_2.
\end{equation}
\end{subequations}
\medbreak

\begin{Theorem} \label{thm:Koszul_NC_M}
    For any finite unitary magma $\Mca$, $\NC\Mca$ is Koszul and the set
    of the normal forms of $\RewTrees_\Mca$ forms a
    Poincaré-Birkhoff-Witt basis of~$\NC\Mca$.
\end{Theorem}
\medbreak

\subsection{Suboperads generated by bubbles}
In this section, we consider suboperads of $\NC\Mca$ generated by finite
sets of $\Mca$-bubbles. We assume here that $\Mca$ is endowed with an
arbitrary total order so that $\Mca = \{x_0, x_1, \dots\}$ with
$x_0 = \Unit_\Mca$.
\medbreak

\subsubsection{Treelike expressions on bubbles}
Let $B$ and $E$ be two subsets of $\Mca$. We denote by
$\Bubbles_\Mca^{B, E}$ the set of all $\Mca$-bubbles $\Pfr$ such that
the bases of $\Pfr$ are labeled on $B$ and all edges of $\Pfr$ are
labeled on $E$. Moreover, we say that $\Mca$ is
\Def{$(E, B)$-quasi-injective} if for all $x, x' \in E$ and
$y, y' \in B$, $x \Op y = x' \Op y' \ne \Unit_\Mca$ implies $x = x'$ and
$y = y'$.
\medbreak

\begin{Lemma} \label{lem:treelike_expression_suboperad_bubbles}
    Let $\Mca$ be a unitary magma, and $B$ and $E$ be two subsets of
    $\Mca$. If $\Mca$ is $(E, B)$-quasi-injective, then any
    $\Mca$-clique admits at most one treelike expression on
    $\Bubbles_\Mca^{B, E}$ of a minimal degree.
\end{Lemma}
\medbreak

\subsubsection{Dimensions}
\label{subsubsec:dimensions_suboperads_triangles}
Let $G$ be a set of $\Mca$-bubbles and
$\Xi := \left\{\xi_{x_0}, \xi_{x_1}, \dots\right\}$ be a set of
noncommutative variables. Given $x_i \in \Mca$, let
$\SeriesBubbles_{x_i}$ be the series of
$\N \langle\langle \Xi \rangle\rangle$ defined by
\begin{equation} \label{equ:series_bubbles}
    \SeriesBubbles_{x_i}\left(\xi_{x_0}, \xi_{x_1}, \dots\right) :=
    \sum_{\substack{
        \Pfr \in \Bubbles_\Mca^G \\
        \Pfr \ne \UnitClique}}
    \enspace
    \prod_{i \in [|\Pfr|]}
    \xi_{\Pfr_i},
\end{equation}
where $\Bubbles_\Mca^G$ is the set of all $\Mca$-bubbles that can be
obtained by partial compositions of elements of $G$. Observe
from~\eqref{equ:series_bubbles} that a noncommutative monomial
$u \in \Xi^{\geq 2}$ appears in $\SeriesBubbles_{x_i}$ with $1$ as
coefficient if and only if there is in the suboperad of $\NC\Mca$
generated by $G$ an $\Mca$-bubble with a base labeled by $x_i$ and with
$u$ as border.
\medbreak

Let also for any $x_i \in \Mca$, the series $\SeriesElements_{x_i}$ of
$\N \langle\langle t \rangle\rangle$ defined by
\begin{equation}
    \SeriesElements_{x_i}(t) :=
    \SeriesBubbles_{x_i}\left(t + \bar \SeriesElements_{x_0}(t),
    t + \bar \SeriesElements_{x_1}(t), \dots\right),
\end{equation}
where for any $x_i \in \Mca$,
\begin{equation} \label{equ:bar_series_elements_based}
    \bar \SeriesElements_{x_i}(t) :=
    \sum_{\substack{
        x_j \in \Mca \\
        x_i \Op x_j \ne \Unit_\Mca
    }}
    \SeriesElements_{x_j}(t).
\end{equation}
\medbreak

\begin{Proposition} \label{prop:suboperads_NC_M_triangles_dimensions}
    Let $\Mca$ be a unitary magma and $\GeneratingSet$ be a finite set
    of $\Mca$-bubbles such that, by denoting by $B$ (resp. $E$) the set
    of the labels of the bases (resp. edges) of the elements of
    $\GeneratingSet$, $\Mca$ is $(E, B)$-quasi-injective. Then, the
    Hilbert series $\HilbSeries_{(\NC\Mca)^\GeneratingSet}(t)$ of the
    suboperad of $\NC\Mca$ generated by $\GeneratingSet$ satisfies
    \begin{equation} \label{equ:suboperads_NC_M_triangles_dimensions}
        \HilbSeries_{(\NC\Mca)^\GeneratingSet}(t) =
        t +
        \sum_{x_i \in \Mca} \SeriesElements_{x_i}(t).
    \end{equation}
\end{Proposition}
\medbreak

As a side remark,
Proposition~\ref{prop:suboperads_NC_M_triangles_dimensions} can be
proved by using the notion of bubble decompositions of operads developed
in Chapter~\ref{chap:enveloping}. This result provides a practical
method to compute the dimensions of some suboperads
$(\NC\Mca)^\GeneratingSet$ of $\NC\Mca$ by describing the
series~\eqref{equ:series_bubbles} of the bubbles of
$\Bubbles_\Mca^\GeneratingSet$. This result implies also, when
$\GeneratingSet$ satisfies the requirement of
Proposition~\ref{prop:suboperads_NC_M_triangles_dimensions}, that the
Hilbert series of $(\NC\Mca)^\GeneratingSet$ is algebraic.
\medbreak

\subsubsection{First example: a cubic suboperad}
Consider the suboperad of $\NC\Ebb_2$ generated by
\begin{equation}
    \GeneratingSet := \left\{
        \TriangleXEX{\Esf_1}{}{\Esf_1},
        \TriangleXEX{\Esf_2}{}{\Esf_2}
    \right\}.
\end{equation}
Computer experiments show that the generators of
$(\NC\Ebb_2)^\GeneratingSet$ are not subjected to any quadratic relation
but are subjected to the four cubic nontrivial relations
\begin{subequations}
\begin{equation} \label{equ:rel_suboperad_1_1}
    \TriangleXEX{\Esf_1}{}{\Esf_1}
    \circ_2
    \left(
    \TriangleXEX{\Esf_1}{}{\Esf_1}
    \circ_2
    \TriangleXEX{\Esf_1}{}{\Esf_1}
    \right)
    =
    \TriangleXEX{\Esf_1}{}{\Esf_1}
    \circ_2
    \left(
    \TriangleXEX{\Esf_2}{}{\Esf_2}
    \circ_2
    \TriangleXEX{\Esf_1}{}{\Esf_1}
    \right),
\end{equation}
\begin{equation} \label{equ:rel_suboperad_1_2}
    \TriangleXEX{\Esf_1}{}{\Esf_1}
    \circ_2
    \left(
    \TriangleXEX{\Esf_1}{}{\Esf_1}
    \circ_2
    \TriangleXEX{\Esf_2}{}{\Esf_2}
    \right)
    =
    \TriangleXEX{\Esf_1}{}{\Esf_1}
    \circ_2
    \left(
    \TriangleXEX{\Esf_2}{}{\Esf_2}
    \circ_2
    \TriangleXEX{\Esf_2}{}{\Esf_2}
    \right),
\end{equation}
\begin{equation} \label{equ:rel_suboperad_1_3}
    \TriangleXEX{\Esf_2}{}{\Esf_2}
    \circ_2
    \left(
    \TriangleXEX{\Esf_1}{}{\Esf_1}
    \circ_2
    \TriangleXEX{\Esf_1}{}{\Esf_1}
    \right)
    =
    \TriangleXEX{\Esf_2}{}{\Esf_2}
    \circ_2
    \left(
    \TriangleXEX{\Esf_2}{}{\Esf_2}
    \circ_2
    \TriangleXEX{\Esf_1}{}{\Esf_1}
    \right),
\end{equation}
\begin{equation} \label{equ:rel_suboperad_1_4}
    \TriangleXEX{\Esf_2}{}{\Esf_2}
    \circ_2
    \left(
    \TriangleXEX{\Esf_1}{}{\Esf_1}
    \circ_2
    \TriangleXEX{\Esf_2}{}{\Esf_2}
    \right)
    =
    \TriangleXEX{\Esf_2}{}{\Esf_2}
    \circ_2
    \left(
    \TriangleXEX{\Esf_2}{}{\Esf_2}
    \circ_2
    \TriangleXEX{\Esf_2}{}{\Esf_2}
    \right).
\end{equation}
\end{subequations}
Hence, $(\NC\Ebb_2)^\GeneratingSet$ is not a quadratic operad. Moreover,
it is possible to prove that this operad does not admit any other
nontrivial relations between its generators. This can be performed by
defining a rewrite rule on the syntax trees on $\GeneratingSet$,
consisting in rewriting the left patterns
of~\eqref{equ:rel_suboperad_1_1}, \eqref{equ:rel_suboperad_1_2},
\eqref{equ:rel_suboperad_1_3}, and~\eqref{equ:rel_suboperad_1_4} into
their respective right patterns, and by checking that this rewrite rule
admits the required properties (like the ones establishing the
presentation of $\NC\Mca$ by Theorem~\ref{thm:presentation_NC_M}). The
existence of this nonquadratic operad shows that $\NC\Mca$ contains
nonquadratic suboperads even if it is quadratic.
\medbreak
By describing the bubbles of $(\NC\Ebb_2)^\GeneratingSet$,
Proposition~\ref{prop:suboperads_NC_M_triangles_dimensions} leads to
the fact that the Hilbert series of $(\NC\Ebb_2)^\GeneratingSet$
satisfies the algebraic equation
\begin{equation}
    t + (t - 1) \HilbSeries_{(\NC\Ebb_2)^\GeneratingSet}(t) +
    (2t + 1)\HilbSeries_{(\NC\Ebb_2)^\GeneratingSet}(t)^2 = 0.
\end{equation}
The first dimensions of $(\NC\Ebb_2)^\GeneratingSet$ are
\begin{equation}
    1, 2, 8, 36, 180, 956, 5300, 30316,
\end{equation}
and form Sequence~\OEIS{A129148} of~\cite{Slo}.
\medbreak

\subsubsection{Second example: a suboperad of Motzkin paths}
Consider the suboperad of $\NC\Dbb_0$ generated by
\begin{equation}
    \GeneratingSet := \left\{
        \TriangleEEE{}{}{},
        \SquareMotz
    \right\}.
\end{equation}
Computer experiments show that the generators of
$(\NC\Dbb_0)^\GeneratingSet$ are subjected to four quadratic nontrivial
relations
\begin{subequations}
\begin{equation} \label{equ:rel_suboperad_2_1}
    \TriangleEEE{}{}{} \circ_1 \TriangleEEE{}{}{}
    =
    \TriangleEEE{}{}{} \circ_2 \TriangleEEE{}{}{}\,,
\end{equation}
\begin{equation} \label{equ:rel_suboperad_2_2}
    \SquareMotz \circ_1 \TriangleEEE{}{}{}
    =
    \TriangleEEE{}{}{} \circ_2 \SquareMotz\,,
\end{equation}
\begin{equation} \label{equ:rel_suboperad_2_3}
    \TriangleEEE{}{}{} \circ_1 \SquareMotz
    =
    \SquareMotz \circ_2 \TriangleEEE{}{}{}\,,
\end{equation}
\begin{equation} \label{equ:rel_suboperad_2_4}
    \SquareMotz \circ_1 \SquareMotz
    =
    \SquareMotz \circ_3 \SquareMotz\,.
\end{equation}
\end{subequations}
It is possible to prove that this operad does not admit any other
nontrivial relations between its generators. This can be performed by
defining a rewrite rule on the syntax trees on $\GeneratingSet$,
consisting in rewriting the left patterns
of~\eqref{equ:rel_suboperad_2_1}, \eqref{equ:rel_suboperad_2_2},
\eqref{equ:rel_suboperad_2_3}, and~\eqref{equ:rel_suboperad_2_4} into
their respective right patterns, and by checking that this rewrite rule
admits the required properties (like the ones establishing the
presentation of $\NC\Mca$ by Theorem~\ref{thm:presentation_NC_M}).
\medbreak

By describing the bubbles of $(\NC\Dbb_0)^\GeneratingSet$,
Proposition~\ref{prop:suboperads_NC_M_triangles_dimensions} leads to
the fact that the Hilbert series of $(\NC\Dbb_0)^\GeneratingSet$
satisfies the algebraic equation
\begin{equation}
    t + (t - 1) \HilbSeries_{(\NC\Dbb_0)^\GeneratingSet}(t)
    + t\HilbSeries_{(\NC\Dbb_0)^\GeneratingSet}(t)^2 = 0.
\end{equation}
The first dimensions of $(\NC\Dbb_0)^\GeneratingSet$ are
\begin{equation}
    1, 1, 2, 4, 9, 21, 51, 127,
\end{equation}
and form Sequence~\OEIS{A001006} of~\cite{Slo}. The operad
$(\NC\Dbb_0)^\GeneratingSet$ has the same presentation by generators and
relations (and thus, the same Hilbert series) as the operad $\Motz$
defined in Section~\ref{subsubsec:operad_Motz} of
Chapter~\ref{chap:monoids}, involving Motzkin paths. Hence,
$(\NC\Dbb_0)^\GeneratingSet$ and $\Motz$ are two isomorphic operads.
Note in passing that these two operads are not isomorphic to the operad
$\Motzkin\Dbb_0$ constructed in
Section~\ref{subsubsec:Motzkin_configurations} and involving Motzkin
configurations. Indeed, the sequence of the dimensions of this last
operad is a shifted version of the one of $(\NC\Dbb_0)^\GeneratingSet$
and~$\Motz$.
\medbreak

\subsection{Algebras over the noncrossing clique operads}
We begin by briefly describing $\NC\Mca$-algebras in terms of relations
between their operations and the free $\NC\Mca$-algebras over one
generator. We continue this section by providing two ways to construct
(non-necessarily free) $\NC\Mca$-algebras. The first one takes as input
an associative algebra endowed with endofunctions satisfying some
conditions, and the second one takes as input a monoid.
\medbreak

\subsubsection{Relations}
From the presentation of $\NC\Mca$ established by
Theorem~\ref{thm:presentation_NC_M}, an $\NC\Mca$-al\-gebra is a vector
space $\Alg$ endowed with binary linear operations
\begin{equation}
    \TriangleOp{\Pfr_0}{\Pfr_1}{\Pfr_2} :
    \Alg \otimes \Alg \to \Alg,
    \qquad
    \Pfr \in \Triangles_\Mca,
\end{equation}
satisfying, for all $a_1, a_2, a_3 \in \Alg$, the relations
\begin{subequations}
\begin{equation} \label{equ:relation_NC_M_algebras_1}
    \left(a_1 \TriangleOp{\Qfr_0}{\Qfr_1}{\Qfr_2} a_2\right)
    \TriangleOp{\Pfr_0}{\Pfr_1}{\Pfr_2} a_3
    =
    \left(a_1 \TriangleOp{\Rfr_0}{\Qfr_1}{\Qfr_2} a_2\right)
    \TriangleOp{\Pfr_0}{\Rfr_1}{\Pfr_2} a_3,
    \qquad
    \mbox{if } \Pfr_1 \Op \Qfr_0 = \Rfr_1 \Op \Rfr_0 \ne \Unit_\Mca,
\end{equation}
\begin{equation} \label{equ:relation_NC_M_algebras_2}
    \left(a_1 \TriangleOp{\Qfr_0}{\Qfr_1}{\Qfr_2} a_2\right)
    \TriangleOp{\Pfr_0}{\Pfr_1}{\Pfr_2} a_3
    =
    a_1 \TriangleOp{\Pfr_0}{\Qfr_1}{\Rfr_2}
    \left(a_2 \TriangleOp{\Rfr_0}{\Qfr_2}{\Pfr_2} a_3 \right),
    \qquad
    \mbox{if } \Pfr_1 \Op \Qfr_0 = \Rfr_2 \Op \Rfr_0 = \Unit_\Mca,
\end{equation}
\begin{equation} \label{equ:relation_NC_M_algebras_3}
    a_1 \TriangleOp{\Pfr_0}{\Pfr_1}{\Pfr_2}
    \left( a_2 \TriangleOp{\Qfr_0}{\Qfr_1}{\Qfr_2} a_3 \right)
    =
    a_1 \TriangleOp{\Pfr_0}{\Pfr_1}{\Rfr_2}
    \left( a_2 \TriangleOp{\Rfr_0}{\Qfr_1}{\Qfr_2} a_3 \right),
    \qquad
    \mbox{if } \Pfr_2 \Op \Qfr_0 = \Rfr_2 \Op \Rfr_0 \ne \Unit_\Mca,
\end{equation}
\end{subequations}
where $\Pfr$, $\Qfr$, and $\Rfr$ are $\Mca$-triangles. Remark that
$\Mca$ has to be finite because Theorem~\ref{thm:presentation_NC_M}
requires this property as premise.
\medbreak

\subsubsection{Free algebras over one generator}
From the realization of $\NC\Mca$ coming from its definition as a
suboperad of $\Cli\Mca$, the free $\NC\Mca$-algebra over one generator
is the linear span $\NC\Mca$ of all noncrossing $\Mca$-cliques endowed
with the linear operations
\begin{equation}
    \TriangleOp{\Pfr_O}{\Pfr_1}{\Pfr_2} :
    \NC\Mca(n) \otimes \NC\Mca(m) \to \NC\Mca(n + m),
    \qquad
    \Pfr \in \Triangles_\Mca,
    n, m \geq 1,
\end{equation}
defined, for any noncrossing $\Mca$-cliques $\Qfr$ and $\Rfr$, by
\begin{equation} \label{equ:free_algebra_MT_M_product}
    \Qfr \TriangleOp{\Pfr_0}{\Pfr_1}{\Pfr_2} \Rfr
    :=
    \left(\Triangle{\Pfr_0}{\Pfr_1}{\Pfr_2}
    \circ_2 \Rfr\right) \circ_1 \Qfr.
\end{equation}
In terms of $\Mca$-Schröder trees (see
Section~\ref{subsubsec:M_Schroder_trees}),
\eqref{equ:free_algebra_MT_M_product} is the $\Mca$-Schröder tree
obtained by grafting the $\Mca$-Schröder trees of $\Qfr$ and $\Rfr$
respectively as left and right children of a binary corolla having its
edge adjacent to the root labeled by $\Pfr_0$, its first edge labeled by
$\Pfr_1 \Op \Qfr_0$, and second edge labeled by $\Pfr_2 \Op \Rfr_0$,
and by contracting each of these two edges when labeled by~$\Unit_\Mca$.
For instance, in the free $\NC\N_3$-algebra, we have
\begin{subequations}
\begin{equation}
\,.
\end{equation}
\end{subequations}
\medbreak

\subsubsection{From associative algebras}
Let $\Alg$ be an associative algebra with associative product denoted by
$\OpAssoc$, and
\begin{equation} \label{equ:compatible_set}
    \omega_x : \Alg \to \Alg, \qquad x \in \Mca,
\end{equation}
be a family of linear maps, not necessarily associative algebra
morphisms, indexed by the elements of $\Mca$. We say that $\Alg$
together with this family~\eqref{equ:compatible_set} of maps is a
\Def{$\Mca$-compatible algebra} if
\begin{equation} \label{equ:compatible_magma_on_algebra}
    \omega_x \circ \omega_y = \omega_{x \Op y},
\end{equation}
for all $x, y \in \Mca$. Observe
that~\eqref{equ:compatible_magma_on_algebra} implies in particular that
$\omega_{\Unit_\Mca} = \Identity_\Alg$ where $\Identity_\Alg$ is the
identity map on $\Alg$. This notion of $\Mca$-compatible algebras
is very similar to the notion of $\Mca$-compatible algebras where
$\Mca$ is a monoid, developed in Section~\ref{subsec:main_properties_T}
of Chapter~\ref{chap:monoids}. Let us now use $\Mca$-compatible
associative algebras to construct $\NC\Mca$-algebras.
\medbreak

\begin{Theorem} \label{thm:NC_M_algebras}
    Let $\Mca$ be a finite unitary magma and $\Alg$ be an
    $\Mca$-compatible associative algebra. The vector space $\Alg$
    endowed with the binary linear operations
    \begin{equation} \label{equ:NC_M_algebras}
        \TriangleOp{\Pfr_0}{\Pfr_1}{\Pfr_2} :
        \Alg \otimes \Alg \to \Alg,
        \qquad
        \Pfr \in \Triangles_\Mca,
    \end{equation}
    defined for each $\Mca$-triangle $\Pfr$ and any $a_1, a_2 \in \Alg$
    by
    \begin{equation} \label{equ:NC_M_algebras_def}
        a_1 \TriangleOp{\Pfr_0}{\Pfr_1}{\Pfr_2} a_2
        :=
        \omega_{\Pfr_0}\left(\omega_{\Pfr_1}\left(a_1\right)
        \OpAssoc \omega_{\Pfr_2}\left(a_2\right)\right),
    \end{equation}
    is an $\NC\Mca$-algebra.
\end{Theorem}
\medbreak

By Theorem~\ref{thm:NC_M_algebras}, $\Alg$ has the structure of an
$\NC\Mca$-algebra. Hence, there is a left action~$\cdot$ of the operad
$\NC\Mca$ on the tensor algebra of $\Alg$ of the form
\begin{equation}
    \cdot : \NC\Mca(n) \otimes \Alg^{\otimes n} \to \Alg,
    \qquad n \geq 1,
\end{equation}
whose definition comes from the ones of the
operations~\eqref{equ:NC_M_algebras} and
Relation~\eqref{equ:algebra_over_operad} of Chapter~\ref{chap:algebra}.
We describe here an algorithm to compute the action of any element of
$\NC\Mca$ of arity $n$ on tensors $a_1 \otimes \dots \otimes a_n$ of
$\Alg^{\otimes n}$. First, if $\Bfr$ is an $\Mca$-bubble of arity $n$,
\begin{equation} \label{equ:NC_M_algebras_action_bubbles}
    \Bfr \cdot \left(a_1 \otimes \dots \otimes a_n\right)
    = \omega_{\Bfr_0}\left(\prod_{i \in [n]}
    \omega_{\Bfr_i}\left(a_i\right) \right),
\end{equation}
where the product
of~\eqref{equ:NC_M_algebras_action_bubbles} denotes the iterated
version of the associative product $\OpAssoc$ of~$\Alg$. When $\Pfr$ is
a noncrossing $\Mca$-clique of arity $n$, $\Pfr$ acts recursively
on $a_1 \otimes \dots \otimes a_n$ as follows. One has
\begin{equation}
    \Pfr \cdot a_1 = a_1
\end{equation}
when $\Pfr = \UnitClique$, and
\begin{equation} \label{equ:NC_M_algebras_action_cliques}
    \Pfr \cdot \left(a_1 \otimes \dots \otimes a_n\right) =
    \Bfr \cdot \left(
        \left(\Rfr_1 \cdot \left(a_1 \otimes \dots
                \otimes a_{|\Rfr_1|}\right)\right)
        \otimes \dots \otimes
        \left(\Rfr_k \cdot
            \left(a_{|\Rfr_1| + \dots + |\Rfr_{k - 1}| + 1}
            \otimes \dots  \otimes a_n\right)\right)
    \right),
\end{equation}
where, by setting $\Tfr$ as the bubble tree $\BubbleTree(\Pfr)$ of
$\Pfr$ (see Section~\ref{subsubsec:treelike_bubbles}), $\Bfr$ and
$\Rfr_1$, \dots, $\Rfr_k$ are the unique $\Mca$-bubble and noncrossing
$\Mca$-cliques such that
\begin{math}
    \Tfr = \Corolla{\Bfr}
    \circ [\BubbleTree(\Rfr_1), \dots, \BubbleTree(\Rfr_k)].
\end{math}
\medbreak

Here are few examples of the construction provided by
Theorem~\ref{thm:NC_M_algebras}.
\begin{description}[fullwidth]
    \item[Noncommutative polynomials and selected concatenation]
    Consider the unitary magma $\Sbb_\ell$ of all subsets of $[\ell]$
    with the union as product. Let $A := \{a_j : j \in [\ell]\}$ be an
    alphabet of noncommutative letters. We define on the associative
    algebra $\K \langle A \rangle$ of polynomials on $A$ the linear maps
    \begin{equation}
        \omega_S : \K \langle A \rangle \to \K \langle A \rangle,
        \qquad S \in \Sbb_\ell,
    \end{equation}
    as follows. For any $u \in A^*$ and $S \in \Sbb_\ell$, we set
    \begin{equation}
        \omega_S(u) :=
        \begin{cases}
            u & \mbox{if } |u|_{a_j} \geq 1 \mbox{ for all } j \in S, \\
            0 & \mbox{otherwise}.
        \end{cases}
    \end{equation}
    Since, for all $u \in A^*$, $\omega_{\emptyset}(u) = u$ and
    $(\omega_S \circ \omega_{S'})(u) = \omega_{S \cup S'}(u)$, and
    $\emptyset$ is the unit of $\Sbb_\ell$, we obtain from
    Theorem~\ref{thm:NC_M_algebras} that the
    operations~\eqref{equ:NC_M_algebras} endow $ \K \langle A \rangle$
    with an $\NC\Sbb_\ell$-algebra structure. For instance, when
    $\ell := 3$, one has
    \begin{subequations}
    \begin{equation}
        \left(a_1 + a_1 a_3 + a_2 a_2\right)
        \;

        \;
        \left(1 + a_3 + a_2 a_1\right)
        =
        2 \; a_1 a_3 + a_1 a_3 a_3 +
        a_1 a_3 a_2 a_1.
    \end{equation}
    \end{subequations}
    Besides, to compute the action
    \begin{equation} \label{equ:example_action_S_clique}
        %
        \cdot
        (f \otimes f \otimes f \otimes f \otimes f \otimes f \otimes f
        \otimes f)
    \end{equation}
    where $f := a_1 + a_2 + a_3$, we use the above algorithm
    and~\eqref{equ:NC_M_algebras_action_bubbles}
    and~\eqref{equ:NC_M_algebras_action_cliques}. By presenting the
    computation on the bubble tree of the noncrossing $\Sbb_3$-clique
    of~\eqref{equ:example_action_S_clique}, we obtain
    \begin{equation}
        %
            };
            \draw(0)edge[Edge]node[CliqueLabel]
                {\begin{math}
                    \textcolor{Col1}{f}
                \end{math}}(4);
            \draw(1)edge[Edge]node[CliqueLabel]
                {\begin{math}
                    \textcolor{Col1}{f}
                \end{math}}(2);
            \draw(10)edge[Edge]node[CliqueLabel]
                {\begin{math}
                    \textcolor{Col1}{f}
                \end{math}}(11);
            \draw(11)edge[Edge]node[CliqueLabel]
                {\begin{math}
                    \textcolor{Col1}{(a_1 + a_2 + a_3)^2}
                \end{math}}(9);
            \draw(12)edge[Edge]node[CliqueLabel]
                {\begin{math}
                    \textcolor{Col1}{f}
                \end{math}}(11);
            \draw(2)edge[Edge]node[CliqueLabel]
                {\begin{math}
                    \textcolor{Col1}{a_1 a_2}
                \end{math}}(4);
            \draw(3)edge[Edge]node[CliqueLabel]
                {\begin{math}
                    \textcolor{Col1}{f}
                \end{math}}(2);
            \draw(4)edge[Edge]node[CliqueLabel]
                {\begin{math}
                    \textcolor{Col1}
                    {(a_1 + a_2 + a_3) a_1 a_2 a_2 a_1 a_3}
                \end{math}}(9);
            \draw(5)edge[Edge]node[CliqueLabel]
                {\begin{math}
                    \textcolor{Col1}{f}
                \end{math}}(6);
            \draw(6)edge[Edge]node[CliqueLabel]
                {\begin{math}
                    \textcolor{Col1}{a_2 (a_1 + a_2 + a_3)\qquad}
                \end{math}}(4);
            \draw(7)edge[Edge]node[CliqueLabel]
                {\begin{math}
                    \textcolor{Col1}{f}
                \end{math}}(6);
            \draw(8)edge[Edge]node[CliqueLabel]
                {\begin{math}
                    \textcolor{Col1}{f}
                \end{math}}(4);
            \node(r)at(9.00,2){};
            \draw(r)edge[Edge]node[CliqueLabel]
                {\begin{math}
                    \textcolor{Col1}
                    {(a_1 + a_2 + a_3) a_1 a_2 a_2 a_1 a_3
                    (a_1 a_2 + a_2 a_1)}
                \end{math}}(9);
        \end{tikzpicture}\,,
    \end{equation}
    so that~\eqref{equ:example_action_S_clique} is equal to the
    polynomial
    \begin{math}
        (a_1 + a_2 + a_3) a_1 a_2 a_2 a_1 a_3 (a_1 a_2 + a_2 a_1).
    \end{math}
    \medbreak

    \item[Noncommutative polynomials and constant term product]
    Consider the unitary magma $\Dbb_0$. Let
    $A := \{a_1, a_2, \dots\}$ be an infinite alphabet of noncommutative
    letters. We define on the associative algebra $\K \langle A \rangle$
    of polynomials on $A$ the linear maps
    \begin{equation}
        \omega_\Unit, \omega_0 :
        \K \langle A \rangle \to \K \langle A \rangle,
    \end{equation}
    as follows. For any $u \in A^*$, we set $\omega_\Unit(u) := u$, and
    \begin{equation}
        \omega_0(u) :=
        \begin{cases}
            1 & \mbox{if } u = \epsilon, \\
            0 & \mbox{otherwise}.
        \end{cases}
    \end{equation}
    In other terms, $\omega_0(f)$ is the constant term, denoted by
    $f(0)$, of the polynomial $f \in \K \langle A \rangle$. Since
    $\omega_\Unit$ is the identity map on $\K \langle A \rangle$ and,
    for all $u \in A^*$,
    \begin{equation}
        (\omega_0 \circ \omega_0)(f)
        = (f(0))(0)
        = f(0)
        = \omega_0(f),
    \end{equation}
    we obtain from Theorem~\ref{thm:NC_M_algebras} that the
    operations~\eqref{equ:NC_M_algebras} endow $\K \langle A \rangle$
    with a $\NC\Dbb_0$-algebra structure. For instance, for all
    polynomials $f_1$ and $f_2$ of $\K \langle A \rangle$, we have
    \vspace{-1.75em}
    \begin{multicols}{2}
    \begin{subequations}
    \begin{equation}
        f_1 \TriangleOp{\Unit}{\Unit}{\Unit} f_2 = f_1 f_2,
    \end{equation}
    \begin{equation}
        f_1 \TriangleOp{0}{\Unit}{\Unit} f_2 = (f_1 f_2)(0)
        = f_1(0) \; f_2(0),
    \end{equation}

    \begin{equation} \label{equ:constant_term_product_example_1}
        f_1 \TriangleOp{\Unit}{0}{\Unit} f_2 = f_1(0) \; f_2,
    \end{equation}
    \begin{equation} \label{equ:constant_term_product_example_2}
        f_1 \TriangleOp{\Unit}{\Unit}{0} f_2 = f_1 \; (f_2(0)).
    \end{equation}
    \end{subequations}
    \end{multicols}
    \noindent From~\eqref{equ:constant_term_product_example_1}
    and~\eqref{equ:constant_term_product_example_2}, when
    $f_1(0) = 1 = f_2(0)$,
    \begin{equation}
        f_1
        \; \left(
        \TriangleOp{\Unit}{0}{\Unit} + \TriangleOp{\Unit}{\Unit}{0}
        \right) \;
        f_2
        =
        f_1(0) \; f_2 + f_1 \; (f_2(0))
        = f_1 + f_2.
    \end{equation}
\end{description}
\medbreak

\subsubsection{From monoids}
If $\Mca$ is a monoid, with binary associative operation $\Op$ and unit
$\Unit_\Mca$, we denote by $\K \langle \Mca^* \rangle$ the space of all
noncommutative polynomials on $\Mca$, seen as an alphabet, with
coefficients in $\K$. This space can be endowed with an
$\NC\Mca$-algebra structure as follows.
\medbreak

For any $x \in \Mca$ and any word $w \in \Mca^*$, let
\begin{equation}
    x * w := (x \Op w_1) \dots (x \Op w_{|w|}).
\end{equation}
This operation $*$ is linearly extended on the right on
$\K \langle \Mca^* \rangle$.
\medbreak

\begin{Proposition} \label{prop:NC_M_algebras_monoid_polynomials}
    Let $\Mca$ be a finite monoid. The vector space
    $\K \langle \Mca^* \rangle$ endowed with the binary linear
    operations
    \begin{equation} \label{equ:NC_M_algebras_monoid_polynomials_op}
        \TriangleOp{\Pfr_0}{\Pfr_1}{\Pfr_2} :
        \K \langle \Mca^* \rangle \otimes \K \langle \Mca^* \rangle
        \to \K \langle \Mca^* \rangle,
        \qquad \Pfr \in \Triangles_\Mca,
    \end{equation}
    defined for each $\Mca$-triangle $\Pfr$ and any
    $f_1, f_2 \in \K \langle \Mca^* \rangle$ by
    \begin{equation} \label{equ:NC_M_algebras_monoid_polynomials}
        f_1 \TriangleOp{\Pfr_0}{\Pfr_1}{\Pfr_2} f_2 :=
        \Pfr_0 * \left(\left(\Pfr_1 * f_1\right) \;
        \left(\Pfr_2 * f_2\right)\right),
    \end{equation}
    is an $\NC\Mca$-algebra.
\end{Proposition}
\medbreak

Here are few examples of the construction provided by
Proposition~\ref{prop:NC_M_algebras_monoid_polynomials}.
\begin{description}[fullwidth]
    \item[Words and double shifted concatenation]
    Consider the monoid $\N_\ell$ for an $\ell \geq 1$. By
    Proposition~\ref{prop:NC_M_algebras_monoid_polynomials}, the
    operations~\eqref{equ:NC_M_algebras_monoid_polynomials_op} endow
    $\K \left\langle \N_\ell^* \right\rangle$ with a structure of an
    $\NC\N_\ell$-algebra. For instance, in
    $\K \left\langle \N_4^* \right\rangle$, one has
    \begin{equation}
        0211 \; \TriangleOp{1}{2}{0} \; 312 = 3100023.
    \end{equation}
    \medbreak

    \item[Words and erasing concatenation]
    Consider the monoid $\Dbb_\ell$ for an $\ell \geq 0$. By
    Proposition~\ref{prop:NC_M_algebras_monoid_polynomials},
    the operations~\eqref{equ:NC_M_algebras_monoid_polynomials_op}
    endow $\K \langle \Dbb_\ell^* \rangle$ with a structure of an
    $\NC\Dbb_\ell$-algebra. For instance, for all words $u$ and $v$
    of $\Dbb_\ell^*$, we have
    \vspace{-1.75em}
    \begin{multicols}{2}
    \begin{subequations}
    \begin{equation}
        u \TriangleOp{\Unit}{\Unit}{\Unit} v = u v,
    \end{equation}
    \begin{equation}
        u \TriangleOp{\Dsf_i}{\Unit}{\Unit} v = (u v)_{\Dsf_i},
    \end{equation}

    \begin{equation}
        u \TriangleOp{0}{\Unit}{\Unit} v = 0^{|u| + |v|},
    \end{equation}
    \begin{equation}
        u \TriangleOp{\Unit}{\Dsf_i}{\Dsf_j} v =
        u_{\Dsf_i} \; v_{\Dsf_j},
    \end{equation}
    \end{subequations}
    \end{multicols}
    \noindent where, for any word $w$ of $\Dbb_\ell^*$ and any element
    $\Dsf_j$ of $\Dbb_\ell$, $j \in [\ell]$, $w_{\Dsf_j}$ is the word
    obtained by replacing each occurrence of $\Unit$ by $\Dsf_j$ and
    each occurrence of $\Dsf_i$, $i \in [\ell]$, by $0$ in~$w$.
\end{description}
\medbreak

\subsection{Koszul dual}\label{subsec:dual_NC_M}
Since by Theorem~\ref{thm:presentation_NC_M}, the operad $\NC\Mca$ is
binary and quadratic, this operad admits a Koszul dual $\NC\Mca^!$. We
end the study of $\NC\Mca$ by collecting the main properties
of~$\NC\Mca^!$.
\medbreak

\subsubsection{Presentation}
Let  $\RelationSpace_{\NC\Mca}^!$ be the subspace of
$\FreeOperad\left(\Triangles_\Mca\right)(3)$ generated by the elements
\begin{subequations}
\begin{equation} \label{equ:relation_1_NC_M_dual}
    \sum_{\substack{
        \Pfr_1, \Qfr_0 \in \Mca \\
        \Pfr_1 \Op \Qfr_0 = \delta
    }}
    \Corolla{\Triangle{\Pfr_0}{\Pfr_1}{\Pfr_2}}
    \circ_1
    \Corolla{\Triangle{\Qfr_0}{\Qfr_1}{\Qfr_2}},
    \qquad
    \Pfr_0, \Pfr_2, \Qfr_1, \Qfr_2 \in \Mca,
    \delta \in \bar{\Mca},
\end{equation}
\begin{equation} \label{equ:relation_2_NC_M_dual}
    \sum_{\substack{
        \Pfr_1, \Qfr_0 \in \Mca \\
        \Pfr_1 \Op \Qfr_0 = \Unit_\Mca
    }}
    \Corolla{\Triangle{\Pfr_0}{\Pfr_1}{\Pfr_2}}
    \circ_1
    \Corolla{\Triangle{\Qfr_0}{\Qfr_1}{\Qfr_2}}
    -
    \Corolla{\Triangle{\Pfr_0}{\Qfr_1}{\Pfr_1}}
    \circ_2
    \Corolla{\Triangle{\Qfr_0}{\Qfr_2}{\Pfr_2}},
    \qquad
    \Pfr_0, \Pfr_2, \Qfr_1, \Qfr_2 \in \Mca,
\end{equation}
\begin{equation} \label{equ:relation_3_NC_M_dual}
    \sum_{\substack{
        \Pfr_2, \Qfr_0 \in \Mca \\
        \Pfr_2 \Op \Qfr_0 = \delta
    }}
    \Corolla{\Triangle{\Pfr_0}{\Pfr_1}{\Pfr_2}}
    \circ_2
    \Corolla{\Triangle{\Qfr_0}{\Qfr_1}{\Qfr_2}},
    \qquad
    \Pfr_0, \Pfr_1, \Qfr_1, \Qfr_2 \in \Mca,
    \delta \in \bar{\Mca},
\end{equation}
\end{subequations}
where $\Pfr$ and $\Qfr$ are $\Mca$-triangles.
\medbreak

\begin{Proposition} \label{prop:presentation_dual_NC_M}
    Let $\Mca$ be a finite unitary magma. Then, the Koszul dual
    $\NC\Mca^!$ of $\NC\Mca$ admits the presentation
    $\left(\Triangles_\Mca, \RelationSpace_{\NC\Mca}^!\right)$.
\end{Proposition}
\medbreak

By Proposition~\ref{prop:presentation_dual_NC_M}, the operad
$\NC\N_2^!$ is generated by
\begin{equation}
    \Triangles_{\N_2} =
    \left\{
    \TriangleEEE{}{}{},
    \TriangleEXE{}{1}{},
    \TriangleEEX{}{}{1},
    \TriangleEXX{}{1}{1},
    \TriangleXEE{1}{}{},
    \TriangleXXE{1}{1}{},
    \TriangleXEX{1}{}{1},
    \Triangle{1}{1}{1}
    \right\},
\end{equation}
and these generators are subjected exactly to the nontrivial relations
\begin{subequations}
\begin{equation}
    \TriangleXEX{a}{}{b_3} \circ_1 \Triangle{1}{b_1}{b_2}
    +
    \Triangle{a}{1}{b_3} \circ_1 \TriangleEXX{}{b_1}{b_2}
    = 0,
    \qquad
    a, b_1, b_2, b_3 \in \N_2,
\end{equation}
\begin{equation}
    \Triangle{a}{1}{b_3} \circ_1 \Triangle{1}{b_1}{b_2}
    +
    \TriangleXEX{a}{}{b_3} \circ_1 \TriangleEXX{}{b_1}{b_2}
    =
    \TriangleXXE{a}{b_1}{} \circ_2 \TriangleEXX{}{b_2}{b_3}
    +
    \Triangle{a}{b_1}{1} \circ_2 \Triangle{1}{b_2}{b_3},
    \qquad
    a, b_1, b_2, b_3 \in \N_2,
\end{equation}
\begin{equation}
    \TriangleXXE{a}{b_1}{} \circ_2 \Triangle{1}{b_2}{b_3}
    +
    \Triangle{a}{b_1}{1} \circ_2 \TriangleEXX{}{b_2}{b_3}
    = 0,
    \qquad
    a, b_1, b_2, b_3 \in \N_2.
\end{equation}
\end{subequations}
\medbreak

\begin{Proposition} \label{prop:dimensions_relations_NC_M_dual}
    Let $\Mca$ be a finite unitary magma. Then, the dimension of the
    space $\RelationSpace_{\NC\Mca}^!$ satisfies
    \begin{equation} \label{equ:dimensions_relations_NC_M_dual}
        \dim \RelationSpace_{\NC\Mca}^! = 2m^5 - m^4,
    \end{equation}
    where $m := \# \Mca$.
\end{Proposition}
\medbreak

Observe that, by Propositions~\ref{prop:dimensions_relations_NC_M}
and~\ref{prop:dimensions_relations_NC_M_dual}, we have
\begin{equation}\begin{split}
    \dim \RelationSpace_{\NC\Mca} + \dim \RelationSpace_{\NC\Mca}^!
        & = 2m^6 - 2m^5 + m^4
          + 2m^5 - m^4 \\
        & = 2m^6 \\
        & = \dim \FreeOperad\left(\Triangles_\Mca\right)(3),
\end{split}\end{equation}
as expected by Koszul duality, where $m := \# \Mca$.
\medbreak

\subsubsection{Dimensions}

\begin{Proposition} \label{prop:Hilbert_series_NC_M_dual}
    Let $\Mca$ be a finite unitary magma. The Hilbert series
    $\HilbSeries_{\NC\Mca^!}(t)$ of $\NC\Mca^!$ satisfies
    \begin{equation} \label{equ:Hilbert_series_NC_M_dual}
        t + (m - 1)t^2
        + \left(2m^2t - 3mt + 2t -1\right)\HilbSeries_{\NC\Mca^!}(t)
        + \left(m^3 - 2m^2 + 2m - 1\right)\HilbSeries_{\NC\Mca^!}(t)^2
        = 0,
    \end{equation}
    where $m := \# \Mca$.
\end{Proposition}
\medbreak

We deduce from Proposition~\ref{prop:Hilbert_series_NC_M_dual} that the
Hilbert series of $\NC\Mca^!$ satisfies
\begin{equation} \label{equ:Hilbert_function_NC_M_dual}
    \HilbSeries_{\NC\Mca^!}(t) =
    \frac{1 - (2m^2 - 3m + 2)t -\sqrt{1 - 2(2m^3 - 2m^2 + m)t +m^2t^2}}
    {2(m^3 - 2m^2 + 2m - 1)},
\end{equation}
where $m := \# \Mca \ne 1$.
\medbreak

\begin{Proposition} \label{prop:dimensions_NC_M_dual}
    Let $\Mca$ be a finite unitary magma. For all $n \geq 2$,
    \begin{equation} \label{equ:dimensions_NC_M_dual}
        \dim \NC\Mca^!(n)
        =
        \sum_{0 \leq k \leq n - 2}
        m^{n + 1}
        (m(m - 1) + 1)^k (m (m - 1))^{n - k - 2}
        \; \Narayana(n, k).
    \end{equation}
\end{Proposition}
\medbreak

We can use Proposition~\ref{prop:dimensions_NC_M_dual} to compute the
first dimensions of $\NC\Mca^!$. For instance, depending on
$m := \# \Mca$, we have the following sequences of dimensions:
\begin{subequations}
\begin{equation}
    1, 1, 1, 1, 1, 1, 1, 1,
    \qquad m = 1,
\end{equation}
\begin{equation}
    1, 8, 80, 992, 13760, 204416, 3180800, 51176960,
    \qquad m = 2,
\end{equation}
\begin{equation}
    1, 27, 1053, 51273, 2795715, 163318599, 9994719033, 632496651597,
    \qquad m = 3,
\end{equation}
\begin{equation}
    1, 64, 6400, 799744, 111923200, 16782082048, 2636161024000,
    428208345579520,
    \qquad m = 4.
\end{equation}
\end{subequations}
The second one is Sequence~\OEIS{A234596} of~\cite{Slo}. The last two
sequences are not listed in~\cite{Slo} at this time. It is worth
observing that the dimensions of $\NC\Mca^!$ when $\# \Mca = 2$ are the
ones of the operad~$\BNC$ of bicolored noncrossing configurations (see
Chapter~\ref{chap:enveloping}).
\medbreak

\subsubsection{Basis} \label{subsubsec:basis_Cli_M_dual}
To describe a basis of $\NC\Mca^!$, let us introduce the following sort
of $\Mca$-deco\-rated cliques. A \Def{dual $\Mca$-clique} is an
$\Mca^2$-clique such that its base and its edges are labeled by pairs
$(a, a) \in \Mca^2$, and all solid diagonals are labeled by pairs
$(a, b) \in \Mca^2$ with $a \ne b$. Observe that a non-solid diagonal of
a dual $\Mca$-clique is labeled by $(\Unit_\Mca, \Unit_\Mca)$. All
definitions about $\Mca$-cliques of
Section~\ref{subsec:decorated_cliques} remain valid for dual
$\Mca$-cliques. For example,
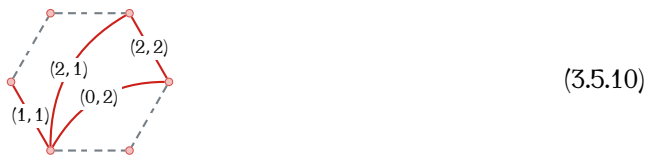
\begin{equation}
    \begin{tikzpicture}[scale=1.05,Centering]
        \node[CliquePoint](1)at(-0.50,-0.87){};
        \node[CliquePoint](2)at(-1.00,-0.00){};
        \node[CliquePoint](3)at(-0.50,0.87){};
        \node[CliquePoint](4)at(0.50,0.87){};
        \node[CliquePoint](5)at(1.00,0.00){};
        \node[CliquePoint](6)at(0.50,-0.87){};
        \draw[CliqueEdge](1)edge[]node[CliqueLabel]
            {\begin{math}(1, 1)\end{math}}(2);
        \draw[CliqueEdge](1)edge[bend left=30]node[CliqueLabel]
            {\begin{math}(0, 2)\end{math}}(5);
        \draw[CliqueEmptyEdge](1)edge[]node[CliqueLabel]{}(6);
        \draw[CliqueEdge](1)edge[bend left=30]node[CliqueLabel]
            {\begin{math}(2, 1)\end{math}}(4);
        \draw[CliqueEmptyEdge](3)edge[]node[CliqueLabel]{}(4);
        \draw[CliqueEdge](4)edge[]node[CliqueLabel]
            {\begin{math}(2, 2)\end{math}}(5);
        \draw[CliqueEmptyEdge](2)edge[]node[CliqueLabel]{}(3);
        \draw[CliqueEmptyEdge](5)edge[]node[CliqueLabel]{}(6);
    \end{tikzpicture}
\end{equation}
is a noncrossing dual $\N_3$-clique.
\medbreak

\begin{Proposition} \label{prop:elements_NC_M_dual}
    Let $\Mca$ be a finite unitary magma. The underlying graded vector
    space of $\NC\Mca^!$ is the linear span of all noncrossing dual
    $\Mca$-cliques.
\end{Proposition}
\medbreak

Proposition~\ref{prop:elements_NC_M_dual} gives a combinatorial
description of the elements of $\NC\Mca^!$. Nevertheless, we do not know
for the time being a partial composition on the linear span of these
elements providing a realization of~$\NC\Mca^!$.
\medbreak

\section{Concrete constructions}\label{sec:concrete_constructions}
The clique construction provides alternative definitions of known
operads. We explore here the cases of the operad $\NCT$ of based
noncrossing trees, the operad $\FF_4$ of formal fractions, the operad
$\BNC$ of bicolored noncrossing configurations, the operads $\MT$
and $\DMT$ of multi-tildes and double multi-tildes, and the gravity
operad~$\Grav$.
\medbreak

\subsection{Rational functions and related operads}
We use here the (noncrossing) clique construction to interpret some few
operads related to the operad $\RatFct$ of rational functions.
\medbreak

\subsubsection{Dendriform and based noncrossing tree operads}
One can build the operad $\NCT$ of based noncrossing trees~\cite{Cha07}
(see for instance Section~\ref{subsubsec:construction_NCT} of
Chapter~\ref{chap:enveloping} where an operad isomorphic to $\NCT$ is
constructed) in the following way.
\medbreak

\begin{Proposition} \label{prop:construction_NCP}
    The suboperad $\Oca_{\NCT}$ of $\Cli\Z$ generated by
    \begin{equation}
        \left\{\TriangleEXE{}{-1}{},\TriangleEEX{}{}{-1}\right\}
    \end{equation}
    is isomorphic to the operad~$\NCT$.
\end{Proposition}
\medbreak

By Theorem~\ref{thm:rat_fct_cliques}, $\Frac_\Identity$ is
an operad morphism from $\Cli\Z$ to $\RatFct$. Hence, the restriction of
$\Frac_\Identity$ on $\Oca_{\NCT}$ is also an operad morphism from
$\Oca_{\NCT}$ to $\RatFct$. Moreover, since
\vspace{-1.75em}
\begin{multicols}{2}
\begin{subequations}
\begin{equation}
    \Frac_\Identity\left(
        \TriangleEXE{}{-1}{}
    \right)
    = \frac{1}{u_1},
\end{equation}

\begin{equation}
    \Frac_\Identity\left(
        \TriangleEEX{}{}{-1}
    \right)
    = \frac{1}{u_2},
\end{equation}
\end{subequations}
\end{multicols}
\noindent
and $\RatFct^{\left\{u_1^{-1}, u_2^{-1}\right\}}$ is isomorphic to the
dendriform operad $\Dendr$~\cite{Lod01,Lod10} (see
Section~\ref{subsubsec:rational_functions_operad} of
Chapter~\ref{chap:algebra}), the map $\Frac_\Identity$ is a surjective
operad morphism from $\Oca_{\NCT}$ to~$\Dendr$.
\medbreak

\subsubsection{Operad of formal fractions} \label{subsubsec:operad_ff}
One can build the suboperad $\FF_4$ of the operad of formal fractions
$\FF$~\cite{CHN16} in the following way.
\medbreak

\begin{Proposition} \label{prop:construction_FF4}
    The suboperad $\Oca_{\FF_4}$ of $\Cli\Z$ generated by
    \begin{equation} \label{equ:generators_FF4}
        \left\{
            \Triangle{-1}{-1}{1}, \Triangle{-1}{1}{-1},
            \Triangle{-1}{1}{1}, \Triangle{1}{-1}{-1}
        \right\}
    \end{equation}
    is isomorphic to the operad~$\FF_4$.
\end{Proposition}
\medbreak

Proposition~\ref{prop:construction_FF4} shows hence that the operad
$\FF_4$ can be built through the construction $\Cli$. Observe also that,
as a consequence of Proposition~\ref{prop:construction_FF4}, all
suboperads of $\FF_4$ defined in~\cite{CHN16} that are generated by a
subset of~\eqref{equ:generators_FF4} can be constructed by the clique
construction.
\medbreak

\subsection{Operad of bicolored noncrossing configurations}
\label{subsec:operad_BNC_C_M}
One can build the operad $\BNC$ of bicolored noncrossing configurations
(see Chapter~\ref{chap:enveloping}) in the following way.
\medbreak

Consider the unitary magma
$\Mca_{\BNC} := \{\Unit, \Asf, \Bsf\}$ wherein operation $\Op$ is
defined by the Cayley table
\begin{small}
\begin{equation}
    \begin{tabular}{c||c|c|c|}
        $\Op$ & \; $\Unit$ \; & \; $\Asf$ \; & \; $\Bsf$ \;
            \\ \hline \hline
        $\Unit$ & $\Unit$ & $\Asf$ & $\Bsf$ \\ \hline
        $\Asf$ & $\Asf$ & $\Asf$ & $\Unit$ \\ \hline
        $\Bsf$ & $\Bsf$ & $\Unit$ & $\Bsf$
    \end{tabular}\, .
\end{equation}
\end{small}
In other words, $\Mca_{\BNC}$ is the unitary magma wherein $\Asf$ and
$\Bsf$ are idempotent, and $\Asf \Op \Bsf = \Unit = \Bsf \Op \Asf$.
Observe that $\Mca_{\BNC}$ is a commutative unitary magma, but, since
\begin{equation}
    (\Bsf \Op \Asf) \Op \Asf = \Unit \Op \Asf = \Asf
    \ne
    \Bsf = \Bsf \Op \Unit = \Bsf \Op (\Asf \Op \Asf),
\end{equation}
the operation $\Op$ is not associative.
\medbreak

Let $\phi_{\BNC} : \BNC \to \NC\Mca_{\BNC}$ be the linear map defined
in the following way. For any bicolored noncrossing configuration $\Cfr$,
$\phi_{\BNC}(\Cfr)$ is the noncrossing $\Mca_{\BNC}$-clique of
$\NC\Mca_{\BNC}$ obtained by replacing all blue arcs of $\Cfr$ by arcs
labeled by $\Asf$, all red diagonals of $\Cfr$ by diagonals labeled by
$\Bsf$, all uncolored edges and bases of $\Cfr$ by edges labeled by
$\Bsf$, and all uncolored diagonals of $\Cfr$ by diagonals labeled by
$\Unit$. For instance,
\begin{equation}
    \phi_{\BNC}\left(
    \begin{tikzpicture}[scale=.7,Centering]
        \node[CliquePoint](0)at(-0.49,-0.86){};
        \node[CliquePoint](1)at(-1.,-0.){};
        \node[CliquePoint](2)at(-0.5,0.87){};
        \node[CliquePoint](3)at(0.5,0.87){};
        \node[CliquePoint](4)at(1.,0.01){};
        \node[CliquePoint](5)at(0.51,-0.86){};
        \draw[CliqueEdgeGray](0)--(1);
        \draw[CliqueEdgeGray](1)--(2);
        \draw[CliqueEdgeGray](2)--(3);
        \draw[CliqueEdgeGray](3)--(4);
        \draw[CliqueEdgeGray](4)--(5);
        \draw[CliqueEdgeGray](5)--(0);
        \draw[CliqueEdgeBlue](1)--(2);
        \draw[CliqueEdgeBlue](1)--(3);
        \draw[CliqueEdgeBlue](0)--(5);
        \draw[CliqueEdgeRed](1)--(4);
    \end{tikzpicture}
    \right)
    \enspace = \enspace
    \begin{tikzpicture}[scale=.7,Centering]
        \node[CliquePoint](0)at(-0.49,-0.86){};
        \node[CliquePoint](1)at(-1.,-0.){};
        \node[CliquePoint](2)at(-0.5,0.87){};
        \node[CliquePoint](3)at(0.5,0.87){};
        \node[CliquePoint](4)at(1.,0.01){};
        \node[CliquePoint](5)at(0.51,-0.86){};
        \draw[CliqueEdge](0)edge[]node[CliqueLabel]
            {\begin{math}\Bsf\end{math}}(1);
        \draw[CliqueEdge](1)edge[]node[CliqueLabel]
            {\begin{math}\Asf\end{math}}(2);
        \draw[CliqueEdge](2)edge[]node[CliqueLabel]
            {\begin{math}\Bsf\end{math}}(3);
        \draw[CliqueEdge](3)edge[]node[CliqueLabel]
            {\begin{math}\Bsf\end{math}}(4);
        \draw[CliqueEdge](4)edge[]node[CliqueLabel]
            {\begin{math}\Bsf\end{math}}(5);
        \draw[CliqueEdge](0)edge[]node[CliqueLabel]
            {\begin{math}\Asf\end{math}}(5);
        \draw[CliqueEdge](1)edge[bend right=30]node[CliqueLabel]
            {\begin{math}\Bsf\end{math}}(4);
        \draw[CliqueEdge](1)edge[bend right=30]node[CliqueLabel]
            {\begin{math}\Asf\end{math}}(3);
    \end{tikzpicture}\,.
\end{equation}
\medbreak

\begin{Proposition} \label{prop:construction_BNC}
    The linear span of $\UnitClique$ together with all noncrossing
    $\Mca_{\BNC}$-cliques without edges nor bases labeled by $\Unit$
    forms a suboperad of $\NC\Mca_{\BNC}$ isomorphic to~$\BNC$.
    Moreover, $\phi_{\BNC}$ is an isomorphism between these two operads.
\end{Proposition}
\medbreak

Proposition~\ref{prop:construction_BNC} shows hence that the operad
$\BNC$ can be built through the noncrossing clique construction.
Moreover, observe that in Section~\ref{subsubsec:symmetries_Bubble}
of Chapter~\ref{chap:enveloping}, an automorphism of $\BNC$ called
\Def{complementary} is considered. The complementary of a bicolored
noncrossing configuration is an involution acting by modifying the
colors of some of its arcs. Under our setting, this automorphism
translates simply as the map $\Cli\theta : \Oca_\BNC \to \Oca_\BNC$
where $\Oca_\BNC$ is the operad isomorphic to $\BNC$ described in the
statement of Proposition~\ref{prop:construction_BNC} and
$\theta : \Mca_{\BNC} \to \Mca_{\BNC}$ is the unitary magma automorphism
of $\Mca_{\BNC}$ satisfying $\theta(\Unit) = \Unit$,
$\theta(\Asf) = \Bsf$, and~$\theta(\Bsf) = \Asf$.
\medbreak

Besides, it is shown in Chapter~\ref{chap:enveloping} that the set of
all bicolored noncrossing configurations of arity $2$ is a minimal
generating set of $\BNC$. Thus, by
Proposition~\ref{prop:construction_BNC}, the set
\begin{equation}
    \left\{
        \Triangle{\Asf}{\Asf}{\Asf},
        \Triangle{\Asf}{\Asf}{\Bsf},
        \Triangle{\Asf}{\Bsf}{\Asf},
        \Triangle{\Asf}{\Bsf}{\Bsf},
        \Triangle{\Bsf}{\Asf}{\Asf},
        \Triangle{\Bsf}{\Asf}{\Bsf},
        \Triangle{\Bsf}{\Bsf}{\Asf},
        \Triangle{\Bsf}{\Bsf}{\Bsf}
    \right\}
\end{equation}
is a minimal generating set of the suboperad $\Oca_\BNC$ of
$\NC\Mca_{\BNC}$ isomorphic to $\BNC$. As a consequence, all the
suboperads of $\BNC$ defined in Chapter~\ref{chap:enveloping} which
are generated by a subset of the set of the generators of $\BNC$ can
be constructed by the noncrossing clique construction. This includes,
among others, the magmatic operad, the free operad on two binary
generators, the operad of noncrossing plants $\NCP$~\cite{Cha07}, the
dipterous operad~\cite{LR03,Zin12}, and the $2$-associative
operad~\cite{LR06,Zin12}.
\medbreak

\subsection{Operads from language theory}
We provide constructions of two operads coming from formal language
theory by using the clique construction.
\medbreak

\subsubsection{Multi-tildes}
One can build the operad $\MT$ of multi-tildes~\cite{LMN13} (see also
Chapter~\ref{chap:languages}) in the following way.
\medbreak

Let $\phi_{\MT} : \MT \to \Cli\Dbb_0$ be the map linearly defined as
follows. For any multi-tilde $(n, \Sfr)$ different from
$(1, \{(1, 1)\})$, $\phi_{\MT}((n, \Sfr))$ is the $\Dbb_0$-clique of
arity $n$ defined, for any $1 \leq x < y \leq n + 1$, by
\begin{equation} \label{equ:isomorphism_MT_Cli_M}
    \phi_{\MT}((n, \Sfr))(x, y) :=
    \begin{cases}
        0 & \mbox{if } (x, y - 1) \in \Sfr, \\
        \Unit & \mbox{otherwise}.
    \end{cases}
\end{equation}
For instance,
\begin{equation}
    \phi_{\MT}((5, \{(1, 5), (2, 4), (4, 5)\}))
    =
    \begin{tikzpicture}[scale=.7,Centering]
        \node[CliquePoint](1)at(-0.50,-0.87){};
        \node[CliquePoint](2)at(-1.00,-0.00){};
        \node[CliquePoint](3)at(-0.50,0.87){};
        \node[CliquePoint](4)at(0.50,0.87){};
        \node[CliquePoint](5)at(1.00,0.00){};
        \node[CliquePoint](6)at(0.50,-0.87){};
        \draw[CliqueEdge](1)edge[]node[CliqueLabel]
            {\begin{math}0\end{math}}(6);
        \draw[CliqueEmptyEdge](1)edge[]node[CliqueLabel]{}(2);
        \draw[CliqueEmptyEdge](2)edge[]node[CliqueLabel]{}(3);
        \draw[CliqueEdge](2)edge[]node[CliqueLabel,near start]
            {\begin{math}0\end{math}}(5);
        \draw[CliqueEmptyEdge](3)edge[]node[CliqueLabel]{}(4);
        \draw[CliqueEdge](4)
            edge[bend right=30]node[CliqueLabel,near start]
            {\begin{math}0\end{math}}(6);
        \draw[CliqueEmptyEdge](4)edge[]node[CliqueLabel]{}(5);
        \draw[CliqueEmptyEdge](5)edge[]node[CliqueLabel]{}(6);
    \end{tikzpicture}\,.
\end{equation}
\medbreak

\begin{Proposition} \label{prop:construction_MT}
    The operad $\Cli\Dbb_0$ is isomorphic to the suboperad
    of $\MT$ consisting in the linear span of all multi-tildes except
    the nontrivial multi-tilde $(1, \{(1, 1)\})$ of arity $1$. Moreover,
    $\phi_{\MT}$ is an isomorphism between these two operads.
\end{Proposition}
\medbreak

\subsubsection{Double multi-tildes}
One can build the operad $\DMT$ of double multi-tildes (see
Chapter~\ref{chap:languages}) in the following way.
\medbreak

Consider the operad $\Cli\Dbb_0^2$ and let
$\phi_{\DMT} : \DMT \to \Cli\Dbb_0^2$ be the map linearly defined as
follows. The image by $\phi_{\DMT}$ of $(1, \emptyset, \emptyset)$ is
the unit of $\Cli\Dbb_0^2$ and, for any double multi-tilde
$(n, \Sfr, \Tfr)$ of arity $n \geq 2$, $\phi_{\DMT}((n, \Sfr, \Tfr))$ is
the $\Dbb_0^2$-clique of arity $n$ defined, for any
$1 \leq x < y \leq n + 1$, by
\begin{equation} \label{equ:isomorphism_DMT_Cli_M}
    \phi_{\DMT}((n, \Sfr, \Tfr))(x, y) :=
    \begin{cases}
        (0, \Unit) & \mbox{if } (x, y - 1) \in \Sfr
            \mbox{ and } (x, y - 1) \notin \Tfr, \\
        (\Unit, 0) & \mbox{if } (x, y - 1) \notin \Sfr
            \mbox{ and } (x, y - 1) \in \Tfr, \\
        (0, 0) & \mbox{if } (x, y - 1) \in \Sfr
            \mbox{ and } (x, y - 1) \in \Tfr, \\
        (\Unit, \Unit) & \mbox{otherwise}.
    \end{cases}
\end{equation}
For instance,
\begin{equation}
    \phi_{\DMT}((4, \{(2, 2), (2, 3)\}, \{(1, 3), (1, 4), (2, 3)\}))
    =
    \begin{tikzpicture}[scale=1.0,Centering]
        \node[CliquePoint](1)at(-0.59,-0.81){};
        \node[CliquePoint](2)at(-0.95,0.31){};
        \node[CliquePoint](3)at(-0.00,1.00){};
        \node[CliquePoint](4)at(0.95,0.31){};
        \node[CliquePoint](5)at(0.59,-0.81){};
        \draw[CliqueEmptyEdge](1)edge[]node[CliqueLabel]{}(2);
        \draw[CliqueEdge](1)edge[bend left=30]node[CliqueLabel]
            {\begin{math}(\Unit, 0)\end{math}}(4);
        \draw[CliqueEdge](1)edge[]node[CliqueLabel]
            {\begin{math}(\Unit, 0)\end{math}}(5);
        \draw[CliqueEdge](2)edge[]node[CliqueLabel]
            {\begin{math}(0, \Unit)\end{math}}(3);
        \draw[CliqueEdge](2)edge[]node[CliqueLabel]
            {\begin{math}(0, 0)\end{math}}(4);
        \draw[CliqueEmptyEdge](3)edge[]node[CliqueLabel]{}(4);
        \draw[CliqueEmptyEdge](4)edge[]node[CliqueLabel]{}(5);
    \end{tikzpicture}\,.
\end{equation}
\medbreak

\begin{Proposition} \label{prop:construction_DMT}
    The operad $\Cli\Dbb_0^2$ is isomorphic to the suboperad of $\DMT$
    consisting in the linear span of all double multi-tildes except the
    three nontrivial double multi-tildes of arity $1$. Moreover,
    $\phi_{\DMT}$ is an isomorphism between these two operads.
\end{Proposition}
\medbreak

\subsection{Gravity operad}
One can build the nonsymmetric version~\cite{AP15} (see
Section~\ref{subsubsec:operad_Grav} of Chapter~\ref{chap:algebra}) of
the operad $\Grav$ of gravity chord diagrams~\cite{Get94} in the
following way.
\medbreak

Let $\phi_{\Grav} : \Grav \to \Cli\Dbb_0$ be the linear map defined in
the following way. For any gravity chord diagram $\Cfr$,
$\phi_{\Grav}(\Cfr)$ is the $\Dbb_0$-clique of $\Cli\Dbb_0$ obtained by
replacing all blue arcs of $\Cfr$ by arcs labeled by $0$ and all
unlabeled arcs by arcs labeled by $\Unit$. For instance,
\begin{equation}
    \phi_{\Grav}\left(
\,.
\end{equation}
\medbreak

Let us say that an $\Mca$-clique $\Pfr$ satisfies the
\Def{gravity condition} if $\Pfr = \UnitClique$, or $\Pfr$ has only
solid edges and bases, and for all crossing diagonals $(x, y)$ and
$(x', y')$ of $\Pfr$ such that $x < x'$,
$\Pfr(x, y) \ne \Unit_\Mca \ne \Pfr(x', y')$
implies~$\Pfr(x', y) = \Unit_\Mca$.
\medbreak

\begin{Proposition} \label{prop:construction_Grav}
    The linear span of all $\Dbb_0$-cliques satisfying the gravity
    condition forms a suboperad of $\Cli\Dbb_0$ isomorphic to $\Grav$.
    Moreover, $\phi_{\Grav}$ is an isomorphism between these two
    operads.
\end{Proposition}
\medbreak

Proposition~\ref{prop:construction_Grav} shows hence that the operad
$\Grav$ can be built through the clique construction. Moreover, as
explained in~\cite{AP15}, $\Grav$ contains the nonsymmetric version of
the $\Lie$ operad, the symmetric operad describing the category of Lie
algebras. This nonsymmetric version of the Lie operad as been introduced
in~\cite{ST09}. Since $\Lie$ is contained in $\Grav$ as the subspace of
all gravity chord diagrams having the maximal number of blue diagonals
for each arity, $\Lie$ can be built through the clique construction as
the suboperad of $\Cli\Dbb_0$ containing all the $\Dbb_0$-cliques that
are images by $\phi_{\Grav}$ of such maximal gravity chord diagrams.
\medbreak

Besides, this alternative construction of $\Grav$ leads to the following
generalization for any unitary magma $\Mca$ of the gravity operad. Let
$\Grav_\Mca$ be the linear span of all $\Mca$-cliques satisfying the
gravity condition. It follows from the definition of the partial
composition of $\Cli\Mca$ that $\Grav_\Mca$ is an operad. Moreover,
observe that when $\Mca$ has nontrivial unit divisors, $\Grav_\Mca$ is
not a free operad.
\medbreak

\section*{Concluding remarks}
This chapter presents and studies the clique construction $\Cli$,
producing operads from unitary magmas. We have seen that $\Cli$ has many
both algebraic and combinatorial properties. Among its most notable
ones, $\Cli\Mca$ admits several quotients involving combinatorial
families of decorated cliques, admits a binary and quadratic suboperad
$\NC\Mca$ which is a Koszul, and contains a lot of already studied and
classic operads. Besides, in the course of this chapter, whose is
already long enough, we have put aside a bunch of questions. Let us
address these here.
\smallbreak

First, we have for the time being no formula to enumerate prime
(resp. white prime) $\Mca$-cliques
(see~\eqref{equ:prime_cliques_numbers_2} (resp.
\eqref{equ:white_prime_cliques_numbers}) for $\#\Mca = 2$). Obtaining
these forms a first combinatorial question.
\smallbreak

When $\Mca$ is a $\Z$-graded unitary magma, a link between $\Cli\Mca$
and the operad of rational functions $\RatFct$ has been
developed in Section~\ref{subsubsec:rational_functions} by means of a
morphism $\Frac_\theta$ between these two operads. We have observed that
$\Frac_\theta$ is not injective (see~\eqref{equ:frac_not_injective_1}
and~\eqref{equ:frac_not_injective_2}). A description of the kernel of
$\Frac_\theta$, even when $\Mca$ is the unitary magma $\Z$, seems not
easy to obtain. Trying to obtain this description is a second
perspective of this work.
\smallbreak

Here is a third perspective. In
Section~\ref{sec:quotients_suboperads_C_M}, we have defined and briefly
studied some quotients and suboperads of $\Cli\Mca$. In particular, we
have considered the quotient $\Deg_1\Mca$ of $\Cli\Mca$, involving
$\Mca$-cliques of degrees at most $1$. As mentioned, $\Deg_1\Dbb_0$ is
an operad defined on the linear span of involutions (except the
nontrivial involution of $\mathfrak{S}_2$). A complete study of this
operad seems worth considering, including a description of a minimal
generating set, a presentation by generators and relations, a
description of its partial composition on the $\BasisH$-basis and on the
$\BasisK$-basis, and a realization of this operad in terms of standard
Young tableaux.
\smallbreak

The last question we develop here concerns the Koszul dual $\NC\Mca^!$
of $\NC\Mca$. Section~\ref{subsec:dual_NC_M} contains results about this
operad, like a description of its presentation and a formula for its
dimensions. We have also established the fact that, as graded vector
spaces, $\NC\Mca^!$ is isomorphic to the linear span of all noncrossing
dual $\Mca$-cliques. To obtain a realization of $\NC\Mca^!$, it is now
enough to endow this last space with an adequate partial composition.
This is the last perspective we address here.
\medbreak

\part{Combinatorial Hopf bialgebras}

\chapter{Hopf bialgebra of packed square matrices} \label{chap:matrices}
The content of this chapter comes from~\cite{CGM15} and is a joint work
with Hayat Cheballah and Rémi Maurice.
\medbreak

\section*{Introduction}
The combinatorial collection of the permutations is naturally endowed
with two operations. One of them, the shifted shuffle product, takes two
permutations as input and put these together by blending their letters.
The other one, the  deconcatenation coproduct, takes one
permutation as input and disassembles it by cutting it into prefixes and
suffixes. These two operations satisfy certain compatibility relations,
resulting in that the linear span of all permutations forms a Hopf
bialgebra~\cite{MR95}, known as the Malvenuto-Reutenauer Hopf bialgebra
or $\FQSym$~\cite{DHT02} (see Section~\ref{subsubsec:FQSym} of
Chapter~\ref{chap:algebra}).
\smallbreak

This Hopf bialgebra plays a central role in algebraic combinatorics for
at least two reasons. On the one hand, $\FQSym$ contains, as Hopf
sub-bialgebras, several structures based on well-known combinatorial
objects as {\em e.g.,} standard Young tableaux~\cite{DHT02}, binary
trees~\cite{HNT05}, and integer compositions~\cite{GKLLRT94}. The
construction of these substructures revisits many algorithms coming
from computer science and combinatorics. Indeed, the insertion of a
letter into a Young tableau (following
Robinson-Schensted~\cite{Sch61,Lot02}) or in a binary search
tree~\cite{Knu98} are algorithms which prove to be as enlightening as
surprising in this algebraic context~\cite{DHT02,HNT05}. On the other
hand, the polynomial realization of $\FQSym$ allows to associate a
polynomial with any permutation~\cite{DHT02} providing a generalization
of symmetric functions, the free quasi-symmetric functions. This
generalization offers alternative ways to prove several properties of
(quasi)symmetric functions.
\smallbreak

It is thus natural to enrich this theory by proposing generalizations
of $\FQSym$. In the last years, several generalizations were proposed
and each of these depends on the way we regard permutations. By
regarding a permutation as a word and allowing repetitions of letters,
Hivert introduced in~\cite{Hiv99} (see~\cite{NT06} for a detailed
study) a Hopf bialgebra $\WQSym$ on packed words. Additionally, by
allowing some jumps for the values of the letters of permutations,
Novelli and Thibon defined in~\cite{NT07} another Hopf bialgebra
$\PQSym$ which involves parking functions. These authors also showed
in~\cite{NT10} that the $k$-colored permutations admit a Hopf bialgebra
structure $\FQSym^{(k)}$. Furthermore, by regarding a permutation
$\sigma$ as a bijection associating the singleton $\{\sigma(i)\}$ with
any singleton $\{i\}$, Aguiar and Orellana constructed~\cite{AO08} a
Hopf bialgebra structure $\UBP$ on uniform block permutations, {\em i.e.},
bijections between set partitions of $[n]$, where each part has the same
cardinality as its image. Finally, by regarding a permutation through
its permutation matrix, Duchamp, Hivert and Thibon introduced
in~\cite{DHT02} a Hopf bialgebra $\MQSym$ which involves some kind of
integer matrices.
\smallbreak

In this chapter we propose a new generalization of $\FQSym$ by regarding
permutations as permutation matrices. For this purpose, we consider the
set of $1$-packed matrices that are square matrices with entries in the
alphabet $\{\Zero, \One\}$ which have at least one $1$ by row and by
column. By equipping these matrices with a product and a coproduct, we
obtain a bigraded Hopf bialgebra, denoted by $\PM{1}$. By only
considering the grading offered by the size (resp. the number of
nonzero entries) of matrices, we obtain a simply graded Hopf bialgebra
denoted by $\PMN{1}$ (resp. $\PML{1}$). Note that since permutation
matrices form a Hopf sub-bialgebra of $\PMN{1}$ (and $\PML{1}$)
isomorphic to $\FQSym$, $\PMN{1}$ (and $\PML{1}$) provides a
generalization of $\FQSym$. Now, by allowing the entries different
from $\Zero$ of a packed matrix to belong to the alphabet $[k]$ where
$k$ is a positive integer, we obtain the notion of a $k$-packed
matrix. The definition of $\PM{1}$ (and $\PMN{1}$ and $\PML{1}$)
obviously extends to these matrices and leads to the Hopf bialgebra
$\PM{k}$ (and $\PMN{k}$ and $\PML{k}$) involving $k$-packed matrices.
Besides, since any $k$-packed matrix is also a $k + 1$-packed matrix,
$(\PM{k})_{k \geq 1}$ is an increasing infinite sequence of Hopf
bialgebras for inclusion. Let us now list some remarkable facts about
these Hopf bialgebras. First, $\FQSym^{(k)}$ embeds into $\PMN{k}$
(and $\PML{k}$), and the dual $\UBP^\star$ of $\UBP$ embeds into
$\PMN{1}$. Besides, as associative algebras, $\PML{1}^\star$ embeds into
$\MQSym$. On the other hand, by considering a bijection between the set
of the alternating sign matrices~\cite{MRR83} and particular
$1$-packed matrices, it appears that the linear span of these
$1$-packed matrices forms a Hopf sub-bialgebra of $\PM{k}$. This Hopf
bialgebra, called $\ASM$, is hence a Hopf bialgebra on alternating
sign matrices. Several statistics defined on alternating sign
matrices through the six-vertex configurations with domain wall
boundary conditions~\cite{Kup96} can be interpreted under this
algebraic point of view.
\smallbreak

Our results are presented as follows. The aims of
Section~\ref{sec:packed_matrices} are to introduce $k$-packed matrices
and the Hopf bialgebra of $k$-packed matrices.
Section~\ref{sec:properties_PM} is devoted to the study
of the algebraic properties of $\PM{k}$. In
Section~\ref{sec:links_packed_matrices},
we describe morphisms between $\PM{k}$ another Hopf bialgebras. We also
provide a general way to construct Hopf sub-bialgebras of $\PM{k}$
analogous to the construction of Hopf sub-bialgebras of $\FQSym$ by
monoid congruences~\cite{Hiv99} (see also~\cite{Gir11}). We end this
chapter by Section~\ref{sec:hopf_algebra_ASM} where we study the Hopf
sub-bialgebra $\ASM$ of $\PMN{1}$.
\medbreak

\section{Hopf algebra of packed matrices} \label{sec:packed_matrices}
We begin this section by defining $k$-packed matrices and by
enumerating them following their sizes and their number of nonzero
entries. Then we introduce the Hopf algebra $\PM{k}$ on the linear span
of the $k$-packed matrices.
\medbreak

\subsection{Packed matrices}
Let us introduce here the most important combinatorial object of this
work.
\medbreak

\subsubsection{First definitions}
Let $k \geq 1$ be an integer. We denote by $\Mca_{k, n, \ell}$ the set
of $n \times n$ matrices with exactly $\ell$ nonzero entries in the
alphabet $A_k := \{\Zero, 1, \dots, k\}$ and by $\NullRows(M)$ (resp.
$\NullColumns(M)$) the set of the indices of the zero rows (resp.
columns) of $M \in \Mca_{k, n, \ell}$. For example, by considering the
matrix
\begin{equation}
    M :=
    \Matrix{
        0 & 1 & 0 & 0 & 1 & 0 \\
        0 & 0 & 0 & 1 & 0 & 1 \\
        0 & 1 & 0 & 0 & 0 & 0 \\
        0 & 0 & 0 & 1 & 1 & 0 \\
        0 & 0 & 0 & 0 & 0 & 0 \\
        0 & 0 & 0 & 0 & 0 & 1},
\end{equation}
we have $\NullRows(M) = \{5\}$ and $\NullColumns(M) = \{1, 3\}$.
\medbreak

A \Def{$k$-packed matrix} $M$ of size $n$ is a matrix of
$\Mca_{k, n, \ell}$ in which each row and each column contains at
least one entry different from $\Zero$, that is to say if the subsets
$\NullRows(M)$ and $\NullColumns(M)$ are empty.
\medbreak

We shall denote in the sequel by $\PackedMatrices_{k, n, \ell}$ the
set of all $k$-packed matrices of size $n$ with exactly $\ell$ nonzero
entries, by $\PackedMatrices_{k, n, -}$ the set of all $k$-packed
matrices of size $n$, by $\PackedMatrices_{k, -, \ell}$ the set of all
$k$-packed matrices with exactly $\ell$ nonzero entries, and by
$\PackedMatrices_k$ the set of all $k$-packed matrices. The $k$-packed
matrix of size $0$ is denoted by $\emptyset$. For instance, the seven
$1$-packed matrices of size $2$ are
\begin{equation}
    \Matrix{1 & 0 \\ 0 & 1},
    \Matrix{0 & 1 \\ 1 & 0},
    \Matrix{1 & 1 \\ 1 & 0},
    \Matrix{1 & 1 \\ 0 & 1},
    \Matrix{1 & 0 \\ 1 & 1},
    \Matrix{0 & 1 \\ 1 & 1},
    \Matrix{1 & 1 \\ 1 & 1},
\end{equation}
and the ten $1$-packed matrices of $\PackedMatrices_{1, -, 3}$ are
\begin{equation}
    \Matrix{1 & 1 \\ 1 & 0},
    \Matrix{1 & 1 \\ 0 & 1},
    \Matrix{1 & 0 \\ 1 & 1},
    \Matrix{0 & 1 \\ 1 & 1},
    \Matrix{1 & 0 & 0 \\ 0 & 1 & 0 \\ 0 & 0 & 1},
    \Matrix{1 & 0 & 0 \\ 0 & 0 & 1 \\ 0 & 1 & 0},
    \Matrix{0 & 1 & 0 \\ 1 & 0 & 0 \\ 0 & 0 & 1},
    \Matrix{0 & 0 & 1 \\ 1 & 0 & 0 \\ 0 & 1 & 0},
    \Matrix{0 & 1 & 0 \\ 0 & 0 & 1 \\ 1 & 0 & 0},
    \Matrix{0 & 0 & 1 \\ 0 & 1 & 0 \\ 1 & 0 & 0}.
\end{equation}
\medbreak

\subsubsection{Operations and decompositions}
Let us now define some operations on packed matrices. We shall denote
by $Z_n^m$ the $n \times m$ null matrix. Given $M_1$ and $M_2$ two
$k$-packed matrices of respective sizes $n_1$ and $n_2$, set
\begin{equation}
    \textcolor{Col1}{M_1} \Over \textcolor{Col4}{M_2} :=
    \left[\begin{array}{c|c}
        \textcolor{Col1}{M_1} & Z_{n_1}^{n_2} \\ \hline
        Z_{n_2}^{n_1} & \textcolor{Col4}{M_2}
    \end{array}\right]
    \qquad \mbox{and} \qquad
    \textcolor{Col1}{M_1} \Under \textcolor{Col4}{M_2} :=
    \left[\begin{array}{c|c}
        Z_{n_1}^{n_2} &  \textcolor{Col1}{M_1} \\ \hline
        \textcolor{Col4}{M_2} & Z_{n_2}^{n_1}
    \end{array}\right].
\end{equation}
Note that these two matrices are $k$-packed matrices of size
$n_1 + n_2$. We shall respectively call $\Over$ and $\Under$ the
\Def{over} and \Def{under} operators. These two operators are
obviously associative.
\medbreak

Given a matrix $M$ whose entries are in $A_k$, the \Def{compression} of
$M$ is the matrix $\Compr(M)$ obtained by deleting in $M$ all null
rows and columns. Let $M$ be a $k$-packed matrix. The tuple
$(M_1, \dots, M_r)$ is a \Def{column decomposition} of $M$, and we
write $M = M_1 \bullet \dots \bullet M_r$, if for all $i \in [r]$ the
$\Compr(M_i)$ are square matrices (and not necessarily column
matrices) and
\begin{equation}
    M = \left[\begin{array}{c|c|c} M_1 & \dots & M_r \end{array}\right].
\end{equation}
Similarly, the tuple $(M_1, \dots, M_r)$ is a \Def{row decomposition}
of $M$, and we write $M = M_1 \DecompC \cdots \DecompC M_r$, if for all
$i \in [r]$ the $\Compr(M_i)$ are square matrices (and not necessarily
row matrices) and
\begin{equation}
    M =
    \left[\begin{array}{c}
        M_1 \\ \hline
        \dots \\ \hline
        M_r
    \end{array}\right].
\end{equation}
\medbreak

For instance, here are a $1$-packed matrix of size $5$, one of its
column decompositions and one of its row decompositions:
\begin{equation}
    \Matrix{
        \Zero & \One & \One & \Zero & \Zero \\
        \Zero & \Zero & \One & \Zero & \Zero \\
        \Zero & \Zero & \Zero & \One & \One \\
        \One & \Zero & \Zero & \Zero & \Zero \\
        \Zero & \Zero & \Zero & \One & \One}
    =
    \Matrix{
        \Zero & \One & \One \\
        \Zero & \Zero & \One \\
        \Zero & \Zero & \Zero \\
        \One & \Zero & \Zero \\
        \Zero & \Zero & \Zero}
    \bullet
    \Matrix{
        \Zero & \Zero \\
        \Zero & \Zero \\
        \One & \One \\
        \Zero & \Zero \\
        \One & \One}
    =
    \Matrix{
        \Zero & \One & \One & \Zero & \Zero \\
        \Zero & \Zero & \One & \Zero & \Zero} \\
    \DecompC
    \Matrix{
        \Zero & \Zero & \Zero & \One & \One \\
        \One & \Zero & \Zero & \Zero & \Zero \\
        \Zero & \Zero & \Zero & \One & \One}.
\end{equation}
\medbreak

These two decompositions have the following property.
\medbreak

\begin{Lemma} \label{lem:packed_matrices_decomposition}
    Let $M$ be a packed matrix and $(M_1,M_2)$ be a column (resp. row)
    decomposition of $M$. Then, there is no integer $i$ such that the
    $i$th rows (resp. columns) of $M_1$ and $M_2$ contain both a
    nonzero entry.
\end{Lemma}
\medbreak

Lemma \ref{lem:packed_matrices_decomposition} provides a sufficient
condition to ensure that a given pair $(M_1, M_2)$ of matrices cannot be
a column (resp. row) decomposition of a matrix $M$. Nevertheless, it
is not a necessary condition. Indeed, let
\begin{equation}
    M :=
    \Matrix{
        \One & \One & \Zero \\
        \Zero & \Zero & \One \\
        \Zero & \Zero & \One}
    \qquad \mbox{and} \qquad
    (M_1, M_2) :=
    \left(
    \Matrix{
        \One & \One \\
        \Zero & \Zero \\
        \Zero & \Zero
    },
    \Matrix{
        \Zero \\
        \One \\
        \One
    }\right).
\end{equation}
Then, even if there is no nonzero entry on the same row in $M_1$ and
$M_2$, $(M_1, M_2)$ is not a column decomposition of~$M$.
\medbreak

\subsection{Enumeration} \label{subsec:enumeration_packed_matrices}
We enumerate here $k$-packed matrices by both their size and their
number  of nonzero entries. We then specialize our enumeration to obtain
formulas enumerating these objects by their size and, separately, by
their number of nonzero entries.
\medbreak

\subsubsection{General enumeration}
Using the sieve principle, we obtain the following enumerative result.
\begin{Proposition} \label{prop:enumeration_packed_matrices}
    For any $k \geq 1$, $n \geq 0$, and $\ell \geq 0$, the number
    $\# \PackedMatrices_{k, n, \ell}$ of $k$-packed matrices of size $n$
    with exactly $\ell$ nonzero entries is
    \begin{equation} \label{eq::Enum_MT}
        \# \PackedMatrices_{k, n, \ell} =
        \sum_{0 \leq i, j \leq n} (-1)^{i + j}
        \binom{n}{i}\binom{n}{j}
        \binom{ij}{\ell}k^\ell.
    \end{equation}
\end{Proposition}
\medbreak

Table~\ref{tab:numbers_packed_matrices} shows the first few values of
$\# \PackedMatrices_{k, n, \ell}$. The enumeration in the case $k = 1$
is Sequence~\OEIS{A055599} of~\cite{Slo}.
\begin{table}[ht]
    \centering
    \subfloat[][Number of $1$-packed matrices.]{
    \begin{tabular}{l|llllllllll}
        & 0 & 1 & 2 & 3 & 4 & 5 & 6 & 7 & 8 & 9 \\ \hline
        0 & 1 \\
        1 &   & 1 \\
        2 &   &  & 2 & 4 & 1 \\
        3 &   &  &   & 6 & 45 & 90 & 78 & 36 & 9 & 1
    \end{tabular}
    \label{subtab:number_1_PM}}
    \bigbreak

    \subfloat[][Number of $2$-packed matrices.]{
    \begin{tabular}{l|llllllllll}
        & 0 & 1 & 2 & 3 & 4 & 5 & 6 & 7 & 8 & 9 \\ \hline
        0 & 1 \\
        1 &   & 2  \\
        2 &   &   & 8 & 32 & 16 \\
        3 &   &   &   & 48 & 720 & 2880 & 4992 & 4608 & 2304 & 512 \\
    \end{tabular}
    \label{subtab:number_2_PM}}
    \bigbreak

    \caption[First cardinalities of sets of packed matrices.]
    {The number of $k$-packed matrices of size $n$ (vertical
    values) with exactly $\ell$ nonzero entries (horizontal values).}
    \label{tab:numbers_packed_matrices}
\end{table}
\medbreak

\subsubsection{Enumeration by size}
Notice that for any $n \geq 0$, since
\begin{equation}
    \PackedMatrices_{k, n, -} =
    \bigsqcup_{n \leq \ell \leq n^2} \PackedMatrices_{k, n, \ell},
\end{equation}
the set $\PackedMatrices_{k, n, -}$ is finite. Hence, by using
Proposition \ref{prop:enumeration_packed_matrices}, we obtain
\begin{equation}
    \# \PackedMatrices_{k, n, -} = \sum_{0 \leq i, j \leq n} (-1)^{i + j}
    \binom{n}{i}\binom{n}{j} (k + 1)^{ij}.
\end{equation}
Sequences $\left(\# \PackedMatrices_{1, n, -}\right)_{n \geq 0}$
and $\left(\# \PackedMatrices_{2, n, -}\right)_{n \geq 0}$
respectively start with
\begin{equation} \label{equ:dimensions_MTN1}
    1, 1, 7, 265, 41503, 24997921, 57366997447,
\end{equation}
and
\begin{equation} \label{equ:dimensions_MTN2}
    1, 2, 56, 16064, 39156608, 813732073472, 147662286695991296.
\end{equation}
These are respectively Sequences~\OEIS{A048291} and~\OEIS{A230879}
of~\cite{Slo}.
\medbreak

\subsubsection{Enumeration by number of nonzero entries}
Similarly, since for any $\ell \geq 0$,
\begin{equation}
    \PackedMatrices_{k, -, \ell} =
    \bigsqcup_{\left\lceil\sqrt{\ell}\right\rceil \leq n \leq \ell}
    \PackedMatrices_{k, n, \ell},
\end{equation}
the set $\PackedMatrices_{k, -, \ell}$ is finite. Hence, by using
Proposition \ref{prop:enumeration_packed_matrices}, we obtain
\begin{equation}
    \# \PackedMatrices_{k, -, \ell} = \sum_{0 \leq i, j \leq n \leq \ell}
    (-1)^{i + j}
    \binom{n}{i}\binom{n}{j}
    \binom{ij}{\ell}k^\ell.
\end{equation}
Sequences $\left(\# \PackedMatrices_{1, -, \ell}\right)_{\ell \geq 0}$
and $\left(\# \PackedMatrices_{2, -, \ell}\right)_{\ell \geq 0}$
respectively start with
\begin{equation} \label{equ:dimensions_MTL1}
    1, 1, 2, 10, 70, 642, 7246, 97052, 1503700,
\end{equation}
and
\begin{equation} \label{equ:dimensions_MTL2}
    1, 2, 8, 80, 1120, 20544, 463744, 12422656, 384947200.
\end{equation}
These are respectively Sequences~\OEIS{A104602} and~\OEIS{A230880}
of~\cite{Slo}.
\medbreak

\subsection{Hopf bialgebra structure}
We are now in position to define a Hopf bialgebra structure on the
linear span of all $k$-packed matrices.
\medbreak

\subsubsection{Bigraded space}
Let, for any $k \geq 1$,
\begin{equation}
    \PM{k} := \bigoplus_{n \geq 0} \; \bigoplus_{\ell \geq 0} \;
    \K \Angle{\PackedMatrices_{k, n, \ell}}
\end{equation}
be the bigraded vector space spanned by the set of all $k$-packed
matrices. The elements $\BasisF_M$, where the $M$ are $k$-packed
matrices, form a basis of $\PM{k}$. We shall call this basis the
\Def{fundamental basis} of $\PM{k}$.
\medbreak

\subsubsection{Product and coproduct}
Given $M_1$ and $M_2$ two $k$-packed matrices of respective sizes $n_1$
and $n_2$, set
\begin{equation}
    M_1 \circ n_2 :=
    \left[\begin{array}{c}
        \textcolor{Col1}{M_1} \\ \hline
        Z_{n_2}^{n_1}
    \end{array}\right]
    \qquad \mbox{and} \qquad
    n_1 \circ M_2 :=
    \left[\begin{array}{c}
        Z_{n_1}^{n_2} \\ \hline
        \textcolor{Col4}{M_2}
    \end{array}\right].
\end{equation}
The \Def{column shifted shuffle} $M_1 \cshuffle M_2$ of $M_1$ and $M_2$
is the set of all matrices obtained by shuffling the columns of
$M_1 \circ n_2$ with the columns of $n_1 \circ M_2$.
\medbreak

Let us endow $\PM{k}$ with a product $\cdot$ linearly defined, for any
$k$-packed matrices $M_1$ and $M_2$, by
\begin{equation} \label{equ:Def_Produit}
    \BasisF_{M_1} \cdot \BasisF_{M_2} :=
    \sum_{M \; \in \; M_1 \cshuffle M_2} \BasisF_M.
\end{equation}
For instance, in $\PM{1}$ one has
\begin{equation}
    \BasisF_{\Matrix{
        \textcolor{Col1}{\Zero} & \textcolor{Col1}{\One} \\
        \textcolor{Col1}{\One} & \textcolor{Col1}{\One}
    }}
    \cdot
    \BasisF_{\Matrix{
        \textcolor{Col4}{\One} & \textcolor{Col4}{\Zero} \\
        \textcolor{Col4}{\Zero} & \textcolor{Col4}{\One}
    }}
    =
    \BasisF_{\Matrix{
        \textcolor{Col1}{\Zero} & \textcolor{Col1}{\One} & \Zero
            & \Zero \\
        \textcolor{Col1}{\One} & \textcolor{Col1}{\One} & \Zero
            & \Zero \\
        \Zero & \Zero & \textcolor{Col4}{\One}
            & \textcolor{Col4}{\Zero} \\
        \Zero & \Zero & \textcolor{Col4}{\Zero}
            & \textcolor{Col4}{\One}
    }}
    +
    \BasisF_{\Matrix{
        \textcolor{Col1}{\Zero} & \Zero & \textcolor{Col1}{\One}
            & \Zero \\
        \textcolor{Col1}{\One} & \Zero & \textcolor{Col1}{\One}
            & \Zero \\
        \Zero & \textcolor{Col4}{\One} & \Zero
            & \textcolor{Col4}{\Zero} \\
        \Zero & \textcolor{Col4}{\Zero} & \Zero
            & \textcolor{Col4}{\One}
    }}
    +
    \BasisF_{\Matrix{
        \textcolor{Col1}{\Zero} & \Zero & \Zero
            & \textcolor{Col1}{\One} \\
        \textcolor{Col1}{\One} & \Zero & \Zero
            & \textcolor{Col1}{\One} \\
        \Zero & \textcolor{Col4}{\One}
            & \textcolor{Col4}{\Zero} & \Zero \\
        \Zero & \textcolor{Col4}{\Zero}
            & \textcolor{Col4}{\One} & \Zero
    }}
    +
    \BasisF_{\Matrix{
        \Zero & \textcolor{Col1}{\Zero} & \textcolor{Col1}{\One}
            & \Zero \\
        \Zero & \textcolor{Col1}{\One} & \textcolor{Col1}{\One}
            & \Zero \\
        \textcolor{Col4}{\One} & \Zero & \Zero
            & \textcolor{Col4}{\Zero} \\
        \textcolor{Col4}{\Zero} & \Zero & \Zero
            & \textcolor{Col4}{\One}
    }}
    +
    \BasisF_{\Matrix{
        \Zero & \textcolor{Col1}{\Zero} & \Zero
            & \textcolor{Col1}{\One} \\
        \Zero & \textcolor{Col1}{\One} & \Zero
            & \textcolor{Col1}{\One} \\
        \textcolor{Col4}{\One} & \Zero
            & \textcolor{Col4}{\Zero} & \Zero \\
        \textcolor{Col4}{\Zero} & \Zero
            & \textcolor{Col4}{\One} & \Zero
    }}
    +
    \BasisF_{\Matrix{
        \Zero & \Zero & \textcolor{Col1}{\Zero}
            & \textcolor{Col1}{\One} \\
        \Zero & \Zero & \textcolor{Col1}{\One}
            & \textcolor{Col1}{\One} \\
        \textcolor{Col4}{\One} & \textcolor{Col4}{\Zero}
            & \Zero & \Zero \\
        \textcolor{Col4}{\Zero} & \textcolor{Col4}{\One}
            & \Zero & \Zero
    }}.
\end{equation}
\medbreak

Moreover, we endow $\PM{k}$ with a coproduct $\Coproduct$ linearly
defined, for any \mbox{$k$-packed} matrix $M$, by
\begin{equation} \label{equ:Def_Coproduit}
    \Coproduct\left(\BasisF_M\right) :=
    \sum_{M = M_1 \bullet M_2}
    \BasisF_{\Compr(M_1)} \otimes \BasisF_{\Compr(M_2)}.
\end{equation}
For instance, in $\PM{1}$ one has
\begin{equation}
    \Coproduct
    \BasisF_{\Matrix{
        \One & \One & \Zero & \Zero \\
        \Zero & \Zero & \Zero & \One \\
        \One & \Zero & \One & \Zero \\
        \Zero & \One & \Zero & \Zero
    }}
    =
    \BasisF_{\Matrix{
        \One & \One & \Zero & \Zero \\
        \Zero & \Zero & \Zero & \One \\
        \One & \Zero & \One & \Zero \\
        \Zero & \One & \Zero & \Zero
    }}
    \otimes
    \BasisF_{\emptyset}
    +
    \BasisF_{\Matrix{
        \One & \One & \Zero \\
        \One & \Zero & \One \\
        \Zero & \One & \Zero
    }}
    \otimes
    \BasisF_{\Matrix{
        \One
    }}
    +
    \BasisF_{\emptyset}
    \otimes
    \BasisF_{\Matrix{
        \One & \One & \Zero & \Zero \\
        \Zero & \Zero & \Zero & \One \\
        \One & \Zero & \One & \Zero \\
        \Zero & \One & \Zero & \Zero
    }}.
\end{equation}
\medbreak

Note that by definition, the product and the coproduct of $\PM{k}$ are
multiplicity free.
\medbreak

\begin{Theorem} \label{thm:PM_Hopf}
    The vector space $\PM{k}$ endowed with the product $\cdot$ and the
    coproduct $\Coproduct$ is a bigraded and connected bialgebra where
    homogeneous components are finite-dimensional.
\end{Theorem}
\medbreak

Notice that since any $k$-packed matrix is also a $k + 1$-packed matrix,
the vector space $\PM{k}$ is included in $\PM{k + 1}$. Hence, and by
Theorem~\ref{thm:PM_Hopf},
\begin{equation}
    \PM{1} \; \hookrightarrow \; \PM{2} \; \hookrightarrow \; \cdots
\end{equation}
is an increasing infinite sequence of bigraded bialgebras for inclusion.
The first few dimensions of $\PM{1}$ and $\PM{2}$ are given by
Table~\ref{tab:numbers_packed_matrices}.
\medbreak

\subsubsection{Antipode}
Since $\PM{k}$ is, by Theorem~\ref{thm:PM_Hopf}, a bigraded and
connected bialgebra, it admits an antipode and hence, is a Hopf
bialgebra. The antipode $\nu$ of $\PM{k}$ satisfies, for any $k$-packed
matrix $M$,
\begin{equation}
    \nu\left(\BasisF_M\right) =
    \sum_{\substack{
        \ell \geq 1 \\ M = M_1 \bullet \dots \bullet M_\ell \\
        M_i \ne \emptyset, \; i \in [\ell]
    }}
    (-1)^\ell \; \BasisF_{\Compr(M_1)}
    \cdot \ldots \cdot
    \BasisF_{\Compr(M_\ell)}.
\end{equation}
For instance, in $\PM{1}$ one has
\begin{equation}
    \begin{split}
    \nu\left(\BasisF_{
    \Matrix{
        \Zero & \One & \One \\
        \One & \Zero & \Zero \\
        \Zero & \One & \Zero
    }}\right) & =
    - \BasisF_{
    \Matrix{
        \Zero & \One & \One \\
        \One & \Zero & \Zero \\
        \Zero & \One & \Zero
    }}
    +
    \BasisF_{
    \Matrix{
        \One
    }} \cdot
    \BasisF_{
    \Matrix{
        \One & \One \\
        \One & \Zero
    }} \\[.5em]
    & =
    \BasisF_{
    \Matrix{
        \One & \Zero & \Zero \\
        \Zero & \One & \One \\
        \Zero & \One & \Zero
    }}
    +
    \BasisF_{
    \Matrix{
        \Zero & \One & \Zero \\
        \One & \Zero & \One \\
        \One & \Zero & \Zero
    }}
    +
    \BasisF_{
    \Matrix{
        \Zero & \Zero & \One \\
        \One & \One & \Zero \\
        \One & \Zero & \Zero
    }}
    -
    \BasisF_{
    \Matrix{
        \Zero & \One & \One \\
        \One & \Zero & \Zero \\
        \Zero & \One & \Zero
    }}.
    \end{split}
\end{equation}
Note besides that $\nu$ is not an involution. Indeed,
\begin{equation}
    \nu^2\left(\BasisF_{
    \Matrix{
        \Zero & \One & \One \\
        \One & \Zero & \Zero \\
        \Zero & \One & \Zero
    }}\right)
    =
    \BasisF_{
    \Matrix{
        \One & \One & \Zero \\
        \One & \Zero & \Zero \\
        \Zero & \Zero & \One
    }}
    +
    \BasisF_{
    \Matrix{
        \One & \Zero & \One \\
        \One & \Zero & \Zero \\
        \Zero & \One & \Zero
    }}
    +
    \BasisF_{
    \Matrix{
        \Zero & \One & \One \\
        \Zero & \One & \Zero \\
        \One & \Zero & \Zero
    }}
    +
    \BasisF_{
    \Matrix{
        \Zero & \One & \One \\
        \One & \Zero & \Zero \\
        \Zero & \One & \Zero
    }}
    -
    \BasisF_{
    \Matrix{
        \One & \Zero & \Zero \\
        \Zero & \One & \One \\
        \Zero & \One & \Zero
    }}
    -
    \BasisF_{
    \Matrix{
        \Zero & \One & \Zero \\
        \One & \Zero & \One \\
        \One & \Zero & \Zero
    }}
    -
    \BasisF_{
    \Matrix{
        \Zero & \Zero & \One \\
        \One & \One & \Zero \\
        \One & \Zero & \Zero
    }}.
\end{equation}
\medbreak

\subsubsection{Two alternative gradings}
Let us now set
\begin{equation}
    \PMN{k} := \bigoplus_{n \geq 0}
    \K \Angle{\PackedMatrices_{k, n, -}}
    \qquad \mbox{and} \qquad
    \PML{k} := \bigoplus_{\ell \geq 0}
    \K \Angle{\PackedMatrices_{k, -, \ell}}
\end{equation}
the vector spaces of $k$-packed matrices respectively graded by the size
and by the number of nonzero entries of matrices. By
Theorem~\ref{thm:PM_Hopf}, and since each homogeneous component of
these vector spaces is finite-dimensional (see
Section~\ref{subsec:enumeration_packed_matrices}), $\PMN{k}$ and
$\PML{k}$ are Hopf bialgebras. Besides,
\begin{equation}
    \PMN{1} \; \hookrightarrow \; \PMN{2} \; \hookrightarrow \; \cdots
    \qquad \mbox{and} \qquad
    \PML{1} \; \hookrightarrow \; \PML{2} \; \hookrightarrow \; \cdots
\end{equation}
are increasing infinite sequences of Hopf bialgebras for inclusion. The
first few dimensions of $\PMN{1}$ and $\PMN{2}$ are given
by~\eqref{equ:dimensions_MTN1} and~\eqref{equ:dimensions_MTN2}, and
the first few dimensions of $\PML{1}$ and $\PML{2}$ are given
by~\eqref{equ:dimensions_MTL1} and~\eqref{equ:dimensions_MTL2}. In the
sequel, we shall denote by $\HilbSeries_{k,n}(t)$ (resp.
$\HilbSeries_{k,\ell}(t)$) the Hilbert series of $\PMN{k}$ (resp.
$\PML{k}$).
\medbreak

\section{Algebraic properties} \label{sec:properties_PM}
A complete study of the algebraic properties of $\PM{k}$ is performed
here. We show that $\PM{k}$ is free as an associative algebra, self-dual,
and admit a bidendriform bialgebra structure.
\medbreak

\subsection{Multiplicative bases and freeness}
\label{subsec:freeness_PM}
To show that $\PM{k}$ is free as an associative algebra, we define two
multiplicative bases of $\PM{k}$. The definitions of these bases rely on
a poset structure on the set of all $k$-packed matrices.
\medbreak

\subsubsection{Poset structure}
\begin{figure}[ht]
    \begin{tikzpicture}[xscale=.5,yscale=.55,font=\footnotesize]
        \filldraw[draw=ColBlack,fill=Col4!40](0,0)rectangle (1,-3);
        \draw[draw=ColBlack,fill=ColWhite](0,-3)rectangle node
            {\begin{math}0\end{math}}(1,-7);
        \draw[draw=ColBlack,fill=ColWhite](1,0)rectangle node
            {\begin{math}0\end{math}}(2,-5);
        \filldraw[draw=ColBlack,fill=Col1!40](1,-5)rectangle(2,-7);
        \node at(.5,.4){\begin{math}i\end{math}};
        \node at(1.5,.4){\begin{math}i\!+\!1\end{math}};
        \draw[draw=ColBlack](-.25,-3)edge[<->] node[anchor=center,left]
            {\begin{math}z_i\end{math}}(-.25,-7);
        \draw[draw=ColBlack](2.25,0)edge[<->]node[anchor=center,right]
            {\begin{math}z_{i + 1}\end{math}}(2.25,-5);
        \draw[draw=ColBlack](-1.25,0)edge[<->]node[anchor=center,left]
            {\begin{math}n\end{math}}(-1.25,-7);
        \node at(4,-3.5){\begin{math}\CoveringRelPM\end{math}};
        \filldraw[draw=ColBlack,fill=Col4!40](7,0)rectangle(8,-3);
        \draw[draw=ColBlack,fill=ColWhite](7,-3)rectangle node
            {\begin{math}0\end{math}}(8,-7);
        \draw[draw=ColBlack,fill=ColWhite](6,0)rectangle node
            {\begin{math}0\end{math}}(7,-5);
        \filldraw[draw=ColBlack,fill=Col1!40](6,-5)rectangle(7,-7);
    \end{tikzpicture}
    \caption[The covering relation of the poset of packed matrices.]
    {The condition for swapping the $i$th and $(i + 1)$st
    columns of a packed matrix according to the relation
    $\CoveringRelPM$. The darker regions contain any entries and the
    white ones, only zeros.}
    \label{fig:covering_relation_packed_matrices}
\end{figure}
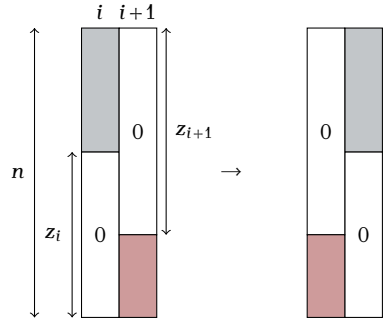

Let $\CoveringRelPM$ be the binary relation on $\PackedMatrices_k$
defined in the following way. If $M_1$ and $M_2$ are two $k$-packed
matrices of size $n$, we have $M_1 \CoveringRelPM M_2$ if there is an
index $i \in [n - 1]$ such that, denoting by $z_i$ the number of
$\Zero$ ending the $i$th column of $M_1$, and by $z_{i + 1}$ the
number of $\Zero$ starting the $(i + 1)$st column of $M_1$, one has
$z_i + z_{i + 1} \geq n$, and $M_2$ is obtained from $M_1$ by
exchanging its $i$th and $(i + 1)$st columns (see
Figure~\ref{fig:covering_relation_packed_matrices}).
\medbreak

\begin{figure}[ht]
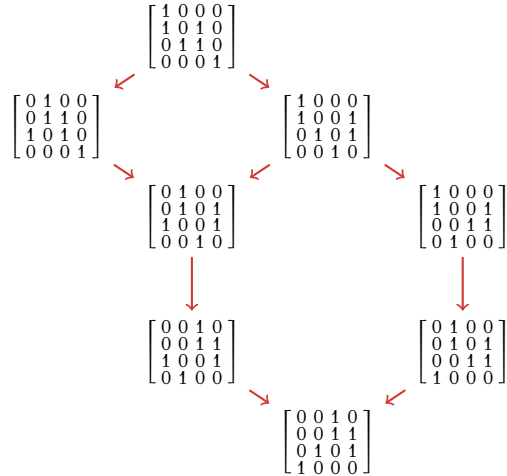


    \caption[An interval of the poset of packed matrices.]
    {The Hasse diagram of an interval for the order $\OrdPM$
    of packed matrices.}
    \label{fig:example_interval_packed_matrices}
\end{figure}

We now endow $\PackedMatrices_k$ with the partial order relation
$\OrdPM$ defined as the reflexive and transitive closure of
$\CoveringRelPM$. Figure~\ref{fig:example_interval_packed_matrices}
shows an interval of this partial order.
\medbreak

Notice that by regarding a permutation $\sigma$ of $\SymmetricGroup_n$
as its \Def{permutation matrix} ({\em i.e.}, the $1$-packed matrix $M$
of size $n$ satisfying $M_{ij} = 1$ if and only if $\sigma_j = i$),
the poset $(\PackedMatrices_{k, n, -}, \OrdPM)$ restricted to
permutation matrices is the right weak order on
permutations~\cite{GR63}.
\medbreak

\subsubsection{Multiplicative bases}
By mimicking definitions of the bases of symmetric functions, for any
$k$-packed matrix $M$, the \Def{elementary elements} $\BasisE_M$ and
the \Def{homogeneous elements} $\BasisH_M$ are respectively defined by
\begin{equation}
    \BasisE_M := \sum_{M \OrdPM M'} \BasisF_{M'}
    \qquad \mbox{and} \qquad
    \BasisH_M := \sum_{M' \OrdPM M} \BasisF_{M'}.
\end{equation}
By triangularity, these two families are bases of $\PM{k}$. For
instance, in $\PM{1}$ one has
\begin{equation}
    \BasisE_{
    \Matrix{
        \One & \Zero & \Zero & \Zero \\
        \One & \Zero & \Zero & \One \\
        \Zero & \Zero & \One & \One \\
        \Zero & \One & \Zero & \Zero
    }}
    =
    \BasisF_{
    \Matrix{
        \One & \Zero & \Zero & \Zero \\
        \One & \Zero & \Zero & \One \\
        \Zero & \Zero & \One & \One \\
        \Zero & \One & \Zero & \Zero
    }}
    +
    \BasisF_{
    \Matrix{
        \Zero & \One & \Zero & \Zero \\
        \Zero & \One & \Zero & \One \\
        \Zero & \Zero & \One & \One \\
        \One & \Zero & \Zero & \Zero
    }}
    +
    \BasisF_{
    \Matrix{
        \Zero & \Zero & \One & \Zero \\
        \Zero & \Zero & \One & \One \\
        \Zero & \One & \Zero & \One \\
        \One & \Zero & \Zero & \Zero
    }},
\end{equation}
and
\begin{equation}
    \BasisH_{
    \Matrix{
        \Zero & \One & \Zero & \Zero \\
        \Zero & \One & \Zero & \One \\
        \One & \Zero & \Zero & \One \\
        \Zero & \Zero & \One & \Zero
    }}
    =
    \BasisF_{
    \Matrix{
        \Zero & \One & \Zero & \Zero \\
        \Zero & \One & \Zero & \One \\
        \One & \Zero & \Zero & \One \\
        \Zero & \Zero & \One & \Zero
    }
    }
    +
    \BasisF_{
    \Matrix{
        \Zero & \One & \Zero & \Zero \\
        \Zero & \One & \One & \Zero \\
        \One & \Zero & \One & \Zero \\
        \Zero & \Zero & \Zero & \One
    }}
    +
    \BasisF_{
    \Matrix{
        \One & \Zero & \Zero & \Zero \\
        \One & \Zero & \Zero & \One \\
        \Zero & \One & \Zero & \One \\
        \Zero & \Zero & \One & \Zero
    }}
    +
    \BasisF_{
    \Matrix{
        \One & \Zero & \Zero & \Zero \\
        \One & \Zero & \One & \Zero \\
        \Zero & \One & \One & \Zero \\
        \Zero & \Zero & \Zero & \One
    }}.
\end{equation}
\medbreak

\begin{Proposition} \label{prop:intervall_product_PM}
    The elements appearing in a product of $\PM{k}$ expressed in the
    fundamental basis form an interval for the $\OrdPM$-partial order.
    More precisely, for any $k$-packed matrices $M_1$ and $M_2$,
    \begin{equation} \label{eq:Produit_intervalle}
        \BasisF_{M_1} \cdot \BasisF_{M_2} =
        \sum_{M_1 \Over M_2 \OrdPM M \OrdPM M_1 \Under M_2} \BasisF_M.
    \end{equation}
\end{Proposition}
\medbreak

\begin{Proposition} \label{prop:multiplicative_bases_PM}
    The product of $\PM{k}$ satisfies, for any $k$-packed matrices $M_1$
    and $M_2$,
    \begin{equation}
        \BasisE_{M_1} \cdot \BasisE_{M_2} = \BasisE_{M_1 \Over M_2}
        \qquad \mbox{and} \qquad
        \BasisH_{M_1} \cdot \BasisH_{M_2} = \BasisH_{M_1 \Under M_2}.
    \end{equation}
\end{Proposition}
\medbreak

\subsubsection{Freeness} \label{subsubsec:freeness_PM}
Given a $k$-packed matrix $M \ne \emptyset$, we say that $M$ is
\Def{connected} (resp. \Def{anti-connected}) if, for all $k$-packed
matrices $M_1$ and $M_2$, $M = M_1 \Over M_2$ (resp.
$M = M_1 \Under M_2$) implies $M_1 = M$ or $M_2 = M$.
\medbreak

\begin{Theorem} \label{thm:freeness_PM}
    The Hopf bialgebra $\PM{k}$ is freely generated as an associative
    algebra by the elements $\BasisE_M$ (resp. $\BasisH_M$) where the
    $M$ are connected (resp. anti-connected) $k$-packed matrices.
\end{Theorem}
\medbreak

Theorem~\ref{thm:freeness_PM} also implies that $\PMN{k}$ and $\PML{k}$
are freely generated by the $\BasisE_M$ (resp. $\BasisH_M$) where the
$M$ are connected (resp. anti-connected) $k$-packed matrices. Hence, the
generating series $\GenSeries_{k,n}(t)$ and $\GenSeries_{k,\ell}(t)$ of
algebraic generators of $\PMN{k}$ and $\PML{k}$ satisfy respectively
\begin{equation}
    \GenSeries_{k,n}(t) = 1 - \frac{1}{\HilbSeries_{k,n}(t)}
    \qquad \mbox{and} \qquad
    \GenSeries_{k,\ell}(t) = 1 - \frac{1}{\HilbSeries_{k,\ell}(t)}.
\end{equation}
The first few numbers of algebraic generators of $\PMN{1}$ and $\PMN{2}$
are respectively
\begin{equation} \label{equ:dimensions_generators_PMN1}
    0, 1, 6, 252, 40944, 24912120, 57316485000
\end{equation}
and
\begin{equation} \label{equ:dimensions_generators_PMN2}
    0, 2, 52, 15848, 39089872, 813573857696, 147659027604370240.
\end{equation}
These are respectively Sequences~\OEIS{A230881} and~\OEIS{A230882}
of~\cite{Slo}. The first few numbers of algebraic generators of
$\PML{1}$ and $\PML{2}$ are respectively
\begin{equation} \label{equ:dimensions_generators_PML1}
   0, 1, 1, 7, 51, 497, 5865, 81305, 1293333
\end{equation}
and
\begin{equation} \label{equ:dimensions_generators_PML2}
    0, 2, 4, 56, 816, 15904, 375360, 10407040, 331093248.
\end{equation}
These are respectively Sequences~\OEIS{A230883} and~\OEIS{A230884}
of~\cite{Slo}.
\medbreak

\subsection{Self-duality} \label{subsec:self_duality_PM}
The product and the coproduct of the dual of $\PM{k}$ are described
here. Moreover, the fact that $\PM{k}$ is a self-dual Hopf bialgebra
is shown.
\medbreak

\subsubsection{Dual Hopf bialgebra}
Let us denote by $\PM{k}^\star$ the bigraded dual vector space of
$\PM{k}$, by $\BasisF^\star_M$, where the $M$ are $k$-packed matrices,
the adjoint basis of the fundamental basis of $\PM{k}$, and by
$\langle - , - \rangle$ the associated duality bracket
(see~\eqref{equ:duality_bracket} of Chapter~\ref{chap:algebra}).
\medbreak

Let $M_1$ and $M_2$ be two $k$-packed matrices of respective sizes $n_1$
and $n_2$. By duality, the product in $\PM{k}^\star$ satisfies
\begin{equation}
    \BasisF^\star_{M_1} \cdot \BasisF^\star_{M_2} =
    \sum_{M \in \PackedMatrices_k}
    \left\langle \Coproduct\left(\BasisF_M\right),
    \BasisF^\star_{M_1} \otimes \BasisF^\star_{M_2} \right\rangle
    \; \BasisF^\star_M.
\end{equation}
Let us set
\begin{equation}
    M_1 \bullet n_2 :=
    \left[\begin{array}{c|c}
        \textcolor{Col1}{M_1} & Z_{n_1}^{n_2}
    \end{array}\right]
    \qquad \mbox{and} \qquad
    n_1 \bullet M_2 :=
    \left[\begin{array}{c|c}
        Z_{n_2}^{n_1} & \textcolor{Col4}{M_2}
    \end{array}\right].
\end{equation}
The \Def{row shifted shuffle} $M_1 \Lshuffle M_2$ of $M_1$ and $M_2$
is the set of all matrices obtained by shuffling the rows of
$M_1 \bullet n_2$ with the rows of $n_1 \bullet M_2$. By a routine
computation, we obtain the following expression for the product of
$\PM{k}^\star$:
\begin{equation}
    \BasisF^\star_{M_1} \cdot \BasisF^\star_{M_2} =
    \sum_{M \in M_1 {\mbox{\tiny \begin{math}\Lshuffle\end{math}}} M_2}
    \BasisF^\star_M.
\end{equation}

For instance, in $\PM{1}^\star$ one has
\begin{equation}
    \BasisF^\star_{\Matrix{
        \textcolor{Col1}{\Zero} & \textcolor{Col1}{\One} \\
        \textcolor{Col1}{\One} & \textcolor{Col1}{\One}
    }}
    \cdot
    \BasisF^\star_{\Matrix{
        \textcolor{Col4}{\One} & \textcolor{Col4}{\Zero} \\
        \textcolor{Col4}{\Zero} & \textcolor{Col4}{\One}
    }}
    =
    \BasisF^\star_{\Matrix{
        \textcolor{Col1}{\Zero} & \textcolor{Col1}{\One} & \Zero
            & \Zero \\
        \textcolor{Col1}{\One} & \textcolor{Col1}{\One} & \Zero
            & \Zero \\
        \Zero & \Zero & \textcolor{Col4}{\One}
            & \textcolor{Col4}{\Zero} \\
        \Zero & \Zero & \textcolor{Col4}{\Zero}
            & \textcolor{Col4}{\One}
    }}
    +
    \BasisF^\star_{\Matrix{
        \textcolor{Col1}{\Zero} & \textcolor{Col1}{\One} & \Zero
            & \Zero \\
        \Zero & \Zero & \textcolor{Col4}{\One}
            & \textcolor{Col4}{\Zero} \\
        \textcolor{Col1}{\One} & \textcolor{Col1}{\One} & \Zero
            & \Zero \\
        \Zero & \Zero & \textcolor{Col4}{\Zero}
            & \textcolor{Col4}{\One}
    }}
    +
    \BasisF^\star_{\Matrix{
        \Zero & \Zero & \textcolor{Col4}{\One}
            & \textcolor{Col4}{\Zero} \\
        \textcolor{Col1}{\Zero} & \textcolor{Col1}{\One} & \Zero
            & \Zero \\
        \textcolor{Col1}{\One} & \textcolor{Col1}{\One} & \Zero
            & \Zero \\
        \Zero & \Zero & \textcolor{Col4}{\Zero}
            & \textcolor{Col4}{\One}
    }}
    +
    \BasisF^\star_{\Matrix{
        \textcolor{Col1}{\Zero} & \textcolor{Col1}{\One} & \Zero
            & \Zero \\
        \Zero & \Zero & \textcolor{Col4}{\One}
            & \textcolor{Col4}{\Zero} \\
        \Zero & \Zero & \textcolor{Col4}{\Zero}
            & \textcolor{Col4}{\One} \\
        \textcolor{Col1}{\One} & \textcolor{Col1}{\One} & \Zero & \Zero
    }}
    +
    \BasisF^\star_{\Matrix{
        \Zero & \Zero & \textcolor{Col4}{\One}
            & \textcolor{Col4}{\Zero} \\
        \textcolor{Col1}{\Zero} & \textcolor{Col1}{\One} & \Zero
            & \Zero \\
        \Zero & \Zero & \textcolor{Col4}{\Zero}
            & \textcolor{Col4}{\One} \\
        \textcolor{Col1}{\One} & \textcolor{Col1}{\One} & \Zero & \Zero
    }}
    +
    \BasisF^\star_{\Matrix{
        \Zero & \Zero & \textcolor{Col4}{\One}
            & \textcolor{Col4}{\Zero} \\
        \Zero & \Zero & \textcolor{Col4}{\Zero}
            & \textcolor{Col4}{\One} \\
        \textcolor{Col1}{\Zero} & \textcolor{Col1}{\One} & \Zero
            & \Zero \\
        \textcolor{Col1}{\One} & \textcolor{Col1}{\One} & \Zero & \Zero
    }}.
\end{equation}
\medbreak

Let $M$ be a $k$-packed matrix. By duality, the coproduct in
$\PM{k}^\star$ satisfies
\begin{equation}
    \Coproduct\left(\BasisF^\star_M\right) =
    \sum_{M_1, M_2 \in \PackedMatrices_k}
    \left\langle \BasisF_{M_1} \cdot \BasisF_{M_2},
    \BasisF^\star_M \right\rangle
    \; \BasisF^\star_{M_1} \otimes \BasisF^\star_{M_2}.
\end{equation}
By a routine computation, we obtain the following expression for the
coproduct of $\PM{k}^\star$:
\begin{equation}
    \Coproduct\left(\BasisF^\star_M\right) =
    \sum_{M = M_1 \DecompC M_2} \BasisF^\star_{\Compr(M_1)}
        \otimes \BasisF^\star_{\Compr(M_2)}.
\end{equation}
For instance, in $\PM{1}^\star$ one has
\begin{equation}
    \Coproduct
    \BasisF^\star_{\Matrix{
        \Zero & \Zero & \One & \Zero \\
        \Zero & \Zero & \Zero & \One \\
        \One & \Zero & \Zero & \Zero \\
        \One & \One & \Zero & \Zero
    }}
    =
    \BasisF^\star_{\Matrix{
        \Zero & \Zero & \One & \Zero \\
        \Zero & \Zero & \Zero & \One \\
        \One & \Zero & \Zero & \Zero \\
        \One & \One & \Zero & \Zero
    }}
    \otimes
    \BasisF^\star_{\emptyset}
    +
    \BasisF^\star_{\Matrix{
        \One & \Zero \\
        \One & \One
    }}
    \otimes
    \BasisF^\star_{\Matrix{
        \One & \Zero \\
        \Zero & \One
    }}
    +
    \BasisF^\star_{\Matrix{
        \Zero & \Zero & \One \\
        \One & \Zero & \Zero \\
        \One & \One & \Zero \\
    }}
    \otimes
    \BasisF^\star_{\Matrix{
        \One
    }}
    +
    \BasisF^\star_{\emptyset}
    \otimes
    \BasisF^\star_{\Matrix{
        \Zero & \Zero & \One & \Zero \\
        \Zero & \Zero & \Zero & \One \\
        \One & \Zero & \Zero & \Zero \\
        \One & \One & \Zero & \Zero
    }}.
\end{equation}
\medbreak

Let us denote by $M^\intercal$ the transpose of $M$.
\medbreak

\begin{Proposition} \label{prop:self_duality_PM}
    The map $\phi : \PM{k} \to \PM{k}^\star$ linearly defined for any
    $k$-packed matrix $M$ by
    \begin{equation}
        \phi\left(\BasisF_M\right) := \BasisF^\star_{M^\intercal}
    \end{equation}
    is a Hopf isomorphism.
\end{Proposition}
\medbreak

Since the transpose of any packed matrix of
$\PackedMatrices_{k, n, \ell}$ also belongs to
$\PackedMatrices_{k, n, \ell}$, Proposition~\ref{prop:self_duality_PM}
also implies that $\PMN{k}$ and $\PML{k}$ are self-dual for the
isomorphism~$\phi$.
\medbreak

\subsubsection{Primitive elements}
For any $k$-packed matrix $M$, define
\begin{equation}
    \BasisW^M :=
    \BasisF^\star_{M_1} \cdot \ldots \cdot \BasisF^\star_{M_r}
\end{equation}
where the $M_i$ are connected packed matrices (see
Section~\ref{subsubsec:freeness_PM}) and
$M = M_1 \Over \dots \Over M_r$. Then, we have
\begin{equation}
    \BasisW^M = \BasisF^\star_M + \sum_{M' \in R} \BasisF^\star_{M'}
\end{equation}
where any matrix $M'$ of $R$ satisfies
$M^\intercal \OrdPM M'^\intercal$ since the product in $\PM{k}^\star$
consists in shifting and shuffling rows of matrices. Thus, by
triangularity, the $\BasisW^M$ form a basis of $\PM{k}^\star$. Moreover,
for any $k$-packed matrices $M_1$ and $M_2$, the product of
$\PM{k}^\star$ can be expressed as
\begin{equation}
    \BasisW^{M_1} \cdot \BasisW^{M_2} = \BasisW^{M_1 \Over M_2}.
\end{equation}
\medbreak

Let us denote by $\BasisV_M$, where the $M$ are $k$-packed matrices, the
adjoint elements of the~$\BasisW^M$.
\begin{Proposition} \label{prop:primitive_elements_PM}
    The elements $\BasisV_M$, where $M$ are connected $k$-packed
    matrices, form a basis of the vector space of primitive elements
    of~$\PM{k}$.
\end{Proposition}
\medbreak

By Proposition~\ref{prop:primitive_elements_PM}, the $\BasisV_M$, where
$M$ are connected $k$-packed matrices, generate the Lie algebra of
primitive elements of $\PM{k}$. The first few dimensions of the Lie
algebras of primitive elements of $\PMN{1}$, $\PMN{2}$, $\PML{1}$,
$\PML{2}$ are respectively given
by~\eqref{equ:dimensions_generators_PMN1},
\eqref{equ:dimensions_generators_PMN2},
\eqref{equ:dimensions_generators_PML1},
and~\eqref{equ:dimensions_generators_PML2}.
\medbreak

\subsection{Bidendriform bialgebra structure}
We show here that $\PM{k}$ admits a bidendriform bialgebra
structure~\cite{Foi07} (see also
Section~\ref{subsubsec:dendriform_algebras} of
Chapter~\ref{chap:algebra}).
\medbreak

\subsubsection{Dendriform algebra structure}
We denote by $\PM{k}^+$ the subspace of $\PM{k}$ restricted on
nonempty matrices. For any nonempty matrix $M$, we shall denote by
$\LastColumn(M)$ its last column. Let us endow $\PM{k}^+$ with two
products $\LDendr$ and $\RDendr$ linearly defined, for any nonempty
$k$-packed matrices $M_1$ and $M_2$ of respective sizes $n_1$ and $n_2$,
by
\begin{equation} \label{equ:left_product_PM}
    \BasisF_{M_1} \LDendr \BasisF_{M_2} :=
    \sum_{\substack{M \in M_1 \cshuffle M_2 \\
    \LastColumn(M) = \LastColumn\left(M_1 \circ n_2\right)}}
    \BasisF_M
\end{equation}
and
\begin{equation} \label{equ:right_product_PM}
    \BasisF_{M_1} \RDendr \BasisF_{M_2} :=
    \sum_{\substack{M \; \in \; M_1 \cshuffle M_2 \\
    \LastColumn(M) = \LastColumn\left(n_1 \circ {M_2}\right)}}
    \BasisF_M.
\end{equation}
\medbreak

In other words, the matrices appearing in a $\LDendr$-product (resp.
$\RDendr$-product) on the fundamental basis involving $M_1$ and $M_2$
are the matrices $M$ obtained by shifting and shuffling the columns
of $M_1$ and $M_2$ such that the last column of $M$ comes from $M_1$
(resp. $M_2$). For example,
\begin{subequations}
\begin{equation}
    \BasisF_{\Matrix{
        \textcolor{Col1}{\Zero} & \textcolor{Col1}{\One} \\
        \textcolor{Col1}{\One} & \textcolor{Col1}{\One}
    }}
    \LDendr
    \BasisF_{\Matrix{
        \textcolor{Col4}{\One} & \textcolor{Col4}{\Zero} \\
        \textcolor{Col4}{\Zero} & \textcolor{Col4}{\One}
    }}
    =
        \BasisF_{\Matrix{
        \textcolor{Col1}{\Zero} & \Zero & \Zero
            & \textcolor{Col1}{\One} \\
        \textcolor{Col1}{\One} & \Zero & \Zero
            & \textcolor{Col1}{\One} \\
        \Zero & \textcolor{Col4}{\One}
            & \textcolor{Col4}{\Zero} & \Zero \\
        \Zero & \textcolor{Col4}{\Zero}
            & \textcolor{Col4}{\One} & \Zero
    }}
    +
    \BasisF_{\Matrix{
        \Zero & \textcolor{Col1}{\Zero} & \Zero
            & \textcolor{Col1}{\One} \\
        \Zero & \textcolor{Col1}{\One} & \Zero
            & \textcolor{Col1}{\One} \\
        \textcolor{Col4}{\One} & \Zero
            & \textcolor{Col4}{\Zero} & \Zero \\
        \textcolor{Col4}{\Zero} & \Zero
            & \textcolor{Col4}{\One} & \Zero
    }}
    +
    \BasisF_{\Matrix{
        \Zero & \Zero & \textcolor{Col1}{\Zero}
            & \textcolor{Col1}{\One} \\
        \Zero & \Zero & \textcolor{Col1}{\One}
            & \textcolor{Col1}{\One} \\
        \textcolor{Col4}{\One} & \textcolor{Col4}{\Zero}
            & \Zero & \Zero \\
        \textcolor{Col4}{\Zero} & \textcolor{Col4}{\One}
            & \Zero & \Zero
    }},
\end{equation}
\begin{equation}
    \BasisF_{\Matrix{
        \textcolor{Col1}{\Zero} & \textcolor{Col1}{\One} \\
        \textcolor{Col1}{\One} & \textcolor{Col1}{\One}
    }}
    \RDendr
    \BasisF_{\Matrix{
        \textcolor{Col4}{\One} & \textcolor{Col4}{\Zero} \\
        \textcolor{Col4}{\Zero} & \textcolor{Col4}{\One}
    }}
    =
    \BasisF_{\Matrix{
        \textcolor{Col1}{\Zero} & \textcolor{Col1}{\One} & \Zero
            & \Zero \\
        \textcolor{Col1}{\One} & \textcolor{Col1}{\One} & \Zero
            & \Zero \\
        \Zero & \Zero & \textcolor{Col4}{\One}
            & \textcolor{Col4}{\Zero} \\
        \Zero & \Zero & \textcolor{Col4}{\Zero}
            & \textcolor{Col4}{\One}
    }}
    +
    \BasisF_{\Matrix{
        \textcolor{Col1}{\Zero} & \Zero & \textcolor{Col1}{\One}
            & \Zero \\
        \textcolor{Col1}{\One} & \Zero & \textcolor{Col1}{\One}
            & \Zero \\
        \Zero & \textcolor{Col4}{\One} & \Zero
            & \textcolor{Col4}{\Zero} \\
        \Zero & \textcolor{Col4}{\Zero} & \Zero
            & \textcolor{Col4}{\One}
    }}
    +
    \BasisF_{\Matrix{
        \Zero & \textcolor{Col1}{\Zero} & \textcolor{Col1}{\One}
            & \Zero \\
        \Zero & \textcolor{Col1}{\One} & \textcolor{Col1}{\One}
            & \Zero \\
        \textcolor{Col4}{\One} & \Zero & \Zero
            & \textcolor{Col4}{\Zero} \\
        \textcolor{Col4}{\Zero} & \Zero & \Zero
            & \textcolor{Col4}{\One}
    }}.
\end{equation}
\end{subequations}
\medbreak

Since the last column of any matrix appearing in the shifted shuffle of
two matrices comes from exactly  of the two operands, for any nonempty
packed matrices $M_1$ and $M_2$, one obviously has
\begin{equation} \label{equ:splitting_product_PM}
    \BasisF_{M_1} \cdot \BasisF_{M_2} =
    \BasisF_{M_1} \LDendr \BasisF_{M_2}
    + \BasisF_{M_1} \RDendr \BasisF_{M_2}.
\end{equation}
\medbreak

\begin{Proposition} \label{prop:dendriform_structure_PM}
    The Hopf algebra $\PM{k}$ admits a dendriform algebra structure for
    the products $\LDendr$ and $\RDendr$.
\end{Proposition}
\medbreak

\subsubsection{Codendriform coalgebra structure}
For any nonempty matrix $M$, we shall denote by $\LastRow(M)$
its last row. Let us endow $\PM{k}$ with two coproducts $\LCoDendr$
and $\RCoDendr$ linearly defined, for any nonempty $k$-packed matrix $M$,
by
\begin{equation}
    \LCoDendr\left(\BasisF_M\right) :=
    \sum_{\substack{
        M = L \bullet R \\ \LastRow(L \bullet r) = \LastRow(M)
    }}
    \BasisF_{\Compr(L)} \otimes \BasisF_{\Compr(R)}
\end{equation}
and
\begin{equation}
    \RCoDendr\left(\BasisF_M\right) :=
    \sum_{\substack{
        M = L \bullet R \\ \LastRow(\ell \bullet R) = \LastRow(M)
    }}
    \BasisF_{\Compr(L)} \otimes \BasisF_{\Compr(R)},
\end{equation}
where $r$ (resp. $\ell$) is the number of columns of $R$ (resp. $L$). In
other words, the pairs of matrices appearing in a $\LCoDendr$-coproduct
(resp. $\RCoDendr$-coproduct) in the fundamental basis are the pairs
$(L, R)$ of packed matrices such that the last row of $L$ (resp. $R$)
comes from the last row of $M$. For example,
\begin{subequations}
\begin{equation}
    \LCoDendr
    \BasisF_{\Matrix{
        \One & \Zero & \Zero & \Zero & \Zero & \Zero \\
        \Zero & \One & \One & \Zero & \Zero & \Zero \\
        \Zero & \Zero & \Zero & \Zero & \Zero & \One \\
        \Zero & \Zero & \Zero & \One & \Zero & \Zero \\
        \Zero & \Zero & \Zero & \One & \One & \Zero \\
        \Zero & \Zero & \textcolor{Col1}{\One} & \Zero & \Zero & \Zero
    }}
    =
    \BasisF_{\Matrix{
        \One & \Zero & \Zero \\
        \Zero & \One & \One \\
        \Zero & \Zero & \textcolor{Col1}{\One}
    }}
    \otimes
    \BasisF_{\Matrix{
        \Zero & \Zero & \One \\
        \One & \Zero & \Zero \\
        \One & \One & \Zero
    }}
    +
    \BasisF_{\Matrix{
        \One & \Zero & \Zero & \Zero & \Zero \\
        \Zero & \One & \One & \Zero & \Zero \\
        \Zero & \Zero & \Zero & \One & \Zero \\
        \Zero & \Zero & \Zero & \One & \One \\
        \Zero & \Zero & \textcolor{Col1}{\One} & \Zero & \Zero
    }}
    \otimes
    \BasisF_{\Matrix{
        \One
    }},
\end{equation}
\begin{equation}
    \RCoDendr
    \BasisF_{\Matrix{
        \One & \Zero & \Zero & \Zero & \Zero & \Zero \\
        \Zero & \One & \One & \Zero & \Zero & \Zero \\
        \Zero & \Zero & \Zero & \Zero & \Zero & \One \\
        \Zero & \Zero & \Zero & \One & \Zero & \Zero \\
        \Zero & \Zero & \Zero & \One & \One & \Zero \\
        \Zero & \Zero & \textcolor{Col4}{\One} & \Zero & \Zero
            & \Zero
    }}
    =
    \BasisF_{\Matrix{
        \One
    }}
    \otimes
    \BasisF_{\Matrix{
        \One & \One & \Zero & \Zero & \Zero \\
        \Zero & \Zero & \Zero & \Zero & \One \\
        \Zero & \Zero & \One & \Zero & \Zero \\
        \Zero & \Zero & \One & \One & \Zero \\
        \Zero & \textcolor{Col4}{\One} & \Zero & \Zero & \Zero
    }}.
\end{equation}
\end{subequations}
\medbreak

Since by Lemma \ref{lem:packed_matrices_decomposition}, one cannot
vertically split a packed matrix by separating two nonzero entries on a
same row, for any nonempty packed matrix $M$, one has
\begin{equation} \label{equ:splitting_coproduct_PM}
    \Coproduct\left(\BasisF_M\right) =
    1 \otimes \BasisF_M + \LCoDendr\left(\BasisF_M\right) +
    \RCoDendr\left(\BasisF_M\right) + \BasisF_M \otimes 1.
\end{equation}
\medbreak

\begin{Proposition} \label{prop:codendriform_structure_PM}
    The Hopf algebra $\PM{k}$ admits a codendriform coalgebra structure
    for the coproducts $\LCoDendr$ and $\RCoDendr$.
\end{Proposition}
\medbreak

\subsubsection{Bidendriform bialgebra structure}

\begin{Theorem} \label{thm:bidendriform_structure_PM}
    The Hopf bialgebra $\PM{k}$ admits a bidendriform bialgebra
    structure for the products $\LDendr$, $\RDendr$ and the coproducts
    $\LCoDendr$, $\RCoDendr$.
\end{Theorem}
\medbreak

Theorem~\ref{thm:bidendriform_structure_PM} also implies that $\PMN{k}$
and $\PML{k}$ admit a bidendriform bialgebra structure.
Following~\cite{Foi07}, the generating series
$\SeriesTotPrimitiveElements_{k, n}(t)$ and
$\SeriesTotPrimitiveElements_{k, \ell}(t)$ of totally primitive elements
of $\PMN{k}$ and $\PML{k}$ satisfy respectively
\begin{equation}
    \SeriesTotPrimitiveElements_{k, n}(t) =
    \frac{\HilbSeries_{k, n}(t) - 1}{\HilbSeries_{k, n}(t)^2}
    \qquad \mbox{ and } \qquad
    \SeriesTotPrimitiveElements_{k, \ell}(t) =
    \frac{\HilbSeries_{k, \ell}(t) - 1}{\HilbSeries_{k, \ell}(t)^2}.
\end{equation}
The first few dimensions of totally primitive elements of $\PMN{1}$ and
$\PMN{2}$ are respectively
\begin{equation}
    0, 1, 5, 240, 40404, 24827208, 57266105928
\end{equation}
and
\begin{equation}
    0, 2, 48, 15640, 39023776, 813415850016, 147655768992433664.
\end{equation}
There are respectively Sequences~\OEIS{A230885} and~\OEIS{A230886}
of~\cite{Slo}. The first few dimensions of totally primitive elements
of $\PML{1}$ and $\PML{2}$ are respectively
\begin{equation}
    0, 1, 0, 5, 36, 381, 4720, 67867, 1109434
\end{equation}
and
\begin{equation}
    0, 2, 0, 40, 576, 12192, 302080, 8686976, 284015104.
\end{equation}
These are respectively Sequences~\OEIS{A230887} and~\OEIS{A230888}
of~\cite{Slo}.
\medbreak

\section{Related Hopf bialgebras} \label{sec:links_packed_matrices}
In this section, we describe links between $\PM{k}$ and some already
known Hopf bialgebras. Next, we provide a method to construct Hopf
sub-bialgebras of~$\PM{k}$.
\medbreak

\subsection{Links with known bialgebras}
We consider here the Hopf bialgebras of $k$-colored permutations, of
uniform block permutations, and of matrix quasi-symmetric functions.
\medbreak

\subsubsection{Hopf bialgebra of colored permutations}
The Hopf bialgebra $\FQSym^{(k)}$ of $k$-colored permutations is
introduced in~\cite{NT10} (see also
Section~\ref{subsubsec:FQSym_colored} of Chapter~\ref{chap:algebra}).
\medbreak

\begin{Proposition} \label{prop:morphism_FQSym_k_PM}
    The map $\alpha_k : \FQSym^{(k)} \to \PMN{k}$ linearly defined,
    for any $k$-colored permutation $(\sigma, c)$ by
    \begin{equation}
        \alpha_k\left(\BasisF_{(\sigma, c)}\right)
        := \BasisF_{M^{(\sigma, c)}}
    \end{equation}
    where $M^{(\sigma, c)}$ is the $k$-packed matrix satisfying
    $M^{(\sigma, c)}_{ij} = c_j \, \delta_{i, \sigma_j}$ is an
    injective Hopf morphism.
\end{Proposition}
\medbreak

In particular, Proposition \ref{prop:morphism_FQSym_k_PM} shows
that $\PMN{1}$ contains $\FQSym$. Notice that the map $\alpha_k$ is
still well-defined on the codomain $\PML{k}$ instead of $\PMN{k}$.
\medbreak

\subsubsection{Hopf bialgebra of uniform block permutations}
The Hopf bialgebra $\UBP$ of uniform block permutations is introduced
in~\cite{AO08} (see also Section~\ref{subsubsec:UBP} of
Chapter~\ref{chap:algebra}).
\medbreak

\begin{Proposition} \label{prop:morphism_PM_UBP}
    The map $\beta : \UBP^\star \to \PMN{1}$
    linearly defined, for any UBP $\pi$ by
    \begin{equation}
        \beta\left(\BasisF^\star_\pi\right) := \BasisF_{M^\pi}
    \end{equation}
    where $M^\pi$ is the $1$-packed matrix satisfying
    \begin{equation}
        M^\pi_{ij} :=
        \begin{cases}
            1 & \mbox{if there is }
                e \in \pi^d \mbox{ such that }
                j \in e \mbox{ and } i \in \pi(e), \\
            0 & \mbox{otherwise}.
        \end{cases}
    \end{equation}
    is an injective Hopf morphism.
\end{Proposition}
\medbreak

For example, if $\pi$ is the UBP defined by
\begin{equation}
    \pi(\{1, 4, 5\}) := \{2, 5, 6\}, \quad
    \pi(\{2\}) := \{1\}, \quad
    \mbox{and} \quad
    \pi(\{3, 6\}) := \{3, 4\},
\end{equation}
we have
\begin{equation}
    \beta\left(\BasisF^\star_\pi\right) =
    \BasisF_{\Matrix{
        \Zero & \One & \Zero & \Zero & \Zero & \Zero \\
        \One & \Zero & \Zero & \One & \One & \Zero \\
        \Zero & \Zero & \One & \Zero & \Zero & \One \\
        \Zero & \Zero & \One & \Zero & \Zero & \One \\
        \One & \Zero & \Zero & \One & \One & \Zero \\
        \One & \Zero & \Zero & \One & \One & \Zero
    }}.
\end{equation}
\medbreak

The existence of this particular morphism $\beta$ exhibited by
Proposition~\ref{prop:morphism_PM_UBP} implies that $\UBP^\star$ is
free (as an associative algebra), cofree (as a coassociative
coalgebra), self-dual, and admits bidendriform bialgebra structure.
\medbreak

Besides, by using same arguments as those used in
Section~\ref{subsec:freeness_PM}, one can build multiplicative bases
of $\UBP^\star$ by setting, for any UBP $\pi$,
\begin{equation}
    \BasisE^\star_{M^\pi} :=
        \sum_{M^\pi \OrdPM M^{\pi'}} \BasisF_{M^{\pi'}}
    \qquad \mbox{and} \qquad
    \BasisH^\star_{M^\pi} :=
        \sum_{M^{\pi'} \OrdPM M^\pi} \BasisF_{M^{\pi'}}.
\end{equation}
This gives another way to prove the freeness of $\UBP^\star$ by using
same arguments as those of Theorem~\ref{thm:freeness_PM}. Hence,
$\UBP^\star$ is freely generated by the elements $\BasisE_{M^\pi}$
(resp. $\BasisH_{M^\pi}$) where the $\pi$ are UBPs such that the $M^\pi$
are connected (resp. anti-connected) $1$-packed matrices. The first few
numbers of algebraic generators of $\UBP^\star$ are
\begin{equation}
    0, 1, 2, 11, 98, 1202, 19052, 375692,  8981392, 255253291,
    8488918198
\end{equation}
and the first few dimensions of totally primitive elements are
\begin{equation}
    0, 1, 1, 7, 72, 962, 16135, 330624, 8117752, 235133003, 7929041828.
\end{equation}
These are Sequence~\OEIS{A230889} and~\OEIS{ A230890} of~\cite{Slo}.
\medbreak

Moreover, since for any UBP $\pi$, there exists a UBP $\pi^{-1}$ such
that the transpose of $M^\pi$ is $M^{\pi^{-1}}$, by
Proposition \ref{prop:self_duality_PM}, the map
$\phi : \UBP^\star \to \UBP$ linearly defined for any UBP $\pi$ by
\begin{equation}
    \phi\left(\BasisF^\star_{M^\pi}\right)
    := \BasisF_{{M^\pi}^\intercal}
\end{equation}
is an isomorphism.
\medbreak

\subsubsection{Algebra of matrix quasi-symmetric functions}
The Hopf algebra of matrix quasi-symmetric functions is introduced
in~\cite{DHT02} (see also \cite{Hiv99} and
Section~\ref{subsubsec:FQSym} of Chapter~\ref{chap:algebra}).
\medbreak

\begin{figure}[ht]
    \begin{tikzpicture}[scale=.45,font=\footnotesize]
        \filldraw[draw=ColBlack,fill=Col4!40](0,0)rectangle (3,-1);
        \draw[draw=ColBlack,fill=ColWhite](3,0)rectangle node
            {\begin{math}0\end{math}}(7,-1);
        \draw[draw=ColBlack,fill=ColWhite](0,-1)rectangle node
            {\begin{math}0\end{math}} (5,-2);
        \filldraw[draw=ColBlack,fill=Col1!40](5,-1)rectangle(7,-2);
        \node at(-.5,-.5){\begin{math}i\end{math}};
        \node at(-.7,-1.5){\begin{math}i\!+\!1\end{math}};
        \draw[draw=ColBlack](3,.25)edge[<->]node[anchor=center,above]
            {\begin{math}z_i\end{math}}(7,.25);
        \draw[draw=ColBlack](0,-2.25)edge[<->]node[anchor=center,below]
            {\begin{math}z_{i + 1}\end{math}}(5,-2.25);
        \draw[draw=ColBlack](0,1.05)edge[<->]node[anchor=center,above]
            {\begin{math}m\end{math}}(7,1.05);
        \node at(9,-1){\begin{math}\rightharpoonup\end{math}};
        \filldraw[draw=ColBlack,fill=Col4!40](11,-.5)
            rectangle(14,-1.5);
        \draw[draw=ColBlack,fill=ColWhite](14,-.5)rectangle node
            {\begin{math}0\end{math}}(16,-1.5);
        \filldraw[draw=ColBlack,fill=Col1!40](16,-.5)
            rectangle(18,-1.5);
    \end{tikzpicture}
    \caption[The covering relation of the poset of matrices.]
    {The condition for overlaying the $i$th and $(i + 1)$st rows of a
    (not necessarily square) packed matrix according to the relation
    $\CoveringRelMQ$. The darker regions contain any entries and the
    white ones, only zeros.}
    \label{fig:covering_relation_matrices}
\end{figure}
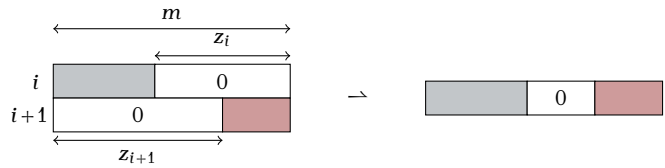

Let us endow the set of matrices indexing $\MQSym$ with a binary
relation $\CoveringRelMQ$ defined in the following way. If $M_1$ and
$M_2$ are two matrices such that $M_1$ has $n$ rows and $m$ columns, we
have $M_1 \CoveringRelMQ M_2$ if there is an index $i \in [n - 1]$
such that, denoting by $z_i$ the number of $0$ which end the $i$th row
of $M_1$, and by $z_{i + 1}$ the number of $0$ which start the
$(i + 1)$st row of $M_1$, one has $z_i + z_{i + 1} \geq m$ and $M_2$
is obtained from $M_1$ by overlaying its $i$th and $(i + 1)$st rows (see
Figure~\ref{fig:covering_relation_matrices}).
\medbreak

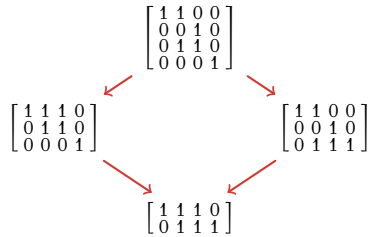
\begin{figure}[ht]
     \begin{tikzpicture}[scale=.6]
         \node(1)at(0,0){\begin{math}
         \Matrix{
             \One & \One & 0   & 0 \\
             0   & 0   & \One & 0 \\
             0   & \One & \One & 0 \\
             0   & 0   &0    & \One
         }\end{math}};
         \node(2)at(-3,-2){\begin{math}
         \Matrix{
             \One & \One & \One & 0 \\
             0   & \One & \One & 0 \\
             0   &0    &0    & \One
         }\end{math}};
         \node(3)at(3,-2){\begin{math}
         \Matrix{
             \One & \One & 0   & 0 \\
             0   & 0   &\One  & 0 \\
             0   &\One  &\One  & \One
         }\end{math}};
         \node(4)at(0,-4){\begin{math}
         \Matrix{
             \One & \One &\One  & 0 \\
             0   &\One  &\One  & \One
         }\end{math}};
         \draw[Arc](1)--(2);
         \draw[Arc](1)--(3);
         \draw[Arc](2)--(4);
         \draw[Arc](3)--(4);
     \end{tikzpicture}
    \caption[An interval of the poset of matrices.]
    {The Hasse diagram of an interval for the order $\OrdMQ$ on (not
    necessarily square) packed matrices.}
    \label{fig:example_interval_matrices}
\end{figure}

We now endow the set of matrices that index $\MQSym$ with the partial
order relation $\OrdMQ$ defined as the reflexive and transitive closure
of $\CoveringRelMQ$. Figure~\ref{fig:example_interval_matrices} shows
an interval of this partial order.
\medbreak

\begin{Proposition} \label{prop:morphism_PM_MQSym}
    The map $\gamma : \PML{1}^\star \to \MQSym$ linearly defined, for
    any $1$-packed matrix $M$ by
    \begin{equation}
        \gamma\left(\BasisF^\star_M\right)
        := \sum_{M \OrdMQ M'} \BasisM_{M'},
    \end{equation}
    is an injective associative algebra morphism.
\end{Proposition}
\medbreak

For instance, one has
\begin{equation}
    \gamma\left(
    \BasisF^\star_{
    \Matrix{
        \One & \One & \Zero & \Zero \\
        \Zero & \Zero & \One & \Zero \\
        \Zero & \One & \One & \Zero \\
        \Zero & \Zero & \Zero & \One
    }}\right) =
    \BasisM_{
    \Matrix{
        \One & \One & 0 & 0 \\
        0 & 0 & \One & 0 \\
        0 & \One & \One & 0 \\
        0 & 0 & 0 & \One
    }}
    +
    \BasisM_{
    \Matrix{
        \One & \One & \One & 0 \\
        0 & \One & \One & 0 \\
        0 & 0 & 0 & \One
    }}
    +
    \BasisM_{
    \Matrix{
        \One & \One & 0 & 0 \\
        0 & 0 & \One & 0 \\
        0 & \One & \One & \One
    }}
    +
    \BasisM_{
    \Matrix{
        \One & \One & \One & 0 \\
        0 & \One & \One & \One
    }}.
\end{equation}
\medbreak

Notice that $\gamma$ is not a Hopf morphism since it is not a
coassociative coalgebra morphism. Indeed, we have
\begin{equation}
    \Coproduct\left( \gamma\left(
    \BasisF^\star_{
    \Matrix{
        \One & \One \\
        \One & \Zero
    }}\right)\right)
    =
    1 \otimes
    \BasisM_{
    \Matrix{
        \One & \One \\
        \One & \Zero
    }}
    +
    \BasisM_{
    \Matrix{
        \One & \One \\
        \One & \Zero
    }}
    \otimes 1,
\end{equation}
but
\begin{equation}
    (\gamma \otimes \gamma)\left( \Coproduct\left(
    \BasisF^\star_{
    \Matrix{
        \One & \One \\
        \One & \Zero
    }}\right)\right)
    =
        1 \otimes
    \BasisM_{
    \Matrix{
        \One & \One \\
        \One & \Zero
    }}
    +
    \BasisM_{
    \Matrix{
        \One & \One
    }}
    \otimes
    \BasisM_{
    \Matrix{
        \One
    }}
    +
    \BasisM_{
    \Matrix{
        \One & \One \\
        \One & \Zero
    }}
    \otimes 1.
\end{equation}
\medbreak

\subsubsection{Diagram of embeddings}
\begin{figure}[ht]
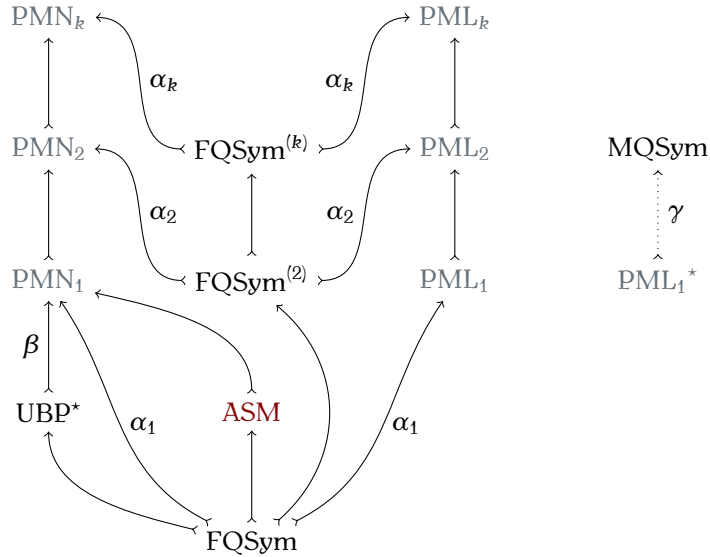

    \centering

    \caption[The diagram of Hopf bialgebras of packed matrices and
    related structures.]
    {The diagram of Hopf bialgebras of packed matrices and related
    structures. Arrows~$\rightarrowtail$ are injective Hopf bialgebra
    morphisms. The dotted arrow is an associative algebra morphism.}
    \label{fig:diagram_Hopf_algebras_packed_matrices}
\end{figure}

The diagram of Figure~\ref{fig:diagram_Hopf_algebras_packed_matrices}
summarizes the relations between known Hopf algebras related to
$\PM{k}$ and, more specifically, to its simply graded versions $\PMN{k}$
and $\PML{k}$. The Hopf bialgebra $\ASM$ is the subject of
Section \ref{sec:hopf_algebra_ASM}.
\medbreak

\subsection{Equivalence relations and Hopf sub-bialgebras}
\label{subsec:congruences_PM}
We provide here a way to construct Hopf sub-bialgebras of $\PM{k}$
analogous to the way using congruences to construct Hopf sub-bialgebras
of $\FQSym$ (see Section~\ref{subsubsec:subbialgebras_FQSym_congruences}
of Chapter~\ref{chap:algebra}).
\medbreak

\subsubsection{The monoid of words of columns}
Let $C_k^*$ be the free monoid generated by the set $C_k$ of all
$n \times 1$-matrices with entries in $A_k$, for all $n \geq 1$. In
other words, the elements of $C_k^*$ are words whose letters are columns
and its product $\bullet$ is the concatenation of such words. When all
the letters of an element $M \in C_k^*$ have, as columns, a same number
of rows, $M$ is a matrix and we shall denote it as such in the sequel.
\medbreak

The alphabet $C_k$ is naturally equipped with the total order $\leq$
where, for any $c_1, c_2 \in C_k$, $c_1 \leq c_2$ if and only if the
bottom to top reading of the column $c_1$ is lexicographically smaller
than the bottom to top reading of $c_2$. For instance,
\begin{equation}
    \Matrix{
        \One \\
        \Zero \\
        \Zero
    }
    \leq
    \Matrix{
        \Zero \\
        \Zero \\
        \One
    }, \qquad
    \Matrix{
        \Zero \\
        \Zero \\
        1 \\
        \Zero \\
        1
    }
    \leq
    \Matrix{
        \Zero \\
        \Zero \\
        1 \\
        1
    }, \qquad
    \Matrix{
        1 \\
        \Zero
    }
    \leq
    \Matrix{
        \Zero \\
        1 \\
        1 \\
        \Zero
    }, \qquad
    \Matrix{
        2 \\
        1 \\
        \Zero
    }
    \leq
    \Matrix{
        1 \\
        2 \\
        \Zero
    }.
\end{equation}
\medbreak

\begin{figure}[ht]
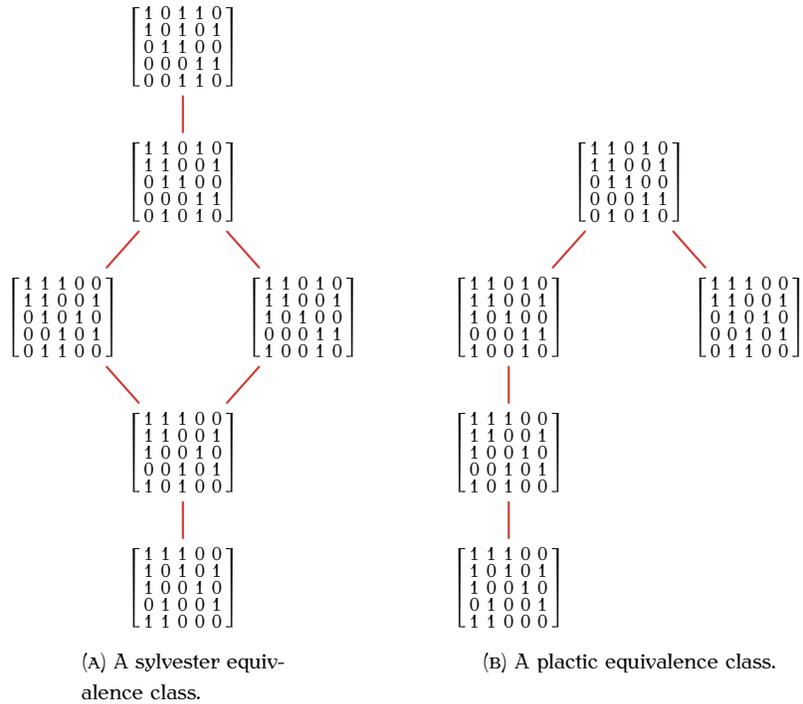

    \centering
    \subfloat[][A sylvester equivalence class.]{

    \label{subfig:sylvester_class_PM}}
    \qquad
    \subfloat[][A plactic equivalence class.]{
    %
    \label{subfig:plactic_class_PM}}
    \caption{Two equivalence classes of packed matrices.}
    \label{fig:exemples_classes_equ_matrices}
\end{figure}

Since $C_k$ is then totally ordered and $C_k^*$ is a free monoid, one
can consider the previous two congruences on $C_k^*$ instead on $A^*$.
For instance, Figure \ref{fig:exemples_classes_equ_matrices} represents
a $\CongrMonoid{S}$-equivalence class and a
$\CongrMonoid{P}$-equivalence class of packed matrices.
\medbreak

The order relation $\leq$ on $C_k$ is compatible with the shifted
shuffle of packed matrices in the following sense. Let $M_1$ and
$M_2$ be two nonempty packed matrices and $M$ be a matrix appearing in
$M_1 \cshuffle M_2$. Then, if $c_1$ (resp. $c_2$) is a column of $M$
coming from $M_1$ (resp. $M_2$), we necessarily have $c_1 \leq c_2$ and
$c_1 \ne c_2$. The obvious analogous property holds for words of $A^*$
and the shifted shuffle of words.
\medbreak

\subsubsection{Properties of equivalence relations}
An equivalence relation $\Equiv$ on $C_k^*$ is \Def{compatible with
the restriction to alphabet intervals} if for any interval $I$ of $C_k$
and for all $u, v \in C_k^*$, $u \Equiv v$ implies
$u_{|I} \Equiv v_{|I}$, where $u_{|I}$ denotes the word obtained by
erasing in $u$ the letters that are not in $I$.
\medbreak

Finally, we say that $\Equiv$ is \Def{compatible with the decompression
process} if for all $u, v \in C_k^*$ such that $u$ and $v$ are matrices,
$u \Equiv v$ if and only if $\Compr(u) \Equiv \Compr(v)$ and $u$ and $v$
have the same commutative image.
\medbreak

\subsubsection{Construction of Hopf sub-bialgebras}
Given an equivalence relation $\Equiv$ on the words of $C_k^*$ and
a $\Equiv$-equivalence class $[M]_\Equiv$ of packed matrices of $C_k^*$,
we consider the elements
\begin{equation} \label{equ:sum_congruence_class}
    \BasisP_{[M]_\Equiv} := \sum_{M' \in [M]_\Equiv} \BasisF_{M'}
\end{equation}
of $\PM{k}$.
\medbreak

One has for instance
\begin{equation}
    \BasisP_{\left[
    \Matrix{
        1 & 1 & 0 & 1 & 0 \\
        1 & 1 & 0 & 0 & 1 \\
        0 & 1 & 1 & 0 & 0 \\
        0 & 0 & 0 & 1 & 1 \\
        0 & 1 & 0 & 1 & 0
    }
    \right]_{\CongrMonoid{P}}} =
    \BasisF_{\Matrix{
        1 & 1 & 0 & 1 & 0 \\
        1 & 1 & 0 & 0 & 1 \\
        0 & 1 & 1 & 0 & 0 \\
        0 & 0 & 0 & 1 & 1 \\
        0 & 1 & 0 & 1 & 0
    }} +
    \BasisF_{\Matrix{
        1 & 1 & 0 & 1 & 0 \\
        1 & 1 & 0 & 0 & 1 \\
        1 & 0 & 1 & 0 & 0 \\
        0 & 0 & 0 & 1 & 1 \\
        1 & 0 & 0 & 1 & 0
    }} +
    \BasisF_{\Matrix{
        1 & 1 & 1 & 0 & 0 \\
        1 & 1 & 0 & 0 & 1 \\
        0 & 1 & 0 & 1 & 0 \\
        0 & 0 & 1 & 0 & 1 \\
        0 & 1 & 1 & 0 & 0
    }} +
    \BasisF_{\Matrix{
        1 & 1 & 1 & 0 & 0 \\
        1 & 1 & 0 & 0 & 1 \\
        1 & 0 & 0 & 1 & 0 \\
        0 & 0 & 1 & 0 & 1 \\
        1 & 0 & 1 & 0 & 0
    }} +
    \BasisF_{\Matrix{
        1 & 1 & 1 & 0 & 0 \\
        1 & 0 & 1 & 0 & 1 \\
        1 & 0 & 0 & 1 & 0 \\
        0 & 1 & 0 & 0 & 1 \\
        1 & 1 & 0 & 0 & 0
    }}.
\end{equation}
\medbreak

In particular, if $\Equiv$ is compatible with the decompression process,
any $\Equiv$-equivalence class of a packed matrix only contains packed
matrices. The family $\BasisP_{[M]_\Equiv}$, where the $[M]_\Equiv$ are
$\Equiv$-equivalence classes of packed matrices, forms then a basis of
a vector subspace of $\PM{k}$ denoted by~$\PM{k}^\Equiv$.
\medbreak

\begin{Theorem} \label{thm:substructure_PM_congruence}
    Let $\Equiv$ be an equivalence relation on the words of $C_k^*$ such
    that $\Equiv$
    \begin{enumerate}[label={\it (\roman*)}]
        \item \label{item:substructure_PM_congruence_congruence}
        is a monoid congruence on $C_k^*$;
        \item \label{item:substructure_PM_congruence_intervals}
        is compatible with the restriction to alphabet intervals;
        \item \label{item:substructure_PM_congruence_decompression}
        is compatible with the decompression process.
    \end{enumerate}
    Then, $\PM{k}^\Equiv$ is a Hopf sub-bialgebra of $\PM{k}$.
\end{Theorem}
\medbreak

Let $\Equiv$ be an equivalence relation on $C_k^*$
satisfying~\ref{item:substructure_PM_congruence_congruence},
\ref{item:substructure_PM_congruence_intervals},
and~\ref{item:substructure_PM_congruence_decompression} of
Theorem~\ref{thm:substructure_PM_congruence}. Note that since $\Equiv$
is compatible with the decompression process, any matrix contained in
a $\Equiv$-equivalence class $[M]_\Equiv$ is obtained by switching
columns of $M$. Then, any $\Equiv$-equivalence class $[M]_\Equiv$ of $k
$-packed matrices only contains matrices whose size and number of
nonzero entries are the same as in $M$. Hence,
Theorem~\ref{thm:substructure_PM_congruence} also implies that the
family~\eqref{equ:sum_congruence_class} forms a basis of Hopf
sub-bialgebras of both $\PMN{k}$ and $\PML{k}$. We respectively denote
these by $\PMN{k}^\Equiv$ and $\PML{k}^\Equiv$.
\medbreak

\subsubsection{Computer experiments}
\label{subsubsec:experiments_congruences_PM}
By Theorem~\ref{thm:substructure_PM_congruence}, the version of
sylvester, plactic, Baxter, Bell, hypoplactic, and total equivalence
relations (see Section~\ref{subsubsec:subbialgebras_FQSym_congruences}
of Chapter~\ref{chap:algebra}) applied to $C_k^*$ lead to bigraded Hopf
sub-bialgebras of $\PM{k}$.
Table~\ref{tab:dimensions_substructures_congruences_PM} shows first few
dimensions of the Hopf subalgebras of $\PMN{1}$ and $\PML{1}$ obtained
from these congruences, computed by computer exploration.
\begin{table}[ht]
    \centering
    \begin{tabular}{l||llllllll}
        Hopf bialgebra & \multicolumn{8}{c}{First dimensions}
            \\ \hline \hline
        $\PMN{1}^{\CongrMonoid{Bx}}$ & 1 & 1 & 7 & 265 & 38051 \\
        $\PMN{1}^{\CongrMonoid{Bl}}$ & 1 & 1 & 7 & 221 & 25789 \\
        $\PMN{1}^{\CongrMonoid{S}}$ & 1 & 1 & 7 & 221 & 24243 \\
        $\PMN{1}^{\CongrMonoid{P}}$ & 1 & 1 & 7 & 177 & 17339 \\
        $\PMN{1}^{\CongrMonoid{H}}$ & 1 & 1 & 7 & 177 & 13887 \\
        $\PMN{1}^{\CongrMonoid{T}}$ & 1 & 1 & 4 & 57 & 2306 \\ \hline
        $\PML{1}^{\CongrMonoid{Bx}}$
            & 1 & 1 & 2 & 10 & 68 & 578 & 5782 & 65745 \\
        $\PML{1}^{\CongrMonoid{Bl}}$
            & 1 & 1 & 2 & 9 & 53 & 390 & 3389 & 33881 \\
        $\PML{1}^{\CongrMonoid{S}}
            $ & 1 & 1 & 2 & 9 & 52 & 364 & 2918 & 26138 \\
        $\PML{1}^{\CongrMonoid{P}}$
            & 1 & 1 & 2 & 8 & 41 & 266 & 1976 & 16569 \\
        $\PML{1}^{\CongrMonoid{H}}$
            & 1 & 1 & 2 & 8 & 39 & 220 & 1396 & 9716 \\
        $\PML{1}^{\CongrMonoid{T}}$
            & 1 & 1 & 1 & 3 & 11 & 43 & 191 & 939
    \end{tabular}
    \bigbreak

    \caption[First dimensions of Hopf sub-bialgebras of packed
    matrices.]
    {First few dimensions of the Hopf sub-bialgebras $\PMN{1}^\Equiv$
    and $\PML{1}^\Equiv$, where $\Equiv$ is successively the Baxter,
    Bell, sylvester, plactic, hypoplactic, and total congruence.}
    \label{tab:dimensions_substructures_congruences_PM}
\end{table}
\medbreak

\section{Alternating sign matrices} \label{sec:hopf_algebra_ASM}
In this last section of the chapter, we construct and study a Hopf
sub-bialgebra of $\PM{1}$ whose bases are indexed by ASMs. We provide
its main properties and investigate how usual statistics on ASMs behave
algebraically inside it.
\medbreak

\subsection{Hopf bialgebra structure}
Let us explain how to encode ASMs by particular $1$-packed matrices.
As a consequence, we obtain a Hopf bialgebra on ASMs.
\medbreak

\subsubsection{From ASMs to $1$-packed matrices}
Let $\delta$ be an ASM~\cite{MRR83} (see also
Section~\ref{subsec:ASMs} of Chapter~\ref{chap:combinatorics}). We
denote by $M^\delta$ the matrix satisfying
\begin{equation}
    M^\delta_{ij} :=
    \begin{cases}
        \One & \mbox{if } \delta_{ij} \in \{\Plus, \Minus\}, \\
        \Zero & \mbox{otherwise}.
    \end{cases}
\end{equation}
\medbreak

For instance, of $\delta$ is the ASM defined by
\begin{equation}
    \delta :=
    \Matrix{
        \Zero & \Plus & \Zero & \Zero & \Zero \\
        \Zero & \Zero & \Plus & \Zero & \Zero \\
        \Plus & \Minus & \Zero & \Zero & \Plus \\
        \Zero & \Plus & \Minus & \Plus & \Zero \\
        \Zero & \Zero & \Plus & \Zero & \Zero
    },
\end{equation}
we obtain
\begin{equation}
    M^\delta =
    \Matrix{
        \Zero & \One & \Zero & \Zero & \Zero \\
        \Zero & \Zero & \One & \Zero & \Zero \\
        \One & \One & \Zero & \Zero & \One \\
        \Zero & \One & \One & \One & \Zero \\
        \Zero & \Zero & \One & \Zero & \Zero
    }.
\end{equation}
\medbreak

It is immediate that $M^\delta$ is a $1$-packed matrix of the same size
as $\delta$. Besides, observe that since the $\Plus$ and the $\Minus$
alternate in an ASM, by starting from  a $1$-packed matrix $M$, there is
at most one ASM $\delta$ such that $M^\delta = M$.
\medbreak

\subsubsection{Hopf bialgebra structure on ASMs}
Let $\ASM$ be the vector space spanned by the set of all ASMs. For any
ASM $\delta$, let us denote by $\BasisF_\delta$ the element
$\BasisF_{M^\delta}$. Due to the above observation, the family
$\BasisF_\delta$, where $\delta$ are ASMs, spans $\ASM$. Moreover, since
the map $\BasisF_{\delta} \mapsto \BasisF_{M^\delta}$ is an injective
morphism from $\ASM$ to $\PM{1}$, this family forms a basis.
\medbreak

The product and the coproduct of $\PM{1}$ induce the product and the
coproduct of $\ASM$. For example, we have
\begin{equation} \label{equ:example_product_ASM}
    \BasisF_{
    \Matrix{
        \textcolor{Col1}{\Zero} & \textcolor{Col1}{\Plus}
            & \textcolor{Col1}{\Zero} \\
        \textcolor{Col1}{\Plus} & \textcolor{Col1}{\Minus}
            & \textcolor{Col1}{\Plus} \\
        \textcolor{Col1}{\Zero} & \textcolor{Col1}{\Plus}
            & \textcolor{Col1}{\Zero}
    }}
    \cdot
    \BasisF_{
    \Matrix{
        \textcolor{Col4}{\Plus}
    }}
    =
    \BasisF_{
    \Matrix{
        \textcolor{Col1}{\Zero} & \textcolor{Col1}{\Plus}
            & \textcolor{Col1}{\Zero} & \Zero \\
        \textcolor{Col1}{\Plus} & \textcolor{Col1}{\Minus}
            & \textcolor{Col1}{\Plus} & \Zero \\
        \textcolor{Col1}{\Zero} & \textcolor{Col1}{\Plus}
            & \textcolor{Col1}{\Zero} & \Zero \\
        \Zero & \Zero & \Zero & \textcolor{Col4}{\Plus}
    }}
    +
    \BasisF_{
    \Matrix{
        \textcolor{Col1}{\Zero} & \textcolor{Col1}{\Plus} & \Zero
            & \textcolor{Col1}{\Zero} \\
        \textcolor{Col1}{\Plus} & \textcolor{Col1}{\Minus} & \Zero
            & \textcolor{Col1}{\Plus} \\
        \textcolor{Col1}{\Zero} & \textcolor{Col1}{\Plus} & \Zero
            & \textcolor{Col1}{\Zero} \\
        \Zero & \Zero & \textcolor{Col4}{\Plus} & \Zero
    }}
    +
    \BasisF_{
    \Matrix{
        \textcolor{Col1}{\Zero} & \Zero & \textcolor{Col1}{\Plus}
            & \textcolor{Col1}{\Zero} \\
        \textcolor{Col1}{\Plus} & \Zero & \textcolor{Col1}{\Minus}
            & \textcolor{Col1}{\Plus} \\
        \textcolor{Col1}{\Zero} & \Zero & \textcolor{Col1}{\Plus}
            & \textcolor{Col1}{\Zero} \\
        \Zero & \textcolor{Col4}{\Plus} & \Zero & \Zero
    }}
    +
    \BasisF_{\Matrix{
        \Zero & \textcolor{Col1}{\Zero} & \textcolor{Col1}{\Plus}
            & \textcolor{Col1}{\Zero} \\
        \Zero & \textcolor{Col1}{\Plus} & \textcolor{Col1}{\Minus}
            & \textcolor{Col1}{\Plus} \\
        \Zero & \textcolor{Col1}{\Zero} & \textcolor{Col1}{\Plus}
            & \textcolor{Col1}{\Zero} \\
        \textcolor{Col4}{\Plus} & \Zero & \Zero & \Zero
    }},
\end{equation}
and
\begin{equation}
    \Coproduct \BasisF_{
    \Matrix{
        \Zero & \Plus & \Zero & \Zero \\
        \Zero & \Zero & \Zero & \Plus \\
        \Plus & \Minus & \Plus & \Zero \\
        \Zero & \Plus & \Zero & \Zero
    }}
    =
    \BasisF_{\emptyset}
    \otimes
    \BasisF_{
    \Matrix{
        \Zero & \Plus & \Zero & \Zero \\
        \Zero & \Zero & \Zero & \Plus \\
        \Plus & \Minus & \Plus & \Zero \\
        \Zero & \Plus & \Zero & \Zero
    }}
    +
    \BasisF_{
    \Matrix{
        \Zero & \Plus & \Zero \\
        \Plus & \Minus & \Plus \\
        \Zero & \Plus & \Zero
    }}
    \otimes
    \BasisF_{
    \Matrix{
        \Plus
    }}
    +
    \BasisF_{
    \Matrix{
        \Zero & \Plus & \Zero & \Zero \\
        \Zero & \Zero & \Zero & \Plus \\
        \Plus & \Minus & \Plus & \Zero \\
        \Zero & \Plus & \Zero & \Zero
    }}
    \otimes
    \BasisF_{\emptyset}.
\end{equation}
\medbreak

\begin{Theorem} \label{thm:ASM_Hopf}
    The vector space $\ASM$, endowed with the product and coproduct
    of $\PM{1}$, forms a free, cofree, and self-dual bigraded Hopf
    bialgebra which admits a bidendriform bialgebra structure.
\end{Theorem}
\medbreak

From now on, we shall see $\ASM$ as a simply graded Hopf bialgebra so
that the degree of any $\BasisF_{\delta}$, where $\delta$ is an ASM,
is the size of $\delta$. The dimensions of $\ASM$ form
Sequence~\OEIS{A005130} of~\cite{Slo} and the first few terms are
\begin{equation}
    1, 1, 2, 7, 42, 429, 7436, 218348, 10850216, 911835460,
    129534272700.
\end{equation}
\medbreak

By using same arguments as those used in
Section~\ref{subsec:freeness_PM}, one can build multiplicative bases
of $\ASM$ by setting, for any ASM $\delta$,
\begin{equation}
    \BasisE_{\delta}
        := \sum_{M^\delta \OrdPM M^{\delta'}} \BasisF_{\delta'}
    \qquad \mbox{and} \qquad
    \BasisH_{\delta}
        := \sum_{M^{\delta'} \OrdPM M^\delta} \BasisF_{\delta'}.
\end{equation}
This gives another way to prove the freeness of $\ASM$ by using same
arguments as those of Theorem~\ref{thm:freeness_PM}. Hence, $\ASM$ is
freely generated by the elements $\BasisE_{\delta}$ (resp.
$\BasisH_{\delta}$) where the $\delta$ are ASMs such that the
$M^\delta$ are connected (resp. anti-connected) $1$-packed matrices.
The first few numbers of algebraic generators of $\ASM$ are
\begin{equation}
    0, 1, 1, 4, 29, 343, 6536, 202890, 10403135, 889855638, 127697994191
\end{equation}
and the first few dimensions of totally primitive elements are
\begin{equation}
    0, 1, 0, 2, 20, 277, 5776, 188900, 9980698, 868571406, 125895356788.
\end{equation}
These are respectively Sequences~\OEIS{A231498} and~\OEIS{A231499}
of~\cite{Slo}.
\medbreak

Moreover, since the transpose of an ASM is also an ASM, by
Proposition~\ref{prop:self_duality_PM}, the map
$\phi : \ASM \to \ASM^\star$ linearly defined for any ASM $\delta$ by
\begin{equation}
    \phi\left(\BasisF_{\delta}\right)
    := \BasisF^\star_{{\delta}^\intercal}
\end{equation}
is an isomorphism.
\medbreak

\subsection{Algebraic interpretation of statistics on ASMs}
\label{subsec:ASM_stats_alg}
We provide algebraic interpretations of very common statistics on
ASMs, whose definitions are recalled in Section~\ref{subsec:ASMs} of
Chapter~\ref{chap:combinatorics}. These algebraic interpretations rely
on the Hopf bialgebra $\ASM$ and morphisms from $\ASM$ to $\K(q)$
(see Section~\ref{subsubsec:rational_functions} of
Chapter~\ref{chap:algebra} for notations about $q$-analogs of integers).
We also study here algebraic quotients of $\ASM$ defined by ideals
involving these statistics.
\medbreak

\subsubsection{Maps from $\ASM$ to $q$-rational functions}
The results presented here are consequences of the following two
combinatorial properties, highlighting compatibility between the
statistics $\NE$, $\SW$, $\SE$, $\NW$, $\OI$, and $\IO$ with the column
shifted shuffle product of ASMs.
\medbreak

\begin{Lemma} \label{lem:statistics_ASM_nonzero}
    Let $\delta$, $\delta_1$, and $\delta_2$ be three ASMs  such that
    $M^\delta \in M^{\delta_1} \cshuffle M^{\delta_2}$. Then, for any
    statistics $\Statistics$ of $\AllStatsB$,
    \begin{equation}
        \Statistics(\delta)
        = \Statistics(\delta_1) + \Statistics(\delta_2).
    \end{equation}
\end{Lemma}
\medbreak

\begin{Lemma} \label{lem:statistics_ASM_zero}
    Let $\delta$, $\delta_1$, and $\delta_2$ be three ASMs  such that
    $M^\delta \in M^{\delta_1} \cshuffle M^{\delta_2}$. Let $m$ be the
    size of $\delta_2$ (resp. $\delta_1$) and
    $\{k_1 < k_2 < \dots < k_m\}$ be the set of the indices of the
    columns of $M^\delta$ coming from $M^{\delta_2}$ (resp.
    $M^{\delta_1}$). Then, for any $\Statistics \in \{\NW, \SE\}$
    (resp. $\Statistics \in \{\SW, \NE\}$),
    \begin{equation}
        \Statistics(\delta) =
        \Statistics(\delta_1) + \Statistics(\delta_2) +
        \sum_{1 \leq j \leq m} (k_j-j).
    \end{equation}
\end{Lemma}
\medbreak

To illustrate Lemmas~\ref{lem:statistics_ASM_nonzero}
and~\ref{lem:statistics_ASM_zero}, we show here the
product~\eqref{equ:example_product_ASM} in $\ASM$, seen on six-vertex
configurations, where the vertices represented by squaresare of kind
$\IO$ while those represented by circles are of kind $\NW$:
\begin{equation}
    \BasisF_{
}
\end{equation}
\medbreak

\begin{Proposition} \label{prop:morphisms_ASM_q_rational}
    The maps $\phi_{s} : \ASM \to \K(q)$ and
    $\phi'_{s'} : \ASM \to \K(q)$ linearly defined, for any
    $s \in \AllStatsB$, $s' \in \AllStatsA$, and any ASM $\delta$ of
    size $n$ by
    \begin{equation}
        \phi_{s}\left(\BasisF_{\delta}\right)
            := \frac{q^{s(\delta)}}{n!}
        \qquad\qquad\mbox{and}\qquad\qquad
        \phi'_{s'}\left(\BasisF_{\delta}\right)
            := \frac{q^{s'(\delta)}}{(n)_q!}
    \end{equation}
    are associative algebra morphisms.
\end{Proposition}
\medbreak

This previous results remain valid in the dual $\ASM^\star$ of $\ASM$.
\medbreak

\begin{Proposition} \label{prop:morphisms_ASM_dual_q_rational}
    The maps $\psi_\Statistics : \ASM^\star \to \K(q)$ and
    $\psi'_{\Statistics'} : \ASM^\star \to \K(q)$ linearly defined,
    for any $\Statistics \in \AllStatsB$, $s' \in \AllStatsA$, and
    any ASM $\delta$ of size $n$ by
    \begin{equation}
        \psi_\Statistics\left(\BasisF^\star_{\delta}\right)
            := \frac{q^{\Statistics(\delta)}}{n!}
        \qquad\qquad\mbox{and}\qquad\qquad
        \psi'_{\Statistics'}\left(\BasisF^\star_{\delta}\right)
            := \frac{q^{\Statistics'(\delta)}}{(n)_q!}
    \end{equation}
    are associative algebra morphisms.
\end{Proposition}
\medbreak

\subsubsection{Equivalence relations on ASMs and associated subspaces
of $\ASM$}
Let $S\subseteq \AllStatsA \cup \AllStatsB$ be a set of statistics and
$\sim_S$ be the equivalence relation on the set of ASMs defined, for
any ASMs $\delta_1$ and $\delta_2$ of the same size, by
\begin{equation}
    \delta_1 \sim_S \delta_2
    \quad \mbox{if and only if} \quad
    s\left(\delta_1\right) = s\left(\delta_2\right)
    \mbox{ for all } s \in S.
\end{equation}
We denote by $\Ideal_S$ the associated vector space spanned by
\begin{equation}
    \left\{
    \BasisF_{\delta_1} - \BasisF_{\delta_2} :
    \delta_1 \sim_S \delta_2
    \right\}.
\end{equation}
\medbreak

\subsubsection{The algebra $\ASM/_{I_\IO}$}
Let us first study the statistics $\IO \in \AllStatsB$.
\medbreak

\begin{Proposition} \label{prop:quotient_ASM_io}
    The quotient $\QASM{\Ideal_\IO}$ is a commutative associative
    algebra.
\end{Proposition}
\medbreak

Note however that $\QASM{\Ideal_\IO}$ does not inherit the structure of
a coalgebra of $\ASM$ because even if
\begin{equation}
    x := \BasisF_{\Matrix{
      0 & \Plus & 0 & 0 \\
      \Plus & \Minus & \Plus & 0 \\
      0 & \Plus & 0 & 0 \\
      0 & 0 & 0 & \Plus
     }} -
     \BasisF_{\Matrix{
      0 & \Plus & 0 & 0 \\
      \Plus & \Minus & 0 & \Plus \\
      0 & \Plus & 0 & 0 \\
      0 & 0 & \Plus & 0
     }}
\end{equation}
is an element of $\Ideal_\IO$, the element
\begin{equation}
    \Coproduct(x) =
        1 \otimes x +
        \BasisF_{\Matrix{
            0 & \Plus & 0 \\
            \Plus & \Minus & \Plus \\
            0 & \Plus & 0
        }} \otimes
        \BasisF_{\Matrix{
            \Plus
        }} +
        x \otimes 1
\end{equation}
is not in $\ASM \otimes \Ideal_\IO + \Ideal_\IO \otimes \ASM$. Hence,
$\Ideal_\IO$ is not a coideal.
\medbreak

\begin{Proposition} \label{prop:dimensions_quotient_ASM_io}
    For any $n \geq 0$, the dimension of the $n$th graded component
    of $\QASM{\Ideal_\IO}$ satisfies
    \begin{equation}
        \dim \QASM{\Ideal_\IO}(n)
        =
        \left\lfloor\frac{n^2}{4}\right\rfloor + 1.
    \end{equation}
\end{Proposition}
\medbreak

The dimensions of $\QASM{\Ideal_\IO}$ form Sequence \OEIS{A033638}
of \cite{Slo} and the first few terms are
\begin{equation}
    1, 1, 1, 2, 3, 5, 7, 10, 13, 17, 21.
\end{equation}
\medbreak

A basic argument on generating series implies that these dimensions
cannot be the ones of a free commutative algebra and hence,
$\QASM{\Ideal_\IO}$ is not free as a commutative associative algebra.
\medbreak

Using the symmetry between the statistics $\IO$ and $\OI$ provided by
Proposition \ref{prop:symmetries_statistics_ASMs} of
Section~\ref{subsec:ASMs} of Chapter~\ref{chap:combinatorics}, we
immediately have $\sim_{\OI} = \sim_{\IO}$ and then,
$\QASM{\Ideal_\OI} = \QASM{\Ideal_\IO}$.
\medbreak

\subsubsection{The algebra $\ASM/_{I_\NW}$}
Let us now study the statistics $\NW \in \AllStatsA$.
\medbreak

\begin{Proposition} \label{prop:quotient_ASM_nw}
    The quotient $\QASM{\Ideal_\NW}$ is a commutative associative
    algebra.
\end{Proposition}
\medbreak

Note however that $\QASM{\Ideal_\NW}$ does not inherit the structure of
a coalgebra of $\ASM$ because even if
\begin{equation}
     x := \BasisF_{\Matrix{
          0 & 0 & 0 & \Plus\\
          \Plus & 0 & 0 & 0\\
          0 & 0 & \Plus & 0\\
          0 & \Plus & 0 & 0
     }} -
     \BasisF_{\Matrix{
          0 & 0 & \Plus & 0\\
          0 & \Plus & 0 & 0\\
          \Plus & 0 & \Minus & \Plus\\
          0 & 0 & \Plus & 0
     }}
\end{equation}
is an element of $\Ideal_\NW$, the element
\begin{equation}
    \Coproduct(x) =
        1 \otimes x +
        \BasisF_{\Matrix{
          \Plus
        }} \otimes
        \BasisF_{\Matrix{
            0 & 0 & \Plus\\
            0 & \Plus & 0\\
            \Plus & 0 & 0
        }} +
        \BasisF_{\Matrix{
            \Plus & 0\\
            0 & \Plus
        }} \otimes
        \BasisF_{\Matrix{
            0 & \Plus\\
            \Plus & 0
        }} +
        \BasisF_{\Matrix{
            \Plus & 0 & 0\\
            0 & 0 & \Plus\\
            0 & \Plus & 0
        }} \otimes
        \BasisF_{\Matrix{
            \Plus
        }} +
        x \otimes 1
\end{equation}
is not in $\ASM \otimes \Ideal_\NW + \Ideal_\NW \otimes \ASM$. Hence,
$\Ideal_\NW$ is not a coideal.
\medbreak

\begin{Proposition} \label{prop:dimensions_quotient_ASM_nw}
    For any $n \geq 0$, the dimension of the $n$th graded component of
    $\QASM{\Ideal_\NW}$ satisfies
    \begin{equation}
        \dim \QASM{\Ideal_\NW}(n) = \binom{n}{2} + 1.
    \end{equation}
\end{Proposition}
\medbreak

The dimensions of $\QASM{\Ideal_\NW}$ form Sequence \OEIS{A152947}
of \cite{Slo} and the first few terms are
\begin{equation}
    1, 1, 2, 4, 7, 11, 16, 22, 29, 37, 46, 56.
\end{equation}
\medbreak

A basic argument on generating series implies that these dimensions
cannot be the ones of a free commutative algebra and hence,
$\QASM{\Ideal_\NW}$ is not free as a commutative associative algebra.
\medbreak

Using the symmetry between the statistics $\NW$ and $\SE$ provided by
Proposition~\ref{prop:symmetries_statistics_ASMs}, we immediately have
$\sim_{\SE} = \sim_{\NW}$ and then,
$\QASM{\Ideal_\SE} = \QASM{\Ideal_\NW}$. Moreover, by using the same
arguments as before, $\QASM{\Ideal_\SW}$ and $\QASM{\Ideal_\NE}$ are
the same commutative algebras.
\medbreak

Note that the map $\theta : \QASM{\Ideal_\NW} \to \QASM{\Ideal_\SW}$
linearly defined for any ASM $\delta$ by
\begin{equation}
    \theta\left(\pi_{\NW}\left(\BasisF_{\delta}\right)\right) :=
    \pi_{\SW}\left(\BasisF_{\overleftarrow{\delta}}\right),
\end{equation}
where $\pi_{\NW}$ (resp. $\pi_{\SW}$) is the canonical projection from
$\ASM$ to $\QASM{\Ideal_\NW}$ (resp. $\QASM{\Ideal_\SW}$) and
$\overleftarrow{\delta}$ is the ASM where, for any $i \in [n]$, the
$i$th column of $\overleftarrow{\delta}$ is the $(n - i + 1)$st column
of $\delta$, is an isomorphism between $\QASM{\Ideal_\NW}$
and~$\QASM{\Ideal_\SW}$.
\medbreak

\subsubsection{The algebra $\ASM/_{I_{\IO, \NW}}$}
Let us finally study the set of statistics $\{\IO, \NW\}$.
\medbreak

\begin{Proposition} \label{prop:quotient_ASM_nw_io}
    The quotient $\QASM{\Ideal_{\IO,\NW}}$ is a commutative associative
    algebra.
\end{Proposition}
\medbreak

Note however that $\QASM{\Ideal_{\IO,\NW}}$ does not inherit the
structure of a coalgebra of $\ASM$ because even if
\begin{equation}
     x :=
     \BasisF_{\Matrix{
          0 & \Plus & 0 & 0\\
          \Plus & \Minus & \Plus & 0\\
          0 & 0 & 0 & \Plus\\
          0 & \Plus & 0 & 0
     }} -
     \BasisF_{\Matrix{
          0 & \Plus & 0 & 0\\
          0 & 0 & \Plus & 0\\
          \Plus & \Minus & 0 & \Plus\\
          0 & \Plus & 0 & 0
     }}
\end{equation}
is an element of $\Ideal_{\IO,\NW}$, the element
\begin{equation}
    \Coproduct(x) =
        1 \otimes x +
        \BasisF_{\Matrix{
            0 & \Plus & 0 \\
            \Plus & \Minus & \Plus \\
            0 & \Plus & 0
        }} \otimes
        \BasisF_{\Matrix{
            \Plus
        }} +
        x \otimes 1
\end{equation}
is not in
\begin{math}
    \ASM \otimes \Ideal_{\IO,\NW} + \Ideal_{\IO,\NW} \otimes \ASM
\end{math}.
Hence, $\Ideal_{\IO,\NW}$ is not a coideal.
\medbreak

By computer exploration, the first few dimensions of
$\QASM{\Ideal_{\IO,\NW}}$ are
\begin{equation}
    1, 1, 2, 5, 13, 31, 66, 127, 225,
\end{equation}
and seems to be Sequence~\OEIS{A116701} of~\cite{Slo}.
\medbreak

A basic argument on generating series implies that these dimensions
cannot be the ones of a free commutative associative algebra and
hence, $\QASM{\Ideal_{\IO,\NW}}$ is not free as a commutative algebra.
\medbreak

\subsubsection{Others quotients of $\ASM$}
Using the symmetries provided by
Proposition~\ref{prop:symmetries_statistics_ASMs},
all the algebras $\QASM{\Ideal_S}$, where $S$ contains two nonsymmetric
statistics, are equal to $\QASM{\Ideal_{\IO,\NW}}$. Moreover, note that
by using the same arguments as before, one can prove that for any
$S \in \AllStatsA \cup \AllStatsB$, $\QASM{\Ideal_S}$ is a commutative
algebra isomorphic to $\QASM{\Ideal_\IO}$, $\QASM{\Ideal_\NW}$,
or~$\QASM{\Ideal_{\IO,\NW}}$.
\medbreak

\section*{Concluding remarks}
The work presented in this chapter contributes to enrich the already
large collection of combinatorial Hopf bialgebras related to $\FQSym$.
Our main contributions are the Hopf bialgebra $\PM{k}$ of $k$-packed
matrices and the Hopf bialgebra $\ASM$ of alternating sign matrices.
\smallbreak

Naturally, our results raise several questions for further research.
First, one can ask for the enumeration of equivalence classes of
$k$-packed matrices for the equivalence relations considered in
Section~\ref{subsubsec:experiments_congruences_PM}. Second, we have
described an injective associative algebra morphism from $\PML{1}^\star$
to $\MQSym$ (see Proposition~\ref{prop:morphism_PM_MQSym}).
Nevertheless, as observed, this morphism is not a Hopf bialgebra
morphism. Then, the question to define a Hopf embedding of
$\PML{1}^\star$ into $\MQSym$ is open. Let us address a last research
direction. Most Hopf bialgebras related to $\FQSym$ have polynomial
realizations, that is, a way to encode their elements as
polynomials~\cite{DHT02}, compatible with the product and the alphabet
doubling (see for instance~\cite{Hiv03}). The question to provide such a
polynomial realization of $\PM{k}$ seems worth studying.
\medbreak


\chapter{From pros to Hopf bialgebras} \label{chap:pros_bialgebras}
The content of this chapter comes from~\cite{BG16} and is a joint work
with Jean-Paul Bultel.
\medbreak

\section*{Introduction}
The theory of operads and the one of Hopf bialgebras have several known
interactions. One of these is a construction~\cite{Vdl04} taking an
operad $\Oca$ as input and producing a Hopf bialgebra $\OvH(\Oca)$ as
output, which is called the natural Hopf bialgebra of $\Oca$. This
construction has been studied in some recent works: in~\cite{CL07}, it
is shown that $\OvH$ can be rephrased in terms of an incidence Hopf
bialgebra of a certain family of posets, and in~\cite{ML14}, a general
formula for its antipode is established. Let us also cite~\cite{Fra08}
in which this construction is considered to study series of trees. The
initial motivation of the work contained in this chapter was to
generalize this $\OvH$ construction with the aim of constructing some
new and interesting Hopf bialgebras. The direction we have chosen is to
start with pros (see~\cite{McL65,Lei04,Mar08}), algebraic structures
which generalize operads in the sense that pros deal with operators with
possibly several outputs (see Section~\ref{subsec:pros} of
Chapter~\ref{chap:algebra}).
\smallbreak

Our main contribution consists in the definition of a new construction
$\PvH$ from pros to bialgebras. Roughly speaking, the construction
$\PvH$ can be described as follows. Given a pro $\Pca$ satisfying some
mild properties, the Hopf bialgebra $\PvH(\Pca)$ has bases indexed by a
particular subset of elements of $\Pca$. The product of $\PvH(\Pca)$ is
the horizontal composition of $\Pca$ and the coproduct of $\PvH(\Pca)$
is defined from the vertical composition of $\Pca$, enabling to separate
a basis element into two smaller parts. The properties satisfied by
$\Pca$ imply, in a nontrivial way, that the product and the coproduct of
$\PvH(\Pca)$ satisfy the required axioms to be a bialgebra. This
construction generalizes $\OvH$ and establishes a new connection between
the theory of pros and the theory of Hopf bialgebras.
\smallbreak

Let us provide some details about our construction $\PvH$. A first
version of this construction is presented, associating a Hopf bialgebra
$\PvH(\Pca)$ with a free pro $\Pca$. The fundamental basis of this Hopf
bialgebra is a set-basis with respect to the product, and the structure
coefficients of the coproduct are nontrivial ({\em i.e.}, they are
possibly different from $0$ and $1$). As an associative algebra,
$\PvH(\Pca)$ is always free. This construction is extended to a class of
non-necessarily free pros. The pros of this class, called stiff pros,
can be described by particular quotients of free pros. These pros arise
somewhat naturally because, under some mild conditions, two well-known
constructions of pros~\cite{Mar08} produce stiff pros. The first one,
$\OvP$, takes as input operads and the second one, $\MvP$, takes as
input monoids. The construction $\OvP$ is used to show that the natural
Hopf bialgebra of an operad can be reformulated as a particular case of
our construction $\PvH$. The Hopf bialgebras that one can construct from
$\PvH$ are very similar to the Connes-Kreimer Hopf
bialgebra~\cite{CK98} in the sense that their coproduct can be computed
by means of admissible cuts in various combinatorial objects. From very
simple stiff pros, it is possible to reconstruct the Hopf bialgebra of
noncommutative symmetric functions $\Sym$~\cite{GKLLRT94} and the
noncommutative Fàa di Bruno Hopf bialgebra $\FdBNC$~\cite{BFK06}.
Besides, we present a way of using $\PvH$ to reconstruct some of the
Hopf bialgebras $\FdBNC_\gamma$, a $\gamma$-deformation of $\FdBNC$
introduced by Foissy~\cite{Foi08}.
\medbreak

This chapter is organized as follows. In
Section~\ref{sec:natural_Hopf_bialgebra_operad}, we recall the natural
Hopf bialgebra construction $\OvH$ of an operad and some background
about the noncommutative Faà di Bruno Hopf bialgebra $\FdBNC$ and its
commutative version $\FdB$. We provide in
Section~\ref{sec:construction_H} the description of our new construction
$\PvH$ and study some of its algebraic and combinatorial properties. We
conclude this chapter by giving some examples of applications of $\PvH$
in Section~\ref{sec:construction_H_examples} from very simple pros. We
hence obtain several Hopf bialgebras, which, respectively, involve
forests of planar rooted trees, some kinds of graphs consisting of nodes
with one parent and several children or several parents and one child
that we call forests of bitrees, heaps of pieces (see~\cite{Vie86} for a
general presentation of these combinatorial objects), and a particular
class of heaps of pieces that we call heaps of friable pieces. All these
Hopf bialgebras depend on a nonnegative integer as parameter~$\gamma$.
\medbreak

\subsubsection*{Note}
This chapter deals only with ns set-operads. For this reason, ``operad''
means ``ns set-operad'' in this chapter. Similarly, ``pro'' means
``set-pro''. Moreover, all the free pros appearing here have generators
with at least one input and one output.
\medbreak

\section{Hopf bialgebras and the natural Hopf bialgebra of an operad}
\label{sec:natural_Hopf_bialgebra_operad}
We recall in this section a  construction associating a combinatorial
Hopf bialgebra with an operad. This construction can be used to define
the Faà di Bruno Hopf bialgebra.
\medbreak

\subsection{Combinatorial Hopf algebras}
Let us start by recalling some definitions and properties about the Faà
di Bruno Hopf bialgebra, the Hopf bialgebra of symmetric functions, and
some of its noncommutative analogs. Here, we assume that the ground
field $\K$ on which all the Hopf bialgebras are defined contains $\R$.
\medbreak

\subsubsection{Faà di Bruno Hopf bialgebra and its deformations}
Let $\FdB$ be the free commutative algebra generated by elements $h_n$,
$n \geq 1$, with $\deg(h_n) = n$. The bases of $\FdB$ are thus indexed
by integer partitions, and the unit is denoted by $h_0$. Alternatively,
$\FdB := \K \Angle{\IntPart}$, where $\IntPart$ is the graded
combinatorial collection of integer partitions defined in
Section~\ref{subsubsec:integer_partitions} of
Chapter~\ref{chap:combinatorics}. This is the \Def{algebra of symmetric
functions}~\cite{Mac15}. There are several ways to endow $\FdB$ with a
coproduct to turn it into a Hopf bialgebra. In~\cite{Foi08}, Foissy
obtains, as a byproduct of his investigation of combinatorial
Dyson-Schwinger equations in the Connes-Kreimer algebra, a one-parameter
family $\Delta_\gamma$, $\gamma \in \R$, of coproducts on $\FdB$,
defined by using alphabet transformations (see~\cite{Mac15}), by
\begin{equation} \label{equ:coproduit_delta_gamma}
    \Delta_\gamma(h_n) :=
    \sum_{0 \leq k \leq n}
    h_k \otimes h_{n - k}((k \gamma + 1) X),
\end{equation}
where, for any $\alpha \in \R$ and $n \in \N$, $h_n(\alpha X)$ is the
coefficient of $t^n$ in the series
\begin{math}
    \left(\sum_{k \geq 0} h_k t^k\right)^\alpha
\end{math}.
In particular,
\begin{equation}
    \Delta_0(h_n) =
    \sum_{0 \leq k \leq n} h_k \otimes h_{n - k}.
\end{equation}
The algebra $\FdB$ with the coproduct $\Delta_0$ is the classical
\Def{Hopf bialgebra of symmetric functions} $\SymCom$~\cite{Mac15}.
Moreover, for all $\gamma \neq 0$, all $\FdB_\gamma$ are isomorphic to
$\FdB_1$, which is known as the \Def{Faà di Bruno bialgebra}~\cite{JR79}.
The coproduct $\Delta_0$ comes from the interpretation of $\FdB$ as the
algebra of polynomial functions on the multiplicative group
\begin{equation}
    G(t) :=
    \left\{1 + \sum_{k \geq 1} a_k t^k : a_k \in \R, k \geq 1\right\}
\end{equation}
of formal power series of constant term $1$, and $\Delta_1$ comes from
its interpretation as the algebra of polynomial functions on the group
$tG(t)$ for the series composition of formal diffeomorphisms of the real
line.
\medbreak

\subsubsection{Noncommutative analogs}
Formal power series in one variable with coefficients in a
noncommutative algebra can be composed (by substitution of the
variable). This operation is not associative, so that they do not form a
group. For example, when $a$ and $b$ belong to a noncommutative algebra,
one has
\begin{equation}
    (t^2 \circ a t) \circ b t = a^2 t^2 \circ b t = a^2 b^2 t^2
\end{equation}
but
\begin{equation}
    t^2 \circ (a t \circ b t) = t^2 \circ a b t = a b a b t^2.
\end{equation}
However, the analogue of the Faà di Bruno Hopf bialgebra still exists in
this noncommutative context and is known as the \Def{noncommutative Faà
di Bruno Hopf bialgebra}. It is investigated in~\cite{BFK06} in view of
applications in quantum field theory. In~\cite{Foi08}, Foissy also
obtains an analogue of the family $\FdB_\gamma$ in this context. Indeed,
considering noncommutative generators $\BasisS_n$ (with
$\deg(\BasisS_n) = n$) instead of the $h_n$, for all $n \geq 1$, leads
to a free noncommutative algebra $\FdBNC$ whose bases are indexed by
integer compositions. This is the \Def{algebra of noncommutative
symmetric functions}~\cite{GKLLRT94}. The addition of the coproduct
$\Delta_\gamma$ defined by
\begin{equation}\label{equ:coproduit_delta_gamma_non_commutatif}
    \Delta_\gamma(\BasisS_n) :=
    \sum_{0 \leq k \leq n}
    \BasisS_k \otimes \BasisS_{n - k}((k \gamma + 1) A),
\end{equation}
where, for any $\alpha\in\R$ and $n\in\N$, $\BasisS_n(\alpha A)$ is the
coefficient of $t^n$ in
\begin{math}
    \left(\sum_{k\geq 0} \BasisS_kt^k\right)^\alpha
\end{math},
forms a noncommutative Hopf bialgebra $\FdBNC_\gamma$. In particular,
\begin{equation}
     \Delta_0(\BasisS_n)=\sum_{k=0}^n \BasisS_k\otimes\BasisS_{n-k},
\end{equation}
where $\BasisS_0$ is the unit. In this way, $\FdBNC$ with the coproduct
$\Delta_0$ is the \Def{Hopf bialgebra of noncommutative symmetric
functions} $\Sym$ \cite{GKLLRT94,KLT97}, and for all $\gamma \neq 0$,
all the $\FdBNC_\gamma$ are isomorphic to $\FdBNC_1$, which is the
\Def{noncommutative Faà di Bruno Hopf bialgebra}.
\medbreak

\subsection{The natural Hopf bialgebra of an operad}
We describe here a construction associating a Hopf bialgebra with an
operad (under some conditions). We then apply this construction to
obtain $\FdB$ from the associative operad.
\medbreak

\subsubsection{The construction}
A slightly different version of the construction we shall present here
is considered in~\cite{Vdl04,CL07,ML14}. Let $\Oca$ be an operad and
denote by $\Oca^+$ the set $\Oca \setminus \{\Unit\}$. The \Def{natural
Hopf bialgebra} of $\Oca$ is the free commutative algebra  $\OvH(\Oca)$
spanned by the $\Tit_x$, where the $x$ are elements of $\Oca^+$. The
bases of $\OvH(\Oca)$ are thus indexed by finite multisets of elements
of $\Oca^+$. Alternatively, $\OvH(\Oca) = \K \Angle{\Multiset(\Oca^+)}$,
where $\Multiset$ is the multiset operation over graded collections
(see Section~\ref{subsubsec:operations_graded_comb_collections} of
Chapter~\ref{chap:combinatorics}). The unit of $\OvH(\Oca)$ is denoted
by $\Tit_{\Unit}$ and the coproduct of $\OvH(\Oca)$ is the unique
associative algebra morphism satisfying, for any element $x$ of
$\Oca^+$,
\begin{equation} \label{equ:coproduct_natural_Hopf_bialgebra_operad}
    \Delta(\Tit_x) =
    \sum_{\substack{y, z_1, \dots, z_\ell \in \Oca \\
    y \circ [z_1, \dots, z_\ell] = x}}
    \Tit_y \otimes \Tit_{z_1} \dots \Tit_{z_\ell}.
\end{equation}
\medbreak

The Hopf bialgebra $\OvH(\Oca)$ can be graded by
$\deg(\Tit_x) := |x| - 1$. Note that with this grading, when
$\Oca(1) = \{\Unit\}$ and when the $\Oca(n)$ are finite for all
$n \geq 1$, $\OvH(\Oca)$ becomes a combinatorial Hopf bialgebra.
\medbreak

\subsubsection{The natural Hopf bialgebra of the associative operad}
Let us consider the associative operad $\As$ (see
Section~\ref{subsubsec:associative_operad} of
Chapter~\ref{chap:algebra}). The set $\As^+$ consists in the elements
$\Afr_n$ with $n \geq 2$. The Hopf bialgebra $\OvH(\As)$ is the linear
span of the elements $\Tit_{\lbag x_1, \dots, x_\ell\rbag}$ where
$x_i \in \As^+$, $i \in [\ell]$. Any multiset
$X := \lbag \Afr_{n_1}, \dots, \Afr_{n_k}\rbag$ of $\As^+$ can be
encoded by a nondecreasing word $u_1^{\alpha_1} \dots u_r^{\alpha_r}$
where $\alpha_i$ is the multiplicity of $\Afr_{i + 1}$ in $X$ for any
$i \in [r]$. For instance, the basis element
$\Tit_{\lbag \Afr_2, \Afr_2, \Afr_4, \Afr_6, \Afr_6, \Afr_7\rbag}$ is
encoded by $\Tit_{655311}$. Moreover the degrees of such basis elements
indexed by words are the sums of their letters. For this reason, the
basis elements of $\OvH(\As)$ are indexed by integer partitions.
Besides, here is an example of a coproduct in $\OvH(\As)$
using~\eqref{equ:coproduct_natural_Hopf_bialgebra_operad}:
\begin{equation}\begin{split}
    \label{equ:example_coproduct_faa_di_bruno}
    \Coproduct(\Tit_2) & =
    \Tit_{\Unit} \otimes \Tit_2 +
    \Tit_1 \otimes (\Tit_\Unit \Tit_1 + \Tit_1 \Tit_\Unit)
    + \Tit_2 \otimes \Tit_\Unit \Tit_\Unit \Tit_\Unit \\
    & =
    \Tit_{\Unit} \otimes \Tit_2 +
    2 \Tit_1 \otimes \Tit_1
    + \Tit_2 \otimes \Tit_{\Unit}.
\end{split}\end{equation}
For instance, the coefficient of $\Tit_1 \otimes \Tit_1$ in
$\Coproduct(\Tit_2)$ is $2$ because there are two ways to factorize
$\Afr_3$ in $\As$ by using the complete composition map where the first
operand is $\Afr_2$:
$\Afr_3 = \Afr_2 \circ [\Afr_1, \Afr_2]$ and
$\Afr_3 = \Afr_2 \circ [\Afr_2, \Afr_1]$. It is known (see for
instance~\cite{ML14}) that $\OvH(\As)$ is isomorphic to~$\FdB_1$.
\medbreak

\section{From pros to combinatorial Hopf algebras}
\label{sec:construction_H}
We introduce in this section the main construction of this work and
review some of its properties. In all this section, $\Pca$ is a free pro
generated by a bigraded set $\GeneratingSet$. We recall that we work
only with generating sets satisfying $\GeneratingSet(p, q) = \emptyset$
when $p = 0$ or $q = 0$. Starting with $\Pca$, our construction produces
a bialgebra $\PvH(\Pca)$ whose bases are indexed by the reduced elements
of $\Pca$. We shall also extend this construction over a class of non
necessarily free pros.
\medbreak

\subsection{The Hopf bialgebra of a free pro}
\label{subsec:free_pro_to_Hopf_bialgebra}
We shall use from now on the notions about prographs introduced in
Section~\ref{subsubsec:prographs} of Chapter~\ref{chap:combinatorics}
and the notions about free pros contained in
Section~\ref{subsubsec:free_pros} of Chapter~\ref{chap:algebra}.
\medbreak

The bases of the vector space
\begin{equation}
    \PvH(\Pca) := \K \Angle{\Reduced(\Pca)}
\end{equation}
are indexed by the reduced elements of $\Pca$. The elements $\BasisS_x$,
$x \in \Reduced(\Pca)$, form thus a basis of $\PvH(\Pca)$, called
\Def{fundamental basis}. We endow $\PvH(\Pca)$ with a product
$\cdot : \PvH(\Pca) \otimes \PvH(\Pca) \to \PvH(\Pca)$ linearly defined,
for any reduced elements $x$ and $y$ of $\Pca$, by
\begin{equation}
    \BasisS_x \cdot \BasisS_y := \BasisS_{x * y},
\end{equation}
and with a coproduct
$\Delta : \PvH(\Pca) \to \PvH(\Pca) \otimes \PvH(\Pca)$ linearly defined,
for any reduced elements $x$ of $\Pca$, by
\begin{equation}
    \Delta\left(\BasisS_x\right) :=
    \sum_{\substack{y, z \in \Pca \\ y \circ z = x}}
    \BasisS_{\Reduced(y)} \otimes \BasisS_{\Reduced(z)}.
\end{equation}
\medbreak

Throughout this section, we shall consider some examples involving
the free pro generated by
\begin{math}
    \GeneratingSet := \GeneratingSet(2, 2) \sqcup \GeneratingSet(3, 1)
\end{math}
where
$\GeneratingSet(2, 2) := \{\Asf\}$ and
$\GeneratingSet(3, 1) := \{\Bsf\}$, denoted by $\AB$. For instance, we
have in $\PvH(\AB)$
\begin{equation}
    \BasisS_{
}
    \otimes
    \BasisS_{\Unit_0}\,.
\end{multline}
\medbreak

As a consequence of Lemma~\ref{lem:square_rule_free_pros}
of Chapter~\ref{chap:algebra}, the coproduct $\Delta$ of $\PvH(\Pca)$
is coassociative. Moreover, this lemma implies that $\Delta$ is a
morphism of associative algebras. Hence, we obtain the following
result.
\medbreak

\begin{Theorem} \label{thm:free_pro_to_Hopf_bialgebra}
    Let $\Pca$ be a free pro. Then, $\PvH(\Pca)$ is a Hopf bialgebra.
\end{Theorem}
\medbreak

\subsection{Properties of the construction}
Let us now study the general properties of the Hopf bialgebras obtained
by the construction~$\PvH$.
\medbreak

\subsubsection{Algebraic generators and freeness}
\begin{Proposition} \label{prop:free_pro_to_Hopf_bialgebra_freeness}
    Let $\Pca$ be a free pro. Then, $\PvH(\Pca)$ is freely generated as
    an associative algebra by the set of all $\BasisS_g$, where the $g$
    are indecomposable and reduced elements of $\Pca$.
\end{Proposition}
\medbreak

\subsubsection{Gradings}
There are several ways to define gradings for $\PvH(\Pca)$ to turn it
into a combinatorial Hopf bialgebra. For this purpose, we say that a map
$\omega : \Reduced(\Pca) \to \N$ is a \Def{grading} of $\Pca$ if it
satisfies the following four properties:
\begin{enumerate}[label = ({\it G{\arabic*})}]
    \item \label{item:property_grading_1}
    for any reduced elements $x$ and $y$ of $\Pca$,
    $\omega(x * y) = \omega(x) + \omega(y)$;
    \item \label{item:property_grading_2}
    for any reduced elements $x$ of $\Pca$ satisfying $x = y \circ z$
    where $y, z \in \Pca$,
    $\omega(x) = \omega(\Reduced(y)) + \omega(\Reduced(z))$;
    \item \label{item:property_grading_3}
    for any $n \geq 0$, the fiber $\omega^{-1}(n)$ is finite;
    \item \label{item:propriete_bonne_graduation_4}
    $\omega^{-1}(0) = \{\Unit_0\}$.
\end{enumerate}
\medbreak

A very generic way to endow $\Pca$ with a grading consists in
providing a map $\omega : \GeneratingSet \to \N \setminus \{0\}$
associating a positive integer with any generator of $\Pca$, namely its
\Def{weight}; the degree $\omega(x)$ of any element $x$ of $\Pca$ being
the sum of the weights of the occurrences of the generators used to
build $x$. For instance, the map $\omega$ defined by $\omega(\Asf) := 3$
and $\omega(\Bsf) := 2$ is a grading of $\AB$ and we have
\begin{equation}
    \omega\left(
    \begin{tikzpicture}[xscale=.3,yscale=.2,Centering]
        \node[Leaf](S1)at(0,0){};
        \node[Leaf](S2)at(2,0){};
        \node[Leaf](S3)at(5,0){};
        \node[Leaf](S4)at(7,0){};
        \node[Operator](N1)at(1,-2){\begin{math}\Asf\end{math}};
        \node[OperatorColorC](N2)at(3,-4){\begin{math}\Bsf\end{math}};
        \node[Operator](N3)at(6,-3){\begin{math}\Asf\end{math}};
        \node[Leaf](E1)at(0,-6){};
        \node[Leaf](E2)at(2,-6){};
        \node[Leaf](E3)at(3,-6){};
        \node[Leaf](E4)at(4,-6){};
        \node[Leaf](E5)at(5,-6){};
        \node[Leaf](E6)at(7,-6){};
        \draw[Edge](N1)--(S1);
        \draw[Edge](N1)--(S2);
        \draw[Edge](N3)--(S3);
        \draw[Edge](N3)--(S4);
        \draw[Edge](N1)--(E1);
        \draw[Edge](N1)--(N2);
        \draw[Edge](N2)--(E2);
        \draw[Edge](N2)--(E3);
        \draw[Edge](N2)--(E4);
        \draw[Edge](N3)--(E5);
        \draw[Edge](N3)--(E6);
    \end{tikzpicture}
    \right)
    = 8.
\end{equation}
\medbreak

\begin{Proposition} \label{prop:free_pro_to_Hopf_bialgebra_graduation}
    Let $\Pca$ be a free pro and $\omega$ be a grading of $\Pca$.
    Then, with the grading
    \begin{equation} \label{equ:graduation}
        \PvH(\Pca) =
        \bigoplus_{n \geq 0}
        \K \Angle{\left\{
            \BasisS_x : x \in \Reduced(\Pca) \mbox{ and } \omega(x) = n
        \right\}},
    \end{equation}
    $\PvH(\Pca)$ is a combinatorial Hopf bialgebra.
\end{Proposition}
\medbreak

\subsubsection{Antipode}
Since the antipode of a combinatorial Hopf bialgebra can be computed by
induction on the degrees, we obtain an expression for the one of
$\PvH(\Pca)$ when $\Pca$ admits a grading. This expression is an
instance of the Takeuchi formula~\cite{Tak71} and is particularly simple
since the product of $\PvH(\Pca)$ is multiplicative.
\medbreak

\begin{Proposition} \label{prop:free_pro_to_Hopf_bialgebra_antipode}
    Let $\Pca$ be a free pro admitting a grading. For any reduced
    element $x$ of $\Pca$ different from $\Unit_0$, the antipode $\nu$
    of $\PvH(\Pca)$ satisfies
    \begin{equation} \label{equ:pro_vers_AHC_antipode}
        \nu(\BasisS_x) =
        \sum_{\substack{x_1, \dots, x_\ell \in \Pca, \ell \geq 1 \\
        x_1 \circ \dots \circ x_\ell = x \\
        \Reduced(x_i) \ne \Unit_0, i \in [\ell]}}
        (-1)^\ell \;
        \BasisS_{\Reduced(x_1 * \dots * x_\ell)}.
    \end{equation}
\end{Proposition}
\medbreak

We have for instance in $\PvH(\AB)$,
\begin{multline}
    \nu\,
    \BasisS_{
}\,.
\end{multline}
\medbreak

\subsubsection{Duality}
When $\Pca$ admits a grading, let us denote by $\PvH(\Pca)^\star$ the
graded dual of $\PvH(\Pca)$. By definition, the dual basis of the
fundamental basis of $\PvH(\Pca)$ consists in the elements
$\BasisS^\star_x$, $x \in \Reduced(\Pca)$.
\medbreak

\begin{Proposition} \label{prop:free_pro_to_Hopf_bialgebra_dual}
    Let $\Pca$ be a free pro admitting a grading. Then, for any
    reduced elements $x$ and $y$ of $\Pca$, the product and the
    coproduct of $\PvH(\Pca)^\star$ satisfy
    \begin{equation} \label{equ:pro_vers_AHC_produit_dual}
        \BasisS^\star_x \cdot \BasisS^\star_y =
        \sum_{\substack{x', y' \in \Pca \\
        x' \circ y' \in \Reduced(\Pca) \\
        \Reduced(x') = x, \Reduced(y') = y}}
        \BasisS^\star_{x' \circ y'}
    \end{equation}
    and
    \begin{equation} \label{equ:pro_vers_AHC_coproduit_dual}
        \Delta\left(\BasisS^\star_x\right) =
        \sum_{\substack{y, z \in \Pca \\ y * z = x}}
        \BasisS^\star_y \otimes \BasisS^\star_z.
    \end{equation}
\end{Proposition}
\medbreak

For instance, we have in $\PvH(\AB)$
\begin{multline}
    \BasisS^\star_{
}
    \otimes
    \BasisS^\star_{\Unit_0}.
\end{multline}
\medbreak

\subsubsection{Quotient bialgebras}

\begin{Proposition}
\label{prop:free_pro_to_Hopf_bialgebra_subset_generators}
    Let $\GeneratingSet$ and $\GeneratingSet'$ be two bigraded sets such
    that $\GeneratingSet' \subseteq \GeneratingSet$. Then, the map
    \begin{math}
        \phi : \PvH(\FreePro(\GeneratingSet))
        \to \PvH\left(\FreePro\left(\GeneratingSet'\right)\right)
    \end{math}
    linearly defined, for any reduced element $x$ of
    $\FreePro(\GeneratingSet)$, by
    \begin{equation}
        \phi(\BasisS_x) :=
        \begin{cases}
            \BasisS_x & \mbox{if }
                x \in \FreePro\left(\GeneratingSet'\right), \\
            0 & \mbox{otherwise},
        \end{cases}
    \end{equation}
    is a surjective bialgebra morphism. Moreover,
    $\PvH\left(\FreePro\left(\GeneratingSet'\right)\right)$ is
    a quotient bialgebra of
    $\PvH\left(\FreePro\left(\GeneratingSet\right)\right)$.
\end{Proposition}
\medbreak

\subsection{The Hopf bialgebra of a stiff pro}
We now extend the construction $\PvH$ to a class a non-necessarily free
pros. Still in this section, $\Pca$ is a free pro.
\medbreak

Let $\equiv$ be a congruence of $\Pca$. For any element $x$ of $\Pca$,
we denote by $[x]_\equiv$ (or by $[x]$ if the context is clear) the
$\equiv$-equivalence class of $x$. We say that $\equiv$ is a
\Def{stiff congruence} if the following three properties are satisfied:
\begin{enumerate}[label = ({\it C{\arabic*}})]
    \item \label{item:stiff_congruences_0}
    for any reduced element $x$ of $\Pca$, the set $[x]$ is finite;
    \item \label{item:stiff_congruences_1}
    for any reduced element $x$ of $\Pca$, $[x]$ contains reduced
    elements only;
    \item \label{item:stiff_congruences_2}
    for any two elements $x$ and $x'$ of $\Pca$ such that $x \equiv x'$,
    the maximal decompositions of $x$ and $x'$ are, respectively of the
    form $(x_1, \dots, x_\ell)$ and $(x'_1, \dots, x'_\ell)$ for some
    $\ell \geq 0$, and for any $i \in [\ell]$, $x_i \equiv x'_i$.
\end{enumerate}
We say that a pro is a \Def{stiff pro} if it is the quotient of a free
pro by a stiff congruence.
\medbreak

For any $\equiv$-equivalence class $[x]$ of reduced elements of $\Pca$,
set
\begin{equation} \label{equ:definition_des_T}
    \BasisT_{[x]} := \sum_{x' \in [x]} \BasisS_{x'}.
\end{equation}
Notice that thanks to \ref{item:stiff_congruences_0}
and~\ref{item:stiff_congruences_1}, $\BasisT_{[x]}$ is a well-defined
element of~$\PvH(\Pca)$.
\medbreak

For instance, if $\Pca$ is the quotient of the free pro generated by
$\GeneratingSet := \GeneratingSet(1, 1) \sqcup \GeneratingSet(2, 2)$
where $\GeneratingSet(1, 1) := \{\Asf\}$ and
$\GeneratingSet(2,2) := \{\Bsf\}$ by the finest congruence $\equiv$
satisfying
\begin{equation}
}\,.
\end{equation}
Moreover, we can observe that $\equiv$ is a stiff congruence.
\medbreak

If $\equiv$ is a stiff congruence of $\Pca$,
\ref{item:stiff_congruences_1} and~\ref{item:stiff_congruences_2} imply
that all the elements of a same $\equiv$-equivalence class $[x]$ have
the same number of factors and are all reduced or all nonreduced. Then,
by extension, we shall say that a $\equiv$-equivalence class $[x]$ of
$\Pca/_\equiv$ is \Def{indecomposable} (resp. \Def{reduced}) if all its
elements are indecomposable (resp. reduced) in $\Pca$. In the same way,
the \Def{wire} of $\Pca/_\equiv$ is the $\equiv$-equivalence class of
the wire of~$\Pca$.
\medbreak

We shall now study how the product and the coproduct of $\PvH(\Pca)$
behave on the $\BasisT_{[x]}$.
\medbreak

\subsubsection{Product}
Let us show that the linear span of the $\BasisT_{[x]}$, where the $[x]$
are $\equiv$-equi\-valence classes of reduced elements of $\Pca$, forms
an associative subalgebra of $\PvH(\Pca)$. The product on the
$\BasisT_{[x]}$ is multiplicative and admits the following simple
description.
\medbreak

\begin{Proposition}
    \label{prop:free_pro_to_Hopf_bialgebra_congruence_product}
    Let $\Pca$ be a free pro and $\equiv$ be a stiff congruence of
    $\Pca$. Then, for any $\equiv$-equivalence classes $[x]$ and $[y]$,
    \begin{equation}
        \BasisT_{[x]} \cdot \BasisT_{[y]} = \BasisT_{[x * y]},
    \end{equation}
    where $x$ (resp. $y$) is any element of $[x]$ (resp. $[y]$).
\end{Proposition}
\medbreak

\subsubsection{Coproduct}
To prove that the linear span of the $\BasisT_{[x]}$, where the $[x]$
are $\equiv$-equi\-valence classes of reduced elements of $\Pca$, forms
a subcoalgebra of $\PvH(\Pca)$ and provides the description of the
coproduct of a $\BasisT_{[x]}$, we need the following notation. For any
element $x$ of~$\Pca$,
\medbreak

\begin{equation} \label{equ:reduced_class_stiff_congruence_pro}
    \Reduced\left([x]\right) :=
    \left\{\Reduced\left(x'\right) : x' \in [x]\right\}.
\end{equation}
\medbreak

\begin{Lemma} \label{lem:stiff_congruence_pro_reduced}
    Let $\Pca$ be a free pro and $\equiv$ be a stiff congruence of
    $\Pca$. For any element $x$ of $\Pca$,
    \begin{equation}
        \Reduced\left([x]\right) = \left[\Reduced(x)\right].
    \end{equation}
\end{Lemma}
\medbreak

\begin{Lemma} \label{lem:stiff_congruence_pro_same_reduced}
    Let $\Pca$ be a free pro, $\equiv$ be a stiff congruence of $\Pca$,
    and $y$ and $z$ be two elements of $\Pca$ such that $y \equiv z$.
    Then, $\Reduced(y) = \Reduced(z)$ implies $y = z$.
\end{Lemma}
\medbreak

The next result is based upon
Lemmas~\ref{lem:stiff_congruence_pro_reduced}
and~\ref{lem:stiff_congruence_pro_same_reduced}.
\medbreak

\begin{Proposition}
\label{prop:free_pro_to_Hopf_bialgebra_congruence_coproduit}
    Let $\Pca$ be a free pro and $\equiv$ be a stiff congruence of
    $\Pca$. Then, for any $\equiv$-equivalence class $[x]$,
    \begin{equation}
        \Delta\left(\BasisT_{[x]}\right) =
        \sum_{\substack{[y], [z] \in \Pca/_\equiv \\
            [y] \circ [z] = [x]}}
            \BasisT_{\Reduced([y])} \otimes \BasisT_{\Reduced([z])}.
    \end{equation}
\end{Proposition}
\medbreak

\subsubsection{Hopf sub-bialgebra}
The description of the product and the coproduct on the $\BasisT_{[x]}$
leads to the following result.
\medbreak

\begin{Theorem} \label{thm:stiff_pro_to_Hopf_bialgebra}
    Let $\Pca$ be a free pro and $\equiv$ be a stiff congruence of
    $\Pca$. Then, the linear span of the $\BasisT_{[x]}$, where the
    $[x]$ are $\equiv$-equivalence classes of reduced elements of
    $\Pca$, forms a Hopf sub-bialgebra of~$\PvH(\Pca)$.
\end{Theorem}
\medbreak

We shall denote, by a slight abuse of notation, by $\PvH(\Pca/_\equiv)$
the sub-bialgebra of $\PvH(\Pca)$ spanned by the $\BasisT_{[x]}$, where
the $[x]$ are $\equiv$-equivalence classes of reduced elements of
$\Pca$. Notice that the construction $\PvH$ as it was presented in
Section~\ref{subsec:free_pro_to_Hopf_bialgebra} is a special case of
this latter when $\equiv$ is the most refined congruence of pros.
\medbreak

Note that this construction of sub-bialgebras of $\PvH(\Pca)$ by taking
an equivalence relation satisfying some precise properties and by
considering the elements obtained by summing over its equivalence
classes is analog to the construction of certain sub-bialgebras of the
Malvenuto-Reutenauer Hopf algebra \cite{MR95}. Indeed, some famous Hopf
algebras are obtained in this way (see Sections~\ref{subsubsec:FQSym}
and~\ref{subsubsec:subbialgebras_FQSym_congruences} of
Chapter~\ref{chap:algebra}).
\medbreak

\subsubsection{The importance of the stiff congruence condition}
Let us now explain why the stiff congruence condition required as a
premise of Theorem~\ref{thm:stiff_pro_to_Hopf_bialgebra} is important by
providing an example of a non-stiff congruence of pros failing to
produce a Hopf bialgebra.
\medbreak

Consider the pro $\Pca$ quotient of the free pro generated by
$\GeneratingSet := \GeneratingSet(1, 1) \sqcup \GeneratingSet(2, 2)$
where $\GeneratingSet(1, 1) := \{\Asf\}$ and
$\GeneratingSet(2, 2) := \{\Bsf\}$ by the finest congruence $\equiv$ satisfying
\begin{equation}

\end{equation}
Here, $\equiv$ is not a stiff congruence since it satisfies
\ref{item:stiff_congruences_1} but not
\ref{item:stiff_congruences_2}.
\medbreak

We have
\begin{equation}
    \BasisT_{
    \left[
    %
}
\end{equation}
but this last element cannot be expressed on the $\BasisT_{[x]}$.
\medbreak

Besides, by a straightforward computation, we have
\begin{multline}
    \Delta \BasisT_{\left[
    %
}\,,
\end{multline}
showing that neither the coproduct is well-defined on
the~$\BasisT_{[x]}$.
\medbreak

\subsubsection{Properties}
By using similar arguments as those used to establish
Proposition~\ref{prop:free_pro_to_Hopf_bialgebra_freeness} together with
the fact that $\equiv$ satisfies \ref{item:stiff_congruences_2} and the
product formula of
Proposition~\ref{prop:free_pro_to_Hopf_bialgebra_congruence_product}, we
obtain that $\PvH(\Pca/_\equiv)$ is freely generated as an algebra by
the $\BasisT_{[x]}$ where the $[x]$ are $\equiv$-equivalence classes of
indecomposable and reduced elements of $\Pca$. Moreover, when $\omega$
is a grading of $\Pca$ so that all elements of a same
$\equiv$-equivalence class have the same degree, the bialgebra
$\PvH(\Pca/_\equiv)$ is graded by the grading inherited the one of
$\PvH(\Pca)$ and forms hence a combinatorial Hopf bialgebra.
\medbreak

\begin{Proposition} \label{prop:stiff_congruences_subbialgebras}
    Let $\Pca$ be a free pro and $\equiv_1$ and $\equiv_2$ be two stiff
    congruences of $\Pca$ such that $\equiv_1$ is finer than $\equiv_2$.
    Then, $\PvH\left(\Pca/_{\equiv_2}\right)$ is a sub-bialgebra of
    $\PvH\left(\Pca/_{\equiv_1}\right)$.
\end{Proposition}
\medbreak

\subsection{Related constructions}
In this section, we first describe two constructions allowing to build
stiff pros. The main interest of these constructions is that the
obtained stiff pros can be placed at the input of the construction
$\PvH$. We next present a way to recover the natural Hopf bialgebra of
an operad through the construction $\PvH$ and the previous constructions
of stiff pros.
\medbreak

\subsubsection{From operads to stiff pros}
\label{subsubsec:operads_to_stiff_pros}
Any operad $\Oca$ gives naturally rise to a pro $\OvP(\Oca)$ whose
elements are sequences of elements of $\Oca$ (see~\cite{Mar08}).
\medbreak

We recall here this construction. Let us set
$\OvP(\Oca) := \sqcup_{p \geq 0} \sqcup_{q \geq 0} \OvP(\Oca)(p, q)$
where
\begin{equation}
    \OvP(\Oca)(p, q) :=
    \left\{
    x_1 \dots x_q : x_i \in \Oca(p_i) \mbox{ for all }  i \in [q]
    \mbox{ and } p_1 + \dots + p_q = p
    \right\}.
\end{equation}
The horizontal composition of $\OvP(\Oca)$ is the concatenation of
sequences, and the vertical composition of $\OvP(\Oca)$ comes directly
from the composition map of $\Oca$. More precisely, for any
$x_1 \dots x_r \in \OvP(\Oca)(q, r)$ and
\begin{math}
    y_{11} \dots y_{1q_1} \dots y_{r1} \dots y_{rq_r}
    \in \OvP(\Oca)(p, q)
\end{math},
we have
\begin{equation} \label{equ:definition_compo_v_construction_r}
    x_1 \dots x_r \circ
    y_{11} \dots y_{1q_1} \dots y_{r1} \dots y_{rq_r}
    := x_1 \circ [y_{11}, \dots, y_{1q_1}] \dots
        x_r \circ [y_{r1}, \dots, y_{rq_r}],
\end{equation}
where for any $i \in [r]$, $x_i \in \Oca(q_i)$ and the occurrences of
$\circ$ in the right-member
of~\eqref{equ:definition_compo_v_construction_r} refer to the total
composition map of~$\Oca$.
\medbreak

For instance, if $\Oca$ is the magmatic operad $\Mag$ (see
Section~\ref{subsubsec:magmatic_operad} of Chapter~\ref{chap:algebra}),
since its elements are binary trees, the elements of the pro
$\OvP(\Oca)$ are forests of binary trees. The horizontal composition of
$\OvP(\Oca)$ is the concatenation of forests, and the vertical
composition $\Ffr_1 \circ \Ffr_2$ in $\OvP(\Oca)$, defined only between
two forests $\Ffr_1$ and $\Ffr_2$ such that the number of leaves of
$\Ffr_1$ is the same as the number of trees in $\Ffr_2$, consists in the
forest obtained by grafting, from left to right, the roots of the trees
of $\Ffr_2$ on the leaves of~$\Ffr_1$.
\medbreak

\begin{Proposition} \label{prop:operad_to_stiff_pro}
    Let $\Oca$ be an operad such that the monoid $(\Oca(1), \circ_1)$
    does not contain any nontrivial subgroup. Then, $\OvP(\Oca)$ is a
    stiff pro.
\end{Proposition}
\medbreak

\subsubsection{From monoids to stiff pros}
\label{subsubsec:monoids_to_stiff_pros}
Any monoid $\Mca$ can be seen as an operad concentrated in arity one.
Then, starting from a monoid $\Mca$, one can construct a pro
$\MvP(\Mca)$ by applying the construction $\OvP$ to $\Mca$ seen as an
operad.
\medbreak

This construction can be rephrased as follows. We have
$\MvP(\Mca) = \sqcup_{p \geq 0} \sqcup_{q \geq 0} \MvP(\Mca)(p, q)$
where
\begin{equation}
    \MvP(\Mca)(p, q) =
    \begin{cases}
        \left\{
            x_1 \dots x_p : x_i \in M \mbox{ for all } i \in [p]
        \right\} &
            \mbox{if } p = q, \\
        \emptyset & \mbox{otherwise}. \\
    \end{cases}
\end{equation}
The horizontal composition of $\MvP(\Mca)$ is the concatenation of
sequences and the vertical composition
$\circ : \MvP(\Mca)(p, p) \times \MvP(\Mca)(p, p) \to \MvP(\Mca)(p, p)$
of $\MvP(\Mca)$ satisfies, for any $x_1 \dots x_p \in \MvP(\Mca)(p, p)$
and $y_1 \dots y_q \in \MvP(\Mca)(q, q)$,
\begin{equation}
    x_1 \dots x_p \circ y_1 \dots y_p =
    (x_1 \Product y_1) \dots (x_p \Product y_p),
\end{equation}
where $\Product$ is the product of~$\Mca$.
\medbreak

For instance, if $\Mca$ is the additive monoid of natural numbers, the
pro $\MvP(\Mca)$ contains all words over $\N$. The horizontal
composition of $\MvP(\Mca)$ is the concatenation of words, and the
vertical composition of $\MvP(\Mca)$, defined only on words with a same
length, is the componentwise addition of their letters.
\medbreak

\begin{Proposition} \label{prop:monoid_to_stiff_pro}
    Let $\Mca$ be a monoid that does not contain any nontrivial
    subgroup. Then, $\MvP(\Mca)$ is a stiff pro.
\end{Proposition}
\medbreak

\subsubsection{The natural Hopf bialgebra of an operad}
We call \Def{abelianization} of a Hopf bialgebra $\Hca$ the Hopf
bialgebra quotient of $\Hca$ by the Hopf bialgebra ideal spanned by the
$x \cdot y - y \cdot x$ for all $x, y \in \Hca$.
\medbreak

Here is the main link between our construction $\PvH$ and the
construction $\OvH$.
\medbreak

\begin{Proposition} \label{prop:natural_Hopf_bialgebra_H}
    Let $\Oca$ be an operad such that the monoid $(\Oca(1), \circ_1)$
    does not contain any nontrivial subgroup. Then, the bialgebra
    $\OvH(\Oca)$ is the abelianization of $\PvH(\OvP(\Oca))$.
\end{Proposition}
\medbreak

\section{Examples of application of the construction}
\label{sec:construction_H_examples}
We conclude this chapter by presenting examples of application of the
construction $\PvH$. The pros considered in this section fit into the
diagram represented by Figure~\ref{fig:diagram_pros} and the obtained
Hopf bialgebras fit into the diagram represented by
Figure~\ref{fig:diagram_Hopf}.
\begin{figure}[ht]
    \centering
    \begin{tikzpicture}[xscale=.55,yscale=.65]
        \node(PRF)at(0,0){\begin{math}\PRF_\gamma\end{math}};
        \node(FBT)at(3,2){\begin{math}\FBT_\gamma\end{math}};
        \node(As)at(3,-2){\begin{math}\As_\gamma\end{math}};
        \node(BAs)at(6,0){\begin{math}\BAs_\gamma\end{math}};
        \node(Heap)at(10,2){\begin{math}\Heap_\gamma\end{math}};
        \node(FHeap)at(10,-2){\begin{math}\FHeap_\gamma\end{math}};
        \draw[Surjection](PRF)--(As);
        \draw[Surjection](FBT)--(BAs);
        \draw[Surjection](Heap)--(FHeap);
        \draw[Injection](PRF)--(FBT);
        \draw[Injection](As)--(BAs);
    \end{tikzpicture}
    \caption[A diagram of pros involving combinatorial objects.]
    {Diagram of pros where arrows $\rightarrowtail$ (resp.
    $\twoheadrightarrow$) are injective (resp. surjective) pro
    morphisms. The parameter $\gamma$ is a positive integer. When
    $\gamma = 0$, $\PRF_0 = \As_0 = \Heap_0 = \FHeap_0$ and
    $\FBT_0 = \BAs_0$.}
    \label{fig:diagram_pros}
\end{figure}
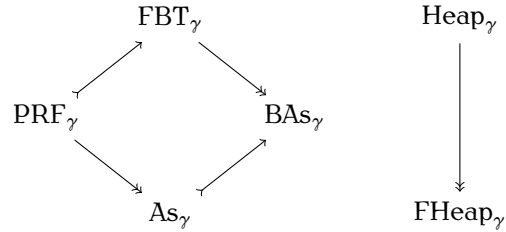
\begin{figure}[ht]
    \centering
    \begin{tikzpicture}[xscale=.55,yscale=.65]
        \node(PRF)at(0,0){\begin{math}\PvH(\PRF_\gamma)\end{math}};
        \node(FBT)at(3,2){\begin{math}\PvH(\FBT_\gamma)\end{math}};
        \node(As)at(3,-2){\begin{math}\PvH(\As_\gamma)\end{math}};
        \node(BAs)at(6,0){\begin{math}\PvH(\BAs_\gamma)\end{math}};
        \node(Heap)at(10,2){\begin{math}\PvH(\Heap_\gamma)\end{math}};
        \node(FHeap)at(10,-2){\begin{math}\PvH(\FHeap_\gamma)\end{math}};
        \draw[Injection](As)--(PRF);
        \draw[Injection](BAs)--(FBT);
        \draw[Injection](FHeap)--(Heap);
        \draw[Surjection](FBT)--(PRF);
        \draw[Surjection](BAs)--(As);
    \end{tikzpicture}
    \caption[A diagram of Hopf bialgebras constructed from pros.]
    {Diagram of combinatorial Hopf bialgebras where arrows
    $\rightarrowtail$ (resp. $\twoheadrightarrow$) are injective (resp.
    surjective) Hopf bialgebras morphisms. The parameter $\gamma$ is a
    positive integer. When $\gamma = 0$,
    $\PvH(\PRF_0) = \PvH(\As_0) = \PvH(\Heap_0) = \PvH(\FHeap_0)$ and
    $\PvH(\FBT_0) = \PvH(\BAs_0)$.}
    \label{fig:diagram_Hopf}
\end{figure}
\medbreak

\subsection{Hopf bialgebra of forests}
We present here the construction of two Hopf bialgebras of forests, one
depending on a nonnegative integer $\gamma$, and with different
gradings. The pro we shall define in this section will intervene in the
next examples.
\medbreak

\subsubsection{Pro of forests with a fixed arity}
Let $\gamma$ be a nonnegative integer and $\PRF_\gamma$ be the free pro
generated by
$\GeneratingSet := \GeneratingSet(\gamma + 1, 1) := \{\Asf\}$, with the
grading $\omega$ defined by $\omega(\Asf) := 1$. Any prograph $x$ of
$\PRF_\gamma$ can be seen as a planar forest of planar rooted trees with
only internal nodes of arity $\gamma + 1$. Since the reduced elements of
$\PRF_\gamma$ have no wire, they are encoded by forests of nonempty
trees.
\medbreak

\subsubsection{Hopf bialgebra}
By Theorem~\ref{thm:free_pro_to_Hopf_bialgebra} and
Proposition~\ref{prop:free_pro_to_Hopf_bialgebra_graduation},
$\PvH(\PRF_\gamma)$ is a combinatorial Hopf bialgebra. By
Proposition~\ref{prop:free_pro_to_Hopf_bialgebra_freeness}, as an
associative algebra, $\PvH(\PRF_\gamma)$ is freely generated by the
$\BasisS_\Tfr$, where the $\Tfr$ are nonempty planar rooted trees with
only internal nodes of arity $\gamma + 1$. Its bases are indexed by
planar forests of such trees where the degree of a basis element
$\BasisS_\Ffr$ is the number of internal nodes of~$\Ffr$.
\medbreak

Notice that the bases of $\PvH(\PRF_0)$ are indexed by forests of linear
trees and that $\PvH(\PRF_0)$ and $\Sym$ are trivially isomorphic as
combinatorial Hopf bialgebras.
\medbreak

\subsubsection{Coproduct}
By definition of the construction $\PvH$, the coproduct of
$\PvH(\PRF_\gamma)$ is given on a generator $\BasisS_\Tfr$ by
\begin{equation} \label{equ:coproduct_forests}
    \Delta(\BasisS_\Tfr) =
    \sum_{\Tfr' \in \AdmissibleCuts(\Tfr)} \BasisS_{\Tfr'} \otimes
    \BasisS_{\Tfr/_{\Tfr'}},
\end{equation}
where $\AdmissibleCuts(\Tfr)$ is the set of \Def{admissible cuts} of
$\Tfr$, that is, the empty tree or the subtrees of $\Tfr$ containing the
root of $\Tfr$ and where $\Tfr/_{\Tfr'}$ denotes the forest consisting
in the maximal subtrees of $\Tfr$ whose roots are leaves of $\Tfr'$, by
respecting the order of these leaves in $\Tfr'$ and by removing the
empty trees. For instance, we have
\begin{multline}
    \Delta \BasisS_{
}
    \otimes
    \BasisS_{\emptyset}.
\end{multline}
\medbreak

This coproduct is similar to the one of the noncommutative
Connes-Kreimer Hopf bialgebra $\CK$~\cite{CK98}. The main difference
between $\PvH(\PRF_\gamma)$ and $\CK$ lies in the fact that in a
coproduct of $\CK$, the admissible cuts can change the arity of some
internal nodes; it is not the case in $\PvH(\PRF_\gamma)$ because for
any $\Tfr' \in \AdmissibleCuts(\Tfr)$, any internal node $u$ of $\Tfr'$
has the same arity as it has in~$\Tfr$.
\medbreak

\subsubsection{Dimensions}
The series of the algebraic generators of
$\PvH(\PRF_\gamma)$ is
\begin{equation}
    \GenSeries(t) :=
    \sum_{n \geq 1} \frac{1}{n\gamma+1}\binom{n(\gamma + 1)}{n} t^n
\end{equation}
since its coefficients are the Fuss-Catalan numbers, counting planar
rooted trees with $n$ internal nodes of arity $\gamma + 1$ (see
Section~\ref{subsubsec:k_ary_trees} of
Chapter~\ref{chap:combinatorics}). Since $\PvH(\PRF_\gamma)$ is free as
an associative algebra, its Hilbert series is
$\HilbSeries_{\PvH(\PRF_\gamma)}(t) := \frac{1}{1 - \GenSeries(t)}$.
\medbreak

The first dimensions of $\PvH(\PRF_1)$ are
\begin{equation} \label{equ:dim_AHC_PRF_1}
    1, 1, 3, 10, 35, 126, 462, 1716, 6435, 24310, 92378,
\end{equation}
and those of $\PvH(\PRF_2)$ are
\begin{equation}
    1, 4, 19, 98, 531, 2974, 17060, 99658, 590563, 3540464, 21430267.
\end{equation}
These two sequences are respectively Sequences~\OEIS{A001700}
and~\OEIS{A047099} of~\cite{Slo}.
\medbreak

\subsubsection{Pro of general forests}
We denote by $\PRF_\infty$ the free pro generated by
\begin{math}
    \GeneratingSet :=
    \sqcup_{n \geq 1} \GeneratingSet(n, 1) :=
    \sqcup_{n \geq 1} \{\Asf_n\}
\end{math}.
Any prograph $x$ of $\PRF_\infty$ can be seen as a planar forest of
planar rooted trees. Since the reduced elements of $\PRF_\infty$ have no
wire, they are encoded by forests of nonempty trees. Observe that for
any nonnegative integer $\gamma$, $\PRF_\gamma$ is a sub-pro
of~$\PRF_\infty$.
\medbreak

\subsubsection{Hopf bialgebra}
By Theorem~\ref{thm:free_pro_to_Hopf_bialgebra}, $\PvH(\PRF_\infty)$ is
a Hopf bialgebra. By
Proposition~\ref{prop:free_pro_to_Hopf_bialgebra_freeness}, as an
associative algebra, $\PvH(\PRF_\infty)$ is freely generated by the
$\BasisS_\Tfr$, where the $\Tfr$ are nonempty planar rooted trees. Its
bases are indexed by planar forests of such trees. Besides, by
Proposition~\ref{prop:free_pro_to_Hopf_bialgebra_subset_generators},
since $\PRF_\gamma$ is generated by a subset of the generators of
$\PRF_\infty$, $\PvH(\PRF_\gamma)$ is a Hopf bialgebra quotient of
$\PvH(\PRF_\infty)$. Moreover, the coproduct of $\PvH(\PRF_\infty)$
satisfies~\eqref{equ:coproduct_forests}.
\medbreak

To turn $\PvH(\PRF_\infty)$ into a combinatorial Hopf bialgebra, we
cannot consider the grading $\omega$ defined by $\omega(\Asf_n) := 1$
because there would be infinitely many elements of degree $1$.
Therefore, we consider on $\PvH(\PRF_\infty)$ the grading $\omega$
defined by $\omega(\Asf_n) := n$. In this way, the degree of a basis
element $\BasisS_\Ffr$ is the number of edges of the forest $\Ffr$. By
Proposition~\ref{prop:free_pro_to_Hopf_bialgebra_graduation},
$\PvH(\PRF_\infty)$ is a combinatorial Hopf bialgebra.
\medbreak

\subsubsection{Dimensions}
The series of the algebraic generators of $\PvH(\PRF_\infty)$ is
\begin{equation}
    \GenSeries(t) :=
    \sum_{n \geq 1} \frac{1}{n + 1} \binom{2n}{n} t^n
\end{equation}
since its coefficients are the Catalan numbers, counting planar rooted
trees with $n$ edges. As $\PvH(\PRF_\infty)$ is free as an
associative algebra, its Hilbert series is
\begin{equation}
    \HilbSeries_{\PvH(\PRF_\infty)}(t) :=
    \frac{1}{1 - \GenSeries(t)} =
    1 + \sum_{n \geq 1} \frac{1}{2}\binom{2n}{n} t^n.
\end{equation}
The dimensions of $\PvH(\PRF_\infty)$ are then the same as the
dimensions of $\PvH(\PRF_1)$ (see~\eqref{equ:dim_AHC_PRF_1}).
\medbreak

\subsection{Faà di Bruno Hopf bialgebra and its deformations}
We shall give here a method to construct the Hopf bialgebras
$\FdBNC_\gamma$ of Foissy~\cite{Foi08} from our construction $\PvH$ in
the case where $\gamma$ is a nonnegative integer.
\medbreak

\subsubsection{Associative pro}
Let $\gamma$ be a nonnegative integer and $\As_\gamma$ be the quotient
of $\PRF_\gamma$ by the finest pro congruence $\equiv$ satisfying
\begin{equation}
\,,
    \qquad k_1 + k_2 = \gamma = \ell_1 + \ell_2.
\end{equation}
\medbreak

We can observe that $\As_\gamma$ is a stiff pro because $\equiv$
satisfies~\ref{item:stiff_congruences_1}
and~\ref{item:stiff_congruences_2}, and that $\As_0 = \PRF_0$. Moreover,
observe that, when $\gamma \geq 1$, there is in $\As_\gamma$ exactly one
indecomposable element of arity $n\gamma + 1$ for any $n \geq 0$. We
denote by $\alpha_n$ this element. We consider on $\As_\gamma$ the
grading $\omega$ inherited the one of $\PRT_\gamma$. This grading
is still well-defined in $\As_\gamma$ since any $\equiv$-equivalence
class contains prographs of a same degree and satisfies, for all
$n \geq 0$, $\omega(\alpha_n) = n$. Any element of $\As_\gamma$ is then
a word $\alpha_{k_1} \dots \alpha_{k_\ell}$ and can be encoded by a word
of nonnegative integers $k_1 \dots k_{\ell}$. Since the reduced elements
of $\As_\gamma$ have no wire, they are encoded by words of positive
integers.
\medbreak

\subsubsection{Hopf bialgebra}
By Theorem~\ref{thm:stiff_pro_to_Hopf_bialgebra} and
Proposition~\ref{prop:free_pro_to_Hopf_bialgebra_graduation},
$\PvH(\As_\gamma)$ is a combinatorial Hopf bialgebra. As an associative
algebra, $\PvH(\As_\gamma)$ is freely generated by the $\BasisT_n$,
$n \geq 1$, and its bases are indexed by words of positive integers
where the degree of a basis element $\BasisT_{k_1 \dots k_\ell}$ is
$k_1 + \dots + k_\ell$.
\medbreak

\subsubsection{Coproduct}
Since any element $\alpha_n$ of $\As_\gamma$ decomposes into
$\alpha_n = x \circ y$ if and only if $x = \alpha_k$ and
$y = \alpha_{i_1} \dots \alpha_{i_{k\gamma + 1}}$ with
$i_1 + \dots + i_{k\gamma + 1} = n - k$, by
Proposition~\ref{prop:free_pro_to_Hopf_bialgebra_congruence_coproduit},
for any $n \geq 1$, the coproduct of $\PvH(\As_\gamma)$ can be
expressed as
\begin{equation} \label{equ:coproduct_analogue_FdBNC_gamma}
    \Delta(\BasisT_n) =
        \sum_{0 \leq k \leq n} \BasisT_k
        \otimes
        \left(\sum_{i_1+\dots+ i_{k\gamma+1} = n - k}
        \BasisT_{i_1}\dots \BasisT_{i_{k\gamma+1}}\right),
\end{equation}
where $\BasisT_0$ is identified with the unit $\BasisT_\epsilon$ of
$\PvH(\As_\gamma)$. For instance, in $\PvH(\As_1)$, we have
\begin{multline}
    \Delta(\BasisT_3) =
        \BasisT_0 \otimes \BasisT_3 +
        \BasisT_1 \otimes
            (\BasisT_0\BasisT_2 + \BasisT_1\BasisT_1 +
            \BasisT_2\BasisT_0) \\
    + \BasisT_2 \otimes
        (\BasisT_0\BasisT_0\BasisT_1 + \BasisT_0\BasisT_1\BasisT_0 +
        \BasisT_1\BasisT_0\BasisT_0)
    + \BasisT_3 \otimes (\BasisT_0\BasisT_0\BasisT_0\BasisT_0) \\
    =
    \BasisT_\epsilon \otimes \BasisT_3 + 2\,\BasisT_1\otimes \BasisT_2 +
    \BasisT_1 \otimes \BasisT_{11} + 3\,\BasisT_2\otimes \BasisT_1 +
    \BasisT_3 \otimes \BasisT_\epsilon,
\end{multline}
and in $\PvH(\As_2)$, we have
\begin{multline}
    \Delta(\BasisT_3) =
    \BasisT_0 \otimes \BasisT_3
    + \BasisT_1 \otimes (\BasisT_0\BasisT_0\BasisT_2 +
            \BasisT_0\BasisT_2\BasisT_0 + \BasisT_2\BasisT_0\BasisT_0 +
            \BasisT_0\BasisT_1\BasisT_1 + \BasisT_1\BasisT_0\BasisT_1 +
            \BasisT_1\BasisT_1\BasisT_0) \\
    + \BasisT_2 \otimes (\BasisT_0\BasisT_0\BasisT_0\BasisT_0\BasisT_1
        + \BasisT_0\BasisT_0\BasisT_0\BasisT_1\BasisT_0
        + \BasisT_0\BasisT_0\BasisT_1\BasisT_0\BasisT_0
        + \BasisT_0\BasisT_1\BasisT_0\BasisT_0\BasisT_0
        + \BasisT_1\BasisT_0\BasisT_0\BasisT_0\BasisT_0) \\
    + \BasisT_3 \otimes
    \BasisT_0\BasisT_0\BasisT_0\BasisT_0\BasisT_0\BasisT_0\BasisT_0 \\
    =
    \BasisT_\epsilon \otimes \BasisT_3 +
    3\,\BasisT_1 \otimes \BasisT_2 +
    3\,\BasisT_1 \otimes \BasisT_{11} +
    5\,\BasisT_2 \otimes \BasisT_1 +
    \BasisT_3 \otimes \BasisT_\epsilon.
\end{multline}
\medbreak

\subsubsection{Deformation of the noncommutative Faà di Bruno Hopf
bialgebra}

\begin{Theorem} \label{thm:construction_FdBNC_gamma}
    For any nonnegative integer $\gamma$, the Hopf bialgebra
    $\PvH(\As_\gamma)$ is the deformation of the noncommutative Faà di
    Bruno Hopf bialgebra $\FdBNC_\gamma$.
\end{Theorem}
\medbreak

\subsection{Hopf bialgebra of forests of bitrees}
To define Hopf bialgebras of forests of bitrees, we need the following
general construction on pros.
\medbreak

\subsubsection{Symmetrization of pros}
If $\GeneratingSet$ is a bigraded collection of the form
\begin{math}
    \GeneratingSet =
    \sqcup_{p \geq 1} \sqcup_{q \geq 1} \GeneratingSet(p, q)
\end{math},
we denote by $\GeneratingSet^\Op$ the bigraded collection defined by
\begin{equation}
    \GeneratingSet^\Op(p, q) :=
    \GeneratingSet(q, p), \qquad p, q \geq 1.
\end{equation}
From a geometrical point of view, any elementary prograph over
$\GeneratingSet^\Op$ is obtained by reversing from bottom to top an
elementary prograph over $\GeneratingSet$. We moreover denote by
\begin{math}
    \Rev : \FreePro(\GeneratingSet^\Op) \to \FreePro(\GeneratingSet)
\end{math}
the bijection sending any prograph $x$ of
$\FreePro(\GeneratingSet^\Op)$ to the prograph $\Rev(x)$ of
$\FreePro(\GeneratingSet)$ obtained by reversing $x$ from bottom to top.
\medbreak

Now, given a pro $\Pca := \FreePro(\GeneratingSet)/_\equiv$, we define
the \Def{symmetrization} $\Sat(\Pca)$ of $\Pca$ as the pro
\begin{equation}
    \Sat(\Pca) :=
    \FreePro\left(
        \GeneratingSet \sqcup \GeneratingSet^\Op\right)/_{\cong},
\end{equation}
where $\cong$ is the finest congruence of
$\FreePro(\GeneratingSet \sqcup \GeneratingSet^\Op)$ satisfying
\begin{equation}
    x \cong y
    \quad \mbox{if }
    (x, y \in \FreePro(\GeneratingSet) \mbox{ and } x \equiv y)
    \enspace \mbox{ or } \enspace
    (x, y \in \FreePro(\GeneratingSet^\Op) \mbox{ and }
        \Rev(x) \equiv \Rev(y)).
\end{equation}
Notice that in this definition, we consider $\FreePro(\GeneratingSet)$
and $\FreePro(\GeneratingSet^\Op)$ as sub-pros of
$\FreePro(\GeneratingSet \sqcup \GeneratingSet^\Op)$ in an obvious way.
Notice also that if $\Pca$ is a free pro $\FreePro(\GeneratingSet)$,
then the congruence $\equiv$ is trivial, so that $\cong$ is also
trivial, and
\begin{math}
    \Sat(\Pca) = \FreePro(\GeneratingSet \sqcup \GeneratingSet^\Op)
\end{math}.
Besides, as another immediate property of this construction, remark
that when $\Pca$ is a stiff pro, the congruence $\cong$
satisfies~\ref{item:stiff_congruences_1}
and~\ref{item:stiff_congruences_2}, and then, $\Sat(\Pca)$ is a stiff
pro.
\medbreak

We shall present here two Hopf bialgebras coming from the construction
$\Sat$ applied to $\PRF_\gamma$ and~$\As_\gamma$.
\medbreak

\subsubsection{Pro of forests of bitrees}
Let $\gamma$ be a nonnegative integer and $\FBT_\gamma$ be the free pro
generated by
\begin{math}
    \GeneratingSet :=
    \GeneratingSet(\gamma + 1, 1) \sqcup \GeneratingSet(1, \gamma + 1)
\end{math}
where $\GeneratingSet(\gamma+1, 1) := \{\Asf\}$ and
$\GeneratingSet(1, \gamma + 1) := \{\Bsf\}$, with the grading $\omega$
defined by $\omega(\Asf) := \omega(\Bsf) := 1$. One has
$\Sat(\PRF_\gamma) = \FBT_\gamma$. Any prograph $x$ of $\FBT_\gamma$ can
be seen as a forest of \Def{$\gamma$-bitrees}, that are labeled planar
trees where internal nodes labeled by $\Asf$ have $\gamma+1$ children
and one parent, and the internal nodes labeled by $\Bsf$ have one child
and $\gamma + 1$ parents. Since the reduced elements of $\FBT_\gamma$
have no wire, they are encoded by forests of nonempty $\gamma$-bitrees.
\medbreak

\subsubsection{Hopf bialgebra}
By Theorem~\ref{thm:free_pro_to_Hopf_bialgebra} and
Proposition~\ref{prop:free_pro_to_Hopf_bialgebra_graduation},
$\PvH(\FBT_\gamma)$ is a combinatorial Hopf bialgebra. By
Proposition~\ref{prop:free_pro_to_Hopf_bialgebra_freeness}, as an
associative algebra, $\PvH(\FBT_\gamma)$ is freely generated by the
$\BasisS_\Tfr$, where the $\Tfr$ are nonempty $\gamma$-bitrees. Its
bases are indexed by planar forests of such bitrees where the degree of
a basis element $\BasisS_\Ffr$ is the total number of internal nodes in
the bitrees of $\Ffr$. Moreover, by
Proposition~\ref{prop:free_pro_to_Hopf_bialgebra_subset_generators},
since $\PRF_\gamma$ is generated by a subset of the generators of
$\FBT_\gamma$, $\PvH(\PRF_\gamma)$ is a quotient bialgebra
of~$\PvH(\FBT_\gamma)$.
\medbreak

\subsubsection{Coproduct}
The coproduct of $\PvH(\FBT_\gamma)$ can be described, like the one of
$\CK$ on forests, by means of admissible cuts on forests of
$\gamma$-bitrees. We have for instance
\begin{multline}
    \Delta \BasisS_{
} \otimes \BasisS_{\emptyset}.
\end{multline}

\subsubsection{Dimensions}
We only know the dimensions of $\PvH(\FBT_\gamma)$ when $\gamma = 0$. In
this case, $0$-bitrees of size $n$ are linear trees and can hence be
seen as words of length $n$ on the alphabet $\{\Asf, \Bsf\}$. Therefore,
as $\PvH(\FBT_\gamma)$ is free as an associative algebra, the bases of
$\PvH(\FBT_0)$ are indexed by multiwords on $\{\Asf, \Bsf\}$ and its
Hilbert series is
\begin{equation}
    \HilbSeries_{\PvH(\FBT_0)}(t) := 1 + \sum_{n \geq 1} 2^{2n - 1} t^n.
\end{equation}
\medbreak

\subsubsection{Pro of biassociative operators and its Hopf bialgebra}
Let $\BAs_\gamma$ be the quotient of $\FBT_\gamma$ by the finest pro
congruence $\equiv$ satisfying
\begin{equation}
\,,
    \qquad k_1 + k_2 = \gamma = \ell_1 + \ell_2.
\end{equation}
\medbreak

We can observe that $\BAs_\gamma$ is a stiff pro because $\equiv$
satisfies~\ref{item:stiff_congruences_1}
and~\ref{item:stiff_congruences_2}. Notice that
$\Sat(\As_\gamma) = \BAs_\gamma$ and $\BAs_0 = \FBT_0$. We consider on
$\BAs_\gamma$ the grading $\omega$ inherited the one of
$\FBT_\gamma$. This grading is still well-defined in $\BAs_\gamma$ since
any $\equiv$-equivalence class contains prographs of a same degree.
Notice that $\BAs_1$ is very similar to the pro governing Hopf
bialgebras (see~\cite{Mar08}). Indeed, it only lacks in $\BAs_1$ the
usual compatibility relation between its two generators. Notice also
that the pro governing bialgebras is not a stiff pro.
\medbreak

By Theorem~\ref{thm:stiff_pro_to_Hopf_bialgebra} and
Proposition~\ref{prop:free_pro_to_Hopf_bialgebra_graduation},
$\PvH(\BAs_\gamma)$ is then a combinatorial Hopf bialgebra. Moreover, we
can observe that $\PvH(\As_\gamma)$ is a quotient Hopf bialgebra
of~$\PvH(\BAs_\gamma)$.
\medbreak

\subsection{Hopf bialgebra of heaps of pieces}
We present here the construction of a Hopf bialgebra depending on a
nonnegative integer $\gamma$, whose bases are indexed by heaps of
pieces.
\medbreak

\subsubsection{Pro of heaps of pieces}
Let $\gamma$ be a nonnegative integer and $\Heap_\gamma$ be the free pro
generated by
\begin{math}
    \GeneratingSet :=
    \GeneratingSet(\gamma + 1, \gamma + 1) :=
    \{\Asf\}
\end{math},
with the grading $\omega$ defined by $\omega(\Asf) := 1$. Any prograph
$x$ of $\Heap_\gamma$ can be seen as a heap of pieces of width
$\gamma + 1$ (see~\cite{Vie86} for some theory about these objects). For
instance, the prograph
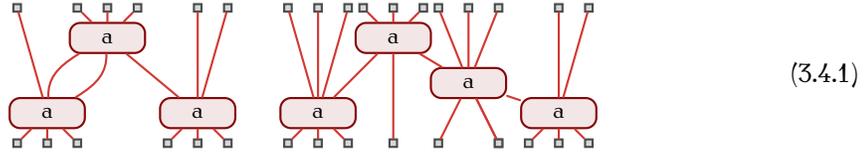
\begin{equation}

\end{equation}
of $\Heap_2$ is encoded by the heap of pieces of width $3$ depicted by
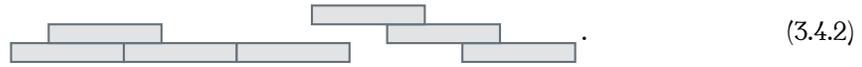
\begin{equation}
    %
\,.
\end{equation}
Notice that $\Heap_0 = \PRF_0$. Besides, since the reduced elements of
$\Heap_\gamma$ have no wire, they are encoded by horizontally connected
heaps of pieces of width $\gamma + 1$.
\medbreak

\subsubsection{Hopf bialgebra}
By Theorem~\ref{thm:free_pro_to_Hopf_bialgebra} and
Proposition~\ref{prop:free_pro_to_Hopf_bialgebra_graduation},
$\PvH(\Heap_\gamma)$ is a combinatorial Hopf bialgebra. By
Proposition~\ref{prop:free_pro_to_Hopf_bialgebra_freeness}, as an
associative algebra, $\PvH(\Heap_\gamma)$ is freely generated by the
$\BasisS_{\LambdaB}$ where the $\LambdaB$ are heaps of pieces that
cannot be obtained by juxtaposing two heaps of pieces. Its bases are
indexed by horizontally connected heaps of pieces of width $\gamma + 1$
where the degree of a basis element $\BasisS_{\LambdaB}$ is the number
of pieces of~$\LambdaB$.
\medbreak

\subsubsection{Coproduct}
The coproduct of $\PvH(\Heap_\gamma)$ can be described, like the one of
$\CK$ on forests, by means of admissible cuts on heaps of pieces.
Indeed, if $\LambdaB$ is a horizontally connected heap of pieces, by
definition of the construction $\PvH$,
\begin{equation}
    \Delta(\BasisS_{\LambdaB}) =
    \sum_{\LambdaB' \in \AdmissibleCuts(\LambdaB)}
    \BasisS_{\LambdaB'} \otimes \BasisS_{\LambdaB/_{\LambdaB'}},
\end{equation}
where $\AdmissibleCuts(\LambdaB)$ is the set of \Def{admissible cuts} of
$\LambdaB$, that is, the set of heaps of pieces obtained by keeping an
upper part of $\LambdaB$ and by readjusting it so that it becomes
horizontally connected and where $\LambdaB/_{\LambdaB'}$ denotes the
heap of pieces obtained by removing from $\LambdaB$ the pieces of
$\LambdaB'$ and by readjusting the remaining pieces so that they form an
horizontally connected heap of pieces. For instance, in $\PvH(\Heap_1)$,
we have
\begin{multline}
    \Delta
    \BasisS_{
}
    \otimes
    \BasisS_{\emptyset}\,.
\end{multline}
\medbreak

\subsubsection{Dimensions}
\begin{Proposition} \label{prop:dimensions_PvH_Heap_gamma}
    For any nonnegative integer $\gamma$, the Hilbert series
    $\HilbSeries_{\PvH(\Heap_\gamma)}(t)$
    of $\PvH(\Heap_\gamma)$ satisfies
    $\HilbSeries_{\PvH(\Heap_\gamma)}(t) = \sum_{n \geq 0} C_{\gamma, n}(t)$, where
    \begin{equation}
        C_{\gamma, n}(t) := P_{\gamma, n}(t)
        - \sum_{k = 0}^{n - 1}
            C_{\gamma, k}(t) P_{\gamma, n - k - 1}(t),
    \end{equation}
    \begin{equation}
        P_{\gamma, n}(t) := \frac{1}{F_{\gamma, n}(t)},
    \end{equation}
    and
    \begin{equation}
        F_{\gamma, n}(t) :=
        \begin{cases}
            1 & \mbox{if } n \leq \gamma, \\
            F_{\gamma, n - 1}(t) - t F_{\gamma, n - \gamma - 1}(t)
                & \mbox{otherwise}.
        \end{cases}
    \end{equation}
\end{Proposition}
\medbreak

By using Proposition~\ref{prop:dimensions_PvH_Heap_gamma}, one can
compute the first dimensions of $\PvH(\Heap_\gamma)$. The first
dimensions of $\PvH(\Heap_1)$ are
\begin{equation} \label{equ:dim_Heap_1}
    1, 1, 4, 18, 85, 411, 2014, 9950, 49417, 246302, 1230623,
\end{equation}
and those of $\PvH(\Heap_2)$ are
\begin{equation} \label{equ:dim_Heap_2}
    1, 1, 6, 42, 313, 2407, 18848, 149271, 1191092, 9553551, 76910632.
\end{equation}
Since by Proposition \ref{prop:free_pro_to_Hopf_bialgebra_freeness},
$\PvH(\Heap_\gamma)$ is free as an associative , the series
$\GenSeries_\gamma(t)$ of its algebraic generators satisfies
\begin{math}
    \GenSeries_\gamma(t) =
    1 - \frac{1}{\HilbSeries_{\PvH(\Heap_\gamma)}(t)}
\end{math}.
The first dimensions of the algebraic generators of $\PvH(\Heap_1)$ are
\begin{equation} \label{equ:dim_generateurs_Heap_1}
    1, 3, 11, 44, 184, 790, 3450, 15242, 67895, 304267, 1369761,
\end{equation}
and those of $\PvH(\Heap_2)$ are
\begin{equation}
    1, 5, 31, 210, 1488, 10826, 80111, 599671, 4525573, 34357725,
    262011295.
\end{equation}
These four integer sequences are respectively
Sequences~\OEIS{A247637}, \OEIS{A247638}, \OEIS{A059715},
and~\OEIS{A247639} of~\cite{Slo}.
\medbreak

\subsection{Hopf bialgebra of heaps of friable pieces}
By considering special quotient of $\Heap_\gamma$, we construct a
Hopf bialgebra structure on the $(\gamma + 1)$st tensor power of the
vector space~$\Sym$.
\medbreak

\subsubsection{Pro of heaps of friable pieces}
Let $\gamma$ be a nonnegative integer and $\FHeap_\gamma$ be the
quotient of $\Heap_\gamma$ by the finest pro congruence $\equiv$
satisfying
\begin{equation}
\,,
    \qquad \ell \in [\gamma].
\end{equation}
For instance, for $\gamma = 2$, the $\equiv$-equivalence class of
\begin{equation}
    %
\end{equation}
contains exactly the prographs
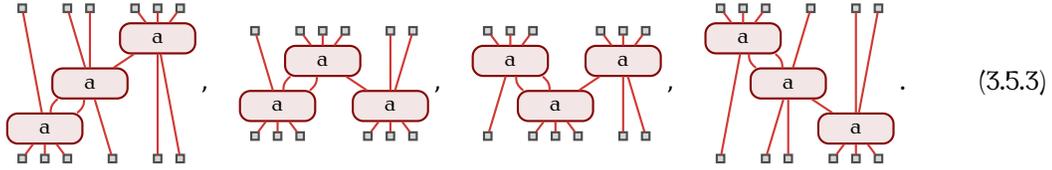
\begin{equation} \label{equ:exemple_classe_FHeap_3}
    %
\,.
\end{equation}
\medbreak

We can observe that $\FHeap_\gamma$ is a stiff pro because $\equiv$
satisfies~\ref{item:stiff_congruences_1}
and~\ref{item:stiff_congruences_2} and $\FHeap_0 = \Heap_0$. We call
$\FHeap_\gamma$ the \Def{pro of heaps of friable pieces} of width
$\gamma + 1$. This terminology is justified by the following
observation. Any piece of width $\gamma + 1$ (depicted by
$\tikz \node[Domino2](0)at(0,0){};$) consists in $\gamma + 1$ small
pieces, called \Def{bursts}, glued together. This forms a
\Def{friable piece} (depicted, for $\gamma = 2$ for instance, by
\begin{math}
\begin{tikzpicture}[Centering]
    \node[DominoFriable](0)at(0,0){};
    \node[DominoFriable](1)at(.5,0){};
    \node[DominoFriable](2)at(1,0){};
\end{tikzpicture}
\end{math}).
The congruence $\equiv$ of $\Heap_\gamma$ can be interpreted by letting
all pieces break under gravity, separating the bursts constituting them.
For instance, the prographs of \eqref{equ:exemple_classe_FHeap_3},
respectively, encoded by the heaps of pieces
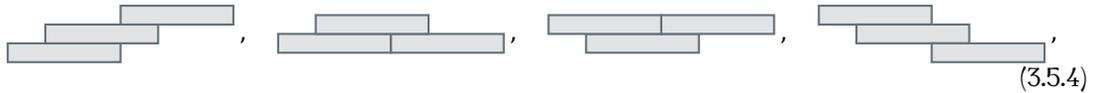
\begin{equation} \label{equ:example_class_Heap_3}
    \begin{tikzpicture}[Centering]
        \node[Domino3](0)at(0,0){};
        \node[Domino3](1)at(.5,.25){};
        \node[Domino3](2)at(1.5,.5){};
    \end{tikzpicture}\,,\quad
    \begin{tikzpicture}[Centering]
        \node[Domino3](0)at(0,0){};
        \node[Domino3](1)at(1.5,0){};
        \node[Domino3](2)at(.5,.25){};
    \end{tikzpicture}\,,\quad
    \begin{tikzpicture}[Centering]
        \node[Domino3](0)at(0,0){};
        \node[Domino3](1)at(-.5,.25){};
        \node[Domino3](2)at(1,.25){};
    \end{tikzpicture}\,,\quad
    \begin{tikzpicture}[Centering]
        \node[Domino3](0)at(0,0){};
        \node[Domino3](1)at(.5,-.25){};
        \node[Domino3](2)at(1.5,-.5){};
    \end{tikzpicture}\,,
\end{equation}
all become the heap of friable pieces
\begin{equation} \label{equ:example_friable_heap}
    \begin{tikzpicture}[Centering]
        \node[DominoFriable](0)at(0,0){};
        \node[DominoFriable](1)at(.5,0){};
        \node[DominoFriable](2)at(1,0){};
        \node[DominoFriable](3)at(1.5,0){};
        \node[DominoFriable](4)at(2,0){};
        \node[DominoFriable](5)at(2.5,0){};
        \node[DominoFriable](6)at(.5,.25){};
        \node[DominoFriable](7)at(1,.25){};
        \node[DominoFriable](8)at(1.5,.25){};
    \end{tikzpicture}
\end{equation}
obtained by replacing each piece of any heap of pieces
of~\eqref{equ:example_class_Heap_3} by friable pieces.
\medbreak

The grading $\omega$ of $\FHeap_\gamma$ is the one inherited the one
of $\Heap_\gamma$. This grading is still well-defined in $\Heap_\gamma$
since any $\equiv$-equivalence class contains prographs of a same
degree. Since the reduced elements of $\FHeap_\gamma$ have no wire,
they are encoded by horizontally connected heaps of friable pieces.
\medbreak

Besides, $\FHeap_\gamma$ admits the following alternative description
using the $\MvP$ construction (see
Section~\ref{subsubsec:monoids_to_stiff_pros}). Indeed,
$\FHeap_\gamma$ is the sub-pro of $\MvP(\N)$ generated by
$1^{\gamma + 1}$, where $\N$ denotes here the additive monoid of
nonnegative integers and $1^{\gamma+1}$ denotes the sequence of
$\gamma + 1$ occurrences of $1 \in \N$. The correspondence between heaps
of friable pieces and words of integers of this second description is
clear since any element $x$ of the sub-pro of $\MvP(\N)$ generated by
$1^{\gamma+1}$ encodes a heap of friable pieces consisting, from left
to right, in columns of $x_i$ bursts for $i \in [n]$, where $n$ is the
length of $x$. For instance, the word $122211$ encodes the heap of
friable pieces of~\eqref{equ:example_friable_heap}.
\medbreak

\subsubsection{Hopf bialgebra}
By Theorem~\ref{thm:stiff_pro_to_Hopf_bialgebra} and
Proposition~\ref{prop:free_pro_to_Hopf_bialgebra_graduation},
$\PvH(\FHeap_\gamma)$ is a combinatorial Hopf sub-bialgebra of
$\PvH(\Heap_\gamma)$. The bases of $\PvH(\FHeap_\gamma)$ are indexed by
horizontally connected heaps of friable pieces of width $\gamma + 1$
where the degree of a basis element $\BasisT_{\LambdaB}$ is the number
of pieces of $\LambdaB$.
\medbreak

\subsubsection{Coproduct}
The coproduct of $\PvH(\FHeap_\gamma)$ can be described with the aid of
the interpretation of $\FHeap_\gamma$ as a sub-pro of $\MvP(\N)$.
Indeed, if $\LambdaB$ is an horizontally connected heap of friable
pieces, by
Proposition~\ref{prop:free_pro_to_Hopf_bialgebra_congruence_coproduit},
\begin{equation}
    \Delta(\BasisT_{\LambdaB}) =
    \sum_{\substack{\LambdaB_1, \LambdaB_2 \in \FHeap_\gamma \\
    \LambdaB = \LambdaB_1 + \LambdaB_2}}
    \BasisT_{\LambdaB'_1} \otimes \BasisT_{\LambdaB'_2},
\end{equation}
where $\LambdaB_1 + \LambdaB_2$ is the heap of friable pieces obtained
by stacking $\LambdaB_2$ onto $\LambdaB_1$ and where $\LambdaB'_1$
(resp. $\LambdaB'_2$) is the readjustment of $\LambdaB_1$ (resp.
$\LambdaB_2$) so that it is horizontally connected. For instance, we
have in $\PvH(\FHeap_1)$
\begin{equation}
    \BasisT_{
}
    \otimes \BasisT_{\emptyset}.
\end{multline}
\medbreak

\subsubsection{Dimensions}

\begin{Proposition} \label{prop:dimensions_PvH_FHeap_gamma}
    For any nonnegative integers $\gamma$ and $n$, the $n$th homogeneous
    component of $\PvH(\FHeap_\gamma)$ has dimension
    $(\gamma + 2)^{n - 1}$.
\end{Proposition}
\medbreak

\subsubsection{Miscellaneous properties}
By the dimensions of $\PvH(\FHeap_\gamma)$ provided by
Proposition~\ref{prop:dimensions_PvH_FHeap_gamma}, as a graded vector
space, $\PvH(\FHeap_\gamma)$ is the $\gamma + 1$st tensor power of the
underlying vector space of $\Sym$. Indeed, the $n$th homogeneous
components of these two spaces have the same dimension. Besides, notice
that since $\FHeap_\gamma$ is by definition a sub-pro of the pro
obtained by applying the construction $\MvP$ to a commutative monoid,
$\PvH(\FHeap_\gamma)$ is cocommutative.
\medbreak

\section*{Concluding remarks}
We have defined a construction $\PvH$ establishing a new link between
the theory of pros and the theory of combinatorial Hopf bialgebras, by
generalizing a former construction from operads to Hopf bialgebras. By
the way, we have exhibited the so-called stiff pros which is the most
general class of pros for which our construction works.
\smallbreak

By using $\PvH$, we have introduced some new and recovered some
already known combinatorial Hopf bialgebras by starting with very
simple pros. Nevertheless, we are very far from having exhausted
the possibilities, and it would not be surprising that $\PvH$ could
reconstruct some other known Hopf bialgebras, maybe in unexpected bases.
\smallbreak

Computing the Hilbert series of a combinatorial Hopf bialgebra is,
usually, a routine work. Nevertheless, in the general case, it is very
difficult to compute the Hilbert series of $\PvH(\Pca)$ when $\Pca$ is a
free pro. Indeed, this computation requires to know, given a free pro
$\Pca$, the series
\begin{equation}
    \HilbSeries_\Pca(t) :=
    \sum_{x \in \Reduced(\Pca)} t^{\deg(x)},
\end{equation}
which seems difficult to explicitly describe in general.
\smallbreak

As another perspective, it is conceivable to go further in the study of
the algebraic structure of the bialgebras obtained by $\PvH$. The
question of the potential autoduality of $\PvH(\Pca)$ depending on some
conditions on the pro $\Pca$ is worth studying. A way to solve this
problem is to provide enough conditions on $\Pca$ to endow $\PvH(\Pca)$
with a bidendriform bialgebra structure~\cite{Foi07} (see also
Section~\ref{subsubsec:dendriform_algebras} of
Chapter~\ref{chap:algebra}). This strategy is based upon the fact that
any bidendriform bialgebra is free and self-dual as a
bialgebra~\cite{Foi07}.
\medbreak

\part{Combinatorics and algorithms}

\chapter{Shuffle of permutations} \label{chap:shuffle}
The content of this chapter comes from~\cite{GV16} and is a joint work
with Stéphane Vialette.
\medbreak

\section*{Introduction}
The shuffle product $\shuffle$ is a well-known operation on words first
defined by Eilenberg and Mac Lane~\cite{EM53}. Given three words $u$,
$v_1$, and $v_2$, $u$ is said to be a shuffle of $v_1$ and $v_2$ if it
can be formed by interleaving the letters from $v_1$ and $v_2$ in a way
that maintains the left-to-right ordering of the letters from each word
(see also Section~\ref{subsubsec:associative_algebras} of
Chapter~\ref{chap:algebra}). Besides purely combinatorial questions, the
shuffle product of words naturally leads to the following computational
problems:
\begin{enumerate}[label={\it (\roman*)}]
    \item \label{item:problem_1}
    Given two words $v_1$ and $v_2$, compute the set of the words
    appearing in the shuffle of $v_1$ with $v_2$;
    \item \label{item:problem_2}
    Given three words $u$, $v_1$, and $v_2$, decide if $u$ appears in
    the shuffle of $v_1$ with $v_2$;
    \item \label{item:problem_3}
    Given a word $u$, decide if there is a word $v$ such that $u$ is
    in the shuffle of $v$ with itself.
\end{enumerate}
\smallbreak

Even if these problems seem similar, they radically differ in terms of
time complexity. Let us now review some facts about these. In what
follows, $n$ denotes the size of $u$ and $m_i$ denotes the size of each
$v_i$. A solution to Problem~\ref{item:problem_1} can be computed in
\begin{math}
  O\left(\binom{m_1 + m_2}{m_1}\right)
\end{math}
time~\cite{Spe86,All00}. Problem~\ref{item:problem_2} is in $\PClass$;
it is indeed a classical textbook exercise to design an efficient
dynamic programming algorithm solving it. It can be tested in
$O\left(n^2 / \log(n)\right)$ time~\cite{LN82}. To the best of our
knowledge, the first $O(n^2)$ time algorithm for this problem appeared
in~\cite{Man83}. This algorithm can easily be extended to check in
polynomial-time whether a word is in the shuffle of any fixed number of
given words. Let us now finally focus on Problem~\ref{item:problem_3}.
It is shown in~\cite{RV13,BS14} that it is $\NPC$ to decide if a word
$u$ is a square with respect to the shuffle, that is a word $u$ with the
property that there exists a word $v$ such that $u$ appears in the
shuffle of $v$ with itself. Hence, Problem~\ref{item:problem_3}
is~$\NPC$.
\smallbreak

This chapter is intended to study a natural generalization of $\shuffle$,
denoted by $\SuperShuffle$, as a shuffle of permutations. Roughly
speaking, given three permutations $\pi$, $\sigma_1$, and $\sigma_2$,
$\pi$ is said to be a shuffle of $\sigma_1$ and $\sigma_2$ if $\pi$
(viewed as a word) appears in the shuffle of two words whose
standardized permutations are respectively $\sigma_1$ and $\sigma_2$.
This shuffle product was first introduced by Vargas~\cite{Var14} under
the name of supershuffle. Our intention in this work is to study this
shuffle product of permutations $\SuperShuffle$ both from a
combinatorial and from a computational point of view by focusing on
square permutations, that are permutations $\pi$ appearing in the
shuffle of a permutation $\sigma$ with itself. Many other shuffle
products on permutations appear in the literature. For instance,
in~\cite{DHT02}, the authors define the convolution product and the
shifted shuffle product (see Section~\ref{subsubsec:FQSym} of
Chapter~\ref{chap:algebra}). It is a simple exercise to prove that,
given three permutations $\pi$, $\sigma_1$, and $\sigma_2$, deciding if
$\pi$ is in the shifted shuffle of $\sigma_1$ and $\sigma_2$ is
in~$\PClass$.
\smallbreak

This chapter is organized as follows. In
Section~\ref{sec:square_elements}, we introduce the general notion of
square elements in algebras. We take as examples the case of the shifted
shuffle product of permutations and the shuffle product of words. We
provide a definition of the supershuffle of permutations in
Section~\ref{sec:supershuffle}, by introducing it from its dual
coproduct $\Delta$, called unshuffling coproduct. Some algebraic and
combinatorial properties of these product and coproduct are reviewed.
Section~\ref{sec:algorithm_supershuffle} is devoted to contain an
algorithmic study of square permutations with respect to
$\SuperShuffle$. The most important result of this work, concerning the
fact that deciding if a permutation is square is $\NPC$, appears here.
\medbreak

\section{Square elements, shuffles, and words}
\label{sec:square_elements}
Before defining and studying the supershuffle product of permutations,
we set here a general algebraic framework about square elements in
algebras.
\medbreak

\subsection{Square elements with respect to a product}
The general notion of square elements in polynomial algebras endowed
with a binary product is introduced here. This notion relies on the
notions of collections (see Section~\ref{sec:collections} of
Chapter~\ref{chap:combinatorics}) and of polynomial algebras (see
Section~\ref{sec:bialgebras} of Chapter~\ref{chap:algebra} for the basic
definitions about these structures).
\medbreak

\subsubsection{General definitions}
Let $C$ be a collection and $\K \Angle{C}$ be an algebra (not
necessarily an associative algebra) endowed with a binary product
$\Product$. An object $x$ of $C$ is a \Def{square} with respect to
$\Product$ if there is an object $y$ of $C$ such that $x$ appears in the
product $y \Product y$. In this case, we say that $y$ is a
\Def{square root} of~$x$. Observe that this notion depends on the basis
of $C$ of $\K \Angle{C}$. Indeed, seen on anther basis $C'$ of
$\K \Angle{C}$, square elements can be different.
\medbreak

By duality, by considering the dual
$({\K \Angle{C}}^\Dual, \Coproduct_\Product)$ of
$(\K \Angle{C}, \Product)$, an element $x$ of $C$ is a square and
$y \in C$ is one of its square root if and only if the tensor
$y \otimes y$ appears in $\Coproduct_\Product(x)$. Indeed, by definition
$x$ is a square and $y$ is one of its square root if and only if the
structure coefficient $\xi^{(y \otimes y, x)}_\Product$ of $\Product$ is
nonzero (see Section~\ref{subsubsec:biproducts} of
Chapter~\ref{chap:algebra}). Hence, this is equivalent to say that the
structure coefficient
$\xi^{(x, y \otimes y)}_{\Coproduct_\Product}$ of $\Coproduct_\Product$
is nonzero.
\medbreak

In the case where $(\K \Angle{C}, \Product)$ is graded, there are two
algorithmic problems related to these concepts. The first one takes as
input an element $x$ of $C$ of size $n$ and consists in deciding whether
$x$ is a square. We call this problem the
\Def{square detection problem $\ProblemSD$}. The second one takes as
input an element $x$ of $C$ of size $n$ and another element $y$ of $C$
and consists in deciding whether $y$ is a square root of $x$. We call
this problem the \Def{square root checking problem $\ProblemRC$}. The
complexity of these two problems is studied with respect to~$n$.
\medbreak

\subsubsection{Square permutations in $\FQSym$}
\label{subsubsec:square_permutations_FQSym}
To give an example of this notion of squares, consider the space
$\FQSym = \K \Angle{\SymmetricGroup}$ of all permutations endowed with
the shifted shuffle product (see Section~\ref{subsubsec:FQSym} of
Chapter~\ref{chap:algebra}). The permutation $316425$ is a square since
it appears in the shifted shuffle of $312$ with itself. The first square
permutations are
\begin{equation}
    \epsilon, \quad
    12, 21, \quad
    1234, 1324, 1342, 3124, 3142, 3412, 2143, 2413, 4213, 2431, 4231,
    4321,
\end{equation}
and the sequence of the number of square permutations begins by
\begin{equation} \label{equ:sequence_square_binary_words}
    1, 0, 2, 0, 12, 0, 120, 0, 1680, 0, 30240, 0, 665280, 0, 17297280,
\end{equation}
forming (after removing the $0$s) Sequence~\OEIS{A001813} of~\cite{Slo}.
\medbreak

To test if a permutation $\sigma$ of an even size $n$ is a square, one
can extract its subword $u$ consisting in the letters in
$\left\{1, 2, \dots, \frac{n}{2}\right\}$, its subword $v$ consisting in
the letters in
$\left\{\frac{n}{2} + 1, \frac{n}{2} + 2, \dots, n\right\}$, and
checking if $\Std(v) = u$, where $\Std$ is the standardization
algorithm (see Section~\ref{subsubsec:permutations} of
Chapter~\ref{chap:combinatorics}). Since all these operations are
obviously polynomial in $n$, $\ProblemSD$ is polynomial.
\medbreak

Besides, to check if a permutation $\nu$ is a square root of a
permutation $\sigma$ of an even size $n$, one can check if
the subword $u$ of $\sigma$ consisting in the letters in
$\left\{1, 2, \dots, \frac{n}{2}\right\}$ of $\sigma$ and if the
standardized of subword $v$ consisting in the letters in
$\left\{\frac{n}{2} + 1, \frac{n}{2} + 2, \dots, n\right\}$ are both
equal to $\nu$. Since these operations are polynomial in $n$,
$\ProblemRC$ is polynomial.
\medbreak

\subsection{Square words for the shuffle product}
We now turn our attention to square words for the shuffle product
and the complexity of $\ProblemSD$. These results come from~\cite{RV13}
and are used as prototype in our upcoming study of square permutations.
\medbreak

\subsubsection{Square words}
Let us consider here the shuffle algebra $(\K \Angle{A^*}, \shuffle)$
where $A$ is a finite alphabet (see
Section~\ref{subsubsec:associative_algebras} of
Chapter~\ref{chap:algebra}). For instance, if $A := \{\Asf, \Bsf\}$, the
word $\Asf \Bsf \Asf \Asf \Bsf \Asf$ is a square since it appears in the
shuffle of $\Asf \Bsf \Asf$ with itself. Contrariwise, the word
$\Asf \Bsf \Bsf \Asf$ is not a square. Observe that having an even
number of occurrences for each letter of $A$ is a necessary condition to
be a square.
\medbreak

\subsubsection{Perfect matchings}
Let us describe a way to decide if a word $u$ is a square. This
comes from~\cite{RV13} and use perfect matchings on words.
A \Def{perfect matching} on a word $u \in A^*$ is a graph
$(V, E)$ such that
\begin{equation}
    V := \left\{(u(i), i) : i \in [|u|]\right\},
\end{equation}
every vertex of $V$ belongs to exactly one edge of $E$, and
$\left\{(u(i), i), (u(j), j)\right\} \in E$ implies $u(i) = u(j)$.
Figure~\ref{fig:example_perfect_matching_word} shows a perfect matching
on a word.
\begin{figure}[ht]
    \begin{tikzpicture}
        [xscale=.5,yscale=.25,node distance=1.5cm,Centering]
        \node[MatchingVertex](U01){\begin{math}(\Bsf, 1)\end{math}};
        \node[MatchingVertex,right of=U01](U02)
            {\begin{math}(\Asf, 2)\end{math}};
        \node[MatchingVertex,right of=U02](U03)
            {\begin{math}(\Bsf, 3)\end{math}};
        \node[MatchingVertexColorA,right of=U03](U04)
            {\begin{math}(\Asf, 4)\end{math}};
        \node[MatchingVertexColorA,right of=U04](U05)
            {\begin{math}(\Bsf, 5)\end{math}};
        \node[MatchingVertexColorA,right of=U05](U06)
            {\begin{math}(\Bsf, 6)\end{math}};
        \node[MatchingVertex,right of=U06](U07)
            {\begin{math}(\Csf, 7)\end{math}};
        \node[MatchingVertexColorA,right of=U07](U08)
            {\begin{math}(\Csf, 8)\end{math}};
        \draw[MatchingEdge]
            (U01.north) .. controls ($ (U01.north) + (0,2.5) $)
            and ($ (U06.north) + (0,2.5) $) .. (U06.north);
        \draw[MatchingEdge]
            (U02.north) .. controls ($ (U02.north) + (0,1.5) $)
            and ($ (U04.north) + (0,1.5) $) .. (U04.north);
        \draw[MatchingEdge]
            (U03.south) .. controls ($ (U03.south) + (0,-1.5) $)
            and ($ (U05.south) + (0,-1.5) $) .. (U05.south);
        \draw[MatchingEdge]
            (U07.north) .. controls ($ (U07.north) + (0,1.5) $)
            and ($ (U08.north) + (0,1.5) $) .. (U08.north);
    \end{tikzpicture}
    \caption{A perfect matching on the word
    $\Bsf\Asf\Bsf\Asf\Bsf\Bsf\Csf\Csf$.}
    \label{fig:example_perfect_matching_word}
\end{figure}
\medbreak

A perfect matching $(V, E)$ is \Def{containment-free} if there are no
edges $\{(u(i), i), (u(j), j)\}$ and $\{(u(i'), i'), (u(j'), j')\}$
of $E$ such that $i < i' < j' < j$. Observe that the perfect matching
of Figure~\ref{fig:example_perfect_matching_word} is not
containment-free.
\medbreak

The criterion of Rizzi and Vialette~\cite{RV13} to recognize square
words is the following.
\medbreak

\begin{Proposition} \label{prop:criterion_recognizing_square_words}
    A word $u \in A^*$ is a square if and only if there exists a
    containment-free perfect matching on~$u$.
\end{Proposition}
\medbreak

Figure~\ref{fig:examples_perfect_matchings_square_words} shows two
examples related to
Proposition~\ref{prop:criterion_recognizing_square_words}.
\begin{figure}[ht]
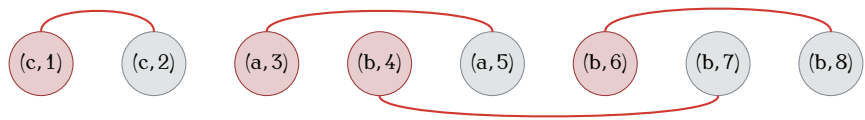
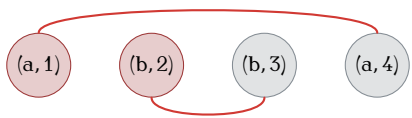

    \subfloat[][A containment-free perfect matching on the word
    $\Csf\Csf\Asf\Bsf\Asf\Bsf\Bsf\Bsf$, showing that it is
    a square. The associated square root is
    $\textcolor{Col1}{\Csf \Asf \Bsf \Bsf}$.]{
    \begin{minipage}[c]{.9\textwidth}
        \centering

    \end{minipage}
    } \\
    \subfloat[][The word $\Asf\Bsf\Bsf\Asf$ is not a square since it
    admits this only perfect matching which is not containment-free.]{
    \begin{minipage}[c]{.9\textwidth}
        \centering
        %
    \end{minipage}
    }
    \caption{Two perfect matchings on words.}
    \label{fig:examples_perfect_matchings_square_words}
\end{figure}
\medbreak

As a consequence of this criterion, given a containment-free perfect
matching on $u$, a square root of $u$ is readable by observing the
subword $u(i_1) u(i_2) \dots u(i_n)$ of $u$ such that the indices $i_1$,
$i_2$, \dots, $i_n$ satisfy $i_1 < i_2 < \dots i_n$ and, for all
$k \in [n]$, $\{(u(i_k), i_k), (u(j_k), j_k)\}$ is an edge of the
perfect matching where~$i_k < j_k$.
\medbreak

\subsubsection{Recognizing square words}
By using the criterion provided by
Proposition~\ref{prop:criterion_recognizing_square_words}, it is
possible to perform a polynomial-time reduction from the longest common
subsequence problem for binary words (which is $\NPC$) to $\ProblemSD$.
This leads to the following result~\cite{RV13}.
\medbreak

\begin{Theorem} \label{thm:hardness_recognizing_square_words}
    In the shuffle algebra $(\K \Angle{A^*}, \shuffle)$ where $A$ is a
    finite alphabet, $\ProblemSD$ is~$\NPC$.
\end{Theorem}
\medbreak

\section{Supershuffle of permutations} \label{sec:supershuffle}
The purpose of this section is to define a shuffle product
$\SuperShuffle$ on permutations, different from the shifted shuffle (see
Section~\ref{subsubsec:square_permutations_FQSym}). Recall that a first
definition of this product was provided by Vargas~\cite{Var14}. To
present an alternative definition of this product adapted to our study,
we shall first define a coproduct denoted by $\Delta$, enabling to
unshuffle permutations. By duality, $\Delta$ implies the definition of
$\SuperShuffle$. The reason why we need to pass through the
definition of $\Delta$ to define $\SuperShuffle$ is justified by the
fact that a lot of properties of $\SuperShuffle$ depend of properties
of $\Delta$, and that this strategy allows to write concise and clear
proofs of them.
\medbreak

\subsection{Unshuffling coalgebra and square permutations}
After defining the unshuffling coproduct of permutations, we define the
supershuffle product. The first properties of this product and of its
square elements are reviewed.
\medbreak

\subsubsection{Unshuffling coproduct}
\label{subsubsec:unshuffling_coproduct}
Let us say that two permutations $\sigma$ and $\nu$ are
\Def{order-iso\-morphic} if $\Std(\sigma) = \Std(\nu)$. We endow the
polynomial space $\FQSym$ with the linear coproduct $\Delta$ defined in
the following way. For any permutation $\pi$, we set
\begin{equation} \label{equ:unshuffling_coproduct}
    \Delta(\pi) =
    \sum_{P_1 \sqcup P_2 = [|\pi|]}
    \Std\left(\pi_{|P_1}\right) \otimes \Std\left(\pi_{|P_2}\right).
\end{equation}
We call $\Delta$ the \Def{unshuffling coproduct of permutations}. For
instance,
\begin{subequations}
\begin{equation}
    \Delta(213) =
    \epsilon \otimes 213 + 2 \cdot 1 \otimes 12 +
    1 \otimes 21 + 2 \cdot 12 \otimes 1 + 21 \otimes 1 +
    213 \otimes \epsilon,
\end{equation}
\begin{equation}
    \Delta(1234) =
    \epsilon \otimes 1234 + 4 \cdot 1 \otimes 123 +
    6 \cdot 12 \otimes 12 + 4 \cdot 123 \otimes 1 +
    1234 \otimes \epsilon,
\end{equation}
\begin{multline} \label{equ:example_unshuffling_coproduct}
    \Delta(1432) =
    \epsilon \otimes 1432 +
    3 \cdot \textcolor{Col1}{1} \otimes \textcolor{Col4}{132}
    + 1 \otimes 321 + 3 \cdot 12 \otimes 21 \\
    + 3 \cdot 21 \otimes 12 + 3 \cdot 132 \otimes 1
    + 321 \otimes 1 + 1432 \otimes \epsilon.
\end{multline}
\end{subequations}
Observe that the coefficient of the tensor
$\textcolor{Col1}{1} \otimes \textcolor{Col4}{132}$ is $3$
in~\eqref{equ:example_unshuffling_coproduct} because there are exactly
three ways to extract from the permutation $1432$ two disjoint subwords
respectively order-isomorphic to the permutations $1$ and~$132$ (that
are $(4, 132)$, $(3, 142)$, and $(2, 143)$).
\medbreak

\subsubsection{Supershuffle product}
Now, by definition of duality, the dual product of $\Delta$, denoted by
$\SuperShuffle$, is a linear binary product on $\FQSym^\Dual$. Since
$\FQSym$ is a graded combinatorial polynomial space,
$\FQSym \simeq \FQSym^\Dual$, so that we shall identify these two
spaces. We call $\SuperShuffle$ the \Def{supershuffle} of permutations.
This product satisfies, for any permutations $\sigma$ and~$\nu$,
\begin{equation}
    \sigma \SuperShuffle \nu =
    \sum_{\pi \in \SymmetricGroup}
    \xi^{(\pi, \sigma \otimes \nu)}_\Delta
    \pi,
\end{equation}
where the coefficients $\xi^{(\pi, \sigma \otimes \nu)}_\Delta$ are the
structure coefficients of $\Delta$. For instance,
\begin{multline} \label{equ:example_supershuffle_product}
    12 \SuperShuffle 21 =
    1243 + 1324 + 2 \cdot 1342 + 2 \cdot 1423
    + 3 \cdot \textcolor{Col1}{1432} + 2134 + 2 \cdot 2314 \\
    + 3 \cdot 2341 + 2413 + 2 \cdot 2431 + 2 \cdot 3124 + 3142 +
    3 \cdot 3214 + 2 \cdot 3241 \\
    + 3421 + 3 \cdot 4123 + 2 \cdot 4132 + 2 \cdot 4213 + 4231 + 4312.
\end{multline}
Observe that the coefficient $3$ of the permutation $1432$
in~\eqref{equ:example_supershuffle_product} comes from the fact that the
coefficient of the tensor $12 \otimes 21$ is $3$
in~\eqref{equ:example_unshuffling_coproduct}.
\medbreak

Intuitively, the supershuffle blends the values and the positions of the
letters of the permutations. One can observe that the empty permutation
$\epsilon$ is a unit for $\SuperShuffle$ and that this product is graded
by the sizes of the permutations.
\medbreak

\subsubsection{Square permutations}
According to Section~\ref{sec:square_elements}, a permutation $\pi$ is
a square with respect to $\SuperShuffle$ if there is a permutation
$\sigma$ such that $\pi$ appears in $\sigma \SuperShuffle \sigma$. In
this case, we say that $\sigma$ is a \Def{square root} of $\pi$.
Equivalently, $\pi$ is a square with $\sigma$ as a square root if and
only if the tensor $\sigma \otimes \sigma$ appears in $\Delta(\pi)$. The
first square permutations are
\begin{multline}
    \epsilon, \quad
    12, 21, \quad
    1234, 1243, 1324, 1342, 1423,
    2134, 2143, 2314, 2413, 2431, \\
    3124, 3142, 3241, 3412, 3421,
    4132, 4213, 4231, 4312, 4321.
\end{multline}
\medbreak

In a more combinatorial way, this is equivalent to saying that there are
two sets $P_1$ and $P_2$ of disjoints indices of letters of $\pi$
satisfying $P_1 \sqcup P_2 = [|\pi|]$ such that the subwords
$\pi_{|P_1}$ and $\pi_{|P_2}$ are order-isomorphic. Computer experiments
give us the first numbers of square permutations with respects to their
size, which are, from size $0$ to $10$,
\begin{equation} \label{equ:sequence_square_permutations}
    1, 0, 2, 0, 20, 0, 504, 0, 21032, 0, 1293418.
\end{equation}
This sequence (after removing the $0$s) is known as
Sequence~\OEIS{A279200} of~\cite{Slo}. We do not have any description
(by a formula, recurrence, or generating series) of these numbers.
\medbreak

\subsection{Square binary words and permutations}
In this section, we shall establish the fact that the square binary
words ({\em i.e.}, square words on the alphabet $\{0, 1\}$ with respect
to the shuffle product) are in one-to-one correspondence with square
permutations avoiding some patterns. This property establishes a link
between the shuffle of binary words and the supershuffle of permutations
and allows to obtain a new description of square binary words.
\medbreak

\subsubsection{From binary words to permutations}
Let $u$ be a binary word of length $n$ with $k$ occurrences of $0$. We
denote by $\BinToPerm$ (Binary word To Permutation) the map sending any
such word $u$ to the permutation obtained by replacing from left to
right each occurrence of $0$ in $u$ by $1$, $2$, \dots, $k$, and from
right to left each occurrence of $1$ in $u$ by $k + 1$, $k + 2$, \dots,
$n$. For instance,
\begin{equation}
    \BinToPerm(\textcolor{Col4}{1}\textcolor{Col1}{0}
    \textcolor{Col4}{11}\textcolor{Col1}{000}
    \textcolor{Col4}{1}\textcolor{Col1}{0})
    =
    \textcolor{Col4}{9}\textcolor{Col1}{1}
    \textcolor{Col4}{87}\textcolor{Col1}{234}
    \textcolor{Col4}{6}\textcolor{Col1}{5}.
\end{equation}
Observe that for any nonempty permutation $\pi$ in the image of
$\BinToPerm$, there is exactly one binary word $u$ such that
$\BinToPerm(u0) = \BinToPerm(u1) = \pi$. In support of this observation,
when $\pi$ has an even size, we denote by $\PermToBin(\pi)$ (Permutation
To Binary word) the word $ua$ such that $|ua|_0$ and $|ua|_1$ are both
even, where $a \in \{0, 1\}$. For instance,
\begin{subequations}
\begin{equation}
    \PermToBin(\textcolor{Col4}{6}
    \textcolor{Col1}{1}\textcolor{Col4}{54}
    \textcolor{Col1}{23})
    =
    \textcolor{Col4}{1}\textcolor{Col1}{0}
    \textcolor{Col4}{11}\textcolor{Col1}{00},
\end{equation}
\begin{equation}
    \PermToBin(\textcolor{Col1}{1}\textcolor{Col4}{4}
    \textcolor{Col1}{2}\textcolor{Col4}{3})
    =
    \textcolor{Col1}{0}\textcolor{Col4}{1}
    \textcolor{Col1}{0}\textcolor{Col4}{1}.
\end{equation}
\end{subequations}
\medbreak

\subsubsection{Link between square binary words and square permutations}
\begin{Proposition} \label{prop:bijection_binary_to_permutations_squares}
  For any $n \geq 0$, the map $\BinToPerm$ restricted to the set of
  square binary words of length $2n$ is a bijection between this last
  set and the set of square permutations of size $2n$ avoiding the
  patterns $213$ and~$231$.
\end{Proposition}
\medbreak

The number of square binary words is (after removing the $0$s)
Sequence~\OEIS{A191755} of~\cite{Slo} beginning by
\begin{equation}
    1, 0, 2, 0, 6, 0, 22, 0, 82, 0, 320, 0, 1268, 0, 5102, 0, 020632.
\end{equation}
According to
Proposition~\ref{prop:bijection_binary_to_permutations_squares}, this
is also the sequence enumerating square permutations avoiding $213$
and~$231$.
Notice that it is conjectured in~\cite{HRS12} that the number of
binary words of length $2n$ is
\begin{equation}
    \binom{2n}{n} \frac{1}{n + 1}2^n -
    \binom{2n - 1}{n + 1}2^{n - 1} +
    O\left(2^{n-2}\right).
\end{equation}
\medbreak

\subsection{Algebraic properties}
The aim of this section is to establish some of properties of the
supershuffle product of permutations $\SuperShuffle$. It is worth to
note that, as we will see, algebraic properties of the unshuffling
coproduct $\Delta$ of permutations defined in
Section~\ref{subsubsec:unshuffling_coproduct} lead to combinatorial
properties of~$\SuperShuffle$.
\medbreak

\subsubsection{Associativity and commutativity}

\begin{Proposition} \label{prop:shuffle_associative_commutative}
    The supershuffle product $\SuperShuffle$ of permutations is
    associative and commutative, that is $(\FQSym, \SuperShuffle)$ is
    an associative commutative algebra.
\end{Proposition}
\begin{proof}
    To prove the associativity of $\SuperShuffle$, we shall prove that
    its dual coproduct $\Delta$ is coassociative. This strategy relies
    on the fact that a product is associative if and only if its dual
    coproduct is coassociative. For any permutation $\pi$, by denoting
    by $\Identity$ the identity map on $\FQSym$, we have
    \begin{equation} \begin{split}
        \label{equ:shuffle_associative_commutative}
        (\Delta \otimes \Identity) \Delta(\pi) & =
        (\Delta \otimes \Identity)
        \sum_{P_1 \sqcup P_2 = [|\pi|]}
        \Std\left(\pi_{|P_1}\right) \otimes \Std\left(\pi_{|P_2}\right) \\
        & =
        \sum_{P_1 \sqcup P_2 = [|\pi|]}
        \Delta\left(\Std\left(\pi_{|P_1}\right)\right)
        \otimes I\left(\Std\left(\pi_{|P_2}\right)\right) \\
        & =
        \sum_{P_1 \sqcup P_2 = [|\pi|]} \;
        \sum_{Q_1 \sqcup Q_2 = [|P_1|]}
        \Std\left(\Std\left(\pi_{|P_1}\right)_{|Q_1}\right)
        \otimes
        \Std\left(\Std\left(\pi_{|P_1}\right)_{|Q_2}\right)
        \otimes \Std\left(\pi_{|P_2}\right) \\
        & =
        \sum_{P_1 \sqcup P_2 \sqcup P_3 = [|\pi|]}
        \Std\left(\pi_{|P_1}\right) \otimes
        \Std\left(\pi_{|P_2}\right) \otimes
        \Std\left(\pi_{|P_3}\right).
    \end{split} \end{equation}
    An analogous computation shows that
    $(\Identity \otimes \Delta) \Delta(\pi)$ is equal to the last member
    of~\eqref{equ:shuffle_associative_commutative}, whence the
    associativity of $\SuperShuffle$.
    \smallbreak

    Finally, to prove the commutativity of $\SuperShuffle$, we shall
    show that $\Delta$ is cocommutative, that is for any permutation
    $\pi$, if in the expansion of $\Delta(\pi)$ there is a tensor
    $\sigma \otimes \nu$ with a coefficient $\lambda$, there is in the
    same expansion the tensor $\nu \otimes \sigma$ with the same
    coefficient $\lambda$. Clearly, a product is commutative if and only
    if its dual coproduct is cocommutative. Now, from the
    definition~\eqref{equ:unshuffling_coproduct} of $\Delta$, one
    observes that if the pair $(P_1, P_2)$ of subsets of $[|\pi|]$
    contributes to the coefficient of
    $\Std\left(\pi_{|P_1}\right) \otimes \Std\left(\pi_{|P_2}\right)$,
    the pair $(P_2, P_1)$ contributes to the coefficient of
    $\Std\left(\pi_{|P_2}\right) \otimes \Std\left(\pi_{|P_1}\right)$.
    This shows that $\Delta$ is cocommutative and hence, that
    $\SuperShuffle$ is commutative.
\end{proof}

Proposition~\ref{prop:shuffle_associative_commutative} implies that
$(\FQSym, \Delta)$ is a coassociative cocommutative coalgebra.
\medbreak

\subsubsection{Endomorphisms}
If $\pi$ is a permutation of $\SymmetricGroup(n)$, we denote by
$\widetilde{u}$ the \Def{mirror image} of $u$, that is the word
$u_{|u|} u_{|u| - 1} \dots u_1$, by $\bar \pi$ the \Def{complement} of
$\pi$, that is the permutation satisfying $\bar \pi(i) = n - \pi(i) + 1$
for all $i \in [n]$, and by $\pi^{-1}$ the \Def{inverse} of $\pi$.
\medbreak

\begin{Proposition} \label{prop:endomorphisms_supershuffle}
    The three linear maps
    \begin{equation}
        \phi_1, \phi_2, \phi_3 : \FQSym \to \FQSym
    \end{equation}
    linearly sending a permutation $\pi$ to, respectively,
    $\widetilde{\pi}$, $\bar \pi$, and $\pi^{-1}$ are endomorphisms of
    the associative algebra $(\FQSym, \SuperShuffle)$.
\end{Proposition}
\medbreak

\subsubsection{Operations preserving square permutations}
We now use the algebraic properties of $\SuperShuffle$ exhibited by
Proposition~\ref{prop:endomorphisms_supershuffle} to obtain
combinatorial properties of square permutations.
\medbreak

\begin{Proposition} \label{prop:square_stability}
    Let $\pi$ be a square permutation and $\sigma$ be a square root of
    $\pi$. Then,
    \begin{enumerate}[label={\it (\roman*)}]
        \item \label{item:square_stability_1}
        the permutation $\widetilde{\pi}$ is a square and
        $\widetilde{\sigma}$ is one of its square roots;
        \item \label{item:square_stability_2}
        the permutation $\bar \pi$ is a square and $\bar \sigma$ is one
        of its square roots;
        \item \label{item:square_stability_3}
        the permutation $\pi^{-1}$ is a square and $\sigma^{-1}$ is one
        of its square roots.
    \end{enumerate}
\end{Proposition}
\begin{proof}
    All statements~\ref{item:square_stability_1},
    \ref{item:square_stability_2}, and~\ref{item:square_stability_3} are
    consequences of Proposition~\ref{prop:endomorphisms_supershuffle}.
    Indeed, since $\pi$ is a square permutation and $\sigma$ is a square
    root of $\pi$, by definition, $\pi$ appears in the product
    $\sigma \SuperShuffle \sigma$. Now, by
    Proposition~\ref{prop:endomorphisms_supershuffle}, for any
    $j \in [3]$, since $\phi_j$ is an endomorphism of associative
    algebras of $\FQSym$, $\phi_j$ commutes with the shuffle product of
    permutations $\SuperShuffle$. Hence, in particular, one has
    \begin{equation}
        \phi_j(\sigma \SuperShuffle \sigma) =
        \phi_j(\sigma) \SuperShuffle \phi_j(\sigma).
    \end{equation}
    Then, since $\pi$ appears in $\sigma \SuperShuffle \sigma$,
    $\phi_j(\pi)$ appears in $\phi_j(\sigma \SuperShuffle \sigma)$ and
    appears also in $\phi_j(\sigma) \SuperShuffle \phi_j(\sigma)$. This
    shows that $\phi_j(\sigma)$ is a square root of $\phi_j(\pi)$ and
    implies~\ref{item:square_stability_1}, \ref{item:square_stability_2},
    and~\ref{item:square_stability_3}.
\end{proof}
\medbreak

Let us make an observation about Wilf-equivalence classes of
permutations restrained on square permutations. Recall that two
permutations $\sigma$ and $\nu$ of the same size are
\Def{Wilf equivalent} if
\begin{math}
    \# \SymmetricGroup(n)^{\{\sigma\}} =
    \# \SymmetricGroup(n)^{\{\nu\}}
\end{math}
for all $n \geq 0$. The well-known~\cite{SS85} fact that there is a
single Wilf-equivalence class of permutations of size $3$ together with
Proposition~\ref{prop:square_stability} imply that $123$ and $321$ are
in the same Wilf-equivalence class of square permutations, and that
$132$, $213$, $231$, and $312$ are in the same Wilf-equivalence class of
square permutations. Computer experiments show us that there are two
Wilf-equivalence classes of square permutations of size $3$. Indeed, the
number of square permutations avoiding $123$ begins by
\begin{equation} \label{equ:sequence_square_123}
    1, 0, 2, 0, 12, 0, 118, 0, 1218, 0, 14272,
\end{equation}
while the number of square permutations avoiding $132$ begins by
\begin{equation} \label{equ:sequence_square_132}
    1, 0, 2, 0, 11, 0, 84, 0, 743, 0, 7108.
\end{equation}
\medbreak

Besides, another consequence of Proposition~\ref{prop:square_stability}
is that it makes sense to enumerate the sets of square permutations
quotiented by the operations of mirror image, complement, and inverse.
The sequence enumerating these sets begins by
\begin{equation} \label{equ:sequence_square_classes}
    1, 0, 1, 0, 6, 0, 81, 0, 2774, 0, 162945.
\end{equation}
\medbreak

These three sequences~\eqref{equ:sequence_square_123},
\eqref{equ:sequence_square_132}, and~\eqref{equ:sequence_square_classes}
are (after removing the $0$s) respectively Sequences~\OEIS{A279201},
\OEIS{A279202}, and~\OEIS{A279203} of~\cite{Slo}. We do not have any
description (by a formula, recurrence, or generating series) of these
numbers.
\medbreak

\section{Algorithmic aspects of square permutations}
\label{sec:algorithm_supershuffle}
We construct here an analog of the criterion to recognize square words
of Rizzi and Valette~\cite{RV13} for the case of square permutations.
This criterion is a center piece to show that $\ProblemSD$ for square
permutation is~$\NPC$.
\medbreak

\subsection{Directed perfect matchings on permutations}
We need a little more complicated combinatorial object than perfect
matching here. We work with directed perfect matchings and some notions
of pattern avoidance.
\medbreak

\subsubsection{Directed graphs}
A \Def{directed graph} is an ordered pair $(V, A)$ where $V$ is a set
whose elements are called \Def{vertices} and $A$ is a set of ordered
pairs of vertices, called \Def{arcs} (from a \Def{source} vertex to a
\Def{sink} vertex). Notice that the aforementioned definition does not
allow a directed graph to have multiple arcs with same source and target
nodes. We shall not allow directed loops (that is, arcs that connect
vertices with themselves). Two arcs are \Def{independent} if they do not
have any common vertex. A directed graph is a \Def{directed matching} if
all its arcs are independent. A directed matching is \Def{perfect} if
every vertex is either a source or a sink.
\medbreak

\subsubsection{Directed perfect matchings}
A \Def{directed perfect matching} on a permutation $\pi$ of an even size
$2n$ is a directed perfect matching $(V, A)$ such that
\begin{equation}
    V := \left\{(\pi(i), i) : i \in [2n] \right\}.
\end{equation}
Figure~\ref{fig:example_directed_perfect_matching_permutation} shows a
directed perfect matching on a permutation. The \Def{word of sources}
(resp. \Def{word of sinks}) of $(V, A)$ is the subword
$\pi(i_1) \pi(i_2) \dots \pi(i_n)$ of $\pi$ where the indices
$i_1 < i_2 < \dots < i_n$ are the sources (resp. sinks) of the arcs
of~$(V, A)$.
\begin{figure}[ht]
    \centering
    \begin{tikzpicture}
        [xscale=0.35,yscale=.3,node distance=1.1cm]
        \node[MatchingVertex](U01){\begin{math}(4, 1)\end{math}};
        \node[MatchingVertexColorA,right of=U01](U02)
            {\begin{math}(1, 2)\end{math}};
        \node[MatchingVertex,right of=U02](U03)
            {\begin{math}(3, 3)\end{math}};
        \node[MatchingVertex,right of=U03](U04)
            {\begin{math}(2, 4)\end{math}};
        \node[MatchingVertexColorA,right of=U04](U05)
            {\begin{math}(8, 5)\end{math}};
        \node[MatchingVertexColorA,right of=U05](U06)
            {\begin{math}(5, 6)\end{math}};
        \node[MatchingVertex,right of=U06](U07)
            {\begin{math}(7, 7)\end{math}};
        \node[MatchingVertexColorA,right of=U07](U08)
            {\begin{math}(6, 8)\end{math}};
        \draw[DirectedMatchingEdge]
            (U01.north) .. controls ($ (U01.north) + (0,4.25) $)
            and ($ (U06.north) + (0,4.25) $) .. (U06.north);
        \draw[DirectedMatchingEdge]
            (U03.north) .. controls ($ (U03.north) + (0,4.25) $)
            and ($ (U08.north) + (0,4.25) $) .. (U08.north);
        \draw[DirectedMatchingEdge]
            (U04.north) .. controls ($ (U04.north) + (0,2.5) $)
            and ($ (U02.north) + (0,2.5) $) .. (U02.north);
        \draw[DirectedMatchingEdge]
            (U07.north) .. controls ($ (U07.north) + (0,2.5) $)
            and ($ (U05.north) + (0,2.5) $) .. (U05.north);
    \end{tikzpicture}
    \caption[A directed perfect matching on the permutation
    $41328576$.]
    {A directed perfect matching $\Mca$ on the permutation
    $41328576$. The word of sources of $\Mca$ is
    $\textcolor{Col1}{4327}$ and its word of sinks is
    $\textcolor{Col4}{1856}$.}
    \label{fig:example_directed_perfect_matching_permutation}
\end{figure}
\medbreak

A \Def{pattern} is a directed perfect matching $([2k], B)$. We say that
a direct perfect matching $(V, A)$ on a permutation $\pi$ admits an
\Def{occurrence} of $([2k], B)$ if
\begin{enumerate}[label={\it (\roman*)}]
    \item there is a map $\phi : [2k] \to V$ such that, for any
    $i, j \in [2k]$, $i < j$ implies that the second coordinate of
    $\phi(i)$ is smaller than the second coordinate of $\phi(j)$;
    \item for any arc $(i, j)$ of $([2k], B)$, $(\phi(i), \phi(j))$ is
    an arc of $(V, A)$.
\end{enumerate}
Observe that this notion of pattern occurrence does not depend on the
permutation $\pi$. Intuitively, $(V, A)$ admits an occurrence of
$([2k], B)$ if $(V, A)$ contains a copy of $([2k], B)$ as a subgraph by
changing some of its labels if necessary and by preserving their order
induced by the second coordinates of the labeling pairs. When $(V, A)$
does not admit any occurrence of $([2k], B)$, we say that $(V, A)$
\Def{avoids} $([2k], B)$. In the sequel, we shall draw patterns as
unlabeled directed graphs. The vertices of the patterns are implicitly
labeled by $1$, $2$, \dots, $2k$ from left to right.
\medbreak

For instance, a directed perfect matching $(V, A)$ on a permutation
$\pi$ admits an occurrence of the pattern $\InclusionRL$ if there are
four vertices $(\pi(i_1), i_1)$, $(\pi(i_2), i_2)$, $(\pi(i_3), i_3)$,
$(\pi(i_4), i_4)$ of $(V, A)$ such that $i_1 < i_2 < i_3 < i_4$, and
$((\pi(i_1), i_1), (\pi(i_4), i_4))$ and
$((\pi(i_3), i_3), (\pi(i_2), i_2))$ are arcs of $(V, A)$. The directed
perfect matching of
Figure~\ref{fig:example_directed_perfect_matching_permutation} admits
hence exactly two occurrences of this pattern: a first one for the arcs
$((4, 1), (5, 6))$ and $((2, 4), (1, 2))$, and a second one for the arcs
$((3, 3), (6, 8))$ and $((7, 7), (8, 5))$.
\medbreak

\subsection{Hardness of recognizing square permutations}
We are now in position to state our criterion to decide if a permutation
is a square.
\medbreak

\subsubsection{Recognizing square permutations}
Let us define two additional properties on directed perfect matchings
on permutations. Let $(V, A)$ be a directed perfect matching on a
permutation $\pi$. We say that $(V, A)$ \Def{satisfies $\PropDPMShape$}
if it avoids all the patterns of the set
\begin{equation}
    \left\{
    \InclusionLL, \InclusionLR, \InclusionRL,
    \InclusionRR, \CrossingLR, \CrossingRL
    \right\}.
\end{equation}
Besides, we say that $(V, A)$ \Def{satisfies $\PropDPMValue$} if for any
two distinct arcs $((\pi(i), i), (\pi(j), j))$ and
$((\pi(i'), i'), (\pi(j'), j'))$ of $(V, A)$, we have $\pi(i) < \pi(i')$
if and only if $\pi(j) < \pi(j')$.
\medbreak

\begin{Proposition}
    \label{prop:criterion_recognizing_square_permutation}
    A permutation $\pi$ is a square if and only if there exists a
    directed perfect matching on $\pi$ satisfying $\PropDPMShape$ and
    $\PropDPMValue$.
\end{Proposition}
\medbreak

Proposition~\ref{prop:criterion_recognizing_square_permutation} is hence
the analogous of
Proposition~\ref{prop:criterion_recognizing_square_words} for the
supershuffle product and associated square permutations. Moreover,
observe that given a square permutation $\pi$ and a directed perfect
matching $(V, A)$ on $\pi$ satisfying $\PropDPMShape$ and
$\PropDPMValue$, one can recover a square root of $\pi$ by considering
the standardized permutation of the word of sources (or, equivalently,
the word of sinks) of~$(V, A)$.
\medbreak

Figure~\ref{fig:examples_directed_perfect_matchings_square_permutations}
an example related to
Proposition~\ref{prop:criterion_recognizing_square_permutation}.
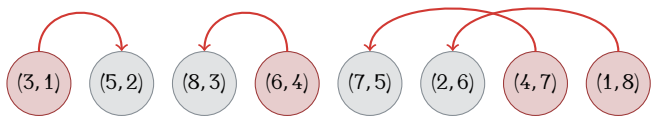
\begin{figure}[ht]
    \centering
    \begin{tikzpicture}
        [xscale=0.35,yscale=.3,node distance=1.1cm,Centering]
        \node[MatchingVertex](U01){\begin{math}(3, 1)\end{math}};
        \node[MatchingVertexColorA,right of=U01](U02)
            {\begin{math}(5, 2)\end{math}};
        \node[MatchingVertexColorA,right of=U02](U03)
            {\begin{math}(8, 3)\end{math}};
        \node[MatchingVertex,right of=U03](U04)
            {\begin{math}(6, 4)\end{math}};
        \node[MatchingVertexColorA,right of=U04](U05)
            {\begin{math}(7, 5)\end{math}};
        \node[MatchingVertexColorA,right of=U05](U06)
            {\begin{math}(2, 6)\end{math}};
        \node[MatchingVertex,right of=U06](U07)
            {\begin{math}(4, 7)\end{math}};
        \node[MatchingVertex,right of=U07](U08)
            {\begin{math}(1, 8)\end{math}};
        \draw[DirectedMatchingEdge]
        (U01.north) .. controls ($ (U01.north) + (0,2.25) $)
            and ($ (U02.north) + (0,2.25) $) .. (U02.north);
        \draw[DirectedMatchingEdge]
        (U04.north) .. controls ($ (U04.north) + (0,2.25) $)
            and ($ (U03.north) + (0,2.25) $) .. (U03.north);
        \draw[DirectedMatchingEdge]
        (U07.north) .. controls ($ (U07.north) + (0,2.5) $)
            and ($ (U05.north) + (0,2.5) $) .. (U05.north);
        \draw[DirectedMatchingEdge]
        (U08.north) .. controls ($ (U08.north) + (0,2.5) $)
            and ($ (U06.north) + (0,2.5) $) .. (U06.north);
    \end{tikzpicture}
    \caption[A directed perfect matching on the square permutation
    $35867241$.]
    {A directed perfect matching on the permutation $\pi := 35867241$
    satisfying $\PropDPMShape$ and $\PropDPMValue$, showing that it is a
    square. It follows also that $\sigma := 2431$ is a square root of
    $\pi$ since $\sigma$ is the standardized of both the word of sources
    $\textcolor{Col1}{3641}$ and the word of sinks
    $\textcolor{Col4}{5872}$ of the directed perfect matching.}
    \label{fig:examples_directed_perfect_matchings_square_permutations}
\end{figure}
\medbreak

\subsubsection{Hardness}
Here is the main algorithmic result of this chapter.
\medbreak

\begin{Theorem} \label{thm:hardness_recognizing_square_permutations}
    In the supershuffle algebra $(\FQSym, \SuperShuffle)$, $\ProblemSD$
    is $\NPC$.
\end{Theorem}
\medbreak

Recall that the \Def{pattern involvement problem} consists, given two
permutations $\pi$ and $\sigma$, in deciding if $\pi$ admits an
occurrence of $\sigma$. This problem is known to be $\NPC$~\cite{BBL98}.
Theorem~\ref{thm:hardness_recognizing_square_permutations} can be shown
by performing a polynomial-time reduction from the pattern involvement
problem to~$\ProblemSD$.
\medbreak

\section*{Concluding remarks}
There are a number of further directions of investigation in this
general subject. They cover several areas: algorithmic, combinatorics,
and algebra. Let us mention several ---not necessarily all new--- open
problems that are, in our opinion, the most interesting.
\smallbreak

First ones are enumerative questions. We have computed few first terms
of some integer sequences, like~\eqref{equ:sequence_square_permutations}
for the number of square permutations,
\eqref{equ:sequence_square_classes} for the number of square
permutations quotiented by their natural symmetries,
or~\eqref{equ:sequence_square_binary_words} for the number of square
permutations avoiding the patterns $213$ and $231$ (equivalently, by
Proposition~\ref{prop:bijection_binary_to_permutations_squares}, this is
also the number of square binary words~\cite{HRS12}). We can ask about
formulas to compute these numbers.
\smallbreak

Second ones are algorithmic questions. One can first ask the difficulty
of deciding whether a permutation avoiding $213$ and $231$ is a square
(see~\cite{HRS12,RV13,BS14} for the point of view of square binary
words). Besides, one can ask about the hardness of $\ProblemRC$ in the
context of the supershuffle. In other terms, the problem consists, given
two permutations $\pi$ and $\sigma$, in deciding if $\sigma$ is a square
root of~$\pi$.
\smallbreak

Finally, in a more algebraic flavor, we can ask about the properties of
the associative algebra $(\FQSym, \SuperShuffle)$, continuing the work
of Vargas~\cite{Var14}. This includes, among others, the description of
a generating family, the definitions of multiplicative bases, and
determining whether this algebra is free as an associative algebra.
\medbreak


\chapter{Bud generating systems} \label{chap:buds}
The content of this chapter comes from~\cite{Gir16d}.
\medbreak

\section*{Introduction}
Coming from theoretical computer science and formal language theory,
formal grammars~\cite{Har78,HMU06} are powerful tools having many
applications in several fields of mathematics. A formal grammar is a
device which describes---more or less concisely and with more or less
restrictions---a set of words, called language. There are several
variations in the definitions of formal grammars and some sorts of them
are classified by the Chomsky-Schützenberger hierarchy~\cite{Cho59,CS63}
according to four different categories, taking into account their
expressive power. In an increasing order of power, there is the class of
Type-3 grammars known as regular grammars, the class of Type-2 grammars
known as context-free grammars, the class of Type-1 grammars known as
context-sensitive grammars, and the class of Type-0 grammars known as
unrestricted grammars. One of the most striking similarities between all
these variations of formal grammars is that they work by constructing
words by applying rewrite rules~\cite{BN98} (see also
Section~\ref{subsec:rewrite_rules} of Chapter~\ref{chap:combinatorics}).
Indeed, a word of the language described by a formal grammar is obtained
by considering a starting word and by iteratively altering some of its
factors in accordance with the production rules of the grammar.
\smallbreak

Similar mechanisms and ideas are translatable into the world of trees,
instead only of those of words. Grammars of trees~\cite{CDGJLLTT07} are
hence the natural counterpart of formal grammars to describe sets of
trees, and here also, there exist many very different types of grammars.
One can cite for instance tree grammars, regular tree
grammars~\cite{GS84}, and synchronous grammars~\cite{Gir12e}, which are
devices providing a way to describe sets of various kinds of treelike
structures. Here also, one of the common points between these grammars
is that they work by applying rewrite rules on trees. In this framework,
trees are constructed by growing from the root to the leaves by
replacing some subtrees by other ones. Like free monoids are algebraic
structures involving words, free ns operads are algebraic structures
involving planar rooted trees. Since monoids are the underlying
structures for most of the generating systems on words, it is natural to
ask whether ns operads can be thought as underlying structures of
generating systems on trees.
\smallbreak

The initial spark of this work has been caused by the following simple
observation. The partial composition $x \circ_i y$ of two elements $x$
and $y$ of a ns operad $\Oca$ can be regarded as the application of a
rewrite rule on $x$ to obtain a new element of $\Oca$---the rewrite rule
being encoded essentially by $y$. This leads to the idea consisting in
considering a ns operad $\Oca$ to define grammars generating some
subsets of $\Oca$. In this way, according to the nature of the elements
of $\Oca$, this provides a way to define grammars which generate objects
different than words (as in the case of formal grammars) and than trees
(as in the case of grammars of trees). We rely in this work on ns
colored operads (see Section~\ref{subsubsec:colored_operads} of
Chapter~\ref{chap:algebra}). Ns colored operads are the suitable devices
to our aim of defining a new kind of grammars since the restrictions
provided by the colors allow a precise control on how the rewrite rules
can be applied.
\smallbreak

Thus, we introduce in this work a new kind of grammars, the bud
generating systems. They are defined mainly from a ground ns operad
$\Oca$, a set $\CFr$ of colors, and a set $\RFr$ of production rules. A
bud generating system describes a subset of $\Bud_\CFr(\Oca)$---the ns
colored operad obtained by augmenting the elements of $\Oca$ with input
and output colors taken from $\CFr$. The generation of an element works
by iteratively altering an element $x$ of $\Bud_\CFr(\Oca)$ by composing
it, if possible, with an element $y$ of $\RFr$. In this context, the
colors play the role analogous of the one of nonterminal symbols in the
formal grammars and in the grammars of trees. Any bud generating system
$\Bca$ specifies two sets of objects: its language $\Lang(\Bca)$ and its
synchronous language $\SyncLang(\Bca)$. Thereby, bud generating systems
can be used to describe sets of combinatorial objects. For instance,
they can be used to describe sets of Motkzin paths with some
constraints, sets of Schröder trees with some constraints, the set of
$\{2, 3\}$-perfect trees~\cite{MPRS79,CLRS09} and some of its
generalizations, and the set of balanced binary trees~\cite{AVL62}. One
remarkable fact is that bud generating systems can emulate both
context-free grammars and regular tree grammars, and allow to see both
of these in a unified manner. In the first case, context-free grammars
are emulated by bud generating systems with the associative operad $\As$
as ground ns operad and in the second case, regular tree grammars are
emulated by bud generating systems with a free ns operad
$\FreeOperad(\GeneratingSet)$ as ground ns operad, where
$\GeneratingSet$ is a precise set of generators.
\smallbreak

A very normal combinatorial question consists, given a bud generating
system $\Bca$, in computing the generating series
$\GenS_{\Lang(\Bca)}(t)$ and $\GenS_{\SyncLang(\Bca)}(t)$, respectively
counting the elements of the language and of the synchronous language of
$\Bca$ with respect to the arity of the elements. To achieve this
objective, we consider a new generalization of formal power series,
namely series on ns colored operads. Series on ns operads and operads
satisfying some precise properties have been
considered~\cite{Cha02,Cha08,Cha09} (see also~\cite{Vdl04,Fra08,LN13}).
In this work, we consider series on ns colored operads which are, in
some sense, generalizations of these notions of series. Any bud
generating system $\Bca$ leads to the definition of three series on ns
colored operads: its {hook generating series $\Hook(\Bca)$, its
syntactic generating series $\Synt(\Bca)$, and its synchronous
generating series $\Sync(\Bca)$. The hook generating series allows to
define analogues of the hook-length statistics of binary
trees~\cite{Knu98} for objects belonging to the language of $\Bca$,
possibly different than trees. The syntactic (resp. synchronous)
generating series bring functional equations and recurrence
formulas to compute the coefficients of $\GenS_{\Lang(\Bca)}(t)$
and~$\GenS_{\SyncLang(\Bca)}(t)$. The definitions of these three series
rely on particular operations on series on ns colored operads: a pre-Lie
product, an associative product, and their respective Kleene stars.
\smallbreak

This chapter is organized as follows.
Section~\ref{sec:series_colored_operads} begins by introduction
the construction $\Bud_\CFr$. Then, we provide elementary definitions
about series on ns colored operads, and define a pre-Lie product and an
associative product on these series. Next,
Section~\ref{sec:bud_generating_systems} is concerned with the
definition of bud generating systems and to study their first
properties. This chapter ends with
Section~\ref{sec:series_bud_generating_systems} wherein we use the
definitions and the results of both previous sections to consider bud
generating systems as devices to define statistics on combinatorial
objects or to enumerate families of combinatorial objects.
\medbreak

\subsubsection*{Note}
This chapter deals only with ns set-operads and ns colored set-operads.
For these reasons, ``operad'' means ``ns set-operad'' and ``colored
operad'' means ``ns colored set-operad''.
\medbreak

\section{Colored operads and formal power series}
\label{sec:series_colored_operads}
The purpose of this section is twofold. First, we present a very
natural construction $\Bud_\CFr$ taking as input a monochrome operad
$\Oca$ and outputting a colored operad by augmenting the outputs and
inputs of the elements of $\Oca$ with colors of $\CFr$. This construction
will be used to define bud generating systems in the next section.
Second, we consider series on colored operads and define two products
on them. These products are considered in the last section of this
chapter for enumerative goals.
\medbreak

\subsection{Bud operads}
Let us first present a simple construction producing colored operads
from operads.
\medbreak

\subsubsection{Sets of colors}
In all this chapter, we consider that $\CFr$ has cardinal $k \geq 1$ and
that the colors of $\CFr$ are arbitrarily indexed so that
$\CFr = \{c_1, \dots, c_k\}$.
\medbreak

\subsubsection{From monochrome operads to colored operads}
\label{subsubsec:operad_to_colored_operads}
If $\Oca$ is a monochrome operad and $\CFr$ is a finite set of colors,
we denote by $\Bud_\CFr(\Oca)$ the $\CFr$-colored collection (see
Section~\ref{subsubsec:colored_collections} of
Chapter~\ref{chap:combinatorics}) by
\begin{equation}
    \Bud_\CFr(\Oca)(n) := \CFr \times \Oca(n) \times \CFr^n,
    \qquad n \geq 1,
\end{equation}
and for all $(a, x, u) \in \Bud_\CFr(\Oca)$, $\Out((a, x, u)) := a$ and
$\In((a, x, u)) := u$. We endow $\Bud_\CFr(\Oca)$ with the partially
defined partial composition $\circ_i$ satisfying, for all triples
$(a, x, u)$ and $(b, y, v)$ of $\Bud_\CFr(\Oca)$ and $i \in [|x|]$ such
that $ \Out((b, y, v)) = \In_i((a, x, u))$,
\begin{equation}
    (a, x, u) \circ_i (b, y, v) := (a, x \circ_i y, u \mapsfrom_i v),
\end{equation}
where $u \mapsfrom_i v$ is the word obtained by replacing the $i$th
letter of $u$ by $v$. Besides, if $\Oca_1$ and $\Oca_2$ are two operads
and $\phi : \Oca_1 \to \Oca_2$ is an operad morphism, we denote by
$\Bud_\CFr(\phi)$ the map
\begin{equation}
    \Bud_\CFr(\phi) : \Bud_\CFr(\Oca_1) \to \Bud_\CFr(\Oca_2)
\end{equation}
defined by
\begin{equation}
    \Bud_\CFr(\phi)((a, x, u)) := (a, \phi(x), u).
\end{equation}
\medbreak

\begin{Proposition} \label{prop:functor_bud_operads}
    For any set of colors $\CFr$, the construction
    \begin{math}
        (\Oca, \phi) \mapsto
        \left(\Bud_\CFr(\Oca), \Bud_\CFr(\phi)\right)
    \end{math}
    is a functor from the category of monochrome operads to the category
    of $\CFr$-colored operads.
\end{Proposition}
\medbreak

Proposition~\ref{prop:functor_bud_operads} shows that $\Bud_\CFr$ is a
functorial construction producing colored operads from monochrome ones.
We call $\Bud_\CFr(\Oca)$ the \Def{$\CFr$-bud operad} of~$\Oca$.
\medbreak

When $\CFr$ is a singleton, $\Bud_\CFr(\Oca)$ is by definition a
monochrome operad isomorphic to~$\Oca$. For this reason, in this case,
we identify $\Bud_\CFr(\Oca)$ with $\Oca$.
\medbreak

As a side observation, remark that in general, the bud operad
$\Bud_\CFr(\Oca)$ of a free operad $\Oca$ is not a free $\CFr$-colored
operad. Indeed, consider for instance the bud operad
$\Bud_{\{1, 2\}}(\Oca)$, where $\Oca := \FreeOperad(C)$ and $C$ is the
monochrome collection defined by $C := C(1) := \{\Asf\}$. Then, a
minimal generating set of $\Bud_{\{1, 2\}}(\Oca)$ is
\begin{equation}
    \left\{
    \left(1, \FreeCorollaOne{\Asf}, 1\right),
    \left(1, \FreeCorollaOne{\Asf}, 2\right),
    \left(2, \FreeCorollaOne{\Asf}, 1\right),
    \left(2, \FreeCorollaOne{\Asf}, 2\right)
    \right\}.
\end{equation}
These elements are subjected to the nontrivial relations
\begin{equation}
    \left(d, \FreeCorollaOne{\Asf}, 1\right)
    \circ_1
    \left(1, \FreeCorollaOne{\Asf}, e\right)
    =
    \left(d,
    \begin{tikzpicture}[xscale=.13,yscale=.15,Centering]
        \node(0)at(1.00,-4.5){};
        \node[NodeST](1)at(1.00,0.00){\begin{math}\Asf\end{math}};
        \node[NodeST](2)at(1.00,-2){\begin{math}\Asf\end{math}};
        \draw[Edge](1)--(2);
        \draw[Edge](0)--(2);
        \node(r)at(1.00,2.5){};
        \draw[Edge](r)--(1);
    \end{tikzpicture}\,,
    e\right)
    =
    \left(d, \FreeCorollaOne{\Asf}, 2\right)
    \circ_1
    \left(2, \FreeCorollaOne{\Asf}, e\right),
\end{equation}
where $d, e \in \{1, 2\}$, implying that $\Bud_{\{1, 2\}}(\Oca)$ is not
free.
\medbreak

\subsubsection{Bud operad of the associative operad}
Let us consider the $\CFr$-bud operads of the associative operad $\As$
(see Section~\ref{subsubsec:associative_operad} of
Chapter~\ref{chap:algebra}). For any set of colors $\CFr$, the bud
operad $\Bud_\CFr(\As)$ is the set of all triples
\begin{equation}
    (a, \Afr_n, u_1 \dots u_n)
\end{equation}
where $a \in \CFr$ and $u_1, \dots, u_n \in \CFr$. For
$\CFr := \{1, 2, 3\}$, one has for instance the partial composition
\begin{equation}
    \left(2, \Afr_4, \textcolor{Col1}{3{\bf 1}12}\right)
    \circ_2
    \left(1, \Afr_3, \textcolor{Col4}{233}\right) =
    \left(2, \Afr_6, \textcolor{Col1}{3}\textcolor{Col4}{233}
        \textcolor{Col1}{12}\right).
\end{equation}
\medbreak

The associative operad and its bud operads will play an important role
in the sequel. For this reason, to gain readability, we shall simply
denote by $(a, u)$ any element $\left(a, \Afr_{|u|}, u\right)$ of
$\Bud_\CFr(\As)$ without any loss of information.
\medbreak

\subsubsection{Pruning map}
Here and in the sequel, we use the fact that any monochrome operad
$\Oca$ can be seen as a $\CFr$-colored operad where all output and input
colors of its elements are equal to $c_1$, where $c_1$ is the first
color of $\CFr$. Let
\begin{equation}
    \Prune : \Bud_\CFr(\Oca) \to \Oca
\end{equation}
be the map defined, for any $(a, x, u) \in \Bud_\CFr(\Oca)$, by
\begin{equation}
    \Prune((a, x, u)) := x.
\end{equation}
We call $\Prune$ the \Def{pruning map} on $\Bud_\CFr(\Oca)$. Observe
that $\Prune$ is not a morphism of $\CFr$-colored operads since it is
not a $\CFr$-colored collection morphism.
\medbreak

\subsection{The space of series on colored operads}
We work here with series on colored operads. We explain how to encode
usual noncommutative multivariate series and series on monoids by series
on colored operads.
\medbreak

\subsubsection{First definitions}
For any $\CFr$-colored operad $\Cca$, a \Def{$\Cca$-series} is a series
on $\Cca$ seen as a collection (see Section~\ref{subsubsec:C_series} of
Chapter~\ref{chap:algebra}). In other terms, a $\Cca$-series is an
element of~$\K \AAngle{\Cca}$.
\medbreak

For any combinatorial $\CFr$-colored collection $C$, we denote by
$\Sbf_C$ the generating series of $C$, seen as a graded collection.
\medbreak

The $\Cca$-series $\Ubf$ defined by
\begin{equation} \label{equ:definition_series_u}
    \Ubf := \sum_{a \in \CFr} \Unit_a
\end{equation}
is the \Def{series of colored units} of $\Cca$ and will play a special
role in the sequel. Since $\CFr$ is finite, $\Ubf$ is a polynomial.
\medbreak

Observe that $\Cca$-series are defined here on fields $\K$ instead of on
the much more general structures of semirings, as it is the case for
series on monoids~\cite{Sak09}. We choose to tolerate this loss of
generality because this considerably simplifies the theory. Furthermore,
we shall use in the sequel $\Cca$-series as devices for combinatorial
enumeration, so that it is sufficient to pick $\K$ as the field
$\Q(q_0, q_1, q_2, \dots)$ of rational functions in an infinite number
of commuting parameters with rational coefficients. The parameters
$q_i$, $i \in \N$ intervene in the enumeration of colored subcollections
of $\Cca$ with respect to several statistics.
\medbreak

\subsubsection{Functorial construction}
If $\Cca_1$ and $\Cca_2$ are two $\CFr$-colored operads and
$\phi : \Cca_1 \to \Cca_2$ is a morphism of colored operads,
$\K \AAngle{\phi}$ is the map
\begin{equation}
    \K \AAngle{\phi} : \K \AAngle{\Cca_1} \to \K \AAngle{\Cca_2}
\end{equation}
defined, for any $\Fbf \in \K \AAngle{\Cca_1}$ and $y \in \Cca_2$, by
\begin{equation} \label{equ:series_morphism}
    \Angle{y, \K \AAngle{\phi}(\Fbf)} :=
    \sum_{\substack{x \in \Cca_1 \\ \phi(x) = y}} \Angle{x, \Fbf}.
\end{equation}
Equivalently, $\K \AAngle{\phi}$ can be defined, by using the sum
notation of series (see Section~\ref{subsubsec:polynomial_spaces} of
Chapter~\ref{chap:algebra}), by
\begin{equation} \label{equ:series_morphism_extended}
    \K \AAngle{\phi}(\Fbf) :=
    \sum_{x \in \Cca_1} \Angle{x, \Fbf} \phi(x).
\end{equation}
\medbreak

Observe first that $\K \AAngle{\Fbf}$ is a linear map. Moreover, notice
that~\eqref{equ:series_morphism} could be undefined for arbitrary
colored operads $\Cca_1$ and $\Cca_2$, and an arbitrary morphism of
colored operads $\phi$. However, when all fibers of $\phi$ are finite,
for any $y \in \Cca_2$, the right member of~\eqref{equ:series_morphism}
is well-defined since the sum has a finite number of terms. Moreover,
since any morphism from a combinatorial colored operad has finite
fibers, one has the following result.
\medbreak

\begin{Proposition} \label{prop:functor_colored_operads_vector_spaces}
    The construction
    $(\Cca, \phi) \mapsto (\K \AAngle{\Cca}, \K \AAngle{\phi})$ is a
    functor from the category of combinatorial $\CFr$-colored operads to
    the category of $\K$-vector spaces.
\end{Proposition}
\medbreak

\subsubsection{Noncommutative multivariate series}
\label{subsubsec:noncommutative_multivariate_series}
For any finite alphabet $A$ of noncommutative letters, recall that
$\K \AAngle{A^*}$ is the set of noncommutative series on
$A$~\cite{Eil74,SS78,BR10}.
\medbreak

Let us explain how to encode any series $\Sbf \in \K \AAngle{A^*}$ by a
series on a particular colored operad. Let $\CFr_A$ be the set of colors
$A \sqcup \{\lozenge\}$ where $\Dummy$ is a virtual letter which is not
in $A$, and $\Cca_A$ be the $\CFr_A$-colored subcollection of
$\Bud_{\CFr_A}(\As)$ consisting in arity one in the colored units of
$\Bud_{\CFr_A}(\As)$ and in arity $n \geq 2$ in the elements of the form
\begin{equation} \label{equ:elements_operad_multivariate_series}
    \left(\Dummy, a_1 \dots a_{n - 1} \Dummy\right),
    \qquad a_1 \dots a_{n - 1} \in A^{n - 1}.
\end{equation}
Since the partial composition of any two elements of the
form~\eqref{equ:elements_operad_multivariate_series} is in $\Cca_A$,
$\Cca_A$ is a colored suboperad of $\Bud_{\CFr_A}(\As)$. Then, any
series
\begin{equation}
    \Sbf := \sum_{u \in A^*} \Angle{u, \Sbf} u
\end{equation}
of $\K \AAngle{A^*}$ is encoded by the series
\begin{equation} \label{equ:colored_series_from_multivariate_series}
    \Multi(\Sbf) := \sum_{u \in A^*}
    \Angle{u, \Sbf} \left(\Dummy, u \Dummy\right).
\end{equation}
of $\K \AAngle{\Cca_A}$. We shall explain a little further how the usual
noncommutative product of series of $\K \AAngle{A^*}$ can be translated
on $\Bud_{\CFr_A}(\As)$-series of the
form~\eqref{equ:colored_series_from_multivariate_series}.
\medbreak

\subsubsection{Series on monoids} \label{subsubsec:series_on_monoids}
For any monoid $\Mca$, recall that $\K \AAngle{\Mca}$ if the set of
noncommutative series on $\Mca$~\cite{Sak09}. Noncommutative
multivariate series are particular cases of series on monoids since any
noncommutative multivariate series of $\K \AAngle{A^*}$ can be seen as
an $A^*$-series, where $A^*$ is the free monoid on~$A$.
\medbreak

Let us explain how to encode any series $\Sbf \in \K \AAngle{\Mca}$ by a
series on a particular colored operad. Let $\Oca_\Mca$ be the
monochrome collection concentrated in arity one where
$\Oca_{\Mca}(1) := \Mca$. We define the map
$\circ_1 : \Oca_\Mca(1) \times \Oca_\Mca(1) \to \Oca_\Mca(1)$ for all
$x, y \in \Oca_\Mca$ by $x \circ_1 y := x \Product y$ where $\Product$
is the operation of $\Mca$. Since $\Product$ is associative and admits a
unit, $\circ_1$ satisfy all relations of operads, so that $\Oca_\Mca$ is
a monochrome operad. Then, any series
\begin{equation}
    \Sbf := \sum_{u \in \Mca} \Angle{u, \Sbf} u
\end{equation}
of $\K \AAngle{\Mca}$ is encoded by the series
\begin{equation} \label{equ:colored_series_from_series_on_monoids}
    \Mono(\Sbf) := \sum_{u \in \Oca_\Mca} \Angle{u, \Sbf} u
\end{equation}
of $\K \AAngle{\Oca_\Mca}$. We shall explain a little further how the
usual product of series of $\K \AAngle{\Mca}$, called Cauchy product
in~\cite{Sak09}, can be translated on series of the
form~\eqref{equ:colored_series_from_series_on_monoids}.
\medbreak

Observe that when $\Mca$ is a free monoid $A^*$ where $A$ is a finite
alphabet of noncommutative letters, we then have two ways to encode a
series $\Sbf$ of $\K \AAngle{A^*}$. Indeed, $\Sbf$ can be encoded as the
series $\Multi(\Sbf)$ of $\K \AAngle{\Cca_A}$ of the
form~\eqref{equ:colored_series_from_multivariate_series}, or as the
series $\Mono(\Sbf)$ of $\K \AAngle{\Oca_{\Abb^*}}$ of the
form~\eqref{equ:colored_series_from_series_on_monoids}. Remark that
the first way to encode $\Sbf$ is preferable since $\Cca_A$ is a
combinatorial operad while $\Oca_{A^*}$ is not.
\medbreak

\subsubsection{Series of colors} \label{subsubsec:series_of_colors}
Let
\begin{equation}
    \Colors : \Cca \to \Bud_\CFr(\As)
\end{equation}
be the morphism of colored operads defined for any $x \in \Cca$ by
\begin{equation}
    \Colors(x) := \left(\Out(x), \In(x)\right).
\end{equation}
By a slight abuse of notation, we denote by
\begin{equation}
    \Colors : \K \AAngle{\Cca} \to \K \AAngle{\Bud_\CFr(\As)}
\end{equation}
the map sending any series $\Fbf$ of $\K \AAngle{\Cca}$ to
$\K \AAngle{\Colors}(\Fbf)$, called \Def{series of colors} of $\Fbf$.
By~\eqref{equ:series_morphism_extended},
\begin{equation} \label{equ:series_of_colors}
    \Colors(\Fbf) =
    \sum_{x \in \Cca} \Angle{x, \Fbf} (\Out(x), \In(x)).
\end{equation}
\medbreak

Intuitively, the series $\Colors(\Fbf)$ can be seen as a version of
$\Fbf$ wherein only the colors of the elements of its support are taken
into account.
\medbreak

\subsubsection{Series of color types}
\label{subsubsec:series_of_color_types}
The \Def{$\CFr$-type} of a word $u \in \CFr^+$ is the word $\Type(u)$ of
$\N^k$ defined by
\begin{equation}
    \Type(u) := |u|_{c_1} \dots |u|_{c_k}.
\end{equation}
By extension, we shall call \Def{$\CFr$-type} any word of $\N^k$ with at
least a nonzero letter and we denote by $\Types_\CFr$ the set of all
$\CFr$-types. The \Def{degree} $\deg(\alpha)$ of
$\alpha \in \Types_\CFr$ is the sum of the letters of $\alpha$. We
denote by $\CFr^\alpha$ the word
\begin{equation}
    \CFr^\alpha := c_1^{\alpha_1} \dots c_k^{\alpha_k}.
\end{equation}
\medbreak

Let $\Zbb_\CFr := \{\VarZ_{c_1}, \dots, \VarZ_{c_k}\}$ be an alphabet
of commutative letters. For any type $\alpha$, we denote by
$\Zbb_\CFr^\alpha$ the monomial
\begin{equation}
    \Zbb_\CFr^\alpha :=
    \VarZ_{c_1}^{\alpha_1} \dots \VarZ_{c_k}^{\alpha_k}
\end{equation}
of $\K \Angle{\Multiset(\Zbb_\CFr)}$.
Moreover, for any two types $\alpha$ and $\beta$, the \Def{sum}
$\alpha \dot{+} \beta$ of $\alpha$ and $\beta$ is the type satisfying
$(\alpha \dot{+} \beta)(i) := \alpha(i) + \beta(i)$ for all $i \in [k]$.
Observe that with this notation,
\begin{math}
    \Zbb_\CFr^\alpha \Zbb_\CFr^\beta = \Zbb_\CFr^{\alpha \dot{+} \beta}
\end{math}.
\medbreak

Now, set $\Xbb_\CFr$ and $\Ybb_\CFr$ respectively as the two alphabets
of commutative letters $\{\VarX_{c_1}, \dots, \VarX_{c_k}\}$
and $\{\VarY_{c_1}, \dots, \VarY_{c_k}\}$. We can see these two
alphabets as graded collections where each letter is of size $1$.
Consider the map
\begin{equation}
    \ColorTypes :
    \K \AAngle{\Cca} \to
    \K \AAngle{\Multiset\left(\Xbb_\CFr + \Ybb_\CFr\right)},
\end{equation}
defined for all $\alpha, \beta \in \Tca_\CFr$ by
\begin{equation} \label{equ:definition_color_types}
    \Angle{\Xbb_\CFr^\alpha \Ybb_\CFr^\beta, \ColorTypes(\Fbf)}
    :=
    \sum_{\substack{
        (a, u) \in \Bud_\CFr(\As) \\
        \Type(a) = \alpha \\
        \Type(u) =  \beta
    }}
    \Angle{(a, u), \Colors(\Fbf)}.
\end{equation}
By the definition of the map $\Colors$,
\begin{equation}
    \ColorTypes(\Fbf) =
    \sum_{x \in \Cca}
    \Angle{x, \Fbf}\enspace
    \Xbb_\CFr^{\Type(\Out(x))}\enspace
    \Ybb_\CFr^{\Type(\In(x))}.
\end{equation}
\medbreak

Observe that for all $\alpha, \beta \in \Types_\CFr$ such that
$\deg(\alpha) \ne 1$, the coefficients of
$\Xbb_\CFr^\alpha \Ybb_\CFr^\beta$ in $\ColorTypes(\Fbf)$ are zero. In
intuitive terms, the series $\ColorTypes(\Fbf)$, called
\Def{series of color types} of $\Fbf$, can be seen as a version of
$\Colors(\Fbf)$ wherein only the output colors and the types of the
input colors of the elements of its support are taken into account, the
variables of $\Xbb_\CFr$ encoding output colors and the variables of
$\Ybb_\CFr$ encoding input colors. In the sequel, we are concerned by
the computation of the coefficients of $\ColorTypes(\Fbf)$ for some
$\Cca$-series~$\Fbf$.
\medbreak

\subsubsection{Pruned series} \label{subsubsec:pruned_series}
Let $\Oca$ be a monochrome operad, $\Bud_\CFr(\Oca)$ be a bud operad,
and $\Fbf$ be a $\Bud_\CFr(\Oca)$-series. Since $\CFr$ is finite, the
series $\K \AAngle{\Prune}(\Fbf)$ is well-defined and, by a slight abuse
of notation, we denote by
\begin{equation}
    \Prune : \K \AAngle{\Bud_\CFr(\Oca)} \to \K \AAngle{\Oca}
\end{equation}
the map sending any series $\Fbf$ of $\K \AAngle{\Bud_\CFr(\Oca)}$ to
$\K \AAngle{\Prune}(\Fbf)$, called \Def{pruned series} of $\Fbf$.
By~\eqref{equ:series_morphism_extended},
\begin{equation} \label{equ:pruned_series}
    \Prune(\Fbf) =
    \sum_{(a, x, u) \in \Bud_\CFr(\Oca)} \Angle{(a, x, u), \Fbf} x.
\end{equation}
\medbreak

Intuitively, the series $\Prune(\Fbf)$ can be seen as a version of
$\Fbf$ wherein the colors of the elements of its support are forgotten.
Besides, $\Fbf$ is said \Def{faithful} if all coefficients of
$\Prune(\Fbf)$ are equal to $0$ or to~$1$.
\medbreak

\subsubsection{Example: series of trees}
\label{subsubsec:example_series_trees}
Let $\Cca$ be the free $\CFr$-colored operad over $C$ where
$\CFr := \{1, 2\}$ and $C$ is the $\CFr$-colored collection defined by
$C := C(2) \sqcup C(3)$ with $C(2) := \{\Asf\}$, $C(3) := \{\Bsf\}$,
$\Out(\Asf) := 1$, $\Out(\Bsf) := 2$, $\In(\Asf) := 21$, and
$\In(\Bsf) := 121$. Let $\Fbf_\Asf$ (resp. $\Fbf_\Bsf$) be the series of
$\K \AAngle{\Cca}$ where for any syntax tree $\Tfr$ of $\Cca$,
$\Angle{\Tfr, \Fbf_\Asf}$ (resp. $\Angle{\Tfr, \Fbf_\Bsf}$) is the
number of internal nodes of $\Tfr$ labeled by $\Asf$ (resp. $\Bsf$). The
series $\Fbf_\Asf$ and $\Fbf_\Bsf$ are of the form
\begin{subequations}
\begin{multline}
    \Fbf_\Asf =

    \enspace + \cdots.
\end{equation}
\end{subequations}
The sum $\Fbf_\Asf + \Fbf_\Bsf$ is the series wherein the coefficient of
any syntax tree $\Tfr$ of $\Cca$ is its degree. Let also $\Fbf_{|_1}$
(resp. $\Fbf_{|_2}$) be the series of $\K \AAngle{\Cca}$ where for any
syntax tree $\Tfr$ of $\Cca$, $\Angle{\Tfr, \Fbf_{|_1}}$ (resp.
$\Angle{\Tfr, \Fbf_{|_2}}$) is the number of inputs colors $1$ (resp.
$2$) of $\Tfr$. The sum $\Fbf_{|_1} + \Fbf_{|_2}$ is the series wherein
the coefficient of any syntax tree $\Tfr$ of $\Cca$ is its arity.
Moreover, the series $\Fbf_\Asf + \Fbf_\Bsf + \Fbf_{|_1} + \Fbf_{|_2}$
is the series wherein the coefficient of any syntax tree $\Tfr$ of
$\Cca$ is its total number of nodes.
\medbreak

The series of colors of $\Fbf_\Asf$ is of the form
\begin{equation}
    \Colors(\Fbf_\Asf) = (1, 21) + 2\, (1, 221) + 3\, (1, 2221)
    + (2, 2121) + (2, 1221) + (1, 1211) + \cdots,
\end{equation}
and the series of color types of $\Fbf$ is of the form
\begin{equation}
    \ColorTypes(\Fbf_\Asf) =
    \VarX_1\VarY_1\VarY_2 + 2\VarX_1\VarY_1\VarY_2^2
    + \VarX_1\VarY_1^3\VarY_2
    + 3\VarX_1\VarY_1\VarY_2^3
    + 2\VarX_2\VarY_1^2\VarY_2^2 + \cdots.
\end{equation}
\medbreak

\subsection{Pre-Lie product on series}
We are now in position to define a binary operation $\PreLieProduct$ on
the space of the $\Cca$-series. As we shall see, this operation is
partially defined, nonunitary, noncommutative, and nonassociative.
\medbreak

\subsubsection{Pre-Lie product}
Given two $\Cca$-series $\Fbf, \Gbf \in \K \AAngle{\Cca}$, the
\Def{pre-Lie product} of $\Fbf$ and $\Gbf$ is the $\Cca$-series
$\Fbf \PreLieProduct \Gbf$ defined, for any $x \in \Cca$, by
\begin{equation} \label{equ:pre_lie_product}
    \Angle{x, \Fbf \PreLieProduct \Gbf}
    := \sum_{\substack{
            y, z \in \Cca \\
            i \in [|y|] \\
            x = y \circ_i z}}
    \Angle{y, \Fbf} \Angle{z, \Gbf}.
\end{equation}
Observe that $\Fbf \PreLieProduct \Gbf$ could be undefined for arbitrary
$\Cca$-series $\Fbf$ and $\Gbf$ on an arbitrary colored operad $\Cca$.
Besides, notice from~\eqref{equ:pre_lie_product} that $\PreLieProduct$
is bilinear and that $\Ubf$ (defined in~\eqref{equ:definition_series_u})
is a left unit of~$\PreLieProduct$. However, since
\begin{equation}
    \Fbf \PreLieProduct \Ubf =
    \sum_{x \in \Cca} |x| \Angle{x, \Fbf} x,
\end{equation}
the $\Cca$-series $\Ubf$ is not a right unit of $\PreLieProduct$. This
product is also nonassociative in the general case since we have, for
instance in~$\K \AAngle{\As}$,
\begin{equation}
    (\Afr_2 \PreLieProduct \Afr_2) \PreLieProduct \Afr_2 =
    6 \Afr_4
    \ne
    4 \Afr_4 =
    \Afr_2 \PreLieProduct (\Afr_2 \PreLieProduct \Afr_2).
\end{equation}
Nevertheless, it satisfies the pre-Lie relation
(see~\eqref{equ:relation_pre_Lie} of
Section~\ref{subsubsec:pre_lie_algebras} of Chapter~\ref{chap:algebra}).
\medbreak

\begin{Proposition} \label{prop:functor_series_pre_lie_algebras}
    The construction
    \begin{math}
        (\Cca, \phi) \mapsto
        ((\K \AAngle{\Cca}, \PreLieProduct), \K \AAngle{\phi})
    \end{math}
    is a functor from the category of combinatorial $\CFr$-colored
    operads to the category of pre-Lie algebras.
\end{Proposition}
\begin{proof}
    Let $\Cca$ be a combinatorial $\CFr$-colored operad. First of all,
    since $\Cca$ is combinatorial, the pre-Lie product of any two
    $\Cca$-series $\Fbf$ and $\Gbf$ is well-defined due to the fact that
    the sum~\eqref{equ:pre_lie_product} has a finite number of terms.
    Let $\Fbf$, $\Gbf$, and $\Hbf$ be three $\Cca$-series and
    $x \in \Cca$. We denote by $\lambda(\Fbf, \Gbf, \Hbf)$ the
    coefficient of $x$ in
    \begin{math}
        (\Fbf \PreLieProduct \Gbf) \PreLieProduct \Hbf -
        \Fbf \PreLieProduct (\Gbf \PreLieProduct \Hbf)
    \end{math}.
    We have
    \begin{equation}\begin{split}
        \label{equ:proof_functor_series_pre_lie_algebra}
        \lambda(\Fbf, \Gbf, \Hbf) & =
        \sum_{\substack{
            y, z, t \in \Cca \\
            i, j \in \N \\
            x = (y \circ_i z) \circ_j t}}
        \Angle{y, \Fbf} \Angle{z, \Gbf} \Angle{t, \Hbf}
        -
        \sum_{\substack{
            y, z, t \in \Cca \\
            i, j \in \N \\
            x = y \circ_i (z \circ_j t)}}
        \Angle{y, \Fbf} \Angle{z, \Gbf} \Angle{t, \Hbf} \\
        & =
        \sum_{\substack{
            y, z, t \in \Cca \\
            i > j \in \N \\
            x = (y \circ_i z) \circ_j t}}
        \Angle{y, \Fbf} \Angle{z, \Gbf} \Angle{t, \Hbf} \\
        & =
        \sum_{\substack{
            y, z, t \in \Cca \\
            i > j \in \N \\
            x = (y \circ_i t) \circ_j z}}
        \Angle{y, \Fbf} \Angle{z, \Gbf} \Angle{t, \Hbf} \\
        & = \lambda(\Fbf, \Hbf, \Gbf).
    \end{split}\end{equation}
    The second and the last equality
    of~\eqref{equ:proof_functor_series_pre_lie_algebra} come from
    Relation~\eqref{equ:operad_axiom_1} of
    Section~\ref{subsubsec:ns_operads} of
    Chapter~\ref{chap:algebra} of operads and the third equality is a
    consequence of Relation~\eqref{equ:operad_axiom_2} of
    Section~\ref{subsubsec:ns_operads} of Chapter~\ref{chap:algebra} of
    operads. Therefore, since
    by~\eqref{equ:proof_functor_series_pre_lie_algebra},
    $\lambda(\Fbf, \Gbf, \Hbf)$ is symmetric in $\Gbf$ and $\Hbf$,
    the series
    \begin{math}
        (\Fbf \PreLieProduct \Gbf) \PreLieProduct \Hbf -
        \Fbf \PreLieProduct (\Gbf \PreLieProduct \Hbf)
    \end{math}
    and
    \begin{math}
        (\Fbf \PreLieProduct \Hbf) \PreLieProduct \Gbf -
        \Fbf \PreLieProduct (\Hbf \PreLieProduct \Gbf)
    \end{math}
    are equal. This shows that $\K \AAngle{\Cca}$, endowed with the
    product $\PreLieProduct$, is a pre-Lie algebra. Finally, by using
    the fact that by
    Proposition~\ref{prop:functor_colored_operads_vector_spaces},
    $(\Cca, \phi) \mapsto (\K \AAngle{\Cca}, \K \AAngle{\phi})$ is
    functorial, we obtain that $\K \AAngle{\phi}$ is a morphism of
    pre-Lie algebras. Hence, the statement of the proposition holds.
\end{proof}
\medbreak

Proposition~\ref{prop:functor_series_pre_lie_algebras} shows that
$\PreLieProduct$ is a pre-Lie product. This product $\PreLieProduct$ is
a generalization of a pre-Lie product defined in~\cite{Cha08}, endowing
the linear span of the underlying monochrome collection of a monochrome
operad with a pre-Lie algebra structure.
\medbreak

\subsubsection{Noncommutative multivariate series and series on monoids}
The pre-Lie product $\PreLieProduct$ on $\Cca$-series provides also a
generalization of the usual product of noncommutative multivariate
series. Indeed, consider the method described in
Section~\ref{subsubsec:noncommutative_multivariate_series} to encode
noncommutative multivariate series on an alphabet $A$ as series on
the colored operad $\Bud_{\CFr_A}(\As)$. For any
$\Sbf, \Tbf \in \K \AAngle{A^*}$ and $u \in A^*$, we have
\begin{equation}\begin{split}
    \Angle{\left(\Dummy, u \Dummy\right),
       \Multi(\Sbf) \PreLieProduct \Multi(\Tbf)} & =
    \sum_{\substack{
        v, w \in A^* \\
        i \in \N \\
        \left(\Dummy, u \Dummy\right)
        = \left(\Dummy, v \Dummy\right) \circ_i
        \left(\Dummy, w \Dummy\right)
    }}
    \Angle{\left(\Dummy, v \Dummy\right), \Multi(\Sbf)}
    \enspace
    \Angle{\left(\Dummy, w \Dummy\right),
    \Multi(\Tbf)} \\
    & = \sum_{\substack{v, w \in A^* \\ u = vw}}
    \Angle{v, \Sbf} \Angle{w, \Tbf} \\
    & = \Angle{u, \Sbf \Tbf},
\end{split}\end{equation}
so that $\Multi(\Sbf \Tbf) = \Multi(\Sbf) \PreLieProduct \Multi(\Tbf)$.
\medbreak

Moreover, through the method presented in
Section~\ref{subsubsec:series_on_monoids} to encode series on a monoid
$\Mca$ as series on the colored operad $\K \AAngle{\Oca_\Mca}$, we have
for any $\Sbf, \Tbf \in \K \AAngle{\Mca}$ and $u \in \Mca$,
\begin{equation}\begin{split}
    \Angle{u, \Mono(\Sbf) \PreLieProduct \Mono(\Tbf)} & =
    \sum_{\substack{
        v, w \in \Oca_\Mca \\
        i \in \N \\
        u = v \circ_i w
    }}
    \Angle{v, \Mono(\Sbf)} \enspace \Angle{w, \Mono(\Tbf)} \\
    & =
    \sum_{\substack{v, w \in \Mca \\ u = v \Product w}}
    \Angle{v, \Sbf} \Angle{w, \Tbf} \\
    & = \Angle{u, \Sbf \Tbf},
\end{split}\end{equation}
where $\Product$ is the operation of $\Mca$, so that
$\Mono(\Sbf \Tbf) = \Mono(\Sbf) \PreLieProduct \Mono(\Tbf)$. Hence, the
pre-Lie product of series on colored operads is a generalization of the
Cauchy product of series on monoids~\cite{Sak09}.
\medbreak

\subsubsection{Pre-Lie star product}
For any $\Cca$-series $\Fbf \in \K \AAngle{\Cca}$ and any $\ell \geq 0$,
let $\Fbf^{\PreLieProduct_\ell}$ be the $\Cca$-series recursively
defined by
\begin{equation}
    \Fbf^{\PreLieProduct_\ell} :=
    \begin{cases}
        \Ubf & \mbox{if } \ell = 0, \\
        \Fbf^{\PreLieProduct_{\ell - 1}} \PreLieProduct \Fbf
            & \mbox{otherwise}.
    \end{cases}
\end{equation}
Immediately from this definition and the definition of the pre-Lie
product $\PreLieProduct$, the coefficients of
$\Fbf^{\PreLieProduct_\ell}$, $\ell \geq 0$, satisfies for any
$x \in \Cca$,
\begin{equation} \label{equ:induction_pre_lie_powers}
    \Angle{x, \Fbf^{\PreLieProduct_\ell}} =
    \begin{cases}
        \delta_{x, \Unit_{\Out(x)}} & \mbox{if } \ell = 0, \\
        \sum\limits_{\substack{
            y, z \in \Cca \\
            i \in [|y|] \\
            x = y \circ_i z}}
        \Angle{y, \Fbf^{\PreLieProduct_{\ell - 1}}} \Angle{z, \Fbf}
        & \mbox{otherwise}.
    \end{cases}
\end{equation}
\medbreak

\begin{Lemma} \label{lem:pre_lie_powers}
    Let $\Cca$ be a combinatorial $\CFr$-colored operad and $\Fbf$ be a
    series of $\K \AAngle{\Cca}$. Then, the coefficients of
    $\Fbf^{\PreLieProduct_{\ell + 1}}$, $\ell \geq 0$, satisfy for any
    $x \in \Cca$,
    \begin{equation} \label{equ:pre_lie_powers}
        \Angle{x, \Fbf^{\PreLieProduct_{\ell + 1}}} =
        \sum_{\substack{
            y_1, \dots, y_{\ell + 1} \in \Cca \\
            i_1, \dots, i_\ell \in \N \\
            x = ( \dots ( y_1 \circ_{i_1} y_2 )
                \circ_{i_2} \dots ) \circ_{i_\ell} y_{\ell + 1}
        }}
        \prod_{j \in [\ell + 1]} \Angle{y_j, \Fbf}.
    \end{equation}
\end{Lemma}
\medbreak

The \Def{$\PreLieProduct$-star} of $\Fbf$ is the series
\begin{equation}\begin{split}
    \Fbf^{\PreLieProduct_*} & :=
    \sum_{\ell \geq 0} \Fbf^{\PreLieProduct_\ell} \\
    & = \Ubf + \Fbf + \Fbf^{\PreLieProduct_2} +
        \Fbf^{\PreLieProduct_3} + \Fbf^{\PreLieProduct_4} + \cdots \\
    & = \Ubf + \Fbf + \Fbf \PreLieProduct \Fbf +
        (\Fbf \PreLieProduct \Fbf) \PreLieProduct \Fbf +
        ((\Fbf \PreLieProduct \Fbf) \PreLieProduct \Fbf) \PreLieProduct
        \Fbf + \cdots.
\end{split}\end{equation}
Observe that $\Fbf^{\PreLieProduct_*}$ could be undefined for an
arbitrary $\Cca$-series $\Fbf$.
\medbreak

In what follows, we shall use the notion of finite factorization
introduced in Section~\ref{subsubsec:left_expressions_hook_length} of
Chapter~\ref{chap:algebra}. More precisely, in this context of colored
operads, we say that a subset $S$ of $\Cca(1)$ finitely factorizes
$\Cca(1)$ if any element of $\Cca(1)$ admits finitely many
factorizations on $S$ with respect to the operation~$\circ_1$.
\medbreak

\begin{Proposition} \label{prop:coefficients_pre_lie_star}
    Let $\Cca$ be a $\CFr$-colored operad and $\Fbf$ be a series of
    $\K \AAngle{\Cca}$. Then, if $\Cca$ is combinatorial and
    $\Support(\Fbf)(1)$ finitely factorizes $\Cca(1)$,
    $\Fbf^{\PreLieProduct_*}$ is a well-defined series. Moreover, in
    this case, for any $x \in \Cca$, the coefficient of $x$ in
    $\Fbf^{\PreLieProduct_*}$ is
    \begin{equation} \label{equ:coefficients_pre_lie_star}
        \Angle{x, \Fbf^{\PreLieProduct_*}} =
        \delta_{x, \Unit_{\Out(x)}} +
        \sum_{\substack{
            y, z \in \Cca \\
            i \in [|y|] \\
            x = y \circ_i z}}
        \Angle{y, \Fbf^{\PreLieProduct_*}} \Angle{z, \Fbf}.
    \end{equation}
\end{Proposition}
\medbreak

Proposition~\ref{prop:coefficients_pre_lie_star} gives hence a way,
given a $\Cca$-series $\Fbf$ satisfying the constraints stated, to
compute recursively the coefficients of its $\PreLieProduct$-star
$\Fbf^{\PreLieProduct_*}$ by using~\eqref{equ:coefficients_pre_lie_star}.
\medbreak

\begin{Proposition} \label{prop:equation_pre_lie_star}
    Let $\Cca$ be a combinatorial $\CFr$-colored operad and $\Fbf$ be a
    series of $\K \AAngle{\Cca}$ such that $\Support(\Fbf)(1)$ finitely
    factorizes $\Cca(1)$. Then, the equation
    \begin{equation} \label{equ:equation_pre_lie_star}
        \Xbf - \Xbf \PreLieProduct \Fbf = \Ubf
    \end{equation}
    admits the unique solution $\Xbf = \Fbf^{\PreLieProduct_*}$.
\end{Proposition}
\medbreak

\subsection{Composition product on series}
We define here a binary operation $\Compo$ on the space of
$\Cca$-series. As we shall see, this operation is partially defined,
unitary, noncommutative, and associative.
\medbreak

\subsubsection{Composition product}
Given two $\Cca$-series $\Fbf, \Gbf \in \K \AAngle{\Cca}$, the
\Def{composition product} of $\Fbf$ and $\Gbf$ is the $\Cca$-series
$\Fbf \Compo \Gbf$ defined, for any $x \in \Cca$, by
\begin{equation} \label{equ:associative_product}
    \Angle{x, \Fbf \Compo \Gbf}
    :=
    \sum_{\substack{
        y, z_1, \dots, z_{|y|} \in \Cca \\
        x = y \circ \left[z_1, \dots, z_{|y|}\right]
    }}
    \Angle{y, \Fbf}
    \prod_{i \in [|y|]} \Angle{z_i, \Gbf}.
\end{equation}
Observe that $\Fbf \Compo \Gbf$ could be undefined for arbitrary
$\Cca$-series $\Fbf$ and $\Gbf$ on an arbitrary colored operad $\Cca$.
Besides, notice from~\eqref{equ:associative_product} that $\Compo$ is
linear on the left and that the series $\Ubf$ is the left and right unit
of $\Compo$. However, this product is not linear on the right since we
have, for instance in~$\K \AAngle{\As}$,
\begin{equation}
    \Afr_2 \Compo (\Afr_2 + \Afr_3) =
    \Afr_4 + 2 \Afr_5 + \Afr_6
    \ne \Afr_4 + \Afr_6 =
    \Afr_2 \Compo \Afr_2 + \Afr_2 \Compo \Afr_3.
\end{equation}
\medbreak

\begin{Proposition} \label{prop:functor_series_monoids}
    The construction
    \begin{math}
        (\Cca, \phi) \mapsto
        ((\K \AAngle{\Cca}, \Compo), \K \AAngle{\phi})
    \end{math}
    is a functor from the category of combinatorial $\CFr$-colored
    operads to the category of monoids.
\end{Proposition}
\medbreak

Proposition~\ref{prop:functor_series_monoids} shows that $\Compo$ is an
associative product. This product $\Compo$ is a generalization of the
composition product of series on operads of~\cite{Cha02,Cha09} (see
also~\cite{Vdl04,Fra08,Cha08,LV12,LN13}).
\medbreak

\subsubsection{Composition star product}
For any $\Cca$-series $\Fbf \in \K \AAngle{\Cca}$ and any $\ell \geq 0$,
let $\Fbf^{\Compo_\ell}$ be the series recursively defined by
\begin{equation} \label{equ:star_composition_product}
    \Fbf^{\Compo_\ell} :=
    \prod_{i \in [\ell]} \Fbf,
\end{equation}
where the product of~\eqref{equ:star_composition_product} denotes the
iterated version of $\Compo$. Since by
Proposition~\ref{prop:functor_series_monoids}, $\Compo$ is associative,
this definition is consistent. Immediately from this definition and the
definition of the composition product $\Compo$, the coefficient of
$\Fbf^{\Compo_\ell}$, $\ell \geq 0$, satisfies for any $x \in \Cca$,
\begin{equation} \label{equ:induction_composition_powers}
    \Angle{x, \Fbf^{\Compo_\ell}} =
    \begin{cases}
        \delta_{x, \Unit_{\Out(x)}} & \mbox{if } \ell = 0, \\
        \sum\limits_{
        \substack{
            y, z_1, \dots, z_{|y|} \in \Cca \\
            x = y \circ \left[z_1, \dots, z_{|y|}\right]}}
        \Angle{y, \Fbf^{\Compo_{\ell - 1}}}
        \prod\limits_{i \in [|y|]} \Angle{z_i, \Gbf}
        & \mbox{otherwise}.
    \end{cases}
\end{equation}
\medbreak

\begin{Lemma} \label{lem:composition_powers}
    Let $\Cca$ be a combinatorial $\CFr$-colored operad and $\Fbf$ be a
    series of $\K \AAngle{\Cca}$. Then, the coefficients of
    $\Fbf^{\Compo_{\ell + 1}}$, $\ell \geq 0$, satisfy for any
    $x \in \Cca$,
    \begin{equation} \label{equ:composition_powers}
        \Angle{x, \Fbf^{\Compo_{\ell + 1}}} =
        \sum_{\substack{
            \Tfr \in \FreeColoredOperad_\Perfect(\Cca) \\
            \Height(\Tfr) = \ell + 1 \\
            \Eval(\Tfr) = x
        }}
        \enspace
        \prod_{v \in \TreeLanguage_\Node(\Tfr)} \Angle{\Tfr(v), \Fbf}.
    \end{equation}
\end{Lemma}
\medbreak

Recall that the notations $\Height(\Tfr)$ and
$\TreeLanguage_\Node(\Tfr)$ appearing in the statement of
Lemma~\ref{lem:composition_powers} stand respectively for the height of
$\Tfr$ and for the set of the internal nodes of $\Tfr$ (see
Section~\ref{subsubsec:definitions_trees} of
Chapter~\ref{chap:combinatorics}). Moreover, the notation
$\FreeColoredOperad_\Perfect(\Cca)$ denotes the set of all perfect
colored $\Cca$-syntax trees (we recall that a tree $\Tfr$ is perfect if
all the maximal paths have the same length).
\medbreak

The \Def{$\Compo$-star} of $\Fbf$ is the series
\begin{equation}\begin{split}
    \Fbf^{\Compo_*} & := \sum_{\ell \geq 0} \Fbf^{\Compo_\ell} \\
    & = \Ubf + \Fbf + \Fbf^{\Compo_2} + \Fbf^{\Compo_3} +
        \Fbf^{\Compo_4} + \cdots \\
    & = \Ubf + \Fbf + \Fbf \Compo \Fbf + \Fbf \Compo \Fbf \Compo \Fbf +
        \Fbf \Compo \Fbf \Compo \Fbf \Compo \Fbf + \cdots.
\end{split}\end{equation}
Observe that $\Fbf^{\Compo_*}$ could be undefined for an arbitrary
$\Cca$-series~$\Fbf$.
\medbreak

\begin{Proposition} \label{prop:coefficients_composition_star}
    Let $\Cca$ be a $\CFr$-colored operad and $\Fbf$ be a series of
    $\K \AAngle{\Cca}$. Then, if $\Cca$ is combinatorial and
    $\Support(\Fbf)(1)$ finitely factorizes $\Cca(1)$, $\Fbf^{\Compo_*}$
    is a well-defined series. Moreover, in this case, for any
    $x \in \Cca$, the coefficient of $x$ in $\Fbf^{\Compo_*}$ is
    \begin{equation} \label{equ:coefficients_composition_star}
        \Angle{x, \Fbf^{\Compo_*}} =
        \delta_{x, \Unit_{\Out(x)}} +
        \sum_{\substack{
            y, z_1, \dots, z_{|y|} \in \Cca \\
            x = y \circ \left[z_1, \dots, z_{|y|}\right]
        }}
        \Angle{y, \Fbf^{\Compo_*}}
        \prod_{i \in [|y|]} \Angle{z_i, \Fbf}.
    \end{equation}
\end{Proposition}
\medbreak

Proposition~\ref{prop:coefficients_composition_star} gives hence a way,
given a $\Cca$-series $\Fbf$ satisfying the constraints stated, to
compute recursively the coefficients of its $\Compo$-star
$\Fbf^{\Compo_*}$ by using~\eqref{equ:coefficients_composition_star}.
\medbreak

\begin{Proposition} \label{prop:equation_composition_star}
    Let $\Cca$ be a combinatorial $\Cca$-colored operad and $\Fbf$ be a
    series of $\K \AAngle{\Cca}$ such that $\Support(\Fbf)(1)$
    finitely factorizes $\Cca(1)$. Then, the equation
    \begin{equation} \label{equ:equation_composition_star}
        \Xbf - \Xbf \Compo \Fbf = \Ubf
    \end{equation}
    admits the unique solution $\Xbf = \Fbf^{\Compo_*}$.
\end{Proposition}
\medbreak

\subsubsection{Invertible elements}
For any $\Cca$-series $\Fbf \in \K \AAngle{\Cca}$, the
\Def{$\Compo$-inverse} of $\Fbf$ is the series $\Fbf^{\Compo_{-1}}$
whose coefficients are defined for any $x \in \Cca$ by
\begin{equation} \label{equ:equation_composition_inverse}
    \Angle{x, \Fbf^{\Compo_{-1}}}
    :=
    \frac{\delta_{x, \Unit_{\Out(x)}}}{\Angle{\Unit_{\Out(x)}, \Fbf}}
     -
    \frac{1}{\Angle{\Unit_{\Out(x)}, \Fbf}}
    \sum_{\substack{
        y, z_1, \dots, z_{|y|} \in \Cca \\
        y \ne \Unit_{\Out(x)} \\
        x = y \circ \left[z_1, \dots, z_{|y|}\right]}}
    \Angle{y, \Fbf}
    \prod_{i \in [|y|]} \Angle{z_i, \Fbf^{\Compo_{-1}}}.
\end{equation}
Observe that $\Fbf^{\Compo_{-1}}$ could be undefined for an arbitrary
$\Cca$-series $\Fbf$.
\medbreak

\begin{Proposition} \label{prop:composition_inverse_expression}
    Let $\Cca$ be a combinatorial colored $\CFr$-operad and $\Fbf$ be a
    series of $\K \AAngle{\Cca}$ such that
    $\Support(\Fbf) = \{\Unit_a : a \in \CFr\} \sqcup S$ where $S$ is
    a $\CFr$-colored subcollection of $\Cca$ such that $S(1)$ finitely
    factorizes $\Cca(1)$. Then, $\Fbf^{\Compo_{-1}}$ is a well-defined
    series and the coefficients of $\Fbf^{\Compo_{-1}}$ satisfy for any
    $x \in \Cca$,
    \begin{equation} \label{equ:composition_inverse_expression}
        \Angle{x, \Fbf^{\Compo_{-1}}} =
        \frac{1}{\Angle{\Unit_{\Out(x)}, \Fbf}}
        \enspace
        \sum_{\substack{
            \Tfr \in \FreeColoredOperad(S) \\
            \Eval(\Tfr) = x
        }}
        (-1)^{\deg(\Tfr)}
        \prod_{v \in \TreeLanguage_\Node(\Tfr)}
        \enspace
        \frac{\Angle{\Tfr(v), \Fbf}}
            {\prod\limits_{j \in [|v|]} \Angle{\Unit_{\In_j(v)}, \Fbf}}.
    \end{equation}
\end{Proposition}
\medbreak

\begin{Proposition} \label{prop:equation_composition_inverse}
    Let $\Cca$ be a combinatorial colored $\CFr$-operad and $\Fbf$ be a
    series of $\K \AAngle{\Cca}$ such that
    $\Support(\Fbf) = \{\Unit_a : a \in \CFr\} \sqcup S$ where $S$ is
    a $\CFr$-colored subcollection of $\Cca$ such that $S(1)$
    finitely factorizes $\Cca(1)$. Then, the equations
    \begin{equation} \label{equ:equation_composition_inverse_1}
        \Fbf \Compo \Xbf = \Ubf
    \end{equation}
    and
    \begin{equation} \label{equ:equation_composition_inverse_2}
        \Xbf \Compo \Fbf = \Ubf
    \end{equation}
    admit both the unique solution $\Xbf = \Fbf^{\Compo_{-1}}$.
\end{Proposition}
\medbreak

Proposition~\ref{prop:equation_composition_inverse} shows that the
$\Compo$-inverse $\Fbf^{\Compo_{-1}}$ of a series $\Fbf$ satisfying the
constraints stated is the inverse of $\Fbf$ for the composition product.
Moreover, $\Fbf^{\Compo_{-1}}$ can be computed recursively by
using~\eqref{equ:equation_composition_inverse} or directly by
using~\eqref{equ:composition_inverse_expression}.
\medbreak

\begin{Proposition} \label{prop:series_composition_group}
    Let $\Cca$ be a combinatorial colored $\CFr$-operad. Then, the
    subset of $\K \AAngle{\Cca}$ consisting in all series $\Fbf$ such
    that $\Support(\Fbf) = \{\Unit_a : a \in \CFr\} \sqcup S$ where $S$
    is a $\CFr$-colored subcollection of $\Cca$ such that $S(1)$
    finitely factorizes $\Cca(1)$ forms a group for the composition
    product~$\Compo$.
\end{Proposition}
\medbreak

The group obtained from $\Cca$ of the $\Cca$-series satisfying the
conditions of Proposition~\ref{prop:series_composition_group} is a
generalization of the groups constructed from operads
of~\cite{Cha02,Cha09} (see also~\cite{Vdl04,Fra08,Cha08,LV12,LN13}).
\medbreak

\section{Bud generating systems and combinatorial generation}
\label{sec:bud_generating_systems}
In this section, we introduce bud generating systems. A bud generating
system relies on a monochrome operad $\Oca$, a set of colors $\CFr$,
and the bud operad $\Bud_\CFr(\Oca)$. The principal interest of these
objects is that they allow to specify sets of objects of
$\Bud_\CFr(\Oca)$. We shall also establish some first properties of bud
generating systems by showing that they can emulate context-free
grammars, regular tree grammars, and synchronous grammars.
\medbreak

\subsection{Bud generating systems}
We introduce here the main definitions and the main tools about bud
generating systems.
\medbreak

\subsubsection{Bud generating systems}
A \Def{bud generating system} is a tuple
$\Bca := (\Oca, \CFr, \RFr, I, T)$ where $\Oca$ is an operad called
\Def{ground operad}, $\CFr$ is a finite set of colors, $\RFr$ is a
finite $\CFr$-colored subcollection of $\Bud_\CFr(\Oca)$ called
\Def{set of rules}, $I$ is a subset of $\CFr$ called \Def{set of initial
colors}, and $T$ is a subset of $\CFr$ called \Def{set of terminal
colors}.
\medbreak

A \Def{monochrome bud generating system} is a bud generating system
whose set $\CFr$ of colors contains a single color, and whose sets of
initial and terminal colors are equal to $\CFr$. In this case, as
explained in Section~\ref{subsubsec:operad_to_colored_operads},
$\Bud_\CFr(\Oca)$ and $\Oca$ are identified. These particular
generating systems are hence simply denoted by pairs $(\Oca, \RFr)$.
\medbreak

Let us explain how bud generating systems specify, in two different
ways, two $\CFr$-colored subcollections of $\Bud_\CFr(\Oca)$. In
what follows, $\Bca := (\Oca, \CFr, \RFr, I, T)$ is a bud generating
system.
\medbreak

\subsubsection{Generation} \label{subsubsec:generation}
We say that $x_2 \in \Bud_\CFr(\Oca)$ is \Def{derivable in one step}
from $x_1 \in \Bud_\CFr(\Oca)$ if there is a rule $r \in \RFr$ and an
integer $i$ such that such that $x_2 = x_1 \circ_i r$. We denote this
property by $x_1 \Deriv x_2$. When $x_1, x_2 \in \Bud_\CFr(\Oca)$ are
such that $x_1 = x_2$ or there are
$y_1, \dots, y_{\ell - 1} \in \Bud_\CFr(\Oca)$, $\ell \geq 1$,
satisfying
\begin{equation}
    x_1 \Deriv y_1 \Deriv \cdots \Deriv y_{\ell - 1} \Deriv x_2,
\end{equation}
we say that $x_2$ is \Def{derivable} from $x_1$. Moreover, $\Bca$
\Def{generates} $x \in \Bud_\CFr(\Oca)$ if there is a color $a$ of
$I$ such that $x$ is derivable from $\Unit_a$ and all colors of $\In(x)$
are in $T$. The \Def{language} $\Lang(\Bca)$ of $\Bca$ is the set of all
the elements of $\Bud_\CFr(\Oca)$ generated by $\Bca$. Finally, $\Bca$
is \Def{faithful} if the characteristic series of $\Lang(\Bca)$ is
faithful (see Section~\ref{subsubsec:pruned_series}). Observe that all
monochrome bud generating systems are faithful.
\medbreak

The \Def{derivation graph} of $\Bca$ is the directed multigraph
$\DerivGraph(\Bca)$ with the set of elements derivable from $\Unit_a$,
$a \in I$, as set of vertices. In $\DerivGraph(\Bca)$, for any
$x_1, x_2 \in \Lang(\Bca)$ such that $x_1 \Deriv x_2$, there are $\ell$
edges from $x_1$ to $x_2$, where $\ell$ is the number of pairs
$(i, r) \in \N \times \RFr$ such that $x_2 = x_1 \circ_i r$.
\medbreak

\subsubsection{A bud generating system for Motzkin paths}
\label{subsubsec:example_bmotz}
Let us consider the operad $\Motz$ on Motzkin paths introduced in
Section~\ref{subsubsec:operad_Motz} of Chapter~\ref{chap:monoids},
seen as a set-operad.
Let the bud generating system
$\BMotz := (\Motz, \{1, 2\}, \RFr, \{1\}, \{1, 2\})$ where
\begin{equation}
    \RFr := \left\{
    \left(1, \MotzHoriz, 22\right),
    \left(1, \MotzPeak, 111\right)
    \right\}.
\end{equation}
Figure~\ref{fig:derivations_bmotz} shows a sequence of derivations in
$\BMotz$ and Figure~\ref{fig:derivation_graph_bmotz} shows the
derivation graph of~$\BMotz$.
\begin{figure}[ht]
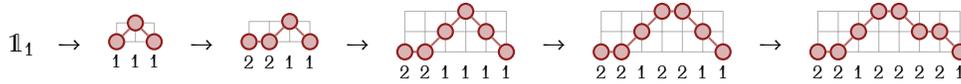

    \centering
    \begin{equation*}
        \Unit_1
        \enspace \Deriv \enspace

    \end{equation*}
    \caption[Derivations in a bud generating system on Motzkin
    paths.]
    {A sequence of derivations in $\BMotz$. The input colors of the
    elements of $\Bud_{\{1, 2\}}(\Motz)$ are depicted below the paths.
    The output color of all these elements is~$1$.}
    \label{fig:derivations_bmotz}
\end{figure}
\begin{figure}[ht]
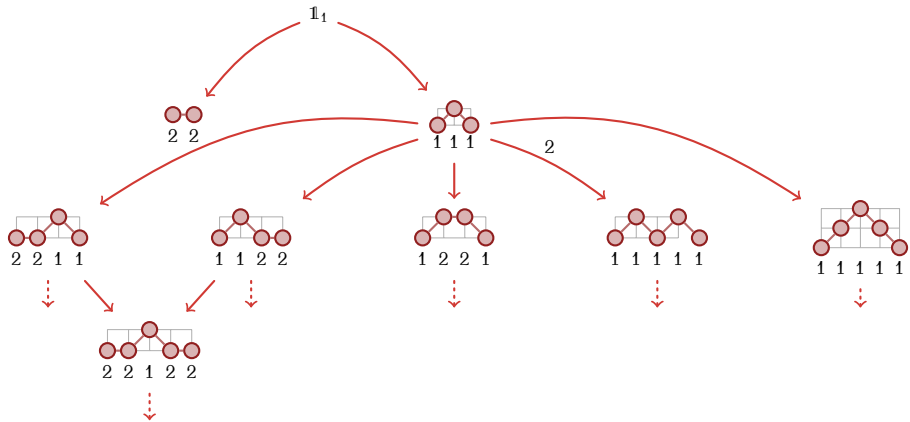

    \centering
    \begin{equation*}
        %
\end{math}};
            \draw[EdgeRew](Unit)
                edge[bend right=12]
                node[anchor=south,font=\scriptsize,color=ColBlack]
                {}(00);
            \draw[EdgeRew](Unit)
                edge[bend left=12]
                node[anchor=south,font=\scriptsize,color=ColBlack]
                {}(010);
            \draw[EdgeRew](010)
                edge[bend right=13]
                node[anchor=south,font=\scriptsize,color=ColBlack]
                {}(0010);
            \draw[EdgeRew](010)
                edge[]
                node[anchor=west,font=\scriptsize,color=ColBlack]
                {}(0110);
            \draw[EdgeRew](010)
                edge[bend right=7]
                node[anchor=south,font=\scriptsize,color=ColBlack]
                {}(0100);
            \draw[EdgeRew](010)
                edge[bend left=7]
                node[anchor=south,font=\scriptsize,color=ColBlack]
                {\begin{math}2\end{math}}(01010);
            \draw[EdgeRew](010)
                edge[bend left=13]
                node[anchor=south,font=\scriptsize,color=ColBlack]
                {}(01210);
            \draw[EdgeRew](0010)
                edge[]
                node[anchor=south,font=\scriptsize,color=ColBlack]
                {}(00100);
            \draw[EdgeRew](0100)
                edge[]
                node[anchor=south,font=\scriptsize,color=ColBlack]
                {}(00100);
            \draw[EdgeRew,dotted,shorten >=6mm](0010)--(-4,-3);
            \draw[EdgeRew,dotted,shorten >=6mm](0110)--(2,-3);
            \draw[EdgeRew,dotted,shorten >=6mm](0100)--(-1,-3);
            \draw[EdgeRew,dotted,shorten >=6mm](01010)--(5,-3);
            \draw[EdgeRew,dotted,shorten >=6mm](01210)--(8,-3);
            \draw[EdgeRew,dotted,shorten >=6mm](00100)--(-2.5,-4);
        \end{tikzpicture}
    \end{equation*}
    \caption[The derivation graph of a bud generating system of
    Motzkin paths.]
    {The derivation graph of $\BMotz$. The input colors of the
    elements of $\Bud_{\{1, 2\}}(\Motz)$ are depicted below the paths.
    The output color of all these elements is~$1$.}
    \label{fig:derivation_graph_bmotz}
\end{figure}
\medbreak

Let $L_{\BMotz}$ be the set of Motzkin paths with no consecutive
horizontal steps.
\medbreak

\begin{Proposition} \label{prop:properties_BMotz}
    The bud generating system $\BMotz$ satisfies the following
    properties.
    \begin{enumerate}[label = ({\it {\roman*})}]
        \item \label{item:properties_BMotz_faithful}
        It is faithful.
        \item \label{item:properties_BMotz_language}
        The restriction of the pruning map $\Prune$ on the domain
        $\Lang(\BMotz)$
        is a bijection between $\Lang(\BMotz)$ and $L_{\BMotz}$.
        \item \label{item:properties_BASchr_finitely_factorizing}
        The set of rules $\RFr(1)$ finitely factorizes
        $\Bud_{\{1, 2\}}(\Motz)(1)$.
    \end{enumerate}
\end{Proposition}
\medbreak

Properties~\ref{item:properties_BMotz_faithful}
and~\ref{item:properties_BMotz_language} of
Proposition~\ref{prop:properties_BMotz} together say that the sequence
enumerating the elements of $\Lang(\BMotz)$ with respect to their arity
is the one enumerating the Motzkin paths with no consecutive horizontal
steps. This sequence is Sequence~\OEIS{A104545} of~\cite{Slo}, starting
by
\begin{equation}
    1, 1, 1, 3, 5, 11, 25, 55, 129, 303, 721, 1743, 4241, 10415,
    25761, 64095.
\end{equation}
\medbreak

\subsubsection{Synchronous generation}
\label{subsubsec:synchronous_generation}
We say that $x_2 \in \Bud_\CFr(\Oca)$ is \Def{synchronously derivable
in one step} from $x_1 \in \Bud_\CFr(\Oca)$ if there are rules $r_1$,
\dots, $r_{|x_1|}$ of $\RFr$ such that
$x_2 = x_1 \circ \left[r_1, \dots, r_{|x_1|}\right]$. We denote this
property by $x_1 \SyncDeriv x_2$. When $x_1, x_2 \in \Bud_\CFr(\Oca)$
are such that $x_1 = x_2$ or there are
\begin{math}
    y_1, \dots, y_{\ell - 1} \in \Bud_\CFr(\Oca)$, $\ell \geq 1,
\end{math}
satisfying
\begin{equation}
    x_1 \SyncDeriv y_1 \SyncDeriv \cdots
    \SyncDeriv y_{\ell - 1} \SyncDeriv x_2,
\end{equation}
we say that $x_2$ is \Def{synchronously derivable} from $x_1$. Moreover,
$\Bca$ \Def{synchronously generates} $x \in \Bud_\CFr(\Oca)$ if there
is a color $a$ of $I$ such that $x$ is synchronously derivable from
$\Unit_a$ and all colors of $\In(x)$ are in $T$. The \Def{synchronous
language} $\SyncLang(\Bca)$ of $\Bca$ is the set of all the elements of
$\Bud_\CFr(\Oca)$ synchronously generated by $\Bca$. Finally, we say
that $\Bca$ is \Def{synchronously faithful} if the characteristic series
of $\SyncLang(\Bca)$ is faithful (see
Section~\ref{subsubsec:pruned_series}). Observe that all monochrome
bud generating systems are synchronously faithful.
\medbreak

The \Def{synchronous derivation graph} of $\Bca$ is the directed
multigraph $\SyncDerivGraph(\Bca)$ with the set of elements
synchronously derivable from $\Unit_a$, $a \in I$, as set of vertices.
In $\SyncDerivGraph(\Bca)$, for any $x_1, x_2 \in \SyncLang(\Bca)$ such
that $x_1 \SyncDeriv x_2$, there are $\ell$ edges from $x_1$ to $x_2$,
where $\ell$ is the number of tuples
$(r_1, \dots, r_{|x_1|}) \in \RFr^{|x_1|}$ such that
$x_2 = x_1 \circ \left[r_1, \dots, r_{|x_1|}\right]$.
\medbreak

\subsubsection{A bud generating system for balanced binary trees}
\label{subsubsec:example_bbaltree}
Let us consider the magmatic operad $\Mag$ (see
Section~\ref{subsubsec:magmatic_operad} of Chapter~\ref{chap:algebra}).
Let the bud generating system
$\BBalTree := (\Mag, \{1, 2\}, \RFr, \{1\}, \{1\})$ where
\begin{equation}
    \RFr := \left\{
    \left(1, \CorollaTwo{}, 11\right),
    \left(1, \CorollaTwo{}, 12\right),
    \left(1, \CorollaTwo{}, 21\right),
    \left(2, \LeafPic, 1\right)\right\}.
\end{equation}
Figure~\ref{fig:derivations_bal_tree} shows a sequence of synchronous
derivations in $\BBalTree$.
\begin{figure}[ht]
    \centering
    \begin{equation*}

    \end{equation*}
    \caption[Derivations in a bud generating system on binary trees.]
    {A sequence of synchronous derivations in $\BBalTree$. The input
    colors of the elements of $\Bud_{\{1, 2\}}(\Mag)$ are depicted below
    the leaves. The output color of all these elements is $1$. Since all
    input colors of the last tree are $1$, this tree is
    in~$\SyncLang(\BBalTree)$.}
    \label{fig:derivations_bal_tree}
\end{figure}
\medbreak

Recall that if $\Tfr$ is a binary tree, the height of $\Tfr$ is the
length of a longest path connecting the root of $\Tfr$ to one of its
leaves (see Section~\ref{subsubsec:definitions_trees} of
Chapter~\ref{chap:combinatorics}). A \Def{balanced binary
tree}~\cite{AVL62} is a binary tree $\Tfr$ wherein, for any internal
node $u$ of $\Tfr$, the difference between the height of the left
subtree and the height of the right subtree of $u$ is $-1$, $0$, or~$1$.
\medbreak

\begin{Proposition} \label{prop:properties_BBalTree}
    The bud generating system $\BBalTree$ satisfies the following
    properties.
    \begin{enumerate}[label = ({\it {\roman*})}]
        \item \label{item:properties_BBalTree_faithful}
        It is synchronously faithful.
        \item \label{item:properties_BBalTree_language}
        The restriction of the pruning map $\Prune$ on the domain
        $\SyncLang(\BBalTree)$
        is a bijection between $\SyncLang(\BBalTree)$ and the set
        of balanced binary trees.
        \item \label{item:properties_BBalTree_finitely_factorizing}
        The set of rules $\RFr(1)$ finitely factorizes
        $\Bud_{\{1, 2\}}(\Mag)(1)$.
    \end{enumerate}
\end{Proposition}
\medbreak

Property~\ref{item:properties_BBalTree_language} of
Proposition~\ref{prop:properties_BBalTree} is based upon
combinatorial properties of a synchronous grammar $\Gca$ of balanced
binary trees defined in~\cite{Gir12e} and satisfying
$\SG(\Gca) = \BBalTree$ (see
Section~\ref{subsubsec:synchronous_grammars} and
Proposition~\ref{prop:emulation_synchronous_grammars}).
Besides, Properties~\ref{item:properties_BBalTree_faithful}
and~\ref{item:properties_BBalTree_language} of
Proposition~\ref{prop:properties_BBalTree} together imply that the
sequence enumerating the elements of $\SyncLang(\BBalTree)$ with respect
to their arity is the one enumerating the balanced binary trees. This
sequence in Sequence~\OEIS{A006265} of~\cite{Slo}, starting by
\begin{equation} \label{equ:first_terms_balanced_binary_trees}
    1, 1, 2, 1, 4, 6, 4, 17, 32, 44, 60, 70,
    184, 476, 872, 1553, 2720, 4288, 6312, 9004.
\end{equation}
\medbreak

\subsection{First properties}
We state now two properties about the languages and the synchronous
languages of bud generating systems.
\medbreak

\begin{Lemma} \label{lem:language_treelike_expressions}
    Let $\Bca := (\Oca, \CFr, \RFr, I, T)$ be a bud generating system.
    Then, for any $x \in \Bud_\CFr(\Oca)$, $x$ belongs to $\Lang(\Bca)$
    if and only if $x$ admits an $\RFr$-treelike expression with output
    color in $I$ and all input colors in~$T$.
\end{Lemma}
\begin{proof}
    Assume that $x$ belongs to $\Lang(\Bca)$. Then, by definition of the
    derivation relation $\Deriv$, $x$ admits an $\RFr$-left expression.
    Lemma~\ref{lem:left_expressions_linear_extensions} of
    Chapter~\ref{chap:algebra} implies in particular that $x$ admits an
    $\RFr$-treelike expression $\Tfr$. Moreover, since $\Tfr$ is a
    treelike expression for $x$, $\Tfr$ has the same output and input
    colors as those of $x$. Hence, because $x$ belongs to $\Lang(\Bca)$,
    its output color is in $I$ and all its input colors are in $T$.
    Thus, $\Tfr$ satisfies the required properties.
    \smallbreak

    Conversely, assume that $x$ is an element of $\Bud_\CFr(\Oca)$
    admitting an $\RFr$-treelike expression $\Tfr$ with output color in
    $I$ and all input colors in $T$.
    Lemma~\ref{lem:left_expressions_linear_extensions} of
    Chapter~\ref{chap:algebra} implies in particular that $x$ admits an
    $\RFr$-left expression. Hence, by definition of the derivation
    relation $\Deriv$, $x$ is derivable from $\Unit_{\Out(x)}$ and all
    its input colors are in $T$. Therefore, $x$ belongs
    to~$\Lang(\Bca)$.
\end{proof}
\medbreak

\begin{Lemma} \label{lem:synchronous_language_treelike_expressions}
    Let $\Bca := (\Oca, \CFr, \RFr, I, T)$ be a bud generating system.
    Then, for any $x \in \Bud_\CFr(\Oca)$, $x$ belongs to
    $\SyncLang(\Bca)$ if and only if $x$ admits an $\RFr$-treelike
    expression with output color in $I$ and all input colors in $T$ and
    which is a perfect tree.
\end{Lemma}
\medbreak

\begin{Proposition} \label{prop:language_sub_operad_generated}
    Let $\Bca := (\Oca, \CFr, \RFr, I, T)$ be a bud generating system.
    Then, the language of $\Bca$ satisfies
    \begin{equation}
        \Lang(\Bca) =
        \left\{x \in \Bud_\CFr(\Oca)^\RFr :
            \Out(x) \in I \mbox{ and } \In(x) \in T^+\right\}.
    \end{equation}
\end{Proposition}
\medbreak

\begin{Proposition} \label{prop:synchronous_language_subset_language}
    Let $\Bca := (\Oca, \CFr, \RFr, I, T)$ be a bud generating system.
    Then, the synchronous language of $\Bca$ is a subset of the
    language of~$\Bca$.
\end{Proposition}
\medbreak

\subsection{Links with other generating systems}
\label{subsec:other_generating_systems}
Context-free grammars, regular tree grammars, and synchronous grammars
are already existing generating systems describing sets of words for the
first, and sets of trees for the last two. We show here that any of
these grammars can be emulated by bud generating systems.
\medbreak

\subsubsection{Context-free grammars}
Recall that a \Def{context-free grammar}~\cite{Har78,HMU06} is a tuple
$\Gca := (V, T, P, s)$ where $V$ is a finite alphabet of
\Def{variables}, $T$ is a finite alphabet of \Def{terminal symbols}, $P$
is a finite subset of $V \times (V \sqcup T)^*$ called \Def{set of
productions}, and $s$ is a variable of $V$ called \Def{start symbol}. If
$x_1$ and $x_2$ are two words of $(V \sqcup T)^*$, $x_2$ is
\Def{derivable in one step} from $x_1$ if $x_1$ is of the form
$x_1 = u a v$ and $x_2$ is of the form $x_2 = u w v$ where
$u, v \in (V \sqcup T)^*$ and $(a, w)$ is a production of $P$. This
property is denoted by $x_1 \Deriv x_2$, so that $\Deriv$ is a binary
relation on $(V \sqcup T)^*$. The reflexive and transitive closure of
$\Deriv$ is the \Def{derivation relation}. A word $x \in T^*$ is
\Def{generated} by $\Gca$ if $x$ is derivable from the word~$s$. The
\Def{language} of $\Gca$ is the set of all words generated by $\Gca$. We
say that $\Gca$ is \Def{proper} if, for any $(a, w) \in P$, $w$ is not
the empty word.
\medbreak

If $\Gca := (V, T, P, s)$ is a proper context-free grammar, we denote by
$\CFG(\Gca)$ the bud generating system
\begin{equation}
    \CFG(\Gca) := (\As, V \sqcup T, \RFr, \{s\}, T)
\end{equation}
wherein $\RFr$ is the set of rules
\begin{equation}
    \RFr :=
    \left\{(a, u) \in \Bud_{V \sqcup T}(\As) : (a, u) \in P\right\}.
\end{equation}
\medbreak

\begin{Proposition} \label{prop:emulation_grammars}
    Let $\Gca$ be a proper context-free grammar. Then, the restriction
    of the map $\In$, sending any $(a, u) \in \Bud_{V \sqcup T}(\As)$
    to $u$, on the domain $\Lang(\CFG(\Gca))$ is a bijection between
    $\Lang(\CFG(\Gca))$ and the language of~$\Gca$.
\end{Proposition}
\begin{proof}
    Let us denote by $V$ the set of variables, by $T$ the set of
    terminal symbols, by $P$ the set of productions, and by $s$ the
    start symbol of $\Gca$.
    \smallbreak

    Let $\left(a, x\right) \in \Bud_{V \sqcup T}(\As)$, $\ell \geq 1$,
    and $y_1, \dots, y_{\ell - 1} \in (V \sqcup T)^*$. Then, by
    definition of $\CFG$, there are in $\CFG(\Gca)$ the derivations
    \begin{equation}
        \Unit_s \Deriv \left(s, y_1\right) \Deriv \cdots
        \Deriv \left(s, y_{\ell - 1}\right)
        \Deriv \left(a, x\right)
    \end{equation}
    if and only if $a = s$ and there are in $\Gca$ the derivations
    \begin{equation}
        s \Deriv y_1 \Deriv \cdots \Deriv y_{\ell - 1} \Deriv x.
    \end{equation}
    Then, $\left(a, x\right)$ belongs to $\Lang(\CFG(\Gca))$ if and
    only if $a = s$ and $x$ belongs to the language of $\Gca$. The fact
    that $\In\left(\left(s, x\right)\right) = x$ completes the proof.
\end{proof}
\medbreak

\subsubsection{Regular tree grammars}
Let $V$ be a finite graded collection of \Def{variables} and $T$ be a
finite graded collection of \Def{terminal symbols}. For any $n \geq 0$
and $a \in T(n)$ (resp. $a \in V(n)$), the arity $|a|$ of $a$ is $n$. We
moreover impose that all the elements of $V$ are of arity $0$. The tuple
$(V, T)$ is called a \Def{signature}.
\medbreak

A \Def{$(V, T)$-tree} is an element of
$\Bud_{V \sqcup T(0)}(\FreeOperad(T \setminus T(0)))$, where
$T \setminus T(0)$ is seen as a monochrome collection. In other words, a
$(V, T)$-tree is a planar rooted $\Tfr$ tree such that, for any
$n \geq 1$, any internal node of $\Tfr$ having $n$ children is labeled
by an element of arity $n$ of $T$, and the output and all leaves of
$\Tfr$ are labeled on $V \sqcup T(0)$.
\medbreak

A \Def{regular tree grammar}~\cite{GS84,CDGJLLTT07} is a tuple
$\Gca := (V, T, P, s)$ where $(V, T)$ is a signature, $P$ is a set of
pairs of the form $(v, \Sfr)$ called \Def{productions} where $v \in V$
and $\Sfr$ is a $(V, T)$-tree, and $s$ is a variable of $V$ called
\Def{start symbol}. If $\Tfr_1$ and $\Tfr_2$ are two $(V, T)$-trees,
$\Tfr_2$ is \Def{derivable in one step} from $\Tfr_1$ if $\Tfr_1$ has a
leaf $y$ labeled by $a$ and the tree obtained by replacing $y$ by the
root of $\Sfr$ in $\Tfr_1$ is $\Tfr_2$, provided that $(a, \Sfr)$ is a
production of $P$. This property is denoted by $\Tfr_1 \Deriv \Tfr_2$,
so that $\Deriv$ is a binary relation on the set of all $(V, T)$-trees.
The reflexive and transitive closure of $\Deriv$ is the derivation
relation. A $(V, T)$-tree $\Tfr$ is \Def{generated} by $\Gca$ if $\Tfr$
is derivable from the tree $\Unit_s$ consisting in one leaf labeled by
$s$ and all leaves of $\Tfr$ are labeled on $T(0)$. The \Def{language}
of $\Gca$ is the set of all $(V, T)$-trees generated by~$\Gca$.
\medbreak

If $\Gca := (V, T, P, s)$ is a regular tree grammar, we denote by
$\RTG(\Gca)$ the bud generating system
\begin{equation}
    \RTG(\Gca) :=
    \left(\FreeOperad(T \setminus T(0)), V \sqcup T(0),
    \RFr, \{s\}, T(0)\right)
\end{equation}
wherein $\RFr$ is the set of rules
\begin{equation}
    \RFr :=
    \left\{
    (a, \Tfr, u) \in \Bud_{V \sqcup T(0)}(\FreeOperad(T \setminus T(0)))
    : (a, \Tfr_{a, u}) \in P \right\},
\end{equation}
where, for any $\Tfr \in \FreeOperad(T \setminus T(0))$,
$a \in V \sqcup T(0)$, and $u \in (V \sqcup T(0))^{|\Tfr|}$,
$\Tfr_{a, u}$ is the $(V, T)$-tree obtained by labeling the output of
$\Tfr$ by $a$ and by labeling from left to right the leaves of $\Tfr$ by
the letters of~$u$.
\medbreak

\begin{Proposition} \label{prop:emulation_tree_grammars}
    Let $\Gca$ be a regular tree grammar. Then, the map
    $\phi : \Lang(\RTG(\Gca)) \to L$ defined by
    $\phi((a, \Tfr, u)) := \Tfr_{a, u}$ is a bijection between the
    language of $\RTG(\Gca)$ and the language $L$ of~$\Gca$.
\end{Proposition}
\medbreak

\subsubsection{Synchronous grammars}
\label{subsubsec:synchronous_grammars}
In this section, we shall denote by $\Tree$ the monochrome operad
$\FreeOperad(C)$ where $C$ is the monochrome collection
$C := \sqcup_{n \geq 1} C(n)$ where $C(n) := \{\Asf_n\}$. The elements
of this operad are planar rooted trees where internal nodes have an
arbitrary arity. Observe that since $C(1) := \{\Asf_1\}$,
$\FreeOperad(C)(1)$ is an infinite set, so that $\FreeOperad(C)$ is not
combinatorial.
\medbreak

Let $B$ be a finite alphabet. A \Def{$B$-bud tree} is an element of
$\Bud_B(\Tree)$. In other words, a $B$-bud tree is a planar rooted
tree $\Tfr$ such that the output and all leaves of $\Tfr$ are labeled on
$B$. The leaves of a $B$-bud tree are indexed from $1$ from left to
right.
\medbreak

A \Def{synchronous grammar}~\cite{Gir12e} is a tuple $\Gca := (B, a, R)$
where $B$ is a finite alphabet of \Def{bud labels}, $a$ is an element of
$B$ called \Def{axiom}, and $R$ is a finite set of pairs of the form
$(b, \Sfr)$ called \Def{substitution rules} where $b \in B$ and $\Sfr$
is a $B$-bud tree. If $\Tfr_1$ and $\Tfr_2$ are two $B$-bud trees such
that $\Tfr_1$ is of arity $n$, $\Tfr_2$ is \Def{derivable in one step}
from $\Tfr_1$ if there are substitution rules
$(b_1, \Sfr_1), \dots, (b_n, \Sfr_n)$ of $R$ such that for all
$i \in [n]$, the $i$th leaf of $\Tfr_1$ is labeled by $b_i$ and $\Tfr_2$
is obtained by replacing the $i$th leaf of $\Tfr_1$ by $\Sfr_i$ for all
$i \in [n]$. This property is denoted by $\Tfr_1 \SyncDeriv \Tfr_2$,
so that $\SyncDeriv$ is a binary relation on the set of all $B$-bud
trees. The reflexive and transitive closure of $\SyncDeriv$ is the
derivation relation. A $B$-bud tree $\Tfr$ is \Def{generated} by $\Gca$
if $\Tfr$ is derivable from the tree $\Unit_a$ consisting is one leaf
labeled by $a$. The \Def{language} of $\Gca$ is the set of all
$B$-bud trees generated by~$\Gca$.
\medbreak

If $\Gca := (B, a, R)$ is a synchronous grammar, we denote by
$\SG(\Gca)$ the bud generating system
\begin{equation}
    \SG(\Gca) := (\Tree, B, \RFr, \{a\}, B)
\end{equation}
wherein $\RFr$ is the set of rules
\begin{equation}
    \RFr :=
    \left\{
    (b, \Tfr, u) \in \Bud_B(\Tree) : \left(b, \Tfr_{b, u}\right) \in R
    \right\},
\end{equation}
where, for any $\Tfr \in \Bud_B(\Tree)$, $b \in B$, and $u \in B^+$,
$\Tfr_{b, u}$ is the $B$-bud tree obtained by labeling the output of
$\Tfr$ by $b$ and by labeling from left to right the leaves of $\Tfr$ by
the letters of~$u$.
\medbreak

\begin{Proposition} \label{prop:emulation_synchronous_grammars}
    Let $\Gca$ be a synchronous grammar. Then, the map
    $\phi : \SyncLang(\SG(\Gca)) \to L$ defined by
    $\phi((b, \Tfr, u)) := \Tfr_{b, u}$ is a bijection between the
    synchronous language of $\SG(\Gca)$ and the language $L$ of~$\Gca$.
\end{Proposition}
\medbreak

\section{Series on colored operads and bud generating systems}
\label{sec:series_bud_generating_systems}
In this section, we explain how to use bud generating systems as tools
to enumerate families of combinatorial objects. For this purpose, we
will define and consider three series on colored operads extracted from
bud generating systems. Each of these series brings information about
the languages or the synchronous languages of bud generating systems.
One of a key issues is, given a bud generating system $\Bca$, to count
arity by arity the elements of the language or the synchronous language
of $\Bca$. In other terms, this amounts to compute the generating series
$\GenS_{\Lang(\Bca)}$ or $\GenS_{\SyncLang(\Bca)}$. As we shall see,
these generating series can be computed from the series of colored
operads extracted from~$\Bca$.
\medbreak

\subsection{General definitions}
Let us list some notations used in this section. In what follows,
$\Bca := (\Oca, \CFr, \RFr, I, T)$ is a bud generating system such that
$\Oca$ is a combinatorial monochrome operad and, as before, $\CFr$ is
a set of colors of the form $\CFr = \{c_1, \dots, c_k\}$.
\medbreak

\subsubsection{Characteristic series}
We shall denote by $\Rbf$ the characteristic series of $\RFr$, by $\Ibf$
the series
\begin{equation}
    \Ibf := \sum_{a \in I} \Unit_a,
\end{equation}
and by $\Tbf$ the series
\begin{equation}
    \Tbf := \sum_{a \in T} \Unit_a.
\end{equation}
\medbreak

\begin{Lemma} \label{lem:initial_terminal_composition}
    Let $\Bca := (\Oca, \CFr, \RFr, I, T)$ be a bud generating system
    and $\Fbf$ be a $\Bud_\CFr(\Oca)$-series. Then, for all
    $x \in \Bud_\CFr(\Oca)$,
    \begin{equation}
        \Angle{x, \Ibf \Compo \Fbf \Compo \Tbf} =
        \begin{cases}
            \Angle{x, \Fbf}
                & \mbox{if } \Out(x) \in I \mbox{ and }
                \In(x) \in T^+, \\
            0 & \mbox{otherwise}.
        \end{cases}
    \end{equation}
\end{Lemma}
\medbreak

\subsubsection{Polynomials}
For all colors $a \in \CFr$ and types $\alpha \in \Types_\CFr$, let
\begin{equation}
    \MultOutIn_{a, \alpha} :=
    \# \left\{r \in \RFr :
    (\Out(r), \Type(\In(r))) = (a, \alpha)\right\}.
\end{equation}
For any $a \in \CFr$, let $\Gbf_a(\VarY_{c_1}, \dots, \VarY_{c_k})$ be
the series of $\K \AAngle{\Multiset(\Ybb_\CFr)}$ defined by
\begin{equation}
    \Gbf_a(\VarY_{c_1}, \dots, \VarY_{c_k}) :=
    \sum_{\gamma \in \Types_\CFr}
    \MultOutIn_{a, \gamma} \enspace \Ybb_\CFr^\gamma.
\end{equation}
Notice that
\begin{equation}
    \Gbf_a(\VarY_{c_1}, \dots, \VarY_{c_k}) =
    \sum_{\substack{
        r \in \RFr \\
        \Out(r) = a
    }}
    \Ybb_\CFr^{\Type(\In(r))}
\end{equation}
and that, since $\RFr$ is finite, this series is a polynomial.
\medbreak

\subsubsection{Maps}
In the sequel, we shall use maps $\phi : \CFr \times \Types_\CFr \to \N$
such that $\phi(a, \gamma) \ne 0$ for a finite number of pairs
$(a, \gamma) \in \CFr \times \Types_\CFr$, to express in a concise
manner some recurrence relations for the coefficients of series on
colored operads. We shall consider the two following notations. If
$\phi$ is such a map and $a \in \CFr$, we define $\phi^{(a)}$ as the
natural number
\begin{equation}
    \phi^{(a)} :=
    \sum_{\substack{
        b \in \CFr \\
        \gamma \in \Types_\CFr
    }}
    \phi(b, \gamma) \gamma_a
\end{equation}
and $\phi_a$ as the finite multiset
\begin{equation}
    \phi_a := \lbag \phi(a, \gamma) : \gamma \in \Types_\CFr \rbag.
\end{equation}
\medbreak

\subsection{Hook generating series}
We call \Def{hook generating series} of $\Bca$ the
$\Bud_\CFr(\Oca)$-series $\Hook(\Bca)$ defined by
\begin{equation} \label{equ:definition_hook_series}
    \Hook(\Bca) :=
    \Ibf \Compo \Rbf^{\PreLieProduct_*} \Compo \Tbf.
\end{equation}
Observe that~\eqref{equ:definition_hook_series} could be undefined
for an arbitrary set of rules $\RFr$ of $\Bca$. Nevertheless, when
$\Rbf$ satisfies the conditions of
Proposition~\ref{prop:coefficients_pre_lie_star}, that is, when $\Oca$
is a combinatorial operad and $\RFr(1)$ finitely factorizes
$\Bud_\CFr(\Oca)(1)$, $\Hook(\Bca)$ is well-defined.
\medbreak

\subsubsection{Expression}
The aim of the following is to provide an expression to compute the
coefficients of~$\Hook(\Bca)$.
\medbreak

\begin{Lemma} \label{lem:pre_lie_star_hook_series}
    Let $\Bca := (\Oca, \CFr, \RFr, I, T)$ be a bud generating system
    such that $\Oca$ is a combinatorial operad and $\RFr(1)$ finitely
    factorizes $\Bud_\CFr(\Oca)(1)$. Then, for any
    $x \in \Bud_\CFr(\Oca)$,
    \begin{equation} \label{equ:pre_lie_star_hook_series}
        \Angle{x, \Rbf^{\PreLieProduct_*}} =
        \delta_{x, \Unit_{\Out(x)}} +
        \sum_{\substack{
            y \in \Bud_\CFr(\Oca) \\
            z \in \RFr \\
            i \in [|y|] \\
            x = y \circ_i z}}
        \Angle{y, \Rbf^{\PreLieProduct_*}}.
    \end{equation}
\end{Lemma}
\medbreak

\begin{Proposition} \label{prop:hook_series_derivation_graph}
    Let $\Bca := (\Oca, \CFr, \RFr, I, T)$ be a bud generating system
    such that $\Oca$ is a combinatorial operad and $\RFr(1)$ finitely
    factorizes $\Bud_\CFr(\Oca)(1)$. Then, for any
    $x \in \Bud_\CFr(\Oca)$ such that $\Out(x) \in I$, the coefficient
    $\Angle{x, \Rbf^{\PreLieProduct_*}}$ is the number of multipaths
    from $\Unit_{\Out(x)}$ to $x$ in the derivation graph of~$\Bca$.
\end{Proposition}
\begin{proof}
    First, since $\RFr(1)$ finitely factorizes $\Bud_\CFr(\Oca)(1)$, by
    Proposition~\ref{prop:coefficients_pre_lie_star},
    $\Rbf^{\PreLieProduct_*}$ is a well-defined series. If $x = \Unit_a$
    for a $a \in I$, since
    $\Angle{\Unit_a, \Rbf^{\PreLieProduct_*}} = 1$, the statement of the
    proposition holds. Let us now assume that $x$ is different from a
    colored unit and let us denote by $\lambda_x$ the number of
    multipaths from $\Unit_{\Out(x)}$ to $x$ in the derivation graph
    $\DerivGraph(\Bca)$ of $\Bca$. By definition of $\DerivGraph(\Bca)$,
    by denoting by $\mu_{y, x}$ the number of edges from
    $y \in \Bud_\CFr(\Oca)$ to $x$ in $\DerivGraph(\Bca)$, we have
    \begin{equation}\begin{split}
        \label{equ:hook_series_derivation_graph}
        \lambda_x
        & = \sum_{y \in \Bud_\CFr(\Oca)} \mu_{y, x}\, \lambda_y \\
        & =
        \sum_{y \in \Bud_\CFr(\Oca)}
        \# \left\{(i, r) \in \N \times \RFr : x = y \circ_i r\right\}\,
        \lambda_y \\
        & =
        \sum_{\substack{
            y \in \Bud_\CFr(\Oca) \\
            i \in [|y|] \\
            r \in \RFr \\
            x = y \circ_i r
        }}
        \lambda_y.
    \end{split}\end{equation}
    We observe that Relation~\eqref{equ:hook_series_derivation_graph}
    satisfied by the $\lambda_x$ is the same as
    Relation~\eqref{equ:pre_lie_star_hook_series} of in the statement of
    Lemma~\ref{lem:pre_lie_star_hook_series} satisfied by the
    $\Angle{x, \Rbf^{\PreLieProduct_*}}$. This implies the statement of
    the proposition.
\end{proof}
\medbreak

\begin{Theorem} \label{thm:hook_series}
    Let $\Bca := (\Oca, \CFr, \RFr, I, T)$ be a bud generating system
    such that $\Oca$ is a combinatorial operad and $\RFr(1)$ finitely
    factorizes $\Bud_\CFr(\Oca)(1)$. Then, the hook generating series of
    $\Bca$ satisfies
    \begin{equation} \label{equ:hook_series}
        \Hook(\Bca)
        = \sum_{\substack{
            \Tfr \in \FreeColoredOperad(\RFr) \\
            \Out(\Tfr) \in I \\
            \In(\Tfr) \in T^+
        }}
        \enspace
        \frac{\deg(\Tfr)!}
            {\prod\limits_{v \in \TreeLanguage_\Node(\Tfr)}
            \deg(\Tfr_v)}
        \enspace
        \Eval(\Tfr).
    \end{equation}
\end{Theorem}
\begin{proof}
    By definition of $\Lang(\Bca)$ and $\DerivGraph(\Bca)$, any
    $x \in \Lang(\Bca)$ can be reached from $\Unit_{\Out(x)}$ by a
    multipath
    \begin{equation}
        \Unit_{\Out(x)} \Deriv y_1 \Deriv y_2 \Deriv \cdots
        \Deriv y_{\ell - 1} \Deriv x
    \end{equation}
    in $\DerivGraph(\Bca)$, where $y_1, \dots, y_{\ell - 1}$ are
    elements of $\Bud_\CFr(\Oca)$ and $\Unit_{\Out(x)} \in I$. Hence,
    by definition of $\Deriv$, $x$ admits an $\RFr$-left expression
    \begin{equation}
        x = (\dots ((\Unit_{\Out(x)} \circ_1 r_1) \circ_{i_1} r_2)
        \circ_{i_2} \dots) \circ_{i_{\ell - 1}} r_\ell
    \end{equation}
    where for any $j \in [\ell]$, $r_j \in \RFr$, and for any
    $j \in [\ell - 1]$,
    \begin{equation}
        y_j = (\dots ((\Unit_{\Out(x)} \circ_1 r_1) \circ_{i_1} r_2)
        \circ_{i_2} \dots) \circ_{i_{j - 1}} r_j
    \end{equation}
    and $i_{j} \in [|y_j|]$. This shows that the set of all multipaths
    from $\Unit_{\Out(x)}$ to $x$ in $\DerivGraph(\Bca)$ is in
    one-to-one correspondence with the set of all $\RFr$-left
    expressions for $x$. Now, observe that since $\RFr(1)$ finitely
    factorizes $\Bud_\CFr(\Oca)(1)$, by
    Proposition~\ref{prop:coefficients_pre_lie_star}
    $\Rbf^{\PreLieProduct_*}$ is a well-defined series. By
    Proposition~\ref{prop:hook_series_derivation_graph},
    Lemmas~\ref{lem:finitely_factorizing_finite_treelike_expressions}
    and~\ref{lem:left_expressions_linear_extensions},
    and~\eqref{equ:number_left_expressions} of
    Chapter~\ref{chap:algebra}, we obtain that
    \begin{equation} \label{equ:hook_series_demo}
        \Angle{x, \Rbf^{\PreLieProduct_*}} =
        \sum_{\substack{
            \Tfr \in \FreeColoredOperad(\RFr) \\
            \Eval(\Tfr) = x
        }}
        \frac{\deg(\Tfr)!}
            {\prod\limits_{v \in \TreeLanguage_\Node(\Tfr)}
            \deg(\Tfr_v)}.
    \end{equation}
    Finally, by Lemma~\ref{lem:initial_terminal_composition}, for any
    $x \in \Bud_\CFr(\Oca)$ such that $\Out(x) \in I$ and
    $\In(x) \in T^+$, we have
    $\Angle{x, \Hook(\Bca)} = \Angle{x, \Rbf^{\PreLieProduct_*}}$. This
    shows that the right member of~\eqref{equ:hook_series} is equal
    to~$\Hook(\Bca)$.
\end{proof}
\medbreak

An alternative way to understand $\Hook(\Bca)$ hence offered by
Theorem~\ref{thm:hook_series} consists is seeing the coefficient
$\Angle{x, \Hook(\Bca)}$, $x \in \Bud_\CFr(\Oca)$, as the number of
$\RFr$-left expressions of~$x$.
\medbreak

\subsubsection{Support}
The following result establishes a link between the hook generating
series of $\Bca$ and its language.
\medbreak

\begin{Proposition} \label{prop:support_hook_series_language}
    Let $\Bca := (\Oca, \CFr, \RFr, I, T)$ be a bud generating system
    such that $\Oca$ is a combinatorial operad and $\RFr(1)$ finitely
    factorizes $\Bud_\CFr(\Oca)(1)$. Then, the support of the hook
    generating series of $\Bca$ is the language of~$\Bca$.
\end{Proposition}
\medbreak

\subsubsection{Analogs of the hook statistics}
\label{subsubsec:hook_bud_generating_systems}
Bud generating systems lead to the definition of analogues of the
hook-length statistics~\cite{Knu98} for combinatorial objects
possibly different than trees in the following way. Let $\Oca$ be a
monochrome operad, $\GeneratingSet$ be a generating set of $\Oca$, and
$\HookSystem_{\Oca, \GeneratingSet} := (\Oca, \GeneratingSet)$ be a
monochrome bud generating system depending on $\Oca$ and
$\GeneratingSet$, called \Def{hook bud generating system}. Since
$\GeneratingSet$ is a generating set of $\Oca$, by
Propositions~\ref{prop:language_sub_operad_generated}
and~\ref{prop:support_hook_series_language}, the support of
$\Hook\left(\HookSystem_{\Oca, \GeneratingSet}\right)$ is equal to
$\Lang\left(\HookSystem_{\Oca, \GeneratingSet}\right)$. We define the
\Def{hook-length coefficient} of any element $x$ of $\Oca$ as the
coefficient
$\Angle{x, \Hook\left(\HookSystem_{\Oca, \GeneratingSet}\right)}$.
\medbreak

Let us consider the hook bud generating system
$\HookSystem_{\Mag, \GeneratingSet}$ where $\Mag$ is the magmatic operad
(whose definition is recalled in Section~\ref{subsubsec:magmatic_operad}
of Chapter~\ref{chap:algebra}) and
\begin{equation}
    \GeneratingSet := \left\{\CorollaTwo{}\right\}.
\end{equation}
This bud generating system leads to the definition of a statistics on
binary trees, provided by the coefficients of the hook generating series
$\Hook\left(\HookSystem_{\Mag, G}\right)$ which begins by
\begin{multline}
    \Hook\left(\HookSystem_{\Mag, G}\right) =

    \enspace + \cdots.
\end{multline}
Theorem~\ref{thm:hook_series} implies that for any binary tree $\Tfr$,
the coefficient
$\Angle{\Tfr, \Hook\left(\HookSystem_{\Mag, \GeneratingSet}\right)}$
can be obtained by the usual hook-length formula of binary trees.
Alternatively, the coefficient
$\Angle{\Tfr, \Hook\left(\HookSystem_{\Mag, \GeneratingSet}\right)}$ is
the cardinal of the sylvester class~\cite{HNT05} of permutations encoded
by~$\Tfr$. This explains the name of hook generating series for
$\Hook(\Bca)$, when $\Bca$ is a bud generating system.
\medbreak

Consider now a second example of a hook generating system involving the
operad $\Motz$ of Motkzin paths (see Section~\ref{subsubsec:operad_Motz}
of Chapter~\ref{chap:monoids}) seen as a set-operad. From its
definition,
\begin{equation}
    \GeneratingSet := \left\{\MotzHoriz, \MotzPeak\right\}
\end{equation}
is a generating set of $\Motz$. Hence,
$\HookSystem_{\Motz, \GeneratingSet}$ is a hook generating system. This
leads to the definition of a statistics on Motzkin paths, provided by
the coefficients of the hook generating series
$\Hook\left(\HookSystem_{\Motz, \GeneratingSet}\right)$ of
$\HookSystem_{\Motz, \GeneratingSet}$ which begins by
\begin{multline}
    \Hook\left(\HookSystem_{\Motz, \GeneratingSet}\right) =

    \enspace + \cdots.
\end{multline}
\medbreak

\subsection{Syntactic generating series}
\label{subsec:syntactic_generating_series}
We call \Def{syntactic generating series} of $\Bca$ the
$\Bud_\CFr(\Oca)$-series $\Synt(\Bca)$ defined by
\begin{equation} \label{equ:definition_syntactic_series}
    \Synt(\Bca) :=
    \Ibf \Compo (\Ubf - \Rbf)^{\Compo_{-1}} \Compo \Tbf.
\end{equation}
Observe that~\eqref{equ:definition_syntactic_series} could be undefined
for an arbitrary set of rules $\RFr$ of $\Bca$. Nevertheless, when
$\Ubf - \Rbf$ satisfies the conditions of
Proposition~\ref{prop:composition_inverse_expression}, $\Synt(\Bca)$ is
well-defined. Remark that this condition is satisfied whenever
$\Oca$ is combinatorial and $\RFr(1)$ finitely
factorizes~$\Bud_\CFr(\Oca)(1)$.
\medbreak

\subsubsection{Expression}
The aim of this section is to provide an expression to compute the
coefficients of~$\Synt(\Bca)$.
\medbreak

\begin{Lemma} \label{lem:inverse_composition_syntactic_series}
    Let $\Bca := (\Oca, \CFr, \RFr, I, T)$ be a bud generating system
    such that $\Oca$ is a combinatorial operad and $\RFr(1)$ finitely
    factorizes $\Bud_\CFr(\Oca)(1)$. Then, for any
    $x \in \Bud_\CFr(\Oca)$,
    \begin{equation} \label{equ:inverse_composition_syntactic_series}
        \Angle{x, (\Ubf - \Rbf)^{\Compo_{-1}}} =
        \delta_{x, \Unit_{\Out(x)}} +
        \sum_{\substack{
            y \in \RFr \\
            z_1, \dots, z_{|y|} \in \Bud_\CFr(\Oca) \\
            x = y \circ \left[z_1, \dots, z_{|y|}\right]
        }}
        \enspace
        \prod_{i \in [|y|]} \Angle{z_i, (\Ubf - \Rbf)^{\Compo_{-1}}}.
    \end{equation}
\end{Lemma}
\medbreak

\begin{Theorem} \label{thm:syntactic_series}
    Let $\Bca := (\Oca, \CFr, \RFr, I, T)$ be a bud generating system
    such that $\Oca$ is a combinatorial operad and $\RFr(1)$ finitely
    factorizes $\Bud_\CFr(\Oca)(1)$. Then, the syntactic generating
    series of $\Bca$ satisfies
    \begin{equation} \label{equ:syntactic_series}
        \Synt(\Bca)
        =
        \sum_{\substack{
            \Tfr \in \FreeColoredOperad(\RFr) \\
            \Out(\Tfr) \in I \\
            \In(\Tfr) \in T^+
        }}
        \Eval(\Tfr).
    \end{equation}
\end{Theorem}
\medbreak

Theorem~\ref{thm:syntactic_series} explains the name of syntactic
generating series for $\Synt(\Bca)$ because this series can be expressed
following~\eqref{equ:syntactic_series} as a sum of evaluations of syntax
trees. An alternative way to see $\Synt(\Bca)$ is that for any
$x \in \Bud_\CFr(\Oca)$, the coefficient $\Angle{x, \Synt(\Bca)}$ is the
number of $\RFr$-treelike expressions for~$x$.
\medbreak

\subsubsection{Support and unambiguity}
The following result establishes a link between the syntactic generating
series of $\Bca$ and its language.
\medbreak

\begin{Proposition} \label{prop:support_syntactic_series_language}
    Let $\Bca := (\Oca, \CFr, \RFr, I, T)$ be a bud generating system
    such that $\Oca$ is a combinatorial operad and $\RFr(1)$ finitely
    factorizes $\Bud_\CFr(\Oca)(1)$. Then, the support of the syntactic
    generating series of $\Bca$ is the language of~$\Bca$.
\end{Proposition}
\medbreak

We rely now on syntactic generating series to define a property of bud
generating systems. We say that $\Bca$ is \Def{unambiguous} if all
coefficients of $\Synt(\Bca)$ are equal to $0$ or to $1$. This property
is important from a combinatorial point of view. Indeed, by definition
of the series of colors $\Colors$ (see
Section~\ref{subsubsec:series_of_colors}) and
Proposition~\ref{prop:support_syntactic_series_language}, when $\Bca$ is
unambiguous, the coefficient of $(a, u) \in \Bud_\CFr(\As)$ in the
series $\Colors(\Synt(\Bca))$ is the number of elements $x$ of
$\Lang(\Bca)$ such that $(\Out(x), \In(x)) = (a, u)$.
\medbreak

For instance, consider the bud generating system $\BMotz$ introduced in
Section~\ref{subsubsec:example_bmotz}. Observe that since the Motzkin
path
\begin{math}

\end{math}
of $\Motz(5)$ admits exactly the two $\RFr$-treelike expressions
\vspace{-1.75em}
\begin{multicols}{2}
\begin{subequations}
\begin{equation}
    %
,
    11111\right),
    \Synt(\BMotz)}
    = 2.
\end{math}
Hence $\BMotz$ is not unambiguous.
\medbreak

As a side remark, observe that Theorem~\ref{thm:syntactic_series}
implies in particular that for any bud generating system of the form
$\Bca := (\Oca, \CFr, \RFr, \CFr, \CFr)$, if $\Synt(\Bca)$ is
unambiguous, then the colored suboperad of $\Bud_\CFr(\Oca)$ generated
by $\RFr$ is free. The converse property does not hold.
\medbreak

\subsubsection{Series of color types}
The purpose of this section is to describe the coefficients of
$\ColorTypes(\Synt(\Bca))$, the series of color types of the syntactic
series of $\Bca$, in the particular case when $\Bca$ is unambiguous.
We shall give two descriptions: a first one involving a system of
equations of series of $\K \AAngle{\Multiset(\Ybb_\CFr)}$, and a second
one involving a recurrence relation on the coefficients of a series
of~$\K \AAngle{\Multiset(\Xbb_\CFr + \Ybb_\CFr)}$.
\medbreak

\begin{Lemma} \label{lem:colt_synt_coefficients_description}
    Let $\Bca := (\Oca, \CFr, \RFr, I, T)$ be an unambiguous bud
    generating system such that $\Oca$ is a combinatorial operad and
    $\RFr(1)$ finitely factorizes $\Bud_\CFr(\Oca)(1)$. Then, for all
    colors $a \in I$ and all types $\alpha \in \Types_\CFr$ such that
    $\CFr^\alpha \in T^+$, the coefficients
    $\Angle{\VarX_a \Ybb_\CFr^\alpha, \ColorTypes(\Synt(\Bca))}$ count
    the number of elements $x$ of $\Lang(\Bca)$ such that
    $(\Out(x), \Type(\In(x))) = (a, \alpha)$.
\end{Lemma}
\medbreak

\begin{Proposition} \label{prop:functional_equation_synt}
    Let $\Bca := (\Oca, \CFr, \RFr, I, T)$ be an unambiguous bud
    generating system such that $\Oca$ is a combinatorial operad and
    $\RFr(1)$ finitely factorizes $\Bud_\CFr(\Oca)(1)$. For all
    $a \in \CFr$, let $\Fbf_a(\VarY_{c_1}, \dots, \VarY_{c_k})$ be the
    series of $\K \AAngle{\Multiset(\Ybb_\CFr)}$ satisfying
    \begin{equation}
        \Fbf_a\left(\VarY_{c_1}, \dots, \VarY_{c_k}\right) =
        \VarY_a +
        \Gbf_a\left(\Fbf_{c_1}\left(\VarY_{c_1}, \dots,
            \VarY_{c_k}\right),
        \dots,
        \Fbf_{c_k}\left(\VarY_{c_1}, \dots, \VarY_{c_k}\right)\right).
    \end{equation}
    Then, for any color $a \in I$ and any type $\alpha \in \Types_\CFr$
    such that $\CFr^\alpha \in T^+$, the coefficients
    $\Angle{\VarX_a \Ybb_\CFr^\alpha, \ColorTypes(\Synt(\Bca))}$ and
    $\Angle{\Ybb_\CFr^\alpha, \Fbf_a}$ are equal.
\end{Proposition}
\medbreak

\begin{Theorem} \label{thm:series_color_types_synt}
    Let $\Bca := (\Oca, \CFr, \RFr, I, T)$ be an unambiguous bud
    generating system such that $\Oca$ is a combinatorial operad and
    $\RFr(1)$ finitely factorizes $\Bud_\CFr(\Oca)(1)$. Let $\Fbf$ be
    the series of $\K \AAngle{\Multiset(\Xbb_\CFr + \Ybb_\CFr)}$
    satisfying, for any $a \in \CFr$ and any type
    $\alpha \in \Types_\CFr$,
    \begin{equation} \label{equ:series_color_types_synt}
        \Angle{\VarX_a \Ybb_\CFr^\alpha, \Fbf}
        =
        \delta_{\alpha, \Type(a)} +
        \sum_{\substack{
            \phi : \CFr \times \Types_\CFr \to \N \\
            \alpha = \phi^{(c_1)} \dots \phi^{(c_k)}
        }}
        \MultOutIn_{a, \sum \phi_{c_1} \dots \sum \phi_{c_k}}
        \left(\prod_{b \in \CFr} \phi_b !\right)
        \left(
        \prod_{\substack{
            b \in \CFr \\
            \gamma \in \Types_\CFr
        }}
        \Angle{\VarX_b \Ybb_\CFr^\gamma, \Fbf}^{\phi(b, \gamma)}
        \right).
    \end{equation}
    Then, for any color $a \in I$ and any type $\alpha \in \Types_\CFr$
    such that $\CFr^\alpha \in T^+$, the coefficients
    $\Angle{\VarX_a \Ybb_\CFr^\alpha, \ColorTypes(\Synt(\Bca))}$ and
    $\Angle{\VarX_a \Ybb_\CFr^\alpha, \Fbf}$ are equal.
\end{Theorem}
\medbreak

\subsubsection{Generating series of languages}
When $\Bca$ is a bud generating system satisfying the conditions of
Proposition~\ref{prop:functional_equation_synt}, the generating
series of the language of $\Bca$ satisfies
\begin{equation} \label{equ:generating_series_synt_functional_equation}
    \GenS_{\Lang(\Bca)} = \sum_{a \in I} \Fbf_a^{T},
\end{equation}
where $\Fbf_a^{T}$ is the specialization of the series
$\Fbf_a(\VarY_{c_1}, \dots, \VarY_{c_k})$ at $\VarY_b := t$ for
all $b \in T$ and at $\VarY_c := 0$ for all $c \in \CFr \setminus T$.
Therefore, the resolution of the system of equations given by
Proposition~\ref{prop:functional_equation_synt} provides a way to
compute the coefficients of~$\GenS_{\Lang(\Bca)}$.
\medbreak

\begin{Theorem} \label{thm:algebraic_generating_series_languages}
    Let $\Bca := (\Oca, \CFr, \RFr, I, T)$ be an unambiguous bud
    generating system such that $\Oca$ is a combinatorial operad and
    $\RFr(1)$ finitely factorizes $\Bud_\CFr(\Oca)(1)$. Then, the
    generating series $\GenS_{\Lang(\Bca)}$ of the language of $\Bca$ is
    algebraic.
\end{Theorem}
\medbreak

When $\Bca$ is a bud generating system satisfying the conditions of
Theorem~\ref{thm:series_color_types_synt} (which are the same as the
ones required by Proposition~\ref{prop:functional_equation_synt}), one
has for any $n \geq 1$,
\begin{equation} \label{equ:generating_series_synt_recurrence}
    \Angle{t^n, \GenS_{\Lang(\Bca)}}
    =
    \sum_{a \in I}
    \sum_{\substack{
        \alpha \in \Types_\CFr \\
        \alpha_i = 0, c_i \in \CFr \setminus T
    }}
    \Angle{\VarX_a \Ybb_\CFr^\alpha, \Fbf}.
\end{equation}
Therefore, this provides an alternative and recursive way to compute the
coefficients of $\GenS_{\Lang(\Bca)}$, different from the one of
Proposition~\ref{prop:functional_equation_synt}.
\medbreak

\subsection{Synchronous generating series}
\label{subsec:synchronous_generating_series}
We call \Def{synchronous generating series} of $\Bca$ the
$\Bud_\CFr(\Oca)$-series $\Sync(\Bca)$ defined by
\begin{equation} \label{equ:definition_synchronous_series}
    \Sync(\Bca) := \Ibf \Compo \Rbf^{\Compo_*} \Compo \Tbf.
\end{equation}
Observe that~\eqref{equ:definition_synchronous_series} could be
undefined for an arbitrary set of rules $\RFr$ of $\Bca$. Nevertheless,
when $\Rbf$ satisfies the conditions of
Proposition~\ref{prop:coefficients_composition_star}, that is, when
$\Oca$ is a combinatorial operad and $\RFr(1)$ finitely factorizes
$\Bud_\CFr(\Oca)(1)$, $\Sync(\Bca)$ is well-defined.
\medbreak

\subsubsection{Expression}
The aim of this section is to provide an expression to compute the
coefficients of~$\Sync(\Bca)$.
\medbreak

\begin{Lemma} \label{lem:composition_star_synchronous_series}
    Let $\Bca := (\Oca, \CFr, \RFr, I, T)$ be a bud generating system
    such that $\Oca$ is a combinatorial operad and $\RFr(1)$
    finitely factorizes $\Bud_\CFr(\Oca)(1)$. Then, for any
    $x \in \Bud_\CFr(\Oca)$,
    \begin{equation} \label{equ:composition_star_synchronous_series}
        \Angle{x, \Rbf^{\Compo_*}} =
        \delta_{x, \Unit_{\Out(x)}} +
        \sum_{\substack{
            y \in \Bud_\CFr(\Oca) \\
            z_1, \dots, z_{|y|} \in \RFr \\
            x = y \circ \left[z_1, \dots, z_{|y|}\right]
        }}
        \Angle{y, \Rbf^{\Compo_*}}.
    \end{equation}
\end{Lemma}
\medbreak

\begin{Theorem} \label{thm:synchronous_series}
    Let $\Bca := (\Oca, \CFr, \RFr, I, T)$ be a bud generating system
    such that $\Oca$ is a combinatorial operad and $\RFr(1)$ finitely
    factorizes $\Bud_\CFr(\Oca)(1)$. Then, the synchronous generating
    series of $\Bca$ satisfies
    \begin{equation} \label{equ:synchronous_series}
        \Sync(\Bca)
        = \sum_{
            \substack{\Tfr \in \FreeColoredOperad_{\Perfect}(\RFr) \\
            \Out(\Tfr) \in I \\
            \In(\Tfr) \in T^+}}
        \Eval(\Tfr).
    \end{equation}
\end{Theorem}
\medbreak

Theorem~\ref{thm:synchronous_series} implies that for any
$x \in \Bud_\CFr(\Oca)$, the coefficient of $\Angle{x, \Sync(\Bca)}$ is
the number of $\RFr$-treelike expressions for $x$ which are perfect
trees.
\medbreak

\subsubsection{Support and unambiguity}
The following result establishes a link between the synchronous
generating series of $\Bca$ and its synchronous language.
\medbreak

\begin{Proposition} \label{prop:support_synchronous_series_language}
    Let $\Bca := (\Oca, \CFr, \RFr, I, T)$ be a bud generating system
    such that $\Oca$ is a combinatorial operad and $\RFr(1)$ finitely
    factorizes $\Bud_\CFr(\Oca)(1)$. Then, the support of the
    synchronous generating series of $\Bca$ is the synchronous language
    of~$\Bca$.
\end{Proposition}
\medbreak

We rely now on synchronous generating series to define a property of
bud generating systems. We say that $\Bca$ is \Def{synchronously
unambiguous} if all coefficients of $\Sync(\Bca)$ are equal to $0$ or to
$1$. This property is important from a combinatorial point of view.
Indeed, by definition of the series of colors $\Colors$ (see
Section~\ref{subsubsec:series_of_colors}) and
Proposition~\ref{prop:support_synchronous_series_language}, when $\Bca$
is synchronously unambiguous, the coefficient of
$(a, u) \in \Bud_\CFr(\As)$ in the series $\Colors(\Sync(\Bca))$ is the
number of elements $x$ of $\SyncLang(\Bca)$ such that
$(\Out(x), \In(x)) = (a, u)$.
\medbreak

For instance, the bud generating system $\BBalTree$ introduced in
Section~\ref{subsubsec:example_bbaltree} is synchronously unambiguous.
\medbreak

\subsubsection{Series of color types}
The purpose of this section is to describe the coefficients of
$\ColorTypes(\Sync(\Bca))$, the series of color types of the
synchronous series of $\Bca$, in the particular case when $\Bca$ is
unambiguous. We shall give two descriptions: a first one involving
a system of functional equations of series of
$\K \AAngle{\Multiset(\Ybb_\CFr)}$, and a second one involving a
recurrence relation on the coefficients of a series of
$\K \AAngle{\Multiset(\Xbb_\CFr + \Ybb_\CFr)}$.
\medbreak

\begin{Lemma} \label{lem:colt_sync_coefficients_description}
    Let $\Bca := (\Oca, \CFr, \RFr, I, T)$ be a synchronously
    unambiguous bud generating system such that $\Oca$ is a
    combinatorial operad and $\RFr(1)$ finitely factorizes
    $\Bud_\CFr(\Oca)(1)$. Then, for all colors $a \in I$ and all types
    $\alpha \in \Types_\CFr$ such that $\CFr^\alpha \in T^+$, the
    coefficients
    $\Angle{\VarX_a \Ybb_\CFr^\alpha, \ColorTypes(\Sync(\Bca))}$ count
    the number of elements $x$ of $\SyncLang(\Bca)$ such that
    $(\Out(x), \Type(\In(x))) = (a, \alpha)$.
\end{Lemma}
\medbreak

\begin{Proposition} \label{prop:functional_equation_sync}
    Let $\Bca := (\Oca, \CFr, \RFr, I, T)$ be a synchronously
    unambiguous bud generating system such that $\Oca$ is a
    combinatorial operad and $\RFr(1)$ finitely factorizes
    $\Bud_\CFr(\Oca)(1)$. For all $a \in \CFr$, let
    $\Fbf_a(\VarY_{c_1}, \dots, \VarY_{c_k})$ be the series of
    $\K \AAngle{\Multiset(\Ybb_\CFr)}$ satisfying
    \begin{equation}
        \Fbf_a\left(\VarY_{c_1}, \dots, \VarY_{c_k}\right) =
        \VarY_a +
        \Fbf_a\left(\Gbf_{c_1}\left(\VarY_{c_1}, \dots,
        \VarY_{c_k}\right),
        \dots,
        \Gbf_{c_k}\left(\VarY_{c_1}, \dots, \VarY_{c_k}\right)\right).
    \end{equation}
    Then, for any color $a \in I$ and any type $\alpha \in \Types_\CFr$
    such that $\CFr^\alpha \in T^+$, the coefficients
    $\Angle{\VarX_a \Ybb_\CFr^\alpha, \ColorTypes(\Sync(\Bca))}$
    and  $\Angle{\Ybb_\CFr^\alpha, \Fbf_a}$ are equal.
\end{Proposition}
\medbreak

\begin{Theorem} \label{thm:series_color_types_sync}
    Let $\Bca := (\Oca, \CFr, \RFr, I, T)$ be a synchronously
    unambiguous bud generating system such that $\Oca$ is a
    combinatorial operad and $\RFr(1)$ finitely factorizes
    $\Bud_\CFr(\Oca)(1)$. Let $\Fbf$ be the series of
    $\K \AAngle{\Multiset(\Xbb_\CFr + \Ybb_\CFr)}$ satisfying, for any
    $a \in \CFr$ and any type $\alpha \in \Types_\CFr$,
    \begin{equation} \label{equ:series_color_types_sync}
        \Angle{\VarX_a \Ybb_\CFr^\alpha, \Fbf}
        =
        \delta_{\alpha, \Type(a)} +
        \sum_{\substack{
            \phi : \CFr \times \Types_\CFr \to \N \\
            \alpha = \phi^{(c_1)} \dots \phi^{(c_k)}
        }}
        \left(
            \prod_{b \in \CFr}
            \phi_b !
        \right)
        \left(
            \prod_{\substack{
                b \in \CFr \\
                \gamma \in \Types_\CFr
            }}
            \MultOutIn_{b, \gamma}^{\phi(b, \gamma)}
        \right)
        \Angle{\VarX_a
        \prod_{b \in \CFr} \VarY_b^{\sum \phi_b},
        \Fbf}.
    \end{equation}
    Then, for any color $a \in I$ and any type $\alpha \in \Types_\CFr$
    such that $\CFr^\alpha \in T^+$, the coefficients
    $\Angle{\VarX_a \Ybb_\CFr^\alpha, \ColorTypes(\Sync(\Bca))}$ and
    $\Angle{\VarX_a \Ybb_\CFr^\alpha, \Fbf}$ are equal.
\end{Theorem}
\medbreak

\subsubsection{Generating series of synchronous languages}
\label{subsubsec:generating_series_synchronous_languages}
When $\Bca$ is a bud generating system satisfying the conditions of
Proposition~\ref{prop:functional_equation_sync}, the generating series
of the synchronous language of $\Bca$ satisfies
\begin{equation} \label{equ:generating_series_sync_functional_equation}
    \GenS_{\SyncLang(\Bca)} = \sum_{a \in I} \Fbf_a^{T},
\end{equation}
where $\Fbf_a^{T}$ is the specialization of the series
$\Fbf_a(\VarY_{c_1}, \dots, \VarY_{c_k})$ at $\VarY_b := t$ for all
$b \in T$ and at $\VarY_c := 0$ for all $c \in \CFr \setminus T$.
Therefore, the resolution of the system of equations given by
Proposition~\ref{prop:functional_equation_sync} provides a way to
compute the coefficients of $\GenS_{\SyncLang(\Bca)}$. This resolution
can be made in most cases by iteration~\cite{BLL98,FS09}.
\medbreak

Moreover, when $\Gca$ is a synchronous grammar~\cite{Gir12e} (see also
Section~\ref{subsubsec:synchronous_grammars} for a description of these
grammars) and when $\SG(\Gca) = \Bca$, the system of functional
equations provided by Proposition~\ref{prop:functional_equation_sync}
and~\eqref{equ:generating_series_sync_functional_equation} for
$\GenS_{\SyncLang(\Bca)}$ is the same as the one which can be extracted
from~$\Gca$.
\medbreak

When $\Bca$ is a bud generating system satisfying the conditions of
Theorem~\ref{thm:series_color_types_sync} (which are the same as the
ones required by Proposition~\ref{prop:functional_equation_sync}), one
has for any $n \geq 1$,
\begin{equation} \label{equ:generating_series_sync_recurrence}
    \Angle{t^n, \GenS_{\SyncLang(\Bca)}}
    =
    \sum_{a \in I}
    \sum_{\substack{
        \alpha \in \Types_\CFr \\
        \alpha_i = 0, c_i \in \CFr \setminus T
    }}
    \Angle{\VarX_a \Ybb_\CFr^\alpha, \Fbf}.
\end{equation}
Therefore, this provides an alternative and recursive way to compute the
coefficients of $\GenS_{\SyncLang(\Bca)}$, different from the one of
Proposition~\ref{prop:functional_equation_sync}.
\medbreak

\subsubsection{Example: enumeration of balanced binary trees}
\label{subsubsec:colt_sync_bbaltree}
Let us consider the bud generating system $\BBalTree$ introduced in
Section~\ref{subsubsec:example_bbaltree}. We have
\begin{equation}
    \MultOutIn_{a, \alpha} =
    \begin{cases}
        1 & \mbox{if } (a, \alpha) = (1, 20), \\
        2 & \mbox{if } (a, \alpha) = (1, 11), \\
        1 & \mbox{if } (a, \alpha) = (2, 10), \\
        0 & \mbox{otherwise},
    \end{cases}
\end{equation}
and
\begin{subequations}
    \begin{equation}
        \Gbf_1(\VarY_1, \VarY_2) = \VarY_1^2 + 2 \VarY_1 \VarY_2,
    \end{equation}
    \begin{equation}
        \Gbf_2(\VarY_1, \VarY_2) = \VarY_1.
    \end{equation}
\end{subequations}
\medbreak

Since by Proposition~\ref{prop:properties_BBalTree}, $\BBalTree$
satisfies the conditions of
Proposition~\ref{prop:functional_equation_sync}, by this last
proposition and~\eqref{equ:generating_series_sync_functional_equation},
the generating series $\GenS_{\SyncLang(\BBalTree)}$ of
$\SyncLang(\BBalTree)$ satisfies
$\GenS_{\SyncLang(\BBalTree)} = \Fbf_1(t, 0)$ where
\begin{equation}
    \Fbf_1(\VarY_1, \VarY_2) =
    \VarY_1 + \Fbf_1\left(\VarY_1^2 + 2 \VarY_1 \VarY_2, \VarY_1\right).
\end{equation}
This functional equation for the generating series of balanced binary
trees is the one obtained in~\cite{BLL88,BLL98,Knu98,Gir12e} by
different methods. As announced in
Section~\ref{subsubsec:generating_series_synchronous_languages}, the
coefficients of $\Fbf_1$ (and hence, those of
$\GenS_{\SyncLang(\BBalTree)}$) can be computed by iteration. This
consists in defining, for any $\ell \geq 0$, the polynomials
$\Fbf_1^{(\ell)}(\VarY_1, \VarY_2)$ as
\begin{equation}
    \label{equ:iteration_functional_equation_balanced_binary_trees}
    \Fbf_1^{(\ell)}(\VarY_1, \VarY_2) :=
    \begin{cases}
        \VarY_1 & \mbox{if } \ell = 0, \\
        \VarY_1 +
        \Fbf_1^{(\ell - 1)}\left(\VarY_1^2 + 2 \VarY_1 \VarY_2,
            \VarY_1\right) & \mbox{otherwise}.
    \end{cases}
\end{equation}
Since
\begin{equation}
    \Fbf_1(\VarY_1, \VarY_2) =
    \lim_{\ell \to \infty} \Fbf_1^{(\ell)}(\VarY_1, \VarY_2),
\end{equation}
Equation~\eqref{equ:iteration_functional_equation_balanced_binary_trees}
provides a way to compute the coefficients of
$\Fbf_1(\VarY_1, \VarY_2)$. First polynomials
$\Fbf_1^{(\ell)}(\VarY_1, \VarY_2)$ are
\begin{subequations}
\begin{equation}
    \Fbf_1^{(0)}(\VarY_1, \VarY_2) = \VarY_1,
\end{equation}
\begin{equation}
    \Fbf_1^{(1)}(\VarY_1, \VarY_2) = \VarY_1 + \VarY_1^2
    + 2\VarY_1\VarY_2,
\end{equation}
\begin{equation}
    \Fbf_1^{(2)}(\VarY_1, \VarY_2) = \VarY_1 + \VarY_1^2
    + 2\VarY_1\VarY_2 + 2\VarY_1^3 + 4\VarY_1^2\VarY_2 + \VarY_1^4
    + 4\VarY_1^3\VarY_2 + 4\VarY_1^2\VarY_2^2,
\end{equation}
\begin{multline}\begin{split}
    \Fbf_1^{(3)}(\VarY_1, \VarY_2) = \VarY_1 + \VarY_1^2
    + 2\VarY_1\VarY_2 + 2\VarY_1^3 + 4\VarY_1^2\VarY_2 + \VarY_1^4
    + 4\VarY_1^3\VarY_2 + 4\VarY_1^2\VarY_2^2 + 4\VarY_1^5 \\
    + 16\VarY_1^4\VarY_2
    + 16\VarY_1^3\VarY_2^2
    + 6\VarY_1^6 + 28\VarY_1^5\VarY_2 + 40\VarY_1^4\VarY_2^2
    + 16\VarY_1^3\VarY_2^3 + 4\VarY_1^7 + 24\VarY_1^6\VarY_2 \\
    + 48\VarY_1^5\VarY_2^2
    + 32\VarY_1^4\VarY_2^3 + \VarY_1^8
    + 8\VarY_1^7\VarY_2
    + 24\VarY_1^6\VarY_2^2 + 32\VarY_1^5\VarY_2^3
    + 16\VarY_1^4\VarY_2^4.
\end{split}\end{multline}
\end{subequations}
\medbreak

Besides, let us recall that $\BBalTree$ is synchronously unambiguous and
satisfies the properties stated by
Proposition~\ref{prop:properties_BBalTree}. Hence, $\BBalTree$
satisfies the conditions of
Theorem~\ref{thm:series_color_types_sync}. By this last theorem
and~\eqref{equ:generating_series_sync_recurrence},
$\GenS_{\SyncLang(\BBalTree)}$ satisfies, for any $n \geq 1$,
\begin{equation}
    \Angle{t^n, \GenS_{\SyncLang(\BBalTree)}}
    = \Angle{\VarY_1^n \VarY_2^0, \Fbf},
\end{equation}
where $\Fbf$ is the series satisfying, for any type
$\alpha \in \Types_{\{1, 2\}}$, the recursive formula
\begin{equation}
    \label{equ:recursive_enumeration_balanced_binary_trees_raw}
    \Angle{\VarY_1^{\alpha_1} \VarY_2^{\alpha_2}, \Fbf}
    = \delta_{\alpha, (1, 0)}
    +
    \sum_{\substack{
        d_1, d_2, d_3 \in \N \\
        \alpha_1 = 2 d_1 + d_2 + d_3 \\
        \alpha_2 = d_2
    }}
    \binom{d_1 + \alpha_2}{d_1}
    2^{d_2}
    \Angle{\VarY_1^{d_1 + d_2} \VarY_2^{d_3}, \Fbf}.
\end{equation}
This recursive formula offers an efficient way to compute the number
of balanced binary trees of a given size.
\medbreak

\section*{Concluding remarks}
We have presented in this chapter a framework for the generation of
combinatorial objects by using colored operads. The described devices
for combinatorial generation, called bud generating systems, are
generalizations of context-free grammars~\cite{Har78,HMU06} generating
words, of regular tree grammars~\cite{GS84,CDGJLLTT07} generating planar
rooted trees, and of synchronous grammars~\cite{Gir12e} generating some
treelike structures. We have provided tools to enumerate the objects of
the languages of bud generating systems or to define new statistics on
these by using formal power series on colored operads and several
products on these. There are many ways to extend this work. Here follow
some few further research directions.
\smallbreak

First, the notion of rationality and recognizability in usual formal
power series~\cite{Sch61,Sch63,Eil74,BR88}, in series on
monoids~\cite{Sak09}, and in series of trees~\cite{BR82} are fundamental.
For instance, a series $\Sbf \in \K \AAngle{\Mca}$ on a monoid $\Mca$ is
rational if it belongs to the closure of the set $\K \Angle{\Mca}$ of
polynomials on $\Mca$ with respect to the addition, the multiplication,
and the Kleene star operations. Equivalently, $\Sbf$ is rational if
there exists a $\K$-weighted automaton accepting it. The equivalence
between these two properties for the rationality property is remarkable.
We ask here for the definition of an analogous and consistent notion of
rationality for series on a colored operad $\Cca$. By consistent, we
mean a property of rationality for $\Cca$-series which can be defined
both by a closure property of the set $\K \Angle{\Cca}$ of the
polynomials on $\Cca$ with respect to some operations, and, at the same
time, by an acceptance property involving a notion of a $\K$-weighted
automaton on~$\Cca$. The analogous question about the definition of a
notion of recognizable series on colored operads also seems worth
studying.
\smallbreak

A second research direction fits mostly in the contexts of computer
science and compression theory. A straight-line grammar (see for
instance~\cite{ZL78,SS82,Ryt04}) is a context-free grammar with a
singleton as language. There exists also the analogous natural
counterpart for regular tree grammars~\cite{LM06}. One of the main
interests of straight-line grammars is that they offer a way to compress
a word (resp. a tree) by encoding it by a context-free grammar (resp. a
regular tree grammar). A word $u$ can potentially be represented by a
context-free grammar (as the unique element of its language) with less
memory than the direct representation of $u$, provided that $u$ is made
of several repeating factors. The analogous definition for bud
generating systems could potentially be used to compress a large variety
of combinatorial objects. Indeed, given a suitable monochrome operad
$\Oca$ defined on the objects we want to compress, we can encode an
object $x$ of $\Oca$ by a bud generating system $\Bca$ with $\Oca$ as
ground operad and such that the language (or the synchronous language)
of $\Bca$ is a singleton $\{y\}$ and $\Prune(y) = x$. Hence, we can hope
to obtain a new and efficient method to compress arbitrary combinatorial
objects.
\smallbreak

Let us finally describe a third extension of this work. Pros (see
Section~\ref{subsec:pros} of Chapter~\ref{chap:algebra}) are algebraic
structures which naturally generalize operads. Indeed, a pro is a set of
operators with several inputs and several outputs, unlike in operads
where operators have only one output (see for
instance~\cite{McL65,Mar08}). It seems fruitful to translate the main
definitions and constructions of this work (as {\em e.g.,} bud operads,
bud generating systems, series on colored operads, pre-Lie and
composition products of series, star operations, {\em etc.}) with pros
instead of operads. We can expect to obtain an even more general class
of grammars and obtain a more general framework for combinatorial
generation.
\medbreak


\chapter{Operads and regular languages} \label{chap:languages}
The content of this chapter comes from~\cite{GLMN16} and is a joint work with Jean-Gabriel Luque, Ludovic Mignot, and Florent Nicart.
\medbreak

\section*{Introduction}
Regular languages form an important class of languages, defined as the
ones that can be generated by Type-3 grammars of the
Chomsky-Schützenberger hierarchy~\cite{Cho59,CS63}. One of the most
surprising property of regular languages is that they can be described
by nonequivalent different ways, for instance by regular grammars,
automata, or regular expressions. These tools are nonequivalent in terms
of spatial complexity: a same family of regular languages can be
represented for example by automata with a linear number of states but
by regular expressions with an exponential number of
symbols~\cite{EZ76}.
\smallbreak

Multi-tildes~\cite{CCM11} are operators acting on languages and
introduced in order to increase the expressiveness of regular
expressions (that is, describing regular languages with the smallest
possible spatial complexity). These operators allow intuitively to jump
forward in a regular expression. Besides, multi-tildes come with a very
natural notion of composition, and it appears that this composition
endows the graded set of all the multi-tildes with a structure of a ns
set-operad~\cite{LMN13}. This establishes an unexpected link between the
theories of formal languages and of ns operads.
\smallbreak

In~\cite{LMN13}, the ns operads $\MT$ of the multi-tildes and $\Poset$
of the pseudo-transitive multi-tildes have been defined. The first one
is the ns operad aforementioned of multi-tildes and the second one is a
quotient operad of $\MT$ involving posets. The set of all the languages
over a finite alphabet is endowed with the structure of an $\MT$-monoid,
and also of a $\Poset$-monoid. The first structure is nonfaithful while
the second is faithful (in the sense that two different elements of
$\Poset$ act differently on languages). The ns operad $\Poset$ provides
hence a new way to express languages with optimality. Moreover, any
finite language can be expressed by the action of an element of $\Poset$
on languages that are empty or consisting only in one word of
length~$1$.
\smallbreak

The purpose of the present work is to generalize these constructions of
ns operads to regular languages (and not only on finite ones). The main
idea for this is to extend the notion of multi-tildes to double
multi-tildes. These are operators acting on languages and allow
intuitively to jump forward or backward in a regular expression. In this
generalization also, double multi-tildes are endowed with a natural
notion of composition and form a ns set-operad $\DMT$. This operad acts
on the set all the languages over a finite alphabet, and provides a way
to express any regular language by the action of an element of $\DMT$ on
languages that are empty or consisting only in one word of length $1$.
In this context, we also construct a quotient operad $\Qoset$ of $\DMT$
which plays the same role for regular languages as $\Poset$ plays for
finite languages. Indeed, the set of all regular languages forms a
faithful $\Qoset$-monoid.
\smallbreak

All the four ns set-operads considered in this chapter can be
constructed in a very similar way. For this reason, we provide an
abstraction for their construction through a functorial construction
$\PrecompToOp$, producing a ns set-operad from a precomposition. These
last structures are kinds of representations of a particular monoid. We
provide, by using precompositions and $\PrecompToOp$, alternative
constructions for the already known operads $\MT$ and $\Poset$, and
interpret our construction of the new operads $\DMT$ and~$\Qoset$.
\smallbreak

This chapter is organized as follows. Section~\ref{sec:precompositions}
contains the definition of the category of the precompositions and of
the functor $\PrecompToOp$. In
Section~\ref{sec:operads_from_precompositions}, we provide alternative
constructions of $\MT$ and $\Poset$, and define $\DMT$ and $\Qoset$. In
Section~\ref{sec:action_on_languages}, we study actions of $\DMT$ and
$\Qoset$ on languages.
\medbreak

\subsubsection*{Note}
This chapter deals only with ns set-operads. For this reason, ``operad''
means ``ns set-operad''.
\medbreak

\section{Breaking operads via precompositions}
\label{sec:precompositions}
The objective of this section is to introduce new algebraic objects, the
precompositions. These objects are a kind of representation of a certain
monoid denoted by $\MonoidPrecomp$ whose elements can be described in
terms of infinite matrices. We present here a functor from the category
of precompositions to the category of operads. We shall use this functor
in the sequel to reconstruct some already known operads and to construct
new ones.
\medbreak

\subsection{Monoids of infinite matrices}
\label{subsec:monoids_infinite_matrices}
We introduce here an associative algebra $\InfMat$ of infinite matrices
whose entries are indexed on $\Z^2$ and a quotient $\InfMatN$ of
$\InfMat$ of infinite matrices whose entries are indexed on $\N^2$.
Moreover, two respective subalgebras $\InfMatP$ and $\InfMatNP$ of
$\InfMat$ and $\InfMatN$ are described. The purpose of this section is
to give a realization and a presentation of $\MonoidPrecomp$, a monoid
defined by seeing $\InfMatNP$ as a monoid.
\medbreak

\subsubsection{A first algebra of infinite matrices}
We consider the vector space $\InfMat$ of all infinite matrices
$(A_{i j})_{i, j \in \Z}$ with a finite number of nonzero diagonals
whose entries belong to $\K$. A typical element $A$ of $\InfMat$ is a
finite linear combination of elements
\begin{equation}
    D^{(k, \lambda)} :=
    \sum_{i \in \Z} \lambda_i E^{(i+k,i)}
\end{equation}
where $\lambda=(\lambda_i)_{i\in\Z}$ and $E^{(k,\ell)}$ is the matrix
such that
\begin{equation}
    E^{(k,\ell)}_{i, j} =
    \begin{cases}
        1 & \mbox{if } (i, j) = (k, \ell), \\
        0 & \mbox{otherwise}.
    \end{cases}
\end{equation}
By observing that
\begin{equation}\begin{split}
    D^{(k, \lambda)} D^{(k', \lambda')}
    & =
    \left(\sum_{i \in \Z} \lambda_i E^{(i + k, i)}\right)
    \left(\sum_{i \in \Z}\lambda'_i E^{(i + k', i)}\right) \\
    & =
    \sum_{i \in \Z} \lambda_{i + k'} \lambda'_i E^{(i + k + k', i)} \\
    & =
    D^{(k + k, \lambda \Product_{k'} \lambda')},
\end{split}\end{equation}
where
\begin{math}
    \lambda \Product_{k'} \lambda'
    :=
    \left(\lambda_{i + k'} \lambda'_i\right)_{i \in \Z}
\end{math},
we deduce that $\InfMat$ is stable for the product of infinite matrices.
Moreover, the unit of $\InfMat$ is
\begin{equation}
    \Unit := D^{(0, (\dots, 1, 1, \dots))}
    =
    \sum_{i \in \Z} E^{(i, i)}.
\end{equation}
This leads to the following result.
\medbreak

\begin{Proposition} \label{prop:infinite_matrices_algebra}
    The space $\InfMat$ is a unitary associative algebra.
\end{Proposition}
\medbreak

Notice that when $\K$ is the field of complex numbers, the algebraic
structure of $\InfMat$ is very rich and has many connections with the
study of infinite Lie algebras (see {\em e.g.},~\cite{Kac90}).
\medbreak

\subsubsection{A first monoid of infinite matrices}
Here, for our purpose, we consider only the structure of monoid of
$\InfMat$. In particular, we define the submonoid $\InfMatP$ of
$\InfMat$ generated by the matrices
\begin{equation}
    M^{(i, n)} :=
    \sum_{j \leq i} E^{(j, j)}
    + \sum_{i < j} E^{(j + n - 1, j)}
\end{equation}
for each $i \in \Z$ and each $n \geq 1$. With these notations we have
\begin{equation}
    M^{(i, 1)} = \Unit,
\end{equation}
for any $i \in \Z$.
\medbreak

\begin{Proposition} \label{prop:presentation_monoid_infinite_matrices}
    The monoid $\InfMatP$ is isomorphic to the monoid
    $\MonoidPrecompBar$ generated by the symbols
    $\left\{\Asf_i^n : i \in \Z, n \geq 1\right\}$ subjected to the
    relations
    \begin{subequations}
    \begin{equation} \label{equ:monoid_infinite_matrices_1}
        \Asf_i^1 = \Unit,
        \qquad i \in \Z,
    \end{equation}
    \begin{equation} \label{equ:monoid_infinite_matrices_2}
        \Asf_i^n \; \Asf_j^m
        =
        \Asf_{j + n - 1}^m \; \Asf_i^n,
        \qquad i \leq j,
    \end{equation}
    \begin{equation} \label{equ:monoid_infinite_matrices_3}
        \Asf_{i + j}^n \; \Asf_i^m
        =
        \Asf_i^{n + m - 1},
        \qquad
        0 \leq j < m,
    \end{equation}
    \end{subequations}
    where $\Unit$ is the unit.
\end{Proposition}
\medbreak

\subsubsection{A second algebra of infinite matrices}
Let $\InfMatN$ be the vector space of all infinite matrices
$(A_{i j})_{i, j \in \N_{\geq 1}}$ with a finite number of nonzero
diagonals whose entries belong to $\K$. An analogous result as the one
stated by Proposition~\ref{prop:infinite_matrices_algebra} shows that
$\InfMatN$ is a unitary associative algebra. Moreover, there is a
surjective monoid morphism from $\InfMat$ to $\InfMatN$ sending any
matrix $(A_{i j})_{i, j \in \Z}$ to $(A_{i j})_{i, j \in \N_{\geq 1}}$.
Hence, $\InfMatN$ is a quotient monoid of $\InfMat$.
\medbreak

\subsubsection{A second monoid of infinite matrices}
Let us define the submonoid $\InfMatNP$ of $\InfMatN$ generated by the
matrices
\begin{equation}
    N^{(i, n)} :=
    \sum_{1 \leq j \leq i} E^{(j, j)}
    + \sum_{i < j} E^{(j + n - 1, j)}
\end{equation}
for each $n, i \geq 1$.
\medbreak

\begin{Proposition}
\label{prop:presentation_monoid_infinite_matrices_positive}
    The monoid $\InfMatNP$ is isomorphic to the monoid $\MonoidPrecomp$
    which is the quotient of the monoid $\MonoidPrecompBar$ satisfying
    the extra relations
    \begin{equation}
        \Asf_i^n = \Asf_0^n,
        \qquad i \leq 0.
    \end{equation}
\end{Proposition}
\medbreak

\subsection{Precompositions and operads} \label{subsec:precompositions}
Before describing a functor from the category of precompositions to the
category of operads, we define this last category. All this relies on
the presentation of the monoid $\MonoidPrecomp$ introduced in
Section~\ref{subsec:monoids_infinite_matrices}.
\medbreak

\subsubsection{Precompositions}
Let $(\Sca, \Product)$ be a commutative monoid endowed with a filtration
\begin{equation}
    \Sca = \bigcup_{n \geq 1} \Sca_n
\end{equation}
with
\begin{equation}
    \Sca_1\subseteq \Sca_2 \subseteq \dots \subseteq \Sca_n
    \subseteq \cdots
\end{equation}
and such that  each $\Sca_n$ is a submonoid of $\Sca$. We will denote by
$\Unit_\Sca$ the unit of $\Sca$.
\medbreak

A \Def{precomposition} is a monoid morphism
\begin{equation}
    \phi : \MonoidPrecomp \to \End(\Sca),
\end{equation}
where $\End(\Sca)$ denotes the set of all monoid endomorphisms
of $\Sca$, satisfying
\begin{subequations}
\begin{equation}
    \phi\left(\Asf_i^n\right) : \Sca_m \to \Sca_{n + m - 1},
    \qquad m \geq 1,
\end{equation}
\begin{equation}
    \phi\left(\Asf_i^n\right)|_{\Sca_m} = \Identity_{\Sca_m},
    \qquad i \geq m + 1,
\end{equation}
\end{subequations}
where $\phi\left(\Asf_i^n\right)|_{\Sca_m}$ denotes the restriction of
the map $\phi\left(\Asf_i^n\right)$ to the domain $\Sca_m$, and
$\Identity_{\Sca_m}$ is the identity map on~$\Sca_m$.
\medbreak

For simplicity, we denote by $\bar{\phi}_i^n$ the map
$\phi\left(\Asf_i^n\right)$. Observe that the maps $\bar{\phi}_i^k$ have
the following intuitive meaning. If $s$ is an element of $\Sca_m$, one
can see $s$ as an element having any number of inputs non smaller than
$m$. Under this point of view, $\bar{\phi}_i^n(s)$ is an element of
$\Sca_{n + m - 1}$ obtained by replacing in $s$ its $i$th input by
$n - 1$ new ones. Axioms~\eqref{equ:monoid_infinite_matrices_1},
\eqref{equ:monoid_infinite_matrices_2},
and~\eqref{equ:monoid_infinite_matrices_3} can be understood in the
light of this interpretation.
\medbreak

Now, let $\phi : \MonoidPrecomp \to \End(\Sca)$ and
$\phi' : \MonoidPrecomp \to \End(\Sca')$ be two precompositions. A map
$\alpha : \Sca \to \Sca'$ is a \Def{precomposition morphism} from $\phi$
to $\phi'$ if $\alpha$ is a monoid morphism and satisfies
\begin{subequations}
\begin{equation}
    \alpha : \Sca_n \to \Sca'_n,
    \qquad n \geq 1,
\end{equation}
\begin{equation}
    \bar{\phi'}_i^n(\alpha(s)) = \alpha(\bar{\phi}_i^n(s)),
    \qquad s \in \Sca.
\end{equation}
\end{subequations}
\medbreak

Let $\Precomp$ be the category wherein objects are all precompositions
and arrows are precomposition morphisms.
\medbreak

If $\phi : \MonoidPrecomp \to \End(\Sca)$ is a precomposition, we
define on $\Sca$ the binary products $\circ_i^{(n)}$ by
\begin{equation}
    s \circ_i^{(n)} t
    :=
    \bar{\phi}_i^n(s) \Product \bar{\phi}_0^i(t),
\end{equation}
where $\Product$ is the binary commutative and associative operation
of~$\Sca$.
\medbreak

\subsubsection{From precompositions to operads}
Let $\phi : \MonoidPrecomp \to \End(\Sca)$ be a precomposition. From
the commutative monoid $\Sca$, we define the set
\begin{equation}
    \Sbb := \bigsqcup_{n \geq 1} \Sbb_n
\end{equation}
where
\begin{equation}
    \Sbb_n :=
    \left\{(n, \Sfr) : \Sfr \in \Sca_n\right\}.
\end{equation}
Hence, $\Sbb$ is the set of all the elements of $\Sca$ endowed with an
arity. Let the partial compositions maps
\begin{equation}
    \circ_i : \Sbb_n \times \Sbb_m \to \Sbb_{n + m - 1}
\end{equation}
defined, for any $(n, \Sfr) \in \Sbb_n$ and $(m, \Tfr) \in \Sbb_m$, by
\begin{equation}
    (n, \Sfr) \circ_i (m, \Tfr)
    :=
    \left(n + m - 1, \Sfr \circ_i^{(m)} \Tfr\right).
\end{equation}
We denote by $\PrecompToOp(\phi)$ the set $\Sbb$ endowed with the maps
$\circ_i$ thus defined.
\medbreak

\begin{Theorem} \label{thm:precompositions_to_operads}
    The construction $\PrecompToOp$ is a functor from the category of
    precompositions to the category of operads.
\end{Theorem}
\medbreak

\subsubsection{Quotients of precompositions}
Let $\phi : \MonoidPrecomp \to \End(\Sca)$ be a precomposition and
$\gamma : \Sca \to \Sca$ be a monoid morphism such that
\begin{subequations}
\begin{equation}
    \gamma : \Sca_n \to \Sca_n,
    \qquad n \geq 1,
\end{equation}
\begin{equation}
    \gamma \circ \gamma = \gamma,
\end{equation}
\begin{equation}
    \bar{\phi}_i^n \circ \gamma = \gamma \circ \bar{\phi}_i^n.
\end{equation}
\end{subequations}
We call such a morphism $\gamma$ a \Def{compatible morphism}.
\medbreak

On the operad $\PrecompToOp(\phi)$, we denote, by a slight abuse of
notation, by
\begin{equation}
    \gamma : \Sbb \to \Sbb
\end{equation}
the map satisfying
\begin{equation}
    \gamma\left((n, \Sfr)\right) = (n, \gamma(\Sfr))
\end{equation}
for any $\Sfr \in \Sca_n$.
\medbreak

Let $\Congr_\gamma$ be the equivalence relation on $\Sca$ satisfying
$\Sfr \Congr_\gamma \Tfr$ if $\gamma(\Sfr) = \gamma(\Tfr)$ for any
$\Sfr, \Tfr \in \Sca$. Since $\gamma$ is a monoid morphism,
$\Congr_\gamma$ is a monoid congruence and hence,
$\Sca/_{\Congr_\gamma}$ is a quotient monoid of~$\Sca$.
\medbreak

On the operad $\PrecompToOp(\phi)$, we denote by a slight abuse of
notation by $\Congr_\gamma$ the equivalence relation satisfying
$(n, \Sfr) \Congr_\gamma (n, \Tfr)$ if $\Sfr \Congr_\gamma \Tfr$ for any
$\Sfr, \Tfr \in \Sca_n$.
\medbreak

Moreover, from the precomposition $\phi$ and the map $\gamma$, one
defines the precomposition
\begin{equation}
    \phi_\gamma : \MonoidPrecomp \to
    \End\left(\Sca/_{\Congr_\gamma}\right)
\end{equation}
defined for any $\Congr_\gamma$-equivalence class $[s]_{\Congr_\gamma}$
by
\begin{equation}
    {\bar{\phi_\gamma}}_i^n\left([s]_{\Congr_\gamma}\right)
    :=
    \left[\bar{\phi}_i^n(s)\right]_{\Congr_\gamma}.
\end{equation}
\medbreak

\begin{Theorem} \label{thm:precompositions_to_operads_quotients}
    Let $\phi : \MonoidPrecomp \to \End(\Sca)$ be a precomposition and
    $\gamma$ be a compatible morphism. Then, the operads
    $\PrecompToOp(\phi)/_{\Congr_\gamma}$ and
    $\PrecompToOp(\phi_\gamma)$ are isomorphic.
\end{Theorem}
\medbreak

\section{Constructing operads from precompositions}
\label{sec:operads_from_precompositions}
We apply here the construction $\PrecompToOp$ introduced in the previous
section to provide alternative constructions of the operads $\MT$ and
$\Poset$ introduced in~\cite{LMN13}, and to construct two new operads
$\DMT$ and $\Qoset$. These operads fit into the diagram
of Figure~\ref{fig:diagram_operads_languages}.
\begin{figure}[ht]
    \centering
    \begin{tikzpicture}[xscale=.75,yscale=.4,Centering]
        \node(DMT)at(0,0){\begin{math}\DMT\end{math}};
        \node(Qoset)at(-2,-3){\begin{math}\Qoset\end{math}};
        \node(MT)at(2,-3){\begin{math}\MT\end{math}};
        \node(Poset)at(0,-6){\begin{math}\Poset\end{math}};
        \draw[Surjection](DMT)--(Qoset);
        \draw[Surjection](MT)--(Poset);
        \draw[Injection](MT)--(DMT);
        \draw[Injection](Poset)--(Qoset);
    \end{tikzpicture}
    \caption[A diagram of operads having links with language theory.]
    {Diagram of operads where arrows $\rightarrowtail$ (resp.
    $\twoheadrightarrow$) are injective (resp. surjective) operad
    morphisms.}
    \label{fig:diagram_operads_languages}
\end{figure}
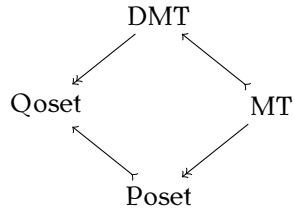
\medbreak

\subsection{Alternative constructions}
We begin by using the methods exposed in
Section~\ref{subsec:precompositions} to construct the operads $\MT$ of
multi-tildes and $\Poset$ of posets.
\medbreak

\subsubsection{Operad of multi-tildes} \label{subsubsec:operad_MT}
Multi-tildes are operators introduced in~\cite{CCM11} in the context of
formal language theory as a convenient way to express regular languages.
Let, for any $n \geq 1$, $P_n$ be the set
\begin{equation}
    P_n :=  \left\{(x, y) \in [n]^2 : x \leq y\right\}.
\end{equation}
A \Def{multi-tilde} is a pair $(n, \Sfr)$ where $n$ is a positive
integer and $\Sfr$ is a subset of $P_n$. The \Def{arity} of the
multi-tilde $(n, \Sfr)$ is $n$. The \Def{binary relation} of $(n, \Sfr)$
is the binary relation $\Rca_{(n, \Sfr)}$ on $[n + 1]$ satisfying
$x \,\Rca_{(n, \Sfr)}\, y$ if $x = y$ or $(x, y - 1) \in \Sfr$.
\medbreak

As shown in~\cite{LMN13}, the graded (by the arity) collection of all
multi-tildes admits a very natural structure of an operad. This operad,
denoted by $\MT$, is defined as follows. The partial composition
$(n, \Sfr) \circ_i (m, \Tfr)$, $i \in [n]$, of two multi-tildes
$(n, \Sfr)$ and $(m, \Tfr)$ is defined by
\begin{equation}
    (n, \Sfr) \circ_i (m, \Tfr) :=
    \left(n + m - 1,
    \left\{\Shift_i^m(x, y) : (x, y) \in \Sfr\right\}
    \cup
    \left\{\Shift_0^i(x, y) : (x, y) \in \Tfr\right\}\right),
\end{equation}
where
\begin{equation}
    \Shift_j^p(x, y) :=
    \begin{cases}
        (x, y) & \mbox{if } y \leq i - 1, \\
        (x, y + p - 1) & \mbox{if } x \leq i \leq y, \\
        (x + p - 1, y + p - 1) & \mbox{otherwise}.
    \end{cases}
\end{equation}
For instance, one has
\begin{subequations}
\begin{equation} \label{equ:example_composition_MT_1}
    (5, \{(1, 5), (2, 4), (4, 5)\}) \circ_4 (6, \{(2, 2), (4, 6)\}) \\
    = (10, \{(1, 10), (2, 9), (4, 10), (5, 5), (7, 9)\}),
\end{equation}
\begin{equation} \label{equ:example_composition_MT_2}
    (5, \{(1, 5), (2, 4), (4, 5)\}) \circ_5 (6, \{(2, 2), (4, 6)\}) \\
    = (10, \{(1, 10), (2, 4), (4, 10), (6, 6), (8, 10)\}).
\end{equation}
\end{subequations}
Observe that the multi-tilde $(1, \emptyset)$ is the unit of~$\MT$.
Since for any $n \geq 1$, $\# \MT(n) = \# \Set(P_n)$ (where $\Set(P_n)$
denote the set of all subsets of $P_n$),
\begin{equation}
    \# \MT(n) = 2^{\binom{n + 1}{2}}.
\end{equation}
Hence, the first dimensions of $\MT$ are
\begin{equation}
    2, 8, 64, 1024, 32768, 2097152, 268435456, 68719476736,
\end{equation}
and form Sequence~\OEIS{A006125} of~\cite{Slo}. Observe that the
dimensions of $\MT$ are very similar to the dimensions of the operads
obtained by the clique construction applied to a unitary magma having
exactly two elements (see
Section~\ref{subsubsec:dimensions_clique_operads} of
Chapter~\ref{chap:cliques}).
\medbreak

Let us provide a construction of $\MT$ through the functor
$\PrecompToOp$. Let $\Sca_n$ be the set of all the subsets of $P_n$. By
observing that $\Sca_n \subseteq \Sca_{n + 1}$, let
$\Sca := \cup_{n \geq 1} \Sca_n$. The pair $(\Sca, \cup)$ is a
commutative monoid whose unit is $\Unit_\Sca := \emptyset$ and belongs
to $\Sca_1$. Observe also that $\Sca$ is, as a monoid, generated by the
set $\{\{(x, y)\} : x \leq y\}$.
\medbreak

Let $\phi : \MonoidPrecomp \to \End(\Sca)$ be the precomposition such
that each morphism $\bar{\phi}_i^n$ is defined by its values on the
generators of $\Sca$ by
\begin{equation}
    \bar{\phi}_i^n(\{x, y\}) :=
    \begin{cases}
        \{(x, y)\} & \mbox{if } y \leq i - 1, \\
        \{(x, y + n - 1)\} & \mbox{if } x \leq i \leq y, \\
        \{(x + n - 1, y + n - 1)\} & \mbox{otherwise}.
    \end{cases}
\end{equation}
One can check that $\phi$ is a precomposition.
\medbreak

\begin{Proposition} \label{prop:precomposition_MT}
    The operads $\MT$ and $\PrecompToOp(\phi)$ are isomorphic.
\end{Proposition}
\medbreak

\subsubsection{Operads of posets}
In~\cite{LMN13}, an operad $\Poset$ defined as the quotient of $\MT$ by
the operad congruence $\Congr$ is considered, where for any multi-tildes
$(n, \Sfr)$ and $(n, \Tfr)$, one sets $(n, \Sfr) \Congr (n, \Tfr)$ if
the binary relations $\Rca_{(n, \Sfr)}$ and $\Rca_{(n, \Tfr)}$ have the
same reflexive and transitive closure. Since any $\equiv$-equivalence
class contains exactly one reflexive, transitive, and antisymmetric
relation, $\Poset$ is an operad on the set of all posets. More
precisely, the elements of $\Poset(n)$ are posets on $[n + 1]$ admitting
$(1, 2, \dots, n + 1)$ as a linear extension. For instance, one has
\begin{subequations}
\begin{equation} \label{equ:example_congruence_poset_1}
    (4, \{(1, 3), (2, 2), (3, 4)\})
    \Congr
    (4, \{
    (1, 3), (2, 2),
    \textcolor{Col4}{(2, 5)},
    (3, 4)
    \}),
\end{equation}
\begin{equation} \label{equ:example_congruence_poset_2}
    (4, \{(1, 1), (2, 3), (4, 4)\})
    \Congr
    (4, \{
    (1, 1), \textcolor{Col4}{(1, 3)},
    \textcolor{Col4}{(1, 4)}, (2, 3),
    \textcolor{Col4}{(2, 4)}, (4, 4)\}),
\end{equation}
\end{subequations}
and
\begin{equation}
    \left[
    (4, \{(1, 3), (2, 2), (3, 4)\})
    \right]_\Congr
    \circ_2
    \left[
    (4, \{(1, 1), (2, 3), (4, 4)\})
    \right]_\Congr
    =
    \left[
    (7, \{(1, 6), (2, 3), (3, 4), (5, 5), (6, 7)\})
    \right]_\Congr.
\end{equation}
The first dimensions of $\Poset$ are
\begin{equation}
    2, 7, 40, 357, 4824, 96428, 2800472, 116473461,
\end{equation}
and form Sequence~\OEIS{A006455} of~\cite{Slo}.
\medbreak

Let us consider the monoid $\Sca$ and the precomposition $\phi$ of
Section~\ref{subsubsec:operad_MT}. Let $\gamma : \Sca \to \Sca$ be the
map sending any set $\Sfr \in P_n$ to the set $\Sfr' \in P_n$ such that
$\Rca_{(n, \Sfr')}$ is the reflexive and transitive closure of
$\Rca_{(n, \Sfr)}$. For instance, the second components of the left
members of~\eqref{equ:example_congruence_poset_1}
and~\eqref{equ:example_congruence_poset_2} show elements of $P_4$
and the second components of their respective right members are their
reflexive and transitive closures. A multi-tilde $(n, \Sfr)$ is
\Def{pseudo-transitive} if $\Sfr$ belongs to the image of $\gamma$. One
can check that $\gamma$ is a compatible morphism. Hence, we can consider
the operad $\PrecompToOp(\phi_\gamma)$ which is, by
Theorem~\ref{thm:precompositions_to_operads_quotients}, isomorphic to
the operad~$\MT/_{\Congr_\gamma}$.
\medbreak

\begin{Proposition} \label{prop:precomposition_Poset}
    The operads $\Poset$ and $\PrecompToOp(\phi_\gamma)$
    are isomorphic.
\end{Proposition}
\medbreak

\subsection{New operads}
We now generalize the concept of multi-tildes to double multi-tildes.
In terms of operators on languages (see
Section~\ref{sec:action_on_languages}), multi-tildes can be seen as
operators allowing to jump forward in a regular expression and double
multi-tildes as operators allowing to jump both forward or backward in a
regular expression. The interest of this extension relies on the fact
that, while multi-tildes can emulate the sum and the concatenation,
double multi-tildes can emulate in addition to this the Kleene star of
regular expressions and their languages. We construct in this section an
operad $\DMT$ of double multi-tildes and a quotient $\Qoset$ of $\DMT$
of quasiorders.
\medbreak

\subsubsection{Operad of double multi-tildes}
\label{subsubsec:operad_DMT}
Let $\DMT$ be the \Def{operad of double multi-tildes} defined as
\begin{equation}
    \DMT := \MT \Hadamard \MT,
\end{equation}
where $\Hadamard$ is the Hadamard product of operads. An element of
arity $n$ of $\DMT$ is, by definition of $\Hadamard$, a pair
$((n, \Sfr), (n, \Tfr))$ where $(n, \Sfr)$ and $(n, \Tfr)$ are
multi-tildes. For simplicity, this element is simply denoted by
$(n, \Sfr, \Tfr)$ and is called a \Def{double multi-tilde}.
The \Def{binary relation} of $(n, \Sfr, \Tfr)$ is the binary relation
$\Rca_{(n, \Sfr, \Tfr)}$ on $[n + 1]$ satisfying
$x \,\Rca_{(n, \Sfr, \Tfr)}\, y$ if $x = y$ or $(x, y - 1) \in \Sfr$
or $(y - 1, x) \in \Tfr$. Since any multi-tilde $(n, \Sfr)$ can be seen
as a double multi-tilde $(n, \Sfr, \emptyset)$, $\MT$ is a suboperad of
$\DMT$. Observe that the double multi-tilde $(1, \emptyset, \emptyset)$
is the unit of $\DMT$. Since for any $n \geq 1$,
$\# \DMT(n) = \left(\# \Set(P_n)\right)^2$,
\begin{equation}
    \# \DMT(n) = 4^{\binom{n + 1}{2}}.
\end{equation}
Hence, the first dimensions of $\DMT$ are
\begin{multline}
    4, 64, 4096, 1048576, 1073741824, 4398046511104,
    72057594037927936, \\
    4722366482869645213696,
\end{multline}
and form Sequence~\OEIS{A053763} of~\cite{Slo}. Observe that the
dimensions of $\DMT$ are very similar to the dimensions of the operads
obtained by the clique construction applied to a unitary magma having
exactly four elements (see
Section~\ref{subsubsec:dimensions_clique_operads} of
Chapter~\ref{chap:cliques}).
\medbreak

Let us provide a construction of $\DMT$ through the function
$\PrecompToOp$. Let $\Dca_n$ be the set of all the pairs $(\Sfr, \Tfr)$,
where $\Sfr$ and $\Tfr$ are subsets of $P_n$. By observing that
$\Dca_n \subseteq \Dca_{n + 1}$, let $\Dca := \cup_{n \geq 1} \Dca_n$.
We endow $\Dca$ with the product $\cup$ defined by
\begin{math}
    (\Sfr, \Tfr) \cup (\Sfr', \Tfr') :=
    (\Sfr \cup \Sfr', \Tfr \cup \Tfr')
\end{math}
for all $(\Sfr, \Tfr), (\Sfr', \Tfr') \in \Dca$. The pair $(\Dca, \cup)$
is a commutative monoid whose unit is
$\Unit_\Dca := (\emptyset, \emptyset)$ and belongs to $\Dca_1$. Observe
also that $\Dca$ is, as a monoid, generated by the set
\begin{equation}
    \{(\{(x, y)\}, \emptyset) : x \leq y\}
    \cup
    \{(\emptyset, \{(x, y)\}) : x \leq y\}.
\end{equation}
\medbreak

Let $\psi : \MonoidPrecomp \to \End(\Dca)$ be the precomposition such
that each morphism $\bar{\psi}_i^n$ is defined by its values on the
generators of $\Dca$ by
\begin{subequations}
\begin{equation}
    \bar{\psi}_i^n((\{(x, y)\}, \emptyset)) :=
    \left(\bar{\phi}_i^n(\{(x, y)\}), \emptyset\right),
\end{equation}
\begin{equation}
    \bar{\psi}_i^n((\emptyset, \{(x, y)\})) :=
    \left(\emptyset, \bar{\phi}_i^n(\{(x, y)\})\right),
\end{equation}
\end{subequations}
where the $\bar{\phi}_i^n$ are the morphisms associated with the
precomposition $\phi$ of Section~\ref{subsubsec:operad_MT}. One can
check that $\psi$ is a precomposition.
\medbreak

\begin{Proposition} \label{prop:precomposition_DMT}
    The operads $\DMT$ and $\PrecompToOp(\psi)$ are isomorphic.
\end{Proposition}
\medbreak

\subsubsection{Operad of quasiorders} \label{subsubsec:operad_Qoset}
We construct here an operad $\Qoset$ which is to $\DMT$ what $\Poset$ is
to~$\MT$.
\medbreak

We define the operad $\Qoset$ as the quotient of $\DMT$ by the operad
congruence $\Congr$ defined as follows. For any double multi-tildes
$(n, \Sfr, \Tfr)$ and $(n', \Sfr', \Tfr')$, one sets
$(n, \Sfr, \Tfr) \Congr (n', \Sfr', \Tfr')$ if the binary relations
$\Rca_{(n, \Sfr, \Tfr)}$ and $\Rca_{(n', \Sfr', \Tfr')}$ have the same
reflexive and transitive closure. For instance, one has
\begin{equation} \label{equ:example_congruence_qoset}
    (4, \{(1, 2), (3, 3)\}, \{(1, 3)\})
    \Congr
    (4, \{(1, 2), \textcolor{Col4}{(1, 3)}, (3, 3)\}, \{(1, 3)\}).
\end{equation}
It is straightforward to check that $\Congr$ is a congruence of $\DMT$,
so that $\Qoset := \DMT/_\Congr$ is an operad. Moreover, since any
$\Congr$-equivalence class contains exactly one reflexive and transitive
relation, $\Qoset$ is an operad on the set of all quasiorders on
$[n + 1]$. The first dimensions of $\Qoset$ are
\begin{equation}
    4, 29, 355, 6942, 209527, 9535241, 642779354, 63260289423,
\end{equation}
and form Sequence~\OEIS{A000798} of~\cite{Slo}.
\medbreak

Let us consider the monoid $\Dca$ and the precomposition $\psi$ of
Section~\ref{subsubsec:operad_DMT}. Let $\gamma : \Dca \to \Dca$ be the
map sending any pair of $(\Sfr, \Tfr)$ of $\Dca_n^2$ to the pair
$(\Sfr', \Tfr')$ of $\Dca_n^2$ such that $\Rca_{(n, \Sfr', \Tfr')}$ is
the reflexive and transitive closure of $\Rca_{(n, \Sfr, \Tfr)}$. For
instance, the pair consisting in the second and third components of the
left member of~\eqref{equ:example_congruence_qoset} shows elements of
$P_4^2$ and the second and third components of the right member is its
reflexive and transitive closure. A double multi-tilde $(n, \Sfr, \Tfr)$
is \Def{pseudo-transitive} if $(\Sfr, \Tfr)$ belongs to the image of
$\gamma$. One can check that $\gamma$ is a compatible morphism. Hence,
we can consider the operad $\PrecompToOp(\psi_\gamma)$ which is, by
Theorem~\ref{thm:precompositions_to_operads_quotients}, isomorphic to
the operad $\DMT/_{\Congr_\gamma}$.
\medbreak

\section{Links with language theory} \label{sec:action_on_languages}
The main motivation for the introduction of the four operads $\MT$,
$\Poset$, $\DMT$, and $\Qoset$ (the first two in~\cite{LMN13} and the
last two here) relies on the fact that they act on languages. In more
precise terms, the set of all languages over a finite alphabet $A$ is
endowed with $\Oca$-monoid structures, where $\Oca$ is one of the four
aforementioned operads. We describe these structures in this section.
\medbreak

\subsection{Action of multi-tildes and double multi-tildes}
The action of $\DMT$ on languages on a finite alphabet can be described
in terms of automata. This leads to the construction of a
$\DMT$-monoid. All this justifies the role of $\DMT$ in formal language
theory since this operad provides a concise way to express languages.
\medbreak

\subsubsection{Automata and regular languages}
An \Def{automaton} is a tuple $(A, Q, \delta, i, t)$ where $A$ is a
\Def{ground alphabet}, $Q$ is a finite set, called \Def{set of states},
$\delta : Q \times A \sqcup \{\epsilon\} \to \Set(Q)$ is a
\Def{transition map}, $i$ is a state of $Q$ called \Def{initial state},
and $t$ is a state of $Q$ called \Def{terminal state}. We consider here
very particular automata (also known as $\epsilon$-automata). We use
the main definitions of the theory (see for instance~\cite{Sak09}), like
the notion of language recognized by an automaton, regular languages,
regular expressions, {\em etc.}
\medbreak

From now on, $\Abb$ is the infinite alphabet
$\Abb := \{\Asf_1, \Asf_2, \dots\}$ and $A$ is any finite alphabet.
\medbreak

\subsubsection{From double multi-tildes to automata}
Let $(n, \Sfr, \Tfr)$ be a double multi-tilde of arity~$n$ of $\DMT$ and
$\Aca_{(n, \Sfr, \Tfr)}$ be the automaton
$(\Abb, Q, \delta, q_1, q_{n + 1})$ defined by
\begin{equation}
    Q := \left\{q_j : j \in [n + 1]\right\},
\end{equation}
\begin{equation}
    \delta(q_j, \Asf_i) := \{q_{j + 1}\},
\end{equation}
\begin{equation}
    \delta(q_j, \epsilon) :=
    \left\{q_{k + 1} : (j, k) \in \Sfr\right\}
    \cup
    \left\{q_k : (k, j - 1) \in \Sfr\right\}.
\end{equation}
For instance, by considering the double multi-tilde
\begin{equation}
    {(n, \Sfr, \Tfr)} :=
    (6, \{(1, 3), (2, 2), (3, 4) \}, \{(2, 2), (2, 3), (4, 5)\}),
\end{equation}
the transition map of the automaton
\begin{math}
    \Aca_{(n, \Sfr, \Tfr)}
    :=
    (\Abb, \{q_1, \dots, q_7\}, \delta, q_1, q_7)
\end{math}
satisfies
\begin{multline}
    \delta(q_1, \Asf_1) := \{q_2\}, \quad
    \delta(q_2, \Asf_2) := \{q_3\}, \quad
    \delta(q_3, \Asf_3) := \{q_4\}, \\
    \delta(q_4, \Asf_4) := \{q_5\}, \quad
    \delta(q_5, \Asf_5) := \{q_6\}, \quad
    \delta(q_6, \Asf_6) := \{q_7\},
\end{multline}
and
\begin{multline}
    \delta(q_1, \epsilon) := \{q_4\}, \quad
    \delta(q_2, \epsilon) := \{q_3\}, \quad
    \delta(q_3, \epsilon) := \{q_2, q_5\}, \\
    \delta(q_4, \epsilon) := \{q_2\}, \quad
    \delta(q_5, \epsilon) := \emptyset, \quad
    \delta(q_6, \epsilon) := \{q_4\}.
\end{multline}
This automaton can be depicted as
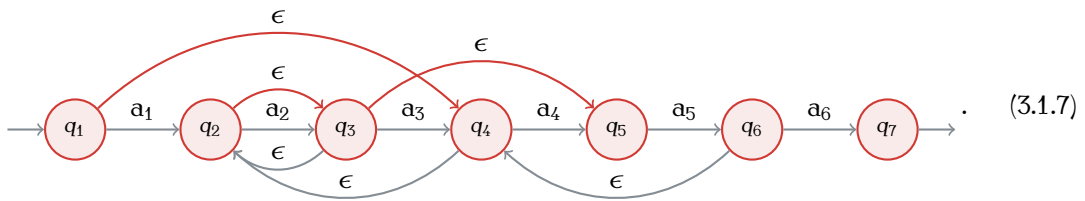
\begin{equation}
      \begin{tikzpicture}[xscale=.9,Centering]
            \node[AutState](q1)at(0,0){\begin{math}q_1\end{math}};
            \node[AutState](q2)at(2,0){\begin{math}q_2\end{math}};
            \node[AutState](q3)at(4,0){\begin{math}q_3\end{math}};
            \node[AutState](q4)at(6,0){\begin{math}q_4\end{math}};
            \node[AutState](q5)at(8,0){\begin{math}q_5\end{math}};
            \node[AutState](q6)at(10,0){\begin{math}q_6\end{math}};
            \node[AutState](q7)at(12,0){\begin{math}q_7\end{math}};
            \draw[AutArc](-1,0)--(q1);
            \draw[AutArc](q7)--(13,0);
            \draw[AutArc](q1)edge[]node[above]
                {\begin{math}\Asf_1\end{math}}(q2);
            \draw[AutArc](q2)edge[]node[above]
                {\begin{math}\Asf_2\end{math}}(q3);
            \draw[AutArc](q3)edge[]node[above]
                {\begin{math}\Asf_3\end{math}}(q4);
            \draw[AutArc](q4)edge[]node[above]
                {\begin{math}\Asf_4\end{math}}(q5);
            \draw[AutArc](q5)edge[]node[above]
                {\begin{math}\Asf_5\end{math}}(q6);
            \draw[AutArc](q6)edge[]node[above]
                {\begin{math}\Asf_6\end{math}}(q7);
            \draw[AutArcColorA](q1)edge[bend left=40]node[above]
                {\begin{math}\epsilon\end{math}}(q4);
            \draw[AutArcColorA](q2)edge[bend left=40]node[above]
                {\begin{math}\epsilon\end{math}}(q3);
            \draw[AutArcColorA](q3)edge[bend left=40]node[above]
                {\begin{math}\epsilon\end{math}}(q5);
            \draw[AutArcColorB](q3)edge[bend left=40]node[above]
                {\begin{math}\epsilon\end{math}}(q2);
            \draw[AutArcColorB](q4)edge[bend left=40]node[above]
                {\begin{math}\epsilon\end{math}}(q2);
            \draw[AutArcColorB](q6)edge[bend left=40]node[above]
                {\begin{math}\epsilon\end{math}}(q4);
      \end{tikzpicture}\,.
\end{equation}
\medbreak

\subsubsection{Action of double multi-tildes on languages}
\label{subsubsec:action_DMT_languages}
Let $\Lca_A$ be the set of all the languages on the alphabet $A$,
$\Rca_A$ be the subset of $\Lca_A$ consisting in all regular languages,
and $\Fca_A$ be the subset of $\Rca_A$ consisting in all finite
languages. Let, for any $n \geq 1$, the action
\begin{equation}
    \Action :
    \DMT(n) \times {\Lca_A}^{\times n} \to
    \Lca_A
\end{equation}
defined, for any double multi-tilde $(n, \Sfr, \Tfr)$ and any languages
$\Lfr_1$, \dots, $\Lfr_n$ of $\Lca_A$ by
\begin{equation}
    (n, \Sfr, \Tfr)
    \Action
    (\Lfr_1, \dots, \Lfr_n)
    :=
    \LangDescr\left(
        \Aca_{(n, \Sfr, \Tfr)}(\Lfr_1, \dots, \Lfr_n)\right),
\end{equation}
where $\Aca_{(n, \Sfr, \Tfr)}(\Lfr_1, \dots, \Lfr_n)$ is the automaton
obtained by replacing each letter $\Asf_i$ by the language $\Lfr_i$,
$i \in [n]$, and for any automaton $\Aca$, $\LangDescr(\Aca)$ denotes
the language described by~$\Aca$. For instance, by using regular
expressions,
\begin{equation}\begin{split}
    (n, \{(1, 4)\}, \{(2, 3), (3, 4)\})
    \Action
    (\Lfr_1, \Lfr_2, \Lfr_3)
    & =
    \LangDescr\left(
    \begin{tikzpicture}[xscale=.8,Centering]
        \node[AutState](q1)at(0,0){\begin{math}q_1\end{math}};
        \node[AutState](q2)at(2,0){\begin{math}q_2\end{math}};
        \node[AutState](q3)at(4,0){\begin{math}q_3\end{math}};
        \node[AutState](q4)at(6,0){\begin{math}q_4\end{math}};
        \draw[AutArc](-1,0)--(q1);
        \draw[AutArc](q4)--(7,0);
        \draw[AutArc](q1)edge[]node[above]
            {\begin{math}\Lfr_1\end{math}}(q2);
        \draw[AutArc](q2)edge[]node[above]
            {\begin{math}\Lfr_2\end{math}}(q3);
        \draw[AutArc](q3)edge[]node[above]
            {\begin{math}\Lfr_3\end{math}}(q4);
        \draw[AutArcColorA](q1)edge[bend left=30]node[above]
            {\begin{math}\epsilon\end{math}}(q4);
        \draw[AutArcColorB](q3)edge[bend left=40]node[above]
            {\begin{math}\epsilon\end{math}}(q2);
        \draw[AutArcColorB](q4)edge[bend left=40]node[above]
            {\begin{math}\epsilon\end{math}}(q3);
    \end{tikzpicture}
    \right) \\
    & =
    \epsilon +
    (\epsilon + \Lfr_1 \Lfr_2)
    (\Lfr_2 + \Lfr_3)^* \Lfr_3.
\end{split}\end{equation}
\medbreak

We can observe that if $(n, \Sfr, \Tfr)$ and $(n', \Sfr', \Tfr')$ are
two double multi-tildes,
$\Aca_{(n, \Sfr, \Tfr) \circ_i (n', \Sfr', \Tfr')}$ is the automaton
obtained by replacing the transition labeled by $\Asf_i$ connecting the
states $q_i$ and $q_{i + 1}$ of $\Aca_{(n, \Sfr, \Tfr)}$ by
$\Aca_{(n', \Sfr', \Tfr')}$, and by relabeling adequately its states and
transitions. Replacing in this way a transition by an automaton is
possible since $\Aca_{(n', \Sfr', \Tfr')}$ has exactly one initial and
one terminal state. From this observation, one has the following result.
\medbreak

\begin{Theorem} \label{thm:DMT_monoid_languages}
    The actions $\Action$ endow the set $\Lca_A$ of all languages on $A$
    with a structure of a $\DMT$-monoid. Moreover, the restrictions of
    the $\Action$ to the set $\Rca_A$ of all regular languages on $A$
    endow $\Rca_A$ with a structure of a $\DMT$-monoid.
\end{Theorem}
\medbreak

\subsubsection{Expressions for regular languages}

\begin{Proposition} \label{prop:DMT_operations_languages}
    Let $\Lfr$ and $\Lfr'$ be two languages of $\Lca_A$ satisfying
    \begin{equation}
        \Lfr = (n, \Sfr, \Tfr) \Action (\alpha_1, \dots, \alpha_n)
    \end{equation}
    and
    \begin{equation}
        \Lfr' = \left(n', \Sfr', \Tfr'\right) \Action
        \left(\alpha'_1, \dots, \alpha'_{n'}\right)
    \end{equation}
    for some $n, n' \geq 1$, $(n, \Sfr, \Tfr) \in \DMT(n)$,
    $\left(n', \Sfr', \Tfr'\right) \in \DMT(n')$,
    $\alpha_i, \alpha'_i \in \{\{a\} \in A\} \cup \{\emptyset\}$. Then,
    \begin{subequations}
    \begin{equation}
        \Lfr + \Lfr' =
        \left( (3, \{(1, 3), (2, 4)\}, \emptyset) \circ
        \left[(n, \Sfr, \Tfr), (1, \emptyset, \emptyset),
        \left(n', \Sfr', \Tfr'\right)\right]
        \right) \Action
        \left(\alpha_1, \dots, \alpha_n, \emptyset,
        \alpha'_1, \dots, \alpha'_{n'}\right),
    \end{equation}
    \begin{equation}
        \Lfr\Lfr' =
        \left((3, \{(2, 3)\}, \emptyset) \circ
        \left[(n, \Sfr, \Tfr), (1, \emptyset, \emptyset),
        \left(n', \Sfr', \Tfr'\right)\right]
        \right)
        \Action
        \left(\alpha_1, \dots, \alpha_n, \emptyset,
        \alpha'_1, \dots, \alpha'_{n'}\right),
    \end{equation}
    \begin{equation}
        \Lfr^* =
        \left((1, \{(1, 2)\}, \{(2, 1)\}) \circ_1 (n, \Sfr, \Tfr)
        \right)
        \Action (\alpha_1, \dots, \alpha_n).
    \end{equation}
    \end{subequations}
\end{Proposition}
\medbreak

In the statement of Proposition~\ref{prop:DMT_operations_languages}, the
symbols $\circ$ denote the complete composition maps of~$\DMT$.
\medbreak

As a consequence of Proposition~\ref{prop:DMT_operations_languages}, one
obtains the following result.
\medbreak

\begin{Proposition} \label{prop:DMT_description_regular_languages}
    Any regular language $\Lfr$ of $\Rca_A$ can be expressed as
    \begin{equation}
        \Lfr = (n, \Sfr, \Tfr) \Action (\alpha_1, \dots, \alpha_n),
    \end{equation}
    for some $n \geq 1$, $(n, \Sfr, \Tfr) \in \DMT(n)$, and
    $\alpha_i \in \{\{a\} \in A\} \cup \{\emptyset\}$.
\end{Proposition}
\medbreak

\subsubsection{Expressions for finite regular languages}
Since $\MT$ can be seen as a suboperad of $\DMT$ consisting in the
double multi-tildes of the form $(n, \Sfr, \emptyset)$, the action
$\Action$ of $\DMT$ on $\Lca_A$ described in
Section~\ref{subsubsec:action_DMT_languages} can be restricted on $\MT$.
In this way, we recover a result of~\cite{LMN13}.
\medbreak

\begin{Proposition} \label{prop:DMT_description_finite_languages}
    Any finite language $\Lfr$ of $\Fca_A$ can be expressed as
    \begin{equation}
        \Lfr = (n, \Sfr, \emptyset) \Action (\alpha_1, \dots, \alpha_n),
    \end{equation}
    for some $n \geq 1$, $(n, \Sfr, \emptyset) \in \DMT(n)$, and
    $\alpha_i \in \{\{a\} \in A\} \cup \{\emptyset\}$.
\end{Proposition}
\medbreak

\subsection{Action of pseudo-transitive double multi-tildes}
The quotient $\Qoset$ of $\DMT$ inherits the action of $\DMT$ on
languages. The main interest to consider the associated $\Qoset$-monoid
instead of the $\DMT$-monoid is that this last one is nonfaithful while
the first is. Hence, $\Qoset$ is an operad providing optimal operators
to describe regular languages.
\medbreak

\subsubsection{A nonfaithful action of $\DMT$ on languages}
Observe that the description of languages by the action of a double
multi-tilde on languages (for instance in the ways provided by
Propositions~\ref{prop:DMT_description_regular_languages}
and~\ref{prop:DMT_description_finite_languages}) is not optimal since a
language $\Lfr$ can be described from different double multi-tildes of
the same arity. Indeed,
\begin{equation}
    (2, \{(1, 2), (2, 3)\}, \emptyset)
    \Action (\{a\}, \{b\})
    =
    (a + \epsilon) (b + \epsilon)
    =
    (2, \{(1, 2), (2, 3), (1, 3)\}, \emptyset)
    \Action (\{a\}, \{b\}).
\end{equation}
In other words, the $\DMT$-monoid consisting in all languages endowed
with the action $\Action$ is a nonfaithful $\DMT$-monoid.
\medbreak

\subsubsection{A faithful action of $\Qoset$ on languages}
Since $\Qoset$ is a quotient operad of $\DMT$, the actions $\Action$ of
$\DMT$ on $\Lca_A$ defined in
Section~\ref{subsubsec:action_DMT_languages} are still well-defined on
$\Qoset$. More precisely, for any $\Congr_\gamma$-equivalence class
$\left[(n, \Sfr, \Tfr)\right]_{\Congr_\gamma}$ of double multi-tildes,
where $\Congr_\gamma$ is the operad congruence introduced in
Section~\ref{subsubsec:operad_Qoset}, and any languages $\Lfr_1$, \dots,
$\Lfr_n$ of $\Lca_A$,
\begin{equation} \label{equ:action_Qoset_languages}
    \left[(n, \Sfr, \Tfr)\right]_{\Congr_\gamma}
    \Action
    (\Lfr_1, \dots, \Lfr_n)
    :=
    (n, \Sfr, \Tfr)
    \Action
    (\Lfr_1, \dots, \Lfr_n),
\end{equation}
where the symbol $\Action$ of the right member
of~\eqref{equ:action_Qoset_languages} denote the actions of $\DMT$ on
languages.
\medbreak

Now, contrariwise to the action of $\DMT$ on languages, the action of
$\Qoset$ is optimal in the following sense.
\medbreak

\begin{Theorem} \label{thm:Qoset_faithful_monoid_languages}
    If $A$ has at least two letters, the actions $\Action$ endow the set
    $\Lca_A$ of all languages on $A$ with a structure of a faithful
    $\Qoset$-monoid.
\end{Theorem}
\medbreak

\subsubsection{Operations of $\Qoset$ as operators on languages}
Let us examine all the actions of the elements of $\Qoset$ up to
arity $2$. We denote each element of $\Qoset$ by pseudo-transitive
double multi-tildes.
\medbreak

In arity $1$,
\begin{subequations}
\begin{equation}
    (1, \emptyset, \emptyset) \Action \{a\} = \epsilon,
\end{equation}
\begin{equation}
    (1, \{(1, 1)\}, \emptyset) \Action \{a\} = \epsilon + a,
\end{equation}
\begin{equation}
    (1, \emptyset, \{(1, 1)\}) \Action \{a\} = a^+,
\end{equation}
\begin{equation}
    (1, \{(1, 1)\}, \{(1, 1)\}) \Action \{a\} = a^*.
\end{equation}
\end{subequations}
In arity $2$,
\begin{equation}
    (2, \emptyset, \emptyset) \Action (\{a\}, \{b\}) = ab,
\end{equation}
\begin{subequations}
\begin{equation}
    (2, \{(1, 1)\}, \emptyset) \Action (\{a\}, \{b\}) = b + ab,
\end{equation}
\begin{equation}
    (2, \{(1, 2)\}, \emptyset) \Action (\{a\}, \{b\}) = \epsilon + ab,
\end{equation}
\begin{equation}
    (2, \{(2, 2)\}, \emptyset) \Action (\{a\}, \{b\}) = a + ab,
\end{equation}
\begin{equation}
    (2, \{(1, 1), (1, 2)\}, \emptyset) \Action (\{a\}, \{b\}) =
    \epsilon + b + ab,
\end{equation}
\begin{equation}
    (2, \{(1, 2), (2, 2)\}, \emptyset) \Action (\{a\}, \{b\})
    = \epsilon + a + ab,
\end{equation}
\begin{equation}
    (2, \{(1, 1), (1, 2), (2, 2)\}, \emptyset) \Action (\{a\}, \{b\}) =
    \epsilon + a + b + ab,
\end{equation}
\begin{equation}
    (2, \emptyset, \{(1, 1)\}) \Action (\{a\}, \{b\}) = a^+b,
\end{equation}
\begin{equation}
    (2, \emptyset, \{(1, 2)\}) \Action (\{a\}, \{b\}) = (ab)^+,
\end{equation}
\begin{equation}
    (2, \emptyset, \{(2, 2)\}) \Action (\{a\}, \{b\}) = ab^+,
\end{equation}
\begin{equation}
    (2, \emptyset, \{(1, 1), (1, 2)\}) \Action (\{a\}, \{b\})
    = {(a^+b)}^+,
\end{equation}
\begin{equation}
    (2, \emptyset, \{(1, 2), (2, 2)\}) \Action (\{a\}, \{b\})
    = {(ab^+)}^+,
\end{equation}
\begin{equation}
    (2, \emptyset, \{(1, 1), (1, 2), (2, 2) \}) \Action (\{a\}, \{b\})
    = {(a^+b^+)}^+,
\end{equation}
\end{subequations}
\begin{subequations}
\begin{equation}
    (2, \{(1, 1)\}, \{(1, 1)\}) \Action (\{a\}, \{b\}) = a^*b
\end{equation}
\begin{equation}
    (2, \{(1, 1)\}, \{(2, 2)\}) \Action (\{a\}, \{b\}) =
    (\epsilon + a)b^+,
\end{equation}
\begin{equation}
    (2, \{(1, 2)\}, \{(1, 2)\}) \Action (\{a\}, \{b\}) =
    (ab)^*,
\end{equation}
\begin{equation}
    (2, \{(2, 2)\}, \{(1, 1)\}) \Action (\{a\}, \{b\}) =
    a^+(\epsilon + b),
\end{equation}
\begin{equation}
    (2, \{(2, 2)\}, \{(2, 2)\}) \Action (\{a\}, \{b\}) = ab^*,
\end{equation}
\end{subequations}
\begin{subequations}
\begin{equation}
    (2, \{(1, 1)\}, \{(1, 2), (2, 2)\}) \Action (\{a\}, \{b\}) =
    {((\epsilon + a)b^+)}^+,
\end{equation}
\begin{equation}
    (2, \{(2, 2)\}, \{(1, 1), (1, 2)\}) \Action (\{a\}, \{b\}) =
    {(a^+ (\epsilon + b))}^+,
\end{equation}
\begin{equation}
    (2, \{(1, 1)\}, \{(1, 1), (1, 2), (2, 2)\}) \Action (\{a\}, \{b\})
    = {(a^*b^+)}^+,
\end{equation}
\begin{equation}
    (2, \{(2, 2)\}, \{(1, 1), (1, 2), (2, 2)\}) \Action (\{a\}, \{b\})
    = {(a^+b^*)}^+,
\end{equation}
\begin{equation}
    (2, \{(1, 2), (2, 2)\}, \{(1, 1)\}) \Action (\{a\}, \{b\}) =
    \epsilon + a^+(\epsilon + b),
\end{equation}
\begin{equation}
    (2, \{(1, 1), (1, 2)\}, \{(2, 2)\}) \Action (\{a\}, \{b\}) =
    \epsilon + (\epsilon + a)b^+,
\end{equation}
\begin{equation}
    (2, \{(1, 1), (1, 2), (2, 2) \}, \{(1, 1)\}) \Action (\{a\}, \{b\})
    = a^* (\epsilon + b),
\end{equation}
\begin{equation}
    (2, \{(1, 1), (1, 2), (2, 2) \}, \{(2, 2)\}) \Action (\{a\}, \{b\})
    = (\epsilon + a)b^*,
\end{equation}
\end{subequations}
\begin{subequations}
\begin{equation}
    (2, \{(1, 1), (1, 2)\}, \{(1, 2), (2, 2) \}) \Action (\{a\}, \{b\})
    = {((\epsilon + a)b^+)}^*,
\end{equation}
\begin{equation}
    (2, \{(1, 1), (1, 2), (2, 2) \}, \{(1, 1), (1, 2), (2, 2)\})
    \Action (\{a\}, \{b\}) = (a + b)^*,
\end{equation}
\begin{equation}
    (2, \{(1, 2), (2, 2)\}, \{(1, 1), (1, 2) \}) \Action (\{a\}, \{b\})
    = {(a^+ (\epsilon + b))}^*.
\end{equation}
\end{subequations}
\medbreak

\section*{Concluding remarks}
The work presented in this chapter provides two kinds of results. The
first one consists in a general construction of operads through a
functor $\PrecompToOp$, producing an operad from a precomposition. The
second one consists in the application of this construction to obtain
new operads or alternative descriptions of already existing ones. In
this context we have introduced and worked with operads acting on formal
languages.
\smallbreak

The operad $\Qoset$, quotient of the operad $\DMT$ of double
multi-tildes, acts faithfully on the set of all regular languages over a
finite alphabet. Since $\Qoset$ is a combinatorial operad, it offers
countable operations for denoting regular languages. The expressions
thus obtained to define regular languages lead to the definition of
several measures for their complexity. For instance, if $\Lfr$ is a
regular language, one can define $w_1(\Lfr)$ (resp. $w_2(\Lfr)$) as the
minimal arity (resp. number of pairs) of the element of $\Qoset$
required to express $\Lfr$ (see
Proposition~\ref{prop:DMT_description_regular_languages}
and Theorem~\ref{thm:Qoset_faithful_monoid_languages}). Intuitively,
$w_1$ and $w_2$ can be respectively interpreted as functions measuring
the width and the height of a language. The first one, $w_1$, is indeed
the minimal number of occurrences of symbols or $\emptyset$ in the
expression of $\Lfr$. The measure $w_2$ expresses the minimal complexity
of an operator involved for denoting the languages. These measures
deserve to be investigated; in particular a parallel with the size of a
minimal  automaton (in terms of states or transitions) should be
established.
\smallbreak

Another perspective is the extension of the conversion methods from
automata to expressions by using double multi-tildes. These conversions
were studied in~\cite{CCM10} and in~\cite{CCM12}. By slightly modifying
the action of the operads, we aim to extend these algorithms of
conversions. Conversely, it seems worth designing algorithms producing
automata from expressions (like {\em e.g.}, position
functions~\cite{Glu61} or expression derivatives~\cite{Brz64,Ant96}).
\smallbreak

A last perspective is the following. By the Alexandroff
correspondence~\cite{Ale37}, quasiorders on finite sets are in bijection
with finite topologies. The question here consists in investigating if
the action of the operad of quasiorders $\Qoset$ on languages has a
topological interpretation.
\medbreak

\backmatter


\chapter*{Conclusion}
In all this dissertation, our main philosophy is to design operations on
combinatorial objects in order to construct algebraic structures on
them. By studying algebraically these structures, we hope to grab
combinatorial  properties on the objects. All this provides a tool to
tackle problems coming from enumerative combinatorics or computer
science.
\smallbreak

We expose in this work numerous constructions inputting simple algebraic
structures (like magmas, monoids, or posets) and outputting more
complicated ones (like Hopf bialgebras, operads, and pros) and involving
combinatorial objects. Therefore, our main contribution is to provide
metatools, in the sense that our constructions can be used to endow
combinatorial collections with algebraic structures.
\smallbreak

Each chapter ends with a section named ``Concluding remarks'' raising
some contextual open questions. For this reason we will not mention
these here. Let us instead speak about the general and cross sectional
ideas and directions for future research.
\medbreak

\SkipTocEntry\section*{About constructions of operads}
A first general direction consists in using the constructions $\T$ (see
Chapters~\ref{chap:monoids} and~\ref{chap:polydendr}), $\As$ (see
Chapter~\ref{chap:posets}), $\Cli$ (see
Chapter~\ref{chap:cliques}), and $\PrecompToOp$ (see
Chapter~\ref{chap:languages}) to define even more operads. As we have
seen, these constructions lead to the definitions of many interesting
operads, involving a large range of combinatorial objects and of partial
composition operations and algorithms. We think that we are far to have
exhausted the subject and that many other operads deserving to be
studied can be obtained.
\smallbreak

A next logical continuation is to develop more connections between
combinatorial algebraic structures and properties of their underlying
combinatorial objects. We have pointed out, mostly in
Chapter~\ref{chap:buds}, that ns colored operads lead to a
generalization of usual formal power series. By using the operations
of operads, we obtain a bunch a natural operations on such generalized
series, forming new tools for enumerative prospects. This axis consists
in constructing new operads on various kind of objects (like integer
partitions, Young tableaux, planar maps, {\em etc.}) and use these and
their formal power series to discover enumerative properties.
\medbreak

\SkipTocEntry\section*{About pros and their combinatorics}
In Chapter~\ref{chap:pros_bialgebras}, a link between the theory of pros
and the one of Hopf bialgebras has been highlighted through a
construction $\PvH$ associating a Hopf bialgebra with some pros. These
last objects are generalizations of operads and their underlying
combinatorics is much less developed and understood than those of
operads. For instance, while free operads are very well-known
structures, free pros are not well described. Some particular phenomena
occur when one considers free pros on generators having null input
and/or output arities. The question here is to provide a good
combinatorial realization of free pros. The described realization in
terms of prographs (see Chapters~\ref{chap:combinatorics}
and~\ref{chap:algebra}) allows only generators with at least one input
and one output. This research axis contains also questions like a
description of the Hilbert series of free pros.
\smallbreak

Pros provide also an interesting framework to work with all the
symmetric groups at the same time. Indeed, the pro $\K \Angle{\Per}$ of
all permutations (see Chapter~\ref{chap:algebra}) encapsulates the
composition operation of all permutations and hence, contains all groups
$\SymmetricGroup_n$, $n \in \N$. This pro, by computing its presentation
by generators and relations, leads to the known presentation of the
symmetric groups in terms of elementary transpositions. Here we ask for
a general construction associating a pro with any sequence $W_n$,
$n \in \N$, of Coxeter groups, analog to what is $\K \Angle{\Per}$ for
$\SymmetricGroup_n$, $n \in \N$. This could lead to combinatorial
realizations of some Coxeter groups. The same, but more general and hard
question, consisting in encapsulating a sequence $M_n$, $n \in \N$, of
monoids in pros holds also.
\smallbreak

Here is a last theme about pros we would like to expose. As said before,
(colored) operads are promising devices to generalize usual formal power
series. Since pros are in some sense generalizations of operads, series
on pros would be an even more powerful generalization of such series.
Additionally, they could be very interesting devices for enumeration.
Due to the richness of the structure, series on pros come with a lot of
different products. At least, a generalization of the pre-Lie and
composition product on series on colored ns operads (see
Chapter~\ref{chap:buds}) can be considered.
\medbreak

\SkipTocEntry\section*{About biproducts and their algorithmic}
As exposed in Chapter~\ref{chap:shuffle}, combinatorial (bi)algebraic
structures are good supports to ask questions of analysis of algorithms.
To be more precise, given a combinatorial space $\K \Angle{C}$ endowed
with a biproduct $\Biproduct$ of arity $p$, the question of the
complexity of the computation of
$\Biproduct(x_1 \otimes \dots \otimes x_p)$, where $x_1$, \dots, $x_p$
are objects of $C$ seems in general unexplored and open. The analysis
may be performed with respect to the sum of the sizes of $x_1$, \dots,
$x_p$. This could lead to a hierarchy of biproducts depending on their
complexity. Some biproducts can have different complexity on different
bases of $\K \Angle{C}$. This research direction mixes in a balanced way
algebraic combinatorics and computer science.
\medbreak

\bibliographystyle{alpha}
\bibliography{Bibliography}

\newcommand{\etalchar}[1]{$^{#1}$}
\begin{thebibliography}{BMFPR11}

\bibitem[Agu00]{Agu00}
M.~Aguiar.
\newblock {Pre-Poisson algebras}.
\newblock {\em Lett. Math. Phys.}, 54(4):263--277, 2000.

\bibitem[AJ08]{AS08}
N.~Alon and Spencer J.
\newblock {\em {The Probabilistic Method}}.
\newblock Wiley-Blackwell, 3rd edition, 2008.

\bibitem[AL04]{AL04}
M.~Aguiar and J.-L Loday.
\newblock {Quadri-algebras}.
\newblock {\em J. Pure Appl. Algebra}, 191(3):205--221, 2004.

\bibitem[Ale37]{Ale37}
P.~Alexandroff.
\newblock {Diskrete R{\"a}ume}.
\newblock {\em Mat. Sb.}, 2(44)(3):501--519, 1937.

\bibitem[All00]{All00}
C.~Allauzen.
\newblock {Calcul efficace du shuffle de $k$ mots}.
\newblock Technical report, Institut Gaspard Monge, Université de Marne la
  Vallée, 2000.

\bibitem[Ant96]{Ant96}
V.~Antimirov.
\newblock {Partial derivatives of regular expressions and finite automaton
  constructions}.
\newblock {\em Theor. Comput. Sci.}, 155:291--319, 1996.

\bibitem[AO08]{AO08}
M.~Aguiar and R.~C. Orellana.
\newblock {The Hopf algebra of uniform block permutations}.
\newblock {\em J. Algebr. Comb.}, 28(1):115--138, 2008.

\bibitem[AP15]{AP15}
J.~Alm and D.~Petersen.
\newblock {Brown's dihedral moduli space and freedom of the gravity operad}.
\newblock {\em To appear in Annales Scientifiques de l'ENS}, 2015.
\newblock \Arxiv{1509.09274v2}.

\bibitem[AVL62]{AVL62}
G.M. Adelson-Velsky and E.~M. Landis.
\newblock {An algorithm for the organization of information}.
\newblock {\em Soviet Mathematics Doklady}, 3:1259--1263, 1962.

\bibitem[BBL98]{BBL98}
P.~Bose, J.~F. Buss, and A.~Lubiw.
\newblock {Pattern Matching for Permutations}.
\newblock {\em Inform. Process. Lett.}, 65(5):277--283, 1998.

\bibitem[BF03]{BF03}
C.~Brouder and A.~Frabetti.
\newblock {QED Hopf algebras on planar binary trees}.
\newblock {\em J. Algebra}, no. 1:298--322, 2003.

\bibitem[BFK06]{BFK06}
C.~Brouder, A.~Frabetti, and C.~Krattenthaler.
\newblock {Non-commutative Hopf algebra of formal diffeomorphisms}.
\newblock {\em Adv. Math.}, 200(2):479--524, 2006.

\bibitem[BG16]{BG16}
J.-P. Bultel and S.~Giraudo.
\newblock {Combinatorial Hopf algebras from PROs}.
\newblock {\em J. Algebr. Comb.}, 44(2):455--493, 2016.

\bibitem[Bj{\"o}84]{Bjo84}
A~Bj{\"o}rner.
\newblock {Orderings of Coxeter groups}.
\newblock In {\em Combinatorics and algebra}, volume~34 of {\em Contemp.
  Math.}, pages 175--195. Amer. Math. Soc., Providence, RI, 1984.

\bibitem[BLL88]{BLL88}
F.~Bergeron, G.~Labelle, and P.~Leroux.
\newblock {Functional Equations for Data Structures}.
\newblock {\em Lect. Notes Comput. Sc.}, 294:73--80, 1988.

\bibitem[BLL98]{BLL98}
F.~Bergeron, G.~Labelle, and P.~Leroux.
\newblock {\em {Combinatorial species and tree-like structures}}, volume~67 of
  {\em Encyclopedia of Mathematics and its Applications}.
\newblock Cambridge University Press, Cambridge, 1998.

\bibitem[BLL13]{BLL13}
F.~Bergeron, G.~Labelle, and P.~Leroux.
\newblock {\em {Introduction to the Theory of Species of Structures}}.
\newblock Université du Québec à Montréal, 2013.

\bibitem[BM03]{BM03}
C.~Berger and I.~Moerdijk.
\newblock {Axiomatic homotopy theory for operads}.
\newblock {\em Comment. Math. Helv.}, 78(4):805--831, 2003.

\bibitem[BMFPR11]{BFP11}
M.~Bousquet-Mélou, É. Fusy, and L.-F. Préville-Ratelle.
\newblock {The number of intervals in the $m$-Tamari lattices}.
\newblock {\em Electron. J. Comb.}, 18(2), 2011.
\newblock Paper 31.

\bibitem[BN98]{BN98}
F.~Baader and T.~Nipkow.
\newblock {\em {Term rewriting and all that}}.
\newblock Cambridge University Press, Cambridge, 1998.

\bibitem[Boz01]{Boz01}
S.~Bozapalidis.
\newblock {Context-free series on trees}.
\newblock {\em Inform. Comput.}, 169(2):186--229, 2001.

\bibitem[BPR12]{BP12}
F.~Bergeron and L.-F. Préville-Ratelle.
\newblock {Higher trivariate diagonal harmonics via generalized Tamari posets}.
\newblock {\em J. Comb.}, 3(3):317--341, 2012.

\bibitem[BR82]{BR82}
J.~Berstel and C.~Reutenauer.
\newblock {Recognizable formal power series on trees}.
\newblock {\em Theor. Comput. Sci.}, 18(2):115--148, 1982.

\bibitem[BR88]{BR88}
J.~Berstel and C.~Reutenauer.
\newblock {\em {Rational series and their languages}}, volume~12 of {\em EATCS
  Monographs on Theoretical Computer Science}.
\newblock Springer-Verlag, Berlin, 1988.

\bibitem[BR10]{BR10}
J.~Berstel and C.~Reutenauer.
\newblock {\em {Noncommutative rational series with applications}}, volume 137.
\newblock Cambridge University Press, 2010.

\bibitem[Brz64]{Brz64}
J.~A. Brzozowski.
\newblock {Derivatives of regular expressions}.
\newblock {\em J. ACM}, 11(4):481--494, 1964.

\bibitem[BS14]{BS14}
S.~Buss and M.~Soltys.
\newblock {Unshuffling a Square is NP-Hard}.
\newblock {\em J. Comput. Syst. Sci.}, 80(4):766--776, 2014.

\bibitem[BV73]{BV73}
J.~M. Boardman and R.~M. Vogt.
\newblock {Homotopy invariant algebraic structures on topological spaces}.
\newblock {\em Lect. Notes Math.}, 347, 1973.

\bibitem[Car07]{Car07}
P.~Cartier.
\newblock {A primer of Hopf algebras}.
\newblock In {\em Frontiers in number theory, physics, and geometry. {II}},
  pages 537--615. Springer, Berlin, 2007.

\bibitem[Cay57]{Cay1857}
A.~Cayley.
\newblock {On the Theory of the Analytical Forms called Trees}.
\newblock {\em Phil. Mag.}, 13:172--176, 1857.

\bibitem[CCM10]{CCM10}
P.~Caron, J.-M. Champarnaud, and L~Mignot.
\newblock {Acyclic automata and small expressions using multi-tilde-bar
  operators}.
\newblock {\em Theor. Comput. Sci.}, 411(38--39):3423--3435, 2010.

\bibitem[CCM11]{CCM11}
P.~Caron, J.-C. Champarnaud, and L.~Mignot.
\newblock {Multi-Bar and Multi-Tilde Regular Operators}.
\newblock {\em J. Autom. Lang. Comb.}, 16(1):11--26, 2011.

\bibitem[CCM12]{CCM12}
P.~Caron, J.-M. Champarnaud, and L~Mignot.
\newblock {Multi-tilde-bar expressions and their automata}.
\newblock {\em Acta Inform.}, 49(6):413--436, 2012.

\bibitem[CDG{\etalchar{+}}07]{CDGJLLTT07}
H.~Comon, M.~Dauchet, R.~Gilleron, F.~Jacquemard, C.~Löding, D.~Lugiez,
  S.~Tison, and M.~Tommasi.
\newblock {\em {Tree Automata Techniques and Applications}}.
\newblock \url{http://www.grappa.univ-lille3.fr/tata}, 2007.

\bibitem[CG14]{CG14}
F.~Chapoton and S.~Giraudo.
\newblock {Enveloping operads and bicoloured noncrossing configurations}.
\newblock {\em Exp. Math.}, 23(3):332--349, 2014.

\bibitem[CGM15]{CGM15}
H.~Cheballah, S.~Giraudo, and R.~Maurice.
\newblock {Hopf algebra structure on packed square matrices}.
\newblock {\em J. Comb. Theory A}, 133:139--182, 2015.

\bibitem[Cha02]{Cha02}
F.~Chapoton.
\newblock {Un théorème de Cartier-Milnor-Moore-Quillen pour les bigèbres
  dendriformes et les algèbres braces}.
\newblock {\em J. Pure Appl. Algebra}, 168(1):1--18, 2002.

\bibitem[Cha05]{Cha05}
F.~Chapoton.
\newblock {On some anticyclic operads}.
\newblock {\em Algebr. Geom. Topol.}, 5:53--69, 2005.

\bibitem[Cha06]{Cha06}
F.~Chapoton.
\newblock {Sur le nombre d'intervalles dans les treillis de Tamari}.
\newblock {\em Sém. Lothar. Combin.}, 55, 2006.

\bibitem[Cha07]{Cha07}
F.~Chapoton.
\newblock {The anticyclic operad of moulds}.
\newblock {\em Int. Math. Res. Notices}, 20:Art. ID rnm078, 36, 2007.

\bibitem[Cha08]{Cha08}
F.~Chapoton.
\newblock {Operads and algebraic combinatorics of trees}.
\newblock {\em Sém. Lothar. Combin.}, 58, 2008.

\bibitem[Cha09]{Cha09}
F.~Chapoton.
\newblock {A rooted-trees $q$-series lifting a one-parameter family of Lie
  idempotents}.
\newblock {\em Algebra \& Number Theory}, 3(6):611--636, 2009.

\bibitem[Cha14]{Cha14}
F.~Chapoton.
\newblock {Flows on rooted trees and the Menous-Novelli-Thibon idempotents}.
\newblock {\em Math. Scand.}, 115(1), 2014.

\bibitem[CHN16]{CHN16}
F.~Chapoton, F.~Hivert, and J.-C. Novelli.
\newblock {A set-operad of formal fractions and dendriform-like sub-operads}.
\newblock {\em J. Algebra}, 465:322--355, 2016.

\bibitem[Cho59]{Cho59}
N.~Chomsky.
\newblock {On certain formal properties of grammars}.
\newblock {\em Inform. Control}, 2:137--167, 1959.

\bibitem[CK98]{CK98}
A.~Connes and D.~Kreimer.
\newblock {Hopf algebras, renormalization and noncommutative geometry}.
\newblock {\em Commun. Math. Phys.}, 199(1):203--242, 1998.

\bibitem[CL01]{CL01}
F.~Chapoton and M.~Livernet.
\newblock {Pre-Lie algebras and the rooted trees operad}.
\newblock {\em Int. Math. Res. Notices}, 8:395--408, 2001.

\bibitem[CL07]{CL07}
F.~Chapoton and M.~Livernet.
\newblock {Relating two Hopf algebras built from an operad}.
\newblock {\em Int. Math. Res. Notices}, 2007(24):Art. ID rnm131, 27 pp, 2007.

\bibitem[CLRS09]{CLRS09}
T.H. Cormen, C.~E. Leiserson, R.L. Rivest, and C.~Stein.
\newblock {\em {Introduction to algorithms}}.
\newblock The MIT Press, 3rd edition, 2009.

\bibitem[Com72]{Com72}
L.~Comtet.
\newblock {Sur les coefficients de l'inverse de la série formelle $\sum n!
  t^n$}.
\newblock {\em C. R. Acad. Sci. Paris}, A 275:569--572, 1972.

\bibitem[Cox34]{Cox34}
H.~S.~M. Coxeter.
\newblock {Discrete groups generated by reflections}.
\newblock {\em Ann. Math.}, 35(3):588--621, 1934.

\bibitem[CP92]{CP92}
V.~Capoyleas and J.~Pach.
\newblock {A Tur\'an-type theorem on chords of a convex polygon}.
\newblock {\em J. Comb. Theory B}, 56(1):9--15, 1992.

\bibitem[CP17]{CP17}
G.~Châtel and V.~Pilaud.
\newblock {Cambrian Hopf algebras}.
\newblock {\em Adv. Math.}, 311:598--633, 2017.

\bibitem[CS63]{CS63}
N.~Chomsky and M.~P. Sch{\"u}tzenberger.
\newblock {The algebraic theory of context-free languages}.
\newblock In {\em Computer programming and formal systems}, pages 118--161.
  North-Holland, Amsterdam, 1963.

\bibitem[CSZ15]{CSZ15}
C.~Ceballos, F.~Santos, and G.~Ziegler.
\newblock {Many non-equivalent realizations of the associahedron}.
\newblock {\em Combinatorica}, 35(5):513--551, 2015.

\bibitem[Cur12]{Cur12}
P.-L. Curien.
\newblock Operads, clones, and distributive laws.
\newblock In {\em Operads and universal algebra}, volume~9 of {\em Nankai Ser.
  Pure Appl. Math. Theoret. Phys.}, pages 25--49. World Sci. Publ., Hackensack,
  NJ, 2012.

\bibitem[DG94]{DG94}
S.~Dulucq and O.~Guibert.
\newblock {Mots de piles, tableaux standards et permutations de Baxter}.
\newblock {\em Formal Power Series and Algebraic Combinatorics}, 1994.

\bibitem[DHNT11]{DHNT11}
G.~Duchamp, F.~Hivert, J.-C. Novelli, and J.-Y. Thibon.
\newblock {Noncommutative symmetric functions VII: free quasi-symmetric
  functions revisited}.
\newblock {\em Ann. Comb.}, 15(4):655--673, 2011.

\bibitem[DHT02]{DHT02}
G.~Duchamp, F.~Hivert, and J.-Y. Thibon.
\newblock {Noncommutative Symmetric Functions VI: Free Quasi- Symmetric
  Functions and Related Algebras}.
\newblock {\em Int. J. of Algebr. Comput.}, 12:671--717, 2002.

\bibitem[DK10]{DK10}
V.~Dotsenko and A.~Khoroshkin.
\newblock {Gr{\"o}bner bases for operads}.
\newblock {\em Duke Math. J.}, 153(2):363--396, 2010.

\bibitem[DLRS10]{DRS10}
J.~A. De~Loera, J.~Rambau, and F.~Santos.
\newblock {\em {Triangulations}}, volume~25 of {\em Algorithms and Computation
  in Mathematics}.
\newblock Springer-Verlag, Berlin, 2010.
\newblock Structures for algorithms and applications.

\bibitem[DM47]{DM47}
A.~Dvoretzky and Th. Motzkin.
\newblock {A problem of arrangements}.
\newblock {\em Duke Math. J.}, 14(2):305--313, 1947.

\bibitem[EFM09]{EM09}
K.~Ebrahimi-Fard and D.~Manchon.
\newblock Dendriform equations.
\newblock {\em J. Algebra}, 322(11):4053--4079, 2009.

\bibitem[EFM14]{EM14}
K.~Ebrahimi-Fard and D.~Manchon.
\newblock {The Magnus expansion, trees and Knuth's rotation correspondence}.
\newblock {\em Found. Comput. Math.}, 14(1):1--25, 2014.

\bibitem[EFMP08]{EMP08}
K.~Ebrahimi-Fard, D.~Manchon, and F.~Patras.
\newblock {New identities in dendriform algebras}.
\newblock {\em J. Algebra}, 320(2):708--727, 2008.

\bibitem[Eil74]{Eil74}
S.~Eilenberg.
\newblock {\em {Automata, languages, and machines. Vol. A}}.
\newblock Academic Press, New York, 1974.
\newblock Pure and Applied Mathematics, Vol. 58.

\bibitem[EML53]{EM53}
S.~Eilenberg and S.~Mac~Lane.
\newblock {On the groups of {$H(\Pi,n)$}. {I}}.
\newblock {\em Ann. Math.}, 58(2):55--106, 1953.

\bibitem[EZ76]{EZ76}
A.~Ehrenfeucht and H.-P. Zeiger.
\newblock {Complexity Measures for Regular Expressions}.
\newblock {\em J. Comput. Syst. Sci.}, 12(2):134--146, 1976.

\bibitem[FFM16]{FFM16}
F.~Fauvet, L.~Foissy, and D.~Manchon.
\newblock {Operads of finite posets}.
\newblock {\em \Arxiv{1604.08149v2}}, 2016.

\bibitem[FH91]{FH91}
W.~Fulton and J.~Harris.
\newblock {\em {Representation Theory: A First Course}}.
\newblock Springer-Verlag, New York, 1991.

\bibitem[FN99]{FN99}
P.~Flajolet and M.~Noy.
\newblock {Analytic combinatorics of non-crossing configurations}.
\newblock {\em Discrete Math.}, 204(1-3):203--229, 1999.

\bibitem[FNT14]{FNT14}
L.~Foissy, J.-C. Novelli, and J.-Y. Thibon.
\newblock {Polynomial realizations of some combinatorial Hopf algebras}.
\newblock {\em J. Noncommut. Geom.}, 8(1):141--162, 2014.

\bibitem[Foi02a]{Foi02a}
L.~Foissy.
\newblock {Les algèbres de Hopf des arbres enracinés décorés. I}.
\newblock {\em B. Sci. Math.}, 126(3):193--239, 2002.

\bibitem[Foi02b]{Foi02b}
L.~Foissy.
\newblock {Les algèbres de Hopf des arbres enracinés décorés. II}.
\newblock {\em B. Sci. Math.}, 126(4):249--288, 2002.

\bibitem[Foi07]{Foi07}
L.~Foissy.
\newblock {Bidendriform bialgebras, trees, and free quasi-symmetric functions}.
\newblock {\em J. Pure Appl. Algebra}, 209(2):439--459, 2007.

\bibitem[Foi08]{Foi08}
L.~Foissy.
\newblock {Faà di Bruno subalgebras of the Hopf algebra of planar trees from
  combinatorial Dyson-Schwinger equations}.
\newblock {\em Adv. Math.}, 218(1):136--162, 2008.

\bibitem[Foi12]{Foi12}
L.~Foissy.
\newblock {Ordered forests and parking functions}.
\newblock {\em Int. Math. Res. Notices}, 2012(7):1603--1633, 2012.

\bibitem[Foi15]{Foi15}
L.~Foissy.
\newblock {Examples of Com-PreLie Hopf algebras}.
\newblock {\em \Arxiv{1501.06375v1}}, 2015.

\bibitem[Fra08]{Fra08}
A.~Frabetti.
\newblock {Groups of tree-expanded series}.
\newblock {\em J. Algebra}, 319(1):377--413, 2008.

\bibitem[FS09]{FS09}
P.~Flajolet and R.~Sedgewick.
\newblock {\em {Analytic Combinatorics}}.
\newblock Cambridge University Press, 2009.

\bibitem[Ful97]{Ful97}
W.~Fulton.
\newblock {\em {Young Tableaux, with Applications to Representation Theory and
  Geometry}}.
\newblock London Mathematical Society Student Texts. Cambridge University
  Press, 1997.

\bibitem[GBV88]{GV88}
D.~Gouyou-Beauchamps and X.~Viennot.
\newblock {Equivalence of the two-dimensional directed animal problem to a
  one-dimensional path problem}.
\newblock {\em Adv. Appl. Math.}, 9(3):334--357, 1988.

\bibitem[Ger63]{Ger63}
M.~Gerstenhaber.
\newblock {The cohomology structure of an associative ring}.
\newblock {\em Ann. Math.}, 78:267--288, 1963.

\bibitem[Get94]{Get94}
E.~Getzler.
\newblock {Two-dimensional topological gravity and equivariant cohomology}.
\newblock {\em Commun. Math. Phys.}, 163(3):473--489, 1994.

\bibitem[Gir11]{Gir11}
S.~Giraudo.
\newblock {\em {Combinatoire algébrique des arbres}}.
\newblock PhD thesis, Université de Marne la Vallée, 2011.

\bibitem[Gir12a]{Gir12a}
S.~Giraudo.
\newblock {Algebraic and combinatorial structures on pairs of twin binary
  trees}.
\newblock {\em J. Algebra}, 360:115--157, 2012.

\bibitem[Gir12b]{Gir12d}
S.~Giraudo.
\newblock {Algebraic and combinatorial structures on pairs of twin binary
  trees}.
\newblock {\em Journal of Algebra}, 360:115--157, 2012.

\bibitem[Gir12c]{Gir12b}
S.~Giraudo.
\newblock {Constructing combinatorial operads from monoids}.
\newblock {\em Formal Power Series and Algebraic Combinatorics}, pages
  229--240, 2012.

\bibitem[Gir12d]{Gir12c}
S.~Giraudo.
\newblock {Construction d'opérades ensemblistes à partir de monoïdes}.
\newblock {\em C. R. Math.}, 350(11--12):549--552, 2012.

\bibitem[Gir12e]{Gir12e}
S.~Giraudo.
\newblock {Intervals of balanced binary trees in the Tamari lattice}.
\newblock {\em Theor. Comput. Sci.}, 420:1--27, 2012.

\bibitem[Gir15]{Gir15a}
S.~Giraudo.
\newblock {Combinatorial operads from monoids}.
\newblock {\em J. Algebr. Comb.}, 41(2):493--538, 2015.

\bibitem[Gir16a]{Gir16d}
S.~Giraudo.
\newblock {Colored operads, series on colored operads, and combinatorial
  generating systems}.
\newblock {\em \Arxiv{1605.04697v1}}, 2016.

\bibitem[Gir16b]{Gir16c}
S.~Giraudo.
\newblock {Operads from posets and Koszul duality}.
\newblock {\em Eur. J. Combin.}, 56C:1--32, 2016.

\bibitem[Gir16c]{Gir16a}
S.~Giraudo.
\newblock {Pluriassociative algebras I: The pluriassociative operad}.
\newblock {\em Adv. Appli. Math.}, 77:1--42, 2016.

\bibitem[Gir16d]{Gir16b}
S.~Giraudo.
\newblock {Pluriassociative algebras II: The polydendriform operad and related
  operads}.
\newblock {\em Adv. Appli. Math.}, 77:43--85, 2016.

\bibitem[Gir17a]{Gir17b}
S.~Giraudo.
\newblock {Combalgebraic structures on decorated cliques}.
\newblock {\em Formal Power Series and Algebraic Combinatorics}, 2017.

\bibitem[Gir17b]{Gir17a}
S.~Giraudo.
\newblock {Operads of decorated cliques}.
\newblock {\em \Arxiv{1702.00217v1}}, 2017.

\bibitem[GK94]{GK94}
V.~Ginzburg and M.~M. Kapranov.
\newblock {Koszul duality for operads}.
\newblock {\em Duke Math. J.}, 76(1):203--272, 1994.

\bibitem[GK95]{GK95}
E.~Getzler and M.~M. Kapranov.
\newblock Cyclic operads and cyclic homology.
\newblock In {\em Geometry, topology, \& physics}, Conf. Proc. Lecture Notes
  Geom. Topology, IV, pages 167--201. Int. Press, Cambridge, MA, 1995.

\bibitem[GKL{\etalchar{+}}95]{GKLLRT94}
I.M. Gelfand, D.~Krob, A.~Lascoux, B.~Leclerc, V.S. Retakh, and J.-Y. Thibon.
\newblock {Noncommutative symmetric functions I}.
\newblock {\em Adv. Math.}, 112, 1995.

\bibitem[GL01]{GL01}
J.~Gordon and M.~Liebeck.
\newblock {\em {Representations and Characters of Groups}}.
\newblock Cambridge University Press, 2nd edition, 2001.

\bibitem[GLMN16]{GLMN16}
S.~Giraudo, J.-G. Luque, L.~Mignot, and F.~Nicart.
\newblock {Operads, quasiorders, and regular languages}.
\newblock {\em Adv. Appl. Math.}, 75:56--93, 2016.

\bibitem[Glu61]{Glu61}
V.~M. Glushkov.
\newblock {The abstract theory of automata}.
\newblock {\em Russ. Math. surv.}, 16:1--53, 1961.

\bibitem[GR63]{GR63}
G.~Th. Guilbaud and P.~Rosenstiehl.
\newblock {Analyse algébrique d'un scrutin}.
\newblock {\em Mathématiques et sciences humaines}, tome 4:9--33, 1963.

\bibitem[GR16]{GR16}
D.~Grinberg and V.~Reiner.
\newblock {Hopf Algebras in Combinatorics}.
\newblock {\em \Arxiv{1409.8356v4}}, 2016.

\bibitem[GS84]{GS84}
F.~Gècseg and M.~Steinby.
\newblock {\em {Tree Automata}}.
\newblock Akadémia Kiadó, 1984.
\newblock 2nd edition (2015) available at
  \href{http://arxiv.org/abs/1509.06233}{arXiv:1509.06233 [cs.FL]}.

\bibitem[GV16]{GV16}
S.~Giraudo and S.~Vialette.
\newblock {Unshuffling Permutations}.
\newblock {\em Latin American Theoretical Informatics Symposium, LNCS},
  9644:509--521, 2016.

\bibitem[Har78]{Har78}
M.~A. Harrison.
\newblock {\em {Introduction to formal language theory}}.
\newblock Addison-Wesley Publishing Co., Reading, Mass., 1978.

\bibitem[Hiv99]{Hiv99}
F.~Hivert.
\newblock {\em {Combinatoire des fonctions quasi-symétriques}}.
\newblock PhD thesis, Université de Marne la Vallée, 1999.

\bibitem[Hiv03]{Hiv03}
F.~Hivert.
\newblock {An introduction to Combinatorial Hopf Algebras}.
\newblock {\em IOS Press}, 2003.

\bibitem[HMU06]{HMU06}
J.~E. Hopcroft, R.~Matwani, and J.~D. Ullman.
\newblock {\em {Introduction to Automata Theory, Languages, and Computation}}.
\newblock Pearson, 3rd edition, 2006.

\bibitem[HN07]{HN07}
F.~Hivert and J.~Nzeutchap.
\newblock {Dual graded graphs in combinatorial Hopf algebras}, 2007.
\newblock unpublished.

\bibitem[HNT05]{HNT05}
F.~Hivert, J.-C. Novelli, and J.-Y. Thibon.
\newblock {The Algebra of Binary Search Trees}.
\newblock {\em Theor. Comput. Sc.}, 339(1):129--165, 2005.

\bibitem[Hof10]{Hof10}
E.~Hoffbeck.
\newblock {A Poincar\'e-Birkhoff-Witt criterion for Koszul operads}.
\newblock {\em Manuscripta Math.}, 131(1-2):87--110, 2010.

\bibitem[HRS12]{HRS12}
D.~Henshall, N.~Rampersad, and J.~Shallit.
\newblock {Shuffling and Unshuffling}.
\newblock {\em Bulletin of the EATCS}, 107:131--142, 2012.

\bibitem[HT72]{HT72}
S.~Huang and D.~Tamari.
\newblock {Problems of associativity: A simple proof for the lattice property
  of systems ordered by a semi-associative law}.
\newblock {\em J. Comb. Theory. A}, 13:7--13, 1972.

\bibitem[Joh63]{Joh63}
S.~Johnson.
\newblock {Generation of permutations by adjacent transposition}.
\newblock {\em Math. Comput.}, 17:282--285, 1963.

\bibitem[Joy81]{Joy81}
A.~Joyal.
\newblock {Une théorie combinatoire des séries formelles}.
\newblock {\em Adv. Math.}, 42(1):1--82, 1981.

\bibitem[JR79]{JR79}
S.~A. Joni and G.-C. Rota.
\newblock {Coalgebras and bialgebras in combinatorics}.
\newblock {\em Stud. Appl. Math.}, 61:93--139, 1979.

\bibitem[Kac90]{Kac90}
V.~G. Kac.
\newblock {\em {Infinite-Dimensional Lie Algebras}}.
\newblock Cambridge University Press, 3rd edition, 1990.

\bibitem[KB70]{KB70}
D.~Knuth and P.~Bendix.
\newblock {Simple word problems in universal algebras}.
\newblock In {\em {Computational Problems in Abstract Algebra}}, pages
  263--297. Pergamon, Oxford, 1970.

\bibitem[KLT97]{KLT97}
D.~Krob, B.~Leclerc, and J.-Y. Thibon.
\newblock {Noncommutative symmetric functions II. Transformations of
  alphabets}.
\newblock {\em Int. J. Algebr. Comput.}, 7(2):181--264, 1997.

\bibitem[Knu97]{Knu97}
D.~Knuth.
\newblock {\em {The Art of Computer Programming, volume 1: Fundamental
  Algorithms}}.
\newblock Addison Wesley Longman, 3rd edition, 1997.

\bibitem[Knu98]{Knu98}
D.~Knuth.
\newblock {\em The art of computer programming, volume 3: Sorting and
  searching}.
\newblock Addison Wesley Longman, 1998.

\bibitem[Koc09]{Koc09}
J.~Kock.
\newblock {\em {Notes on Polynomial Functors}}.
\newblock Unpublished, 2009.
\newblock {\tt \url{http://mat.uab.es/~kock/cat/%
  notes-on-polynomial-functors.html}}.

\bibitem[KT97]{KT97}
D.~Krob and J.-Y. Thibon.
\newblock {Noncommutative symmetric functions IV: Quantum linear groups and
  Hecke algebras at $q = 0$}.
\newblock {\em Journal Algebr. Comb.}, 6:339--376, 1997.

\bibitem[KT99]{KT99}
D.~Krob and J.-Y. Thibon.
\newblock {Noncommutative symmetric functions V: A degenerate version of
  $U_q(gl_N)$}.
\newblock {\em Int. J. Algebr. Comput.}, 3-4:405--430, 1999.

\bibitem[Kup96]{Kup96}
G.~Kuperberg.
\newblock {Another proof of the alternating-sign matrix conjecture}.
\newblock {\em Int. Math. Res. Notices}, 3:139--150, 1996.

\bibitem[Lab81]{Lab81}
G.~Labelle.
\newblock {Une nouvelle démonstration combinatoire des formules d'inversion de
  Lagrange}.
\newblock {\em Adv. Math.}, 42(3):217--247, 1981.

\bibitem[Laf03]{Laf03}
Y.~Lafont.
\newblock {Towards an algebraic theory of Boolean circuits}.
\newblock {\em J. Pure Appl. Algebra}, 184(2-3):257--310, 2003.

\bibitem[Laf11]{Laf11}
Y.~Lafont.
\newblock {Diagram rewriting and operads}.
\newblock {\em Sémin. Congr.}, 26:163--179, 2011.

\bibitem[Las84]{Las84}
A.~Lascoux.
\newblock {Fonctions symétriques}.
\newblock {\em Sém. Lothar. Combin.}, 8:37--58, 1984.

\bibitem[Lei04]{Lei04}
T.~Leinster.
\newblock {\em Higher operads, higher categories}, volume 298 of {\em London
  Mathematical Society Lecture Note Series}.
\newblock Cambridge University Press, Cambridge, 2004.

\bibitem[Ler04]{Ler04}
P.~Leroux.
\newblock {Ennea-algebras}.
\newblock {\em J. Algebra}, 281(1):287--302, 2004.

\bibitem[Ler07]{Ler07}
P.~Leroux.
\newblock {A simple symmetry generating operads related to rooted planar
  {$m$}-ary trees and polygonal numbers}.
\newblock {\em J. Integer Seq.}, 10(4):Article 07.4.7, 23, 2007.

\bibitem[Ler11]{Ler11}
P.~Leroux.
\newblock L-algebras, triplicial-algebras, within an equivalence of categories
  motivated by graphs.
\newblock {\em Comm. Algebra}, 39(8):2661--2689, 2011.

\bibitem[Liv06]{Liv06}
M.~Livernet.
\newblock {A rigidity theorem for pre-Lie algebras}.
\newblock {\em J. Pure Appl. Algebra}, 207(1):1--18, 2006.

\bibitem[LM06]{LM06}
M.~Lohrey and S.~Maneth.
\newblock {The complexity of tree automata and XPath on grammar-compressed
  trees}.
\newblock {\em Theor. Comput. Sci.}, 363(2):196--210, 2006.

\bibitem[LMN13]{LMN13}
J.-G. Luque, L.~Mignot, and F.~Nicart.
\newblock {Some Combinatorial Operators in Language Theory}.
\newblock {\em J. Autom. Lang. Comb.}, 18(1):27--52, 2013.

\bibitem[LN13]{LN13}
J.-L. Loday and N.~M. Nikolov.
\newblock {Operadic construction of the renormalization group}.
\newblock In {\em Lie theory and its applications in physics}, volume~36 of
  {\em Springer Proc. Math. Stat.}, pages 191--211. Springer, Tokyo, 2013.

\bibitem[Lod93]{Lod93}
J.-L. Loday.
\newblock {Une version non commutative des algèbres de Lie~: les algèbres de
  Leibniz}.
\newblock In {\em R.{C}.{P}.\ 25, {V}ol.\ 44}, volume~41 of {\em Prépubl.
  Inst. Rech. Math. Av.}, pages 127--151. Univ. Louis Pasteur, Strasbourg,
  1993.

\bibitem[Lod95]{Lod95}
J.-L. Loday.
\newblock {Cup-product for Leibniz cohomology and dual Leibniz algebras}.
\newblock {\em Math. Scand.}, 77(2), 1995.

\bibitem[Lod96]{Lod96}
J.-L. Loday.
\newblock {La renaissance des opérades}.
\newblock {\em Séminaire Bourbaki}, 37(792):47--74, 1996.

\bibitem[Lod01]{Lod01}
J.-L. Loday.
\newblock {Dialgebras}.
\newblock In {\em Dialgebras and related operads}, volume 1763 of {\em Lecture
  Notes in Math.}, pages 7--66. Springer, Berlin, 2001.

\bibitem[Lod02]{Lod02}
J.-L. Loday.
\newblock Arithmetree.
\newblock {\em J. Algebra}, 258(1):275--309, 2002.
\newblock Special issue in celebration of Claudio Procesi's 60th birthday.

\bibitem[Lod06]{Lod06}
J.-L. Loday.
\newblock {Completing the operadic butterfly}.
\newblock {\em Georgian Math. J.}, 13(4):741--749, 2006.

\bibitem[Lod08]{Lod08}
J.-L. Loday.
\newblock {Generalized bialgebras and triples of operads}.
\newblock {\em Astérisque}, 320:x+116, 2008.

\bibitem[Lod10]{Lod10}
J.-L. Loday.
\newblock {On the operad of associative algebras with derivation}.
\newblock {\em Georgian Math. J.}, 17(2):347--372, 2010.

\bibitem[Lot02]{Lot02}
M.~Lothaire.
\newblock {\em {Algebraic combinatorics on words}}.
\newblock Encyclopedia of mathematics and its applications. Cambridge
  university press, New York, 2002.

\bibitem[LPRR15]{LPR15}
D.~López, L.-F. Préville-Ratelle, and M.~Ronco.
\newblock {Algebraic structures defined on $m$-Dyck paths}.
\newblock {\em \Arxiv{1508.01252v2}}, 2015.

\bibitem[LR34]{LR34}
D.~E. Littlewood and A.~R. Richardson.
\newblock {Group characters and algebra}.
\newblock {\em Philos. T. Soc. A}, 233:99--141, 1934.

\bibitem[LR98]{LR98}
J.-L. Loday and M.~Ronco.
\newblock {Hopf Algebra of the Planar Binary Trees}.
\newblock {\em Adv. Math.}, 139:293--309, 1998.

\bibitem[LR02]{LR02}
J.-L. Loday and M.~Ronco.
\newblock {Order Structure on the Algebra of Permutations and of Planar Binary
  Trees}.
\newblock {\em J. Algebr. Comb.}, 15(3):253--270, 2002.

\bibitem[LR03]{LR03}
J.-L. Loday and M.~Ronco.
\newblock {Algèbres de Hopf colibres}.
\newblock {\em C. R. Math.}, 337(3):153--158, 2003.

\bibitem[LR04]{LR04}
J.-L. Loday and M.~Ronco.
\newblock {Trialgebras and families of polytopes}.
\newblock In {\em Homotopy theory: relations with algebraic geometry, group
  cohomology, and algebraic {$K$}-theory}, volume 346 of {\em Contemp. Math.},
  pages 369--398. Amer. Math. Soc., Providence, RI, 2004.

\bibitem[LR06]{LR06}
J.-L. Loday and M.~Ronco.
\newblock On the structure of cofree {H}opf algebras.
\newblock {\em J. Reine Angew. Math.}, 592:123--155, 2006.

\bibitem[LR12]{LR12}
S.~Law and N.~Reading.
\newblock {The Hopf algebra of diagonal rectangulations}.
\newblock {\em J. Comb. Theory A}, 119(3):788--824, 2012.

\bibitem[LS81]{LS81}
A.~Lascoux and M.-P. Sch{\"u}tzenberger.
\newblock {Le monoïde plaxique}.
\newblock In {\em Noncommutative structures in algebra and geometric
  combinatorics ({N}aples, 1978)}, volume 109 of {\em Quad. ``Ricerca Sci.''},
  pages 129--156. CNR, Rome, 1981.

\bibitem[LV12]{LV12}
J.-L. Loday and B.~Vallette.
\newblock {\em {Algebraic Operads}}, volume 346 of {\em Grundlehren der
  mathematischen Wissenschaften}.
\newblock Springer, Heidelberg, 2012.

\bibitem[Mac15]{Mac15}
I.~G. Macdonald.
\newblock {\em {Symmetric Functions and Hall Polynomials}}.
\newblock Oxford Classic Texts in the Physical Sciences. Oxford University
  Press, second edition, 2015.

\bibitem[Man83]{Man83}
A.~Mansfield.
\newblock {On the computational complexity of a merge recognition problem}.
\newblock {\em Discrete Appl. Math.}, 5:119--122, 1983.

\bibitem[Man11]{Man11}
D.~Manchon.
\newblock {A short survey on pre-Lie algebras}.
\newblock In {\em Noncommutative geometry and physics: renormalisation,
  motives, index theory}, ESI Lect. Math. Phys., pages 89--102. Eur. Math.
  Soc., Zürich, 2011.

\bibitem[Mar08]{Mar08}
M.~Markl.
\newblock {Operads and PROPs}.
\newblock In {\em Handbook of Algebra}, volume~5, pages 87--140.
  Elsevier/North-Holland, Amsterdam, 2008.

\bibitem[May72]{May72}
J.~P. May.
\newblock {\em {The geometry of iterated loop spaces}}.
\newblock Springer-Verlag, Berlin-New York, 1972.
\newblock Lectures Notes in Mathematics, Vol. 271.

\bibitem[ML65]{McL65}
S.~Mac~Lane.
\newblock {Categorical algebra}.
\newblock {\em Bull. Amer. Math. Soc.}, 71:40--106, 1965.

\bibitem[ML14]{ML14}
M.~Méndez and J.~Liendo.
\newblock {An antipode formula for the natural Hopf algebra of a set operad}.
\newblock {\em Adv. Appl. Math.}, 53:112--140, 2014.

\bibitem[MN93]{MN93}
M.~Méndez and O.~Nava.
\newblock {Colored species, $c$-monoids, and plethysm I}.
\newblock {\em J. Comb. Theory A}, 64(1):102--129, 1993.

\bibitem[Mot48]{Mot48}
Th. Motzkin.
\newblock {Relations between hypersurface cross ratios, and a combinatorial
  formula for partitions of a polyion, for permanent preponderance, and for
  non-associative products}.
\newblock {\em Bull. Amer. Math. Soc.}, 54:352--360, 1948.

\bibitem[MPRS79]{MPRS79}
R.~E. Miller, N.~Pippenger, A.~L. Rosenberg, and L.~Snyder.
\newblock {Optimal 2,3-Trees}.
\newblock {\em SIAM J. Comput.}, 8:42--59, 1979.

\bibitem[MR95]{MR95}
C.~Malvenuto and C.~Reutenauer.
\newblock {Duality between quasi-symmetric functions and the Solomon descent
  algebra}.
\newblock {\em J. Algebra}, 177(3):967--982, 1995.

\bibitem[MRR83]{MRR83}
W.H. Mills, D.P. Robbins, and H.~Jr. Rumsey.
\newblock {Alternating sign matrices and descending plane partitions}.
\newblock {\em J. Comb. Theory A.}, 34:340--359, 1983.

\bibitem[MY91]{MY91}
M.~Méndez and J.~Yang.
\newblock {M{\"o}bius Species}.
\newblock {\em Adv. Math.}, 85(1):83--128, 1991.

\bibitem[Mé15]{Men15}
M.~A. Méndez.
\newblock {\em {Set operads in combinatorics and computer science}}.
\newblock SpringerBriefs in Mathematics. Springer, Cham, 2015.

\bibitem[Nar55]{Nar55}
T.V. Narayana.
\newblock {Sur les treillis formés par les partitions d'un entier et leurs
  applications à la théorie des probabilités}.
\newblock {\em C. R. Acad. Sci. Paris}, 240:1188--1189, 1955.

\bibitem[New42]{New42}
M.~H.~A. Newman.
\newblock {On theories with a combinatorial definition of ``equivalence.''}.
\newblock {\em Ann. Math.}, 43(2):223--243, 1942.

\bibitem[Nov14]{Nov14}
J.-C. Novelli.
\newblock {$m$-dendriform algebras}.
\newblock {\em \Arxiv{1406.1616v1}}, 2014.

\bibitem[NT06]{NT06}
J.-C. Novelli and J.-Y. Thibon.
\newblock {Polynomial realizations of some trialgebras}.
\newblock {\em Formal Power Series and Algebraic Combinatorics}, 2006.

\bibitem[NT07]{NT07}
J.-C. Novelli and J.-Y. Thibon.
\newblock {Hopf algebras and dendriform structures arising from parking
  functions}.
\newblock {\em Fund. Math.}, 193:189--241, 2007.

\bibitem[NT10]{NT10}
J.-C. Novelli and J.-Y. Thibon.
\newblock {Free quasi-symmetric functions and descent algebras for wreath
  products, and noncommutative multi-symmetric functions}.
\newblock {\em Discrete Math.}, 310(24):3584--3606, 2010.

\bibitem[NT13]{NT13}
J.-C. Novelli and J.-Y. Thibon.
\newblock {Duplicial algebras and Lagrange inversion}.
\newblock {\em \Arxiv{1209.5959v2}}, 2013.

\bibitem[NT14]{NT14}
J.-C. Novelli and J.-Y. Thibon.
\newblock {Hopf Algebras of m-permutations, $(m + 1)$-ary trees, and
  $m$-parking functions}.
\newblock {\em \Arxiv{1403.5962v2}}, 2014.

\bibitem[PR85]{PR85}
A.~Proskurowski and F.~Ruskey.
\newblock {Binary tree gray codes}.
\newblock {\em J. Algorithms}, 6(2):225--238, 1985.

\bibitem[PR95]{PR95}
S.~Poirier and C.~Reutenauer.
\newblock {Algèbres de Hopf de tableaux}.
\newblock {\em Ann. Sci. Math. Québec}, 19:79--90, 1995.

\bibitem[Rea05]{Rea05}
N.~Reading.
\newblock {Lattice congruences, fans and Hopf algebras}.
\newblock {\em J. Comb. Theory A}, 110(2):237--273, 2005.

\bibitem[Ree58]{Ree58}
R.~Ree.
\newblock {Lie elements and an algebra associated with shuffles}.
\newblock {\em Ann. Math.}, 68(2):210--220, 1958.

\bibitem[Reu93]{Ret93}
C.~Reutenauer.
\newblock {\em {Free Lie algebras}}.
\newblock Clarendon Press, 1993.

\bibitem[Rey07]{Rey07}
M.~Rey.
\newblock {Algebraic constructions on set partitions}.
\newblock {\em Formal Power Series and Algebraic Combinatorics}, 2007.

\bibitem[Rot64]{Rot64}
G.-C. Rota.
\newblock {On the foundations of combinatorial theory. I. Theory of Möbius
  functions}.
\newblock {\em Z. Wahrscheinlichkeit}, 2:340--368, 1964.

\bibitem[Rus04]{Rus03}
F.~Ruskey.
\newblock {\em {Combinatorial Generation}}.
\newblock Preliminary working draft, 2004.

\bibitem[RV13]{RV13}
R.~Rizzi and S.~Vialette.
\newblock {On Recognizing Words That Are Squares for the Shuffle Product}.
\newblock In A.A Bulatov and A.M. Shur, editors, {\em 8th International
  Computer Science Symposium in Russia, CSR 2013, Ekaterinburg, Russia}, pages
  235--245, 2013.

\bibitem[Ryt04]{Ryt04}
W.~Rytter.
\newblock {Grammar compression, Z-encodings, and string algorithms with
  implicit input}.
\newblock In {\em Automata, languages and programming}, volume 3142 of {\em
  Lect. Notes Comput. Sci.}, pages 15--27. Springer, Berlin, 2004.

\bibitem[Ré85]{Rem85}
J.-L. Rémy.
\newblock {Un procédé itératif de dénombrement d'arbres binaires et son
  application à leur génération aléatoire}.
\newblock {\em RAIRO-Theor. Inf. Appl.}, 19(2):179--195, 1985.

\bibitem[Sak09]{Sak09}
J.~Sakarovitch.
\newblock {\em {Elements of Automata Theory}}.
\newblock Cambridge University Press, 2009.

\bibitem[Sch61]{Sch61}
C.~Schensted.
\newblock {Longest increasing and decreasing subsequences}.
\newblock {\em Canadian J. Math.}, 13:179--191, 1961.

\bibitem[Sch63]{Sch63}
M.~P. Sch{\"u}tzenberger.
\newblock {Certain elementary families of automata}.
\newblock In {\em {Proc. Sympos. Math. Theory of Automata (New York, 1962)}},
  pages 139--153. Polytechnic Press of Polytechnic Inst. of Brooklyn, Brooklyn,
  New York, 1963.

\bibitem[Slo]{Slo}
N.~J.~A. Sloane.
\newblock {The On-Line Encyclopedia of Integer Sequences}.
\newblock \url{https://oeis.org/}.

\bibitem[Soi10]{Soi10}
A.~Soifer.
\newblock {\em {Ramsey Theory: Yesterday, Today, and Tomorrow}}.
\newblock Birkhäuser, 2010.

\bibitem[Spe86]{Spe86}
J.-C. Spehner.
\newblock {Le calcul rapide des melanges de deux mots}.
\newblock {\em Theor. Comput. Sci.}, 47:181--203, 1986.

\bibitem[SS78]{SS78}
A.~Salomaa and M.~Soittola.
\newblock {\em {Automata-theoretic aspects of formal power series}}.
\newblock Springer-Verlag, New York-Heidelberg, 1978.
\newblock Texts and Monographs in Computer Science.

\bibitem[SS82]{SS82}
J.~A. Storer and T.~G. Szymanski.
\newblock {Data compression via textual substitution}.
\newblock {\em J. Assoc. Comput. Mach.}, 29(4):928--951, 1982.

\bibitem[SS85]{SS85}
R.~Simion and F.~W. Schmidt.
\newblock {Restricted permutations}.
\newblock {\em Eur. J. Combin.}, 6(4):383--406, 1985.

\bibitem[ST09]{ST09}
P.~Salvatore and R.~Tauraso.
\newblock {The operad Lie is free}.
\newblock {\em J. Pure Appl. Algebra}, 213(2):224--230, 2009.

\bibitem[Sta99]{Sta99}
R.~P. Stanley.
\newblock {\em {Enumerative Combinatorics}}, volume~2.
\newblock Cambridge University Press, 1999.

\bibitem[Sta11]{Sta11}
R.~P. Stanley.
\newblock {\em {Enumerative Combinatorics}}, volume~1.
\newblock Cambridge University Press, 2nd edition edition, 2011.

\bibitem[Ste64]{Ste64}
H.~Steinhaus.
\newblock {\em {One hundred problems in elementary mathematics}}.
\newblock Basic Books, Inc., Publishers, New York, 1964.

\bibitem[Tak71]{Tak71}
M.~Takeuchi.
\newblock {Free Hopf algebras generated by coalgebras}.
\newblock {\em J. Math. Soc. Jpn}, 23(4):561--582, 1971.

\bibitem[Tam62]{Tam62}
D.~Tamari.
\newblock {The algebra of bracketings and their enumeration}.
\newblock {\em Nieuw Arch. Wisk.}, 10(3):131--146, 1962.

\bibitem[Tro62]{Tro62}
H.~Trotter.
\newblock {Algorithm 115: Perm}.
\newblock {\em Commun. ACM}, 5(8):434--435, 1962.

\bibitem[Val07]{Val07}
B.~Vallette.
\newblock {Homology of generalized partition posets}.
\newblock {\em J. Pure Appl. Algebra}, 208(2):699--725, 2007.

\bibitem[Var14]{Var14}
Y.~Vargas.
\newblock {Hopf algebra of permutation pattern functions}.
\newblock {\em Formal Power Series and Algebraic Combinatorics}, pages
  839--850, 2014.

\bibitem[vdL04]{Vdl04}
P.~van~der Laan.
\newblock {\em {Operads. Hopf algebras and coloured Koszul duality}}.
\newblock PhD thesis, Universiteit Utrecht, 2004.

\bibitem[Vie86]{Vie86}
X.~Viennot.
\newblock Heaps of pieces. {I}. {B}asic definitions and combinatorial lemmas.
\newblock In {\em Combinatoire énumérative}, volume 1234 of {\em Lect. Notes
  Math.}, pages 321--350. Springer, 1986.

\bibitem[Vin63]{Vin63}
È.~B. Vinberg.
\newblock {The theory of homogeneous convex cones}.
\newblock {\em Trudy Moskov. Mat. Ob\v s\v c.}, 12:303--358, 1963.

\bibitem[vLN82]{LN82}
J.~van Leeuwen and M.~Nivat.
\newblock {Efficient Recognition of Rational Relations}.
\newblock {\em Inform. Process. Lett.}, 14(1):34--38, 1982.

\bibitem[Yau16]{Yau16}
D.~Yau.
\newblock {\em {Colored Operads}}.
\newblock Graduate Studies in Mathematics. American Mathematical Society, 2016.

\bibitem[YO69]{YO69}
T.~Yanagimoto and M.~Okamoto.
\newblock {Partial orderings of permutations and monotonicity of a rank
  correlation statistic}.
\newblock {\em Ann. I. Stat. Math.}, 21:489--506, 1969.

\bibitem[Zei96]{Zei96}
D~Zeilberger.
\newblock {Proof of the alternating sign matrix conjecture}.
\newblock {\em Electron. J. Comb.}, 3(2), 1996.
\newblock Paper 13.

\bibitem[Zin12]{Zin12}
G.~W. Zinbiel.
\newblock {Encyclopedia of types of algebras 2010}.
\newblock In {\em Operads and universal algebra}, volume~9 of {\em Nankai Ser.
  Pure Appl. Math. Theoret. Phys.}, pages 217--297. World Sci. Publ.,
  Hackensack, NJ, 2012.

\bibitem[ZL78]{ZL78}
J.~Ziv and A.~Lempel.
\newblock {Compression of individual sequences via variable-rate coding}.
\newblock {\em IEEE T. Inform. Theory}, 24(5):530--536, 1978.

\end{thebibliography}

\cleardoublepage

\BackCover

\end{document}